\documentclass[9pt]{amsart}
\usepackage{amsmath}
\usepackage{amsxtra}
\usepackage{amscd}
\usepackage{amsthm}
\usepackage{amsfonts}
\usepackage{amssymb}
\usepackage{eucal}
\usepackage[all]{xy}
\usepackage{graphicx}
\usepackage[usenames]{color}
\usepackage[hidelinks]{hyperref}

\newtheorem{cor}[subsubsection]{Corollary}

\newtheorem{lem}[subsubsection]{Lemma}
\newtheorem{prop}[subsubsection]{Proposition}

\newtheorem{thmconstr}[subsubsection]{Theorem-Construction}
\newtheorem{lemconstr}[subsubsection]{Lemma-Construction}

\newtheorem{thm}[subsubsection]{Theorem}

\newtheorem{defn}[subsubsection]{Definition}
\newtheorem{quest}[subsubsection]{Question}

\theoremstyle{remark}
\newtheorem{rem}[subsubsection]{Remark}

\theoremstyle{definition}

\theoremstyle{remark}

\newcommand{\thmref}[1]{Theorem~\ref{#1}}

\newcommand{\secref}[1]{Sect.~\ref{#1}}
\newcommand{\lemref}[1]{Lemma~\ref{#1}}
\newcommand{\propref}[1]{Proposition~\ref{#1}}
\newcommand{\corref}[1]{Corollary~\ref{#1}}

\newcommand{\remref}[1]{Remark~\ref{#1}}

\numberwithin{equation}{section}

\newcommand{\nc}{\newcommand}
\nc{\renc}{\renewcommand}
\nc{\ssec}{\subsection}
\nc{\sssec}{\subsubsection}
\nc{\on}{\operatorname}

\nc{\ips}{{\iota_P^{(S)}}}
\nc{\ipms}{{\iota_{P^-}^{(S)}}}
\nc{\sfpps}{{\sfp_P^{(S)}}}
\nc{\sfppms}{{\sfp_{P^-}^{(S)}}}

\nc\ol{\overline}
\nc\ul{\underline}
\nc\wt{\widetilde}
\nc\tboxtimes{\wt{\boxtimes}}
\nc\tstar{\wt{\star}}
\nc{\alp}{\alpha}

\nc{\ZZ}{{\mathbb Z}}
\nc{\NN}{{\mathbb N}}
\nc{\OO}{{\mathbb O}}
\renc{\SS}{{\mathbb S}}
\nc{\DD}{{\mathbb D}}
\nc{\GG}{{\mathbb G}}

\nc{\Fq}{{\mathbb F}_q}
\nc{\Fqb}{\ol{\mathbb F}_q}
\nc{\Ql}{{\mathbb Q}_\ell}
\nc{\Qlb}{{\ol{\mathbb Q}_\ell}}
\nc{\id}{\text{id}}
\nc\X{\mathcal X}

\nc{\red}{\on{red}}
\nc{\Ho}{\on{Ho}}
\nc{\Hom}{\on{Hom}}
\nc{\coHom}{\ul{\on{coHom}}}
\nc{\coMaps}{{\bf{coMaps}}}
\nc{\coef}{\on{coef}}
\nc{\Lie}{\on{Lie}}
\nc{\Loc}{\on{Loc}}
\nc{\coLoc}{\on{coLoc}}
\nc{\Pic}{\on{Pic}}
\nc{\Bun}{\on{Bun}}
\nc{\IC}{\on{IC}}
\nc{\Aut}{\on{Aut}}
\nc{\rk}{\on{rk}}
\nc{\Sh}{\on{Sh}}
\nc{\Perv}{\on{Perv}}
\nc{\pos}{{\on{pos}}}
\nc{\Conv}{\on{Conv}}
\nc{\Sph}{{\on{Sph}}}
\nc{\Sym}{\on{Sym}}
\nc{\BunBb}{\overline{\Bun}_B}
\nc{\BunNb}{\overline{\Bun}_N}
\nc{\BunTb}{\overline{\Bun}_T}
\nc{\BunBbm}{\overline{\Bun}_{B^-}}
\nc{\BunBbel}{\overline{\Bun}_{B,el}}
\nc{\BunBbmel}{\overline{\Bun}_{B^-,el}}
\nc{\Buno}{\overset{o}{\Bun}}
\nc{\BunPb}{{\overline{\Bun}_P}}
\nc{\BunBM}{\Bun_{B(M)}}
\nc{\BunBMb}{\overline{\Bun}_{B(M)}}
\nc{\BunPbw}{{\widetilde{\Bun}_P}}
\nc{\BunBP}{\widetilde{\Bun}_{B,P}}
\nc{\GUb}{\overline{G/U}}
\nc{\GUPb}{\overline{G/U(P)}}

\nc{\Hhom}{\underline{\on{Hom}}}
\nc\syminfty{\on{Sym}^{\infty}}
\nc\lal{\ol{\lambda}}
\nc\xl{\ol{x}}
\nc\thl{\ol{\theta}}
\nc\nul{\ol{\nu}}
\nc\mul{\ol{\mu}}
\nc\Sum\Sigma
\nc{\oX}{\overset{o}{X}{}}
\nc{\hl}{\overset{\leftarrow}h{}}
\nc{\hr}{\overset{\rightarrow}h{}}
\nc{\M}{{\mathcal M}}
\nc{\N}{{\mathcal N}}
\nc{\F}{{\mathcal F}}
\nc{\D}{{\mathcal D}}
\nc{\Q}{{\mathcal Q}}
\nc{\Y}{{\mathcal Y}}
\nc{\G}{{\mathcal G}}
\nc{\E}{{\mathcal E}}
\nc{\CalC}{{\mathcal C}}
\nc\Dh{\widehat{\D}}

\nc{\C}{{\mathcal C}}
\nc{\K}{{\mathcal K}}
\renewcommand{\H}{{\mathcal H}}

\nc{\T}{{\mathcal T}}
\nc{\V}{{\mathcal V}}
\renc{\P}{{\mathcal P}}
\nc{\A}{{\mathcal A}}
\nc{\B}{{\mathcal B}}
\nc{\U}{{\mathcal U}}

\nc{\Gr}{{\on{Gr}}}

\nc{\frn}{{\check{\mathfrak u}(P)}}

\nc{\fC}{\mathfrak C}
\nc{\fT}{\mathfrak T}
\nc{\p}{\mathfrak p}
\nc{\q}{\mathfrak q}
\nc\f{{\mathfrak f}}

\nc{\qo}{{\mathfrak q}}
\nc{\po}{{\mathfrak p}}
\nc{\s}{{\mathfrak s}}
\nc\w{\text{w}}

\renewcommand{\mod}{{\on{-mod}}}
\newcommand{\comod}{{\on{-comod}}}

\nc\mathi\iota
\nc\Spec{\on{Spec}}
\nc\Proj{\on{Proj}}
\nc\Mod{\on{Mod}}
\nc{\tw}{\widetilde{\mathfrak t}}
\nc{\pw}{\widetilde{\mathfrak p}}
\nc{\qw}{\widetilde{\mathfrak q}}
\nc{\jw}{\widetilde j}

\nc{\grb}{\overline{\Gr}}
\nc{\I}{\mathcal I}
\renewcommand{\i}{\mathfrak i}

\nc{\lambdach}{{\check\lambda}}
\nc{\Lambdach}{{\check\Lambda}{}}
\nc{\much}{{\check\mu}}
\nc{\omegach}{{\check\omega}}
\nc{\nuch}{{\check\nu}}
\nc{\etach}{{\check\eta}}
\nc{\alphach}{{\check\alpha}}
\nc{\oblvtach}{{\check\oblvta}}
\nc{\rhoch}{{\check\rho}}
\nc{\ch}{{\check h}}

\nc{\Hb}{\overline{\H}}

\nc{\BA}{{\mathbb{A}}}
\nc{\BC}{{\mathbb{C}}}
\nc{\BE}{{\mathbb{E}}}
\nc{\BF}{{\mathbb{F}}}
\nc{\BG}{{\mathbb{G}}}
\nc{\BL}{{\mathbb{L}}}
\nc{\BM}{{\mathbb{M}}}
\nc{\BO}{{\mathbb{O}}}
\nc{\BD}{{\mathbb{D}}}
\nc{\BN}{{\mathbb{N}}}
\nc{\BP}{{\mathbb{P}}}
\nc{\BQ}{{\mathbb{Q}}}
\nc{\BR}{{\mathbb{R}}}
\nc{\BV}{{\mathbb{V}}}
\nc{\BZ}{{\mathbb{Z}}}
\nc{\BS}{{\mathbb{S}}}
\nc{\Deep}{{\bf{deep}}}
\nc{\deep}{deep}

\nc{\CA}{{\mathcal{A}}}
\nc{\CB}{{\mathcal{B}}}

\nc{\CE}{{\mathcal{E}}}
\nc{\CF}{{\mathcal{F}}}
\nc{\CH}{{\mathcal{H}}}

\nc{\CL}{{\mathcal{L}}}
\nc{\CC}{{\mathcal{C}}}
\nc{\CG}{{\mathcal{G}}}
\nc{\CalD}{{\mathcal{D}}}
\nc{\CM}{{\mathcal{M}}}
\nc{\CN}{{\mathcal{N}}}
\nc{\CK}{{\mathcal{K}}}
\nc{\CO}{{\mathcal{O}}}
\nc{\CP}{{\mathcal{P}}}
\nc{\CQ}{{\mathcal{Q}}}
\nc{\CR}{{\mathcal{R}}}
\nc{\CS}{{\mathcal{S}}}
\nc{\CT}{{\mathcal{T}}}
\nc{\CU}{{\mathcal{U}}}
\nc{\CV}{{\mathcal{V}}}
\nc{\CW}{{\mathcal{W}}}
\nc{\CX}{{\mathcal{X}}}
\nc{\CY}{{\mathcal{Y}}}
\nc{\CZ}{{\mathcal{Z}}}
\nc{\CI}{{\mathcal{I}}}

\nc{\csM}{{\check{\mathcal A}}{}}
\nc{\oM}{{\overset{\circ}{\mathcal M}}{}}
\nc{\obM}{{\overset{\circ}{\mathbf M}}{}}
\nc{\oCA}{{\overset{\circ}{\mathcal A}}{}}
\nc{\obA}{{\overset{\circ}{\mathbf A}}{}}
\nc{\ooM}{{\overset{\circ}{M}}{}}
\nc{\osM}{{\overset{\circ}{\mathsf M}}{}}
\nc{\vM}{{\overset{\bullet}{\mathcal M}}{}}
\nc{\nM}{{\underset{\bullet}{\mathcal M}}{}}
\nc{\oD}{{\overset{\circ}{\mathcal D}}{}}
\nc{\obD}{{\overset{\circ}{\mathbf D}}{}}
\nc{\oA}{{\overset{\circ}{A}}{}}
\nc{\op}{{\overset{\bullet}{\mathbf p}}{}}
\nc{\cp}{{\overset{\circ}{\mathbf p}}{}}
\nc{\oU}{{\overset{\bullet}{\mathcal U}}{}}
\nc{\oZ}{{\overset{\circ}{\mathcal Z}}{}}
\nc{\ofZ}{{\overset{\circ}{\mathfrak Z}}{}}
\nc{\oF}{{\overset{\circ}{\fF}}}

\nc{\fa}{{\mathfrak{a}}}
\nc{\ofa}{\overset{\circ}{\mathfrak{a}}}
\nc{\fb}{{\mathfrak{b}}}
\nc{\fd}{{\mathfrak{d}}}
\nc{\ff}{{\mathfrak{f}}}
\nc{\fg}{{\mathfrak{g}}}
\nc{\fgl}{{\mathfrak{gl}}}
\nc{\fh}{{\mathfrak{h}}}
\nc{\fj}{{\mathfrak{j}}}
\nc{\fl}{{\mathfrak{l}}}
\nc{\fm}{{\mathfrak{m}}}
\nc{\ofm}{\overset{\circ}{\mathfrak{m}}}
\nc{\fn}{{\mathfrak{n}}}
\nc{\fu}{{\mathfrak{u}}}
\nc{\fp}{{\mathfrak{p}}}
\nc{\fr}{{\mathfrak{r}}}
\nc{\fs}{{\mathfrak{s}}}
\nc{\ft}{{\mathfrak{t}}}
\nc{\oft}{\overset{\circ}{\mathfrak{t}}}
\nc{\fz}{{\mathfrak{z}}}
\nc{\fsl}{{\mathfrak{sl}}}
\nc{\hsl}{{\widehat{\mathfrak{sl}}}}
\nc{\hgl}{{\widehat{\mathfrak{gl}}}}
\nc{\hg}{{\widehat{\mathfrak{g}}}}
\nc{\hm}{{\widehat{\mathfrak{m}}}}
\nc{\chg}{{\widehat{\mathfrak{g}}}{}^\vee}
\nc{\hn}{{\widehat{\mathfrak{n}}}}
\nc{\chn}{{\widehat{\mathfrak{n}}}{}^\vee}

\nc{\fA}{{\mathfrak{A}}}
\nc{\fB}{{\mathfrak{B}}}
\nc{\fD}{{\mathfrak{D}}}
\nc{\fE}{{\mathfrak{E}}}
\nc{\fF}{{\mathfrak{F}}}
\nc{\fG}{{\mathfrak{G}}}
\nc{\fK}{{\mathfrak{K}}}
\nc{\fL}{{\mathfrak{L}}}
\nc{\fM}{{\mathfrak{M}}}
\nc{\fN}{{\mathfrak{N}}}
\nc{\fP}{{\mathfrak{P}}}
\nc{\fU}{{\mathfrak{U}}}
\nc{\fV}{{\mathfrak{V}}}
\nc{\fZ}{{\mathfrak{Z}}}

\nc{\ba}{{\mathbf{a}}}
\nc{\bb}{{\mathbf{b}}}
\nc{\bc}{{\mathbf{c}}}
\nc{\bd}{{\mathbf{d}}}
\nc{\bbf}{{\mathbf{f}}}
\nc{\be}{{\mathbf{e}}}
\nc{\bi}{{\mathbf{i}}}
\nc{\bj}{{\mathbf{j}}}
\nc{\bh}{{\mathbf{h}}}
\nc{\bm}{{\mathbf{m}}}
\nc{\bn}{{\mathbf{n}}}
\nc{\bo}{{\mathbf{o}}}
\nc{\bp}{{\mathbf{p}}}
\nc{\bq}{{\mathbf{q}}}
\nc{\bu}{{\mathbf{u}}}
\nc{\bv}{{\mathbf{v}}}
\nc{\bx}{{\mathbf{x}}}
\nc{\bs}{{\mathbf{s}}}
\nc{\by}{{\mathbf{y}}}
\nc{\bw}{{\mathbf{w}}}
\nc{\bA}{{\mathbf{A}}}
\nc{\bK}{{\mathbf{K}}}
\nc{\bB}{{\mathbf{B}}}
\nc{\bC}{{\mathbf{C}}}
\nc{\bG}{{\mathbf{G}}}
\nc{\bD}{{\mathbf{D}}}
\nc{\bE}{{\mathbf{E}}}
\nc{\bH}{{{\mathbf{H}}}}
\nc{\bL}{{\mathbf{L}}}
\nc{\bM}{{\mathbf{M}}}
\nc{\bN}{{\mathbf{N}}}
\nc{\bO}{{\mathbf{O}}}
\nc{\bQ}{{\mathbf{Q}}}
\nc{\bV}{{\mathbf{V}}}
\nc{\bW}{{\mathbf{W}}}
\nc{\bX}{{\mathbf{X}}}
\nc{\bZ}{{\mathbf{Z}}}
\nc{\bS}{{\mathbf{S}}}

\nc{\sA}{{\mathsf{A}}}
\nc{\sB}{{\mathsf{B}}}
\nc{\sC}{{\mathsf{C}}}
\nc{\sD}{{\mathsf{D}}}
\nc{\sF}{{\mathsf{F}}}
\nc{\sG}{{\mathsf{G}}}
\nc{\sH}{{\mathsf{H}}}
\nc{\sK}{{\mathsf{K}}}
\nc{\sM}{{\mathsf{M}}}
\nc{\sN}{{\mathsf{N}}}
\nc{\sO}{{\mathsf{O}}}
\nc{\sW}{{\mathsf{W}}}
\nc{\sQ}{{\mathsf{Q}}}
\nc{\sP}{{\mathsf{P}}}
\nc{\sR}{{\mathsf{R}}}
\nc{\sT}{{\mathsf{T}}}
\nc{\sZ}{{\mathsf{Z}}}
\nc{\sfi}{{\mathsf{i}}}
\nc{\sfj}{{\mathsf{j}}}
\nc{\sfp}{{\mathsf{p}}}
\nc{\sfq}{{\mathsf{q}}}
\nc{\sfs}{{\mathsf{s}}}
\nc{\sft}{{\mathsf{t}}}
\nc{\sr}{{\mathsf{r}}}
\nc{\bk}{{\mathsf{k}}}
\nc{\sa}{{\mathsf{s}}}
\nc{\sg}{{\mathsf{g}}}
\nc{\sn}{{\mathsf{n}}}
\nc{\so}{{\mathsf{ob}}}
\nc{\sh}{{\mathsf{h}}}
\nc{\sff}{{\mathsf{f}}}
\nc{\sfb}{{\mathsf{b}}}
\nc{\sfc}{{\mathsf{c}}}
\nc{\sfe}{{\mathsf{e}}}
\nc{\sd}{{\mathsf{d}}}

\nc{\BK}{{\bar{K}}}

\nc{\tA}{{\widetilde{\mathbf{A}}}}
\nc{\tB}{{\widetilde{\mathcal{B}}}}
\nc{\tg}{{\widetilde{\mathfrak{g}}}}
\nc{\tG}{{\widetilde{G}}}
\nc{\TM}{{\widetilde{\mathbb{M}}}{}}
\nc{\tO}{{\widetilde{\mathsf{O}}}{}}
\nc{\tU}{{\widetilde{\mathfrak{U}}}{}}
\nc{\TZ}{{\tilde{Z}}}
\nc{\tx}{{\tilde{x}}}
\nc{\tbv}{{\tilde{\bv}}}
\nc{\tfP}{{\widetilde{\mathfrak{P}}}{}}
\nc{\tz}{{\tilde{\zeta}}}
\nc{\tmu}{{\tilde{\mu}}}

\nc{\urho}{\underline{\rho}}
\nc{\uB}{\underline{B}}
\nc{\uC}{{\underline{\mathbb{C}}}}
\nc{\ui}{\underline{i}}
\nc{\uj}{\underline{j}}
\nc{\ofP}{{\overline{\mathfrak{P}}}}
\nc{\oB}{{\overline{\mathcal{B}}}}
\nc{\og}{{\overline{\mathfrak{g}}}}
\nc{\oI}{{\overline{I}}}

\nc{\eps}{\varepsilon}
\nc{\hrho}{{\hat{\rho}}}

\nc{\one}{{\mathbf{1}}}
\nc{\two}{{\mathbf{t}}}

\nc{\Rep}{{\mathop{\operatorname{\rm Rep}}}}
\nc{\Tot}{{\mathop{\operatorname{\rm Tot}}}}
\nc{\Ker}{{\mathop{\operatorname{\rm Ker}}}}
\nc{\im}{{\mathop{\operatorname{\rm Im}}}}
\nc{\Hilb}{{\mathop{\operatorname{\rm Hilb}}}}
\nc{\End}{{\mathop{\operatorname{\rm End}}}}
\nc{\Ext}{{\mathop{\operatorname{\rm Ext}}}}
\nc{\CHom}{{\mathop{\operatorname{{\mathcal{H}}\it om}}}}
\nc{\CEnd}{{\mathop{\operatorname{{\mathcal{E}}\it nd}}}}
\nc{\GL}{{\mathop{\operatorname{\rm GL}}}}
\nc{\gr}{{\mathop{\operatorname{\rm gr}}}}
\nc{\HN}{{\mathop{\operatorname{\rm HN}}}}
\nc{\Id}{{\mathop{\operatorname{\rm Id}}}}
\nc{\de}{{\mathop{\operatorname{\rm def}}}}
\nc{\length}{{\mathop{\operatorname{\rm length}}}}
\nc{\supp}{{\mathop{\operatorname{\rm supp}}}}

\nc{\Cliff}{{\mathsf{Cliff}}}
\nc{\Fl}{\on{Fl}}
\nc{\Fib}{{\mathsf{Fib}}}
\nc{\Coh}{{\on{Coh}}}
\nc{\QCoh}{{\on{QCoh}}}
\nc{\IndCoh}{{\on{IndCoh}}}
\nc{\FCoh}{{\mathsf{FCoh}}}

\nc{\reg}{{\on{reg}}}
\nc{\mer}{{\on{mer}}}
\nc{\mf}{{\on{mon-free}}}

\nc{\cplus}{{\mathbf{C}_+}}
\nc{\cminus}{{\mathbf{C}_-}}
\nc{\cthree}{{\mathbf{C}_\bullet}}
\nc{\Qbar}{{\bar{Q}}}
\nc\Eis{\on{Eis}}
\nc\Eisb{\ol\Eis{}}
\nc\Eisr{\on{Eis}^{rat}{}}
\nc\wh{\widehat}
\nc{\Def}{\on{Def_{\check{\fb}}(E)}}
\nc{\barZ}{\overline{Z}{}}
\nc{\barbarZ}{\overline{\barZ}{}}
\nc{\barpi}{\overline\pi}
\nc{\barbarpi}{\overline\barpi}
\nc{\barpip}{\overline\pi{}^+}
\nc{\barpim}{\overline\pi{}^-}

\nc{\fq}{\mathfrak q}

\nc{\fqb}{\ol{\sfq}{}}
\nc{\fpb}{\ol{\sfp}{}}
\nc{\fpr}{{\mathsf{pair}^{rat}}{}}
\nc{\fqr}{{\sfq^{rat}}{}}

\nc{\hattimes}{\wh\otimes}

\nc{\bOmega}{{\overline{\Omega(\check \fn)}}}

\nc{\seq}[1]{\stackrel{#1}{\sim}}

\nc{\cT}{{\check{T}}}
\nc{\cG}{{\check{G}}}
\nc{\cM}{{\check{M}}}
\nc{\cB}{{\check{B}}}
\nc{\cP}{{\check{P}}}

\nc{\ct}{{\check{\mathfrak t}}}
\nc{\cg}{{\check{\fg}}}
\nc{\cb}{{\check{\fb}}}
\nc{\cn}{{\check{\fn}}}

\nc{\cLambda}{{\check\Lambda}}

\nc{\cla}{{\check\lambda}}
\nc{\cmu}{{\check\mu}}
\nc{\cnu}{{\check\nu}}
\nc{\ceta}{{\check\eta}}

\nc{\DefbE}{{\on{Def}_{\cB}(E_\cT)}}

\nc{\imathb}{{\ol{\imath}}}
\nc{\rlr}{\overset{\longrightarrow}{\underset{\longrightarrow}\longleftarrow}}

\nc{\oBun}{\overset{\circ}\Bun}
\nc{\LS}{\on{LS}}
\nc{\BunBbb}{\ol{\ol{Bun}}_B}
\nc{\BunBr}{\Bun_B^{rat}}
\nc{\BunBrsg}{\Bun_B^{rat,\on{s.g.}}}
\nc{\BunBrp}{\Bun_B^{rat,polar}}
\nc{\BunBrpbg}{\Bun_B^{rat,polar,\on{b.g.}}}
\nc{\BunBrpsg}{\Bun_B^{rat,polar,\on{s.g.}}}
\nc{\BunTrp}{\Bun_T^{rat,polar}}
\nc{\BunTrpbg}{\Bun_T^{rat,polar,\on{b.g.}}}
\nc{\BunTrpsg}{\Bun_T^{rat,polar,\on{s.g.}}}
\nc{\BunNr}{\Bun_N^{rat}}
\nc{\BunNre}{\Bun_N^{enh,rat}}
\nc{\BunTr}{\Bun_T^{rat}}
\nc{\Vect}{\on{Vect}}
\nc{\Whit}{{\on{Whit}}}
\nc{\CTb}{\ol{\on{CT}}}
\nc{\Ran}{{\on{Ran}}}
\nc{\CTr}{\on{CT}^{rat}{}}
\nc\jmathr{\jmath^{rat}{}}
\nc{\ux}{\underline{x}}
\nc{\clambda}{{\check\lambda}}
\nc{\calpha}{{\check\alpha}}
\nc{\ind}{{\mathbf{ind}}}
\nc{\coinv}{{\mathbf{coinv}}}
\nc{\oblv}{{\mathbf{oblv}}}
\nc{\free}{{\mathbf{free}}}
\nc{\ox}{{\overline{x}}}
\nc{\cLa}{\check{\Lambda}}
\nc{\StinftyCat}{\on{DGCat}}
\nc{\inftyCat}{\infty\on{-Cat}}
\nc{\inftygroup}{\infty\on{-Grpd}}
\nc{\Dmod}{\on{D-mod}}
\nc{\CMaps}{{\mathcal Maps}}
\nc{\Maps}{\on{Maps}}
\nc{\Sch}{\on{Sch}}
\nc{\affSch}{\on{Sch}^{\on{aff}}}
\nc{\IndSch}{\on{indSch}}
\nc{\dr}{{\on{dR}}}
\nc{\oCF}{\overset{\circ}\CF}
\nc{\oCY}{\overset{\circ}\CY}
\nc{\opi}{\overset{\circ}\pi}
\nc{\leqG}{\underset{G}\leq}
\nc{\leqM}{\underset{M}\leq}
\nc{\leqGad}{\underset{G_{ad}}\leq}
\nc{\leqMad}{\underset{M_{ad}}\leq}
\nc{\Tr}{\on{Tr}}
\nc{\Frob}{{\on{Frob}}}
\nc{\DGCat}{\on{DGCat}}
\nc{\tDGCat}{2\on{-DGCat}_{\on{u.g.}}}
\nc{\ev}{\on{ev}}
\nc{\mmod}{\!\on{-}\!\mathbf{mod}}
\nc{\commod}{\on{-}\mathbf{comod}}
\nc{\sotimes}{\overset{!}\otimes}
\nc{\arrowtimes}{\overset{\to}\otimes}
\nc{\Shv}{\on{Shv}}
\nc{\Spc}{\on{Spc}}
\nc{\Res}{\on{Res}}
\nc{\bDelta}{{\mathbf{\Delta}}}
\nc{\bMaps}{{\mathbf{Maps}}}
\nc{\cD}{\mathcal D}
\nc{\ocD}{\overset{\circ}\cD}
\nc{\ppart}{(\!(t)\!)}
\nc{\qqart}{[\![t]\!]}
\nc{\oCU}{\overset{\circ}{\CU}}
\nc{\Exc}{{\mathcal{E}xc}}
\nc{\Sht}{\on{Sht}}
\nc{\Nilp}{{\on{Nilp}}}
\nc{\Drinf}{\on{Drinf}}
\nc{\Sing}{\on{Sing}}
\nc{\IndLisse}{\Lisse}
\nc{\Shvl}{\on{Shv}_{\on{lisse}}} 
\nc{\Lisse}{\on{Lisse}}
\nc{\Mir}{\on{Mir}}
\nc{\fSet}{\on{fSet}}
\nc{\qLisse}{\on{QLisse}}
\nc{\Ev}{\on{Ev}}
\nc{\Sat}{\on{Sat}}
\nc{\Se}{\on{Se}}
\nc{\coSht}{\on{co-Sht}}
\nc{\coCK}{\on{co-}\!\CK}
\nc{\FLE}{\on{FLE}}
\nc{\BRST}{\on{BRST}}
\nc{\KL}{\on{KL}}
\nc{\KM}{\on{KM}}
\nc{\crit}{{\on{crit}}}
\nc{\Op}{{\on{Op}}}
\nc{\MOp}{\on{MOp}}
\nc{\Wak}{\on{Wak}}
\nc{\Av}{\on{Av}}
\nc{\semiinf}{{\frac{\infty}{2}}}
\nc{\DS}{\on{DS}}
\nc{\dR}{{\on{dR}}}
\nc{\Poinc}{{\on{Poinc}}}
\nc{\Vac}{{\on{Vac}}}
\renc{\det}{\on{det}}
\nc{\oG}{\overset{\circ}{G}}
\nc{\fCDO}{\mathfrak{CDO}}
\nc{\CDO}{\on{CDO}}
\nc{\Sect}{\on{Sect}}

\begin{document}


\vskip1cm

\title[Proof of the geometric Langlands conjecture II]{Proof of the geometric Langlands conjecture II: \\
Kac-Moody localization and the FLE}



\author[Arinkin, Beraldo, Chen, Faergeman, Gaitsgory, Lin, Raskin, Rozenblyum]
{D.~Arinkin, D.~Beraldo, L.~Chen, J.~F\ae{}rgeman, \\ D.~Gaitsgory, K.~Lin, S.~Raskin and N.~Rozenblyum} 

\date{\today}


\dedicatory{To Joseph Bernstein}


\date{\today}

\maketitle

\bigskip

\bigskip


\tableofcontents

\section*{Introduction}

\ssec{Overview}

This paper is the second in a series of five that together prove the 
geometric Langlands conjecture. In this paper, we study the interaction
between Kac-Moody localization and the global geometric Langlands functor of
\cite{GLC1}. We do so following the methodology of Beilinson-Drinfeld, using chiral (a.k.a, factorization) 
homology.

\sssec{}

The main result of this paper, which appears in the main body as \thmref{t:Langlands critical compat}, says: 

\begin{thm}\label{t:main intro}
There is a commutative diagram
 $$
 \CD
 \Dmod_{\crit}(\Bun_G)  @>{\BL_G}>>  \IndCoh_{\Nilp}(\LS_{\cG}) \\
 @A{\Loc_{G,\crit}\otimes \fl}AA @AA{\Poinc^{\on{spec}}_{\cG,*}}A \\
 \KL(G)_{\crit,\Ran} @>{\FLE_{G,\crit}}>> \IndCoh^*(\Op^{\mf}_\cG)_\Ran. 
 \endCD
 $$
\end{thm}

The terms appearing in the above diagram warrant further discussion.
We will do so at more length later in this introduction, but here is a brief
synopsis: 

\begin{itemize}

\item $\Dmod_{\crit}(\Bun_G)$ is the category of 
\emph{critically twisted} D-modules on $\Bun_G$, as considered originally
in \cite{BD1};

\item $\IndCoh_{\Nilp}(\LS_{\cG})$ was defined in \cite{AG} as the spectral
category in geometric Langlands;
 
\item The functor $\BL_G$ is the \emph{Langlands functor} as constructed
in \cite{GLC1};

\item $\KL(G)_{\crit,\Ran}$ is the Ran space\footnote{The Ran space of $X$ parameterizes finite
collection of points of $X$.} version of the Kazhdan-Lusztig category at the critical level
(i.e., Kac-Moody modules at the critical level integrable with respect to the arc group $\fL^+(G)$); 

\item $\Op^{\mf}_\cG$ is the factorization space parametrizing
local systems on the formal disc equipped with an oper structure on the punctured disc;

\item $\FLE_{G,\crit}$ is the \emph{fundamental local equivalence} 
at critical level. This equivalence of factorization
categories appears in \thmref{t:critical FLE} and extends the pointwise
equivalence of \cite{FG2}. It is the main theorem of Part I of this paper;

\item $\Loc_{G,\crit}$ is the functor of critical level Kac-Moody localization;

\item The Poincar\'e series functor, denoted functor $\Poinc^{\on{spec}}_{\cG,*}$, is given at 
each finite set $\ul{x} \in \Ran$ of points in $X$ by 
pull-push along the correspondence
\[
\Op^{\mf}_{\cG,\ul{x}} \leftarrow 
\Op^{\mf,\on{glob}}_{\cG,\ul{x}} \rightarrow \LS_{\cG},
\]

\noindent where the middle term parametrizes local systems on the global
curve $X$ equipped with an oper structure on $X-\ul{x}$;

\item $\fl$ is a cohomologically shifted 1-dimensional vector space that
can be ignored at first approximation.
Using notation defined in the paper, it is 
$\fl^{\otimes \frac{1}{2}}_{G,N_{\rho(\omega_X)}}\otimes 
\fl^{\otimes -1}_{N_{\rho(\omega_X)}}[-\delta_{N_{\rho(\omega_X)}}]$.

\end{itemize}

\sssec{}

The above theorem has had folklore status in the subject. Its main 
ingredients were discussed at the 2014 conference 
``Towards the proof of the geometric Langlands conjecture." 

\medskip

However, some of the key technical aspects have not been addressed in the existing literature. 
This is most notably true for the category $\IndCoh^*(\Op^{\mf}_\cG)_\Ran$, i.e., the category
of ind-coherent sheaves on the \emph{Ran space} version of the space of monodromy-free
opers (see \secref{sss:mf intro} below). 

\sssec{}

The role of \thmref{t:main intro} in the geometric Langlands
program is as follows:

\medskip

The functor $\Loc_{G,\crit}$ is not surjective, but neither
is it so far from being surjective (see \cite[Prop. 10.1.6]{Ga1}). 
Therefore, understanding the interaction
between $\BL_G$ and Kac-Moody localization 
plays a crucial role in understanding $\Dmod_{\crit}(\Bun_G)$ in spectral terms. 
See \remref{r:spectral action} for an example where this idea is applied.

\sssec{What goes into the proof?}

As indicated above, we first need to define the various categories.

\medskip 

Second, we need to construct the functors appearing in the commutative
diagram. Perhaps the most interesting is $\FLE_{G,\crit}$, which is the
subject of Part I of this paper. We discuss it further below. 

\medskip 

Finally, we need to prove the diagram commutes. Ultimately, we do this by
expressing both circuits in terms of chiral homology for the 
critical level $\CW$-algebra and appealing to the Feigin-Frenkel isomorphism.

\sssec{}

As was mentioned already, this paper builds on the ideas of Beilinson and Drinfeld. 

\medskip

In their seminal works \cite{BD1} and \cite{BD2}, Beilinson and Drinfeld
introduced the theory of \emph{chiral algebras} -- which are equivalent to
\emph{factorization algebras} and, suitably understood, vertex algebras --
and of \emph{chiral homology} as a tool for studying 
interactions between categories of \emph{local} nature,
such as sheaves on the affine Grassmannian, and categories of
\emph{global} nature, such as sheaves on $\Bun_G$.

\medskip

The functors appearing in \thmref{t:main intro} are of local-to-global nature, 
and may be viewed as generalizations of the functor of chiral homology. It
is in this sense that one can view the present work as a continuation of
\cite{BD1,BD2}. 

\sssec{}

In writing this text, we found that we needed to refine many foundational
parts of the original work of Beilinson-Drinfeld. This ultimately accounts for 
the length of the present work. 

\medskip

A significant part of these refinements has to do with the fact that we (have to) work
with $\infty$-categories (whereas in \cite{BD1,BD2} one mostly works with abelian categories). 

\ssec{What is done in this paper?}

We now highlight what we think are the most important contributions of this paper.

\sssec{Monodromy-free opers} \label{sss:mf intro}

First, as an algebro-geometric object, $\Op^{\mf}_{\cG}$ parametrizes
a point $\ul{x} \in \Ran$, a local system $\sigma$ on the formal disc $\cD_{\ul{x}}$
at $x$, and an oper structure on the restriction $\sigma|_{\cD^{\times}_{\ul{x}}}$
of $\sigma$ to the punctured disc.

\medskip

When we work over a fixed point $x\in X$, the corresponding space $\Op^{\mf}_{\cG,x}$
was introduced and studied in \cite{FG2}. However, the Ran space version presents a host
of new challenges. 

\medskip

This space has infinite type, so it is not immediately obvious how
to define (ind-)coherent sheaves on it. We explain the relevant geometry needed
to make sense of $\IndCoh^*(\Op^{\mf}_{\cG})$
in \secref{s:opers}. 

\medskip

In \secref{ss:indcoh via fact}, we 
show that $\IndCoh^*(\Op^{\mf}_{\cG})$ can \emph{almost} be realized
as a category of factorization modules. More precisely, in 
\secref{sss:R G Op} 
we define a factorization algebra $R_{\cG,\Op} \in \Rep(\cG)$ and
prove in \propref{p:IndCoh Op via fact almost} that its category of factorization modules
is equivalent to $\IndCoh^*(\Op^{\mf}_{\cG})$ modulo homological convergence issues
(more precisely: the corresponding bounded below categories are equivalent). 

\sssec{The critical FLE}

This result appears as \thmref{t:critical FLE}. It asserts that we have
a t-exact equivalence of factorization categories
\[
\FLE_{G,\crit}:\KL(G)_{\crit} \to \IndCoh^*(\Op^{\mf}_\cG).
\]

\medskip

One can view this equivalence as a (substantially amplified) categorical incarnation 
of the Feigin-Frenkel isomorphism. 

\medskip

The idea of the proof is as follows:

\medskip 

First, we construct the 
functor $\FLE_{G,\crit}$. The ingredients are Feigin's Drinfeld-Sokolov
functor and Beilinson-Drinfeld's birth of opers construction.

\medskip 

Second, we prove that $\FLE_{G,\crit}$ preserves compact objects
(in the sense suitable for factorization categories). This expresses a
finiteness property of Drinfeld-Sokolov reductions that is not 
immediate using classical VOA methods; 
we show that it is immediate from the categorical construction of
$\CW$-algebras from \cite{Ra2}
(i.e., the \emph{affine Skryabin theorem}).

\medskip 

Thanks to the preservation of compactness mentioned above, we are reduced to 
proving that $\FLE_{G,\crit}$ is a \emph{pointwise} equivalence.
This is a theorem of \cite{FG2}. We actually reprove this theorem 
here to illustrate a more modern point of view on 
studying Kac-Moody representations using categorical tools. 

\medskip

Namely, we show that in general, for a category $\bC$ with an 
$\fL(G)$-action, the \emph{tempered quotient}\footnote{See \secref{ss:temp}, where this notion is defined.} 
$\Sph(\bC)_{\on{temp}}$ of
$\Sph(\bC) := \bC^{\fL^+(G)}$ can be algorithmically recovered
from $\Whit(\bC)$, the Whittaker model (i.e., $\fL(N)$-invariants against a non-degenerate
character) of $\bC$. Heuristically, $\Whit(\bC)$ should live as a sheaf over 
$\LS_{\cG}(\cD^{\times})$ and its sections over $\LS_{\cG}(\cD)$ should
recover $\Sph(\bC)_{\on{temp}}$; we give a precise assertion of this
type in \propref{p:Sph to unr Whit}. The key input for this result is the pointwise version of 
derived Satake.

\medskip

We then show that $\KL(G)_{\crit} = \Sph(\hg\mod_{\crit,x})$ equals its
tempered quotient. However, by \cite{Ra2}, 
$$\Whit(\hg\mod_{\crit,x}) \simeq \IndCoh^*(\Op^\mer_\cG),$$
and hence we obtain the FLE from the previous paragraph.

\medskip

To summarize: we deduce the FLE at critical level as an 
essentially formal consequence of derived Satake and affine Skryabin.

\begin{rem}

We should add that a particular case of the pointwise abelian category
version of the FLE was established already in \cite{BD1}:

\medskip

Namely, in {\it loc. cit.} it was shown that the subcategory of $(\KL(G)_{\crit,x})^\heartsuit$
consisting of modules with \emph{regular central characters} is freely generated
over $(\QCoh(\Op^\reg_{\cG,x}))^\heartsuit$ by the vacuum module. 

\medskip

Note, however, that a parallel statement would be \emph{false} at the derived level;
this observation is what led the authors of \cite{FG2} to considering the ind-scheme
of monodromy-free opers,

\end{rem}

\sssec{The FLE and duality}

The category $\KL(G)_{\crit}$ is canonically self-dual by a construction
with semi-infinite cohomology, see \secref{sss:KL self-duality crit}.

\medskip 

In \secref{s:opers}, we show that $\IndCoh^*(\Op^{\mf}_\cG)$ 
is canonically self-dual,
which we express as an equivalence
$$\Theta_{\Op^\mf_\cG}:\IndCoh^!(\Op^{\mf}_{\cG}) \simeq \IndCoh^*(\Op^{\mf}_{\cG}).$$
This equivalence comes from a similar equivalence 
$\Theta_{\Op^\mer_\cG}$ using all opers
in place of monodromy-free opers;
the latter should be thought of as a critical limit of
the semi-infinite cohomology for $\CW$-algebras considered
by Dhillon in \cite{Dh}. 

\medskip 

In \secref{s:FLE and duality}, we prove that these two self-duality 
constructions match under the FLE. As indicated above, this 
result should be considered as a compatibility
between the FLE and two flavors of semi-infinite cohomology. 

\sssec{The formalism of local-to-global functors}

We develop axiomatics in \secref{s:unitality}. In some part, the constructions
here abstract Beilinson-Drinfeld's construction of chiral homology.

\medskip 

There is a separate introduction to this material in Sect. 12.0, so we 
describe the material briefly:

\medskip 

One often finds the following situation: there is a local (factorization)
category $\ul\bC^{\on{loc}}$, a global category $\bC^{\on{glob}}$, and 
a \emph{local-to-global} functor $$\sF:\bC^{\on{loc}}_{\Ran} \to \bC^{\on{glob}}.$$

Here are examples we have in mind:

\begin{itemize}

\item For a chiral algebra $\CA$, take\footnote{Technically,
$\CA\mod^{\on{fact}}$ only forms a \emph{lax} factorization category.
In fact, the material of \secref{s:unitality} does not assume any sort
of factorization, just the existence of suitable categories over the
(unital) Ran space.}
 $\ul\bC^{\on{loc}} = \CA\mod^{\on{fact}}$, $\bC^{\on{glob}} = \Vect$,
 and $\sF = \ul{\on{C}}^{\on{fact}}_\cdot(X,\CA,-)$, i.e., the functor of
 chiral homology;
 
 \medskip
 
\item Take $\ul\bC^{\on{loc}} = \KL(G)_{\kappa}$, $\bC^{\on{glob}} = 
\Dmod_{\kappa}(\Bun_G)$, and $\sF = \Loc_{G,\kappa}$;

\medskip

\item Take $\ul\bC^{\on{loc}} = \Whit_{\kappa}(G)$ the category of 
$\kappa$-twisted Whittaker D-modules on the affine Grassmannian,   
$\bC^{\on{glob}} = \Dmod_{\kappa}(\Bun_G)$, and 
$\sF = \Poinc_{G,!}$ (or $\Poinc_{G,*}$);

\medskip

\item Take $\ul\bC^{\on{loc}} = \Rep(\cG)$,    
$\bC^{\on{glob}} = \QCoh(\LS_{\cG})$, and 
$\sF = \Loc^{\on{spec}}$ the spectral localization functor;

\medskip

\item Take $\ul\bC^{\on{loc}} = \IndCoh^*(\Op^{\mf}_{\cG})$, 
$\bC^{\on{glob}} = \IndCoh(\LS_{\cG})$, and $\sF = \ul{\Poinc}^{\on{spec}}_{\cG,*}$
(or $\ul{\Poinc}^{\on{spec}}_{\cG,!}$).

\end{itemize}

The last two examples can be considered as $\kappa\to\infty$ limits of 
the second and third examples. 

\medskip 

One key feature of each of the above constructions is that they 
are \emph{unital}, in the sense that our factorization categories
are themselves unital and $\sF$ (canonically) commutes with vacuum insertion.

\medskip 

The following construction plays a key role: 

\medskip

One starts with a non-unital (but \emph{lax} unital) functor
$\sF_0$ and it turns out that there exists a procedure that canonically produces from it 
a \emph{strictly unital} local-to-global
functor $\sF$.

\medskip 

For example, for a chiral algebra
$\CA$, we could take $\sF_0$ as 
$$\CA\mod_{\Ran}^{\on{fact}} \xrightarrow{\oblv} \Dmod(\Ran) 
\xrightarrow{\on{C}^\cdot_c(\Ran,-)} \Vect$$ with second functor that of compactly supported 
de Rham cochains. The resulting functor $\sF$ is that of chiral homology. I.e., chiral homology
can be thought of as a universal procedure that forces 
$\sF_0$ to commute with vacuum insertion.

\medskip 

We discuss this passage from lax to strict unital globlization functors 
at length in \secref{s:unitality}. 

\begin{rem} 
In favorable cases, Betti analogues of local-to-global functors have
TQFT interpretations. Namely, given a 4d TQFT $Z$ and a boundary condition 
$\CB$ for it, we obtain $\bC^{\on{glob}}$ as $Z(X)$;
$\ul\bC^{\on{loc}}$ as the evaluation of $Z$ on a closed disc $\cD^{\on{Betti}}$,
putting $\CB$ on the boundary $\partial \cD^{\on{Betti}}$;
and for $\ul{x} \in \Ran$, we suture the disc $\cD_{\ul{x}}^{\on{Betti}}$
around $\ul{x}$ into $X \setminus (\cD_{\ul{x}}^{\on{Betti}} \setminus \partial 
\cD_{\ul{x}}^{\on{Betti}})$ to obtain $\sF$. 
See \remref{r:unitality qft} for a related discussion.
\end{rem} 

\sssec{Localization}

In \secref{s:loc}, we construct and study the Kac-Moody localization functor,
which appears in \thmref{t:main intro}. We do this in a loop group equivariant
way, which has the effect of making $\Loc_G$ respect the Hecke actions. 

\medskip 

We reproduce some results from \cite{CF}. 
For the purposes of \thmref{t:main intro}, the
most important outcome is \thmref{t:int loc over BunN kappa chi}, which
says we have a commutative diagram
$$
\CD
\Dmod_{\kappa}(\Bun_G)  @>{\on{coeff}_G^{\on{Vac,glob}}}>> \Vect \\
@A{\Loc_{G,\kappa} \otimes \fl_{N_{\rho(\omega_X)}}[\delta_{N_{\rho(\omega_X)}}{]}}AA @AA{\ul{\on{C}}^{\on{fact}}_{\cdot}(X{,}\CW{,}-)}A \\
\KL(G)_{\kappa,\Ran} @>{\DS^{\on{enh}}}>> \CW_{\kappa}\mod^{\on{fact}}. 
\endCD
$$
%

\noindent Here the bottom horizontal arrow is Drinfeld-Sokolov reduction, the right vertical 
arrow is chiral homology, and the top horizontal arrow is the functor of
\emph{vacuum Whittaker coefficient}.

\medskip

The similarity between this result and \thmref{t:main intro} is plain.
In fact, the Langlands functor $\BL_G$ is \emph{characterized} using
$\on{coeff}_G^{\on{Vac,glob}}$, so this result is quite close to 
\thmref{t:main intro}.

\sssec{Localization at the critical level and Hecke eigensheaves}\label{sss:crit loc hecke intro}

In \secref{s:Hecke Loc}, we prove the \emph{Hecke eigenproperty} of 
localization at the critical level. This result extends 
one of the main theorems 
of \cite{BD1}; there Beilinson-Drinfeld considered the vacuum representation,
but our result allows consideration of arbitrary objects of $\KL(G)_{\crit}$.

\medskip 

More precisely, Feigin-Frenkel duality (or the FLE) allows us to consider 
$\KL(G)_{\crit,\Ran}$ as a module category for $\QCoh(\Op^{\mf}_{\cG})_{\Ran}$.

\medskip 

Now let $\Op^{\mf,\on{glob}}_{\cG}$ be the space over $\Ran$ parametrizing
$\ul{x} \in \Ran$, a $\cG$-local system $\sigma$ on $X$, 
and an oper structure on $\sigma|_{X\setminus \ul{x}}$. There is 
an evident map $\Op^{\mf,\on{glob}}_{\cG} \to \Op^{\mf}_{\cG}$.

\medskip 

We prove in \corref{c:Hecke crit Loc ult} 
that $\Loc_{G,\crit}$ factors through a functor
\[
\Loc_{G,\crit}^{\Op}:
\KL(G)_{\crit,\Ran} \underset{\QCoh(\Op^{\mf}_{\cG})}{\otimes} 
\QCoh(\Op^{\mf,\on{glob}}_{\cG}) \to \Dmod_{\crit}(\Bun_G) 
\]

\noindent that is $\Rep(\cG)_{\Ran}$-linear with respect to the following
actions:

\begin{itemize}

\item $\Rep(\cG)_{\Ran}$ acts on $\Dmod_{\crit}(\Bun_G)$ through the Hecke action. 

\item $\Rep(\cG)_{\Ran}$ acts on 
$$\KL(G)_{\crit,\Ran} \underset{\QCoh(\Op^{\mf}_{\cG})_{\Ran}}{\otimes} 
\QCoh(\Op^{\mf,\on{glob}}_{\cG})_{\Ran}$$
by 
$$\Rep(\cG)_{\Ran} \overset{\Loc^{\on{spec}}_\cG}\longrightarrow \QCoh(\LS_\cG)$$ and pullback along
the tautological map $\Op^{\mf,\on{glob}}_{\cG} \to \LS_\cG$.

\end{itemize}

Applying the FLE, we can rewrite $\Loc_{G,\crit}^{\Op}$
as a functor
\[
\IndCoh^*(\Op^{\mf}_{\cG})_{\Ran} \underset{\QCoh(\Op^{\mf}_{\cG})}{\otimes} 
\QCoh(\Op^{\mf,\on{glob}}_{\cG}) \to \Dmod_{\crit}(\Bun_G).
\]

\noindent We have canonical functors 
$$\QCoh(\Op^{\mf}_{\cG})_{\Ran} \longrightarrow
\IndCoh^!(\Op^{\mf}_{\cG})_{\Ran} \xrightarrow{\Theta_{\Op^{\mf}_{\cG}}} 
\IndCoh^*(\Op^{\mf}_{\cG})_{\Ran}$$ that we can compose with the above
to obtain a functor
\[
\QCoh(\Op^{\mf,\on{glob}}_{\cG}) \to \Dmod_{\crit}(\Bun_G).
\]

Our Hecke property implies that this functor sends
the skyscraper sheaf at an oper $\chi \in \Op^{\mf,\on{glob}}_{\cG}$ to 
an eigensheaf for the local system underlying $\chi$. 
When $\chi$ is a \emph{regular} oper on $X$, this is the main construction
of \cite{BD1}.\footnote{Actually, even in this case, our notion of eigensheaf
is somewhat more homotopically robust than the one from \cite{BD1}.}

\begin{rem}

Beilinson-Drinfeld show that their eigensheaves are non-zero by computing
their characteristic cycles. This does not directly apply in the monodromy-free
setting, but one can use \thmref{t:main intro} to verify that they are non-zero
by calculating the Whittaker coefficients of these localized D-modules.

\end{rem}

\begin{rem}\label{r:spectral action}

The above eigen-property for Kac-Moody localization was used in 
Drinfeld-Gaitsgory's proof of the
 \emph{spectral action} of $\QCoh(\LS_{\cG})$ on 
$\Dmod_{\crit}(\Bun_G)$, cf. \cite{Ga1} Theorem 4.5.2. 
The proof presented there is based on Kac-Moody localization,
with \cite[Theorem 10.3.4]{Ga1} essentially asserting the Hecke property discussed above. 
In this sense, the present work fills an important gap in the literature.

\end{rem}

\ssec{Structure of this paper}

\sssec{}

Overall, the paper proceeds as follows. In Part I, we formulate and prove the critical level 
FLE, which is the main local result of this paper. We also consider the
interactions of the FLE with various duality functors.

\medskip 

Part II considers the vertical local-to-global functors from 
\thmref{t:main intro} as well as the Hecke actions. Some of our main results
here reproduce results and arguments from \cite{CF}. We put particular emphasis
on the role of \emph{unital} structures, building on ideas from \cite{BD2},
\cite{Ra6}, \cite{Ro2} and \cite{Ga4}.

\medskip 

The proof of \thmref{t:main intro} relies on chiral homology and 
some mild variants thereof. Because \cite{BD2} largely considered 
abelian categories of chiral modules, we need some extensions of their 
ideas to the derived setting; these appear in the appendices to this paper. 

\sssec{}

In more granular detail, the paper is structured as follows. 

\medskip 

Part I deals with the local theory:

\medskip 

\secref{s:Sat} reviews various factorization categories on the geometric side 
associated with the affine Grassmannian of $G$,
as well as their spectral counterparts. 
The key points here are the geometric Casselman-Shalika formula
(\thmref{t:geom CS}) and (derived) geometric Satake equivalence (\thmref{t:geom Satake}). This material is largely taken from \cite{CR}.

\medskip 

In \secref{s:KL} we discuss the Kazhdan-Lusztig category and 
quantum Drinfeld-Sokolov reduction for its modules.

\medskip

In \secref{s:opers}, we consider ind-coherent sheaves on various spaces of
local opers. We reinterpret these objects using factorization
module categories in \secref{ss:indcoh via fact}.

\medskip

In \secref{s:fact alg} we explain how various factorization categories of interest can be 
expressed as factorization module categories over factorization algebras. We use this to
define a \emph{categorical} action of the Feigin-Frenkel center on Kac-Moody modules
at the critical level. 

\medskip 

\secref{s:FF} reviews the Feigin-Frenkel isomorphism at the critical level and 
Beilinson-Drinfeld's \emph{birth of opers}, i.e., the local interplay of
Hecke symmetries and Feigin-Frenkel duality. 

\medskip

In \secref{s:FLE}, we formulate the FLE at critical level. We reduce
the statement to the pointwise assertion, which was proved in 
\cite{FG2}. We also prove some important compatibilities for the FLE here.

\medskip

As discussed above, we also provide another proof of the pointwise
assertion in the present paper, by combining general considerations about derived Satake with 
Feigin and Frenkel's duality
for $\CW$-algebras.

\medskip

In \secref{s:FLE and duality}, we study the interaction between 
the FLE and natural duality functors between the DG categories appearing
in it.

\sssec{}

We now turn to Part II, which deals with local-to-global constructions.  

\medskip 

\secref{s:coeff} reviews the definition of 
Whittaker coefficient functors and establishes some conventions related to them.

\medskip

In \secref{s:loc}, we introduce the Kac-Moody localization functor $\Loc_G$.

\medskip

In \secref{s:unitality}, we discuss axiomatics for local-to-global functors for 
factorization categories (or even just categories over $\Ran$). 
The key construction produces \emph{unital} local-to-global functors from 
\emph{lax unital} such functors, abstracting the construction of chiral homology.

\medskip

\secref{s:prop Loc} considers the interaction between Kac-Moody localization
and restriction/inflation along group homomorphisms $H \to G$. This material
appeared previously in \cite{CF}. 

\medskip

In \secref{s:loc Hecke gpd}, we discuss an alternative construction 
of Kac-Moody localization from the eighth author's thesis. The idea is to
realize D-modules on $\Bun_G$ as quasi-coherent sheaves equipped with
infinitesimal Hecke equivariance structures at every point of $\Ran$.

\medskip

In \secref{s:coeff Loc}, we apply the material from \secref{s:prop Loc} to 
calculate Whittaker coefficients of Kac-Moody representations in terms
of chiral homology for the critical level $\CW$-algebra.

\medskip

\secref{s:Hecke Loc} constructs a Hecke equivariance structure for $\Loc_G$.
This can be considered as an extension of the main
construction of \cite[Sect. 7]{BD1}, allowing more choices of characters for the
Feigin-Frenkel center while adding homotopy coherence to \emph{loc. cit}.  

\medskip
 
\secref{s:ins vac reg} contains the proof of \thmref{t:ins vac reg}, 
a technical point from \secref{s:Hecke Loc}.

\medskip

\secref{s:spectral Poinc} constructs the functor $\Poinc^{\on{spec}}_{\cG,*}$
and its relative $\Poinc^{\on{spec}}_{\cG,!}$.

\medskip

Finally, \secref{s:Langlands functor} considers the interaction 
between the Langlands functor $\BL_G$ and the constructions of
earlier sections. Most importantly, we conclude the proof of \thmref{t:main intro}
here. We also prove that $\BL_G$ is compatible with the factorizable
derived Satake equivalence here, to be used in the sequel to this paper. 

\sssec{}

This paper relies on a lot of foundational material, a big part (but not all) of which has 
not been previously written down. This material is developed in the Appendix\footnote{The
Appendix is coauthored by J.~Campbell, L.~Chen, D.~Gaitsgory, K.~Lin, S.~Raskin and N.~Rozenblyum.} 
to this paper.

\medskip

In Appendix \ref{s:IndCoh inf type} we develop the IndCoh theory for algebro-geometric objects
that are \emph{not} of finite type. As it turns out, there are two versions, denoted $\IndCoh^!(-)$
and $\IndCoh^*(-)$, respectively, which in good situations are mutually dual. We introduce the property
of schemes, called \emph{placidity}, which guarantees that these categories behave particularly
well. In addition, we introduce another category, denoted $\QCoh_{\on{co}}(-)$,
useful in many situations, and which is a \emph{pre-dual} of $\QCoh_{\on{co}}$. 

\medskip

In Appendix \ref{s:fact} we discuss the pattern of factorization. We introduce factorization spaces,
and construct examples of such (e.g., loops or arcs into a given target, or various spaces
attached to the formal disc). We introduce factorization algebras and modules, and various
operations between them. One of the central notions in this paper is that of 
\emph{factorization category}. We show how various categories of algebro-geometric 
or representation-theoretic nature acquire this structure (notably, $\IndCoh^*(-)$ of 
monodromy-free opers and the category of Kac-Moody representations). 

\medskip

In Appendix \ref{s:unit} we discuss the phenomenon of \emph{unitality}. We introduce 
categorical prestacks, D-modules and sheaves of categories on them. Our main example
is the unital Ran space. We introduce unital and counital factorization spaces,
and their common generalization, called ``unital-in-correspondences" factorization spaces; 
it is this latter notion that plays the most important role. We introduce unital factorization 
algebras and categories. We emphasize that some phenomena (such as restriction of 
module categories) work differently in unital and non-unital settings, and it is the
former that are responsible for some of the fundamental constructions in this paper. 

\medskip

In Appendix \ref{s:chiral mods} we prove one of the general fundamental theorems that describe
a category of algebro-geometric nature as modules over a factorization algebra. Namely,
we show that (at the level of bounded below categories), the category 
$\on{QCoh}_{\on{co}}(\fL_\nabla(\CY))$ (here $\CY$ is an affine D-scheme, and $\fL_\nabla(\CY))$
is the space of its horizontal sections on the punctured disc) identifies with factorization modules 
over the (commutative) factorization algebra of regular functions on $\CY$. The equivalence
at the level of abelian categories is nearly evident. However, at the derived level, it is quite
non-trivial, and requires that $\CY$ be \emph{of finite presentation in the D-sense}. 

\medskip

In Appendix \ref{s:spec Sph fact} we explain how to make sense of the \emph{spectral spherical
category}, i.e., the category of ind-coherent sheaves on the local spectral Hecke stack
$\on{Hecke}_\cG^{\on{spec,loc}}$, as a factorization category. 
The problem is that $\on{Hecke}_\cG^{\on{spec,loc}}$ does not
quite fit into the paradigm of \secref{s:IndCoh inf type}, in which we can make sense of
$\IndCoh(-)$ by an algorithmic procedure. Yet, we give an algebro-geometric definition of
$\IndCoh^*(\on{Hecke}_\cG^{\on{spec,loc}})$, and then compare it with a representation-theoretic
one of \cite{CR}. 
\medskip

In Appendix \ref{s:hor sect D sch} we recap (essentially, following \cite{BD2}) the relation between
the scheme of horizontal sections of an affine D-scheme $\CY$ and the factorization homology
of the factorization algebra of regular functions on $\CY$. 

\medskip

In Appendix \ref{s:from LS} we describe a procedure that attaches factorization module categories over $\Rep(\cG)$ to module
categories over $\QCoh(\LS^\mer_\cG)$, and show that this functor is fully faithful on a certain subcategory.

\medskip

In Appendix \ref{s:indep} we recast some of the material from \secref{s:unitality} using
the notion of the ``independent" category, attached to a crystal of categories on
the unital Ran space. We then discuss various notions of action of a factorization
monoidal category on a DG category.

\medskip

Appendix \ref{s:add unit colax} contains some complementary material to \secref{s:unitality}: we
give an interpretation to the functor of the integrated insertion of the unit in terms of left-lax
functors between crystals of categories over the unital Ran space. 

\medskip

Appendix \ref{s:device} is homotopy-theoretic. Here we introduce a device that allows us
to construct us monoidal actions from Sects. \ref{ss:action of center} and \ref{ss:z on KL} 
up to \emph{coherent} homotopy. These monoidal actions play a key role in the definition
of the FLE functor. 

\ssec{Conventions and notation: generalities}

\sssec{The players}

Throughout the paper we work over a fixed algebraically closed field $k$ of characteristic $0$. Thus,
all algebro-geometric objects are defined over $k$. 

\medskip

In particular, $X$ is a smooth projective curve over $k$, $G$ is a reductive group over $k$,
and $\cG$ is the Langlands dual of $G$. 

\sssec{Categories}

When we say ``category", we mean an $\infty$-category. Conventions pertaining to the $\infty$-categorical 
language are borrowed from \cite[Chapter 1, Sect. 1]{GaRo3}.

\sssec{}

Conventions pertaining to DG categories follow those in \cite[Chapter 1, Sect. 10]{GaRo3}. 
Unless explicitly stated otherwise, a DG category $\bC$ is assumed cocomplete (i.e., to contain arbitrary direct sums). 
(An exception would be, e.g., the category of compact objects in a given $\bC$, denoted $\bC^c$.)

\medskip

Unless explicitly stated otherwise, given a pair of DG categories $\bC_1$ and $\bC_2$, by a functor $F:\bC_1\to \bC_2$ 
we will always understand a \emph{continuous} functor, i.e., one that commutes with arbitrary direct sums (equivalently, colimits). 

\sssec{}

Given a DG category $\bC$ with a t-structure, we will use \emph{cohomological} conventions. I.e., $\bC^{\leq 0}$
will denote the subcategory of \emph{connective} objects. We will denote by $\bC^\heartsuit$ the heart of the t-structure. 



\sssec{} 

Conventions adopted in this paper regarding higher algebra and derived algebraic geometry 
follow closely those of \cite{AGKRRV}.

\sssec{Factorization}

Conventions and notation pertaining to the Ran space and \emph{factorization} are explained in 
\secref{s:fact}. 

\medskip

There are several pieces of notation associated with factorization categories: 

\medskip

Given a factorization
category $\bC$, we will denote by $\ul\bC$ the corresponding sheaf of categories over $\Ran$,
by $\bC_\Ran$ its category of global sections, and for $\CZ\to \Ran$ by $\bC_\CZ$ the
category of sections of the pullback of $\ul\bC$ to $\CZ$. In particular, for a $k$-point $\ul{x}\in \Ran$,
we will denote by $\bC_{\ul{x}}$ the fiber of $\ul\bC$ at $\ul{x}$. 

\medskip

Given a pair of factorization categories $\bC_1$ and $\bC_2$ and a functor $\Phi$ between,
we will distinguish between a property of this functor (such as admitting an adjoint or being an
equivalence) taking place at the \emph{pointwise} or \emph{factorization} level.

\medskip

The former means that the given property holds for the corresponding functor
$$\Phi:\bC_{1,\ul{x}}\to \bC_{2,\ul{x}}$$
for any $k$-point $\ul{x}$ of the Ran space. The latter means that the given property holds
for
$$\Phi:\bC_{1,\CZ}\to \bC_{2,\CZ}$$
for any prestack $\CZ\to \Ran$ (equivalently, one can take $\CZ$ to be $\Ran$ itself).

\ssec{Acknowledgements} 

As should be clear from what we said above, the majority of the second part of this paper 
can be traced back to the ideas of A.~Beilinson and V.~Drinfeld recorded 
in \cite{BD1} and \cite{BD2}.

\medskip

The FLE as presented in Part I relies crucially on the Feigin-Frenkel isomorphism, as a passage
between $G$ and $\cG$, see \secref{s:FF}.  

\medskip

A crucial role in local and local-to-global constructions is played by the concept of \emph{factorization}. Its appearance in representation
theory was pioneered by M.~Finkelberg, I.~Mirkovi\'c and V.~Schechtman, and it was further subsequently elucidated by A.~Beilinson and J.~Lurie. 

\medskip

Separate thanks are due to J.~Lurie for enabling representation theorists to work within Higher Algebra. The mathematics
developed in this paper would not be possible if one worked ``up to homotopy". 

\medskip

The fifth and seventh authors wish to thank IH\'ES, where a significant part of this paper was written, for creating an excellent
working environment. 

\medskip

The work of D.G. was supported by NSF grant DMS-2005475. 
The work of S.R. was supported by a Sloan Research Fellowship and NSF grants DMS-2101984 and DMS-2416129 while this work was in preparation.

\newpage

\centerline{\bf Part I: Local Theory}

\bigskip

This Part is mainly dedicated to the proof of a key local result: the critical FLE. It says that the
Kazhdan-Lusztig category at the critical level (for $G$) is equivalent to the category of ind-coherent
sheaves on the space of monodromy-free opers on the punctured disc (for $\cG$). 

\medskip

The FLE involves crossing the Langlands bridge. I.e., at some point, we will need to know something
about the relationship between $G$ and $\cG$. In fact, there are exactly two sources of
such results (as long as we stay at the critical level for $G$ and level $\infty$ for $\cG$):
one is the geometric Casselman-Shalika formula (\thmref{t:geom CS}), and the other is the Feigin-Frenkel
isomorphism (\thmref{t:FF}). The compatibility between the two is encapsulated by 
\thmref{t:birth}. The other results of local Langlands nature, including the FLE, are ultimately deduced from 
one (or a combination) of these two. 

\medskip

Once the FLE is proved, we will use it in Part II to establish a certain global compatibility of the Langlands functor,
which will play a key role in subsequent papers in this series. This property will essentially say that the Langlands
functor is compatible with the Beilinson-Drinfeld construction of eigensheaves via Kac-Moody localization and 
opers. 

\bigskip

\section{Geometric Satake and Casselman-Shalika formula: recollections} \label{s:Sat}

In this section we will review the constructions of categories of \emph{geometric nature}
associated, on the geometric side, to spaces of maps
$$\cD\to G \text{ and } \cD^\times\to G,$$
and (twisted) D-modules on these spaces, 
and on the spectral side to spaces of maps
$$\cD_\dR\to \cG \text{ and } \cD^\times_\dR\to \cG$$
and ind-coherent sheaves on these spaces. 

\medskip

Thus, the main players are:

\begin{itemize}

\item The category $\Whit(G)$ of Whittaker D-modules on the affine Grassmannian;

\smallskip

\item Its spectral counterpart $\QCoh(\LS^{\on{reg}}_\cG)\simeq \Rep(\cG)$; 

\smallskip

\item The equivalence $\Whit(G)\simeq \Rep(\cG)$, which we call the geometric Casselman-Shalika
formula (\thmref{t:geom CS});

\smallskip

\item The local spherical category $\Sph_G$;

\smallskip

\item Its spectral counterpart $\Sph^{\on{spec}}_\cG$;

\smallskip

\item The (derived) geometric Satake equivalence $\Sat_G:\Sph_G\simeq \Sph^{\on{spec}}_\cG$ (\thmref{t:geom Satake}).

\end{itemize} 

When dealing with these objects there is one major trouble and three ``annoyances", all of which will be introduced in this section, and that will plague us throughout the paper:

\medskip

\begin{enumerate}

\item The trouble is that the local algebro-geometric objects on the spectral side are \emph{not} of finite type 
(once we consider their factorization versions), so the $\IndCoh(-)$ categories associated to them need extra work
to define;

\smallskip

\item This paper is concerned with the \emph{classical} geometric Langlands. However, ``classical" for $G$ means
the \emph{critical} level. This means that the categories on the geometric side\footnote{We avoid using the word ``automorphic"
in the local context, as automorphy refers to the global situation.} will consist not of D-modules, but of
critically or half-twisted D-modules. As a result, throughout the paper, we will have to watch carefully what happens
with these twistings as we move between different spaces.

\smallskip

\item Ultimately, on the geometric side, the object we need to consider is not the constant group-scheme on $X$
with fiber $G$, but rather its twist by the $T$-torsor $\rho(\omega_X)$. This twist is analogous to the usual $\rho$-shift
in the representation theory of the finite-dimensional $G$. Thus, all spaces associated with $G$ will undergo the corresponding 
twist. 

\smallskip

\item 
Both categories $\Sph_G$ and $\Sph^{\on{spec}}_\cG$ are endowed with anti-involutions, denoted $\sigma$ and $\sigma^{\on{spec}}$.
A source of constant headache throughout this paper is that these anti-involutions are compatible under $\Sat_G$, 
\emph{up to the Chevalley involution} on $G$, denoted $\tau_G$. This can be seen as a vestige (in a rather precise sense) of
the fact that the square of the usual Fourier transform is not the identity, but rather is given by the action of $-1$. 

\end{enumerate} 

\ssec{The critical twist} 

\sssec{} \label{sss:omega X} 

We choose once and for all a square root $\omega^{\otimes \frac{1}{2}}_X$ of the canonical line bundle $\omega_X$ on $X$. 

\medskip

\noindent{\it Warning:} In this series of papers, $\omega_X$ denotes the canonical line bundle on $X$, and \emph{not} the dualizing
sheaf on $X$, which is the $[1]$ shift of that. (So, properly, we should have used $\Omega^1_X$, rather than $\omega_X$.) This
deviates from the convention, according to which, for a prestack $\CY$ we denote by $\omega_\CY$ its dualizing sheaf. So
the only exception for this rule is when $\CY$ is the curve $X$ itself. 

\sssec{}

Consider the affine Grassmannian $\Gr_G$ as a factorization space over $X$, equipped with an action
of the (factorization) group indscheme $\fL(G)$.

\medskip

We refer the reader to Sects. \ref{sss:aff grass} and \ref{sss:just loops}, respectively, where the definition of these objects is recalled,
and to \secref{sss:fact spaces}, where the general theory of factorization spaces is set up. 

\sssec{}

Let $\det_{\Gr_G}$ denote the determinant (factorization) line bundle on $\Gr_G$. 

\begin{rem} \label{r:Pfaff}

According to \cite[Sect. 4]{BD1}, the choice 
of $\omega^{\otimes \frac{1}{2}}_X$ gives rise to a square root of $\det_{\Gr_G}$, as a line bundle over $\Gr_{G,\Ran}$. 
However, this square root is \emph{incompatible} with factorization.\footnote{More precisely, this square root exists as a factorization
$\BZ/2\BZ$-graded line bundle, where the grading over the connected component $\Gr_G^\lambda$ of $\Gr_G$ 
(here $\lambda\in \Lambda_{G,G}=\pi_0(\Gr_G)$) equals $\langle \lambda,2\rhoch\rangle\on{mod}2$.}

\end{rem}

\sssec{}

For a line bundle $\CL$ on a space $\CY$ and an integer $n$, let $\CL^{\frac{1}{n}}$ denote the \'etale 
$\mu_n$-gerbe of $n$th roots of $\CL$. 

\medskip

Recall now that given a $\mu_n$-gerbe $\CG$ on a space $\CY$, we can consider the $\CG$-twisted category of D-modules
on $\CY$, to be denoted 
$$\Dmod_\CG(\CY).$$

\medskip

Thus, for $(\CY,\CL,n)$ as above we can consider the corresponding category 
$$\Dmod_{\CL^{\frac{1}{n}}}(\CY).$$

\sssec{} \label{sss:1/2 Gr}

Consider the $\mu_2$-gerbe $\det^{\frac{1}{2}}_{\Gr_G}$. 

\medskip

We will use the short-hand notation 
$$\Dmod_{\frac{1}{2}}(\Gr_G)$$
to denote the (factorization) category 
$$\Dmod_{\det^{\frac{1}{2}}_{\Gr_G}}(\Gr_G)$$
of $\det^{\frac{1}{2}}_{\Gr_G}$-twisted D-modules on $\Gr_G$. 

\begin{rem} \label{r:Pfaff bis}
According to Remark \ref{r:Pfaff}, a choice of $\omega^{\otimes \frac{1}{2}}_X$ gives rise to a trivialization of the gerbe
$\det^{\frac{1}{2}}_{\Gr_G}$. However, this trivialization is incompatible with factorization.

\medskip

For that reason, henceforth, we will avoid using it. 

\end{rem}

\sssec{} \label{sss:dR twistings from line bundles}

Recall that for a space $\CY$, we can consider \emph{de Rham} twistings on $\CY$ (see, \cite[Sect. 6]{GaRo2}). 
These are by definition $\CO^\times$-gerbes on $\CY_\dR$, equipped with a trivialization of their pullback to $\CY$. 

\medskip

Given a de Rham twisting $\CT$, we can consider the corresponding twisted category of D-modules
$$\Dmod_\CT(\CY),$$
see \cite[Sect. 7]{GaRo2}.

\medskip

Recall also that to a line bundle $\CL$ on $\CY$, we can associate a de Rham twisting, which in this paper we denote by 
$\on{dlog}(\CL)$ (the corresponding $\CO^\times$-gerbe on $\CY_\dr$ is trivial, but the trivialization of its pullback to $\CY$
differs from the tautological one by tensoring with $\CL$).

\medskip

Note that tensoring
by $\CL$ defines an equivalence 
\begin{equation} \label{e:tw line bundle}
\Dmod(\CY)\to \Dmod_{\on{dlog}(\CL)}(\CY). 
\end{equation} 

\medskip

Finally, recall (see \cite[Corollary 6.4.5]{GaRo2}) that the space of de Rham twistings on a given space $\CY$ carries a natural a $k$-linear structure.
Thus, for $c\in k$, we have a well-defined twisting $c\cdot \on{dlog}(\CL)$, and the corresponding category
$$\Dmod_{c\cdot\on{dlog}(\CL)}(\CY).$$

\sssec{}

Let $(\CY,\CL,n)$ be as above. Note that for $c=n\in \BZ\subset k$, 
we have 
$$n\cdot \on{dlog}(\CL)=\on{dlog}(\CL^{\otimes n}).$$

In particular, we have a canonical identification
of the corresponding twisted categories of D-modules:
\begin{equation} \label{e:etale vs dR twistings}
\Dmod_{\CL^{\frac{1}{n}}}(\CY) \overset{\sim}\to  \Dmod_{\frac{1}{n}\cdot \on{dlog}(\CL)}(\CY). 
\end{equation}

\medskip 

For example, when $n=1$, the identification \eqref{e:etale vs dR twistings} is the identification
of \eqref{e:tw line bundle}.

\sssec{} \label{sss:crit}

We will use the short-hand notation 
$$\Dmod_{\on{crit}}(\Gr_G)$$
for the (factorization) category 
$$\Dmod_{\frac{1}{2}\cdot\on{dlog}(\det_{\Gr_G})}(\Gr_G).$$

\sssec{} \label{sss:crit vs 1/2}

Applying \eqref{e:etale vs dR twistings} to $\CY=\Gr_G$ and $\CL=\det_{\Gr_G}$, 
we obtain a canonical equivalence of (factorization) categories
$$\Dmod_{\frac{1}{2}}(\Gr_G)\simeq \Dmod_{\on{crit}}(\Gr_G).$$ 

\begin{rem}
According to Remark \ref{r:Pfaff}, we can also identify 
$$\Dmod_{\frac{1}{2}}(\Gr_{G,\Ran})\simeq \Dmod(\Gr_{G,\Ran}),$$
or equivalently
$$\Dmod_{\on{crit}}(\Gr_{G,\Ran})\simeq \Dmod(\Gr_{G,\Ran}),$$
as plain categories, but these identifications are incompatible with the factorization
structures.  

\end{rem}

\begin{rem} \label{r:1/2 vs crit}

We distinguish $\Dmod_{\on{crit}}(\Gr_G)$ and $\Dmod_{\frac{1}{2}}(\Gr_G)$ notationally for the following two reasons:

\smallskip

\noindent(1) The \'etale gerbe-twisted version makes sense not just in the context of D-modules, but also
in other sheaf-theoretic contexts (e.g., Betti, $\ell$-adic).

\smallskip

\noindent(2) The category $\Dmod_{\on{crit}}(\Gr_G)$ comes equipped with a natural forgetful functor to $\IndCoh(\Gr_G)$,
while for a general \'etale gerbe, the gerbe-twisted category of D-modules does not carry such a functor. 

\smallskip

Thus, the distinction between gerbes and twistings becomes relevant when discussing connections 
between D-modules and modules over Lie algebras, as we often do in this paper. We use the $\Dmod_{\frac{1}{2}}(\Gr_G)$
(or $\Dmod_{\frac{1}{2}}(\Bun_G)$) to evoke the sheaf-theoretic geometry of these spaces, 
while  $\Dmod_\crit(\Gr_G)$ (or $\Dmod_\crit(\Bun_G)$) evokes the connection to Kac-Moody representation theory
at the critical level. 

\end{rem} 

\sssec{}

We can also consider the corresponding \emph{multiplicative} factorization $\mu_2$-gerbe 
on $\fL(G)$, equipped with a multiplicative trivialization of its restriction to $\fL^+(G)$. 

\medskip

Since the group indscheme $\fL(N)$ is contractible, the restriction of the above gerbe to it also admits
a canonical multiplicative trivialization.

\medskip

In particular, if $H$ is a factorization subgroup of either $\fL^+(G)$ or $\fL(N)$, it makes sense to consider the (factorization) category
$$\Dmod_{\frac{1}{2}}(\Gr_G)^H$$
of $H$-equivariant D-modules.

\ssec{A geometric twisting construction}  \label{ss:twist by G-bundle}

\sssec{} \label{sss:twist by G-bundle}

Let $H$ be a group mapping to $G$, and let $\CP_H$ be an $H$-torsor over $X$. Taking sections over the
formal disc, $\CP_H$ gives rise to a factorization torsor over $\fL^+(H)$; by a slight abuse of notation, 
we will denote this $\fL^+(H)$-factorization torsor by the same symbol $\CP_H$. 

\medskip

Given a space $\CY$ over $X$, equipped with an action of $\fL^+(H)$, we can form a twist,
to be denoted $\CY_{\CP_H}$, i.e.,
$$\CY_{\CP_H}:=(\CP_H\times \CY)/\fL^+(H).$$

 If $\CY$ was endowed with a factorization structure compatible with the $\fL^+(H)$-action, then so 
is $\CY_{\CP_H}$.

\sssec{}

The space $\CY_{\CP_H}$ is acted on by the adjoint twist $\fL^+(H)_{\CP_H}$ of $\fL^+(H)$. 

\medskip

Note that
we have a canonical isomorphism
\begin{equation} \label{e:alpha}
\CY/\fL^+(H) \simeq \CY_{\CP_H}/\fL^+(H)_{\CP_H}.
\end{equation} 

\sssec{}  \label{sss:remove twist by G-bundle}

We will denote by the subscript $\CP_H$ the various categories of D-modules associated with
the above geometric objects, such as
$$\Dmod(\CY) \rightsquigarrow \Dmod(\CY)_{\CP_H} \text{ and }
\Dmod(\CY)^{\fL^+(H)} \rightsquigarrow (\Dmod(\CY)^{\fL^+(H)})_{\CP_H}.$$

\medskip

Note, however, that thanks to the identification \eqref{e:alpha}, 
the category $(\Dmod(\CY)^{\fL^+(H)})_{\CP_H}$ is canonically equivalent to 
the original category $\Dmod(\CY)^{\fL^+(H)}$. We will denote this equivalence by 
$$\alpha_{\CP_H,\on{taut}}:\Dmod(\CY)^{\fL^+(H)}\overset{\sim}\to (\Dmod(\CY)^{\fL^+(H)})_{\CP_H}.$$ 

\sssec{} \label{sss:rho twist}

A typical example of the above situation that we will consider is when $H=T$, and the $T$-bundle is
$\rho(\omega_X)$, i.e., the bundle induced from $\omega^{\otimes \frac{1}{2}}_X$ by means of
$$2\rho:\BG_m\to T.$$

\ssec{The Whittaker category on the affine Grassmannian} \label{ss:Whit}

\sssec{}

We apply the construction of \secref{sss:rho twist} to $\CY:=\Gr_G$,
viewed as a scheme acted on by 
$\fL^+(T)\subset \fL^+(G)$, and the group indscheme $\fL(N)$. 

\medskip

Thus, we can form the (factorization) space $\Gr_{G,\rho(\omega_X)}$, which is acted on by $\fL(G)_{\rho(\omega_X)}$,
and in particular $\fL(N)_{\rho(\omega_X)}$.

\sssec{} \label{sss:can char}

The group indscheme $\fL(N)_{\rho(\omega_X)}$ is equipped with a homomorphism 
\begin{equation} \label{e:character on N}
\chi:\fL(N)_{\rho(\omega_X)}\to \BG_a,
\end{equation} 
equal to the composition
$$\fL(N)_{\rho(\omega_X)} \to \fL(N/[N,N])_{\rho(\omega_X)}\simeq \underset{I}\Pi\, \fL(\BG_a)_{\omega_X} \overset{\on{Res}}\to
\underset{I}\Pi\, \BG_a\overset{\chi_0}\to \BG_a,$$
where:

\begin{itemize}

\item $I$ is the set of vertices of the Dynkin diagram of $G$;

\item $\fL(\BG_a)_{\omega_X}$ is the twist formed with respect to the $\fL^+(\BG_m)$-action on $\fL(\BG_a)$;

\item $\Res:\fL(\BG_a)_{\omega_X} \to \BG_a$ is the canonical residue map;

\item $\chi_0$ is a non-degenerate character (i.e., a character non-trivial along each factor). 

\end{itemize}

\sssec{}

Let $\bC$ be a category acted on by $\fL(G)_{\rho(\omega_X)}$ at the critical level.\footnote{The discussion here is applicable both when 
we work over a fixed point $\ul{x}\in \Ran$ and in the factorization setting.} Denote:
$$\Whit^!(\bC):=\bC^{\fL(N)_{\rho(\omega_X)},\chi} \text{ and } \Whit_*(\bC):=\bC_{\fL(N)_{\rho(\omega_X)},\chi},$$
where we impose equivariance against the pullback of
$$\on{exp}\in \Dmod(\BG_a)$$
by means of $\chi$ (see \cite{Ra2} for more details). Our normalization for $\on{exp}$ is that it is a character sheaf in the *-sense, i.e.,
$$\on{add}^*(\on{exp})\simeq \on{exp}\boxtimes \on{exp}.$$

\medskip

Note that 
\begin{equation} \label{e:dual Whit}
(\Whit_*(\bC))^\vee \simeq \Whit^!(\bC^\vee),
\end{equation}
up to replacing $\chi_0$ by its inverse, where\footnote{In the next formula $\on{Funct}(-,-)$ stands for colimit-preserving functors. We will
always use this convention when talking about functors between cocomplete categories, unless explicitly specified otherwise.}
$$(-)^\vee:=\on{Funct}((-),\Vect).$$ 

\sssec{}

Although the assignments 
$$\bC\rightsquigarrow \Whit^!(\bC) \text{ and } \bC\rightsquigarrow \Whit_*(\bC)$$
involve the group \emph{ind}-scheme $\fL(N)_{\rho(\omega_X)}$, they behave nicely on the 2-category of $\fL(G)_{\rho(\omega_X)}$-module
categories (see \cite{Ra2}).

\medskip

Namely, they both commute with limits and colimits. Combined with \eqref{e:dual Whit}, this implies that if $\bC$ is dualizable,
then so are $\Whit^!(\bC)$ and $\Whit_*(\bC)$. 

\medskip

However, more is true.  

\sssec{} 

Let $\omega^{\on{ren}}_{\fL(N)_{\rho(\omega_X)}}\in \Dmod(\fL(N)_{\rho(\omega_X)})$ be \emph{the renormalized} dualizing sheaf on $\fL(N)_{\rho(\omega_X)}$,
defined to be the *-pullback of the dualizing sheaf along the projection
$$\fL(N)_{\rho(\omega_X)}\to \fL(N)_{\rho(\omega_X)}/\fL^+(N)_{\rho(\omega_X)}.$$

Consider the object
$$\omega^{\on{ren}}_{\fL(N)_{\rho(\omega_X)},\chi}:=\omega^{\on{ren}}_{\fL(N)_{\rho(\omega_X)}}\overset{*}\otimes \chi^*(\on{exp})\in \Dmod(\fL(N)_{\rho(\omega_X)}).$$

\sssec{} \label{sss:ren Whit aver}

Let $\bC$ be as above. The operation of *-convolution with $\omega^{\on{ren}}_{\fL(N)_{\rho(\omega_X)},\chi}$
is an endofunctor of $\bC$ (as a plain DG category), and this endofunctor factors as
$$\bC\twoheadrightarrow \Whit_*(\bC)\to \Whit^!(\bC)\hookrightarrow \bC.$$

Denote the resulting functor $\Whit_*(\bC)\to \Whit^!(\bC)$ by 
$$\Theta_{\Whit(\bC)}:\Whit_*(\bC)\to \Whit^!(\bC).$$

The following fundamental result was established in \cite{Ra2}:

\begin{thm} \label{t:Whit self-dual gen}
The functor $\Theta_{\Whit(\bC)}$ is an equivalence.
\end{thm} 

\begin{rem}
The proof of \thmref{t:Whit self-dual}, as recorded in \cite{Ra2}, is given for a fixed formal disc,
but the same argument applies to prove a version of this theorem over $\Ran$. 
\end{rem} 

\sssec{}

We apply the above discussion to 
$$\bC:=\Dmod_{\frac{1}{2}}(\Gr_{G,\rho(\omega_X)}).$$
Thus we obtain the (factorization) categories
$$\Whit^!(\Dmod_{\frac{1}{2}}(\Gr_{G,\rho(\omega_X)})) \text{ and } \Whit_*(\Dmod_{\frac{1}{2}}(\Gr_{G,\rho(\omega_X)})).$$

\medskip

We will use for them short-hand notations
$$\Whit^!(G) \text{ and } \Whit_*(G),$$
respectively. 

\begin{rem}  \label{r:indep char}
The categories $\Whit^!(G)$ and $\Whit_*(G)$ are canonically independent of the choice of $\chi_0$:

\medskip

Indeed, given two non-degenerate characters $\chi^1_0$ and $\chi^2_0$, there exists an element $t\in T$
that conjugates $\chi^1_0$ to $\chi^2_0$. Translation by $t$ on $\Gr_{G,\rho(\omega_X)}$
defines then an equivalence between the corresponding Whittaker categories.

\medskip

The choice of $t$ is unique up to an element $z\in Z_G$. However, the translation action of $z$ on
$\Gr_{G,\rho(\omega_X)}$ is trivial.

\end{rem}

\sssec{}

By \eqref{e:dual Whit}, the categories $\Whit^!(G)$ and $\Whit_*(G)$ are naturally mutually dual, up to replacing $\chi_0$ by its inverse. 
Note, however, that due to Remark \ref{r:indep char}, they are actually mutually dual.

\medskip

Furthermore, as is shown in \cite{Ga6}, both $\Whit^!(G)$ and $\Whit_*(G)$ are compactly generated (see \secref{sss:comp gen fact} for what compact
generation means in the factorization setting). 

\sssec{}

Let 
$$\Theta_{\Whit(G)}:\Whit_*(G)\to \Whit^!(G)$$
denote the functor from \secref{sss:ren Whit aver}.

\medskip

As a particular case of \thmref{t:Whit self-dual gen}, we obtain:

\begin{thm} \label{t:Whit self-dual}
The functor $\Theta_{\Whit(G)}$ is an equivalence (of factorization categories).
\end{thm} 

\sssec{}

The factorization categories $\Whit^!(G)$ and $\Whit_*(G)$ are unital (see \secref{sss:fact cat untl}) for what this means. 
Here is the explicit description of their
factorization units:

\medskip

The factorization unit $\one_{\Whit_*(G)}\in \Whit_*(G)$ is the object, denoted $\on{Vac}_{\Whit_*(G)}$, equal to the 
projection along 
$$\Dmod_{\frac{1}{2}}(\Gr_{G,\rho(\omega_X)})\to \Whit_*(G)$$
of $\delta_{1_{\Gr_{G,\rho(\omega_X)}}}\in \Dmod_{\frac{1}{2}}(\Gr_{G,\rho(\omega_X)})$, the latter being the factorization unit
$\one_{\Dmod_{\frac{1}{2}}(\Gr_{G,\rho(\omega_X)})}$ for $\Dmod_{\frac{1}{2}}(\Gr_{G,\rho(\omega_X)})$ itself. 

\sssec{} \label{sss:Whit clean}

The factorization unit $\one_{\Whit^!(G)}\in \Whit^!(G)$ is the object, denoted $\on{Vac}_{\Whit^!(G)}$, equal to 
the *-direct image along the locally-closed embedding
$$\fL(N)_{\rho(\omega_X)}/\fL^+(N)_{\rho(\omega_X)}\hookrightarrow \Gr_{G,\rho(\omega_X)}$$
of 
$$\omega_{\fL(N)_{\rho(\omega_X)}/\fL^+(N)_{\rho(\omega_X)}}\otimes \chi^*(\on{exp})\in
\Dmod(\fL(N)_{\rho(\omega_X)}/\fL^+(N)_{\rho(\omega_X)})^{\fL(N)_{\rho(\omega_X)},\chi}.$$

Note that the above *-extension is \emph{clean}, i.e., receives an isomorphism from the !-extension. 

\medskip

This implies that the functor co-represented by $\on{Vac}_{\Whit^!(G)}$ identifies with the functor of
!-fiber at the unit point $1_{\Gr_{G,\rho(\omega_X)}}\in \Gr_{G,\rho(\omega_X)}$, restricted to 
$$\Whit^!(G)\subset \Dmod_{\frac{1}{2}}(\Gr_{G,\rho(\omega_X)}).$$

\ssec{The geometric Casselman-Shalika formula}

\sssec{}

The following is the statement of the geometric Casselman-Shalika formula 
(see \cite[Theorem 6.36.1]{Ra3}\footnote{The original result in this direction is the main theorem of \cite{FGV}.}):

\begin{thm} \label{t:geom CS}
There exists a canonically defined equivalence of factorization categories:
$$\on{CS}_G:\Whit^!(G)\to \Rep(\cG).$$
\end{thm}

\begin{rem}
In the course of the proof of \thmref{t:geom CS} one uses the \emph{naive}
(i.e., non-derived) geometric Satake to construct a functor
$$\Rep(\cG)\to \Whit^!(G),$$
and one shows that it is an equivalence, see Remark \ref{r:CS again}. 
\end{rem}  
 
\sssec{}

The functor $\on{CS}_G$ is normalized so that it sends the standard object
$$\Delta^\lambda\in \Whit^!(G), \quad \lambda\in \Lambda_G^+,$$
corresponding to the $\fL(N)_{\rho(\omega_X)}$-orbit
$$S^\lambda:=\fL(N)_{\rho(\omega_X)}\cdot t^\lambda$$
to the highest weight module
$$V^{\lambda}\in \Rep(\cG).$$
(In the above formula, $t$ denotes the uniformizer on $\cD$.)

\begin{rem} \label{r:choice CS}

By fixing the above normalization for $\on{CS}_G$ we made a choice. We could have made a 
different choice by applying the Chevalley involution $\tau_G$ on $G$, or equivalently, on $\cG$. 

\medskip

The normalization for $\on{CS}_G$ ultimately forces how we normalize the Langlands functor $\BL_G$. 
Our particular choice for $\on{CS}_G$ is dictated by the following: we want the Langlands functor $\BL_G$ for
$G$ and its Levi subgroups to be compatible via the Eisenstein and Constant Term functors. 

\medskip

If we composed $\BL_G$ with $\tau_G$, the resulting functor would intertwine the Eisenstein/Constant Term functors on the
geometric side with similar functors on the spectral side, but with respect to the \emph{opposite parabolic}. 

\medskip

Note also that when $G$ is a torus $T$, in our normalization, the Langlands functor $\BL_T$ is as in \cite[Sect. 1.5.2]{GLC1}. 

\end{rem}

\begin{rem}
For the validity of \thmref{t:geom CS} at the factorization level, it is crucial that in the definition of $\Whit^!(G)$
we use the twisted category $\Dmod_{\frac{1}{2}}(\Gr_G)$, rather than the untwisted one, i.e., $\Dmod(\Gr_G)$.
\end{rem} 

\sssec{} \label{sss:FLE infty}

The following is a basic pattern of how the equivalence $\on{CS}_G$ interacts with duality.

\medskip

Let us denote by
$$\FLE_{\cG,\infty}:\Rep(\cG)\to \Whit_*(G)$$
the functor equal to $\on{CS}_G^\vee$, with respect to the canonical dualities:
$$\Whit_*(G)=(\Whit^!(G))^\vee \text{ and } \Rep(\cG)^\vee \simeq \Rep(\cG).$$

\begin{rem}
The notation $\FLE_{\cG,\infty}$ stems from the fact that the above functor is indeed
the limiting value of the (positive level) $\FLE$ equivalence. 
\end{rem}

\sssec{Example} 

Note, in particular, that the functor $\FLE_{\cG,\infty}$ sends 
$$V^{-w_0(\lambda)}\in \Rep(\cG)\, \mapsto \nabla^{\lambda}\in \Whit_*(G),$$
where for $\mu\in \Lambda^+$ we denote by 
$$\nabla^\mu\in \Whit_*(G)$$
the object dual to $\Delta^\mu\in \Whit^!(G)$, i.e.,
$$\langle \CF,\nabla^\mu\rangle= \CHom_{\Whit^!(G)}(\Delta^\mu,\CF), \quad \CF\in \Whit^!(G),$$
where 
$$\langle -,-\rangle:\Whit^!(G)\otimes \Whit_*(G)\to \Vect$$
is the canonical pairing.

\sssec{} \label{sss:Cartan inv}

Note that the Whittaker category is canonically attached to the pair $(G,B)$. Hence,
the group of \emph{outer automorphisms} of $G$ (i.e., the group of automorphisms of the 
polarized\footnote{By a polarization of a root datum we mean a choice of the subset of positive roots.} 
root datum of $G$) acts on both versions of the Whittaker category.

\medskip

Let $\tau_G$ be the Chevalley involution, viewed as an outer automorphism of $G$. 
The corresponding automorphism of the polarized root datum acts as $\lambda\mapsto -w_0(\lambda)$.

%

\sssec{}

We have:

\begin{lem} \label{l:CS and duality}
The composition
$$\Rep(\cG) \overset{\FLE_{\cG,\infty}}\longrightarrow \Whit_*(G)\overset{\Theta_{\Whit(G)}}\longrightarrow \Whit^!(G)$$
identifies canonically with
$$\tau_G\circ (\on{CS}_G)^{-1}.$$
\end{lem}

\begin{rem}

As a reality check, note that both functors in \lemref{l:CS and duality} send 
$$V^{-w_0(\lambda)}\in \Rep(\cG)\,\, \mapsto \,\, \Delta^{\lambda}\in \Whit^!(G).$$

Indeed, the functor $\Theta_{\Whit(G)}$ is easily seen to send $\Delta^\lambda$ to $\nabla^\lambda$. 

\medskip

The proof of \lemref{l:CS and duality} follows easily from the construction of $\on{CS}_G$ via
naive geometric Satake.

\end{rem}

\ssec{The spherical category}

\sssec{}

We denote by $\Sph^{\on{non-ren}}_G$ the (factorization) monoidal category 
$$\Dmod_{\frac{1}{2}}(\fL^+(G)\backslash \fL(G)/\fL^+(G)).$$

\medskip

We have a naturally defined right action of $\Sph^{\on{non-ren}}_G$ on $\Dmod_{\frac{1}{2}}(\Gr_G)$, compatible
with the left action of $\fL(G)$. 

\sssec{}

We let $\Sph_G$ denote its \emph{renormalized} version, which is defined as the ind-completion of the
full subcategory in $\Sph^{\on{non-ren}}_G$ consisting of objects whose image under (either of) the forgetful functors
$$\Dmod_{\frac{1}{2}}(\fL^+(G)\backslash \fL(G))\leftarrow \Dmod_{\frac{1}{2}}(\fL^+(G)\backslash \fL(G)/\fL^+(G))\to
\Dmod_{\frac{1}{2}}(\fL(G)/\fL^+(G))$$
is compact (see \cite[Proposition 6.3.2]{CR} for more details).

\medskip

By construction, the monoidal (and also factorization) unit
$$\one_{\Sph_G}\simeq \delta_{1_{\Gr_G}}\in \Sph_G$$
is compact.

\sssec{}

We have an adjoint pair of functors 
$$\on{ren}:\Sph^{\on{non-ren}}_G\rightleftarrows \Sph_G:\on{non-ren},$$
with $\on{ren}$ being fully faithful and $\on{non-ren}$ monoidal. This makes
$\Sph^{\on{non-ren}}_G$ into a monoidal colocalization of $\Sph_G$. 

\medskip

In particular, we have a right action of $\Sph_G$ on $\Dmod_{\frac{1}{2}}(\Gr_G)$, compatible
with the left action of $\fL(G)$ and factorization. 

\sssec{} \label{sss:left to right Sph}

Inversion on the group $\fL(G)$ defines an anti-involution, denoted $\sigma$, of 
$\Sph_G$. We will refer to it as the ``flip" anti-involution.

\medskip

Henceforth, we will use $\sigma$ to pass between left and right module
categories over $\Sph_G$. In light of this, we will not necessarily distinguish between
left and right actions of $\Sph_G$. 

\sssec{} \label{sss:duals in Sph}

The fact that $\Gr_G$ is ind-proper implies that the composition of the involution $\sigma$ with Verdier duality (on compact objects)
defines an equivalence
\begin{equation} \label{e:SphG is rigid}
\Sph_G^\vee \simeq \Sph_G,
\end{equation} 
which identifies both with right and left monoidal dualization.

\medskip

Combined with the fact that the unit in $\Sph_G$ is compact, we obtain that $\Sph_G$ is \emph{rigid} as a 
monoidal category.\footnote{Being a monoidal
colocalization of a rigid category, $\Sph^{\on{non-ren}}_G$ is semi-rigid (see \cite[Appendix C]{AGKRRV}).}

%
%
%
%
%
%
%

%
%

\sssec{}

Recall the setting of \secref{ss:twist by G-bundle}. For any $G$-bundle $\CP_G$ 
on $X$, we can form the twisted version 
$$\Sph_{G,\CP_G}$$ 
of $\CP_G$. 

\medskip

We have a naturally defined action of $\Sph_{G,\CP_G}$ on $\Dmod_{\frac{1}{2}}(\Gr_{G,\CP_G})$, compatible
with the left action of $\fL(G)_{\CP_G}$ and factorization.

\sssec{}

In particular, we have a natural action of $\Sph_{G,\rho(\omega_X)}$ on $\Whit^!(G)$ and $\Whit_*(G)$. 

\medskip

These actions are compatible both with the duality
\begin{equation} \label{e:Whi !* duality}
(\Whit^!(G))^\vee \simeq \Whit_*(G)
\end{equation}
(see \secref{sss:left to right Sph}) and the functor $\Theta_{\Whit(G)}$. 

\sssec{}

Note, however, that according to \secref{sss:remove twist by G-bundle}, we can identify\footnote{In the formula below we consider $\fL(G)$ as acted on
by $\fL^+(G)\times \fL^+(G)$.}
$$\Sph_G\overset{\alpha_{\rho(\omega_X),\on{taut}}}\simeq \Sph_{G,\rho(\omega_X)},$$
and thus we can regard $\Whit^!(G)$ and $\Whit_*(G)$ as acted on by $\Sph_G$ itself. 

\ssec{The spectral spherical category}

\sssec{}

Consider the \emph{local spectral Hecke stack} 
$$\on{Hecke}_\cG^{\on{spec,loc}}:=\LS^{\on{reg}}_\cG\underset{\LS^{\on{mer}}_\cG}\times \LS^{\on{reg}}_\cG.$$
as a factorization space.

\medskip

In the above formula $\LS^{\on{reg}}_\cG$ (resp., $\LS^{\on{mer}}_\cG$) is the factorization space that attaches
to $\ul{x}\in \Ran$ the stack $\LS_\cG(\cD_{\ul{x}})$ (resp., $\LS_\cG(\cD^\times_{\ul{x}})$)
of $\cG$-local systems on the formal multi-disc $\cD_{\ul{x}}$ 
(resp., the punctured multi-disc $\cD^\times_{\ul{x}}:=\cD_{\ul{x}}-\ul{x}$), see \secref{sss:LS loc}.

\sssec{}

The fiber $\on{Hecke}_{\cG,\ul{x}}^{\on{spec,loc}}$ of $\on{Hecke}_\cG^{\on{spec,loc}}$ over a given point $\ul{x}\in \Ran$ is the stack 
\begin{equation} \label{e:Hecke spec local x}
\on{Hecke}_{\cG,\ul{x}}^{\on{spec,loc}}:=\LS^\reg_{\cG,\ul{x}}\underset{\LS^\mer_{\cG,\ul{x}}}\times \LS^\reg_{\cG,\ul{x}}.
\end{equation} 

The stack \eqref{e:Hecke spec local x} is locally of finite type. In fact, its is isomorphic of the product of copies of
\begin{equation} \label{e:local Hecke spec fin dim}
\on{pt}/\cG\underset{\cg^\wedge_0/\on{Ad}(\cG)}\times \on{pt}/\cG
\end{equation}
for each distinct point that comprises $\ul{x}$. 

\sssec{}

Hence, it makes sense to consider the category 
$$\Sph^{\on{spec}}_{\cG,\ul{x}}:=\IndCoh(\on{Hecke}_{\cG,\ul{x}}^{\on{spec,loc}}).$$

We endow $\Sph^{\on{spec}}_{\cG,\ul{x}}$ with a monoidal structure via 
*-pull and *-push along the standard convolution diagram. 

\sssec{} \label{sss:Sph spec tech}

As we let $\ul{x}$ move along $\Ran$ (or $X^n$ for a fixed integer $n$), the resulting prestack 
is no longer locally almost of finite type, so the category of ind-coherent sheaves on it is not a priori-defined. 

\medskip

In fact, $\on{Hecke}_\cG^{\on{spec,loc}}$ violates the condition of being (locally almost) of finite type so badly,
that we do not really know how to define the corresponding category $\IndCoh^*(\on{Hecke}_\cG^{\on{spec,loc}})$
algorithmically. 

\medskip

We refer the reader to \secref{s:spec Sph fact}, where the definition is given (and is compared to another
working definition, adopted in \cite{CR}).

\medskip

Accordingly, the proofs of all the statements that involve $\Sph^{\on{spec}}_\cG$ are also delegated to
\secref{s:spec Sph fact}. In the main body of the text, we will supply prototypes of the
corresponding proofs for the pointwise version $\Sph^{\on{spec}}_{\cG,\ul{x}}$. 

\sssec{}

The pointwise version of the spectral Hecke category $\Sph^{\on{spec}}_{\cG}$, i.e.,  $\Sph^{\on{spec}}_{\cG,\ul{x}}$, is equipped with a tautological action 
on $\QCoh(\LS^\reg_{\cG,\ul{x}})$. 

\medskip 

This construction persists in the factorization setting, i.e., we have an action of the monoidal factorization category 
$\Sph^{\on{spec}}_\cG$ on 
\begin{equation} \label{e:Rep fact}
\Rep(\cG)\simeq \QCoh(\LS^{\on{reg}}_\cG),
\end{equation} 
viewed as a plain\footnote{As opposed to a (symmetric) monoidal factorization category.}
factorization category (see \secref{sss:QCOh LS loc} where the equivalence \eqref{e:Rep fact}
is established). 


\sssec{}

By construction, the category $\Sph^{\on{spec}}_\cG$ receives a monoidal functor, denoted 
$$\on{nv}:\Rep(\cG)\to \Sph^{\on{spec}}_\cG,$$
to be thought of\footnote{And literally so over a fixed point of $\Ran$.} as the direct image functor along
$$\LS^{\on{reg}}_\cG\to \on{Hecke}_\cG^{\on{spec,loc}},$$
where we now view the categories appearing in 
\eqref{e:Rep fact} as (symmetric) \emph{monoidal} factorization categories. 

\sssec{}

The flip of two factors defines an anti-involution on $\Sph^{\on{spec}}_\cG$ to be 
denoted $\sigma^{\on{spec}}$. 

\medskip

We will use $\sigma^{\on{spec}}$ to pass between left and right $\Sph^{\on{spec}}_\cG$-module
categories. 

\medskip

Note that we have a commutative diagram 
\begin{equation} \label{e:sigma spec}
\CD
\Rep(\cG) @>{\on{nv}}>> \Sph_\cG^{\on{spec}} \\
@V{\on{Id}}VV @VV{\sigma^{\on{spec}}}V \\
\Rep(\cG) @>{\on{nv}}>> \Sph_\cG^{\on{spec}},
\endCD
\end{equation}  
where $\on{Id}$ makes sense as an anti-involution of $\Rep(\cG)$, since this category is symmetric monoidal.

\ssec{Geometric Satake equivalence}

\sssec{}

The following is the statement of the factorization version of the \emph{derived} geometric Satake equivalence (see \cite{CR}):

\begin{thm} \label{t:geom Satake}
There exists a unique equivalence of monoidal factorization categories
$$\Sat_G:\Sph_G\to \Sph_\cG^{\on{spec}},$$
compatible with the actions of $\Sph_G$ on $\Whit(G)$ and
$\Sph_\cG^{\on{spec}}$ on $\Rep(\cG)$ via the equivalence
$$\on{CS}_G:\Whit^!(G)\simeq \Rep(\cG).$$
\end{thm} 

The construction of the functor $\Sat_G$ will be recalled in \secref{ss:FLE and Sat}. 

\begin{rem}
In this series of papers we will refer to $\Sat_G$ just as ``geometric Satake equivalence", omitting the word
``derived". What is more commonly referred to as ``geometric Satake" is not an equivalence, but a functor 
in one direction, which we will refer to as ``naive Satake" and denote by $\Sat_G^{\on{nv}}$, see 
\secref{sss:naive Sat}.
\end{rem} 

\sssec{Example}

Unwinding the construction, we obtain that $\Sat_G$ sends the object in $\Sph_G$ corresponding to the double coset of
the point $t^\lambda$ (for $\lambda\in \Lambda^+$) to the object
$$\on{nv}(V^{-w_0(\lambda)})\in \Sph_\cG^{\on{spec}}.$$

The above object in $\Sph_G$ is what is usually denoted by
$$\on{IC}_{\ol\Gr^{\lambda}_G},$$
the intersection cohomology sheaf on the closure of the $\fL^+(G)$-orbit $\Gr^{\lambda}$ of
$t^{\lambda}$.

\begin{rem}
As in the case of \thmref{t:geom CS}, for the validity of \thmref{t:geom Satake} at the factorization level, it is crucial
that we work with the twisted category
$$\Dmod_{\frac{1}{2}}(\fL^+(G)\backslash \fL(G)/\fL^+(G))$$
rather than with $\Dmod(\fL^+(G)\backslash \fL(G)/\fL^+(G))$.
\end{rem} 

\sssec{} \label{sss:naive Sat}

In what follows, we will denote by $\Sat_G^{\on{nv}}$ the functor
$$\Rep(\cG)\overset{\on{nv}}\to \Sph_\cG^{\on{spec}} \overset{\Sat^{-1}}\to \Sph_G.$$

\begin{rem}

One thing that is naive about the naive Satake functor is its direction.
Our conventions are that functors in geometric Langlands go from the 
geometric side to the spectral side.
However, the naive Satake functor produces sheaves from representations.

\end{rem}

\begin{rem} \label{r:CS again}

Note, for example that the functor
$$\Rep(\cG)\overset{\Sat_G^{\on{nv}}}\longrightarrow \Sph_G 
\overset{-\star \on{Vac}_{\Whit^!(G)}}\longrightarrow \Whit^!(G)$$
is $\on{CS}_G^{-1}$. 

\medskip

The functor 
$$\Rep(\cG)\overset{\Sat_G^{\on{nv}}}\longrightarrow \Sph_G \overset{\sigma}\to \Sph_G 
\overset{-\star \on{Vac}_{\Whit_*(G)}}\longrightarrow \Whit_*(G)$$
is $\FLE_{\cG,\infty}$. 

\end{rem}

\ssec{The curse of \texorpdfstring{$\sigma$}{sigma} and \texorpdfstring{$\tau$}{tau}} \label{ss:curse}

\sssec{}

The following statement results from the uniqueness assertion in \thmref{t:geom Satake}
combined with \lemref{l:CS and duality}:

\begin{cor} \label{c:sigma and Sat}
The following diagram of anti-equivalences commutes:
$$
\CD
\Sph_G  @>{\Sat_G}>> \Sph_\cG^{\on{spec}} \\
@V{\sigma}VV  \\
\Sph_G & & @VV{\sigma^{\on{spec}}}V \\
@V{\tau_G}VV \\ 
\Sph_G  @>{\Sat_G}>> \Sph_\cG^{\on{spec}}.
\endCD
$$
\end{cor}  

\sssec{} \label{sss:Sat -1 nv tau}

Denote by $\Sat_{G,\tau}$ the (factorization) equivalence
$$\Sph_G  \overset{\tau}\to \Sph_G\overset{\Sat_G}\longrightarrow \Sph_\cG^{\on{spec}}.$$

Denote by $\Sat^{\on{nv}}_{G,\tau}$ the functor
$$\tau_G\circ \Sat^{\on{nv}}_G, \quad \Rep(\cG)\to \Sph_G.$$

%
%

\sssec{}

As another corollary of \lemref{l:CS and duality} we obtain:

\begin{cor} \label{c:geom Sat tau}
The equivalence
$$\Rep(\cG)\overset{\FLE_{\cG,\infty}}\simeq \Whit_*(G)$$
is compatible with the actions of 
$\Sph_G$ and $\Sph_\cG^{\on{spec}}$ via $\Sat_{G,\tau}$.
\end{cor} 

\sssec{Warning} \label{sss:curse}

As has been mentioned above, we will use $\sigma$ (resp., $\sigma^{\on{spec}}$) to pass between
left and right module categories for $\Sph_G$ (resp., $\Sph^{\on{spec}}_\cG$).

\medskip

Note, however, that due to \corref{c:sigma and Sat}, this procedure is compatible with the geometric
Satake equivalence \emph{up to} the Chevalley involution.

\medskip

In practice, this will manifest itself as follows. Let $\bC_1$ and $\bC_2$ (resp., $\bC_1^{\on{spec}}$ and $\bC_2^{\on{spec}}$)
be left module categories over $\Sph_G$ (resp., $\Sph^{\on{spec}}_\cG$). Thanks to the above left-right passage, we can form
the tensor products
$$\bC_1\underset{\Sph_G}\otimes \bC_2 \text{ and } \bC^{\on{spec}}_1\underset{\Sph^{\on{spec}}_\cG}\otimes \bC^{\on{spec}}_2.$$

Suppose that we have a given a functor
$$F_1:\bC_1\to \bC_1^{\on{spec}},$$
which is compatible with the actions via
\begin{equation} \label{e:Sat curse}
\Sph_G\overset{\Sat_G}\simeq \Sph^{\on{spec}}_\cG
\end{equation} 
and a functor 
$$F_2:\bC_2\to \bC_2^{\on{spec}},$$
which is compatible with the actions via
\begin{equation} \label{e:Sat tau curse}
\Sph_G\overset{\Sat_{G,\tau}}\simeq \Sph^{\on{spec}}_\cG.
\end{equation} 

In this case, we obtain a functor
\begin{equation} \label{e:curse}
F_1\otimes F_2:\bC_1\underset{\Sph_G}\otimes \bC_2  \simeq \bC^{\on{spec}}_1\underset{\Sph^{\on{spec}}_\cG}\otimes \bC^{\on{spec}}_2.
\end{equation} 

\sssec{Warning} \label{sss:curse duality}

Similarly, let $\bC$ and $\bC'$ be left module categories over $\Sph_G$ and $\Sph_\cG^{\on{spec}}$,
respectively. Let us view $\bC^\vee$ (resp., $\bC'{}^\vee$) again as a left module, using $\sigma$
(resp., $\sigma^{\on{spec}}$).

\medskip

Let $\bC\simeq \bC'$ be an equivalence compatible with the actions via \eqref{e:Sat curse}. 
Then the induced equivalence
$$\bC^\vee\simeq \bC'{}^\vee$$
is compatible with the actions via \eqref{e:Sat tau curse}.

\section{Kac-Moody modules and the Kazhdan-Lusztig category} \label{s:KL}

In this section we study the local representation-theoretic category on the geometric side,
which we will later connect to the global category $\Dmod_{\frac{1}{2}}(\Bun_G)$ by a 
local-to-global procedure. 

\medskip

The category in question is the Kazhdan-Lusztig category at the critical level,
$$\KL(G)_\crit:=\hg\mod_\crit^{\fL^+(G)}.$$

We will need the following aspects of the theory associated with $\KL(G)_\crit$:

\begin{itemize}

\item Self-duality;



%

\item The functor of Drinfeld-Sokolov reduction.

\end{itemize}

\ssec{Definition and basic properties}

\sssec{} 

Let $\kappa$ be a level for $\fg$. We consider 
$$\hg\mod_\kappa,$$
the category of Kac-Moody modules at level $\kappa$. This category carries a natural action of $\fL(G)$ at level $\kappa$. 

\medskip

The definition of this category at a fixed point $\ul{x}\in \Ran$ is given in \cite{Ra5}. The factorization version is defined in 
\secref{ss:KM fact}. 

\sssec{}

Let 
$$\KL(G)_\kappa:=\hg\mod_\kappa^{\fL^+(G)},$$
denote the corresponding category of spherical objects.

\medskip

We have an adjunction
$$\oblv_{\fL^+(G)}:\KL(G)_\kappa\rightleftarrows \hg\mod_\kappa:\Av^{\fL^+(G)}_*.$$

\sssec{}

We have a monadic adjunction
\begin{equation} \label{e:KL Rep adj}
\ind_{\fL^+(G)}^{(\hg,\fL^+(G))_\kappa}:\Rep(\fL^+(G))\rightleftarrows \KL(G)_\kappa:\oblv_{\fL^+(G)}^{(\hg,\fL^+(G))_\kappa}. 
\end{equation} 

In particular, $\KL(G)_\kappa$ is compactly generated by the image of compact generators of $\Rep(\fL^+(G))$,
where the latter is \emph{by definition} the ind-completion of the small category consisting of finite-dimensional representations.



\sssec{} \label{sss:KM fact alg}

Let $\on{Vac}(G)_\kappa$ denote the factorization unit in $\KL(G)_\kappa$. By by a slight abuse of notation, we will denote by the 
same symbol $\on{Vac}(G)_\kappa$ its image under the (strictly unital) factorization functor $\oblv_{\fL^+(G)}$. 

\medskip

We let $\BV_{\fg,\kappa}$ denote the image of $\on{Vac}(G)_\kappa$ under the tautological forgetful functor
$$\hg\mod_\kappa\to \Vect.$$

The latter is the usual factorization (a.k.a. chiral) algebra attached to the pair $(\fg,\kappa)$. 

\sssec{}

Our primary interest in this paper is the case when $\kappa=-\frac{1}{2}\cdot \on{Kil}$,
where $\on{Kil}$ is the Killing form. 
We will denote the corresponding level by the symbol $\on{crit}$. 

\medskip

By construction, the category $\KL(G)_\crit$ carries a monoidal action of $\Sph_G$. 

\ssec{Duality} \label{ss:KL duality} 

\sssec{}

For a given level $\kappa$, denote 
$$\kappa':=-\kappa+2\cdot\on{crit}.$$

(In particular, $\on{crit}'=\on{crit}$.)

\sssec{} \label{sss:KM duality}

It is known that the categories 
$$\hg\mod_\kappa \text{ and } \hg\mod_{\kappa'}$$
are canonically dual to one another, in a way compatible with the (unital) factorization structure and the $\fL(G)$-action,
see \cite[Sect. 9.16.1]{Ra5}. 

\medskip

The counit of the duality is the functor
$$\hg\mod_\kappa\otimes \hg\mod_{\kappa'}\overset{\otimes}\to \hg\mod_{-\!\on{Kil}}\to \Vect,$$
where the second arrow is the functor of semi-infinite cohomology. 

\medskip

By \secref{sss:duality lax untl}, the above functor pairing has a structure of 
(lax unital\footnote{See \secref{sss:lax vs strict unital} for what ``lax unital" means.}) 
factorization functor. 

\sssec{} \label{sss:KL self-duality kappa}

The above duality induces a duality between 
\begin{equation} \label{e:KL self-duality kappa}
(\KL(G)_\kappa)^\vee \simeq \KL(G)_{\kappa'},
\end{equation} 
so that 
$$\left(\oblv_{\fL^+(G)}\right)^\vee\simeq \Av^{\fL^+(G)}_* \text{ and }
\left(\Av^{\fL^+(G)}_*\right)^\vee\simeq \oblv_{\fL^+(G)}.$$

\medskip

The unit of the duality \eqref{e:KL self-duality kappa} is the object 
$$\fCDO(G)_{\kappa,\kappa'}\in \KL(G)_\kappa\otimes \KL(G)_{\kappa'}.$$

\medskip

Under this duality and the canonical self-duality of $\Rep(\fL^+(G))$, we have 
$$\left(\ind_{\fL^+(G)}^{(\hg,\fL^+(G))_\kappa}\right)^\vee\simeq \oblv_{\fL^+(G)}^{(\hg,\fL^+(G))_\kappa} \text{ and }
\left(\oblv_{\fL^+(G)}^{(\hg,\fL^+(G))_\kappa}\right)^\vee\simeq \ind_{\fL^+(G)}^{(\hg,\fL^+(G))_\kappa}.$$

\sssec{} \label{sss:KL self-duality crit}

Specializing to the critical level, we obtain a canonical self-duality
\begin{equation} \label{e:KM crit self dual}
(\hg\mod_{\on{crit}})^\vee \simeq \hg\mod_{\on{crit}},
\end{equation} 
compatible with the $\fL(G)$-actions, and
\begin{equation} \label{e:KL crit self dual}
(\KL(G)_{\on{crit}})^\vee \simeq \KL(G)_{\on{crit}},
\end{equation} 
compatible with the $\Sph_G$-actions.

\ssec{The functor of Drinfeld-Sokolov reduction} \label{ss:DS}

\sssec{}

The duality in \secref{sss:KM duality} is applicable to \emph{any} finite-dimensional Lie algebra
(where the role of the level $2\cdot \crit$ is played by the Tate extension). In particular, for a unipotent Lie algebra $\fn'$,
the corresponding category $\fL(\fn')$ is canonically self-dual, in a way compatible with the $\fL(N')$-action. 

\medskip

This construction is functorial. Hence, if a group $H$ acts on $\fn'$, and $\CP_H$ is an $H$-torsor on $X$,
we obtain a canonical self-duality on $\fL(\fn')_{\CP_H}\mod$.

\medskip

In the particular case $N'=N$, $H=T$ and $\CP_H=\rho(\omega_X)$, 
we obtain a canonical duality on 
$$\fL(\fn)_{\rho(\omega_X)}\mod.$$

\sssec{}  

The character $\chi$ (see \secref{sss:can char}) on the group $\fL(N)_{\rho(\omega_X)}$ gives rise to a chracter
(denoted by the same symbol)
$$\fL(\fn)_{\rho(\omega_X)}\to k.$$

We can regard this character as an object
$$k_\chi \in \fL(\fn)_{\rho(\omega_X)}\mod.$$

The factorization property of $\chi$ equips $k_\chi$ with a structure of \emph{factorization algebra} in 
$\fL(\fn)_{\rho(\omega_X)}\mod$, where the latter is regarded as a lax factorization category. 

\medskip

The fact that
$\chi|_{\fL(\fn)^+_{\rho(\omega_X)}}=0$ implies that this factorization algebra is naturally unital.

\sssec{} \label{sss:twisted BRST}

We define the functor 
\begin{equation} \label{e:pre DS}
\BRST_{\fL(\fn)_{\rho(\omega_X)},\chi}:\fL(\fn)_{\rho(\omega_X)}\mod\to \Vect
\end{equation}
to be given by 
$$\fL(\fn)_{\rho(\omega_X)}\mod\overset{\on{Id}\otimes k_\chi}\longrightarrow 
\fL(\fn)_{\rho(\omega_X)}\mod\otimes \fL(\fn)_{\rho(\omega_X)}\mod\to \Vect,$$
where the second arrow is the counit of the self-duality on $\fL(\fn)_{\rho(\omega_X)}\mod$.

\medskip

By the above, the functor \eqref{e:pre DS} has a natural (lax unital) factorization structure. 

\sssec{}  \label{sss:DS}

Precomposing with 
$$\hg\mod_{\kappa,\rho(\omega_X)}\to \fL(\fn)_{\rho(\omega_X)}\mod,$$
we obtain a functor of Drinfeld-Sokolov reduction, which we denote by 
\begin{equation} \label{e:initial DS}
\DS:\hg\mod_{\kappa,\rho(\omega_X)}\to \Vect.
\end{equation}

The functor \eqref{e:initial DS} inherits a (lax unital) factorization structure. 

\sssec{}

It follows from the construction that the functor $\DS$ factors as 
\begin{equation} 
\hg\mod_{\kappa,\rho(\omega_X)}\to 
\Whit_*(\hg\mod_{\kappa,\rho(\omega_X)})\to \Vect,
\end{equation} 

\medskip

We denote the resulting functor 
$$\Whit_*(\hg\mod_{\kappa,\rho(\omega_X)})\to \Vect$$
by 
\begin{equation} \label{e:DS factor}
\ol{\DS}:\Whit_*(\hg\mod_{\kappa,\rho(\omega_X)})\to \Vect.
\end{equation} 

\begin{rem}
The category $\Whit_*(\hg\mod_{\kappa,\rho(\omega_X)})$ and the functor $\DS$ are canonically
independent of the choice of the character $\chi_0$ of $\fn/[\fn,\fn]$. This happens by the same mechanism
as in Remark \ref{r:indep char}: the action of the center of the derived group of $G$ on $\hg\mod_\kappa$
is trivial.
\end{rem}

\sssec{}

In the sequel, we will use the following assertion (see \cite[Theorem 3.2.2]{FG2}): 

\begin{lem} \label{l:DS on KL exact}
For a fixed $\ul{x}\in \Ran$, the functor 
$$\KL(G)_{\kappa,\rho(\omega_X),\ul{x}} \to \hg\mod_{\kappa,\ul{x}}\overset{\DS}\to \Vect$$
is t-exact.
\end{lem} 

\begin{rem} 
One can show that the assertion of \lemref{l:DS on KL exact} holds also in the factorization setting:
the proof of \cite[Corollary 7.2.2]{Ra2} given at the end of Sect. B.3 of \emph{loc. cit.} adapts to the factorization setting.
\end{rem}

\section{Ind-coherent sheaves on monodromy-free opers} \label{s:opers}

In this section we study the local counterpart of the Kazhdan-Lusztig category on the spectral side:
this is the category 
$$\IndCoh^*(\Op^\mf_\cG)$$
of ind-coherent sheaves on the space of monodromy-free $\cG$-opers on the punctured disc. This category
will be related to the global spectral category (in this case $\QCoh(\LS_\cG(X))$) 
by a local-to-global procedure. 

\medskip

We study $\IndCoh^*(\Op^\mf_\cG)$ along with its cousins, the factorization categories
$$\IndCoh^*(\Op^{\on{reg}}_\cG) \text{ and } \IndCoh^*(\Op^{\on{mer}}_\cG).$$

\medskip

Since the geometric objects involved are \emph{not} locally of finite type, the definition
of $\IndCoh(-)$ on them is not automatic. However, thankfully, these objects turn out to be
\emph{(ind)-placid} (see \secref{sss:placid} for what this means), so we have well-behaved
categories $\IndCoh^*(-)$ and $\IndCoh^!(-)$ attached to them. 

%
%
%
%
%
%
%

\ssec{Monodromy-free opers}

\sssec{}

Let $Y$ be a D-scheme over $X$. Recall the notion of $\cG$-oper on $Y$. This is a datum of a triple
$$(\CP_\cB,\epsilon,\nabla),$$ where

\begin{itemize}

\item $\CP_\cB$ is a $\cB$-bundle on $\CY$;

\item $\epsilon$ is the identification between the induced $\cT$-bundle $\CP_\cT$ with $\rhoch(\omega_X)|_Y:=2\rhoch(\omega^{\otimes \frac{1}{2}})|_Y$;

\item $\nabla$ is a connection \emph{along} $X$ on the induced $\cG$-bundle $\CP_\cG$.

\end{itemize}

These data are supposed to satisfy the following compatibility condition:

\medskip

The \emph{incompatibilty} of $\CP_\cB$ and $\nabla$, which is an element
$$\nabla\,\on{mod}\, \cb\in (\cg/\cb)_{\CP_\CG}\otimes \omega_X|_Y$$
belongs to 
$$\on{Fil}_{-1}(\cg/\cb)_{\CP_\CG}\otimes \omega_X|_Y\subset (\cg/\cb)_{\CP_\CG}\otimes \omega_X|_Y$$
(here $\on{Fil}_{-1}(\cg/\cb)\subset \cg/\cb$ is the bottom piece of the principal filtration), 
and its evaluation by means of every negative simple root $-\alpha_i$ of $\cg$
$$-\alpha_i(\nabla\,\on{mod}\, \cb) \in -\alpha_i(\CP_\cT)\otimes \omega_X|_Y \overset{\epsilon}\simeq 
-\alpha_i(\rhoch(\omega_X))|_Y\otimes \omega_X|_Y \simeq \CO_Y$$
is the unit section. 

\sssec{}

A priori, opers form a D-prestack over $X$. However, one shows (see, e.g., \secref{sss:Op is a torsor}) that it is actually an affine
D-scheme over $X$. 

\sssec{}

We will denote the D-scheme of $\cG$-opers by $\Op_\cG$. Its fiber over a given point $x\in X$
is the scheme $\Op_\cG(\cD_x)$ of $\cG$-opers on the formal disc $\cD_x$ around $x$. 

\medskip

We will denote by the symbol 
$$\Op^{\on{reg}}_\cG:=\fL^+_\nabla(\Op_\cG)$$ 
the corresponding factorization (affine) scheme (see \secref{sss:forming arcs gen}), i.e., its fiber $\Op^\reg_{\cG,\ul{x}}$
over a given $\ul{x}\in \Ran$
is the scheme $\Op_\cG(\cD_{\ul{x}})$ of $\cG$-opers on the formal multi-disc $\cD_{\ul{x}}$ around $\ul{x}$. 

\medskip

We let 
$$\Op^{\on{mer}}_\cG:=\fL_\nabla(\Op_\cG)$$ denote the factorization ind-scheme of $\cG$-opers on the formal punctured disc
(see \secref{sss:forming loops aff}). Its
fiber $\Op^{\on{mer}}_{\cG,\ul{x}}$ over a given $\ul{x}\in \Ran$ is the ind-scheme $\Op_\cG(\cD^\times_{\ul{x}})$ of 
$\cG$-opers on the \emph{punctured} 
multi-disc $\cD^\times_{\ul{x}}$. 

\sssec{} \label{sss:oper bundle is fixed}

We recall the following basic fact about opers:

\medskip

Once the ambient curve $X$ is fixed, 
we can assume that the $G$-bundle underlying an oper (on $X$ itself, a multi-disc in $X$, or a punctured multi-disc in $X$)
is induced from a fixed $B$-bundle, to be denoted $\CP_\cB^{\Op}$ (see \cite[Proposition 3.1.10(iii)]{BD1}). 

\medskip

In what follows we will denote by $\CP_\cG^{\Op}$ the induced $\cG$-bundle. 

\sssec{}

By construction, we have a map
$$\Op^{\on{reg}}_\cG\to \LS^{\on{reg}}_\cG,$$
to be denoted $\fr^{\on{reg}}$.

\medskip

Note now that thanks to \secref{sss:oper bundle is fixed} we also have 
a map\footnote{We alert the reader to the fact that $\LS^{\on{mer}}_\cG$ is \emph{not}
the space of loops into $\on{pt}/\cG$, see \secref{sss:LS punctured bad}.}
$$\Op^{\on{mer}}_\cG \to \LS^{\on{mer}}_\cG,$$
to be denoted $\fr$. 

\medskip 

We have a commutative \emph{but non-Cartesian} diagram
$$
\CD
\Op^{\on{reg}}_\cG @>{\iota}>> \Op^{\on{mer}}_\cG \\
@V{\fr^{\on{reg}}}VV @VV{\fr}V \\
\LS^{\on{reg}}_\cG@>>> \LS^{\on{mer}}_\cG.
\endCD
$$

\sssec{}

We define the factorization ind-scheme of \emph{monodromy-free opers} as the fiber product
$$\Op^{\on{mon-free}}_\cG:=\LS^{\on{reg}}_\cG\underset{\LS^{\on{mer}}_\cG}\times \Op^{\on{mer}}_\cG.$$

I.e., for a fixed $\ul{x}\in \Ran$, the fiber $\Op^{\on{mon-free}}_{\cG,\ul{x}}$ is the fiber product
$$\Op^\mf_{\cG,\ul{x}}:=\LS^\reg_{\cG,\ul{x}} \underset{\LS^\mer_{\cG,\ul{x}}}\times \Op^\mer_{\cG,\ul{x}}.$$

Denote by $\iota^\mf$ and $\iota^{+,\mf}$ the resulting maps
$$\Op^\reg \overset{\iota^{+,\mf}}\longrightarrow
\Op^{\on{mon-free}}_\cG\overset{\iota^\mf}\to \Op^{\on{mer}}_\cG.$$

%
%
%
%
%
%
%

\sssec{} \label{sss:Op is a torsor}

Recall also that the D-scheme $\Op_\cG$ is acted on simply transitively by the D-scheme $\on{Jets}(\fa(\cg)_{\omega_X})$
(see \secref{sss:jet construction})
of jets into $\fa(\cg)_{\omega_X}$, where $\fa(\cg)\subset \cg$ is the centralizer of a regular nilpotent element, and the
twist by $\omega_X$ is performed with respect to the canonical $\BG_m$-action on $\fa(\cg)$ (see, e.g., \cite[Sect. 3.1.9]{BD1}). 

\medskip

From here we obtain that $\Op^\reg_\cG$ (resp., $\Op^\mer_\cG$) is acted on simply-transitively by
$\fL^+(\fa(\cg)_{\omega_X})$ (resp., $\fL(\fa(\cg)_{\omega_X})$).

\sssec{}

Recall now the notion of \emph{formal smoothness} of a prestack (see, e.g., \cite[Sect. 8.1]{GaRo1}). This notion has
an evident relative variant. 

\medskip

We record the following (well-known) assertion: 

\begin{lem} \label{l:Op form smooth over LS}
The morphism $\fr:\Op^\mer_\cG\to \LS^\mer_\cG$ is formally smooth.
\end{lem}

\begin{proof} 

We will show that for any classical affine scheme $S$ and a map $\sigma:S\to \Op^\mer_\cG$
relative procotangent sheaf 
\begin{equation} \label{e:cotan Op/LS}
T^*_\sigma(\Op^\mer_\cG/\LS^\mer_\cG)\in \on{Pro}(\QCoh(S)^{\leq 0})
\end{equation}
is a \emph{Tate} vector bundle. This would imply the assertion of the lemma by \cite[Proposition 8.2.2]{GaRo1}. 

\medskip

In order to simplify the notation we will assume that $S=\on{pt}$ (and in particular, we work over a fixed 
point $\ul{x}\in \Ran$). However, we will perform the analysis in such a way that it will be clear that it works
in families. 

\medskip

Let $\CP^\Op_\cB$ be as in \secref{sss:oper bundle is fixed}. We can represent the tangent space to $\Op^\mer_{\cG,\ul{x}}$ at $\sigma$
is 
$$\on{coFib}\left((\cn\otimes \CO_{\cD^\times_{\ul{x}}})_{\CP^\Op_\cB}\overset{\nabla_\sigma}\to (\cb\otimes \omega_{\cD^\times_{\ul{x}}})_{\CP^\Op_\cB}\right),$$
where:

\begin{itemize}

\item $\CO_{\cD^\times_{\ul{x}}}$ is the ring of functions on the punctured multi-disc $\cD^\times_{\ul{x}}$;

\item $\omega_{\cD^\times_{\ul{x}}}$ is the space of 1-forms on $\cD^\times_{\ul{x}}$;

\item The notation $(-)_{\CP^\Op_\cB}$ indicates the twist by $\CP^\Op_\cB$, viewed as a $\cB$-bundle on $\cD^\times_{\ul{x}}$;

\item $\nabla_\sigma$ is the connection defined by $\sigma$.

\end{itemize}

The tangent space to $\LS^\mer_{\cG,\ul{x}}$ at the image of $\sigma$ is 
\begin{equation} \label{e:tang Op}
\on{coFib}\left((\cg\otimes \CO_{\cD^\times_{\ul{x}}})_{\CP^\Op_\cG}\overset{\nabla_\sigma}\to (\cg\otimes \omega_{\cD^\times_{\ul{x}}})_{\CP^\Op_\cG}\right).
\end{equation}

Hence, the relative tangent space along $\fr$ is
$$\on{Fib}\left((\cg/\cn)\otimes \CO_{\cD^\times_{\ul{x}}})_{\CP^\Op_\cB}\overset{\nabla_\sigma}\to ((\cg/\cb)\otimes \omega_{\cD^\times_{\ul{x}}})_{\CP^\Op_\cB}\right).$$

Thus, choosing a non-degenerate invariant form on $\cg$, we can identify the \emph{cotangent} space at $\sigma$ again with\footnote{
The fact that the tangent space and the relative cotangent space to opers are isomorphic is no coincidence: it reflects the interaction
of the Poisson structure on $\Op^\mer_{\cG,\ul{x}}$ and the symplectic structure on $\LS^\mer_{\cG,\ul{x}}$. In fact, the morphism $\fr$
is Lagrangian.}
\begin{equation} \label{e:cotang Op/LS}
\on{coFib}\left((\cn\otimes \CO_{\cD^\times_{\ul{x}}})_{\CP^\Op_\cB}\overset{\nabla_\sigma}\to (\cb\otimes \omega_{\cD^\times_{\ul{x}}})_{\CP^\Op_\cB}\right).
\end{equation}

Thus, we need to show that \eqref{e:cotang Op/LS} is indeed a Tate vector bundle (in degree $0$). Consider the principal
filtration on $\cg$, and the induced filtrations on $\cn$ and $\cb$. The map in \eqref{e:cotang Op/LS} sends 
$$\on{Fil}_i(\cn\otimes \CO_{\cD^\times_{\ul{x}}})_{\CP^\Op_\cB}\overset{\nabla_\sigma}\to \on{Fil}_{i-1}(\cb\otimes \omega_{\cD^\times_{\ul{x}}})_{\CP^\Op_\cB}.$$

It suffices to show that for every $i$, 
\begin{equation} \label{e:cotang Op/LS gr i}
\on{coFib}\left(
\on{gr}_i(\cn\otimes \CO_{\cD^\times_{\ul{x}}})_{\CP^\Op_\cB}\overset{\nabla_\sigma}\to \on{gr}_{i-1}(\cb\otimes \omega_{\cD^\times_{\ul{x}}})_{\CP^\Op_\cB}\right)
\end{equation}
is a Tate vector bundle (in degree $0$).

\medskip

However, the latter is evident. In fact, the maps in \eqref{e:cotang Op/LS gr i} are independent of $\sigma$, 
and the assertion follows fro the fact that the maps
$$\on{gr}_i(\cn) \overset{\on{ad}_f}\to \on{gr}_{i-1}(\cb)$$
are injective (where $f$ is a negative principal nilpotent, fixed as an element in $\cg/\cb$).

\end{proof} 

\sssec{}

As a corollary of lemmas \ref{l:Op form smooth over LS} and \ref{l:LS formall smooth}, we obtain:

\begin{cor} \label{c:mf formally smooth}
The ind-scheme $\Op^\mf_{\cG,\ul{x}}$ is formally smooth.\footnote{See \secref{sss:local prop fact} for what
formal smoothness means in the factorization setting.}
\end{cor} 

\ssec{The \texorpdfstring{$\IndCoh^*$}{IndCoh*} categories} \label{ss:IndCoh op}

In this subsection, for expositional purposes, we will work over a fixed point $\ul{x}\in \Ran$.
However, the entire discussion works when $\ul{x}$ forms a family over $\Ran$.

\sssec{}

First, since $\Op^\reg_{\cG,\ul{x}}$ (resp., $\Op^\mer_{\cG,\ul{x}}$, $\Op^\mf_{\cG,\ul{x}}$)
is a scheme (resp., ind-scheme), we have a priori defined categories
$\IndCoh^*(-)$ and $\IndCoh^!(-)$ attached to them, see Sects. \ref{sss:IndCoh non-placid !} and \ref{sss:IndCoh non-placid *}
(the Ran space version is discussed in Sects. \ref{sss:IndCoh ! Ran} and \ref{sss:IndCoh * Ran}).


\sssec{} \label{sss:Op pro-smooth}

Using \secref{sss:Op is a torsor}, we can write 
$$\Op^\reg_{\cG,\ul{x}}\simeq \underset{\bL}{\on{lim}}\, \Op^\reg_{\cG,\ul{x}}/\bL,$$
where $\bL$ runs over the filtered poset of lattices in $\fL^+(\fa(\cg)_{\omega_X})_{\ul{x}}$, viewed as a Tate
vector space. 

\medskip

This exhibits $\Op^\reg_{\cG,\ul{x}}$ as a limit of smooth schemes with smooth transition maps. 

\sssec{}

In particular, we obtain that 
$\Op^\reg_{\cG,\ul{x}}$ is \emph{placid} (see \secref{sss:placid} for what this means), so that the categories 
$$\IndCoh_*(\Op^\reg_{\cG,\ul{x}}) \text{ and } \IndCoh^!(\Op^\reg_{\cG,\ul{x}})$$
are well-behaved; in particular, they are both compactly generated and are mutually dual. 

\medskip

Note, however, that since $\Op^\reg_{\cG,\ul{x}}$ is pro-smooth, the coarsening functor
$$\Psi_{\Op^\reg_{\cG,\ul{x}}}: \IndCoh^*(\Op^\reg_{\cG,\ul{x}})\to \QCoh(\Op^\reg_{\cG,\ul{x}})$$
is an equivalence, as is the functor
$$\Upsilon_{\Op^\reg_{\cG,\ul{x}}}: \QCoh(\Op^\reg_{\cG,\ul{x}})\to \IndCoh^!(\Op^\reg_{\cG,\ul{x}}).$$

\sssec{} \label{sss:Op placid}

We can identify
$$\Op^\mer_{\cG,\ul{x}} \simeq (\Op^\reg_{\cG,\ul{x}}\times \fL(\fa(\cg)_{\omega_X})_{\ul{x}})/\fL^+(\fa(\cg)_{\omega_X})_{\ul{x}}.$$

In particular, we have a pro-smooth projection
$$\Op^\mer_{\cG,\ul{x}}\to \Op^\mer_{\cG,\ul{x}}/\fL^+(\fa(\cg)_{\omega_X})_{\ul{x}}\simeq 
\fL(\fa(\cg)_{\omega_X})_{\ul{x}}/\fL^+(\fa(\cg)_{\omega_X})_{\ul{x}}.$$

This exhibits $\Op^\mer_{\cG,\ul{x}}$ as an \emph{ind-placid} ind-scheme (see \secref{sss:ind-placid} for what this means). 

\sssec{}

This ensures that the categories 
$$\IndCoh^*(\Op^\mer_{\cG,\ul{x}}) \text{ and } \IndCoh^!(\Op^\mer_{\cG,\ul{x}})$$
are well-behaved; in particular, they are both compactly generated and are mutually dual. 

\sssec{} \label{sss:mf placid}

Note that the map
$$\LS^\reg_{\cG,\ul{x}}\to \LS^\mer_{\cG,\ul{x}}$$
is an ind-closed embedding, locally almost of finite presentation (see \lemref{l:LS reg to LS mer}). 

\medskip

This implies that 
$$\Op^\mf_{\cG,\ul{x}}\overset{\iota^\mf}\to \Op^\mer_{\cG,\ul{x}}$$
is also an ind-closed embedding locally almost of finite presentation. In particular, since $\Op^\mer_{\cG,\ul{x}}$
is ind-placid, we obtain that $\Op^\mf_{\cG,\ul{x}}$ is also ind-placid
(see \corref{c:placid inherited}). 

\medskip

Hence, we obtain that the categories 
$$\IndCoh^*(\Op^\mf_{\cG,\ul{x}}) \text{ and } \IndCoh^!(\Op^\mf_{\cG,\ul{x}})$$
are also well-behaved; in particular, they are both compactly generated and are mutually dual. 

\sssec{} \label{sss:Op mon-free expl} 

We will now show how one 
can explicitly exhibit $\Op^\mf_{\cG,\ul{x}}$ as a colimit of placid schemes. 

\medskip

To simplify the notation, 
we will assume that $\ul{x}$ consists of a single point $x$. Henceforth, we will omit $x$
from the subscript, so we will write $\cD$ instead of $\cD_{\ul{x}}$. 

\medskip

Recall (see \secref{sss:oper bundle is fixed}) that opers can be thought of as connections on a fixed $\cG$-bundle $\CP^{\Op}_\cG$. 
Trivializing this bundle on $\cD$, we will think of opers as connection forms, to be denoted $\alpha$. 

\medskip

Thus, we can write $\Op_\cG^\mf(\cD^\times)$ as 
\begin{equation} \label{e:mon-free leq n bis}
\{(\alpha\in \Op_\cG(\cD^\times),g\in \fL(\cG))\,|\, g\cdot \alpha \in \cg\otimes \omega_{\cD^\times_{\ul{x}}}\}/\fL^+(\cG)\subset \Op_\cG(\cD^\times)\times \Gr_\cG,
\end{equation} 

Therefore, we can write $\Op_\cG^\mf(\cD^\times)$ as\footnote{In the next formula, the symbol $``\on{colim}"$ (i.e., with quotes) refers to the fact that
we are forming an ind-scheme, rather than taking the colimit in the category of schemes.}
$$\underset{Y}{``\on{colim}"}\, \Op_\cG^\mf(\cD^\times)\underset{\Gr_\cG}\times Y,$$
where $Y$ runs over the filtered poset of closed subschemes of $\Gr_\cG$.

\medskip

Let us show that each $\Op_\cG^\mf(\cD^\times)\underset{\Gr_\cG}\times Y$ is 
a limit of schemes almost of finite type with smooth transition maps. 

\medskip


Let us consider $\Op_\cG(\cD^\times)$ as acted on by $\fL(\fa(\cg)_{\omega_X})$. Now, it is clear that for a fixed closed subscheme $Y\subset \Gr_\cG$,
there is an action of \emph{any} small enough lattice $\bL\subset \fL^+(\fa(\cg)_{\omega_X})$ on 
$$\Op_\cG^\mf(\cD^\times)\underset{\Gr_\cG}\times Y$$
via
$$(\alpha,g)\mapsto (\alpha+\alpha_0,g), \quad \alpha_0\in \fL^+(\fa(\cg)_{\omega_X}).$$

Further, it easy to see that for any such $\bL$, the quotient of \eqref{e:mon-free leq n bis} by it is locally almost of finite type. 

\ssec{Properties of the \texorpdfstring{$(\iota^\mf)^\IndCoh_*$}{iota} functor}

\sssec{}

The map $\iota^\mf$ gives rise to the $\IndCoh$-pushforward functor
\begin{equation} \label{e:mon-free to all cons}
(\iota^\mf)^\IndCoh_*:\IndCoh^*(\Op^\mf_{\cG,\ul{x}})\to \IndCoh^*(\Op^\mer_{\cG,\ul{x}}).
\end{equation} 

The functor $(\iota^\mf)^\IndCoh_*$ is t-exact with respect to the natural t-structures. 

\medskip

Note that since $\iota^\mf$ is a closed embedding locally almost of finite presentation, the right adjoint $(\iota^\mf)^!$ of $(\iota^\mf)^\IndCoh_*$ is continuous, 
in particular, the functor $(\iota^\mf)^\IndCoh_*$ preserves compactness, see \secref{sss:closed emb placid}. 

\sssec{} \label{sss:! dual of *}

By the same logic, the functor
$$(\iota^\mf)^!: \IndCoh^!(\Op^\mer_{\cG,\ul{x}})\to \IndCoh^!(\Op^\mf_{\cG,\ul{x}})$$
admit a left adjoint, to be denoted also by
$$(\iota^\mf)^\IndCoh_*:\IndCoh^!(\Op^\mf_{\cG,\ul{x}})\to \IndCoh^!(\Op^\mer_{\cG,\ul{x}}).$$

(We allow ourselves to use the same symbol $(\iota^\mf)^\IndCoh_*$ in both instances, as the two are unlikely to be confused.)

\sssec{}

The adjoint pairs
$$(\iota^\mf)^\IndCoh_*:\IndCoh^*(\Op^\mf_{\cG,\ul{x}})\rightleftarrows \IndCoh^*(\Op^\mer_{\cG,\ul{x}}):(\iota^\mf)^!$$
and 
$$(\iota^\mf)^\IndCoh_*:\IndCoh^!(\Op^\mf_{\cG,\ul{x}})\rightleftarrows \IndCoh^!(\Op^\mer_{\cG,\ul{x}}):(\iota^\mf)^!$$
are dual to one another with respect to the identifications
$$\IndCoh^*(\Op^\mf_{\cG,\ul{x}})^\vee \simeq \IndCoh^!(\Op^\mf_{\cG,\ul{x}})$$
 and 
$$\IndCoh^*(\Op^\mer_{\cG,\ul{x}})^\vee \simeq \IndCoh^!(\Op^\mer_{\cG,\ul{x}}).$$

\sssec{}

We will prove:

\begin{prop} \label{p:mon-free to all} \hfill

\smallskip

\noindent{\em(a)} 
The functor \eqref{e:mon-free to all cons} is conservative.

\smallskip

\noindent{\em(b)} An object of $\IndCoh^*(\Op^\mf_{\cG,\ul{x}})$ is compact if (and only if)
its image under $(\iota^\mf)^\IndCoh_*$ is compact. 

\end{prop}

The rest of this subsection is devoted to the proof of this proposition.



\sssec{}

Let 
$$(\Op^\mer_{\cG,\ul{x}})^\wedge_{\on{mon-free}}$$
denote the formal completion of $\Op^\mer_{\cG,\ul{x}}$ along $\Op^\mf_{\cG,\ul{x}}$.

\medskip

Since $\iota^\mf$ is (locally) almost of finite presentation, so is the embedding
$$(\iota^\mf)^\wedge:(\Op^\mer_{\cG,\ul{x}})^\wedge_{\on{mon-free}}\hookrightarrow \Op^\mer_{\cG,\ul{x}}.$$

In particular, $(\Op^\mer_{\cG,\ul{x}})^\wedge_{\on{mon-free}}$ is ind-placid, and we have a well-behaved
category 
$$\IndCoh^*((\Op^\mer_{\cG,\ul{x}})^\wedge_{\on{mon-free}}).$$


\sssec{}

The functor 
$$((\iota^\mf)^\wedge)^\IndCoh_*:\IndCoh^*((\Op^\mer_{\cG,\ul{x}})^\wedge_{\on{mon-free}})\to \IndCoh^*(\Op^\mer_{\cG,\ul{x}})$$ gives rise to an equivalence
$$\IndCoh^*((\Op^\mer_{\cG,\ul{x}})^\wedge_{\on{mon-free}})\overset{\sim}\to 
\IndCoh^*(\Op^\mer_{\cG,\ul{x}})_{\on{mon-free}},$$
where
$$\IndCoh^*(\Op^\mer_{\cG,\ul{x}})_{\on{mon-free}}\subset \IndCoh^*(\Op^\mer_{\cG,\ul{x}})$$
is the full subcategory of objects with set-theoretic support on $\Op^\mf_{\cG,\ul{x}}$.

\medskip

Furthermore, the functor $((\iota^\mf)^\wedge)^\IndCoh_*$ admits a continuous right adjoint, to be denoted $((\iota^\mf)^\wedge)^!$,
so $((\iota^\mf)^\wedge)^\IndCoh_*$ preserves compactness. 

\sssec{}

The functor $(\iota^\mf)^\IndCoh_*$ of \eqref{e:mon-free to all cons} factors as 
$$\IndCoh^*(\Op^\mf_{\cG,\ul{x}})\to
\IndCoh^*((\Op^\mer_{\cG,\ul{x}})^\wedge_{\on{mon-free}})\overset{((\iota^\mf)^\wedge)^\IndCoh_*}\longrightarrow 
\IndCoh^*(\Op^\mer_{\cG,\ul{x}}),$$
and in order to prove \propref{p:mon-free to all}, it is enough to establish the corresponding properties for 
the above functor
\begin{equation} \label{e:mon-free to compl cons}
({}'\!\iota^\mf)^\IndCoh_*:\IndCoh^*(\Op^\mf_{\cG,\ul{x}})\to
\IndCoh^*((\Op^\mer_{\cG,\ul{x}})^\wedge_{\on{mon-free}}).
\end{equation}

\sssec{}

Recall that according to \corref{c:mf formally smooth}, the ind-scheme $\Op^\mf_{\cG,\ul{x}}$ is formally smooth. 
Hence, its embedding into any nilpotent thickening\footnote{Assuming all ind-schemes involved are $\aleph_0$.}
admits a retraction. This implies that the embedding
\begin{equation} \label{e:mon-free to compl geom}
\Op^\mf_{\cG,\ul{x}}\hookrightarrow (\Op^\mer_{\cG,\ul{x}})^\wedge_{\on{mon-free}}
\end{equation}
admits a retraction
\begin{equation} \label{e:mon-free to compl retract}
(\Op^\mer_{\cG,\ul{x}})^\wedge_{\on{mon-free}}\to \Op^\mf_{\cG,\ul{x}}.
\end{equation}

\sssec{}

The existence of the retraction \eqref{e:mon-free to compl retract} readily implies that $({}'\!\iota^\mf)^\IndCoh_*$ is conservative: 

\medskip

Indeed, the functor of $\IndCoh$-pushforward along \eqref{e:mon-free to compl retract} is a left inverse of $({}'\!\iota^\mf)^\IndCoh_*$.

\medskip

This proves point (a) of \propref{p:mon-free to all}. 

\begin{rem} \label{r:Op mf comonadic}

Note that the above argument implies that the functor $({}'\!\iota^\mf)^\IndCoh_*$ of \eqref{e:mon-free to compl cons} is \emph{co-monadic}. Indeed,
according to what we just proved, it is conservative, and it admits a right adjoint.
Hence, it remains to check that it preserves totalizations of $({}'\!\iota^\mf)^\IndCoh_*$-split cosimplicial objects. 

\medskip

However, since $({}'\!\iota^\mf)^\IndCoh_*$ admits a left inverse, such cosimplicial objects are themselves split, and hence
their totalizations are preserved by any functor.

\end{rem} 

\sssec{}

Let $\CF\in \IndCoh^*(\Op^\mf_{\cG,\ul{x}})$ be an object, such that $(\iota^\mf)^\IndCoh_*(\CF)$ is compact.
Let us show that $\CF$ is itself compact. 

\medskip

Since $(\iota^\mf)^\IndCoh_*(\CF)$ is compact, it is cohomologically bounded, i.e., there exists an integer $n$ such that the map
$$(\iota^\mf)^\IndCoh_*(\CF)\to \tau^{\geq -n}((\iota^\mf)^\IndCoh_*(\CF))$$
is an isomorphism. 

\medskip

Since the functor $(\iota^\mf)^\IndCoh_*$ is t-exact, we obtain that 
$$(\iota^\mf)^\IndCoh_*(\CF)\to (\iota^\mf)^\IndCoh_*(\tau^{\geq -n}(\CF))$$
is an isomorphism. 

\medskip

However, since we already know that $(\iota^\mf)^\IndCoh_*$ is conservative, this implies that 
$$\CF\to \tau^{\geq -n}(\CF)$$
is an isomorphism, i.e., $\CF$ is itself cohomologically bounded. 

\medskip

Hence, it remains to check that the individual cohomologies $H^i(\CF)$ of $\CF$ are coherent.
However, 
$$(\iota^\mf)^\IndCoh_*(H^i(\CF))\simeq H^i((\iota^\mf)^\IndCoh_*(\CF)),$$
and it is easy to see that an object $\CF^0\in  \IndCoh^*(\Op^\mf_{\cG,\ul{x}})^\heartsuit$
is coherent if and only if 
$$(\iota^\mf)^\IndCoh_*(\CF^0)\in \IndCoh^*(\Op^\mer_{\cG,\ul{x}})^\heartsuit$$
is coherent (this is true for any closed embedding almost of finite presentation between ind-placid ind-schemes). 

\qed[\propref{p:mon-free to all}]

\begin{rem} 
The implication ``$(\iota^\mf)^\IndCoh_*(\CF)$ is compact" $\Rightarrow$ ``$\CF$ is compact" can also
be proved using the retraction \eqref{e:mon-free to compl retract}: since this map is ind-finite, 
the functor IndCoh-pushforward along it preserves compactness.
\end{rem} 

\begin{rem}
Note that the existence of a retraction implies that the ind-scheme
$\Op^\mf_{\cG,\ul{x}}$ is \emph{classical}.
Indeed, the ind-scheme $\Op^\mer_{\cG,\ul{x}}$ is classical, and hence so is its formal
completion $(\Op^\mer_{\cG,\ul{x}})^\wedge_{\on{mon-free}}$ (the latter follows from
\cite[Proposition 6.8.2]{GaRo1}: we reduce to the Noetherian situation using placidity). 

\end{rem}

\begin{rem}

The contents of this subsection apply ``as-is" when $\ul{x}$ forms family over $\Ran$.
In particular, the functor $(\iota^\mf)^\IndCoh_*$ has a natural factorization structure.

\medskip

Moreover, when viewed as a functor between unital factorization categories, $(\iota^\mf)^\IndCoh_*$
has a natural lax unital factorization structure. 

\medskip

We claim, however, that this unital structure is actually
strict. Indeed, this follows from the fact that $(\iota^\mf)^\IndCoh_*$ sends
$$\CO^\reg_\cG=\one_{\IndCoh^*(\Op^\mf_\cG)}\to \one_{\IndCoh^*(\Op^\mer_\cG)}=\CO^\reg_\cG,$$
see \lemref{l:unit determines strict}.


%
%
%

\end{rem} 

%
%
%
%

\ssec{A direct product decomposition} \label{ss:form neighb of mon-free}

\sssec{}

We will show that we actually have a (non-canonical) isomorphism
\begin{equation} \label{e:formal neighb as prod}
(\Op^\mer_{\cG,\ul{x}})^\wedge_{\on{mon-free}}\simeq 
\Op^\mf_{\cG,\ul{x}} \times (\cg^\wedge_0)^n, \quad |\ul{x}|=n
\end{equation} 
(here $\cg^\wedge_0$ is the formal completion of $\cg$ at $0$), 
so that \eqref{e:mon-free to compl geom} identifies with the base change of
$$0\to \cg^\wedge_0.$$

\begin{rem}

The material in this subsection is specific to the situation over a given $\ul{x}\in \Ran$. I.e., we do not know what
how to even formulate the corresponding statement over $\Ran$ (or even $X^n$ for $n\geq 2$). 

\end{rem}

\sssec{}

With no restriction of generality, we can assume that $\ul{x}$ consists of a single point $x$. Henceforth in this proof,
we will drop the subscript ``$x$" and simply write $\cD$ instead of $\cD_x$.

\medskip

We have
$$\Op_\cG(\cD^\times)^\wedge_{\on{mon-free}} 
\simeq \Op_\cG(\cD^\times)  \underset{\LS_\cG(\cD^\times)}\times\LS_\cG(\cD^\times)^\wedge_{\on{reg}},$$
where $\LS_\cG(\cD^\times)^\wedge_{\on{reg}}$ is the formal completion of $\LS_\cG(\cD^\times)$ along 
$\LS_\cG(\cD)$.

\medskip

Note that we can identify 
$$\LS_\cG(\cD^\times)^\wedge_{\on{reg}}\simeq \cg^\wedge_0/\on{Ad}(\cG),$$
so that
$$
\Op^\mf_\cG(\cD)
\simeq \Op_\cG(\cD^\times)^\wedge_{\on{mon-free}} \underset{\cg^\wedge_0/\on{Ad}(\cG)}\times \on{pt}/\cG.$$

\sssec{}

Note that the $\cG$-bundle on $\Op_\cG^\mf(\cD^\times)$ corresponding to the map
$$\Op_\cG^\mf(\cD^\times)\overset{\fr}\to \LS_\cG(\cD)\simeq \on{pt}/\cG$$
can be (non-canonically) trivialized, see \secref{sss:oper bundle is fixed}. I.e., the map
$$\Op_\cG^\mf(\cD^\times) \to \Op_\cG(\cD^\times)^\wedge_{\on{mon-free}} 
\overset{\fr}\to
\LS_\cG(\cD^\times)^\wedge_{\on{reg}} \simeq \cg^\wedge_0/\on{Ad}(\cG)\to \on{pt}/\cG$$
factors though a map
$$\Op_\cG^\mf(\cD^\times) \to \on{pt}.$$

Hence, so does the map
$$\Op_\cG(\cD^\times)^\wedge_{\on{mon-free}}  \overset{\fr}\to
\LS_\cG(\cD^\times)^\wedge_{\on{reg}} \simeq \cg^\wedge_0/\on{Ad}(\cG)\to \on{pt}/\cG.$$

Hence, the map
$$\Op_\cG(\cD^\times)^\wedge_{\on{mon-free}} \overset{\fr}\to
\LS_\cG(\cD^\times)^\wedge_{\on{reg}} \simeq \cg^\wedge_0/\on{Ad}(\cG)$$
can be (non-canonically) lifted to a map 
\begin{equation} \label{e:form compl to residue}
\Op_\cG(\cD^\times)^\wedge_{\on{mon-free}} \to \cg^\wedge_0.
\end{equation} 

\sssec{}

Combining the maps \eqref{e:mon-free to compl retract} and \eqref{e:form compl to residue}, we obtain a map
\begin{equation} \label{e:formal neighb as prod constr}
\Op_\cG(\cD^\times)^\wedge_{\on{mon-free}}\to 
\Op_\cG^\mf(\cD^\times) \times \cg^\wedge_0,
\end{equation} 
such that if we base change both sides with respect to $0\to \cg^\wedge_0$, we obtain the identity map on 
$\Op_\cG^\mf(\cD^\times)$.

\medskip

Thus, the map \eqref{e:formal neighb as prod constr} becomes an isomorphism after a base change
by a nil-isomorphism. This implies that the map \eqref{e:formal neighb as prod constr} is itself an isomorphism.

\ssec{The action of \texorpdfstring{$\IndCoh^!$}{IndCoh!}}

As in \secref{ss:IndCoh op}, for expositional purposes, we will work over a fixed point $\ul{x}\in \Ran$.
However, the entire discussion works when $\ul{x}$ forms a family over $\Ran$. 

\sssec{}

Recall that the category $\IndCoh^*(\CY)$ of an ind-scheme $\CY$ is naturally acted on by $\IndCoh^!(\CY)$. For a morphism 
$f:\CY_1\to \CY_2$, the corresponding functor $f_*:\IndCoh^*(\CY_1)\to \IndCoh^*(\CY_2)$ is $\IndCoh^!(\CY_2)$-linear,
where $\IndCoh^!(\CY_2)$ acts on $\IndCoh^*(\CY_1)$ via $f^!:\IndCoh^!(\CY_2)\to \IndCoh^!(\CY_1)$, see \secref{sss:! acts on *}.  

\medskip

In particular, we obtain that the category $\IndCoh^*(\Op^\mer_{\cG,\ul{x}})$ (resp., $\IndCoh^*(\Op^\mf_{\cG,\ul{x}})$)
is acted on by $\IndCoh^!(\Op^\mer_{\cG,\ul{x}})$ (resp., $\IndCoh^!(\Op^\mf_{\cG,\ul{x}})$), and the functor
$$(\iota^\mf)^\IndCoh_*:\IndCoh^*(\Op^\mf_{\cG,\ul{x}})\to \IndCoh^*(\Op^\mer_{\cG,\ul{x}})$$
is linear with respect to $\IndCoh^!(\Op^\mer_{\cG,\ul{x}})$.

\sssec{}

Being the right adjoint of a $\IndCoh^!(\Op^\mer_{\cG,\ul{x}})$-linear functor, the functor 
$$(\iota^\mf)^!:\IndCoh^*(\Op^\mer_{\cG,\ul{x}})\to \IndCoh^*(\Op^\mf_{\cG,\ul{x}})$$
is right-lax $\IndCoh^!(\Op^\mer_{\cG,\ul{x}})$-linear .

\medskip

However, it is easy to see that this right-lax $\IndCoh^!(\Op^\mer_{\cG,\ul{x}})$-linearity structure is actually strict, 
i.e., the adjunction
\begin{equation} \label{e:iota adj linear}
(\iota^\mf)^\IndCoh_*:\IndCoh^*(\Op^\mf_{\cG,\ul{x}})\rightleftarrows \IndCoh^*(\Op^\mer_{\cG,\ul{x}}):(\iota^\mf)^!
\end{equation}
takes place in the 2-category of $\IndCoh^!(\Op^\mer_{\cG,\ul{x}})$-module categories.

\sssec{}

Since $\IndCoh^!(\Op^\mer_{\cG,\ul{x}})$ (resp., $\IndCoh^!(\Op^\mf_{\cG,\ul{x}})$) is
\emph{symmetric} monoidal, we can view the dual categories
$$\IndCoh^*(\Op^\mer_{\cG,\ul{x}})^\vee  \text{ and } \IndCoh^*(\Op^\mf_{\cG,\ul{x}})^\vee$$
as modules over $\IndCoh^!(\Op^\mer_{\cG,\ul{x}})$ and $\IndCoh^!(\Op^\mf_{\cG,\ul{x}})$,
respectively.

\medskip

Note that the identifications
$$\IndCoh^*(\Op^\mer_{\cG,\ul{x}})^\vee \simeq \IndCoh^!(\Op^\mer_{\cG,\ul{x}})$$
and 
$$\IndCoh^*(\Op^\mf_{\cG,\ul{x}})^\vee \simeq \IndCoh^!(\Op^\mf_{\cG,\ul{x}})$$
are compatible with the 
$\IndCoh^!(\Op^\mer_{\cG,\ul{x}})$- and $\IndCoh^!(\Op^\mf_{\cG,\ul{x}})$-actions, respectively. 

\medskip

Recall (see \secref{sss:! dual of *}) that the dual of the adjunction \eqref{e:iota adj linear} identifies with
\begin{equation} \label{e:iota ! adj linear}
(\iota^\mf)^\IndCoh_*:\IndCoh^!(\Op^\mf_{\cG,\ul{x}})\rightleftarrows \IndCoh^!(\Op^\mer_{\cG,\ul{x}}):(\iota^\mf)^!
\end{equation}

It is easy to see that the resulting $\IndCoh^!(\Op^\mer_{\cG,\ul{x}})$-linear structure on \eqref{e:iota ! adj linear} arising from
the $\IndCoh^!(\Op^\mer_{\cG,\ul{x}})$-linear structure on \eqref{e:iota adj linear} is the natural $\IndCoh^!(\Op^\mer_{\cG,\ul{x}})$-linear structure
on 
$$(\iota^\mf)^!: \IndCoh^!(\Op^\mer_{\cG,\ul{x}})\to \IndCoh^!(\Op^\mf_{\cG,\ul{x}}).$$

\sssec{}

Being a $\IndCoh^!(\Op^\mer_{\cG,\ul{x}})$-linear, the functor  
$$(\iota^\mf)^!:\IndCoh^*(\Op^\mer_{\cG,\ul{x}})\to \IndCoh^*(\Op^\mf_{\cG,\ul{x}})$$
induces a functor
\begin{equation} \label{e:ten prod category}
\IndCoh^!(\Op^\mf_{\cG,\ul{x}})\underset{\IndCoh^!(\Op^\mer_{\cG,\ul{x}})}\otimes 
\IndCoh^*(\Op^\mer_{\cG,\ul{x}})\to \IndCoh^*(\Op^\mf_{\cG,\ul{x}}).
\end{equation} 

\medskip

We will prove:

\begin{prop} \label{p:ten prod property}
The functor \eqref{e:ten prod category} is an equivalence.
\end{prop}

The proof will be given in \secref{sss:proof of ten prod property}.

\sssec{}

Note that in addition to the functor \eqref{e:ten prod category}, we have a tautologically defined functor
\begin{equation} \label{e:hom category}
\IndCoh^*(\Op^\mf_{\cG,\ul{x}})
\to \on{Funct}_{\IndCoh^!(\Op^\mer_{\cG,\ul{x}})}(\IndCoh^!(\Op^\mf_{\cG,\ul{x}}),\IndCoh^*(\Op^\mer_{\cG,\ul{x}})).
\end{equation}

We claim:

\begin{lem} \label{l:hom category}
The functor \eqref{e:hom category} is an equivalence.
\end{lem}

\begin{proof}

Recall that the category $\IndCoh^*(\Op^\mer_{\cG,\ul{x}})$ identifies with the dual of $\IndCoh^!(\Op^\mer_{\cG,\ul{x}})$,
and this identification is compatible with the $\IndCoh^!(\Op^\mer_{\cG,\ul{x}})$-module structures. Therefore, for any
$\IndCoh^!(\Op^\mer_{\cG,\ul{x}})$-module category $\bC$, we have
$$\on{Funct}_{\IndCoh^!(\Op^\mer_{\cG,\ul{x}})}(\bC,\IndCoh^*(\Op^\mer_{\cG,\ul{x}}))\simeq \on{Funct}(\bC,\Vect),$$
where $\on{Funct}(-,-)$ refers to colimit-preserving functors. 

\medskip

Applying this to $\bC=\IndCoh^!(\Op^\mf_{\cG,\ul{x}})$, we obtain that the right-hand side on \eqref{e:hom category}
identifies with
$$\IndCoh^!(\Op^\mf_{\cG,\ul{x}})^\vee\simeq \IndCoh^*(\Op^\mf_{\cG,\ul{x}}).$$

It is easy to see, however, that the endomorphism of $\IndCoh^*(\Op^\mf_{\cG,\ul{x}})$, induced by
\eqref{e:hom category} and the above identfication, is the identity functor.

\end{proof} 

\ssec{Action of the spherical category}

The discussion in this subsection will be specific to the situation when $\ul{x}\in \Ran$ is fixed. The generalization in the factorization
setting will be discussed in \secref{ss:Sph spec acts on mon-free}.

\sssec{} \label{sss:act sph mf}

Let us write $\Op^\mf_{\cG,\ul{x}}$ as
$$(\Op^\mer_{\cG,\ul{x}})^\wedge_{\on{mon-free}}\underset{(\LS^\mer_{\cG,\ul{x}})^\wedge_{\on{reg}}}\times
\LS^\reg_{\cG,\ul{x}}\simeq (\Op^\mer_{\cG,\ul{x}})^\wedge_{\on{mon-free}}\underset{(\cg^\wedge_0/\on{Ad}(\cG))^{|\ul{x}|}}\times (\on{pt}/\cG)^{|\ul{x}|}.$$

From this presentation it is clear that
$$\Sph^{\on{spec}}_{\cG,\ul{x}}\simeq \IndCoh(\LS^\reg_{\cG,\ul{x}}\underset{\LS^\mer_{\cG,\ul{x}}}\times \LS^\reg_{\cG,\ul{x}})\simeq
\IndCoh(\on{pt}/\cG\underset{\cg^\wedge_0/\on{Ad}(\cG)}\times \on{pt}/\cG)^{\otimes |\ul{x}|}$$
acts on both 
$$\IndCoh^*(\Op^\mf_{\cG,\ul{x}}) \text{ and } \IndCoh^!(\Op^\mf_{\cG,\ul{x}}).$$

Moreover, these actions are $\IndCoh^!((\Op^\mer_{\cG,\ul{x}})^\wedge_{\on{mon-free}})$-linear, and hence
$\IndCoh^!(\Op^\mer_{\cG,\ul{x}})$-linear. 

\sssec{} \label{sss:Sph and duality}

The identification
$$\IndCoh^*(\Op^\mf_{\cG,\ul{x}})^\vee \simeq \IndCoh^!(\Op^\mf_{\cG,\ul{x}})$$
is compatible with the structure of $(\Sph^{\on{spec}}_{\cG,\ul{x}},\IndCoh^!(\Op^\mer_{\cG,\ul{x}}))$-bimodule on the two sides. 

\sssec{} \label{sss:sph via ten}

It follows from the constructions that the functor \eqref{e:ten prod category} (resp., \eqref{e:hom category}) respects the $\Sph^{\on{spec}}_{\cG,\ul{x}}$-actions on the two sides, 
where the action on the left-hand side of \eqref{e:ten prod category} (resp., right-hand side of \eqref{e:hom category})
is via the $\IndCoh^!(\Op^\mf_{\cG,\ul{x}})$-factor.

\sssec{}

Furthermore, it again follows from the construction that the functor
$$\QCoh(\LS^\reg_{\cG,\ul{x}})\otimes \IndCoh^*(\Op^\mf_{\cG,\ul{x}})\overset{\otimes}
\to \IndCoh^*(\Op^\mf_{\cG,\ul{x}}) \overset{(\iota^\mf)^\IndCoh_*}\to \IndCoh^*(\Op^\mer_{\cG,\ul{x}})_{\on{mon-free}}$$
canonically factors via a functor 
\begin{equation} \label{e:recover compl via Sph}
\QCoh(\LS^\reg_{\cG,\ul{x}})\underset{\Sph^{\on{spec}}_{\cG,\ul{x}}}\otimes 
\IndCoh^*(\Op^\mf_{\cG,\ul{x}})\to \IndCoh^*(\Op^\mer_{\cG,\ul{x}})_{\on{mon-free}}.
\end{equation}

We claim:

\begin{prop} \label{p:recover compl via Sph}
The functor \eqref{e:recover compl via Sph} is an equivalence.
\end{prop}

\begin{proof}

Since the monoidal categories $\Sph^{\on{spec}}_{\cG,\ul{x}}$ and $\QCoh(\LS^\reg_{\cG,\ul{x}})$ are rigid, the projection functor
\begin{multline} \label{e:monadic 1}
\IndCoh^*(\Op^\mf_{\cG,\ul{x}})\simeq 
\QCoh(\LS^\reg_{\cG,\ul{x}})\underset{\QCoh(\LS^\reg_{\cG,\ul{x}})}\otimes \IndCoh^*(\Op^\mf_{\cG,\ul{x}})\to \\
\to \QCoh(\LS^\reg_{\cG,\ul{x}})\underset{\Sph^{\on{spec}}_{\cG,\ul{x}}}\otimes 
\IndCoh^*(\Op^\mf_{\cG,\ul{x}})
\end{multline}
admits a continuous right adjoint. Moreover, the corresponding adjunction is monadic. 

\medskip

The functor 
$$((\iota^\mf)^\wedge)^!|_{\IndCoh^*(\Op^\mer_{\cG,\ul{x}})_{\on{mon-free}}}:
\IndCoh^*(\Op^\mer_{\cG,\ul{x}})_{\on{mon-free}}\to \IndCoh^*(\Op^\mf_{\cG,\ul{x}})$$
is also monadic. 

\medskip

Hence, we need to show that the functor \eqref{e:recover compl via Sph} induces an isomorphism between the two monads
acting on $\IndCoh^*(\Op^\mf_{\cG,\ul{x}})$ as plain endofunctors. 

\medskip

It is easy to see that the composition of \eqref{e:monadic 1} with \eqref{e:recover compl via Sph} is the functor
$$(\iota^\mf)^\IndCoh_*:\IndCoh^*(\Op^\mf_{\cG,\ul{x}})\to \IndCoh^*(\Op^\mer_{\cG,\ul{x}})_{\on{mon-free}}.$$

This gives rise to a map between the two monads. Let us show that this map is indeed an isomorphism of the underlying endofunctors. 

\medskip

The monad corresponding to \eqref{e:monadic 1} is given by the action on $\IndCoh^*(\Op^\mf_{\cG,\ul{x}})$
by the (algebra) object in $\Sph^{\on{spec}}_{\cG,\ul{x}}$, equal to the !-pullback of
$$\CO_{\LS^\reg_{\cG,\ul{x}}}\in \QCoh(\LS^\reg_{\cG,\ul{x}})\overset{\Psi_{\LS^\reg_{\cG,\ul{x}}}}\simeq \IndCoh(\LS^\reg_{\cG,\ul{x}})$$
along the first projection
$$p_1:\LS^\reg_{\cG,\ul{x}}\underset{\LS^\mer_{\cG,\ul{x}}}\times \LS^\reg_{\cG,\ul{x}}\to \LS^\reg_{\cG,\ul{x}}.$$

The monad corresponding to $(\iota^\mf)^\IndCoh_*$ is given by !-pull followed by *-push along
$$\Op^\mf_{\cG,\ul{x}}\leftarrow 
\Op^\mf_{\cG,\ul{x}}\underset{(\Op^\mer_{\cG,\ul{x}})^\wedge_{\on{mon-free}}}\times \Op^\mf_{\cG,\ul{x}}\to
\Op^\mf_{\cG,\ul{x}}.$$

However, it is easy to see that this functor is given by the action of the same object in $\Sph^{\on{spec}}_{\cG,\ul{x}}$. 

\end{proof}

\ssec{Self-duality for \texorpdfstring{$\IndCoh$}{IndCoh} on opers}  \label{ss:duality on Opers}

\sssec{} \label{sss:Theta op mer}

First, we claim that there is a canonically defined equivalence 
\begin{equation} \label{e:self-duality Op}
\Theta_{\Op^\mer_\cG}:\IndCoh^!(\Op^\mer_{\cG,\ul{x}})\to \IndCoh^*(\Op^\mer_{\cG,\ul{x}}),
\end{equation} 
compatible with the monoidal action of $\IndCoh^!(\Op^\mer_{\cG,\ul{x}})$ on both sides. 

\medskip

By $\IndCoh^!(\Op^\mer_{\cG,\ul{x}})$-linearity, the datum of a functor \eqref{e:self-duality Op} 
is equivalent to a choice of an object in $\IndCoh^*(\Op^\mer_{\cG,\ul{x}})$.

\medskip

The corresponding object, to be denoted 
$$\omega^{*,\on{fake}}_{\Op^\mer_{\cG,\ul{x}}}\in \IndCoh^*(\Op^\mer_{\cG,\ul{x}}),$$
is constructed as follows.

\sssec{}  \label{sss:omega fake mer}

Consider $\Op^\mer_{\cG,\ul{x}}$ as equipped with an action of $\fL^+(\fa(\cg)_{\omega_X})_{\ul{x}}$, and note that the quotient
$$\Op^\mer_{\cG,\ul{x}}/\fL^+(\fa(\cg)_{\omega_X})_{\ul{x}}$$ is an ind-scheme of ind-finite type. In particular, we have a well-defined category 
$$\IndCoh(\Op^\mer_{\cG,\ul{x}}/\fL^+(\fa(\cg)_{\omega_X})_{\ul{x}}),$$ and an object
$$\omega_{\Op^\mer_{\cG,\ul{x}}/\fL^+(\fa(\cg)_{\omega_X})_{\ul{x}}}\in \IndCoh(\Op^\mer_{\cG,\ul{x}}/\fL^+(\fa(\cg)_{\omega_X})_{\ul{x}}).$$

The operation of *-pullback along
$$\Op^\mer_{\cG,\ul{x}}\to \Op^\mer_{\cG,\ul{x}}/\fL^+(\fa(\cg)_{\omega_X})_{\ul{x}}$$
is a well-defined functor
\begin{equation} \label{e:mod out pos transl}
\IndCoh(\Op^\mer_{\cG,\ul{x}}/\fL^+(\fa(\cg)_{\omega_X})_{\ul{x}})\to \IndCoh^*(\Op^\mer_{\cG,\ul{x}}).
\end{equation} 

We let $\omega^{*,\on{fake}}_{\Op^\mer_{\cG,\ul{x}}}$ be the image of $\omega_{\Op^\mer_{\cG,\ul{x}}/\fL^+(\fa(\cg)_{\omega_X})_{\ul{x}}}$
under \eqref{e:mod out pos transl}. 

\sssec{}

We claim:

\begin{lem} \label{l:Op duality all Opers}
The functor $\Theta_{\Op^\mer_\cG}$ of \eqref{e:self-duality Op}, defined by $\omega^{*,\on{fake}}_{\Op^\mer_{\cG,\ul{x}}}$, 
is an equivalence.
\end{lem} 

\begin{proof}

We can write
$$\IndCoh^!(\Op^\mer_{\cG,\ul{x}}) \simeq \underset{\bL}{\on{colim}}\, \IndCoh(\Op^\mer_{\cG,\ul{x}}/\bL)$$
(where the transition functors are given by !-pullback)
and 
$$\IndCoh^*(\Op^\mer_{\cG,\ul{x}}) \simeq \underset{\bL}{\on{colim}}\, \IndCoh(\Op^\mer_{\cG,\ul{x}}/\bL)$$
(where the transition functors are given by *-pullback), and where the $\bL$'s run over the poset of lattices in $\fL^+(\fa(\cg)_{\omega_X})_{\ul{x}}$. 

\medskip

The functor $\Theta_{\Op^\mer_\cG}$ is given by the compatible family of (endo)functors
$$\IndCoh(\Op^\mer_{\cG,\ul{x}}/\bL)\to \IndCoh(\Op^\mer_{\cG,\ul{x}}/\bL),$$ 
each given by tensoring by the graded line 
$$\det(\fL^+(\fa(\cg)_{\omega_X})_{\ul{x}}/\bL)[-\dim(\fL^+(\fa(\cg)_{\omega_X})_{\ul{x}}/\bL)].$$

Since all these functors are equivalences, so is their colimit.

\end{proof} 

\begin{rem}
We can combine the functor $\Theta_{\Op^\mer_\cG}$ of \eqref{e:self-duality Op} with the identification
$$\IndCoh^*(\Op^\mer_{\cG,\ul{x}})^\vee \simeq \IndCoh^!(\Op^\mer_{\cG,\ul{x}})$$
and thus view it as the datum of self-duality on 
$\IndCoh^*(\Op^\reg_{\cG,\ul{x}})$.
\end{rem}

\begin{rem}
By the same token, we can define the functor
$$\Theta_{\Op^\reg_\cG}:\IndCoh^!(\Op^\reg_{\cG,\ul{x}})\to \IndCoh^*(\Op^\reg_{\cG,\ul{x}})\simeq \QCoh(\Op^\reg_{\cG,\ul{x}})$$
and show that it is an equivalence. 

\medskip

Note that the corresponding object 
$$\omega^{*,\on{fake}}_{\Op^\reg_{\cG,\ul{x}}}\in \IndCoh^*(\Op^\reg_{\cG,\ul{x}})\simeq \QCoh(\Op^\reg_{\cG,\ul{x}})$$
is $\CO_{\Op^\reg_{\cG,\ul{x}}}$. 

\end{rem} 

\sssec{} \label{sss:omega fake mf}

Let 
$$\omega^{*,\on{fake}}_{\Op^\mf_{\cG,\ul{x}}}\in \IndCoh^*(\Op^\mf_{\cG,\ul{x}})$$
be defined by
$$\omega^{*,\on{fake}}_{\Op^\mf_{\cG,\ul{x}}}:=(\iota^\mf)^!(\omega^{*,\on{fake}}_{\Op^\mer_{\cG,\ul{x}}}).$$

Let
\begin{equation} \label{e:self-duality mon free Op}
\Theta_{\Op^\mf_\cG}:\IndCoh^!(\Op^\mf_{\cG,\ul{x}})\to \IndCoh^*(\Op^\mf_{\cG,\ul{x}}),
\end{equation} 
be the $\IndCoh^!(\Op^\mf_{\cG,\ul{x}})$-linear functor, corresponding to $\omega^{*,\on{fake}}_{\Op^\mf_{\cG,\ul{x}}}$.

\sssec{}

Note that the functor $\Theta_{\Op^\mf_\cG}$ is rigged so that it makes the diagram 
\begin{equation} \label{e:duality and pullback}
\CD
\IndCoh^!(\Op^\mer_{\cG,\ul{x}})  @>{(\iota^\mf)^!}>> \IndCoh^!(\Op^\mf_{\cG,\ul{x}})   \\
@V{\Theta_{\Op^\mer_\cG}}VV @VV{\Theta_{\Op^\mf_\cG}}V \\
\IndCoh^*(\Op^\mer_{\cG,\ul{x}})  @>{(\iota^\mf)^!}>> \IndCoh^*(\Op^\mf_{\cG,\ul{x}}) 
\endCD
\end{equation} 
commute. 

%

\sssec{}

We claim:

\begin{prop} \label{p:Op duality mon-free Opers}
The functor $\Theta_{\Op^\mf_\cG}$ of \eqref{e:self-duality mon free Op} is an equivalence.
\end{prop} 

\begin{proof} 

Write 
$$\Op^\mf_{\cG,\ul{x}}=\underset{i}{``\on{colim}"}\, Y^0_i \text{ and }
\Op^\mer_{\cG,\ul{x}}=\underset{i}{``\on{colim}"}\, Y_i,$$
where $Y^0_i$ and $Y_i$ are schemes, and the map $\iota^\mf$ is given by a compatible family of maps
$$Y^0_i\overset{\iota^\mf_i}\to Y_i$$
almost of finite presentation. Moreover, we can choose $Y_i$ so that its map to $\Op^\mer_{\cG,\ul{x}}$
is almost of finite presentation. In this case both $Y_i$ and $Y^0_i$ are placid. 

\medskip

Set
$$\omega^{*,\on{fake}}_{Y^0_i}:=(\iota^\mf_i)^!(\omega^{*,\on{fake}}_{Y_i}),$$
where $\omega^{\on{fake},*}_{Y_i}$ is the !-restriction of $\omega^{*,\on{fake}}_{\Op^\mer_{\cG,\ul{x}}}$
along $Y_i\to \Op^\mer_{\cG,\ul{x}}$. Let 
$$\Theta_{Y^0_i}:\IndCoh^!(Y^0_i)\to \IndCoh^*(Y^0_i)$$
be the $\IndCoh^!(Y^0_i)$-linear functor, defined by $\omega^{*,\on{fake}}_{Y^0_i}$. We will show that each 
$\Theta_{Y^0_i}$ is an equivalence.

\medskip

Write
$$Y_i=\underset{\alpha}{\on{lim}}\, Y_{i,\alpha},$$
where $Y_\alpha$ are schemes almost of finite type, with smooth transition maps. 

\medskip

Since $\iota^\mf_i$ is almost of finite presentation (up to truncation)\footnote{The issue of truncation is taken over by passing to the limit.}, 
we can find an index $\alpha$ such that $Y^0_i$ fits into a Cartesian diagram
$$
\CD
Y^0_i @>{\iota^\mf_i}>> Y_i \\
@V{\pi^0_\alpha}VV @VV{\pi_\alpha}V \\ 
Y^0_{i,\alpha} @>{(\iota^\mf)^0_i}>> Y_{i,\alpha}.
\endCD
$$

Furthermore, up to enlarging $\alpha$, by the construction of $\omega^{*,\on{fake}}_{\Op^\mer_{\cG,\ul{x}}}$,
we can assume that
$$\omega^{\on{fake},*}_{Y_i}\simeq \pi_\alpha^*(\omega_{Y_{i,\alpha}}\otimes \CL_{i,\alpha}),$$
where $\CL_{i,\alpha}$ is a (comologically graded) line bundle on $Y_{i,\alpha}$. Hence,
$$\omega^{*,\on{fake}}_{Y^0_i}\simeq (\pi^0_\alpha)^*(\omega_{Y^0_{i,\alpha}}\otimes \CL^0_{i,\alpha}),$$
where 
$$\CL^0_{i,\alpha}:=((\iota^\mf)^0_i)^*(\CL_{i,\alpha}).$$

\medskip

Since the category of indices $\alpha$ is filtered, we can write 
$$Y^0_i\simeq \underset{\beta\geq \alpha}{\on{lim}}\, Y^0_{i,\beta}, \quad Y^0_{i,\beta}:=Y^0_{i,\alpha}\underset{Y_{i,\alpha}}\times Y_{i,\beta},$$
so that
$$\IndCoh^*(Y^0_i)\simeq \underset{\beta\geq \alpha}{\on{colim}}\, \IndCoh(Y^0_{i,\beta})$$
under *-pullbacks, and 
$$\IndCoh^!(Y^0_i)\simeq \underset{\beta\geq \alpha}{\on{colim}}\, \IndCoh(Y^0_{i,\beta})$$
under !-pullbacks. 

\medskip

For each $\beta$, let $\CL^0_{i,\beta}$ be the (canonically defined) line bundle on $Y^0_{i,\beta}$ so that
$$\omega_{Y^0_{i,\beta}}\otimes \CL^0_{i,\beta}\simeq (\pi^0_{\beta,\alpha})^*(\omega_{Y^0_{i,\alpha}}\otimes \CL^0_{i,\alpha}), \quad
\pi^0_{\beta,\alpha}: Y^0_{i,\beta}\to Y^0_{i,\alpha}.$$

We obtain that the functor $\Theta_{Y^0_i}$ is given by the compatible system of (endo)functors
$$\IndCoh(Y^0_{i,\beta})\overset{\CL^0_{i,\beta}\otimes -}\to \IndCoh(Y^0_{i,\beta}),$$
which are all equivalences. 

\end{proof} 

\begin{rem}

One could make the above proof more explicit by using the presentation of $\Op^\mf_{\cG,\ul{x}}$ as in 
\secref{sss:Op mon-free expl}.

\end{rem}

\begin{rem}
Note that combined with the identification
$$\IndCoh^*(\Op^\mf_{\cG,\ul{x}})^\vee \simeq \IndCoh^!(\Op^\mf_{\cG,\ul{x}}),$$
functor $\Theta_{\Op^\mf_\cG}$ of \eqref{e:self-duality mon free Op} can be viewed as the datum of self-duality on 
$\IndCoh^*(\Op^\mf_{\cG,\ul{x}})$.
\end{rem} 

\sssec{}

As a first consequence, passing to the right adjoint functors along the horizontal arrows in \eqref{e:duality and pullback}, and 
knowing that the vertical arrows are equivalences, we obtain another commutative diagram
$$
\CD
\IndCoh^!(\Op^\mf_{\cG,\ul{x}})   @>{(\iota^\mf)^\IndCoh_*}>> \IndCoh^!(\Op^\mer_{\cG,\ul{x}})  \\
@V{\Theta_{\Op^\mf_\cG}}VV @VV{\Theta_{\Op^\mer_\cG}}V \\
\IndCoh^*(\Op^\mf_{\cG,\ul{x}})   @>{(\iota^\mf)^\IndCoh_*}>> \IndCoh^*(\Op^\mer_{\cG,\ul{x}}). 
\endCD
$$

\sssec{Proof of \propref{p:ten prod property}} \label{sss:proof of ten prod property}

From \eqref{e:duality and pullback} we obtain a commutative diagram
$$
\CD
\IndCoh^!(\Op^\mf_{\cG,\ul{x}})\underset{\IndCoh^!(\Op^\mer_{\cG,\ul{x}})}\otimes
\IndCoh^!(\Op^\mer_{\cG,\ul{x}})  @>{(\iota^\mf)^!}>> \IndCoh^!(\Op^\mf_{\cG,\ul{x}})   \\
@V{\on{Id}\otimes \Theta_{\Op^\mer_\cG}}VV @VV{\Theta_{\Op^\mf_\cG}}V \\
\IndCoh^!(\Op^\mf_{\cG,\ul{x}})\underset{\IndCoh^!(\Op^\mer_{\cG,\ul{x}})}\otimes 
\IndCoh^*(\Op^\mer_{\cG,\ul{x}})  
@>{(\iota^\mf)^!}>> \IndCoh^*(\Op^\mf_{\cG,\ul{x}}). 
\endCD
$$

The vertical arrows in this diagram are equivalences by \lemref{l:Op duality all Opers} and \propref{p:Op duality mon-free Opers}, respectively. 
Since the top horizontal arrow is an equivalence, we obtain that so is the bottom horizontal arrow. 

\qed[\propref{p:ten prod property}]

\sssec{} \label{sss:Theta mf via duality}

Let us observe now that once we know \propref{p:ten prod property}, we could view the construction of $\Theta_{\Op^\mf_\cG}$ differently:

\medskip

We start with the equivalence
$$\IndCoh^!(\Op^\mf_{\cG,\ul{x}})\underset{\IndCoh^!(\Op^\mer_{\cG,\ul{x}})}\otimes 
\IndCoh^*(\Op^\mer_{\cG,\ul{x}}) \simeq \IndCoh^*(\Op^\mf_{\cG,\ul{x}})$$
and pass to dual categories. We obtain 
\begin{equation} \label{e:Theta from duality 0}
\on{Funct}_{\IndCoh^!(\Op^\mer_{\cG,\ul{x}})}(\IndCoh^!(\Op^\mf_{\cG,\ul{x}}),\IndCoh^!(\Op^\mer_{\cG,\ul{x}}))
\simeq \IndCoh^!(\Op^\mf_{\cG,\ul{x}}).
\end{equation} 

Applying $\Theta_{\Op^\mer_\cG}$, we replace the left-hand side in \eqref{e:Theta from duality 0} by
$$\on{Funct}_{\IndCoh^!(\Op^\mer_{\cG,\ul{x}})}(\IndCoh^!(\Op^\mf_{\cG,\ul{x}}),\IndCoh^*(\Op^\mer_{\cG,\ul{x}})),$$
and applying \lemref{l:hom category}, we rewrite it further as $\IndCoh^*(\Op^\mf_{\cG,\ul{x}})$. 

\medskip

Thus, we can interpret \eqref{e:Theta from duality 0} as an equivalence
\begin{equation} \label{e:Theta from duality 1}
\IndCoh^*(\Op^\mf_{\cG,\ul{x}})\simeq \IndCoh^!(\Op^\mf_{\cG,\ul{x}}).
\end{equation} 

It is easy to see, however, that \eqref{e:Theta from duality 1} equals the (inverse of the) equivalence $\Theta_{\Op^\mf_\cG}$
constructed above.

\sssec{}

The following results from the definition of the $\Sph^{\on{spec}}_{\cG,x}$-action on $\IndCoh^*(\Op^\mf_{\cG,\ul{x}})$ and 
$\IndCoh^*(\Op^\mf_{\cG,\ul{x}})$ in \secref{sss:act sph mf}:

\begin{lem} \label{l:Theta and Sph}
The equivalence \eqref{e:Theta from duality 1} is compatible with the $\Sph^{\on{spec}}_{\cG,x}$-actions.
\end{lem} 

\sssec{} All the preceding discussion in this subsection applies also in the factorization setting. 

\ssec{Relation to quasi-coherent sheaves}

We now consider the relationship between ind-coherent and quasi-coherent sheaves on $\Op_{\cG}^\mf$. 
We observe two pleasant categorical properties over a point and ask if they extend to the factorization setting.

\sssec{}

For any prestack $\CY$ we have a canonically defined (symmetric monoidal) functor
$$\Upsilon_\CY:\QCoh(\CY)\to \IndCoh^!(\CY).$$

\medskip

It is known that if $\CY$ is a formally smooth ind-scheme locally almost of finite type,
then $\Upsilon_\CY$ is an equivalence (see \cite[Theorem 10.1.1]{GaRo1})\footnote{This result was originally proved by J.~Lurie.}. 

\sssec{}

Consider the functors
\begin{equation}  \label{e:Ups Op}
\Upsilon_{\Op^\mer_{\cG,\ul{x}}}: \QCoh(\Op^\mer_{\cG,\ul{x}})\to \IndCoh^!(\Op^\mer_{\cG,\ul{x}})
\end{equation}
and 
\begin{equation}  \label{e:Ups mon-free Op}
\Upsilon_{\Op^\mf_{\cG,\ul{x}}}: \QCoh(\Op^\mf_{\cG,\ul{x}})\to \IndCoh^!(\Op^\mf_{\cG,\ul{x}}),
\end{equation}
respectively. 

\sssec{}

First, we claim:

\begin{lem}  \label{l:Ups Op}
The functor $\Upsilon_{\Op^\mer_{\cG,\ul{x}}}$ is an equivalence.
\end{lem}

\begin{proof} 

Repeats that of \lemref{l:Op duality all Opers}.

\end{proof}

\begin{rem}

Both the statement and the proof of \lemref{l:Ups Op} carry over to the factorization setting.

\end{rem} 

\sssec{}

We now claim:

\begin{prop}  \label{p:Ups mon-free Op}
For a fixed $\ul{x}\in \Ran$, the functor $\Upsilon_{\Op^\mf_{\cG,\ul{x}}}$ is an equivalence.
\end{prop}

\begin{proof}

We will use the direct product decomposition of \secref{ss:form neighb of mon-free}. 

\medskip

First, it is easy to see that the fact that $\Upsilon_{\Op^\mer_{\cG,\ul{x}}}$ is an equivalence implies that the functor
$$\Upsilon_{(\Op^\mer_{\cG,\ul{x}})^\wedge_{\on{mon-free}}}:
\QCoh((\Op^\mer_{\cG,\ul{x}})^\wedge_{\on{mon-free}})\to \IndCoh^!((\Op^\mer_{\cG,\ul{x}})^\wedge_{\on{mon-free}})$$
is also an equivalence. 
 
\medskip

Since the category $\QCoh(\cg^\wedge_0/\on{Ad}(\cG))$ is dualizable, the $\boxtimes$ functor
$$\QCoh(\cg^\wedge_0/\on{Ad}(\cG))\otimes  \QCoh(\Op^\mf_{\cG,\ul{x}})\to
\QCoh(\cg^\wedge_0/\on{Ad}(\cG)\times  \Op^\mf_{\cG,\ul{x}})$$
is an equivalence.

\medskip

The $\boxtimes$ functor 
$$\IndCoh(\cg^\wedge_0/\on{Ad}(\cG))\otimes  \IndCoh^!(\Op^\mf_{\cG,\ul{x}})\to
\IndCoh^!(\cg^\wedge_0/\on{Ad}(\cG)\times  \Op^\mf_{\cG,\ul{x}})$$
is an equivalence tautologically. 

\medskip

Since the functor 
$$\Upsilon_{\cg^\wedge_0/\on{Ad}(\cG)}:\QCoh(\cg^\wedge_0/\on{Ad}(\cG))\to \IndCoh(\cg^\wedge_0/\on{Ad}(\cG))$$
is an equivalence, in order to prove that $\Upsilon_{\Op^\mf_{\cG,\ul{x}}}$ is an equivalence,
it suffices to show that
$$\Upsilon_{\cg^\wedge_0/\on{Ad}(\cG)\times  \Op^\mf_{\cG,\ul{x}}}:
\QCoh(\cg^\wedge_0/\on{Ad}(\cG)\times  \Op^\mf_{\cG,\ul{x}})
\to \IndCoh^!(\cg^\wedge_0/\on{Ad}(\cG)\times  \Op^\mf_{\cG,\ul{x}})$$
is an equivalence. 

\medskip

However, this follows from the fact that $\Upsilon_{(\Op^\mer_{\cG,\ul{x}})^\wedge_{\on{mon-free}}}$ is an equivalence,
combined with the existence of an isomorphism
$$\cg^\wedge_0/\on{Ad}(\cG)\times  \Op^\mf_{\cG,\ul{x}}\simeq
(\Op^\mer_{\cG,\ul{x}})^\wedge_{\on{mon-free}}.$$

\end{proof}

\sssec{}

The proof of \propref{p:Ups mon-free Op} given above is specific to the situation when $\ul{x}\in \Ran$ is fixed. Yet, we propose:

\begin{quest} \label{q:Ups fact}
Is the functor
$$\Upsilon_{\Op^{\on{mon-free}}_\cG}: \QCoh(\Op^{\on{mon-free}}_\cG)\to \IndCoh^!(\Op^{\on{mon-free}}_\cG)$$
a factorization equivalence?
\end{quest} 

\sssec{}

Note that we can write
$$\Op^\mf_{\cG,\ul{x}}\simeq (\Op^\mer_{\cG,\ul{x}})^\wedge_{\on{mon-free}}\underset{(\LS^\mer_{\cG,\ul{x}})^\wedge_{\on{reg}}}\times
\LS^\reg_{\cG,\ul{x}}.$$

Hence, the functor of !-pullback along
$$\Op^\mf_{\cG,\ul{x}}\to (\Op^\mer_{\cG,\ul{x}})^\wedge_{\on{mon-free}}$$
gives rise to a functor
\begin{equation} \label{e:IndCoh mf as QCoh ten prod}
\QCoh(\LS^\reg_{\cG,\ul{x}})\underset{\QCoh((\LS^\mer_{\cG,\ul{x}})^\wedge_{\on{reg}})}\otimes \IndCoh^*((\Op^\mer_{\cG,\ul{x}})^\wedge_{\on{mon-free}})\to
\IndCoh^*(\Op^\mf_{\cG,\ul{x}}).
\end{equation}

We claim:

\begin{prop} \label{p:IndCoh mf as QCoh ten prod}
The functor \eqref{e:IndCoh mf as QCoh ten prod} is an equivalence.
\end{prop}

\begin{proof}

Given that the functors
$$\Upsilon_{\Op^\mf_{\cG,\ul{x}}} \text{ and }
\Upsilon_{(\Op^\mer_{\cG,\ul{x}})^\wedge_{\on{mon-free}}},$$
as well as
$$\Theta_{\Op_\cG^\mf} \text{ and } \Theta_{\Op_\cG^\mer}$$
are equivalences, in order to prove \propref{p:Ups mon-free Op}, it suffices to show that the functor
$$\QCoh(\LS^\reg_{\cG,\ul{x}})\underset{\QCoh((\LS^\mer_{\cG,\ul{x}})^\wedge_{\on{reg}})}\otimes \QCoh((\Op^\mer_{\cG,\ul{x}})^\wedge_{\on{mon-free}})\to
\QCoh(\Op^\mf_{\cG,\ul{x}})$$
is an equivalence.

\medskip

However, this follows from the fact that the prestack
$$(\LS^\mer_{\cG,\ul{x}})^\wedge_{\on{reg}}\simeq \cg^\wedge_0/\on{Ad}(\cG)$$
is \emph{passable} (see \cite[Chapter 3, Proposition 3.5.3]{GaRo3}).

\end{proof}

\sssec{}

As in the case of \propref{p:Ups mon-free Op}, the assertion of \propref{p:IndCoh mf as QCoh ten prod} is specific to 
the situation when $\ul{x}\in \Ran$ is fixed. Parallel to Question \ref{q:Ups fact}, we propose:

\begin{quest}
Is the functor 
$$\QCoh(\LS_\cG)\underset{\QCoh((\LS^{\on{mer}}_\cG)^\wedge_{\on{reg}})}\otimes \IndCoh^*(\Op^{\on{mer}}_\cG)^\wedge_{\on{mon-free}}\to
\IndCoh^*(\Op^{\on{mon-free}}_\cG)$$
a factorization equivalence? 
\end{quest}

%
%

\section{Digression: \texorpdfstring{$\IndCoh^*$}{IndCoh*} via factorization algebras} \label{s:fact alg}

In this section we discuss the approach to factorization categories arising in the local Langlands theory,
on both the geometric and spectral sides, as \emph{factorization modules} over
\emph{factorization algebras}\footnote{The factorization algebras in question may be
either plain ones (i.e., in $\Vect$) or in some simpler or better understood factorization
categories.}.

\medskip

This approach is most efficient when we want to cross the Langlands bridge, i.e., map 
a category on the geometric side and a category on the spectral side to one another. 
Indeed, it is often possible to compare the corresponding factorization algebras directly
(a prominent example of this is the Feigin-Frenkel isomorphism, see \thmref{t:FF}).

\medskip

However, this approach comes with a caveat: typically, the given representation-theoretic or
algebro-geometric category will \emph{not} be exactly equivalent to the corresponding category of 
factorization modules. Rather, the two will be related by a \emph{renormalization procedure}.
Most often, this will be manifested by the fact that both sides will be endowed with t-structures,
and the corresponding \emph{eventually coconnective} subcategories will be equivalent on the
nose. 

\medskip

We apply these ideas to construct a \emph{categorical} action of the Feigin-Frenkel center on Kac-Moody
modules at the critical level; see \secref{ss:action of center}. 

\ssec{Factorization algebras and modules}

\sssec{}

Let $\CA$ be a unital factorization algebra. To it we can attach a \emph{lax} factorization category 
$$\CA\mod^{\on{fact}}$$
of unital $\CA$-factorization modules (see \secref{sss:lax fact}).

\medskip

It comes equipped with a conservative forgetful functor
$$\oblv_\CA:\CA\mod^{\on{fact}}\to \Vect.$$

\sssec{}

By definition, the value of $\CA\mod^{\on{fact}}$ over a given $\ul{x}\in \Ran$ is the category
$$\CA\mod^{\on{fact}}_{\ul{x}}$$
of factorization $\CA$-modules at $\ul{x}$. 

\medskip

In general, we cannot say much about homological properties of the category $\CA\mod^{\on{fact}}_{\ul{x}}$.
In particular, we do not know whether it is compactly generated. 

\sssec{}

This is also reflected by the following phenomenon:

\medskip

For a pair of disjoint points $\ul{x}_1$ and $\ul{x}_2$, the lax factorization structure on $\CA\mod^{\on{fact}}$
recovers the naturally defined functor 
$$\CA\mod^{\on{fact}}_{\ul{x}_1}\otimes \CA\mod^{\on{fact}}_{\ul{x}_2}\to
\CA\mod^{\on{fact}}_{\ul{x}_1\sqcup \ul{x}_2}.$$

However, it is not clear whether this functor is an equivalence. (If it were, and if this were true in families over
$\ul{x}_1,\ul{x}_2$ moving over $\Ran$, this would mean that the lax factorization structure on $\CA\mod^{\on{fact}}$
is strict.)

\sssec{}

Assume for a moment that $\CA$ is connective, i.e., $\oblv^l(\CA_X)\in \QCoh(X)$ is connective. 
Then the category $\CA\mod^{\on{fact}}$ carries a (uniquely defined) t-structure
(see \secref{sss:t-str fact} for what this means in the factorization setting), for which the functor $\oblv_\CA$ is t-exact,
see \secref{sss:A mod is lax}. 

\medskip

In addition, it follows from the definition that $\CA\mod^{\on{fact}}$ is \emph{left-complete} in its t-structure.

\begin{rem}

The left-completeness of $\CA\mod^{\on{fact}}_{\ul{x}}$ is an indication of its failure of compact generation: 

\medskip

Let $\bC$ be a category, equipped with a t-structure, in which it is left-complete. 
Then every object $\bc\in \bC^c$ is of bounded projective dimension, i.e., the functor $\CHom_\bC(\bc,-)$ is 
of bounded cohomological amplitude. However, typically, 
the category $\CA\mod^{\on{fact}}_{\ul{x}}$ does not contain any objects of bounded projective dimension.

\end{rem}

\sssec{} \label{sss:fact alg from fact cat}

Here is how factorization algebras and modules will typically arise in this paper. Let $\bA$ be a (unital) factorization category equipped
with a (lax unital) factorization functor
$$F:\bA\to \Vect.$$

Then $F(\one_\bA)$ is a factorization algebra (in $\Vect$). Moreover, the functor $F$ naturally upgrades to a (factorization) functor, denoted
$$F^{\on{enh}}:\bA\to F(\one_\bA)\mod^{\on{fact}},$$
see \lemref{l:enhancement modules}.

\begin{rem}
The factorization structure on $F^{\on{enh}}$ means for example that for disjoint points $\ul{x}_1,\ul{x}_2\in \Ran$, the diagram 
$$
\CD
\bA_{\ul{x}_1}\otimes \bA_{\ul{x}_2} @>{\sim}>> \bA_{\ul{x}_1\sqcup \ul{x}_2} \\
@V{F_{\ul{x}_1}\otimes F_{\ul{x}_2}}VV @VV{F_{\ul{x}_1\sqcup \ul{x}_2}}V \\
F(\one_\bA)\mod^{\on{fact}}_{\ul{x}_1}\otimes F(\one_\bA)\mod^{\on{fact}}_{\ul{x}_2} @>>> F(\one_\bA)\mod^{\on{fact}}_{\ul{x}_1\sqcup \ul{x}_2}
\endCD
$$
commutes (even though the bottom horizontal arrow is not in general an equivalence).
\end{rem} 

\sssec{} \label{sss:fact alg from fact cat t}

If in the situation of \secref{sss:fact alg from fact cat}, the category $\bC$ is equipped with a t-structure so that $\one_\bC$ lies in the heart,
and the functor $F$ is t-exact, we obtain that $F(\one_\bC)$ is a \emph{classical} factorization algebra, so that the category $F(\one_\bC)\mod^{\on{fact}}$
carries a t-structure. 

\medskip

In this case, the functor $F^{\on{enh}}$ is obviously t-exact. 

\sssec{} \label{sss:fact alg from fact cat gen}

More generally, if $F:\bC_1\to \bC_2$ is a lax unital factorization functor between unital factorization categories, the object
$$F(\one_{\bC_1})\in \bC_2$$
has a natural structure of factorization algebra, and the functor $F$ upgrades to a functor
$$F^{\on{enh}}:\bC_1\to F(\one_{\bC_1})\mod^{\on{fact}}(\bC_2).$$

\ssec{Kac-Moody modules as factorization modules}

Here is a typical example of the paradigm described in Sects. \ref{sss:fact alg from fact cat}-\ref{sss:fact alg from fact cat t}. 

\sssec{} 

Consider the tautological forgetful functor
$$\oblv_{\hg}:\hg\mod_\kappa\to \Vect.$$

Recall that $\BV_{\fg,\kappa}$ denotes the factorization algebra $\oblv_\hg(\on{Vac}(G)_\kappa)$, so that $\oblv_\hg$
upgrades to a (t-exact) functor:
\begin{equation} \label{e:enh oblv KM}
\oblv_{\hg}^{\on{enh}}:\hg\mod_\kappa\to \BV_{\fg,\kappa}\mod^{\on{fact}}.
\end{equation}

\sssec{}

We have the following basic observation:

\begin{lem} \label{l:plus part KM} \hfill

\smallskip

\noindent{\em(a)} 
The functor $\oblv_{\hg}^{\on{enh}}$ of \eqref{e:enh oblv KM} induces an equivalence between the eventually
coconnective subcategories of the two sides.

\smallskip

\noindent{\em(b)} The essential image of the subcategory of compact objects of $\hg\mod_\kappa$ under is $\oblv_{\hg}^{\on{enh}}$
is contained in 
$(\BV_{\fg,\kappa}\mod^{\on{fact}})^{>-\infty}$.

\end{lem}
%

\sssec{}

The rest of this subsection is devoted to the proof of \lemref{l:plus part KM}. We will prove a pointwise version for 
$\ul{x}=x\in \Ran$. The factorization version is just a variant of this in families.

\medskip

Let $\BV^{\on{ch}}_{\fg,\kappa}$ be the \emph{chiral algebra} corresponding to $\BV_{\fg,\kappa}$. I.e., as a D-module on $X$,
$$\BV^{\on{ch}}_{\fg,\kappa}=\BV_{\fg,\kappa,X}[-1].$$

First, by \cite[Proposition 3.4.19]{BD1} (see \cite{FraG} for the derived version), we have an equivalence
$$\BV_{\fg,\kappa}\mod_x^{\on{fact}}\simeq \BV^{\on{ch}}_{\fg,\kappa}\mod^{\on{ch}}_x,$$
which commutes with the forgetful functors of both sides into $\Vect$, and hence preserves the 
t-structures on the two sides.

\sssec{}

Let $L_{\fg,\kappa}$ be the Lie-* algebra 
$$\omega_X\oplus (\fg\otimes \on{D}_X)$$
of \cite[Sect. 2.5.9]{BD1}. Then by \cite[Proposition 3.7.17]{BD1} (which applies as-is in the derived setting), we have
$$\BV^{\on{ch}}_{\fg,\kappa}\mod^{\on{ch}}_x\simeq L_{\fg,\kappa}\mod^{\on{ch}}_x,$$
which commutes with the forgetful functors of both sides into $\Vect$, and hence preserves the 
t-structures on the two sides. 

\sssec{}

By the construction of $\hg\mod_{\kappa,x}$ in \cite{Ra5} (or, equivalently, in \cite[Sect. 23.1]{FG6}),
we have an embedding
\begin{equation} \label{e:KM to Lie*}
(\hg\mod_{\kappa,x})^c\hookrightarrow (L_{\fg,\kappa}\mod^{\on{ch}}_x)^{>-\infty},
\end{equation}
which satisfies the conditions of \cite[Sect. 22.1.4]{FG6} (see \secref{sss:ass of top Lie} for an explanation of why this happens). 

\medskip

Hence,
$$(\hg\mod_{\kappa,x})^{>-\infty}\to (L_{\fg,\kappa}\mod^{\on{ch}}_x)^{>-\infty}$$
is an equivalence by \cite[Proposition 22.1.5 and 22.2.1]{FG6}.

\qed[\lemref{l:plus part KM}]

\ssec{The case of commutative factorization algebras} \label{ss:com fact}

\sssec{}

Let $\CY$ be an affine D-scheme over $X$, i.e., $\CY=\Spec_X(A)$, where
$A\in \on{ComAlg}(\Dmod(X))$ with $\oblv^l(A)\in \QCoh(X)^{\leq 0}$.

\medskip

Let $\CA\in \on{ComAlg}(\on{FactAlg}^{\on{untl}}(X))$ denote the corresponding commutative factorization 
algebra, i.e., $\CA:=\on{Fact}(A)$, see \secref{sss:com fact vs Dmod com}.

\sssec{} \label{sss:fact schemes ass to aff}

Let $\fL^+_\nabla(\CY)$ denote the affine factorization scheme corresponding to $\CY$ (see \secref{sss:forming arcs gen}), i.e.,
the fiber $\fL^+_\nabla(\CY)_{\ul{x}}$ of $\fL^+_\nabla(\CY)$ at $\ul{x}\in \Ran$ is the space
\begin{equation} \label{e:sect arcs}
\on{Sect}_\nabla(\cD_{\ul{x}},\CY).
\end{equation} 

According to \secref{sss:com fact and arcs aff}, for $\CZ\to \Ran$,
$$\fL^+_\nabla(\CY)_\CZ=\Spec_\CZ(\CA_\CZ).$$

In particular, for $\ul{x}=\{x_1,...,x_n\}$, we have 
$$\fL^+_\nabla(\CY)_{\ul{x}}\simeq \underset{i=1,...,n}\Pi\, \CY_{x_i}, \quad \ul{x}=\{x_1,...,x_n\}.$$

(Note that this is compatible with \eqref{e:sect arcs}, since for a singleton $\ul{x}=\{x\}$, we have
$\on{Sect}_\nabla(\cD_x,\CY)\simeq \CY_x$, in agreement with \eqref{e:sect arcs}). 

\sssec{}

Consider the corresponding factorization category $\QCoh(\fL^+_\nabla(\CY))$ (see \secref{sss:QCoh fact}),
so that for $\ul{x}\in \Ran$, we have 
$$\QCoh(\fL^+_\nabla(\CY))_{\ul{x}}:=\QCoh(\fL^+_\nabla(\CY)_{\ul{x}}).$$

The factorization category $\QCoh(\fL^+_\nabla(\CY))$ is unital, with the structure sheaf $\CO_{\fL^+_\nabla(\CY)}$
being the factorization unit.

\sssec{}

The functor of global sections $\Gamma(\fL^+_\nabla(\CY),-)$ sends
$$\CO_{\fL^+_\nabla(\CY)}\mapsto \CA$$ 
and induces an equivalence
$$\QCoh(\fL^+_\nabla(\CY))\simeq \CA\mod^{\on{com}}.$$

\sssec{}

Let $\fL_\nabla(\CY)$ denote the factorization D-ind-scheme that attaches to a point $\ul{x}\in \Ran$ the space
$$\fL_\nabla(\CY)_{\ul{x}}:=\on{Sect}_\nabla(\cD^\times_{\ul{x}},\CY),$$
see \secref{sss:forming loops aff}.

\medskip

We consider the corresponding lax factorization category 
$\QCoh(\fL_\nabla(\CY))$ (see \secref{sss:QCoh fact}),
so that for $\ul{x}\in \Ran$, we have 
$$\QCoh(\fL_\nabla(\CY))_{\ul{x}}:=\QCoh(\fL_\nabla(\CY)_{\ul{x}}).$$

For a general $\CY$, this category may be quite ill-behaved (basically, because the category
of quasi-coherent sheaves on an ind-scheme may be quite unwieldy); in particular, it is not
clear whether $\QCoh(\fL_\nabla(\CY))$ is unital. 

\sssec{}

Recall now that we can also consider the factorization category 
$$\QCoh_{\on{co}}(\fL_\nabla(\CY)),$$
see \secref{sss:QCoh co fact}. 

\medskip

Its factorization unit is the direct image of the structure sheaf on $\fL^+_\nabla(\CY)$
along the tautological closed embedding $\fL^+_\nabla(\CY)\overset{\iota}\to \fL_\nabla(\CY)$. 
By a slight abuse of notation, we will denote it by the same symbol $\CO_{\fL^+_\nabla(\CY)}$. 

\sssec{} \label{sss:QCohco to fact mod}

The operation of taking global sections is a (t-exact\footnote{See \secref{sss:t on QCoh co} for the definition
of the t-structure on $\QCoh_{\on{co}}$ of an ind-scheme.}) 
factorization functor 
$$\Gamma(\fL_\nabla(\CY),-):\QCoh_{\on{co}}(\fL_\nabla(\CY))\to \Vect.$$

Hence, by Sects. \ref{sss:fact alg from fact cat}-\ref{sss:fact alg from fact cat t}, the functor $\Gamma(\fL_\nabla(\CY),-)$
upgrades to a (t-exact) lax unital factorization functor
\begin{equation} \label{e:from QCoh star to fact}
\Gamma(\fL_\nabla(\CY),-)^{\on{enh}}:\QCoh_{\on{co}}(\fL_\nabla(\CY))\to \CA\mod^{\on{fact}}.
\end{equation}

\sssec{}

We have the following basic assertion: 

\begin{thm} \label{t:from QCoh star to fact} 
Assume that $\CY$ is almost finitely presented in the D-sense\footnote{See \secref{sss:afp D} for what this means.}. 
Then the functor \eqref{e:from QCoh star to fact} induces an equivalence between the eventually coconnective subcategories
of the two sides.
\end{thm} 

The proof of this theorem will be given in \secref{ss:proof of fact over com}. We note that the assertion of the theorem would be 
\emph{false} without the finite presentation hypothesis, see \secref{ss:QCoh* non-fp}.

\sssec{} \label{sss:IndCoh^* to QCoh_*}

Recall that if $\CZ$ is an ind-scheme, we have a well-defined (t-exact) functor 
$$\Psi_\CZ:\IndCoh^*(\CZ)\to \QCoh_{\on{co}}(\CZ),$$
which induces an equivalence
$$\IndCoh^*(\CZ)^{>-\infty}\overset{\sim}\to \QCoh_{\on{co}}(\CZ)^{>-\infty},$$
see \lemref{l:Psi IndCoh* IndSch}.

\medskip

Furthermore, if $\CZ$ is ind-placid, $\Psi_\CZ$ gives rise to an equivalence between 
$\IndCoh^*(\CZ)^c$ and the subcategory of almost compact objects
in $\QCoh_{\on{co}}(\CZ)^{>-\infty}$. 

\medskip

Note also that the composition 
$$\Gamma(\CZ,-)\circ \Psi_\CZ:\IndCoh^*(\CZ)\to \Vect$$
is the functor $\Gamma^\IndCoh(\CZ,-)$ of IndCoh-global sections.

\sssec{} \label{sss:fact com placid}

Assume for a moment that $\fL_\nabla(\CY)$ is ind-placid. 
Applying \secref{sss:IndCoh^* to QCoh_*}, we obtain a factorization functor
$$\Gamma^{\IndCoh}(\fL_\nabla(\CY),-)\simeq \Gamma(\fL_\nabla(\CY),-)\circ \Psi_{\fL_\nabla(\CY)},\quad \IndCoh^*(\fL_\nabla(\CY))\to \Vect$$
and its enhancement
\begin{equation} \label{e:from IndCoh * to fact}
\Gamma^{\IndCoh}(\fL_\nabla(\CY),-)^{\on{enh}}:\IndCoh^*(\fL_\nabla(\CY))\to  \CO_\CY\mod^{\on{fact}}.
\end{equation}

Combining with \thmref{t:from QCoh star to fact} we obtain:

\begin{cor} \label{c:from QCoh * to fact placid}  \hfill

\smallskip

\noindent{\em(a)} The functor $\Gamma^{\IndCoh}(\fL_\nabla(\CY),-)^{\on{enh}}$ of \eqref{e:from IndCoh * to fact} is t-exact and 
induces an equivalence between the eventually coconnective subcategories of the two sides.

\smallskip

\noindent{\em(b)} The essential image of the subcategory of compact objects in 
$\IndCoh^*(\fL_\nabla(\CY))$ under the functor $\Gamma^{\IndCoh}(\fL_\nabla(\CY),-)^{\on{enh}}$ is
contained in $(\CO_\CY\mod^{\on{fact}})^{> -\infty}$. 
\end{cor}

\begin{rem}  \label{r:recover IndCoh}

Point (b) in \corref{c:from QCoh * to fact placid} can be strengthened as follows: the essential image of 
$\IndCoh^*(\fL_\nabla(\CY))^c$ under $\Gamma^{\IndCoh}(\fL_\nabla(\CY),-)^{\on{enh}}$ equals the
category of \emph{almost compact} objects\footnote{Recall that
an object $\bc$ in a DG category $\bC$ equipped with a t-structure (assumed compatible with filtered colimits)
is said to be \emph{almost compact} if the functor $\CHom_\bC(\bc,-)$ commutes with filtered colimits 
on $\bC^{\geq -n}$ for all $n$.} in $(\CO_\CY\mod^{\on{fact}})^{> -\infty}$.

\medskip

Note that from \corref{c:from QCoh * to fact placid} allows us to recover the
(factorization) category $\IndCoh^*(\fL_\nabla(\CY))$ from $\CO_\CY\mod^{\on{fact}}$ equipped with its t-structure:

\medskip

Namely, 
we can identify $\IndCoh^*(\fL_\nabla(\CY))$ with the ind-completion of the category of almost compact objects in 
$(\CO_\CY\mod^{\on{fact}})^{\geq -\infty}$.

\end{rem}

%
%
%
%
%

%

\ssec{Recovering \texorpdfstring{$\IndCoh^*$}{IndCoh*} of opers as factorization modules} \label{ss:indcoh via fact} 

\sssec{}

The setup of \secref{sss:fact com placid} is directly applicable to the case when $\CY=\Op_\cG$
(the D-afp assumption is satisfied by \secref{sss:Op is a torsor}), so that
$$\fL^+_\nabla(\CY)=\Op_\cG^\reg \text{ and } \fL_\nabla(\CY)=\Op_\cG^\mer.$$

\medskip

By a slight abuse of notation, we will denote by $\CO_{\Op_\cG^\reg}$ (rather than 
$\Gamma(\Op_\cG^\reg,\CO_{\Op_\cG^\reg})$) the corresponding 
factorization algebra in $\Vect$. 

\medskip

In particular, we obtain that the functor
$$\Gamma^\IndCoh(\Op^{\on{mer}}_\cG,-):\IndCoh^*(\Op^{\on{mer}}_\cG)\to \Vect$$
upgrades to a (t-exact) functor 
\begin{equation} \label{e:global sect opers punct enh}
\Gamma^\IndCoh(\Op^{\on{mer}}_\cG,-)^{\on{enh}}:\IndCoh^*(\Op^{\on{mer}}_\cG)\to \CO_{\Op^\reg_\cG}\mod^{\on{fact}},
\end{equation}
and we have:

\begin{cor} \label{c:IndCoh* Op bdd below} \hfill

\smallskip

\noindent{\em(a)}
The functor \eqref{e:global sect opers punct enh} induces an equivalence between the corresponding eventually coconnective
(a.k.a. bounded below) subcategories.

\smallskip

\noindent{\em(b)}
The essential image of the subcategory of compact objects in $\IndCoh^*(\Op^\mer_{\cG,\ul{x}})$ under the functor 
\eqref{e:global sect opers punct enh} is contained in 
$(\CO_{\Op^\reg_\cG}\mod^{\on{fact}}_{\ul{x}})^{>-\infty}$.
\end{cor}


\medskip

We now consider the case of \emph{monodromy-free} opers. 

\sssec{}

Direct image along the projection 
$$\fr^\reg:\Op^\reg_\cG\to \LS^{\on{reg}}_\cG$$
defines a t-exact (lax unital factorization) functor
\begin{equation} \label{e:Op to LS reg}
\IndCoh^*(\Op^\reg_\cG)\simeq \QCoh(\Op^\reg_\cG)\overset{\fr^\reg_*}\to
\QCoh(\LS^{\on{reg}}_\cG)\simeq \Rep(\cG).
\end{equation} 

Denote
\begin{equation} \label{e:R G Op}
R_{\cG,\Op}:=\fr^\reg_*(\CO_{\Op^{\on{reg}}_\cG}).
\end{equation} 

This is naturally a commutative factorization algebra in $\Rep(\cG)$.

\sssec{} \label{sss:R G Op}

Explicitly, 
$$R_{\cG,\Op}:=\left(\Gamma(\Op^{\on{reg}}_\cG,-)\otimes \on{Id}\right)\circ \left((\fr^\reg)^*\otimes \on{Id})\right)(R_\cG),$$
where 
$$R_\cG\in \Rep(\cG)\otimes \Rep(\cG)$$
is the regular representation.

\sssec{}

Consider now the factorization functor
\begin{equation} \label{e:Op to LS}
\fr^\IndCoh_*:\IndCoh^*(\Op^{\on{mon-free}}_\cG)\to \QCoh(\LS^{\on{reg}}_\cG)\simeq \Rep(\cG).
\end{equation} 

Note that we can interpret $R_{\cG,\Op}$ also as 
$$\fr^\IndCoh_*(\CO_{\Op^{\on{reg}}_\cG}),$$
where by a slight abuse of notation we view $\CO_{\Op^{\on{reg}}_\cG}$ as an object of 
$\IndCoh^*(\Op^{\on{mon-free}}_\cG)$ using
$$\QCoh(\Op^\reg_\cG)\overset{\Psi_{\Op^\reg_\cG}}\simeq \IndCoh^*(\Op^\reg_\cG)
\overset{(\iota^{+,\mf})^\IndCoh_*}\longrightarrow \IndCoh^*(\Op^\mf_\cG).$$

\sssec{} \label{sss:indcoh via fact} 

The functor \eqref{e:Op to LS} naturally upgrades to a t-exact factorization functor
\begin{equation} \label{e:Op to LS enh}
(\fr^{\IndCoh}_*)^{\on{enh}}:\IndCoh^*(\Op^{\on{mon-free}}_\cG)\to 
R_{\cG,\Op}\mod^{\on{fact}}(\Rep(\cG)).
\end{equation} 

\medskip

%

We will prove:

\begin{prop} \label{p:IndCoh Op via fact almost} \hfill

\smallskip

\noindent{\em(a)}
The functor \eqref{e:Op to LS enh} induces an equivalence between the eventually coconnective
subcategories of the two sides. 

\smallskip

\noindent{\em(b)}
The essential image of the subcategory $\IndCoh^*(\Op^{\on{mon-free}}_\cG)^c$ under \eqref{e:Op to LS enh} is contained in 
$(R_{\cG,\Op}\mod^{\on{fact}}(\Rep(\cG)))^{>-\infty}$. 
\end{prop} 

\ssec{Proof of \propref{p:IndCoh Op via fact almost}} 

We will provide a general framework, of which \propref{p:IndCoh Op via fact almost} is a particular case. 
The assertion is local, so we can assume that $X$ is affine. 

\sssec{} \label{sss:assump i}

Let $\CY$ be an affine D-scheme, and consider $\fL^+_\nabla(\CY):=\CT^+$ as a factorization space. Let
$\CT$ be a factorization space, equipped with a map
$$\iota:\CT^+\to \CT$$
that extends to a \emph{unital-in-correspondences structure} on $\CT$ relative to $\CT^+$ (see \secref{sss:rel unital}).

\medskip

Consider the factorization category 
$$\QCoh_{\on{co}}(\CT),$$
see \secref{sss:QCoh co fact}.

\medskip

Note its factorization unit is given by
$$\iota_*(\CO_{\CT^+}).$$

We will assume:

\begin{itemize}

\item(i) $\CT$ is an ind-placid ind-scheme.\footnote{See \secref{sss:local prop fact} for what this means.}

\end{itemize}

\sssec{} \label{sss:assump ii}

Let us be given a D-prestack $\CY_0$, equipped with a map
$$f:\CY\to \CY_0,$$
and an extension of the map
$$\CT^+=\fL^+_\nabla(\CY) \overset{\fL^+(f)}\to \fL^+_\nabla(\CY_0)=:\CT^+_0$$
to a map
$$\CT \overset{\fL(f)}\to \CT^+_0.$$

We will assume: 

\begin{itemize}

\item(ii) $\CT^+_0$ has an affine diagonal. 

\end{itemize}

\medskip

Note that assumptions (i) and (ii) imply in particular that the map $\fL(f)$ is ind-schematic. 

\sssec{}

Consider the functor
\begin{equation} \label{e:f fact}
\QCoh_{\on{co}}(\CT) \overset{\fL(f)_*}\to
\QCoh_{\on{co}}(\CT^+_0) \overset{\Omega_{\fL_\nabla^+(\CY_0)}}\to \QCoh(\CT^+_0),
\end{equation}
where $\Omega_{\fL_\nabla^+(\CY_0)}$ is as in \eqref{e:from QCoh co to QCoh}. By a slight abuse of notation, we
will denote the composite functor in \eqref{e:f fact} by the same symbol $\fL(f)_*$. 

\medskip

The functor $\fL(f)_*$ of \eqref{e:f fact} upgrades to a (factorization) functor
\begin{multline} \label{e:f fact enh}
(\fL(f)_*)^{\on{enh}}:\QCoh_{\on{co}}(\CT)\to \\
\to \fL(f)_*\circ \iota_*(\CO_{\CT^+})\mod^{\on{fact}}(\QCoh(\CT^+_0))=
\fL^+(f)_*(\CO_{\CT^+})\mod^{\on{fact}}(\QCoh(\CT^+_0)). 
.\end{multline}

\sssec{} \label{sss:rel fact flat}

We now make an additional assumptions: 

\begin{itemize}

\item(iii) The prestack $\CY_0$ admits a map 
$$g:\wt\CY_0\to \CY_0,$$
where $\wt\CY_0$ is an affine D-scheme, such that the map
$$\wt\CT^+_0:=\fL^+_\nabla(\wt\CY_0)\overset{\fL^+(g)}\to \fL^+_\nabla(\CY_0)=:\CT^+_0$$
is an fpqc cover.\footnote{See \secref{sss:local prop fact} for what this means.}

\smallskip

\item(iv) For $\wt\CY:=\CY\underset{\CY_0}\times \wt\CY^0$, the resulting map
$$\wt\CT^+:=\CT^+\underset{\CT^+_0}\times \wt\CT^+_0 \overset{\iota\times \on{id}}\to \CT\underset{\CT^+_0}\times \wt\CT^+_0=:\wt\CT$$
identifies with
$$\fL^+_\nabla(\wt\CY)\to  \fL_\nabla(\wt\CY)\underset{\fL_\nabla(\wt\CY_0)}\times \fL^+_\nabla(\wt\CY_0).$$

\end{itemize} 

\bigskip

Note that assumption (iii) implies, in particular, that the category $\QCoh(\CT^+_0)$
has a well-behaved t-structure: it is characterized uniquely by the property that the functor 
$$(\fL^+(g))^*:\QCoh(\CT^+_0)\to \QCoh(\wt\CT^+_0)$$
is t-exact. 

\medskip

It follows from assumption (iv) and base change that the functor $\fL(f)_*$ of \eqref{e:f fact} is t-exact
(see \secref{sss:QCoh co fact}, where the t-structure on the left-hand side is defined). Hence, $(\fL(f)_*)^{\on{enh}}$ is also t-exact. 

\sssec{}

Finally, we make the following assumption:

\begin{itemize}

\item(v) The map 
$$\wt\CY:=\CY\underset{\CY_0}\times \wt\CY_0\overset{\wt{f}}\to \wt\CY_0$$
is D-afp (see \secref{ss:afp D} for what this means).

\end{itemize}

\sssec{}

We claim:

\begin{cor} \label{c:from QCoh * to fact rel}
Under the above assumptions, the functor \eqref{e:f fact enh} induces an equivalence between 
the eventually coconnective subcategories of the two sides.
\end{cor}

\begin{proof}

Let $\wt\CY_0^\bullet$ be the \v{C}ech nerve of the map $g$. Denote
$$\wt\CY^\bullet:=\CY\underset{\CY_0}\times \wt\CY_0^\bullet$$ and 
$$\wt\CT_0^{+,\bullet}:=\fL^+_\nabla(\wt\CY_0^\bullet), \quad 
\wt\CT^{+,\bullet}:=\CT^+\underset{\CT^+_0}\times \wt\CT_0^{+,\bullet} \simeq \fL_\nabla^+(\wt\CY^\bullet)$$
and
$$\wt\CT^\bullet:=\CT^+\underset{\CT^+_0}\times \wt\CT_0^{+,\bullet} \simeq 
\fL_\nabla(\wt\CY^\bullet)\underset{\fL_\nabla(\wt\CY_0^\bullet)}\times \fL^+_\nabla(\wt\CY_0^\bullet).$$

Consider the resulting maps:
$$\wt\CT^\bullet \overset{\fL(\wt{f})^\bullet}\to \wt\CT_0^{+,\bullet} \text{ and }
\wt\CT^{+,\bullet} \overset{\fL^+(\wt{f})^{\bullet}}\to \wt\CT_0^{+,\bullet}.$$

First, by fpqc descent we have a t-exact equivalence
$$\QCoh(\CT^+_0) \simeq \on{Tot}(\QCoh(\wt\CT_0^{+,\bullet})),$$
from which we obtain a t-exact equivalence 
$$\fL^+(f)_*(\CO_{\CT^+})\mod^{\on{fact}}(\QCoh(\CT^+_0)) \simeq
\on{Tot}\left(\fL^+(\wt{f}^\bullet)_*(\CO_{\wt\CT^{+,\bullet}})\mod^{\on{fact}}(\QCoh(\wt\CT_0^{+,\bullet}))\right)$$
and hence
$$\fL^+(f)_*(\CO_{\CT^+})\mod^{\on{fact}}(\QCoh(\CT^+_0))^{>-\infty} \simeq
\on{Tot}\left(\fL^+(\wt{f}^\bullet)_*(\CO_{\wt\CT^{+,\bullet}})\mod^{\on{fact}}(\QCoh(\wt\CT_0^{+,\bullet}))^{>-\infty}\right).$$

\medskip

Next, by assumption (i) in \secref{sss:assump i} and \propref{p:QCoh co descent}, the functor
$$\QCoh_{\on{co}}(\CT)^{>-\infty} \to \on{Tot}\left(\QCoh_{\on{co}}(\wt\CT^\bullet)^{>-\infty}\right)$$
is also an equivalence.

\medskip 

Finally, by assumption (v) and a relative version of \thmref{t:from QCoh star to fact}, the functor
\begin{multline*}
\QCoh_{\on{co}}(\wt\CT^\bullet)^{>-\infty} \simeq
\QCoh_{\on{co}}\left(\fL_\nabla(\wt\CY^\bullet)\underset{\fL_\nabla(\wt\CY_0^\bullet)}\times \fL^+_\nabla(\wt\CY_0^\bullet)\right){}^{>-\infty} \to \\
\to \fL^+(\wt{f}^\bullet)_*(\CO_{\wt\CT^{+,\bullet}})\mod^{\on{fact}}(\QCoh(\wt\CT_0^{+,\bullet}))^{>-\infty}\simeq
\fL^+(\wt{f}^\bullet)_*(\CO_{\fL_\nabla^+(\wt\CY^\bullet)})\mod^{\on{fact}}(\QCoh(\fL^+_\nabla(\wt\CY_0^\bullet)))^{>-\infty}
\end{multline*}
is a term-wise equivalence.

\medskip

Combining, we obtain that \eqref{e:f fact enh} is also an equivalence, as required.

\end{proof} 

\sssec{}

Precomposing the equivalence \eqref{e:f fact enh} with the equivalence
$$\IndCoh^*(\CT)^{>-\infty} \overset{\Psi_{\CT}}\longrightarrow 
\QCoh_{\on{co}}(\CT)^{>-\infty}$$
of \lemref{l:Psi IndCoh* IndSch}, we obtain that under assumptions (i)-(iv) above, the functor
$$((\fL(f))^\IndCoh_*)^{\on{enh}}:=(\fL(f)_*)^{\on{enh}}\circ \Psi_{\CT}$$
induces an equivalence
\begin{equation} \label{e:from QCoh star to fact rel}
\IndCoh^*(\CT)^{>-\infty} \to 
((\fL(f))^\IndCoh_*)^{\on{enh}}\mod^{\on{fact}}(\QCoh(\CT^+_0))^{>-\infty}.
\end{equation} 

\sssec{}

We apply the above to 
$$\CY=\Op_\cG, \quad \CY_0=\on{pt}/\cG \text{ and } \CT:=\Op^\mf_\cG:=\Op^\mer_\cG\underset{\LS_\cG^\mer}\times \LS_\cG^\reg.$$

\medskip

Hence, in order to deduce the assertion of \propref{p:IndCoh Op via fact almost}, we have
to show that conditions (i)-{v) above hold.

\sssec{}

Condition (i) says that $\Op_\cG^\mf$ is placid; this has been established in \secref{sss:mf placid}.

\medskip

We take $\wt\CY_0=\on{pt}$ with the tautological map $\on{pt} \to \LS^\reg_\cG$.

\medskip
 
Condition (iii) is the content of \lemref{l:pt to pt/H}. Condition (ii) also follows from \lemref{l:pt to pt/H}, since the property of a map being
affine can be checked fpqc-locally, and
$$\LS^\reg_\cG\underset{\LS^\reg_\cG\times \LS^\reg_\cG}\times (\on{pt}\times \on{pt})\simeq
\on{pt}\underset{\LS^\reg_\cG}\times \on{pt}\simeq \fL_\nabla^+(\cG).$$

\medskip 

Condition (iv) is automatic from the construction. Finally, condition (v) is the content of the next lemma: 

\begin{lem} 
The affine D-scheme $\Op_\cG\underset{\on{pt}/\cG}\times \on{pt}$ is D-afp.
\end{lem}

\begin{proof}

We have: 
$$\Op_\cG\underset{\on{pt}/\cG}\times \on{pt}\simeq \Op_\cG\underset{\on{Jets}(\cg\otimes \omega_X)}\times \on{Jets}(\cG),$$
where:

\begin{itemize}

\item The map $\Op_\cG\to \on{Jets}(\cg\otimes \omega_X)$ is well-defined (Zariski-locally on $X$) 
thanks to \secref{sss:oper bundle is fixed};

\medskip

\item The map $\on{Jets}(\cG)\to \on{Jets}(\cg\otimes \omega_X)$ is given by the gauge action
on the trivial connection.

\end{itemize}

This makes the assertion of the lemma manifest, as 
$$\Op_\cG,\,\, \on{Jets}(\cg\otimes \omega_X) \text{ and } \on{Jets}(\cG)$$
are all D-afp. 

\end{proof} 

\qed[\propref{p:IndCoh Op via fact almost}]

\ssec{Action of the center on Kac-Moody modules} \label{ss:action of center}

\sssec{}

Let $\fz_\fg$ be the (classical) center of $\BV_{\fg,\crit}$, viewed as a plain factorization (chiral) algebra
(i.e., a factorization algebra in $\Vect$). 

\medskip

By construction, $\fz_\fg$ is a commutative factorization algebra. It acts as such on $\on{Vac}(G)_\crit\in \KL(G)_\crit$. 
In particular, we obtain a map of factorization algebras in $\KL(G)_\kappa$: 
\begin{equation} \label{e:action of ff}
\fz_\fg\otimes \on{Vac}(G)_\crit\to \on{Vac}(G)_\crit.
\end{equation}

\sssec{} \label{sss:fact sch z}

We will denote by the symbols 
$\Spec(\fz_\fg)$ and $``\Spec"(\fZ_\fg)$ the corresponding factorization scheme and ind-scheme, respectively, see \secref{ss:com fact}. 


\medskip

In this subsection we will construct the an action of $\IndCoh^!(``\Spec"(\fZ_\fg))$ on $\hg\mod_\crit$, compatible with factorization.

\begin{rem}

The existence of such an action at the level of abelian categories is essentially evident: the topological algebra
of global functions on $``\Spec"(\fZ_\fg)$ maps to the center of the completed universal enveloping algebra
of $U(\hg_\crit)$. 

\medskip

At the derived level (for a fixed point $\ul{x}\in \Ran$) the construction of such an action was carried out
in \cite[Sect. 23.2-23.4]{FG6} and \cite[Sect. 11]{Ra5} in the language of topological associative algebras. 

\medskip

The methods of {\it loc. cit.} could be adapted to the factorization setting. However, below we present a different
construction. Even though it looks more complicated (at least more abstract), its advantage is that it is compatible with the construction
of \secref{ss:z on KL}, where we do not know how to make other methods work.

\medskip

The construction of the action presented below has another advantage in that it is manifestly compatible
with the action of $\fL(G)$, see \secref{ss:center acts on Whit coinv} below. 

\end{rem}

\sssec{} \label{sss:action via coaction}

Let $\CZ$ be an ind-placid ind-scheme. The categories $\IndCoh^!(\CZ)$ and $\IndCoh^*(\CZ)$ are each mutually dual, with the pairing 
given by 
$$\IndCoh^!(\CZ)\otimes \IndCoh^*(\CZ)\overset{\sotimes}\to \IndCoh^*(\CZ) \overset{\Gamma^{\IndCoh}(\CZ,-)}\longrightarrow \Vect.$$

\medskip

Hence, we can view $\IndCoh^*(\CZ)$ as a comonoidal category. Moreover, the datum of an action of $\IndCoh^!(\CZ)$ on a
(factorization) category $\bC$ is is equivalent to the datum of a coaction of $\IndCoh^*(\CZ)$ on $\bC$. 

\sssec{}

As we shall see shortly (see \thmref{t:FF} and \secref{sss:Op placid}), the factorization ind-scheme $``\Spec"(\fZ_\fg)$ is ind-placid. 

\medskip

In particular,
the category $\IndCoh^*(``\Spec"(\fZ_\fg))$ is well-defined (see \secref{sss:IndCoh * Ran}). Moreover, it is compactly generated
and identifies with the dual of $\IndCoh^!(``\Spec"(\fZ_\fg))$. 

\medskip

The contents of \secref{sss:action via coaction} apply also in the factorization context. Hence, our task will be to define a coaction
of $\IndCoh^*(``\Spec"(\fZ_\fg))$, viewed as a comonoidal category, on $\hg\mod_\crit$. 

\sssec{} \label{sss:coaction functor first term}

We will first explain how to construct the coaction functor
\begin{equation} \label{e:coaction functor first term}
\on{coact}:\hg\mod_\crit\to \IndCoh^*(``\Spec"(\fZ_\fg))\otimes \hg\mod_\crit.
\end{equation} 

\medskip

Since $\on{Vac}(G)_\crit$ is the factorization unit in $\hg\mod_\crit$, we can write
$$\hg\mod_\crit\simeq \on{Vac}(G)_\crit\mod^{\on{fact}}(\hg\mod_\crit)$$
and
$$(\fz_\fg\otimes \on{Vac}(G)_\crit)\mod^{\on{fact}}(\hg\mod_\crit)\simeq 
\fz_\fg\mod^{\on{fact}}(\hg\mod_\crit).$$

\medskip

Now, restriction along the map \eqref{e:action of ff} gives rise to a t-exact functor
$$\on{Vac}(G)_\crit\mod^{\on{fact}}(\hg\mod_\crit)\to (\fz_\fg\otimes \on{Vac}(G)_\crit)\mod^{\on{fact}}(\hg\mod_\crit),$$
i.e., a functor
$$\hg\mod_\crit\to \fz_\fg\mod^{\on{fact}}(\hg\mod_\crit).$$

In particular, since the compact generators of $\hg\mod_\crit$ are eventually coconnective, we obtain a functor
\begin{equation} \label{e:coaction functor first term 0}
(\hg\mod_\crit)^c\to \left(\fz_\fg\mod^{\on{fact}}(\hg\mod_\crit)\right)^{>-\infty}.
\end{equation} 

\sssec{}

Note now that by combining \corref{c:R mod in A} and \corref{c:IndCoh* Op bdd below}(a), we obtain:

\begin{cor} \label{c:z-mod KM}
The functor
\begin{equation} \label{e:coaction functor first term comp}
\IndCoh^*(``\Spec"(\fZ_\fg))\otimes \hg\mod_\crit\to \fz_\fg\mod^{\on{fact}}(\hg\mod_\crit)
\end{equation}
is t-exact and induces an equivalence between the eventually coconnective subcategories of the two sides.
\end{cor}

\medskip

Hence, \eqref{e:coaction functor first term 0} can be thought of as a functor
\begin{equation} \label{e:coaction functor first term 1}
(\hg\mod_\crit)^c\to \left(\IndCoh^*(``\Spec"(\fZ_\fg))\otimes \hg\mod_\crit\right)^{>-\infty}\hookrightarrow 
\IndCoh^*(``\Spec"(\fZ_\fg))\otimes \hg\mod_\crit.
\end{equation} 

Ind-extending, from \eqref{e:coaction functor first term 1}, we obtain the desired functor \eqref{e:coaction functor first term}. 

\sssec{}

Our next task is to extend the functor \eqref{e:coaction functor first term} to a coaction of 
$\IndCoh^*(``\Spec"(\fZ_\fg))$ on $\hg\mod_\crit$. In doing so we will have to overcome two
hurdles:

\smallskip

\noindent{(i)} Homological-algebraic, which has to do with inverting the functor \eqref{e:coaction functor first term}
on the eventually coconnective subcategories.

\smallskip

\noindent{(ii)} Homotopic-algebraic, which has to do with equipping the functor  \eqref{e:coaction functor first term}
with a homotopy-coherent associativity datum.

\medskip

We will deal with (i) in the rest of this subsection, and with (ii) in \secref{s:device}. 

\sssec{} \label{sss:define coaction 1}

First, proceeding as in \secref{sss:coaction functor first term}, for an integer $n$, we define an $n$-ry operation
\begin{equation} \label{e:coaction functor term n}
\on{coact}_n:\hg\mod_\crit\to \IndCoh^*(``\Spec"(\fZ_\fg))^{\otimes n}\otimes \hg\mod_\crit,
\end{equation} 
so that the composition with 
\begin{equation} \label{e:coaction functor term n oblv}
\IndCoh^*(``\Spec"(\fZ_\fg))^{\otimes n}\otimes \hg\mod_\crit\to \fz^{\otimes n}_\fg\mod^{\on{fact}}(\hg\mod_\crit)
\end{equation}
is the restriction functor along the action map
\begin{equation} \label{e:action of ff n}
\fz_\fg^{\otimes n}\otimes \on{Vac}(G)_\crit\to \on{Vac}(G)_\crit.
\end{equation} 

\sssec{}

We claim:

\begin{lem} \label{l:coaction functor term n}
The functor \eqref{e:coaction functor term n} is t-exact.
\end{lem}

\begin{proof}

By construction, the composition of \eqref{e:coaction functor term n} with \eqref{e:coaction functor term n oblv}
is t-exact. Since the functor \eqref{e:coaction functor term n oblv} induces an equivalence on eventually eventually coconnective 
subcategories (see \corref{c:z-mod KM}),
it suffices to show that \eqref{e:coaction functor term n} has a bounded cohomological amplitude (over each $X^I$). By factorization,
this reduces to the case when $I$ is a singleton, and by evaluating at field-valued points of $X$, we reduce to the pointwise
case. The latter was established in \cite[Sect. 11.13]{Ra5}. 

\end{proof}

\sssec{} \label{sss:define coaction 2} 

From \lemref{l:coaction functor term n} we obtain:

\begin{cor} \label{c:coaction homotopy}
The functor \eqref{e:coaction functor first term} satisfies associativity \emph{at the homotopy level}.
\end{cor} 

\begin{proof}

We need to show:

\begin{itemize}

\item For every $n=n_1+n_2$, the diagram 
{\small$$
\CD
\hg\mod_\crit @>{\on{coact}_{n_2}}>> \IndCoh^*(``\Spec"(\fZ_\fg))^{\otimes n_2}\otimes \hg\mod_\crit \\
@V{\on{coact}_n}VV @VV{\on{Id}\otimes \on{coact}_{n_1}}V \\ 
\IndCoh^*(``\Spec"(\fZ_\fg))^{\otimes n}\otimes \hg\mod_\crit @>{\sim}>> 
\IndCoh^*(``\Spec"(\fZ_\fg))^{\otimes n_2}\otimes \IndCoh^*(``\Spec"(\fZ_\fg))^{\otimes n_1}\otimes \hg\mod_\crit
\endCD
$$}
commutes;

\medskip

\item 
For any $n$, the diagram
$$
\CD 
\hg\mod_\crit @>{\on{coact}}>> \IndCoh^*(``\Spec"(\fZ_\fg))\otimes \hg\mod_\crit  \\
@V{\on{coact}_n}VV @VV{\otimes \on{comult}_n\otimes \on{Id}}V  \\
\IndCoh^*(``\Spec"(\fZ_\fg))^{\otimes n}\otimes \hg\mod_\crit  @>{\sim}>> 
\IndCoh^*(``\Spec"(\fZ_\fg))^{\otimes n}\otimes \hg\mod_\crit  
\endCD
$$
commutes.

\end{itemize} 

In both cases, it suffices to show that the natural transformation in question is an isomorphism
when evaluated on compact objects. In particular, it suffices to show that it is an isomorphism
when evaluated on eventually coconnected subcategories. 

\medskip

We know that the natural transformation becomes an isomorphism after composing with the functor
\eqref{e:coaction functor term n oblv}. Since the functor \eqref{e:coaction functor term n oblv} is 
an equivalence on eventually coconnected subcategories (see \corref{c:z-mod KM}), it suffices to show that all the functors
involved have cohomological amplitude bounded on the left. However, this follows from 
\lemref{l:coaction functor term n} (for $\on{coact}_n$), while the functor 
$$\on{comult}_n:\IndCoh^*(``\Spec"(\fZ_\fg))\to
\IndCoh^*(``\Spec"(\fZ_\fg))^{\otimes n}\simeq \IndCoh^*(``\Spec"(\fZ_\fg)^n)$$
is t-exact (see \corref{c:dir im IndCoh* t exact}).

\end{proof}

\ssec{Action of the center and the loop group action} \label{ss:action of center and loop} \label{ss:center acts on Whit coinv}

\sssec{}

Our current goal is to show that  the $\IndCoh^!(``\Spec"(\fZ_\fg))$-module structure on $\hg\mod_\crit$
is compatible with the action of $\fL(G)$ on $\hg\mod_\crit$. 

\medskip

By the construction of the module structure, we need to show that the each of the categories
$$(\fz_\fg^{\otimes n}\otimes \on{Vac}(G)_\crit)\mod^{\on{fact}}(\hg\mod_\crit)$$
carries an action of $\fL(G)$, such that:

\begin{itemize}

\item It is compatible with the functor
\begin{multline*}
\IndCoh^!(``\Spec"(\fZ_\fg))^{\otimes n}\otimes \hg\mod_\crit \overset{\text{\eqref{e:coaction functor term n oblv}}}\longrightarrow \\
\to \fz_\fg^{\otimes n}\mod^{\on{fact}}(\hg\mod_\crit)\simeq (\fz_\fg^{\otimes n}\otimes \on{Vac}(G)_\crit)\mod^{\on{fact}}(\hg\mod_\crit);
\end{multline*}
where the $\fL(G)$ on the left-hand side is via the $\hg\mod_\crit$-factor;

\smallskip

\item The restriction functors 
$$(\fz_\fg^{\otimes n_1}\otimes \on{Vac}(G)_\crit)\mod^{\on{fact}}(\hg\mod_\crit)\to 
(\fz_\fg^{\otimes n_2}\otimes \on{Vac}(G)_\crit)\mod^{\on{fact}}(\hg\mod_\crit)$$
along the maps
$$\fz_\fg^{\otimes n_2}\otimes \on{Vac}(G)_\crit\to \fz_\fg^{\otimes n_1}\otimes \on{Vac}(G)_\crit$$
that encode the $\fz_\fg$-action on $\on{Vac}(G)_\crit$ carry a natural $\fL(G)$-equivariant structure.

\end{itemize}

\sssec{}

In order to do so, it suffices to show that for \emph{any} factorization algebra $\CA\in \KL(G)_\crit$, the lax factorization category
$$\CA\mod^{\on{fact}}(\hg\mod_\crit)$$
carries an action of $\fL(G)$, compatible with the forgetful functor
$$\oblv_\CA:\CA\mod^{\on{fact}}(\hg\mod_\crit)\to \hg\mod_\crit,$$
and this construction is functorial with respect to the functors
$$\Res_\phi:\CA_2\mod^{\on{fact}}(\hg\mod_\crit)\to \CA_1\mod^{\on{fact}}(\hg\mod_\crit),$$
corresponding to homomorphisms
$$\phi:\CA_1\to \CA_2$$
of factorization algebras in $\KL(G)_\crit$. 

\sssec{}

Let $\oblv_{\fL^+(G)}$ denote the forgetful functor. For $\CZ\to \Ran$, consider
\begin{equation} \label{e:restr along Sph}
\Res_{\oblv_{\fL^+(G)}}(\hg\mod_\crit^{\on{fact}_\CZ})\in \KL(G)_\crit\mmod^{\on{fact}}_\CZ.
\end{equation} 

By \lemref{l:modules for fact alg restr untl}, for $\CA\in \on{FactAlg}^{\on{untl}}(X,\KL(G)_\crit)$, we have
$$\oblv_{\fL^+(G)}(\CA)\mmod^{\on{fact}}(\hg\mod_\crit)_\CZ\simeq 
\CA\mmod^{\on{fact}}\left(\Res_{\oblv_{\fL^+(G)}}(\hg\mod_\crit^{\on{fact}_\CZ})\right)_\CZ.$$

Hence, it suffices to show that \eqref{e:restr along Sph} carries an action of $\fL(G)_\CZ$. 
We will show this in the following general framework. 

\sssec{} \label{sss:sph vs LG act}

Let $\bA$ be a factorization category, and let $\bC$ be a factorization module category 
over $\bA$ at $\CZ\to \Ran$. We assume that $\bA$ and $\bC$ carry compatible actions 
of $\fL(G)$ at some level $\kappa$. 

\medskip

Set $\bA_0:=\bA^{\fL^+(G)}$; denote by $\oblv_{\fL^+(G)}$ the forgetful functor
$$\bA_0\to \bA.$$

\medskip

Consider
$$\Res_{\oblv_{\fL^+(G)}}(\bC)\in \bA_0\mmod^{\on{fact}}_\CZ.$$

We claim that $\Res_{\oblv_{\fL^+(G)}}(\bC)$, viewed as a factorization module category over $\bA_0$, 
carries an action of $\fL(G)_\CZ$, such that the induced action on
$$\Res_{\oblv_{\fL^+(G)}}(\bC)_\CZ\simeq \bC_\CZ$$
is the original $\fL(G)_\CZ$-action on $\bC_\CZ$.

\medskip

We sketch the construction of this action below; a more detailed exposition will be given in \cite{CFGY}.

\sssec{}

Let $\Gr^{\on{level}_\CZ}_G$ be the factorization $\Gr_G$-module space from \secref{sss:Gr G level}.
It is equipped with a compatible action of $\fL(G)$ in the left and a commuting $\fL(G)_\CZ$-action
on the right. 

\medskip

Let $\fL(G)^{\on{fact}_\CZ}$ be the vacuum factorization module space over $\fL(G)$ at $\CZ$ (see \secref{sss:vac fact space});
it carries a compatible action of $\fL(G)\times \fL(G)$.

\medskip

We have a naturally defined projection
$$\pi_\CZ:\fL(G)^{\on{fact}_\CZ}\to \Gr^{\on{level}_\CZ}_G,$$
(see \secref{sss:Gr G disc}) which gives rise to the pullback functor 
$$\pi_\CZ^*:\Dmod(\Gr^{\on{level}_\CZ}_G)\to \Dmod(\fL(G)^{\on{fact}_\CZ}),$$
compatible\footnote{See \secref{sss:factorization restriction} for what this means.} 
with the factorization functor
$$\pi^*:\Dmod(\Gr_G)\to \Dmod(\fL(G)),$$
given by pullback along the projection
$$\pi:\fL(G)\to \Gr_G.$$

\medskip

In particular, $\pi_\CZ^*$ gives rise to a functor
\begin{equation} \label{e:restr Gr}
\Dmod(\Gr^{\on{level}_\CZ}_G)\to \Res_{\pi^*}(\Dmod(\fL(G)^{\on{fact}_\CZ}),
\end{equation}
as factorization module categories over $\Dmod(\Gr_G)$,
see \secref{sss:univ property restr cat}. 

\sssec{}

We claim:

\begin{lem} \label{l:restr Gr}
The functor \eqref{e:restr Gr} is an equivalence.
\end{lem}

\begin{proof}

Follows from \lemref{l:fact res crit}.

\end{proof}

\sssec{} \label{sss:from LG to Gr_G}

We now return to the setting of \secref{sss:sph vs LG act}. Note that
$$\bC\simeq  \Dmod(\fL(G)^{\on{fact}_\CZ})\underset{\fL(G)}\otimes \bC$$
as factorization module categories over 
$$\bA\simeq \Dmod(\fL(G))\underset{\fL(G)}\otimes \bA$$ at $\CZ$, where we use the action of
$\fL(G)$ on itself on the \emph{left} to form the tensor product.

\medskip

Note that the functor $\pi_\CZ^*$ and hence the equivalence \eqref{e:restr Gr} are compatible with the
actions of $\fL(G)$ on the left. Hence, from \lemref{l:restr Gr} we obtain:

\begin{cor} \label{c:restr Gr}
There is a canonical equivalence
$$\Res_{\oblv_{\fL^+(G)}}(\bC) \simeq \Dmod(\Gr^{\on{level}_\CZ}_G)\underset{\fL(G)}\otimes \bC,$$
as factorization categories over
$$\bA_0\simeq \Dmod(\Gr_G)\underset{\fL(G)}\otimes \bA.$$
\end{cor}

\sssec{}

Now, the action of $\fL(G)_\CZ$-action on $\Gr^{\on{level}_\CZ}_G$ on the right gives rise to an 
action of $\fL(G)_\CZ$ on $\Dmod(\Gr^{\on{level}_\CZ}_G)\underset{\fL(G)}\otimes \bC$, 
commuting with the factorization module structure over $\Dmod(\Gr_G)\underset{\fL(G)}\otimes \bA$.

\medskip

Applying \corref{c:restr Gr}, we produce the sought-for $\fL(G)_\CZ$-action on $\Res_{\oblv_{\fL^+(G)}}(\bC)$.

\sssec{}

As a consequence of the compatibility of the $\IndCoh^!(``\Spec"(\fZ_\fg))$ and $\fL(G)$-actions, we obtain:

\begin{cor} \label{c:z action on KL} \hfill

\smallskip

\noindent{\em(a)}
The category $\KL(G)_\crit$ carries an action of $\IndCoh^!(``\Spec"(\fZ_\fg))$ compatible with the action of $\IndCoh^!(``\Spec"(\fZ_\fg))$ on $\hg\mod_\crit$ 
and the forgetful functor
$$\KL(G)_\crit\to \hg\mod_\crit.$$

\smallskip

\noindent{\em(b)} The action of $\IndCoh^!(``\Spec"(\fZ_\fg))$ on $\KL(G)_\crit$ is compatible with the action of $\Sph_G$. 

\end{cor}

\begin{rem}
Note that we could have equivalently defined the action of $\IndCoh^!(``\Spec"(\fZ_\fg))$ on $\KL(G)_\crit$ directly, 
by repeating the procedure in \secref{ss:action of center}, replacing $\hg\mod_\crit$ by $\KL(G)_\crit$.  An analog of 
\lemref{l:coaction functor term n}
follows from the original variant of this lemma, since the corresponding functors
$$\IndCoh^*(``\Spec"(\fZ_\fg))^{\otimes n}\otimes \KL(G)_\crit\to \IndCoh^*(``\Spec"(\fZ_\fg))^{\otimes n}\otimes \hg\mod_\crit$$
are conservative.
\end{rem}

\sssec{}

Note that since $\fz_\fg\subset \BV_{\fg,\crit}$ is invariant under the adjoint action, we can view it also as the (classical)
center of the twisted version $\BV_{\fg,\crit,\rho(\omega_X)}$. 

\medskip

In particular, we can regard $\fz_\fg$ as acting on $\on{Vac}(G)_{\crit,\rho(\omega_X)}$ as an object of $\hg\mod_{\crit,\rho(\omega_X)}$
(or $\KL(G)_{\crit,\rho(\omega_X)}$). 

\medskip

In particular, we obtain an action of $\IndCoh^!(``\Spec"(\fZ_\fg))$ on $\hg\mod_{\crit,\rho(\omega)}$, compatible with the action of
$\fL(G)_{\rho(\omega_X)}$. The conclusion of \lemref{c:z action on KL} renders automatically to the present twisted context.

\medskip

In addition, we have: 

\begin{cor} \label{c:z action on Wht(KM)} \hfill

\smallskip

\noindent{\em(a)}
The category $\Whit_*(\hg\mod_{\crit,\rho(\omega)})$ carries a unique action of $\IndCoh^!(``\Spec"(\fZ_\fg))$, compatible with the action of 
$\IndCoh^!(``\Spec"(\fZ_\fg))$ on $\hg\mod_{\crit,\rho(\omega)}$ and the projection
$$\hg\mod_{\crit,\rho(\omega)}\to \Whit_*(\hg\mod_{\crit,\rho(\omega)}).$$ 

\smallskip

\noindent{\em(b)} The category $\Whit^!(\hg\mod_{\crit,\rho(\omega)})$ carries a unique action of $\IndCoh^!(``\Spec"(\fZ_\fg))$, compatible with the action of 
$\IndCoh^!(``\Spec"(\fZ_\fg))$ on $\hg\mod_{\crit,\rho(\omega)}$ and the embedding
$$\Whit^!(\hg\mod_{\crit,\rho(\omega)})\hookrightarrow \hg\mod_{\crit,\rho(\omega)}.$$ 

\smallskip

\noindent{\em(c)} The functor
$$\Theta_{\Whit(\hg\mod_{\crit,\rho(\omega)})}:\Whit_*(\hg\mod_{\crit,\rho(\omega)})\to \Whit^!(\hg\mod_{\crit,\rho(\omega)})$$
carries a natural $\IndCoh^!(``\Spec"(\fZ_\fg))$-linear structure.
\end{cor}

\ssec{The enhanced functor of Drinfeld-Sokolov reduction at the critical level} \label{ss:enh DS}

In this subsection we will study the functor
\begin{equation} \label{e:DS crit}
\DS:\hg\mod_{\crit,\rho(\omega_X)}\to \Vect
\end{equation}
of \eqref{e:initial DS}. 

\sssec{}

Consider the factorization unit 
$$\on{Vac}(G)_{\crit,\rho(\omega_X)}\in \hg\mod_{\crit,\rho(\omega_X)}.$$

\medskip

The functor \eqref{e:DS crit} has a natural lax unital factorization structure (see \secref{sss:DS}). In particular, the object
$$\DS(\on{Vac}(G)_{\crit,\rho(\omega_X)})$$
is naturally a factorization algebra (in $\Vect$). 

\sssec{}

The action map
$$\fz_\fg\otimes \on{Vac}(G)_{\crit,\rho(\omega_X)}\to \on{Vac}(G)_{\crit,\rho(\omega_X)}$$
gives rise to a map
$$\fz_\fg\otimes \DS(\on{Vac}(G)_{\crit,\rho(\omega_X)})\to \DS(\on{Vac}(G)_{\crit,\rho(\omega_X)})$$
as factorization algebras.

\medskip

Pre-composing with the unit for $\DS(\on{Vac}(G)_{\crit,\rho(\omega_X)})$, we obtain a map of factorization algebras 
\begin{equation} \label{e:center to DS of Vac}
\fz_\fg\to \DS(\on{Vac}(G)_{\crit,\rho(\omega_X)}).
\end{equation}

We have the following fundamental result, see \cite{FF}:

\begin{thm} \label{t:center as DS}
The map \eqref{e:center to DS of Vac} is an isomorphism.
\end{thm} 

\sssec{}

By \secref{sss:fact alg from fact cat}, the functor $\DS$ of \eqref{e:DS crit} naturally lifts to a functor
$$\hg\mod_{\crit,\rho(\omega_X)}\to \DS(\on{Vac}(G)_{\crit,\rho(\omega_X)})\mod^{\on{fact}}.$$

Restricting along \eqref{e:center to DS of Vac}, we can view it as a functor, to be denoted
$$\DS^{\on{enh}}:\hg\mod_{\crit,\rho(\omega_X)}\to \fz_\fg\mod^{\on{fact}}.$$

\sssec{}

Consider the functor
$$\Gamma^{\IndCoh}(``\Spec"(\fZ_\fg),-):\IndCoh^*(``\Spec"(\fZ_\fg))\to \fz_\fg\mod^{\on{fact}}.$$

We claim:

\begin{prop} \label{p:DS enh coarse}
There exists a uniquely defined (continuous) functor $$\DS^{\on{enh,rfnd}}:\hg\mod_{\crit,\rho(\omega_X)}\to \IndCoh^*(``\Spec"(\fZ_\fg)),$$
satisfying
\begin{itemize}

\item There exists an isomorphism
\begin{equation} \label{e:DS enh}
\DS^{\on{enh}}\simeq \Gamma^\IndCoh(``\Spec"(\fZ_\fg),-)^{\on{enh}}\circ \DS^{\on{enh,rfnd}};
\end{equation}

\item $\DS^{\on{enh,rfnd}}$ sends compact objects in $\hg\mod_{\crit,\rho(\omega_X)}$ to eventually
coconnective (i.e., bounded below) objects in $\IndCoh^*(``\Spec"(\fZ_\fg))$. 

\end{itemize}

Furthermore, $\DS^{\on{enh,rfnd}}$ carries a uniquely defined factorization structure, so that \eqref{e:DS enh}
is an isomorphism of factorization functors.

\end{prop}

\begin{proof}

It is enough to show that the restriction of $\DS^{\on{enh}}$ to the subcategory
$$(\hg\mod_{\crit,\rho(\omega_X)})^c\subset \hg\mod_{\crit,\rho(\omega_X)}$$
can be uniquely lifted to a functor 
$$(\hg\mod_{\crit,\rho(\omega_X)})^c\to \IndCoh^*(``\Spec"(\fZ_\fg))^{>-\infty}.$$

However, this follows from \corref{c:from QCoh * to fact placid}(a), using the fact that the initial functor
$\DS$ sends 
$$(\hg\mod_{\crit,\rho(\omega_X)})^c\to \Vect^{>-\infty}.$$

\end{proof} 

\sssec{}

Recall now that the functor $\DS$ factors via a functor
$$\ol{\DS}:\Whit_*(\hg\mod_{\crit,\rho(\omega_X)})\to \Vect.$$

It follows formally that the functor $\DS^{\on{enh}}$ also factors via a functor, denoted
$$\ol\DS^{\on{enh}}:\Whit_*(\hg\mod_{\crit,\rho(\omega_X)})\to \fz_\fg\mod^{\on{fact}}.$$

We now quote the following fundamental result of \cite{Ra2}:

\begin{thm} \label{t:Ra-W}
The functor $\DS^{\on{enh,rfnd}}$ factors via a functor 
\begin{equation} \label{e:DS enh bar}
\ol\DS^{\on{enh,rfnd}}:\Whit_*(\hg\mod_{\crit,\rho(\omega_X)})\to 
\IndCoh^*(``\Spec"(\fZ_\fg)),
\end{equation}
and the resulting functor $\ol\DS^{\on{enh,rfnd}}$ is an equivalence of factorization categories.
\end{thm} 

Note that by construction
$$\Gamma^\IndCoh(``\Spec"(\fZ_\fg),-)\circ \ol\DS^{\on{enh,rfnd}}\simeq \ol\DS^{\on{enh}}.$$

\sssec{} \label{sss:DS and action of center}

According to \secref{ss:action of center} (applied to the twist $\hg\mod_{\crit,\rho(\omega_X)}$ instead of the original $\hg\mod_\crit$), 
the category $\hg\mod_{\crit,\rho(\omega_X)}$ carries an action of the monoidal category $\IndCoh^!(``\Spec"(\fZ_\fg))$. 

\medskip

By construction, the functor $\DS^{\on{enh,rfnd}}$ intertwines the above $\IndCoh^!(``\Spec"(\fZ_\fg))$-action
on $\hg\mod_{\crit,\rho(\omega_X)}$ and the natural $\IndCoh^!(``\Spec"(\fZ_\fg))$-action on $\IndCoh^*(``\Spec"(\fZ_\fg))$.

\medskip

According to \corref{c:z action on Wht(KM)}(a), the $\IndCoh^!(``\Spec"(\fZ_\fg))$-action
on $\hg\mod_{\crit,\rho(\omega_X)}$ descends to an (a priori, uniquely defined) action on 
$\Whit_*(\hg\mod_{\crit,\rho(\omega_X)})$.

\medskip

It follows formally that the functor $\ol\DS^{\on{enh,rfnd}}$ intertwines this action and the
$\IndCoh^!(``\Spec"(\fZ_\fg))$-action on $\IndCoh^*(``\Spec"(\fZ_\fg))$.

\section{The Feigin-Frenkel isomorphism and its applications} \label{s:FF}

In this section we review the Feigin-Frenkel isomorphism, which provides a bridge between
Kac-Moody representations and opers. 

\medskip

Using the Feigin-Frenkel isomorphism, we construct an action of $\IndCoh^!(\Op^\mf_\cG)$
on $\KL(G)_\crit$, which is a key ingredient of the critical FLE functor, studied in the next
section.

\ssec{The Feigin-Frenkel isomorphism}

\sssec{}

We quote the following fundamental result of Feigin and Frenkel (\cite{FF}):

\begin{thm} \label{t:FF}
There exists a canonically defined isomorphism of factorization algebras 
$$\fz_\fg\overset{\on{FF}_G}\simeq \CO_{\Op^\reg_\cG}.$$
\end{thm}

Below we will complement \thmref{t:FF} by an assertion that describes how it interacts 
with geometric Satake, see \thmref{t:birth}.

\sssec{}

As in Remark \ref{r:choice CS}, one has a choice in the normalization of the functor $\on{FF}_G$;
the two different choices differ by the Chevalley involution.

\medskip

However, the choices for the normalizations of $\on{CS}_G$ and $\on{FF}_G$ force one another
via the compatibility in \thmref{t:FLE and Sat}. So, once $\on{CS}_G$ is fixed, $\on{FF}_G$ is 
unambiguous.

\sssec{Example}

Let $G=T$ be a torus. Then $\fz_\fg$ is the commutative factorization algebra associated with the
commutative algebra object
$$\Sym^!(\ft\otimes \on{D}_X[1])\in \on{ComAlg}(\Dmod(X)),$$
i.e.,
$$\fz_\fg=\on{Fact}(\Sym^!(\ft\otimes \on{D}_X[1])),$$
see \secref{sss:com fact vs Dmod com} for the notation. 

\medskip

When we think of $\Dmod(X)$ as ``left D-modules", the above object is 
$$\Sym_{\CO_X}(\ft\otimes \on{D}_X\otimes \omega^{\otimes -1}_X)\in \on{ComAlg}(\Dmod^l(X)).$$

The affine D-scheme $\Op_\cT$ identifies with the scheme of jets $\on{Jets}(\ct\otimes \omega_X)$ (see \secref{sss:jet construction}). 
Under our normalization, the isomorphism $\on{FF}_T$ is the tautological identification
\begin{equation} \label{e:naive FF for torus}
\Spec_X(\Sym_{\CO_X}(\ft\otimes \on{D}_X\otimes \omega^{\otimes -1}_X))\simeq J(\ct\otimes \omega_X).
\end{equation}

At the level of fibers at a given point $x\in X$, the identification \eqref{e:naive FF for torus} is the isomorphism
$$\Spec(\Sym(\ft\otimes \CK_x/\CO_x))\simeq \ct\otimes \omega_{\cD_x},$$
corresponding to the canonical identification of pro-finite dimensional vector spaces
$$(\ft\otimes \CK_x/\CO_x)^\vee \simeq \ct\otimes \omega_{\cD_x}.$$

\sssec{}

For a general $G$, the isomorphism $\on{FF}_G$ is normalized so that the resulting isomorphism 
at the associated graded level
$$\on{Jets}(\fa(\cg)_{\omega_X})\simeq \Spec_X(\on{gr}(\fz_{\fg,X})) \overset{\on{FF}_G}\simeq \Op^{\on{cl}}_\cG\simeq \on{Jets}(\fa(\cg)_{\omega_X}),$$
is the identity map, where $\Op^{\on{cl}}_\cG$ denotes the D-scheme of \emph{classical} $\cG$-opers.

\ssec{The ``birth" of opers}

In this subsection we will formulate \thmref{t:birth}, which in \cite[Sect. 5.3]{BD1} was called ``the birth of opers",
that explains how the isomorphism $\on{FF}_G$ interacts with geometric Satake. 

\sssec{} 

Consider the (symmetric monoidal) functor
\begin{equation} \label{e:Rep G to z}
\Rep(\cG)\overset{\fr^*}\to \QCoh(\Op_\cG^\reg)\simeq \CO_{\Op^\reg_\cG}\mod^{\on{com}}\overset{\on{FF}_G}\simeq \fz_\fg\mod^{\on{com}}.
\end{equation}

In particular, the functor \eqref{e:Rep G to z} allows us to view $\fz_\fg\mod^{\on{com}}$ as a $\Rep(\cG)$-module category,
in a way compatible with factorization.

\sssec{} \label{sss:Rep acts on KM}

Let us view $\KL(G)_\crit$ as a module category over $\Rep(\cG)$ via
$$\Rep(\cG)\overset{\Sat_{G}^{\on{nv}}}\longrightarrow \Sph_G$$
and the $\Sph_G$-action on $\KL(G)_\crit$. This structure is also compatible with factorization.

\sssec{}

Finally, note that the action of $\fz_\fg$ on $\on{Vac}(G)_\crit$ gives rise to a factorization functor 
\begin{equation} \label{e:z to KL}
\fz_\fg\mod^{\on{com}}\to \KL(G)_\crit, \quad -\underset{\fz_\fg}\otimes \on{Vac}(G)_\crit.
\end{equation} 

\sssec{}

We claim (see \cite[Theorem 5.5.3]{BD1}):

\begin{thm} \label{t:birth}
The functor \eqref{e:z to KL} admits a lift to a functor between $\Rep(\cG)$-module categories.
This structure is compatible with factorization.
\end{thm} 

\begin{rem}
Concretely, \thmref{t:birth} says that the object $\on{Vac}(G)_\crit$ satisfies the \emph{Hecke property}
with respect to the action of $\Rep(\cG)$ on $\KL(G)_\crit$: i.e., we have
$$\Sat_{G}^{\on{nv}}(V)\star \on{Vac}(G)_\crit \simeq \on{Vac}(G)_\crit \underset{\fz_\fg}\otimes (\on{FF}_G\circ (\fr^\reg)^*(V)),\quad V\in \Rep(\cG),$$
where in the right-hand side we denoted by $\on{FF}_G$ the equivalence
$$\fz_\fg\mod^{\on{com}}\simeq \CO_{\Op^\reg_\cG}\mod^{\on{com}}=\QCoh(\Op_\cG^\reg),$$
induced by the isomorphism of algebras $\on{FF}_G$.

\medskip

The above isomorphisms are compatible with tensor products of the $V$'s in the natural sense. 
\end{rem} 

\sssec{}

From now on we will identify
$$\Op_\cG^\reg\simeq \Spec(\fz_\fg) \text{ and } \Op_\cG^{\on{mer}}\simeq ``\Spec"(\fZ_\fg)$$
using $\on{FF}_G$.

\sssec{}

In particular, we will view the category $\hg\mod_\crit$ as acted on by $\IndCoh^!(\Op_\cG^{\on{mer}})$.
Similarly, we will view the functors $\ol\DS^{\on{enh,rfnd}}$ (resp., $\DS^{\on{enh,rfnd}}$) as $\IndCoh^!(\Op_\cG^{\on{mer}})$-linear functors from 
$\Whit_*(\hg\mod_{\crit,\rho(\omega_X)})$ (resp., $\hg\mod_{\crit,\rho(\omega_X)}$) to
$\IndCoh^*(\Op_\cG^{\on{mer}})$. 

\sssec{} \label{sss:R Op acts on Vac}

Here is one particular application of \thmref{t:birth} that will be used in the sequel. Recall the commutative algebra
(factorization) object
$$R_{\cG,\Op}\in \Rep(\cG),$$
\secref{e:R G Op}.

\medskip

We claim that it acts on $\on{Vac}(G)_\crit\in \KL(G)_\crit$, when we consider $\KL(G)_\crit$ as a 
$\Rep(\cG)$-module category as in \secref{sss:Rep acts on KM}.

\medskip

Indeed, by \thmref{t:birth}, in order to construct this structure, it suffices to construct an action of $R_{\cG,\Op}$
on 
$$\CO_{\Op^\reg_\cG}\in \QCoh(\Op^\reg_\cG),$$
when we consider $\QCoh(\Op^\reg_\cG)$ as a $\Rep(\cG)$-module category via $(\fr^\reg)^*$. However, the latter structure 
comes from the map of commutative algebras in $\QCoh(\Op^\reg_\cG)$
$$(\fr^\reg)^*(R_{\cG,\Op})\to \CO_{\Op^\reg_\cG},$$
given by the counit of the $((\fr^\reg)^*,\fr^\reg_*)$-adjunction.

\ssec{The Kazhdan-Lusztig category at the critical level and monodromy-free opers} \label{ss:z on KL}

In this subsection we will show that the $\IndCoh^!(\Op_\cG^{\on{mer}})$-action on $\KL(G)_\crit$ factors through an action of
$\IndCoh^!(\Op_\cG^\mf)$.

\medskip

The construction will emulate the construction of the $\IndCoh^!(\Op_\cG^{\on{mer}})$-action on $\hg\mod_\crit$
in \secref{ss:action of center}, with a ``decoration" by $\Rep(\cG)$. 

\begin{rem}
The construction of such an action (along with its properties discussed in \secref{ss:properies act z on KL})
at fixed $\ul{x}\in \Ran$ was the subject of the paper \cite{FG5}. 
However, we do not know how to adapt the methods of \emph{loc. cit.} to the factorization setting. 
\end{rem} 

\sssec{}

Using the duality between $\IndCoh^!(\Op_\cG^\mf)$ and $\IndCoh^*(\Op_\cG^\mf)$, it suffices to construct
a coaction of $\IndCoh^*(\Op_\cG^\mf)$, viewed as a comonoidal factorization category, on $\KL(G)_\crit$. 

\medskip

As in \secref{sss:coaction functor first term}, we first explain how to construct the coaction functor
\begin{equation} \label{e:co-action on KL}
\KL(G)_\crit\to \IndCoh^*(\Op_\cG^\mf)\otimes \KL(G)_\crit.
\end{equation}

We will then upgrade this to the datum of coaction. 

\sssec{} \label{sss:mf on KL 1}

We start with the monoidal action of $\Rep(\cG)$ on $\KL(G)_\crit$ as in \secref{sss:Rep acts on KM}. Since $\Rep(\cG)$ is rigid,
the right adjoint to the action is a (lax unital) factorization functor 
\begin{equation} \label{e:co-action of Rep}
\on{coact}_{\Rep(\cG),\KL(G)_\crit}:\KL(G)_\crit\to \Rep(\cG)\otimes \KL(G)_\crit.
\end{equation}

This functor upgrades to a factorization functor
\begin{multline} \label{e:co-action of Rep enh}
\on{coact}_{\Rep(\cG),\KL(G)_\crit}^{\on{enh}}:\KL(G)_\crit\to \\
\to \on{coact}_{\Rep(\cG),\KL(G)_\crit}(\on{Vac}(G)_\crit)\mod^{\on{fact}}(\Rep(\cG)\otimes \KL(G)_\crit).
\end{multline}

\sssec{}

Recall now (see \secref{sss:R Op acts on Vac}) that $R_{\cG,\Op}$, viewed as an associative algebra object in $\Rep(\cG)$ acts
on $\on{Vac}(G)_\crit$, compatibly with factorization. By adjunction, we obtain a map of factorization algebras 
$$R_{\cG,\Op}\otimes \on{Vac}(G)_\crit\to \on{coact}_{\Rep(\cG),\KL(G)_\crit}(\on{Vac}(G)_\crit).$$

Restriction along the above map defines a functor
\begin{multline} \label{e:co-action on KL 1}
\on{coact}_{\Rep(\cG),\KL(G)_\crit}(\on{Vac}(G)_\crit)\mod^{\on{fact}}(\Rep(\cG)\otimes \KL(G)_\crit)\to \\
\to (R_{\cG,\Op}\otimes \on{Vac}(G)_\crit)\mod^{\on{fact}}(\Rep(\cG)\otimes \KL(G)_\crit)=R_{\cG,\Op}\mod^{\on{fact}}(\Rep(\cG)\otimes \KL(G)_\crit)
\end{multline}

Composing \eqref{e:co-action of Rep enh} and \eqref{e:co-action on KL 1} we obtain a functor
\begin{equation} \label{e:co-action on KL 2}
\KL(G)_\crit   \to R_{\cG,\Op}\mod^{\on{fact}}(\Rep(\cG)\otimes \KL(G)_\crit).
\end{equation} 

\medskip

The functor \eqref{e:co-action of Rep} is left t-exact. Hence, since the compact generators of $\KL(G)_\crit$ are eventually coconnective,
the functor \eqref{e:co-action on KL 2} gives rise to a functor
\begin{equation} \label{e:co-action on KL 3}
(\KL(G)_\crit)^c   \to (R_{\cG,\Op}\mod^{\on{fact}}(\Rep(\cG)\otimes \KL(G)_\crit))^{>-\infty}.
\end{equation} 

\sssec{} \label{sss:mf on KL 2}

Consider the functor 
\begin{multline} \label{e:R-mod KL}
\IndCoh^*(\Op_\cG^\mf)\otimes \KL(G)_\crit\to 
R_{\cG,\Op}\mod^{\on{fact}}(\Rep(\cG))\otimes \KL(G)_\crit\to \\
\to R_{\cG,\Op}\mod^{\on{fact}}(\Rep(\cG)\otimes \KL(G)_\crit)
\end{multline}

As in \corref{c:z-mod KM}, by combining \corref{c:R mod in A} and \propref{p:IndCoh Op via fact almost}(a), we obtain: 

\begin{lem} \label{l:R-mod KL}
The functor \eqref{e:R-mod KL} induces an equivalence between the eventually coconnective subcategories of the two sides.
\end{lem}

\medskip

Hence, we can view \eqref{e:co-action on KL 3} as a functor
\begin{equation} \label{e:co-action on KL 4}
(\KL(G)_\crit)^c  \to \left(\IndCoh^*(\Op_\cG^\mf)\otimes \KL(G)_\crit\right)^{>-\infty}\hookrightarrow
\IndCoh^*(\Op_\cG^\mf)\otimes \KL(G)_\crit.
\end{equation} 

Ind-extending \eqref{e:co-action on KL 4} we obtain the sought-for functor \eqref{e:co-action on KL}. 

\sssec{} \label{sss:define coaction KM}

Our next goal is upgrade \eqref{e:co-action on KL} to a datum of coaction of 
$\IndCoh^*(\Op_\cG^\mf)$ on $\KL(G)_\crit$. We do so by mimicking the strategy 
in Sects. \ref{sss:define coaction 1}-\ref{sss:define coaction 2}.

\medskip

First, we generalize the construction in Sects. \ref{sss:mf on KL 1}-\ref{sss:mf on KL 2} above and define the $n$-ry operation
\begin{equation} \label{e:co-action on KL n}
\KL(G)_\crit\to \IndCoh^*(\Op_\cG^\mf)^{\otimes n}\otimes \KL(G)_\crit.
\end{equation}

\sssec{}

We have the following analog of \lemref{l:coaction functor term n}:

\begin{lem} \label{l:coaction functor term KL n}
The functor \eqref{e:co-action on KL n} is t-exact.
\end{lem}

\begin{proof}

Note that the functor
\begin{equation} \label{e:co-action on KL n bis}
\IndCoh^*(\Op_\cG^\mf)^{\otimes n}\otimes \KL(G)_\crit \overset{(\iota^\mf)^\IndCoh_*\otimes \oblv_{\fL^+(G)}}\longrightarrow
\IndCoh^*(\Op_\cG^\mer)^{\otimes n}\otimes \hg\mod_\crit
\end{equation}
is t-exact and conservative.

\medskip

Hence, it is enough to show that the composition of \eqref{e:co-action on KL n} with \eqref{e:co-action on KL n bis}
is t-exact. 

\medskip

Note, however, that by construction, we have a commutative diagram
$$
\CD
\KL(G)_\crit @>{\text{\eqref{e:co-action on KL n}}}>>   \IndCoh^*(\Op_\cG^\mf)^{\otimes n}\otimes \KL(G)_\crit \\
@V{\oblv_{\fL^+(G)}}VV @VV{(\iota^\mf)^\IndCoh_*\otimes \oblv_{\fL^+(G)}}V \\
\hg\mod_\crit @>{\text{\eqref{e:coaction functor term n oblv}}}>> \IndCoh^*(\Op_\cG^\mer)^{\otimes n}\otimes \hg\mod_\crit. 
\endCD
$$

Since $\oblv_{\fL^+(G)}$ is t-exact, the assertion follows from that of \lemref{l:coaction functor term n}.

\end{proof} 

\sssec{} 

As in \corref{c:coaction homotopy}, from \lemref{l:coaction functor term KL n}, we obtain that the coaction
functor \eqref{e:co-action on KL} is associative \emph{at the homotopy level}. We will equip it with a structure
of \emph{coherent homotopy} in \secref{ss:coaction functor term KL coh hom}.

\ssec{Properties of the \texorpdfstring{$\IndCoh^!(\Op_\cG^\mf)$}{IndCoh!}-action on \texorpdfstring{$\KL(G)_\crit$}{KL}} \label{ss:properies act z on KL}

In this subsection we will discuss those properties of the $\IndCoh^!(\Op_\cG^\mf)$-action on $\KL(G)_\crit$
that will be used in the sequel. 

\sssec{}

First, unwinding the construction, we obtain:

\begin{lem} \label{l:Rep G act on KM}
The action of $\Rep(\cG)$ on $\KL(G)_\crit$ given by
\begin{equation} \label{e:r^*}
\Rep(\cG)\simeq \QCoh(\LS_\cG^\reg)\overset{\fr^*}\to \QCoh(\Op_\cG^\mf) 
\overset{-\otimes \omega_{\Op_\cG^\mf}}\longrightarrow  \IndCoh^!(\Op_\cG^\mf)
\end{equation} 
and the above action of $\IndCoh^!(\Op_\cG^\mf)$ on $\KL(G)_\crit$ identifies canonically with the
action given by
$$\Rep(\cG)\overset{\Sat_{G}^{\on{nv}}}\to \Sph_G$$
and the $\Sph_G$-action on $\KL(G)_\crit$.
\end{lem} 

\sssec{} \label{sss:z act on KL and KM}

Further, comparing with the construction of the action of $\IndCoh^!(``\Spec"(\fZ_\fg))$ on $\KL(G)_\crit$ given by \corref{c:z action on KL},
we obtain that this action coincides with the precomposition of the above action of $\IndCoh^!(\Op_\cG^\mf)$
on $\KL(G)_\crit$ with
$$\iota^!:\IndCoh^!(\Op_\cG^{\on{mer}})\to \IndCoh^!(\Op_\cG^\mf).$$

\sssec{}

We now claim:

\begin{cor} \label{c:Hecke action t-exact}
For a fixed $\ul{x}\in \Ran$, the action functor of $\Rep(\cG)_{\ul{x}}$ on $\KL(G)_{\crit,\ul{x}}$ 
via $\Sat_{G}^{\on{nv}}$ and the $\Sph_{G,\ul{x}}$-action on $\KL(G)_{\crit,\ul{x}}$ is t-exact.
\end{cor}

\begin{proof}

By \lemref{l:Rep G act on KM}, we need to show that for $V\in \Rep(\cG)^\heartsuit_{\ul{x}}$, the functor
$$\KL(G)_{\crit,\ul{x}} \overset{V\otimes \on{Id}}\to 
\Rep(\cG)_{\ul{x}}\otimes \KL(G)_{\crit,\ul{x}} \to \IndCoh^!(\Op^\mf_{\cG,\ul{x}})) \otimes \KL(G)_{\crit,\ul{x}} \to \KL(G)_{\crit,\ul{x}}$$
is t-exact, where $\Rep(\cG)\to  \IndCoh^!(\Op^\mf_{\cG,\ul{x}})$ is the functor \eqref{e:r^*}.

\medskip

Hence, it suffices to show that if $\CE^\mf \in \QCoh(\Op^\mf_{\cG,\ul{x}})$ is a \emph{vector bundle}, then its action on 
$\KL(G)_{\crit,\ul{x}}$ via 
$$\QCoh(\Op^\mf_{\cG,\ul{x}}) \overset{-\otimes \omega_{\Op^\mf_{\cG,\ul{x}}}}\longrightarrow  \IndCoh^!(\Op^\mf_{\cG,\ul{x}})$$
and the $\IndCoh^!(\Op^\mf_{\cG,\ul{x}})$-action on $\KL(G)_{\crit,\ul{x}}$ is t-exact.

\medskip

With no restriction of generality, we can assume that $\CE^\mf$ is the restriction of a vector bundle $\CE^\mer$ over $\Op^\mer_{\cG,\ul{x}}$. Hence,
by \secref{sss:z act on KL and KM}, it suffices to show that for a vector bundle $\CE^\mer$ on  $\Op^\mer_{\cG,\ul{x}}$, its action on $\KL(G)_{\crit,\ul{x}}$ via
$$\QCoh(\Op^\mer_{\cG,\ul{x}}) \overset{-\otimes \omega_{\Op^\mer_{\cG,\ul{x}}}}\longrightarrow  \IndCoh^!(\Op^\mer_{\cG,\ul{x}})$$
and the action of $\IndCoh^!(\Op^\mer_{\cG,\ul{x}})$ given by \corref{c:z action on KL} is t-exact. 

\medskip

Since the forgetful functor
$$\KL(G)_{\crit,\ul{x}}\to \hg\mod_{\crit,\ul{x}}$$
is t-exact and conservative, it suffices to show that the action of $\CE^\mer$ on $\hg\mod_{\crit,\ul{x}}$ is t-exact.

\medskip

Unwinding the construction, it suffices to show that the composition
\begin{multline*}
\hg\mod_{\crit,\ul{x}} \to \hg\mod_{\crit,\ul{x}}\otimes \IndCoh^*(\Op^\mer_{\cG,\ul{x}})
\overset{\on{Id}\otimes (\CE^\mer\otimes -)}\longrightarrow \\
\to \hg\mod_{\crit,\ul{x}}\otimes \IndCoh^*(\Op^\mer_{\cG,\ul{x}})
\overset{\on{Id}\otimes \Gamma^\IndCoh(\Op^\mer_{\cG,\ul{x}},-)}\longrightarrow \hg\mod_{\crit,\ul{x}}
\end{multline*}
is t-exact.

\medskip

In the above composition, the first and the third arrows are t-exact. Hence, it suffices to show that the functor
$$\hg\mod_{\crit,\ul{x}}\otimes \IndCoh^*(\Op^\mer_{\cG,\ul{x}})
\overset{\on{Id}\otimes (\CE^\mer\otimes -)}\longrightarrow \hg\mod_{\crit,\ul{x}}\otimes \IndCoh^*(\Op^\mer_{\cG,\ul{x}})$$
is t-exact.

\medskip

However, this easily follows from the fact that the functor
$$\IndCoh^*(\Op^\mer_{\cG,\ul{x}})\overset{(\CE^\mer\otimes -)}\longrightarrow  \IndCoh^*(\Op^\mer_{\cG,\ul{x}})$$
is t-exact. 

\end{proof}

\section{The critical FLE} \label{s:FLE}

In this section we prove the main result of Part I, namely, the \emph{critical FLE}, \thmref{t:critical FLE}, which says that there
exists a canonical equivalence of factorization categories 
\begin{equation} \label{e:critical FLE functor}
\FLE_{G,\crit}:\KL(G)_\crit\overset{\sim}\to \IndCoh^*(\Op^{\on{mon-free}}_\cG),
\end{equation}

The functor in one direction in \eqref{e:critical FLE functor} is a variation on the theme of the functor $\ol\DS^{\on{enh,rfnd}}$ from
\secref{ss:enh DS}. Essentially $\FLE_{G,\crit}$ is obtained by base changing $\ol\DS^{\on{enh,rfnd}}$ along the map from $\Op^\mf_\cG$
to $\Op^\mer_\cG$. 

\ssec{Construction of the critical FLE functor}

\sssec{}

Let $\bC$ be a category equipped with a $\fL(G)_{\rho(\omega_X)}$-action at the critical level, in a way compatible with factorization. 
Consider the functor
\begin{equation} \label{e:pre-FLE -1}
\Sph(\bC):=\bC^{\fL(G)^+_{\rho(\omega_X)}}\to \bC\to \bC_{\fL(N)_{\rho(\omega_X)},\chi}=:\Whit_*(\bC).
\end{equation}

We apply this to $\bC=\hg\mod_{\crit,\rho(\omega_X)}$. Consider the resulting (factorization) functor
\begin{equation} \label{e:pre-FLE 0}
\KL(G)_{\crit,\rho(\omega_X)} \to \Whit_*(\hg\mod_{\crit,\rho(\omega_X)})
\end{equation}

Composing, we obtain a functor
\begin{equation} \label{e:FLE 1}
\KL(G)_{\crit,\rho(\omega_X)} \to \Whit_*(\hg\mod_{\crit,\rho(\omega_X)})
\overset{\ol\DS^{\on{enh,rfnd}}}\longrightarrow \IndCoh^*(\Op_\cG^{\on{mer}}).
\end{equation}

\sssec{} \label{sss:construct FLE}

We regard $\IndCoh^*(\Op_\cG^{\on{mer}})$ as equipped with a natural action of $\IndCoh^!(\Op_\cG^{\on{mer}})$. We regard 
$\KL(G)_{\crit,\rho(\omega_X)}$ as acted on by $\IndCoh^!(\Op_\cG^\mf)$. 

\medskip

By Sects. \ref{sss:z act on KL and KM} and 
\ref{sss:DS and action of center} and \corref{c:z action on KL}(a), 
the functor \eqref{e:FLE 1} is compatible with the $\IndCoh^!(\Op_\cG^{\on{mer}})$-actions on the two sides.
Furthermore, by \secref{sss:z act on KL and KM}, the $\IndCoh^!(\Op_\cG^{\on{mer}})$-action on $\KL(G)_{\crit,\rho(\omega_X)}$ factors through an action of 
$\IndCoh^!(\Op_\cG^\mf)$.

\medskip

Hence, the functor \eqref{e:FLE 1} gives rise to a (factorization) functor
\begin{equation} \label{e:FLE 2}
\KL(G)_{\crit,\rho(\omega_X)} \to \on{Funct}_{\IndCoh^!(\Op_\cG^{\on{mer}})}(\IndCoh^!(\Op_\cG^\mf),\IndCoh^*(\Op_\cG^{\on{mer}})).
\end{equation}

Finally, recall that by \lemref{l:hom category}, we have a canonical identification 
$$\IndCoh^*(\Op_\cG^\mf)\simeq  \on{Funct}_{\IndCoh^!(\Op_\cG^{\on{mer}})}(\IndCoh^!(\Op_\cG^\mf),\IndCoh^*(\Op_\cG^{\on{mer}})).$$

\medskip

Thus, we can interpret \eqref{e:FLE 2} as a $\IndCoh^!(\Op_\cG^\mf)$-linear functor
\begin{equation}  \label{e:critical FLE omega}
\KL(G)_{\crit,\rho(\omega_X)} \to \IndCoh^*(\Op_\cG^\mf)
\end{equation}

\sssec{}

Precomposing \eqref{e:critical FLE omega}
with $\KL(G)_{\crit} \overset{\alpha_{\rho(\omega_X),\on{taut}}} \simeq \KL(G)_{\crit,\rho(\omega_X)}$, we obtain a functor 

\begin{equation} \label{e:critical FLE}
\FLE_{G,\crit}:\KL(G)_{\crit} \to \IndCoh^*(\Op_\cG^\mf)
\end{equation}

The functor \eqref{e:critical FLE} is the critical FLE functor. The main result of Part I of this paper reads:

\begin{thm} \label{t:critical FLE}
The functor $\FLE_{G,\crit}$ is an equivalence of factorization categories.
\end{thm}

This theorem will be proved in the course of this and the next two sections.

\begin{rem}
A pointwise version of \thmref{t:critical FLE}, formulated below as \thmref{t:critical FLE pointwise}, was originally
proved in \cite{FG2}. 

\medskip

The proof of the factorization version requires substantially new ideas. 

\medskip

Progress towards  \thmref{t:critical FLE} has been made earlier in the papers \cite{FLMM1,FLMM2}.

\end{rem}

\sssec{}

Unwinding the definitions, we observe that the functor $\FLE_{G,\crit}$ carries a natural lax unital
structure (as a functor between unital factorization categories). In particular, we obtain a canonical homomorphism
\begin{equation} \label{e:FLE is unital}
\CO_{\Op^\reg_\cG}=\one_{\IndCoh^*(\Op_\cG^\mf)}\to \FLE_{G,\crit}(\one_{\KL(G)_{\crit}})=\on{Vac}(G)_\crit
\end{equation} 
as factorization algebras in $\IndCoh^*(\Op_\cG^\mf)$.

\medskip

However, we claim:

\begin{lem} \label{l:FLE is unital}
The map \eqref{e:FLE is unital} is an isomorphism.
\end{lem}

\begin{proof}

By \propref{p:mon-free to all}(a), it suffices to show that the map \eqref{e:FLE is unital} becomes an isomorphism
after applying the functor $(\iota^\mf)^\IndCoh_*$. Since the latter is also a strict unital 
factorization functor, the resulting homomorphism identifies with a homomorphism of factorization algebras 
$$\CO_{\Op^\reg_\cG}=\one_{\IndCoh^*(\Op_\cG^\reg)}\to \DS^{\on{enh,rfnd}}(\on{Vac}(G)_{\crit,\rho(\omega_X)}),$$
corresponding to the lax unital functor
$$\KL(G)_{\crit,\rho(\omega_X)} \to \hg\mod_{\crit,\rho(\omega_X)} \overset{\DS^{\on{enh,rfnd}}}\to  \IndCoh^*(\Op_\cG^\mer).$$

However, the latter isomorphism is the content of \thmref{t:center as DS}.

\end{proof} 

\sssec{}

Combining Lemmas \ref{l:FLE is unital} and \ref{l:unit determines strict}, we obtain:

\begin{cor} \label{c:FLE is unital}
The functor $\FLE_{G,\crit}$ is strictly unital.
\end{cor} 

\ssec{Reduction to the pointwise version}

\sssec{}

Fix a point $x\in X$. The pointwise version of \thmref{t:critical FLE} reads:

\begin{thm} \label{t:critical FLE pointwise}
The functor $\FLE_{G,\crit}$ induces an equivalence
$$\KL(G)_{\crit,x}\overset{\sim}\to  \IndCoh^*(\Op^\mf_{\cG,x}).$$
\end{thm}

Obviously, \thmref{t:critical FLE} implies \thmref{t:critical FLE pointwise}. However, in this subsection
we will show that the coverse implication also takes place. 

\medskip

In its turn, \thmref{t:critical FLE pointwise} is known: it is the main result of the paper \cite{FG2}. We will, however,
supply a different proof, in which we deduce it from \thmref{t:Ra-W}, see \secref{ss:proof of local FLE}.

\sssec{}

A key step in proving the implication  \thmref{t:critical FLE pointwise} $\Rightarrow$ \thmref{t:critical FLE} 
is the following:

\begin{prop} \label{p:FLE pres comp}
The functor $\FLE_{G,\crit}$ preserves compactness.\footnote{See \secref{sss:pres compac fact} for what it means 
for a factorization functor to preserve compactness.}
\end{prop}

\begin{proof}

By \propref{p:mon-free to all}(b), it suffices to show that the composite functor
$$\KL(G)_{\crit} \overset{\FLE_{G,\crit}}\longrightarrow \IndCoh^*(\Op_\cG^\mf) \overset{(\iota^\mf)^{\IndCoh}_*}\longrightarrow 
\IndCoh^*(\Op_\cG^{\on{mer}})$$
preserves compactness.

\medskip

I.e., it suffices to show that \eqref{e:FLE 1} preserves compactness. Since $\ol\DS^{\on{enh,rfnd}}$ is an equivalence, it suffices to show
that the functor \eqref{e:pre-FLE 0} preserves compactness. However, we claim that this is true more generally. 

\medskip

Namely, we claim that the functor \eqref{e:pre-FLE -1} admits a continuous right adjoint (and hence, preserves compactness). 
Indeed, the right adjoint in question is given 
by\footnote{In \lemref{l:Sph and Whit adj} we will give another description of this right adjoint.} 
convolution with the vacuum object (i.e., the factorization unit) 
$$\on{Vac}_{\Whit^!(G)}\in \Whit^!(G)\simeq  \Dmod_{\frac{1}{2}}(\fL(G)_{\rho(\omega_X)})^{(\fL(N)_{\rho(\omega_X)},\chi),\fL^+(G)_{\rho(\omega_X)}}.$$

\end{proof}

\sssec{}

Given \propref{p:FLE pres comp} and \thmref{t:critical FLE pointwise}, we will deduce \thmref{t:critical FLE} using the following principle:

\medskip

Let $F:\bC^1\to \bC^2$ be a factorization functor between factorization categories. 
Assume that $\bC^1$ is compactly generated, and assume that $F$ preserves compactness. 

\begin{prop} \label{p:ptw => fact}
If the induced functor $F_x:\bC^1_x\to \bC^2_x$ is an equivalence for any field-valued point $x$, then the original functor $F$ is
also an equivalence.
\end{prop}

\begin{proof}

The assumption that $F$ preserves compactness implies that its right adjoint $F^R$ is also equipped
with a factorization structure. We need to show that the unit and the counit of the $(F,F^R)$-adjunction
are isomorphisms. 

\medskip

The latter assertion can be checked strata-wise on $\Ran$. I.e., we have to show that for every $n$, the corresponding
functor
$$F_{\overset{\circ}X{}^{(n)}}:\bC^1_{\overset{\circ}X{}^{(n)}}\to \bC^2_{\overset{\circ}X{}^{(n)}}$$
is an equivalence. 

\medskip

By factorization, the latter statement reduces to the case $n=1$, i.e., we have to show that
$$F_X:\bC^1_X\to \bC^2_X$$
is an equivalence.

\medskip

The latter fact can be also checked after base-changing to field-valued points.  

\end{proof} 

\qed[\thmref{t:critical FLE}]

\ssec{The inverse of the critical FLE functor}

In this subsection we will assume the statement of \thmref{t:critical FLE}, which was proved modulo
\thmref{t:critical FLE pointwise}. 

\sssec{} \label{sss:inverse mech}

Let $\bC$ be a category, acted on by $\fL(G)_{\rho(\omega_X)}$. Note that in addition to the functor
\begin{equation} \label{e:Sph to Whit* again}
\Sph(\bC)\to \bC\to \Whit_*(\bC),
\end{equation}
one can consider the functor
\begin{equation} \label{e:Sph to Whit^! again}
\Whit^!(\bC)\to \bC\overset{\on{Av}^{\fL^+(G)_{\rho(\omega_X)}}_*}\to \Sph(\bC).
\end{equation} 

We have the following elementary assertion:

\begin{lem} \label{l:Sph and Whit adj}
The composite
\begin{equation} \label{e:Sph to Whit^! Theta}
\Whit_*(\bC) \overset{\Theta_{\Whit(\bC)}}\to \Whit^!(\bC) \overset{\text{\eqref{e:Sph to Whit^! again}}}\longrightarrow  \Sph(\bC)
\end{equation} 
identifies canonically with the right adjoint of \eqref{e:Sph to Whit* again}.
\end{lem}

\begin{rem}
This lemma is embedded into the machinery developed in \cite{Ra2}. We supply a proof
for completeness.
\end{rem} 

\begin{proof}

We need to check that for $\CF_\Sph\in \Sph(\bC)$ and $\CF\in \bC$, we have
$$\CHom_{\Whit_*(\bC)}(\ol\CF_\Sph,\ol\CF)\simeq \CHom_{\bC}(\CF_\Sph,\Theta_{\Whit(\bC)}(\ol\CF)),$$
where $\ol\CF_\Sph,\ol\CF$ denotes the image of $\CF_\Sph$ and $\CF$ along
$\bC\to  \Whit_*(\bC)$. 

\medskip

Unwinding the definitions, we reduce the assertion to the case when $\bC:=\Dmod_{\frac{1}{2}}(\Gr_{G,\rho(\omega_X)})$
and $\CF_\Sph=\delta_{1,\Gr_{G,\rho(\omega_X)}}$. 

\medskip

Applying the definition of $\Whit_*(G)$, we calculate
\begin{equation} \label{e:Hom from av delta}
\CHom_{\Whit_*(G)}(\ol\delta_{1,\Gr_{G,\rho(\omega_X)}},\ol\CF)=
\underset{\alpha}{\on{colim}}\, \CHom_{\Dmod_{\frac{1}{2}}(\Gr_{G,\rho(\omega_X)})_{N^\alpha,\chi}}
((\delta_{1,\Gr_{G,\rho(\omega_X)}})^\alpha,\CF^\alpha),
\end{equation} 
where:

\begin{itemize}

\item $N^\alpha$ is a filtered family of group subschemes that comprise $\fL(N)_{\rho(\omega_X)}$;

\item $(\delta_{1,\Gr_{G,\rho(\omega_X)}})^\alpha$ and $\CF^\alpha$ denote the projections of the corresponding
objects along $$\Dmod_{\frac{1}{2}}(\Gr_{G,\rho(\omega_X)})\to \Dmod_{\frac{1}{2}}(\Gr_{G,\rho(\omega_X)})_{N^\alpha,\chi}.$$

\end{itemize}

We have
$$
\CHom_{\Dmod_{\frac{1}{2}}(\Gr_{G,\rho(\omega_X)})_{N^\alpha}}
((\delta_{1,\Gr_{G,\rho(\omega_X)}})^\alpha,\CF^\alpha)\simeq 
\CHom_{\Dmod_{\frac{1}{2}}(\Gr_{G,\rho(\omega_X)})}(\on{Av}_*^{N^\alpha,\chi}(\delta_{1,\Gr_{G,\rho(\omega_X)}}),\CF),$$
which we further rewrite
as
$$\on{C}^\cdot\left(\Gr_{G,\rho(\omega_X)}, \BD\left(\on{Av}_*^{N^\alpha,\chi}(\delta_{1,\Gr_{G,\rho(\omega_X)}})\right)\sotimes \CF\right) 
\simeq
\on{C}^\cdot(\Gr_{G,\rho(\omega_X)}, \on{Av}_!^{N^\alpha,\chi}(\delta_{1,\Gr_{G,\rho(\omega_X)}})\sotimes \CF).$$

Hence, we can rewrite \eqref{e:Hom from av delta} as 
\begin{multline} \label{e:Hom from av delta 1}
\on{C}^\cdot\left(\Gr_{G,\rho(\omega_X)}, \underset{\alpha}{\on{colim}}\, (\on{Av}_!^{N^\alpha,\chi}(\delta_{1,\Gr_{G,\rho(\omega_X)}}))
\sotimes \CF\right)\simeq \\
\simeq \on{C}^\cdot\left(\Gr_{G,\rho(\omega_X)}, 
\on{Av}_!^{\fL(N)_{\rho(\omega_X)},\chi}(\delta_{1,\Gr_{G,\rho(\omega_X)}})\sotimes \CF\right).
\end{multline}

Now, the cleanness property from \secref{sss:Whit clean} implies that the natural map
$$\on{Av}_!^{\fL(N)_{\rho(\omega_X)},\chi}(\delta_{1,\Gr_{G,\rho(\omega_X)}})\to
\on{Av}_{*,\on{ren}}^{\fL(N)_{\rho(\omega_X)},\chi}(\delta_{1,\Gr_{G,\rho(\omega_X)}})$$
is an isomorphism, where $\on{Av}_{*,\on{ren}}^{\fL(N)_{\rho(\omega_X)},\chi}$ is the functor of $*$-convolution
with $\omega^{\on{ren}}_{\fL(N)_{\rho(\omega_X)},\chi}$. 

\medskip 

Hence, we further rewrite \eqref{e:Hom from av delta 1} as 
\begin{multline*} 
\on{C}^\cdot\left(\Gr_{G,\rho(\omega_X)}, \on{Av}_{*,\on{ren}}^{\fL(N)_{\rho(\omega_X)},\chi}(\delta_{1,\Gr_{G,\rho(\omega_X)}})
\sotimes \CF\right) \simeq \\
\simeq \on{C}^\cdot\left(\Gr_{G,\rho(\omega_X)}, \delta_{1,\Gr_{G,\rho(\omega_X)}}\sotimes 
\on{Av}_{*,\on{ren}}^{\fL(N)_{\rho(\omega_X)},\chi}(\CF)\right) \simeq \\
\simeq \CHom_{\Dmod_{\frac{1}{2}}(\Gr_{G,\rho(\omega_X)})}( \delta_{1,\Gr_{G,\rho(\omega_X)}},
\Theta_{\Whit(G)}(\ol\CF)),
\end{multline*}
as desired.

\end{proof}

\sssec{}

We will now use \lemref{l:Sph and Whit adj} to give an explicit description of the inverse of the functor $\FLE_{G,\crit}$. 

\medskip

Consider the functor \eqref{e:Sph to Whit^! Theta} for $\hg\mod_{\crit,\rho(\omega_X)}$
\begin{equation} \label{e:Whit KM to KL}
\Whit_*(\hg\mod_{\crit,\rho(\omega_X)})\to \KL(G)_{\crit,\rho(\omega_X)}.
\end{equation}

By the same logic as in \secref{sss:construct FLE}, the functor \eqref{e:Whit KM to KL} gives rise to a functor
\begin{equation} \label{e:pre-inverse FLE 0}
\IndCoh^!(\Op_\cG^\mf)\underset{\IndCoh^!(\Op_\cG^\mer)}\otimes \Whit_*(\hg\mod_{\crit,\rho(\omega_X)})\to \KL(G)_{\crit,\rho(\omega_X)}.
\end{equation}

Combining with the equivalence $\ol\DS^{\on{enh,rfnd}}$, we obtain a functor
\begin{equation} \label{e:pre-inverse FLE 1}
\IndCoh^!(\Op_\cG^\mf)\underset{\IndCoh^!(\Op_\cG^\mer)}\otimes \IndCoh^*(\Op_\cG^\mer) \to \KL(G)_{\crit,\rho(\omega_X)}.
\end{equation}

Combining with the equivalence \eqref{e:ten prod category}, from \eqref{e:pre-inverse FLE 1} we obtain a functor
\begin{equation} \label{e:pre-inverse FLE}
\IndCoh^*(\Op_\cG^\mf)\to \KL(G)_{\crit,\rho(\omega_X)}.
\end{equation} 

We will prove:

\begin{prop} \label{p:inverse FLE}
The functor \eqref{e:pre-inverse FLE} is the inverse of \eqref{e:critical FLE omega}.
\end{prop}

The rest of this subsection is devoted to the proof of \propref{p:inverse FLE}. 

\sssec{}

We need to show that the composition
\begin{multline} \label{e:pre-inverse FLE 2}
\IndCoh^!(\Op_\cG^\mf)\underset{\IndCoh^!(\Op_\cG^\mer)}\otimes \Whit_*(\hg\mod_{\crit,\rho(\omega_X)})
\overset{\on{Id}\otimes \ol\DS^{\on{enh,rfnd}}}\longrightarrow \\
\to \IndCoh^!(\Op_\cG^\mf)\underset{\IndCoh^!(\Op_\cG^\mer)}\otimes \IndCoh^*(\Op_\cG^\mer) \to 
\IndCoh^*(\Op_\cG^\mf)
\end{multline} 
is isomorphic to 
\begin{multline} \label{e:pre-inverse FLE 3}
\IndCoh^!(\Op_\cG^\mf)\underset{\IndCoh^!(\Op_\cG^\mer)}\otimes \Whit_*(\hg\mod_{\crit,\rho(\omega_X)})
\overset{\text{\eqref{e:pre-inverse FLE 1}}}\longrightarrow \\
\to \KL(G)_{\crit,\rho(\omega_X)}\overset{\text{\eqref{e:critical FLE omega}}}\longrightarrow
\IndCoh^*(\Op_\cG^\mf).
\end{multline} 

Both functors are $\IndCoh^!(\Op_\cG^\mf)$-linear. Hence, by the
$$\IndCoh^!(\Op_\cG^\mer)\mmod \rightleftarrows \IndCoh^!(\Op_\cG^\mf)\mmod$$
adjunction, it suffices to show that the functors
\begin{equation} \label{e:pre-inverse FLE 2'}
\Whit_*(\hg\mod_{\crit,\rho(\omega_X)}) \overset{\ol\DS^{\on{enh,rfnd}}}\longrightarrow \IndCoh^*(\Op_\cG^\mer) \overset{(\iota^\mf)^!}\to \IndCoh^*(\Op_\cG^\mf)
\end{equation} 
and 
$$\Whit_*(\hg\mod_{\crit,\rho(\omega_X)}) \overset{\text{\eqref{e:pre-inverse FLE 0}}}\longrightarrow 
\KL(G)_{\crit,\rho(\omega_X)}\overset{\text{\eqref{e:critical FLE omega}}}\longrightarrow \IndCoh^*(\Op_\cG^\mf)$$
are isomorphic as $\IndCoh^!(\Op_\cG^\mer)$-linear functors, i.e., that the diagram
$$
\CD
\Whit_*(\hg\mod_{\crit,\rho(\omega_X)})  @>{\text{\eqref{e:pre-inverse FLE 0}}}>>  \KL(G)_{\crit,\rho(\omega_X)} \\
@V{\ol\DS^{\on{enh,rfnd}}}VV @VV{\text{\eqref{e:critical FLE omega}}}V \\
\IndCoh^*(\Op_\cG^\mer)  @>{(\iota^\mf)^!}>> \IndCoh^*(\Op_\cG^\mf)
\endCD
$$
commutes, in a way compatible with the $\IndCoh^!(\Op_\cG^\mer)$-actions. 

\sssec{}

Since $\ol\DS^{\on{enh,rfnd}}$ and \eqref{e:critical FLE omega} are both equivalences, it suffices
to show that the diagram obtained by passing to left adjoints along the horizontal arrows, i.e.,
$$
\CD
\Whit_*(\hg\mod_{\crit,\rho(\omega_X)})  @<<<  \KL(G)_{\crit,\rho(\omega_X)} \\
@V{\ol\DS^{\on{enh,rfnd}}}VV @VV{\text{\eqref{e:critical FLE omega}}}V \\
\IndCoh^*(\Op_\cG^\mer)  @<{(\iota^\mf)^\IndCoh_*}<< \IndCoh^*(\Op_\cG^\mf),
\endCD
$$
commutes, in a way compatible with the $\IndCoh^!(\Op_\cG^\mer)$-actions. 
 
\sssec{}
 
However, according to \lemref{l:Sph and Whit adj}, the top vertical arrow in the latter diagram 
is the functor \eqref{e:pre-FLE 0}, and the corresponding diagram commutes by construction.

\qed[\propref{p:inverse FLE}]

\ssec{Compatibility of \texorpdfstring{$\FLE_{G,\crit}$}{FLEG} and \texorpdfstring{$\FLE_{\cG,\infty}$}{FLEGinfty}}

\sssec{}

Note that by construction, the functor $\FLE_{G,\crit}$ makes the following diagram commute
$$
\CD
\KL(G)_{\crit}  @>{\alpha_{\rho(\omega_X),\on{taut}}}>{\sim}> \KL(G)_{\crit,\rho(\omega_X)} @>>> \Whit_*(\hg\mod_{\crit,\rho(\omega_X)}) \\
@V{\FLE_{G,\crit}}VV & & @VV{\ol\DS^{\on{enh,rfnd}}}V \\
\IndCoh^*(\Op_\cG^\mf) & @>{(\iota^\mf)^{\IndCoh}_*}>> & \IndCoh^*(\Op_\cG^{\on{mer}}).
\endCD
$$

\sssec{}

Note that the functor \eqref{e:pre-FLE -1} can be expanded to a functor
\begin{equation} \label{e:pairing 1 prel abs}
\Whit_*(G)\underset{\Sph_G}\otimes \Sph(\bC)\to \Whit_*(\bC).
\end{equation} 

\medskip

Applying this to $\bC=\hg\mod_{\crit,\rho(\omega_X)}$, we obtain a functor
\begin{equation} \label{e:pairing 1 prel}
\Whit_*(G)\underset{\Sph_G}\otimes \KL(G)_{\crit,\rho(\omega_X)} \to \Whit_*(\hg\mod_{\crit,\rho(\omega_X)}).
\end{equation}

Composing with $\KL(G)_{\crit} \overset{\alpha_{\rho(\omega_X),\on{taut}}}\simeq \KL(G)_{\crit,\rho(\omega_X)}$ and
$\ol\DS^{\on{enh,rfnd}}$, we obtain a (factorization) functor
\begin{equation} \label{e:pairing 1}
\Whit_*(G)\underset{\Sph_G}\otimes \KL(G)_\crit\to \IndCoh^*(\Op_\cG^{\on{mer}}).
\end{equation}

\sssec{}

Similarly, the functor
$$(\iota^\mf)^{\IndCoh}_*:\IndCoh^*(\Op_\cG^\mf)\to \IndCoh^*(\Op_\cG^{\on{mer}})$$
can be expanded to a (factorization) functor
\begin{equation} \label{e:pairing 2}
\Rep(\cG)\underset{\Sph_\cG^{\on{spec}}}\otimes \IndCoh^*(\Op_\cG^\mf)\to \IndCoh^*(\Op_\cG^{\on{mer}}),
\end{equation}

\sssec{}

We claim:

\begin{thm} \label{t:FLE and Sat} \hfill

\smallskip

\noindent{\em(a)} 
The functor $\FLE_{G,\crit}$ can be canonically endowed with the datum of compatibility with the
$\Sph_G$-action on $\KL(G)_{\crit}$ and the $\Sph_\cG^{\on{spec}}$-action on $\IndCoh^*(\Op_\cG^\mf)$,
where we identify $$\Sph_G\simeq \Sph_\cG^{\on{spec}}$$ via $\Sat_G$. 

\smallskip

\noindent{\em(b)} Under the identification of point (a), the diagram
\begin{equation} \label{e:FLE and Sat diag}
\CD
\Whit_*(G)\underset{\Sph_G}\otimes \KL(G)_\crit @>{\text{\eqref{e:pairing 1}}}>> \IndCoh^*(\Op_\cG^{\on{mer}}) \\
@V{\FLE_{\cG,\infty}^{-1}\otimes \FLE_{G,\crit}}V{\sim}V @V{\sim}V{\on{Id}}V \\
\Rep(\cG)\underset{\Sph_\cG^{\on{spec}}}\otimes \IndCoh^*(\Op_\cG^\mf) @>{\text{\eqref{e:pairing 2}}}>> \IndCoh^*(\Op_\cG^{\on{mer}})
\endCD
\end{equation} 
canonically commutes.

\end{thm} 

This theorem will be proved in the factorization setting in \secref{ss:FLE and Sat}. 

\sssec{} \label{sss:local pairing G}

Let us denote by $\sP_G^{\on{loc,enh,rfnd}}$ the precomposition of \eqref{e:pairing 1} with the projection
$$\Whit_*(G)\otimes \KL(G)_\crit\to \Whit_*(G)\underset{\Sph_G}\otimes \KL(G)_\crit.$$

Explicitly, it is given by 
\begin{multline} \label{e:pairing G one}
\Whit_*(G) \otimes \KL(G)_\crit \overset{\on{Id}\otimes \alpha_{\rho(\omega_X),\on{taut}}}\longrightarrow 
\Whit_*(G)\otimes \KL(G)_{\crit,\rho(\omega_X)}\to \\
\to \Whit_*(G)\underset{\Sph_G}\otimes \KL(G)_{\crit,\rho(\omega_X)}
\overset{\text{\eqref{e:pairing 1 prel}}}\longrightarrow
\Whit_*(\hg\mod_{\crit,\rho(\omega_X)})\overset{\ol\DS^{\on{enh,rfnd}}}\to \IndCoh^*(\Op^{\on{mer}}_\cG).
\end{multline} 

Let $\sP^{\on{loc}}_G$ and $\sP^{\on{loc},\on{enh}}_G$
denote the compositions of $\sP^{\on{loc},\on{enh,rfnd}}_\cG$ with the forgetful functors
\begin{equation} \label{e:Gamma Op again}
\Gamma^{\IndCoh}(\Op^{\on{mer}}_\cG,-):\IndCoh^*(\Op^{\on{mer}}_\cG) \to \Vect
\end{equation}
and
\begin{equation} \label{e:Gamma Op enh again}
\Gamma^{\IndCoh}(\Op^{\on{mer}}_\cG,-)^{\on{enh}}:\IndCoh^*(\Op^{\on{mer}}_\cG) \to \CO_{\Op^\reg_\cG}\mod^{\on{fact}},
\end{equation}
respectively. (These two functors are obtained by  replacing the last arrow in \eqref{e:pairing G one} by $\ol\DS$
and $\ol\DS^{\on{enh}}$, respectively.)

\sssec{} \label{sss:PGc pairing}

Let us denote by $\sP_\cG^{\on{loc,enh,rfnd}}$ the precomposition of \eqref{e:pairing 2} with the projection
$$\Rep(\cG) \otimes \IndCoh^*(\Op_\cG^\mf)\to
\Rep(\cG)\underset{\Sph_\cG^{\on{spec}}}\otimes \IndCoh^*(\Op_\cG^\mf).$$

\medskip

Explicitly, it is given by 
\begin{multline} \label{e:pairing Gc one}
\Rep(\cG) \otimes \IndCoh^*(\Op^{\on{mon-free}}_\cG) \overset{\fr^*\otimes \on{Id}}\to \\
\to \QCoh(\Op^{\on{mon-free}}_\cG)\otimes \IndCoh^*(\Op^{\on{mon-free}}_\cG)\overset{\otimes}
\to \IndCoh^*(\Op^{\on{mon-free}}_\cG) \overset{(\iota^\mf)^\IndCoh_*}\to\\ \to  \IndCoh^*(\Op^{\on{mer}}_\cG).
\end{multline} 

Let $\sP^{\on{loc}}_\cG$ and $\sP^{\on{loc},\on{enh}}_\cG$
denote the compositions of $\sP^{\on{loc},\on{enh,rfnd}}_\cG$ with the forgetful functors \eqref{e:Gamma Op again}
and \eqref{e:Gamma Op enh again}, respectively.  

\sssec{}

From \thmref{t:FLE and Sat} we immediately obtain: 

\begin{cor} \label{c:two pairings fine}
The functors $\sP^{\on{loc},\on{enh,rfnd}}_G$ and $\sP^{\on{loc},\on{enh,rfnd}}_\cG$ match under the equivalences
$$\KL(G)_\crit \overset{\FLE_{G,\crit}}\simeq \IndCoh^*(\Op^{\on{mon-free}}_\cG) \text{ and }
\Rep(\cG) \overset{\FLE_{\cG,\infty}}\simeq \Whit_*(G).$$
\end{cor} 

And hence: 

\begin{cor} \label{c:two pairings coarse}
The functors $\sP^{\on{loc}}_G$ and $\sP^{\on{loc}}_\cG$ (resp., $\sP^{\on{loc},\on{enh}}_G$ and $\sP^{\on{loc},\on{enh}}_\cG$)
match under the equivalences
$$\KL(G)_\crit \overset{\FLE_{G,\crit}}\simeq \IndCoh^*(\Op^{\on{mon-free}}_\cG) \text{ and }
\Rep(\cG) \overset{\FLE_{\cG,\infty}}\simeq \Whit_*(G).$$
\end{cor} 

\ssec{The functor \texorpdfstring{$\on{pre-FLE}_{G,\crit}$}{preFLE}}

\sssec{}

By the construction of the functor $\FLE_{G,\crit}$ we have the following explicit descriptions 
of its compositions with various forgetful functors out of $\IndCoh^*(\Op^{\on{mon-free}}_\cG)$:

\begin{itemize}

\item The composition with the functor 
$$\IndCoh^*(\Op^{\on{mon-free}}_\cG)\overset{(\iota^\mf)^\IndCoh_*}\to \IndCoh^*(\Op^{\on{mer}}_\cG)$$
is the functor 
$$\KL(G)_\crit \overset{\alpha_{\rho(\omega_X),\on{taut}}}\longrightarrow \KL(G)_{\crit,\rho(\omega_X)}\to
\hg\mod_{\crit,\rho(\omega_X)} \overset{\DS^{\on{enh,rfnd}}}\to \IndCoh^*(\Op^{\on{mer}}_\cG);$$

\item The composition with the functor 
$$\IndCoh^*(\Op^{\on{mon-free}}_\cG)\overset{(\iota^\mf)^\IndCoh_*}\to \IndCoh^*(\Op^{\on{mer}}_\cG)
\overset{\Gamma^\IndCoh(\Op^{\on{mer}}_\cG,-)^{\on{enh}}}\longrightarrow \CO_{\Op^\reg_\cG}\mod^{\on{fact}}$$
is the functor 
$$\KL(G)_\crit \overset{\alpha_{\rho(\omega_X),\on{taut}}}\longrightarrow \KL(G)_{\crit,\rho(\omega_X)}\to
\hg\mod_{\crit,\rho(\omega_X)} \overset{\DS^{\on{enh}}}\to \CO_{\Op^\reg_\cG}\mod^{\on{fact}};$$

\item The composition with the functor
$$\Gamma^\IndCoh(\Op^{\on{mon-free}}_\cG,-): \IndCoh^*(\Op^{\on{mon-free}}_\cG)\to \Vect$$
is the functor
\begin{equation} \label{e:DS on KL again}
\KL(G)_\crit \overset{\alpha_{\rho(\omega_X),\on{taut}}}\longrightarrow \KL(G)_{\crit,\rho(\omega_X)}\to
\hg\mod_{\crit,\rho(\omega_X)} \overset{\DS}\to \Vect;
\end{equation} 

\end{itemize}

\medskip

In this subsection we will describe explicitly the composition of $\FLE_{G,\crit}$ with the functor
$$\fr^{\IndCoh}_*:\IndCoh^*(\Op^{\on{mon-free}}_\cG)\to \Rep(\cG).$$

\sssec{}

Define the (factorization) functor
\begin{equation} \label{e:pre-FLE}
\on{pre-FLE}_{G,\crit}:\KL(G)_\crit \to \Rep(\cG)
\end{equation}
as the composition
\begin{multline} \label{e:pre-FLE 1}
\KL(G)_\crit \overset{\alpha_{\rho(\omega_X),\on{taut}}}\longrightarrow \KL(G)_{\crit,\rho(\omega_X)}\to 
\Whit^!(G)\otimes \Whit_*(\hg\mod_{\crit,\rho(\omega_X)})\overset{\on{CS}_G\otimes \on{Id}}\longrightarrow \\
\to \Rep(\cG) \otimes \Whit_*(\hg\mod_{\crit,\rho(\omega_X)}) \overset{\on{Id}\otimes \ol\DS}\longrightarrow \Rep(\cG),
\end{multline}
where the second arrow is obtained by duality from the pairing 
$$\Whit_*(G)\otimes \KL(G)_{\crit,\rho(\omega_X)}\to \Whit_*(\hg\mod_{\crit,\rho(\omega_X)}).$$

\sssec{}

We claim:

\begin{prop} \label{p:pre-FLE}
The functor $\on{pre-FLE}_{G,\crit}$ identifies canonically with $\fr^\IndCoh_*\circ \FLE_{G,\crit}$. 
\end{prop}

The rest of this subsection is devoted to the proof of this proposition.

\sssec{}

The next assertion results from the construction of the functor $\FLE_{G,\crit}$ and \lemref{l:Rep G act on KM}:

\begin{lem} \label{l:pre-FLE}
The functor $\fr^\IndCoh_*\circ \FLE_{G,\crit}$ identifies with the composition
\begin{multline} \label{e:pre-FLE 2}
\KL(G)_\crit \to \Rep(\cG)\otimes \KL(G)_\crit  \overset{\on{Id}\otimes \alpha_{\rho(\omega_X),\on{taut}}}\longrightarrow \\
\to \Rep(\cG)\otimes \KL(G)_{\crit,\rho(\omega_X)}\to 
\Rep(\cG)\otimes \hg\mod_{\crit,\rho(\omega_X)}
\overset{\on{Id}\otimes \DS}\to \Rep(\cG),
\end{multline}
where the first arrow is the functor, right adjoint to the action of $\Rep(\cG)$ on $\KL(G)_\crit$,
given by 
$$\Rep(\cG)\overset{\Sat_{G}^{\on{nv}}}\longrightarrow \Sph_G$$
and the $\Sph_G$-action on $\KL(G)_\crit$. 
\end{lem}

Hence, in order to prove \propref{p:pre-FLE}, it suffices to establish an isomorphism between
\eqref{e:pre-FLE 1} and \eqref{e:pre-FLE 2}. 

\sssec{}

We rewrite the functor in \eqref{e:pre-FLE 2} as
\begin{multline*}
\KL(G)_\crit \overset{\alpha_{\rho(\omega_X),\on{taut}}}\longrightarrow \KL(G)_{\crit,\rho(\omega_X)}\to \\
\to \Rep(\cG)\otimes \KL(G)_{\crit,\rho(\omega_X)} \to \Rep(\cG)\otimes \Whit_*(\hg\mod_{\crit,\rho(\omega_X)}) 
\overset{\on{Id}\otimes \ol\DS}\longrightarrow \Rep(\cG).
\end{multline*}

Hence, in order to prove an isomorphism between \eqref{e:pre-FLE 1} and \eqref{e:pre-FLE 2}, it suffices to show 
that the functors 
$$\KL(G)_{\crit,\rho(\omega_X)}
\to \Rep(\cG)\otimes \KL(G)_{\crit,\rho(\omega_X)} \to \Rep(\cG)\otimes \Whit_*(\hg\mod_{\crit,\rho(\omega_X)})$$ 
and
$$\KL(G)_{\crit,\rho(\omega_X)}\to 
\Whit^!(G)\otimes \Whit_*(\hg\mod_{\crit,\rho(\omega_X)})\overset{\on{CS}_G\otimes \on{Id}}\longrightarrow 
\Rep(\cG) \otimes \Whit_*(\hg\mod_{\crit,\rho(\omega_X)})$$
are canonically isomorphic. 

\sssec{}

By duality, this amounts to showing that the functors 
\begin{multline*}
\Rep(\cG)\otimes \KL(G)_{\crit,\rho(\omega_X)} \overset{\Sat_{G}^{\on{nv}}\otimes \on{Id}}\longrightarrow \\
\to \Sph_G\otimes \KL(G)_{\crit,\rho(\omega_X)} \to \KL(G)_{\crit,\rho(\omega_X)} \to \Whit_*(\hg\mod_{\crit,\rho(\omega_X)})
\end{multline*}
and
$$\Rep(\cG)\otimes \KL(G)_{\crit,\rho(\omega_X)} \overset{\FLE_{\cG,\infty}\otimes \on{Id}}\longrightarrow 
\Whit_*(G)\otimes \KL(G)_{\crit,\rho(\omega_X)}\to \Whit_*(\hg\mod_{\crit,\rho(\omega_X)})$$
are canonically identified. 

\sssec{}

This follows by combining the following observations:

\medskip

\begin{itemize}

\item The functor $\Whit_*(G)\otimes \KL(G)_{\crit,\rho(\omega_X)}\to \Whit_*(\hg\mod_{\crit,\rho(\omega_X)})$ factors as
$$\Whit_*(G)\otimes \KL(G)_{\crit,\rho(\omega_X)}\to \Whit_*(G)\underset{\Sph_G}\otimes \KL(G)_{\crit,\rho(\omega_X)}\to \Whit_*(\hg\mod_{\crit,\rho(\omega_X)});$$

\smallskip

\item The functor $\KL(G)_{\crit,\rho(\omega_X)} \to \Whit_*(\hg\mod_{\crit,\rho(\omega_X)})$ identifies with
$$\KL(G)_{\crit,\rho(\omega_X)} \overset{\on{Vac}_{\Whit_*(G)}\otimes \on{Id}}\longrightarrow 
\Whit_*(G)\otimes \KL(G)_{\crit,\rho(\omega_X)}\to \Whit_*(\hg\mod_{\crit,\rho(\omega_X)}),$$
where $\one_{\Whit_*(G)}\in \Whit_*(G)$ is the vacuum object;

\smallskip

\item The functor $\FLE_{\cG,\infty}$ identifies with
$$\Rep(\cG) \overset{\Sat_{G}^{\on{nv}}}\longrightarrow \Sph_G \overset{\on{Vac}_{\Whit_*(G)}\star -}\longrightarrow \Whit_*(G)$$
(see Remark \ref{r:CS again}).

\end{itemize}

\qed[\propref{p:pre-FLE}]

\ssec{An alternative construction of the critical FLE functor}

\sssec{}

We start by observing:

\begin{lem} \label{l:pre FLE unit}
There exists a canonical identification of factorization algebras in $\Rep(\cG)$ 
\begin{equation} \label{e:pre FLE unit}
\on{pre-FLE}_{G,\crit}(\Vac(G)_\crit) \simeq R_{\cG,\Op}.
\end{equation} 
\end{lem} 

\begin{proof}

Follows from \propref{p:pre-FLE} and the fact that the functor $\FLE_{G,\crit}$ is unital
(see \corref{c:FLE is unital}). 

\end{proof}

\begin{rem}
One can establish the isomorphism \eqref{e:pre FLE unit} directly, i.e., without
appealing to the functor $\FLE_{G,\crit}$. 
\end{rem} 

\sssec{}

By \secref{sss:fact alg from fact cat gen} and \lemref{l:pre FLE unit}, the functor $\on{pre-FLE}_{G,\crit}$ upgrades to a functor
$$\on{pre-FLE}_{G,\crit}^{\on{enh}}: \KL(G)_\crit\to R_{\cG,\Op}\mod^{\on{fact}}(\Rep(\cG)).$$

Note that by construction we have a canonical isomorphism 
$$\on{pre-FLE}_{G,\crit}^{\on{enh}} \simeq \Gamma^\IndCoh(\Op_\cG^\mf,-)^{\on{enh}}\circ \FLE_{G,\crit}.$$

Because of this isomorphism we will also use the notation
$$\FLE_{G,\crit}^{\on{coarse}}:=\on{pre-FLE}_{G,\crit}^{\on{enh}}.$$

\sssec{}

We claim now that the functor $\FLE_{G,\crit}$ can be uniquely recovered from $\FLE_{G,\crit}^{\on{coarse}}$. 

\medskip

Namely, 
by \propref{p:IndCoh Op via fact almost}(a), it suffices to show that the functor $\FLE_{G,\crit}^{\on{coarse}}$ sends compact
objects in $\KL(G)_\crit$ to eventually coconnective objects in $R_{\cG,\Op}\mod^{\on{fact}}(\Rep(\cG))$.

\sssec{}

Since the forgetful functor 
$$R_{\cG,\Op}\mod^{\on{fact}}(\Rep(\cG))\to \Rep(\cG)$$
is t-exact and conservative, it suffices to show that the functor $\on{pre-FLE}_{G,\crit}$ sends  
compact objects in $\KL(G)_\crit$ to eventually coconnective objects in $\Rep(\cG)$. 

\medskip

Since the compact generators of $\KL(G)_\crit$ are eventually coconnective, it suffices to prove the following:

\begin{lem} \label{l:pre FLE t exact}
For a fixed $\ul{x}\in \Ran$, the functor 
$$\on{pre-FLE}_{G,\crit}:\KL(G)_{\crit,\ul{x}}\to \Rep(\cG)_{\ul{x}}$$
is t-exact. 
\end{lem}

\begin{proof} 

We rewrite the functor $\on{pre-FLE}_{G,\crit}$ as \eqref{e:pre-FLE 2}, or equivalently
\begin{multline}
\KL(G)_\crit \overset{R_\cG\otimes \on{Id}}\to \Rep(\cG)\otimes \Rep(\cG) \otimes \KL(G)_\crit
\overset{\on{Id} \otimes \Sat_{G}^{\on{nv}}\otimes \on{Id}}\longrightarrow  
\Rep(\cG)\otimes \Sph_G \otimes \KL(G)_\crit \overset{\on{Id} \otimes (-\star-)}\longrightarrow \\
\to \Rep(\cG)\otimes  \KL(G)_\crit \overset{\on{Id}\otimes \alpha_{\rho(\omega_X),\on{taut}}}
\longrightarrow  \Rep(\cG)\otimes \KL(G)_{\crit,\rho(\omega_X)} \overset{\on{Id}\otimes \DS}\to  \Rep(\cG)
\end{multline}

In this composition, the first arrow is tautologically exact, and the second arrow is t-exact by 
\corref{c:Hecke action t-exact}. Hence, the assertion follows from \lemref{l:DS on KL exact}. 

\end{proof}

\sssec{}

Note that as a corollary of \lemref{l:pre FLE t exact} and \propref{p:pre-FLE}, we obtain:

\begin{cor} \label{c:FLE t exact}
For a fixed $\ul{x}\in \Ran$, the functor
$$\FLE_{G,\crit}:\KL(G)_{\crit,\ul{x}}\to \IndCoh^*(\Op^\mf_{\cG,\ul{x}})$$
is t-exact.
\end{cor}

\section{Proof of the pointwise version of the critical FLE} \label{s:FLE pointwise} 

In this section we will give a proof of the pointwise version of the critical FLE by deducing it
from \thmref{t:Ra-W}.

\medskip

The idea of the proof is the that the critical FLE is essentially the \emph{base change} of the equivalence
of \thmref{t:Ra-W} along $\LS^\reg_{\cG,x}\to \LS^\mer_{\cG,x}$. In fact, such an equivalence is a general
phenomenon for categories acted on by $\fL(G)_x$, given a \emph{temperedness} condition (see \propref{p:Sph to unr Whit}
for a precise statement). 

\medskip

The reason this proof only works for the pointwise version is that it is only in this case that we have a good grip
on the base change operation alluded to above. 

\medskip

We note that a completely different proof of the pointwise FLE was given in the paper \cite{FG2}. 

\ssec{Temperedness} \label{ss:temp}

\sssec{} \label{sss:temp spec}

Let 
\begin{equation} \label{e:embed temp}
\Sph_{\cG,\on{temp},x}^{\on{spec}}\hookrightarrow \Sph_{\cG,x}^{\on{spec}}
\end{equation}
be the \emph{tempered subcategory}. 

\medskip

By definition, this is the essential image of
$$\Xi_{\on{Hecke}_{\cG,x}^{\on{spec,loc}}}:\QCoh(\on{Hecke}_{\cG,x}^{\on{spec,loc}})\to \IndCoh(\on{Hecke}_{\cG,x}^{\on{spec,loc}}).$$

\medskip

The embedding \eqref{e:embed temp} admits a right adjoint, namely,
$$\Psi_{\on{Hecke}_{\cG,x}^{\on{spec,loc}}}:\IndCoh(\on{Hecke}_{\cG,x}^{\on{spec,loc}})\to \QCoh(\on{Hecke}_{\cG,x}^{\on{spec,loc}}),$$
whose kernel is a monoidal
ideal. 

\medskip

This allows us to view $\Sph_{\cG,\on{temp},x}^{\on{spec}}$ as a monoidal colocalization
of $\Sph_{\cG,x}^{\on{spec}}$.

\begin{rem} \label{r:no Ran temp}

The definition of $\Sph_{\cG,\on{temp},x}^{\on{spec}}$ is specific to the pointwise version. 
We do not know how to define it in the factorization setting. The reason for this is the following:

\medskip

Although we can define $\Sph_{\cG,\on{temp}}^{\on{spec}}:=\IndCoh^*(\on{Hecke}^{\on{spec,loc}}_\cG)$
in the factorization setting, we do not have the $(\Xi,\Psi)$-adjunction. The latter is a feature of a locally of
finite type situation, which we \emph{are} in at a fixed $\ul{x}\in \Ran$, but \emph{not} when we are allowed to vary over $\Ran$
in families.

\end{rem}

\sssec{}

Let us regard 
$$\QCoh(\LS^\reg_{\cG,x})\simeq \Rep(\cG)$$ as a bimodule with respect to $\QCoh((\LS^\mer_{\cG,x})^\wedge_\reg)$
and $\Sph_{\cG,x}^{\on{spec}}$. 

\medskip

Note, however, that the $\Sph_{\cG,x}^{\on{spec}}$-action on $\QCoh(\LS^\reg_{\cG,x})$
factors via $\Sph_{\cG,\on{temp},x}^{\on{spec}}$: indeed, the action is
given by t-exact functors and the t-structure on $\Rep(\cG)$ is separared.

\sssec{}

Consider the corresponding functor
\begin{equation} \label{e:temp quot via action}
\Sph_{\cG,\on{temp},x}^{\on{spec}}\to 
\on{Funct}_{\QCoh((\LS^\mer_{\cG,x})^\wedge_\reg)}(\QCoh(\LS^\reg_{\cG,x}), \QCoh(\LS^\reg_{\cG,x})).
\end{equation}

\medskip

The following results from the definitions:

\begin{lem} \label{l:temp quot via action}
The functor \eqref{e:temp quot via action} is an equivalence. 
\end{lem} 

\sssec{} \label{sss:temp}

Let $\bC$ be a module category over $\Sph_{\cG,x}^{\on{spec}}$. Denote 
$$\bC_{\on{temp}}:=\Sph_{\cG,\on{temp},x}^{\on{spec}}\underset{\Sph_{\cG,x}^{\on{spec}}}\otimes \bC.$$

The adjunction
$$\Sph_{\cG,\on{temp},x}^{\on{spec}}\rightleftarrows \Sph_{\cG,x}^{\on{spec}}$$
gives rise to an adjunction
\begin{equation} \label{e:temp adj}
\bC_{\on{temp}} \rightleftarrows  \bC,
\end{equation} 
making $\bC_{\on{temp}}$ into a colocalization of $\bC$. 

\medskip

We let $\on{temp}_\bC$ denote the comonad on $\bC$ corresponding to the adjunction
\eqref{e:temp adj}

\sssec{}

Let us regard the two sides of \eqref{e:temp quot via action} as right modules with respect to $\Sph_{\cG,x}^{\on{spec}}$,
where:

\begin{itemize}

\item $\Sph_{\cG,x}^{\on{spec}}$ acts on $\Sph_{\cG,\on{temp},x}^{\on{spec}}$ by right multiplication;

\item $\Sph_{\cG,x}^{\on{spec}}$ acts on $\on{Funct}_{\QCoh((\LS^\mer_{\cG,x})^\wedge_\reg)}(\QCoh(\LS^\reg_{\cG,x}), \QCoh(\LS^\reg_{\cG,x}))$
via the target.

\end{itemize}

Tensoring \eqref{e:temp quot via action} over $\Sph_{\cG,x}^{\on{spec}}$ with $\bC$
we obtain a functor
\begin{equation} \label{e:temp quot via action C}
\bC_{\on{temp}}\to 
\on{Funct}_{\QCoh((\LS^\mer_{\cG,x})^\wedge_\reg)}\left(\QCoh(\LS^\reg_{\cG,x}), 
\QCoh(\LS^\reg_{\cG,x})\underset{\Sph_{\cG,x}^{\on{spec}}}\otimes \bC\right)
\end{equation}

From \lemref{l:temp quot via action} and the fact that $\Sph_{\cG,x}^{\on{spec}}$ is rigid we obtain:

\begin{cor} \label{c:temp quot via action C} 
The functor \eqref{e:temp quot via action C} is an equivalence.
\end{cor} 

\sssec{}

We will say that $\bC$ is \emph{tempered}
if the action of $\Sph_{\cG,x}^{\on{spec}}$ on $\bC$ factors via $\Sph_{\cG,\on{temp},x}^{\on{spec}}$.

\medskip

This is equivalent to the condition that the functors \eqref{e:temp adj} are mutually inverse 
equivalences.

\begin{rem}
Note that \corref{c:temp quot via action C} says that if $\bC$ is tempered, then it can be recovered
from the $\QCoh((\LS^\mer_{\cG,x})^\wedge_\reg)$-module category
$$\QCoh(\LS^\reg_{\cG,x})\underset{\Sph_{\cG,x}^{\on{spec}}}\otimes \bC$$
by applying the functor
\begin{multline*}
\on{Funct}_{\QCoh((\LS^\mer_{\cG,x})^\wedge_\reg)}\left(\QCoh(\LS^\reg_{\cG,x}),-\right): \\
\QCoh((\LS^\mer_{\cG,x})^\wedge_\reg)\mmod\to \Sph_{\cG,\on{temp},x}^{\on{spec}}\mmod \hookrightarrow \Sph_{\cG,x}^{\on{spec}}\mmod.
\end{multline*}

\end{rem} 

\sssec{}

From \corref{c:temp quot via action C} we obtain:

\begin{cor} \label{c:temp quot via action C 1}
The functor
$$\bC\mapsto \QCoh(\LS^\reg_{\cG,x})\underset{\Sph_{\cG,x}^{\on{spec}}}\otimes \bC,\quad 
\Sph_{\cG,x}^{\on{spec}}\mmod\to \DGCat$$
is conservative, when restricted to the subcategory
$$\Sph_{\cG,\on{temp},x}^{\on{spec}}\mmod \subset \Sph_{\cG,x}^{\on{spec}}\mmod.$$
\end{cor}

\sssec{} 

Let \begin{equation} \label{e:project on temp}
\Sph_{G,x}\to \Sph_{G,\on{temp},x}
\end{equation}
be the colocalization corresponding to the colocalization
$$\Sph_{\cG,x}^{\on{spec}}\to \Sph_{\cG,\on{temp},x}^{\on{spec}}$$
(we can use either $\Sat_G$ or $\Sat_{G,\tau}$ to identify $\Sph_G$ with $\Sph_\cG^{\on{spec}}$; the resulting
colocalizations are the same). 

\medskip

The definitions and results from this subsection render automatically to the setting of $\Sph_{G,x}$-module categories. 

\ssec{Proof of \thmref{t:critical FLE pointwise}} \label{ss:proof of local FLE} 

We are now ready to prove the pointwise version of the FLE.

\sssec{}

Recall (see \thmref{t:FLE and Sat}) that the functor 
\begin{equation} \label{e:FLE x again}
\FLE_{G,\crit}:\KL(G)_{\crit,x}\to \IndCoh^*(\Op^\mf_{\cG,x})
\end{equation} 
intertwines the actions of $\Sph_{G,x}$ on the left-hand side with
the $\Sph^{\on{spec}}_{G,x}$-action on the right-hand side (via $\on{Sat}_G$),
and makes the diagram 
\begin{equation} \label{e:FLE compat x}
\CD
\Whit_*(G)_x\underset{\Sph_{G,x}}\otimes \KL(G)_{\crit,x} @>{\FLE_{\cG,\infty}^{-1}\otimes \FLE_{G,\crit}}>> 
\Rep(\cG) \underset{\Sph_{\cG,x}^{\on{spec}}}\otimes \IndCoh^*(\Op_{\cG,x}^\mf) \\
@V\text{\eqref{e:pairing 1}}VV @VV{\text{\eqref{e:pairing 2}}}V \\
\IndCoh^*(\Op_{\cG,x}^{\on{mer}})  @>{\on{Id}}>>   \IndCoh^*(\Op_{\cG,x}^{\on{mer}})
\endCD
\end{equation} 
commute. 

\medskip

Note that \propref{p:recover compl via Sph} says that the right vertical arrow in \eqref{e:FLE compat x} is fully faithful 
with essential image equal to
$$\IndCoh^*(\Op_{\cG,x}^{\on{mer}})_\mf\subset  \IndCoh^*(\Op_{\cG,x}^{\on{mer}}).$$

\medskip

We will prove:

\begin{prop} \label{p:image of pairing}
The left vertical arrow in \eqref{e:FLE compat x} is fully faithful 
with essential image equal to
$$\IndCoh^*(\Op_{\cG,x}^{\on{mer}})_\mf\subset  \IndCoh^*(\Op_{\cG,x}^{\on{mer}}).$$
\end{prop}

\sssec{}

We now claim:

\begin{prop} \label{p:opers is temp}
The $\Sph^{\on{spec}}_{\cG,x}$-module category $\IndCoh^*(\Op^\mf_{\cG,x})$ is tempered.
\end{prop}

\begin{proof} 

Recall the functor 
\begin{equation} \label{e:IndCoh mf as QCoh ten prod again}
\QCoh(\LS^\reg_{\cG,\ul{x}})\underset{\QCoh((\LS^\mer_{\cG,\ul{x}})^\wedge_{\on{reg}})}\otimes
\IndCoh^*((\Op^\mer_{\cG,\ul{x}})^\wedge_{\on{mon-free}})\to
\IndCoh^*(\Op^\mf_{\cG,\ul{x}}).
\end{equation}

It intertwines the actions of $\Sph^{\on{spec}}_{\cG,x}$, where the action on the left-hand side
is via the first factor. 

\medskip

Now, according to \propref{p:IndCoh mf as QCoh ten prod}, the functor \eqref{e:IndCoh mf as QCoh ten prod again}
is an equivalence. Hence, it is enough to show that the action of $\Sph^{\on{spec}}_{\cG,x}$ on the left-hand side of
\eqref{e:IndCoh mf as QCoh ten prod again} factors through $\Sph^{\on{spec}}_{\cG,\on{temp},x}$. 

\medskip

However, this follows from the fact that the action of $\Sph^{\on{spec}}_{\cG,\on{temp},x}$ on 
$$\QCoh(\LS^\reg_{\cG,\ul{x}})\simeq \Rep(\cG)$$
factors through $\Sph^{\on{spec}}_{\cG,\on{temp},x}$. 

\end{proof}

\sssec{}

Finally, we claim: 

\begin{prop} \label{p:KL is temp}
The $\Sph_{G,x}$-module category $\KL(G)_{\kappa,x}$ is tempered.
\end{prop}

\medskip

We now observe that the combination of Propositions \ref{p:image of pairing}, \ref{p:opers is temp} and \ref{p:KL is temp},
together with \corref{c:temp quot via action C 1}, implies that \eqref{e:FLE x again} is an equivalence. 

\qed[Pointwise FLE]

\medskip

The rest of this section is devoted to the proof of Propositions \ref{p:image of pairing} and \ref{p:KL is temp}.

\ssec{Proof of \propref{p:image of pairing}} 

\sssec{} \label{sss:sph vs Whit sph}

Let $\bC$ be a category equipped with a $\fL(G)_x$-action at the critical level. 
Consider now the category 
$$\Sph(\bC):=\bC^{\fL^+(G)_x}.$$ 

Interpreting $\bC^{\fL^+(G)_x}$ as
$$\on{Funct}_{\fL(G)_x}(\Dmod_{\frac{1}{2}}(\Gr_{G,x}),\bC),$$
we obtain that $\Sph(\bC)$ carries a natural action of
$$\Sph_{G,x}\simeq \End_{\fL(G)_x}\left(\Dmod_{\frac{1}{2}}(\Gr_{G,x})\right).$$

\sssec{}

Denote 
$$\bC^{\on{Sph-gen}}:\Dmod_{\frac{1}{2}}(\Gr_{G,x})\underset{\Sph_{G,x}}\otimes \Sph(\bC).$$

We have a tautological functor 
\begin{equation} \label{e:Sph gen to all}
\bC^{\on{Sph-gen}}\to \bC
\end{equation} 
commuting with the $\fL(G)_x$-action.

\medskip

The following is a standard result that results from the ind-properness of the affine Grassmannian:

\begin{lem}  \label{l:Sph gen to all}
The functor \eqref{e:Sph gen to all} is fully faithful and admits a continuous right 
adjoint.\footnote{This right adjoint automatically commutes with the $\fL(G)_x$-action, essentially
because $\fL(G)_x$ is a group.}
\end{lem}

\sssec{}

We shall say that $\bC$ is \emph{spherically generated} if the functor \eqref{e:Sph gen to all} is an
equivalence.

\medskip

This is equivalent to requiring that $\bC$ is generated, as a category acted 
on by $\fL(G)_x$, by the essential image of the forgetful functor
$$\Sph(\bC)\to \bC.$$

\sssec{}

The above definitions apply when we replace $\fL(G)_x$ by $\fL(G)_{\rho(\omega_X),x}$. Applying
to both sides of \eqref{e:Sph gen to all} the functor $\Whit_*(-)$, we obtain a functor
\begin{equation} \label{e:pairing 1 prel abs again}
\Whit_*(G)\underset{\Sph_{G,x}}\otimes \Sph(\bC) \to \Whit_*(\bC),
\end{equation} 
i.e., the functor \eqref{e:pairing 1 prel abs}.

\medskip

From \lemref{l:Sph gen to all} we obtain:

\begin{cor} \label{c:Sph gen to all Whit}
The functor \eqref{e:pairing 1 prel abs again} is fully faithful.
\end{cor}

\sssec{}

We apply \corref{c:Sph gen to all Whit} to
$$\bC:=\hg\mod_{\crit,\rho(\omega_X),x}.$$

Hence, we obtain that the left vertical arrow in \eqref{e:FLE compat x}
is fully faithful. Thus, to complete the proof of \propref{p:image of pairing}, it suffices to show that the essential 
image of \eqref{e:pairing 1} is contained in and generates 
$$\IndCoh^*(\Op^\mer_{\cG,x})_\mf\subset \IndCoh^*(\Op^\mer_{\cG,x}).$$

\sssec{}

Note that the essential image of
$$\KL(G)_{\crit,x} \overset{\one_{\Whit_*(G)}\otimes \on{Id}}\longrightarrow 
\Whit_*(G)_x\underset{\Sph_{G,x}}\otimes \KL(G)_{\crit,x}$$
generates the target. 

\medskip

Indeed, this follows by interpreting the above functor as 
\begin{multline*} 
\KL(G)_{\crit,x}  \simeq \Rep(\cG)\underset{\Rep(\cG),\on{Sat}_G^{\on{nv}}}\otimes \KL(G)_{\crit,x}
\overset{\FLE_{\cG,\infty}\otimes \on{Id}}\simeq \\
\simeq \Whit_*(G)_x\underset{\on{Sat}_{G,\tau}^{\on{nv}},\Rep(\cG),\on{Sat}_G^{\on{nv}}}\otimes \KL(G)_{\crit,x}\to 
 \Whit_*(G)_x\underset{\on{Sat}_{G,\tau}^{-1},\Sph^{\on{spec}}_{\cG,x},\on{Sat}_G^{-1}}\otimes \KL(G)_{\crit,x} \simeq \\
\simeq \Whit_*(G)_x\underset{\Sph_{G,x}}\otimes \KL(G)_{\crit,x}.
\end{multline*}

\medskip

Since the essential image of $\FLE_{G,x}$ is contained in $\IndCoh^*(\Op^\mer_{\cG,x})_\mf$, we obtain that
so is the essential image of \eqref{e:pairing 1}. 

\sssec{}

Hence, it remains to prove the following: 

\begin{lem} \label{l:DS sph gen}
The essential image $\KL(G)_{\crit,\rho(\omega_X),x}$ under
$$\KL(G)_{\crit,\rho(\omega_X),x}\to \hg\mod_{\crit,\rho(\omega_X),x} \overset{\DS^{\on{enh,rfnd}}}\longrightarrow  \IndCoh^*(\Op^\mf_{\cG,x})$$
generates $\IndCoh^*(\Op^\mer_{\cG,x})_\mf$. 
\end{lem}

\begin{proof}

We prove the lemma by matching the generators.

\medskip

The compact generators of $\KL(G)_{\crit,\rho(\omega_X),x}$ are the Weyl modules
\begin{equation} \label{e:Weyl}
\BV^\lambda_\crit:=\ind_{\fL^+(G)}^{(\hg,\fL^+(G))_\crit}(V^\clambda),
\end{equation}
where $V^\clambda\in \Rep(G)$ is the irreducible representation of highest weight $\lambda$.

\medskip

According to \cite{FG3}, the image of $\BV^\lambda_\crit$ under $\DS^{\on{enh,rfnd}}$ is the structure sheaf 
$\CO_{\Op^{\lambda-\!\reg}_{\cG,x}}$ of the subscheme
$$\Op^{\lambda-\!\reg}_{\cG,x}\subset \Op^\mer_{\cG,x}$$
of $\clambda$-opers.

\medskip

We have
$$^{\on{red}}(\Op^\mf_{\cG,x})=\underset{\lambda}\sqcup\, \Op^{\lambda-\!\reg}_{\cG,x},$$
from which it is clear that the objects
$$\CO_{\Op^{\lambda-\!\reg}_{\cG,x}}\in \IndCoh^*(\Op^\mer_{\cG,x})$$
are the (compact) generators of $\IndCoh^*(\Op^\mer_{\cG,x})_\mf$.

\end{proof} 

\ssec{Proof of \propref{p:KL is temp}} \label{ss:KL is temp}

\sssec{}

We need to show that the functor
\begin{equation} \label{e:KL is temp}
\KL(G)_{\crit,x}:=\Sph(\hg\mod_{\crit,,x})\to  \Sph(\hg\mod_{\crit,x})_{\on{temp}}=:
(\KL(G)_{\crit,x})_{\on{temp}}
\end{equation} 
is an equivalence.

\medskip

First, note that by combining \corref{c:temp quot via action C} with Propositions \ref{p:image of pairing} and \ref{p:opers is temp} 
we obtain that the functor $\FLE_{G,\crit}$ factors as 
$$\KL(G)_{\crit,x} \overset{\text{\eqref{e:KL is temp}}}\to 
(\KL(G)_{\crit,x})_{\on{temp}} \overset{\FLE_{G,\crit,\on{temp}}}\to \IndCoh^*(\Op^\mf_{\cG,x}),$$
where $\FLE_{G,\crit,\on{temp}}$ is an equivalence.

\medskip

Combined with \propref{p:FLE pres comp}, we obtain that the functor \eqref{e:KL is temp} preserves
compactness. Since \eqref{e:KL is temp} is a colocalization, we obtain that \eqref{e:KL is temp}
restricts to a colocalization on compacts. Hence, it is sufficient to prove
that \eqref{e:KL is temp} is conservative on $(\KL(G)_{\crit,x})^c$.  

\medskip

To prove the latter, it is sufficient to prove that the functor $\FLE_{G,\crit}$ is conservative on $(\KL(G)_{\crit,x})^c$. 

\sssec{}

Since $(\KL(G)_{\crit,x})^c\subset (\KL(G)_{\crit,x})^b$
and since $\FLE_{G,\crit}$ is t-exact (by \corref{c:FLE t exact}), its suffices to show that 
$\FLE_{G,\crit}$ is conservative on $(\KL(G)_{\crit,x})^\heartsuit$. 

\medskip

Using the fact that $\FLE_{G,\crit,\on{temp}}$ is an equivalence, we obtain that it is enough to show
that the functor \eqref{e:KL is temp} is conservative on $(\KL(G)_{\crit,x})^\heartsuit$. 

\sssec{}

Let $\on{temp}_{\KL(G)_{\crit,,x}}$ denote the temperization functor (see \secref{sss:temp}). We will prove:

\begin{lem} \label{l:temp KL}
The functor $\on{temp}_{\KL(G)_{\crit,x}}$ is right t-exact, and the counit map
$$\on{temp}_{\KL(G)_{\crit,x}}\to \on{Id}$$
induces an isomorphism on $H^0$ when applied to objects in $(\KL(G)_{\crit,x})^\heartsuit$.
\end{lem}

The lemma immediately implies the conservativity of 
\eqref{e:KL is temp} on $(\KL(G)_{\crit,x})^\heartsuit$. 

\qed[\propref{p:KL is temp}]

\sssec{Proof of \lemref{l:temp KL}}

The assertion of the lemma holds more generally for a $\Sph^{\on{spec}}_{\cG,x}$-module category, equipped with a t-structure,
such that the action functor is t-exact. The corresponding property for $\KL(G)_{\crit,x}\simeq \KL(G)_{\crit,x}$ 
is guaranteed by \corref{c:Hecke action t-exact}.

\medskip

Consider the temperization functor $\on{temp}_{\Sph^{\on{spec}}_{\cG,x}}$ on $\Sph^{\on{spec}}_{\cG,x}$ itself, i.e., the composition of
\begin{multline*}
\Sph^{\on{spec}}_{\cG,x}=\IndCoh(\on{Hecke}_{\cG,x}^{\on{spec,loc}})\overset{\Psi_{\on{Hecke}_{\cG,x}^{\on{spec,loc}}}}\longrightarrow \\
\to \QCoh(\on{Hecke}_{\cG,x}^{\on{spec,loc}})\overset{\Xi_{\on{Hecke}_{\cG,x}^{\on{spec,loc}}}}\longrightarrow \IndCoh(\on{Hecke}_{\cG,x}^{\on{spec,loc}})=\Sph^{\on{spec}}_{\cG,x}.
\end{multline*}

\medskip

It suffices to show that the object
$$\on{temp}_{\Sph^{\on{spec}}_{\cG,x}}(\one_{\Sph^{\on{spec}}_{\cG,x}})$$
lives in cohomological degrees $\leq 0$ and that its $0$th cohomology maps isomorphically to $\one_{\Sph^{\on{spec}}_{\cG,x}}$.

\medskip

However, this is a general property of the composition $\Xi_Y\circ \Psi_Y$ on an eventually coconnective stack locally almost of finite type. 
Namely, for every $\CF\in \IndCoh(\CY)$, the counit of the adjuction
$$\Xi_Y\circ \Psi_Y(\CF)\to \CF$$
induces an isomorphism on the truncation $\tau^{\geq -n}$ for \emph{any} $n$. 

\qed[\lemref{l:temp KL}]

\ssec{Spherical vs Whittaker}

This subsection is not logically necessary for the sequel, but it carries an ideological significance.
Here we explain how to realize the pointwise $\FLE_{G,\crit}$ functor as the base change of the functor $\ol\DS^{\on{enh,rfnd}}$
along $\LS_{\cG,x}^\reg\to \LS^\mer_{\cG,x}$. 

\sssec{} \label{sss:Whit ten over LS}

Let $\bC$ be a category equipped with a $\fL(G)_{\rho(\omega_X),x}$-action, and assume that 
$\bC$ is spherically generated. 

\medskip

Consider the corresponding category $$\Whit_*(\bC):=\bC_{\fL(N)_{\rho(\omega_X),x},\chi}.$$

We claim that $\Whit_*(\bC)$ has a natural structure of module category over $\QCoh((\LS^\mer_{\cG,x})^\wedge_\reg)$. 

\sssec{}

Indeed, we have
\begin{equation} \label{e:Whi Sph gen}
\Whit_*(\bC)\simeq \Whit_*(G)_x\underset{\Sph_{G,x}}\otimes \Sph(\bC),
\end{equation}
so it is enough to show that $\Whit_*(G)_x$ carries a structure of $\QCoh((\LS^\mer_{\cG,x})^\wedge_\reg)$-module
category in a way that commutes with the $\Sph_{G,x}$-action.

\medskip

We identify
$$\Whit_*(G)_x\overset{\FLE_{\cG,\infty}}\simeq \Rep(\cG)\simeq \QCoh(\LS^\reg_{\cG,x}),$$
where the $\Sph_{G,x}$-action on $\Whit_*(G)_x$ corresponds to the natural $\Sph^{\on{spec}}_{\cG,x}$-action 
on $\QCoh(\LS^\reg_{\cG,x})$ via $\Sat_{G,\tau}$ (see \corref{c:geom Sat tau}). 

\medskip

The desired $\QCoh((\LS^\mer_{\cG,x})^\wedge_\reg)$-module structure on $\Whit(G)_x$ comes from the natural action of
$\QCoh((\LS^\mer_{\cG,x})^\wedge_\reg)$-action on $\QCoh(\LS^\reg_{\cG,x})$, 
which naturally commutes with the 
$\Sph^{\on{spec}}_{\cG,x}$-action. 

\sssec{} \label{sss:abs FLE}

For $\bC$ as above, consider the category 
\begin{equation} \label{e:Funct unr to Whit}
\on{Funct}_{\QCoh((\LS^\mer_{\cG,x})^\wedge_\reg)}(\QCoh(\LS^\reg_{\cG,x}),\Whit_*(\bC)).
\end{equation} 

Recall the functor \eqref{e:pre-FLE -1}
\begin{equation} \label{e:pre-FLE -1 again}
\Sph(\bC)\to \Whit_*(\bC).
\end{equation} 

\sssec{}

We claim:

\begin{prop} \label{p:Sph to unr Whit}
The functor \eqref{e:pre-FLE -1 again} lifts to a functor
\begin{equation} \label{e:Sph to unr Whit}
\Sph(\bC)\to \on{Funct}_{\QCoh((\LS^\mer_{\cG,x})^\wedge_\reg)}(\QCoh(\LS^\reg_{\cG,x}),\Whit_*(\bC)).
\end{equation}
Moreover, the functor \eqref{e:Sph to unr Whit} factors as
\begin{equation} \label{e:Sph to unr Whit bis}
\Sph(\bC)\to \Sph(\bC)_{\on{temp}} \overset{\sim}\to
\on{Funct}_{\QCoh((\LS^\mer_{\cG,x})^\wedge_\reg)}(\QCoh(\LS^\reg_{\cG,x}),\Whit_*(\bC)),
\end{equation}
where the second arrow is an equivalence. 
\end{prop} 

\begin{proof}

By \eqref{e:Whi Sph gen}, we have
\begin{multline} \label{e:Funct unr to Whit 1}
\on{Funct}_{\QCoh((\LS^\mer_{\cG,x})^\wedge_\reg)}(\QCoh(\LS^\reg_{\cG,x}),\Whit_*(\bC))
\simeq  \\
\simeq \on{Funct}_{\QCoh((\LS^\mer_{\cG,x})^\wedge_\reg)}(\QCoh(\LS^\reg_{\cG,x}),\Whit_*(G)_x)\underset{\Sph_{G,x}}\otimes \Sph(\bC)
\overset{\FLE_{G,\infty}}\simeq \\
\simeq \on{Funct}_{\QCoh((\LS^\mer_{\cG,x})^\wedge_\reg)}(\QCoh(\LS^\reg_{\cG,x}),\QCoh(\LS^\reg_{\cG,x}))
\underset{\Sph_{G,x}}\otimes \Sph(\bC)\simeq \Sph^{\on{spec}}_{\cG,\on{temp}.x}\underset{\Sph_{G,x}}\otimes \Sph(\bC).
\end{multline}

\end{proof} 

%
%
%
%
%
%
%
%
%

\sssec{}

Take $\bC:=\hg\mod_{\crit,\rho(\omega_X),x}^{\on{Sph-gen}}$. So we can regard 
$$\Whit_*(\hg\mod_{\crit,\rho(\omega_X),x}^{\on{Sph-gen}})$$
as acted on by $\QCoh((\LS^\mer_{\cG,x})^\wedge_\reg)$ via the recipe of \secref{sss:Whit ten over LS}. 

\medskip

Let us regard 
$$\IndCoh^*(\Op^\mer_{\cG,x})_\mf\simeq \IndCoh^*(\Op^\mer_{\cG,x})^\wedge_\mf$$
as acted on by $\QCoh((\LS^\mer_{\cG,x})^\wedge_\reg)$ via pullback along
$$\fr:(\Op^\mer_{\cG,x})^\wedge_\mf\to (\LS^\mer_{\cG,x})^\wedge_\reg.$$

Recall now that according to \propref{p:image of pairing}, the functor $\ol\DS^{\on{enh,rfnd}}$ gives rise to an equivalence
\begin{equation} \label{e:DS mf}
\Whit_*(\hg\mod_{\crit,\rho(\omega_X),x}^{\on{Sph-gen}}) \overset{\ol\DS^{\on{enh,rfnd}}}\longrightarrow  
\IndCoh^*(\Op^\mer_{\cG,x})_\mf.
\end{equation}

We claim:

\begin{prop} \label{p:structure over LS compat}
The equivalence \eqref{e:DS mf} is compatible with the above actions of $\QCoh((\LS^\mer_{\cG,x})^\wedge_\reg)$.
\end{prop}

\begin{proof}

Follows from the fact that the functor \eqref{e:DS mf} can be recovered from $\FLE_{G,\crit}$ by applying to both sides
the functor
$$\QCoh(\LS^\reg_{G,x})\underset{\Sph^{\on{spec}}_{\cG,x}}\otimes -.$$

\end{proof} 

\begin{rem} \label{r:structure over LS}

The above proof of \propref{p:structure over LS compat} relies relies on \thmref{t:FLE and Sat} as an essential
ingredient. 

\medskip

In \secref{s:from LS} we will give a different proof of \propref{p:structure over LS compat}, by showing that the 
functor \eqref{e:DS mf} is
compatible with the structure of \emph{factorization module category} over $\Rep(\cG)$. 

\end{rem} 
 
\sssec{} \label{sss:reinterp local FLE}

From \propref{p:structure over LS compat} we obtain that we can regard the FLE functor
$$\FLE_{G,\crit}:\KL(G)_{\crit,x}\to \IndCoh^*(\Op^\mf_{\cG,x})$$
as obtained from the functor \eqref{e:Sph to unr Whit} for
$$\bC:=\hg\mod_{\crit,\rho(\omega_X),x}^{\on{Sph-gen}}$$
by precomposing with 
$$\KL(G)_{\crit,x} \overset{\alpha_{\rho(\omega_X),\on{taut}}}\longrightarrow \KL(G)_{\crit,\rho(\omega_X),x}=\Sph(\hg\mod_{\crit,\rho(\omega_X),x})$$
and post-composing with 
\begin{multline} \label{e:FLE fund}
\on{Funct}_{\QCoh((\LS^\mer_{\cG,x})^\wedge_\reg)}\left(\QCoh(\LS^\reg_{\cG,x}),\Whit_*(\hg\mod_{\crit,\rho(\omega_X),x}^{\on{Sph-gen}})\right)
\overset{\text{\propref{p:structure over LS compat}}}\simeq  \\
\simeq \on{Funct}_{\QCoh((\LS^\mer_{\cG,x})^\wedge_\reg)}\left(\QCoh(\LS^\reg_{\cG,x}),\IndCoh^*(\Op^\mer_{\cG,x})_\mf\right) \simeq
\IndCoh^*(\Op^\mf_{\cG,x}),
\end{multline}
where the last equivalence is obtained as follows:
\begin{multline*} 
\on{Funct}_{\QCoh((\LS^\mer_{\cG,x})^\wedge_\reg)}\left(\QCoh(\LS^\reg_{\cG,x}),\IndCoh^*(\Op^\mer_{\cG,x})_\mf\right)  \simeq \\
\on{Funct}_{\IndCoh^!(\Op^\mer_{\cG,x})_\mf}
\left(\IndCoh^!(\Op^\mer_{\cG,x})_\mf\underset{\QCoh((\LS^\mer_{\cG,x})^\wedge_\reg)}\otimes 
\QCoh(\LS^\reg_{\cG,x}),\IndCoh^*(\Op^\mer_{\cG,x})_\mf\right) \\
\overset{\text{\propref{p:IndCoh mf as QCoh ten prod}}}\simeq
\on{Funct}_{\IndCoh^!(\Op^\mer_{\cG,x})_\mf}
\left(\IndCoh^!(\Op^\mf_{\cG,x}),\IndCoh^*(\Op^\mer_{\cG,x})_\mf\right) \simeq \\
\simeq \on{Funct}_{\IndCoh^!(\Op^\mer_{\cG,x})}
\left(\IndCoh^!(\Op^\mf_{\cG,x}),\IndCoh^*(\Op^\mer_{\cG,x})\right)\overset{\text{\lemref{l:hom category}}}\simeq \IndCoh^*(\Op^\mf_{\cG,x}).
\end{multline*} 

\begin{rem}

The proof of the pointwise FLE given in \secref{ss:proof of local FLE} relied on \propref{p:structure over LS compat},
and hence on \thmref{t:FLE and Sat} ingredient.
 
\medskip

As was mentioned in Remark \ref{r:structure over LS}, we will supply a different proof of 
\propref{p:structure over LS compat}.

\medskip

This allows us to give a proof of the pointwise FLE, avoiding \thmref{t:FLE and Sat}: 

\medskip

Then one interprets the pointwise $\FLE_{G,\crit}$ functor as in \secref{sss:reinterp local FLE} above. The assertion
of the pointwise FLE follows now by combining Propositions \ref{p:KL is temp} and \ref{p:Sph to unr Whit}. 

\end{rem}

\sssec{}

We remark also that one can deduce the pointwise version of \thmref{t:FLE and Sat} from \propref{p:structure over LS compat}
by the following argument:

\medskip

First, we note that in the context of \secref{sss:abs FLE}, the functor
$$\Sph(\bC)\to \on{Funct}_{\QCoh((\LS^\mer_{\cG,x})^\wedge_\reg)}(\QCoh(\LS^\reg_{\cG,x}),\Whit_*(\bC))$$
intertwines the action of $\Sph_{G,x}$ on $\Sph(\bC)$ and the action of $\Sph^{\on{spec}}_{\cG,x}$ on 
$$\on{Funct}_{\QCoh((\LS^\mer_{\cG,x})^\wedge_\reg)}(\QCoh(\LS^\reg_{\cG,x}),\Whit_*(\bC))$$
via the source.

\medskip

Unwinding the construction, we obtain that the following diagram commutes
$$
\CD
\Whit_*(G) \underset{\Sph_{G,x}}\otimes \Sph(\bC) @>{\text{\eqref{e:Whi Sph gen}}}>{\sim}>  \Whit_*(\bC) \\
@V{\FLE^{-1}_{\cG,\infty} \otimes \text{\eqref{e:Sph to unr Whit}}}VV @VV{\on{Id}}V \\
\QCoh(\LS^\reg_{\cG,x})\underset{\Sph^{\on{spec}}_{\cG,x}}\otimes 
\on{Funct}_{\QCoh((\LS^\mer_{\cG,x})^\wedge_\reg)}(\QCoh(\LS^\reg_{\cG,x}),\Whit_*(\bC)) @>>> \Whit_*(\bC),
\endCD
$$
where the bottom horizontal arrow is the natural evaluation functor. 

\medskip

This proves the desired compatibility, since the functor \eqref{e:FLE fund} is compatible
with the $\Sph^{\on{spec}}_{\cG,x}$-actions.

\section{Compatibility of the critical FLE with duality} \label{s:FLE and duality}

In this section we show that the FLE equivalence is compatible with the natural self-dualities of the two sides. 

\medskip

The proof proceeds along the following steps:

\begin{enumerate} 

\item We show that the self-duality on $\hg\mod_\crit$ is compatible with the $\IndCoh^!(``\Spec"(\fZ_\fg))$-action.
This boils down to a particular property of the factorization algebra $\CDO(G)_{\crit,\crit}$, given by 
\lemref{l:crit CDO};

\item We show that the equivalence $\ol\DS^{\on{enh,rfnd}}:\Whit_*(\hg\mod_{\crit,\rho(\omega_X)})\to \IndCoh^*(\Op_\cG^\mer)$
is compatible with the self-dualities of the two sides. We deduce this from the $\IndCoh^!(\Op^\mer_\cG)$-linearity 
(guaranteed by the previous point), combined with a general uniqueness statement;

\smallskip

\item Finally, we establish the compatibility of $\FLE_{G,\crit}$ with the self-dualities by essentially base-changing it from
the previous point along $\Op_\cG^\mf\to \Op_\cG^\mer$.

\end{enumerate}

As we highlight below, our methods are robust enough to show that established compatibility is automatically compatible with the actions of
$$\Sph_G\simeq \Sph^{\on{spec}}_\cG$$
on the two sides.

%
%

\ssec{Statement of the compatibility theorem} \label{ss:FLE and duality}

\sssec{}

Recall that according to \eqref{e:KL crit self dual}, we have a canonical identification
\begin{equation} \label{e:KL self-dual again}
(\KL(G)_\crit)^\vee \simeq \KL(G)_\crit.
\end{equation}

By construction, the equivalence \eqref{e:KL self-dual again} respects the actions of $\Sph_G$. 

\sssec{}

In addition, we have an equivalence
\begin{equation} \label{e:Op mf self-dual again}
\IndCoh^*(\Op_\cG^\mf)^\vee\simeq
\IndCoh^!(\Op_\cG^\mf)\overset{\Theta_{\Op^\mf_\cG}}\simeq 
\IndCoh^*(\Op_\cG^\mf).
\end{equation}

This equivalence respects the actions of $\Sph_\cG^{\on{spec}}$ (see Sects. \ref{sss:Sph and duality}, \ref{sss:Theta mf via duality} and 
\lemref{l:Theta and Sph}). 

\sssec{}

The goal of this section is to prove the following:

\begin{thm} \label{t:FLE and duality}
With respect to the identifications \eqref{e:KL self-dual again} and \eqref{e:Op mf self-dual again},
the (factorization) functor
$$(\FLE_{G,\crit})^\vee:\IndCoh^*(\Op_\cG^\mf)\to \KL(G)_\crit$$
identifies with 
$$\tau_G\circ (\FLE_{G,\crit})^{-1}.$$
Moreover, this identification of functors respects the compatibility with the actions of 
$$\Sph_G\overset{\Sat_{G,\tau}}\simeq \Sph_\cG^{\on{spec}}.$$
\end{thm}

One can rephrase \thmref{t:FLE and duality} as a commutative diagram of factorization categories
\begin{equation} \label{e:self duality opers mf}
\CD
\IndCoh^*(\Op_\cG^\mf)^\vee @>{\text{\eqref{e:Op mf self-dual again}}}>{\sim}> \IndCoh^*(\Op_\cG^\mf) \\
@V{(\FLE_{G,\crit})^\vee}V{\sim}V @A{\sim}A{\tau_G\circ \FLE_{G,\crit}}A \\
(\KL(G)_\crit)^\vee @>{\sim}>{\text{\eqref{e:KL self-dual again}}}> \KL(G)_\crit. 
\endCD
\end{equation}

\begin{rem}
Note the similarity between the statement of \thmref{t:FLE and duality} and \lemref{l:CS and duality}:
in both cases a non-tautological self-equivalence of the Whittaker side makes the FLE inverse to
its dual, up to the Chevalley involution.
\end{rem}

\begin{rem}
Note again that the appearance of the Chevalley involution in \thmref{t:FLE and duality} is in line with
the curse in \secref{sss:curse duality}.
\end{rem}

\sssec{}

The following assertion will not be needed in the sequel, but it provides a concrete perspective on
what \thmref{t:FLE and duality} really says. 

\medskip

Recall that the unit of the self-duality \eqref{e:KL self-dual again} is the object
$$\fCDO_{\crit,\crit}\in \KL(G)_\crit\otimes \KL(G)_\crit.$$

The unit of the self-duality of $\IndCoh^*(\Op_\cG^\mf)$ is 
$$(\Delta_{\Op_\cG^\mf})^\IndCoh_*(\omega^{*,\on{fake}}_{\Op^{\on{mon-free}}_\cG}),$$
where
$$\omega^{*,\on{fake}}_{\Op^{\on{mon-free}}_\cG}\in \IndCoh^*(\Op_\cG^\mf)$$
is as in \secref{sss:omega fake mf}, and $\Delta_{\Op_\cG^\mf}$ is the diagonal map of
$\Op^{\on{mon-free}}_\cG$. 

\medskip

Thus, from \thmref{t:FLE and duality}, we obtain:

\begin{cor}  \label{c:FLE CDO}
We have a canonical isomorphism (of factorization algebra objects)
$$(\FLE_{G,\crit}\otimes \FLE_{G,\crit})(\fCDO_{\crit,\crit})\simeq 
(\on{Id}\otimes \tau_G)\circ(\Delta_{\Op_\cG^\mf})^\IndCoh_*(\omega^{*,\on{fake}}_{\Op^{\on{mon-free}}_\cG}).$$
\end{cor} 

\begin{rem}
Note that, on the one hand, the statement of \corref{c:FLE CDO} is actually equivalent to \thmref{t:FLE and duality}
(without the $\Sph_G\simeq \Sph^{\on{spec}}_\cG$ compatibility). 

\medskip

One the other hand, the two factorization algebras appearing in \corref{c:FLE CDO} are \emph{classical}
(i.e., the corresponding chiral algebras lie in $\Dmod(X)^\heartsuit$), and one can actually prove the existence
of an isomorphism between them directly.

\end{rem}

\sssec{}

Applying the forgetful functor to $\Vect$, from \corref{c:FLE CDO} we obtain:

\begin{cor}
We have a canonical isomorphism of factorization algebras
\begin{equation} \label{e:FLE CDO}
(\DS\otimes \DS)\circ (\alpha_{\rho(\omega_X),\on{taut}}\otimes \alpha_{\rho(\omega_X),\on{taut}})(\fCDO_{\crit,\crit})
\simeq \Gamma^\IndCoh(\Op_\cG^\mf,\omega^{*,\on{fake}}_{\Op^{\on{mon-free}}_\cG}).
\end{equation} 
\end{cor}

\begin{rem}
We note that the factorization algebra 
$$\CB_G:=(\DS\otimes \DS)\circ (\alpha_{\rho(\omega_X),\on{taut}}\otimes \alpha_{\rho(\omega_X),\on{taut}})(\fCDO_{\crit,\crit})$$
was studied in \cite{FG2}. 

\medskip

One can view \eqref{e:FLE CDO} as an extension of the Feigin-Frenkel isomorphism $\on{FF}_G$: indeed
according to \lemref{l:crit CDO} below, the factorization algebra $\CB_\fg$ receives a homomorphism from
$$(\DS\otimes \DS)(\BV_{\fg,\crit,\rho(\omega_X)}\underset{\fz_\fg}\otimes \BV_{\fg,\crit,\rho(\omega_X)})\simeq \fz_\fg,$$
while
$$\CB^{\on{spec}}_\cG:= \Gamma^\IndCoh(\Op_\cG^\mf,\omega^{*,\on{fake}}_{\Op^{\on{mon-free}}_\cG})$$
receives a homomorphism from
$$\Gamma(\Op^\reg_\cG,\CO_{\Op^\reg_\cG})=:\CO_{\Op^\reg_\cG}.$$

\end{rem} 

%
%

\ssec{Feigin-Frenkel center and self-duality}

\sssec{}

Recall the duality identification 
\begin{equation} \label{e:KM self-dual again}
(\hg\mod_\crit)^\vee \simeq \hg\mod_\crit.
\end{equation}
of \eqref{e:KM crit self dual}.

\medskip

Recall also that $\hg\mod_\crit$ carries a canonical action of $\IndCoh^!(``\Spec"(\fZ_\fg))$,
see \secref{ss:action of center}. Since the category $\IndCoh^!(``\Spec"(\fZ_\fg))$ is \emph{symmetric}
monoidal, this action induces an action of $\IndCoh^!(``\Spec"(\fZ_\fg))$ on $(\hg\mod_\crit)^\vee$.

\sssec{}

The goal of this subsection is to prove the following:

\begin{thmconstr} \label{t:center and duality}
The identification \eqref{e:KM self-dual again} carries a canonical structure of
compatibility with the $\IndCoh^!(``\Spec"(\fZ_\fg))$-actions, up to the automorphism of 
$``\Spec"(\fZ_\fg)$ induced by the Chevalley involution $\tau_G$. 
\end{thmconstr}

The rest of this subsection is devoted to the proof of this theorem.

\sssec{}

Let
$$\on{u}_{\hg\mod_\crit,\hg\mod_\crit}\in \hg\mod_\crit\otimes \hg\mod_\crit$$
be the unit of the self-duality. The statement of \thmref{t:center and duality} is equivalent to the assertion
that $\on{u}_{\hg\mod_\crit,\hg\mod_\crit}$ can be lifted to an object of the category 
$$\on{Funct}_{\IndCoh^!(``\Spec"(\fZ_\fg))\otimes \IndCoh^!(``\Spec"(\fZ_\fg))}(\IndCoh^!(``\Spec"(\fZ_\fg)),\hg\mod_\crit\otimes \hg\mod_\crit)$$
where $\IndCoh^!(``\Spec"(\fZ_\fg))\otimes \IndCoh^!(``\Spec"(\fZ_\fg))$ acts on $\IndCoh^!(``\Spec"(\fZ_\fg))$ via $\tau_G$ on one of the factors. 

\sssec{}

Consider the (factorization) category 
$$\Dmod(\fL(G))_{\crit,\crit}$$
of critically twisted D-modules on the loop group. We have a naturally defined (factorization) functor
\begin{equation} \label{e:sect loop group}
\Gamma^\IndCoh(\fL(G),-):\Dmod(\fL(G))_{\crit,\crit}\to \hg\mod_\crit\otimes \hg\mod_\crit.
\end{equation}

We have:
$$\on{u}_{\hg\mod_\crit,\hg\mod_\crit}\simeq \Gamma^\IndCoh(\fL(G),\delta_{1,\fL(G)}),$$
where $\delta_{1,\fL(G)}\in \Dmod(\fL(G))_{\crit,\crit}$ is the $\delta$-function at the origin. 

\medskip 

The required property of $\on{u}_{\hg\mod_\crit,\hg\mod_\crit}$ follows from the next general assertion:

\begin{prop} \label{p:Dmods LG and center}
The functor \eqref{e:sect loop group} factors as
\begin{multline*}
\Dmod(\fL(G))_{\crit,\crit}\to \\
\to \on{Funct}_{\IndCoh^!(``\Spec"(\fZ_\fg))\otimes \IndCoh^!(``\Spec"(\fZ_\fg))}(\IndCoh^!(``\Spec"(\fZ_\fg)),\hg\mod_\crit\otimes \hg\mod_\crit)\to \\
\to \hg\mod_\crit\otimes \hg\mod_\crit.
\end{multline*}
\end{prop}

\qed[\thmref{t:center and duality}]

\ssec{Proof of \propref{p:Dmods LG and center}}

\sssec{}

First, we record the initial input, from which we will deduce \propref{p:Dmods LG and center}. Recall
the (factorization algebra) object
$$\fCDO(G)_{\crit,\crit}\in \KL(G)_\crit\otimes \KL(G)_\crit.$$

By a slight abuse of notation we will denote by the same symbol the image of $\fCDO(G)_{\crit,\crit}$ under the
forgetful functor
$$\KL(G)_\crit\otimes \KL(G)_\crit\to \hg\mod_\crit\otimes \hg\mod_\crit.$$

\medskip

Let us denote by $\CDO(G)_{\crit,\crit}$ the image of $\fCDO(G)_{\crit,\crit}$ along the further forgetful functor
$$\oblv_{\hg\times \hg}:\hg\mod_\crit\otimes \hg\mod_\crit\to \Vect.$$

\sssec{}

The unit of $\fCDO(G)_{\crit,\crit}$ as a factorization algebra in $\KL(G)_\crit\otimes \KL(G)_\crit$ is a map (of factorization algebras)
\begin{equation} \label{e:Vac to CDO}
\on{Vac}(G)_\crit\otimes \on{Vac}(G)_\crit\to \fCDO(G)_{\crit,\crit},
\end{equation}
which gives rise to a map
\begin{equation} \label{e:Vac to CDO bis}
\BV_{\fg,\crit}\otimes \BV_{\fg,\crit}\to \CDO(G)_{\crit,\crit}.
\end{equation}

The following was established in \cite[Theorem 5.4]{FG1}:

\begin{lem} \label{l:crit CDO}
The map \eqref{e:Vac to CDO bis} factors as
$$\BV_{\fg,\crit}\otimes \BV_{\fg,\crit}\to \BV_{\fg,\crit}\underset{\fz_\fg}\otimes \BV_{\fg,\crit}\to \CDO(G)_{\crit,\crit},$$
where the $\fz_\fg$-action on one of the tensor factors is twisted by $\tau_G$.
\end{lem}

Since the factorization algebras involved lie in the heart of the t-structure, from \lemref{l:crit CDO} we obtain:

\begin{cor}  \label{c:crit CDO}
The map \eqref{e:Vac to CDO} factors as
$$\on{Vac}(G)_\crit\otimes \on{Vac}(G)_\crit\to \on{Vac}(G)_\crit\underset{\fz_\fg}\otimes \on{Vac}(G)_\crit\to \fCDO(G)_{\crit,\crit},$$
where the $\fz_\fg$-action on the right tensor factor is twisted by $\tau_G$.
\end{cor}

We will now show how to use \corref{c:crit CDO} to prove \propref{p:Dmods LG and center}. 

\sssec{}

Note that $\Dmod(\fL(G))_{\crit,\crit}$ is equipped with a t-structure (see \secref{sss:t-str fact} for what this means 
in the factorization setting), so that the object $\delta_{\fL^+(G)\subset \fL(G)}\in \Dmod(\fL(G))_{\crit,\crit}$
lies in the heart.

\medskip

The category $\Dmod(\fL(G))_{\crit,\crit}$ is compactly generated by objects that lie in $\Dmod(\fL(G))_{\crit,\crit}^{>-\infty}$.
Hence, in order to prove \propref{p:Dmods LG and center}, it suffices to show that the restriction of the functor
$\Gamma^\IndCoh(\fL(G),-)$ to $\Dmod(\fL(G))_{\crit,\crit}^{>-\infty}$ factors as
\begin{multline*}
\Dmod(\fL(G))^{>-\infty}_{\crit,\crit}\to \\
\to \on{Funct}_{\IndCoh^!(``\Spec"(\fZ_\fg))\otimes \IndCoh^!(``\Spec"(\fZ_\fg))}(\IndCoh^!(``\Spec"(\fZ_\fg)),\hg\mod_\crit\otimes \hg\mod_\crit)\to \\
\to \hg\mod_\crit\otimes \hg\mod_\crit.
\end{multline*}

\sssec{}

Consider the (factorization) functor $\Gamma^\IndCoh(\fL(G),-)$. It sends the factorization unit
$$\one_{\Dmod(\fL(G))_{\crit,\crit}}\simeq \delta_{\fL^+(G)\subset \fL(G)}\in \Gamma^\IndCoh(\fL(G),-)$$
to 
$$\fCDO(G)_{\crit,\crit}\in \hg\mod_\crit\otimes \hg\mod_\crit.$$

In particular, by \secref{sss:fact alg from fact cat gen}, it upgrades to a functor
$$\Gamma^\IndCoh(\fL(G),-)^{\on{enh}}: 
\Dmod(\fL(G))_{\crit,\crit}\to \fCDO\mod^{\on{fact}}(\hg\mod_\crit\otimes \hg\mod_\crit).$$

\sssec{}

Applying \corref{c:crit CDO}, we obtain that the functor $\Gamma^\IndCoh(\fL(G),-)$ factors as
\begin{multline*}
\Dmod(\fL(G))_{\crit,\crit}\to \\
\to (\on{Vac}(G)_\crit\underset{\fz_\fg}\otimes \on{Vac}(G)_\crit)\mod^{\on{fact}}(\hg\mod_\crit\otimes \hg\mod_\crit)\to
\hg\mod_\crit\otimes \hg\mod_\crit.
\end{multline*}

In particular, the restriction of the functor $\Gamma^\IndCoh(\fL(G),-)$ to $\Dmod(\fL(G))_{\crit,\crit}^{>-\infty}$ factors via 
the forgetful functor
\begin{multline} \label{e:over z}
\left((\on{Vac}(G)_\crit\underset{\fz_\fg}\otimes \on{Vac}(G)_\crit)\mod^{\on{fact}}(\hg\mod_\crit\otimes \hg\mod_\crit)\right)^{>-\infty}\to \\
\to \biggl(\hg\mod_\crit\otimes \hg\mod_\crit\biggr)^{>-\infty}.
\end{multline}

\sssec{}

However, unwinding the construction of the $\IndCoh^!(``\Spec"(\fZ_\fg))$-action on $\hg\mod_\crit$, we obtain that the functor
\eqref{e:over z}
factors as
\begin{multline*} 
\left((\on{Vac}(G)_\crit\underset{\fz_\fg}\otimes \on{Vac}(G)_\crit)\mod^{\on{fact}}(\hg\mod_\crit\otimes \hg\mod_\crit)\right)^{>-\infty}\to \\
\to \on{Funct}_{\IndCoh^!(``\Spec"(\fZ_\fg))\otimes \IndCoh^!(``\Spec"(\fZ_\fg))}(\IndCoh^!(``\Spec"(\fZ_\fg)),\hg\mod_\crit\otimes \hg\mod_\crit)\to \\
\to \hg\mod_\crit\otimes \hg\mod_\crit.
\end{multline*}

\qed[\propref{p:Dmods LG and center}]

\ssec{Self-duality on opers via Kac-Moody}

\sssec{}

Recall the identification \eqref{e:dual Whit}. Applying this to the category $\hg\mod_{\crit,\rho(\omega)}$ and using the identification
\eqref{e:KM self-dual again}, we obtain an identification
\begin{equation} \label{e:dual Whit KM}
\Whit_*(\hg\mod_{\crit,\rho(\omega)})^\vee\simeq  \Whit^!(\hg\mod_{\crit,\rho(\omega)}),
\end{equation}
which fits into the commutative diagram
$$
\CD
\Whit_*(\hg\mod_{\crit,\rho(\omega)})^\vee @>{\text{\eqref{e:dual Whit KM}}}>> \Whit^!(\hg\mod_{\crit,\rho(\omega)}) \\
@VVV @VVV \\
(\hg\mod_{\crit,\rho(\omega)})^\vee @>{\text{\eqref{e:KM self-dual again}}}>> \hg\mod_{\crit,\rho(\omega)},
\endCD
$$
where the left vertical arrow is the dual of the projection
$$\hg\mod_{\crit,\rho(\omega)}\to \Whit_*(\hg\mod_{\crit,\rho(\omega)}).$$

\medskip

Since the vertical arrows in the above diagram are fully faithful, from 
\thmref{t:center and duality}\footnote{Here we apply \thmref{t:center and duality} in the twisted setting, i.e., to 
$\hg\mod_{\crit,\rho(\omega)}$ instead of $\hg\mod_\crit$.} (combined with \corref{c:z action on Wht(KM)}(a,b)), we obtain:

\begin{cor} \label{c:dual Whit KM z}
The functor \eqref{e:dual Whit KM} is equipped with a natural $\IndCoh^!(``\Spec"(\fZ_\fg))$-linear structure, up to the automorphism of 
$``\Spec"(\fZ_\fg)$, induced by the Chevalley involution $\tau_G$. 
\end{cor}

\sssec{} \label{sss:dual Whit KM z}

Recall now identification of \thmref{t:Whit self-dual gen}. Applying this to the category $\hg\mod_{\crit,\rho(\omega)}$, we obtain an identification
\begin{equation} \label{e:Theta Whit KM}
\Theta_{\Whit(\hg\mod_{\crit,\rho(\omega)})}:\Whit_*(\hg\mod_{\crit,\rho(\omega)})\simeq  \Whit^!(\hg\mod_{\crit,\rho(\omega)}).
\end{equation} 

Concatenating \eqref{e:dual Whit KM} with \eqref{e:Theta Whit KM} we obtain an identification
\begin{equation} \label{e:self-dual Whit KM}
\Whit_*(\hg\mod_{\crit,\rho(\omega)})^\vee\simeq  \Whit_*(\hg\mod_{\crit,\rho(\omega)}).
\end{equation} 

Combining with Corollaries \ref{c:dual Whit KM z} and \ref{c:z action on Wht(KM)}(c), we obtain that the functor \eqref{e:self-dual Whit KM} is also
endowed with an $\IndCoh^!(``\Spec"(\fZ_\fg))$-linear structure, up to the automorphism of 
$``\Spec"(\fZ_\fg)$ induced by the Chevalley involution $\tau_G$. 

\sssec{}

Recall now that we have an identification
\begin{equation} \label{e:Op self-dual mer again}
\IndCoh^*(\Op_\cG^\mer)^\vee\simeq
\IndCoh^!(\Op_\cG^\mer)\overset{\Theta_{\Op^\mer_\cG}}\simeq 
\IndCoh^*(\Op_\cG^\mer).
\end{equation}

We will prove:

\begin{thm} \label{t:self duality opers mer}
With respect to the identifications \eqref{e:self-dual Whit KM} and \eqref{e:Op self-dual mer again}, the functor dual to
$$\ol\DS^{\on{enh,rfnd}}:\Whit_*(\hg\mod_{\crit,\rho(\omega)})\overset{\sim}\to \IndCoh^*(\Op_\cG^\mer)$$
identifies with $\tau_G\circ (\ol\DS^{\on{enh,rfnd}})^{-1}$, compatibly with the actions of 
$$\IndCoh^!(``\Spec"(\fZ_\fg))\overset{\on{FF}_G}\simeq \IndCoh^!(\Op_\cG^\mer).$$
\end{thm} 

One can rephrase \thmref{t:FLE and duality} as a commutative diagram
\begin{equation} \label{e:self duality opers mer}
\CD
\IndCoh^*(\Op_\cG^\mer)^\vee @>{\text{\eqref{e:Op self-dual mer again}}}>{\sim}> \IndCoh^*(\Op_\cG^\mer) \\
@V{(\ol\DS^{\on{enh,rfnd}})^\vee}V{\sim}V @A{\sim}A{\tau_G\circ \ol\DS^{\on{enh,rfnd}}}A \\
\Whit_*(\hg\mod_{\crit,\rho(\omega)})^\vee @>{\sim}>{\text{\eqref{e:self-dual Whit KM}}}> \Whit_*(\hg\mod_{\crit,\rho(\omega)}).
\endCD
\end{equation}

\begin{rem}
As we shall see below, \thmref{t:self duality opers mer} is actually easy. However, it can be seen as a particular case of a 
conjecture, proposed by G.~Dhillon, which says that at any level $\kappa$, the self-dualities of the (renormalized) categories
of factorization modules 
$$\CW_{\fg,\kappa}\mod^{\on{fact}}\simeq \CW_{\cg,\check\kappa}\mod^{\on{fact}}$$
that come from the identifications
$$\CW_{\fg,\kappa}\mod^{\on{fact}}=\Whit_*(\hg\mod_\kappa) \text{ and } 
\CW_{\cg,\check\kappa}\mod^{\on{fact}}=\Whit_*(\wh{\cg}\mod_{\check\kappa})$$
and \thmref{t:Whit self-dual gen}, respectively, agree. 

\medskip

For non-critical $\kappa$, this conjecture is completely open. What makes it tractable at the critical
level is precisely the interpretation of $\CW_{\fg,\crit}$ as the Feigin-Frenkel center.
\end{rem}

\ssec{Proof of \thmref{t:self duality opers mer}}

\sssec{}

Consider the $\IndCoh^!(\Op_\cG^\mer)$-linear self-equivalence of $\IndCoh^*(\Op_\cG^\mer)$ obtained by going clockwise along the edges of 
\eqref{e:self duality opers mer}. We need to show that this functor is isomorphic to the identity. 

\medskip

Using the equivalence 
$\Theta_{\Op^\mer_\cG}$, and further, the equivalence
$$\Upsilon_{\Op^\mer_\cG}:\QCoh(\Op^\mer_\cG)\to \IndCoh^!(\Op^\mer_\cG)$$
of \propref{p:Ups mon-free Op}, we can transform the above $\IndCoh^!(\Op_\cG^\mer)$-linear self-equivalence 
of the category $\IndCoh^*(\Op_\cG^\mer)$
into a $\QCoh(\Op^\mer_\cG)$-linear self-equivalence of $\QCoh(\Op^\mer_\cG)$.

\medskip

Such a self-equivalence is given by a (graded) line bundle, to be denoted $\CL_{\Op^\mer_\cG}$, 
on $\Op^\mer_\cG$. Since all the equivalences in sight are compatible
with factorization, this line bundle has a natural factorization structure. 

\medskip

We will now show that any such (graded) line bundle is automatically trivial. 

\sssec{} \label{sss:triviality on Ran}

The question is local, 
so let us choose a $\cG$-oper $\sigma$ on $X$. 

\medskip

The datum of $\sigma$ gives rise
to a section 
$$\sigma_\Ran:\Ran\to \Op^\reg_{\cG,\Ran}\to \Op^\mer_{\cG,\Ran},$$
compatible with factorization.

\medskip

Set $\CL_{\Ran,\sigma}:=(\sigma_\Ran)^*(\CL_{\Op^\mer_\cG})$. This is a 
factorization line bundle on $\Ran$. We claim that it is automatically trivial. 

\sssec{}

Indeed, write 
$$\Ran\simeq \underset{I}{\on{colim}}\, X^I_\dr,$$
where the index $I$ runs over the category
of non-empty finite sets and surjective maps.

\medskip

Set $\CL_{X^I_\dr,\sigma}:=\CL_{\Ran,\sigma}|_{X^I_\dr}$. 
The collection of local systems 
\begin{equation} \label{e:loc sys X I}
I\rightsquigarrow \CL_{X^I_\dr,\sigma}
\end{equation} 
is compatible with the factorization structure. 

\medskip

In particular,
$$\CL_{X^I_\dr,\sigma}|_{\oX^I}\simeq (\CL_{X_\dr,\sigma})^{\boxtimes I}|_{\oX^I},$$
where $\oX^I\subset X^I$ is the complement of the diagonal divisor.

\medskip

Hence, 
\begin{equation} \label{e:L I dr}
\CL_{X^I_\dr,\sigma} \simeq (\CL_{X_\dr,\sigma})^{\boxtimes I},
\end{equation} 
compatibly with the factorization structure. 

\medskip

Consider \eqref{e:L I dr} for $I=\{1,2\}$. Restricting to the diagonal $X\to X\times X$, we obtain 
$$\CL_{X_\dr,\sigma} \simeq (\CL_{X_\dr,\sigma})^{\otimes 2}.$$

Hence, $\CL_{X_\dr,\sigma}$ is trivial. By \eqref{e:L I dr}, this trivializes the system \eqref{e:loc sys X I}. 

\sssec{}

For a fixed $I$, consider the fiber products
$$\Op^\mer_{\cG,X^I_\dr}:=X^I_\dr\underset{\Ran}\times  \Op^\mer_{\cG,\Ran} \text{ and }
\Op^\mer_{\cG,X^I}:=X^I\underset{\Ran}\times  \Op^\mer_{\cG,\Ran}$$
and the line bundles 
$$\CL_{\Op^\mer_{\cG,X^I_\dr}}:=\CL_{\Op^\mer_{\cG,\Ran}}|_{\Op^\reg_{\cG,X^I_\dr}} \text{ and }
\CL_{\Op^\mer_{\cG,X^I}}:=\CL_{\Op^\mer_{\cG,\Ran}}|_{\Op^\reg_{\cG,X^I}}.$$

The map $\Op^\mer_{\cG,X^I}\to X^I$ is a Zariski-locally trivial fibration with 
ind-pro-affine spaces as fibers. Since $X^I$ is smooth, 
we obtain that $\CL_{\Op^\reg_{\cG,X^I}}$ descends to a canonically defined
line bundle $\CL_{X^I}$ on $X^I$. Moreover, the collection
$$I\rightsquigarrow \CL_{X^I}$$
has a natural factorization structure.

\medskip

Furthermore, by construction, we have 
$$\CL_{X^I}\simeq \CL_{X^I_\dr,\sigma}|_{X^I},$$
compatibly with factorization.

\medskip

Hence, by \secref{sss:triviality on Ran}, we obtain that the system
\begin{equation} \label{e:line bundles on XI}
I\mapsto \CL_{\Op^\mer_{\cG,X^I}}
\end{equation}
is canonically trivial, compatibly with factorization.

\medskip

It remains to show that this trivialization descends to a trivialization of the system
\begin{equation} \label{e:local systems on XI}
I\mapsto \CL_{\Op^\mer_{\cG,X^I_\dr}}.
\end{equation}

\sssec{}

Note a priori, the obstruction to triviality is given by a function on $\CL_{\Op^\mer_{\cG,X^I}}$
with values in the pullback of the sheaf of 1-forms on $X^I$; denote this function by $\alpha_{X^I}$.

\medskip

Note that by factorization, 
$$\alpha_{X^I}|_{\oX^I}\simeq \alpha_X^{\boxtimes I}.$$

Hence, it is enough to show that $\alpha_X=0$. 

\sssec{} \label{sss:conn form triv on reg}

We will first show that $\alpha_{X}|_{\Op^\reg_{\cG,X}}$ is trivial. 

\medskip

Recall that the factorization scheme $\Op_\cG^\reg$ is counital (see \secref{sss:counital fact spaces} for what this means).
In particular, there exist canonical projections
$$p_i:\Op^\reg_{\cG,X^2}\to \Op^\reg_{\cG,X}, \quad i=1,2$$
covering the two projections $p_i:X^2\to X$. 

\medskip

We claim that 
$$\alpha_{X^2}=p_1^*(\alpha_X)+p_2^*(\alpha_X).$$

Indeed, the equality takes place over $\oX^2$, by factorization, and hence over the
entire $X^2$ by density.

\sssec{} \label{sss:conn form pullback one factor}
 
Recall now that $\Op^\mer_\cG$ has a unital-in-correspondences structure relative to $\Op^\reg_\cG$ 
(see \secref{sss:rel unital} for what this means). We claim that the connection forms
$\alpha_I$ are compatible with this structure in the following sense:

\medskip

For an injection of finite sets $\phi:I_1\to I_2$, let 
$$\Op^\mer_{\cG,X^{I_1}} \overset{\on{pr}^\Op_{\on{small}}}\longleftarrow \Op^{\mer\rightsquigarrow \reg}_{\cG,X^\phi} 
\overset{\on{pr}^\Op_{\on{big}}}\longrightarrow 
\Op^\mer_{\cG,X^{I_2}}$$
be the correspondence covering
$$X^{I_1} \overset{\Delta_\phi}\longleftarrow X^{I_2} \overset{\on{id}}\longrightarrow X^{I_2}.$$

We claim that
\begin{equation} \label{e:source and target form}
(\on{pr}^\Op_{\on{small}})^*(\alpha_{X^{I_1}})=(\on{pr}^\Op_{\on{big}})^*(\alpha_{X^{I_2}}).
\end{equation} 

Indeed, write $I_2=I_1\sqcup J$, and let
$$(X^{I_1}\times X^J)_{\on{disj}}\subset X^{I_1}\times X^J$$
be the corresponding open subset. 

\medskip

It suffices that the equality \eqref{e:source and target form} takes place over 
$$\Op^{\mer\rightsquigarrow\reg}_{\cG,X^\phi}\underset{X^{I_2}}\times (X^{I_1}\times X^J)_{\on{disj}}.$$

\medskip

We have
\begin{equation} \label{e:source and target form 1}
\Op^{\mer\rightsquigarrow\reg}_{\cG,X^\phi}\underset{X^{I_2}}\times (X^{I_1}\times X^J)_{\on{disj}}\simeq
(\Op^\mer_{\cG,X^{I_1}}  \times \Op^\reg_{\cG,X^J})\underset{X^{I_1}\times X^J}\times (X^{I_1}\times X^J)_{\on{disj}}
\end{equation} 
and 
$$\Op^\mer_{\cG,X^{I_2}}\underset{X^{I_2}}\times (X^{I_1}\times X^J)_{\on{disj}})\simeq 
(\Op^\mer_{\cG,X^{I_1}}  \times \Op^\mer_{\cG,X^J})\underset{X^{I_1}\times X^J}\times (X^{I_1}\times X^J)_{\on{disj}},$$
where the map $\on{pr}^\Op_{\on{small}}$ identifies with projection on the first factor in \eqref{e:source and target form 1}, 
and the map $\on{pr}^\Op_{\on{big}}$ is the inclusion
$$\Op^\reg_{\cG,X^J}\to  \Op^\mer_{\cG,X^J}$$
along the second factor.

\medskip

By factorization, we obtain that 
$$(\on{pr}^\Op_{\on{big}})^*(\alpha_{X^{I_2}})|_{\Op^{\mer\rightsquigarrow\reg}_{\cG,X^\phi}\underset{X^{I_2}}\times (X^{I_1}\times X^J)_{\on{disj}}}$$
equals the sum of $(\on{pr}^\Op_{\on{small}})^*(\alpha_{X^{I_1}})$ and the pullback of $\alpha_{X^J}|_{\Op^\reg_{\cG,X^J}}$
along the projection of \eqref{e:source and target form 1} on the second factor.

\medskip

However, the latter is zero by \secref{sss:conn form triv on reg}. 

\sssec{}

We are now ready to show that $\alpha_X=0$. In doing so we will mimic the argument in \cite[Proposition 3.4.7]{BD2}.

\medskip

Write $\Omega^1_{X^2}$ as 
$$\omega_X\boxtimes \CO_X \oplus \CO_X\boxtimes \omega_X.$$

We will show that the restriction of $\alpha_{X^2}$ to 
$$\Op^\mer_{\cG,X}\simeq X\underset{X^2}\times \Op^\mer_{\cG,X^2},$$
viewed as a function on $\Op^\mer_{\cG,X}$ with values in the pullback of $\Omega^1_{X^2}$,
takes values both in the pullback of $\omega_X\boxtimes \CO_X$ \emph{ and} in the pullback of
$\CO_X\boxtimes \omega_X$. This would implies that this restriction is $0$, and hence also that $\alpha_X=0$.

\medskip

We will in fact show that the restriction of $\alpha_{X^2}$ to 
\begin{equation} \label{e:source and target form 2}
(X^2)^\wedge \underset{X^2}\times \Op^\mer_{\cG,X^2}
\end{equation} 
is $0$, where $(X^2)^\wedge$ is the formal completion of the diagonal in $X^2$. 

\sssec{}

By symmetry, it suffices to show that the restriction of $\alpha_{X^2}$ to 
takes values in the pullback of $\omega_X\boxtimes \CO_X$. 

\medskip

Consider the inclusion $I_1:=\{1\}\overset{\phi}\hookrightarrow \{1,2\}:=I_2$, and the corresponding map
$$\Op^{\mer\rightsquigarrow\reg}_{\cG,X^\phi} \overset{\on{pr}^\Op_{\on{big}}}\longrightarrow \Op^\mer_{\cG,X^2}.$$

This map is an isomorphism over $X\subset X^2$, and hence, induces an isomorphism
$$(X^2)^\wedge \underset{X^2}\times \Op^{\mer\rightsquigarrow\reg}_{\cG,X^\phi}\overset{\sim}\to 
(X^2)^\wedge \underset{X^2}\times \Op^\mer_{\cG,X^2}.$$

Hence, it is enough to prove that the pullback of $\alpha_{X^2}$ along $\on{target}$
takes values in the pullback of $\omega_X\boxtimes \CO_X$. However, this has been 
established in \secref{sss:conn form pullback one factor}.

\qed[\thmref{t:self duality opers mer}]

\ssec{Proof of \thmref{t:FLE and duality}}

By a slight abuse of notation we will use the symbol $\FLE_{G,\crit}$ for the functor \eqref{e:critical FLE omega}. 

\sssec{}

Consider the following diagram 
{\tiny
$$
\CD
(\KL(G)_{\crit,\rho(\omega_X)})^\vee @<{\tau_G\circ (\FLE_{G,\crit})^\vee}<< \IndCoh^*(\Op_\cG^\mf)^\vee \\
@VV{\sim}V \\
(\IndCoh^!(\Op^\mf_\cG)\underset{\IndCoh^!(\Op^\mer_\cG)}\otimes \Whit^!(\hg\mod_{\crit,\rho(\omega_X)}))^\vee & & @VV{\sim}V \\
@VV{\sim}V  \\
(\IndCoh^!(\Op^\mf_\cG)\underset{\IndCoh^!(\Op^\mer_\cG)}\otimes \Whit_*(\hg\mod_{\crit,\rho(\omega_X)}))^\vee
@<{(\on{Id}\otimes \tau_G\circ \ol\DS^{\on{enh,rfnd}})^\vee}<<
(\IndCoh^!(\Op^\mf_\cG)\underset{\IndCoh^!(\Op^\mer_\cG)}\otimes \IndCoh^*(\Op_\cG^\mer))^\vee \\
@V{\sim}VV @VV{\sim}V \\
\on{Funct}_{\IndCoh^!(\Op^\mer_\cG)}(\IndCoh^!(\Op^\mf_\cG),\Whit_*(\hg\mod_{\crit,\rho(\omega_X)})^\vee)
@<{\tau_G\circ (\ol\DS^{\on{enh,rfnd}})^\vee}<<
\on{Funct}_{\IndCoh^!(\Op^\mer_\cG)}(\IndCoh^!(\Op^\mf_\cG),\IndCoh^*(\Op_\cG^\mer)^\vee) \\
@V{\sim}VV @VV{\sim}V \\
\on{Funct}_{\IndCoh^!(\Op^\mer_\cG)}(\IndCoh^!(\Op^\mf_\cG),\Whit^!(\hg\mod_{\crit,\rho(\omega_X)})) & & 
\on{Funct}_{\IndCoh^!(\Op^\mer_\cG)}(\IndCoh^!(\Op^\mf_\cG),\IndCoh^!(\Op_\cG^\mer)) \\
@A{\sim}A{\Theta_{\Whit(\hg\mod_{\crit,\rho(\omega_X)})}}A @V{\Theta_{\Op^\mer_\cG}}V{\sim}V  \\
\on{Funct}_{\IndCoh^!(\Op^\mer_\cG)}(\IndCoh^!(\Op^\mf_\cG),\Whit_*(\hg\mod_{\crit,\rho(\omega_X)}))  @>{\ol\DS^{\on{enh,rfnd}}}>>
\on{Funct}_{\IndCoh^!(\Op^\mer_\cG)}(\IndCoh^!(\Op^\mf_\cG),\IndCoh_*(\Op_\cG^\mer)) \\
@A{\sim}AA @AA{\sim}A \\
\KL(G)_{\crit,\rho(\omega_X)} @>{\FLE_{G,\crit}}>> \IndCoh_*(\Op_\cG^\mf),
\endCD
$$}
in which the upper vertical arrows are the duals of
$$
\IndCoh^!(\Op^\mf_\cG)\underset{\IndCoh^!(\Op^\mer_\cG)}\otimes \Whit^!(\hg\mod_{\crit,\rho(\omega_X)})
\overset{\text{\eqref{e:pre-inverse FLE 0}}}\longrightarrow \KL(G)_{\crit,\rho(\omega_X)}$$
and
$$\IndCoh^!(\Op^\mf_\cG)\underset{\IndCoh^!(\Op^\mer_\cG)}\otimes \IndCoh^*(\Op_\cG^\mer)\to \IndCoh^*(\Op^\mf_\cG),$$
respectively. 

\medskip

We will show:

\begin{itemize}

\item The left vertical composition is the identification \eqref{e:KL self-dual again};

\item The right vertical composition is the identification \eqref{e:Op mf self-dual again}; 

\item All inner squares commute.

\end{itemize}

This will establish the commutativity of \eqref{e:self duality opers mf}.  The compatibility of this isomorphism 
with the actions of $$\Sph_G\overset{\Sat_{G,\tau}}\simeq \Sph_\cG^{\on{spec}}.$$
is automatic from the construction.

\sssec{The left vertical composition}

We need to establish the commutativity of the following diagram
$$
\CD
(\KL(G)_{\crit,\rho(\omega_X)})^\vee @>>>  (\IndCoh^!(\Op^\mf_\cG)\underset{\IndCoh^!(\Op^\mer_\cG)}\otimes \Whit^!(\hg\mod_{\crit,\rho(\omega_X)}))^\vee  \\
& & @VV{\sim}V \\
@V{\text{\eqref{e:KL self-dual again}}}VV  \on{Funct}_{\IndCoh^!(\Op^\mer_\cG)}(\IndCoh^!(\Op^\mf_\cG),\Whit^!(\hg\mod_{\crit,\rho(\omega_X)})^\vee) \\
& &  @VV{\sim}V \\
\KL(G)_{\crit,\rho(\omega_X)} @>>> \on{Funct}_{\IndCoh^!(\Op^\mer_\cG)}(\IndCoh^!(\Op^\mf_\cG),\Whit_*(\hg\mod_{\crit,\rho(\omega_X)})).
\endCD
$$

However, this is just the fact that in the context of \secref{sss:inverse mech}, the dual of the functor \eqref{e:Sph to Whit* again} for $\bC$
is the functor \eqref{e:Sph to Whit^! again} for $\bC^\vee$. 

\sssec{The right vertical composition}

The identification of the right vertical composition follows from \secref{sss:Theta mf via duality}.

\sssec{The top square}

Passing to the dual functors, we need to establish the commutativity of

{\tiny
$$
\CD
\KL(G)_{\crit,\rho(\omega_X)}) @>{\FLE_{G,\crit}}>> \IndCoh^*(\Op_\cG^\mf) \\
@AA{\sim}A \\
\IndCoh^!(\Op^\mf_\cG)\underset{\IndCoh^!(\Op^\mer_\cG)}\otimes \Whit^!(\hg\mod_{\crit,\rho(\omega_X)})  & & @AA{\sim}A \\
@AA{\sim}A  \\
\IndCoh^!(\Op^\mf_\cG)\underset{\IndCoh^!(\Op^\mer_\cG)}\otimes \Whit_*(\hg\mod_{\crit,\rho(\omega_X)})
@>{\on{Id}\otimes  \ol\DS^{\on{enh,rfnd}}}>{\sim}>
\IndCoh^!(\Op^\mf_\cG)\underset{\IndCoh^!(\Op^\mer_\cG)}\otimes \IndCoh^*(\Op_\cG^\mer)).
\endCD
$$}

However, this is the content of \propref{p:inverse FLE}. 

\sssec{The 2nd square from the top}

This square commutes tautologically.

\sssec{The bottom square}

The commutation follows from the definition of the functor $\FLE_{G,\crit}$.

\sssec{}

Finally, it remains to show that the 3rd square from the top commutes.\footnote{This is the only non-tautological
point in the proof.}

\medskip

However, the required commutation is given by  \thmref{t:self duality opers mer}. 

\qed[\thmref{t:FLE and duality}]

\newpage

\centerline{\bf Part II: Local-to-global constructions}

\section{The coefficient and Poincar\'e functor(s)} \label{s:coeff}

This section begins by introducing our main object of study: the critically twisted category of D-modules on $\Bun_G$.
In this section we will mostly think of its incarnation as $\Dmod_{\frac{1}{2}}(\Bun_G)$, see Remark \ref{r:1/2 vs crit},
as the main characters in this section are sheaf-theoretic in nature. 

\medskip

The focus of this section is Poincar\'e and Whattaker coefficient functors. In fact, there are two Poincar\'e functors
$$\on{Poinc}_{G,!}:\Whit^!(G)_\Ran\to \Dmod_{\frac{1}{2}}(\Bun_G) \text{ and } 
\on{Poinc}_{G,*}:\Whit_*(G)_\Ran\to \Dmod_{\frac{1}{2},\on{co}}(\Bun_G),$$
where $\Dmod_{\frac{1}{2},\on{co}}(\Bun_G)$ is the dual category of $\Dmod_{\frac{1}{2}}(\Bun_G)$. These
two functors are Verdier-conjugate: the dual functor of $\on{Poinc}_{G,*}$ is isomorphic to the right adjoint of $\on{Poinc}_{G,!}$;
this is the functor
$$\on{coeff}_G:\Dmod_{\frac{1}{2}}(\Bun_G)\to \Whit^!(G)_\Ran.$$


\medskip

One can also give a global interpretation of the above functors, where instead of the affine Grassmannian,
one uses the \emph{twisted Drinfeld compactification}
$$\overline\Bun_{N,\rho(\omega_X)}\to \Bun_G.$$
This is how the global geometric Whittaker model had been mostly approached so far (see, e.g., \cite{Ga1}).
The two approaches are, however, equivalent (see \cite{Ga6}).

\medskip

For the purposes of this paper, we will only explicitly need the global interpretation of the 
\emph{vacuum} cases of the above functors, see \secref{ss:glob Whit}.

\ssec{Twisted D-modules on \texorpdfstring{$\Bun_G$}{BunG}}

\sssec{}

Let $\det_{\Bun_G}$ be the determinant line bundle on $\Bun_G$, normalized so that it sends a $G$-bundle $\CP_G$ 
to 
$$\det\left(\Gamma(X,\fg_{\CP_G})\right) \otimes \det\left(\Gamma(X,\fg_{\CP^0_G})\right)^{\otimes -1},$$
where $\CP^0_G$ is the trivial bundle.

\sssec{}

Note that we have
$$\pi^*(\det_{\Bun_G})\simeq \det_{\Gr_{G,\Ran}},$$
where $\pi$ denotes the projection
\begin{equation} \label{e:Gr to Bun}
\Gr_{G,\Ran}\to \Bun_G.
\end{equation} 

\sssec{} \label{sss:crit 1/2 can}

Note also that up to the (constant) line $\det\left(\Gamma(X,\fg_{\CP^0_G})\right)$, the line bundle $\det_{\Bun_G}$
identifies with the canonical line bundle on $\Bun_G$.

\sssec{} \label{sss:crit on Bun_G}

Let $\crit$ be half of the de Rham twisting defined by $\det_{\Bun_G}$, i.e.,
$$\crit=\frac{1}{2}\cdot \on{dlog}(\det_{\Bun_G}).$$

We will denote by 
$$\Dmod_\crit(\Bun_G)$$
the corresponding category of twisted D-modules.

\medskip 

Note that by \secref{sss:crit 1/2 can}, the critical twisting on $\Bun_G$ is canonically isomorphic to
the half-canonical twisting. 

\sssec{} 

As in \secref{sss:crit vs 1/2}, we obtain a \emph{canonical} identification
\begin{equation} \label{e:1/2 vs crit}
\Dmod_{\frac{1}{2}}(\Bun_G) \overset{\sim}\to \Dmod_\crit(\Bun_G),
\end{equation}
where $\Dmod_{\frac{1}{2}}(\Bun_G)$ is the short-hand for
$$\Dmod_{\det^{\frac{1}{2}}_{\Bun_G}}(\Gr_G),$$
cf. \secref{sss:1/2 Gr}.

\begin{rem} \label{r:critical is int}

According to \cite[Sect. 4]{BD1}, the choice of $\omega_X^{\otimes \frac{1}{2}}$ gives rise to a choice of the square root of 
$\det_{\Bun_G}$ as a line bundle. This allows us to identify $\Dmod_\crit(\Bun_G)$ (or equivalently $\Dmod_{\frac{1}{2}}(\Bun_G)$)
with the untwisted category
$\Dmod(\Bun_G)$. 

\medskip

However, we will avoid using this identification.

\end{rem}

\sssec{}

Pullback along $\pi$ defines functors
$$\pi^!:\Dmod_\crit(\Bun_G)\to \Dmod_\crit(\Gr_{G,\Ran})$$
and 
$$\pi^!:\Dmod_{\frac{1}{2}}(\Bun_G)\to \Dmod_{\frac{1}{2}}(\Gr_{G,\Ran}),$$
so that the diagram 
$$
\CD
 \Dmod_{\frac{1}{2}}(\Gr_{G,\Ran}) @>>>  \Dmod_\crit(\Gr_{G,\Ran}) \\
 @A{\pi^!}AA @AA{\pi^!}A \\ 
 \Dmod_{\frac{1}{2}}(\Bun_G) @>>> \Dmod_\crit(\Bun_G) 
\endCD
$$
commutes. 

\ssec{Restricting to (twists of) \texorpdfstring{$\Bun_N$}{BunN}}

\sssec{} \label{sss:crit to BunN}

Let $\CP_T$ be any $T$-bundle. Consider the stack
\begin{equation} \label{e:twisted Bun N}
\Bun_{N,\CP_T}\simeq \Bun_B\underset{\Bun_T}\times \on{pt},
\end{equation}
where $\on{pt}\to \Bun_T$ is the point $\CP_T$.

\medskip

Denote by $\fp$ the map
$$\Bun_{N,\CP_T}\to \Bun_G.$$
Note that the pullback of $\det_{\Bun_G}$ along this map is canonically constant. 
Denote the resulting line by
$$\fl_{G,N_{\CP_T}},$$
see \cite[Sect. 1.3.1]{GLC1}. 

\sssec{}

Note that on the one hand, we obtain an identification
\begin{equation} \label{e:triv gerbe Bun N P T 1}
\Dmod_{\frac{1}{2}\cdot\on{dlog}(\fl_{G,N_{\CP_T}})}(\Bun_{N,\CP_T})  \overset{\sim}\to \Dmod(\Bun_{N,\CP_T}), 
\end{equation} 
since the $\on{dlog}$ map over $\on{pt}$ is trivial.

\medskip

On the other hand, recall (see \cite[Proposition 1.3.3]{GLC1}) that the line bundle $\fl_{G,N_{\CP_T}}$ admits a canonical square root, to be denoted
$\fl^{\otimes \frac{1}{2}}_{G,N_{\CP_T}}$. Hence, we
obtain \emph{another} identification
\begin{equation} \label{e:triv gerbe Bun N P T 2}
\Dmod_{\frac{1}{2}\cdot\on{dlog}(\fl_{G,N_{\CP_T}})}(\Bun_{N,\CP_T}) \overset{\text{\eqref{e:etale vs dR twistings}}}\simeq 
\Dmod_{\fl^{\frac{1}{2}}_{G,N_{\CP_T}}}(\Bun_{N,\CP_T}) \simeq \Dmod(\Bun_{N,\CP_T}).
\end{equation} 

The discrepancy between the two trivializations is given by tensoring by the line $\fl^{\otimes \frac{1}{2}}_{G,N_{\CP_T}}$.

\sssec{}

We will denote by $\fp_\crit^!$ the functor
\begin{multline*}
\Dmod_\crit(\Bun_G) \to \Dmod_{\frac{1}{2}\cdot \on{dlog}(\det_{\Bun_G}|_{\Bun_{N,\CP_T}})}(\Bun_{N,\CP_T})= 
\Dmod_{\frac{1}{2}\cdot\on{dlog}(\fl_{G,N_{\CP_T}})}(\Bun_{N,\CP_T}) \to \\
\overset{\text{\eqref{e:triv gerbe Bun N P T 1}}}\longrightarrow \Dmod(\Bun_{N,\CP_T}).
\end{multline*}

We will denote by $\fp^!$ the functor 
\begin{multline*}
\Dmod_\crit(\Bun_G) \to \Dmod_{\frac{1}{2}\cdot \on{dlog}(\det_{\Bun_G}|_{\Bun_{N,\CP_T}})}(\Bun_{N,\CP_T})= 
\Dmod_{\frac{1}{2}\cdot\on{dlog}(\fl_{G,N_{\CP_T}})}(\Bun_{N,\CP_T}) \to \\
\overset{\text{\eqref{e:triv gerbe Bun N P T 2}}}\longrightarrow \Dmod(\Bun_{N,\CP_T}).
\end{multline*}

Thus, we obtain a commutative diagram
\begin{equation} \label{e:tw pullback to BunN diag N P T}
\CD
\Dmod(\Bun_{N,\CP_T}) @>{\otimes \fl^{\otimes \frac{1}{2}}_{G,N_{\CP_T}}}>> \Dmod(\Bun_{N,\CP_T})  \\
@A{\fp^!}AA  @AA{\fp^!_\crit}A  \\ 
\Dmod_{\frac{1}{2}}(\Bun_G) @>{\text{\eqref{e:1/2 vs crit}}}>> \Dmod_\crit(\Bun_G). 
\endCD
\end{equation}

\sssec{} \label{sss:pull back to 1/2 rho omega} 

In this section we will take $\CP_T:=\rho(\omega_X)$. Consider the corresponding line 
\begin{equation} \label{e:fl G N}
\fl^{\otimes \frac{1}{2}}_{G,N_{\rho(\omega_X)}}.
\end{equation}

\medskip

In this case, \eqref{e:tw pullback to BunN diag N P T} reads as 

\begin{equation} \label{e:tw pullback to BunN diag rho omega}
\CD
\Dmod(\Bun_{N,\rho(\omega_X)}) @>{\otimes \fl^{\otimes \frac{1}{2}}_{G,N_{\rho(\omega_X)}}}>> \Dmod(\Bun_{N,\rho(\omega_X)})  \\
@A{\fp^!}AA  @AA{\fp^!_\crit}A  \\ 
\Dmod_{\frac{1}{2}}(\Bun_G) @>{\text{\eqref{e:1/2 vs crit}}}>> \Dmod_\crit(\Bun_G). 
\endCD
\end{equation}

\ssec{The coefficient functor}

In this subsection we will recall the definition of the functor of Whittaker coefficient(s). 

\sssec{}

Consider the $\rho(\omega_X)$-twisted version of the map \eqref{e:Gr to Bun}
$$\Gr_{G,\rho(\omega_X),\Ran}\to \Bun_G.$$

\medskip

Due to the trivialization of the $\mu_2$-gerbe $\fl_{G,N_{\rho(\omega_X)}}^{\frac{1}{2}}$, we have
$$\pi^*(\det^{\frac{1}{2}}_{\Bun_G})\simeq \det^{\frac{1}{2}}_{\Gr_{G,\rho(\omega_X)}}$$
as $\mu_2$-gerbes on $\Gr_{G,\rho(\omega_X),\Ran}$.

\medskip

Hence, $\pi$ gives rise to a well-defined functor
$$\pi^!:\Dmod_{\frac{1}{2}}(\Bun_G) \to \Dmod_{\frac{1}{2}}(\Gr_{G,\rho(\omega_X),\Ran}).$$

\sssec{} 

In this subsection we will the functor of Whittaker coefficient(s), denoted $\on{coeff}_G$, which maps
$$\Dmod_{\frac{1}{2}}(\Bun_G) \to \Whit^!(G)_\Ran.$$
and is defined as follows.

\medskip

To simplify the notation, we will work over a fixed point $\ul{x}\in \Ran$. So we need to define
the functor
$$\on{coeff}_{G,\ul{x}}:\Dmod_{\frac{1}{2}}(\Bun_G) \to \Whit^!(G)_{\ul{x}}.$$

Denote by
$$\pi_{\ul{x}}:\Gr_{G,\rho(\omega_X),\ul{x}}\to \Bun_G$$
the restriction of $\pi$ along
\begin{equation} \label{e:Gr x to Gr Ran}
\Gr_{G,\rho(\omega_X),\ul{x}}\to \Gr_{G,\rho(\omega_X),\Ran}.
\end{equation}

\sssec{} \label{sss:N alpha}

Write $\fL(N)_{\rho(\omega_X),\ul{x}}$ as a filtered union of subschemes $N^\alpha$. For every $\alpha$, 
consider the functor
$$\Av^{(N^\alpha,\chi)}_*:\Dmod_{\frac{1}{2}}(\Gr_{G,\rho(\omega_X),\ul{x}})\to
\Dmod_{\frac{1}{2}}(\Gr_{G,\rho(\omega_X),\ul{x}})^{N^\alpha,\chi} \hookrightarrow \Dmod_{\frac{1}{2}}(\Gr_{G,\rho(\omega_X),\ul{x}}).$$

For $N^{\alpha}\subset N^{\alpha'}$, we have a canonically defined natural transformation
\begin{equation} \label{e:passage Whit av}
\Av^{(N^{\alpha'},\chi)}_*\to \Av^{(N^{\alpha},\chi)}_*.
\end{equation} 

We have the following (elementary) observation:

\begin{lem} \label{l:passage Whit av}
The natural transformation
$$\Av^{(N^{\alpha'},\chi)}_*\circ \pi^!_{\ul{x}}\to \Av^{(N^{\alpha},\chi)}_*\circ \pi^!_{\ul{x}},$$
induced by \eqref{e:passage Whit av}, is an isomorphism when $N^{\alpha}$ 
is large enough.
\end{lem}

\begin{proof}

Let 
$$\on{Sect}(X-\ul{x},N_{\rho(\omega_X)})\subset \on{Sect}(X-\ul{x},G_{\rho(\omega_X)})$$ be the group ind-schemes of sections of 
$$N_{\rho(\omega_X)}\subset G_{\rho(\omega_X)}$$ over $X-\ul{x}$.  Laurent expansion defines maps
$$\on{Sect}(X-\ul{x},N_{\rho(\omega_X)})\to \fL(N)_{\rho(\omega_X),\ul{x}} \text{ and }
\on{Sect}(X-\ul{x},G_{\rho(\omega_X)})\to \fL(G)_{\rho(\omega_X),\ul{x}}.$$

\medskip

Note that the restriction of $\chi$ to $\on{Sect}(X-\ul{x},N_{\rho(\omega_X)})$ is trivial. 

\medskip

For $\CF\in \Dmod_{\frac{1}{2}}(\Bun_G)$, the pullback 
$$\pi^!_{\ul{x}}(\CF)\in \Dmod_{\frac{1}{2}}(\Gr_{G,\rho(\omega_X),\ul{x}})$$
is $\on{Sect}(X-\ul{x},G_{\rho(\omega_X)})$-equivariant, and hence 
$\on{Sect}(X-\ul{x},N_{\rho(\omega_X)})$-equivariant.

\medskip

Hence, the map in the lemma is an isomorphism any time 
\begin{equation} \label{e:alpha large}
N^\alpha\cdot \on{Sect}(X-\ul{x},N_{\rho(\omega_X)})=\fL(N)_{\rho(\omega_X),\ul{x}}.
\end{equation} 

\end{proof} 

\sssec{}

By \lemref{e:passage Whit av}, for $N^\alpha$ large enough, the functor 
\begin{equation} \label{e:Whit prel}
\Av^{(N^{\alpha},\chi)}_*\circ \pi_{\ul{x}}^!
\end{equation}
does not depend on the choice of $N^\alpha$. In particular, its essential image is contained in
$$\underset{N^\alpha}\cap\, \Dmod_{\frac{1}{2}}(\Gr_{G,\rho(\omega_X),\ul{x}})^{N^\alpha,\chi} =
\Dmod_{\frac{1}{2}}(\Gr_{G,\rho(\omega_X),\ul{x}})^{\fL(N)_{\rho(\omega_X),\ul{x}},\chi} =\Whit^!(G)_{\ul{x}}.$$

Thus, we let $\on{coeff}_{G,\ul{x}}$ be the functor \eqref{e:Whit prel} for $N^\alpha$ large enough. 

\sssec{} \label{sss:coeff and Hecke}

By construction, the functor $\on{coeff}_{G,\ul{x}}$ is compatible with the action of $\Sph_{G,\ul{x}}$. 

\medskip

When working over the Ran space,
we consider the functor, to be denoted $\on{coeff}_{G,\Ran}$,
\begin{multline*}
\Dmod_{\frac{1}{2}}(\Bun_G\times \Ran) \simeq 
\Dmod_{\frac{1}{2}}(\Bun_G)\otimes \Dmod(\Ran) \overset{\on{coeff}_G\otimes \on{Id}}\to \\
\to \Whit^!(G)_\Ran\otimes \Dmod(\Ran) \overset{\sotimes}\to \Whit^!(G)_\Ran.
\end{multline*}

This functor is compatible with the natural action of $\Sph_{G,\Ran}$ on the two sides. 

\sssec{}

We let 
$$\on{coeff}_G:\Dmod_{\frac{1}{2}}(\Bun_G)\to \Whit^!(G)_\Ran$$
denote the composition
$$\Dmod_{\frac{1}{2}}(\Bun_G)\overset{\on{Id}\otimes \omega_\Ran}\longrightarrow 
\Dmod_{\frac{1}{2}}(\Bun_G\times \Ran)\overset{\on{coeff}_{G,\Ran}}\longrightarrow \Whit^!(G)_\Ran.$$

\sssec{} \label{sss:coeff unital}

The functors $\on{coeff}_{G,\ul{x}}$ have the following property: 

\medskip

For $\ul{x}\subseteq \ul{x}'$ consider the natural embedding
$$\on{incl}_{\ul{x}\subseteq \ul{x}'}:\Gr_{G,\rho(\omega_X),\ul{x}}\hookrightarrow \Gr_{G,\rho(\omega_X),\ul{x}'}.$$

The functor $(\on{incl}_{\ul{x}\subseteq \ul{x}'})^!$ maps $\Whit(G)_{\ul{x}'}\to \Whit(G)_{\ul{x}}$, and we have 
\begin{equation} \label{e:coeff unital}
\on{coeff}_{G,\ul{x}}\simeq \on{incl}_{\ul{x}\subseteq \ul{x}'}^!\circ \on{coeff}_{G,\ul{x}'}.
\end{equation}

The isomorphism \eqref{e:coeff unital} expresses the \emph{unital} structure on the functor $\on{coeff}_G$,
to be discussed in \secref{s:unitality}. 


\sssec{} \label{sss:vac coeff}

Let 
$$\on{coeff}_G^{\on{Vac}}:\Dmod_{\frac{1}{2}}(\Bun_G)\to \Vect$$ 
denote the composition of $\on{coeff}_{G,\ul{x}}$ with the functor
$$\Whit^!(G)_{\ul{x}}\hookrightarrow \Dmod_{\frac{1}{2}}(\Gr_{G,\rho(\omega_X),\ul{x}})\to \Vect,$$
where the second arrow is the functor of !-fiber at the unit point.

\medskip

By \eqref{e:coeff unital}, the above definition of $\on{coeff}_G^{\on{Vac}}$ is canonically independent of
the choice of $\ul{x}$.

\medskip

Equivalently, $\on{coeff}_G^{\on{Vac}}$ is the unique functor $\Dmod_{\frac{1}{2}}(\Bun_G) \to \Vect$ so
that the diagram commutes
$$
\CD
\Whit^!(G)_{\Ran} @>>> \Dmod_{\frac{1}{2}}(\Gr_{G,\rho(\omega_X),\Ran}) @>{(\one_{\Gr_{G,\rho(\omega_X),\Ran}})^!}>> \Dmod(\Ran) \\
@A{\on{coeff}_G}AA & & @AAA \\
\Dmod_{\frac{1}{2}}(\Bun_G)  & @>{\on{coeff}_G^{\on{Vac}}}>> & \Vect.
\endCD
$$
(In the above diagram the right vertical arrow is the !-pullback along $\Ran\to \on{pt}$, which is fully faithful by
the contractibility of the Ran space.)

\sssec{}

As in Remark \ref{r:indep char}, both the category $\Whit^!(G)_{\Ran}$ and the functor $\on{coeff}_G$ are canonically independent
of the choice of a non-degenerate character $\chi_0:N\to \BG_a$.

\ssec{The !-Poincar\'e functor}

\sssec{}

We start again by working with a fixed $\ul{x}\in \Ran$. We claim:

\begin{prop} \label{p:Poinc exist}
The functor $\on{coeff}_{G,\ul{x}}$ admits a left adjoint, to be denoted $\on{Poinc}_{G,!,\ul{x}}$.
\end{prop}

\begin{rem}

In fact, as we work over a fixed point $\ul{x}\in \Ran$, all objects in 
$\Whit^!(G)_{\ul{x}}$ are ind-holonomic, which implies the assertion of the proposition.
Below we give a different argument, which works also when $\ul{x}$ is allowed to move
in families over $\Ran$, see \secref{sss:Poinc ! Ran}. 
\end{rem} 

\begin{proof}[Proof of \propref{p:Poinc exist}]

Consider the  \emph{partially defined}\footnote{The issue here is that the ``lower-!" functors are not necessarily defined on non-holonomic objects.} 
functor
$$(\pi_{\ul{x}})_!:\Dmod_{\frac{1}{2}}(\Gr_{G,\rho(\omega_X),\ul{x}})\to \Dmod_{\frac{1}{2}}(\Bun_G),$$
left adjoint to $\pi^!_{\ul{x}}$.

\medskip

The assertion of the proposition is equivalent to the fact that $(\pi_{\ul{x}})_!$ \emph{is} defined on 
$\Whit^!(G)_{\Ran}\hookrightarrow \Dmod_{\frac{1}{2}}(\Gr_{G,\rho(\omega_X),\ul{x}})$. 

\medskip

First, it is easy to see that if $(\pi_{\ul{x}})_!$ is defined on some object $\CF\in \Dmod_{\frac{1}{2}}(\Gr_{G,\rho(\omega_X),\ul{x}})$
and $\CS$ is an object of $\Sph_{G,\ul{x}}$, then $(\pi_{\ul{x}})_!$ is defined on $\CS\star \CF$, and in fact
$$(\pi_{\ul{x}})_!(\CS\star \CF)\simeq \CS\star (\pi_{\ul{x}})_!(\CF).$$

This follows from the properness of the convolution diagram that defines the $\Sph_{G,\ul{x}}$-action on 
$\Dmod_{\frac{1}{2}}(\Gr_{G,\rho(\omega_X),\ul{x}})$ and $\Dmod_{\frac{1}{2}}(\Bun_G)$.

\medskip

Next, we observe that $(\pi_{\ul{x}})_!$ is defined on the vacuum object
$$\on{Vac}_{\Whit^!(G),\ul{x}}\in \Whit^!(G)_{\ul{x}}\subset \Dmod_{\frac{1}{2}}(\Gr_{G,\rho(\omega_X),\ul{x}}).$$
Indeed, this follows from the fact that $\on{Vac}_{\Whit^!(G),\ul{x}}$ is ind-holonomic, and the !-direct image 
functor is defined on holonomic D-modules.

\medskip

Finally, we claim that any object of $\Whit^!(G)_{\ul{x}}$ can be obtained as a convolution of an object of $\Sph_{G,\ul{x}}$
with $\on{Vac}_{\Whit^!(G),\ul{x}}$. In fact, by Remark \ref{r:CS again}, the functor
$$\Rep(\cG)_{\ul{x}}\overset{\Sat^{\on{nv}}}\longrightarrow \Sph_{\ul{x}} \overset{-\star \on{Vac}_{\Whit^!(G),\ul{x}}}\longrightarrow  \Whit^!(G)_{\ul{x}}$$
is an equivalence.

\end{proof} 

\sssec{}

The above proof shows that the functor $\on{Poinc}_{G,!,\ul{x}}$ is also compatible 
with the action of $\Sph_{G,\ul{x}}$. 

\medskip

Note, however, that this also follows a priori from the compatibility of 
$\on{coeff}_{G,x}$ with $\Sph_{G,x}$-actions and the observation that $\Sph_{G,x}$ is rigid as a 
monoidal category.\footnote{In fact, this was implicitly used in the proof of \propref{p:Poinc exist}: the properness of the convolution diagram
is the reason for the rigidity of $\Sph_{G,x}$.} 

\sssec{}

For a pair of points $\ul{x},\ul{x}'$ of $\Ran$ with $\ul{x}\subseteq \ul{x}'$, let
$$\on{ins.vac}_{\ul{x}\subseteq \ul{x}'}:\Whit(G)_{\ul{x}}\to \Whit(G)_{\ul{x}'}$$
be the functor left adjoint to 
$$(\on{incl}_{\ul{x}\subseteq \ul{x}'})^!:\Whit(G)_{\ul{x}'}\to \Whit(G)_{\ul{x}}.$$

Explicitly, if $\ul{x}'=\ul{x}\sqcup \ul{x}''$, so that
$$\Whit(G)_{\ul{x}'}\simeq \Whit(G)_{\ul{x}}\otimes \Whit(G)_{\ul{x''}},$$
then 
$$\on{ins.vac}_{\ul{x}\subseteq \ul{x}'}\simeq \on{Id}\otimes \on{Vac}_{\Whit^!(G),\ul{x}''}.$$

By adjunction, from \eqref{e:coeff unital}, we obtain:
\begin{equation} \label{e:Poinc unital x}
\on{Poinc}_{G,!,\ul{x}'}\circ \on{ins.vac}_{\ul{x}\subseteq \ul{x}'} \simeq \on{Poinc}_{G,!,\ul{x}}.
\end{equation} 

In \secref{sss:unitality fam} we will formulate a version of this isomorphism when the points $\ul{x}$
and $\ul{x}'$ move in families over the Ran space. 

\sssec{} \label{sss:Poinc ! Ran}

By the same token, the functor 
$$\on{coeff}_{G,\Ran}: \Dmod_{\frac{1}{2}}(\Bun_G\times \Ran)\to \Whit^!(G)_\Ran$$
admits a left adjoint, to be denoted 
$$\on{Poinc}_{G,!,\Ran}:\Whit^!(G)_\Ran\to \Dmod_{\frac{1}{2}}(\Bun_G\times \Ran).$$

\medskip

The functor $\on{Poinc}_{G,!,\Ran}$ is compatible with the actions of $\Sph_{G,\Ran}$ on the two sides. 

\sssec{}

The functor 
$$\on{coeff}_{G}: \Dmod_{\frac{1}{2}}(\Bun_G)\to \Whit^!(G)_\Ran$$
admits a left adjoint, to be denoted 
$$\on{Poinc}_{G,!}:\Whit^!(G)_\Ran\to  \Dmod_{\frac{1}{2}}(\Bun_G),$$ 
and given by
restricting the partially defined functor
$$\Dmod_{\frac{1}{2}}(\Gr_{G,\rho(\omega_X),\Ran})\overset{(\pi_{\ul{x}})_!}\to \Dmod_{\frac{1}{2}}(\Bun_G)$$
to 
$$\Whit^!(G)_\Ran\hookrightarrow \Dmod_{\frac{1}{2}}(\Gr_{G,\rho(\omega_X),\Ran}).$$

Explicitly, the $\on{Poinc}_{G,!}$ identifies with the composition
$$\Whit^!(G)_\Ran \overset{\on{Poinc}_{G,!,\Ran}}\longrightarrow \Dmod_{\frac{1}{2}}(\Bun_G\times \Ran)\to
\Dmod_{\frac{1}{2}}(\Bun_G),$$
where the second arrow is the functor of !-direct image.

\medskip

The functor $\on{Poinc}_{G,!,\ul{x}}$ is obtained from $\on{Poinc}_{G,!}$ by restriction along \eqref{e:Gr x to Gr Ran}. 

\sssec{}

It follows formally from \secref{sss:vac coeff} that the object
$$\on{Poinc}_{G,!,\ul{x}}(\one_{\Whit^!(G),\ul{x}})\in  \Dmod_{\frac{1}{2}}(\Bun_G)$$
is canonically independent of the choice of $\ul{x}$.

\medskip

We will denote it by
$$\on{Poinc}^{\on{Vac}}_{G,!}\in \Dmod_{\frac{1}{2}}(\Bun_G).$$

\medskip

We also have
$$\on{Poinc}^{\on{Vac}}_{G,!} \simeq \on{Poinc}_{G,!}(\on{Vac}_{\Whit^!(G),\Ran}),$$
where
$$\on{Vac}_{\Whit^!(G),\Ran}\in \Whit^!(G)_{\Ran}$$
is the factorization unit spread over the Ran space. 

\sssec{} \label{sss:recover Poinc from Vac}

As we saw in the proof of \propref{p:Poinc exist}, we can recover the functor $\on{Poinc}_{G,!,\ul{x}}$ from the object 
$\on{Poinc}^{\on{Vac}}_{G,!}$ using the Hecke action. 

\medskip

By adjunction, we obtain that the functor $\on{coeff}_{G,\ul{x}}$ can be uniquely recovered from the knowledge of 
$\on{coeff}^{\on{Vac}}_G$ and the action of $\Rep(G)_{\ul{x}}$ on $\Dmod_{\frac{1}{2}}(\Bun_G)$ via
$$\Rep(G)_{\ul{x}}\overset{\Sat_G^{\on{nv}}}\longrightarrow \Sph_{G,\ul{x}}.$$

\medskip

The same applies to the functors $(\on{Poinc}^{\on{Vac}}_{G,!,\Ran},\on{coeff}^{\on{Vac}}_{G,\Ran})$.

\ssec{The *-Poincar\'e functor}

\sssec{}

Recall that along with the category $\Dmod(\Bun_G)$, one can consider its version $\Dmod_{\on{co}}(\Bun_G)$,
and similarly for gerbe-twisted versions $\Dmod_\CG(\Bun_G)$. 

\medskip

In the untwisted case, we have the identification 
$$(\Dmod(\Bun_G))^\vee \simeq \Dmod_{\on{co}}(\Bun_G).$$

In the twisted case, this becomes
\begin{equation} \label{e:dual of half BunG}
(\Dmod_\CG(\Bun_G))^\vee \simeq \Dmod_{\CG^{\otimes -1},\on{co}}(\Bun_G).
\end{equation}

For $\CG=\det_{\Bun_G}^{\frac{1}{2}}$, the identification \eqref{e:dual of half BunG} becomes a self-duality 
\begin{equation} \label{e:dual of half twisted BunG}
(\Dmod_{\frac{1}{2}}(\Bun_G))^\vee \simeq \Dmod_{\frac{1}{2},\on{co}}(\Bun_G).
\end{equation}

\sssec{}

Let
$$\on{Poinc}_{G,*}:\Whit_*(G)_{\Ran}\to \Dmod_{\frac{1}{2},\on{co}}(\Bun_G)$$
be the functor \emph{dual} to $\on{coeff}_G$. 

\medskip

Let 
$$\on{Poinc}_{G,*,\ul{x}}:\Whit_*(G)_{\ul{x}}\to \Dmod_{\frac{1}{2},\on{co}}(\Bun_G)$$
be the functor dual to $\on{coeff}_{G,\ul{x}}$. It is easy to see that the functor
$\on{Poinc}_{G,*,\ul{x}}$ is obtained from $\on{Poinc}_{G,*}$ by restriction along 
\eqref{e:Gr x to Gr Ran}.

\medskip

The functor $\on{Poinc}_{G,*,\ul{x}}$ is also compatible with the action of $\Sph_{G,\ul{x}}$. 

\sssec{}

Let 
$$\on{Poinc}^{\on{Vac}}_{G,*} \in \Dmod_{\frac{1}{2},\on{co}}(\Bun_G)$$
be the image under $\on{Poinc}_{G,*,\ul{x}}$ of the factorization unit
$$\one_{\Whit_*(G),\ul{x}}\in \Whit_*(G)_{\ul{x}}$$
at some/any $\ul{x}\in \Ran$ or, equivalently, of
$$\one_{\Whit_*(G),\Ran}\in \Whit_*(G)_{\Ran}$$
under $\on{Poinc}_{G,*}$.

\medskip

By definition, the pairing with $\on{Poinc}^{\on{Vac}}_{G,*}$, viewed as a functor
$$\Dmod_{\frac{1}{2}}(\Bun_G)\to \Vect,$$
is the functor $\on{coeff}_G^{\on{Vac}}$. 

\sssec{}

The functor $\on{Poinc}_{G,*,\ul{x}}$ can be explicitly described as follows. For $N^\alpha$ as in 
\secref{sss:N alpha}, consider the composition
$$\Dmod_{\frac{1}{2}}(\Gr_{G,\rho(\omega_X),\ul{x}}) \overset{\Av^{(N^\alpha,\chi)}_*}\longrightarrow
\Dmod_{\frac{1}{2}}(\Gr_{G,\rho(\omega_X),\ul{x}}) \overset{(\pi_{\ul{x}})_*}\to \Dmod_{\frac{1}{2},\on{co}}(\Bun_G).$$

For $N^{\alpha}\subset N^{\alpha'}$, we have a canonically defined natural transformation
\begin{equation} \label{e:passage Whit av *}
(\pi_{\ul{x}})_*\circ \Av^{(N^{\alpha'},\chi)}_*\to (\pi_{\ul{x}})_*\circ \Av^{(N^{\alpha},\chi)}_*,
\end{equation} 
and it follows from \lemref{l:passage Whit av} that the maps \eqref{e:passage Whit av *}
are isomorphisms for $N^\alpha$ large enough.

\medskip

It follows formally that the resulting functor $(\pi_{\ul{x}})_*\circ \Av^{(N^{\alpha},\chi)}_*$ 
$$\Dmod_{\frac{1}{2}}(\Gr_{G,\rho(\omega_X),\ul{x}}) \to \Dmod_{\frac{1}{2},\on{co}}(\Bun_G),$$
for some/any $\alpha$ that is large enough, factors via the projection
$$\Dmod_{\frac{1}{2}}(\Gr_{G,\rho(\omega_X),\ul{x}}) \to 
\left(\Dmod_{\frac{1}{2}}(\Gr_{G,\rho(\omega_X),\ul{x}})\right)_{\fL(N)_{\rho(\omega_X),\ul{x}},\chi}
\to\Dmod_{\frac{1}{2},\on{co}}(\Bun_G).$$

The resulting functor
$$\Whit_*(G)_{\ul{x}}:=
\left(\Dmod_{\frac{1}{2}}(\Gr_{G,\rho(\omega_X),\ul{x}})\right)_{\fL(N)_{\rho(\omega_X),\ul{x}},\chi}
\to\Dmod_{\frac{1}{2},\on{co}}(\Bun_G)$$
is the functor $\on{Poinc}_{G,*,\ul{x}}$. 

\ssec{Coefficient and Poincar\'e functors: global interpretation} \label{ss:glob Whit}

\sssec{}

Consider the stack $\Bun_{N,\rho(\omega_X)}$ and the map 
$$\fp:\Bun_{N,\rho(\omega_X)}\to \Bun_G.$$

Recall the functor
\begin{equation} \label{e:pullback to BunN half}
\fp^!:\Dmod_{\frac{1}{2}}(\Bun_G)\to \Dmod(\Bun_{N,\rho(\omega_X)}),
\end{equation}
see \secref{sss:pull back to 1/2 rho omega}. 

\sssec{}

The character $\chi$ has a global counterpart, which is a map
$$\chi^{\on{glob}}:\Bun_{N,\rho(\omega_X)}\to \BG_a.$$

Namely, it is the composition
$$\Bun_{N,\rho(\omega_X)} \to \Bun_{N/[N,N],\rho(\omega_X)} \simeq 
\underset{i}\Pi\, \Bun_{(\BG_a)_{\omega_X}} \to \underset{i}\Pi\, H^1(X,\omega_X)\simeq 
\underset{i}\Pi\, \BG_a \overset{\chi_0}\to \BG_a,$$
where:

\begin{itemize}

\item $(\BG_a)_{\omega_X}$ is the twist of the constant group-scheme with fiber $\BG_a$ using the $\BG_m$-action
on $\BG_a$ and the line bundle $\omega_X$, viewed as a $\BG_m$-torsor\footnote{I.e., it is the total space of $\omega_X$ as
a line bundle.};

\smallskip

\item $\Bun_{(\BG_a)_{\omega_X}} \to H^1(X,\omega_X)$ is the map that records the class of a torsor.

\end{itemize}

\sssec{}

We let
$$\on{coeff}_G^{\on{Vac,glob}}:\Dmod_{\frac{1}{2}}(\Bun_G)\to \Vect$$
denote the functor 
$$\Dmod_{\frac{1}{2}}(\Bun_G) \overset{\fp^!}\to 
\Dmod(\Bun_{N,\rho(\omega_X)}) \overset{-\overset{*}\otimes (\chi^{\on{glob}})^*(\on{exp})}\longrightarrow 
\Dmod(\Bun_{N,\rho(\omega_X)}) \overset{\on{C}^\cdot_\dr(\Bun_{N,\rho(\omega_X)},-)}\longrightarrow \Vect.$$

\sssec{} \label{sss:Poinc Vac glob}

Let
$$\on{Poinc}^{\on{Vac,glob}}_{G,!} \in \Dmod_{\frac{1}{2}}(\Bun_G)$$
be the object
$$\fp_!\circ (\chi^{\on{glob}})^*(\on{exp}).$$
I.e., the functor 
$$\Vect \to \Dmod_{\frac{1}{2}}(\Bun_G), \quad k\mapsto \on{Poinc}^{\on{Vac,glob}}_{G,!}$$
is the left adjoint of $\on{coeff}_G^{\on{Vac,glob}}$. 

\medskip

Let 
$$\on{Poinc}^{\on{Vac,glob}}_{G,*} \in \Dmod_{\frac{1}{2},\on{co}}(\Bun_G)$$
be the object
$$\fp_*\circ (\chi^{\on{glob}})^*(\on{exp}).$$

\sssec{} \label{sss:delta N}

Denote
$$\delta_{N_{\rho(\omega_X)}}:=\dim(\Bun_{N,\rho(\omega_X)}).$$

\medskip

We have
$$\BD^{\on{Verdier}}(\on{Poinc}^{\on{Vac,glob}}_{G,!})=\on{Poinc}^{\on{Vac,glob}}_{G,*}[2\delta_{N_{\rho(\omega_X)}}],$$
where $\BD^{\on{Verdier}}$ is the usual Verdier dualization functor
$$(\Dmod_{\frac{1}{2}}(\Bun_G)^c)^{\on{op}}\to (\Dmod_{\frac{1}{2},\on{co}}(\Bun_G))^c.$$

In other words, the object $\on{Poinc}^{\on{Vac,glob}}_{G,*}[2\delta_{N_{\rho(\omega_X)}}]$, viewed as a functor
$$\Vect\to \Dmod_{\frac{1}{2},\on{co}}(\Bun_G),$$
is the dual of $\on{coeff}_G^{\on{Vac,glob}}$. 

\sssec{}

We claim:

\begin{lem} \label{l:Vac coeff loc vs glob} \hfill

\smallskip

\noindent{\em(i)} $\on{coeff}_G^{\on{Vac,glob}}\simeq \on{coeff}_G^{\on{Vac}}[2\delta_{N_{\rho(\omega_X)}}]$;

\smallskip

\noindent{\em(ii)} $\on{Poinc}^{\on{Vac,glob}}_{G,!}\simeq \on{Poinc}^{\on{Vac}}_{G,!}[-2\delta_{N_{\rho(\omega_X)}}]$;

\smallskip

\noindent{\em(iii)} $\on{Poinc}^{\on{Vac,glob}}_{G,*}\simeq \on{Poinc}^{\on{Vac}}_{G,*}$.

\end{lem} 

\begin{proof}

The three statements are logically equivalent. We will prove point (iii). 

\medskip

Pick $\ul{x}\in \Ran$, and consider the map
$$N^\alpha/\fL^+(N)_{\rho(\omega_X),\ul{x}}
\hookrightarrow \fL(N)_{\rho(\omega_X),\ul{x}}/\fL^+(N)_{\rho(\omega_X),\ul{x}}
\overset{\text{act on the unit}}\to \Gr_{G,\rho(\omega_X),\ul{x}}\to \Bun_G$$
for $N^\alpha\supset \fL^+(N)_{\rho(\omega_X),\ul{x}}$ as in \secref{sss:N alpha}. 

\medskip

By definition, the object $\on{Poinc}^{\on{Vac,glob}}_{G,*}$ is the direct image along this map
of $(\chi|_{N^\alpha/\fL^+(N)_{\rho(\omega_X),\ul{x}}}^*(\on{exp})$, where:

\begin{itemize} 

\item By a slight abuse of notation, we regard $\chi$ as a map $\fL(N)_{\rho(\omega_X),\ul{x}}/\fL^+(N)_{\rho(\omega_X),\ul{x}}\to \BG_a$;

\item The index $\alpha$ is taken to be large enough so that \eqref{e:alpha large} holds.

\end{itemize}

Note that, however, the map 
$$\fL(N)_{\rho(\omega_X),\ul{x}}/\fL^+(N)_{\rho(\omega_X),\ul{x}}
\overset{\text{act on the unit}}\to \Gr_{G,\rho(\omega_X),\ul{x}}\overset{\pi}\to \Bun_G$$
factors as
$$\fL(N)_{\rho(\omega_X),\ul{x}}/\fL^+(N)_{\rho(\omega_X),\ul{x}}\to \Bun_{N,\rho(\omega_X)}\to \Bun_G.$$

Hence, it suffices to show that if \eqref{e:alpha large} holds, then 
\begin{equation} \label{e:chi local vs global}
(\pi_{N,\rho(\omega_X)}|_{N^\alpha/\fL^+(N)_{\rho(\omega_X),\ul{x}}})_*\circ
(\chi|_{N^\alpha/\fL^+(N)_{\rho(\omega_X),\ul{x}}}^*(\on{exp})\simeq (\chi^{\on{glob}})^*(\on{exp}).
\end{equation}

Note that the map 
$$\chi:\fL(N)_{\rho(\omega_X),\ul{x}}/\fL^+(N)_{\rho(\omega_X),\ul{x}}\to \BG_a$$
identifies with
$$\fL(N)_{\rho(\omega_X),\ul{x}}/\fL^+(N)_{\rho(\omega_X),\ul{x}}\overset{\pi_{N,\rho(\omega_X)}}\to \Bun_{N,\rho(\omega_X)} 
\overset{\chi^{\on{glob}}}\to \BG_a.$$

This implies \eqref{e:chi local vs global} by the projection formula, since if \eqref{e:alpha large} holds, the map
$$\pi_{N,\rho(\omega_X)}|_{N^\alpha/\fL^+(N)_{\rho(\omega_X),\ul{x}}}:
N^\alpha/\fL^+(N)_{\rho(\omega_X),\ul{x}}\to \Bun_{N,\rho(\omega_X)}$$
is smooth with homologically contractible fibers; in fact, fibers are isomorphic to 
$$N^\alpha\cap \on{Sect}(X-\ul{x},N_{\rho(\omega_X)}).$$

\end{proof}

\section{The localization functor} \label{s:loc}

A fundamental insight of Beilinson-Drinfeld in \cite{BD1} is that the localization functor
\begin{equation} \label{e:usual localization preamble}
\Loc_{G,\kappa}:\KL(G)_{\kappa,\Ran} \to \Dmod_\kappa(\Bun_G).
\end{equation}
may be used as a key local-to-global tool in geometric Langlands theory. 

\medskip

There are multiple (equivalent) ways to set this up, and
in this section we will describe one of them.\footnote{The reader who is willing to take the existence of
$\Loc_{G,\kappa}$ and its basic properties of faith may choose to skip directly to \secref{s:prop Loc}.}

\medskip

For the sake of completeness we will define $\Loc_{G,\kappa}$ for any level $\kappa$. We will
specialize to the critical value of $\kappa$ starting from \secref{s:coeff Loc}. 

\ssec{The de Rham twisting on \texorpdfstring{$\Bun_G$}{BunG} corresponding to a level}

In this subsection we will show how a level $\kappa$ gives rise to a de Rham twisting $\CT^{\on{glob}}_{\kappa}$ on
$\Bun_G$. 

\sssec{} \label{sss:mult line bundle local}

Let $\fL(G)^\wedge$ denote the formal completion of $\fL(G)$ along $\fL^+(G)$, viewed as a factorization
group ind-scheme. 

\medskip

Note that a level $\kappa$ may be thought of as a central extension of $\fL(G)^\wedge$ equipped with a splitting
over $\fL^+(G)$. Equivalently, we can think of it as a factorization line bundle $\CL^{\on{loc}}_\kappa$ on the groupoid
$$\on{pt}/\fL^+(G)\overset{\hl^{\on{loc},\wedge}}\leftarrow \on{Hecke}_G^{\on{loc},\wedge} \overset{\hr^{\on{loc},\wedge}}\to \on{pt}/\fL^+(G),$$
compatible with the groupoid structure, where:

\smallskip

\begin{itemize}

\item $\on{Hecke}_G^{\on{loc}}$ is the local Hecke stack, i.e., $\fL^+(G)\backslash \fL(G)/\fL^+(G)$, viewed
as a groupoid on $\on{pt}/\fL^+(G)$;

\medskip

\item $\on{Hecke}_G^{\on{loc},\wedge}$ is the formal completion of $\on{Hecke}_G^{\on{loc}}$ along the unit
section $\on{pt}/\fL^+(G)\to \on{Hecke}_G^{\on{loc}}$.

\end{itemize}

\sssec{}

Let 
$$\Bun_G\times \Ran \overset{\hl^{\on{glob}}}\leftarrow \on{Hecke}_{G,\Ran}^{\on{glob}} \overset{\hr^{\on{glob}}}\to \Bun_G\times \Ran$$
be the global Hecke groupoid, and let $\on{Hecke}_{G,\Ran}^{\on{glob},\wedge}$ denote its formal
completion along the unit section
$$\Bun_G\times \Ran\to \on{Hecke}_{G,\Ran}^{\on{glob}}.$$

Note that we have a map of groupoids
\begin{equation} \label{e:log glob Hecke diagram}
\CD
\Bun_G\times \Ran @<{\hl^{\on{glob}}}<< \on{Hecke}_{G,\Ran}^{\on{glob}} @>{\hr^{\on{glob}}}>> \Bun_G\times \Ran \\
@V{\on{ev}_\Ran}VV @V{\on{ev}_\Ran}VV @VV{\on{ev}_\Ran}V \\
(\on{pt}/\fL^+(G))_\Ran @<{\hl^{\on{loc}}}<< \on{Hecke}_{G,\Ran}^{\on{loc}}  @>{\hr^{\on{loc}}}>> (\on{pt}/\fL^+(G))_\Ran,
\endCD
\end{equation} 
in which both squares are Cartesian, where we denote by $\on{ev}_\Ran$ the ``global-to-local" map given by retsriction
to the parameterized multi-disc. 

\medskip

Taking the formal completion along the unit sections yields the diagram
\begin{equation} \label{e:log glob Hecke diagram compl}
\CD
\Bun_G\times \Ran @<{\hl^{\on{glob},\wedge}}<< \on{Hecke}_{G,\Ran}^{\on{glob},\wedge} @>{\hr^{\on{glob},\wedge}}>> \Bun_G\times \Ran \\
@V{\on{ev}_\Ran}VV @V{\on{ev}_\Ran}VV @VV{\on{ev}_\Ran}V \\
(\on{pt}/\fL^+(G))_\Ran @<{\hl^{\on{loc},\wedge}}<< \on{Hecke}_{G,\Ran}^{\on{loc},\wedge}  @>{\hr^{\on{loc},\wedge}}>> (\on{pt}/\fL^+(G))_\Ran. 
\endCD
\end{equation} 

\sssec{} \label{sss:mult line bundle global}

The pullback of the line bundle $\CL^{\on{loc}}_{\kappa,\Ran}$ on $\on{Hecke}_{G,\Ran}^{\on{loc},\wedge}$ along
$\on{ev}_\Ran$ gives rise to a line bundle, to be denoted $\CL^{\on{glob}}_{\kappa,\Ran}$, on 
$\on{Hecke}_{G,\Ran}^{\on{glob},\wedge}$ that is multiplicative with respect to  
the groupoid structure. 

\medskip

Consider the prestack quotient
$$(\Bun_G\times \Ran)/\on{Hecke}_{G,\Ran}^{\on{loc},\wedge}.$$

The datum of $\CL^{\on{glob}}_{\kappa,\Ran}$ is equivalent to that of a $\CO^\times$-gerbe
on $(\Bun_G\times \Ran)/\on{Hecke}_{G,\Ran}^{\on{loc},\wedge}$, to be denoted $\CG^{/\on{Hecke}}_{\kappa,\Ran}$, 
equipped with a trivialization of its pullback to $\Bun_G\times \Ran$. 

\sssec{}

Note now that we have a tautological map
$$(\Bun_G\times \Ran)/\on{Hecke}_{G,\Ran}^{\on{loc},\wedge}\to (\Bun_G)_\dr\times \Ran.$$

Consider the composite map
\begin{equation} \label{e:inf Hecke to inf}
(\Bun_G\times \Ran)/\on{Hecke}_{G,\Ran}^{\on{loc},\wedge}\to (\Bun_G)_\dr
\end{equation} 

We have the following fundamental assertion (see \cite[Theorem 4.5.3]{Ro2}):

\begin{thm} \label{t:Nick 1}
The functor of pullback along \eqref{e:inf Hecke to inf} on $\IndCoh(-)$ 
is fully faithful.
\end{thm} 

\begin{cor} \label{c:Nick 1}
The functor of pullback along \eqref{e:inf Hecke to inf} is fully faithful on\footnote{For a prestack $\CY$, we denote by $\on{Perf}(\CY)$
the category of dualizable objects in $\QCoh(\CY)$, i.e., the objects whose pullback to any affine scheme is perfect.} 
$\on{Perf}(-)$.
\end{cor}

\sssec{}

The construction of the twising $\CT_\kappa$ on $\Bun_G$ corresponding to $\kappa$ is provided by the
following assertion:

\begin{cor}  \label{c:existence of the twising}
There exists a uniquely defined de Rham twisting $\CT^{\on{glob}}_\kappa$ on $\Bun_G$, such that:

\begin{itemize}

\item The pullback of the underlying $\CO^\times$-gerbe $\CG^{\on{glob}}_\kappa$ on $(\Bun_G)_\dr$ along 
\eqref{e:inf Hecke to inf} identifies with $\CG^{/\on{Hecke}}_{\kappa,\Ran}$;

\item The trivialization of $\CG^{\on{glob}}_\kappa|_{\Bun_G}$ is compatible with the trivialization
of $\CG^{/\on{Hecke}}_{\kappa,\Ran}|_{\Bun_G\times \Ran}$.

\end{itemize}

\end{cor}

\begin{proof}

According to \cite[Sect. 6.3]{GaRo2}, using the exponential isomorphism
$$(\BG_a)^\wedge\overset{\on{exp}}\to (\BG_m)^\wedge,$$
we can think of a twisting on a prestack $\CY$ as a point in
$\Maps(\CO_{\CY_\dr},\CO_{\CY_\dr}[2])$ equipped with a trivialization 
of its pullback to $\CY$.

\medskip

By the same logic, we can think of $\CG^{/\on{Hecke}}_{\kappa,\Ran}$ as a point
of 
$$\Maps(\CO_{(\Bun_G\times \Ran)/\on{Hecke}_{G,\Ran}^{\on{loc},\wedge}},
\CO_{(\Bun_G\times \Ran)/\on{Hecke}_{G,\Ran}^{\on{loc},\wedge}}[2])$$
equipped with a trivialization of its pullback to $\Bun_G\times \Ran$.

\medskip

The assertion of the corollary follows now from \corref{c:Nick 1},
combined with the fact that the functor
$$\Vect\to \QCoh(\Ran), \quad k\mapsto \CO_\Ran$$
is fully faithful. 

\end{proof}

\sssec{}

In what follows we will denote the category of $\CT_\kappa^{\on{glob}}$-twisted D-modules on $\Bun_G$ by
$$\Dmod_\kappa(\Bun_G).$$

\sssec{}

Take $\kappa={2\cdot\crit}$. We claim:

\begin{prop} \label{p:Tate}
The resulting twisting $\CT_{2\cdot \crit}$ on $\Bun_G$ identifies canonically with 
$\on{dlog}(\det_{\Bun_G})$.
\end{prop} 

\begin{rem}

This proposition implies that our notations for $\Dmod_\crit(\Bun_G)$ (see \secref{sss:crit on Bun_G}) are consistent. 

\end{rem}

\begin{proof} 

Unwinding the construction, we need to show that the multiplicative line bundle $\CL^{\on{glob}}_{2\cdot \crit,\Ran}$
on $\on{Hecke}_{G,\Ran}^{\on{loc},\wedge}$ identifies with the restriction of the multiplicative line bundle
on $\on{Hecke}_{G,\Ran}^{\on{loc}}$ given by
$$\hr^*(\det_{\Bun_G})\otimes \hl^*(\det_{\Bun_G})^{\otimes -1}.$$

Recall that that the multiplicative line bundle $\CL^{\on{loc}}_{2\cdot \crit}$ on $\on{Hecke}_G^{\on{loc},\wedge}$
is itself obtained as the restriction of the \emph{inverse} of the Tate line bundle $\CL^{\on{loc}}_{\on{Tate}}$
on $\on{Hecke}_G^{\on{loc}}$, constructed as follows:

\medskip

The line bundle $\CL^{\on{loc}}_{\on{Tate}}$ associates to a pair of $G$-bundles $\CP^1_G$ and $\CP^2_G$ on
$\cD_{\ul{x}}$ equipped with an identification 
$$\CP^1_G|_{\cD^\times_{\ul{x}}}\simeq \CP^2_G|_{\cD^\times_{\ul{x}}}$$
the relative determinant of the two lattices
$$\Gamma(\cD_{\ul{x}},\fg_{\CP^1_G}) \subset \Gamma(\cD^\times_{\ul{x}},\fg_{\CP^1_G}) =
\Gamma(\cD^\times_{\ul{x}},\fg_{\CP^2_G}) \supset \Gamma(\cD_{\ul{x}},\fg_{\CP^2_G}),$$
i.e.,
\begin{equation} \label{e:det loc}
\det(\Gamma(\cD_{\ul{x}},\fg_{\CP^1_G})/\bL)\otimes \det(\Gamma(\cD_{\ul{x}},\fg_{\CP^2_G})/\bL)^{-1}
\end{equation}
for some/any lattice $\bL$ contained in both. 

\medskip

Given a pair of $G$-bundles $\CP^1_G$ and $\CP^2_G$ on $X$ equipped with an identification 
$$\CP^1_G|_{X-\ul{x}}\simeq \CP^2_G|_{X-\ul{x}},$$
the fiber of $\hr^*(\det_{\Bun_G})\otimes \hl^*(\det_{\Bun_G})^{\otimes -1}$ at the corresponding point of 
$\on{Hecke}_{G,\Ran}^{\on{loc}}$ is given by
\begin{equation} \label{e:det glob}
\det(\Gamma(X,\fg_{\CP^2_G}))\otimes \det(\Gamma(X,\fg_{\CP^1_G}))^{\otimes -1}.
\end{equation}

We claim that the lines \eqref{e:det loc} and \eqref{e:det glob} are indeed canonically inverse to one another.

\medskip

Indeed, we can rewrite \eqref{e:det loc} as
$$\det\left(\on{Fib}\left(\Gamma(\cD_{\ul{x}},\fg_{\CP^1_G})\oplus \fg^{\on{out}}\to \fg^{\on{mer}}\right)\right)
\otimes \det\left(\on{Fib}\left(\Gamma(\cD_{\ul{x}},\fg_{\CP^2_G})\oplus \fg^{\on{out}}\to \fg^{\on{mer}}\right)\right)^{\otimes -1},$$
where:

\begin{itemize}

\item $\Gamma(\cD^\times_{\ul{x}},\fg_{\CP^1_G})=:\fg^\mer:=\Gamma(\cD^\times_{\ul{x}},\fg_{\CP^2_G})$;

\item $\Gamma(X-\ul{x},\fg_{\CP^1_G})=:\fg^{\on{out}}:=\Gamma(X-\ul{x},\fg_{\CP^2_G})$.

\end{itemize}

However,
$$\on{Fib}\left(\Gamma(\cD_{\ul{x}},\fg_{\CP^i_G})\oplus \fg^{\on{out}}\to \fg^{\on{mer}}\right)\simeq
\Gamma(X,\fg_{\CP^i_G}).$$

\end{proof} 

\ssec{The functor \texorpdfstring{$\Gamma_G$}{Gamma}}  \label{ss:functor Gamma}

Our approach to the construction of the localization functor is by defining it as the left adjoint 
of the functor $\Gamma_G$ of \emph{global sections} (not quite literally, though, see \secref{ss:defn Loc}).

\medskip

In this subsection we introduce the functor $\Gamma_G$. 

\medskip

To simplify the notation, for most of this subsection we fix a point $\ul{x}\in \Ran$.

\sssec{}

For an integer $n$, consider the stack $\Bun_G^{\on{level}_{n\cdot \ul{x}}}$
of $G$-bundles with structure of level $n$ at $\ul{x}$. 

\medskip

Consider the corresponding category 
$$\Dmod_{\kappa,\on{co}}(\Bun^{\on{level}_{n\cdot \ul{x}}}_G),$$
see \cite[Sect. 1.4.25]{Ra5}.

\medskip

We endow it with the forgetful functor
\begin{equation} \label{e:ren oblv l}
\oblv^{l,\on{ren}}_{\kappa}:\Dmod_{\kappa,\on{co}}(\Bun^{\on{level}_{n\cdot \ul{x}}}_G)\to
\QCoh_{\on{co}}(\Bun^{\on{level}_{n\cdot \ul{x}}}_G),
\end{equation} 
which is the composition with the usual \emph{left} forgetful functor
$$\oblv^l_{\kappa}:\Dmod_{\kappa,\on{co}}(\Bun^{\on{level}_{n\cdot \ul{x}}}_G)\to
\QCoh_{\on{co}}(\Bun^{\on{level}_{n\cdot \ul{x}}}_G)$$
(see \corref{c:QCoh co alg stacks}), followed by the cohomological shift $[2n\dim(\fg)]$. 

\medskip

Note that for $n_1\leq n_2$, we have a commutative diagram
$$
\CD
\Dmod_{\kappa,\on{co}}(\Bun^{\on{level}_{n_2\cdot \ul{x}}}_G) @>{\oblv^{l,\on{ren}}_{\kappa}}>> \QCoh_{\on{co}}(\Bun^{\on{level}_{n_2\cdot \ul{x}}}_G) \\
@AAA @AAA \\
\Dmod_{\kappa,\on{co}}(\Bun^{\on{level}_{n_1\cdot \ul{x}}}_G) @>\oblv^{l,\on{ren}}_{\kappa}>> \QCoh_{\on{co}}(\Bun^{\on{level}_{n_1\cdot \ul{x}}}_G),
\endCD$$
where:

\begin{itemize}

\item The left vertical arrow is the functor of *-pullback on twisted D-modules;

\item The right vertical arrow is the functor of *-pullback on $\QCoh_{\on{co}}$.

\end{itemize}

Note that by definition, for $n=0$, we have $\oblv^{l,\on{ren}}_{\kappa}=\oblv^l_{\kappa}$. 

\sssec{}

Consider the stack (in fact, a scheme) 
$$\Bun_G^{\on{level}_{\ul{x}}}=\underset{n}{\on{lim}}\, \Bun_G^{\on{level}_{n\cdot \ul{x}}}$$
of $G$-bundles with full level structure at $\ul{x}$. 

\medskip

Define
$$\Dmod_{\kappa,\on{co}}(\Bun^{\on{level}_{\ul{x}}}_G):=
\underset{n}{\on{colim}}\, \Dmod_{\kappa,\on{co}}(\Bun_G^{\on{level}_{n\cdot \ul{x}}}),$$
where the colimit is formed using the *-pullback functors
$$\Dmod_\kappa(\Bun_G^{\on{level}_{n_1\cdot \ul{x}}})\to \Dmod_\kappa(\Bun_G^{\on{level}_{n_2\cdot \ul{x}}}), \quad n_2\geq n_1.$$

\medskip

The functors \eqref{e:ren oblv l} combine to a functor
$$\oblv^{l,\on{ren}}_{\kappa}:\Dmod_{\kappa,\on{co}}(\Bun^{\on{level}_{\ul{x}}}_G)\to 
\QCoh_{\on{co}}(\Bun_G^{\on{level}_{\ul{x}}}).$$

\medskip

Consider the composite functor, to be denoted $\Gamma^{\on{ren}}_\kappa(\Bun^{\on{level}_{\ul{x}}}_G,-)$,
\begin{equation} \label{e:global sections level}
\Dmod_{\kappa,\on{co}}(\Bun^{\on{level}_{\ul{x}}}_G)\overset{\oblv^{l,\on{ren}}_{\kappa}}\longrightarrow 
\QCoh_{\on{co}}(\Bun_G^{\on{level}_{\ul{x}}}) \overset{\Gamma(\Bun_G^{\on{level}_{\ul{x}}},-)}\longrightarrow \Vect;
\end{equation} 

\sssec{}

According to \cite[Sect. 1.4.25]{Ra5}, the category $\Dmod_{\kappa,\on{co}}(\Bun^{\on{level}_{\ul{x}}}_G)$ 
carries a \emph{strong} action of $\fL(G)_{\ul{x}}$ at level $\kappa$. Furthermore, the functor
$\Gamma^{\on{ren}}_\kappa(\Bun^{\on{level}_{\ul{x}}}_G,-)$ of \eqref{e:global sections level} is
\emph{weakly} $\fL(G)_{\ul{x}}$-equivariant.

\medskip

Hence, by the \emph{universal property of} $\hg\mod_{\kappa,\ul{x}}$ (see \cite[Sect. 1.4.25]{Ra5}), 
the functor $\Gamma^{\on{ren}}_\kappa(\Bun^{\on{level}_{\ul{x}}}_G,-)$ upgrades to a functor
$$\Gamma^{\on{ren}}_\kappa(\Bun^{\on{level}_{\ul{x}}}_G,-)^{\on{enh}}:
\Dmod_{\kappa,\on{co}}(\Bun^{\on{level}_{\ul{x}}}_G)\to \hg\mod_{\kappa,\ul{x}},$$
\emph{strongly} compatible with the $\fL(G)_{\ul{x}}$-actions.

\sssec{}

In particular, the functor $\Gamma^{\on{ren}}_\kappa(\Bun^{\on{level}_{\ul{x}}}_G,-)^{\on{enh}}$ gives rise
to a functor, to be denoted $\Gamma_{G,\kappa,\ul{x}}$:
$$\Dmod_{\kappa,\on{co}}(\Bun_G)\simeq \left(\Dmod_{\kappa,\on{co}}(\Bun^{\on{level}_{\ul{x}}}_G)\right)^{\fL^+(G)_{\ul{x}}}\to
(\hg\mod_{\kappa,\ul{x}})^{\fL^+(G)_{\ul{x}}}=\KL(G)_{\kappa,\ul{x}}.$$

\sssec{}

By a similar token, letting $\ul{x}$ vary over the Ran space, we obtain a functor
$$\Gamma_{G,\kappa}:\Dmod_{\kappa,\on{co}}(\Bun_G)\to \KL(G)_{\kappa,\Ran}.$$

Furthermore, we can consider a $\Dmod(\Ran)$-linear functor 
$$\Gamma_{G,\kappa,\Ran}:\Dmod_{\kappa,\on{co}}(\Bun_G)\otimes \Dmod(\Ran)\to \KL(G)_{\kappa,\Ran},$$
so that $\Gamma_{G,\kappa}$ is the composition
$$\Dmod_{\kappa,\on{co}}(\Bun_G) \overset{\on{Id}\otimes \omega_\Ran}\longrightarrow
\Dmod_{\kappa,\on{co}}(\Bun_G)\otimes \Dmod(\Ran)\overset{\Gamma_{G,\kappa,\Ran}}\longrightarrow \KL(G)_{\kappa,\Ran}.$$

\sssec{} \label{sss:Loc crit init}

Let us specialize for a moment to the case when $\kappa=\crit$. In this case, both sides of
$$\Gamma_{G,\crit,\ul{x}}:\Dmod_{\crit,\on{co}}(\Bun_G)\to \KL(G)_{\crit,\ul{x}}$$
are acted on by $\Sph_{G,\ul{x}}$, and it follows from the construction that the functor $\Gamma_{G,\crit,\ul{x}}$
is compatible with these actions. 

\medskip

Similarly, the functor $\Gamma_{G,\kappa,\Ran}$ is compatible with the action of $\Sph_{G,\Ran}$.

\sssec{}

Note that by construction, we have a commutative diagram
\begin{equation} \label{e:Gamma G and oblv x}
\CD
\Dmod_{\kappa,\on{co}}(\Bun_G) @>{\Gamma_{G,\kappa,\ul{x}}}>> \KL(G)_{\kappa,\ul{x}} \\
@V{\oblv^l_\kappa}VV @VV{\oblv^{(\hg,\fL^+(G))_\kappa}_{\fL^+(G)}}V \\
\QCoh_{\on{co}}(\Bun_G) @>>{\Gamma^{\QCoh}_{G,\kappa,\ul{x}}}> \Rep(\fL^+(G)_{\ul{x}}),
\endCD
\end{equation} 
where $\Gamma^{\QCoh}_{G,\kappa,\ul{x}}$ is the functor of pushforward along 
$$\on{ev}_{\ul{x}}:\Bun_G\to \on{pt}/\fL^+(G)_{\ul{x}}.$$

\begin{rem} \label{r:ignore ren Rep L+G for _*}

Note that $\Rep(\fL^+(G)_{\ul{x}})$ is the renormalized version of $\QCoh(\on{pt}/\fL^+(G)_{\ul{x}})$
(see \secref{sss:Rep L^+G}); 
however, this difference is immaterial for the definition of the functor$\Gamma^{\QCoh}_{G,\kappa,\ul{x}}$:  

\medskip

We have:
$$\QCoh_{\on{co}}(\Bun_G)\simeq \underset{U}{\on{colim}}\, \QCoh(U),$$
where:

\begin{itemize}

\item The index $U$ runs over the filtered posets of quasi-compact open substacks of $\Bun_G$;

\item For $U_1\overset{j_{1,2}}\hookrightarrow U_2$, the transition functor $\QCoh(U_1)\to \QCoh(U_2)$
is given by $(j_{1,2})_*$.

\end{itemize}

\smallskip

A functor out of $\QCoh_{\on{co}}(\Bun_G)$ amounts to a compatible collection of functors out of $\QCoh(U)$.
Thus, in order to define 
$$\Gamma^{\QCoh}_{G,\kappa,\ul{x}}:\QCoh_{\on{co}}(\Bun_G)\to \Rep(\fL^+(G)_{\ul{x}}),$$
we need to define the functors
$$\Gamma^{\QCoh}_{G,\kappa,\ul{x},U}:\QCoh(U)\to \Rep(\fL^+(G)_{\ul{x}}).$$

\medskip

The sought-for functor $\Gamma^{\QCoh}_{G,\kappa,\ul{x},U}$ are defined as the ind-extension of the functor
\begin{multline*}
\QCoh(U)^c\hookrightarrow \QCoh(U)^{>-\infty} \overset{(\on{ev}_{\ul{x},U})_*}\longrightarrow 
\QCoh(\on{pt}/\fL^+(G)_{\ul{x}})^{>-\infty}\simeq \\
\simeq \Rep(\fL^+(G)_{\ul{x}})^{>-\infty}\hookrightarrow
\Rep(\fL^+(G)_{\ul{x}}),
\end{multline*}
where $\on{ev}_{\ul{x},U}$ denotes the restriction of $\on{ev}_{\ul{x}}$ to $U$. 

\end{rem}

\sssec{}

Similarly, we have a commutative diagram 
\begin{equation} \label{e:Gamma G and oblv}
\CD
\Dmod_{\kappa,\on{co}}(\Bun_G)\otimes \Dmod(\Ran) @>{\Gamma_{G,\kappa,\Ran}}>> \KL(G)_{\kappa,\Ran} \\
@V{\oblv^l_\kappa}VV @VV{\oblv^{(\hg,\fL^+(G))_\kappa}_{\fL^+(G)}}V \\
\QCoh_{\on{co}}(\Bun_G) \otimes \Dmod(\Ran) @>>{\Gamma^{\QCoh}_{G,\kappa,\Ran}}> \Rep(\fL^+(G))_\Ran.
\endCD
\end{equation}

\sssec{} \label{sss:ins vac x x'}

Let now $\ul{x}$ and $\ul{x}'$ be two points of $\Ran$ with $\ul{x}\subseteq \ul{x}'$. Consider the
functor
$$\on{ins.vac}_{\ul{x}\subseteq \ul{x}'}:\KL(G)_{\kappa,\ul{x}} \to \KL(G)_{\kappa,\ul{x}'}$$
obtained by inserting the vacuum objects at the points $\ul{x}'-\ul{x}$. 

\medskip

Consider its right adjoint, $(\on{ins.vac}_{\ul{x}\subseteq \ul{x}'})^R$. Explicitly, the functor $(\on{ins.vac}_{\ul{x}\subseteq \ul{x}'})^R$
is given by
\begin{multline*}
\KL(G)_{\kappa,\ul{x}'}\simeq \KL(G)_{\kappa,\ul{x}}\otimes \KL(G)_{\kappa,\ul{x}'-\ul{x}}
\overset{\on{Id}\otimes \oblv^{(\hg,\fL^+(G))_\kappa}_{\fL^+(G)}}\longrightarrow \\
\to \KL(G)_{\kappa,\ul{x}}\otimes \Rep(\fL^+(G)_{\ul{x}'-\ul{x}}) \overset{\on{Id}\otimes \on{inv}^{\fL^+(G)_{\ul{x}'-\ul{x}}}}\longrightarrow
\KL(G)_{\kappa,\ul{x}}.
\end{multline*} 

It follows from the commutation of \eqref{e:Gamma G and oblv x} that we have a canonical isomorphism
\begin{equation} \label{e:unitality Gamma x}
\Gamma_{G,\kappa,\ul{x}} \simeq (\on{ins.vac}_{\ul{x}\subseteq \ul{x}'})^R\circ \Gamma_{G,\kappa,\ul{x}'}.
\end{equation} 

This isomorphism expresses the \emph{unital} structure on the assignment
$$\ul{x}\rightsquigarrow \Gamma_{G,\crit,\ul{x}},$$
to be discussed in \secref{s:unitality}.

\sssec{} \label{sss:inv of Gamma}

In particular, for any $\ul{x}$, the functor
$$\Dmod_{\kappa,\on{co}}(\Bun_G) \overset{\Gamma_{G,\kappa,\ul{x}}}\longrightarrow \KL(G)_{\kappa,\ul{x}} 
\overset{\on{inv}^{\fL^+(G)_{\ul{x}}}}\longrightarrow \Vect$$
identifies with $\Gamma(\Bun_G,\oblv_\kappa^l(-))$.

\ssec{Localization functor as a left adjoint} \label{ss:defn Loc}

As was mentioned previously, we construct the localization functor $\Loc_{G,\kappa}$
to be \emph{essentially} the left adjoint of $\Gamma_G$. However, there is a caveat:
this adjunction takes place over quasi-compact open substacks $U\subset \Bun_G$, 
and we obtain the corresponding functors $\Loc_{G,\kappa,U}$. We then obtain the
sought-for functor $\Loc_{G,\kappa}$ by passing to the limit. 

\sssec{}

Let 
\begin{equation} \label{e:qc U}
U\overset{j_U}\hookrightarrow \Bun_G
\end{equation}
be a quasi-compact open substack. Consider the corresponding functor
$$j_{*,\on{co}}:\Dmod_\kappa(U)\to \Dmod_{\kappa,\on{co}}(\Bun_G).$$

\medskip

Denote
$$\Gamma_{G,\kappa,\ul{x},U}:=\Gamma_{G,\kappa,\ul{x}}\circ j_{*,\on{co}},
\quad \Dmod_\kappa(U)\to \KL(G)_{\kappa,\ul{x}}.$$

\sssec{}

We claim:

\begin{lem} \label{l:Loc U} 
The functor $\Gamma_{G,\kappa,\ul{x},U}$
admits a left adjoint. 
\end{lem} 

\begin{proof}

Since the essential image of
$$\ind_{\fL^+(G)}^{(\hg,\fL^+(G))_\kappa}:\Rep(\fL^+(G)_{\ul{x}})\to \KL(G)_{\kappa,\ul{x}}$$
generates the target category, it suffices to show that the composite functor
$$\oblv^{(\hg,\fL^+(G))_\kappa}_{\fL^+(G)}\circ \Gamma_{G,\kappa,\ul{x},U},\quad \Dmod_\kappa(U)\to  \Rep(\fL^+(G)_{\ul{x}})$$
admits a left adjoint.

\medskip 

The above functor identifies with
$$(\on{ev}_{\ul{x}}\circ j)_*\circ \oblv^l_\kappa.$$

In this composition, both arrows admit left adjoints: the left adjoint of $\oblv^l_\kappa$ is $\ind^l_\kappa$, and 
the left adjoint of $(\on{ev}_{\ul{x}}\circ j)_*$ is $(\on{ev}_{\ul{x}}\circ j)^*$. 

\end{proof}

\sssec{}

Let us denote the left adjoint in \lemref{l:Loc U} by $\Loc_{G,\kappa,\ul{x},U}$. 

\medskip

For an inclusion between quasi-compact open substacks
$$U_1\overset{j_{1,2}}\hookrightarrow U_2,$$
we have
$$(j_2)_{*,\on{co}}\circ (j_{1,2})_*\simeq (j_1)_{*,\on{co}}.$$

Hence, we obtain a canonical identification
$$\Loc_{G,\kappa,\ul{x},U_1}\simeq j_{1,2}^*\circ \Loc_{G,\kappa,\ul{x},U_2}.$$

Therefore, the system of functors 
$$U\rightsquigarrow \{\Loc_{G,\kappa,\ul{x},U}\}$$
gives rise to a functor
\begin{equation} \label{e:Loc as right adj}
\Loc_{G,\kappa,\ul{x}}:\KL(G)_{\kappa,\ul{x}}\to \Dmod_\kappa(\Bun_G),
\end{equation} 
so that for every \eqref{e:qc U}, we have
$$j^*\circ \Loc_{G,\kappa,\ul{x}}\simeq \Loc_{G,\kappa,\ul{x},U}.$$

\medskip

The functor \eqref{e:Loc as right adj} is the sought-for localization functor. 

\sssec{}

The entire preceding discussion generalizes to the case when $\ul{x}$ is allowed to move in families 
over $\Ran$. In particular, we obtain a functor
$$\Loc_{G,\kappa,\Ran}:\KL(G)_{\kappa,\Ran}\to \Dmod_\kappa(\Bun_G\times \Ran).$$

\medskip

Let $\Loc_{G,\kappa}$ denote the composition
$$\KL(G)_{\kappa,\Ran}\overset{\Loc_{G,\kappa,\Ran}}\longrightarrow \Dmod_\kappa(\Bun_G\times \Ran)\to \Dmod_\kappa(\Bun_G),$$
where the second arrow is the functor of !-pushforward.

\begin{rem}

The above construction of the localization functor is essentially equivalent to the one from \cite[Sect. 4.1]{CF}. 

\end{rem} 

\sssec{}

Properties of the functor $\Gamma_{G,\kappa',\Ran}$ induce corresponding properties
of the functor $\Loc_{G,\kappa,\Ran}$. We will now list some of them. 

\sssec{}

By adjunction, for every quasi-compact open as in \eqref{e:qc U}, from diagram \eqref{e:Gamma G and oblv x} 
we obtain a commutative diagram:
\begin{equation} \label{e:localization induction diagram x U}
\CD
\QCoh(U) @>{\ind^l_\kappa}>> \Dmod_\kappa(U) \\
@A{j^*\circ \on{Loc}^{\QCoh}_{G,\ul{x}}}AA @AA{\Loc_{G,\kappa,\ul{x},U}}A \\
\Rep(\fL^+(G)_{\ul{x}}) @>>{\ind_{\fL^+(G)}^{(\hg,\fL^+(G))_\kappa}}> \KL(G)_{\kappa,\ul{x}},
\endCD
\end{equation} 
where:
$$\on{Loc}^{\QCoh}_{G,\ul{x}}:\Rep(\fL^+(G)_{\ul{x}})\to \QCoh(\Bun_G)$$ is the functor of 
pullback along $\on{ev}_{\ul{x}}:\Bun_G\to \on{pt}/\fL^+(G)_{\ul{x}}$. 

\begin{rem} \label{r:ignore ren Rep L+G for ^*}

As in Remark \ref{r:ignore ren Rep L+G for _*}, the difference between $\Rep(\fL^+(G)_{\ul{x}})$ and
$\QCoh(\on{pt}/\fL^+(G)_{\ul{x}})$ does not play a role in the definition of the functor $\on{Loc}^{\QCoh}_{G,\ul{x}}$.

\medskip

Namely, $\on{Loc}^{\QCoh}_{G,\ul{x}}$ is defined as the ind-extension of
$$\Rep(\fL^+(G)_{\ul{x}})^c \hookrightarrow  \QCoh(\on{pt}/\fL^+(G)_{\ul{x}}) \overset{(\on{ev}_{\ul{x}})^*}\hookrightarrow
\QCoh(\Bun_G).$$

Note also that for a quasi-compact $U$, the functor
$$\on{Loc}^{\QCoh}_{G,\ul{x},U}:=j^*\circ \on{Loc}^{\QCoh}_{G,\ul{x}}$$
is the left adjoint of the functor $\Gamma^{\QCoh}_{G,\kappa,\ul{x},U}$. 

\end{rem}

\sssec{}

Passing to the limit over \eqref{e:qc U}, we obtain a commutative diagram 
\begin{equation} \label{e:localization induction diagram x}
\CD
\QCoh(\Bun_G) @>{\ind^l_\kappa}>> \Dmod_\kappa(\Bun_G) \\
@A{\on{Loc}^{\QCoh}_{G,\ul{x}}}AA @AA{\Loc_{G,\kappa,\ul{x}}}A \\
\Rep(\fL^+(G)_{\ul{x}}) @>>{\ind_{\fL^+(G)}^{(\hg,\fL^+(G))_\kappa}}> \KL(G)_{\kappa,\ul{x}}.
\endCD
\end{equation} 

Similarly, we have 
\begin{equation} \label{e:localization induction diagram Ran}
\CD
\QCoh(\Bun_G)\otimes \Dmod(\Ran) @>{\ind^l_\kappa}>> \Dmod_\kappa(\Bun_G) \otimes \Dmod(\Ran) \\
@A{\on{Loc}^{\QCoh}_{G,\Ran}}AA @AA{\Loc_{G,\kappa,\Ran}}A \\
\Rep(\fL^+(G))_\Ran @>>{\ind_{\fL^+(G)}^{(\hg,\fL^+(G))_\kappa}}> \KL(G)_{\kappa,\Ran}
\endCD
\end{equation} 
and
\begin{equation} \label{e:localization induction diagram}
\CD
\QCoh(\Bun_G) @>{\ind^l_\kappa}>> \Dmod_\kappa(\Bun_G) \\
@A{\on{Loc}^{\QCoh}_G}AA @AA{\Loc_{G,\kappa}}A \\
\Rep(\fL^+(G))_\Ran @>{\ind_{\fL^+(G)}^{(\hg,\fL^+(G))_\kappa}}>> \KL(G)_{\kappa,\Ran},
\endCD
\end{equation} 
where $\on{Loc}^{\QCoh}_G$ is the functor of pull-push along
$$\Bun_G\leftarrow \Bun_G\times \Ran \to (\on{pt}/\fL^+(G))_\Ran.$$

\sssec{}

Let now $\ul{x}$ and $\ul{x}'$ be two points of $\Ran$ with $\ul{x}\subseteq \ul{x}'$. Recall the functor
\begin{equation} \label{e:ins vac x x'}
\on{ins.vac}_{\ul{x}\subseteq \ul{x}'}:\KL(G)_{\kappa,\ul{x}} \to \KL(G)_{\kappa,\ul{x}'}.
\end{equation} 

\medskip

For every \eqref{e:qc U}, from \eqref{e:unitality Gamma x} we obtain a canonical isomorphism
\begin{equation} \label{e:unitality Loc x U}
\Loc_{G,\kappa,\ul{x}',U} \circ \on{ins.vac}_{\ul{x}\subseteq \ul{x}'}\simeq \Loc_{G,\kappa,\ul{x},U}.
\end{equation}

Passing to the limit over \eqref{e:qc U}, we obtain a canonical isomorphism
\begin{equation} \label{e:unitality Loc x}
\Loc_{G,\kappa,\ul{x}'} \circ \on{ins.vac}_{\ul{x}\subseteq \ul{x}'}\simeq \Loc_{G,\kappa,\ul{x}}.
\end{equation}

\medskip

In \secref{sss:unitality fam} we will formulate a version of \eqref{e:unitality Loc x} when the points 
$\ul{x}$ and $\ul{x}'$ move in families over the Ran space. 

\sssec{} \label{sss:Loc of Vac}

It follows from \eqref{e:localization induction diagram x} that  for any $\ul{x}$, we have 
$$\Loc_{G,\kappa,\ul{x}}(\on{Vac}(G)_{\kappa,\ul{x}})\simeq \ind^l_\kappa(\CO_{\Bun_G}).$$

Equivalently,
\begin{equation} \label{e:Loc of Vac}
\Loc_{G,\kappa,\Ran}(\on{Vac}(G)_{\kappa,\Ran})\simeq \ind^l_\kappa(\CO_{\Bun_G})\boxtimes \omega_\Ran.
\end{equation}

Note that
$$\ind^l_\kappa(\CO_{\Bun_G})\simeq \on{D}_{\Bun_G,\kappa},$$
where
$$\on{D}_{\Bun_G,\kappa}\in \Dmod_\kappa(\Bun_G)$$
is the D-module of differential operators, viewed a twisted \emph{left} D-module. 

\sssec{}

Let us specialize for a moment to the case when $\kappa=\crit$. Then from 
\secref{sss:Loc crit init} we obtain that for every $U$, the functor $\Loc_{G,\kappa,\ul{x},U}$
is compatible with the action of $\Sph_{G,\ul{x}}$. Hence, so is the functor $\Loc_{G,\kappa,\ul{x}}$.

\medskip

Similarly, the functor $\Loc_{G,\kappa,\Ran}$ is compatible with the action of $\Sph_{G,\Ran}$.

\sssec{}

For future reference we note that the entire discussion in this subsection applies to to the case of
an infinite level structure at a given $\ul{x}_0\in \Ran$. I.e., for a given quasi-compact $U\subset \Bun_G$
and 
$$U^{\on{level}_{\ul{x}_0}}:=U\underset{\Bun_G}\times \Bun^{\on{level}_{\ul{x}_0}},$$
we can consider the left adjoint 
$$\Loc_{G,\kappa,\ul{x}_0,U}:\hg\mod_{\kappa,\ul{x}_0}\to \Dmod_{\kappa,\on{co}}(U^{\on{level}_{\ul{x}_0}}_G)$$
of
$$\Gamma(U^{\on{level}_{\ul{x}_0}},-)^{\on{enh}}:\Dmod_{\kappa,\on{co}}(U^{\on{level}_{\ul{x}_0}}_G)\to 
\hg\mod_{\kappa,\ul{x}_0}.$$

The functors $\Loc_{G,\kappa,\ul{x}_0,U}$ glue to a functor
 $$\Loc_{G,\kappa,\ul{x}_0,U}:\hg\mod_{\kappa,\ul{x}_0}\to \Dmod_{\kappa,\on{co}}(\Bun^{\on{level}_{\ul{x}_0}}_G).$$
 
 \medskip
 
 For the Ran space version, one should consider 
 $$\Ran_{\ul{x}_0}:=\{\ul{x}_0\}\underset{\Ran}\times \Ran^{\subseteq},$$
 where:
 
 \begin{itemize}
 
 \item $\Ran^{\subseteq}$ is as in \secref{sss:Ran subset};
 
 \item The fiber product is formed using the map $\on{pr}_{\on{small}}:\Ran^{\subseteq}\to \Ran$.
 
 \end{itemize}

\ssec{The fiber of the localization functor}

\sssec{}

Fix a $k$-point $\CP_G\in \Bun_G$. The goal of this subsection is to describe the functor
\begin{equation} \label{e:fiber of localization}
\Loc_{G,\kappa,\ul{x}} \overset{\Loc_{G,\kappa,\ul{x}}}\longrightarrow 
\Dmod_\kappa(\Bun_G) \overset{!\on{-fiber\,at}\,\CP_G}\longrightarrow \Vect.
\end{equation}

\sssec{}

Consider the Lie algebra $\Gamma(X-\ul{x},\fg_{\CP_G})$. Laurent expansion defines a map
$$\Gamma(X-\ul{x},\fg_{\CP_G})\to \fL(\fg_{\CP_G})_{\ul{x}},$$
and recall that the Kac-Moody extension
$$0\to k\to \hg_{\kappa,\CP_G,\ul{x}}\to \fL(\fg_{\CP_G})_{\ul{x}}\to 0$$
admits a canonical splitting over $\Gamma(X-\ul{x},\fg_{\CP_G})$. Hence, we have a well-defined restriction functor
\begin{equation} \label{e:KM to out}
\hg\mod_{\kappa,\CP_G,\ul{x}}\to \Gamma(X-\ul{x},\fg_{\CP_G})\mod.
\end{equation}

\sssec{}

Composing with
$$\KL(G)_{\kappa,\ul{x}}\overset{\alpha_{\CP_G,\on{taut}}}\longrightarrow 
\KL(G)_{\kappa,\CP_G,\ul{x}}\to \hg\mod_{\kappa,\CP_G,\ul{x}}$$
we obtain a functor 
\begin{equation} \label{e:from KL to out}
\KL(G)_{\kappa,\ul{x}}\to \Gamma(X-\ul{x},\fg_{\CP_G})\mod.
\end{equation}

\sssec{}

We claim:

\begin{prop} \label{p:fiber of localization}
The functor \eqref{e:fiber of localization} identifies canonically with the composition of \eqref{e:from KL to out}
and the functor of $\Gamma(X-\ul{x},\fg_{\CP_G})$-coinvariants
$$\Gamma(X-\ul{x},\fg_{\CP_G})\mod\to \Vect.$$
\end{prop}

\begin{rem} \label{r:loc f.d}
An analog of \propref{p:fiber of localization} for the localization functor in the finite-dimensional situation
is obvious:

\medskip

Let $\fh$ be a (discrete\footnote{As opposed to Tate.}) 
Lie algebra and let $\CY$ be a smooth variety equipped with an action of $\fh$
by vector fields. Then the corresponding localization functor 
$$\Loc_{\fh,\CY}:\fh\mod\to \Dmod(\CY),$$
left adjoint to 
$$\Gamma(\CY,\oblv^l(-))^{\on{enh}}:\Dmod(\CY)\to \fh\mod,$$
is given by
\begin{equation} \label{e:naive Loc}
\on{D}_\CY\underset{U(\fh)}\otimes -.
\end{equation} 

If $y\in \CY$ is a point for which the action map $\fh\to T_y(\CY)$ is surjective, the composition
$$\fh\mod\overset{\Loc_{\fh,\CY}}\to \Dmod(\CY)\overset{!\on{-fiber\,at}\,y}\longrightarrow \Vect$$
identifies with the functor of coinvariants with respect to
$$\on{Stab}_y(\fh)\subset \fh.$$

This follows from the fact that the *-fiber at $y$ of $\on{D}_\CY$ (as an object of $\QCoh(\CY)$ via left
multiplication) identifies, as a $\fh$-module, with
$$\ind^\fh_{\on{Stab}_y(\fh)}(k).$$

\end{rem}

\sssec{}

One can prove \propref{p:fiber of localization} directly by emulating the argument in Remark \ref{r:loc f.d}. 

\medskip

In fact, such an assertion is valid for $\Bun_G^{\on{level}_{\ul{x}}}$ replaced by a pro-scheme $\CY$ equipped with
an action of $\fL(G)^\wedge_{\ul{x}}$ (the formal completion of $\fL(G)_{\ul{x}}$ along $\fL^+(G)_{\ul{x}}$), 
such that $\CY/\fL^+(G)_{\ul{x}}$ is locally of finite type, and a point $y\in \CY$ at which the action is
infinitesimally transitive, i.e.,
$$\fL(\fg)_{\ul{x}}\to T_y(\CY)$$
is surjective.

\medskip

We will, however, supply a different argument, specific to the case of $\Bun_G$, see  
\secref{ss:proof of fiber of localization}.

\sssec{}

As an immediate corollary of \propref{p:fiber of localization} we obtain:

\begin{cor} \label{c:Loc is right-exact}
The functor 
$$\Loc_{G,\kappa,\ul{x}}:\KL(G)_{\kappa,\ul{x}}\to \Dmod_\kappa(\Bun_G)$$
is right t-exact, when $\Dmod_\kappa(\Bun_G)$ is equipped with the \emph{left}
t-structure, i.e., one for which the functor $\oblv^l_\kappa$ is t-exact. 
\end{cor}

\begin{proof}

We need to show that the composite functor
$$\KL(G)_{\kappa,\ul{x}}\overset{\Loc_{G,\kappa,\ul{x}}}\longrightarrow 
\Dmod_\kappa(\Bun_G) \overset{\oblv^l_\kappa}\to \QCoh(\Bun_G)$$
is right t-exact. 

\medskip

In order to prove that, it suffices to show that the composition of the above
functor with the functor of *-fiber at any field-valued point of $\Bun_G$
is right t-exact. 

\medskip

By base change, we can assume that the point in question is rational. In this case,
the corresponding functor identifies with the functor \eqref{e:fiber of localization}.

\end{proof} 

\begin{cor}
The functor $\Loc_{G,\kappa,\ul{x}}$ annihilates infinitely connective objects
(i.e., objects that belong to $(\KL(G)_{\kappa,\ul{x}})^{<-n}$ for any $n$).
\end{cor}

\begin{proof}

Follows from the fact that the t-structure on $\Dmod_\kappa(\Bun_G)$ is
separated.

\end{proof}

\ssec{Localization functor as the dual} \label{ss:Loc as dual}

\sssec{}

Let $\kappa'$ be the reflected level, i.e.,
$$\kappa':=-\kappa+{2\cdot\crit} .$$

We claim that that we have a canonical duality
\begin{equation} \label{e:global duality}
(\Dmod_{\kappa',\on{co}}(\Bun_G))^\vee \simeq \Dmod_{\kappa}(\Bun_G)
\end{equation}
for which the dual of the functor
$$\oblv^l_{\kappa'}:\Dmod_{\kappa',\on{co}}(\Bun_G)\to \QCoh_{\on{co}}(\Bun_G)$$
is the functor
$$\ind^l_\kappa:\QCoh(\Bun_G)\to \Dmod_\kappa(\Bun_G),$$ 
with respect to the identification\footnote{We warn the reader that the category $\QCoh(\Bun_G)$
is not dualizable.} 
$$\on{Funct}_{\on{cont}}(\QCoh_{\on{co}}(\Bun_G),\Vect)\simeq \QCoh(\Bun_G).$$

\sssec{} \label{sss:shift duality}

Indeed, we start with the identification
$$(\Dmod_{\kappa',\on{co}}(\Bun_G))^\vee \simeq \Dmod_{-\kappa'}(\Bun_G),$$
given by Verdier duality, and compose it with the functor
\begin{multline}  \label{e:tw can bundle}
\Dmod_{-\kappa'}(\Bun_G) \overset{\otimes K_{\Bun_G}}\longrightarrow 
\Dmod_{-\kappa'+\on{dlog}(K_{\Bun_G})}(\Bun_G)\simeq \Dmod_{-\kappa'+{2\cdot\crit} }(\Bun_G)=\\
=\Dmod_{\kappa}(\Bun_G)\overset{[\dim(\Bun_G)]}\longrightarrow \Dmod_{\kappa}(\Bun_G),
\end{multline}
where:

\begin{itemize}

\item $K_{\Bun_G}$ is the canonical line bundle on $\Bun_G$, so that $K_{\Bun_G}[\dim(\Bun_G)]\simeq \omega_{\Bun_G}$;

\item We have used the identification $\on{dlog}(K_{\Bun_G})=\on{dlog}(\det_{\Bun_G})={2\cdot\crit}$ from \secref{sss:crit 1/2 can}.

\end{itemize} 

\sssec{}

We have the following assertion:

\begin{prop} \label{p:Loc as dual x} With respect to the identifications \eqref{e:global duality} and 
\begin{equation} \label{e:KL duality kappa}
(\KL(G)_{\kappa,\ul{x}})^\vee\simeq \KL(G)_{\kappa',\ul{x}}
\end{equation} 
of \eqref{e:KL self-duality kappa}, the functor
$$\Loc_{G,\kappa,\ul{x}}:\KL(G)_{\kappa,\ul{x}} \to \Dmod_\kappa(\Bun_G)$$
identifies canonically with the dual of
$$\Gamma_{G,\kappa',\ul{x}}:\Dmod_{\kappa',\on{co}}(\Bun_G)\to \KL(G)_{\kappa',\ul{x}}.$$
The induced identification
\begin{multline*} 
\ind^l_\kappa \circ (\on{ev}_{\ul{x}})^*\simeq 
\Loc_{G,\kappa,\ul{x}}\circ \ind_{\fL^+(G)}^{(\hg,\fL^+(G))_\kappa} \simeq 
(\Gamma_{G,\kappa',\ul{x}})^\vee \circ (\oblv_{\fL^+(G)}^{(\hg,\fL^+(G))_{\kappa'}})^\vee\simeq \\
\simeq \left(\oblv_{\fL^+(G)}^{(\hg,\fL^+(G))_{\kappa'}}\circ \Gamma_{G,\kappa',\ul{x}}\right)^\vee
\simeq ((\on{ev}_{\ul{x}})_*\circ \oblv^l_{\kappa'})^\vee\simeq (\oblv^l_{\kappa'})^\vee\circ 
((\on{ev}_{\ul{x}})_*)^\vee\simeq \\
\simeq \ind^l_\kappa \circ ((\on{ev}_{\ul{x}})_*)^\vee \simeq \ind^l_\kappa \circ (\on{ev}_{\ul{x}})^*
\end{multline*}
is the identity map. 

\end{prop}

This assertion is proved in \cite[Theorem 4.0.5(2)]{CF}.\footnote{In {\it loc. cit.} the dual functor to $\Loc_{G,\kappa}$ is
denoted $\Loc_{\on{co}}$.}

\begin{rem}

The proof of \propref{p:Loc as dual x} in \cite{CF} essentially emulates the following finite-dimensional 
phenomenon. 

\medskip

Let $\CY$ and $\fh$ be as in Remark \ref{r:loc f.d}. On the one hand, we can consider 
the adjoint pair
$$\Loc_{\fh,\CY}:\fh\mod \leftrightarrows \Dmod(\CY):\Gamma(\CY,\oblv^l(-))^{\on{enh}}.$$

On the other hand, consider the canonical line bundle $K_\CY$ as a line bundle acted on by $\fh$,
and consider the corresponding functor
$$\Gamma(\CY,K_\CY\otimes \oblv^l(-))^{\on{enh}}:\Dmod(\CY)\to \fh\mod.$$

Let 
$$\Loc_{\fh,\CY,K_\CY}: \fh\mod\to \Dmod(\CY)$$
denote the left adjoint of $\Gamma(\CY,K_\CY\otimes \oblv^l(-))^{\on{enh}}$.

\medskip

Then the functors $\Loc_{\fh,\CY,K_\CY}[\dim(\CY)]$ and $\Gamma(\CY,\oblv^l(-))^{\on{enh}}$
are mutually dual in terms of the Verdier duality identification
$$\Dmod(\CY)^\vee \simeq \Dmod(\CY).$$

This follows from the expression for $\Loc_{\fh,\CY}$ given by formula \eqref{e:naive Loc},
and a similar formula for $\Loc_{\fh,\CY,K_\CY}$. 

\end{rem} 

\sssec{}

The assertion of \propref{p:Loc as dual x} admits an immediate generalization when $\ul{x}$ moves in families
over the Ran space: 

\begin{prop} \label{p:Loc as dual}
With respect to the identifications \eqref{e:global duality} and 
$$(\KL(G)_{\kappa,\Ran})^\vee\simeq \KL(G)_{\kappa',\Ran}$$

\smallskip

\noindent{\em(a)} The functor
$$\Loc_{G,\kappa,\Ran}:\KL(G)_{\kappa,\Ran} \to \Dmod_\kappa(\Bun_G)\otimes \Dmod(\Ran)$$
identifies canonically with the dual of
$$\Gamma_{G,\kappa',\Ran}:\Dmod_{\kappa',\on{co}}(\Bun_G)\otimes \Dmod(\Ran)\to \KL(G)_{\kappa',\Ran}.$$

The induced identification
\begin{multline*} 
\ind^l_\kappa \circ (\on{ev}_\Ran)^*\simeq 
\Loc_{G,\kappa,\Ran}\circ \ind_{\fL^+(G)}^{(\hg,\fL^+(G))_\kappa} \simeq 
(\Gamma_{G,\kappa',\Ran})^\vee \circ (\oblv_{\fL^+(G)}^{(\hg,\fL^+(G))_{\kappa'}})^\vee\simeq \\
\simeq \left(\oblv_{\fL^+(G)}^{(\hg,\fL^+(G))_{\kappa'}}\circ \Gamma_{G,\kappa',\Ran}\right)^\vee
\simeq ((\on{ev}_\Ran)_*\circ \oblv^l_{\kappa'})^\vee\simeq (\oblv^l_{\kappa'})^\vee\circ 
((\on{ev}_\Ran)_*)^\vee\simeq \\
\simeq \ind^l_\kappa \circ ((\on{ev}_\Ran)_*)^\vee \simeq \ind^l_\kappa \circ (\on{ev}_\Ran)^*
\end{multline*}
is the identity map. 

\smallskip

\noindent{\em(b)} The functor
$$\Loc_{G,\kappa}:\KL(G)_{\kappa,\Ran} \to \Dmod_\kappa(\Bun_G)\otimes \Dmod(\Ran)$$
identifies canonically with the dual of
$$\Gamma_{G,\kappa'}:\Dmod_{\kappa',\on{co}}(\Bun_G)\to \KL(G)_{\kappa',\Ran}.$$
\end{prop}

For the proof, see \cite[Theorem 4.0.5(2)]{CF}. 

\sssec{} \label{sss:Loc of Vac via dual}

Note that by combining \secref{sss:inv of Gamma} and \propref{p:Loc as dual} with the fact 
that the functors
$$\KL(G)_{\kappa'}\overset{\on{inv}_{\fL^+(G)}}\longrightarrow \Vect \text{ and } 
\Vect\overset{k\mapsto \on{Vac}(G)_\kappa}\longrightarrow \KL(G)_\kappa$$
and 
$$\Dmod_{\kappa',\on{co}}\overset{\Gamma(\Bun_G,-)\circ \oblv^l_{\kappa}}\longrightarrow \Vect \text{ and }
\Vect\overset{k\mapsto \ind^l_\kappa(\CO_{\Bun_G})}\longrightarrow \Dmod_{\kappa}(\Bun_G)$$
are mutually dual, we obtain an identification
\begin{equation} \label{e:Loc of Vac via dual}
\Loc_{G,\kappa,\Ran}(\on{Vac}_{G,\Ran})\simeq \ind^l_\kappa(\CO_{\Bun_G})\boxtimes \omega_\Ran.
\end{equation}

However, it follows formally that the
identification \eqref{e:Loc of Vac via dual} is the same as that in \eqref{e:Loc of Vac}. 

\section{Digression: local-to-global functors and unitality} \label{s:unitality}

In this section we will introduce a general framework that formalizes the unital property of the functors 
$$\Poinc_{G,!,\Ran},\,\, \Poinc_{G,*,\Ran} \text{ and }\Loc_{G,\kappa,\Ran}.$$

\medskip

The unital property says, roughly speaking, that the insertion of the vacuum\footnote{In the main body of this section, we use the 
word ``unit" instead of ``vacuum".} does not change the value
of the functor (see \secref{sss:global unitality example}).

\medskip

A key phenomenon that we will observe is the following: insertion of the vacuum along the entire Ran space
improves the unital property of the functor, see \secref{ss:integrated ins}. The functor of factorization
homology and its generalizations are particular cases of this construction, see \secref{ss:ch homology}. 

\ssec*{11.0. What is this section about?}

As this section deals with some abstract material, a general introduction is in order.

\sssec{}

In this section, we study the general formalism of local-to-global functors. The (\emph{local}) source of such a functor 
is a crystal of categories $\ul\bC^{\on{loc}}$ over the Ran space, while its (\emph{global}) target is a single category $\bC^{\on{glob}}$. 
Roughly, for a space $\CZ$ equipped with a map $\ul{x}:\CZ\to\Ran$, the value of $\ul\bC^{\on{loc}}$ on $\CZ$ is a 
category $\bC^{\on{loc}}_\CZ=\bC^{\on{loc}}_{\CZ,\ul{x}}$, and a local-to-global functor $\ul\sF$ is a 
compatible collection of functors
\[\sF_\CZ=\sF_{\CZ,\ul{x}}:\bC^{\on{loc}}_{\CZ,\ul{x}}\to \bC^{\on{glob}}\otimes\Dmod(\CZ)\]
for every $\CZ$ and $\ul{x}$. (Here and below, the words ``compatible collection" mean ``collection equipped with higher coherence data".) 

\sssec{}

Next, we introduce the notion of a unital crystal of categories over the Ran space (see \secref{ss:unitality loc}). 
Informally, a unital structure on a sheaf $\ul\bC^{\on{loc}}$ is a compatible collection of functors
\[\on{ins.unit}_{\ul{x_1}\subseteq \ul{x}_2}: \bC^{\on{loc}}_{\CZ,\ul{x}_1}\to \bC^{\on{loc}}_{\CZ,\ul{x}_2}\]
for every space $\CZ$ and two maps $\ul{x}_1,\ul{x}_2:\CZ\to \Ran$ such that $\ul{x}_1\subseteq \ul{x}_2$. 
Here we view $\CZ$-points of the Ran space as $\CZ$-families of finite subsets of $X$.

\sssec{}

Suppose now that $\ul\bC^{\on{loc}}$ is a unital crystal of categories, and $\ul\sF$ is a local-to-global functor from $\ul\bC^{\on{loc}}$ to a 
category $\bC^{\on{glob}}$. We then introduce the notion of a unital structure on $\ul\sF$; informally, it is a compatible collection of 
natural transformations
\begin{equation}\label{e:unital functor}
\sF_{\CZ,\ul{x}_1}\to \sF_{\CZ,\ul{x}_2}\circ \on{ins.unit}_{\ul{x_1}\subseteq \ul{x}_2}.
\end{equation}
for every $\CZ$, $\ul{x}_1$, and $\ul{x}_2$ as above. 

\medskip

In fact, we have two notions: a (\emph{strict}) unital structure, where the transformations \eqref{e:unital functor} 
are required to be isomorphisms, and a \emph{lax} unital structure, where \eqref{e:unital functor} can be arbitrary transformations. 

\medskip

Accordingly, we obtain two categories of unital local-to-global functors: the category of (strictly) 
unital local-to-global functors, and the category of lax unital  local-to-global functors; the former is a full subcategory of the latter. 
We denote the categories by

\begin{equation} \label{e:strict untl into lax unital-intro}
\on{Funct}^{\on{loc}\to \on{glob},\on{untl}}(\ul\bC^{\on{loc}},\bC^{\on{glob}}) \subset
\on{Funct}^{\on{loc}\to \on{glob},\on{lax-untl}}(\ul\bC^{\on{loc}},\bC^{\on{glob}}).
\end{equation}

\sssec{}

The main subject of this section is a construction on local-to-global functors, which we call \emph{the integrated insertion of the unit}. It can be defined as the left adjoint of the embedding \eqref{e:strict untl into lax unital-intro}:

\[{\int\on{ins.unit}}:\on{Funct}^{\on{loc}\to \on{glob},\on{lax-untl}}(\ul\bC^{\on{loc}},\bC^{\on{glob}})\to
\on{Funct}^{\on{loc}\to \on{glob},\on{untl}}(\ul\bC^{\on{loc}},\bC^{\on{glob}}).
\]
However, the functor admits a geometric description. Remarkably, the description makes sense for all
(i.e., not necessarily lax unital ) local-to-global functors. 

\sssec{}

We will use this formalism in \secref{s:prop Loc}
in the context of compatibility between certain natural constructions and local-to-global functors. We will see that, in three different situations, the compatibility is only lax at the start, but composition with the functor $\int\on{ins.unit}$ makes it strict.

\ssec{Setup for local-to-global functors} \label{ss:local-to-global set up}

\sssec{} \label{sss:local to global functor abs}

Let $\ul\bC^{\on{loc}}$ be a crystal of categories over $\Ran$ (see \secref{ss:cat over Ran}). Let $\bC^{\on{glob}}$
be a DG category, and let us be given a functor
$$\ul\sF:\ul\bC^{\on{loc}} \to \bC^{\on{glob}} \otimes \ul\Dmod(\Ran),$$
where $\ul\Dmod(\Ran)$ is the unit crystal of categories over $\Ran$. 

\medskip

Thus, for every space $\CZ$ mapping to $\Ran$, we have a category $\bC^{\on{loc}}_\CZ$, tensored over 
$\Dmod(\CZ)$ and a functor
\begin{equation} \label{e:F Z}
\sF_\CZ:\bC^{\on{loc}}_\CZ\to  \bC^{\on{glob}}\otimes \Dmod(\CZ).
\end{equation} 

\begin{rem}

In the above procedure, we associate to $\CZ\to \Ran$ the category of \emph{cristalline} sections
of $\ul\bC^{\on{loc}}$ over $\CZ$, i.e., the category of sections of $\ul\bC^{\on{loc}}$ over $\CZ_\dR$,
cf. \secref{sss:cryst term}. 

\end{rem}

\sssec{} \label{e:F Z int}

Assume for a moment that $\CZ$ is pseudo-proper, so that the functor
$$\on{C}^\cdot_c(\CZ,-):\Dmod(\CZ)\to \Vect$$
left adjoint to $k\mapsto \omega_\CZ$ is defined (see \secref{c:integration ps-proper}). In this case we will denote by $\sF_{\int_\CZ}$ the composition
$$\bC^{\on{loc}}_\CZ\overset{\sF_\CZ}\to \Dmod(\CZ)\otimes \bC^{\on{glob}} 
\overset{\on{C}^\cdot_c(\CZ,-)\otimes \on{Id}}\longrightarrow \bC^{\on{glob}}.$$

\sssec{}

In particular, for $\CZ=\Ran$ and the identity map, we obtain the category 
$\bC^{\on{loc}}_\Ran$ and a functor
$$\sF_\Ran:\bC^{\on{loc}}_\Ran\to \bC^{\on{glob}}\otimes \Dmod(\Ran).$$

\medskip

We will also use the symbol $\sF: \bC^{\on{loc}}_\Ran\to \bC^{\on{glob}}$
for the functor $\sF_{\int_\Ran}$. 

\sssec{} \label{sss:F recovers uF}

Note that the datum of $\sF$ recovers that of $\ul\sF$. Namely, the functor $\sF_\CZ$ identifies with
$$\bC^{\on{loc}}_\CZ \to \bC^{\on{loc}}_{\Ran\times \CZ} \simeq \bC^{\on{loc}}_\Ran\otimes \Dmod(\CZ)
\overset{\sF\otimes \on{Id}}\longrightarrow  \bC^{\on{glob}}\otimes \Dmod(\CZ),$$
where:

\begin{itemize}

\item $\Ran\times \CZ$ is viewed as a space over $\Ran$ via the projection on the first factor;

\smallskip

\item The arrow $\bC^{\on{loc}}_\CZ \to \bC^{\on{loc}}_{\Ran\times \CZ}$ is the !-pushforward 
along the graph $\CZ\to \Ran\times \CZ$ of the original map $\CZ\to \Ran$.

\end{itemize}

\sssec{}

For a general pseudo-proper $\CZ$, the functor $\sF_{\int_\CZ}$
factors as
$$\bC^{\on{loc}}_\CZ\to \bC^{\on{loc}}_\Ran \overset{\sF}\to \bC^{\on{glob}},$$
where the first arrow is the functor of !-puhsforward (see \corref{c:!-dir image ps-proper}).

\medskip

For $\CZ=\on{pt}$ and the map $\CZ\to \Ran$ given by $\ul{x}\in \Ran$, we obtain the
category denoted $\bC^{\on{loc}}_{\ul{x}}$ and a functor
$$\sF_{\ul{x}}:\bC^{\on{loc}}_{\ul{x}}\to \bC^{\on{glob}}.$$

\sssec{} \label{sss:ex local to global}

The main examples of the above are when $\ul\bC^{\on{loc}}$ is one of the 
factorization categories 
$$\Whit^!(G),\,\, \Whit_*(G),\,\,\KL(G)_\kappa.$$

In each of these cases, the corresponding global category is
$$\Dmod_{\frac{1}{2}}(\Bun_G),\,\, \Dmod_{\frac{1}{2},\on{co}}(\Bun_G), \text{ and } \Dmod_\kappa(\Bun_G),$$
and the functor $\sF_\Ran$ is 
$$\Poinc_{G,!,\Ran},\,\, \Poinc_{G,*,\Ran} \text{ and }\Loc_{G,\kappa,\Ran},$$
respectively.

\sssec{} \label{sss:category of local-to-global functors}

For given $\ul\bC^{\on{loc}}$ and $\bC^{\on{glob}}$, we can consider the totality of functors $\ul\sF$ as
above as a category, denoted
$$\on{Funct}^{\on{loc}\to \on{glob}}(\ul\bC^{\on{loc}},\bC^{\on{glob}}).$$

By \secref{sss:F recovers uF}, this is the same as just the category 
$$\on{Funct}_{\on{cont}}(\bC^{\on{loc}}_\Ran,\bC^{\on{glob}}).$$

\ssec{The local unital structure} \label{ss:unitality loc} 

\sssec{} \label{sss:unitality loc} 

Let $\ul\bC^{\on{loc}}$ be a crystal of categories over $\Ran$. By a local unital structure on $\ul\bC^{\on{loc}}$:
we mean an extension $\ul\bC^{\on{loc,untl}}$ of $\ul\bC^{\on{loc}}$ to a crystal of categories over
$\Ran^{\on{untl}}$ (see \secref{sss:shvs-of-cats categ} for what this means). 

\medskip

An example of such a structure is provided by a unital lax factorization category. 

\medskip

Let us explain what the unital structure means in concrete terms.  

\sssec{} \label{sss:Ran subset}

Let $\Ran^{\subseteq}$ be the moduli space
of pairs 
$$(\ul{x},\ul{x}'\,|\,\ul{x}\subseteq \ul{x}'),$$
see \secref{sss:Ran subseteq}. 

\medskip

We have the maps
$$\on{pr}_{\on{small}},\on{pr}_{\on{big}}:\Ran^{\subseteq}\rightrightarrows \Ran$$
that remember $\ul{x}$ and $\ul{x}'$, respectively. 

\medskip

Let $\on{diag}$ denote the diagonal map
$$\Ran\to \Ran^{\subseteq}.$$

Note that
$$\on{pr}_{\on{small}}\circ \on{diag}\simeq \on{Id} \simeq \on{pr}_{\on{diag}}\circ \on{diag}.$$

\sssec{}

Denote
$$\Ran^{\subseteq^2}:=\Ran^{\subseteq}\underset{\on{pr}_{\on{small}},\Ran,\on{pr}_{\on{big}}}\times \Ran^{\subseteq}.$$

In addition to the two projections
$$\on{pr}_{\on{small}^2},\on{pr}_{\on{big}^2}:\Ran^{\subseteq^2}\rightrightarrows \Ran^{\subseteq},$$
we have a map 
$$\on{pr}_{\on{comp}}:\Ran^{\subseteq^2}\to \Ran^{\subseteq}$$
that sends
$$(\ul{x},\ul{x}',\ul{x}''\,|\,\ul{x}\subseteq \ul{x}'\subseteq \ul{x}'')\mapsto (\ul{x},\ul{x}'').$$

\sssec{} 

Note that $\Ran^{\subseteq}$ is the prestack of morphisms of $\Ran^{\on{untl}}$. Hence, 
at the level of 1-morphisms, an extension of $\ul\bC^{\on{loc}}$ to $\ul\bC^{\on{loc,untl}}$ amounts
to a functor
$$(\on{pr}_{\on{small}})^*(\ul\bC^{\on{loc}})\to (\on{pr}_{\on{big}})^*(\ul\bC^{\on{loc}})$$
as crystals of categories over $\Ran^{\subseteq}$, or equivalently, to a functor 
\begin{equation} \label{e:ins unit pre}
\ul\bC^{\on{loc}}\to (\on{pr}_{\on{small}})_*\circ (\on{pr}_{\on{big}})^*(\ul\bC^{\on{loc}})
\end{equation} 
as crystals of categories over $\Ran$. 

\medskip

In the above formula: 
\begin{itemize}

\item $(\on{pr}_{\on{big}})^*$ (resp., $(\on{pr}_{\on{small}})^*$) 
is the functor of pullback along $\on{pr}_{\on{big}}$ (resp., $\on{pr}_{\on{small}}$) 
from crystals of categories over $\Ran$ to crystals of categories over $\Ran^{\subseteq}$;

\smallskip

\item $(\on{pr}_{\on{small}})_*$ is the functor of pushforward along $\on{pr}_{\on{small}}$ from sheaves of
categories over $\Ran^{\subseteq}$ to crystals of categories over $\Ran$.

\end{itemize} 

We refer the reader to \secref{ss:dir mage sheaf of cat}, where the operation of pushforward
for crystals of categories is reviewed. 

\sssec{} \label{sss:unital insertion}

Denote the functor \eqref{e:ins unit pre} by\footnote{In the formula below $\on{ins.unit}$ is the abbreviation
of ``insert unit".}
$$\on{ins.unit}:\ul\bC^{\on{loc}}\to (\on{pr}_{\on{small}})_*\circ (\on{pr}_{\on{big}})^*(\ul\bC^{\on{loc}}).$$

\medskip

The functor $\on{ins.unit}$ has an associativity structure explained in \secref{sss:ins unit ass}. The full
datum of the upgrade
$$\ul\bC^{\on{loc}}\rightsquigarrow \ul\bC^{\on{loc,untl}}$$
is encoded by $\on{ins.unit}$, together with the associativity structure satisfying a homotopy-coherent system
of compatibilities. 

\medskip

We will now explain the concrete meaning of the functor $\on{ins.unit}$. 

\sssec{} \label{sss:diag Z}

%

For $\CZ\to \Ran$, denote 
$$\CZ^{\subseteq}:=\CZ\underset{\Ran}\times \Ran^{\subseteq},$$ where in the formation of the fiber product the map 
$\Ran^{\subseteq}\to \Ran$ is $\on{pr}_{\on{small}}$, see \secref{sss:Z Ran}. 

\medskip

Denote by $\on{pr}_{\on{small},\CZ}$ the map
$$\CZ^{\subseteq}\to \CZ,$$
and by $\on{pr}_{\on{big}}$ the projection 
$$\CZ^{\subseteq}\to \Ran^{\subseteq} \overset{\on{pr}_{\on{big}}}\longrightarrow \Ran.$$

We view $\CZ^{\subseteq}$ as mapping to $\Ran$ via $\on{pr}_{\on{big}}$.
The map $\on{diag}$ induces a map 
$$\on{diag}_\CZ:\CZ\to \CZ^{\subseteq}.$$

\sssec{}

The map $\on{diag_\CZ}$ gives rise to a functor
$$(\on{diag}_\CZ)^!: \bC^{\on{loc}}_{\CZ^{\subseteq}}\to \bC^{\on{loc}}_{\CZ}.$$

Since $\on{diag_\CZ}$ is pseudo-proper, the functor $(\on{diag}_\CZ)^!$ admits a left 
adjoint, to be denoted $(\on{diag}_\CZ)_!$ (see \corref{c:!-dir image ps-proper}). Thus, we have an adjoint pair: 
\begin{equation} \label{e:diag adj Z}
(\on{diag}_\CZ)_!:\bC^{\on{loc}}_{\CZ}\rightleftarrows \bC^{\on{loc}}_{\CZ^{\subseteq}}:(\on{diag}_\CZ)^!.
\end{equation} 

\sssec{}

The functor $\on{ins.unit}$ assigns to $\CZ$ a $\Dmod(\CZ)$-linear functor
\begin{equation} \label{e:insert unit Z}
\on{ins.unit}_\CZ:\bC^{\on{loc}}_\CZ \to \bC^{\on{loc}}_{\CZ^{\subseteq}}.
\end{equation} 

Note also that 
\begin{equation} \label{e:insert unit diag Z}
(\on{diag}_\CZ)^!\circ \on{ins.unit}_\CZ\simeq \on{Id}
\end{equation} 
as endofunctors of $\bC^{\on{loc}}_\CZ$. 

\sssec{} \label{sss:ins unit ass}

Denote
$$\CZ^{\subseteq^2}:=(\CZ^{\subseteq})^{\subseteq}\simeq 
\CZ\underset{\Ran,\on{pr}_{\on{small}^2}}\times \Ran^{\subseteq^2}\simeq 
\CZ^{\subseteq}\underset{\Ran,\on{pr}_{\on{small}}}\times \Ran^{\subseteq}.$$

We view $\CZ^{\subseteq^2}$ as mapping to $\Ran$ via $\on{pr}_{\on{big}^2}$.
The map $\on{pr}_{\on{comp}}:\Ran^{\subseteq^2}\to \Ran^{\subseteq}$ gives rise
to a map 
$$\on{pr}_{\on{comp},\CZ}:\CZ^{\subseteq^2}\to \CZ^{\subseteq}$$
as spaces over $\Ran$.

\medskip

The associativity property of the functor $\on{ins.unit}$ is encoded by the following
diagram
$$
\CD
\bC^{\on{loc}}_\CZ @>{\on{ins.unit}_\CZ}>> \bC^{\on{loc}}_{\CZ^{\subseteq}} \\
@V{\on{ins.unit}_\CZ}VV @VV{\on{pr}_{\on{comp},\CZ}^!}V \\
\bC^{\on{loc}}_{\CZ^{\subseteq}} @>>{\on{ins.unit}_{\CZ^{\subseteq}}}> \bC^{\on{loc}}_{\CZ^{\subseteq^2}}.
\endCD
$$

\sssec{Example}

At the pointwise level, the datum of \eqref{e:insert unit Z} is a system of functors
$$\on{ins.unit}_{\ul{x}_1\subseteq \ul{x}_2}:\bC^{\on{loc}}_{\ul{x}_1}\to \bC^{\on{loc}}_{\ul{x}_2} \text{ for } \ul{x}_1\subseteq \ul{x}_2.$$

When 
$$\ul{x}_2=\ul{x}_1\sqcup \ul{x}',$$
and $\ul\bC^{\on{loc}}$ is a unital lax factorization category $\bC$, the above functor is
$$\bC_{\ul{x}_1}\overset{\on{Id}\otimes \one_{\bC,\ul{x}'}}\longrightarrow 
 \bC_{\ul{x}_1}\otimes \bC_{\ul{x}'}\to \bC_{\ul{x}_2},$$
where the last arrow is given by the lax factorization structure.

\sssec{} \label{sss:ins Z int}

Assume for a moment that $\CZ$ is pseudo-proper. In this case the map
$\on{pr}_{\on{big}}:\CZ^{\subseteq}\to \Ran$
is pseudo-proper, and hence the functor
$$(\on{pr}_{\on{big}})_!:\bC^{\on{loc}}_{\CZ^{\subseteq}}\to \bC^{\on{loc}}_\Ran$$
left adjoint to $(\on{pr}_{\on{big}})^!$ is defined (see \corref{c:!-dir image ps-proper}).

\medskip

We will consider the functor
$$\int_\CZ\on{ins.unit}:\bC^{\on{loc}}_\CZ \to \bC^{\on{loc}}_\Ran$$
equal to the composition
$$\bC^{\on{loc}}_\CZ \overset{\on{ins.unit}_\CZ}\longrightarrow \bC^{\on{loc}}_{\CZ^{\subseteq}}
\overset{(\on{pr}_{\on{big}})_!}\longrightarrow \bC^{\on{loc}}_\Ran.$$

\sssec{} \label{sss:ins Ran int}

In particular, we obtain an endofunctor
\begin{equation} \label{e:ins Ran}
\int_\Ran\on{ins.unit}:\bC^{\on{loc}}_\Ran\to \bC^{\on{loc}}_\Ran.
\end{equation}

Note that the adjunction \eqref{e:insert unit Z} and the identification \eqref{e:insert unit diag Z} give rise
to a natural transformation 
\begin{equation} \label{e:ins unit int Ran}
\on{Id}\to \int_\Ran\on{ins.unit}
\end{equation}
as endofunctors of $\bC^{\on{loc}}_\Ran$. Indeed, \eqref{e:ins unit int Ran} is given by
\begin{multline*}
\on{Id}\simeq (\on{pr}_{\on{big}})_!\circ (\on{diag}_\Ran)_!  \simeq 
(\on{pr}_{\on{big}})_!\circ (\on{diag}_\Ran)_! \circ (\on{diag}_\Ran)^! \circ \on{ins.unit}_\Ran \to \\
\to (\on{pr}_{\on{big}})_!\circ \on{ins.unit}_\Ran=\int_\Ran\on{ins.unit}.
\end{multline*}

\sssec{Inventory of notation}

We briefly summarize the notation related to insertion of the unit. 

\medskip

We denote by $\on{ins.unit}_\CZ$ the functor
$$\on{ins.unit}_\CZ:\bC^{\on{loc}}_\CZ \to \bC^{\on{loc}}_{\CZ^{\subseteq}}$$

\medskip

For $\CZ$ pseudo-proper, we denote by
$$\int_\CZ \on{ins.unit}: \bC^{\on{loc}}_\CZ \to \bC^{\on{loc}}_\Ran$$
the composition of $\on{ins.unit}_\CZ$ with $(\on{pr}_{\on{big}})_!$.

\medskip

In particular, for $\CZ=\Ran$, we have $\int_\Ran \on{ins.unit}$, which is an endofunctor of
$\bC^{\on{loc}}_\Ran$. 

\sssec{}

Yet, in \eqref{e:ins unit int as functor} we will introduce yet another symbol: just $\int \on{ins.unit}$.
It will be an endofunctor of the category 
$$\on{Funct}^{\on{loc}\to \on{glob}}(\ul\bC^{\on{loc}},\bC^{\on{glob}}), \quad \ul\sF\mapsto \ul\sF^{\int \on{ins.unit}}.$$
(see \secref{sss:category of local-to-global functors}), defined when $\ul\bC^{\on{loc}}$ is equipped with a local unital structure. 

\medskip

We will have
$$\sF_{\int_{\CZ^\subseteq}}\circ \on{ins.unit}_\CZ \simeq \sF'_{\int_\CZ}$$
as functors $\bC^{\on{loc}}_\CZ\to \bC^{\on{glob}}$, 
where:

\smallskip

\begin{itemize}

\item $\ul\sF':=\sF^{\int \on{ins.unit}}\in \on{Funct}^{\on{loc}\to \on{glob}}(\ul\bC^{\on{loc}},\bC^{\on{glob}})$;

\smallskip

\item The notation $\sF'_{\int_\CZ}$ is as in \secref{e:F Z int}. 

\end{itemize}

\ssec{A (lax) unital structure on a local-to-global functor}

\sssec{} \label{sss:lax unital}

Let $(\ul\bC^{\on{loc}},\bC^{\on{glob}},\ul\sF)$ be as in \secref{sss:local to global functor abs}. Assume now
that $\ul\bC^{\on{loc}}$ is equipped with a local unital structure.

\medskip

A lax unital  structure on $\ul\sF$ is its upgrade to a \emph{right-lax} functor
\begin{equation} \label{e:lax unital}
\ul\sF^{\on{untl}}:\ul\bC^{\on{loc,untl}}\to \bC^{\on{glob}}\otimes \ul\Dmod(\Ran^{\on{untl}})
\end{equation} 
between crystals of categories over $\Ran^{\on{untl}}$, see \secref{sss:lax vs strict functors} for what this means. 

\sssec{} \label{sss:unitality concrete}

Concretely, a lax unital  structure on $\ul\sF$ means the following. Let $\CZ$ be a space,
and let 
$$\ul{x}_1\overset{\alpha}\to \ul{x}_2$$
be a morphism in the category $\Maps(\CZ,\Ran^{\on{untl}})$.

\medskip

The maps $\ul{x}_i$ give rise to categories $\bC^{\on{loc}}_{\CZ,\ul{x}_i}$ tensored over $\Dmod(\CZ)$, $i=1,2$. 
The datum of $\ul\sF$ gives rise to $\Dmod(\CZ)$-linear functors
$$\sF_{\CZ,\ul{x}_i}:\bC^{\on{loc}}_{\CZ,\ul{x}_i}\to \bC^{\on{glob}}\otimes \Dmod(\CZ).$$

\medskip

The local unital structure on $\ul\bC^{\on{loc}}$ gives rise to a $\Dmod(\CZ)$-linear functor
$$\bC^{\on{loc}}_{\CZ,\ul{x}_1}\overset{\bC^{\on{loc}}_\alpha}\to \bC^{\on{loc}}_{\CZ,\ul{x}_2}.$$

Then the datum of $\ul\sF^{\on{untl}}$ gives rise to a natural transformation
\begin{equation} \label{e:lax unital nat trans}
\sF_{\CZ,\ul{x}_1}\to \sF_{\CZ,\ul{x}_2}\circ \bC^{\on{loc}}_\alpha
\end{equation} 
as functors
$$\bC^{\on{loc}}_{\CZ,\ul{x}_1}\to \bC^{\on{glob}}\otimes \Dmod(\CZ).$$

\sssec{Example} \label{sss:global unitality example} 

Set $\CZ=\on{pt}$, so that $\ul{x}_1\overset{\alpha}\to \ul{x}_2$ corresponds to an inclusion
$$\ul{x}_1\subseteq \ul{x}_2.$$

Then $\ul\sF^{\on{untl}}$ gives rise to a natural transformation
$$\sF_{\ul{x}_1}\to \sF_{\ul{x}_2}\circ \on{ins.unit}_{\ul{x}_1\subseteq \ul{x}_2}.$$

\sssec{} \label{sss:F subset Z}

We can rewrite the datum of natural transformations \eqref{e:lax unital nat trans} as follows:

\medskip

Let $\ul\sF$ be as in \secref{ss:local-to-global set up}. Evaluating $\ul\sF$ on $\CZ^{\subseteq}$, 
we obtain a $\Dmod(\CZ)$-linear functor
$$\sF_{\CZ^\subseteq}:\bC^{\on{loc}}_{\CZ^{\subseteq}}\to \bC^{\on{glob}}\otimes \Dmod(\CZ^{\subseteq}).$$

\medskip

The datum of $\ul\sF^{\on{untl}}$ gives rise to a natural transformation 
\begin{equation} \label{e:unitality transform original}
(\on{Id}\otimes (\on{pr}_{\on{small},\CZ})^!) \circ \sF_\CZ\to 
\sF_{\CZ^\subseteq} \circ \on{ins.unit}_\CZ
\end{equation}
as functors
$$\bC^{\on{loc}}_\CZ\to \bC^{\on{glob}}\otimes \Dmod(\CZ^{\subseteq}).$$

The natural transformation \eqref{e:unitality transform original} encodes the datum of $\ul\sF^{\on{untl}}$ at the level of 1-morphisms.
One can recover the full datum of $\ul\sF^{\on{untl}}$ by imposing a datum of associativity that \eqref{e:unitality transform original}
is supposed to satisfy. 

\sssec{} \label{sss:strictly unital}

We shall say that a lax unital structure on $\ul\sF$ is \emph{strict} if $\ul\sF^{\on{untl}}$ 
is a strict functor between crystals of categories over $\Ran^{\on{untl}}$, see \secref{sss:lax vs strict functors}
for what this means. 

\medskip

By definition, this means that the natural transformations \eqref{e:lax unital nat trans}
are isomorphisms. 

\medskip

In this case we will call $\ul\sF^{\on{untl}}$ a \emph{unital structure} on $\ul\sF$.

\sssec{} \label{sss:Axiom 2 for untl}

Equivalently, a lax unital structure on $\ul\sF$ is strict if the natural transformation 
\eqref{e:unitality transform original} is an isomorphism for any $\CZ$. 

\sssec{} \label{sss:unitality fam} 

Each of the examples from \secref{sss:ex local to global} has a natural unital structure. 

\sssec{} \label{sss:loc-to-glob notation}

We can consider the categories
\begin{equation} \label{e:strict untl into lax unital}
\on{Funct}^{\on{loc}\to \on{glob},\on{untl}}(\ul\bC^{\on{loc}},\bC^{\on{glob}}) \subset 
\on{Funct}^{\on{loc}\to \on{glob},\on{lax-untl}}(\ul\bC^{\on{loc}},\bC^{\on{glob}})
\end{equation}
of local-to-global functors equipped with a unital or lax unital structures,
respectively, with the former being a full subcategory of the latter.

\medskip

Note that we have a forgetful functor
\begin{equation} \label{e:forget lax unital}
\on{Funct}^{\on{loc}\to \on{glob},\on{lax-untl}}(\ul\bC^{\on{loc}},\bC^{\on{glob}})\to
\on{Funct}^{\on{loc}\to \on{glob}}(\ul\bC^{\on{loc}},\bC^{\on{glob}}),
\end{equation}
where $\on{Funct}^{\on{loc}\to \on{glob}}(\ul\bC^{\on{loc}},\bC^{\on{glob}})$ is as in \secref{sss:category of local-to-global functors}. 

\begin{rem} \label{r:unitality qft} 

In \cite[Sect. 1.2.6]{HR}, axioms for an 
algebro-geometric avatar of [1,2]-extended 3d quantum field theories on $X$ were 
considered, although a detailed definition was
not provided.
In the present setting, we can easily spell out the complete axioms: 

\medskip

We should have the data of a unital factorization category 
$\bC$ (viewed as a crystal of categories over $\Ran^{\on{untl}}$),
a local-to-global functor $\sF:\bC_\Ran \to \Vect$, 
and a unital structure on $\ul\sF$. 

\medskip

We refer to {\it loc. cit.} 
for a discussion of why these axiomatics can be geometrically interpreted in terms 
of 3d QFTs. 

\medskip

Moreover, the discussion from \cite[Sect. 1.2.13-15]{HR} suggests
that local-to-global functors
valued in more general global categories $\bC^{\on{glob}}$ 
should generally be interpreted in terms of boundary conditions for 4d QFTs; 
this applies for all the examples we consider here.

\end{rem}

\ssec{Integrated insertion of the unit} \label{ss:integrated ins}

The main construction in this subsection (i.e., the operation $\int \on{ins.unit}$) may be viewed 
as an abstraction of the definition of chiral (a.k.a. factorization) homology in \cite[Sect. 4.2]{BD2}. 

\medskip

As we shall see, the framework introduced above allows us to reproduce this construction 
automatically: it amounts to the left of adjoint to the embedding
$$\on{Funct}^{\on{loc}\to \on{glob},\on{untl}}(\ul\bC^{\on{loc}},\bC^{\on{glob}})\to 
\on{Funct}^{\on{loc}\to \on{glob},\on{lax-untl}}(\ul\bC^{\on{loc}},\bC^{\on{glob}}).$$

\medskip

See also \secref{ss:ch homology}, where the specific example of the functor of 
factorization homology is considered. 

\sssec{}

Let $\ul\bC^{\on{loc}}$ and $\bC^{\on{glob}}$ be as in \secref{ss:local-to-global set up}.
Let 
$$\on{Funct}^{\on{loc}\to \on{glob}}(\ul\bC^{\on{loc}},\bC^{\on{glob}})$$
be the corresponding category of local-to-global functors.

\medskip

Assume now that $\ul\bC^{\on{loc}}$ is equipped with a local unital structure.

\sssec{}

We define an endofunctor 
\begin{equation} \label{e:ins unit int as functor}
\int\on{ins.unit}: \on{Funct}^{\on{loc}\to \on{glob}}(\ul\bC^{\on{loc}},\bC^{\on{glob}})\to \on{Funct}^{\on{loc}\to \on{glob}}(\ul\bC^{\on{loc}},\bC^{\on{glob}})
\end{equation} 
$$\ul\sF\mapsto \ul\sF^{\int \on{ins.unit}}$$
by 
$$\sF^{\int \on{ins.unit}}_\CZ:=
(\on{Id}\otimes (\on{pr}_{\on{small},\CZ})_!)\circ \sF_{\CZ^\subseteq}\circ \on{ins.unit}_\CZ.$$

In other words, $\sF^{\int \on{ins.unit}}_\CZ$ is the composition
$$\bC^{\on{loc}}_\CZ \overset{\on{ins.unit}_\CZ}\longrightarrow \bC^{\on{loc}}_{\CZ^{\subseteq}}
\overset{\sF_{\CZ^\subseteq}}\longrightarrow \bC^{\on{glob}}\otimes \Dmod(\CZ^{\subseteq})
\overset{\on{Id}\otimes (\on{pr}_{\on{small},\CZ})_!}\longrightarrow \bC^{\on{glob}}\otimes \Dmod(\CZ).$$

\sssec{}

Note that we have a natural transformation 
\begin{equation} \label{e:correlator F}
\on{Id}\to \int\on{ins.unit} 
\end{equation} 
so that for a given $\CZ$ the corresponding map 
\begin{equation} \label{e:correlator F Z}
\sF_\CZ\to \sF^{\int \on{ins.unit}}_\CZ
\end{equation} 
is given by 
\begin{multline} \label{e:pre-unitality transform}
\sF_\CZ \simeq \left(\on{Id}\otimes (\on{pr}_{\on{small},\CZ})_!\right) \circ 
\left(\on{Id}\otimes (\on{diag}_\CZ)_!\right) \circ \sF_\CZ 
\simeq \left(\on{Id}\otimes (\on{pr}_{\on{small},\CZ})_!\right) \circ \sF_{\CZ^\subseteq} \circ (\on{diag}_\CZ)_!
\overset{\text{\eqref{e:insert unit diag Z}}}\simeq \\
\simeq \left(\on{Id}\otimes (\on{pr}_{\on{small},\CZ})_!\right) \circ \sF_{\CZ^\subseteq} \circ (\on{diag}_\CZ)_!\circ 
(\on{diag}_\CZ)^!\circ \on{ins.unit}_\CZ \to \\
\to \left(\on{Id}\otimes (\on{pr}_{\on{small},\CZ})_!\right) \circ \sF_{\CZ^\subseteq} \circ \on{ins.unit}_\CZ
=\ul\sF^{\int\on{ins.unit}}_\CZ.
\end{multline} 

\sssec{}

Assume for a moment that $\CZ$ is pseudo-proper. Applying $\on{C}_c^\cdot(\CZ,-)\otimes \on{Id}$ to 
both sides of \eqref{e:correlator F Z}, we obtain a natural transformation
\begin{equation} \label{e:ins Z int F}
\sF_{\int_\CZ} \to \sF\circ \int_\CZ\on{ins.unit},
\end{equation} 
where $\int_\CZ\on{ins.unit}$ is as in \secref{sss:ins Z int}. 

\medskip

Take $\CZ=\Ran$. In this case, the resulting natural transformation \eqref{e:ins Z int F} is
\begin{equation} \label{e:ins Ran int F}
\sF  \to \sF\circ \int_\Ran\on{ins.unit},
\end{equation} 
where $\int_\Ran\on{ins.unit}$ is as in \secref{sss:ins Ran int}. 

\medskip

It is easy to see, however, that \eqref{e:ins Ran int F} equals the natural transformation
obtained by applying $\sF$ to the natural transformation \eqref{e:ins unit int Ran}. 

\sssec{} \label{sss:verify Axiom 1 for unital}

Suppose for a moment that $\ul\sF$ is equipped with a unital structure, i.e., it is 
the image under the forgetful functor \eqref{e:forget lax unital}
of an object $\ul\sF^{\on{untl}}\in \on{Funct}^{\on{loc}\to \on{glob},\on{untl}}(\ul\bC^{\on{loc}},\bC^{\on{glob}})$. 

\medskip

We claim that in this case the map \eqref{e:correlator F Z}
is an isomorphism. Indeed, in this
case, the \emph{isomorphism} \eqref{e:unitality transform original} identifies 
\begin{equation} \label{e:ins of unital}
\sF^{\int\on{ins.unit}}_\CZ\simeq (\on{pr}_{\on{small},\CZ})_! \circ
(\on{pr}_{\on{small},\CZ})^!\circ \sF_\CZ,
\end{equation} 
and the map \eqref{e:correlator F} is the map
\begin{equation} \label{e:ins of unital receives}
\sF_\CZ \simeq (\on{pr}_{\on{small},\CZ})_!\circ (\on{diag}_\CZ)_! \circ (\on{diag}_\CZ)^! \circ
(\on{pr}_{\on{small},\CZ})^!\circ \sF_\CZ\to 
(\on{pr}_{\on{small},\CZ})_! \circ
(\on{pr}_{\on{small},\CZ})^! \circ \sF_\CZ.
\end{equation} 

Now, the contractibility of the Ran space implies that the counit of the 
$((\on{pr}_{\on{small},\CZ})_!,(\on{pr}_{\on{small},\CZ})^!)$-adjunction is an isomorphism. 
Hence, the right-hand side of \eqref{e:ins of unital receives} maps isomorphically to $\sF_\CZ$, and the
composition
$$\sF_\CZ \overset{\text{\eqref{e:ins of unital receives}}}\longrightarrow 
(\on{pr}_{\on{small},\CZ})_! \circ
(\on{pr}_{\on{small},\CZ})^! \circ \sF_\CZ\to \sF_\CZ$$
is the identity map.

\medskip

In particular, in this case the natural transformation \eqref{e:ins Ran int F} is an isomorphism. 

\sssec{} \label{sss:int untl 1}

Our next goal, carried out in Sects. \ref{ss:constr int}-\ref{ss:prop of int}, 
is to perform similar constructions with the same input, but in the unital context,
i.e., working over $\Ran^{\on{untl}}$ rather than over $\Ran$. Namely, we will show that, parallel 
to \eqref{e:ins unit int as functor}, there exists an endofunctor
\begin{equation} \label{e:ins unit int as functor untl}
\int\on{ins.unit}: \on{Funct}^{\on{loc}\to \on{glob},\on{lax-untl}}(\ul\bC^{\on{loc}},\bC^{\on{glob}})\to 
\on{Funct}^{\on{loc}\to \on{glob},\on{lax-untl}}(\ul\bC^{\on{loc}},\bC^{\on{glob}})
\end{equation} 
$$\ul\sF^{\on{untl}}\mapsto \ul\sF^{\on{untl},\int \on{ins.unit}}$$
that makes the diagram
\begin{equation} \label{e:ins unit int as functor untl compat}
\CD
\on{Funct}^{\on{loc}\to \on{glob},\on{lax-untl}}(\ul\bC^{\on{loc}},\bC^{\on{glob}}) @>{\int\on{ins.unit}}>>
\on{Funct}^{\on{loc}\to \on{glob},\on{lax-untl}}(\ul\bC^{\on{loc}},\bC^{\on{glob}})  \\
@VVV @VVV \\
\on{Funct}^{\on{loc}\to \on{glob}}(\ul\bC^{\on{loc}},\bC^{\on{glob}}) @>{\int\on{ins.unit}}>>
\on{Funct}^{\on{loc}\to \on{glob}}(\ul\bC^{\on{loc}},\bC^{\on{glob}})
\endCD
\end{equation} 
commute.

\sssec{} \label{sss:int untl 2}

In addition, the functor \eqref{e:ins unit int as functor untl}
will be equipped with a natural transformation
\begin{equation} \label{e:correlator F untl}
\on{Id}\to \int\on{ins.unit}, 
\end{equation} 
which is compatible with \eqref{e:correlator F} via  \eqref{e:ins unit int as functor untl compat}. 

\medskip

We will also show:

\begin{itemize}

\item The essential
image of \eqref{e:ins unit int as functor untl} belongs to 
$\on{Funct}^{\on{loc}\to \on{glob},\on{untl}}(\ul\bC^{\on{loc}},\bC^{\on{glob}})$. 

\smallskip

\item The natural transformation \eqref{e:correlator F untl} evaluates to an isomorphism
on objects that belong to $\on{Funct}^{\on{loc}\to \on{glob},\on{untl}}(\ul\bC^{\on{loc}},\bC^{\on{glob}})$;

\item The two natural transformations
\begin{equation} \label{e:double insert}
\int\on{ins.unit}\rightrightarrows \int\on{ins.unit} \circ \int\on{ins.unit},
\end{equation}
arising from \eqref{e:correlator F untl} coincide (it follows that they are isomorphisms). 

\end{itemize}

\sssec{} \label{sss:ins left adj}

The above properties combined imply that the functor \eqref{e:ins unit int as functor untl} is the left adjoint of the 
embedding
$$\on{Funct}^{\on{loc}\to \on{glob},\on{untl}}(\ul\bC^{\on{loc}},\bC^{\on{glob}}) \hookrightarrow
\on{Funct}^{\on{loc}\to \on{glob},\on{lax-untl}}(\ul\bC^{\on{loc}},\bC^{\on{glob}}).$$

\begin{rem}

The reason for discussing both versions of $\int\on{ins.unit}$, i.e., \eqref{e:ins unit int as functor untl}
and \eqref{e:ins unit int as functor}, is that the former has a clear categorical meaning
(i.e., it is the left adjoint of the forgetful functor), while the latter is easily computable (a priori,
the functor \eqref{e:ins unit int as functor untl} involves taking cohomology over categorical
prestacks). 

\medskip

However, the
commutation of \eqref{e:ins unit int as functor untl compat} implies that \eqref{e:ins unit int as functor untl}
is computable as well.

\end{rem}

\ssec{Construction of the integrated functor} \label{ss:constr int}

This and the next subsection are devoted to the construction of the functor \eqref{e:ins unit int as functor untl}
and the verification of its properties. 
The reader who is willing to take this on faith may choose
to skip these two subsections. 

\medskip

We are going to present the construction of the functor \eqref{e:ins unit int as functor untl}
in a hands-on manner. See, however, \secref{sss:left adj lax to strict} for its abstract interpretation. 

\sssec{} \label{sss:arrows in cat prestack}

Let $\CY$ be a categorical prestack and let $\CY^\to$ be the categorical prestack of 1-morphisms in $\CY$.
I.e., for an affine scheme $S$, objects of $\Maps(S,\CY^\to)$ are 
$$y_1,y_2\in \Maps(S,\CY), \,\,y_1\to y_2,$$
and morphisms are commutative diagrams
$$
\CD
y_1 @>>> y_2 \\
@VVV @VVV \\
y'_1 @>>> y'_2.
\endCD
$$

We have the projections
$$\on{pr}_{\on{source}},\on{pr}_{\on{target}}:\CY^\to \to \CY,$$
with $\on{pr}_{\on{source}}$ being a Cartesian fibration.

\medskip

Let $\ul\bC$ be a crystal of categories  on $\CY$. Tautologically, we have a (strict) functor
\begin{equation} \label{e:source to target}
(\on{pr}_{\on{source}})^*(\ul\bC)\to (\on{pr}_{\on{target}})^*(\ul\bC),
\end{equation}
as crystals of categories on $\CY^\to$.

\medskip

Recall the construction of the direct image of a crystal of categories , reviewed in \secref{ss:dir mage sheaf of cat}. 
According to \secref{sss:push-pull categorical prestacks}, we have a (strict) functor
$$\ul\bC \to (\on{pr}_{\on{source}})_{*,\on{strict}}\circ (\on{pr}_{\on{source}})^*(\ul\bC).$$

Composing with \eqref{e:source to target} we obtain a (strict) functor
\begin{equation} \label{e:source target pullpush}
\ul\bC \to (\on{pr}_{\on{source}})_{*,\on{strict}}\circ (\on{pr}_{\on{target}})^*(\ul\bC).
\end{equation}

\sssec{}

We apply the construction in \secref{sss:arrows in cat prestack} to $\CY=\Ran^{\on{untl}}$. Denote
$$\Ran^{\on{untl},\to}=:\Ran^{\subseteq,\on{untl}},$$
viewed as a categorical prestacks.

\medskip

Note that the prestack in groupoids underlying $\Ran^{\subseteq,\on{untl}}$ is the prestack $\Ran^{\subseteq}$
introduced in \secref{sss:Ran subset}. We will use the symbols $\on{pr}^{\on{untl}}_{\on{small}}$ and 
$\on{pr}^{\on{untl}}_{\on{big}}$ for the corresponding maps $\on{pr}_{\on{source}}$ and $\on{pr}_{\on{target}}$. 

\medskip

Thus, for $\ul\bC^{\on{loc,untl}}$ as in \secref{ss:unitality loc}, the functor \eqref{e:source target pullpush} is a functor
\begin{equation} \label{e:source target pullpush Ran}
\ul\bC^{\on{loc,untl}} \to (\on{pr}^{\on{untl}}_{\on{small}})_{*,\on{strict}}\circ (\on{pr}^{\on{untl}}_{\on{big}})^*(\ul\bC^{\on{loc,untl}}).
\end{equation} 

\sssec{}

Let now $\ul\sF^{\on{untl}}$ be an object of $\on{Funct}^{\on{loc}\to \on{glob},\on{lax-untl}}(\ul\bC^{\on{loc}},\bC^{\on{glob}})$.
Applying pullback along $\on{pr}^{\on{untl}}_{\on{big}}$, we obtain a \emph{lax} functor
$$(\on{pr}^{\on{untl}}_{\on{big}})^*(\ul\bC^{\on{loc,untl}}) 
\overset{(\on{pr}^{\on{untl}}_{\on{big}})^*(\ul\sF^{\on{untl}})}\longrightarrow
(\on{pr}^{\on{untl}}_{\on{big}})^*(\bC^{\on{glob}} \otimes \ul{\Dmod}(\Ran^{\on{untl}}))$$
as crystals of categories on $\Ran^{\subseteq,\on{untl}}$. 

\medskip

Using 
\secref{sss:pushforward sheaf of cat funct} we obtain a lax functor of crystal of categories  on $\Ran^{\on{untl}}$
$$(\on{pr}^{\on{untl}}_{\on{small}})_{*,\on{lax}}\circ (\on{pr}^{\on{untl}}_{\on{big}})^*(\ul\bC^{\on{loc,untl}}) 
\overset{(\on{pr}^{\on{untl}}_{\on{small}})_*\circ (\on{pr}^{\on{untl}}_{\on{big}})^*(\ul\sF^{\on{untl}})}\longrightarrow
(\on{pr}^{\on{untl}}_{\on{small}})_{*,\on{lax}}\circ (\on{pr}^{\on{untl}}_{\on{big}})^*(\bC^{\on{glob}} \otimes \ul{\Dmod}(\Ran^{\on{untl}})).$$

Combining, we obtain a functor  
\begin{multline} \label{e:ins unit untl 1}
\ul\bC^{\on{loc,untl}}\to (\on{pr}^{\on{untl}}_{\on{small}})_{*,\on{strict}}\circ (\on{pr}^{\on{untl}}_{\on{big}})^*(\ul\bC^{\on{loc,untl}}) 
\to (\on{pr}^{\on{untl}}_{\on{small}})_{*,\on{lax}}\circ (\on{pr}^{\on{untl}}_{\on{big}})^*(\ul\bC^{\on{loc,untl}}) \to \\
\to (\on{pr}^{\on{untl}}_{\on{small}})_{*,\on{lax}}\circ (\on{pr}^{\on{untl}}_{\on{big}})^*(\bC^{\on{glob}} \otimes \ul{\Dmod}(\Ran^{\on{untl}})) 
\simeq \bC^{\on{glob}} \otimes (\on{pr}^{\on{untl}}_{\on{small}})_{*,\on{lax}}(\ul\Dmod(\Ran^{\subseteq,\on{untl}})).
\end{multline}


\sssec{}

Consider the (strict) functor
\begin{equation} \label{e:pullback untl from Ran to arr}
\ul\Dmod(\Ran^{\on{untl}})\overset{\text{\eqref{e:push-pull categorical prestacks}}}\longrightarrow 
(\on{pr}^{\on{untl}}_{\on{small}})_{*,\on{strict}}(\ul\Dmod(\Ran^{\subseteq,\on{untl}}))\to 
(\on{pr}^{\on{untl}}_{\on{small}})_{*,\on{lax}}(\ul\Dmod(\Ran^{\subseteq,\on{untl}})).
\end{equation}

\begin{lem} \label{l:pr small untl}
The functor \eqref{e:pullback untl from Ran to arr} admits a value-wise left adjoint, to be denoted
$(\on{pr}^{\on{untl}}_{\on{small}})_!$. Moreover, this value-wise left adjoint, which is a priori a \emph{left-lax}
functor, is strict. 
\end{lem}

The proof will be given in \secref{sss:proof pr small untl}. 

\sssec{} 

Thus, composing the lax functor \eqref{e:ins unit untl 1} with $(\on{pr}^{\on{untl}}_{\on{small}})_!$, we obtain a
lax functor
\begin{equation} \label{e:ins unit untl}
\ul\bC^{\on{loc,untl}}\to \bC^{\on{glob}} \otimes (\on{pr}^{\on{untl}}_{\on{small}})_{*,\on{lax}}(\ul\Dmod(\Ran^{\subseteq,\on{untl}}))
\overset{\on{Id}\otimes (\on{pr}^{\on{untl}}_{\on{small}})_!}\longrightarrow \bC^{\on{glob}} \otimes \ul\Dmod(\Ran).
\end{equation}


\sssec{}

The functor \eqref{e:ins unit untl} is the sought-for object
$$\ul\sF^{\on{untl},\int \on{ins.unit}}.$$

\ssec{Properties of the integrated functor} \label{ss:prop of int}

We now proceed to establish the properties of the functor \eqref{e:ins unit int as functor untl}. 

\sssec{} \label{sss:functoriality integral}

Note that the construction in \secref{sss:arrows in cat prestack} is functorial in the following sense:

\medskip

If $\ul\Phi:\ul\bC_1^{\on{loc,untl}}\to \ul\bC_2^{\on{loc,untl}}$ is a \emph{strict} functor between sheaves of categories
on $\Ran^{\on{untl}}$, and 
$$\ul\sF^{\on{untl}}:\ul\bC_2^{\on{loc,untl}}\to \ul\Dmod(\Ran^{\on{untl}})$$
is a lax unital  functor, then we have a (tautological) isomorphism
$$(\ul\sF^{\on{untl}}\circ \Phi)^{\int\on{ins.unit}}\simeq \ul\sF^{\on{untl},\int\on{ins.unit}}\circ \Phi.$$

\sssec{}

We first show that
$\int \on{ins.unit}$ acts as identity on objects 
$\ul\sF^{\on{untl}}\in \on{Funct}^{\on{loc}\to \on{glob},\on{untl}}(\ul\bC^{\on{loc}},\bC^{\on{glob}})$. 

\medskip

Indeed, if the functor
$$\ul\sF^{\on{untl}}:\ul\bC^{\on{loc,untl}}\to \bC^{\on{glob}}\otimes \ul\Dmod(\Ran^{\on{untl}})$$
is strict, then by \secref{sss:functoriality integral} above, the functor
\eqref{e:ins unit untl 1} identifies with the composition of $\ul\sF^{\on{untl}}$ with the tensor product
of the identity endofunctor of $\bC^{\on{glob}}$ with
\begin{multline} \label{e:endo Dmod Ran}
\ul\Dmod(\Ran^{\on{untl}}) \overset{\text{\eqref{e:source target pullpush}}}\longrightarrow
(\on{pr}^{\on{untl}}_{\on{small}})_{*,\on{strict}}\circ (\on{pr}^{\on{untl}}_{\on{big}})^*(\ul\Dmod(\Ran^{\on{untl}})) = \\
= (\on{pr}^{\on{untl}}_{\on{small}})_{*,\on{strict}}(\ul\Dmod(\Ran^{\subseteq,\on{untl}}))\overset{(\on{pr}^{\on{untl}}_{\on{small}})_!}\longrightarrow
\ul\Dmod(\Ran^{\on{untl}}).
\end{multline}

We claim that \eqref{e:endo Dmod Ran} is the identity endofunctor of $\ul\Dmod(\Ran^{\on{untl}})$. Indeed, observe that for
$\ul\bC=\ul\Dmod(\Ran^{\on{untl}})$, the functor \eqref{e:source to target} is the identity endofunctor of
$\ul\Dmod(\Ran^{\subseteq,\on{untl}})$, so the composition
\begin{multline*} 
\ul\Dmod(\Ran^{\on{untl}}) \overset{\text{\eqref{e:source target pullpush}}}\longrightarrow
(\on{pr}^{\on{untl}}_{\on{small}})_{*,\on{strict}}\circ (\on{pr}^{\on{untl}}_{\on{big}})^*(\ul\Dmod(\Ran^{\on{untl}})) = \\
= (\on{pr}^{\on{untl}}_{\on{small}})_{*,\on{strict}}(\ul\Dmod(\Ran^{\subseteq,\on{untl}}))
\end{multline*}
is the functor
\begin{multline*} 
\ul\Dmod(\Ran^{\on{untl}}) \overset{\text{\eqref{e:source target pullpush}}}\longrightarrow
(\on{pr}^{\on{untl}}_{\on{small}})_{*,\on{strict}}\circ (\on{pr}^{\on{untl}}_{\on{small}})^*(\ul\Dmod(\Ran^{\on{untl}})) = \\
= (\on{pr}^{\on{untl}}_{\on{small}})_{*,\on{strict}}(\ul\Dmod(\Ran^{\subseteq,\on{untl}})).
\end{multline*}

The assertion follows now from the next lemma:

\begin{lem} \label{l:untl contr}
The functor \eqref{e:pullback untl from Ran to arr} is fully faithful.\footnote{See \secref{sss:ff strict functors} for what this means.}
\end{lem}

The proof will be given in \secref{sss:proof untl contr}.

\sssec{} \label{sss:ins unitality}

Next we show that the essential image of the functor \eqref{e:ins unit int as functor untl} lies in
the subcategory $\on{Funct}^{\on{loc}\to \on{glob},\on{untl}}(\ul\bC^{\on{loc}},\bC^{\on{glob}})$, i.e.,
$\ul\sF^{\on{untl},\int \on{ins.unit}}$ is strict. 

\medskip

Let $S$ be an affine scheme and let 
$\ul{x}_1 \overset{\alpha}\to \ul{x}_2$ be a 1-morphism in $\Maps(S,\Ran^{\on{untl}})$. 
Consider the fiber products
$$S_{\ul{x}_i}^{\subseteq,\on{untl}} :=S\underset{\ul{x}_i,\Ran^{\on{untl}}}\times \Ran^{\subseteq,\on{intl}}, \quad i=1,2,$$
and the resulting map
\begin{equation} \label{e:pullback fibers}
S_{\ul{x}_2}^{\subseteq,\on{untl}} \overset{\alpha^*}\to S_{\ul{x}_1}^{\subseteq,\on{untl}}.
\end{equation}

Consider the diagram
$$
\xy  
(0,0)*+{\ul\bC^{\on{loc,untl}}_{S,\ul{x}_1}}="A";
(0,-30)*+{\ul\bC^{\on{loc,untl}}_{S,\ul{x}_2}}="A'";
(45,0)*+{\Gamma^{\on{lax}}(S_{\ul{x}_1}^{\subseteq,\on{untl}},\ul\bC^{\on{loc,untl}})}="B";
(45,-30)*+{\Gamma^{\on{lax}}(S_{\ul{x}_2}^{\subseteq,\on{untl}},\ul\bC^{\on{loc,untl}})}="B'";
(110,0)*+{\bC^{\on{glob}}\otimes \Dmod(S_{\ul{x}_1}^{\subseteq,\on{untl}})}="C";
(110,-30)*+{\bC^{\on{glob}}\otimes \Dmod(S_{\ul{x}_2}^{\subseteq,\on{untl}}),}="C'";
{\ar@{->} "A";"B"};
{\ar@{->}^{\ul\sF^{\on{untl}}} "B";"C"};
{\ar@{->} "A'";"B'"};
{\ar@{->}^{\ul\sF^{\on{untl}}} "B'";"C'"};
{\ar@{->} "A";"A'"};
{\ar@{->} "B";"B'"};
{\ar@{->}^{\on{Id}\otimes (\alpha^*)^!} "C";"C'"};
{\ar@{=>} "C";"B'"};
\endxy
$$
where:

\begin{itemize}

\item The left vertical arrow is given by the structure of crystal of categories on $\ul\bC^{\on{loc,untl}}$;

\item The middle vertical arrow is \eqref{e:moving between fibers};

\item Both left horizontal arrows are \eqref{e:source target pullpush}; 

\item The left square commutes because the functor \eqref{e:push-pull categorical prestacks}
is strict; 

\item The natural transformation in the right square is given by the lax functor 
structure on $\ul\sF^{\on{untl}}$.

\end{itemize} 

\medskip

Given \lemref{l:pr small untl}, it suffices to show that the above natural transformation
is an isomorphism. 

\medskip

Note, however, that the middle vertical arrow in the above diagram is the functor
\eqref{e:pullback on lax sections} corresponding to the morphism $\alpha^*$ of
\eqref{e:pullback fibers}. This makes the assertion manifest.

\sssec{}

We now construct the natural transformation \eqref{e:correlator F untl}. This is done in the same way as in 
\eqref{e:pre-unitality transform} using the map
$$\on{diag}^{\on{untl}}:\Ran^{\on{untl}}\to \Ran^{\subseteq,\on{untl}},$$
corresponding to the identity morphisms in $\Ran^{\on{untl}}$ as a categorical prestack. 

\sssec{}  \label{sss:double insert}

We now show that the two maps in \eqref{e:double insert} coincide. This amounts to the following
assertion. 

\medskip

Let $S$ be an affine scheme equipped with a map to $\Ran^{\on{untl}}$. Consider 
$$S^{\subseteq,\on{untl}} :=S\underset{\Ran^{\on{untl}},\on{pr}^{\on{untl}}_{\on{small}}}\times \Ran^{\subseteq,\on{intl}} \text{ and } 
S^{\subseteq^2,\on{untl}}:=S^{\subseteq,\on{untl}} 
\underset{\on{pr}^{\on{untl}}_{\on{big}},\Ran^{\on{untl}},\on{pr}^{\on{untl}}_{\on{small}}}\times \Ran^{\subseteq,\on{intl}}.$$

We have the naturally defined maps
$$\on{pr}^{\on{untl}}_{\on{big}^2}: S^{\subseteq^2,\on{untl}}\to \Ran^{\subseteq,\on{intl}} \text{ and }
\on{pr}^{\on{untl}}_{\on{small}^2,S}: S^{\subseteq^2,\on{untl}}\to S^{\subseteq,\on{untl}}.$$

We also have the maps
$$\on{diag}^{\on{untl}}_{S,\on{big}}  \text{ and } \on{diag}^{\on{untl}}_{S,\on{small}} , \quad 
S^{\subseteq,\on{untl}}\rightrightarrows S^{\subseteq^2,\on{untl}}$$
so that
$$\on{pr}^{\on{untl}}_{\on{small}^2,S}\circ \on{diag}^{\on{untl}}_{S,\on{big}}  =\on{Id}, \quad
\on{pr}^{\on{untl}}_{\on{small}^2,S}\circ \on{diag}^{\on{untl}}_{S,\on{small}}  = 
\on{diag}^{\on{untl}}_{S}\circ \on{pr}^{\on{untl}}_{\on{small},S}$$
and
$$
\on{pr}^{\on{untl}}_{\on{big}^2}\circ \on{diag}^{\on{untl}}_{S,\on{small}}=\on{pr}^{\on{untl}}_{\on{big}},\quad
\on{pr}^{\on{untl}}_{\on{big}^2}\circ \on{diag}^{\on{untl}}_{S,\on{big}} =\on{pr}^{\on{untl}}_{\on{big}}.$$

\sssec{}

Let us explain explicitly what these maps are when $S=\on{pt}$, and the map $S\to \Ran^{\on{untl}}$ corresponds 
to a point $\ul{x}\in \Ran$.

\medskip

The categorical prestacks classify $S^{\subseteq,\on{untl}}$ and $S^{\subseteq^2,\on{untl}}$
$$(\ul{x}\subseteq \ul{x}_1) \text{ and } (\ul{x}\subseteq \ul{x}_1 \subseteq \ul{x}_2)$$
respectively, with the morphisms given by inclusions of the $\ul{x}_1$'s and $\ul{x}_2$. 

\medskip

The map $\on{pr}^{\on{untl}}_{\on{small}^2,S}$ sends 
$$(\ul{x}\subseteq \ul{x}_1 \subseteq \ul{x}_2) \mapsto (\ul{x}\subseteq \ul{x}_1);$$
the map $\on{pr}^{\on{untl}}_{\on{big}^2}$ sends 
$$(\ul{x}\subseteq \ul{x}_1 \subseteq \ul{x}_2) \mapsto \ul{x}_2;$$
the map $\on{diag}^{\on{untl}}_{S,\on{small}}$ sends
$$(\ul{x}\subseteq \ul{x}_1) \mapsto (\ul{x}\subseteq \ul{x}\subseteq \ul{x}_1);$$
the map $\on{diag}^{\on{untl}}_{S,\on{big}}$ sends
$$(\ul{x}\subseteq \ul{x}_1) \mapsto (\ul{x}\subseteq \ul{x}_1\subseteq \ul{x}_1).$$

\sssec{}

We obtain the natural transformations
\begin{multline} \label{e:double insert 1}
(\on{pr}^{\on{untl}}_{\on{small},S})_!\circ 
(\on{pr}^{\on{untl}}_{\on{big}})^!\simeq
(\on{pr}^{\on{untl}}_{\on{small},S})_!\circ (\on{pr}^{\on{untl}}_{\on{small}^2,S})_!\circ (\on{diag}^{\on{untl}}_{S,\on{big}})_!\circ 
(\on{pr}^{\on{untl}}_{\on{big}})^!\simeq \\
\simeq (\on{pr}^{\on{untl}}_{\on{small},S})_!\circ 
(\on{pr}^{\on{untl}}_{\on{small}^2,S})_!\circ (\on{diag}^{\on{untl}}_{S,\on{big}})_!\circ 
(\on{diag}^{\on{untl}}_{S,\on{big}})^!\circ (\on{pr}^{\on{untl}}_{\on{big}^2})^!\to 
(\on{pr}^{\on{untl}}_{\on{small},S})_!\circ (\on{pr}^{\on{untl}}_{\on{small}^2,S})_!\circ (\on{pr}^{\on{untl}}_{\on{big}^2})^!
\end{multline}
and 
\begin{multline} \label{e:double insert 2}
(\on{pr}^{\on{untl}}_{\on{small},S})_!\circ (\on{pr}^{\on{untl}}_{\on{big}})^!\simeq
(\on{pr}^{\on{untl}}_{\on{small},S})_! \circ (\on{diag}^{\on{untl}}_{S})_! \circ (\on{pr}^{\on{untl}}_{\on{small},S})_! \circ (\on{pr}^{\on{untl}}_{\on{big}})_!\simeq \\
\simeq (\on{pr}^{\on{untl}}_{\on{small},S})_! \circ (\on{pr}^{\on{untl}}_{\on{small}^2,S})_!\circ 
(\on{diag}^{\on{untl}}_{S,\on{small}} )_! \circ (\on{diag}^{\on{untl}}_{S,\on{small}} )^! \circ (\on{pr}^{\on{untl}}_{\on{big}^2})^! \to
(\on{pr}^{\on{untl}}_{\on{small},S})_! \circ (\on{pr}^{\on{untl}}_{\on{small}^2,S})_!\circ (\on{pr}^{\on{untl}}_{\on{big}^2})^!
\end{multline}
as functors
$$\Dmod(\Ran^{\on{untl}})\to \Dmod(S).$$

We need to show that the natural transformations \eqref{e:double insert 1} and \eqref{e:double insert 2} coincide. 

\medskip

This follows from the following observation: there exists a 1-morphism between the maps
$$\on{diag}^{\on{untl}}_{S,\on{small}} \overset{\alpha}\to \on{diag}^{\on{untl}}_{S,\on{big}},$$
so that the induced map
$$\on{pr}^{\on{untl}}_{\on{big}} \simeq 
\on{pr}^{\on{untl}}_{\on{big}^2}\circ \on{diag}^{\on{untl}}_{S,\on{small}} \overset{\alpha}\to 
\on{pr}^{\on{untl}}_{\on{big}^2}\circ \on{diag}^{\on{untl}}_{S,\on{big}} \simeq \on{pr}^{\on{untl}}_{\on{big}}$$
is the identity map. 

\sssec{}

Finally, we prove the commutativity of \eqref{e:ins unit int as functor untl compat}. 

\medskip

Let $\ul\sF^{\on{untl}}$ be an object of $\on{Funct}^{\on{loc}\to \on{glob},\on{lax-untl}}(\ul\bC^{\on{loc}},\bC^{\on{glob}})$,
and let $\ul\sF$ be the corresponding object of $\on{Funct}^{\on{loc}\to \on{glob}}(\ul\bC^{\on{loc}},\bC^{\on{glob}})$.
We need to establish an isomorphism between $\ul\sF^{\int \on{ins.unit}}$ and 
$$\ul\sF^{\on{untl},\int \on{ins.unit}}|_{\ul\bC^{\on{loc}}}:\ul\bC^{\on{loc}}\to \bC^{\on{glob}}\otimes \ul\Dmod(\Ran).$$

\medskip

Let $S$ be an affine scheme and let us be given an $S$-point of $\Ran$. Denote
$$S^{\subseteq,\on{untl}} :=S\underset{\Ran^{\on{untl}}}\times \Ran^{\subseteq,\on{intl}}
\text{ and } S^{\subseteq} :=S\underset{\Ran}\times \Ran^{\subseteq},$$
so that $S^{\subseteq}$ is the prestack in groupoids underlying $S^{\subseteq,\on{untl}}$.

\medskip

The value of $\ul\sF^{\on{untl},\int \on{ins.unit}}|_{\ul\bC^{\on{loc}}}$ at the above $S$-point of $\Ran$
is given by the composition
\begin{multline} \label{e:untl vs nonuntl int 1}
\bC^{\on{loc}}_S =\ul\bC^{\on{loc,untl}}_S \overset{\text{\eqref{e:source target pullpush}}}\longrightarrow 
\Gamma^{\on{lax}}(S^{\subseteq,\on{untl}},\ul\bC^{\on{loc,untl}}) \overset{\ul\sF^{\on{untl}}}\longrightarrow \\
\to \bC^{\on{glob}}\otimes \Dmod(S^{\subseteq,\on{untl}}) \overset{(\on{pr}^{\on{untl}}_{\on{small},S})_!}\longrightarrow 
\bC^{\on{glob}}\otimes \Dmod(S)
\end{multline}

The value of $\ul\sF^{\int \on{ins.unit}}$ at the above $S$-point of $\Ran$ is given by 
\begin{multline} \label{e:untl vs nonuntl int 2}
\bC^{\on{loc}}_S =\ul\bC^{\on{loc,untl}}_S \overset{\text{\eqref{e:source target pullpush}}}\longrightarrow 
\Gamma^{\on{lax}}(S^{\subseteq,\on{untl}},\ul\bC^{\on{loc,untl}}) \overset{\ul\sF^{\on{untl}}}\longrightarrow \\
\to \bC^{\on{glob}}\otimes \Dmod(S^{\subseteq,\on{untl}}) \to \bC^{\on{glob}}\otimes \Dmod(S^{\subseteq}) 
\overset{(\on{pr}_{\on{small},S})_!}\longrightarrow 
\bC^{\on{glob}}\otimes \Dmod(S).
\end{multline}

\medskip

Let $\sft$ denote the tautological map $S^{\subseteq} \to S^{\subseteq,\on{untl}}$. We have a natural transformation
\begin{equation}  \label{e:untl vs non-untl int}
(\on{pr}_{\on{small},S})_! \circ  \sft^!\simeq  
(\on{pr}^{\on{untl}}_{\on{small},S})_!\circ \sft_!\circ \sft^!\to (\on{pr}^{\on{untl}}_{\on{small},S})_!
\end{equation} 
as functors
$$\Dmod(S^{\subseteq,\on{untl}}) \to \Dmod(S).$$

Thus, to establish an isomorphism between \eqref{e:untl vs nonuntl int 1} and \eqref{e:untl vs nonuntl int 2},
it suffices to prove that the natural transformation \eqref{e:untl vs non-untl int} is an isomorphism.
However, this is a variant of \lemref{l:int over Ran and Ran untl} (with the same proof).

\sssec{}

The compatibility of \eqref{e:correlator F} and \eqref{e:correlator F untl} follows from the commutativity of the diagram
$$
\CD 
(\on{diag}_S)^! \circ  \sft^! @>{\sim}>> (\on{pr}_{\on{small},S})_! \circ (\on{diag}_S)_!\circ (\on{diag}_S)^! \circ  \sft^! @>>> 
(\on{pr}_{\on{small},S})_! \circ  \sft^! \\
@VVV & & @VVV \\
(\on{diag}^{\on{untl}}_S)^! @>{\sim}>> (\on{pr}^{\on{untl}}_{\on{small},S})_! \circ (\on{diag}^{\on{untl}}_S)_!\circ (\on{diag}^{\on{untl}}_S)^!
@>>> (\on{pr}^{\on{untl}}_{\on{small},S})_!,
\endCD
$$ 
as functors $\Dmod(S^{\subseteq,\on{untl}}) \to \Dmod(S)$, where
$$\on{diag}_S:S\to S^{\subseteq} \text{ and } \on{diag}^{\on{untl}}_S:S\to S^{\subseteq,\on{untl}}$$
are the corresponding maps.  

\ssec{Proofs of Lemmas \ref{l:pr small untl}, \lemref{l:untl contr} and \ref{l:non untl to untl cofinal}}

\sssec{Proof of \lemref{l:untl contr}} \label{sss:proof untl contr}

We need to show that, for any affine scheme $S$ equipped with a map $S\to \Ran^{\on{untl}}$, the functor
$$\Dmod(S) \overset{(\on{pr}_{\on{small},S})^!}\longrightarrow \Dmod(S^{\subseteq,\on{untl}})^{\on{strict}}\hookrightarrow
\Dmod(S^{\subseteq,\on{untl}})^{\on{lax}},$$
is fully faithful.

\medskip

It suffices to show that the first arrow, i.e., $(\on{pr}_{\on{small},S})^!$ is fully faithful. However, 
this follows from the fact that morphism $\on{pr}_{\on{small},S}$ admits a left adjoint. 
Namely, it is given by $\on{diag}_S$. 

\qed[\lemref{l:untl contr}]

\sssec{Proof of \lemref{l:pr small untl}} \label{sss:proof pr small untl}

Let $S$ be an affine scheme and let $\ul{x}$ be an $S$-point of $\Ran^{\on{untl}}$. Consider the map 
$$S_{\ul{x}}^{\subseteq,\on{untl}}:=S\underset{\ul{x},\Ran^{\on{untl}}}\times \Ran^{\subseteq,\on{untl}}
\overset{\on{pr}^{\on{untl}}_{\on{small},S,\ul{x}}}\longrightarrow S.$$

The first assertion of the lemma is that the functor
$$\Dmod(S) \overset{(\on{pr}^{\on{untl}}_{\on{small},S,\ul{x}})^!}\to \Dmod(S_{\ul{x}}^{\subseteq,\on{untl}})^{\on{strict}} \to
\Dmod(S_{\ul{x}}^{\subseteq,\on{untl}})$$
admits a left adjoint, to be denoted $(\on{pr}^{\on{untl}}_{\on{small},S,\ul{x}})_!$. This is a particular case of
\corref{c:integration ps-proper}. 

\medskip

The second assertion of the lemma is that for a 1-morphism 
$\ul{x}_1 \overset{\alpha}\to \ul{x}_2$ in $\Maps(S,\Ran^{\on{untl}})$
and the corresponding map 
$$S_{\ul{x}_2}^{\subseteq,\on{untl}} \overset{\alpha^*}\to S_{\ul{x}_1}^{\subseteq,\on{untl}},$$
the natural transformation
$$(\on{pr}^{\on{untl}}_{\on{small},S,\ul{x}_2})_! \circ (\alpha^*)^! \simeq 
(\on{pr}^{\on{untl}}_{\on{small},S,\ul{x}_1})_! \circ (\alpha^*)_! \circ (\alpha^*)^! 
\to (\on{pr}^{\on{untl}}_{\on{small},S,\ul{x}_1})_!$$
is an isomorphism.  

\medskip

This follows from the fact that the map $\alpha^*$ admits a value-wise
left adjoint. Namely, the left adjoint in question attaches to an affine scheme $S'$ with a map
$g_1:S'\to S_{\ul{x}_1}^{\subseteq,\on{untl}}$ the map $g_2:S'\to S_{\ul{x}_2}^{\subseteq,\on{untl}}$
defined as follows: the corresponding map
$$S' \overset{g_2}\to S_{\ul{x}_2}^{\subseteq,\on{untl}} \overset{\on{pr}^{\on{untl}}_{\on{big}}}\longrightarrow \Ran^{\on{untl}}$$
is obtained by applying the map
$$\on{union}: \Ran^{\on{untl}}\times  \Ran^{\on{untl}}\to  \Ran^{\on{untl}}$$
to 
$$S' \overset{g_1}\to S_{\ul{x}_1}^{\subseteq,\on{untl}} \overset{\on{pr}^{\on{untl}}_{\on{big}}}\longrightarrow \Ran^{\on{untl}}$$
and 
$$S'\overset{g_1}\to S_{\ul{x}_1}^{\subseteq,\on{untl}}  \overset{\on{pr}^{\on{untl}}_{\on{small},S,\ul{x}_1}}\longrightarrow S\overset{\ul{x}_2}\to
\Ran^{\on{untl}}.$$

\qed[\lemref{l:pr small untl}]

\ssec{Unitality as a property}

A somewhat surprising fact is that, given a local unital structure on $\ul\bC^{\on{loc}}$, 
a \emph{unital}\footnote{As opposed to lax unital.} structure on (a non-unital local-to-global functor) $\ul\sF$ is actually a \emph{property},
and not an additional piece of structure, as we shall presently explain.

\medskip

The contents of this subsection are not necessary for the sequel. 

\sssec{}

Let $\ul\sF$ be as in \secref{ss:local-to-global set up}. Recall the morphism \eqref{e:pre-unitality transform}
\begin{equation} \label{e:pre-unitality transform again}
\sF_\CZ \to \left(\on{Id}\otimes (\on{pr}_{\on{small},\CZ})_!\right) \circ \sF_{\CZ^\subseteq} \circ \on{ins.unit}_\CZ
\end{equation} 
as functors $\bC^{\on{loc}}_\CZ\to \bC^{\on{glob}}\otimes \Dmod(\CZ)$. 

\begin{defn}
We shall say that $\ul\sF$ satisfies Global Unitality Axiom 1 if the natural transformation \eqref{e:pre-unitality transform} is an isomorphism
(for any $\CZ\to \Ran$).
\end{defn}

\sssec{}

Assume that $\sF$ satisfies Global Unitality Axiom 1. Then inverting \eqref{e:pre-unitality transform again} and applying the
$((\on{pr}_{\on{small},\CZ})_!,(\on{pr}_{\on{small},\CZ})^!)$-adjunction, we obtain a map
\begin{equation} \label{e:unitality transform}
\sF_{\CZ^\subseteq} \circ \on{ins.unit}_\CZ \to 
\left(\on{Id}\otimes (\on{pr}_{\on{small},\CZ})^!\right)(\sF_\CZ)
\end{equation} 
as functors $\bC^{\on{loc}}_\CZ\to \bC^{\on{glob}}\otimes \Dmod(\CZ^{\subseteq})$. 

\begin{defn}
We shall say that $\ul\sF$ satisfies Global Unitality Axiom 2 if the natural transformation \eqref{e:unitality transform} is an isomorphism
(for any $\CZ\to \Ran$).
\end{defn}

\begin{defn}
We shall say that $\ul\sF$ has a global unital property if it satisfies Axioms 1 and 2. 
\end{defn}

\sssec{}

It is clear that if $\ul\sF$ is the image under
$$\on{Funct}^{\on{loc}\to \on{glob},\on{untl}}(\ul\bC^{\on{loc}},\bC^{\on{glob}})\to
\on{Funct}^{\on{loc}\to \on{glob}}(\ul\bC^{\on{loc}},\bC^{\on{glob}})$$
of 
$$\ul\sF^{\on{untl}}\in \on{Funct}^{\on{loc}\to \on{glob},\on{untl}}(\ul\bC^{\on{loc}},\bC^{\on{glob}}),$$
then $\ul\sF$ has a global unitality property, see \secref{sss:verify Axiom 1 for unital}. 

\sssec{}

Vice versa, suppose that 
$$\ul\sF\in \on{Funct}^{\on{loc}\to \on{glob}}(\ul\bC^{\on{loc}},\bC^{\on{glob}})$$
has a global unitality property. 

\medskip

Then the inverse of the isomorphism 
\eqref{e:unitality transform} provides an isomorphism as in \eqref{e:unitality transform original}.

\medskip

More generally, one can show that in this case $\ul\sF$ comes from a uniquely defined object
$\ul\sF^{\on{untl}}\in \on{Funct}^{\on{loc}\to \on{glob},\on{untl}}(\ul\bC^{\on{loc}},\bC^{\on{glob}})$.

\medskip

In other words, we claim: 

\begin{prop} \label{p:unitality as a property}
The composite functor
$$\on{Funct}^{\on{loc}\to \on{glob},\on{untl}}(\ul\bC^{\on{loc}},\bC^{\on{glob}})\to
\on{Funct}^{\on{loc}\to \on{glob},\on{lax-untl}}(\ul\bC^{\on{loc}},\bC^{\on{glob}})\to 
\on{Funct}^{\on{loc}\to \on{glob}}(\ul\bC^{\on{loc}},\bC^{\on{glob}})$$
is fully faithful, and its essential image consists of objects that have 
a global unitality property.
\end{prop} 

We will give two, rather different in spirit, proofs of \propref{p:unitality as a property}: one
in \secref{sss:unitality as a property Take 2}, and another in \secref{s:add unit colax}.

\ssec{Factorization homology} \label{ss:ch homology}

\sssec{} 

In this section we will assume that $\ul\bC^{\on{loc,untl}}$ comes from a unital lax factorization category 
$\bA$, i.e., $\ul\bC^{\on{loc,untl}}=\ul\bA$ in the notations of \secref{sss:fact cat}. Let $\ul\sF$ be a functor
$$\ul\bA\to \bC^{\on{glob}}\otimes \ul\Dmod(\Ran),$$
equipped with a lax unital structure. 

\medskip

Let now $\CA$ be a unital factorization algebra in $\bA$. We will regard $\CA\mod^{\on{fact}}(\bA)$ as 
a lax factorization category (see \secref{sss:fact mod in A untl}), and consider the corresponding crystal of categories
$\CA\ul\mod^{\on{fact}}(\bA)$ on $\Ran$, which naturally extends to a crystal of categories over $\Ran^{\on{untl}}$
(see \secref{sss:untl str on fact mod}). 

\medskip

In the particular case of factorization algebras, we will denote the functor $\on{ins.unit}_\CZ$ of \secref{e:insert unit Z} by
$$\on{ins.vac}_\CZ:\CA\mod^{\on{fact}}(\bA)_\CZ\to \CA\mod^{\on{fact}}(\bA)_{\CZ^{\subseteq}}, \quad \CZ\to \Ran.$$

This notation is meant to emphasize that the unit in $\CA\mod^{\on{fact}}(\bA)$ is the ``vacuum module", i.e., $\CA$,
viewed as a factorization module over itself. 

\sssec{}

Let
$$\ul\oblv_\CA:\CA\mod^{\on{fact}}(\bA)\to \bA$$
be the tautological forgetful functor, viewed as a functor between crystals of categories over $\Ran$. 

\medskip

Note that the unital structure on $\CA$ defines on $\ul\oblv_\CA$ a structure of \emph{lax functor}
between crystals of categories over $\Ran^{\on{untl}}$, see \secref{sss:unital on oblv}. 

\medskip

In particular
$$\ul\sF\circ \ul\oblv_\CA:\CA\ul\mod^{\on{fact}}(\bA)\to \bC^{\on{glob}}\otimes \ul\Dmod(\Ran)$$
also acquires a lax unital structure. 

\sssec{} \label{sss:fact homology defn}

The functor of \emph{factorization homology}
$$\ul{\on{C}}^{\on{fact}}_\cdot(X,\CA,-)^\sF:\CA\ul\mod^{\on{fact}}(\bA)\to \bC^{\on{glob}}\otimes \ul\Dmod(\Ran)$$
is by definition
$$(\ul\sF\circ \ul\oblv_\CA)^{\int \on{ins.unit}},$$
see \eqref{e:ins unit int as functor}. 

\medskip

For $\CZ\to \Ran$ we will denote the corresponding functor
$$\CA\mod^{\on{fact}}(\bA)_\CZ\to \bC^{\on{glob}}\otimes \Dmod(\CZ)$$
by $\on{C}^{\on{fact}}_\cdot(X,\CA,-)^\sF_\CZ$.

\medskip

For $\CZ$ pseudo-proper, we will denote the composition of $\on{C}^{\on{fact}}_\cdot(X,\CA,-)^\sF_\CZ$ with
$\on{Id}\otimes \on{C}^\cdot_c(\CZ,-)$ by $\on{C}^{\on{fact}}_\cdot(X,\CA,-)^\sF_{\int_\CZ}$.

\medskip

For $\CZ=\Ran$ and the identity map we will denote the resulting functor
$$\CA\mod^{\on{fact}}(\bA)_\Ran\to \bC^{\on{glob}}$$
by $\on{C}^{\on{fact}}_\cdot(X,\CA,-)^\sF$.

\sssec{}

Recall the natural transformation \eqref{e:correlator F}. In our case, this is a map
\begin{equation} \label{e:crltr univ}
\ul\sF\circ \ul\oblv_\CA\to \ul{\on{C}}^{\on{fact}}_\cdot(X,\CA,-)^\sF,
\end{equation}
which we will denote by $\ul{\on{Cltr}}_\CA$ and refer to as the ``correlator" map.

\medskip

For a given $\CZ\to \Ran$, this is a map
$$\on{Cltr}_{\CA,\CZ}:\sF\circ \ul\oblv_{\CA,\CZ}\to \on{C}^{\on{fact}}_\cdot(X,\CA,-)^\sF_\CZ.$$

\sssec{Example} \label{sss:usual fact homology}

Let $\bA=\Vect$, $\bC^{\on{glob}}=\Vect$ and $\ul\sF=\on{Id}$. In this case,
$$\ul{\on{C}}^{\on{fact}}_\cdot(X,\CA,-):=\ul{\on{C}}^{\on{fact}}_\cdot(X,\CA,-)^\sF, \quad
\CA\mod^{\on{fact}}(\bA)\to \ul\Dmod(\Ran)$$
is the usual functor of factorization homology.

\sssec{} \label{sss:fact homology is unital}

A key fact for us is that 
according to \secref{sss:int untl 1}, the functor $\ul{\on{C}}^{\on{fact}}_\cdot(X,\CA,-)^\sF$ acquires a natural
lax unital structure. Moreover, by \secref{sss:int untl 2}, this lax unital structure is actually \emph{strict}. 

\sssec{} \label{sss:vac fact hom}

Let us apply the functor $\on{C}^{\on{fact}}_\cdot(X,\CA,-)^\sF_\Ran$ to the object  
$$(\one_{\CA\mod^{\on{fact}}(\bA)})_\Ran=\CA^{\on{fact}}_\Ran\in \CA\mod^{\on{fact}}(\bA)_\Ran$$
(see \secref{sss:fact A mod in A untl}) 
i.e., to $\CA_\Ran$, viewed as a factorization module over itself at $\Ran$.

\medskip

By unitality, the above object is of the form
$$\on{C}^{\on{fact}}_\cdot(X,\CA)\otimes \omega_{\Ran},$$
for a canonically defined object
\begin{equation} \label{e:vac fact hom}
\on{C}^{\on{fact}}_\cdot(X,\CA)\in \bC^{\on{glob}}.
\end{equation}

The object \eqref{e:vac fact hom} is called the \emph{vacuum factorization homology} of $\CA$. 

\medskip

Explicitly,
$$\on{C}^{\on{fact}}_\cdot(X,\CA)\simeq (\on{Id}\otimes \on{C}_c^\cdot(\Ran,-)) \circ \sF(\CA_\Ran).$$

\sssec{} \label{sss:vac fact hom bis}

Note that again by the unitality (\secref{sss:fact homology is unital}) of the functor $\ul{\on{C}}^{\on{fact}}_\cdot(X,\CA,-)$, 
for any $\CZ\to \Ran$ and 
$$\CA^{\on{fact}_\CZ}\in \CA\mod^{\on{fact}}_\CZ,$$
(see Sects. \ref{sss:vacuum module} and \ref{sss:A mod in A} for the notation), we have
\begin{equation} \label{e:ins vacuum again}
\on{C}^{\on{fact}}_\cdot(X,\CA,\CA^{\on{fact}_\CZ})_\CZ\simeq \on{C}^{\on{fact}}_\cdot(X,\CA)\otimes \omega_\CZ.
\end{equation} 

In particular, for any $\ul{x}\in \Ran$, 
\begin{equation} \label{e:ins vacuum again x}
\on{C}^{\on{fact}}_\cdot(X,\CA,\CA^{\on{fact}_{\ul{x}}})_{\ul{x}}\simeq \on{C}^{\on{fact}}_\cdot(X,\CA).
\end{equation} 

\sssec{} \label{sss:expl unitality fact homology}

Let us write out explicitly the proof of the fact that the functor $\ul{\on{C}}^{\on{fact}}_\cdot(X,\CA,-)^\sF$
is unital (we will essentially repeat the argument from \secref{sss:ins unitality}). 

\medskip

Fix $\CZ\to \Ran$ and consider an object $\CM\in \CA\mod^{\on{fact}}(\bA)_\CZ$. From it we produce 
an object
$$\on{ins.unit}_\CZ(\CM)\in \CA\mod^{\on{fact}}(\bA)_{\CZ^\subseteq},$$
and further
$$\on{ins.unit}_{\CZ^{\subseteq}}(\on{ins.unit}_\CZ(\CM))\in \CA\mod^{\on{fact}}(\bA)_{\CZ^{\subseteq^2}},$$
where
$$\CZ^{\subseteq^2}:=(\CZ^{\subseteq})^{\subseteq}.$$

Consider the object
$$\CM':=\sF_{\CZ^{\subseteq^2}}\circ \oblv_{\CA,\CZ^{\subseteq^2}}\left(\on{ins.unit}_{\CZ^{\subseteq}}(\on{ins.unit}_\CZ(\CM))\right)\in 
\bC^{\on{glob}}\otimes \Dmod(\CZ^{\subseteq^2}).$$

Consider now the map
$$\on{diag}_{\CZ,\on{big}}:\CZ^\subseteq\to \CZ^{\subseteq^2}$$
(see \secref{sss:double insert}). 

\medskip

It gives rise to a map
\begin{multline} \label{e:fact hom untl expl}
(\on{pr}_{\on{small},\CZ})_!\circ (\on{diag}_{\CZ,\on{big}})^!(\CM')\simeq 
(\on{pr}_{\on{small},\CZ})_!\circ (\on{pr}_{\on{small},\CZ^\subseteq})_! \circ (\on{diag}_{\CZ,\on{big}})_!
\circ (\on{diag}_{\CZ,\on{big}})^!(\CM')\to \\
\to (\on{pr}_{\on{small},\CZ})_!\circ (\on{pr}_{\on{small},\CZ^\subseteq})_!(\CM'),
\end{multline}
and we wish to show that this map is an isomorphism.

\medskip

Now, the lax unital structure on $\ul\sF$ implies that $\CM'$ is the pullback of an object in 
$$\CM''\in \bC^{\on{glob}}\otimes \Dmod(\CZ^{\subseteq^2,\on{untl}}).$$

As in \lemref{l:non untl to untl cofinal}, one shows that one can replace both sides in 
\eqref{e:fact hom untl expl} by their unital versions, i.e., it is sufficient to show that the corresponding map
\begin{multline*} 
(\on{pr}^{\on{untl}}_{\on{small},\CZ})_!\circ (\on{diag}^{\on{untl}}_{\CZ,\on{big}})^!(\CM'')\simeq 
(\on{pr}^{\on{untl}}_{\on{small},\CZ})_!\circ (\on{pr}^{\on{untl}}_{\on{small},\CZ^\subseteq})_! \circ (\on{diag}^{\on{untl}}_{\CZ,\on{big}})_!
\circ (\on{diag}^{\on{untl}}_{\CZ,\on{big}})^!(\CM'')\to \\
\to (\on{pr}^{\on{untl}}_{\on{small},\CZ})_!\circ (\on{pr}^{\on{untl}}_{\on{small},\CZ^\subseteq})_!(\CM'')
\end{multline*}
is an isomorphism.

\medskip

However, the latter follows from the fact that the map 
$$\on{diag}^{\on{untl}}_{\CZ,\on{big}}:\CZ^{\subseteq,\on{untl}}\to \CZ^{\subseteq^2,\on{untl}}$$
is value-wise cofinal.

\qed

\sssec{}

In what follows we will need the following variant of the isomorphism we just proved:

\medskip

The natural transformation $\ul\sF\to \ul\sF^{\int\, \on{ins.unit}}$ induces a natural transformation
\begin{equation} \label{e:AB nat trans abs}
\ul{\on{C}}^{\on{fact}}_\cdot(X,\CA,-)^\sF \to \ul{\on{C}}^{\on{fact}}_\cdot(X,\CA,-)^{\sF^{\int\, \on{ins.unit}}}.
\end{equation}

We claim:

\begin{lem} \label{l:AB nat trans abs}
The natural transformation \eqref{e:AB nat trans abs} is an isomorphism.
\end{lem}

\begin{proof}

By the construction of \eqref{e:ins unit int as functor}, the unit in $\CA$ gives rise
to a natural transformation
\begin{equation} \label{e:A to B}
\sF^{\int\, \on{ins.unit}}\circ \ul\oblv_\CA\to (\sF\circ \ul\oblv_\CA)^{\int\, \on{ins.unit}}= \ul{\on{C}}^{\on{fact}}_\cdot(X,\CA,-)^\sF.
\end{equation}

Applying $\int \on{ins.unit}$ we obtain a natural transformation
$$\ul{\on{C}}^{\on{fact}}_\cdot(X,\CA,-)^{\sF^{\int\, \on{ins.unit}}}\to \int \on{ins.unit}(\ul{\on{C}}^{\on{fact}}_\cdot(X,\CA,-)^\sF),$$
so that the diagram 
$$
\CD
\ul{\on{C}}^{\on{fact}}_\cdot(X,\CA,-)^{\sF^{\int\, \on{ins.unit}}} @>{\text{\eqref{e:A to B}}}>> 
\int \on{ins.unit}\circ \int \on{ins.unit}(\sF\circ \ul\oblv_\CA)\\
@AAA @A{\sim}A{\int \on{ins.unit}(\ul{\on{Cltr}})}A \\
\ul{\on{C}}^{\on{fact}}_\cdot(X,\CA,-)^\sF @>{=}>>  \int \on{ins.unit}(\sF\circ \ul\oblv_\CA)
\endCD
$$
commutes.

\medskip

Hence, it suffices to check that the map \eqref{e:A to B} is an isomorphism. However, this is done
by the same argument as in \secref{sss:expl unitality fact homology}.

\end{proof}

\sssec{}

Here is a particular case of \lemref{l:AB nat trans abs} that we will need:

\medskip

Let $\CA_1$ and $\CA_2$ be a pair of unital factorization algebras in a unital lax factorization category $\bA_0$,
and let $\phi:\CA_1\to \CA_2$ be a unital homomorphism. Denote by $\on{res}^\phi$ the resulting functor
$$\CA_2\mod(\bA_0)\to \CA_1\mod(\bA_0).$$

Let us be given a lax unital functor
$$\ul\sF_0:\bA_0\to \bC^{\on{glob}}\otimes \ul\Dmod(\Ran).$$

For a given $\CZ\to \Ran$ we have a 
natural transformation 
\begin{multline} \label{e:AB}
\on{C}^{\on{fact}}_\cdot(X;\CA_2,-)^{\sF_0}_\CZ =
(\on{pr}_{\on{small},\CZ})_! \circ \oblv_{\CA_2,\CZ^{\subseteq}}\circ \on{ins.vac}_{\CZ,\CA_2} \simeq \\
\simeq (\on{pr}_{\on{small},\CZ})_! \circ \oblv_{\CA_1,\CZ^{\subseteq}}\circ \on{res}^\phi \circ \on{ins.vac}_{\CZ,\CA_2} \simeq \\
\simeq (\on{pr}_{\on{small},\CZ})_! \circ (\on{pr}_{\on{small},\CZ^{\subseteq}})_! \circ 
(\on{diag}_{\CZ^{\subseteq}})_!\circ \oblv_{\CA_1,\CZ^{\subseteq}}\circ \on{res}^\phi \circ \on{ins.vac}_{\CZ,\CA_2} \simeq \\
\simeq 
(\on{pr}_{\on{small},\CZ})_! \circ (\on{pr}_{\on{small},\CZ^{\subseteq}})_! \circ 
\oblv_{\CA_1,(\CZ^{\subseteq})^{\subseteq}}\circ (\on{diag}_{\CZ^{\subseteq}})_!\circ
\on{res}^\phi \circ \on{ins.vac}_{\CZ,\CA_2} \simeq \\
\simeq 
(\on{pr}_{\on{small},\CZ})_! \circ (\on{pr}_{\on{small},\CZ^{\subseteq}})_! \circ 
\oblv_{\CA_1,(\CZ^{\subseteq})^{\subseteq}}\circ (\on{diag}_{\CZ^{\subseteq}})_!\circ
(\on{diag}_{\CZ^{\subseteq}})^!\circ \on{ins.vac}_{\CZ^{\subseteq},\CA_1} \circ \on{res}^\phi \circ \on{ins.vac}_{\CZ,\CA_2} \to \\
\to 
(\on{pr}_{\on{small},\CZ})_! \circ (\on{pr}_{\on{small},\CZ^{\subseteq}})_! \circ 
\oblv_{\CA_1,(\CZ^{\subseteq})^{\subseteq}}\circ \on{ins.vac}_{\CZ^{\subseteq},\CA_1} \circ \on{res}^\phi \circ \on{ins.vac}_{\CZ,\CA_2} \simeq \\
\simeq (\on{pr}_{\on{small},\CZ})_! \circ  \on{C}^{\on{fact}}_\cdot(X;\CA_1,-))^{\ul\sF_0}_{\CZ^{\subseteq}} \circ \on{res}^\phi \circ \on{ins.vac}_{\CZ,\CA_2},
\end{multline} 

We claim:

\begin{cor} \label{c:A B lemma}
The natural transformation \eqref{e:AB} is an isomorphism. 
\end{cor}

\begin{proof}

Take $\bA:=\CA_1\mod(\bA_0)$ and $\CA=\CA_2$, viewed as a unital factorization algebra in $\CA_1\mod(\bA_0)$.
Take $\ul\sF=\sF_0\circ \ul\oblv_{\CA_1}$. 

\medskip

Then the assertion follows from \lemref{l:AB nat trans abs}, where we note that the natural transformation
\eqref{e:AB} is \eqref{e:AB nat trans abs}. 

\end{proof}

\sssec{}

Here is another application of \lemref{l:AB nat trans abs}. Let
$\Phi:\bA_1\to \bA$
be a lax unital factorization functor between lax factorization categories. Denote by $\ul\Phi$
the corresponding right-lax functor
$$\ul\bA_1\to \ul\bA$$
as crystals of categories on $\Ran^{\on{untl}}$. 

\medskip

Let
us be given a lax unital local-to-global functor
$$\ul\sF:\bA\to \bC^{\on{glob}}\otimes \ul\Dmod(\Ran).$$

Note that $\ul\sF_1:=\ul\sF\circ \ul\Phi$ also acquires a lax unital structure. By functoriality, the natural
transformation
$$\ul\sF_1=\ul\sF\circ \ul\Phi\to \ul\sF^{\int \on{ins.unit}}\circ \ul\Phi$$
gives rise to a natural transformation
\begin{equation} \label{e:two fact categ}
\ul\sF_1^{\int \on{ins.unit}}\to (\ul\sF^{\int \on{ins.unit}}\circ \ul\Phi)^{\int \on{ins.unit}}.
\end{equation}

\begin{cor} \label{c:A B categ}
The natural transformation \eqref{e:two fact categ} is an isomorphism. 
\end{cor}

\begin{proof}

By \secref{sss:fact alg from fact cat}, the functor $\Phi$ factors as
$$\bA_1\overset{\Phi^{\on{enh}}}\longrightarrow \CA\mod^{\on{fact}}(\bA)\overset{\ul\oblv_\CA}\to \bA,$$
where $\CA=\Phi(\one_{\bA_1})$.

\medskip

Since the functor $\Phi^{\on{enh}}$ is strictly unital, the assertion of the corollary reduces to the case when
$\bA_1=\CA\mod(\bA)$ and $\Phi=\oblv_\CA$. However, in the latter case, the map \eqref{e:two fact categ}
is the map \eqref{e:AB nat trans abs}.

\end{proof}

\section{Properties of the localization functor}  \label{s:prop Loc}

In this section we study the composition of the localization functors with three constructions of global nature: 

\begin{itemize}

\item The forgetful functor 
$\Dmod_\kappa(\Bun_G)\to \QCoh(\Bun_G)$;

\smallskip

\item The pullback functor 
$\Dmod_\kappa(\Bun_G)\to \Dmod_\kappa(\Bun_{G'})$ corresponding to a group homomorphism $G'\to G$;

\item For a unipotent group-scheme $N'$, the functor of de Rham
cohomology  $\Dmod(\Bun_{N'})\to \Vect$.

\end{itemize}

The pattern in the three composite functors mentioned above is that they can all be expressed via a local operation,
followed by another localization functor:

\begin{itemize}

\item Restriction $\KL(G)_{\kappa,\Ran}\to \Rep(\fL^+(G))$, followed by $\CO$-module localization 
$\Rep(\fL^+(G))_{\Ran}\to \QCoh(\Bun_G)$;

\smallskip

\item Restriction $\KL(G)_{\kappa,\Ran}\to \KL(G')_{\kappa,\Ran}$, followed by
$$\Loc_{G',\kappa}:\KL(G')_{\kappa,\Ran} \to \Dmod_\kappa(\Bun_{G'});$$

\item The functor of BRST reduction $\KL(N')_\Ran\to \Vect$.

\end{itemize}

However, there is a caveat, common to all three of these situations: in order for the local operation to reproduce
the global one, we need to precompose the former with an endofunctor of the source given by $\on{ins.vac}_{\int_\Ran}$,
see \eqref{e:ins Ran}. 

\ssec{Localization and the forgetful functor}

\sssec{}

Note that by adjunction, the commutative diagram \eqref{e:localization induction diagram} gives rise to a 
natural transformation
\begin{equation} \label{e:localization oblv diagram}
\xy
(0,0)*+{\QCoh(\Bun_G)}="A";
(60,0)*+{\Dmod_\kappa(\Bun_G)}="B";
(0,-30)*+{\Rep(\fL^+(G))_\Ran}="C";
(60,-30)*+{\KL(G)_{\kappa,\Ran}.}="D";
{\ar@{->}^{\Loc_G^{\on{\QCoh}}} "C";"A"};
{\ar@{->}_{\Loc_{G,\kappa}} "D";"B"};
{\ar@{->}_{\oblv^l_\kappa} "B";"A"};
{\ar@{->}^{\oblv^{(\hg,\fL^+(G))_\kappa}_{\fL^+(G)}} "D";"C"};
{\ar@{->} "D";"C"};
{\ar@{=>} "C";"B"};
\endxy
\end{equation}

\medskip

The natural transformation in \eqref{e:localization oblv diagram} is \emph{not} an isomorphism (unless $G=1$):
namely, evaluate both circuits on 
$$\on{Vac}(G)_{\kappa,x}\in \KL(G)_{\kappa,x}\hookrightarrow \KL(G)_{\kappa,\Ran}$$
for some $x\in X$.

\medskip

We will now draw another diagram, in which a natural transformation will be an isomorphism, which encodes 
another basic property of the localization functor. 

\sssec{} \label{sss:localization oblv 1}

Being unital factorization categories, both $\KL(G)_\kappa$ and $\Rep(\fL^+(G))$, viewed as 
crystals of categories over $\Ran$, carry local unital structures (see \secref{sss:unitality loc} for what this means). 
Furthermore, the local-to-global functors\footnote{In the formulas below, the underline has the meaning from \secref{sss:local to global functor abs}.} 
\begin{equation} \label{e:Loc untl}
\ul\Loc_{G,\kappa}:\ul\KL(G)_\kappa\to \Dmod_\kappa(\Bun_G)\otimes \ul\Dmod(\Ran)
\end{equation}
and 
\begin{equation} \label{e:Loc QCoh untl}
\ul\Loc^{\QCoh}_G:\ul\Rep(\fL^+(G))\to \QCoh(\Bun_G)\otimes \ul\Dmod(\Ran)
\end{equation}
both carry naturally defined (strict) unital structures (see \secref{sss:lax unital} for what this means). 

\medskip

Note, however, that the restriction functor
$$\ul\oblv^{(\hg,\fL^+(G))_\kappa}_{\fL^+(G)}:\ul\KL(G)_\kappa\to \ul\Rep(\fL^+(G))$$
is merely right-lax, (as is the case for any factorization functor that is lax unital as opposed to strictly unital). 

\sssec{} \label{sss:localization oblv 2}

The natural transformation 
$$\Loc_G^{\QCoh}\circ \oblv^{(\hg,\fL^+(G))_\kappa}_{\fL^+(G)}\to \oblv_\kappa^l\circ \Loc_{G,\kappa}$$
from diagram \eqref{e:localization oblv diagram} can be viewed as a natural transformation
\begin{equation} \label{e:localization oblv sheaves}
\ul{\Loc}_G^{\QCoh}\circ \ul\oblv^{(\hg,\fL^+(G))_\kappa}_{\fL^+(G)}\to \oblv_\kappa^l\circ \ul\Loc_{G,\kappa}
\end{equation} 
between functors
$$\ul\KL(G)_\kappa\rightrightarrows \QCoh(\Bun_G)\otimes \ul\Dmod(\Ran)$$
between crystals of categories over $\Ran$ (see \secref{sss:F recovers uF}).

\medskip

The two sides in \eqref{e:localization oblv sheaves} are lax unital local-to-global functors, and 
the map \eqref{e:localization oblv sheaves} is compatible with the lax unital structures.

\sssec{} \label{sss:localization oblv 3}

Note now that the right-hand side in \eqref{e:localization oblv sheaves} is strictly unital
(because $\ul\Loc_{G,\kappa}$ is unital). Hence
by \secref{sss:ins left adj}, the map \eqref{e:localization oblv sheaves} gives rise to a 
map\footnote{In the formula below, for the factorization category $\KL(G)_\kappa$ we use the notation $\on{ins.vac}$ instead of $\on{ins.unit}$.}
\begin{equation} \label{e:localization oblv sheaves int}
\left(\ul{\Loc}_G^{\QCoh}\circ \ul\oblv^{(\hg,\fL^+(G))_\kappa}_{\fL^+(G)}\right)^{\int \on{ins.vac}}\to 
\oblv_\kappa^l\circ \ul\Loc_{G,\kappa}.
\end{equation} 

We claim:

\begin{thm} \label{t:localization and oblv}
The natural transformation \eqref{e:localization oblv sheaves int} is an isomorphism.
\end{thm} 

\begin{rem}
One can interpret \thmref{t:localization and oblv} as follows: the natural transformation \eqref{e:localization oblv sheaves}
fails to be an isomorphism because the right-hand side is unital (i.e., insertion of vacuum does not change the value
of the functor), but the left-hand side is only lax unital. But once we correct this by applying $\int \on{ins.vac}$,
the corresponding map becomes an isomorphism.
\end{rem}

\sssec{}

We will now reformulate \thmref{t:localization and oblv} is concrete terms, which do 
not explicitly involve categorical prestacks:

\medskip

Let  $\int \on{ins.vac}$ be the endofunctor of $\KL(G)_{\kappa,\Ran}$ from \secref{sss:ins Ran int}.
We have:

\begin{thm} \label{t:localization and oblv orig}
The natural transformation in \eqref{e:localization oblv diagram} becomes an isomorphism after precomposing with
$\int \on{ins.vac}$. 
\end{thm} 

Note that Theorems \ref{t:localization and oblv} and \ref{t:localization and oblv orig} are logically
equivalent. This follows from \secref{sss:F recovers uF} and the commutativity of \eqref{e:ins unit int as functor untl compat}. 

\medskip

\thmref{t:localization and oblv orig} will be proved in \secref{ss:proof localization and oblv}\footnote{We supply a proof for completness.
An equivalent statement appears in \cite[Lemma 4.4.16 and  Variant 4.4.17]{CF}}. 

\sssec{}

Note that since $\ul\Loc_{G,\kappa}$ is strictly unital, the map
$$\Loc_{G,\kappa}\to  \Loc_{G,\kappa} \circ \int_\Ran\on{ins.vac}$$
is an isomorphism (see \secref{sss:verify Axiom 1 for unital}). Hence, \thmref{t:localization and oblv orig} implies:

\begin{cor} \label{c:localization and oblv}
We have a commutative diagram
$$
\CD
\Dmod_\kappa(\Bun_G) & @>{\oblv^l_\kappa}>> & \QCoh(\Bun_G) \\
@A{\Loc_{G,\kappa}}AA & & @AA{\Loc_G^{\on{\QCoh}}}A \\
\KL(G)_{\kappa,\Ran} @>{\int \on{ins.vac}}>> \KL(G)_{\kappa,\Ran} 
@>{\oblv^{(\hg,\fL^+(G))_\kappa}_{\fL^+(G)}}>> \Rep(\fL^+(G))_{\Ran}. 
\endCD
$$
\end{cor}


\sssec{}

Consider the functor 
$$\ul\Rep(\fL^+(G))\overset{\ul\Loc_G^{\on{\QCoh}}} \longrightarrow \QCoh(\Bun_G)\otimes \ul\Dmod(\Ran).$$

\medskip

By \secref{sss:fact homology defn}, for a factorization algebra $\CA\in \Rep(\fL^+(G))$, we can consider the 
local-to-global functor 
$$\ul{\on{C}}^{\on{fact}}_\cdot(X,\CA,-)^{\Loc_G^{\on{\QCoh}}}:\CA\mod^{\on{fact}}(\Rep(\fL^+(G)))\to \QCoh(\Bun_G)\otimes \ul\Dmod(\Ran).$$

Moreover, by \secref{sss:fact homology is unital}, the functor $\ul{\on{C}}^{\on{fact}}_\cdot(X,\CA,-)^{\Loc_G^{\on{\QCoh}}}$
is \emph{strictly} unital. 

\sssec{} \label{sss:BVG}

Denote 
$$\BV_{G,\kappa}:=\oblv^{(\hg,\fL^+(G))_\kappa}_{\fL^+(G)}(\on{Vac}(G)_\kappa),$$
viewed as a factorization algebra in $\Rep(\fL^+(G))$.

\medskip

The functor $\oblv^{(\hg,\fL^+(G))_\kappa}_{\fL^+(G)}$ upgrades to a strictly unital factorization functor 
$$(\oblv^{(\hg,\fL^+(G))_\kappa}_{\fL^+(G)})^{\on{enh}}:\KL(G)_\kappa\to \BV_{G,\kappa}\mod^{\on{fact}}(\Rep(\fL^+(G))).$$

Note that whereas $\ul\oblv^{(\hg,\fL^+(G))_\kappa}_{\fL^+(G)}$ was right-lax, when viewed as a functor between sheaves 
of categories over $\Ran^{\on{untl}}$, the functor $(\ul\oblv^{(\hg,\fL^+(G))_\kappa}_{\fL^+(G)})^{\on{enh}}$ is strict.

\sssec{}

With these notations, from \thmref{t:localization and oblv}, we obtain: 
\begin{cor} \label{c:localization and oblv fact hom}
We have a commutative diagram
\begin{equation} \label{e:localization and oblv fact hom}
\CD
\Dmod_\kappa(\Bun_G)\otimes \ul\Dmod(\Ran) @>{\oblv^l_\kappa\otimes \on{Id}}>>
\QCoh(\Bun_G)\otimes \ul\Dmod(\Ran) \\
@A{\ul{\Loc}_{G,\kappa}}AA  @AA{\ul{\on{C}}^{\on{fact}}_\cdot(X,\BV_{G,\kappa},-)^{\Loc_G^{\on{\QCoh}}}}A     \\
\ul\KL(G)_\kappa @>>{(\ul\oblv^{(\hg,\fL^+(G))_\kappa}_{\fL^+(G)})^{\on{enh}}}> \BV_{G,\kappa}\ul\mod^{\on{fact}}(\Rep(\fL^+(G))).
\endCD
\end{equation} 
\end{cor}

Integrating over $\Ran$, we obtain:

\begin{cor} \label{c:localization and oblv fact hom int}
We have a commutative diagram
\begin{equation} \label{e:localization and oblv fact hom int}
\CD
\Dmod_\kappa(\Bun_G) @>{\oblv^l_\kappa}>>
\QCoh(\Bun_G) \\
@A{\Loc_{G,\kappa}}AA  @AA{\on{C}^{\on{fact}}_\cdot(X,\BV_{G,\kappa},-)^{\Loc_G^{\on{\QCoh}}}}A     \\
\KL(G)_{\kappa,\Ran} @>>{(\oblv^{(\hg,\fL^+(G))_\kappa}_{\fL^+(G)})^{\on{enh}}}> \BV_{G,\kappa}\mod^{\on{fact}}(\Rep(\fL^+(G)))_\Ran.
\endCD
\end{equation} 
\end{cor}

\ssec{Proof of \propref{p:fiber of localization}}  \label{ss:proof of fiber of localization}

In this subsection we will use \thmref{t:localization and oblv orig} (or rather \corref{c:localization and oblv fact hom int})
to deduce \propref{p:fiber of localization}. 

\medskip

Let $\CP_G$ be a $k$-point of $\Bun_G$. Applying \corref{c:localization and oblv fact hom int}, 
we need to construct an isomorphism between 
\begin{equation} \label{e:corrected O-Loc LHS}
\KL(G)_{\kappa,\ul{x}} \overset{\oblv^{(\hg,\fL^+(G))_\kappa}_{\fL^+(G)}}\longrightarrow 
\Rep(\fL^+(G))_{\ul{x}}\overset{\on{C}^{\on{fact}}_\cdot(X,\BV_{G,\kappa},-)^{\Loc_G^{\on{\QCoh}}}_{\ul{x}}}\longrightarrow
\QCoh(\Bun_G) \overset{*\text{-fiber at }\CP_G}\longrightarrow \Vect
\end{equation}
and 
\begin{equation} \label{e:corrected O-Loc RHS}
\KL(G)_{\kappa,\ul{x}}\overset{\alpha_{\CP_G,\on{taut}}}\longrightarrow 
\KL(G)_{\kappa,\CP_G,\ul{x}}\to \hg\mod_{\kappa,\CP_G,\ul{x}} \overset{\text{\eqref{e:KM to out}}}\longrightarrow
\Gamma(X-\ul{x},\fg_{\CP_G})\mod \overset{\on{coinv}_{\Gamma(X-\ul{x},\fg_{\CP_G})}}\longrightarrow \Vect.
\end{equation}

\sssec{} \label{sss:twist of vacuum}

First, we note that the functor 
$$\ul\Rep(\fL^+(G))\overset{\ul\Loc_G^{\on{\QCoh}}} \longrightarrow \QCoh(\Bun_G)\otimes \ul\Dmod(\Ran) 
\overset{(*\text{-fiber at }\CP_G)\otimes \on{Id}}\longrightarrow  \ul\Dmod(\Ran)$$
is a $\CP_G$-twisted version of the forgeftul functor. Denote it by $\oblv_{\fL^+(G),\CP_G}$. We can view it
as a factorization functor
$$\oblv_{\fL^+(G),\CP_G}:\Rep(\fL^+(G))\to \Vect.$$

\medskip

Denote
$$\BV_{\fg,\kappa,\CP_G}:=\oblv_{\fL^+(G),\CP_G}\circ \oblv^{(\hg,\fL^+(G))_\kappa}_{\fL^+(G)}(\on{Vac}(G)_\kappa).$$

This is a factorization algebra in $\Vect$, which is a $\CP_G$-twisted version of the vacuum representation $\BV_{\fg,\kappa}$.
We can rewrite the functor 
$$\ul\Rep(\fL^+(G))\overset{\ul{\on{C}}^{\on{fact}}_\cdot(X,\BV_{G,\kappa},-)^{\Loc_G^{\on{\QCoh}}}}\longrightarrow
\QCoh(\Bun_G)\otimes \ul\Dmod(\Ran) \overset{*\text{-fiber at }\CP_G\otimes \on{Id}}\longrightarrow \ul\Dmod(\Ran)$$
as
$$\ul\Rep(\fL^+(G))\overset{\ul{\on{C}}^{\on{fact}}_\cdot(X,\BV_{\fg,\kappa,\CP_G},-)^{\oblv_{\fL^+(G),\CP_G}}}\longrightarrow \ul\Dmod(\Ran).$$

\sssec{}

Consider the forgetful (factorization) functor
$$\oblv_{\hg,\CP_G}:\hg\mod_{\kappa,\CP_G,\ul{x}}\to \Vect.$$

Note that
$$\oblv_{\hg,\CP_G}\circ \oblv_{\fL^+(G)}\circ \alpha_{\CP_G,\on{taut}} 
\simeq \oblv_{\fL^+(G),\CP_G}\circ \oblv^{(\hg,\fL^+(G))_\kappa}_{\fL^+(G)}$$
as factorization functors
$$\KL(G)_\kappa\to \Vect.$$

In particular,
$$\oblv_{\hg,\CP_G}(\on{Vac}(G)_{\kappa,\CP_G})\simeq \BV_{\fg,\kappa,\CP_G},$$
and we obtain that the functor $\oblv_{\hg,\CP_G}$ enhances to a functor
\begin{equation} \label{e:to naive 1}
\oblv_{\hg,\CP_G}^{\on{enh}}:\hg\mod_{\kappa,\CP_G,\ul{x}}\to  \BV_{\fg,\kappa,\CP_G}\mod^{\on{fact}},
\end{equation}
see \secref{sss:fact alg from fact cat}. 

\sssec{}

We obtain that we can rewrite \eqref{e:corrected O-Loc LHS} as
\begin{multline} \label{e:corrected O-Loc LHS 1}
\KL(G)_{\kappa,\ul{x}}\overset{\alpha_{\CP_G,\on{taut}}}\longrightarrow 
\KL(G)_{\kappa,\CP_G,\ul{x}} \overset{\oblv_{\fL^+(G)}}\longrightarrow \hg\mod_{\kappa,\CP_G,\ul{x}}
\overset{\oblv_{\hg,\CP_G}^{\on{enh}}}\longrightarrow \\
\to \BV_{\fg,\kappa,\CP_G}\mod^{\on{fact}}_{\ul{x}}
\overset{\ul{\on{C}}^{\on{fact}}_\cdot(X,\BV_{\fg,\kappa,\CP_G},-)}\longrightarrow \Vect.
\end{multline}

\medskip

Thus, we obtain that it suffices to show that the composition
\begin{equation} \label{e:corrected O-Loc LHS 2}
\hg\mod_{\kappa,\CP_G,\ul{x}}
\overset{\oblv_{\hg,\CP_G}^{\on{enh}}}\longrightarrow 
\BV_{\fg,\kappa,\CP_G}\mod^{\on{fact}}_{\ul{x}}
\overset{\ul{\on{C}}^{\on{fact}}_\cdot(X,\BV_{\fg,\kappa,\CP_G},-)}\longrightarrow \Vect
\end{equation}
identifies with 
\begin{equation} \label{e:corrected O-Loc RHS 2}
\hg\mod_{\kappa,\CP_G,\ul{x}}\overset{\text{\eqref{e:KM to out}}}\longrightarrow
\Gamma(X-\ul{x},\fg_{\CP_G})\mod \overset{\on{coinv}_{\Gamma(X-\ul{x},\fg_{\CP_G})}}\longrightarrow \Vect.
\end{equation}

\medskip

We will relate \eqref{e:corrected O-Loc LHS 2} to the functor of $\Gamma(X-\ul{x},\fg_{\CP_G})$-coinvariants 
using the calculation of factorization homology of (twisted) chiral envelopes of Lie-* algebras performed in \cite{BD2}.

\sssec{}

Observe that the chiral algebra corresponding to the factorization algebra $\BV_{\fg,\kappa,\CP_G}$
identifies with ($\kappa$-twisted) the chiral envelope 
$$U^{\on{ch}}(L_{\fg,\CP_G})_\kappa,$$
where $L_{\fg,\kappa,\CP_G}$ is the $\CP_G$-twist of the Lie-algebra 
$$\omega_X\oplus (\fg\otimes \on{D}_X).$$

\medskip

Moreover, we have a canonical equivalence
$$L_{\fg,\kappa,\CP_G}\mod^{\on{ch}}_{\ul{x}} \simeq U^{\on{ch}}(L_{\fg,\CP_G})_\kappa\mod^{\on{ch}}_{\ul{x}}\simeq 
\BV_{\fg,\kappa,\CP_G}\mod^{\on{fact}}_{\ul{x}}$$
and a functor

\begin{equation} \label{e:restr to out Lie alg}
L_{\fg,\kappa,\CP_G}\mod^{\on{ch}}_{\ul{x}} \to \Gamma(X-\ul{x},\fg_{\CP_G})\mod.
\end{equation} 

Under these identifications, the functor
$$\hg\mod_{\kappa,\CP_G,\ul{x}}\overset{\oblv_{\hg,\CP_G}^{\on{enh}}}\longrightarrow 
\BV_{\fg,\kappa,\CP_G}\mod^{\on{fact}}_{\ul{x}} \simeq L_{\fg,\kappa,\CP_G}\mod^{\on{ch}}_{\ul{x}}
\overset{\text{\eqref{e:restr to out Lie alg}}}\longrightarrow 
\Gamma(X-\ul{x},\fg_{\CP_G})\mod$$
identifies with the functor \eqref{e:KM to out}. 

\sssec{}\label{sss:proof fiber of localization end}

Thus, it remains to show that the functor
$$L_{\fg,\kappa,\CP_G}\mod^{\on{ch}}_{\ul{x}} \simeq 
\BV_{\fg,\kappa,\CP_G}\mod^{\on{fact}}_{\ul{x}} \overset{\on{C}^{\on{fact}}_\cdot(X,\BV_{\fg,\kappa,\CP_G},-)_{\ul{x}}}\longrightarrow \Vect$$
identifies canonically with
$$L_{\fg,\kappa,\CP_G}\mod^{\on{ch}}_{\ul{x}} 
\overset{\text{\eqref{e:restr to out Lie alg}}}\longrightarrow 
\Gamma(X-\ul{x},\fg_{\CP_G})\mod
\overset{\on{coinv}_{\Gamma(X-\ul{x},\fg_{\CP_G})}}\longrightarrow \Vect.$$

\medskip

However, the latter is the assertion of \cite[Proposition 4.8.2]{BD2}

\begin{rem} An alternative proof of the latter assertion can be found in
\cite[Corollary 6.4.4]{FraG} in the special case when the coefficient module is the vacuum representation. 

\medskip

However, the method from \cite{FraG} easily adapts to the present setting and can also be used to reprove 
\cite[Proposition 4.8.2]{BD2} in the generality in which we are using it.

\end{rem} 

\qed[\propref{p:fiber of localization}]

\ssec{Localization and restriction}

\sssec{}

Let $\phi:G'\to G$ be a group homomorphism. We restrict the level $\kappa$ to $G'$ and consider the corresponding Kazhdan-Lusztig category
$$\KL(G')_\kappa:= \wh\fg'\mod_\kappa^{\fL^+(G')}$$
and the localization functor
$$\Loc_{G',\kappa}:\KL(G')_{\kappa,\Ran}\to \Dmod_\kappa(\Bun_{G'}).$$

\sssec{}

The map $\phi$ gives rise to (factorization) restriction functors
$$\Rep(\fL^+(G))\overset{\on{res}^\phi}\longrightarrow \Rep(\fL^+(G')) \text{ and }
\KL(G)_\kappa\overset{\on{res}^\phi}\longrightarrow \KL(G')_\kappa$$
so that the diagram 
$$
\CD
\KL(G)_\kappa @>{\on{res}^\phi}>> \KL(G')_\kappa \\
@V{\oblv^{(\hg,\fL^+(G))_\kappa}_{\fL^+(G)}}VV @VV{\oblv^{(\hg',\fL^+(G'))_\kappa}_{\fL^+(G')}}V \\
\Rep(\fL^+(G)) @>{\on{res}^\phi}>>  \Rep(\fL^+(G'))
\endCD
$$
commutes.

\medskip

In addition, the map $\phi$ gives rise to a map
$$\phi^{\on{glob}}:\Bun_{G'}\to \Bun_G,$$
which is compatible with the twistings and thus gives rise to a functor
$$(\phi^{\on{glob}})^!_\kappa:\Dmod_\kappa(\Bun_G)\to \Dmod_\kappa(\Bun_{G'}),$$
which makes the diagram
$$
\CD
\QCoh(\Bun_G) @>{(\phi^{\on{glob}})^*}>> \QCoh(\Bun_{G'})  \\
@A{\oblv_\kappa^l}AA    @AA{\oblv_\kappa^l}A \\
\Dmod_\kappa(\Bun_G) @>{(\phi^{\on{glob}})^!_\kappa}>> \Dmod_\kappa(\Bun_{G'})
\endCD
$$
commute. 

\sssec{}

Let $U\subset \Bun_G$ and $U'\subset \Bun_{G'}$ be a pair of quasi-compact substacks so that
$\phi^{\on{glob}}$ maps $U'\to U$.

\medskip

Consider the corresponding functors
$$\Gamma_{G,\kappa,U}: \Dmod_\kappa(U)\to  \KL(G)_{\kappa,\Ran} \text{ and }
\Gamma_{G',\kappa,U}: \Dmod_\kappa(U')\to  \KL(G')_{\kappa,\Ran}.$$

By construction, we have a natural transformation
\begin{equation} \label{e:Gamma restr diag U}
\xy
(0,0)*+{\Dmod_\kappa(U)}="A";
(40,0)*+{\Dmod_\kappa(U')}="B";
(0,-30)*+{\KL(G)_{\kappa,\Ran}}="C";
(40,-30)*+{\KL(G')_{\kappa,\Ran}.}="D";
{\ar@{<-}^{\Gamma_{G,\kappa,U}} "C";"A"};
{\ar@{<-}_{\Gamma_{G',\kappa,U}} "D";"B"};
{\ar@{->}_{\on{res}^\phi} "C";"D"};
{\ar@{->}^{(\phi^{\on{glob}})^!_\kappa} "A";"B"};
{\ar@{=>} "C";"B"};
\endxy
\end{equation}

By adjunction, we obtain a natural transformation

\begin{equation} \label{e:Loc restr diag U}
\xy
(0,0)*+{\Dmod_\kappa(U)}="A";
(40,0)*+{\Dmod_\kappa(U')}="B";
(0,-30)*+{\KL(G)_{\kappa,\Ran}}="C";
(40,-30)*+{\KL(G')_{\kappa,\Ran}.}="D";
{\ar@{->}^{\Loc_{G,\kappa,U}} "C";"A"};
{\ar@{->}_{\Loc_{G',\kappa,U'}} "D";"B"};
{\ar@{->}_{\on{res}^\phi} "C";"D"};
{\ar@{->}^{(\phi^{\on{glob}})^!_\kappa} "A";"B"};
{\ar@{=>} "D";"A"};
\endxy
\end{equation}

Passing to the limit over $U$, from \eqref{e:Loc restr diag U}, we obtain a natural transformation
\begin{equation} \label{e:Loc restr diag}
\xy
(0,0)*+{\Dmod_\kappa(\Bun_G)}="A";
(40,0)*+{\Dmod_\kappa(\Bun_{G'})}="B";
(0,-30)*+{\KL(G)_{\kappa,\Ran}}="C";
(40,-30)*+{\KL(G')_{\kappa,\Ran}.}="D";
{\ar@{->}^{\Loc_{G,\kappa}} "C";"A"};
{\ar@{->}_{\Loc_{G',\kappa}} "D";"B"};
{\ar@{->}_{\on{res}^\phi} "C";"D"};
{\ar@{->}^{(\phi^{\on{glob}})^!_\kappa} "A";"B"};
{\ar@{=>} "D";"A"};
\endxy
\end{equation}

\sssec{}

The natural transformation in \eqref{e:Loc restr diag} is not an isomorphism (unless $\phi$ itself is). We will
now draw another diagram, in which the natural transformation is an isomorphism, and which expresses the composition
$$(\phi^{\on{glob}})^!_\kappa\circ \Loc_{G,\kappa}$$
via $\Loc_{G',\kappa}$.

\medskip

This will be completely parallel to Sects. \ref{sss:localization oblv 1}-\ref{sss:localization oblv 3}. 

\sssec{}

The natural transformation
$$\Loc_{G',\kappa}\circ \on{res}^\phi\to (\phi^{\on{glob}})^!_\kappa\circ \Loc_{G,\kappa}$$
in \eqref{e:Loc restr diag} can be viewed as a natural transformation 
\begin{equation} \label{e:localization restr sheaves}
\ul\Loc_{G',\kappa}\circ \ul{\on{res}}^\phi\to (\phi^{\on{glob}})^!_\kappa\circ \ul\Loc_{G,\kappa}
\end{equation} 
between functors
$$\ul\KL(G)_\kappa\rightrightarrows \Dmod_\kappa(\Bun_{G'})\otimes \ul\Dmod(\Ran)$$
between crystals of categories over $\Ran$ (see \secref{sss:F recovers uF}).

\medskip

The two sides in \eqref{e:localization restr sheaves} are lax unital local-to-global functors, and 
the map \eqref{e:localization oblv sheaves} is compatible with the lax unital structures.

\medskip

Note now that the right-hand side in \eqref{e:localization restr sheaves} is strictly unital. Hence
by \secref{sss:ins left adj}, the map \eqref{e:localization restr sheaves} gives rise to a 
map
\begin{equation} \label{e:localization restr sheaves int}
\left(\ul\Loc_{G',\kappa}\circ \ul{\on{res}}^\phi\right)^{\int \on{ins.vac}}\to (\phi^{\on{glob}})^!_\kappa\circ \ul\Loc_{G,\kappa}
\end{equation} 
as unital local-to-global functors. 

\medskip

We claim:

\begin{prop} \label{p:localization restr ult}
The natural transformation \eqref{e:localization restr sheaves int} is an isomorphism.
\end{prop} 

Parallel to \secref{p:localization restr ult}, we will now reformulate \propref{p:localization restr ult} is several 
(equivalent) ways. 

\medskip

\begin{prop} \label{p:localization restr orig}
The natural transformation in \eqref{e:Loc restr diag} becomes an isomorphism after precomposing with
$\int \on{ins.vac}$. 
\end{prop} 

\begin{cor} \label{c:localization and restr}
We have a commutative diagram
$$
\CD
\Dmod_\kappa(\Bun_G) & @>{(\phi^{\on{glob}})^!_\kappa}>> & \Dmod_\kappa(\Bun_{G'}) \\
@A{\Loc_{G,\kappa}}AA & & @AA{\Loc_{G',\kappa}}A \\
\KL(G)_{\kappa,\Ran} @>{\int \on{ins.vac}}>> \KL(G)_{\kappa,\Ran} @>{\on{res}^\phi}>> \KL(G')_{\kappa,\Ran} 
\endCD
$$
\end{cor} 

\sssec{}

Consider the factorization algebra
$$\on{Vac}(G|_{G'})_\kappa:=\on{res}^\phi(\on{Vac}(G)_\kappa)\in \KL(G')_\kappa.$$

The functor $\on{res}^\phi$ upgrades to a factorization functor
$$(\on{res}^\phi)^{\on{enh}}:\KL(G)_\kappa\to \on{Vac}(G|_{G'})_\kappa\mod^{\on{fact}}(\KL(G')_\kappa).$$

\medskip

Consider the functor
$$\ul{\on{C}}^{\on{fact}}_\cdot(X,\on{Vac}(G|_{G'})_\kappa,-)^{\Loc_{G',\kappa}}:
\on{Vac}(G|_{G'})_\kappa\mod^{\on{fact}}(\KL(G')_\kappa)\to \Dmod_\kappa(\Bun_{G'}),$$
see \secref{sss:fact homology defn}. 

\medskip

\begin{cor} \label{c:localization and restr fact hom}
We have a commutative diagram
\begin{equation} \label{e:localization and restr fact hom}
\CD
\Dmod_\kappa(\Bun_G)\otimes \ul\Dmod(\Ran) @>{(\phi^{\on{glob}})^!_\kappa\otimes \on{Id}}>>
\Dmod_\kappa(\Bun_{G'}) \otimes \ul\Dmod(\Ran) \\
@A{\ul{\Loc}_{G,\kappa}}AA  @AA{\ul{\on{C}}^{\on{fact}}_\cdot(X,\on{Vac}(G|_{G'})_\kappa,-)^{\Loc_{G',\kappa}}}A     \\
\KL(G)_\kappa @>>{(\on{res}^\phi)^{\on{enh}}}> \on{Vac}(G|_{G'})_\kappa\mod^{\on{fact}}(\KL(G')_\kappa).
\endCD
\end{equation} 
\end{cor}

Integrating over $\Ran$, we obtain:

\begin{cor} \label{c:localization and restr fact hom Int}
We have a commutative diagram
\begin{equation} \label{e:localization and restr fact hom Int}
\CD
\Dmod_\kappa(\Bun_G) @>{(\phi^{\on{glob}})^!_\kappa}>>
\Dmod_\kappa(\Bun_{G'}) \\
@A{\Loc_{G,\kappa}}AA  @AA{\ul{\on{C}}^{\on{fact}}_\cdot(X,\on{Vac}(G|_{G'})_\kappa,-)}A     \\
\KL(G)_{\kappa,\Ran} @>>{(\on{res}^\phi)^{\on{enh}}}> \on{Vac}(G|_{G'})_\kappa\mod^{\on{fact}}(\KL(G')_\kappa)_\Ran.
\endCD
\end{equation} 
\end{cor}

\ssec{Proof of \propref{p:localization restr ult}}

\sssec{}

Since the functor
\begin{equation} \label{e:oblv l G'}
\oblv_{\kappa,G'}^l: \Dmod_\kappa(\Bun_{G'})\to \QCoh(\Bun_{G'})
\end{equation} 
is conservative, it is sufficient to prove that the natural transformation in \eqref{e:localization restr sheaves int} becomes
an isomorphism after composing with 
$$\oblv_{\kappa,G'}^l\otimes \on{Id}:\Dmod_\kappa(\Bun_{G'})\otimes \ul\Dmod(\Ran^{\on{untl}})\to \QCoh(\Bun_{G'})\otimes \ul\Dmod(\Ran^{\on{untl}}).$$

\begin{rem}
The idea of the proof is the following: the diagram 
$$
\CD
\QCoh(\Bun_G)@>{(\phi^{\on{glob}})^*}>> \QCoh(\Bun_{G'}) \\
@A{\Loc^{\QCoh}_G}AA @AA{\Loc^{\QCoh}_{G'}}A \\
\Rep(\fL^+(G))_\Ran @>{\on{res}^\phi}>> \Rep(\fL^+(G'))_\Ran 
\endCD
$$
commutes tautologically.

\medskip

Combining this observation with \corref{c:localization and oblv}, we can express both 
$$\oblv^l_\kappa\circ (\phi^{\on{glob}})^!_\kappa \circ \Loc_{G,\kappa} \text{ and } 
\oblv^l_\kappa\circ  \Loc_{G',\kappa} \circ \on{res}^\phi$$
in terms of 
$$\Loc^{\QCoh}_{G'} \circ \on{res}^\phi \circ \oblv^{(\hg,\fL^+(G))_\kappa}_{\fL^+(G)}.$$

The two sides will not match on the nose, but the difference will be accounted for
by \corref{c:A B categ}.

\end{rem}

\sssec{}

By construction, we have a commutative diagram of lax unital local-to global functors
$$
\CD
(\oblv^l_{G',\kappa}\otimes \on{Id}) \circ \ul{\Loc}_{G',\kappa}\circ \ul{\on{res}}^\phi 
@>{(\oblv_{\kappa,G'}^l\otimes \on{Id})\text{\eqref{e:localization restr sheaves}}}>> 
(\oblv^l_{G',\kappa}\otimes \on{Id}) \circ (\phi^{\on{glob}})^!_\kappa \circ \ul{\Loc}_{G,\kappa} \\
@A{\text{\eqref{e:localization oblv sheaves}}\circ \ul{\on{res}}^\phi }AA \\
\ul{\Loc}^{\QCoh}_{G'}\circ \ul\oblv^{(\hg',\fL^+(G'))_\kappa}_{\fL^+(G')} \circ \ul{\on{res}}^\phi \\
@A{\sim}AA @AA{\sim}A \\
\ul{\Loc}^{\QCoh}_{G'}\circ \ul{\on{res}}^\phi \circ \ul\oblv^{(\hg,\fL^+(G))_\kappa}_{\fL^+(G)}  \\
@A{\sim}AA \\
(\phi^{\on{glob}})^* \circ \ul{\Loc}^{\QCoh}_{G} \circ  \ul\oblv^{(\hg,\fL^+(G))_\kappa}_{\fL^+(G)} 
@>{(\phi^{\on{glob}})^*\text{\eqref{e:localization oblv sheaves}}}>> (\phi^{\on{glob}})^* \circ \oblv^l_{G,\kappa}\circ \ul{\Loc}_{G,\kappa}.
\endCD
$$

By adjunction, we obtain a diagram
$$
\CD
(\oblv^l_{G',\kappa}\otimes \on{Id}) \circ \left(\ul{\Loc}_{G',\kappa}\circ \ul{\on{res}}^\phi\right)^{\int \on{ins.vac}} 
@>{(\oblv_{\kappa,G'}^l\otimes \on{Id})\text{\eqref{e:localization restr sheaves int}}}>>
(\oblv^l_{G',\kappa}\otimes \on{Id}) \circ (\phi^{\on{glob}})^!_\kappa \circ \ul{\Loc}_{G,\kappa} \\ 
@A{\sim}AA @AA{=}A \\
\left((\oblv^l_{G',\kappa}\otimes \on{Id}) \circ \ul{\Loc}_{G',\kappa}\circ \ul{\on{res}}^\phi\right)^{\int \on{ins.vac}}
 @>>> 
(\oblv^l_{G',\kappa}\otimes \on{Id}) \circ (\phi^{\on{glob}})^!_\kappa \circ \ul{\Loc}_{G,\kappa} \\
@A{(\text{\eqref{e:localization oblv sheaves}}\circ \on{res}^\phi )^{\int \on{ins.vac}}}AA \\
\left(\ul{\Loc}^{\QCoh}_{G'}\circ \ul\oblv^{(\hg',\fL^+(G'))_\kappa}_{\fL^+(G')} \circ \ul{\on{res}}^\phi\right)^{\int \on{ins.vac}} \\
@A{\sim}AA @AA{\sim}A \\
\left(\ul{\Loc}^{\QCoh}_{G'}\circ \ul{\on{res}}^\phi \circ \ul\oblv^{(\hg,\fL^+(G))_\kappa}_{\fL^+(G)}\right)^{\int \on{ins.vac}}  \\
@A{\sim}AA \\
\left((\phi^{\on{glob}})^* \circ \ul{\Loc}^{\QCoh}_{G} \circ  \ul\oblv^{(\hg,\fL^+(G))_\kappa}_{\fL^+(G)}\right)^{\int \on{ins.vac}}
@>>> (\phi^{\on{glob}})^* \circ \oblv^l_{G,\kappa}\circ \ul{\Loc}_{G,\kappa} \\
@A{\sim}AA @AA{=}A \\
(\phi^{\on{glob}})^*\circ \left(\ul{\Loc}^{\QCoh}_{G} \circ  \ul\oblv^{(\hg,\fL^+(G))_\kappa}_{\fL^+(G)}\right)^{\int \on{ins.vac}}
@>{(\phi^{\on{glob}})^*\text{\eqref{e:localization oblv sheaves int}}}>{\sim}> (\phi^{\on{glob}})^* \circ \oblv^l_{G,\kappa}\circ \ul{\Loc}_{G,\kappa},
\endCD
$$
where the bottom arrow is an isomorphism by \thmref{t:localization and oblv} (for $G$). 

\sssec{}

Hence, to prove that $(\oblv_{\kappa,G'}^l\otimes \on{Id})\text{\eqref{e:localization restr sheaves int}}$ is an isomorphism,
it suffices to show that the map
\begin{multline}  \label{e:localization restr 1}
\left(\ul{\Loc}^{\QCoh}_{G'}\circ \ul\oblv^{(\hg',\fL^+(G'))_\kappa}_{\fL^+(G')} \circ \ul{\on{res}}^\phi\right)^{\int \on{ins.vac}} 
\overset{(\text{\eqref{e:localization oblv sheaves}}\circ \ul{\on{res}}^\phi )^{\int \on{ins.vac}}}\longrightarrow \\
\to \left((\oblv^l_{G',\kappa}\otimes \on{Id}) \circ \ul{\Loc}_{G',\kappa}\circ \ul{\on{res}}^\phi\right)^{\int \on{ins.vac}}
\end{multline}
induced by
$$\ul{\Loc}^{\QCoh}_{G'}\circ \ul\oblv^{(\hg',\fL^+(G'))_\kappa}_{\fL^+(G')} \circ \ul{\on{res}}^\phi
\overset{\text{\eqref{e:localization oblv sheaves}}\circ \ul{\on{res}}^\phi}\longrightarrow
(\oblv^l_{G',\kappa}\otimes \on{Id}) \circ \ul{\Loc}_{G',\kappa}\circ \ul{\on{res}}^\phi,$$
is an isomorphism.

\medskip

Note that the map
$$\left(\ul{\Loc}^{\QCoh}_{G'}\circ \ul\oblv^{(\hg',\fL^+(G'))_\kappa}_{\fL^+(G')}\right)^{\int \on{ins.vac}} 
\overset{\text{\eqref{e:localization oblv sheaves}}^{\int \on{ins.vac}}}\longrightarrow
\left((\oblv^l_{G',\kappa}\otimes \on{Id}) \circ \ul{\Loc}_{G',\kappa}\right)^{\int \on{ins.vac}},$$
induced by 
$$\ul{\Loc}^{\QCoh}_{G'}\circ \ul\oblv^{(\hg',\fL^+(G'))_\kappa}_{\fL^+(G')}
\overset{\text{\eqref{e:localization oblv sheaves}}}\longrightarrow
(\oblv^l_{G',\kappa}\otimes \on{Id}) \circ \ul{\Loc}_{G',\kappa},$$
is an isomorphism, by \thmref{t:localization and oblv} (for $G'$). 

\medskip

This implies that \eqref{e:localization restr 1} is an isomorphism by \corref{c:A B categ}.

\qed[\propref{p:localization restr ult}]

%

\ssec{Localization for unipotent group-schemes} \label{ss:loc unip}

Let $N'$ be a unipotent group-scheme over $X$.  We will make the following technical assumption:
$N'$ admits a filtration by normal group-schemes with abelian subquotients. 

\sssec{}

Consider the factorization category
$$\KL(N'):= \fL(\fn')\mod^{\fL^+(N')}.$$

Note that the critical level for $N'$ is zero. In particular, we have the self-dualites
\begin{equation} \label{e:KM unip self-dual}
(\fL(\fn')\mod)^\vee \simeq \fL(\fn')\mod.
\end{equation} 
\begin{equation} \label{e:KL unip self-dual}
\KL(N')^\vee \simeq \KL(N').
\end{equation} 

Both dualities take place in the sense of unital factorization categories, see \secref{sss:duality lax untl}. 

\sssec{}

Note that since $\fL^+(N')$ is pro-unipotent, the forgetful functor
$$\oblv_{\fL^+(N')}:\KL(N')\to \fL(\fn')\mod$$
is fully faithful.

\medskip

Recall also that the right adjoint $$\on{Av}^{\fL^+(N')}_*:\fL(\fn')\mod\to \KL(N')$$
of $\oblv_{\fL^+(N')}$ identifies also with the dual of $\oblv_{\fL^+(N')}$ with respect
to the self-dualities \eqref{e:KM unip self-dual} and \eqref{e:KL unip self-dual}.

\sssec{}

Consider the (factorization) functor of semi-infinite cohomology with respect to $\fL(\fn')$:
$$\BRST_{\fn'}:\fL(\fn')\mod\to \Vect.$$

In terms of the duality \eqref{e:KM unip self-dual}, the functor $\BRST_{\fn'}$ is given by
$$\langle k,-\rangle_{\fL(\fn')\mod},$$
where:

\begin{itemize}

\item $k\in \KL_{N'}\subset \fL(\fn')\mod$ is the \emph{trivial} representation;

\medskip

\item $\langle -,-\rangle_{\fL(\fn')\mod}$ denotes the pairing $\fL(\fn')\mod\otimes \fL(\fn')\mod\to \Vect$,
corresponding to \eqref{e:KM unip self-dual}.

\end{itemize}

Since $k$ upgrades to an object of $\on{FactAlg}^{\on{untl}}(X,\fL(\fn')\mod)$, the functor $\BRST_{\fn'}$
has a natural lax unital factorization structure. 

\sssec{}

The value of $\BRST_{\fn'}$ on $\on{Vac}(N')$ is
\begin{multline} \label{e:BRST of Vac 0}
\langle k,\on{Vac}(N')\rangle_{\fL(\fn')\mod}\simeq 
\langle \on{Av}^{\fL^+(N')}_*(k),\on{Vac}(N')\rangle_{\KL(N')}\simeq \\
\simeq \langle k,\on{Vac}(N')\rangle_{\KL(N')}\simeq 
\langle k,\one_{\Rep(\fL^+(N'))}\rangle_{\Rep(\fL^+(N'))},
\end{multline}
where:

\begin{itemize}

\item $\langle -,-\rangle_{\KL(N')}$ is the pairing $\KL(N')\otimes \KL(N')\to \Vect$
corresponding to \eqref{e:KL unip self-dual};

\smallskip

\item $\one_{\Rep(\fL^+(N'))}\in \Rep(\fL^+(N'))$ is the \emph{trivial} representation;

\smallskip

\item $\langle -,-\rangle_{\Rep(\fL^+(N'))}$ denotes pairing $\Rep(\fL^+(N'))\otimes \Rep(\fL^+(N'))\to \Vect$,
corresponding to the natural self-duality of $\Rep(\fL^+(N'))$, i.e.,
$$\langle V_1,V_2\rangle_{\Rep(\fL^+(N'))}=\CHom_{\Rep(\fL^+(N'))}(\one_{\Rep(\fL^+(N'))},V_1\otimes V_2).$$

\end{itemize}

Denote
$$\Omega(\fn'):=\on{inv}_{\fL^+(N')}(\one_{\Rep(\fL^+(N'))})\simeq \on{C}^\cdot(\fL^+(\fn'))\in \on{FactAlg}^{\on{untl}}(X).$$

\begin{rem} \label{r:Omega as Chev}

Note also that $\Omega(\fn')$ is the commutative factorization algebra canonically isomorphic to
$$\on{Fact}(\on{C}^\cdot_{\on{chev}}(L_{\fn'})),$$
where:

\begin{itemize}

\item $L_{\fn'}$ is the Lie-* algebra $\fn'\otimes \on{D}_X$;

\medskip

\item For a Lie-* algebra $L$, we denote by 
$\on{C}^\cdot_{\on{chev}}(L)\in \on{ComAlg}(\Dmod(X))$ its cohomological Chevalley complex, see 
\cite[Sect. 1.4.10]{BD2}.

\end{itemize}

\end{rem}

\sssec{}

Hence, from \eqref{e:BRST of Vac 0}, we obtain
\begin{equation} \label{e:BRST of Vac}
\BRST_{\fn'}(\on{Vac}(N'))\simeq \Omega(\fn'),
\end{equation}
as factorization algebras (in $\Vect$).

\medskip

In particular, the functor $\BRST_{\fn'}$ enhances to a functor
$$\BRST^{\on{enh}}_{\fn'}:\fL(\fn')\mod\to \Omega(\fn')\mod^{\on{fact}}.$$

By a slight abuse of notation, we will denote by the same symbols $\BRST_{\fn'}$ and $\BRST^{\on{enh}}_{\fn'}$
the restrictions of the above functors along
$$\KL(N')\to \fL(\fn')\mod.$$

\sssec{} \label{sss:l N' line}

Let $\delta_{N'}$ denote the integer $\dim(\Bun_{N'})$.

\medskip

Note that the canonical line bundle $K_{\Bun_{N'}}$ of $\Bun_{N'}$, i.e.,
$$\det(T^*(\Bun_{N'})),$$
is canonically constant. Let $\fl_{N'}$ denote the corresponding (ungraded) line. 

\medskip

The material in \secref{s:loc} applies as-is to the group scheme $N'$ over $X$, so we can consider the localization functor
$$\Loc_{N'}:\KL(N')_\Ran\to \Dmod(\Bun_{N'}).$$

\medskip

Note also that since $N'$ is unipotent, the stack $\Bun_{N'}$ is quasi-compact. In particular, there is no difference
between $\Dmod(\Bun_{N'})$ and $\Dmod_{\on{co}}(\Bun_{N'})$. Moreover, $\Bun_{N'}$ is \emph{safe} in the sense 
of \cite[Sect. 10.2]{DrGa1} so the ``constant sheaf"
$$\ul{k}\in \Dmod(\Bun_{N'})$$
is compact, and the functor
$$\on{C}^\cdot_\dr(\Bun_{N'},-)=\CHom_{\Dmod(\Bun_{N'})}(\ul{k},-)$$
is continuous.

\sssec{}

Our next goal is to construct a natural transformation
from the composition
\begin{equation} \label{e:unip localization LHS}
\KL(N')_\Ran \overset{\BRST_{\fn'}}\longrightarrow \Dmod(\Ran)
\overset{\on{C}^\cdot_c(\Ran,-)}\longrightarrow \Vect
\overset{-\otimes \fl_{N'}[\delta_{N'}]}\longrightarrow \Vect
\end{equation} 
to 
\begin{equation} \label{e:unip localization RHS}
\KL(N')_\Ran\overset{\Loc_{N'}}\longrightarrow \Dmod(\Bun_{N'})
\overset{\on{C}^\cdot_\dr(\Bun_{N'},-)}\longrightarrow \Vect,
\end{equation} 
i.e., 
\begin{equation} \label{e:unip localization}
\left(\on{C}^\cdot_c(\Ran,-) \circ \BRST_{\fn'}\right) \otimes \fl_{N'}[\delta_{N'}]\to 
\on{C}^\cdot_\dr(\Bun_{N'},-)\circ \Loc_{N'}.
\end{equation} 

\sssec{}

Let us interpret $\on{C}^\cdot_\dr(\Bun_{N'},-)\circ \Loc_{N'}$ as 
\begin{equation} \label{e:unip localization 1}
\langle \omega_{\Bun_{N'}}, \Loc_{N'}(-)\rangle_{\Bun_{N'}},
\end{equation} 
where 
$$\langle -,-\rangle_{\Bun_{N'}}:\Dmod(\Bun_{N'})\otimes \Dmod(\Bun_{N'})\to \Vect$$
is the Verdier duality pairing.

\medskip

Using \propref{p:Loc as dual}, we rewrite \eqref{e:unip localization 1} as 
\begin{equation} \label{e:unip localization 2}
\langle \Gamma_{N',\Ran}(\omega_{\Bun_{N'}}),-\rangle_{\KL(N')_\Ran} \otimes \fl_{N'}[\delta_{N'}],
\end{equation} 
where 
$$\langle -,-\rangle_{\KL(N')_\Ran}$$
is the self-duality on $\KL_{N',\Ran}$ induced by \eqref{e:KL unip self-dual}. 

\sssec{}

Note now that for any $\ul{x}\in \Ran$
$$\Gamma^{\on{ren}}(\Bun^{\on{level}_{\ul{x}}}_{N'},\omega_{\Bun_{N'}})$$
(see Equation \eqref{e:global sections level} for the definition of $\Gamma^{\on{ren}}(\Bun^{\on{level}_{\ul{x}}}_{N'},-)$)
receives a map from
$$\Gamma(\Bun_{N'},\oblv^l(\omega_{\Bun_{N'}}))\simeq \Gamma(\Bun_{N'},\CO_{\Bun_{N'}}),$$
and hence from 
$$k\to \Gamma(\Bun_{N'},\CO_{\Bun_{N'}}).$$

Furthermore, it is easy to see that the resulting map 
$$k\to \Gamma^{\on{ren}}(\Bun^{\on{level}_{\ul{x}}}_{N'},\omega_{\Bun_{N'}})$$
in $\Vect$ upgrades to a map
$$k\to \Gamma^{\on{ren}}(\Bun^{\on{level}_{\ul{x}}}_{N'},\omega_{\Bun_{N'}})^{\on{enh}}=\Gamma_{N',\ul{x}}(\omega_{\Bun_{N'}}).$$
in $\KL(N')_{\ul{x}}$.

\medskip

Making $\ul{x}$ move in families over $\Ran$, we obtain a map
$$k_\Ran\to \Gamma_{N',\Ran}(\omega_{\Bun_{N'}}).$$

\sssec{}

Thus, we obtain a map
\begin{equation} \label{e:unip localization 3}
\langle k_\Ran,-\rangle_{\KL(N'),\Ran} \otimes \fl_{N'}[\delta_{N'}]\to \on{C}^\cdot_\dr(\Bun_{N'},-)\circ \Loc_{N'}.
\end{equation} 

\medskip

Finally, we note that the functor 
$$\langle k_\Ran,-\rangle_{\KL(N'),\Ran}:\KL(N')_\Ran\to \Vect$$
identifies with
$$\KL(N')_\Ran \overset{\BRST_{\fn'}}\longrightarrow \Dmod(\Ran)
\overset{\on{C}^\cdot_c(\Ran,-)}\longrightarrow \Vect.$$

Combining with \eqref{e:unip localization 3} we obtain the desired map
\eqref{e:unip localization}. 

\sssec{}

We will prove:\footnote{This result is established in \cite[Theorem 4.0.5(4)]{CF}. We will provide proof for completeness.} 

\begin{thm} \label{t:unip localization}
The natural transformation \eqref{e:unip localization} becomes an isomorphism after precomposing with the endofunctor
$$\int \on{ins.vac}:\KL(N')_\Ran\to \KL(N')_\Ran.$$
\end{thm}

\sssec{}

Let us reformulate \thmref{t:unip localization} in more concrete terms. Note that we can rewrite the precomposition 
of \eqref{e:unip localization LHS} with $\int \on{ins.vac}$ as 
\begin{equation} \label{e:unip localization LHS rewrite}
\KL(N')_\Ran \overset{\BRST^{\on{enh}}_{\fn'}}\longrightarrow \Omega(\fn')\mod^{\on{fact}}_\Ran
\overset{\on{C}^{\on{fact}}_\cdot(X,\Omega(\fn'),-)}\longrightarrow \Vect
\overset{-\otimes \fl_{N'}[\delta_{N'}]}\longrightarrow \Vect.
\end{equation} 

Further, since the functor $\Loc_{N'}$ has a unital structure, the precomposition of 
\eqref{e:unip localization RHS} with $\int \on{ins.vac}$ is canonically isomorphic to
\eqref{e:unip localization RHS} itself. 

\medskip

Hence, \thmref{t:unip localization} implies:

\begin{cor} \label{c:unip localization}
There exists a canonical isomorphism between \eqref{e:unip localization LHS rewrite} and
the functor \eqref{e:unip localization RHS}. 
\end{cor}

\sssec{Variant} \label{sss:BRST twisted by char}

Consider the Lie algebra $\Gamma(X,\fn')$, and let $\chi^{\on{Lie}}$ be its character.
The datum of $\chi^{\on{Lie}}$ gives rise to a factorization character 
$$\chi:\fL(\fn')\to \BG_a,$$
trivial in $\fL^+(\fn')$, and a map
$$\chi^{\on{glob}}:\Bun_{N'}\to \BG_a.$$

\medskip

Let $\BRST_{\fn',\chi}$ be the $\chi$-twisted version of the semi-infinite cohomology functor, i.e.
$$\BRST_{\fn',\chi}(-)=\BRST_{\fn'}(-\otimes \chi).$$

Note that
$$\on{Vac}(N')\otimes \chi\simeq \on{Vac}(N').$$

Hence, 
$$\BRST_{\fn',\chi}(\on{Vac}(N'))\simeq \BRST_{\fn'}(\on{Vac}(N'))\simeq \Omega(\fn')$$
as factorization algebras. 

\medskip

Let $\BRST^{\on{enh}}_{\fn',\chi}$ be the enhancement of $\BRST_{\fn',\chi}$
$$\BRST^{\on{enh}}_{\fn',\chi}:\fL(\fn')\mod\to \Omega(\fn')\mod^{\on{fact}}.$$

\sssec{}

As in \eqref{e:unip localization} one constructs a natural transformation
\begin{equation} \label{e:unip localization chi}
\left(\on{C}^\cdot_c(\Ran,-) \circ \BRST_{\fn',\chi}\right) \otimes \fl_{N'}[\delta_{N'}]\to 
\on{C}^\cdot_\dr(\Bun_{N'},-\overset{*}\otimes \chi^*(\on{exp}))\circ \Loc_{N'}.
\end{equation} 

And parallel to \thmref{t:unip localization}, we have:

\begin{thm} \label{t:unip localization chi}
The map \eqref{e:unip localization chi} becomes an isomorphism after precomposing with 
$$\int \on{ins.vac}:\KL(N')_\Ran\to \KL(N')_\Ran.$$
\end{thm}

\begin{cor} \label{c:unip localization chi}
There exists a canonical isomorphism between 
\begin{equation} \label{e:unip localization chi LHS rewrite}
\KL(N')_\Ran \overset{\BRST^{\on{enh}}_{\fn',\chi}}\longrightarrow \Omega(\fn')\mod^{\on{fact}}_\Ran
\overset{\on{C}^{\on{fact}}_\cdot(X,\Omega(\fn'),-)}\longrightarrow \Vect
\overset{-\otimes \fl_{N'}[\delta_{N'}]}\longrightarrow \Vect
\end{equation} 
and 
\begin{equation} \label{e:unip localization chi RHS}
\KL(N')_\Ran\overset{\Loc_{N'}}\longrightarrow \Dmod(\Bun_{N'})\overset{-\overset{*}\otimes \chi^*(\on{exp})}\longrightarrow 
\Dmod(\Bun_{N'})
\overset{\on{C}^\cdot_\dr(\Bun_{N'},-)}\longrightarrow \Vect,
\end{equation} 
\end{cor}

\ssec{Proof of \thmref{t:unip localization}}

\sssec{}

Since $\fL^+(N')$ is pro-unipotent, the category $\Rep(\fL^+(G))$ is generated by 
$\one_{\Rep(\fL^+(N'))}$. And hence the category $\KL(N')$ is generated by $\on{Vac}(N')$. 
Therefore, $\KL(N')_\Ran$ is generated by $\on{Vac}(N')_\Ran$ as a $\Dmod(\Ran)$-module
category. 

\medskip

By unitality, for both sides in \thmref{t:unip localization}, tensoring the source by
an object $\CF\in \Dmod(\Ran)$ has the effect of tensoring the target by $\on{C}^\cdot_c(\Ran,\CF)$.
Hence, it is enough to show that the map in \thmref{t:unip localization} evaluates to an isomorphism on 
$\on{Vac}(N')_\Ran$.

\sssec{}

By construction, the left-hand side is 
$$\on{C}^{\on{fact}}_\cdot(X;\Omega(\fn'))$$
(see \secref{sss:vac fact hom} for the notation).

\medskip

By \eqref{e:Loc of Vac via dual}, we have
$$\Loc_{N'}(\on{Vac}(N)_\Ran)\simeq \ind^l(\CO_{\Bun_{N'}}).$$

We rewrite
$$\ind^l(\CO_{\Bun_{N'}})\simeq \ind^r(\omega_{\Bun_{N'}})\simeq
\ind^r(\CO_{\Bun_{N'}})\otimes \fl_{N'}[\delta_{N'}].$$

\medskip

Hence, the map in \thmref{t:unip localization} becomes a map
\begin{equation} \label{e:glob funct BunN}
\on{C}^{\on{fact}}_\cdot(X;\Omega(\fn'))\to \Gamma(\Bun_{N'},\CO_{\Bun_{N'}}).
\end{equation}

The fact that \eqref{e:glob funct BunN} is an isomorphism is well-known. For completeness,
we will supply a proof in the next subsection.

\sssec{}

The material in the rest of this subsection is not logically necessary, except for the example
considered in \secref{sss:ch homology Omega ab}. 

\medskip

According to Remark \ref{r:Omega as Chev}, the factorization algebra $\Omega(\fn')$ can be thought of as 
the factorization algebra associated to the cohomological Chevalley complex of a Lie-* algebra. Let us
consider $\on{C}^{\on{fact}}_\cdot(X;\Omega(\fn'))$ in this paradigm. 

\medskip

Let $L$ be a Lie-* algebra, whose underlying D-module is classical, finitely generated
and projective (as a D-module). Consider its cohomological Chevalley complex
$$\on{C}_{\on{chev}}^\cdot(L)\in \on{ComAlg}(\Dmod(X)).$$
Denote
$$\Omega(L):=\on{Fact}(\on{C}_{\on{chev}}^\cdot(L))\in \on{ComAlg}(\on{FactAlg}(X)),$$
and consider
$$\on{C}^{\on{fact}}_\cdot(X,\Omega(L))\in \on{ComAlg}(\Vect).$$

Unfortunately, we do not have a good grip on what $\on{C}^{\on{fact}}_\cdot(X,\Omega(L))$ looks
like. 

\sssec{}

Note that $\on{C}_{\on{chev}}^\cdot(L)$ is naturally written as 
\begin{equation} \label{e:chev as limit}
\on{C}_{\on{chev}}^\cdot(L)\simeq \underset{n}{\on{lim}}\, \on{C}_{\on{chev}}^\cdot(L)_n,
\end{equation} 
where $\on{C}_{\on{chev}}^\cdot(L)_n\in \on{ComAlg}(\Dmod(X))$ is the $n$-step cohomological Chevalley complex. 

\medskip

Note, however, that the assumptions on $L$ 
imply that for every $n$, the composite map
$$\tau^{\leq n}(\on{C}^\cdot(L))\to \on{C}^\cdot(L)\to \on{C}^\cdot(L)_n$$
is an isomorphism (here $\tau^{\leq n}$ refers to the left D-module structure, i.e., one for which $\oblv^l$ is t-exact). So, in the formation
$\on{C}_{\on{ch}}^\cdot(L)$ no ``actual completion" is involved. 

\sssec{}

Denote 
$$\Omega(L)_n:=\on{Fact}(\on{C}_{\on{chev}}^\cdot(L)_n)\in \on{ComAlg}(\on{FactAlg}(X)).$$

Set:
$$\on{C}^{\on{fact}}_\cdot(X,\Omega(L))^\wedge:=\underset{n}{\on{lim}}\, \on{C}^{\on{fact}}_\cdot(X,\Omega(L)_n).$$

Unlike $\on{C}^{\on{fact}}_\cdot(X,\Omega(L))$, the algebra $\on{C}^{\on{fact}}_\cdot(X,\Omega(L))^\wedge$ can be 
described explicitly. Namely, according to \cite[Proposition 7.4.1]{BD2}, 
$$\on{C}^{\on{fact}}_\cdot(X,\Omega(L))^\wedge\simeq 
\on{C}_{\on{chev}}^\cdot(\on{C}^\cdot_\dr(X,L))=\on{C}_{\on{chev},\cdot}(\on{C}^\cdot_\dr(X,L))^\vee,$$
where:

\begin{itemize}

\item $\on{C}^\cdot_\dr(X,L)$ is considered as a Lie algebra (in $\Vect$);

\smallskip

\item  $\on{C}_{\on{chev}}^\cdot(-)$ and $\on{C}_{\on{chev},\cdot}(-)$ are the cohomological and homological
Chevalley complexes of a Lie algebra, respectively. 

\end{itemize}

\sssec{}

The map $\to$ in \eqref{e:chev as limit} gives rise to a map
\begin{equation}  \label{e:compl of ch hom Omega}
\on{C}^{\on{fact}}_\cdot(X,\Omega(L)) \to \on{C}^{\on{fact}}_\cdot(X,\Omega(L))^\wedge 
\end{equation} 
but this map is in general \emph{not} an isomorphism. 

\sssec{Example} \label{sss:ch homology Omega ab}

Let $L$ be abelian. Denote $L^\vee:=\BD(L)[-1]$. Then
$$\on{C}_{\on{chev}}^\cdot(L)\simeq \Sym^!(\BD(L)),$$
and hence
$$\on{C}^{\on{fact}}_\cdot(X,\Omega(L))\simeq \Sym(\on{C}^\cdot_\dr(X,\BD(L)))\simeq
\Sym(\on{C}^\cdot_\dr(X,L)^\vee[-1]).$$

By contrast,
$$\on{C}_{\on{chev}}^\cdot(\on{C}^\cdot_\dr(X,L))\simeq \Sym(\on{C}^\cdot_\dr(X,L)[-1])^\vee.$$

So the difference between the two sides in \eqref{e:compl of ch hom Omega} in this case is
that between a polynomial algebra and its completion.

\sssec{}

The map \eqref{e:compl of ch hom Omega} is the best approximation to $\on{C}^{\on{fact}}_\cdot(X,\Omega(L))$
that we have in general. In certain situations, it allows us to recover $\on{C}^{\on{fact}}_\cdot(X,\Omega(L))$
completely. 

\medskip

This happens, for example, if $L$ carries a strictly positive grading. In this case, the map \eqref{e:compl of ch hom Omega}
defines an isomorphism on each graded component. I.e., we have 
$$\on{C}^{\on{fact}}_\cdot(X,\Omega(L))^d\simeq (\on{C}^{\on{fact}}_\cdot(X,\Omega(L))^\wedge)^d\simeq 
\on{C}^\cdot_{\on{chev}}(\on{C}^\cdot_\dr(X,L))^d \simeq 
(\on{C}_{\on{chev},\cdot}(\on{C}^\cdot_\dr(X,L))^{-d})^\vee, \quad
d\in \BZ^{<0}.$$

Note, however, that $\on{C}^{\on{fact}}_\cdot(X,\Omega(L))^\wedge$ is \emph{not} the direct sum of its graded pieces.
Rather,
$$(\on{C}^{\on{fact}}_\cdot(X,\Omega(L))^\wedge)^d=\underset{n}{\on{lim}}\, \on{C}^{\on{fact}}_\cdot(X,\Omega(L)_n)^d.$$

\begin{rem}

Using cohomological truncations on powers of $X$ one can show: 

\smallskip 

\noindent{(i)} $\on{C}^{\on{fact}}_\cdot(X,\Omega(L))$ is coconnective. 

\smallskip 

\noindent{(ii)} Suppose that $L$ is such that $H^0(\on{C}^\cdot_\dr(X,L))=0$. Then $\on{C}^{\on{fact}}_\cdot(X,\Omega(L))$
is classical. 

\end{rem} 

\ssec{Global functions on the moduli space of bundles for a unipotent group-scheme} \label{ss:as good as affine}

\sssec{}

Let $\CY$ be a D-prestack over $X$. Let $A\in \on{ComAlg}(\Dmod(X))$ be the algebra of global functions
on $\CY$, and let $\CA=\on{Fact}(A)\in \on{ComAlg}(\on{FactAlg}(X))$ be the corresponding factorization algebra
(see \secref{sss:com fact vs Dmod com}). 

\medskip

Recall that by \eqref{e:fact hom A on Ran}, the evaluation map $$\on{Sect}_\nabla(X,\CY)\times X \to \CY$$
gives rise to a map
\begin{equation} \label{e:glob funct sect}
\on{C}^{\on{fact}}_\cdot(X,\CA)\to \Gamma(\on{Sect}_\nabla(X,\CY),\CO_{\on{Sect}_\nabla(X,\CY)}),
\end{equation}

\medskip

Unwinding the construction, it is easy to see that the map \eqref{e:glob funct BunN} is the map
\eqref{e:glob funct sect} for $\CY=\on{Jets}(\on{pt}/N')$, where 
$$\Omega(\fn')=\CO_{\on{pt}/\fL^+(N')} \overset{\text{\eqref{e:arcs as fact}}}\simeq \on{Fact}(\CO_{\on{Jets}(\on{pt}/N')})$$
and 
$$\on{Sect}_\nabla(X,\on{Jets}(\on{pt}/N'))\simeq \on{Sect}(X,\on{pt}/N')\simeq \Bun_{N'}.$$

\medskip

Thus, we need to show that \eqref{e:glob funct sect} is an isomorphism in this case.

\begin{rem} 

Recall that \eqref{e:glob funct sect} is an isomorphism for $\CY$ that is affine over $X$,
see \propref{p:hor sect}. So we want to prove that $\CY=\on{Jets}(\on{pt}/N')$ is not too
different from the affine case. 

\medskip

The proof below follows closely that of \propref{p:hor sect}. 

\end{rem} 

\sssec{}

First, we will show that the map \eqref{e:glob funct sect} induces an \emph{isomorphism} 
\begin{multline} \label{e:Maps into Bun N'}
\Maps(\Spec(R),\on{Sect}_\nabla(X,\on{Jets}(\on{pt}/N'))) \to \\
\to \Maps_{\on{ComAlg}(\Vect)}(\Gamma(\on{Sect}_\nabla(X,\on{Jets}(\on{pt}/N')),\CO_{\on{Sect}_\nabla(X,\on{Jets}(\on{pt}/N'))}),R)\to \\
\to \Maps_{\on{ComAlg}(\Vect)}(\on{C}^{\on{fact}}_\cdot(X,\CO_{\on{pt}/\fL^+(N')}),R)
\end{multline} 
for $R\in \on{ComAlg}(\Vect^{\leq 0})$.

\medskip

This essentially follows from the fact that $\on{Jets}(N')$ is pro-unipotent, 
and hence $\on{Jets}(\on{pt}/N')$ is \emph{as good as affine}\footnote{Up to issues of renormalization, which are irrelevant here.} 
(see \secref{sss:as good as affine} below for what this means), i.e., for 
$\CR\in \on{ComAlg}(\Dmod(X)^{\leq 0})$, the map 
\begin{equation} \label{e:express sect}
\Maps_{X,\nabla}(\Spec_X(\CR),\on{Jets}(\on{pt}/N')) \to \Maps_{\on{ComAlg}(\Dmod(X))}(\CO_{\on{Jets}(\on{pt}/N')},\CR)
\end{equation}
is an isomorphism. 

\sssec{}

In more detail, for $R\in \on{ComAlg}(\Vect^{\leq 0})$ we have: 
\begin{multline} \label{e:express sect 1}
\Maps(\Spec(R),\on{Sect}_\nabla(X,\on{Jets}(\on{pt}/N')))=
\Maps_{X,\nabla}(\Spec(R)\times X,\on{Jets}(\on{pt}/N'))) \overset{\text{\eqref{e:express sect}}}\simeq \\
\simeq \Maps_{\on{ComAlg}(\Dmod(X))}(\CO_{\on{Jets}(\on{pt}/N')},R\otimes \CO_X).
\end{multline} 

\medskip

Using \corref{c:ch homology as left adj}, we rewrite the expression in the right-hand side in \eqref{e:express sect 1} as
\begin{equation} \label{e:express sect 2}
\Maps_{\on{ComAlg}(\Vect)}(\on{C}^{\on{fact}}_\cdot(X,\CO_{\on{pt}/\fL^+(N')}),R). 
\end{equation}

\medskip

Combining \eqref{e:express sect 1} and \eqref{e:express sect 2}, we obtain an isomorphism
$$\Maps(\Spec(R),\on{Sect}_\nabla(X,\on{Jets}(\on{pt}/N'))) \simeq 
\Maps_{\on{ComAlg}(\Vect)}(\on{C}^{\on{fact}}_\cdot(X,\CO_{\on{pt}/\fL^+(N')}),R),$$
and unwinding the definitions we obtain that this isomorphism equals the map in \eqref{e:Maps into Bun N'}. 

\sssec{} \label{sss:as good as affine}

Let us call a prestack $\CZ$ \emph{as good as affine} if the functor
$$\Gamma(\CZ,-): \QCoh(\CZ)\to \CO_{\CZ}\mod$$
is an equivalence.

\sssec{} \label{sss:as good as connective}

Let $R$ be an object of $\on{ComAlg}(\Vect)$. Define the prestack $``\Spec(R)"$ by
$$\Maps(\Spec(R'), ``\Spec(R)"):= \Maps_{\on{ComAlg}(\Vect)}(R,R'), \quad R'\in \on{ComAlg}(\Vect^{\leq 0}).$$

Note the formation of $``\Spec(R)"$ is functorial in $R$. In particular, we obtain a map
\begin{multline} \label{e:functions on fake Spec}
R\simeq \Maps_{\on{ComAlg}(\Vect)}(k[t],R) \to \Maps(``\Spec(R)",\Spec(k[t]))\simeq \\
\simeq \Maps(``\Spec(R)",\BA^1)\simeq 
\Gamma(``\Spec(R)", \CO_{``\Spec(R)"}).
\end{multline} 

\medskip

We shall say that $R$ is ``as good as connective" if:

\begin{itemize}

\item The prestack $``\Spec(R)"$ is as good as affine;

\item The map \eqref{e:functions on fake Spec} is an isomorphism.

\end{itemize}

\sssec{}

Thus, given that \eqref{e:Maps into Bun N'} is an isomorphism, we need to show that $\on{C}^{\on{fact}}_\cdot(X,\Omega(\fn'))$
is as good as connective.

\medskip

We will now use the assumption that $N'$ admits a filtration by normal subgroups with abelian quotients. We will argue
by induction on the length of such a filtration.

\sssec{} \label{sss:vector group}

We first consider the base of the induction, i.e., case when $N'$ is a vector group-scheme, i.e., is the total space of a vector bundle 
$\CE$ on $X$. In this case, the computation of $\on{C}^{\on{fact}}_\cdot(X,\Omega(\fn'))$ has been performed in 
\secref{sss:ch homology Omega ab}. 

\medskip

We obtain that $\on{C}^{\on{fact}}_\cdot(X,\Omega(\fn'))$ is (non-canonically) isomorphic to the tensor product
$$\Sym(V_1)\otimes \Sym(V_2[-1]),$$
where $V_1$ and $V_2$ are classical finite-dimensional vector spaces. 

\medskip

It is clear that the tensor product of two algebras both of which are as good as connective is itself 
as good as connective. Hence, it remains to see that algebras of the form $\Sym(V[-1])$, where 
$V$ is a classical finite-dimensional vector space, are as good as connective.

\medskip

However, this is well-known: in this case 
$$``\Spec(\Sym(V[-1]))"\simeq \on{pt}/V^\vee$$
and the assertion is manifest. 

\sssec{}

We now perform the induction step. Thus, we fix a short exact sequence
$$1\to N'_2\to N'\to N'_1\to 1,$$
where $N'_2$ is a vector group-scheme.

\medskip

We observe:

\begin{lem}
Let $R_1\to R$ be a map of commutative algebras in $\Vect$. Assume that:

\begin{itemize}

\item $R_1$ is as good as connective;

\item For any homomorphism $R_1\to R'$ with $R'$ connective, the base change $R'\underset{R_1}\otimes R$
is as good as connective.

\end{itemize}

Then $R$ is as good as connective. 

\end{lem}

\medskip

We apply this lemma to
$$R:=\on{C}^{\on{fact}}_\cdot(X,\Omega(\fn')) \text{ and } R_1:=\on{C}^{\on{fact}}_\cdot(X,\Omega(\fn'_1)).$$

By the induction hypothesis, $\on{C}^{\on{fact}}_\cdot(X,\Omega(\fn'_1))$ is as good as connective.
Hence, it remains to show that for any connective $R'$ and a homomorphism 
$$\on{C}^{\on{fact}}_\cdot(X,\Omega(\fn'_1))\to R',$$
the algebra
$$R'\underset{\on{C}^{\on{fact}}_\cdot(X,\Omega(\fn'_1))}\otimes \on{C}^{\on{fact}}_\cdot(X,\Omega(\fn'))$$
is as good as connective.

\sssec{}

We now apply \lemref{l:rel fact homology}, and hence we can rewrite
$$R'\underset{\on{C}^{\on{fact}}_\cdot(X,\Omega(\fn'_1))}\otimes \on{C}^{\on{fact}}_\cdot(X,\Omega(\fn'))\simeq 
\on{C}^{\on{fact}}_\cdot(X,\Omega(\fn')_{R'}),$$
where 
$$\Omega(\fn')_{R'}:=\Omega(\fn')\underset{\Omega(\fn'_1)}\otimes (R'\otimes \omega_\Ran)\in
\on{ComAlg}(\on{FactAlg}(X)\otimes R'\mod).$$

\sssec{}

Recall that 
$$\Omega(\fn')\simeq \on{Fact}(\on{C}^\cdot_{\on{chev}}(L_{\fn'})) \text{ and } \Omega(\fn'_1)\simeq \on{Fact}(\on{Chev}^\cdot_{\on{chev}}(L_{\fn'_1})),$$
where
$$L_{\fn'}=\fn'\otimes \on{D}_X, \quad L_{\fn'_1}=\fn'_1\otimes \on{D}_X.$$

Hence, we can rewrite
$$\Omega(\fn')_{R'}\simeq
\on{Fact}\left(\on{C}^\cdot_{\on{chev}}(L_{\fn'})\underset{\on{C}^\cdot_{\on{chev}}(L_{\fn'_1})}\otimes (R'\otimes \CO_X)\right).$$

\sssec{}

We can interpret the datum of $\on{C}^{\on{fact}}_\cdot(X,\Omega(\fn'_1))\to R'$ as a map
$$\Spec(R')\to \Bun_{N'_1}.$$

The adjoint action of $N'_1$ on $N'_2$ gives rise to an $R'$-family of twisted forms of $N'_2$,
denoted $N'_{2,R'}$. Consider the corresponding $R'$-family of Lie-* algebras 
$$L_{\fn'_{2,R'}}:=\fn'_{2,R'}\otimes \on{D}_X.$$

We have
$$\on{C}^\cdot_{\on{chev}}(L_{\fn'})\underset{\on{C}^\cdot_{\on{chev}}(L_{\fn'_1})}\otimes (R'\otimes \CO_X)\simeq
\on{C}^\cdot_{\on{chev}}(L_{\fn'_{2,R'}}).$$

Thus, we can consider
$$\Omega(\fn'_{2,R'})\in \on{ComAlg}(\on{FactAlg}(X)\otimes R'\mod)$$
and we obtain:
$$\Omega(\fn')_{R'}\simeq \on{C}^{\on{fact}}_\cdot(X,\Omega(\fn'_{2,R'})).$$

\sssec{}

It remains to show that $\on{C}^{\on{fact}}_\cdot(X,\Omega(\fn'_{2,R'}))$ is as good as connective. 

\medskip

Recall that $\fn'_2$ is abelian. Hence, $N'_{2,R'}$ is $R'$-family of vector group-schemes. Hence, the required
assertion is a relative (over $\Spec(R')$) version of the case considered in \secref{sss:vector group}.

\qed[\thmref{t:unip localization}]

\ssec{Application: integration over (twists of) \texorpdfstring{$\Bun_N$}{BunN} via BRST}

\sssec{}

Let $\CP_T$ be a $T$-bundle on $X$. Consider the unipotent group-scheme $N_{\CP_T}$.
Denote the corresponding moduli stack $\Bun_{N_{\CP_T}}$; note that it identifies with 
$\Bun_{N,\CP_T}$ (see \eqref{e:twisted Bun N}). 

\medskip

The resulting map 
$$\Bun_{N_{\CP_T}}\simeq \Bun_{N,\CP_T} \overset{\fp}\to \Bun_G$$ can be thought of as 
$$\Bun_{N_{\CP_T}}\to \Bun_{G_{\CP_T}}\overset{\alpha_{\CP_T,\on{taut}}^{-1}}\longrightarrow \Bun_G.$$

\sssec{}

Since the restriction of $\kappa$ to $\fn$ is trivial, we obtain that the restriction of the twisting $\CT_\kappa$ 
along $\fp$ is canonically trivial. In particular, we have a well-defined functor 
\begin{equation} \label{e:pullback to BunN CPT kappa}
\fp^!_\kappa:\Dmod_\kappa(\Bun_G)\to \Dmod(\Bun_{N,\CP_T}).
\end{equation}

\sssec{}

Note that the embedding
$$N_{\CP_T}\to G_{\CP_T}$$
gives rise to a map 
$$\fL(\fn_{\CP_T})\to \hg_{\kappa,\CP_T},$$
and this map lifts to the Kac-Moody extension.

\medskip

In particular, we obtain a well-defined restriction functor
$$\KL(G)_\kappa \overset{\alpha_{\CP_T},\on{taut}}\longrightarrow 
\KL(G)_{\kappa,\CP_T} \to \KL(N_{\CP_T}).$$

Denote
$$\Omega(\fn_{\CP_T},\fg)_\kappa:=\BRST_{\fn_{\CP_T}}(\on{Vac}(G)_{\kappa,\CP_T}).$$
This is a factorization algebra, which receives a homomorphism from $\Omega(\fn_{\CP_T})$.

\medskip

Thus, the composition
$$\KL(G)_\kappa \overset{\alpha_{\CP_T},\on{taut}}\longrightarrow 
\KL(G)_{\kappa,\CP_T} \to \KL(N_{\CP_T})\overset{\BRST^{\on{enh}}_{\fn_{\CP_T}}}\longrightarrow 
\Omega(\fn_{\CP_T})\mod^{\on{fact}}$$
further enhances to a (factorization) functor
$$\BRST^{\fg\!\on{-enh}}_{\fn_{\CP_T}}:\KL(G)_\kappa\to \Omega(\fn_{\CP_T},\fg)_\kappa\mod^{\on{fact}}.$$

\sssec{}

We are going to prove: 

\begin{thm} \label{t:int loc over BunN kappa}
The composition 
\begin{equation} \label{e:int loc over BunN kappa 1}
\KL(G)_{\kappa,\Ran} \overset{\Loc_{G,\kappa}}\longrightarrow \Dmod_\kappa(\Bun_G)\overset{\fp^!_\kappa}
\to \Dmod(\Bun_{N,\CP_T}) \overset{\on{C}^\cdot_\dr(\Bun_{N,\CP_T},-)}\longrightarrow \Vect
\end{equation} 
identifies with the functor
\begin{equation} \label{e:int loc over BunN kappa 2}
\KL(G)_{\kappa,\Ran} \overset{\BRST^{\fg\!\on{-enh}}_{\fn_{\CP_T}}}\longrightarrow 
\Omega(\fn_{\CP_T},\fg)\mod^{\on{fact}}_\Ran 
\overset{\on{C}^{\on{fact}}_\cdot(X;\Omega(\fn_{\CP_T},\fg)_\kappa,-)}\longrightarrow 
\Vect \overset{-\otimes \fl_{N_{\CP_T}}[\delta_{N_{\CP_T}}]}\longrightarrow \Vect,
\end{equation} 
where the notations $\delta_{N_{\CP_T}}$ and $\fl_{N_{\CP_T}}$ are as in \secref{sss:l N' line}. 
\end{thm}

%
%

\sssec{Proof of \thmref{t:int loc over BunN kappa}}

First, we rewrite the functor
$$\KL(G)_{\kappa,\Ran} \overset{\Loc_{G,\kappa}}\longrightarrow \Dmod_\kappa(\Bun_G)\overset{\fp^!_\kappa}
\to \Dmod(\Bun_{N,\CP_T})$$
using (a $\CP_T$-twisted version) of \corref{c:localization and restr fact hom Int}.

\medskip

We obtain that it identifies with
\begin{multline*}
\KL(G)_{\kappa,\Ran} \overset{\alpha_{\CP_T},\on{taut}}\longrightarrow 
\KL(G)_{\kappa,\CP_T,\Ran} \overset{\int \on{ins.vac}}\longrightarrow \KL(G)_{\kappa,\CP_T,\Ran}\to \\
\to \KL(N_{\CP_T})_{\Ran}\overset{\Loc_{N_{\CP_T}}}\longrightarrow 
\Dmod(\Bun_{N,\CP_T}).
\end{multline*}

By \thmref{t:unip localization}, the functor
$$\KL(N_{\CP_T})_{\Ran}\overset{\Loc_{N_{\CP_T}}}\longrightarrow 
\Dmod(\Bun_{N,\CP_T}) \overset{\on{C}^\cdot_\dr(\Bun_{N,\CP_T},-)}\longrightarrow \Vect$$
identifies with
$$\KL(N_{\CP_T})_{\Ran}\overset{\BRST^{\on{enh}}_{\fn_{\CF_T}}}\longrightarrow
\Omega(\fn_{\CP_T})\mod^{\on{fact}}_\Ran
\overset{\on{C}^{\on{fact}}_\cdot(X;\Omega(\fn_{\CP_T}),-)} \longrightarrow \Vect 
\overset{-\otimes \fl_{N_{\CP_T}}[\delta_{N_{\CP_T}}]}\longrightarrow \Vect.$$

The assertion of the theorem follows now by applying \corref{c:A B lemma}.

\qed[\thmref{t:int loc over BunN kappa}]

\sssec{}

Note that the same proof applies in the situation twisted by a character. Namely,  
$\chi^{\on{Lie}}$ be a character of $\Gamma(X,\fn_{\CP_T})$ as in \secref{sss:BRST twisted by char}.

\medskip

Denote 
$$\Omega(\fn_{\CP_T},\chi,\fg)_\kappa:=\BRST_{\fn_{\CP_T},\chi}(\on{Vac}(G)_{\kappa,\CP_T}).$$

Consider the corresponding functor
$$\BRST^{\fg\!\on{-enh}}_{\fn_{\CP_T},\chi}:\KL(G)_\kappa\to \Omega(\fn_{\CP_T},\chi,\fg)_\kappa\mod^{\on{fact}}.$$

Then:
\begin{thm} \label{t:int loc over BunN kappa chi}
The composition 
\begin{multline*} 
\KL(G)_{\kappa,\Ran} \overset{\Loc_{G,\kappa}}\longrightarrow \Dmod_\kappa(\Bun_G)\overset{\fp^!_\kappa}
\to \Dmod(\Bun_{N,\CP_T}) \overset{-\overset{*}\otimes \chi^*(\on{exp})}\longrightarrow \\
\to \Dmod(\Bun_{N,\CP_T})
\overset{\on{C}^\cdot_\dr(\Bun_{N,\CP_T},-)}\longrightarrow \Vect
\end{multline*} 
identifies with the functor
$$\KL(G)_{\kappa,\Ran} \overset{\BRST^{\fg\!\on{-enh}}_{\fn_{\CP_T},\chi}}\longrightarrow 
\Omega(\fn_{\CP_T},\chi,\fg)\mod^{\on{fact}}_\Ran 
\overset{\on{C}^{\on{fact}}_\cdot(X;\Omega(\fn_{\CP_T},\chi,\fg)_\kappa,-)}\longrightarrow 
\Vect \overset{-\otimes \fl_{N_{\CP_T}}[\delta_{N_{\CP_T}}]}\longrightarrow \Vect.$$
\end{thm}

\section{Localization via the infinitesimal Hecke groupoid}  \label{s:loc Hecke gpd} 

This goal of this section is to prove \thmref{t:localization and oblv orig}. This will 
be based on the approach to the localization functor via the infinitesimal Hecke
groupoid, which was developed in the unpublished part of the thesis of the eighth
author of this paper. 

\medskip

We will then use these ideas to prove that the functor $\Loc_{G,\kappa}$ is 
\emph{almost} a localization (i.e., its right adjoint is fully faithful). 
Namely, it becomes a localization when we compose
it with restriction to any quasi-compact open substack, see \thmref{t:Loc is loc}.

\ssec{Another take on the functor \texorpdfstring{$\Gamma_G$}{Gamma}}

\sssec{} \label{sss:KL via local inf Hecke}

Recall the local Hecke stack completed along the diagonal, viewed as a groupoid acting on $\on{pt}/\fL^+(G)$:
$$\on{pt}/\fL^+(G)\overset{\hl^{\on{loc},\wedge}}\leftarrow \on{Hecke}_G^{\on{loc},\wedge} \overset{\hr^{\on{loc},\wedge}}\to \on{pt}/\fL^+(G).$$

Recall also that the datum of a level $\kappa$ gives rise to a \emph{multiplicative} line bundle $\CL^{\on{loc}}_\kappa$ on 
$\on{Hecke}_G^{\on{loc},\wedge}$, see \secref{sss:mult line bundle local}. 

\medskip

According to \cite[Sect. 3.3]{CF}, we can identify the category $\KL(G)_\kappa$ with the category 
$$\Rep(\fL^+(G))^{\on{Hecke}^{\on{loc},\wedge}_G,\CL^{\on{loc}}_\kappa}$$
of
$\CL^{\on{loc}}_\kappa$-twisted $\on{Hecke}^{\on{loc},\wedge}_G$-equivariant objects in $\Rep(\fL^+(G))$.

\medskip

I.e., this is the category of $\CM\in \Rep(\fL^+(G))$ equipped with an isomorphism
\begin{equation} \label{e:QCoh equiv}
(\hl^{\on{loc},\wedge})^*(\CM)\simeq  \CL^{\on{loc}}_\kappa\otimes (\hr^{\on{loc},\wedge})^*(\CM)
\end{equation} 
(the isomorphism taking place in $\QCoh(\on{Hecke}^{\on{loc},\wedge}_G)$), 
equipped with a homotopy-coherent system of compatibilities. 

\sssec{}

Let $\on{inf}(\Bun_G)$ denote the infinitesimal groupoid of $\Bun_G$, i.e.,
$$\on{inf}(\Bun_G):=(\Bun_G\times \Bun_G)^\wedge,$$
where $(-)^\wedge$ means formal completion along the diagonal.  Consider the corresponding diagram
$$\Bun_G \overset{\hl^{\on{inf}}}\longleftarrow \on{inf}(\Bun_G)
\overset{\hr^{\on{inf}}}\longrightarrow \Bun_G.$$

\medskip

Note that the datum of the de Rham twisting $\CT_\kappa$ gives rise to a multiplicative line bundle,
to be denoted $\CL^{\on{inf}}_\kappa$ on $\on{inf}(\Bun_G)$. 

\medskip

Let $\CZ$ be a prestack mapping to $\Ran$. Since $\Bun_G$ is eventually coconnective, 
it follows from \cite[Proposition 3.4.3]{GaRo2} that we have a canonical equivalence
$$\Dmod_\kappa(\Bun_G)\otimes \Dmod(\CZ)
\simeq \left(\QCoh(\Bun_G)\otimes \Dmod(\CZ)\right)^{\on{inf}(\Bun_G)\times \CZ,\CL^{\on{inf}}_\kappa}$$
that intertwines the functor 
$$(\oblv^l_\kappa\otimes \on{Id}):\Dmod_\kappa(\Bun_G)\otimes \Dmod(\CZ)\to \QCoh(\Bun_G)\otimes \Dmod(\CZ)$$
with the tautological forgetful functor 
$$\left(\QCoh(\Bun_G)\otimes \Dmod(\CZ)\right)^{\on{inf}(\Bun_G)\times \CZ,\CL^{\on{inf}}_\kappa}\to 
\QCoh(\Bun_G)\otimes \Dmod(\CZ).$$

The same applies when we replace $\Bun_G$ by its open substack $U$. 

\sssec{} \label{sss:inf Hecke to inf}

Note now that we have a tautological map of groupoids
\begin{equation} \label{e:Hecke to inf}
\on{Hecke}^{\on{glob},\wedge}_{G,\CZ}\to \on{inf}(\Bun_G)\times \CZ
\end{equation} 
over $\CZ$. 

\medskip

By construction, the pullback of $\CL^{\on{inf}}_\kappa$ along \eqref{e:Hecke to inf} identifies canonically
with $\CL^{\on{glob}}_\kappa$ as a multiplicative line bundle. 

\medskip

From here we obtain that *-pullback along \eqref{e:Hecke to inf} defines a functor
\begin{multline} \label{e:Dmod to glob inf Hecke equiv}
\Dmod_\kappa(\Bun_G)\otimes \Dmod(\CZ)
\simeq \left(\QCoh(\Bun_G)\otimes \Dmod(\CZ)\right)^{\on{inf}(\Bun_G)\times \CZ,\CL^{\on{inf}}_\kappa}\to \\
\to \left(\QCoh(\Bun_G)\otimes \Dmod(\CZ)\right)^{\on{Hecke}^{\on{glob},\wedge}_{G,\CZ},\CL^{\on{glob}}_\kappa},
\end{multline}
that intertwines the forgetful functor
$$(\oblv^l_\kappa\otimes \on{Id}):\Dmod_\kappa(\Bun_G)\otimes \Dmod(\CZ)\to \QCoh(\Bun_G)\otimes \Dmod(\CZ)$$
with the tautological forgetful functor
$$\left(\QCoh(\Bun_G)\otimes \Dmod(\CZ)\right)^{\on{Hecke}^{\on{glob},\wedge}_{G,\CZ},\CL^{\on{glob}}_\kappa}\to
\QCoh(\Bun_G)\otimes \Dmod(\CZ).$$

Denote the functor \eqref{e:Dmod to glob inf Hecke equiv} by
$$\oblv_{\on{inf}\to \on{Hecke}^\wedge,\CZ}.$$

\sssec{}

Let $U\subset \Bun_G$ be a quasi-compact open substack. Note that 
it makes sense to restrict $\on{Hecke}^{\on{glob},\wedge}_{G,\CZ}$ to $U$:
$$(\hl^{\on{glob},\wedge})^{-1}(U\times \CZ)=:\on{Hecke}_{G,\CZ,U}^{\on{glob},\wedge}=:
(\hr^{\on{glob},\wedge})^{-1}(U\times \CZ).$$

\medskip

The contents of \secref{sss:inf Hecke to inf} apply over $U$ as well, and we obtain a functor, denoted 
$$\oblv_{\on{inf}\to \on{Hecke}^\wedge,\CZ,U}$$
that maps
\begin{multline} \label{e:Dmod to glob inf Hecke equiv U}
\Dmod_\kappa(U)\otimes \Dmod(\CZ)
\simeq \left(\QCoh(U)\otimes \Dmod(\CZ)\right)^{\on{inf}(\Bun_U)\times \CZ,\CL^{\on{inf}}_\kappa}\to \\
\to \left(\QCoh(U)\otimes \Dmod(\CZ)\right)^{\on{Hecke}^{\on{glob},\wedge}_{G,\CZ,U},\CL^{\on{glob}}_\kappa}
\end{multline}
and that intertwines the forgetful functor
$$(\oblv^l_\kappa\otimes \on{Id}):\Dmod_\kappa(U)\otimes \Dmod(\CZ)\to \QCoh(U)\otimes \Dmod(\CZ)$$
with the tautological forgetful functor
$$\left(\QCoh(U)\otimes \Dmod(\CZ)\right)^{\on{Hecke}^{\on{glob},\wedge}_{G,\CZ,U},\CL^{\on{glob}}_\kappa}\to
\QCoh(U)\otimes \Dmod(\CZ).$$

\sssec{} \label{sss:from KM to inf Hecke U}

Assume now that $U$ is quasi-compact. Let $j$ denote its embedding into $\Bun_G$. 
Consider the diagram:

\begin{equation} \label{e:log glob Hecke diagram compl U}
\CD
U\times \CZ @<{\hl^{\on{glob},\wedge}}<< \on{Hecke}_{G,\CZ,U}^{\on{glob},\wedge} @>{\hr^{\on{glob},\wedge}}>> U\times \CZ \\
@V{\on{ev}_{\CZ,U}}VV @VVV @VV{\on{ev}_{\CZ,U}}V \\
(\on{pt}/\fL^+(G))_\CZ @<{\hl^{\on{loc},\wedge}}<< \on{Hecke}_{G,\CZ}^{\on{loc},\wedge}  @>{\hr^{\on{loc},\wedge}}>> (\on{pt}/\fL^+(G))_\CZ,
\endCD
\end{equation} 
where $\on{ev}_{\CZ,U}:=\on{ev}_{\CZ}\circ j$.

\medskip

Since both squares in \eqref{e:log glob Hecke diagram compl U} are Cartesian, we obtain that the 
functor\footnote{Note that as in Remark \ref{r:ignore ren Rep L+G for _*}, the difference between 
$\Rep(\fL^+(G)$ and $\QCoh(\on{pt}/\fL^+(G))$ is immaterial here.}
$$(\on{ev}_{\CZ,U})_*:\QCoh(U)\otimes \Dmod(\CZ)\to \Rep(\on{pt}/\fL^+(G))_\CZ$$
gives rise to a functor 
\begin{equation} \label{e:global inf Hecke equiv to local inf Hecke equiv U}
\left(\QCoh(U)\otimes \Dmod(\CZ)\right)^{\on{Hecke}^{\on{glob},\wedge}_{G,\CZ,U},\CL^{\on{glob}}_\kappa} \to
\left(\Rep(\fL^+(G))_\CZ\right)^{\on{Hecke}^{\on{loc},\wedge}_{G,\CZ},\CL^{\on{loc}}_\kappa}
\overset{\text{\secref{sss:KL via local inf Hecke}}}\simeq \KL(G)_{\kappa,\CZ}.
\end{equation} 

We will denote the functor \eqref{e:global inf Hecke equiv to local inf Hecke equiv U} by
$$(\on{ev}_{\CZ,U})^{\on{Hecke}^\wedge\on{-enh}}_*.$$

\sssec{} 

Composing, we obtain that
the functor
$$(\on{ev}_{\CZ,U})_*\circ (\oblv^l_\kappa \otimes \on{Id}):\Dmod_\kappa(U)\otimes \Dmod(\CZ) \to \Rep(\fL^+(G))_\CZ$$
lifts to a functor
\begin{equation} \label{e:Dmod to local inf Hecke equiv}
\Dmod_\kappa(U)\otimes \Dmod(\CZ) \overset{(\on{ev}_{\CZ,U})^{\on{Hecke}^\wedge\on{-enh}}_* \circ \oblv_{\on{inf}\to \on{Hecke}^\wedge,\CZ,U}}
\longrightarrow \KL(G)_{\kappa,\CZ}.
\end{equation}

\sssec{} \label{sss:alt Gamma U}

It follows from the construction of the identification 
\begin{equation} \label{e:KL via local inf Hecke}
\left(\Rep(\fL^+(G))_\CZ\right)^{\on{Hecke}^{\on{loc},\wedge}_{G,\CZ},\CL^{\on{loc}}_\kappa} \simeq \KL(G)_{\kappa,\CZ}
\end{equation}
that the functor $$(\on{ev}_{\CZ,U})^{\on{Hecke}^\wedge\on{-enh}}_* \circ \oblv_{\on{inf}\to \on{Hecke}^\wedge,\CZ,U}$$
of \eqref{e:Dmod to local inf Hecke equiv} identifies canonically with the functor
$$\Gamma_{G,\kappa,\CZ,U}:=\Gamma_{G,\kappa,\CZ}\circ j_{*,\on{co}},$$
so that the diagram
$$
\CD
\oblv^{(\hg,\fL^+(G))_\kappa}_{\fL^+(G)}\circ (\on{ev}_{\CZ,U})^{\on{Hecke}^\wedge\on{-enh}}_* \circ \oblv_{\on{inf}\to \on{Hecke}^\wedge,\CZ,U}
@>{\sim}>> \oblv^{(\hg,\fL^+(G))_\kappa}_{\fL^+(G)}\circ \Gamma_{G,\kappa,\CZ,U} \\
@V{\sim}VV @VV{\sim}V \\
(\on{ev}_{\CZ,U})_*\circ (\oblv^l_\kappa \otimes \on{Id}) @>{\on{id}}>> (\on{ev}_{\CZ,U})_*\circ (\oblv^l_\kappa \otimes \on{Id}) 
\endCD
$$
commutes. 

\ssec{The functor \texorpdfstring{$\Loc_{G,\kappa}$}{Loc} and the infinitesimal Hecke groupoid}

\sssec{}

Let $U\subset \Bun_G$ be a quasi-compact open substack, and let $\CZ$ be a space mapping to $\Ran$.

\medskip

By the same logic as in \secref{sss:from KM to inf Hecke U}, the functor
$$(\on{ev}_{\CZ,U})^*:(\on{pt}/\fL^+(G))_\Ran\to \QCoh(U)\otimes \Dmod(\CZ)$$
lifts to a functor
\begin{equation} \label{e:almost Loc via inf}
\KL(G)_{\kappa,\CZ} \simeq \left(\Rep(\fL^+(G))_\CZ\right)^{\on{Hecke}^{\on{loc},\wedge}_{G,\CZ},\CL^{\on{loc}}_\kappa} 
\to \left(\QCoh(U)\otimes \Dmod(\CZ)\right)^{\on{Hecke}^{\on{glob},\wedge}_{G,\CZ,U},\CL^{\on{glob}}_\kappa},
\end{equation}
to be denoted 
$$(\on{ev}_{\CZ,U})^{*,\on{Hecke}^\wedge\on{-enh}}.$$

Furthermore, the functors $(\on{ev}_{\CZ,U})^{*,\on{Hecke}^\wedge\on{-enh}}$ and $(\on{ev}_{\CZ,U})^{\on{Hecke}^\wedge\on{-enh}}_*$
are adjoint. 

\sssec{}

Denote by $\sM_\kappa^{\on{loc}}$ the (factorization) monad
$$\oblv^{(\hg,\fL^+(G))_\kappa}_{\fL^+(G)}\circ \ind^{(\hg,\fL^+(G))_\kappa}_{\fL^+(G)}$$
acting on $\Rep(\fL^+(G))$. 

\medskip

Denote by $\sM_{\kappa,\CZ,U}^{\on{glob}}$ the monad acting on $\QCoh(U)\otimes \Dmod(\CZ)$,
corresponding to the forgetful functor
$$\left(\QCoh(U)\otimes \Dmod(\CZ)\right)^{\on{Hecke}^{\on{glob},\wedge}_{G,\CZ,U},\CL^{\on{glob}}_\kappa} \to
\QCoh(U)\otimes \Dmod(\CZ)$$
and its left adjoint. 

\medskip

By adjunction, we obtain that the functor $\on{ev}^*_{\CZ,U}$ intertwines the monads 
$\sM_\kappa^{\on{loc}}$ and $\sM_{\kappa,\CZ,U}^{\on{glob}}$, i.e., we have a commutative diagram
\begin{equation} \label{e:loc vs global monad U}
\CD
\Rep(\fL^+(G))_\CZ @>{\on{ev}^*_{\CZ,U}}>> \QCoh(U)\otimes \Dmod(\CZ) \\
@A{\sM_\kappa^{\on{loc}}}AA  @AA{\sM_{\kappa,\CZ,U}^{\on{glob}}}A \\
\Rep(\fL^+(G))_\CZ @>{\on{ev}^*_{\CZ,U}}>> \QCoh(U)\otimes \Dmod(\CZ).
\endCD
\end{equation} 

\sssec{}

The diagrams \eqref{e:loc vs global monad U} are compatible under inclusions $U_1\subset U_2$.
Passing to the limit over $U$, we obtain a commutative diagram
\begin{equation} \label{e:loc vs global monad}
\CD
\Rep(\fL^+(G))_\CZ @>{\on{ev}^*_{\CZ}}>> \QCoh(\Bun_G)\otimes \Dmod(\CZ) \\
@A{\sM_\kappa^{\on{loc}}}AA  @AA{\sM_{\kappa,\CZ}^{\on{glob}}}A \\
\Rep(\fL^+(G))_\CZ @>{\on{ev}^*_{\CZ}}>> \QCoh(\Bun_G)\otimes \Dmod(\CZ),
\endCD
\end{equation} 
where
$\sM_{\kappa,\CZ}^{\on{glob}}:=\sM_{\kappa,\CZ,U}^{\on{glob}}$ fot $U=\Bun_G$. 

\sssec{} \label{sss:replace monad diagram}

Let $\sM_{\kappa,\CZ}^{\on{inf}}$ denote the monad
$$(\oblv^l_\kappa\circ \ind^k_\kappa)\otimes \on{Id}$$
acting on $\QCoh(\Bun_G)\otimes \Dmod(\CZ)$.

\medskip

The map of groupoids \eqref{e:Hecke to inf}
gives rise to a map of monads
$$\sM_{\kappa,\CZ}^{\on{glob}}\to \sM_{\kappa,\CZ}^{\on{inf}}.$$

In particular, we obtain a diagram
\begin{equation} \label{e:diag two monads}
\xy
(0,0)*+{\QCoh(\Bun_G)\otimes \Dmod(\CZ)}="X";
(0,-20)*+{\QCoh(\Bun_G)\otimes \Dmod(\CZ)}="Y";
(50,0)*+{\QCoh(\Bun_G)\otimes \Dmod(\CZ)}="Z";
(50,-20)*+{\QCoh(\Bun_G)\otimes \Dmod(\CZ)}="W";
{\ar@{->}^{\sM_{\kappa,\CZ}^{\on{glob}}}"Y";"X"};
{\ar@{->}_{\sM_{\kappa,\CZ}^{\on{inf}}}"W";"Z"};
{\ar@{->}^{\on{id}}"X";"Z"};
{\ar@{->}_{\on{id}}"Y";"W"};
{\ar@{=>} "X";"W"}
\endxy
\end{equation} 

Concatenating diagrams \eqref{e:loc vs global monad} and \eqref{e:diag two monads} we obtain a diagram
\begin{equation} \label{e:loc vs inf monad}
\xy
(0,0)*+{\Rep(\fL^+(G))_\CZ}="X";
(0,-20)*+{\Rep(\fL^+(G))_\CZ}="Y";
(60,0)*+{\QCoh(\Bun_G)\otimes \Dmod(\CZ)}="Z";
(60,-20)*+{\QCoh(\Bun_G)\otimes \Dmod(\CZ)}="W";
{\ar@{->}^{\sM_{\kappa}^{\on{loc}}}"Y";"X"};
{\ar@{->}_{\sM_{\kappa,\CZ}^{\on{inf}}}"W";"Z"};
{\ar@{->}^{(\Loc^{\QCoh}_G)_\CZ}"X";"Z"};
{\ar@{->}_{(\Loc^{\QCoh}_G)_\CZ}"Y";"W"};
{\ar@{=>} "X";"W"}
\endxy
\end{equation} 

It follows from \secref{sss:alt Gamma U} that diagram \eqref{e:loc vs inf monad} identifies with the outer diagram in 
\begin{equation} \label{e:loc vs inf monad via Loc}
\xy
(0,0)*+{\Rep(\fL^+(G))_\CZ}="X";
(0,-20)*+{\KL(G)_{\kappa,\CZ}}="A";
(0,-40)*+{\Rep(\fL^+(G))_\CZ}="Y";
(60,0)*+{\QCoh(\Bun_G)\otimes \Dmod(\CZ)}="Z";
(60,-20)*+{\Dmod_\kappa(\Bun_G)\otimes \Dmod(\CZ)}="B";
(60,-40)*+{\QCoh(\Bun_G)\otimes \Dmod(\CZ),}="W";
{\ar@{->}^{(\Loc^{\QCoh}_G)_\CZ}"X";"Z"};
{\ar@{->}_{(\Loc^{\QCoh}_G)_\CZ}"Y";"W"};
{\ar@{->}_{(\Loc_{G,\kappa})_\CZ}"A";"B"};
{\ar@{->}^{\oblv^{(\hg,\fL^+(G))_\kappa}_{\fL^+(G)}}"A";"X"};
{\ar@{->}^{\ind^{(\hg,\fL^+(G))_\kappa}_{\fL^+(G)}}"Y";"A"};
{\ar@{->}_{\oblv^l_\kappa}"B";"Z"};
{\ar@{->}_{\ind^l_\kappa}"W";"B"};
{\ar@{=>} "X";"B"}
\endxy
\end{equation} 
in which the lower square is (the base change along $\CZ\to \Ran$) of 
\eqref{e:localization induction diagram Ran}, and the upper square is obtained from
the lower square by passing to adjoints along the vertical arrows.

\ssec{Proof of \thmref{t:localization and oblv orig}} \label{ss:proof localization and oblv}

We are now ready to prove \thmref{t:localization and oblv orig}. 

\sssec{}

We wish to show that the natural
transformation in
$$
\xy 
(0,0)*+{\Rep(\fL^+(G))_{\Ran^{\subseteq}}}="X";
(0,-20)*+{\KL(G)_{\kappa,\Ran^{\subseteq}}}="A";
(80,0)*+{\QCoh(\Bun_G)\otimes \Dmod(\Ran^{\subseteq})}="Z";
(80,-20)*+{\Dmod_\kappa(\Bun_G)\otimes \Dmod(\Ran^{\subseteq})}="B";
{\ar@{->}^{(\Loc^{\QCoh}_G)_{\Ran^{\subseteq}}}"X";"Z"};
{\ar@{->}_{(\Loc_{G,\kappa})_{\Ran^{\subseteq}}}"A";"B"};
{\ar@{->}^{\oblv^{(\hg,\fL^+(G))_\kappa}_{\fL^+(G)}}"A";"X"};
{\ar@{->}_{\oblv^l_\kappa\otimes \on{Id}}"B";"Z"};
{\ar@{=>} "X";"B"}
\endxy$$
becomes an isomorphism, after we: 

\begin{itemize}

\item We precompose with
$\on{ins.vac}_{\Ran}:\KL(G)_{\kappa,\Ran}\to \KL(G)_{\kappa,\Ran^{\subseteq}}$; 

\item Postcompose with $\on{Id}\otimes \on{C}^\cdot_c(\Ran^{\subseteq},-)$. 

\end{itemize}

Since the essential image of 
$$\ind^{(\hg,\fL^+(G))_\kappa}_{\fL^+(G)}:\Rep(\fL^+(G))_\Ran\to \KL(G)_{\kappa,\Ran}$$
generates the target, it is sufficient to show that the natural transformation becomes an isomorphism
when we further precompose with this functor. 

\medskip

Thus, we obtain the diagram
\begin{equation} \label{e:loc vs inf monad via Loc 1}
\xy 
(0,0)*+{\Rep(\fL^+(G))_{\Ran^\subseteq}}="X";
(0,-20)*+{\KL(G)_{\kappa,\Ran^{\subseteq}}}="A";
(80,0)*+{\QCoh(\Bun_G)\otimes \Dmod(\Ran^{\subseteq})}="Z";
(80,-20)*+{\Dmod_\kappa(\Bun_G)\otimes \Dmod(\Ran^{\subseteq})}="B";
(0,-40)*+{\KL(G)_{\kappa,\Ran}}="C";
(80,-40)*+{\Dmod_\kappa(\Bun_G)\otimes \Dmod(\Ran)}="D";
(0,-60)*+{\Rep(\fL^+(G))_\Ran}="E";
(80,-60)*+{\QCoh(\Bun_G)\otimes \Dmod(\Ran),}="F";
{\ar@{->}^{\on{ins.vac}_{\Ran}}"C";"A"};
{\ar@{->}_{\on{Id}\otimes (\on{pr}_{\on{small}})^!}"D";"B"};
{\ar@{->}_{(\Loc_{G,\kappa})_\Ran}"C";"D"};
{\ar@{->}^{(\Loc^{\QCoh}_G)_{\Ran^\subseteq}}"X";"Z"};
{\ar@{->}_{(\Loc_{G,\kappa})_{\Ran^\subseteq}}"A";"B"};
{\ar@{->}^{\oblv^{(\hg,\fL^+(G))_\kappa}_{\fL^+(G)}}"A";"X"};
{\ar@{->}_{\oblv^l_\kappa \otimes \on{Id}}"B";"Z"};
{\ar@{->}_{(\Loc^{\QCoh}_G)_\Ran}"E";"F"};
{\ar@{->}^{\ind^{(\hg,\fL^+(G))_\kappa}_{\fL^+(G)}}"E";"C"};
{\ar@{->}_{\ind^l_\kappa\otimes \on{Id}}"F";"D"};
{\ar@{=>} "X";"B"}
\endxy
\end{equation}
which we then compose with
$$\on{Id}\otimes \on{C}^\cdot_c(\Ran^{\subseteq},-):
\QCoh(\Bun_G)\otimes \Dmod(\Ran^{\subseteq})\to \QCoh(\Bun_G),$$
and we need to show that the resulting natural transformation commutes. 

\sssec{}

Consider the commutative diagram obtained by concatenating the lower two squares in \eqref{e:loc vs inf monad via Loc 1}.
It is easy to see that it identifies with the outer diagram in 
$$
\CD
\KL(G)_{\kappa,\Ran^{\subseteq}} @>{(\Loc_{G,\kappa})_{\Ran^\subseteq}}>> 
\Dmod_\kappa(\Bun_G)\otimes \Dmod(\Ran^{\subseteq}) \\
@A{\ind^{(\hg,\fL^+(G))_\kappa}_{\fL^+(G)}}AA @AA{\ind^l_\kappa \otimes \on{Id}}A \\
\Rep(\fL^+(G))_{\Ran^{\subseteq}} @>{(\Loc^{\QCoh}_G)_{\Ran^\subseteq}}>> 
\QCoh(\Bun_G)\otimes \Dmod(\Ran^{\subseteq}) \\
@A{\on{ins.vac}_{\Ran}}AA @AA{\on{Id}\otimes (\on{pr}_{\on{small}})^!}A \\
\Rep(\fL^+(G)) _\Ran @>{(\Loc^{\QCoh}_G)_\Ran}>> 
\QCoh(\Bun_G)\otimes \Dmod(\Ran). 
\endCD
$$

Hence, instead of \eqref{e:loc vs inf monad via Loc 1}, we can consider the diagram
\begin{equation} \label{e:loc vs inf monad via Loc 2}
\xy 
(0,0)*+{\Rep(\fL^+(G))_{\Ran^{\subseteq}}}="X";
(0,-20)*+{\KL(G)_{\kappa,\Ran^{\subseteq}}}="A";
(80,0)*+{\QCoh(\Bun_G)\otimes \Dmod(\Ran^{\subseteq})}="Z";
(80,-20)*+{\Dmod_\kappa(\Bun_G)\otimes \Dmod(\Ran^{\subseteq})}="B";
(0,-40)*+{\Rep(\fL^+(G))_{\Ran^{\subseteq}}}="C";
(80,-40)*+{\QCoh(\Bun_G)\otimes \Dmod(\Ran^{\subseteq})}="D";
(0,-60)*+{\Rep(\fL^+(G))_\Ran}="E";
(80,-60)*+{\QCoh(\Bun_G)\otimes \Dmod(\Ran).}="F";
{\ar@{->}^{\ind^{(\hg,\fL^+(G))_\kappa}_{\fL^+(G)}}"C";"A"};
{\ar@{->}_{\ind^l_\kappa \otimes \on{Id}}"D";"B"};
{\ar@{->}_{(\Loc^{\QCoh}_G)_{\Ran^\subseteq}}"C";"D"};
{\ar@{->}^{(\Loc^{\QCoh}_G)_{\Ran^\subseteq}}"X";"Z"};
{\ar@{->}_{(\Loc_{G,\kappa})_{\Ran^\subseteq}}"A";"B"};
{\ar@{->}^{\oblv^{(\hg,\fL^+(G))_\kappa}_{\fL^+(G)}}"A";"X"};
{\ar@{->}_{\oblv^l_\kappa \otimes \on{Id}}"B";"Z"};
{\ar@{->}_{(\Loc^{\QCoh}_G)_\Ran}"E";"F"};
{\ar@{->}^{\on{ins.vac}_{\Ran}}"E";"C"};
{\ar@{->}_{\on{Id}\otimes (\on{pr}_{\on{small}})^!}"F";"D"};
{\ar@{=>} "X";"B"}
\endxy
\end{equation}

\sssec{}

Using \secref{sss:replace monad diagram}, we can replace \eqref{e:loc vs inf monad via Loc 2} by the diagram
\begin{equation} \label{e:loc vs inf monad via Loc 3}
\xy 
(0,0)*+{\Rep(\fL^+(G))_{\Ran^{\subseteq}}}="X";
(65,0)*+{\QCoh(\Bun_G)\otimes \Dmod(\Ran^{\subseteq})}="Z";
(0,-20)*+{\Rep(\fL^+(G))_{\Ran^{\subseteq}}}="Y";
(65,-20)*+{\QCoh(\Bun_G)\otimes \Dmod(\Ran^{\subseteq})}="W";
(0,-40)*+{\Rep(\fL^+(G))_\Ran}="E";
(65,-40)*+{\QCoh(\Bun_G)\otimes \Dmod(\Ran).}="F";
{\ar@{->}^{\sM_{\kappa}^{\on{loc}}}"Y";"X"};
{\ar@{->}_{\sM_{\kappa,\Ran^{\subseteq}}^{\on{inf}}}"W";"Z"};
{\ar@{->}^-{(\Loc^{\QCoh}_G)_{\Ran^\subseteq}}"X";"Z"};
{\ar@{->}_-{(\Loc^{\QCoh}_G)_{\Ran^\subseteq}}"Y";"W"};
{\ar@{=>} "X";"W"};
{\ar@{->}^{\on{ins.vac}_{\Ran}}"E";"Y"};
{\ar@{->}_{\on{Id}\otimes (\on{pr}_{\on{small}})^!}"F";"W"};
{\ar@{->}_-{(\Loc^{\QCoh}_G)_\Ran}"E";"F"};
\endxy
\end{equation}

We further rewrite \eqref{e:loc vs inf monad via Loc 3} as the outer diagram in 

\medskip

{\tiny
\begin{equation} \label{e:loc vs inf monad via Loc 4}
\xy 
(0,0)*+{\Rep(\fL^+(G))_{\Ran^\subseteq}}="X";
(60,0)*+{\QCoh(\Bun_G)\otimes \Dmod(\Ran^{\subseteq})}="Z";
(0,-20)*+{\Rep(\fL^+(G))_{\Ran^\subseteq}}="Y";
(60,-20)*+{\QCoh(\Bun_G)\otimes \Dmod(\Ran^{\subseteq})}="W";
(0,-40)*+{\Rep(\fL^+(G))_\Ran}="E";
(60,-40)*+{\QCoh(\Bun_G)\otimes \Dmod(\Ran)}="F";
(110,0)*+{\QCoh(\Bun_G)\otimes \Dmod(\Ran^{\subseteq})}="Z'";
(110,-20)*+{\QCoh(\Bun_G)\otimes \Dmod(\Ran^{\subseteq})}="W'";
(110,-40)*+{\QCoh(\Bun_G)\otimes \Dmod(\Ran).}="F'";
{\ar@{->}^{\sM_{\kappa}^{\on{loc}}}"Y";"X"};
{\ar@{->}_{\sM_{\kappa,\Ran^{\subseteq}}^{\on{glob}}}"W";"Z"};
{\ar@{->}^-{(\Loc^{\QCoh}_G)_{\Ran^\subseteq}}"X";"Z"};
{\ar@{->}_-{(\Loc^{\QCoh}_G)_{\Ran^\subseteq}}"Y";"W"};
{\ar@{->}^{\on{ins.vac}_{\Ran}}"E";"Y"};
{\ar@{->}_{\on{Id}\otimes (\on{pr}_{\on{small}})^!}"F";"W"};
{\ar@{->}_-{(\Loc^{\QCoh}_G)_\Ran}"E";"F"};
{\ar@{->}_{\sM_{\kappa,\Ran^{\subseteq}}^{\on{inf}}}"W'";"Z'"};
{\ar@{->}_{\on{Id}\otimes (\on{pr}_{\on{small}})^!}"F'";"W'"};
{\ar@{->}^{\on{Id}}"Z";"Z'"};
{\ar@{->}^{\on{Id}}"W";"W'"};
{\ar@{->}^{\on{Id}}"F";"F'"};
{\ar@{=>} "Z";"W'"};
\endxy
\end{equation}}

\sssec{}

The left portion of \eqref{e:loc vs inf monad via Loc 4} commutes. Hence, it is enough to show that the outer diagram in 
{\small
\begin{equation} \label{e:loc vs inf monad via Loc 5}
\xy 
(40,20)*+{\QCoh(\Bun_G)}="P";
(110,20)*+{\QCoh(\Bun_G)}="P'";
(40,0)*+{\QCoh(\Bun_G)\otimes \Dmod(\Ran^{\subseteq})}="Z";
(40,-20)*+{\QCoh(\Bun_G)\otimes \Dmod(\Ran^{\subseteq})}="W";
(40,-40)*+{\QCoh(\Bun_G)\otimes \Dmod(\Ran)}="F";
(110,0)*+{\QCoh(\Bun_G)\otimes \Dmod(\Ran^{\subseteq})}="Z'";
(110,-20)*+{\QCoh(\Bun_G)\otimes \Dmod(\Ran^{\subseteq})}="W'";
(110,-40)*+{\QCoh(\Bun_G)\otimes \Dmod(\Ran).}="F'";
{\ar@{->}^{\sM_{\kappa,\Ran^{\subseteq}}^{\on{glob}}}"W";"Z"};
{\ar@{->}^{\on{Id}\otimes (\on{pr}_{\on{small}})^!}"F";"W"};
{\ar@{->}_{\sM_{\kappa,\Ran^{\subseteq}}^{\on{inf}}}"W'";"Z'"};
{\ar@{->}_{\on{Id}\otimes (\on{pr}_{\on{small}})^!}"F'";"W'"};
{\ar@{->}^{\on{Id}}"Z";"Z'"};
{\ar@{->}^{\on{Id}}"W";"W'"};
{\ar@{->}^{\on{Id}}"F";"F'"};
{\ar@{->}^{\on{Id}}"P";"P'"};
{\ar@{=>} "Z";"W'"};
{\ar@{->}^{\on{Id}\otimes \on{C}^\cdot_c(\Ran^{\subseteq},-)}"Z";"P"};
{\ar@{->}_{\on{Id}\otimes \on{C}^\cdot_c(\Ran^{\subseteq},-)}"Z'";"P'"};
\endxy
\end{equation}}
commutes. 

\medskip

We will show that already the diagram
{\small
\begin{equation} \label{e:loc vs inf monad via Loc 6}
\xy 
(40,20)*+{\QCoh(\Bun_G)\otimes \Dmod(\Ran)}="P";
(110,20)*+{\QCoh(\Bun_G)\otimes \Dmod(\Ran)}="P'";
(40,0)*+{\QCoh(\Bun_G)\otimes \Dmod(\Ran^{\subseteq})}="Z";
(40,-20)*+{\QCoh(\Bun_G)\otimes \Dmod(\Ran^{\subseteq})}="W";
(40,-40)*+{\QCoh(\Bun_G)\otimes \Dmod(\Ran)}="F";
(110,0)*+{\QCoh(\Bun_G)\otimes \Dmod(\Ran^{\subseteq})}="Z'";
(110,-20)*+{\QCoh(\Bun_G)\otimes \Dmod(\Ran^{\subseteq})}="W'";
(110,-40)*+{\QCoh(\Bun_G)\otimes \Dmod(\Ran).}="F'";
{\ar@{->}^{\sM_{\kappa,\Ran^{\subseteq}}^{\on{glob}}}"W";"Z"};
{\ar@{->}^{\on{Id}\otimes (\on{pr}_{\on{small}})^!}"F";"W"};
{\ar@{->}_{\sM_{\kappa,\Ran^{\subseteq}}^{\on{inf}}}"W'";"Z'"};
{\ar@{->}_{\on{Id}\otimes (\on{pr}_{\on{small}})^!}"F'";"W'"};
{\ar@{->}^{\on{Id}}"Z";"Z'"};
{\ar@{->}^{\on{Id}}"W";"W'"};
{\ar@{->}^{\on{Id}}"F";"F'"};
{\ar@{->}^{\on{Id}}"P";"P'"};
{\ar@{=>} "Z";"W'"};
{\ar@{->}^{\on{Id}\otimes (\on{pr}_{\on{small}})_!}"Z";"P"};
{\ar@{->}_{\on{Id}\otimes (\on{pr}_{\on{small}})_!}"Z'";"P'"};
\endxy
\end{equation}}
commutes. 

\medskip

This is a particular case of the next assertion: 

\begin{thm} \label{t:Nick monad}
For any $\CZ\to \Ran$, the natural transformation in the diagram
$$
\xy 
(40,20)*+{\QCoh(\Bun_G)\otimes \Dmod(\CZ)}="P";
(110,20)*+{\QCoh(\Bun_G)\otimes \Dmod(\CZ)}="P'";
(40,0)*+{\QCoh(\Bun_G)\otimes \Dmod(\CZ^{\subseteq})}="Z";
(40,-20)*+{\QCoh(\Bun_G)\otimes \Dmod(\CZ^{\subseteq})}="W";
(40,-40)*+{\QCoh(\Bun_G)\otimes \Dmod(\CZ)}="F";
(110,0)*+{\QCoh(\Bun_G)\otimes \Dmod(\CZ^{\subseteq})}="Z'";
(110,-20)*+{\QCoh(\Bun_G)\otimes \Dmod(\CZ^{\subseteq})}="W'";
(110,-40)*+{\QCoh(\Bun_G)\otimes \Dmod(\CZ).}="F'";
{\ar@{->}^{\sM_{\kappa,\CZ^{\subseteq}}^{\on{glob}}}"W";"Z"};
{\ar@{->}^{\on{Id}\otimes (\on{pr}_{\on{small},\CZ})^!}"F";"W"};
{\ar@{->}_{\sM_{\kappa,\CZ^{\subseteq}}^{\on{inf}}}"W'";"Z'"};
{\ar@{->}_{\on{Id}\otimes (\on{pr}_{\on{small},\CZ})^!}"F'";"W'"};
{\ar@{->}^{\on{Id}}"Z";"Z'"};
{\ar@{->}^{\on{Id}}"W";"W'"};
{\ar@{->}^{\on{Id}}"F";"F'"};
{\ar@{->}^{\on{Id}}"P";"P'"};
{\ar@{=>} "Z";"W'"};
{\ar@{->}^{\on{Id}\otimes (\on{pr}_{\on{small}})_!}"Z";"P"};
{\ar@{->}_{\on{Id}\otimes (\on{pr}_{\on{small}})_!}"Z'";"P'"};
\endxy
$$
induced by the map of monads 
$$\sM_{\kappa,\CZ^{\subseteq}}^{\on{glob}}\to \sM_{\kappa,\CZ^{\subseteq}}^{\on{inf}},$$
is an isomorphism.
\end{thm}

This theorem is a particular case of \cite[Theorem 4.3.6]{Ro2}, combined with Remark 4.5.6 in 
{\it loc. cit.}

\qed[\thmref{t:localization and oblv orig}]

\ssec{Localization is (almost) a localization}

\sssec{}

The goal of the next few subsections is prove the following assertion:

\begin{thm} \label{t:Loc is loc}
Let $U\overset{j}\hookrightarrow \Bun_G$ be a quasi-compact open substack. Then the functor
$$\Loc_{G,\kappa,U}:=j^*\circ \Loc_{G,\kappa}, \quad \KL(G)_{\kappa,\Ran}\to \Dmod_\kappa(U)$$
is a \emph{localization}.\footnote{I.e., its right adjoint is fully faithful.}
\end{thm}

\sssec{}

We first consider a version of \thmref{t:Loc is loc} when instead of $\Dmod_\kappa(\Bun_G)$,
we take $\QCoh(\Bun_G)$.

\medskip 

Consider the local-to-global functor
$$\ul\Loc_G^\QCoh:\ul\Rep(\fL^+(G))\to \QCoh(\Bun_G)\otimes \ul\Dmod(\Ran).$$

\medskip

For an open substack $U\subset \Bun_G$, denote
$$\Loc_{G,U}^\QCoh:=j^*\circ \Loc_G^\QCoh, \quad \Rep(\fL^+(G))_\Ran\to \QCoh(U).$$

We will prove:\footnote{A more general assertion, of which \thmref{t:Loc is loc QCoh} is a particular case, appears in \cite[Proposition C.1.7]{GLC4}.}

\begin{thm} \label{t:Loc is loc QCoh}
If $U$ is quasi-compact, then the functor $\Loc_{G,U}^\QCoh$ is a localization.
\end{thm} 

\sssec{}

Consider the following general paradigm: let 
$$F:\bC_1\to \bC_2$$
be a continuous functor between compactly generated categories. Assume that $F$ preserves compactness.

\medskip

Let
$$F^{\on{fake-op}}:\bC_1^\vee\to \bC_2^\vee$$
be the conjugate functor, i.e., 
$$F^{\on{fake-op}}=(F^R)^\vee.$$

\medskip

In other words, $F^{\on{fake-op}}$ is the ind-extension of
$$(\bC_1^\vee)^c\simeq (\bC_1^c)^{\on{op}} \overset{F^{\on{op}}}\longrightarrow  (\bC_2^c)^{\on{op}} \simeq (\bC_2^\vee)^c.$$

\medskip

Let $\on{u}_{\bC_i}\in \bC_i\otimes \bC^\vee_i$ be the unit of the duality. Note that we have a canonically defined map
\begin{equation} \label{e:unit to unit cat}
(F\otimes F^{\on{fake-op}})(\on{u}_{\bC_1})\to \on{u}_{\bC_2}.
\end{equation}

Namely, 
\begin{multline*}
(F\otimes F^{\on{fake-op}})(\on{u}_{\bC_1})=
(F\otimes \on{Id})\circ (\on{Id}\otimes F^{\on{fake-op}})(\on{u}_{\bC_1})=
(F\otimes \on{Id})\circ (\on{Id}\otimes (F^R)^\vee)(\on{u}_{\bC_1})\simeq  \\
\simeq (F\otimes \on{Id})\circ (F^R\otimes \on{Id})(\on{u}_{\bC_1})=((F\circ F^R)\otimes \on{Id})(\on{u}_{\bC_1})\to
\on{u}_{\bC_1}.
\end{multline*} 

\sssec{}

The following is elementary:

\begin{lem} \label{l:loc crit}
The functor $F$ is a localization if and only if \eqref{e:unit to unit cat} is an isomorphism.
\end{lem}

\sssec{}

We will prove \thmref{t:Loc is loc QCoh} by applying 
\lemref{l:loc crit} to $\bC_1:=\Rep(\fL^+(G))_\Ran$, $\bC_2=\QCoh(U)$ and $F:=\Loc^\QCoh_{G,U}$.

\medskip

First, it is easy to see that $\Loc^\QCoh_{G,U}$ preserves compactness. 

\sssec{}

We have a canonical self-duality 
\begin{equation} \label{e:Rep GO self-dual}
\Rep(\fL^+(G))^\vee\simeq \Rep(\fL^+(G))
\end{equation}
as a factorization category. Its unit object $u_{\Rep(\fL^+(G))}$
is the regular representation
$$R_{\fL^+(G)}\in \Rep(\fL^+(G))\otimes \Rep(\fL^+(G)),$$
which has a natural structure of (commutative) factorization algebra. 

\medskip

From \eqref{e:Rep GO self-dual} we obtain a self-duality
\begin{equation} \label{e:Rep GO self-dual Ran}
(\Rep(\fL^+(G))_\Ran)^\vee\simeq \Rep(\fL^+(G))_\Ran.
\end{equation}

The unit $\on{u}_{\Rep(\fL^+(G))_\Ran}$ of \eqref{e:Rep GO self-dual Ran} is given by the image of 
$$R_{\fL^+(G),\Ran}\in (\Rep(\fL^+(G))\otimes \Rep(\fL^+(G)))_\Ran$$
under 
\begin{multline*}
(\Rep(\fL^+(G))\otimes \Rep(\fL^+(G)))_\Ran\simeq \\
\simeq (\Rep(\fL^+(G))_\Ran\otimes \Rep(\fL^+(G))_\Ran)\underset{\Dmod(\Ran\times \Ran)}\otimes \Dmod(\Ran)
\overset{\on{Id}\otimes (\Delta_\Ran)_!}\longrightarrow  \\
\to (\Rep(\fL^+(G))_\Ran\otimes \Rep(\fL^+(G))_\Ran)\underset{\Dmod(\Ran\times \Ran)}\otimes (\Dmod(\Ran) \otimes \Dmod(\Ran))= \\
=\Rep(\fL^+(G))_\Ran\otimes \Rep(\fL^+(G))_\Ran.
\end{multline*}

\sssec{}

With respect to the canonical self-duality
$$\QCoh(U)^\vee\simeq \QCoh(U),$$
we have
$$(\Loc^\QCoh_{G,U})^{\on{fake-op}}\simeq \Loc^\QCoh_{G,U}.$$

\medskip

From here we obtain
\begin{multline*}
(\Loc^\QCoh_{G,U}\otimes (\Loc^\QCoh_{G,U}))^{^{\on{fake-op}}}(\on{u}_{\Rep(\fL^+(G))_\Ran})\simeq
(\Loc_{G,U}^\QCoh\otimes \Loc_{G,U}^\QCoh)\circ (\Delta_\Ran)_!(R_{\fL^+(G),\Ran})\simeq \\
\simeq \Loc^\QCoh_{G\times G,U\times U}(R_{\fL^+(G),\Ran}),
\end{multline*}
see \secref{sss:almost unital} for the latter isomorphism. 

\medskip

Thus, the map \eqref{e:unit to unit cat} is a map
\begin{equation} \label{e:Loc diag}
\Loc^\QCoh_{G\times G,U\times U}(R_{\fL^+(G),\Ran})\to (\Delta_{\Bun_G})_*(\CO_{\Bun_G}).
\end{equation}

\sssec{}

Hence, by \lemref{l:loc crit}, in order to prove \thmref{t:Loc is loc QCoh}, we have to show that the map
\eqref{e:Loc diag} is an isomorphism.

\medskip

The latter is obtained by repeating verbatim the proof of \cite[Theorem 12.6.3]{AGKRRV}, see also Remark 12.6.6 in {\it loc. cit.} 

\qed[\thmref{t:Loc is loc QCoh}]

\sssec{}

We are now ready to prove \thmref{t:Loc is loc}. We will give two proofs: one in 
\secref{ss:Loc is loc 1} and another in \secref{ss:Loc is loc 2}.

\ssec{First proof of \thmref{t:Loc is loc}}  \label{ss:Loc is loc 1}

\sssec{}

We need to show that the functor
$$\Gamma_{G,\kappa,U}: \Dmod_\kappa(U)\to \KL(G)_{\kappa,\Ran}$$
is fully faithful.

\medskip

For $\CZ\to \Ran$, we interpret the category
$\KL(G)_{\kappa,\CZ}$ as
$$\left(\Rep(\fL^+(G))_\CZ\right)^{\on{Hecke}^{\on{loc},\wedge}_{G,\CZ},\CL^{\on{loc}}_\kappa}.$$

\medskip

We interpret the functor $\Gamma_{G,\kappa,U}$ as a composition
\begin{multline} \label{e:Gamma U}
\Dmod_\kappa(U) \overset{\on{Id}\otimes \omega_\Ran}\longrightarrow \Dmod_\kappa(U)\otimes \Dmod(\CZ)
\overset{\oblv_{\on{inf}\to \on{Hecke}^\wedge,\Ran,U}}\longrightarrow \\
\to \left(\QCoh(U)\otimes \Dmod(\Ran)\right)^{\on{Hecke}^{\on{glob},\wedge}_{G,\Ran,U},\CL^{\on{glob}}_\kappa}
\overset{(\on{ev}_{\Ran,U})^{\on{Hecke}^\wedge\on{-enh}}_*}\longrightarrow \\
\to \left(\Rep(\fL^+(G))_\Ran\right)^{\on{Hecke}^{\on{loc},\wedge}_{G,\Ran},\CL^{\on{loc}}_\kappa},
\end{multline}
where:

\begin{itemize}

\item $\oblv_{\on{inf}\to \on{Hecke}^\wedge,\Ran,U}$ is the functor \eqref{e:Dmod to glob inf Hecke equiv U};

\item $(\on{ev}_{\Ran,U})^{\on{Hecke}^\wedge\on{-enh}}_*$  is the functor \eqref{e:global inf Hecke equiv to local inf Hecke equiv U}.

\end{itemize} 

\medskip

We will show that \eqref{e:Gamma U} is a \emph{retract} of a fully faithful functor\footnote{We say that a functor $F:\bC\to \bD$
is a retract of a functor $F_1:\bC\to \bD_1$ if there exists a retraction $\bD\overset{\phi}\to \bD_1\overset{\psi}\to \bD$ so that
$F_1\simeq \phi\circ F$ and $F\simeq F_1\circ \psi$.} and is therefore itself fully faithful. 

\begin{rem}

Here is the reason the retraction appears:

\medskip

There are two copies of the Ran space that play a role: one is in \thmref{t:Loc is loc QCoh}, and the other is 
\thmref{t:Nick 1}. To achieve fully faithfulness of \eqref{e:Gamma U}, we need both of these copies, 
and that is how that $\Ran$ space becomes ``doubled"  in the guise of $\Ran^\subseteq$. 

\medskip

We then pass from $\Ran^\subseteq$ to just $\Ran$, but at the expense of replacing the original
fully faithful functor \eqref{e:Gamma overkill} by its retract. 

\end{rem}

\sssec{}

Note that the pullbacks of $\on{Hecke}^{\on{loc},\wedge}_{G,\Ran}$ along $\on{pr}_{\on{small}}$ and $\on{pr}_{\on{big}}$
give rise to a well-defined
groupoids on $(\on{pt}/\fL^+(G))_{\Ran^{\subseteq}}$; we will denote them by 
$$\on{Hecke}^{\on{loc},\wedge}_{G,\Ran^\subseteq,\on{small}} \text{ and }
\on{Hecke}^{\on{loc},\wedge}_{G,\Ran^\subseteq,\on{big}},$$
respectively. 

\medskip

The unital structure on $\on{Hecke}^{\on{loc},\wedge}_{G,\Ran}$ gives rise to a map
\begin{equation} \label{e:restr big small groupoids}
\on{Hecke}^{\on{loc},\wedge}_{G,\Ran^\subseteq,\on{small}} \to \on{Hecke}^{\on{loc},\wedge}_{G,\Ran^\subseteq,\on{big}}.
\end{equation}

\medskip

Similar definitions apply to $\on{Hecke}^{\on{glob},\wedge}_{G,\Ran,U}$.

\sssec{}

Consider the functor
\begin{multline} \label{e:Gamma overkill}
\Dmod_\kappa(U) \overset{\on{Id}\otimes \omega_\Ran}\longrightarrow \Dmod_\kappa(U)\otimes \Dmod(\CZ)
\overset{\oblv_{\on{inf}\to \on{Hecke}^\wedge,\Ran,U}}\longrightarrow \\
\to \left(\QCoh(U)\otimes \Dmod(\Ran)\right)^{\on{Hecke}^{\on{glob},\wedge}_{G,\Ran,U},\CL^{\on{glob}}_\kappa}\to \\
\overset{(\on{pr}_{\on{small}})^!}\longrightarrow 
\left(\QCoh(U)\otimes \Dmod(\Ran^\subseteq)\right)^{\on{Hecke}^{\on{glob},\wedge}_{G,\Ran^\subseteq,\on{small},U},\CL^{\on{glob}}_\kappa}\to \\
\overset{(\on{ev}_{\Ran^\subseteq,U})^{\on{Hecke}_{\on{small}}^\wedge\on{-enh}}_*}\longrightarrow 
\left(\Rep(\fL^+(G))_{\Ran^\subseteq}\right)^{\on{Hecke}^{\on{loc},\wedge}_{G,\Ran^\subseteq,\on{small}},\CL^{\on{loc}}_\kappa},
\end{multline} 
where $(\on{ev}_{\Ran^\subseteq,U})^{\on{Hecke}_{\on{small}}^\wedge\on{-enh}}_*$ is the enhancement for the corresponding
groupoids of the functor $\on{ev}_{\CZ^\subseteq,U}$, cf. formula \eqref{e:Dmod to glob inf Hecke equiv U}. 

\sssec{}

We claim that the functor \eqref{e:Gamma overkill} is fully faithful. In fact, we claim that it is a composition of two fully
faithful functors. 

\medskip

Namely, the functor
\begin{multline*} 
\Dmod_\kappa(U) \overset{\on{Id}\otimes \omega_\Ran}\longrightarrow \Dmod_\kappa(U)\otimes \Dmod(\CZ)
\overset{\oblv_{\on{inf}\to \on{Hecke}^\wedge,\Ran,U}}\longrightarrow \\
\to \left(\QCoh(U)\otimes \Dmod(\Ran)\right)^{\on{Hecke}^{\on{glob},\wedge}_{G,\Ran,U},\CL^{\on{glob}}_\kappa}
\end{multline*} 
is fully faithful thanks to \cite[Theorem 4.5.3]{Ro2} (cf. \thmref{t:Nick 1}). 

\medskip

The fact that 
\begin{multline}  \label{e:partial equiv}
\left(\QCoh(U)\otimes \Dmod(\Ran)\right)^{\on{Hecke}^{\on{glob},\wedge}_{G,\Ran,U},\CL^{\on{glob}}_\kappa}
\overset{(\on{pr}_{\on{small}})^!}\longrightarrow \\
\to \left(\QCoh(U)\otimes \Dmod(\Ran^\subseteq)\right)^{\on{Hecke}^{\on{glob},\wedge}_{G,\Ran^\subseteq,\on{small},U},\CL^{\on{glob}}_\kappa}
\overset{(\on{ev}_{\Ran^\subseteq,U})^{\on{Hecke}_{\on{small}}^\wedge\on{-enh}}_*}\longrightarrow \\
\to \left(\Rep(\fL^+(G))_{\Ran^\subseteq}\right)^{\on{Hecke}^{\on{loc},\wedge}_{G,\Ran^\subseteq,\on{small}},\CL^{\on{loc}}_\kappa},
\end{multline}
is fully faithful is a formal consequence of the fact that
\begin{equation}  \label{e:param Loc QCoh}
\QCoh(U)\otimes \Dmod(\Ran)
\overset{(\on{pr}_{\on{small}})^!}\longrightarrow 
\QCoh(U)\otimes \Dmod(\Ran^\subseteq)
\overset{\on{ev}_{\Ran^\subseteq,\on{small},U}}\longrightarrow \Rep(\fL^+(G))_{\Ran^\subseteq}
\end{equation} 
is fully faithful, which is a parameterized version of \thmref{t:Loc is loc QCoh}:

\medskip

Indeed, each side of \eqref{e:partial equiv} is obtained from the corresponding side of \eqref{e:param Loc QCoh}
as modules for the corresponding monad, and the functor \eqref{e:param Loc QCoh} intertwines the actions of 
these monads. 

\sssec{}

We now claim that the functor \eqref{e:Gamma U} is a retract of \eqref{e:Gamma overkill}. 

\medskip

Indeed, the corresponding functor
$$\left(\Rep(\fL^+(G))_\Ran\right)^{\on{Hecke}^{\on{loc},\wedge}_{G,\Ran},\CL^{\on{loc}}_\kappa} \to 
\left(\Rep(\fL^+(G))_{\Ran^\subseteq}\right)^{\on{Hecke}^{\on{loc},\wedge}_{G,\Ran^\subseteq,\on{small}},\CL^{\on{loc}}_\kappa}$$
is 
\begin{multline} \label{e:restr groupoids 1}
\left(\Rep(\fL^+(G))_\Ran\right)^{\on{Hecke}^{\on{loc},\wedge}_{G,\Ran},\CL^{\on{loc}}_\kappa} \to \\
\overset{(\on{pr}_{\on{big}})^!}\longrightarrow 
\left(\Rep(\fL^+(G))_{\Ran^\subseteq}\right)^{\on{Hecke}^{\on{loc},\wedge}_{G,\Ran^\subseteq,\on{big}},\CL^{\on{loc}}_\kappa}\to \\
\to \left(\Rep(\fL^+(G))_{\Ran^\subseteq}\right)^{\on{Hecke}^{\on{loc},\wedge}_{G,\Ran^\subseteq,\on{small}},\CL^{\on{loc}}_\kappa},
\end{multline}
where the second arrow is given by restriction along \eqref{e:restr big small groupoids}. 

\medskip

The functor
$$\left(\Rep(\fL^+(G))_{\Ran^\subseteq}\right)^{\on{Hecke}^{\on{loc},\wedge}_{G,\Ran^\subseteq,\on{small}},\CL^{\on{loc}}_\kappa}\to
\left(\Rep(\fL^+(G))_\Ran\right)^{\on{Hecke}^{\on{loc},\wedge}_{G,\Ran},\CL^{\on{loc}}_\kappa}$$
is given by
\begin{equation} \label{e:restr groupoids 2}
\left(\Rep(\fL^+(G))_{\Ran^\subseteq}\right)^{\on{Hecke}^{\on{loc},\wedge}_{G,\Ran^\subseteq,\on{small}},\CL^{\on{loc}}_\kappa}
\overset{\on{diag}^!}\longrightarrow 
\left(\Rep(\fL^+(G))_\Ran\right)^{\on{Hecke}^{\on{loc},\wedge}_{G,\Ran},\CL^{\on{loc}}_\kappa}.
\end{equation}

It is a straightforward verification that
$$\text{\eqref{e:restr groupoids 1}}\circ \text{\eqref{e:Gamma U}}\simeq \text{\eqref{e:Gamma overkill}} \text{ and }
\text{\eqref{e:restr groupoids 2}}\circ \text{\eqref{e:Gamma overkill}} \simeq \text{\eqref{e:Gamma U}}.$$

\qed[First proof of \thmref{t:Loc is loc}]

\ssec{Second proof of \thmref{t:Loc is loc}} \label{ss:Loc is loc 2}

In this subsection we will apply \lemref{l:loc crit} directly to the functor $\Loc_{G,\kappa,U}$ 
to deduce the assertion of \thmref{t:Loc is loc}. 

\sssec{}

Let $\kappa':=2\cdot \crit-\kappa$. Recall that the unit $\on{u}_{\KL(G)_\kappa}$ of the duality
$$(\KL(G)_\kappa)^\vee\simeq \KL(G)_{\kappa'}$$
is the (factorization algebra) object $\fCDO(G)_{\kappa,\kappa'}$, see \secref{sss:KL self-duality kappa}. 

\medskip

The dual of $\KL(G)_{\kappa,\Ran}$ identifies with $\KL(G)_{\kappa',\Ran}$ with the unit
$\on{u}_{\KL(G)_{\kappa,\Ran}}$ of the duality
being the image of 
$$(\fCDO(G)_{\kappa,\kappa'})_\Ran\in (\KL(G)_\kappa\otimes \KL(G)_{\kappa'})_\Ran$$
under
\begin{multline*}
(\KL(G)_\kappa\otimes \KL(G)_{\kappa'})_\Ran\simeq
(\KL(G)_{\kappa,\Ran}\otimes \KL(G)_{\kappa',\Ran})\underset{\Dmod(\Ran\times \Ran)}\otimes \Dmod(\Ran)
\overset{\on{Id}\otimes (\Delta_\Ran)_!}\longrightarrow \\
\to (\KL(G)_{\kappa,\Ran}\otimes \KL(G)_{\kappa',\Ran})\underset{\Dmod(\Ran\times \Ran)}\otimes (\Dmod(\Ran) \otimes \Dmod(\Ran))=\\
=\KL(G)_{\kappa,\Ran}\otimes \KL(G)_{\kappa',\Ran}.
\end{multline*}

\sssec{} \label{sss:dual DU first}

We identify the category dual to $\Dmod_\kappa(U)$ with $\Dmod_{\kappa'}(U)$ by means of the functor
$$(\Dmod_{\kappa}(U))^\vee\simeq \Dmod_{-\kappa}(U) \overset{\text{\eqref{e:tw can bundle}}}\simeq \Dmod_{\kappa'}(U),$$
where we apply \eqref{e:tw can bundle} to $U$ instead of $\Bun_G$ and $\kappa$ instead of $\kappa'$. 

\medskip

In terms of this identification, the unit $\on{u}_{\Dmod_\kappa(U)}$
of the duality is given by
$$(j_U\times j_U)^*((\Delta_{\Bun_G})_*(\omega_{\Bun_G})'),$$
where 
$$(\Delta_{\Bun_G})_*(\omega_{\Bun_G})'\in \Dmod_\kappa(\Bun_G)\otimes \Dmod_{\kappa'}(\Bun_G)$$
is the image of 
$$(\Delta_{\Bun_G})_*(\omega_{\Bun_G})\in \Dmod_\kappa(\Bun_G)\otimes \Dmod_{-\kappa}(\Bun_G)$$
under $\on{Id}\otimes \text{\eqref{e:tw can bundle}}$.

\sssec{} \label{sss:dual DU last}

According to \propref{p:Loc as dual}, with respect to the above dualities, the functor conjugate to $\Loc_{G,\kappa,U}$ identifies with
$\Loc_{G,\kappa',U}$.

\medskip

We obtain 
\begin{multline*}
(\Loc_{G,\kappa,U}\otimes (\Loc_{G,\kappa,U})^{\on{fake-op}})(\on{u}_{\KL(G)_{\kappa,\Ran}})\simeq
(\Loc_{G,\kappa,U}\otimes \Loc_{G,\kappa',U})\circ (\Delta_\Ran)_!((\fCDO(G)_{\kappa,\kappa'})_\Ran)\simeq \\
\simeq (j_U\times j_U)^*\left(\Loc_{G\times G,(\kappa,\kappa')}(\fCDO(G)_{\kappa,\kappa'})\right),
\end{multline*} 
see \secref{sss:almost unital} for the latter isomorphism. 

\medskip

Thus, the map \eqref{e:unit to unit cat} is a map 
\begin{equation} \label{e:unit to unit KL U}
(j_U\times j_U)^*\left(\Loc_{G\times G,(\kappa,\kappa')}((\fCDO(G)_{\kappa,\kappa'})_\Ran)\right)\to
(j_U\times j_U)^*((\Delta_{\Bun_G})_*(\omega_{\Bun_G})').
\end{equation}

Hence, by \lemref{l:loc crit}, in order to prove \thmref{t:Loc is loc}, we need to show that the map
\eqref{e:unit to unit KL U} is an isomorphism.

\begin{rem}

The maps \eqref{e:unit to unit KL U} are compatible under $U\subset U'$. Hence, once we establish 
the isomorphism \eqref{e:unit to unit KL U}, we will have established an isomorphism
\begin{equation} \label{e:unit to unit KL BunG}
\Loc_{G\times G,(\kappa,\kappa')}((\fCDO(G)_{\kappa,\kappa'})_\Ran)\simeq 
(\Delta_{\Bun_G})_*(\omega_{\Bun_G})'.
\end{equation}

\end{rem}

\sssec{}

We will show that the map 
\begin{equation} \label{e:Loc isom to prove}
(\on{Id}\otimes \oblv^l_{\kappa'})\circ (\Loc_{G,\kappa,U}\otimes \Loc_{G,\kappa',U})((\fCDO(G)_{\kappa,\kappa'})_\Ran)\to 
(\on{Id}\otimes \oblv^l_{\kappa'})((\Delta_{\Bun_G})_*(\omega_{\Bun_G})')
\end{equation}
in $\Dmod_\kappa(U)\otimes \QCoh(U)$, induced by the map \eqref{e:unit to unit cat}, is an isomorphism.

\medskip

This would prove that the map \eqref{e:unit to unit cat} is an isomorphism, since the fuctor 
$\on{Id}\otimes \oblv^l_{\kappa'}$ is conservative. 

\begin{rem} 

The idea of the proof that \eqref{e:Loc isom to prove} is an equivalence is that we will identify this map with the map
obtained by applying $\ind^l_\kappa\otimes \on{Id}$ to the map \eqref{e:Loc diag}. 

\medskip

This is not completely automatic, since the functor the natural transformation
$$ \Loc^\QCoh_{G,U}\circ \oblv^{(\hg,\fL^+(G))_{\kappa'}}_{\fL^+(G)}\to \oblv^l_{\kappa'}\circ \Loc_{G,\kappa',U}$$
is \emph{not} an isomorphism. 

\medskip 

However, the map in question will become an isomorphism when evaluated on $\fCDO(G)_{\kappa,\kappa'}$, since
the latter is a \emph{unital} factorization algebra.

\end{rem}

\sssec{}

Consider the functor
$$\oblv^{(\hg,\fL^+(G))_{\kappa'}}_{\fL^+(G)}:\KL(G)_{\kappa'}\to 
\Rep(\fL^+(G))$$
and the corresponding functor
$$\KL(G)_{\kappa',\Ran}\to \Rep(\fL^+(G))_\Ran,$$
which we denote by the same character 
$\oblv^{(\hg,\fL^+(G))_{\kappa'}}_{\fL^+(G)}$.

\medskip

Consider the object
\begin{multline*}
(\on{Id}\otimes \oblv^{(\hg,\fL^+(G))_{\kappa'}}_{\fL^+(G)})(\on{u}_{\KL(G)_{\kappa,\Ran}})\simeq \\
\simeq (\Delta_\Ran)_!\circ (\on{Id}\otimes \oblv^{(\hg,\fL^+(G))_{\kappa'}}_{\fL^+(G)})((\fCDO(G)_{\kappa,\kappa'})_\Ran)\in 
\KL(G)_{\kappa,\Ran}\otimes \Rep(\fL^+(G))_\Ran.
\end{multline*}

Consider the object
$$(\on{Id}\otimes \oblv^l_{\kappa'})(\on{u}_{\Dmod_\kappa(U)})\simeq 
(\on{Id}\otimes \oblv^l_{\kappa'})\circ (j_U\times j_U)^*((\Delta_{\Bun_G})_*(\omega_{\Bun_G})')\in 
\Dmod_\kappa(U)\otimes \QCoh(U).$$

We claim that there is a canonically defined map
\begin{equation} \label{e:right forget}
(\Loc_{G,\kappa,U}\otimes \Loc^\QCoh_{G,U}) \circ 
(\on{Id}\otimes \oblv^{(\hg,\fL^+(G))_{\kappa'}}_{\fL^+(G)})(\on{u}_{\KL(G)_{\kappa,\Ran}})\to 
(\on{Id}\otimes \oblv^l_{\kappa'})(\on{u}_{\Dmod_\kappa(U)}).
\end{equation} 

Namely, the map \eqref{e:right forget} is the composition

\begin{multline*}
(\Loc_{G,\kappa,U}\otimes \Loc^\QCoh_{G,U}) \circ 
(\on{Id}\otimes \oblv^{(\hg,\fL^+(G))_{\kappa'}}_{\fL^+(G)})(\on{u}_{\KL(G)_{\kappa,\Ran}}) \simeq \\
\simeq 
(\Loc_{G,\kappa,U}\otimes\on{Id}) \circ (\on{Id}\otimes \Loc^\QCoh_{G,U})\circ
(\on{Id}\otimes \oblv^{(\hg,\fL^+(G))_{\kappa'}}_{\fL^+(G)})(\on{u}_{\KL(G)_{\kappa,\Ran}})  \simeq \\
\simeq
(\on{Id}\otimes \Loc^\QCoh_{G,U})\circ
(\on{Id}\otimes \oblv^{(\hg,\fL^+(G))_{\kappa'}}_{\fL^+(G)})\circ (\Loc_{G,\kappa,U}\otimes\on{Id})(\on{u}_{\KL(G)_{\kappa,\Ran}})\simeq \\
\simeq
(\on{Id}\otimes \Loc^\QCoh_{G,U})\circ (\on{Id}\otimes \oblv^{(\hg,\fL^+(G))_{\kappa'}}_{\fL^+(G)})\circ 
(\on{Id}\otimes (\Loc_{G,\kappa,U})^\vee)(\on{u}_{\Dmod_\kappa(U)}) \simeq \\
\simeq (\on{Id}\otimes \Loc^\QCoh_{G,U})\circ (\on{Id}\otimes \oblv^{(\hg,\fL^+(G))_{\kappa'}}_{\fL^+(G)})\circ 
(\on{Id}\otimes (\Loc_{G,\kappa',U})^R)(\on{u}_{\Dmod_\kappa(U)}) = \\
=\left(\on{Id} \otimes (\Loc^\QCoh_{G,U}\circ \oblv^{(\hg,\fL^+(G))_{\kappa'}}_{\fL^+(G)}\circ (\Loc_{G,\kappa',U})^R)\right)(\on{u}_{\Dmod_\kappa(U)})
\to \\
\to \left(\on{Id} \otimes (\oblv^l_{\kappa'}\circ \Loc_{G,\kappa',U}\circ (\Loc_{G,\kappa',U})^R)\right)(\on{u}_{\Dmod_\kappa(U)})\to
(\on{Id}\otimes \oblv^l_{\kappa'})(\on{u}_{\Dmod_\kappa(U)}),
\end{multline*} 
where the next-to-last arrow is induced by the natural transformation \eqref{e:localization oblv sheaves int}. 

\sssec{}

Note that by construction, the map \eqref{e:right forget} can reinterpreted as
\begin{multline} \label{e:right forget again}
(\Loc_{G,\kappa,U}\otimes \Loc^\QCoh_{G,U}) \circ 
(\on{Id}\otimes \oblv^{(\hg,\fL^+(G))_{\kappa'}}_{\fL^+(G)})(\on{u}_{\KL(G)_{\kappa,\Ran}})\simeq \\
\simeq 
(\Loc_{G,\kappa,U}\otimes \Loc^\QCoh_{G,U}) \circ 
(\on{Id}\otimes \oblv^{(\hg,\fL^+(G))_{\kappa'}}_{\fL^+(G)})((\fCDO(G)_{\kappa,\kappa'})_\Ran)\to \\
\to (\on{Id}\otimes \oblv^l_{\kappa'})\circ (\Loc_{G,\kappa,U}\otimes \Loc_{G,\kappa',U})
((\fCDO(G)_{\kappa,\kappa'})_\Ran) \overset{\text{\eqref{e:Loc isom to prove}}}\longrightarrow \\
\to (\on{Id}\otimes \oblv^l_{\kappa'})((\Delta_{\Bun_G})_*(\omega_{\Bun_G})')\simeq 
(\on{Id}\otimes \oblv^l_{\kappa'})(\on{u}_{\Dmod_\kappa(U)}),
\end{multline} 
where the second arrow is induced by \eqref{e:localization oblv sheaves int}. 

\medskip

Note, however, that the fact that $\fCDO(G)_{\kappa,\kappa'}$ is a \emph{unital} factorization algebra,
we obtain that the second arrow in \eqref{e:right forget again} is an isomorphism, see \lemref{l:AB nat trans abs}. 

\medskip

Hence, we obtain that in order to show that \eqref{e:Loc isom to prove} is an isomorphism, it suffices
to show that \eqref{e:right forget} is an isomorphism. 

\sssec{}

Consider now the functor 
$$\ind^{(\hg,\fL^+(G))_\kappa}_{\fL^+(G)}:\Rep(\fL^+(G))\to \KL(G)_\kappa,$$
and the corresponding functor
$$\Rep(\fL^+(G))_\Ran\to \KL(G)_{\kappa,\Ran},$$
which we denote by the same character 
$\ind^{(\hg,\fL^+(G))_\kappa}_{\fL^+(G)}$.

\medskip

Consider the object
\begin{multline*}
(\ind^{(\hg,\fL^+(G))_\kappa}_{\fL^+(G)}\otimes \on{Id})(\on{u}_{\Rep(\fL^+(G))_\Ran})\simeq \\
\simeq (\Delta_\Ran)_!\circ (\ind^{(\hg,\fL^+(G))_\kappa}_{\fL^+(G)}\otimes \on{Id})(R_{\fL^+(G),\Ran})\in 
\KL(G)_{\kappa,\Ran}\otimes \Rep(\fL^+(G))_\Ran.
\end{multline*}

\medskip

Consider the object
$$(\ind^l_\kappa\otimes \on{Id})(\on{u}_{\QCoh(U)})\simeq 
(\ind^l_\kappa\otimes \on{Id})\circ (j_U\times j_U)^*\circ (\Delta_{\Bun_G})_*(\CO_{\Bun_G})\in 
\Dmod_\kappa(U)\otimes \QCoh(U).$$

We claim that there is a canonically defined map
\begin{equation} \label{e:left induce}
(\Loc_{G,\kappa,U}\otimes \Loc^\QCoh_{G,U})\circ (\ind^{(\hg,\fL^+(G))_\kappa}_{\fL^+(G)}\otimes \on{Id})
(\on{u}_{\Rep(\fL^+(G))_\Ran})
\to (\ind^l_\kappa\otimes \on{Id})(\on{u}_{\QCoh(U)}).
\end{equation}

Namely, the map \eqref{e:left induce} is the composition
\begin{multline} \label{e:left induce again}
(\Loc_{G,\kappa,U}\otimes \Loc^\QCoh_{G,U}) \circ 
(\ind^{(\hg,\fL^+(G))_\kappa}_{\fL^+(G)}\otimes \on{Id})(\on{u}_{\Rep(\fL^+(G))_\Ran})\simeq \\
\simeq (\Loc_{G,\kappa,U}\otimes \on{Id})\circ (\on{Id}\otimes  \Loc^\QCoh_{G,U}) \circ 
(\ind^{(\hg,\fL^+(G))_\kappa}_{\fL^+(G)}\otimes \on{Id})(\on{u}_{\Rep(\fL^+(G))_\Ran})\simeq \\
\simeq (\Loc_{G,\kappa,U}\otimes \on{Id})\circ (\ind^{(\hg,\fL^+(G))_\kappa}_{\fL^+(G)}\otimes \on{Id})
\circ (\on{Id}\otimes  \Loc^\QCoh_{G,U})(\on{u}_{\Rep(\fL^+(G))_\Ran})\simeq \\
\simeq  
(\Loc_{G,\kappa,U}\otimes \on{Id})\circ (\ind^{(\hg,\fL^+(G))_\kappa}_{\fL^+(G)}\otimes \on{Id})
((\Loc^\QCoh_{G,U})^\vee\otimes \on{Id})(\on{u}_{\QCoh(U)})\simeq \\
\simeq (\Loc_{G,\kappa,U}\otimes \on{Id})\circ (\ind^{(\hg,\fL^+(G))_\kappa}_{\fL^+(G)}\otimes \on{Id}) \circ 
((\Loc^\QCoh_{G,U})^R\otimes \on{Id})(\on{u}_{\QCoh(U)})\simeq \\
\simeq 
(\ind^l_\kappa\otimes \on{Id})\circ (\Loc^\QCoh_{G,U}\otimes \on{Id}) \circ 
((\Loc^\QCoh_{G,U})^R\otimes \on{Id})(\on{u}_{\QCoh(U)})\to (\ind^l_\kappa\otimes \on{Id})(\on{u}_{\QCoh(U)}).
\end{multline}

Note that the last arrow in \eqref{e:left induce again} is an isomorphism by \thmref{t:Loc is loc QCoh}. Hence, the 
map \eqref{e:left induce} is an isomorphism.

\sssec{}

Note now that we have canonical isomorphisms
$$(\on{Id}\otimes \oblv^{(\hg,\fL^+(G))_{\kappa'}}_{\fL^+(G)})(\on{u}_{\KL(G)_{\kappa,\Ran}})\simeq 
(\ind^{(\hg,\fL^+(G))_\kappa}_{\fL^+(G)}\otimes \on{Id})(\on{u}_{\Rep(\fL^+(G))_\Ran}),$$
and hence
\begin{multline*}
(\Loc_{G,\kappa,U}\otimes \Loc^\QCoh_{G,U}) \circ (\on{Id}\otimes \oblv^{(\hg,\fL^+(G))_{\kappa'}}_{\fL^+(G)})(\on{u}_{\KL(G)_{\kappa,\Ran}})\simeq \\
\simeq (\Loc_{G,\kappa,U}\otimes \Loc^\QCoh_{G,U}) \circ (\ind^{(\hg,\fL^+(G))_\kappa}_{\fL^+(G)}\otimes \on{Id})(\on{u}_{\Rep(\fL^+(G))_\Ran})
\end{multline*}
and
$$(\on{Id}\otimes \oblv^l_{\kappa'})(\on{u}_{\Dmod_\kappa(U)})\simeq (\ind^l_\kappa\otimes \on{Id})(\on{u}_{\QCoh(U)}).$$

\medskip

Now, unwinding and using \propref{p:Loc as dual}(a), we obtain that under these identifications, the map \eqref{e:right forget} equals the map
\eqref{e:left induce}.

\medskip

Hence, we obtain that the map \eqref{e:right forget} is an isomorphism, as required.

\qed[Second proof of \thmref{t:Loc is loc}]

\section{The composition of localization and coefficient functors} \label{s:coeff Loc}

The goal of this section is to give an expression for the composition
\begin{equation} \label{e:coeff Loc}
\KL(G)_{\crit,\Ran} \overset{\Loc_{G,\crit}}\longrightarrow \Dmod_{\crit}(\Bun_G) \simeq 
\Dmod_{\frac{1}{2}}(\Bun_G) 
\overset{\on{coeff}}\longrightarrow \Whit^!(G)_\Ran
\end{equation}
in terms of factorization homology. 

\medskip

This expression (combined with the compatibility of $\FLE_{G,\crit}$ and $\FLE_{\cG,\infty}$ expressed by
\corref{c:two pairings coarse}) will be used in \secref{s:Langlands functor} in order to show that the Langlands functor is compatible with 
the critical localization and the spectral Poincar\'e series functors via the critical FLE. 

\ssec{The vacuum case}

\sssec{}

Let us specialize the setting of \thmref{t:int loc over BunN kappa chi} to the case
$\kappa=\crit$ and $\CP_T=\rho(\omega_X)$. In this case, the integer that we denoted 
$\delta_{N_{\CP_T}}$ is $\delta_{N_{\rho(\omega_X)}}$. Denote the corresponding line $\fl_{N_{\CP_T}}$ by 
\begin{equation} \label{e:fl N}
\fl_{N_{\rho(\omega_X)}}.
\end{equation}

\sssec{}

We obtain: 

\begin{thm} \label{t:vac coeff of Loc prel}
The composition 
\begin{multline*}
\KL(G)_{\crit,\Ran} \overset{\Loc_{G,\crit}}\longrightarrow \Dmod_\crit(\Bun_G)\overset{\fp^!_\crit}\to
\Dmod(\Bun_{N,\rho(\omega_X)}) \overset{-\overset{*}\otimes \chi^*(\on{exp})}\longrightarrow  \\
\to \Dmod(\Bun_{N,\rho(\omega_X)})\overset{\on{C}^\cdot_\dr(\Bun_{N,\rho(\omega_X)},-)}\longrightarrow \Vect
\end{multline*}
identifies with the functor
\begin{multline*}
\KL(G)_{\crit,\Ran} \overset{\alpha_{\rho(\omega_X),\on{taut}}}\longrightarrow 
\KL(G)_{\crit,\rho(\omega_X),\Ran} \overset{\DS^{\on{enh}}}\longrightarrow \\
\to \fz_\fg\mod^{\on{fact}}_{\Ran} \overset{\on{C}^{\on{fact}}_\cdot(X;\fz_\fg,-)}
\longrightarrow 
\Vect\overset{-\otimes \fl_{N_{\rho(\omega_X)}}[\delta_{N_{\rho(\omega_X)}}]}\longrightarrow \Vect.
\end{multline*}
\end{thm}

\sssec{} \label{sss:just Loc}

Denote by $\Loc_G$ the functor
\begin{equation} \label{e:Loc crit}
\KL(G)_{\crit,\Ran} \overset{\Loc_{G,\crit}}\longrightarrow \Dmod_\crit(\Bun_G) 
\overset{\text{\eqref{e:1/2 vs crit}}}\longrightarrow \Dmod_{\frac{1}{2}}(\Bun_G).
\end{equation} 
  
\medskip

From \eqref{e:tw pullback to BunN diag rho omega} we obtain a commutative diagram
$$
\CD
\Vect @>{\on{Id}}>> \Vect \\
@A{\on{C}^\cdot_\dr(\Bun_{N,\rho(\omega_X)},-)\circ (-\overset{*}\otimes \chi^*(\on{exp}))}AA  @AA{-\otimes \fl^{\otimes \frac{1}{2}}_{G,N_{\rho(\omega_X)}}}A \\
\Dmod(\Bun_{N,\rho(\omega_X)}) & & \Vect \\
@A{\fp^!_\crit}AA   @AA{\on{coeff}_G^{\on{Vac,glob}}}A \\ 
\Dmod_\crit(\Bun_G) @>{\text{\eqref{e:1/2 vs crit}}}>>  \Dmod_{\frac{1}{2}}(\Bun_G) \\
@A{\Loc_{G,\crit}}AA @AA{\Loc_G}A \\
\KL(G)_{\crit,\Ran}  @>{\on{Id}}>> \KL(G)_{\crit,\Ran}.
\endCD
$$

Recall also that
$$\on{coeff}_G^{\on{Vac,glob}}\simeq \on{coeff}_G^{\on{Vac}}[2\delta_{N_{\rho(\omega_X)}}].$$

\sssec{}

Hence, \thmref{t:vac coeff of Loc prel} can be restated as:

\begin{thm} \label{t:vac coeff of Loc}
The composition 
$$\KL(G)_{\crit,\Ran} \overset{\Loc_G}\longrightarrow \Dmod_{\frac{1}{2}}(\Bun_G)
\overset{\on{coeff}_G^{\on{Vac}}}\longrightarrow \Vect \overset{-\otimes \fl^{\otimes \frac{1}{2}}_{G,N_{\rho(\omega_X)}}
\otimes \fl^{\otimes -1}_{N_{\rho(\omega_X)}}[\delta_{N_{\rho(\omega_X)}}]}
\longrightarrow \Vect$$
identifies with the functor
$$\KL(G)_{\crit,\Ran} \overset{\alpha_{\rho(\omega_X),\on{taut}}}\longrightarrow 
\KL(G)_{\crit,\rho(\omega_X),\Ran} \overset{\DS^{\on{enh}}}\longrightarrow 
\fz_\fg\mod^{\on{fact}}_{\Ran} \overset{\on{C}^{\on{fact}}_\cdot(X;\fz_\fg,-)}
\longrightarrow \Vect.$$
\end{thm}

\ssec{Composition of coefficient and localization functors: the general case}

We are now ready to state the general theorem, describing the composition of the functors
$$\KL(G)_{\crit,\Ran} \overset{\Loc_G}\longrightarrow \Dmod_{\frac{1}{2}}(\Bun_G)
\overset{-\otimes \fl^{\otimes \frac{1}{2}}_{G,N_{\rho(\omega_X)}}\otimes \fl^{\otimes -1}_{N_{\rho(\omega_X)}}[\delta_{N_{\rho(\omega_X)}}]}
\longrightarrow \Dmod_{\frac{1}{2}}(\Bun_G)$$
and
$$\on{coeff}_G: \Dmod_{\frac{1}{2}}(\Bun_G)\to \Whit^!(G)_{\Ran}.$$

\sssec{}

Recall that the category $\Whit^!(G)_\Ran$ is the dual of $\Whit_*(G)_\Ran$. 
Hence, the description of the above composition is equivalent to describing the pairing 
\begin{multline} \label{e:Loc and coeff pairing}
\KL(G)_{\crit,\Ran} \otimes \Whit_*(G)_{\Ran} 
\overset{\Loc_G\otimes \on{Id}}\longrightarrow \\
\to \Dmod_{\frac{1}{2}}(\Bun_G)\otimes \Whit_*(G)_{\Ran} 
\overset{(-\otimes \fl^{\otimes \frac{1}{2}}_{G,N_{\rho(\omega_X)}}\otimes \fl^{\otimes -1}_{N_{\rho(\omega_X)}}[\delta_{N_{\rho(\omega_X)}}])\otimes \on{Id}}\longrightarrow \\
\to \Dmod_{\frac{1}{2}}(\Bun_G)\otimes \Whit_*(G)_{\Ran} 
\overset{\on{coeff}_G\otimes \on{Id}}\longrightarrow 
\Whit^!(G)_{\Ran}\otimes \Whit_*(G)_{\Ran}\to \Vect. 
\end{multline}

\sssec{}

Recall the factorization functor $\sP^{\on{loc,enh}}_G$, see \secref{sss:local pairing G}. We will think of it as the functor
\begin{multline*} 
\Whit_*(G) \otimes \KL(G)_\crit \overset{\on{Id}\otimes \alpha_{\rho(\omega_X),\on{taut}}}\longrightarrow 
\Whit_*(G)\otimes \KL(G)_{\crit,\rho(\omega_X)}\to \\
\to \Whit_*(\hg\mod_{\crit,\rho(\omega_X)})\overset{\ol\DS^{\on{enh}}}\to 
\fz_\fg\mod^{\on{fact}}.
\end{multline*} 

Let  $\sP^{\on{loc}}_G$ be the composition of $\sP^{\on{loc,enh}}_G$ with
$$\oblv_{\fz_\G}:\fz_\fg\mod^{\on{fact}}\to \Vect.$$

\sssec{}

We will prove:

\begin{thm} \label{t:Loc and coeff}
The functor \eqref{e:Loc and coeff pairing} identifies canonically with
\begin{multline}  \label{e:Loc and coeff pairing 1}
\KL(G)_{\crit,\Ran} \otimes \Whit_*(G)_\Ran 
\overset{\on{ins.vac}_\Ran\otimes \on{ins.unit}_\Ran}\longrightarrow 
\KL(G)_{\crit,\Ran^{\subseteq}} \otimes \Whit_*(G)_{\Ran^{\subseteq}}\to \\
\to \left(\KL(G)_{\crit}\otimes \Whit_*(G)\right)_{\Ran^{\subseteq}\underset{\Ran}\times \Ran^{\subseteq}}
\overset{\sP^{\on{loc,enh}}_G}\longrightarrow \\
\to \fz_\fg\mod^{\on{fact}}_{\Ran^{\subseteq}\underset{\Ran}\times \Ran^{\subseteq}} 
\overset{\on{C}^{\on{fact}}_\cdot(X;\fz_\fg,-)_{\Ran^{\subseteq}\underset{\Ran}\times \Ran^{\subseteq}}}\longrightarrow
\Dmod(\Ran^{\subseteq}\underset{\Ran}\times \Ran^{\subseteq}) \to \\
\overset{\on{C}^\cdot_c(\Ran^{\subseteq}\underset{\Ran}\times \Ran^{\subseteq},-)}
\longrightarrow  \Vect,
\end{multline}
where the fiber product $\Ran^{\subseteq}\underset{\Ran}\times \Ran^{\subseteq}$ is formed using the projections
$\on{pr}_{\on{big}}:\Ran^{\subseteq}\to \Ran$. 
\end{thm}

\begin{rem} \label{r:simplify Loc and coeff}

Note that the functor \eqref{e:Loc and coeff pairing 1}, appearing in \thmref{t:Loc and coeff} can also be rewritten as
\begin{multline}  \label{e:Loc and coeff pairing 1.5}
\KL(G)_{\crit,\Ran} \otimes \Whit_*(G)_\Ran 
\overset{\on{ins.vac}_\Ran\otimes \on{ins.unit}_\Ran}\longrightarrow 
\KL(G)_{\crit,\Ran^{\subseteq}} \otimes \Whit_*(G)_{\Ran^{\subseteq}}\to \\
\to \left(\KL(G)_{\crit}\otimes \Whit_*(G)\right)_{\Ran^{\subseteq}\underset{\Ran}\times \Ran^{\subseteq}}
\overset{\sP^{\on{loc}}_G}\longrightarrow \Dmod(\Ran^{\subseteq}\underset{\Ran}\times \Ran^{\subseteq}) 
\overset{\on{C}^\cdot_c(\Ran^{\subseteq}\underset{\Ran}\times \Ran^{\subseteq})}\longrightarrow \Vect
\end{multline} 
i.e., instead of $\on{C}^{\on{fact}}_\cdot(X;\fz_\fg,-)_{\Ran^{\subseteq}\underset{\Ran}\times \Ran^{\subseteq}}$
we can use the functor $\oblv_{\fz_\fg,\Ran^{\subseteq}\underset{\Ran}\times \Ran^{\subseteq}}$. 

Indeed, since the functor $\sP^{\on{loc,enh}}_G$ is strictly unital, 
a priori, by \thmref{t:Loc and coeff}, the pairing \eqref{e:Loc and coeff pairing 1} is 

\begin{multline}   \label{e:Loc and coeff pairing 1.75}
\KL(G)_{\crit,\Ran} \otimes \Whit_*(G)_\Ran 
\overset{\on{ins.vac}_\Ran\otimes \on{ins.unit}_\Ran}\longrightarrow 
\KL(G)_{\crit,\Ran^{\subseteq}} \otimes \Whit_*(G)_{\Ran^{\subseteq}}\to \\
\to \left(\KL(G)_{\crit}\otimes \Whit_*(G)\right)_{\Ran^{\subseteq}\underset{\Ran}\times \Ran^{\subseteq}}
\overset{\on{ins.unit}_{\Ran^{\subseteq}\underset{\Ran}\times \Ran^{\subseteq}}}\longrightarrow \\
\to \left(\KL(G)_{\crit}\otimes \Whit_*(G)\right)_{(\Ran^{\subseteq}\underset{\Ran}\times \Ran^{\subseteq})^{\subseteq}}
\overset{\sP^{\on{loc}}_G}\longrightarrow 
\Dmod((\Ran^{\subseteq}\underset{\Ran}\times \Ran^{\subseteq})^{\subseteq}) \to \\
\overset{\on{C}^\cdot_c((\Ran^{\subseteq}\underset{\Ran}\times \Ran^{\subseteq})^{\subseteq})}\longrightarrow \Vect.
\end{multline} 

However, the composition of the first three lines in \eqref{e:Loc and coeff pairing 1.75} lies in the essential image
of the functor
$$\sft^!:\Dmod((\Ran^{\subseteq}\underset{\Ran}\times \Ran^{\subseteq})^{\subseteq,\on{untl}})\to 
\Dmod((\Ran^{\subseteq}\underset{\Ran}\times \Ran^{\subseteq})^{\subseteq}).$$

Hence, by the same mechanism as in \secref{sss:expl unitality fact homology}, the expression in 
\eqref{e:Loc and coeff pairing 1.75} is isomorphic to that in \eqref{e:Loc and coeff pairing 1.5}. 

\end{rem}

\ssec{Proof of \thmref{t:Loc and coeff}.}

\sssec{}

We rewrite the functor
\begin{multline} \label{e:Loc and coeff pairing 2}
\KL(G)_{\crit,\Ran} \otimes \Whit_*(G)_{\Ran} 
\overset{\Loc_G\otimes \on{Id}}\longrightarrow \Dmod_{\frac{1}{2}}(\Bun_G)\otimes \Whit_*(G)_{\Ran} 
\overset{\on{coeff}_G\otimes \on{Id}}\longrightarrow \\
\to \Whit^!(G)_{\Ran}\otimes \Whit_*(G)_{\Ran}\to \Vect. 
\end{multline}
as
\begin{multline} \label{e:Loc and coeff pairing 3}
\KL(G)_{\crit,\Ran} \otimes \Whit_*(G)_{\Ran} 
\overset{\Loc_G\otimes \on{Poinc}_{G,*}}\longrightarrow \\
\to \Dmod_{\frac{1}{2}}(\Bun_G)\otimes \Dmod_{\frac{1}{2},\on{co}}(\Bun_G)\overset{\langle -,-\rangle_{\Bun_G}}
\longrightarrow \Vect,
\end{multline}
where the last arrow is the canonical pairing
$$\Dmod_{\frac{1}{2}}(\Bun_G)\otimes \Dmod_{\frac{1}{2},\on{co}}(\Bun_G)\to \Vect.$$

\sssec{}

By the unital property of $\Loc_G$ and $\on{Poinc}_{G,*}$, we can rewrite 
the functor 
$$\KL(G)_{\crit,\Ran} \otimes \Whit_*(G)_{\Ran} 
\overset{\Loc_G\otimes \on{Poinc}_{G,*}}\longrightarrow 
\Dmod_{\frac{1}{2}}(\Bun_G)\otimes \Dmod_{\frac{1}{2},\on{co}}(\Bun_G)$$
as
\begin{multline} \label{e:Loc and coeff pairing 3.5}
\KL(G)_{\crit,\Ran} \otimes \Whit_*(G)_{\Ran}  \overset{\on{ins.vac}_\Ran\otimes \on{ins.unit}_\Ran}\longrightarrow  \\
\to \KL(G)_{\crit,\Ran^{\subseteq}} \otimes \Whit_*(G)_{\Ran^{\subseteq}}
\overset{(\Loc_G)_{\Ran^\subseteq}\otimes (\on{Poinc}_{G,*})_{\Ran^\subseteq}}\longrightarrow \\
\to \Dmod_{\frac{1}{2}}(\Bun_G)\otimes \Dmod_{\frac{1}{2},\on{co}}(\Bun_G) \otimes 
\Dmod(\Ran^{\subseteq})\otimes \Dmod(\Ran^{\subseteq}) \to \\
\overset{\on{Id}\otimes \on{Id}\otimes \on{C}^\cdot_c(\Ran^{\subseteq}\times \Ran^{\subseteq},-)}\longrightarrow 
\Dmod_{\frac{1}{2}}(\Bun_G)\otimes \Dmod_{\frac{1}{2},\on{co}}(\Bun_G).
\end{multline}

\sssec{} \label{sss:Loc and coeff pairing diag}

In this subsection we will use an analog of \secref{sss:almost unital}, which will allow us to replace 
$$\Ran^{\subseteq}\times \Ran^{\subseteq}\rightsquigarrow \Ran^{\subseteq}\underset{\Ran}\times \Ran^{\subseteq}.$$

\medskip

We note that the composition in the first three lines in \eqref{e:Loc and coeff pairing 3.5} actually factors via
the tensor product of $\Dmod_{\frac{1}{2}}(\Bun_G)\otimes \Dmod_{\frac{1}{2},\on{co}}(\Bun_G)$ with 

\begin{multline*}
\Dmod(\Ran\underset{\Ran,\on{pr}_{\on{small}}}\times \Ran^{\subseteq,\on{untl}})\otimes 
\Dmod(\Ran\underset{\Ran,\on{pr}_{\on{small}}}\times \Ran^{\subseteq,\on{untl}}) 
\overset{\on{Id}\otimes \on{Id}\otimes \sft^!\otimes \sft^!}\longrightarrow \\
\to \Dmod(\Ran^{\subseteq})\otimes \Dmod(\Ran^{\subseteq}) 
\end{multline*}

Hence, by the same principle as in \secref{sss:almost unital}, its further composition with 
\begin{multline*} 
\Dmod_{\frac{1}{2}}(\Bun_G)\otimes \Dmod_{\frac{1}{2},\on{co}}(\Bun_G) \otimes 
\Dmod(\Ran^{\subseteq})\otimes \Dmod(\Ran^{\subseteq}) \to \\
\overset{\on{Id}\otimes \on{Id}\otimes \on{C}^\cdot_c(\Ran^{\subseteq}\times \Ran^{\subseteq},-)}\longrightarrow 
\Dmod_{\frac{1}{2}}(\Bun_G)\otimes \Dmod_{\frac{1}{2},\on{co}}(\Bun_G)
\end{multline*}
is isomorphic to its composition with
\begin{multline*} 
\Dmod_{\frac{1}{2}}(\Bun_G)\otimes \Dmod_{\frac{1}{2},\on{co}}(\Bun_G) \otimes 
\Dmod(\Ran^{\subseteq})\otimes \Dmod(\Ran^{\subseteq}) 
\overset{\on{Id}\otimes \on{Id}\otimes \Delta^!_{\Ran^{\subseteq}}}\longrightarrow \\
\to \Dmod_{\frac{1}{2}}(\Bun_G)\otimes \Dmod_{\frac{1}{2},\on{co}}(\Bun_G) \otimes 
\Dmod(\Ran^{\subseteq}\underset{\Ran}\times \Ran^{\subseteq}) \overset{\on{Id}\otimes \on{Id}\otimes 
\on{C}^\cdot_c(\Ran^{\subseteq}\underset{\Ran}\times \Ran^{\subseteq},-)}\longrightarrow \\
\to \Dmod_{\frac{1}{2}}(\Bun_G)\otimes \Dmod_{\frac{1}{2},\on{co}}(\Bun_G). 
\end{multline*}

\medskip

Hence, we can rewrite \eqref{e:Loc and coeff pairing 3.5} as 

\begin{multline*} 
\KL(G)_{\crit,\Ran} \otimes \Whit_*(G)_{\Ran}  \overset{\on{ins.vac}_\Ran\otimes \on{ins.unit}_\Ran}\longrightarrow  \\
\to \KL(G)_{\crit,\Ran^{\subseteq}} \otimes \Whit_*(G)_{\Ran^{\subseteq}}
\overset{(\Loc_G)_{\Ran^\subseteq}\otimes (\on{Poinc}_{G,*})_{\Ran^\subseteq}}\longrightarrow \\
\to \Dmod_{\frac{1}{2}}(\Bun_G)\otimes \Dmod_{\frac{1}{2},\on{co}}(\Bun_G) \otimes 
\Dmod(\Ran^{\subseteq})\otimes \Dmod(\Ran^{\subseteq}) \to \\
\overset{\on{Id}\otimes \on{Id}\otimes (\Delta_\Ran)^!}\longrightarrow 
\Dmod_{\frac{1}{2}}(\Bun_G)\otimes \Dmod_{\frac{1}{2},\on{co}}(\Bun_G) \otimes \Dmod(\Ran^{\subseteq}\underset{\Ran}\times \Ran^{\subseteq}) \to \\
\overset{\on{Id}\otimes \on{Id}\otimes \on{C}^\cdot_c(\Ran^{\subseteq},\underset{\Ran}\times \Ran^{\subseteq},-)}\longrightarrow 
\Dmod_{\frac{1}{2}}(\Bun_G)\otimes \Dmod_{\frac{1}{2},\on{co}}(\Bun_G),
\end{multline*}
which is the same as 
\begin{multline*} 
\KL(G)_{\crit,\Ran} \otimes \Whit_*(G)_{\Ran}  \overset{\on{ins.vac}_\Ran\otimes \on{ins.unit}_\Ran}\longrightarrow \\
\to \KL(G)_{\crit,\Ran^{\subseteq}} \otimes \Whit_*(G)_{\Ran^{\subseteq}}
\to \left(\KL(G)_{\crit}\otimes \Whit_*(G)\right)_{\Ran^{\subseteq}\underset{\Ran}\times \Ran^{\subseteq}} \to \\
\overset{(\Loc_G\otimes \on{Poinc}_{G,*})_{\Ran^{\subseteq}\underset{\Ran}\times \Ran^{\subseteq}}}\longrightarrow 
\Dmod_{\frac{1}{2}}(\Bun_G)\otimes \Dmod_{\frac{1}{2},\on{co}}(\Bun_G)\otimes 
\Dmod(\Ran^{\subseteq}\underset{\Ran}\times \Ran^{\subseteq}) \to \\
\overset{\on{Id}\otimes \on{Id} \otimes\on{C}^\cdot_c(\Ran^{\subseteq}\underset{\Ran}\times \Ran^{\subseteq},-)}\longrightarrow 
\Dmod_{\frac{1}{2}}(\Bun_G)\otimes \Dmod_{\frac{1}{2},\on{co}}(\Bun_G).
\end{multline*}

\sssec{}

Thus, we obtain that we can rewrite \eqref{e:Loc and coeff pairing 3} as
\begin{multline} \label{e:Loc and coeff pairing 4}
\KL(G)_{\crit,\Ran} \otimes \Whit_*(G)_{\Ran}  \overset{\on{ins.vac}\otimes \on{ins.unit}}\longrightarrow 
\KL(G)_{\crit,\Ran^{\subseteq}} \otimes \Whit_*(G)_{\Ran^{\subseteq}}\to \\
\to \left(\KL(G)_{\crit}\otimes \Whit_*(G)\right)_{\Ran^{\subseteq}\underset{\Ran}\times \Ran^{\subseteq}} \to \\
\overset{(\Loc_G\otimes \on{Poinc}_{G,*})_{\Ran^{\subseteq}\underset{\Ran}\times \Ran^{\subseteq}}}\longrightarrow 
\Dmod_{\frac{1}{2}}(\Bun_G)\otimes \Dmod_{\frac{1}{2},\on{co}}(\Bun_G) \otimes \Dmod(\Ran^{\subseteq}\underset{\Ran}\times \Ran^{\subseteq})
\to \\
\overset{\on{Id}\otimes \on{Id}\otimes \on{C}^\cdot_\dr(\Ran^{\subseteq}\underset{\Ran}\times \Ran^{\subseteq},-)}\longrightarrow 
\Dmod_{\frac{1}{2}}(\Bun_G)\otimes \Dmod_{\frac{1}{2},\on{co}}(\Bun_G)\overset{\langle -,-\rangle_{\Bun_G}}
\longrightarrow \Vect.
\end{multline}

Hence, in order to prove the theorem,  it is enough to identify the composition
\begin{multline} \label{e:Loc and coeff pairing 5}
\left(\KL(G)_{\crit}\otimes \Whit_*(G)\right)_{\Ran^{\subseteq}\underset{\Ran}\times \Ran^{\subseteq}} \to \\
\overset{(\Loc_G\otimes \on{Poinc}_{G,*})_{\Ran^{\subseteq}\underset{\Ran}\times \Ran^{\subseteq}}}\longrightarrow 
\to
\Dmod_{\frac{1}{2}}(\Bun_G)\otimes \Dmod_{\frac{1}{2},\on{co}}(\Bun_G)\otimes \Dmod(\Ran^{\subseteq}\underset{\Ran}\times \Ran^{\subseteq})
 \to \\
\overset{\on{Id}\otimes \on{Id}\otimes \on{C}^\cdot_\dr(\Ran^{\subseteq}\underset{\Ran}\times \Ran^{\subseteq},-)}\longrightarrow 
\Dmod_{\frac{1}{2}}(\Bun_G)\otimes \Dmod_{\frac{1}{2},\on{co}}(\Bun_G) \overset{\langle -,-\rangle_{\Bun_G}}
\longrightarrow \\
\to \Vect 
\overset{(-\otimes \fl^{\otimes \frac{1}{2}}_{G,N_{\rho(\omega_X)}}\otimes\fl^{\otimes -1}_{N_{\rho(\omega_X)}})[\delta_{N_{\rho(\omega_X)}}]}
\longrightarrow \Vect
\end{multline}
with
\begin{multline} \label{e:Loc and coeff pairing 6}
\left(\KL(G)_{\crit}\otimes \Whit_*(G)\right)_{\Ran^{\subseteq}\underset{\Ran}\times \Ran^{\subseteq}}
\overset{\sP^{\on{loc,enh}}_G}\longrightarrow \\
\to \fz_\fg\mod^{\on{fact}}_{\Ran^{\subseteq}\underset{\Ran}\times \Ran^{\subseteq}} 
\overset{\on{C}^{\on{fact}}_\cdot(X;\fz_\fg,-)_{\Ran^{\subseteq}\underset{\Ran}\times \Ran^{\subseteq}}}\longrightarrow
\Dmod(\Ran^{\subseteq}\underset{\Ran}\times \Ran^{\subseteq}) \to \\
\overset{\on{C}^\cdot_\dr(\Ran^{\subseteq}\underset{\Ran}\times \Ran^{\subseteq},-)}
\longrightarrow  \Vect.
\end{multline}

\sssec{}

Applying the functor $(\on{pr}_{\on{big}})_!$, we obtain that it suffices to identify
\begin{multline} \label{e:Loc and coeff pairing 5.5}
\left(\KL(G)_{\crit}\otimes \Whit_*(G)\right)_\Ran 
\overset{(\Loc_G\otimes \on{Poinc}_{G,*})_\Ran}\longrightarrow \\
\to \Dmod_{\frac{1}{2}}(\Bun_G)\otimes \Dmod_{\frac{1}{2},\on{co}}(\Bun_G) \otimes \Dmod(\Ran) 
\overset{\on{Id}\otimes \on{Id}\otimes \on{C}^\cdot_c(\Ran,-)}\longrightarrow \\
\to \Dmod_{\frac{1}{2}}(\Bun_G)\otimes \Dmod_{\frac{1}{2},\on{co}}(\Bun_G) 
\overset{\langle -,-\rangle_{\Bun_G}}\longrightarrow  \Vect 
\overset{(-\otimes \fl^{\otimes \frac{1}{2}}_{G,N_{\rho(\omega_X)}}\otimes\fl^{\otimes -1}_{N_{\rho(\omega_X)}})[\delta_{N_{\rho(\omega_X)}}]}
\longrightarrow \Vect
\end{multline}
with 
\begin{equation} \label{e:Loc and coeff pairing 6.5}
\left(\KL(G)_{\crit}\otimes \Whit_*(G)\right)_\Ran 
\overset{\sP^{\on{loc,enh}}_G}\longrightarrow 
\fz_\fg\mod^{\on{fact}}_\Ran \overset{\on{C}^{\on{fact}}_\cdot(X;\fz_\fg,-)}\longrightarrow \Vect.
\end{equation}

\sssec{}

Note that both functors \eqref{e:Loc and coeff pairing 5.5} and \eqref{e:Loc and coeff pairing 6.5} factor naturally via
$$(\KL(G)_{\crit}\otimes \Whit_*(G))_\Ran 
\to (\KL(G)_{\crit}\underset{\Sph_G}\otimes \Whit_*(G))_\Ran,$$
and in particular via
$$(\KL(G)_{\crit}\otimes \Whit_*(G))_\Ran \to (\KL(G)_{\crit}\underset{\Rep(\cG)}\otimes \Whit_*(G))_\Ran,$$
where $\Rep(\cG)$ maps to $\Sph_G$ via $\on{Sat}^{\on{nv}}_G$. 

\medskip

Hence, using the fact that the action of $\Rep(\cG)$ on $\on{Vac}_{\Whit_*(G)}$ defines an equivalence
$$\Rep(\cG)\to \Whit_*(G),$$
we are reduced to
establishing an identification between
\begin{multline} \label{e:Loc and coeff pairing 5.75}
\KL(G)_{\crit,\Ran} \overset{\on{Id}\otimes \on{Vac}_{\Whit_*(G)}}\longrightarrow
(\KL(G)_{\crit}\otimes \Whit_*(G))_\Ran \to \\
\overset{(\Loc_G\otimes \on{Poinc}_{G,*})_\Ran}\longrightarrow \Dmod_{\frac{1}{2}}(\Bun_G)\otimes \Dmod_{\frac{1}{2},\on{co}}(\Bun_G) \otimes \Dmod(\Ran) 
\overset{\on{Id}\otimes \on{Id}\otimes \on{C}^\cdot_c(\Ran,-)}\longrightarrow \\
\to \Dmod_{\frac{1}{2}}(\Bun_G)\otimes \Dmod_{\frac{1}{2},\on{co}}(\Bun_G) 
\overset{\langle -,-\rangle_{\Bun_G}}\longrightarrow  \Vect 
\overset{(-\otimes \fl^{\otimes \frac{1}{2}}_{G,N_{\rho(\omega_X)}}\otimes\fl^{\otimes -1}_{N_{\rho(\omega_X)}})[\delta_{N_{\rho(\omega_X)}}]}
\longrightarrow \Vect
\end{multline}
and 
\begin{multline} \label{e:Loc and coeff pairing 6.75}
\KL(G)_{\crit,\Ran}\overset{\on{Id}\otimes \on{Vac}_{\Whit_*(G)}}\longrightarrow
(\KL(G)_{\crit}\otimes \Whit_*(G))_\Ran
\overset{\sP^{\on{loc,enh}}_G}\longrightarrow \\
\to \fz_\fg\mod^{\on{fact}}_\Ran
\overset{\on{C}^{\on{fact}}_\cdot(X;\fz_\fg,-)}\longrightarrow \Vect.
\end{multline}

\sssec{}

Using the unital property of $\on{Poinc}_{G,*}$ we can identify 
\eqref{e:Loc and coeff pairing 5.75} with
\begin{multline} \label{e:Loc and coeff pairing 7}
\KL(G)_{\crit,\Ran}
\overset{(\Loc_G)_\Ran\otimes \on{Poinc}^{\on{Vac}}_{G,*}}\longrightarrow \\
\to \Dmod_{\frac{1}{2}}(\Bun_G)\otimes \Dmod(\Ran)\otimes 
 \Dmod_{\frac{1}{2},\on{co}}(\Bun_G) 
\overset{\on{Id}\otimes \on{C}^\cdot_c(\Ran,-)\otimes \on{Id}}\longrightarrow \\
\to \Dmod_{\frac{1}{2}}(\Bun_G)\otimes \Dmod_{\frac{1}{2},\on{co}}(\Bun_G) \overset{\langle -,-\rangle_{\Bun_G}}
\longrightarrow \Vect 
\overset{(-\otimes \fl^{\otimes \frac{1}{2}}_{G,N_{\rho(\omega_X)}}\otimes \fl^{\otimes -1}_{N_{\rho(\omega_X)}})[\delta_{N_{\rho(\omega_X)}}]}
\longrightarrow \Vect.
\end{multline}

By definition, we can identify \eqref{e:Loc and coeff pairing 6.75} with 
\begin{equation} \label{e:Loc and coeff pairing 8}
\KL(G)_{\crit,\Ran} \overset{\alpha_{\rho(\omega_X),\on{taut}}}\simeq 
\KL(G)_{\crit,\rho(\omega_X),\Ran}
\overset{\DS^{\on{enh}}}\longrightarrow  \fz_\fg\mod^{\on{fact}}_\Ran
\overset{\on{C}^{\on{fact}}_\cdot(X;\fz_\fg,-)}\longrightarrow  \Vect. 
\end{equation}

So, it remains to identify \eqref{e:Loc and coeff pairing 7} and \eqref{e:Loc and coeff pairing 8}.

\sssec{}

We rewrite \eqref{e:Loc and coeff pairing 7} as
\begin{equation} \label{e:Loc and coeff pairing 9}
\KL(G)_{\crit,\Ran} \overset{\Loc_G}\longrightarrow \Dmod_{\frac{1}{2}}(\Bun_G)\overset{\on{coeff}_G^{\on{Vac}}}\longrightarrow \Vect
\overset{(-\otimes \fl^{\otimes \frac{1}{2}}_{G,N_{\rho(\omega_X)}}\otimes \fl^{\otimes -1}_{N_{\rho(\omega_X)}})[\delta_{N_{\rho(\omega_X)}}]}
\longrightarrow \Vect.
\end{equation}

\medskip

Now the isomorphism between \eqref{e:Loc and coeff pairing 9} and \eqref{e:Loc and coeff pairing 8} 
is the assertion of \thmref{t:vac coeff of Loc}.

\qed[\thmref{t:Loc and coeff}]

\section{The Hecke eigen-property of critical localization} \label{s:Hecke Loc}

The goal of this section is to establish the Hecke eigenproprety of the functor $\Loc_{G,\crit}$,
which lies in the heart of the manuscript \cite{BD1}. What we do will essentially amount to a 
souped-up version of the construction in {\it loc. cit.}

\medskip

The contents of this (and the next) section are not logically necessary either for this paper, nor for
the other papers in the GLC series. However, here we fill the gap in the literature: we supply 
a proof of a \cite[Theorem 10.3.4]{Ga1}, which is used in \cite[Corollary 4.5.5]{Ga1}, while 
the latter is essential for the GLC series. 

\ssec{Statement of the result}

\sssec{} \label{sss:gen opers}

Consider the prestacks
$$\LS^{\on{mer,glob}}_{\cG,\Ran} \text{ and } \Op^{\on{mer,glob}}_{\cG,\Ran}$$
that attach to $\ul{x}\in \Ran$ the spaces
$$\LS_\cG(X-\ul{x}) \text{ and } \Op_\cG(X-\ul{x}),$$
respectively (see \secref{sss:LS open curve} for what we mean by local systems on an open curve). 

\medskip

Set
$$\Op^{\mf,\on{glob}}_{\cG,\Ran}:=(\LS_\cG\times \Ran) \underset{\LS^{\on{mer,glob}}_{\cG,\Ran}}\times  \Op^{\on{mer,glob}}_{\cG,\Ran}.$$

\medskip

Thus, we have a Cartesian diagram
$$
\CD
\Op^{\mf,\on{glob}}_{\cG,\Ran} @>{\iota^{\mf,\on{glob}}}>> \Op^{\on{mer,glob}}_{\cG,\Ran}  \\
@V{\fr^{\on{glob}}}VV @V{\fr^{\on{glob}}}VV \\
\LS_\cG\times \Ran @>>> \LS^{\on{mer,glob}}_{\cG,\Ran}.
\endCD
$$

All of the spaces and maps in the above diagram have natural unital structures (see \secref{sss:unital spaces} for what this means). 

\medskip 

As we shall see shortly, $\Op^{\mf,\on{glob}}_{\cG,\Ran}$ is a relative ind-affine ind-scheme over $\Ran$. 

\sssec{}

Restriction to the formal disc gives rise to the maps
$$\Op^{\on{mer,glob}}_{\cG,\Ran}\to \Op^\mer_{\cG,\Ran} \text{ and } 
\Op^{\mf,\on{glob}}_{\cG,\Ran}\to \Op^\mf_{\cG,\Ran};$$
we will denote both by $\on{ev}_\Ran$. 

\medskip

Note that the following commutative square is Cartesian
\begin{equation} \label{e:loc vs global mf opers}
\CD
\Op^{\mf,\on{glob}}_{\cG,\Ran} @>{\on{ev}_{\Ran}}>>  \Op^\mf_{\cG,\Ran} \\
@V{\iota^{\mf,\on{glob}}}VV @VV{\iota^\mf}V \\
\Op^{\on{mer,glob}}_{\cG,\Ran} @>{\on{ev}_{\Ran}}>>  \Op^\mer_{\cG,\Ran} 
\endCD
\end{equation}

\sssec{}

For a given $\CZ\to \Ran$, we will change the subscript
$$\Ran \rightsquigarrow \CZ$$
for the corresponding base-changed spaces and maps.\footnote{We remind that, according to our conventions
in \secref{sss:cristalline obj}, when discussing crystals over $\Ran$, given $\CZ\to \Ran$, by default we base change
to $\CZ_\dr$ rather than too $\CZ$.} 

\medskip

We can consider the assignment
$$\CZ \mapsto \QCoh(\Op^{\mf,\on{glob}}_\cG)_\CZ:=\QCoh(\Op^{\mf,\on{glob}}_{\cG,\CZ_\dr})$$
as a crystal of categories over $\Ran$, see \secref{sss:QCoh co fact}. 

\sssec{}

Fix $\CZ\to \Ran$, and recall that according to \secref{ss:z on KL}, we have a canonically defined action of the (symmetric) monoidal category 
$\IndCoh^!(\Op^\mf_\cG)_\CZ$ on $\KL(G)_{\crit,\CZ}$.
Composing with
$$\Upsilon_{\Op^\mf_{\cG,\CZ}}:\QCoh(\Op^\mf_\cG)_\CZ\to
\IndCoh^!(\Op^\mf_\cG)_\CZ,$$
we obtain an action of $\QCoh(\Op^\mf_\cG)_\CZ$ on $\KL(G)_{\crit,\CZ}$.

\medskip

Consider the category
\begin{equation} \label{e:loc-to-glob base change opers}
\KL(G)^{\Op_\cG^{\on{glob}}}_{\crit,\CZ}:=\KL(G)_{\crit,\CZ}\underset{\QCoh(\Op^\mf_\cG)_\CZ}\otimes \QCoh(\Op^{\mf,\on{glob}}_\cG)_\CZ.
\end{equation}

Denote by $\on{Id}\otimes (\on{ev}_\CZ)^*$ the resulting functor
$$\KL(G)_{\crit,\CZ}\to  
\KL(G)_{\crit,\CZ}\underset{\QCoh(\Op^\mf_\cG)_\CZ}\otimes \QCoh(\Op^{\mf,\on{glob}}_\cG)_\CZ=\KL(G)^{\Op_\cG^{\on{glob}}}_{\crit,\CZ}.$$

Note that the assignment
\begin{equation} \label{e:KL Op glob}
\CZ\mapsto \KL(G)^{\Op_\cG^{\on{glob}}}_{\crit,\CZ}
\end{equation}
is naturally a crystal of categories over $\Ran$. 

%

\sssec{}

We define an action of $\Rep(\cG)_\CZ$ on the category \eqref{e:loc-to-glob base change opers} as follows.

\medskip

Recall that we have a symmetric monoidal functor
$$\Loc_{\cG,\CZ}^{\on{spec}}:\Rep(\cG)_\CZ\to \QCoh(\LS_\cG)\otimes \Dmod(\CZ),$$
see \secref{sss:spec loc}.

\medskip

Composing with 
$$(\fr^{\on{glob}})^*:\QCoh(\LS_\cG)\otimes \Dmod(\CZ) \to \QCoh(\Op^{\mf,\on{glob}}_\cG)_\CZ$$
we obtain a symmetric monoidal functor
$$(\fr^{\on{glob}})^*\circ \Loc_{\cG,\CZ}^{\on{spec}}:\Rep(\cG)_\CZ\to \QCoh(\Op^{\mf,\on{glob}}_\cG)_\CZ.$$

We let $\Rep(\cG)_\CZ$ act on \eqref{e:loc-to-glob base change opers} by 
$(\fr^{\on{glob}})^*\circ \Loc_{\cG,\CZ}^{\on{spec}}$ via the second factor.

\sssec{}

The main result of this section is the following:

\begin{thmconstr} \label{t:Hecke crit Loc} 
There exists a canonically defined functor
$$\Loc^\Op_{G,\crit,\CZ}:\KL(G)^{\Op^{\on{glob}}}_{\crit,\CZ}\to
\Dmod_\crit(\Bun_G)\otimes \Dmod(\CZ)$$ such that:

\smallskip

\noindent{\em(a)} The functor
$$\Loc_{G,\crit,\CZ}:\KL(G)_{\crit,\CZ}\to \Dmod_\crit(\Bun_G)\otimes \Dmod(\CZ)$$
factors as 
\begin{equation} \label{e:Hecke crit loc}
\KL(G)_{\crit,\CZ} \overset{\on{Id}\otimes (\on{ev}_\CZ)^*}\longrightarrow 
\KL(G)^{\Op^{\on{glob}}}_{\crit,\CZ}
\overset{\Loc^\Op_{G,\crit,\CZ}}\longrightarrow
\Dmod_\crit(\Bun_G)\otimes \Dmod(\CZ).
\end{equation} 

\smallskip

\noindent{\em(b)} The functor $\Loc^\Op_{G,\crit,\CZ}$ intertwines the above action of $\Rep(\cG)_\CZ$ on 
\eqref{e:loc-to-glob base change opers} and the action of $\Rep(\cG)_\CZ$ on $\Dmod_\crit(\Bun_G)\otimes \Dmod(\CZ)$
obtained from $\on{Sat}_G^{\on{nv}}$ and the action of $\Sph_{G,\CZ}$ on $\Dmod_\crit(\Bun_G)\otimes \Dmod(\CZ)$. 

\end{thmconstr} 

\begin{rem} \label{r:Hecke at x}
The Hecke eigen-property of the functor $\Loc^\Op_{G,\crit,\CZ}$ formulated in point (b) of \thmref{t:Hecke crit Loc}
is not quite the full Hecke property one wants: 

\medskip

For example, when $\CZ=\on{pt}$ and $\CZ\to \Ran$ corresponds to $\ul{x}\in \Ran$, the compatibility assertion in point (b)
only talks about the Hecke action at $\ul{x}$, whereas one wants the compatibility with Hecke action over
the entire Ran space. 

\medskip

See \corref{c:Hecke crit Loc ult} for a stronger assertion, which we will deduce from \propref{p:Hecke crit Loc} below.

\end{rem}

\ssec{The key local construction}

In this subsection we will formulate \thmref{t:ins vac reg}, which is a local counterpart of \thmref{t:Hecke crit Loc}. 

\sssec{}

Let 
$$\Op^{\mf\rightsquigarrow\reg}_{\cG,\Ran^{\subseteq}}\to \Ran^{\subseteq}$$
be the relative indscheme that attaches to
$$(\ul{x}\subseteq \ul{x}')\in \Ran^{\subseteq}$$
the space
$$\Op^{\mf\rightsquigarrow\reg}_{\cG,\ul{x}\subseteq \ul{x}'}
:=\Op_\cG(\cD_{\ul{x}'}-\ul{x})\underset{\LS_\cG(\cD_{\ul{x}'}-\ul{x})}\times \LS_\cG(\cD_{\ul{x}'}).$$

\begin{rem}

Note that $\Op^{\mf\rightsquigarrow\reg}_{\cG,\Ran^{\subseteq}}$ is exactly the geometric object that
encodes the unital-in-correspondences structure on $\Op^\mf_\cG$, see \secref{sss:unital on mf}.

\end{rem} 

\sssec{}

We have the projections
\begin{equation} \label{e:groupoid squig}
\Op^\mf_{\cG,\Ran} \overset{\on{pr}^\Op_{\on{small}}}\longleftarrow \Op^{\mf\rightsquigarrow\reg}_{\cG,\Ran^{\subseteq}}
\overset{\on{pr}^\Op_{\on{big}}}\longrightarrow \Op^\mf_{\cG,\Ran}
\end{equation} 
given by restrictions along
$$\cD_{\ul{x}}-\ul{x}\hookrightarrow \cD_{\ul{x}'}-\ul{x}\hookleftarrow \cD_{\ul{x}'}-\ul{x}',$$
respectively.

\sssec{Example}

For $\ul{x}'=\ul{x}\sqcup \ul{x}''$, we have
$$\Op^{\mf\rightsquigarrow\reg}_\cG(\cD_{\ul{x}'}-\ul{x})=\Op^\mf_\cG(\cD_{\ul{x}}-\ul{x})\times \Op_\cG(\cD_{\ul{x}''}).$$

%
%
%

\sssec{}

For $\CZ\to \Ran$, let us denote by $\Op^{\mf\rightsquigarrow\reg}_{\cG,\CZ^{\subseteq}}$ the 
base-change
$$\CZ\underset{\Ran,\on{pr}_{\on{small}}}\times \Op^{\mf\rightsquigarrow\reg}_{\cG,\Ran^{\subseteq}}.$$

\medskip

The assignment 
$$\CZ\mapsto \QCoh(\Op^{\mf\rightsquigarrow\reg}_\cG)_{\CZ^{\subseteq}}$$
is a crystal of symmetric monoidal categories over $\Ran$, which we will denote by
$$\QCoh(\Op^{\mf\rightsquigarrow\reg}_\cG).$$

\medskip

The map $\on{pr}^\Op_{\on{small}}$ gives rise to a symmetric monoidal functor 
$$(\on{pr}^\Op_{\on{small}})^*:\QCoh(\Op^\mf_\cG)\to \QCoh(\Op^{\mf\rightsquigarrow\reg}_\cG).$$

\sssec{}

Denote
\begin{equation} \label{e:KL subset}
\KL(G)^{\Op^{\on{loc}}_\cG}_\crit:=\KL(G)_\crit \underset{\QCoh(\Op^\mf_\cG)}\otimes \QCoh(\Op^{\mf\rightsquigarrow\reg}_\cG).
\end{equation} 

We will regard it as a crystal of categories over $\Ran$. 
Denote by $\on{Id}\otimes (\on{pr}^\Op_{\on{small}})^*$ the functor
$$\KL(G)_\crit\to \KL(G)_\crit \underset{\QCoh(\Op^\mf_\cG)}\otimes \QCoh(\Op^{\mf\rightsquigarrow\reg}_\cG)=
\KL(G)^{\Op^{\on{loc}}_\cG}_\crit.$$

\sssec{}

Recall that $\on{ins.vac.}$ is a functor between crystals of categories over $\Ran$ 
\begin{equation} \label{e:ins vac KL again}
\KL(G)_\crit \to (\on{pr}_{\on{small}})_*\circ (\on{pr}_{\on{big}})^*(\KL(G)_\crit),
\end{equation} 
and actually a (strict) functor between crystals of categories over $\Ran^{\on{untl}}$
$$\KL(G)_\crit \to (\on{pr}^{\on{untl}}_{\on{small}})_{*,\on{strict}}\circ (\on{pr}^{\on{untl}}_{\on{big}})^*(\KL(G)_\crit).$$

\medskip 

The key construction that we will need says the following: 

\begin{thmconstr} \label{t:ins vac reg} 
There exists a canonically defined functor
\begin{equation} \label{e:ins vac reg mf}
\KL(G)^{\Op^{\on{loc}}_\cG}_\crit 
\overset{\on{ins.vac}^{\mf\rightsquigarrow\reg}}\longrightarrow 
(\on{pr}_{\on{small}})_*\circ (\on{pr}_{\on{big}})^*(\KL(G)_\crit),
\end{equation}
such that $\on{ins.vac.}$ factors as 
$$\KL(G)_\crit \overset{\on{Id}\otimes (\on{pr}^\Op_{\on{small}})^*}\longrightarrow \KL(G)^{\Op^{\on{loc}}_\cG}_\crit
\overset{\on{ins.vac}^{\mf\rightsquigarrow\reg}}\longrightarrow 
(\on{pr}_{\on{small}})_*\circ (\on{pr}_{\on{big}})^*(\KL(G)_\crit).$$

Furthermore, the functor $\on{ins.vac}^{\mf\rightsquigarrow\reg}$ is linear with respect to 
$$(\on{pr}_{\on{small}})_*\circ (\on{pr}_{\on{big}})^*(\QCoh(\Op^\mf_\cG)),$$
where:

\begin{itemize}

\item $(\on{pr}_{\on{small}})_*\circ (\on{pr}_{\on{big}})^*(\QCoh(\Op^\mf_\cG))$ acts on 
the left-hand side via $(\on{pr}^\Op_{\on{big}})^*$;

\smallskip

\item $(\on{pr}_{\on{small}})_*\circ (\on{pr}_{\on{big}})^*(\QCoh(\Op^\mf_\cG))$ acts on 
the right-hand side by applying the functor $(\on{pr}_{\on{small}})_*\circ (\on{pr}_{\on{big}})^*$ to the 
$\QCoh(\Op^\mf_\cG)$-action on $\KL(G)_\crit$.

\end{itemize}

\end{thmconstr}

The proof of \thmref{t:ins vac reg} will be given in \secref{s:ins vac reg}. 

\begin{rem}

Note that \thmref{t:ins vac reg} sounds semantically close to \thmref{t:Hecke crit Loc}. And indeed, 
\thmref{t:Hecke crit Loc} will be deduced from \thmref{t:ins vac reg} by a manipulation that involves 
factorization homology.

\medskip

Yet, \thmref{t:ins vac reg} is purely local, while \thmref{t:Hecke crit Loc} is ``local-to-global". 

\medskip

That said, we will need \thmref{t:ins vac reg} \emph{not only} for the proof of \thmref{t:Hecke crit Loc}.
In the next subsection, we will use it to extract a stronger Hecke property for the functor $\Loc^\Op_{G,\crit,\CZ}$. 

\end{rem} 

\sssec{}

Let us explain what \thmref{t:ins vac reg} says at the pointwise level, i.e., for a fixed $(\ul{x}\subseteq \ul{x}')\in \Ran^{\subseteq}$.

\medskip

Write
$$\ul{x}'=\ul{x}\sqcup \ul{x}''.$$

The corresponding functor 
$$\on{ins.vac}_{\ul{x}\subseteq \ul{x}'}: \KL(G)_{\crit,\ul{x}}\to \KL(G)_{\crit,\ul{x}'}\simeq \KL(G)_{\crit,\ul{x}}\otimes \KL(G)_{\crit,\ul{x}''}$$
acts as
$$\CM\mapsto \CM\otimes \on{Vac}(G)_{\crit,\ul{x}''}.$$

The pointwise statement of \thmref{t:ins vac reg} is that this functor can be factored as 
$$\KL(G)_{\crit,\ul{x}} \to \KL(G)_{\crit,\ul{x}} \otimes \QCoh(\Op^\reg_{\cG,\ul{x}''}) \to \KL(G)_{\crit,\ul{x}}\otimes \KL(G)_{\crit,\ul{x}''},$$
where the second arrow is $\QCoh(\Op^\mf_{\cG,\ul{x}''})$-linear. 

\medskip

In other words, we are saying that the object $\on{Vac}(G)_{\crit,\ul{x}''}\in  \KL(G)_{\crit,\ul{x}''}$ naturally lifts to an object of the category
$$\on{Funct}_{\QCoh(\Op^\mf_{\cG,\ul{x}''})}(\QCoh(\Op^\reg_{\cG,\ul{x}''}),\KL(G)_{\crit,\ul{x}''}).$$

This lift is the basic feature of the vacuum object: it says that the structure of factorization $\fz_\fg$-module 
on $\on{Vac}(G)_\crit$ as an object of $\KL(G)_\crit$ is obtained by restriction from a structure of \emph{commutative} 
$\fz_\fg$-module.

\medskip

The proof of \thmref{t:ins vac reg} will amount to spelling out the above construction 
in the factorization setting. 

\begin{rem}

In the proof of \thmref{t:ins vac reg} that we will give, we will avoid using the FLE. We do this for aesthetical reasons:
the construction of the functor \thmref{t:ins vac reg} is more or less tautological if one says the right words. 

\medskip

However, if we use the FLE, there would be almost nothing to prove: the FLE, combined with the equivalence $\Theta_{\Op_\cG^\mf}$, 
allows us to identity the two sides in \eqref{e:ins vac reg mf} with 
$$\IndCoh^!(\Op^{\mf\rightsquigarrow\reg}_\cG)  \text{ and } 
(\on{pr}_{\on{small}})_*\circ (\on{pr}_{\on{big}})^*(\IndCoh^!(\Op_\cG^\mf)),$$
respectively, and the functor in question is obtained by taking direct image along
$$\Op^{\mf\rightsquigarrow\reg}_{\cG,\Ran^{\subseteq}} \overset{\on{pr}^\Op_{\on{big}}}\longrightarrow 
\Ran^{\subseteq}\underset{\on{pr}_{\on{big}},\Ran}\times \Op^\mf_{\cG,\Ran}.$$

That said, the functor $\on{ins.vac}^{\mf\rightsquigarrow\reg}$ that we will construct does reproduce the above functor,
by the nature of the FLE. 

\end{rem}

\begin{rem} \label{r:grpd acts on KL}

Note that the space $\Op^{\mf\rightsquigarrow\reg}_{\cG,\Ran^{\subseteq}}$, equipped with the maps \eqref{e:groupoid squig}, 
has a structure of groupoid\footnote{This groupoid structure is part of the structure of being unital-in-correspondences, see \secref{ss:untl in corr}.}
acting on $\Op^\mf_{\cG,\Ran}$,
with the composition given by
\begin{multline*} 
\Op^{\mf\rightsquigarrow\reg}_\cG(\cD_{\ul{x}'}-\ul{x})\underset{\Op^\mf_\cG(\cD_{\ul{x}'}-\ul{x'})}\times 
\Op^{\mf\rightsquigarrow\reg}_\cG(\cD_{\ul{x}''}-\ul{x}')\simeq \\
\simeq (\Op_\cG(\cD_{\ul{x}'}-\ul{x})\underset{\Op_\cG(\cD_{\ul{x}'}-\ul{x'})}\times 
\Op_\cG(\cD_{\ul{x}''}-\ul{x}'))\underset{(\LS_\cG(\cD_{\ul{x}'}-\ul{x})\underset{\LS_\cG(\cD_{\ul{x}'}-\ul{x'})}\times 
\LS_\cG(\cD_{\ul{x}''}-\ul{x}'))}\times \LS_\cG(\cD_{\ul{x}''})\to \\
\to \Op_\cG(\cD_{\ul{x}''}-\ul{x})\underset{\LS_\cG(\cD_{\ul{x}''}-\ul{x})}\times \LS_\cG(\cD_{\ul{x}''})=
\Op^{\mf\rightsquigarrow\reg}_\cG(\cD_{\ul{x}''}-\ul{x}),
\end{multline*} 
where the middle arrow is given by gluing along 
$$\cD_{\ul{x}'}-\ul{x}\underset{\cD_{\ul{x}'}-\ul{x'}}\sqcup \cD_{\ul{x}''}-\ul{x}'\simeq \cD_{\ul{x}''}-\ul{x}.$$

\medskip

Note that one can interpret \thmref{t:ins vac reg} as saying that the above action of $\Op^{\mf\rightsquigarrow\reg}_{\cG,\Ran^{\subseteq}}$
on $\Op^\mf_{\cG,\Ran}$ can be lifted to an action on $\KL(G)_\crit$ \emph{at the level of 1-morphisms}.

\medskip

In fact, \thmref{t:ins vac reg} has a natural upgrade to a statement that we have a full datum of action of 
$\Op^{\mf\rightsquigarrow\reg}_{\cG,\Ran^{\subseteq}}$ on $\KL(G)_\crit$. This is again an automatic if
we allow ourselves to use the FLE. 

\end{rem}

\ssec{The expanded Hecke eigen-property}

In this subsection we will assume \thmref{t:Hecke crit Loc} and explain that the Hecke property of the functor 
$\Loc^\Op_{G,\crit,\CZ}$ stated in point (b) of the theorem 
implies a stronger property (the difference between the two versions of the Hecke 
property is what was alluded to in Remark \ref{r:Hecke at x}). 

\sssec{}

Let $\CZ\to \Ran$ be as above. Consider the space
\begin{equation} \label{e:global oper added}
\Op^{\mf\rightsquigarrow\reg,\on{glob}}_{\cG,\CZ^{\subseteq}}:=\Op^{\mf,\on{glob}}_{\cG,\CZ}\underset{\CZ,\on{pr}_{\on{small},\CZ}}\times \CZ^{\subseteq}.
\end{equation}

Denote by $\on{pr}^{\Op^{\on{glob}}}_{\on{small},\CZ}$ the projection
$$\Op^{\mf\rightsquigarrow\reg,\on{glob}}_{\cG,\CZ^{\subseteq}}\to \Op^{\mf,\on{glob}}_{\cG,\CZ}.$$

\medskip

Note that we have a naturally defined ind-closed embedding
\begin{equation} \label{e:adding poles to the global oper}
\on{pr}^{\Op^{\on{glob}}}_{\on{big},\CZ}:\Op^{\mf\rightsquigarrow\reg,\on{glob}}_{\cG,\CZ^{\subseteq}}\hookrightarrow
\Op^{\mf,\on{glob}}_{\cG,\CZ^{\subseteq}}.
\end{equation}

Denote by $\on{pr}^{\Op^{\on{glob}}}_{\on{big}}$ the composition of \eqref{e:adding poles to the global oper} with the projection
$$\Op^{\mf,\on{glob}}_{\cG,\CZ^{\subseteq}}\overset{\on{id}\times \on{pr}_{\on{big}}}\longrightarrow \Op^{\mf,\on{glob}}_{\cG,\Ran}.$$

\sssec{Example}
Let us explain what the map \eqref{e:adding poles to the global oper} looks like for $\CZ=\on{pt}$ 
so that $\CZ\to \Ran$ corresponds to $\ul{x}\in \Ran$ (in which case $\CZ^{\subseteq}=\Ran_{\ul{x}}$, so a point
on it corresponds to $\ul{x}\subseteq \ul{x}'$).

\medskip

Then a point of the left-hand (resp., right-hand) side in \eqref{e:adding poles to the global oper} is a local system 
on $X$ with an oper structure away from $\ul{x}$ (resp., $\ul{x}'$), and the map \eqref{e:adding poles to the global oper}
is given by restriction along
$$X-\ul{x}'\subseteq X-\ul{x}.$$

\sssec{}

Consider the map
$$\Op^{\mf\rightsquigarrow\reg,\on{glob}}_{\cG,\CZ^{\subseteq}}=
\Op^{\mf,\on{glob}}_{\cG,\CZ}\underset{\CZ,\on{pr}_{\on{small},\CZ}}\times \CZ^{\subseteq} \overset{\fr^{\on{glob}}\times \on{id}}\longrightarrow 
(\LS_\cG\times \CZ) \underset{\CZ}\times \CZ^{\subseteq}=\LS_\cG\times \CZ^{\subseteq},$$
to be denoted $'\fr^{\on{glob}}$.

\medskip

Thus, we can consider the symmetric monoidal functor 
$$({}'\fr^{\on{glob}})^*\circ \Loc^{\on{spec}}_{\cG,\CZ^{\subseteq}}:
\Rep(\cG)_{\CZ^{\subseteq}}\to \QCoh(\Op^{\mf\rightsquigarrow\reg,\on{glob}}_\cG)_{\CZ^{\subseteq}}.$$

\sssec{} \label{sss:action Rep subset}

We can view $\Op^{\mf\rightsquigarrow\reg,\on{glob}}_{\cG,\CZ^{\subseteq}}$ as mapping to $\Op^\mf_{\cG,\CZ}$ by 
$$\Op^{\mf\rightsquigarrow\reg,\on{glob}}_{\cG,\CZ^{\subseteq}}=
\Op^{\mf,\on{glob}}_{\cG,\CZ}\underset{\CZ,\on{pr}_{\on{small},\CZ}}\times \CZ^{\subseteq}\to
\Op^{\mf,\on{glob}}_{\cG,\CZ} \overset{\on{ev}_\CZ}\to \Op^{\mf}_{\cG,\CZ}.$$

Hence, we can form the category
\begin{equation} \label{e:adding poles KL}
\KL(G)_{\crit,\CZ}\underset{\QCoh(\Op^\mf_\cG)_\CZ}\otimes 
\QCoh(\Op^{\mf\rightsquigarrow\reg,\on{glob}}_\cG)_{\CZ^{\subseteq}}.
\end{equation}

Since 
$$\QCoh(\Op^{\mf,\on{glob}}_\cG)_\CZ\underset{\Dmod(\CZ)}\otimes \Dmod(\CZ^{\subseteq})\to
\QCoh(\Op^{\mf\rightsquigarrow\reg,\on{glob}}_\cG)_{\CZ^{\subseteq}}$$
is an equivalence, we can rewrite \eqref{e:adding poles KL} as
\begin{multline} \label{e:adding poles KL bis}
\KL(G)_{\crit,\CZ}\underset{\QCoh(\Op^\mf_\cG)_\CZ}\otimes 
\QCoh(\Op^{\mf,\on{glob}}_\cG)_\CZ \underset{\Dmod(\CZ)}\otimes \Dmod(\CZ^{\subseteq})= \\
=\KL(G)^{\Op_\cG^{\on{glob}}}_{\crit,\CZ}\underset{\Dmod(\CZ)}\otimes \Dmod(\CZ^{\subseteq}).
\end{multline} 

Thus, we obtain that the category 
$$\KL(G)^{\Op_\cG^{\on{glob}}}_{\crit,\CZ}\underset{\Dmod(\CZ)}\otimes \Dmod(\CZ^{\subseteq})$$
carries a monoidal action of $\Rep(\cG)_{\CZ^{\subseteq}}$. 

\sssec{}

Consider the functor
\begin{multline} \label{e:Loc Op added}
\KL(G)^{\Op_\cG^{\on{glob}}}_{\crit,\CZ} \underset{\Dmod(\CZ)}\otimes \Dmod(\CZ^{\subseteq})  
\overset{\Loc^\Op_{G,\crit,\CZ}\otimes \on{Id}}\longrightarrow \\
\to (\Dmod_\crit(\Bun_G)\otimes \Dmod(\CZ)) \underset{\Dmod(\CZ)}\otimes \Dmod(\CZ^{\subseteq}) 
\simeq \Dmod_\crit(\Bun_G)\otimes \Dmod(\CZ^{\subseteq}).
\end{multline} 

We claim:

\begin{prop}  \label{p:Hecke crit Loc}
The functor \eqref{e:Loc Op added} intertwines the actions of $\Rep(\cG)_{\CZ^{\subseteq}}$
on the two sides, where the action of $\Rep(\cG)_{\CZ^{\subseteq}}$ on the left-hand sides is the
one specified in \secref{sss:action Rep subset}, and on the right-hand side it is
obtained from $\on{Sat}_G^{\on{nv}}$ and the action of $\Sph_{G,\CZ^{\subseteq}}$ on $\Dmod_\crit(\Bun_G)\otimes \Dmod(\CZ^{\subseteq})$. 
\end{prop} 

\ssec{Proof of \propref{p:Hecke crit Loc}}

\sssec{} \label{sss:pr Op big}

Note that we have a Cartesian square
$$
\CD
\Op^{\mf\rightsquigarrow\reg,\on{glob}}_{\cG,\CZ^{\subseteq}} @>{\on{pr}^{\Op^{\on{glob}}}_{\on{big},\CZ}}>> \Op^{\mf,\on{glob}}_{\cG,\CZ^{\subseteq}} \\
@V{\on{ev}_{\CZ^{\subseteq}}}VV @VV{\on{ev}_{\CZ^{\subseteq}}}V \\
\Op^{\mf\rightsquigarrow\reg}_{\cG,\CZ^{\subseteq}} @>>{\on{pr}^\Op_{\on{big},\CZ}}> 
\Op^\mf_{\cG,\CZ^{\subseteq}},
\endCD
$$
where $\on{pr}^\Op_{\on{big},\CZ}$ is the map whose composition with
$$\Op^\mf_{\cG,\CZ^{\subseteq}}\simeq \CZ^{\subseteq}\underset{\on{pr}_{\on{big}},\Ran}\times \Op^\mf_{\cG,\Ran}\to
\Op^\mf_{\cG,\Ran}$$
is the map
$$\Op^\mf_{\cG,\CZ^{\subseteq}}\to \Op^\mf_{\cG,\Ran^{\subseteq}}\overset{\on{pr}^\Op_{\on{big}}}\longrightarrow \Op^\mf_{\cG,\Ran}.$$

\medskip

We claim:

\begin{lem} \label{l:ten prod glob mixed Op}
The functor
$$\QCoh(\Op^{\mf\rightsquigarrow\reg}_\cG)_{\CZ^{\subseteq}}\underset{\QCoh(\Op^\mf_\cG)_{\CZ^{\subseteq}}}\otimes
\QCoh(\Op^{\mf,\on{glob}}_\cG)_{\CZ^{\subseteq}}
\to \QCoh(\Op^{\mf\rightsquigarrow\reg,\on{glob}}_\cG)_{\CZ^{\subseteq}}$$
is an equivalence.
\end{lem}

The lemma will be proved in \secref{ss:base change mf}. 

\sssec{}

Consider the category 
\begin{multline} \label{e:KL Op 1}
\KL(G)^{\Op_\cG^{\on{glob}}}_{\crit,\CZ} \underset{\Dmod(\CZ)}\otimes \Dmod(\CZ^{\subseteq})  = \\
=\KL(G)_{\crit,\CZ}\underset{\QCoh(\Op^\mf_\cG)_{\CZ}}\otimes \QCoh(\Op^{\mf\rightsquigarrow\reg,\on{glob}}_\cG)_\CZ
\underset{\Dmod(\CZ)}\otimes \Dmod(\CZ^{\subseteq})\simeq \\
\simeq 
\KL(G)_{\crit,\CZ}\underset{\QCoh(\Op^\mf_\cG)_{\CZ}}\otimes \QCoh(\Op^{\mf\rightsquigarrow\reg,\on{glob}}_\cG)_{\CZ^{\subseteq}}.
\end{multline} 

By \lemref{l:ten prod glob mixed Op}, we can rewrite it as
$$\KL(G)_{\crit,\CZ}\underset{\QCoh(\Op^\mf_\cG)_{\CZ}}\otimes 
\QCoh(\Op^{\mf\rightsquigarrow\reg}_\cG)_{\CZ^{\subseteq}}\underset{\QCoh(\Op^\mf_\cG)_{\CZ^{\subseteq}}}\otimes
\QCoh(\Op^{\mf,\on{glob}}_\cG)_{\CZ^{\subseteq}},$$
i.e.,
$$\KL^{\Op^{\on{loc}}_\cG}_{\crit,\CZ} \underset{\QCoh(\Op^\mf_\cG)_{\CZ^{\subseteq}}}\otimes
\QCoh(\Op^{\mf,\on{glob}}_\cG)_{\CZ^{\subseteq}},$$
where 
$$\CZ\mapsto \KL^{\Op^{\on{loc}}_\cG}_{\crit,\CZ}$$
is the factorization category \eqref{e:KL subset}.

\medskip

Using the functor $\on{ins.vac}^{\mf\rightsquigarrow\reg}$ we obtain a functor
\begin{multline} \label{e:unital for KL Op}
\KL(G)^{\Op_\cG^{\on{glob}}}_{\crit,\CZ} \underset{\Dmod(\CZ)}\otimes \Dmod(\CZ^{\subseteq})  \simeq \\
\simeq \KL^{\Op^{\on{loc}}_\cG}_{\crit,\CZ} \underset{\QCoh(\Op^\mf_\cG)_{\CZ^{\subseteq}}}\otimes
\QCoh(\Op^{\mf,\on{glob}}_\cG)_{\CZ^{\subseteq}} \to \\
%
%
\overset{\on{ins.vac}^{\mf\rightsquigarrow\reg}_\CZ\otimes \on{Id}}\longrightarrow
\KL(G)_{\crit,\CZ^{\subseteq}} \underset{\QCoh(\Op^\mf_\cG)_{\CZ^{\subseteq}}}\otimes
\QCoh(\Op^{\mf,\on{glob}}_\cG)_{\CZ^{\subseteq}} \simeq \\
\simeq \KL(G)^{\Op_\cG^{\on{glob}}}_{\crit,\CZ^{\subseteq}}. 
\end{multline}

\begin{rem} \label{r:unital for KL Op}
In the spirit of Remark \ref{r:grpd acts on KL}, one can show that the functors \eqref{e:unital for KL Op} upgrade to
a \emph{local unital structure} (see \secref{sss:unitality loc} for what this means) on the crystal of categories over $\Ran$ given by
\eqref{e:KL Op glob}. 
\end{rem} 

\sssec{}

The following property will be embedded into the construction of the assignment
\begin{equation} \label{e:Loc Op added 0}
\CZ\mapsto \Loc^\Op_{G,\crit,\CZ}:
\end{equation}

The composition
\begin{multline} \label{e:Loc Op added 1}
\KL(G)^{\Op_\cG^{\on{glob}}}_{\crit,\CZ}  \overset{\on{Id}\otimes \on{pr}^!_{\on{small}}}\longrightarrow 
\KL(G)^{\Op_\cG^{\on{glob}}}_{\crit,\CZ} \underset{\Dmod(\CZ)}\otimes \Dmod(\CZ^{\subseteq}) 
%
%
%
\overset{\text{\eqref{e:unital for KL Op}}}\longrightarrow \\
\to \KL(G)^{\Op_\cG^{\on{glob}}}_{\crit,\CZ^{\subseteq}} \overset{\Loc^\Op_{G,\crit,\CZ^{\subseteq}}}\longrightarrow 
\Dmod_\crit(\Bun_G)\otimes \Dmod(\CZ^{\subseteq})
\end{multline} 
identifies with the functor
\begin{multline} \label{e:Loc Op added 2}
\KL(G)^{\Op_\cG^{\on{glob}}}_{\crit,\CZ} 
\overset{\Loc^\Op_{G,\crit,\CZ}}\longrightarrow \\
\to \Dmod_\crit(\Bun_G)\otimes \Dmod(\CZ) \overset{\on{Id}\otimes \on{pr}_{\on{small}}^!}\longrightarrow 
\Dmod_\crit(\Bun_G)\otimes \Dmod(\CZ^{\subseteq}).
\end{multline} 

\begin{rem}
In the spirit of Remark \ref{r:unital for KL Op} one can show that the isomorphism between 
\eqref{e:Loc Op added 1} and \eqref{e:Loc Op added 2} upgrades to a \emph{unital structure}
(see \secref{sss:strictly unital} for what this means) on the assignment \eqref{e:Loc Op added 0}. 
\end{rem}

\sssec{}

Since the functors
\begin{multline} \label{e:pre-Loc Op added 1}
\KL(G)^{\Op_\cG^{\on{glob}}}_{\crit,\CZ} \underset{\Dmod(\CZ)}\otimes \Dmod(\CZ^{\subseteq}) 
\overset{\text{\eqref{e:unital for KL Op}}}\longrightarrow \\
\to \KL(G)^{\Op_\cG^{\on{glob}}}_{\crit,\CZ^{\subseteq}} \overset{\Loc^\Op_{G,\crit,\CZ^{\subseteq}}}\longrightarrow 
\Dmod_\crit(\Bun_G)\otimes \Dmod(\CZ^{\subseteq})
\end{multline} 
and \eqref{e:Loc Op added} are $\Dmod(\CZ^{\subseteq})$-linear, and the isomorphism between 
\eqref{e:Loc Op added 1} and \eqref{e:Loc Op added 2} is $\Dmod(\CZ)$-linear, we obtain that 
the functor \eqref{e:pre-Loc Op added 1} identifies with \eqref{e:Loc Op added}.

\medskip

Thus, in order to prove \propref{p:Hecke crit Loc}, it suffices to construct the datum of compatibility with the $\Rep(\cG)_{\CZ^{\subseteq}}$-action 
for the functor \eqref{e:pre-Loc Op added 1}. 

\sssec{}

We will show that each of the two arrows in  \eqref{e:pre-Loc Op added 1} is compatible with the $\Rep(\cG)_{\CZ^{\subseteq}}$-action. 

\medskip

For the second arrow, 
 this follows from \thmref{t:Hecke crit Loc}(b). 

\sssec{}

For the first arrow, unwinding the definition of the functor \eqref{e:unital for KL Op}, 
we have to show that the functor 
\begin{multline*} 
\QCoh(\Op^{\mf\rightsquigarrow\reg}_\cG)_{\CZ^{\subseteq}}\underset{\QCoh(\Op^\mf_\cG)_{\CZ^{\subseteq}}}\otimes
\QCoh(\Op^{\mf,\on{glob}}_\cG)_{\CZ^{\subseteq}}\to \\
\to \QCoh(\Op^{\mf\rightsquigarrow\reg,\on{glob}}_\cG)_{\CZ^{\subseteq}}
\end{multline*} 
is $\Rep(\cG)_{\CZ^{\subseteq}}$-linear, where:

\smallskip

\begin{itemize} 
 
\item $\Rep(\cG)_{\CZ^{\subseteq}}$ acts on the right-hand side via $\Loc^{\on{spec}}_{\cG,\CZ^{\subseteq}}$ followed by pullback along
\begin{equation} \label{e:r Op 1}
'\fr^{\on{glob}}:\Op^{\mf\rightsquigarrow\reg,\on{glob}}_{\cG,\CZ^{\subseteq}}\to \LS_\cG\times \CZ^{\subseteq};
\end{equation} 

\item $\Rep(\cG)_{\CZ^{\subseteq}}$ acts on the left-hand side on the second factor via $\Loc^{\on{spec}}_{\cG,\CZ^{\subseteq}}$ 
followed by pullback along
$$\fr^{\on{glob}}:\Op^{\mf,\on{glob}}_{\cG,\CZ^{\subseteq}} \to \LS_\cG\times \CZ^{\subseteq}.$$

\end{itemize}

The required compatibility follows from the fact that
$$'\fr^{\on{glob}}=\fr^{\on{glob}}\circ \on{pr}^{\Op^{\on{glob}}}_{\on{big},\CZ}$$
as maps
$$\Op^{\mf\rightsquigarrow\reg,\on{glob}}_{\cG,\CZ^{\subseteq}}\rightrightarrows \LS_\cG\times \CZ^{\subseteq}.$$

\ssec{The \emph{integrated} Hecke eigen-property}

In this subsection we will apply the paradigm of \secref{ss:loc act param}, and deduce 
an ultimate form of compatibility of the functor $\Loc^\Op_{G,\crit,\CZ}$ with the Hecke action.

\sssec{}

Consider $\Rep(\cG)$ as a \emph{unital} monoidal factorization category. Note that we consider the Hecke action 
as a local Ran-unital action of $\Rep(\cG)$ on $\Dmod_\crit(\Bun_G)$ (see \secref{sss:local actions} for what this means), i.e., 
$$\Dmod_\crit(\Bun_G)\in \Rep(\cG)^{\on{loc}}\mmod.$$

\medskip

In particular, the action of $\Rep(\cG)_\Ran$, endowed with the \emph{convolution} monoidal structure 
(i.e., the monoidal category $(\Rep(\cG)_\Ran)^\star$, see \secref{sss:conv non-unital}), on $\Dmod_\crit(\Bun_G)$ factors via
$$(\Rep(\cG)_\Ran)^\star \twoheadrightarrow \Rep(\cG)_{\Ran^{\on{untl}},\on{indep}}.$$

\medskip

For any $\CZ\to \Ran$, we will consider 
$$\Dmod_\crit(\Bun_G)\otimes \Dmod(\CZ)$$
as a module over 
$$\Rep(\cG)_{\Ran^{\on{untl}},\on{indep}}\otimes \Dmod(\CZ).$$

\sssec{}

Note also that the spectral localization functor can be viewed as a (strictly) unital (symmetric) monoidal
functor
$$\ul\Loc^{\on{spec}}_\cG:\ul\Rep(\cG)\to \QCoh(\LS_\cG)\otimes \ul\Dmod(\Ran^{\on{untl}}),$$
where $\ul\Rep(\cG)$ is the crystal of categories over $\Ran^{\on{untl}}$ corresponding to 
$\Rep(\cG)$, viewed as a factorization category. 

\medskip

In particular, the functor 
$$\Loc_\cG^{\on{spec}}:\Rep(\cG)_\Ran\to \QCoh(\LS_\cG)$$
factors as
$$\Rep(\cG)_\Ran\twoheadrightarrow \Rep(\cG)_{\Ran^{\on{untl}},\on{indep}}\to \QCoh(\LS_\cG).$$

Given $\CZ\to \Ran$, we will consider the category 
$$\KL(G)^{\Op_\cG^{\on{glob}}}_{\crit,\CZ}=\KL(G)_{\crit,\CZ}\underset{\QCoh(\Op^\mf_\cG)_\CZ}\otimes \QCoh(\Op^{\mf,\on{glob}}_\cG)_\CZ$$
as acted on by 
$$\Rep(\cG)_{\Ran^{\on{untl}},\on{indep}}\otimes \Dmod(\CZ)$$
via the projection
$$\Op^{\mf,\on{glob}}_{\cG,\CZ}\to \LS_\cG\times \CZ$$
and the action on the second factor. 

\sssec{}

We claim:

\begin{cor} \label{c:Hecke crit Loc}
The functor $$\Loc^\Op_{G,\crit,\CZ}:\KL(G)^{\Op^{\on{glob}}}_{\crit,\CZ}\to
\Dmod_\crit(\Bun_G)\otimes \Dmod(\CZ)$$ 
is compatible with the action of $\Rep(\cG)_{\Ran^{\on{untl}},\on{indep}}\otimes \Dmod(\CZ)$.
\end{cor}

\begin{proof}

First, it is easy to see that we can assume that $\CZ$ is an affine scheme $S$. By \corref{c:get rid of parameters}, 
we can consider both sides, i.e., 
$$\Dmod_\crit(\Bun_G)\otimes \Dmod(S) \text{ and } \KL(G)^{\Op_\cG^{\on{glob}}}_{\crit,S},$$
as objects of 
$$\Rep(\cG)_S^{\on{loc,untl}}\mmod$$
(see \secref{sss:local actions params} for the notation), and we have to show that $\Loc^\Op_{G,\crit,S}$ extends to a map
inside this category (see \corref{c:get rid of parameters}).

\medskip

Now, \propref{p:Hecke crit Loc} implies that the functor $\Loc^\Op_{G,\crit,S}$ gives rise to a functor 
between the images of these two objects under the forgetful functor
\begin{equation} \label{e:forget untl Rep}
\Rep(\cG)_S^{\on{loc,untl}}\mmod\to \Rep(\cG)_S^{\on{loc}}\mmod.
\end{equation} 

Now, the assertion follows from the fact that the functor \eqref{e:forget untl Rep} is fully
faithful, see \secref{sss:local actions params}.

\end{proof}

\begin{cor} \label{c:Hecke crit Loc bis}
The functor $$\Loc^\Op_{G,\crit,\CZ}:\KL(G)^{\Op^{\on{glob}}}_{\crit,\CZ}\to
\Dmod_\crit(\Bun_G)\otimes \Dmod(\CZ)$$ 
is compatible with the action of $(\Rep(\cG)_\Ran)^\star$.
\end{cor}

\sssec{}

Let now $\CZ$ be pseudo-proper, and consider the functor
$$(\on{Id}\otimes \on{C}^\cdot_c(\CZ,-)): \Dmod_\crit(\Bun_G)\otimes \Dmod(\CZ)\to \Dmod_\crit(\Bun_G).$$

This functor is obviously compatible with the action of $(\Rep(\cG)_\Ran)^\star$.

\medskip

Denote
$$\Loc^\Op_{G,\crit,\int_\CZ}:=(\on{Id}\otimes \on{C}^\cdot_c(\CZ,-))\circ \Loc^\Op_{G,\crit,\CZ},\quad 
\KL(G)^{\Op^{\on{glob}}}_{\crit,\CZ}\to \Dmod_\crit(\Bun_G).$$

Hence, from \corref{c:Hecke crit Loc bis} we obtain:

\begin{cor} \label{c:Hecke crit Loc bis bis}
The functor $\Loc^\Op_{G,\crit,\int_\CZ}$
is compatible with the action of $(\Rep(\cG)_\Ran)^\star$ 
on the two sides.
\end{cor}

\sssec{}

Finally, we take $\CZ=\Ran$. Denote the corresponding functor 
$$\Loc^\Op_{G,\crit,\int_\Ran}:\KL(G)^{\Op^{\on{glob}}}_{\crit,\Ran} \to \Dmod_\crit(\Bun_G)$$
by $\Loc^\Op_{G,\crit}$. 

\medskip

We obtain:

\begin{cor} \label{c:Hecke crit Loc ult}
The the functor 
$\Loc^\Op_{G,\crit}$ is compatible with the action of $(\Rep(\cG)_\Ran)^\star$ (with the convolution monoidal structure)
on the two sides.
\end{cor}

This corollary is the ultimate form of the compatibility between the (globalized) critical localization
functor and the (non-derived) Hecke action. 

\sssec{}

A variant of \corref{c:Hecke crit Loc bis} was at the core of the construction of Hecke eigensheaves in \cite{BD1}. In our language,
this construction can be reformulated as follows.

\medskip

Take $\CZ=\on{pt}$, so that $\CZ\to \Ran$ corresponds to $\ul{x}\in \Ran$.

\medskip

Recall that according to \thmref{t:critical FLE}, \propref{p:Ups mon-free Op} and \propref{p:Op duality mon-free Opers}, we can identify
\begin{multline} 
\KL(G)_{\crit,\ul{x}} \overset{\FLE_{G,\crit}}\simeq \IndCoh^*(\Op^\mf_{\cG,\ul{x}}) \overset{\Theta_{\Op^\mf_\cG}}\simeq \\
\simeq \IndCoh^!(\Op^\mf_{\cG,\ul{x}})\overset{\Upsilon_{\Op^\mf_{\cG,\ul{x}}}}\simeq
\QCoh(\Op^\mf_{\cG,\ul{x}})
\end{multline} 
as $\QCoh(\Op^\mf_{\cG,\ul{x}})$-module categories. 

\medskip

Hence, we can identify the category $\KL(G)^{\Op^{\on{glob}}}_{\crit,\CZ}$ 
with $\QCoh(\Op^\mf_\cG(X-\ul{x}))$, as a module category over $\QCoh(\Op^\mf_\cG(X-\ul{x}))$. 
Since $\Op^\mf_\cG(X-\ul{x})$ is locally almost of finite type and formally smooth, by \cite[Theorem 10.1.1]{GaRo1}, the functor
$$\Upsilon_{\Op^\mf_\cG(X-\ul{x})}:\QCoh(\Op^\mf_\cG(X-\ul{x}))\to \IndCoh(\Op^\mf_\cG(X-\ul{x}))$$
is an equivalence. Hence, we can further identify $\KL(G)^{\Op^{\on{glob}}}_{\crit,\CZ}$ with 
$\IndCoh(\Op^\mf_\cG(X-\ul{x}))$ as a $\QCoh(\Op^\mf_\cG(X-\ul{x}))$-module category. 

\medskip

Hence, we can view $\Loc^\Op_{G,\crit,\ul{x}}$ as a functor
$$\IndCoh(\Op^\mf_\cG(X-\ul{x}))\to \Dmod_\crit(\Bun_G)$$
that intertwines the action of $(\Rep(\cG)_\Ran)^\star$ on $\IndCoh(\Op^\mf_\cG(X-\ul{x}))$, given by 
$$(\fr^{\on{glob}})^*\circ \Loc_\cG^{\on{spec}}:(\Rep(\cG)_\Ran)^\star\to \QCoh(\Op^\mf_\cG(X-\ul{x}))$$
and the action of $\QCoh(\Op^\mf_\cG(X-\ul{x}))$ on $\IndCoh(\Op^\mf_\cG(X-\ul{x}))$, 
and the $(\Rep(\cG)_\Ran)^\star$-action on $\Dmod_\crit(\Bun_G)$. 

\medskip

In particular, for a $k$-point $\sigma\in \Op^\mf_\cG(X-\ul{x})$ and the corresponding sky-scraper sheaf
$$k_\sigma\in \IndCoh(\Op^\mf_\cG(X-\ul{x})),$$
the object
$$\Loc^\Op_{G,\crit,\ul{x}}(k_\sigma)\in \Dmod_\crit(\Bun_G)$$
is a Hecke eigensheaf with eigenvalue $\fr^{\on{glob}}(\sigma)\in \LS_\cG$.

\begin{rem}
The case that is actually considered in \cite{BD1} is when $\sigma$ is a regular oper. In our language,
this corresponds replacing $\KL(G)^{\Op^{\on{glob}}}_{\crit,\CZ}$ with
$$\KL(G)^\reg_{\crit,\CZ}\underset{\QCoh(\Op^\reg_\cG)_\CZ}\otimes (\QCoh(\Op_\cG(X))\otimes \Dmod(\CZ)),$$
where 
$$\KL(G)^\reg_{\crit,\CZ}:=\on{Funct}_{\QCoh(\Op^\mf_\cG)_\CZ}(\QCoh(\Op^\reg_\cG)_\CZ,\KL(G)_{\crit,\CZ}).$$
\end{rem}

\ssec{Proof of \thmref{t:Hecke crit Loc}} \label{ss:proof Hecke crit Loc}

\sssec{}

In the proof of \thmref{t:Hecke crit Loc}, for expositional purposes we will assume that $\CZ=\on{pt}$, so
its map to $\Ran$ corresponds to $\ul{x}\in \Ran$. We will denote the corresponding 
space $\CZ^{\subseteq}$ by $\Ran_{\ul{x}}$. 

\sssec{}

Before we launch the proof, let us explain the idea that lies behind it. Let $\CA$ be a factorization algebra,
and let $\CA^{\on{ch}}$ be the corresponding chiral algebra, see \secref{sss:chiral algebras}. Let $\fz$ be 
a commutative factorization algebra, such that $\fz^{\on{ch}}$ maps to the center of $\CA^{\on{ch}}$. 

\medskip

Let $\CM$ be an object of
$$\CA\mod^{\on{fact}}_{\ul{x}}=\CA^{\on{ch}}\mod^{\on{ch}}_{\ul{x}}$$
equipped with a \emph{commutative} action of $\fz^{\on{ch}}$, which is compatible
with the action of $\CA^{\on{ch}}$. 

\medskip

We claim that in this case $\on{C}^{\on{fact}}_\cdot(X,\CA,\CM)_{\ul{x}}$ carries an action of 
$\on{C}^{\on{fact}}_\cdot(X,\fz)$. Let us construct the action morphism. 

\medskip

We can interpret the given structure on $\CM$ as a map of modules 
\begin{equation} \label{e:action of z on M}
\fz_x\otimes \CM\to \CM
\end{equation} 
compatible with a map of the chiral algebras
$$\fz^{\on{ch}}\otimes \CA^{\on{ch}}\to \CA^{\on{ch}}.$$

By the functoriality of factorization homology, we obtain a map
\begin{multline*}
\on{C}^{\on{fact}}_\cdot(X,\fz)\otimes \on{C}^{\on{fact}}_\cdot(X,\CA,\CM)_{\ul{x}}\simeq 
\on{C}^{\on{fact}}_\cdot(X,\fz,\fz_x)\otimes \on{C}^{\on{fact}}_\cdot(X,\CA,\CM)_{\ul{x}}\overset{\sim}\to \\
\to \on{C}^{\on{fact}}_\cdot(X,\fz\otimes \CA,\fz_x\otimes \CM)_{\ul{x}} \to \on{C}^{\on{fact}}_\cdot(X,\CA,\CM)_{\ul{x}},
\end{multline*}
which is the required action map. 

\medskip

The resulting action of $\on{C}^{\on{fact}}_\cdot(X,\fz)$ on $\on{C}^{\on{fact}}_\cdot(X,\CA,\CM)_{\ul{x}}$ has the following
property: 

\medskip

The action of $\fz_x$ on $\on{C}^{\on{fact}}_\cdot(X,\CA,\CM)_{\ul{x}}$ obtained from the homomorphism
$$\fz_x\to \on{C}^{\on{fact}}_\cdot(X,\fz,\fz_x)\simeq \on{C}^{\on{fact}}_\cdot(X,\fz)$$
equals the action obtained from the $\fz_x$-action on $\CM$ by endomorphisms of the chiral
$\CA^{\on{ch}}$-module structure. 

\medskip

It is a souped-up version of this construction that will be used in \secref{sss:R glob acting} below.  See also
Remark \ref{r:indep com}. 

\sssec{} \label{sss:Y for opers}

Write 
$$\Op^\mf_{\cG,\ul{x}} \simeq \underset{R}{``\on{colim}"}\, \Spec(R),$$
where $Y:=\Spec(R)\to \Op^\mf_{\cG,\ul{x}}$ are closed embeddings almost of finite presentation. 

\medskip

We claim:

\begin{lem} \label{l:approx mod}
For any $\IndCoh^!(\Op^\mf_{\cG,\ul{x}})$-module category $\bC$, the naturally defined functor
$$\underset{R}{\on{colim}}\, \on{Funct}_{\IndCoh^!(\Op^\mf_{\cG,\ul{x}})}(\IndCoh^!(Y),\bC)\to \bC$$
is an equivalence.
\end{lem}

This lemma will be proved in \secref{ss:base change mf}. 

\sssec{}

For every $Y=\Spec(R)$ as above, denote
$$\KL(G)_{\crit,\ul{x},Y}:=
\on{Funct}_{\IndCoh^!(\Op^\mf_{\cG,\ul{x}})}(\IndCoh^!(Y),\KL(G)_{\crit,\ul{x}}).$$
We consider it as a $\QCoh(Y)$-linear category via
$$\Upsilon_Y:\QCoh(Y)\to \IndCoh^!(Y).$$

\sssec{} 

Denote
$$Y^{\on{glob}}:=Y\underset{\Op^\mf_{\cG,\ul{x}}}\times \Op^\mf_\cG(X-\ul{x}).$$

Denote by 
\begin{equation} \label{e:ev Y x}
\on{ev}_{\ul{x}}:Y^{\on{glob}}\to Y
\end{equation}
the evaluation map.

\medskip

Denote
$$\KL(G)^{\Op^{\on{glob}}_\cG}_{\crit,\ul{x},Y}:=\KL(G)_{\crit,\ul{x},Y}\underset{\QCoh(Y)}\otimes 
\QCoh(Y^{\on{glob}}).$$

Denote by
$$\on{Id}\otimes \on{ev}^*_{\ul{x}}:\KL(G)_{\crit,\ul{x},Y}\to \KL(G)^{\Op^{\on{glob}}_\cG}_{\crit,\ul{x},Y}$$
the corresponding functor.

\sssec{} \label{sss:Loc Op Y}

By \lemref{l:approx mod}, in order to construct the 
functor $\Loc^\Op_{G,\crit,\ul{x}}$, it suffices to construct a compatible family of functors 
$$\Loc^\Op_{G,\crit,\ul{x},Y}:\KL(G)^{\Op^{\on{glob}}_\cG}_{\crit,\ul{x},Y} \to \Dmod_\crit(\Bun_G),$$
such that:

\medskip

\noindent{(a)} The functor
\begin{equation} \label{e:Loc Y}
\KL(G)_{\crit,\ul{x},Y}\to \KL(G)_{\crit,\ul{x}} \overset{\Loc_{G,\kappa,\ul{x}}}\longrightarrow \Dmod_\crit(\Bun_G)
\end{equation} 
factors as
\begin{equation} \label{e:factor Loc via glob Y}
\KL(G)_{\crit,\ul{x},Y} \overset{\on{Id}\otimes \on{ev}^*_{\ul{x}}}\longrightarrow  \KL(G)^{\Op^{\on{glob}}_\cG}_{\crit,\ul{x},Y} 
\overset{\Loc^\Op_{G,\crit,\ul{x},Y}}\longrightarrow \Dmod_\crit(\Bun_G);
\end{equation}

\medskip

\noindent{(b)} The functor $\Loc^\Op_{G,\crit,\ul{x},Y}$ is $\Rep(\cG)_{\ul{x}}$-linear.

\sssec{} \label{sss:R glob acting}

Let $R^{\on{glob}}$ denote the algebra of functions on the (affine) scheme $Y^{\on{glob}}$. 
The closed embedding \eqref{e:ev Y x} gives rise to a homomorphism 
\begin{equation} \label{e:R to R glob}
R\to R^{\on{glob}}.
\end{equation} 

\medskip

The datum of $\Loc^\Op_{G,\crit,\ul{x},Y}$ together with the factorization in point (a) above 
is equivalent to the datum of action of $R^{\on{glob}}$ on the functor 
\eqref{e:Loc Y}, such that the action of $R$ obtained by precomposing with \eqref{e:R to R glob} is identified
with the action coming from the action of $R$ on the identity endofunctor of $\KL(G)_{\crit,\ul{x},Y}$.

\sssec{} \label{sss:mer alg acting}

By unitality, we can rewrite \eqref{e:Loc Y} as
\begin{multline} \label{e:Loc Y 1}
\KL(G)_{\crit,\ul{x},Y}\to \KL(G)_{\crit,\ul{x}} \overset{\on{ins.vac}_{\ul{x}}}\longrightarrow 
\KL(G)_{\crit,\Ran_{\ul{x}}} \overset{(\Loc_{G,\kappa})_{\Ran_{\ul{x}}}}\longrightarrow \\
\to \Dmod_\crit(\Bun_G)\otimes \Dmod(\Ran_{\ul{x}}) \overset{\on{Id}\otimes \on{C}^\cdot_c(\Ran_{\ul{x}},-)}
\longrightarrow \Dmod_\crit(\Bun_G).
\end{multline}

Let $\CO_{\Op_\cG^{\reg,Y}}$ be the object in $\on{ComAlg}(\CO_{\Op_\cG^\reg}\mod^{\on{fact}}_x)$ from \secref{sss:modified vacuum}. 
Let $(\CO_{\Op_\cG^{\reg,Y}})_{\Ran_{\ul{x}}}$ be the underlying object in $\on{ComAlg}(\Dmod(\Ran_{\ul{x}}))$.

\medskip

Let us view $\KL(G)_{\crit,\Ran_{\ul{x}}}$ as a category tensored over $\Dmod(\Ran_{\ul{x}})$. In particular, it makes sense to talk
about algebras in $\Dmod(\Ran_{\ul{x}})$ acting on objects in $\KL(G)_{\crit,\Ran_{\ul{x}}}$ or on functors with values in 
$\KL(G)_{\crit,\Ran_{\ul{x}}}$.

\medskip
 
We will show that $(\CO_{\Op_\cG^{\reg,Y}})_{\Ran_{\ul{x}}}$ acts on the functor 
\begin{equation} \label{e:Loc Y 1.5}
\KL(G)_{\crit,\ul{x},Y}\to \KL(G)_{\crit,\ul{x}} \overset{\on{ins.vac}_{\ul{x}}}\longrightarrow \KL(G)_{\crit,\Ran_{\ul{x}}}.
\end{equation}

Moreover, when we evaluate the natural transformation 
$$(\on{diag}_{\ul{x}})_!\simeq (\on{diag}_{\ul{x}})_!\circ (\on{diag}_{\ul{x}})^!\circ \on{ins.vac}_{\ul{x}}\to \on{ins.vac}_{\ul{x}}$$
on $\KL(G)_{\crit,\ul{x},Y}$, the action of $(\CO_{\Op_\cG^{\reg,Y}})_{\Ran_{\ul{x}}}$ on the right-hand side will be compatible with the
action of
\begin{equation} \label{e:R at x}
R\simeq (\on{diag}_{\ul{x}})^!((\CO_{\Op_\cG^{\reg,Y}})_{\Ran_{\ul{x}}})
\end{equation}
on the left-hand side. 

\sssec{} \label{sss:constr action of glob 1}

Assume for a moment the existence of such an action of $(\CO_{\Op_\cG^{\reg,Y}})_{\Ran_{\ul{x}}}$ on \eqref{e:Loc Y 1.5}, and let us produce 
from it an action of $R^{\on{glob}}$ on \eqref{e:Loc Y}. 

\medskip

First, by functoriality, we obtain an action of $(\CO_{\Op_\cG^{\reg,Y}})_{\Ran_{\ul{x}}}$ on the functor 
\begin{multline} \label{e:Loc Y 1.75}
\KL(G)_{\crit,\ul{x},Y}\to \KL(G)_{\crit,\ul{x}} \overset{\on{ins.vac}_{\ul{x}}}\longrightarrow 
\KL(G)_{\crit,\Ran_{\ul{x}}} \overset{(\Loc_{G,\kappa})_{\Ran_{\ul{x}}}}\longrightarrow \\
\to \Dmod_\crit(\Bun_G)\otimes \Dmod(\Ran_{\ul{x}}).
\end{multline}

\sssec{} \label{sss:constr action of glob 2}

Note that $(\CO_{\Op_\cG^{\reg,Y}})_{\Ran_{\ul{x}}}\in \Dmod(\Ran_{\ul{x}})$ belongs to the essential image of the restriction functor
$$\sft^!:\Dmod(\Ran^{\on{untl}}_{\ul{x}})\to \Dmod(\Ran_{\ul{x}}).$$
In particular, it belongs to the subcategory 
$$\Dmod(\Ran_{\ul{x}})^{\on{almost-untl}}\subset \Dmod(\Ran_{\ul{x}})$$
(see \secref{sss:almost unital}).

\medskip

In particular, $\on{C}^\cdot_c(\Ran_{\ul{x}},(\CO_{\Op_\cG^{\reg,Y}})_{\Ran_{\ul{x}}})$ acquires a structure of (commutative) algebra. 

\sssec{} \label{sss:constr action of glob 3}

The essential image of 
$$\KL(G)_{\crit,\ul{x}} \overset{\on{ins.vac}_{\ul{x}}}\longrightarrow \KL(G)_{\crit,\Ran_{\ul{x}}}$$
belongs to the essential image of 
$$\sft^!:\KL(G)_{\crit,\Ran^{\on{untl}}_{\ul{x}}}\to \KL(G)_{\crit,\Ran_{\ul{x}}}.$$

Hence, the essential image of \eqref{e:Loc Y 1.75} belongs to the essential image of
$$(\on{Id}\otimes \sft^!): \Dmod_\crit(\Bun_G)\otimes \Dmod(\Ran^{\on{untl}}_{\ul{x}})\to  \Dmod_\crit(\Bun_G)\otimes \Dmod(\Ran_{\ul{x}}).$$
In particular, it belongs to the subcategory 
$$\Dmod_\crit(\Bun_G)\otimes \Dmod(\Ran_{\ul{x}})^{\on{almost-untl}}\subset \Dmod_\crit(\Bun_G)\otimes \Dmod(\Ran_{\ul{x}}).$$

\sssec{} \label{sss:constr action of glob 4}

Hence, by \secref{sss:almost unital}, the action of $(\CO_{\Op_\cG^{\reg,Y}})_{\Ran_{\ul{x}}}$ on \eqref{e:Loc Y 1.75} gives rise to an action 
of the commutative algebra
$$\on{C}^\cdot_c(\Ran_{\ul{x}},(\CO_{\Op_\cG^{\reg,Y}})_{\Ran_{\ul{x}}})$$
on the functor \eqref{e:Loc Y 1}.

\sssec{} \label{sss:constr action of glob 5}

Finally, by Lemmas \ref{l:fact hom R} and \ref{l:mixed with poles}, we identify
$$\on{C}^\cdot_c(\Ran_{\ul{x}},(\CO_{\Op_\cG^{\reg,Y}})_{\Ran_{\ul{x}}})\simeq \on{C}^{\on{fact}}_\cdot(X,\CO_{\Op^\reg},R)\simeq R^{\on{glob}}$$
as (commutative) algebras.

\medskip

This gives the sought-for action of $R^{\on{glob}}$ on \eqref{e:Loc Y}. The compatbility with the $R$-action on 
the identity endofunctor of $\KL(G)_{\crit,\ul{x},Y}$ follows from the fact that the (iso)morphism
\begin{multline*}
\Loc_{G,\kappa,\ul{x}}\simeq \Loc_{G,\kappa,\Ran_{\ul{x}}}\circ (\on{diag}_{\ul{x}})_! \simeq \\
\simeq \Loc_{G,\kappa,\Ran_{\ul{x}}}\circ (\on{diag}_{\ul{x}})_!\circ (\on{diag}_{\ul{x}})^!\circ \on{ins.vac}_{\ul{x}}\to
\Loc_{G,\kappa,\Ran_{\ul{x}}}\circ  \on{ins.vac}_{\ul{x}}
\end{multline*}
intertwines the $R$-action on the left-hand side and the $\on{C}^{\on{fact}}_\cdot(X,\CO_{\Op^\reg},R)$-action on the right-hand side.

\begin{rem} \label{r:indep com}
A simplified version of the above argument proves the following statement: 

\medskip

Let $\CA$ be a commutative factorization algebra. Consider the factorization category $\CA\mod^{\on{com}}$ 
as a unital sheaf of categories on $\Ran$. Consider the local-to-global functor
$$\CA\mod^{\on{com}}\to \ul\Dmod(\Ran^{\on{untl}})$$
given by
$$(\CZ\to \Ran) \rightsquigarrow \on{C}^{\on{fact}}_\cdot(X,\CA,-)_\CZ.$$

Then:

\smallskip

\noindent(a) The above functor is acted on by $\on{C}^{\on{fact}}_\cdot(X,\CA)$; in particular, 
upgrades to a (strictly unital) local-to-global functor 
\begin{equation} \label{e:fact hom of com univ}
\CA\mod^{\on{com}}\to \on{C}^{\on{fact}}_\cdot(X,\CA)\mod\otimes \ul\Dmod(\Ran^{\on{untl}}).
\end{equation} 

\smallskip

\noindent(b) The functor \eqref{e:fact hom of com univ} is \emph{universal} among strictly unital local-to-global functors.

\medskip

In the language of \secref{ss:indep}, point (b) can reformulated as saying that \eqref{e:fact hom of com univ}
induces an equivalence
$$\CA\mod^{\on{com}}_{\Ran,\on{indep}}\to \on{C}^{\on{fact}}_\cdot(X,\CA)\mod.$$

\end{rem} 

\ssec{Construction of the algebra action} \label{ss:mer alg acting}

In this subsection we will construct the sought-for action of $(\CO_{\Op_\cG^{\reg,Y}})_{\Ran_{\ul{x}}}$ on the functor
\eqref{e:Loc Y 1.5}. 

\sssec{}

By \thmref{t:ins vac reg}, the functor \eqref{e:Loc Y 1.5}
factors as 
\begin{multline}  \label{e:Loc Y 2}
\KL(G)_{\crit,\ul{x},Y} \to  \KL(G)_{\crit,\ul{x}} \overset{\on{Id}\otimes (\on{pr}^\Op_{\on{small}})^*}\longrightarrow \\
\to \KL(G)_{\crit,\ul{x}} \underset{\QCoh(\Op^\mf_{\cG,\ul{x}})}\otimes \QCoh(\Op^{\mf\rightsquigarrow\reg}_{\cG,\Ran_{\ul{x}}}) 
\overset{\on{ins.vac}^{\mf\rightsquigarrow\reg}_{\ul{x}}}\longrightarrow \KL(G)_{\crit,\Ran_{\ul{x}}},
\end{multline} 
while the functor
$$\KL(G)_{\crit,\ul{x},Y} \to  \KL(G)_{\crit,\ul{x}} \overset{\on{Id}\otimes (\on{pr}^\Op_{\on{small}})^*}\longrightarrow 
\KL(G)_{\crit,\ul{x}} \underset{\QCoh(\Op^\mf_{\cG,\ul{x}})}\otimes \QCoh(\Op^{\mf\rightsquigarrow\reg}_{\cG,\Ran_{\ul{x}}})$$
which appears in \eqref{e:Loc Y 2}, factors naturally as
\begin{multline*}
\KL(G)_{\crit,\ul{x},Y}  \to 
\KL(G)_{\crit,\ul{x},Y} \underset{\QCoh(\Op^\mf_{\cG,\ul{x}})}\otimes \QCoh(\Op^{\mf\rightsquigarrow\reg}_{\cG,\Ran_{\ul{x}}})\to \\
\to \KL(G)_{\crit,\ul{x}}  \underset{\QCoh(\Op^\mf_{\cG,\ul{x}})}\otimes \QCoh(\Op^{\mf\rightsquigarrow\reg}_{\cG,\Ran_{\ul{x}}}).
%
\end{multline*}

\sssec{}

We rewrite $\KL(G)_{\crit,\ul{x},Y} \underset{\QCoh(\Op^\mf_{\cG,\ul{x}})}\otimes \QCoh(\Op^{\mf\rightsquigarrow\reg}_{\cG,\Ran_{\ul{x}}})$
tautologically as
$$\KL(G)_{\crit,\ul{x},Y} \underset{\QCoh(Y)}\otimes 
\Bigl(\QCoh(Y)\underset{\QCoh(\Op^\mf_{\cG,\ul{x}})}\otimes \QCoh(\Op^{\mf\rightsquigarrow\reg}_{\cG,\Ran_{\ul{x}}})\Bigr).$$

\medskip

We now claim:

\begin{lem} \label{l:base change poles at x}
The naturally defined functor
$$
\QCoh(Y)\underset{\QCoh(\Op^\mf_{\cG,\ul{x}})}\otimes \QCoh(\Op^{\mf\rightsquigarrow\reg}_{\cG,\Ran_{\ul{x}}})\to
\QCoh\Bigl(Y\underset{\Op^\mf_{\cG,\ul{x}}}\times \Op^{\mf\rightsquigarrow\reg}_{\cG,\Ran_{\ul{x}}}\Bigr)$$
is an equivalence. 
\end{lem}

The lemma will be proved in \secref{ss:base change mf}.

\sssec{}

Applying \lemref{l:base change poles at x}, we obtain that the functor \eqref{e:Loc Y 1.5} can be factored as
\begin{multline}  \label{e:Loc Y 3}
\KL(G)_{\crit,\ul{x},Y}  
\to \KL(G)_{\crit,\ul{x},Y} \underset{\QCoh(Y)}\otimes 
\QCoh\Bigl(Y\underset{\Op^\mf_{\cG,\ul{x}}}\times \Op^{\mf\rightsquigarrow\reg}_{\cG,\Ran_{\ul{x}}}\Bigr) \to \\
\to \KL(G)_{\crit,\ul{x}} \underset{\QCoh(\Op^\mf_{\cG,\ul{x}})}\otimes 
\QCoh(\Op^{\mf\rightsquigarrow\reg}_{\cG,\Ran_{\ul{x}}})\overset{\on{ins.vac}^{\mf\rightsquigarrow\reg}_{\ul{x}}}\longrightarrow 
\KL(G)_{\crit,\Ran_{\ul{x}}}. 
\end{multline}

Hence, it is enough to construct an action of $(\CO_{\Op_\cG^{\reg,Y}})_{\Ran_{\ul{x}}}$ on the composition
\begin{multline}  \label{e:Loc Y 4}
\KL(G)_{\crit,\ul{x},Y} \underset{\QCoh(Y)}\otimes 
\QCoh\Bigl(Y\underset{\Op^\mf_{\cG,\ul{x}}}\times \Op^{\mf\rightsquigarrow\reg}_{\cG,\Ran_{\ul{x}}}\Bigr) \to \\
\to \KL(G)_{\crit,\ul{x}} \underset{\QCoh(\Op^\mf_{\cG,\ul{x}})}\otimes 
\QCoh(\Op^{\mf\rightsquigarrow\reg}_{\cG,\Ran_{\ul{x}}})\overset{\on{ins.vac}^{\mf\rightsquigarrow\reg}_{\ul{x}}}\longrightarrow 
\KL(G)_{\crit,\Ran_{\ul{x}}}. 
\end{multline}

\sssec{}

Consider the category 
\begin{equation} \label{e:KM IndCoh on opers with poles}
\KL(G)_{\crit,\ul{x},Y} \underset{\QCoh(Y)}\otimes 
\QCoh\Bigl(Y\underset{\Op^\mf_{\cG,\ul{x}}}\times \Op^{\mf\rightsquigarrow\reg}_{\cG,\Ran_{\ul{x}}}\Bigr) 
\end{equation} 
as tensored over 
\begin{equation} \label{e:IndCoh on opers with poles}
\QCoh\Bigl(Y\underset{\Op^\mf_{\cG,\ul{x}}}\times \Op^{\mf\rightsquigarrow\reg}_{\cG,\Ran_{\ul{x}}}\Bigr) 
\end{equation} 
and hence over $\Dmod(\Ran_{\ul{x}})$. 

\medskip

We note that, according to \lemref{l:mixed with poles},
$$Y\underset{\Op^\mf_{\cG,\ul{x}}}\times \Op^{\mf\rightsquigarrow\reg}_{\cG,\Ran_{\ul{x}}}\simeq
\Spec_{\Ran_{\ul{x}}}((\CO_{\Op_\cG^{\reg,Y}})_{\Ran_{\ul{x}}}).$$

Hence, $(\CO_{\Op_\cG^{\reg,Y}})_{\Ran_{\ul{x}}}$ maps (isomorphically) to endomorphisms of the monoidal unit 
$$\CO_{Y\underset{\Op^\mf_{\cG,\ul{x}}}\times \Op^{\mf\rightsquigarrow\reg}_{\cG,\Ran_{\ul{x}}}}\in 
\QCoh\Bigl(Y\underset{\Op^\mf_{\cG,\ul{x}}}\times \Op^{\mf\rightsquigarrow\reg}_{\cG,\Ran_{\ul{x}}}\Bigr).$$

Hence, it acts by endomorphisms
of the identity functor on \eqref{e:KM IndCoh on opers with poles}. 

\sssec{} \label{sss:mer alg acting end}

The functors in the composition \eqref{e:Loc Y 4} are $\Dmod(\Ran_{\ul{x}})$-linear. This produces the sought-for
action of $(\CO_{\Op_\cG^{\reg,Y}})_{\Ran_{\ul{x}}}$ on \eqref{e:Loc Y 4}. 

\medskip

The compatibility of this action with \eqref{e:R at x} follows from the construction.

\ssec{Verification of the Hecke property: reduction to a local statement}

The goal of this and the next subsections is to verify property (b) from \secref{sss:Loc Op Y}. 

\medskip

For expositional reasons, we will fix an object $V\in \Rep(\cG)_{\ul{x}}$ and show that the functor $\Loc^\Op_{G,\crit,\ul{x},Y}$
intertwines the actions of $V$ on the two sides.

\medskip

We will reduce the local-to-global assertion we are after to a purely local one, namely, \eqref{e:Hecke property Y 2}. 

\sssec{}

Let $\CV$ denote the object of $\QCoh(\Op^\mf_{\cG,\ul{x}})$ equal to $(\fr^{\on{reg}})^*(V)$,
where we identify
$$\Rep(\cG)_{\ul{x}}\simeq \QCoh(\LS^\reg_{\cG,\ul{x}}).$$

Let 
$$\CV^{\on{glob}}:=(\fr^{\on{glob}})^*\circ \Loc^{\on{spec}}_{\cG,\ul{x}}(V)\in \QCoh(\Op_\cG^\mf(X-\ul{x})).$$

Note that we have
$$\CV^{\on{glob}}\simeq \on{ev}_{\ul{x}}^*(\CV),$$
where
$$\on{ev}_{\ul{x}}:\Op_\cG^\mf(X-\ul{x})\to \Op^\mf_{\cG,\ul{x}}.$$

\medskip

Denote by
$$\CV_Y \text{ and } \CV^{\on{glob}}_Y$$
the restrictions of $\CV$ and $\CV^{\on{glob}}$ to
$$Y\hookrightarrow \Op^\mf_{\cG,\ul{x}} \text{ and }
Y^{\on{glob}}\hookrightarrow \Op_\cG^\mf(X-\ul{x}),$$
respectively. 

\medskip

By a slight abuse of notation, we will denote by the same symbols $\CV_Y$ and $\CV^{\on{glob}}_Y$
global sections of the corresponding vector bundles, viewed as modules over $R$ and $R^{\on{glob}}$, respectively. 

\sssec{}

Let $\CM$ be an object of $\KL(G)_{\kappa,\ul{x},Y}$, and consider the object
$$\Loc_{G,\kappa,\ul{x}}(\CM)\in \Dmod_\crit(\Bun_G).$$

The construction in Sects. \ref{sss:constr action of glob 1}-\ref{sss:constr action of glob 5}
endows $\Loc_{G,\kappa,\ul{x}}(\CM)$ with an action of $R^{\on{glob}}$.

\medskip

Let 
$$\on{H}_V: \Dmod_\crit(\Bun_G)\to  \Dmod_\crit(\Bun_G)$$
be the Hecke endofunctor corresponding to $V$. 

\medskip

On the one hand, by functoriality, the object 
$$\on{H}_V(\Loc_{G,\kappa,\ul{x}}(\CM))\in \Dmod_\crit(\Bun_G)$$ 
acquires an action of $R^{\on{glob}}$.

\medskip

On the other hand, we can consider
$$\CV^{\on{glob}}_Y\underset{R^{\on{glob}}}\otimes \Loc_{G,\kappa,\ul{x}}(\CM) \in R^{\on{glob}}\mod(\Dmod_\crit(\Bun_G)).$$

\sssec{}

The statement of (b) in \secref{sss:Loc Op Y} is that we have a canonical isomorphism 
\begin{equation} \label{e:Hecke property Y}
\on{H}_V(\Loc_{G,\kappa,\ul{x}}(\CM))\simeq \CV^{\on{glob}}_Y\underset{R^{\on{glob}}}\otimes \Loc_{G,\kappa,\ul{x}}(\CM)
\end{equation} 
as objects of $R^{\on{glob}}\mod(\Dmod_\crit(\Bun_G))$. 

\medskip

Thus, our goal is to establish \eqref{e:Hecke property Y}.

\sssec{}

First, we rewrite the right-hand side in \eqref{e:Hecke property Y}. Namely,
$$\CV^{\on{glob}}_Y\underset{R^{\on{glob}}}\otimes \Loc_{G,\kappa,\ul{x}}(\CM)
\simeq \Loc_{G,\kappa,\ul{x}}(\CV_Y\underset{R}\otimes \CM),$$
where:

\begin{itemize}

\item We regard $\CV_Y\underset{R}\otimes \CM$ as an object $\KL(G)_{\kappa,\ul{x},Y}$;

\item $\Loc_{G,\kappa,\ul{x}}(-)$ acquires an action of $R^{\on{glob}}$ via the construction in 
Sects. \ref{sss:constr action of glob 1}-\ref{sss:constr action of glob 5}.

\end{itemize}

\sssec{} 

We will now rewrite the left-hand side in \eqref{e:Hecke property Y}.

\medskip

Consider the category $\KL(G)_{\kappa,\Ran_{\ul{x}}}$. It carries an action of $\Sph_{G,\Ran_{\ul{x}}}$.
We will denote the action functor by
$$\CF,\CM'\mapsto \CF\cdot \CM'.$$

\medskip

The unital structure on $\Sph_G$ gives rise to a (monoidal) functor
$$\on{ins.unit}_{\ul{x}}:\Sph_{G,\ul{x}}\to \Sph^{\sotimes}_{G,\Ran_{\ul{x}}}$$
(see \secref{sss:pointwise monoidal on Ran} for the notation). 

\medskip

In particular, we obtain a monoidal action of $\Sph_{G,\ul{x}}$ on $\KL(G)_{\kappa,\Ran_{\ul{x}}}$.

\sssec{} \label{sss:Hecke property Y LHS 1}

Recall that the object
$$\on{ins.vac}_{\ul{x}}(\CM)\in \KL(G)_{\kappa,\Ran_{\ul{x}}}$$
carries an action of the algebra object $(\CO_{\Op_\cG^{\reg,Y}})_{\Ran_{\ul{x}}}\in \Dmod(\Ran_{\ul{x}})$.  By functoriality, we obtain that
\begin{equation} \label{e:ins Hecke}
\on{ins.unit}_{\ul{x}}(\on{Sat}^{\on{nv}}_G(V))\cdot \on{ins.vac}_{\ul{x}}(\CM)\in \KL(G)_{\kappa,\Ran_{\ul{x}}}
\end{equation} 
also carries an action of $(\CO_{\Op_\cG^{\reg,Y}})_{\Ran_{\ul{x}}}$. 

\medskip

Further, the object 
\begin{equation} \label{e:ins Hecke 1}
\Loc_{G,\crit,\Ran_{\ul{x}}}\Bigl(\on{ins.unit}_{\ul{x}}(\on{Sat}^{\on{nv}}_G(V))\cdot \on{ins.vac}_{\ul{x}}(\CM)\Bigr)\in 
\Dmod_\crit(\Bun_G)\otimes \Dmod(\Ran_{\ul{x}})
\end{equation} 
also carries an action of $(\CO_{\Op_\cG^{\reg,Y}})_{\Ran_{\ul{x}}}$.

\sssec{} \label{sss:Hecke property Y LHS 2}

The object \eqref{e:ins Hecke} belongs to the essential image of the restriction functor
$$\sft^!:\KL(G)_{\kappa,\Ran^{\on{untl}}_{\ul{x}}}\to \KL(G)_{\kappa,\Ran_{\ul{x}}}.$$

Hence, the object \eqref{e:ins Hecke 1} belongs to the essential image of
$$(\on{Id}\otimes \sft^!): \Dmod_\crit(\Bun_G)\otimes \Dmod(\Ran^{\on{untl}}_{\ul{x}})\to \Dmod_\crit(\Bun_G)\otimes \Dmod(\Ran_{\ul{x}}).$$
In particular, it belongs to
$$\Dmod_\crit(\Bun_G)\otimes \Dmod(\Ran_{\ul{x}})^{\on{almost-untl}} \subset \Dmod_\crit(\Bun_G)\otimes \Dmod(\Ran_{\ul{x}}).$$

\medskip

Hence, by \secref{sss:almost unital}, we obtain that the object
\begin{equation} \label{e:Hecke property Y 1}
(\on{Id}\otimes \on{C}^\cdot_c(\Ran_{\ul{x}},-))\circ 
\Loc_{G,\crit,\Ran_{\ul{x}}}\Bigl(\on{ins.unit}_{\ul{x}}(\on{Sat}^{\on{nv}}_G(V))\cdot \on{ins.vac}_{\ul{x}}(\CM)\Bigr)
\end{equation}
acquires an action of 
$$\on{C}^\cdot_c(\Ran_{\ul{x}},(\CO_{\Op_\cG^{\reg,Y}})_{\Ran_{\ul{x}}})\simeq \on{C}^{\on{fact}}_\cdot(X,\CO_{\Op^\reg},R)\simeq R^{\on{glob}}.$$

\sssec{} \label{sss:Hecke property Y LHS 4}

Recall now that  the functor
$$\Loc_{G,\crit,\Ran_{\ul{x}}}:\KL(G)_{\kappa,\Ran_{\ul{x}}}\to \Dmod_\crit(\Bun_G)\otimes \Dmod(\Ran_{\ul{x}})$$
is $\Sph^{\sotimes}_{G,\Ran_{\ul{x}}}$-linear. 

\medskip

Note also that the functor 
$$\on{Id}\otimes \on{C}^\cdot_c(\Ran_{\ul{x}},-):\Dmod_\crit(\Bun_G)\otimes \Dmod(\Ran_{\ul{x}})\to \Dmod_\crit(\Bun_G)$$
is $\Sph_{G,\ul{x}}$-linear, where $\Sph_{G,\ul{x}}$ acts on the left-hand side via $\on{ins.unit}_{\ul{x}}$.

\medskip

Combining, we obtain that the functor
$$(\on{Id}\otimes \on{C}^\cdot_c(\Ran_{\ul{x}},-))\circ \Loc_{G,\crit,\Ran_{\ul{x}}}, \quad
\KL(G)_{\kappa,\Ran_{\ul{x}}}\to \Dmod_\crit(\Bun_G)$$
is $\Sph_{G,\ul{x}}$-linear.

\sssec{}

Combining Sects. \ref{sss:Hecke property Y LHS 1}-\ref{sss:Hecke property Y LHS 2} with \secref{sss:Hecke property Y LHS 4}, 
we obtain that 
$$\on{H}_V(\Loc_{G,\kappa,\ul{x}}(\CM))\in R^{\on{glob}}\mod(\Dmod_\crit(\Bun_G))$$
identifies with the object \eqref{e:Hecke property Y 1}, with the $R^{\on{glob}}$-action specified in 
\ref{sss:Hecke property Y LHS 1}-\ref{sss:Hecke property Y LHS 2}.  

\sssec{}

Hence, we obtain that  in order to prove \eqref{e:Hecke property Y}, it suffices to establish an isomorphism
\begin{equation} \label{e:Hecke property Y 2}
\on{ins.unit}_{\ul{x}}(\on{Sat}^{\on{nv}}_G(V))\cdot \on{ins.vac}_{\ul{x}}(\CM) \simeq
\on{ins.vac}_{\ul{x}}(\CV_Y\underset{R}\otimes \CM)
\end{equation}
as $(\CO_{\Op_\cG^{\reg,Y}})_{\Ran_{\ul{x}}}$-modules in $\KL(G)_{\kappa,\Ran_{\ul{x}}}$, where:

\medskip

\begin{itemize}

\item We regard $\CV_Y\underset{R}\otimes \CM$ as an object $\KL(G)_{\kappa,\ul{x},Y}$;

\item $\on{ins.vac}_{\ul{x}}(-)$ on each side acquires an action of $(\CO_{\Op_\cG^{\reg,Y}})_{\Ran_{\ul{x}}}$ via the construction in 
Sect. \ref{ss:mer alg acting}.

\end{itemize}

\begin{rem}
As was mentioned in the preamble to this subsection, we fixed objects $\CM\in \KL(G)_{\kappa,\ul{x},Y}$
and $V\in \Rep(\cG)_{\ul{x}}$ for expositional reasons. The actual assertion behind \eqref{e:Hecke property Y 2},
and one that we actually prove in \secref{ss:verify local Hecke}, is that that $\on{ins.vac}_{\ul{x}}$, viewed
as a functor
$$\KL(G)_{\kappa,\ul{x},Y}\to ((\CO_{\Op_\cG^{\reg,Y}})_{\Ran_{\ul{x}}})\mod(\KL(G)_{\kappa,\Ran_{\ul{x}}}),$$
is $\Rep(\cG)_{\ul{x}}$-linear. 
\end{rem}

\ssec{Verification of the Hecke property at the local level} \label{ss:verify local Hecke}

In this subsection we will construct the identification \eqref{e:Hecke property Y 2} and thereby complete the
verification of point (b) in \secref{sss:Loc Op Y}.  

\sssec{}

Let us apply \lemref{l:Rep G act on KM}. It implies that we can rewrite the left-hand side in
\eqref{e:Hecke property Y 2} as
\begin{equation} \label{e:Hecke property Y 3}
\fr^*(\on{ins.unit}_{\ul{x}}(V))\otimes  \on{ins.vac}_{\ul{x}}(\CM),
\end{equation}
where:

\begin{itemize}

\item $\on{ins.unit}_{\ul{x}}:\Rep(\cG)_{\ul{x}}\to \Rep(\cG)_{\Ran_{\ul{x}}}$ is the unital structure on $\Rep(\cG)$;

\smallskip

\item We identify $\Rep(\cG)_{\Ran_{\ul{x}}}\simeq \QCoh(\LS^\reg_{\cG,\Ran_{\ul{x}}})$;

\smallskip

\item $\fr$ denotes the map $\Op^\mf_{\cG,\Ran_{\ul{x}}}\to \LS^\reg_{\cG,\Ran_{\ul{x}}}$;

\smallskip

\item $\otimes$ refers to the action of $\QCoh(\Op^\mf_{\cG,\Ran_{\ul{x}}})$ on 
$\KL(G)_{\kappa,\Ran_{\ul{x}}}$.

\end{itemize}

\medskip

The action of $(\CO_{\Op_\cG^{\reg,Y}})_{\Ran_{\ul{x}}}$ on \eqref{e:Hecke property Y 3} is obtained by functoriality
from the $(\CO_{\Op_\cG^{\reg,Y}})_{\Ran_{\ul{x}}}$-action on $\on{ins.vac}_{\ul{x}}(\CM)$.

\sssec{}

We now note that the linearity with respect to
$$(\on{pr}_{\on{small}})_*\circ (\on{pr}_{\on{big}})^*(\QCoh(\Op^\mf_\cG))$$
in \thmref{t:ins vac reg}
implies that the functor 
\begin{multline}  \label{e:Loc Y 4 again}
\KL(G)_{\crit,\ul{x},Y} \underset{\QCoh(Y)}\otimes 
\QCoh\Bigl(Y\underset{\Op^\mf_{\cG,\ul{x}}}\times \Op^{\mf\rightsquigarrow\reg}_{\cG,\Ran_{\ul{x}}}\Bigr) \to \\
\to \KL(G)_{\crit,\ul{x}} \underset{\QCoh(\Op^\mf_{\cG,\ul{x}})}\otimes 
\QCoh(\Op^{\mf\rightsquigarrow\reg}_{\cG,\Ran_{\ul{x}}})\overset{\on{ins.vac}^{\mf\rightsquigarrow\reg}_{\ul{x}}}\longrightarrow 
\KL(G)_{\crit,\Ran_{\ul{x}}}
\end{multline}
of \eqref{e:Loc Y 4} is linear with respect to $\QCoh(\Op^\mf_{\cG,\Ran_{\ul{x}}})$,
which acts on 
$$\KL(G)_{\crit,\ul{x},Y} \underset{\QCoh(Y)}\otimes 
\QCoh\Bigl(Y\underset{\Op^\mf_{\cG,\ul{x}}}\times \Op^{\mf\rightsquigarrow\reg}_{\cG,\Ran_{\ul{x}}}\Bigr)$$
via the pullback along
$$Y\underset{\Op^\mf_{\cG,\ul{x}}}\times \Op^{\mf\rightsquigarrow\reg}_{\cG,\Ran_{\ul{x}}}\overset{p_2}\to
\Op^{\mf\rightsquigarrow\reg}_{\cG,\Ran_{\ul{x}}} \overset{\on{pr}^\Op_{\on{big},\ul{x}}}\longrightarrow \Op^\mf_{\cG,\Ran_{\ul{x}}},$$
where $\on{pr}^\Op_{\on{big},\ul{x}}$ is the map from \secref{sss:pr Op big}. 

\sssec{}

Hence, we can rewrite \eqref{e:Hecke property Y 3} as the value on $\CM\in \KL(G)_{\kappa,\ul{x},Y}$ of 
the functor
\begin{multline} \label{e:Loc Y 5}
\KL(G)_{\crit,\ul{x},Y}  
\to \KL(G)_{\crit,\ul{x},Y} \underset{\QCoh(Y)}\otimes 
\QCoh\Bigl(Y\underset{\Op^\mf_{\cG,\ul{x}}}\times \Op^{\mf\rightsquigarrow\reg}_{\cG,\Ran_{\ul{x}}}\Bigr)\to \\
\overset{\on{Id}\otimes (-\otimes (\fr\circ \on{pr}^\Op_{\on{big},\ul{x}}\circ p_2)^*(\on{ins.unit}_{\ul{x}}(V)))}\longrightarrow 
 \KL(G)_{\crit,\ul{x},Y} \underset{\QCoh(Y)}\otimes 
 \QCoh\Bigl(Y\underset{\Op^\mf_{\cG,\ul{x}}}\times \Op^{\mf\rightsquigarrow\reg}_{\cG,\Ran_{\ul{x}}}\Bigr) \to \\
\overset{\text{\eqref{e:Loc Y 4}}}\longrightarrow \KL(G)_{\kappa,\Ran_{\ul{x}}}.
\end{multline}

In terms of this identification, the action of $(\CO_{\Op_\cG^{\reg,Y}})_{\Ran_{\ul{x}}}$ on \eqref{e:Hecke property Y 3} is obtained from
the action of $(\CO_{\Op_\cG^{\reg,Y}})_{\Ran_{\ul{x}}}$ on the functor \eqref{e:Loc Y 4 again} (which is the same as \eqref{e:Loc Y 4}) 
from \secref{sss:mer alg acting end}. 

\sssec{}

Thus, we obtain that in order to prove \eqref{e:Hecke property Y 2}, it suffices to establish
an isomorphism between the value on $\CM$ of the functor 
\begin{multline} \label{e:Loc Y 6}
\KL(G)_{\crit,\ul{x},Y}  
\to \KL(G)_{\crit,\ul{x},Y} \underset{\IndCoh^!(Y)}\otimes 
\QCoh\Bigl(Y\underset{\Op^\mf_{\cG,\ul{x}}}\times \Op^{\mf\rightsquigarrow\reg}_{\cG,\Ran_{\ul{x}}}\Bigr) \to \\
\overset{\on{Id}\otimes (-\otimes (\fr\circ \on{pr}^\Op_{\on{big},\ul{x}}\circ p_2)^*(\on{ins.unit}_{\ul{x}}(V)))}\longrightarrow 
\KL(G)_{\crit,\ul{x},Y} \underset{\QCoh(Y)}\otimes 
\QCoh\Bigl(Y\underset{\Op^\mf_{\cG,\ul{x}}}\times \Op^{\mf\rightsquigarrow\reg}_{\cG,\Ran_{\ul{x}}}\Bigr) 
\end{multline} 
and the value on $\CV_Y\underset{R}\otimes \CM$ of the functor
\begin{equation} \label{e:Loc Y 7}
\KL(G)_{\crit,\ul{x},Y}  
\to \KL(G)_{\crit,\ul{x},Y} \underset{\QCoh(Y)}\otimes 
\QCoh\Bigl(Y\underset{\Op^\mf_{\cG,\ul{x}}}\times \Op^{\mf\rightsquigarrow\reg}_{\cG,\Ran_{\ul{x}}}\Bigr).
\end{equation}

\sssec{}

In order to do that, it suffices to establish an isomorphism between the following vector bundles on 
$$Y\underset{\Op^\mf_{\cG,\ul{x}}}\times \Op^{\mf\rightsquigarrow\reg}_{\cG,\Ran_{\ul{x}}}:$$

\medskip

\begin{itemize}

\item The pullback of $\on{ins.unit}_{\ul{x}}(V)\in \Rep(\cG)_{\Ran_{\ul{x}}}\simeq \QCoh(\LS^\reg_{\cG,\Ran_{\ul{x}}})$
along
$$Y\underset{\Op^\mf_{\cG,\ul{x}}}\times \Op^{\mf\rightsquigarrow\reg}_{\cG,\Ran_{\ul{x}}}\overset{p_2}\to
\Op^{\mf\rightsquigarrow\reg}_{\cG,\Ran_{\ul{x}}}\overset{\on{pr}^\Op_{\on{big},\ul{x}}}\longrightarrow
 \Op^\mf_{\cG,\Ran_{\ul{x}}} \overset{\fr}\to \LS^\reg_{\cG,\Ran_{\ul{x}}};$$

\item The pullback of $\CV_Y$ along
$$Y\underset{\Op^\mf_{\cG,\ul{x}}}\times \Op^{\mf\rightsquigarrow\reg}_{\cG,\Ran_{\ul{x}}}\to Y.$$

\end{itemize}

\sssec{}

Note that for $(\ul{x}\subseteq \ul{x}')\in \Ran_{\ul{x}}$, restriction along
$$\cD_{\ul{x}}\hookrightarrow \cD_{\ul{x}'}$$
gives rise to a map 
\begin{equation} \label{e:LS restr to disc}
\LS^\reg_{\cG,\Ran_{\ul{x}}} \to \LS^\reg_{\cG,\ul{x}},
\end{equation}
so that
$$\on{ins.unit}_{\ul{x}}\simeq \text{\eqref{e:LS restr to disc}}^*.$$

\medskip

The required isomorphism of vector bundles follows from the commutative diagram
$$
\CD
Y\underset{\Op^\mf_{\cG,\ul{x}}}\times \Op^{\mf\rightsquigarrow\reg}_{\cG,\Ran_{\ul{x}}} @>{\on{id}\times p_1}>> Y  @>>>  \Op^\mf_{\cG,\ul{x}} \\
@V{p_2}VV  \\
\Op^{\mf\rightsquigarrow\reg}_{\cG,\Ran_{\ul{x}}} & & & & @VV{\fr}V \\
@V{\on{pr}^\Op_{\on{big},\ul{x}}}VV  \\
\Op^\mf_{\cG,\Ran_{\ul{x}}}  @>{\fr}>> \LS^\reg_{\cG,\Ran_{\ul{x}}} @>{\text{\eqref{e:LS restr to disc}}}>> \LS^\reg_{\cG,\ul{x}}. 
\endCD
$$

\ssec{Proofs of Lemmas \ref{l:approx mod} 
and \ref{l:base change poles at x}} \label{ss:base change mf}

\sssec{}

We first prove \lemref{l:approx mod}. In fact, the assertion holds for $\Op_\cG^\mf(\cD^\times_{\ul{x}})$
replaced by an arbitrary ind-placid ind-scheme $Z$.

\medskip

Namely, for $Y_1\overset{\iota_{1,2}}\hookrightarrow Y_2$, the functor
$$\on{Funct}_{\IndCoh^!(Z)}(\IndCoh^!(Y_1),\bC)\to
\on{Funct}_{\IndCoh^!(Z)}(\IndCoh^!(Y_2),\bC)$$
admits a right adjoint, given by precomposition with $(\iota_{1,2})^\IndCoh_*$.

\medskip

Hence, we can rewrite the colimit
$$\underset{Y}{\on{colim}}\, \on{Funct}_{\IndCoh^!(Z)}(\IndCoh^!(Y),\bC)$$
as a limit with respect to the above right adjoint functors. 

\medskip 

The latter limit is the same as
$$\on{Funct}_{\IndCoh^!(Z)}\left(\underset{Y}{\on{colim}}\, \IndCoh^!(Y),\bC\right).$$

We now use the fact that the functor
$$\underset{Y}{\on{colim}}\, \IndCoh^!(Y)\to \IndCoh(Z)$$
is an equivalence.

\qed[\lemref{l:approx mod}]

\sssec{}

The rest of this subsection is devoted to the proof of \lemref{l:base change poles at x}. Consider the Cartesian diagram
$$
\CD
\Op^{\mf\rightsquigarrow\reg}_{\cG,\Ran_{\ul{x}}} @>>> \Op^{\mer\rightsquigarrow\reg}_{\cG,\Ran_{\ul{x}}} \\
@V{\on{pr}^\Op_{\on{small},\ul{x}}}VV @VV{\on{pr}^\Op_{\on{small},\ul{x}}}V \\
\Op^\mf_{\cG,\ul{x}} @>{\iota^\mf}>> \Op^\mer_{\cG,\ul{x}},
\endCD
$$
where $\Op^{\mer\rightsquigarrow\reg}_{\cG,\Ran_{\ul{x}}}$ is as in \secref{sss:Y mer to reg}. 

\sssec{}

We will prove:

\begin{lem} \label{l:base change poles at x better} \hfill

\smallskip

\noindent{\em(a)} The category $\QCoh(\Op^{\mer\rightsquigarrow\reg}_{\cG,\Ran_{\ul{x}}})$
is dualizable as a $\QCoh(\Op^\mer_{\cG,\ul{x}})$-module. 

\smallskip

\noindent{\em(b)} For any affine $Y\to \Op^\mer_{\cG,\ul{x}}$, the functor
$$\QCoh(Y)\underset{\QCoh(\Op^\mer_{\cG,\ul{x}})}\otimes 
\QCoh(\Op^{\mer\rightsquigarrow\reg}_{\cG,\Ran_{\ul{x}}})\to \QCoh(Y\underset{\Op^\mer_{\cG,\ul{x}}}\times
\Op^{\mer\rightsquigarrow\reg}_{\cG,\Ran_{\ul{x}}})$$
is an equivalence.

\end{lem}

It is easy to see that \lemref{l:base change poles at x better} implies \eqref{l:base change poles at x} by passage 
to the limit.

\medskip

Thus, the rest of this subsection is devoted to the proof of \lemref{l:base change poles at x better}. 

\sssec{}

First, a standard limit-colimit procedure reduces the assertion to the case when we replace $\Ran_{\ul{x}}$ by 
$$\CZ:=(X^I_{\ul{x}})_\dr,$$
where:

\begin{itemize}

\item We think of $\ul{x}$ as a point of $X^J$ for some  finite set $J$;

\item $I$ is a finite set with a map $J\to I$;

\item $X^I_{\ul{x}}:=X^I\underset{X^J,\ul{x}}\times \on{pt}$.

\end{itemize}

\medskip

Further, we can assume that $X$ is affine and admits an \'etale map to $\BA^1$. 

\medskip 

Second, we can replace the D-scheme $\Op_\cG$ by $\on{Jets}(\CE)$,
where $\CE$ is the total space of a vector bundle on $X$ (see \secref{sss:Op is a torsor}). 

\sssec{}

Note that for $X$ as above, we have (non-canonical) isomorphisms
$$\fL^{\mer\rightsquigarrow\reg}_{\CZ} \simeq \fL^+(\CE)_{\CZ}\times 
(\fL(\CE)_{\ul{x}}/\fL^+(\CE)_{\ul{x}})$$
and 
$$\fL(\CE)_{\ul{x}}\simeq \fL^+(\CE)_{\ul{x}}\times (\fL(\CE)_{\ul{x}}/\fL^+(\CE)_{\ul{x}}),$$
so that the projection
$$\fL^{\mer\rightsquigarrow\reg}_{X^I} \to \fL(\CE)_{\ul{x}}$$ 
corresponds to the projection
$$\fL^+(\CE)_{\CZ}\to \fL^+(\CE)_{\ul{x}}.$$

We can therefore identify
$$\QCoh(\fL^{\mer\rightsquigarrow\reg}_{\CZ}) \simeq \QCoh(\fL^+(\CE)_{\CZ})\otimes 
\QCoh(\fL(\CE)_{\ul{x}}/\fL^+(\CE)_{\ul{x}})$$
and
$$\QCoh(\fL(\CE)_{\ul{x}})\simeq \QCoh(\fL^+(\CE)_{\ul{x}}) \otimes \QCoh(\fL(\CE)_{\ul{x}}/\fL^+(\CE)_{\ul{x}}).$$

\sssec{}

To prove point (a), it suffices to show that $\QCoh(\fL^+(\CE)_{\CZ})$ is dualizable as a module over
$\QCoh(\fL(\CE)_{\ul{x}})$. However, this is obvious, since both geometric objects are affine schemes.

\sssec{}

To prove point (b), it is sufficient to do so for a cofinal family of $Y$'s. Hence, we can assume that $Y$
is invariant under translations with respect to $\fL^+(\CE)_{\ul{x}}$. Hence, we can identify
$$Y\underset{\fL(\CE)_{\ul{x}}}\times \fL^{\mer\rightsquigarrow\reg}_{\CZ} \simeq
\fL^+(\CE)_{\ul{x}}\times (Y/\fL^+(\CE)_{\ul{x}}).$$

We have
$$\QCoh(Y\underset{\fL(\CE)_{\ul{x}}}\times \fL^{\mer\rightsquigarrow\reg}_{\CZ})\simeq 
\QCoh(\fL^+(\CE)_{\CZ})\otimes \QCoh(Y/\fL^+(\CE)_{\ul{x}}),$$
which makes the assertion of point (b) manifest. 

\qed[\lemref{l:base change poles at x better}]

\section{Proof of \thmref{t:ins vac reg}} \label{s:ins vac reg}

\ssec{Reformulation and strategy}

\sssec{} \label{sss:ins vac reg pt}

For expositional purposes we will let $\CZ=\on{pt}$, so that $\CZ\to \Ran$ corresponds to $\ul{x}\in \Ran$.

\medskip

Hence, our goal is to construct a functor 
$$\QCoh(\Op^{\mf\rightsquigarrow\reg}_{\cG,\Ran_{\ul{x}}})\underset{\QCoh(\Op^\mf_{\cG,\ul{x}})}\otimes \KL(G)_{\crit,\ul{x}}  
\overset{\on{ins.vac}^{\mf\rightsquigarrow\reg}_{\ul{x}}}\longrightarrow \KL(G)_{\kappa,\Ran_{\ul{x}}}$$
such that 
$$\on{ins.vac}_{\ul{x}}:\KL(G)_{\crit,\ul{x}} \to \KL(G)_{\kappa,\Ran_{\ul{x}}}$$
factors as 
\begin{multline}
\KL(G)_{\crit,\ul{x}} \overset{(\on{pr}^\Op_{\on{small},\ul{x}})^*\otimes \on{Id}}\longrightarrow 
\QCoh(\Op^{\mf\rightsquigarrow\reg}_{\cG,\Ran_{\ul{x}}})\underset{\QCoh(\Op^\mf_{\cG,\ul{x}})}\otimes \KL(G)_{\crit,\ul{x}} \to  \\
\overset{\on{ins.vac}^{\mf\rightsquigarrow\reg}_{\ul{x}}}\longrightarrow \KL(G)_{\kappa,\Ran_{\ul{x}}},
\end{multline}
and such that the functor $\on{ins.vac}^{\mf\rightsquigarrow\reg}_{\ul{x}}$ is $\QCoh(\Op^\mf_{\cG,\Ran_{\ul{x}}})$-linear via
$$(\on{pr}^\Op_{\on{big},\ul{x}})^*:\QCoh(\Op^\mf_{\cG,\Ran_{\ul{x}}})\to \QCoh(\Op^{\mf\rightsquigarrow\reg}_{\cG,\Ran_{\ul{x}}}).$$

\sssec{}

Precomposing with the functor
\begin{equation} \label{e:Ups mf reg}
\Upsilon_{\Op^{\mf\rightsquigarrow\reg}_{\cG,\Ran_{\ul{x}}}}:
\QCoh(\Op^{\mf\rightsquigarrow\reg}_{\cG,\Ran_{\ul{x}}})\to \IndCoh^!(\Op^{\mf\rightsquigarrow\reg}_{\cG,\Ran_{\ul{x}}}),
\end{equation}
we obtain that it suffices to carry out the construction in \secref{sss:ins vac reg pt} above for $\QCoh(-)$ replaced\footnote{One can show
(using \lemref{l:base change poles at x better} combined with a parallel statement for $\IndCoh^!$)
that the functor \eqref{e:Ups mf reg} is actually an equivalence.} by $\IndCoh^!(-)$. 

\medskip

I.e., from now on our goal will be to construct a functor 
\begin{equation} \label{e:ins vac reg mf IndCoh}
\IndCoh^!(\Op^{\mf\rightsquigarrow\reg}_{\cG,\Ran_{\ul{x}}})\underset{\IndCoh^!(\Op^\mf_{\cG,\ul{x}})}\otimes \KL(G)_{\crit,\ul{x}}  
\overset{\on{ins.vac}^{\mf\rightsquigarrow\reg}_{\ul{x}}}\longrightarrow \KL(G)_{\kappa,\Ran_{\ul{x}}}
\end{equation}
such that $\on{ins.vac}_{\ul{x}}$ factors as 
\begin{multline}
\KL(G)_{\crit,\ul{x}} \overset{(\on{pr}^\Op_{\on{small},\ul{x}})^!\otimes \on{Id}}\longrightarrow 
\IndCoh^!(\Op^{\mf\rightsquigarrow\reg}_{\cG,\Ran_{\ul{x}}})\underset{\IndCoh^!(\Op^\mf_{\cG,\ul{x}})}\otimes \KL(G)_{\crit,\ul{x}} \to  \\
\overset{\on{ins.vac}^{\mf\rightsquigarrow\reg}_{\ul{x}}}\longrightarrow \KL(G)_{\kappa,\Ran_{\ul{x}}},
\end{multline}
and such that the functor $\on{ins.vac}^{\mf\rightsquigarrow\reg}_{\ul{x}}$ is $\IndCoh^!(\Op^\mf_{\cG,\Ran_{\ul{x}}})$-linear via
$$(\on{pr}^\Op_{\on{big},\ul{x}})^!:\IndCoh^!(\Op^\mf_{\cG,\Ran_{\ul{x}}})\to \IndCoh^!(\Op^{\mf\rightsquigarrow\reg}_{\cG,\Ran_{\ul{x}}}).$$

\sssec{}

Consider the factorization functor
$${\mathcal Vac}(G)_\crit:\QCoh(\Op^\reg_\cG)\simeq \fz\mod^{\on{cl}}\to \KL(G)_\kappa,$$
given by 
$$\CF\mapsto \CF\underset{\fz}\otimes \on{Vac}(G)_\crit.$$ 

\sssec{}

Recall that according to the conventions in \secref{sss:vac fact mod cat}, for a factorization category $\bA$, we denote by
$\bA^{\on{fact}_{\ul{x}}}$ the vacuum object of $\bA\mmod^{\on{fact}}_{\ul{x}}$, i.e., the object whose
underlying category is $\bA_x$.

\medskip

Thus, we can consider
$$(\KL(G)_\kappa)^{\on{fact}_{\ul{x}}}\in \KL(G)_\kappa\mmod^{\on{fact}}_{\ul{x}}.$$

Consider the object 
$$\Res_{{\mathcal Vac}(G)_\crit}((\KL(G)_\kappa)^{\on{fact}_{\ul{x}}})\in \QCoh(\Op^\reg_\cG)\mmod^{\on{fact}}_{\ul{x}}.$$

\sssec{Example}

For $\ul{x}\sqcup \ul{x}''=\ul{x'}\in \Ran_{\ul{x}'}$, the fiber of $\Res_{{\mathcal Vac}(G)_\crit}((\KL(G)_\kappa)^{\on{fact}_{\ul{x}}})$
at $\ul{x}'$ is
$$\KL(G)_{\kappa,\ul{x}}\otimes \QCoh(\Op^\reg_{\cG,\ul{x}''}).$$

\sssec{} \label{sss:indCoh mf reg Ran x}

Consider the assignment
$$(\CZ \to \Ran_{\ul{x}}) \mapsto \IndCoh^!(\Op^{\mf\rightsquigarrow\reg}_{\cG,\Ran_{\ul{x}}}\underset{\Ran_{\ul{x}}}\times \CZ_\dr)$$
as a crystal of (symmetric) monoidal categories over $\Ran_{\ul{x}}$, denote it by $\IndCoh^!(\Op^{\mf\rightsquigarrow\reg}_\cG)^{\on{fact}_{\ul{x}}}$.

\medskip

Note that it extends naturally to a crystal of (symmetric) monoidal categories over $\Ran^{\on{untl}}_{\ul{x}}$, in which the monoidal
structure is given by \emph{strict} functors between sheaves of categories. 

\sssec{} \label{sss:goal construct action}

The key step will be to construct an action of $$\IndCoh^!(\Op^{\mf\rightsquigarrow\reg}_\cG)^{\on{fact}_{\ul{x}}}$$ on 
$$\Res_{{\mathcal Vac}(G)_\crit}((\KL(G)_\kappa)^{\on{fact}_{\ul{x}}})$$ (as crystals of categories over $\Ran_{\ul{x}}$). Furthermore, we will
show this action extends to a \emph{strict} action as crystals of categories over $\Ran^{\on{untl}}_{\ul{x}}$.

\sssec{} \label{sss:ins x action}

The above strict compatibility means in particular that the functor
$$\on{ins.vac}_{\ul{x}}:\KL(G)_{\kappa,\ul{x}}\to \Res_{{\mathcal Vac}(G)_\crit}((\KL(G)_\kappa)^{\on{fact}_{\ul{x}}})_{\Ran_{\ul{x}}}$$
intertwines the $\IndCoh^!(\Op^\mf_{\cG,\ul{x}})$-action on $\KL(G)_{\kappa,\ul{x}}$ and the 
$\IndCoh^!(\Op^{\mf\rightsquigarrow\reg}_{\cG,\Ran_{\ul{x}}})$-action on 
$\Res_{{\mathcal Vac}(G)_\crit}((\KL(G)_\kappa)^{\on{fact}_{\ul{x}}})_{\Ran_{\ul{x}}}$
via the functor
$$\on{ins.unit}_{\ul{x}}:\IndCoh^!(\Op^\mf_{\cG,\ul{x}})\to \IndCoh^!(\Op^{\mf\rightsquigarrow\reg}_{\cG,\Ran_{\ul{x}}}),$$
while the latter is the functor of !-pullback along
$$\on{pr}^\Op_{\on{small},\ul{x}}:\Op^{\mf\rightsquigarrow\reg}_{\cG,\Ran_{\ul{x}}}\to \Op^\mf_{\cG,\ul{x}}.$$

Hence, we obtain a functor
\begin{equation} \label{e:act IndCoh on Res}
\IndCoh^!(\Op^{\mf\rightsquigarrow\reg}_{\cG,\Ran_{\ul{x}}})\underset{\IndCoh^!(\Op^\mf_{\cG,\ul{x}})}\otimes 
\KL(G)_{\kappa,\ul{x}}\to \Res_{{\mathcal Vac}(G)_\crit}((\KL(G)_\kappa)^{\on{fact}_{\ul{x}}})_{\Ran_{\ul{x}}}.
\end{equation}

\sssec{}

Composing with the tautological functor 
$$\Res_{{\mathcal Vac}(G)_\crit}((\KL(G)_\kappa)^{\on{fact}_{\ul{x}}})_{\Ran_{\ul{x}}}\to \KL(G)_{\kappa,\Ran_{\ul{x}}},$$
we obtain the sought-for functor
$$\on{ins.vac}^{\mf\rightsquigarrow\reg}_{\ul{x}}:
\IndCoh^!(\Op^{\mf\rightsquigarrow\reg}_{\cG,\Ran_{\ul{x}}})\underset{\IndCoh^!(\Op^\mf_{\cG,\ul{x}})}\otimes 
\KL(G)_{\kappa,\ul{x}}\to \KL(G)_{\kappa,\Ran_{\ul{x}}}$$
of \eqref{e:ins vac reg mf IndCoh}.

\ssec{The acting agents}

In this subsection we will interpret the category $\IndCoh^!(\Op^{\mf\rightsquigarrow\reg}_{\cG,\Ran_{\ul{x}}})$,
in terms of \emph{factorization restriction}. 

\sssec{}

Recall the map of factorization spaces
$$\iota^{+,\mf}:\Op^\reg_\cG\to \Op^\mf_\cG.$$

Consider the corresponding factorization functors
\begin{equation} \label{e:reg to mf fact *}
\QCoh(\Op^\reg_\cG)\simeq \IndCoh^*(\Op^\reg_\cG)\overset{(\iota^{+,\mf})^\IndCoh_*}\longrightarrow 
\IndCoh^*(\Op^\mf_\cG)
\end{equation} 
and 
\begin{equation} \label{e:reg to mf fact !}
\IndCoh^!(\Op^\reg_\cG)\overset{(\iota^{+,\mf})^\IndCoh_*}\longrightarrow 
\IndCoh^!(\Op^\mf_\cG).
\end{equation} 

\sssec{}

We start with
$$\IndCoh^*(\Op^\mf_\cG)^{\on{fact}_{\ul{x}}}\in \IndCoh^*(\Op^\mf_\cG)\mmod^{\text{fact}}_{\ul{x}}$$
and 
$$\IndCoh^!(\Op^\mf_\cG)^{\on{fact}_{\ul{x}}}\in \IndCoh^!(\Op^\mf_\cG)\mmod^{\text{fact}}_{\ul{x}}$$
and consider the resulting objects 
\begin{equation} \label{e:reg to mf fact * module}
\Res_{(\iota^{+,\mf})^\IndCoh_*}(\IndCoh^*(\Op^\mf_\cG)^{\on{fact}_{\ul{x}}})\in 
\QCoh(\Op^\reg_\cG)\mmod^{\on{fact}}_{\ul{x}}
\end{equation} 
and 
\begin{equation} \label{e:reg to mf fact ! module}
\Res_{(\iota^{+,\mf})^\IndCoh_*}(\IndCoh^!(\Op^\mf_\cG)^{\on{fact}_{\ul{x}}})\in 
\IndCoh^!(\Op^\reg_\cG)\mmod^{\on{fact}}_{\ul{x}}.
\end{equation} 

\medskip

Since the factorization functors \eqref{e:reg to mf fact *} and \eqref{e:reg to mf fact !} are unital,
the module categories \eqref{e:reg to mf fact * module} and \eqref{e:reg to mf fact ! module} 
have natural unital structures, see \secref{sss:unital restr categ}. 

\sssec{}

Let us consider the assignments 
$$(\CZ \to \Ran_{\ul{x}}) \mapsto \IndCoh^*(\Op^{\mf\rightsquigarrow\reg}_{\cG,\Ran_{\ul{x}}}\underset{\Ran_{\ul{x}}}\times \CZ_\dr)$$
and 
$$(\CZ \to \Ran_{\ul{x}}) \mapsto \IndCoh^!(\Op^{\mf\rightsquigarrow\reg}_{\cG,\Ran_{\ul{x}}}\underset{\Ran_{\ul{x}}}\times \CZ_\dr)$$
as crystals of categories over $\Ran_{\ul{x}}$ (the latter is the crystal of categories that we have introduced in \secref{sss:indCoh mf reg Ran x}).

\medskip

They have natural structures of unital module categories (at $\ul{x}$) over
$$\IndCoh^*(\Op^\reg_\cG) \text{ and } \IndCoh^!(\Op^\reg_\cG),$$
respectively. We will denote them by 
$$\IndCoh^*(\Op^{\mf\rightsquigarrow\reg}_\cG)^{\on{fact}_{\ul{x}}} \text{ and } 
\IndCoh^!(\Op^{\mf\rightsquigarrow\reg}_\cG)^{\on{fact}_{\ul{x}}},$$
respectively.

\sssec{}

We will regard $\IndCoh^*(\Op^{\mf\rightsquigarrow\reg}_\cG)^{\on{fact}_{\ul{x}}}$ 
(resp., $\IndCoh^!(\Op^{\mf\rightsquigarrow\reg}_\cG)^{\on{fact}_{\ul{x}}}$) as equipped with a comonoidal (resp., monoidal)
structure given by $\IndCoh^*$-pushforward (resp., $!$-pullback) along the diagonal morphism.

\medskip

We note that when we view $\IndCoh^*(\Op^{\mf\rightsquigarrow\reg}_\cG)^{\on{fact}_{\ul{x}}}$ as a crystal of categories 
over $\Ran^{\on{untl}}_{\ul{x}}$, its comonoidal structure is given by right-lax functors. 

\medskip

By contrast, when we view $\IndCoh^!(\Op^{\mf\rightsquigarrow\reg}_\cG)^{\on{fact}_{\ul{x}}}$ as a crystal of categories over $\Ran^{\on{untl}}_{\ul{x}}$,
its monoidal structure is given by \emph{strict} functors.

\sssec{}

The map $\iota^{+,\mf}$ gives rise to a map 
$$\iota^{+,\mf\to \reg}:\Op^{\mf\rightsquigarrow\reg}_{\cG,\Ran_{\ul{x}}}\to \Op^\mf_{\cG,\Ran_{\ul{x}}}.$$
(Note that $\iota^{+,\mf}$ is the same as the map $\on{pr}^\Op_{\on{big},\ul{x}}$). 

\medskip

We obtain functors of (unital) module categories
\begin{equation} \label{e:reg to mf fact * module pre comparison}
(\iota^{+,\mf\to \reg})^\IndCoh_*:\IndCoh^*(\Op^{\mf\rightsquigarrow\reg}_\cG)^{\on{fact}_{\ul{x}}} \to \IndCoh^*(\Op^\mf_\cG)^{\on{fact}_{\ul{x}}}
\end{equation} 
and 
\begin{equation} \label{e:reg to mf fact ! module pre comparison}
(\iota^{+,\mf\to \reg})^\IndCoh_*:\IndCoh^!(\Op^{\mf\rightsquigarrow\reg}_\cG)^{\on{fact}_{\ul{x}}} \to \IndCoh^!(\Op^\mf_\cG)^{\on{fact}_{\ul{x}}},
\end{equation} 
compatible with the functors \eqref{e:reg to mf fact *} and \eqref{e:reg to mf fact !}, respectively.

\sssec{}

By \secref{sss:univ property restr cat}, the functors \eqref{e:reg to mf fact * module pre comparison} and 
\eqref{e:reg to mf fact ! module pre comparison}
give rise to maps
\begin{multline} \label{e:reg to mf fact * module comparison}
(\iota^{+,\mf})^\IndCoh_*: \IndCoh^*(\Op^{\mf\rightsquigarrow\reg}_\cG)^{\on{fact}_{\ul{x}}} \to \\
\Res_{(\iota^{+,\mf})^\IndCoh_*}(\IndCoh^*(\Op^\mf_\cG)^{\on{fact}_{\ul{x}}})
\end{multline} 
and
\begin{multline} \label{e:reg to mf fact ! module comparison}
(\iota^{+,\mf})^\IndCoh_*: \IndCoh^!(\Op^{\mf\rightsquigarrow\reg}_\cG)^{\on{fact}_{\ul{x}}} \to \\
\to \Res_{(\iota^{+,\mf})^\IndCoh_*}(\IndCoh^!(\Op^\mf_\cG)^{\on{fact}_{\ul{x}}})
\end{multline} 
in 
$$\QCoh(\Op^\reg_\cG)\mmod^{\on{fact}}_{\ul{x}} \text{ and } \IndCoh^!(\Op^\reg_\cG)\mmod^{\on{fact}}_{\ul{x}},$$
respectively. Moreover, the maps \eqref{e:reg to mf fact * module comparison} and \eqref{e:reg to mf fact ! module comparison}
are compatible with the unital structures, see \lemref{l:unital restr categ univ}.

\sssec{}

The following is a variant of \lemref{l:iota Y mer to reg restr}: 

\begin{lem} \label{l:reg to mf fact module comparison}
The functors \eqref{e:reg to mf fact * module comparison} and \eqref{e:reg to mf fact ! module comparison} 
are equivalences.
\end{lem}

\begin{proof} 

We prove the assertion for $\IndCoh^*$. The case of $\IndCoh^!$ is analogous. 

\medskip

By \lemref{l:fact res crit}, it suffices to check that:

\smallskip

\noindent{(i)} The functor \eqref{e:reg to mf fact *} admits a right adjoint (as a functor between sheaves of categories);

\smallskip

\noindent{(ii)} The functor \eqref{e:reg to mf fact * module pre comparison} admits a right adjoint (as a functor between sheaves of categories);

\smallskip

\noindent{(iii)} The functor \eqref{e:reg to mf fact * module pre comparison} induces an equivalence between the fibers of the two
sides at $\ul{x}\in \Ran_{\ul{x}}$.

\medskip

We note that point (iii) holds tautologically.

\medskip 

Points (i) and (ii) are also automatic: the right adjoints in question are given by the funtors
$(\iota^{+,\mf})^!$ and $(\iota^{+,\mf\to \reg})^!$, respectively. 

\end{proof}

\ssec{Construction of the action}

\sssec{}

Recall the object 
$$\Res_{{\mathcal Vac}(G)_\crit}((\KL(G)_\kappa)^{\on{fact}_{\ul{x}}})\in \QCoh(\Op^\reg_\cG)\mmod^{\on{fact}}_{\ul{x}}.$$

\medskip

We claim that $\Res_{{\mathcal Vac}(G)_\crit}((\KL(G)_\kappa)^{\on{fact}_{\ul{x}}})$, viewed as a sheaf of categories
over $\Ran_{\ul{x}}$, carries a canonically defined 
action of $\Res_{(\iota^{+,\mf})^\IndCoh_*}(\IndCoh^!(\Op^\mf_\cG)^{\on{fact}_{\ul{x}}})$.

\sssec{}

By duality, a datum of such an action is equivalent to the datum of a coaction of the sheaf of comonoidal
categories 
$\Res_{(\iota^{+,\mf})^\IndCoh_*}(\IndCoh^*(\Op^\mf_\cG)^{\on{fact}_{\ul{x}}})$.

\medskip

We will construct the corresponding coaction functor
\begin{multline} \label{e:coact mm->reg on}
\Res_{{\mathcal Vac}(G)_\crit}((\KL(G)_\kappa)^{\on{fact}_{\ul{x}}})\to \\
\to \Res_{(\iota^{+,\mf})^\IndCoh_*}(\IndCoh^*(\Op^\mf_\cG)^{\on{fact}_{\ul{x}}})\otimes
\Res_{{\mathcal Vac}(G)_\crit}((\KL(G)_\kappa)^{\on{fact}_{\ul{x}}}).
\end{multline} 

The full datum of coaction is defined similar to \secref{ss:z on KL}, using the device from \secref{s:device}.

\sssec{} \label{sss:coact mm->reg on}

We interpret the right-hand side in \eqref{e:coact mm->reg on} as the restriction of
$$(\IndCoh^*(\Op^\mf_\cG)\otimes \KL(G)_\kappa)^{\on{fact}_{\ul{x}}}\in
(\IndCoh^*(\Op^\mf_\cG)\otimes \KL(G)_\kappa)\mmod^{\on{fact}}_{\ul{x}}$$
along the factorization functor
$$((\iota^{+,\mf})^\IndCoh_*\otimes {\mathcal Vac}(G)_\crit):
\QCoh(\Op^\reg_\cG)\otimes \QCoh(\Op^\reg_\cG)\to
\IndCoh^*(\Op^\mf_\cG)\otimes \KL(G)_\kappa.$$

\medskip

The functor \eqref{e:coact mm->reg on} is given by the procedure of restriction from \secref{sss:restr pairs} along the diagram
$$
\CD
\KL(G)_\kappa @>>> \IndCoh^*(\Op^\mf_\cG)\otimes \KL(G)_\kappa \\
@A{{\mathcal Vac}(G)_\crit}AA @AA{(\iota^{+,\mf})^\IndCoh_*\otimes {\mathcal Vac}(G)_\crit}A \\
\QCoh(\Op^\reg_\cG) @>>> \QCoh(\Op^\reg_\cG)\otimes \QCoh(\Op^\reg_\cG),
\endCD
$$
where:

\begin{itemize}

\item The top horizontal arrow is the coaction of $\IndCoh^*(\Op^\mf_\cG)$ on $\KL(G)_\kappa$;

\smallskip

\item The bottom horizontal arrow is the comonoidal operation, i.e., the
functor of direct image along the diagonal map.

\end{itemize}

\sssec{Example}

Here is what the above action (resp., coaction) does at the pointwise level. Write
$$\ul{x}'=\ul{x}\sqcup \ul{x}'',$$
so that 
$$\Res_{{\mathcal Vac}(G)_\crit}((\KL(G)_\kappa)^{\on{fact}_{\ul{x}}})_{\ul{x}'}\simeq 
\KL(G)_{\kappa,\ul{x}}\otimes \QCoh(\Op^\reg_{\cG,\ul{x}''}),$$
$$\Res_{(\iota^{+,\mf})^\IndCoh_*}(\IndCoh^*(\Op^\mf_\cG)^{\on{fact}_{\ul{x}}})_{\ul{x'}}\simeq
\IndCoh^*(\Op^\mf_{\cG,\ul{x}}) \otimes \QCoh(\Op^\reg_{\cG,\ul{x}''}),$$
$$\Res_{(\iota^{+,\mf})^\IndCoh_*}(\IndCoh^!(\Op^\mf_\cG)^{\on{fact}_{\ul{x}}})_{\ul{x'}}\simeq
\IndCoh^!(\Op^\mf_{\cG,\ul{x}}) \otimes \IndCoh^!(\Op^\reg_{\cG,\ul{x}''}).$$

The coaction of $\Res_{(\iota^{+,\mf})^\IndCoh_*}(\IndCoh^*(\Op^\mf_\cG)^{\on{fact}_{\ul{x}}})_{\ul{x'}}$
on $\Res_{{\mathcal Vac}(G)_\crit}((\KL(G)_\kappa)^{\on{fact}_{\ul{x}}})_{\ul{x}'}$ acts as the tensor product of

\begin{itemize}

\item The coaction of $\IndCoh^*(\Op^\mf_{\cG,\ul{x}})$ on $\KL(G)_{\kappa,\ul{x}}$, and 

\item The functor of direct image along the diagonal map $\QCoh(\Op^\reg_{\cG,\ul{x}''})\to 
\QCoh(\Op^\reg_{\cG,\ul{x}''})\otimes \QCoh(\Op^\reg_{\cG,\ul{x}''})$.

\end{itemize}

The action of $\Res_{(\iota^{+,\mf})^\IndCoh_*}(\IndCoh^!(\Op^\mf_\cG)^{\on{fact}_{\ul{x}}})_{\ul{x'}}$
on $\Res_{{\mathcal Vac}(G)_\crit}((\KL(G)_\kappa)^{\on{fact}_{\ul{x}}})_{\ul{x}'}$ acts as the tensor product of

\begin{itemize}

\item The action of $\IndCoh^!(\Op^\mf_{\cG,\ul{x}})$ on $\KL(G)_{\kappa,\ul{x}}$, and 

\item The canonical action of $\IndCoh^!(\Op^\reg_{\cG,\ul{x}''})$ on 
$\IndCoh^*(\Op^\reg_{\cG,\ul{x}''})\simeq \QCoh(\Op^\reg_{\cG,\ul{x}''})$.

\end{itemize}

\ssec{The unital structure on the action functor}

\sssec{}

Let us regard 
$$\KL(G)_\kappa \text{ and } \IndCoh^*(\Op^\mf_\cG)$$
as crystals of categories on $\Ran$, equipped with a \emph{unital} structure (see \secref{sss:unitality loc} for what this means).

\medskip

The coaction of $\IndCoh^*(\Op^\mf_\cG)$ on $\KL(G)_\kappa$ has the following feature: it extends to a coaction 
in the 2-category of crystals of categories on $\Ran^{\on{untl}}$ with \emph{right-lax} functors.

\medskip

This follows from the construction of this coaction in \secref{ss:z on KL}, using the following observation:

\medskip

For a map of factorization algebras $\CA_1\to \CA_2$ in a given factorization category $\bA$, the restriction 
functor
$$\CA_2\mod(\bA)\to \CA_1\mod(\bA),$$
viewed as a functor between crystals of categories on $\Ran$, admits a natural extension to a right-lax
functor between crystals of categories on $\Ran^{\on{untl}}$. 
 
\sssec{}

It follows from the construction in \secref{sss:coact mm->reg on} and \secref{sss:unital restr categ pairs} that
the functor \eqref{e:coact mm->reg on} extends to a right-lax functor between crystals of categories on $\Ran_{\ul{x}}^{\on{untl}}$. 

\medskip

By a similar token, we obtain that 
the full datum of coaction of the comonoidal category $\Res_{(\iota^{+,\mf})^\IndCoh_*}(\IndCoh^*(\Op^\mf_\cG)^{\on{fact}_{\ul{x}}})$
on $\Res_{{\mathcal Vac}(G)_\crit}((\KL(G)_\kappa)^{\on{fact}_{\ul{x}}})$ extends to a coaction 
in the 2-category of crystals of categories on $\Ran^{\on{untl}}_{\ul{x}}$ with right-lax functors.

\sssec{}

Combining with \lemref{l:reg to mf fact module comparison}, we obtain that $\Res_{{\mathcal Vac}(G)_\crit}((\KL(G)_\kappa)^{\on{fact}_{\ul{x}}})$
carries a coaction of $\IndCoh^*(\Op^{\mf\rightsquigarrow\reg}_\cG)^{\on{fact}_{\ul{x}}}$. 

\medskip

Furthermore, this coaction 
extends to a coaction  in the 2-category of crystals of categories on $\Ran^{\on{untl}}_{\ul{x}}$ with right-lax functors.

\sssec{} \label{sss:counit right-lax}

We note now that the unital structures on 
$$\IndCoh^!(\Op^{\mf\rightsquigarrow\reg}_\cG)^{\on{fact}_{\ul{x}}} \text{ and }
\IndCoh^*(\Op^{\mf\rightsquigarrow\reg}_\cG)^{\on{fact}_{\ul{x}}}$$
have the following feature:

\medskip

The counit of the duality
$$\IndCoh^!(\Op^{\mf\rightsquigarrow\reg}_\cG)^{\on{fact}_{\ul{x}}} \otimes \IndCoh^!(\Op^{\mf\rightsquigarrow\reg}_\cG)^{\on{fact}_{\ul{x}}} \to
\ul\Dmod(\Ran_{\ul{x}})$$
extends to a right-lax functor between crystals of categories on $\Ran^{\on{untl}}_{\ul{x}}$.

\sssec{}

The coaction of $\IndCoh^*(\Op^{\mf\rightsquigarrow\reg}_\cG)^{\on{fact}_{\ul{x}}}$ on 
$\Res_{{\mathcal Vac}(G)_\crit}((\KL(G)_\kappa)^{\on{fact}_{\ul{x}}})$
gives rise to an action of $\IndCoh^!(\Op^{\mf\rightsquigarrow\reg}_\cG)^{\on{fact}_{\ul{x}}}$ on 
$\Res_{{\mathcal Vac}(G)_\crit}((\KL(G)_\kappa)^{\on{fact}_{\ul{x}}})$.

\medskip

Moreover, by \secref{sss:counit right-lax}, this action extends to an action in the 2-category of sheaves of 
categories on $\Ran^{\on{untl}}$ with \emph{right-lax} functors.

\sssec{}

We now claim:

\begin{lem} \label{l:action is strictly unital}
The right-lax functors that define the 
action of $\IndCoh^!(\Op^{\mf\rightsquigarrow\reg}_\cG)^{\on{fact}_{\ul{x}}}$ on $\Res_{{\mathcal Vac}(G)_\crit}((\KL(G)_\kappa)^{\on{fact}_{\ul{x}}})$
as crystals of categories on $\Ran^{\on{untl}}$ are \emph{strict}. 
\end{lem}

\begin{proof}

We need to show the following: for $\ul{x}\subseteq \ul{x}_1\subseteq \ul{x}_2$, the natural transformation from
\begin{multline} \label{e:action is strictly unital 1}
\IndCoh^!(\Op^{\mf\rightsquigarrow\reg}_{\cG,\ul{x}\subseteq \ul{x}_1})\otimes 
\Res_{{\mathcal Vac}(G)_\crit}((\KL(G)_\kappa)^{\on{fact}_{\ul{x}}})_{\ul{x}\subseteq \ul{x}_1}\overset{\text{action}}\to \\
\to \Res_{{\mathcal Vac}(G)_\crit}((\KL(G)_\kappa)^{\on{fact}_{\ul{x}}})_{\ul{x}\subseteq \ul{x}_1} \overset{\on{ins.unit}_{\ul{x}_1\subseteq \ul{x}_2}}\longrightarrow
\Res_{{\mathcal Vac}(G)_\crit}((\KL(G)_\kappa)^{\on{fact}_{\ul{x}}})_{\ul{x}\subseteq \ul{x}_2}
\end{multline} 
to 
\begin{multline} \label{e:action is strictly unital 2}
\IndCoh^!(\Op^{\mf\rightsquigarrow\reg}_{\cG,\ul{x}\subseteq \ul{x}_1})\otimes 
\Res_{{\mathcal Vac}(G)_\crit}((\KL(G)_\kappa)^{\on{fact}_{\ul{x}}})_{\ul{x}\subseteq \ul{x}_1} 
\overset{\on{ins.unit}_{\ul{x}_1\subseteq \ul{x}_2}\otimes \on{ins.unit}_{\ul{x}_1\subseteq \ul{x}_2}}\longrightarrow \\
\to \IndCoh^!(\Op^{\mf\rightsquigarrow\reg}_{\cG,\ul{x}\subseteq \ul{x}_2})\otimes 
\Res_{{\mathcal Vac}(G)_\crit}((\KL(G)_\kappa)^{\on{fact}_{\ul{x}}})_{\ul{x}\subseteq \ul{x}_2} \overset{\text{action}}\to \\
\to \Res_{{\mathcal Vac}(G)_\crit}((\KL(G)_\kappa)^{\on{fact}_{\ul{x}}})_{\ul{x}\subseteq \ul{x}_2}
\end{multline} 
is an isomorphism.

\medskip

The question of a functor between crystals of categories (in this case, over $\Ran_{\ul{x}}$) being an isomorphism 
can be checked strata-wise. So we can assume that
$$\ul{x}_2=\ul{x}_1\sqcup \ul{x}'.$$

We identify 
\begin{equation} \label{e:action is strictly unital 3}
\Res_{{\mathcal Vac}(G)_\crit}((\KL(G)_\kappa)^{\on{fact}_{\ul{x}}})_{\ul{x}\subseteq \ul{x}_2} \simeq
\Res_{{\mathcal Vac}(G)_\crit}((\KL(G)_\kappa)^{\on{fact}_{\ul{x}}})_{\ul{x}\subseteq \ul{x}_1} \otimes \QCoh(\Op^\reg_{\cG,\ul{x}''})
\end{equation} 
and
\begin{multline} \label{e:action is strictly unital 4}
\IndCoh^!(\Op^{\mf\rightsquigarrow\reg}_{\cG,\ul{x}\subseteq \ul{x}_2})\otimes 
\Res_{{\mathcal Vac}(G)_\crit}((\KL(G)_\kappa)^{\on{fact}_{\ul{x}}})_{\ul{x}\subseteq \ul{x}_2} \simeq \\
\simeq 
\IndCoh^!(\Op^{\mf\rightsquigarrow\reg}_{\cG,\ul{x}\subseteq \ul{x}_1})\otimes 
\Res_{{\mathcal Vac}(G)_\crit}((\KL(G)_\kappa)^{\on{fact}_{\ul{x}}})_{\ul{x}\subseteq \ul{x}_1} \bigotimes \\ 
\bigotimes \IndCoh^!(\Op^\reg_{\cG,\ul{x}''})\otimes \QCoh(\Op^\reg_{\cG,\ul{x}''}),
\end{multline}
so that the functor 
$$\Res_{{\mathcal Vac}(G)_\crit}((\KL(G)_\kappa)^{\on{fact}_{\ul{x}}})_{\ul{x}\subseteq \ul{x}_1} \overset{\on{ins.unit}_{\ul{x}_1\subseteq \ul{x}_2}}\longrightarrow
\Res_{{\mathcal Vac}(G)_\crit}((\KL(G)_\kappa)^{\on{fact}_{\ul{x}}})_{\ul{x}\subseteq \ul{x}_2}$$
identifies with 
$$\on{Id}\otimes \CO_{\Op^\reg_{\cG,\ul{x}''}}$$
and the functor
\begin{multline*} 
\IndCoh^!(\Op^{\mf\rightsquigarrow\reg}_{\cG,\ul{x}\subseteq \ul{x}_1})\otimes 
\Res_{{\mathcal Vac}(G)_\crit}((\KL(G)_\kappa)^{\on{fact}_{\ul{x}}})_{\ul{x}\subseteq \ul{x}_1} 
\overset{\on{ins.unit}_{\ul{x}_1\subseteq \ul{x}_2}\otimes \on{ins.unit}_{\ul{x}_1\subseteq \ul{x}_2}}\longrightarrow \\
\to \IndCoh^!(\Op^{\mf\rightsquigarrow\reg}_{\cG,\ul{x}\subseteq \ul{x}_2})\otimes 
\Res_{{\mathcal Vac}(G)_\crit}((\KL(G)_\kappa)^{\on{fact}_{\ul{x}}})_{\ul{x}\subseteq \ul{x}_2}
\end{multline*}
identifies with 
$$\on{Id}\otimes \omega_{\Op^\reg_{\cG,\ul{x}''}}\otimes \CO_{\Op^\reg_{\cG,\ul{x}''}}.$$

\medskip

In terms of \eqref{e:action is strictly unital 3} and \eqref{e:action is strictly unital 4}, the functor
$$\IndCoh^!(\Op^{\mf\rightsquigarrow\reg}_{\cG,\ul{x}\subseteq \ul{x}_2})\otimes 
\Res_{{\mathcal Vac}(G)_\crit}((\KL(G)_\kappa)^{\on{fact}_{\ul{x}}})_{\ul{x}\subseteq \ul{x}_2} \overset{\text{action}}\to 
\Res_{{\mathcal Vac}(G)_\crit}((\KL(G)_\kappa)^{\on{fact}_{\ul{x}}})_{\ul{x}\subseteq \ul{x}_2}$$
is the tensor product of
$$\IndCoh^!(\Op^{\mf\rightsquigarrow\reg}_{\cG,\ul{x}\subseteq \ul{x}_1})\otimes 
\Res_{{\mathcal Vac}(G)_\crit}((\KL(G)_\kappa)^{\on{fact}_{\ul{x}}})_{\ul{x}\subseteq \ul{x}_1} \overset{\text{action}}\to 
\Res_{{\mathcal Vac}(G)_\crit}((\KL(G)_\kappa)^{\on{fact}_{\ul{x}}})_{\ul{x}\subseteq \ul{x}_1}$$
along the first factor and 
\begin{multline*} 
\IndCoh^!(\Op^\reg_{\cG,\ul{x}''})\otimes \QCoh(\Op^\reg_{\cG,\ul{x}''}) \simeq 
\IndCoh^!(\Op^\reg_{\cG,\ul{x}''})\otimes \IndCoh^*(\Op^\reg_{\cG,\ul{x}''}) \overset{\sotimes}\to \\
\to \IndCoh^*(\Op^\reg_{\cG,\ul{x}''}) \simeq \QCoh(\Op^\reg_{\cG,\ul{x}''})
\end{multline*}
along the second factor. 

\medskip

Hence, we obtain that, in terms of \eqref{e:action is strictly unital 3} and \eqref{e:action is strictly unital 4},
both functors \eqref{e:action is strictly unital 1} and \eqref{e:action is strictly unital 2} are identified with
\begin{multline}  \label{e:action is strictly unital 5}
\IndCoh^!(\Op^{\mf\rightsquigarrow\reg}_{\cG,\ul{x}\subseteq \ul{x}_1})\otimes 
\Res_{{\mathcal Vac}(G)_\crit}((\KL(G)_\kappa)^{\on{fact}_{\ul{x}}})_{\ul{x}\subseteq \ul{x}_1} \overset{\text{action}}\to \\
\to \Res_{{\mathcal Vac}(G)_\crit}((\KL(G)_\kappa)^{\on{fact}_{\ul{x}}})_{\ul{x}\subseteq \ul{x}_1} 
\overset{\on{Id} \otimes \omega_{\Op^\reg_{\cG,\ul{x}''}}}\longrightarrow \\
\to \Res_{{\mathcal Vac}(G)_\crit}((\KL(G)_\kappa)^{\on{fact}_{\ul{x}}})_{\ul{x}\subseteq \ul{x}_1} \otimes \IndCoh^!(\Op^\reg_{\cG,\ul{x}''}).
\end{multline}

Unwinding the construction, we obtain that the endomorphism of \eqref{e:action is strictly unital 5} defined
by the structure of right-lax functor on the action map is the identity map.

\end{proof} 

\ssec{End of the construction}

\sssec{}

Thus, we have carried out the construction announced in \secref{sss:goal construct action}. In particular, we obtain a functor
\begin{equation} \label{e:ins vac reg mf pre} 
\IndCoh^!(\Op^{\mf\rightsquigarrow\reg}_{\cG,\Ran_{\ul{x}}})\underset{\IndCoh^!(\Op^\mf_{\cG,\ul{x}})}\otimes 
\KL(G)_{\crit,\ul{x}}\to \Res_{{\mathcal Vac}(G)_\crit}((\KL(G)_\kappa)^{\on{fact}_{\ul{x}}})_{\Ran_{\ul{x}}},
\end{equation} 
and its composition of the functor \eqref{e:act IndCoh on Res} with 
\begin{equation} \label{e:from res taut}
\Res_{{\mathcal Vac}(G)_\crit}((\KL(G)_\kappa)^{\on{fact}_{\ul{x}}})_{\Ran_{\ul{x}}}\to \KL(G)_{\kappa,\Ran_{\ul{x}}}
\end{equation}
produces the sought-for functor 
\begin{equation} \label{e:act IndCoh on Res again}
\on{ins.vac}^{\mf\rightsquigarrow\reg}_{\ul{x}}:
\IndCoh^!(\Op^{\mf\rightsquigarrow\reg}_{\cG,\Ran_{\ul{x}}})\underset{\IndCoh^!(\Op^\mf_{\cG,\ul{x}})}\otimes 
\KL(G)_{\kappa,\ul{x}}\to \KL(G)_{\kappa,\Ran_{\ul{x}}}. 
\end{equation}

\sssec{}

It remains to show that the functor \eqref{e:act IndCoh on Res again} is $\IndCoh^!(\Op^{\mf\rightsquigarrow\reg}_{\cG,\Ran_{\ul{x}}})$-linear. 

\medskip

The functor \eqref{e:ins vac reg mf pre} is $\IndCoh^!(\Op^{\mf\rightsquigarrow\reg}_{\cG,\Ran_{\ul{x}}})$-linear by construction.
So it remains to show that the functor \eqref{e:from res taut} is also $\IndCoh^!(\Op^{\mf\rightsquigarrow\reg}_{\cG,\Ran_{\ul{x}}})$-linear. 

\sssec{}

Unwinding the construction of the coaction of $\Res_{(\iota^{+,\mf})^\IndCoh_*}(\IndCoh^*(\Op^\mf_\cG)^{\on{fact}_{\ul{x}}})$
on $\Res_{{\mathcal Vac}(G)_\crit}((\KL(G)_\kappa)^{\on{fact}_{\ul{x}}})$ in \secref{sss:coact mm->reg on}, we obtain that the functor
\eqref{e:from res taut}
intertwines the coaction of $\Res_{(\iota^{+,\mf})^\IndCoh_*}(\IndCoh^*(\Op^\mf_\cG)^{\on{fact}_{\ul{x}}})$ on the left-hand side
with the coaction of $\IndCoh^*(\Op^{\mf\rightsquigarrow\reg}_{\cG,\Ran_{\ul{x}}})$ on the right hand side via the tautological forgetful functor
$$\Res_{(\iota^{+,\mf})^\IndCoh_*}(\IndCoh^*(\Op^\mf_\cG)^{\on{fact}_{\ul{x}}})\to 
\IndCoh^*(\Op^{\mf\rightsquigarrow\reg}_{\cG,\Ran_{\ul{x}}}).$$

\medskip

Hence, the functor \eqref{e:from res taut} intertwines the 
coaction of $\IndCoh^*(\Op^\mf_{\cG,\Ran_{\ul{x}}})$ on the left-hand side of \eqref{e:from res taut} with the action of 
$\IndCoh^*(\Op^{\mf\rightsquigarrow\reg}_{\cG,\Ran_{\ul{x}}})$ on the right-hand side of \eqref{e:from res taut} via 
$(\on{pr}^\Op_{\on{big},\ul{x}})^\IndCoh_*$.

\medskip

Passing to the dual of the acting agents, we obtain that functor \eqref{e:from res taut} intertwines the 
action of $\IndCoh^!(\Op^\mf_{\cG,\Ran_{\ul{x}}})$ on the left-hand side of \eqref{e:from res taut} with the action of 
$\IndCoh^!(\Op^{\mf\rightsquigarrow\reg}_{\cG,\Ran_{\ul{x}}})$ on the right-hand side of \eqref{e:from res taut} via 
$(\on{pr}^\Op_{\on{big},\ul{x}})^!$, as desired. 

\qed[\thmref{t:ins vac reg}]

\sssec{}

The next assertion is not needed for the sequel; we mention it for the sake of completeness:

\begin{prop}  \label{p:ins vac reg mf pre QCoh}
The functor 
\begin{equation} \label{e:ins vac reg mf pre QCoh}
\QCoh(\Op^{\mf\rightsquigarrow\reg}_{\cG,\Ran_{\ul{x}}})\underset{\QCoh(\Op^\mf_{\cG,\ul{x}})}\otimes 
\KL(G)_{\crit,\ul{x}}\to \Res_{{\mathcal Vac}(G)_\crit}((\KL(G)_\kappa)^{\on{fact}_{\ul{x}}})_{\Ran_{\ul{x}}},
\end{equation}
induced by \eqref{e:ins vac reg mf pre}, is an equivalence.
\end{prop} 

The rest of this subsection is devoted to the proof of this proposition.

\sssec{}

The functor \eqref{e:ins vac reg mf pre QCoh} comes from a morphism in the 2-category $\QCoh(\Op^\reg_\cG)\mmod^{\on{fact}}_{\ul{x}}$:
\begin{equation} \label{e:ins vac reg mf pre QCoh 1}
\QCoh(\Op^{\mf\rightsquigarrow\reg}_\cG)^{\on{fact}_x}\underset{\QCoh(\Op^\mf_{\cG,\ul{x}})}\otimes 
\KL(G)_{\crit,\ul{x}}\to \Res_{{\mathcal Vac}(G)_\crit}((\KL(G)_\kappa)^{\on{fact}_{\ul{x}}})
\end{equation}

\medskip

Hence, by \lemref{l:fact res crit}, in order to check that \eqref{e:ins vac reg mf pre QCoh} is an equivalence, it suffices to check that:

\medskip

\noindent{(i)} The composition of \eqref{e:ins vac reg mf pre QCoh 1} with the tautological functor
$$\Res_{{\mathcal Vac}(G)_\crit}((\KL(G)_\kappa)^{\on{fact}_{\ul{x}}})\to
(\KL(G)_\kappa)^{\on{fact}_{\ul{x}}}$$
admits a right adjoint (as a functor between sheaves of categories);

\medskip

\noindent{(ii)} The functor ${\mathcal Vac}(G)_\crit$ admits a right adjoint (as a functor between sheaves of categories);

\medskip

\noindent{(iii)} The functor \eqref{e:ins vac reg mf pre QCoh 1} induces an equivalence between the fibers of the two
sides at $\ul{x}\in \Ran_{\ul{x}}$.

\sssec{}

Point (iii) above is immediate. Point (ii) follows from the fact that the functor ${\mathcal Vac}(G)_\crit$ preserves compactness.

\medskip

Hence, it remains to show that the composition in point (i) preserves compactness (and the left-hand side is compactly generated). 

\sssec{}

Let $Y$ be as in \secref{sss:Y for opers}.  

\medskip

It is easy to see that the corresponding category 
$\KL(G)_{\crit,\ul{x},Y}$ is compactly generated. By \lemref{l:approx mod}, it suffices to show that the composition
\begin{multline} \label{e:ins vac reg mf pre QCoh 2}
\QCoh(\Op^{\mf\rightsquigarrow\reg}_\cG)^{\on{fact}_x}\underset{\QCoh(\Op^\mf_{\cG,\ul{x}})}\otimes \KL(G)_{\crit,\ul{x},Y}\to \\
\to 
\QCoh(\Op^{\mf\rightsquigarrow\reg}_\cG)^{\on{fact}_x}\underset{\QCoh(\Op^\mf_{\cG,\ul{x}})}\otimes \KL(G)_{\crit,\ul{x}}\to \\
\to \Res_{{\mathcal Vac}(G)_\crit}((\KL(G)_\kappa)^{\on{fact}_{\ul{x}}})\to (\KL(G)_\kappa)^{\on{fact}_{\ul{x}}}
\end{multline}
preserves compactness. 

\medskip

We rewrite the left-hand side in \eqref{e:ins vac reg mf pre QCoh 2} as
$$\QCoh(\Op^{\mf\rightsquigarrow\reg}_\cG)^{\on{fact}_x}\underset{\QCoh(\Op^\mf_{\cG,\ul{x}})}\otimes \QCoh(Y)\underset{\QCoh(Y)}\otimes 
\KL(G)_{\crit,\ul{x},Y},$$
and further, by \lemref{l:base change poles at x} as
\begin{equation} \label{e:ins vac reg mf pre QCoh 3}
\QCoh(\Op^{\mf\rightsquigarrow\reg}_\cG\underset{\Op^\mf_{\cG,\ul{x}}}\times Y)^{\on{fact}_x}\underset{\QCoh(Y)}\otimes 
\KL(G)_{\crit,\ul{x},Y}.
\end{equation} 

\sssec{}

Now, since $\Op^{\mf\rightsquigarrow\reg}_\cG\underset{\Op^\mf_{\cG,\ul{x}}}\times Y$ is a relative affine \emph{scheme}
(as opposed to ind-scheme) over $\Ran_{\ul{x}}$, the category \eqref{e:ins vac reg mf pre QCoh 3} is compactly generated 
by the essential image of compact objects along the functor 
$$\KL(G)_{\crit,\ul{x},Y}\to \QCoh(\Op^{\mf\rightsquigarrow\reg}_\cG\underset{\Op^\mf_{\cG,\ul{x}}}\times Y)^{\on{fact}_x}\underset{\QCoh(Y)}\otimes 
\KL(G)_{\crit,\ul{x},Y},$$
given by tensoring with the structure sheaf along the first factor.

\medskip

Hence, it suffices to show that the functor 
\begin{multline} \label{e:ins vac reg mf pre QCoh 4}
\KL(G)_{\crit,\ul{x},Y}\to \QCoh(\Op^{\mf\rightsquigarrow\reg}_\cG)^{\on{fact}_x}\underset{\QCoh(\Op^\mf_{\cG,\ul{x}})}\otimes \KL(G)_{\crit,\ul{x},Y}\to \\
\to 
\QCoh(\Op^{\mf\rightsquigarrow\reg}_\cG)^{\on{fact}_x}\underset{\QCoh(\Op^\mf_{\cG,\ul{x}})}\otimes \KL(G)_{\crit,\ul{x}}\to \\
\to \Res_{{\mathcal Vac}(G)_\crit}((\KL(G)_\kappa)^{\on{fact}_{\ul{x}}})\to (\KL(G)_\kappa)^{\on{fact}_{\ul{x}}}
\end{multline}
preserves compactness.

\sssec{}

The forgetful functor $$\KL(G)_{\crit,\ul{x},Y}\to \KL(G)_{\crit,\ul{x}}$$ preserves compactness (indeed, it admits a continuous
right adjoint). Hence, it suffices to show that the functor
$$
\KL(G)_{\crit,\ul{x}}\to \QCoh(\Op^{\mf\rightsquigarrow\reg}_\cG)^{\on{fact}_x}\underset{\QCoh(\Op^\mf_{\cG,\ul{x}})}\otimes 
\KL(G)_{\crit,\ul{x}}\to (\KL(G)_\kappa)^{\on{fact}_{\ul{x}}}$$
preserves compactness.

\medskip

However, by construction, the latter functor is 
$$\on{ins.vac}_{\ul{x}}:\KL(G)_{\crit,\ul{x}}\to (\KL(G)_\kappa)^{\on{fact}_{\ul{x}}},$$
and the assertion follows.

\qed[\propref{p:ins vac reg mf pre QCoh}]

\section{Spectral Poincar\'e functor(s)} \label{s:spectral Poinc}

In this section we start dealing with the local-to-global constructions on the spectral side,
i.e., when the recipient category is $\IndCoh(\LS_\cG)$.

\medskip

We introduce two versions of the spectral Poincar\'e functor:
$$\IndCoh^!(\Op^{\on{mon-free}}_\cG)_\Ran \overset{\Poinc^{\on{spec}}_{\cG,!}}\longrightarrow
\IndCoh_\Nilp(\LS_\cG(X))$$
and 
$$\IndCoh^*(\Op^{\on{mon-free}}_\cG)_\Ran \overset{\Poinc^{\on{spec}}_{\cG,*}}\longrightarrow
\IndCoh_\Nilp(\LS_\cG(X)).$$

However, we show (\thmref{t:Poinc spec * vs !}) that they are intertwined by the ``self-duality" functor
$$\Theta_{\Op^\mf_\cG}:\IndCoh^!(\Op^{\on{mon-free}}_\cG)\to 
\IndCoh^*(\Op^{\on{mon-free}}_\cG),$$
up to tensoring by a graded line.

\medskip

Next we recall the definition of the \emph{spectral localization and global sections functors}
$$\Loc_\cG^{\on{spec}}:\Rep(\cG)_\Ran\rightleftarrows \IndCoh_\Nilp(\LS_\cG):\Gamma^{\on{spec},\IndCoh}_\cG.$$

Finally, we give the expression for the composition
$$\IndCoh^*(\Op^{\on{mon-free}}_\cG)_\Ran \overset{\Poinc^{\on{spec}}_{\cG,*}}\longrightarrow
\IndCoh_\Nilp(\LS_\cG)\overset{\Gamma^{\on{spec},\IndCoh}_\cG}\longrightarrow \Rep(\cG)_\Ran$$
via factorization homology, which exactly matches the composition \eqref{e:coeff Loc} under
$\FLE_{G,\crit}$ and $\FLE_{\cG,\infty}$.

\ssec{Ind-coherent sheaves on local vs. global opers}

\sssec{}

For $\CZ\to \Ran$ consider the morphism 
$$\Op^{\mer,\on{glob}}_{\cG,\CZ} \overset{\on{ev}_\CZ}\to \Op^{\mer}_{\cG,\CZ}.$$

Note that the prestack $\Op^{\mer,\on{glob}}_{\cG,\CZ}$ is locally almost of finite type, so we have a well-defined 
category $\IndCoh(\Op^{\mer,\on{glob}}_\cG)_\CZ:=\IndCoh(\Op^{\mer,\on{glob}}_{\cG,\CZ})$.

\medskip

Consider the pair of \emph{mutually dual} functors
\begin{equation} \label{e:local to global Op _* on *}
(\on{ev}_\CZ)^\IndCoh_*:\IndCoh(\Op^{\mer,\on{glob}}_\cG)_\CZ\to \IndCoh^*(\Op^{\mer}_\cG)_\CZ
\end{equation} 
and 
\begin{equation} \label{e:local to global Op ^! on !}
(\on{ev}_\CZ)^!:\IndCoh^!(\Op^{\mer}_\cG)_\CZ\to \IndCoh(\Op^{\mer,\on{glob}}_\cG)_\CZ.
\end{equation} 

\sssec{}

We claim:

\begin{lem} \label{l:local to global Op ^! on !}
The functor \eqref{e:local to global Op ^! on !} preserves compactness.
\end{lem} 

\sssec{} \label{sss:lattice notation}

Before we prove \lemref{l:local to global Op ^! on !}, we need to introduce some notation. 
For expositional purposes, we will assume that $\CZ=\on{pt}$, so that $\CZ\to \Ran$
corresponds to $\ul{x}\in \Ran$. 

\medskip

Let $\bV$ denote the Tate vector space 
$$\Gamma(\cD^\times_{\ul{x}},\fa(\cg)_{\omega_X}).$$

Let $\bL_0\subset \bV$ denote the standard lattice, i.e., 
$$\bL_0:=\Gamma(\cD_{\ul{x}},\fa(\cg)_{\omega_X}).$$

Recall that, according to \secref{sss:Op is a torsor}, we have a simply-transitive action of 
$\bL_0$ on $\Op^\reg_{\cG,\ul{x}}$, so that
$$\bV\overset{\bL_0}\times \Op^\reg_{\cG,\ul{x}}\simeq \Op^\mer_{\cG,\ul{x}}.$$

\medskip

For a lattice $\bL\supset \bL_0$ denote
$$\bL\overset{\bL_0}\times \Op^\reg_{\cG,\ul{x}}=:\Op^\bL_{\cG,\ul{x}}\subset \Op^\mer_{\cG,\ul{x}}.$$

\medskip

Denote
$$\bV^{\on{glob}}:=\Gamma(X-\ul{x},\fa(\cg)_{\omega_X}).$$

We have a simply-transitive action of $\bV^{\on{glob}}$ on $\Op^{\mer,\on{glob}}_{\cG,\ul{x}}$, compatible
with the embedding $\on{ev}_{\ul{x}}$ via
$$\bV^{\on{glob}}\hookrightarrow \bV.$$

For $\bL \supset \bL_0$ set
$$\Op^{\bL,\on{glob}}_{\cG,\ul{x}}:=\Op^{\mer,\on{glob}}_{\cG,\ul{x}}\underset{\Op^\mer_{\cG,\ul{x}}}\times \Op^\bL_{\cG,\ul{x}}.$$

\sssec{Proof of \lemref{l:local to global Op ^! on !}}

We have
\begin{equation} \label{e:IndCoh mer as colimit}
\IndCoh^!(\Op^{\mer}_\cG)_\CZ\simeq \underset{\bL\subset \bL_0}{\on{colim}}\, \IndCoh^!(\Op^\bL_\cG)_\CZ,
\end{equation} 
where the colimit is taken with respect to the $\IndCoh$-pushforward functors.

\medskip

Hence, it is enough to show that the composition
$$\IndCoh^!(\Op^\bL_\cG)_\CZ\to \IndCoh^!(\Op^\mer_\cG)_\CZ \overset{\on{ev}^!_{\ul{x}}}\longrightarrow 
\IndCoh(\Op^{\mer,\on{glob}}_\cG)_\CZ$$
preserves compactness.

\medskip

The Cartesian diagram
\begin{equation} \label{e:cart diag L Op}
\CD
\Op^{\bL,\on{glob}}_{\cG,\ul{x}} @>{\on{ev}_{\ul{x}}}>> \Op^{\bL}_{\cG,\ul{x}} \\
@VVV @VVV \\
\Op^{\mer,\on{glob}}_{\cG,\ul{x}} @>>{\on{ev}_{\ul{x}}}> \Op^\mer_{\cG,\ul{x}} 
\endCD
\end{equation} 
gives rise to a commutative diagram
$$
\CD
\IndCoh^!(\Op^{\bL,\on{glob}}_{\cG,\ul{x}}) @<{(\on{ev}_{\ul{x}})^!}<< \IndCoh^!(\Op^{\bL}_{\cG,\ul{x}}) \\
@V{\on{pushforward}}VV @VV{\on{pushforward}}V \\
\IndCoh^!(\Op^{\mer,\on{glob}}_{\cG,\ul{x}}) @<<{(\on{ev}_{\ul{x}})^!}< \IndCoh^!(\Op^\mer_{\cG,\ul{x}}),
\endCD
$$
see \secref{sss:BeckChev placid ind}.

\medskip

Hence, it suffices to show that the functor
$$(\on{ev}_{\ul{x}})^!:\IndCoh^!(\Op^{\bL}_{\cG,\ul{x}}) \to \IndCoh^!(\Op^{\bL,\on{glob}}_{\cG,\ul{x}})$$
preserves compactness.

\medskip

Write
$$\IndCoh^!(\Op^{\bL}_{\cG,\ul{x}}) \simeq \underset{\bL'\subset \bL_0}{\on{colim}}\, \IndCoh(\Op^\bL_{\cG,\CZ}/\bL'),$$
where the colimit is taken with respect to the !-pullback functors.

\medskip

Hence, it suffices to show that the !-pullback functors along
$$\Op^{\bL}_{\cG,\ul{x}}\to \Op^{\bL,\on{glob}}_{\cG,\ul{x}}\to \Op^\bL_{\cG,\CZ}/\bL'$$
preserve compactness. 

\medskip

However, the latter is obvious, since the above morphism goes between two smooth schemes.

\qed[\lemref{l:local to global Op ^! on !}]

\sssec{}

As an immediate corollary of \lemref{l:local to global Op ^! on !}, we obtain:

\begin{cor} \label{c:local to global Op _* on *}
The functor \eqref{e:local to global Op _* on *} admits a left adjoint,
to be denoted $\on{ev}_\CZ^{*,\IndCoh}$.
\end{cor}

\sssec{}

Note that we have a tautological commutative diagram 
\begin{equation} \label{e:IndCoh* vs QCohco on mer}
\CD
\IndCoh(\Op^{\mer,\on{glob}}_\cG)_\CZ @>{(\on{ev}_\CZ)^\IndCoh_*}>>  \IndCoh^*(\Op^{\mer}_\cG)_\CZ \\
@V{\Psi_{\Op^{\mer,\on{glob}}_{\cG,\CZ}}}VV @VV{\Psi_{\Op^\mer_{\cG,\CZ}}}V \\
\QCoh_{\on{co}}(\Op^{\mer,\on{glob}}_\cG)_\CZ @>>{(\on{ev}_\CZ)_*}>  \QCoh_{\on{co}}(\Op^{\mer}_\cG)_\CZ,
\endCD
\end{equation} 
see \secref{sss:Psi for IndCoh* indsch}. 

\medskip

Since the morphism $\on{ev}_\CZ$ is schematic, the functor 
$$(\on{ev}_\CZ)_*:\QCoh_{\on{co}}(\Op^{\mer,\on{glob}}_\cG)_\CZ \to \QCoh_{\on{co}}(\Op^{\mer}_\cG)_\CZ$$
admits a left adjoint, denoted $(\on{ev}_\CZ)^*$, see \secref{sss:functoriality of QCohco}.

\medskip

Passing to left adjoints along the horizontal arrows in \eqref{e:IndCoh* vs QCohco on mer}, we obtain a diagram
\begin{equation} \label{e:IndCoh* vs QCohco on mer adj}
\xy
(0,-20)*+{\QCoh_{\on{co}}(\Op^{\mer,\on{glob}}_\cG)_\CZ}="X";
(0,0)*+{\IndCoh(\Op^{\mer,\on{glob}}_\cG)_\CZ }="Y";
(50,-20)*+{\QCoh_{\on{co}}(\Op^{\mer}_\cG)_\CZ.}="X'";
(50,0)*+{\IndCoh^*(\Op^{\mer}_\cG)_\CZ }="Y'";
{\ar@{<-}^{\Psi_{\Op^{\mer,\on{glob}}_{\cG,\CZ}}} "X";"Y"};
{\ar@{<-}_{\Psi_{\Op^\mer_{\cG,\CZ}}} "X'";"Y'"};
{\ar@{<-}_{(\on{ev}_\CZ)^*} "X";"X'"};
{\ar@{<-}^-{(\on{ev}_\CZ)^{*,\IndCoh}} "Y";"Y'"};
{\ar@{<=} "Y";"X'"}
\endxy
\end{equation} 

We claim:

\begin{lem} \label{l:IndCoh* vs QCohco on mer adj}
The natural transformation in \eqref{e:IndCoh* vs QCohco on mer adj} is an isomorphism.
\end{lem}

\begin{proof}

With no restriction of generality, we can assume that $\CZ=X^I$; in particular, it is smooth. 

\medskip

We claim that the vertical arrows in \eqref{e:IndCoh* vs QCohco on mer} are in fact equivalences.
Indeed, we write
$$\IndCoh^*(\Op^{\mer}_\cG)_\CZ\simeq \underset{\bL\subset \bL_0}{\on{colim}}\, \IndCoh^*(\Op^\bL_\cG)_\CZ$$
and 
$$\IndCoh(\Op^{\mer,\on{glob}}_\cG)_\CZ\simeq \underset{\bL\subset \bL_0}{\on{colim}}\, \IndCoh(\Op^{\bL,\on{glob}}_\cG)_\CZ$$
where both colimits are formed with respect to the pushforward functors.

\medskip

Hence it enough to show that the functors
$$\Psi_{\Op^\bL_{\cG,\CZ}}:\IndCoh^*(\Op^\bL_\cG)_\CZ\to \QCoh(\Op^\bL_\cG)_\CZ$$
and
$$\Psi_{\Op^{\bL,\on{glob}}_{\cG,\CZ}}:\IndCoh(\Op^{\bL,\on{glob}}_\cG)_\CZ\to \QCoh(\Op^{\bL,\on{glob}}_\cG)_\CZ$$
are equivalences. 

\medskip

However, this follows from the fact that $\Op^{\bL,\on{glob}}_{\cG,\CZ}$ (resp., $\Op^\bL_{\cG,\CZ}$)
is smooth (resp., pro-smooth).

\end{proof}

\ssec{Interaction with self-duality}

\sssec{}

Recall now the functor 
$$\Theta_{\Op^\mer_{\cG,\CZ}}:\IndCoh^!(\Op^\mer_\cG)_\CZ\to \IndCoh^*(\Op^{\mer}_\cG)_\CZ,$$
see \secref{sss:Theta op mer}.

\sssec{}  \label{sss:Kost line}

Denote by $\fl_{\on{Kost}(\cG)}$ the (non-graded) line
$$\det (\Gamma(X,\fa(\cg)_{\omega_X})).$$

Set
$$\delta_G:=\dim(\Bun_G)=(g-1)\cdot \dim(G).$$

\sssec{}

We claim:

\begin{prop} \label{p:ev and duality} 
There exists a commutative diagram
$$
\CD
\IndCoh^!(\Op^\mer_\cG)_\CZ @>{\Theta_{\Op^\mer_{\cG,\CZ}}}>> \IndCoh^*(\Op^\mer_\cG)_\CZ \\
@V{(\on{ev}_\CZ)^!}VV @VV{(\on{ev}_\CZ)^{*,\IndCoh}}V \\
\IndCoh(\Op^{\mer,\on{glob}}_\cG)_\CZ @>>{-\otimes \fl_{\on{Kost}(\cG)}[-\delta_G]}> \IndCoh(\Op^{\mer,\on{glob}}_\cG)_\CZ.
\endCD
$$
\end{prop} 

The rest of this subsection is devoted to proof of \propref{p:ev and duality}.

\sssec{}

For expositional purposes, we will assume that $\CZ=\on{pt}$, so that $\CZ\to \Ran$
corresponds to $\ul{x}\in \Ran$. We will use the notation from \secref{sss:lattice notation}. 

\sssec{}

By the definition of the functor $\Theta_{\Op^\mer_{\cG,\ul{x}}}$, we need to establish an isomorphism
of the following two objects in $\IndCoh(\Op^{\mer,\on{glob}}_{\cG,\ul{x}})$:
$$\omega_{\Op^{\mer,\on{glob}}_{\cG,\ul{x}}}\otimes \fl_{\on{Kost}(\cG)}[-\delta_G]\simeq
(\on{ev}_{\ul{x}})^{*,\IndCoh}(\omega^{*,\on{fake}}_{\Op^\mer_{\cG,\ul{x}}}),$$
where 
$$\omega^{*,\on{fake}}_{\Op^\mer_{\cG,\ul{x}}}\in  \IndCoh^*(\Op^{\mer}_{\cG,\ul{x}})$$
is as in \secref{sss:omega fake mer}.

\sssec{}

In terms of the presentation 
$$\IndCoh^*(\Op^\mer_{\cG,\ul{x}})\simeq \underset{\bL\subset \bL_0}{\on{colim}}\, \IndCoh^*(\Op^\bL_{\cG,\ul{x}})
\overset{\Psi_{\Op^\bL_{\cG,\ul{x}}}}\simeq \QCoh(\Op^\bL_{\cG,\ul{x}}),$$
the object $\omega^{*,\on{fake}}_{\Op^\mer_{\cG,\ul{x}}}$ is, by construction, the colimit of the images of
$$\CO_{\Op^\bL_{\cG,\ul{x}}}\otimes 
\det(\bL/\bL_0)^{\otimes -1}[\dim(\bL/\bL_0)].$$

In terms of the presentation
$$\IndCoh(\Op^{\on{glob}}_{\cG,\ul{x}})\simeq 
\underset{\bL\subset \bL_0}{\on{colim}}\,\IndCoh^*(\Op^{\bL,\on{glob}}_{\cG,\ul{x}})
\overset{\Psi_{\Op^{\bL,\on{glob}}_{\cG,\ul{x}}}}\simeq \QCoh(\Op^{\bL,\on{glob}}_{\cG,\ul{x}}),$$
the object $\omega_{\Op^{\on{glob}}_{\cG,\ul{x}}}$ is, tautologically, the colimit of the images of
$$\Psi_{\Op^{\bL,\on{glob}}_{\cG,\ul{x}}}(\omega_{\Op^{\bL,\on{glob}}_{\cG,\ul{x}}}).$$

\sssec{}

The Cartesian diagram \eqref{e:cart diag L Op} gives rise to a commutative diagram
$$
\CD
\IndCoh(\Op^{\bL,\on{glob}}_{\cG,\ul{x}}) @>{(\on{ev}_{\ul{x}})^\IndCoh_*}>> \IndCoh^*(\Op^{\bL}_{\cG,\ul{x}})  \\
@A{!\text{-pullback}}AA @AA{!\text{-pullback}}A \\
\IndCoh(\Op^{\on{glob}}_{\cG,\ul{x}}) @>>{(\on{ev}_{\ul{x}})^\IndCoh_*}> \IndCoh^*(\Op^\mer_{\cG,\ul{x}}),
\endCD
$$
see \secref{sss:BeckChev placid ind}.

\medskip

Passing to the left adjoints, we obtain a commutative diagram 
\begin{equation} \label{e:ev and lattice diag}
\CD
\IndCoh(\Op^{\bL,\on{glob}}_{\cG,\ul{x}}) @<{(\on{ev}_{\ul{x}})^{*,\IndCoh}}<< \IndCoh^*(\Op^{\bL}_{\cG,\ul{x}})  \\
@V{*\text{-pushforward}}VV @VV{*\text{-pushforward}}V \\
\IndCoh(\Op^{\on{glob}}_{\cG,\ul{x}}) @<<{(\on{ev}_{\ul{x}})^{*,\IndCoh}}< \IndCoh^*(\Op^\mer_{\cG,\ul{x}}).
\endCD
\end{equation} 

As in \lemref{l:IndCoh* vs QCohco on mer adj}, we also have a commutative diagram
$$
\CD
\QCoh(\Op^{\bL,\on{glob}}_{\cG,\ul{x}}) @<{(\on{ev}_{\ul{x}})^*}<< \QCoh(\Op^{\bL}_{\cG,\ul{x}})   \\
@A{\Psi_{\Op^{\bL,\on{glob}}_{\cG,\ul{x}}}}AA @AA{\Psi_{\Op^\bL_{\cG,\ul{x}}}}A \\
\IndCoh(\Op^{\bL,\on{glob}}_{\cG,\ul{x}}) @<<{(\on{ev}_{\ul{x}})^{*,\IndCoh}}< \IndCoh^*(\Op^{\bL}_{\cG,\ul{x}}).
\endCD
$$

Hence, it is enough to construct a compatible collection of identifications 
$$\Psi_{\Op^{\bL,\on{glob}}_{\cG,\ul{x}}}(\omega_{\Op^{\bL,\on{glob}}_{\cG,\ul{x}}})
\otimes \fl_{\on{Kost}(\cG)}[-\delta_G]\simeq  \CO_{\Op^{\bL,\on{glob}}_{\cG,\ul{x}}}\otimes 
\det(\bL/\bL_0)^{\otimes -1}[\dim(\bL/\bL_0)],$$
taking place in $\QCoh(\Op^{\bL,\on{glob}}_{\cG,\ul{x}})$. 

\sssec{}

Note now that $\Op^{\bL,\on{glob}}_{\cG,\ul{x}}$ is an affine space with respect to
$$\bL^{\on{glob}}:=\bV^{\on{glob}}\cap \bL.$$

Hence, 
$$\Psi_{\Op^{\bL,\on{glob}}_{\cG,\ul{x}}}(\omega_{\Op^{\bL,\on{glob}}_{\cG,\ul{x}}})\simeq
\CO_{\Op^{\bL,\on{glob}}_{\cG,\ul{x}}}\otimes \det(\bL^{\on{glob}})^{\otimes -1}[\dim(\bL^{\on{glob}})].$$

\sssec{}

Thus, it remains to establish a compatible collection of isomorphisms between the lines
\begin{equation} \label{e:identify lattice lines}
\det(\bL^{\on{glob}})^{\otimes -1} \otimes \fl_{\on{Kost}(\cG)}[\dim(\bL^{\on{glob}})-\delta_G] \simeq
\det(\bL/\bL_0)^{\otimes -1}[\dim(\bL/\bL_0)].
\end{equation}

However, this follows from the fact that
$$\fl_{\on{Kost}(\cG)}\simeq \det(\bV^{\on{glob}}\cap \bL_0) \text{ and }
\delta_G=\dim(\Gamma(X,\fa(\cg)_{\omega_X}))=\dim(\bV^{\on{glob}}\cap \bL_0).$$

\ssec{Ind-coherent sheaves on local vs. global \emph{monodromy-free} opers}

\sssec{}

For $\CZ\to \Ran$ recall the (relative ind-scheme) $\Op^{\mf,\on{glob}}_{\cG,\CZ}$, which fits into the Cartesian square 
\begin{equation} \label{e:loc vs global mf opers Z}
\CD
\Op^{\mf,\on{glob}}_{\cG,\CZ} @>{\on{ev}_\CZ}>>  \Op^\mf_{\cG,\CZ} \\
@V{\iota^{\mf,\on{glob}}}VV @VV{\iota^\mf}V \\
\Op^{\on{mer,glob}}_{\cG,\CZ} @>>{\on{ev}_\CZ}>  \Op^\mer_{\cG,\CZ} 
\endCD
\end{equation}

\medskip

Consider the morphism: 
$$\Op^{\mf,\on{glob}}_{\cG,\CZ} \overset{\on{ev}_\CZ}\to \Op^\mf_{\cG,\CZ}$$
and the resulting pair of mutually dual functors
\begin{equation} \label{e:local to global Op _* on * mf}
(\on{ev}_\CZ)^\IndCoh_*:\IndCoh(\Op^{\mf,\on{glob}}_\cG)_\CZ\to \IndCoh^*(\Op^\mf_\cG)_\CZ
\end{equation} 
and 
\begin{equation} \label{e:local to global Op ^! on ! mf}
(\on{ev}_\CZ)^!:\IndCoh^!(\Op^\mf_\cG)_\CZ\to \IndCoh(\Op^{\mf,\on{glob}}_\cG)_\CZ.
\end{equation} 

\sssec{}

We claim:

\begin{lem} \label{l:local to global Op ^! on ! mf}
The functor \eqref{e:local to global Op ^! on ! mf} preserves compactness.
\end{lem} 

\begin{proof}

From \eqref{e:loc vs global mf opers Z} we obtain a commutative square
\begin{equation} \label{e:ev and iota}
\CD
\IndCoh(\Op^{\mf,\on{glob}}_\cG)_\CZ @<{\on{ev}_\CZ^!}<< \IndCoh^!(\Op^\mf_\cG)_\CZ  \\
@V{(\iota^{\mf,\on{glob}})^\IndCoh_*}VV @VV{(\iota^\mf)^\IndCoh_*}V  \\
\IndCoh(\Op^{\mer,\on{glob}}_\cG)_\CZ @<<{\on{ev}_\CZ^!}< \IndCoh^!(\Op^\mer_\cG)_\CZ,
\endCD
\end{equation} 
see \secref{sss:BeckChev placid ind}.

\medskip

As in \propref{p:mon-free to all}(b), one shows that an object in $\IndCoh(\Op^{\mf,\on{glob}}_\cG)_\CZ$ is compact
if and only if its image in $\IndCoh(\Op^{\mer,\on{glob}}_\cG)_\CZ$ under $(\iota)^{\on{glob}})^\IndCoh_*$
is compact. 

\medskip

Hence, it suffices to show that the clockwise circuit in \eqref{e:ev and iota} preserves compactness.

\medskip

For the functor $(\iota^\mf)^\IndCoh_*$ this is evident (since the morphism $\iota^\mf$ is of finite presentation).
For the bottom horizontal arrow in \eqref{e:ev and iota} this follows from 
\lemref{l:local to global Op ^! on !}.

\end{proof} 

\sssec{}

As a corollary of \lemref{l:local to global Op ^! on ! mf} we obtain:

\begin{cor} \label{c:local to global Op _* on * mf}
The functor \eqref{e:local to global Op _* on * mf} admits a left adjoint,
to be denoted $\on{ev}_\CZ^{*,\IndCoh}$.
\end{cor}

\sssec{}

The fact that \eqref{e:loc vs global mf opers Z} is Cartesian implies that the diagram
\begin{equation} \label{e:ev and iota 1}
\CD
\IndCoh(\Op^{\mf,\on{glob}}_\cG)_\CZ @>{(\on{ev}_\CZ)^\IndCoh_*}>> \IndCoh^*(\Op^\mf_\cG)_\CZ  \\
@A{(\iota^{\mf,\on{glob}})^!}AA @AA{(\iota^\mf)^!}A  \\
\IndCoh(\Op^{\mer,\on{glob}}_\cG)_\CZ @>>{(\on{ev}_\CZ)^\IndCoh_*}> \IndCoh^*(\Op^\mer_\cG)_\CZ
\endCD
\end{equation}
commutes, see \secref{sss:BeckChev placid ind}.

\medskip

By passing to left adjoints along all arrows in \eqref{e:ev and iota 1} we obtain that the diagram 
\begin{equation} \label{e:ev and iota 2}
\CD
\IndCoh(\Op^{\mf,\on{glob}}_\cG)_\CZ @<{(\on{ev}_\CZ)^{*,\IndCoh}}<< \IndCoh^*(\Op^\mf_\cG)_\CZ  \\
@V{(\iota^{\mf,\on{glob}})^\IndCoh_*}VV @VV{(\iota^\mf)^\IndCoh_*}V  \\
\IndCoh(\Op^{\mer,\on{glob}}_\cG)_\CZ @<<{(\on{ev}_\CZ)^{*,\IndCoh}}< \IndCoh^*(\Op^\mer_\cG)_\CZ
\endCD
\end{equation}
commutes as well. 

\medskip

However, passing left adjoints only along the horizontal arrows in \eqref{e:ev and iota 1} we obtain a diagram 
\begin{equation} \label{e:ev and iota 3}
\xy
(0,0)*+{\IndCoh(\Op^{\mf,\on{glob}}_\cG)_\CZ}="X";
(50,0)*+{\IndCoh^*(\Op^\mf_\cG)_\CZ}="Y";
(0,-20)*+{\IndCoh(\Op^{\mer,\on{glob}}_\cG)_\CZ }="Z";
(50,-20)*+{\IndCoh^*(\Op^\mer_\cG)_\CZ.}="W";
{\ar@{<-}^{(\on{ev}_\CZ)^{*,\IndCoh}} "X";"Y"};
{\ar@{<-}_{(\on{ev}_\CZ)^{*,\IndCoh}} "Z";"W"};
{\ar@{<-}_{(\iota^{\mf,\on{glob}})^!} "X";"Z"};
{\ar@{<-}^{(\iota^\mf)^!} "Y";"W"};
{\ar@{<=} "Z";"Y"}
\endxy
\end{equation} 

We claim:

\begin{lem} \label{l:ev and iota}
The natural transformation in \eqref{e:ev and iota 3} is an isomorphism.
\end{lem}

\begin{proof}

For expositional purposes we will assume that $\CZ=\on{pt}$, so that $\CZ\to \Ran$
corresponds to $\ul{x}\in \Ran$. We will use the notations from \secref{sss:lattice notation}. 

\medskip

For $\bL\supset \bL_0$ denote
$$\Op^{\bL\cap\mf}_{\cG,\ul{x}}:=\Op^{\mf}_{\cG,\ul{x}}\underset{\Op^\mer_{\cG,\ul{x}}}\times \Op^\bL_{\cG,\ul{x}},$$
$$\Op^{\bL\cap\mf,\on{glob}}_{\cG,\ul{x}}:=\Op^{\mf,\on{glob}}_{\cG,\ul{x}}\underset{\Op^\mer_{\cG,\ul{x}}}\times \Op^\bL_{\cG,\ul{x}},$$

Using \eqref{e:ev and lattice diag} and a similar diagram for ``$\mer$" replaced my ``$\mf$", it suffices to show that the 
natural transformation in the diagram 
$$
\xy
(0,0)*+{\IndCoh(\Op^{\bL\cap\mf,\on{glob}}_{\cG,\ul{x}})}="X";
(70,0)*+{\IndCoh^*(\Op^{\bL\cap\mf}_{\cG,\ul{x}})}="Y";
(0,-20)*+{\IndCoh(\Op^{\bL,\on{glob}}_{\cG,\ul{x}}) }="Z";
(70,-20)*+{\IndCoh^*(\Op^\bL_{\cG,\ul{x}}).}="W";
{\ar@{<-}^-{(\on{ev}_{\ul{x}})^{*,\IndCoh}} "X";"Y"};
{\ar@{<-}_-{(\on{ev}_{\ul{x}})^{*,\IndCoh}} "Z";"W"};
{\ar@{<-}_{(\iota^{\mf,\on{glob}})^!} "X";"Z"};
{\ar@{<-}^{(\iota^\mf)^!} "Y";"W"};
{\ar@{<=} "Z";"Y"}
\endxy
$$
is an isomorphism. 

\medskip

According to \secref{sss:Op mon-free expl}, for a small enough lattice $\bL'\subset \bL$, we have a well-defined action of $\bL'$
on $\Op^{\bL\cap\mf}_{\cG,\ul{x}}$ by translations, and the quotient $\Op^{\bL\cap\mf}_{\cG,\ul{x}}/\bL'$ is a prestack locally 
almost of finite type. 

\medskip

We have a commutative diagram
$$
\CD
\IndCoh^*(\Op^{\bL\cap\mf}_{\cG,\ul{x}}) @<{*\on{-pullback}}<< \IndCoh(\Op^{\bL\cap\mf}_{\cG,\ul{x}}/\bL')  \\
@A{(\iota^\mf)^!}AA @AA{(\iota^\mf/\bL')^!}A \\
\IndCoh^*(\Op^\bL_{\cG,\ul{x}}) @<{*\on{-pullback}}<< \IndCoh^*(\Op^\bL_{\cG,\ul{x}}/\bL'). 
\endCD
$$

Hence, it suffices to show that the natural transformation in the diagram 
$$
\xy
(0,0)*+{\IndCoh(\Op^{\bL\cap\mf,\on{glob}}_{\cG,\ul{x}})}="X";
(70,0)*+{\IndCoh^*(\Op^{\bL\cap\mf}_{\cG,\ul{x}}/\bL')}="Y";
(0,-20)*+{\IndCoh(\Op^{\bL,\on{glob}}_{\cG,\ul{x}}) }="Z";
(70,-20)*+{\IndCoh^*(\Op^\bL_{\cG,\ul{x}}/\bL').}="W";
{\ar@{<-}^-{*\on{-pullback}} "X";"Y"};
{\ar@{<-}_-{*\on{-pullback}} "Z";"W"};
{\ar@{<-}_{(\iota^{\mf,\on{glob}})^!} "X";"Z"};
{\ar@{<-}^{(\iota^\mf/\bL')^!} "Y";"W"};
{\ar@{<=} "Z";"Y"}
\endxy
$$
is an isomorphism. 

\medskip

However, this follows from the fact that the diagram
$$
\CD
\Op^{\bL\cap\mf,\on{glob}}_{\cG,\ul{x}} @>>> \Op^{\bL\cap\mf}_{\cG,\ul{x}}/\bL' \\
@VVV @VVV \\
\Op^{\bL,\on{glob}}_{\cG,\ul{x}} @>>> \Op^\bL_{\cG,\ul{x}}/\bL'
\endCD
$$
is Cartesian, combined with the fact that $\Op^\bL_{\cG,\ul{x}}/\bL'$ is a smooth scheme.

\end{proof} 

\sssec{}

Recall the functor 
$$\Theta_{\Op^\mf_{\cG,\CZ}}:\IndCoh^!(\Op^\mf_\cG)_\CZ\to \IndCoh^*(\Op^\mf_\cG)_\CZ,$$
see \secref{sss:omega fake mf}.

\medskip

We claim: 

\begin{prop} \label{p:ev and duality mf} 
There exists a commutative diagram
$$
\CD
\IndCoh^!(\Op^\mf_\cG)_\CZ @>{\Theta_{\Op^\mf_{\cG,\CZ}}}>> \IndCoh^*(\Op^\mf_\cG)_\CZ \\
@V{(\on{ev}_\CZ)^!}VV @VV{(\on{ev}_\CZ)^{*,\IndCoh}}V \\
\IndCoh(\Op^{\mf,\on{glob}}_\cG)_\CZ @>>{-\otimes \fl_{\on{Kost}(\cG)}[-\delta_G]}> \IndCoh(\Op^{\mf,\on{glob}}_\cG)_\CZ.
\endCD
$$
\end{prop} 

\begin{proof}

Both circuits of the diagram are $\IndCoh^!(\Op^\mf_\cG)_\CZ$-linear functors. Hence, it suffices to identify
the objects that correspond to the image of the unit. 

\medskip

I.e., we wish to identify 
\begin{equation} \label{e:ev and duality mf}
\omega_{\Op^{\mf,\on{glob}}_{\cG,\CZ}}\otimes \fl_{\on{Kost}(\cG)}[-\delta_G]\simeq
(\on{ev}_\CZ)^{*,\IndCoh}(\omega^{*,\on{fake}}_{\Op^\mf_{\cG,\CZ}}).
\end{equation} 

We start with 
$$\omega_{\Op^{\mer,\on{glob}}_{\cG,\CZ}}\otimes \fl_{\on{Kost}(\cG)}[-\delta_G]\simeq
(\on{ev}_\CZ)^{*,\IndCoh}(\omega^{*,\on{fake}}_{\Op^\mer_{\cG,\CZ}}),$$
given by \propref{p:ev and duality} and apply the functor $(\iota^{\mf,\on{glob}})^!$.

\medskip

The left-hand side gives the left-hand side of \eqref{e:ev and duality mf}. The right-hand side 
gives the right-hand side of \eqref{e:ev and duality mf} thanks to \lemref{l:ev and iota}. 

\end{proof} 

\ssec{Two versions of the spectral Poincar\'e functor}

\sssec{} \label{sss:Poinc spec !}

For $\CZ\to \Ran$, we define the spectral !-Poincar\'e functor
$$\Poinc^{\on{spec}}_{\cG,!,\CZ}:\IndCoh^!(\Op^\mf_\cG)_\CZ\to \IndCoh(\LS_\cG)\otimes \Dmod(\CZ)$$
as
$$\IndCoh^!(\Op^\mf_\cG)_\CZ\overset{\on{ev}_\CZ^!}\longrightarrow
\IndCoh(\Op^{\mf,\on{glob}}_\cG)_\CZ\overset{(\fr^{\on{glob}})^\IndCoh_*}\to \IndCoh(\LS_\cG)\otimes \Dmod(\CZ).$$

\sssec{} \label{sss:Poinc spec *}

We define the spectral *-Poincar\'e functor
$$\Poinc^{\on{spec}}_{\cG,*,\CZ}:\IndCoh^*(\Op^\mf_\cG)_\CZ\to \IndCoh(\LS_\cG)\otimes \Dmod(\CZ)$$
as
$$\IndCoh^*(\Op^\mf_\cG)_\CZ\overset{\on{ev}_\CZ^{*,\IndCoh}}\longrightarrow
\IndCoh(\Op^{\mf,\on{glob}}_\cG)_\CZ\overset{(\fr^{\on{glob}})^\IndCoh_*}\to \IndCoh(\LS_\cG)\otimes \Dmod(\CZ).$$

\sssec{}

Both 
\begin{equation}  \label{e:Poinc spec ! and *}
\CZ \rightsquigarrow \Poinc^{\on{spec}}_{\cG,!,\CZ} \text{ and } 
\CZ \rightsquigarrow \Poinc^{\on{spec}}_{\cG,*,\CZ}
\end{equation} 
are naturally local-to-global functors in the sense of \secref{sss:local to global functor abs},
to be denoted 
$$\ul{\Poinc}^{\on{spec}}_{\cG,!} \text{ and } \ul{\Poinc}^{\on{spec}}_{\cG,*},$$
respectively. 

\sssec{} \label{sss:untl Poinc}

Furthermore, the assignments \eqref{e:Poinc spec ! and *} have natural \emph{unital} structures,
in the sense of \secref{sss:strictly unital}. Let us spell it out explicitly for the *-version (!-version is analogous).

\medskip

The local unital structure on the source crystal of categories , i.e., $\IndCoh^*(\Op^\mf_\cG)$ assigns to
$$(\ul{x}\subseteq \ul{x}')\in \Ran$$
the functor
$$\IndCoh^*(\Op^\mf_{\cG,\ul{x}})\to \IndCoh^*(\Op^\mf_{\cG,\ul{x}'})$$ given by *-pull followed by *-push along the diagram
\begin{equation} \label{e:untl Poinc 1}
\xy
(0,0)*+{\Op_\cG(\cD_{\ul{x}}-\ul{x})\underset{\LS_\cG(\cD_{\ul{x}}-\ul{x})}\times \LS_\cG(\cD_{\ul{x}})}="A";
(60,0)*+{\Op_\cG(\cD_{\ul{x}'}-\ul{x}')\underset{\LS_\cG(\cD_{\ul{x}'}-\ul{x}')}\times \LS_\cG(\cD_{\ul{x}'})}="B";
(30,30)*+{\Op_\cG(\cD_{\ul{x}'}-\ul{x})\underset{\LS_\cG(\cD_{\ul{x}'}-\ul{x})}\times \LS_\cG(\cD_{\ul{x}'})}="C";
(0,-20)*+{\Op^\mf_{\cG,\ul{x}}}="D";
(60,-20)*+{\Op^\mf_{\cG,\ul{x}'},}="E";
{\ar@{->} "C";"A"};
{\ar@{->} "C";"B"};
{\ar@{=} "A";"D"};
{\ar@{=} "B";"E"};
\endxy
\end{equation} 
in which the slanted arrows are given by restriction along the inclusions
$$(\cD_{\ul{x}}-\ul{x})\to (\cD_{\ul{x}'}-\ul{x}) \leftarrow (\cD_{\ul{x}'}-\ul{x}'),$$
respectively. 

\sssec{}

Consider the diagram 
$$
\CD
\Op^\mf_{\cG,\ul{x}'} @<{\on{ev}_{\ul{x}'}}<< \Op^{\mf,\on{glob}}_{\cG,\ul{x}'} @>{\fr^{\on{glob}}_{\ul{x'}}}>> \LS_\cG \\
@AAA @AAA @AA{\on{id}}A \\
\Op_\cG(\cD_{\ul{x}'}-\ul{x})\underset{\LS_\cG(\cD_{\ul{x}}-\ul{x})}\times \LS_\cG(\cD_{\ul{x}}) @<<< 
\Op^{\mf,\on{glob}}_{\cG,\ul{x}} @>{\fr^{\on{glob}}_{\ul{x}}}>> \LS_\cG \\
@VVV @VV{\on{ev}_{\ul{x}}}V \\
\Op^\mf_{\cG,\ul{x}} @<{\on{id}}<< \Op^\mf_{\cG,\ul{x}},
\endCD
$$
in which the square
$$
\CD
\Op^\mf_{\cG,\ul{x}'} @<{\on{ev}_{\ul{x}'}}<< \Op^{\mf,\on{glob}}_{\cG,\ul{x}'} \\
@AAA @AAA  \\
\Op_\cG(\cD_{\ul{x}'}-\ul{x})\underset{\LS_\cG(\cD_{\ul{x}}-\ul{x})}\times \LS_\cG(\cD_{\ul{x}}) @<<< 
\Op^{\mf,\on{glob}}_{\cG,\ul{x}} 
\endCD
$$
is Cartesian.

\medskip

By \secref{sss:BeckChev placid ind}, we have a commutative diagram
$$
\CD
\IndCoh^*(\Op^\mf_{\cG,\ul{x}'}) @>{*\on{-pushforward}}>> \IndCoh(\Op^{\mf,\on{glob}}_{\cG,\ul{x}'}) \\
@V{!\on{-pullback}}VV @VV{!\on{-pullback}}V \\ 
\IndCoh^*(\Op_\cG(\cD_{\ul{x}'}-\ul{x})\underset{\LS_\cG(\cD_{\ul{x}}-\ul{x})}\times \LS_\cG(\cD_{\ul{x}})) @>{*\on{-pushforward}}>>
\IndCoh(\Op^{\mf,\on{glob}}_{\cG,\ul{x}}).
\endCD
$$

Passing to left adjoints, we obtain a commutative diagram
$$
\CD
\IndCoh^*(\Op^\mf_{\cG,\ul{x}'}) @>{*\on{-pullback}}>> \IndCoh(\Op^{\mf,\on{glob}}_{\cG,\ul{x}'}) \\
@A{*\on{-pushforward}}AA @AA{*\on{-pushforward}}A \\ 
\IndCoh^*(\Op_\cG(\cD_{\ul{x}'}-\ul{x})\underset{\LS_\cG(\cD_{\ul{x}}-\ul{x})}\times \LS_\cG(\cD_{\ul{x}})) @>{*\on{-pullback}}>>
\IndCoh(\Op^{\mf,\on{glob}}_{\cG,\ul{x}}).
\endCD
$$

\medskip

Now, the unital structure on $\ul{\Poinc}^{\on{spec}}_{\cG,*}$ is encoded by the following diagram:
$$
\CD
\IndCoh^*(\Op^\mf_{\cG,\ul{x}'}) @>{\on{ev}^{*,\IndCoh}_{\ul{x}'}}>> 
\IndCoh(\Op^{\mf,\on{glob}}_{\cG,\ul{x}'}) @>{(\fr^{\on{glob}}_{\ul{x'}})^\IndCoh_*}>> \IndCoh(\LS_\cG) \\
@A{*\on{-pushforward}}AA @A{*\on{-pushforward}}AA @AA{\on{id}}A \\
\IndCoh^*(\Op_\cG(\cD_{\ul{x}'}-\ul{x})\underset{\LS_\cG(\cD_{\ul{x}'}-\ul{x})}\times \LS_\cG(\cD_{\ul{x}'})) @>{*\on{-pullback}}>>
\IndCoh(\Op^{\mf,\on{glob}}_{\cG,\ul{x}}) @>{(\fr^{\on{glob}}_{\ul{x}})^\IndCoh_*}>> \IndCoh(\LS_\cG) \\
@A{*\on{-pullback}}AA @AA{\on{ev}^*_{\ul{x}}}A \\
\IndCoh^*(\Op^\mf_{\cG,\ul{x}}) @>{\on{Id}}>> \IndCoh^*(\Op^\mf_{\cG,\ul{x}}).
\endCD
$$

\sssec{}

Note that from \propref{p:ev and duality mf} we obtain:

\begin{thm} \label{t:Poinc spec * vs !}
There is an isomorphism of local-to-global functors
$$\ul{\Poinc}^{\on{spec}}_{\cG,!}\otimes \fl_{\on{Kost}(\cG)}[-\delta_G]\simeq 
\ul{\Poinc}^{\on{spec}}_{\cG,*}\circ \Theta_{\Op^\mf_\cG}.$$
\end{thm} 

\begin{rem}
The reason that we discuss both $\ul{\Poinc}^{\on{spec}}_{\cG,!}$ and $\ul{\Poinc}^{\on{spec}}_{\cG,*}$
(despite the fact that, thanks to \thmref{t:Poinc spec * vs !}, they are easily expressible one through another) 
is that the !-version is naturally compatible with Eisenstein series (which we will exploit in the sequel to
this paper), and the *-version is naturally compatible with the functor
$$\Gamma^\IndCoh(\LS_\cG,-):\IndCoh(\LS_\cG)\to \Vect,$$
which we will use in the next section.
\end{rem} 

\ssec{Action of the spectral spherical category and temperedness}

In this subsection we will work with $\ul{\Poinc}^{\on{spec}}_{\cG,*}$, but a parallel discussion is
applicable to $\ul{\Poinc}^{\on{spec}}_{\cG,!}$.

\sssec{} \label{sss:Sph spec action x}

For a fixed $\ul{x}\in \Ran$ we have a naturally defined action of $\Sph^{\on{spec}}_{\cG,\ul{x}}$ on $\IndCoh(\LS_\cG)$. 
Namely, it is given by *-pull followed by *-push along the following diagram
$$
\CD
\LS_\cG  @<<< \on{Hecke}_{\cG,\ul{x}}^{\on{spec,glob}} @>>> \LS_\cG  \\
@V{\on{ev}_{\ul{x}}}VV @V{\on{ev}_{\ul{x}}}VV @VV{\on{ev}_{\ul{x}}}V \\
\LS^\reg_{\cG,\ul{x}} @<<< \on{Hecke}_{\cG,\ul{x}}^{\on{spec,loc}}  @>>> \LS^\reg_{\cG,\ul{x}},
\endCD
$$
in which both squares are Cartesian.

\begin{rem} 
In \secref{sss:Hecke spec action Ran} we will consider a Ran version of this action. This involves
some technical difficulties, inherent to the definition of $\Sph^{\on{spec}}_\cG$ as a factorization 
category, see \secref{sss:Sph spec tech}.
\end{rem} 

\sssec{}

We have a natural action of $\Sph^{\on{spec}}_{\cG,\ul{x}}$ on $\IndCoh^*(\Op^\mf_{\cG,\Ran_{\ul{x}}})$ given by pull-push along the 
following diagram
$$
\CD
\Op^\mf_{\cG,\ul{x}'} @<<< \on{Hecke}_{\cG,\ul{x}}^{\on{spec},\Op^\mf_{\cG,\ul{x}'}} @>>> \Op^\mf_{\cG,\ul{x}'}  \\
@VVV @VVV @VVV \\
\LS^\reg_{\cG,\ul{x}} @<<< \on{Hecke}_{\cG,\ul{x}}^{\on{spec,loc}}  @>>> \LS^\reg_{\cG,\ul{x}},
\endCD
$$
in which both squares are Cartesian, where:

\begin{itemize}

\item $\ul{x}\subseteq \ul{x}'$; 

\item $\on{Hecke}_{\cG,\ul{x}}^{\on{spec},\Op^\mf_{\cG,\ul{x}'}}:=
\Op^\mer_{\cG,\ul{x}'}\underset{\LS^\mer_{\cG,\ul{x}'}}\times 
(\LS^\reg_{\cG,\ul{x}'}\underset{\LS_{\cG,\ul{x}\subseteq \ul{x}'}^{\mer\rightsquigarrow\reg}}\times \LS^\reg_{\cG,\ul{x}'})$;

\item $\LS_{\cG,\ul{x}\subseteq \ul{x}'}^{\mer\rightsquigarrow\reg}:=\LS_\cG(\cD_{\ul{x}'}-\ul{x})$.

\end{itemize}

\sssec{}

Consider the functor
$$\Poinc^{\on{spec}}_{\cG,*,\Ran_{\ul{x}}}:\IndCoh^*(\Op^\mf_{\cG,\Ran_{\ul{x}}})\to \IndCoh(\LS_\cG)\otimes \Dmod(\Ran_{\ul{x}}).$$

The following results by unwinding the constructions:

\begin{lem} \label{l:Poinc intertwines Sph}
The functor $\Poinc^{\on{spec}}_{\cG,*,\Ran_{\ul{x}}}$ intertwines the actions of $\Sph^{\on{spec}}_{\cG,\ul{x}}$ on
the two sides.
\end{lem}

\begin{rem}
One can define an action of $\Sph^{\on{spec}}_\cG$ on $\IndCoh^*(\Op^\mf_\cG)$ as a factorization category.
Furthermore, one can show that that the functor $\ul{\Poinc}^{\on{spec}}_{\cG,*}$ is compatible with the action of 
$\IndCoh^*(\Op^\mf_\cG)$, where the latter is thought of as a crystal of monoidal categories over $\Ran$.

\medskip

Moreover, the above action and compatibility are in turn compatible with the unital structures. This will
be performed in \secref{ss:Sph spec acts on mon-free}. 

\end{rem} 

\sssec{} \label{sss:QCoh as temp}

Recall (see \cite[Sect. 12.8.2]{AG}) that the subcategory
$$\QCoh(\LS_\cG) \subset \IndCoh(\LS_\cG)$$
can be singled out by the \emph{temperedness} condition:

\medskip

Namely, it is the maximal subcategory on which for some/any $\ul{x}\in \Ran$, 
the action of $\Sph^{\on{spec}}_{\cG,\ul{x}}$ on $\IndCoh(\LS_\cG)$ factors via 
$$\Sph^{\on{spec}}_{\cG,\ul{x}}\twoheadrightarrow \Sph^{\on{spec}}_{\cG,\on{temp},\ul{x}},$$
see \secref{sss:temp spec}. 

\sssec{}

A basic property of the spectral Poincar\'e functor is the following:

\begin{prop} \label{p:Poinc to QCoh}
The essential image of the functor 
$$\Poinc^{\on{spec}}_{\cG,*}:\IndCoh^*(\Op^\mf_\cG)_\Ran \to \IndCoh(\LS_\cG)$$
lies in
$$\QCoh(\LS_\cG) \subset \IndCoh(\LS_\cG).$$
\end{prop}

\begin{proof}

Choose some/any $\ul{x}\in \Ran$. By the unital property of $\ul{\Poinc}^{\on{spec}}_{\cG,*}$, it suffices to show that
the functor
$$\Poinc^{\on{spec}}_{\cG,*,\Ran_{\ul{x}}}:\IndCoh^*(\Op^\mf_{\cG,\Ran_{\ul{x}}})\to \IndCoh(\LS_\cG)\otimes \Dmod(\Ran_{\ul{x}})$$
takes values in 
$$\QCoh(\LS_\cG)\otimes \Dmod(\Ran_{\ul{x}})\subset \IndCoh(\LS_\cG)\otimes \Dmod(\Ran_{\ul{x}}).$$

By \secref{sss:QCoh as temp} and \lemref{l:Poinc intertwines Sph}, it suffices to show that the $\Sph^{\on{spec}}_{\cG,\ul{x}}$-action 
on $\IndCoh^*(\Op^\mf_{\cG,\Ran_{\ul{x}}})$ factors via $\Sph^{\on{spec}}_{\cG,\on{temp},\ul{x}}$.

\medskip

The latter assertion can be checked strata-wise, so it is enough to show that for a fixed $\ul{x}\subseteq \ul{x}'$, the action of
$\Sph^{\on{spec}}_{\cG,\ul{x}}$ on $\IndCoh^*(\Op^\mf_{\cG,\Ran_{\ul{x}'}})$ factors via $\Sph^{\on{spec}}_{\cG,\on{temp},\ul{x}}$.

\medskip

Write $\ul{x}'=\ul{x}\sqcup \ul{x}''$. In terms of the factorization 
$$\IndCoh^*(\Op^\mf_{\cG,\Ran_{\ul{x}'}})\simeq \IndCoh^*(\Op^\mf_{\cG,\Ran_{\ul{x}}})\otimes \IndCoh^*(\Op^\mf_{\cG,\Ran_{\ul{x}''}}),$$
the action of  $\Sph^{\on{spec}}_{\cG,\ul{x}}$ on $\IndCoh^*(\Op^\mf_{\cG,\Ran_{\ul{x}'}})$ is via the first factor.

\medskip

The required assertion follows now from \propref{p:opers is temp}.

\end{proof} 

\begin{rem}
Note that thanks to \thmref{t:Poinc spec * vs !}, we obtain that an assertion parallel to \propref{p:Poinc to QCoh}
holds for the functor $\Poinc^{\on{spec}}_{\cG,!}$. (Alternatively, one can prove it by the same
argument.) 
\end{rem}

\ssec{The spectral localization and global sections functors}

\sssec{} \label{sss:spec loc}

The spectral localization functor $\ul\Loc^{\on{spec}}_\cG$, i.e., the collectiion functors
$$\Loc_{\cG,\CZ}^{\on{spec}}:\Rep(\cG)_\CZ\to \QCoh(\LS_\cG)\otimes \Dmod(\CZ)$$
for $\CZ\to \Ran$, is defined as pullback along 
\begin{equation} \label{e:spectral localization}
\LS_\cG\times \CZ\to \LS^{\on{reg}}_{\cG,\CZ},
\end{equation} 
where we identify
$$\Rep(\cG)_\CZ\simeq \QCoh(\LS^{\on{reg}}_\cG)_\CZ.$$

\medskip

The functor $\ul\Loc_\cG^{\on{spec}}$ possesses a natural unital structure (see \secref{sss:strictly unital} for what this means).

\sssec{}

The functor 
$$\Loc_\cG^{\on{spec}}:\Rep(\cG)_\Ran\to \QCoh(\LS_\cG)$$
admits a right adjoint, denoted
$$\Gamma^{\on{spec}}_\cG: \QCoh(\LS_\cG)\to \Rep(\cG)_\Ran.$$

Explicitly, for a given $\ul{x}\in \Ran$, the corresponding functor
$$\Gamma^{\on{spec}}_{\cG,\ul{x}}:\QCoh(\LS_\cG)\to \Rep(\cG)_{\ul{x}}$$
is given by *-direct image along
$$\LS_\cG\to \LS^\reg_{\cG,\ul{x}}.$$

\sssec{}

Note also that the categories $\QCoh(\LS_\cG)$ and $\Rep(\cG)_\Ran$ are both canonically self-dual, and 
with respect to these dualities, we have
$$(\Loc_\cG^{\on{spec}})^\vee\simeq \Gamma^{\on{spec}}_\cG.$$

\sssec{} \label{sss:Gamma spec IndCoh}

By a slight abuse of notation we will denote by the same symbol $\Loc_\cG^{\on{spec}}$ the composite functor
$$\Rep(\cG)_\Ran\overset{\Loc_\cG^{\on{spec}}}\longrightarrow \QCoh(\LS_\cG)
\overset{\Xi_{\LS_\cG}}\hookrightarrow \IndCoh(\LS_\cG).$$

We will denote by $\Gamma^{\on{spec},\IndCoh}_\cG$ the functor
$$\IndCoh(\LS_\cG)\overset{\Psi_{\LS_\cG}}\twoheadrightarrow \QCoh(\LS_\cG)\overset{\Gamma^{\on{spec}}_\cG}\longrightarrow
\Rep(\cG)_\Ran.$$

The functors
$$\Loc_\cG^{\on{spec}}:\Rep(\cG)_\Ran\leftrightarrows \IndCoh(\LS_\cG):\Gamma^{\on{spec},\IndCoh}_\cG$$
also form an adjoint pair.

\sssec{}

Note that the category $\IndCoh(\LS_\cG)$ is also self-dual by means of \emph{Serre} duality.
Under this duality and the standard self-duality of $\QCoh(\LS_\cG)$, we have
$$\Psi_{\LS_\cG}^\vee \simeq \Upsilon_{\LS_\cG}.$$

However, note that $\LS_\cG$ is quasi-smooth and Calabi-Yau:
$$\omega_{\LS_\cG}=\CO_{\LS_\cG}\otimes \det(\Lie(Z_G))^{\otimes (2-2g)}[2\delta_G].$$

Hence, we have
$$\Upsilon_{\LS_\cG}\simeq \Xi_{\LS_\cG}\otimes \det(\Lie(Z_G))^{\otimes (2-2g)}[2\delta_G].$$

\sssec{}

Hence, we obtain that with respect to the self-duality of $\Rep(\cG)_\Ran$ and the Serre duality of $\IndCoh(\LS_\cG)$,
we have
$$(\Loc_\cG^{\on{spec}})^\vee \simeq \Gamma^{\on{spec},\IndCoh}_\cG\otimes \det(\Lie(Z_G))^{\otimes (2-2g)}[2\delta_G].$$

\ssec{Composing spectral Poincar\'e and global sections functors}

\sssec{}

Our current goal is to study the composite functor
\begin{equation} \label{e:spec Poinc and Gamma}
\IndCoh^*(\Op^{\on{mon-free}}_\cG)_\Ran \overset{\Poinc^{\on{spec}}_{\cG,*}}\longrightarrow 
\QCoh(\LS_\cG)\overset{\Gamma^{\on{spec}}_\cG}\longrightarrow \Rep(\cG)_\Ran.
\end{equation} 

Applying the canonical self-duality of $\Rep(\cG)_\Ran$, the datum of the functor \eqref{e:spec Poinc and Gamma}
is equivalent to the datum of the pairing
$$\IndCoh^*(\Op^{\on{mon-free}}_\cG)_\Ran\otimes \Rep(\cG)_\Ran\to \Vect,$$
given by
\begin{multline}  \label{e:spec Poinc and Gamma pairing}
\IndCoh^*(\Op^{\on{mon-free}}_\cG)_\Ran \otimes \Rep(\cG)_\Ran
\overset{\Poinc^{\on{spec}}_{\cG,*}\otimes \on{Id}}\longrightarrow 
\QCoh(\LS_\cG)\otimes \Rep(\cG)_\Ran \to \\
\overset{\Gamma^{\on{spec}}_\cG\otimes \on{Id}} \longrightarrow \Rep(\cG)_\Ran\otimes \Rep(\cG)_\Ran \to \Vect.
\end{multline} 

We will prove (cf. \thmref{t:Loc and coeff}):

\begin{thm} \label{t:spec Poinc and Gamma}
The functor \eqref{e:spec Poinc and Gamma pairing} identifies canonically with
\begin{multline} \label{e:spec Poinc and Gamma pairing 1}
\IndCoh^*(\Op^{\on{mon-free}}_\cG)_\Ran \otimes \Rep(\cG)_\Ran \overset{\on{ins.unit}_\Ran\otimes \on{ins.unit}_\Ran}\longrightarrow \\
\to \IndCoh^*(\Op^{\on{mon-free}}_\cG)_{\Ran^{\subseteq}} \otimes \Rep(\cG)_{\Ran^{\subseteq}}\to \\
\to \left(\IndCoh^*(\Op^{\on{mon-free}}_\cG)\otimes \Rep(\cG)\right)_{\Ran^{\subseteq}\underset{\Ran}\times \Ran^{\subseteq}}
\overset{\sP^{\on{loc,enh}}_\cG}\longrightarrow \\
\to \CO_{\Op^\reg_\cG}\mod^{\on{fact}}_{\Ran^{\subseteq}\underset{\Ran}\times \Ran^{\subseteq}} 
\overset{\on{C}^{\on{fact}}_\cdot(X;\CO_{\Op^{\on{reg}}_\cG},-)_{\Ran^{\subseteq}\underset{\Ran}\times \Ran^{\subseteq}}}\longrightarrow
\Dmod(\Ran^{\subseteq}\underset{\Ran}\times \Ran^{\subseteq}) \to \\
\overset{\on{C}^\cdot_c(\Ran^{\subseteq}\underset{\Ran}\times \Ran^{\subseteq},-)}
\longrightarrow  \Vect,
\end{multline}
where $\sP^{\on{loc,enh}}_\cG$ is the functor introduced in \secref{sss:PGc pairing}. 
\end{thm}

\begin{rem}

Note that the functor \eqref{e:spec Poinc and Gamma pairing 1}, appearing in \thmref{t:spec Poinc and Gamma} can also be rewritten as
\begin{multline*} 
\IndCoh^*(\Op^{\on{mon-free}}_\cG)_\Ran \otimes \Rep(\cG)_\Ran \overset{\on{ins.unit}_\Ran\otimes \on{ins.unit}_\Ran}\to \\
\to \IndCoh^*(\Op^{\on{mon-free}}_\cG)_{\Ran^{\subseteq}} \otimes \Rep(\cG)_{\Ran^{\subseteq}}\to \\
\left(\IndCoh^*(\Op^{\on{mon-free}}_\cG)\otimes \Rep(\cG)\right)_{\Ran^{\subseteq}\underset{\Ran}\times \Ran^{\subseteq}}
\overset{\sP^{\on{loc}}_\cG}\longrightarrow
\Dmod(\Ran^{\subseteq}\underset{\Ran}\times \Ran^{\subseteq}) 
\overset{\on{C}^\cdot_c(\Ran^{\subseteq}\underset{\Ran}\times \Ran^{\subseteq})}\longrightarrow \Vect,
\end{multline*} 
i.e., instead of $\on{C}^{\on{fact}}_\cdot(X;\CO_{\Op^{\on{reg}}_\cG},-)_{\Ran^{\subseteq}\underset{\Ran}\times \Ran^{\subseteq}}$
we can use the functor $\oblv_{\CO_{\Op^{\on{reg}}_\cG},\Ran^{\subseteq}\underset{\Ran}\times \Ran^{\subseteq}}$.
This follows by the same manipulation as in Remark \ref{r:simplify Loc and coeff}

\end{rem}

\sssec{}

The rest of this subsection is devoted to the proof of \thmref{t:spec Poinc and Gamma}.

\medskip

First, using the (non-derived) Satake action, as in the proof of \thmref{t:Loc and coeff}, we obtain that the assertion of the theorem 
is equivalent to that of the following: 

\begin{thm} \label{t:spec Poinc and Gamma vac}
The functor
$$\IndCoh^*(\Op^{\on{mon-free}}_\cG)_\Ran 
\overset{\Poinc^{\on{spec}}_{\cG,*}}\longrightarrow \QCoh(\LS_\cG) \overset{\Gamma(\LS_\cG,-)}\longrightarrow \Vect$$
identifies canonically with
\begin{multline} 
\IndCoh^*(\Op^{\on{mon-free}}_\cG)_\Ran \overset{(\iota^\mf)^\IndCoh_*}\to  \IndCoh^*(\Op^{\on{mer}}_\cG)_\Ran 
\overset{\Gamma^\IndCoh(\Op^{\on{mer}}_\cG,-)^{\on{enh}}}\longrightarrow \\
\to \CO_{\Op^\reg_\cG}\mod^{\on{fact}}_\Ran \overset{\on{C}^{\on{fact}}_\cdot(X;\CO_{\Op^{\on{reg}}_\cG},-)}\longrightarrow \Vect.$$
\end{multline}
\end{thm} 

\sssec{}

For expositional purposes, will replace the situation over $\Ran$ by one with a fixed $\ul{x}\in \Ran$.
So, we want to show that the composition
\begin{equation} \label{e:spec Poinc and Gamma 1}
\IndCoh^*(\Op^\mf_{\cG,\ul{x}})
\overset{\Poinc^{\on{spec}}_{\cG,*,\ul{x}}}\longrightarrow \QCoh(\LS_\cG) \overset{\Gamma(\LS_\cG,-)}\longrightarrow \Vect
\end{equation} 
identifies canonically with
\begin{multline} \label{e:spec Poinc and Gamma 2}
\IndCoh^*(\Op^\mf_{\cG,\ul{x}}) \overset{(\iota^\mf)^\IndCoh_*}\to 
\IndCoh^*(\Op^\mer_{\cG,\ul{x}}) \overset{\Gamma^\IndCoh(\Op^{\on{mer}}_\cG,-)^{\on{enh}}}\longrightarrow \\
\to \CO_{\Op^\reg_\cG}\mod^{\on{fact}}_{\ul{x}} 
\overset{\on{C}_\cdot^{\on{ch}}(X,\CO_{\Op^{\on{reg}}_\cG},-)_{\ul{x}}}\longrightarrow \Vect.
\end{multline} 

\sssec{}

The functor \eqref{e:spec Poinc and Gamma 1} can be tautologically rewritten as the composition
$$\IndCoh^*(\Op^\mf_{\cG,\ul{x}})\overset{\on{ev}^{*,\IndCoh}_{\ul{x}}}\longrightarrow
\IndCoh^*(\Op^{\mf,\on{glob}}_{\cG,\ul{x}})\overset{\Gamma^\IndCoh(\Op^{\mf,\on{glob}}_{\cG,\ul{x}},-)}\longrightarrow \Vect,$$
and further as
\begin{multline} \label{e:spec Poinc and Gamma 3}
\IndCoh^*(\Op^\mf_{\cG,\ul{x}})\overset{\on{ev}^{*,\IndCoh}_{\ul{x}}}\longrightarrow
\IndCoh^*(\Op^{\mf,\on{glob}}_{\cG,\ul{x}}) \overset{(\iota^{\mf,\on{glob}})^\IndCoh_*}\longrightarrow \\
\to \IndCoh^*(\Op^{\mer,\on{glob}}_{\cG,\ul{x}})\overset{\Gamma^\IndCoh(\Op^{\mer,\on{glob}}_{\cG,\ul{x}},-)}\longrightarrow \Vect.
\end{multline}

Applying \eqref{e:ev and iota 2}, we rewrite \eqref{e:spec Poinc and Gamma 3} as
\begin{multline} \label{e:spec Poinc and Gamma 4}
\IndCoh^*(\Op^\mf_{\cG,\ul{x}}) \overset{(\iota^\mf)^\IndCoh_*}\longrightarrow \IndCoh^*(\Op^\mer_{\cG,\ul{x}}) 
\overset{\on{ev}^{*,\IndCoh}_{\ul{x}}}\longrightarrow \\
\to \IndCoh^*(\Op^{\mer,\on{glob}}_{\cG,\ul{x}})
\overset{\Gamma^\IndCoh(\Op^{\mer,\on{glob}}_{\cG,\ul{x}},-)}\longrightarrow \Vect.
\end{multline}

\medskip

Thus, it suffices to establish an isomorphism between 
\begin{equation} \label{e:spec Poinc and Gamma 5}
\IndCoh^*(\Op^\mer_{\cG,\ul{x}}) 
\overset{\on{ev}^{*,\IndCoh}_{\ul{x}}}\longrightarrow \IndCoh^*(\Op^{\mer,\on{glob}}_{\cG,\ul{x}})
\overset{\Gamma^\IndCoh(\Op^{\mer,\on{glob}}_{\cG,\ul{x}},-)}\longrightarrow \Vect
\end{equation}
and 
\begin{equation} \label{e:spec Poinc and Gamma 6}
\IndCoh^*(\Op^\mer_{\cG,\ul{x}}) \overset{\Gamma^\IndCoh(\Op^{\on{mer}}_\cG,-)^{\on{enh}}}\longrightarrow 
\CO_{\Op^\reg_\cG}\mod^{\on{fact}}_{\ul{x}} 
\overset{\on{C}_\cdot^{\on{ch}}(X,\CO_{\Op^{\on{reg}}_\cG},-)_{\ul{x}}}\longrightarrow \Vect.
\end{equation} 

\sssec{}

By \lemref{l:IndCoh* vs QCohco on mer adj}, we can rewrite \eqref{e:spec Poinc and Gamma 5} as
\begin{equation} \label{e:spec Poinc and Gamma 7}
\IndCoh^*(\Op^\mer_{\cG,\ul{x}}) \overset{\Psi_{\Op^\mer_{\cG,\ul{x}}}}\longrightarrow 
\QCoh_{\on{co}}(\Op^\mer_{\cG,\ul{x}}) \overset{\on{ev}^*_{\ul{x}}}\longrightarrow \QCoh_{\on{co}}(\Op^{\mer,\on{glob}}_{\cG,\ul{x}})
\overset{\Gamma(\Op^{\mer,\on{glob}}_{\cG,\ul{x}},-)}\longrightarrow \Vect,
\end{equation}
while \eqref{e:spec Poinc and Gamma 6} is by definition
\begin{equation} \label{e:spec Poinc and Gamma 8}
\IndCoh^*(\Op^\mer_{\cG,\ul{x}})  \overset{\Psi_{\Op^\mer_{\cG,\ul{x}}}}\longrightarrow 
\QCoh_{\on{co}}(\Op^\mer_{\cG,\ul{x}}) 
\overset{\Gamma(\Op^{\on{mer}}_\cG,-)^{\on{enh}}}\longrightarrow 
\CO_{\Op^\reg_\cG}\mod^{\on{fact}}_{\ul{x}} 
\overset{\on{C}_\cdot^{\on{ch}}(X,\CO_{\Op^{\on{reg}}_\cG},-)_{\ul{x}}}\longrightarrow \Vect.
\end{equation} 

Hence, it suffices to establish an isomorphism between 
\begin{equation} \label{e:spec Poinc and Gamma 9}
\QCoh_{\on{co}}(\Op^\mer_{\cG,\ul{x}}) \overset{\on{ev}^*_{\ul{x}}}\longrightarrow \QCoh_{\on{co}}(\Op^{\mer,\on{glob}}_{\cG,\ul{x}})
\overset{\Gamma(\Op^{\mer,\on{glob}}_{\cG,\ul{x}},-)}\longrightarrow \Vect
\end{equation}
and
\begin{equation} \label{e:spec Poinc and Gamma 10}
\QCoh_{\on{co}}(\Op^\mer_{\cG,\ul{x}}) 
\overset{\Gamma(\Op^{\on{mer}}_\cG,-)^{\on{enh}}}\longrightarrow 
\CO_{\Op^\reg_\cG}\mod^{\on{fact}}_{\ul{x}} 
\overset{\on{C}_\cdot^{\on{ch}}(X,\CO_{\Op^{\on{reg}}_\cG},-)_{\ul{x}}}\longrightarrow \Vect.
\end{equation} 

However, the latter is the statement of \propref{p:global section modules}. 

\section{The Langlands functor} \label{s:Langlands functor}

In this section we recall the construction of the Langlands functor, and establish the following of its properties:

\begin{itemize}

\item Compatibility with the functors $\on{coeff}_G$ and $\Gamma^{\on{spec}}_\cG$;

\item Compatibility with the actions of $\Sph_G$ and $\Sph^{\on{spec}}_\cG$;

\item Compatibility with the functors $\Loc_G$ and $\on{Poinc}^{\on{spec}}_{\cG,*}$.

\end{itemize}

\ssec{Recollections on the Langlands functor--the coarse version} \label{ss:L functor coarse}

In this and the next subsections we recall the construction of the \emph{coarse version of the} Langlands functor 
\begin{equation}  \label{e:Langlands functor}
\BL_{G,\on{coarse}}:\Dmod_{\frac{1}{2}}(\Bun_G)\to \QCoh(\LS_\cG).
\end{equation} 

\sssec{}

We consider $(\Rep(\cG)_\Ran)^\star$ as a monoidal category (see \secref{sss:conv non-unital}), and 
\begin{equation} \label{e:Loc spec again}
\Loc^{\on{spec}}_\cG:(\Rep(\cG)_\Ran)^\star\to \QCoh(\LS_\cG)
\end{equation} 
as a monoidal functor. Recall that $\Loc^{\on{spec}}_\cG$ is a \emph{localization},
i.e., its right adjoint is fully faithful (the proof is given, e.g., in \cite[Corollary C.1.8 and Sect. C.1.9]{GLC4}).

\sssec{}

We consider $\Dmod_{\frac{1}{2}}(\Bun_G)$ as acted on by $(\Rep(\cG)_\Ran)^\star$ via
the action of $(\Sph_{G,\Ran})^\star$ on $\Dmod_{\frac{1}{2}}(\Bun_G)$ (see \secref{sss:action of conv}) and 
$$\on{Sat}_G^{\on{nv}}:(\Rep(\cG)_\Ran)^\star\to (\Sph_{G,\Ran})^\star.$$

\medskip

According to \cite[Corollary 4.5.5]{Ga1}, the action of $(\Rep(\cG)_\Ran)^\star$ on 
$\Dmod_{\frac{1}{2}}(\Bun_G)$ factors through \eqref{e:Loc spec again}, so we obtain an
action on $\Dmod_{\frac{1}{2}}(\Bun_G)$ of $\QCoh(\LS_\cG)$.

\sssec{}

The coarse Langlands functor
$$\BL_{G,\on{coarse}}:\Dmod_{\frac{1}{2}}(\Bun_G)\to \QCoh(\LS_\cG),$$
as constructed in \cite[Sect. 1.4]{GLC1}, is uniquely characterized by the following two properties:

\medskip

\begin{itemize}

\item The functor $\BL_{G,\on{coarse}}$ is $\QCoh(\LS_\cG)$-linear;

\smallskip

\item The diagram
\begin{equation} \label{e:coeff comp L vac}
\CD
\Vect @>{\on{Id}}>> \Vect \\
@A{\on{coeff}_G^{\on{Vac,glob}}}AA @AA{\Gamma(\LS_\cG,-)}A \\
\Dmod_{\frac{1}{2}}(\Bun_G) @>{\BL_{G,\on{coarse}}}>>  \QCoh(\LS_\cG)
\endCD
\end{equation} 
commutes.

\end{itemize}

\sssec{}

Note that since $\Loc^{\on{spec}}_\cG$ is a localization, the second property can be equivalently formulated
as linearity with respect to $(\Rep(\cG)_\Ran)^\star$.

\medskip

By \corref{c:indep action from conv}, we obtain that the functor
$$\BL_{G,\on{coarse}}\otimes \on{Id}: \Dmod_{\frac{1}{2}}(\Bun_G) \otimes \Dmod(\Ran)\to \QCoh(\LS_\cG)\otimes \Dmod(\Ran)$$
intertwines the actions of $(\Rep(\cG)_\Ran)^{\sotimes}$ on both sides.

\medskip

In fact, \corref{c:indep action from conv} implies that the functor 
$$\BL_{G,\on{coarse}}\otimes \on{Id}: \Dmod_{\frac{1}{2}}(\Bun_G) \otimes 
\ul\Dmod(\Ran^{\on{untl}})\to \QCoh(\LS_\cG)\otimes \ul\Dmod(\Ran^{\on{untl}})$$
intertwines the actions of $\Rep(\cG)$, viewed as a crystal of monoidal categories over $\Ran^{\on{untl}}$.

\sssec{}

We now claim:

\begin{prop} \label{p:L coarse Whit compat} 
The following diagram commutes:
\begin{equation} \label{e:coeff comp L coarse}
\CD
\Whit^!(G)_\Ran @>{\on{CS}_G}>> \Rep(\cG)_\Ran \\
@A{\on{coeff}_G[2\delta_{N_{\rho(\omega_X)}}]}AA @AA{\Gamma^{\on{spec}}_\cG}A \\
\Dmod_{\frac{1}{2}}(\Bun_G) @>{\BL_{G,\on{coarse}}}>> \QCoh(\LS_\cG),
\endCD
\end{equation} 
where $\delta_{N_{\rho(\omega_X)}}$ is as in \secref{sss:delta N}.
\end{prop} 

\begin{proof}

It suffices to construct the datum of commutativity for the diagram 
\begin{equation} \label{e:coeff comp L coarse Ran}
\CD
\Whit^!(G)_\Ran @>{\on{CS}_G}>> \Rep(\cG)_\Ran \\
@A{\on{coeff}_{G,\Ran}[2\delta_{N_{\rho(\omega_X)}}]}AA @AA{\Gamma^{\on{spec}}_{\cG,\Ran}}A \\ 
\Dmod_{\frac{1}{2}}(\Bun_G)\otimes \Dmod(\Ran) @>{\BL_{G,\on{coarse}}\otimes \on{Id}}>> \QCoh(\LS_\cG)\otimes \Dmod(\Ran). 
\endCD
\end{equation} 

Consider the categories appearing in \eqref{e:coeff comp L coarse Ran} as equipped with an action of 
$(\Rep(\cG)_\Ran)^{\sotimes}$. We will construct a datum of commutativity of \eqref{e:coeff comp L coarse Ran} 
as $(\Rep(\cG)_\Ran)^{\sotimes}$-module categories.

\medskip

Note, however, that the upper right corner, i.e., 
$$\Rep(\cG)_\Ran=(\Rep(\cG)_\Ran)^{\sotimes},$$ 
is \emph{co-free}, when viewed as a module over itself. Hence, the datum of commutativity of 
\eqref{e:coeff comp L coarse Ran} as $(\Rep(\cG)_\Ran)^{\sotimes}$-module categories, is equivalent to
the datum of commutativity of the \emph{outer diagram} in 
\begin{equation} \label{e:coeff comp L coarse Ran vac}
\CD
\Dmod(\Ran) @>{\on{Id}}>> \Dmod(\Ran) \\
@AAA @AAA \\
\Whit^!(G)_\Ran @>{\on{CS}_G}>> \Rep(\cG)_\Ran \\
@A{\on{coeff}_{G,\Ran}[2\delta_{N_{\rho(\omega_X)}}]}AA @AA{\Gamma^{\on{spec}}_{\cG,\Ran}}A \\ 
\Dmod_{\frac{1}{2}}(\Bun_G)\otimes \Dmod(\Ran) @>{\BL_{G,\on{coarse}}\otimes \on{Id}}>> \QCoh(\LS_\cG)\otimes \Dmod(\Ran),
\endCD
\end{equation} 
as $\Dmod(\Ran) $-linear categories, 
where:

\begin{itemize}

\item The upper right vertical arrow is the factorization functor $\on{inv}_\cG:\Rep(\cG)\to \Vect$;

\smallskip

\item The  upper left vertical arrow is the factorization functor 
$$\Whit^!(G)_\Ran \to \Dmod_{\frac{1}{2}}(\Gr_{G,\rho(\omega_X)})\to \Vect,$$
where the second arrow is the functor of !-fiber at the unit.

\end{itemize}

\medskip

In its turn, the datum of commutativity of the outer diagram in \eqref{e:coeff comp L coarse Ran vac} is equivalent to the datum
of commutativity of 
\begin{equation} \label{e:coeff comp L coarse Ran vac again}
\CD
\Dmod(\Ran) @>{\on{Id}}>> \Dmod(\Ran) \\
@AAA @AAA \\
\Whit^!(G)_\Ran & & \Rep(\cG)_\Ran \\
@A{\on{coeff}_{G,\Ran}[2\delta_{N_{\rho(\omega_X)}}]}AA @AA{\Gamma^{\on{spec}}_{\cG,\Ran}}A \\  
\Dmod_{\frac{1}{2}}(\Bun_G)\otimes \Dmod(\Ran) & &  \QCoh(\LS_\cG)\otimes \Dmod(\Ran) \\
@A{\on{Id}\otimes \omega_{\Ran}}AA @AA{\on{Id}\otimes \omega_{\Ran}}A \\
\Dmod_{\frac{1}{2}}(\Bun_G) @>{\BL_{G,\on{coarse}}}>> \QCoh(\LS_\cG)
\endCD
\end{equation} 
just as DG categories.

\medskip

Note, however, that the composite left vertical arrow in \eqref{e:coeff comp L coarse Ran vac again} is the functor
$$\Dmod_{\frac{1}{2}}(\Bun_G)  \overset{\on{coeff}_G^{\on{Vac}}[2\delta_{N_{\rho(\omega_X)}}]}\longrightarrow 
\Vect \overset{\omega_{\Ran}}\longrightarrow \Dmod(\Ran)$$
and the composite right vertical arrow in \eqref{e:coeff comp L coarse Ran vac again} is the functor
$$\QCoh(\LS_\cG)\overset{\Gamma(\LS_\cG,-)}\longrightarrow \Vect \overset{\omega_{\Ran}}\longrightarrow \Dmod(\Ran).$$

Now, the required commutativity is supplied by \eqref{e:coeff comp L vac}, combined with \lemref{l:Vac coeff loc vs glob}.  

\end{proof} 

\ssec{Compatibility with the full spherical action}

\sssec{} \label{sss:Hecke spec action Ran}

As was mentioned in \secref{sss:Sph spec action x}, for a fixed $\ul{x}\in \Ran$, the category 
$\IndCoh(\LS_\cG)$ carries an action of $\Sph^{\on{spec}}_{\cG,\ul{x}}$. 

\medskip

In \secref{sss:Sph spec action on LS} we will extend this to an action of 
$(\Sph^{\on{spec}}_{\cG,\Ran})^{\sotimes}$ on $\IndCoh(\LS_\cG)\otimes \Dmod(\Ran)$. 
In fact, we have an action of $\Sph^{\on{spec}}_\cG$, viewed as a crystal of monoidal categories 
over $\Ran^{\on{untl}}$ on $\IndCoh(\LS_\cG)\otimes \ul\Dmod(\Ran^{\on{untl}})$. 

\medskip

By \secref{sss:action of conv}, this gives rise to an action of 
$$(\Sph^{\on{spec}}_{\cG,\Ran})^\star\twoheadrightarrow \Sph^{\on{spec}}_{\cG,\Ran^{\on{untl}},\on{indep}}$$
on $\IndCoh(\LS_\cG)$.

\sssec{}

A basic feature of the above action is that the functor
$$\Gamma^{\on{spec},\on{IndCoh}}_{\cG,\Ran}:\IndCoh(\LS_\cG)\otimes \Dmod(\Ran) \to \Rep(\cG)_\Ran$$
is compatible with the actions of $(\Sph^{\on{spec}}_{\cG,\Ran})^{\sotimes}$ on the two sides. 

\sssec{}

We claim:

\begin{lem} \hfill \label{l:Sph spec acts on QCoh}

\smallskip

\noindent{\em(a)} 
The action of $(\Sph^{\on{spec}}_{\cG,\Ran})^{\sotimes}$ on $\IndCoh(\LS_\cG)\otimes \Dmod(\Ran)$
preserves the subcategory 
$$\QCoh(\LS_\cG)\otimes \Dmod(\Ran)\overset{\Xi_{\LS_\cG}\otimes \on{Id}}\hookrightarrow 
\IndCoh(\LS_\cG)\otimes \Dmod(\Ran).$$

\smallskip

\noindent{\em(b)} The resulting action of $(\Sph^{\on{spec}}_{\cG,\Ran})^{\sotimes}$ 
on
$\QCoh(\LS_\cG)\otimes \Dmod(\Ran)$ is compatible with the projection
$$\IndCoh(\LS_\cG)\otimes \Dmod(\Ran) \overset{\Psi_{\LS_\cG}\otimes \on{Id}}\twoheadrightarrow 
\QCoh(\LS_\cG)\otimes \Dmod(\Ran).$$

\end{lem} 

\begin{proof}

To prove point (a), it suffices to show that the generators of $(\Sph^{\on{spec}}_{\cG,\Ran})^{\sotimes}$ preserve 
the subcategory $\QCoh(\LS_\cG)\otimes \Dmod(\Ran)$. We take these generators to be the essential image of the factorization functor
$$\on{nv}:\Rep(\cG)\to \Sph^{\on{spec}}_\cG.$$

This makes the assertion evident: the resulting action is the natural action of $(\Rep(\cG)_\Ran)^{\sotimes}$ on
$\QCoh(\LS_\cG)\otimes \Dmod(\Ran)$.

\medskip

Point (b) of the lemma follows similarly.

\end{proof}

\begin{rem}
Note that for a fixed $\ul{x}\in \Ran$, the action of $\Sph^{\on{spec}}_{\cG,\ul{x}}$ on $\QCoh(\LS_\cG)$ factors
through the quotient 
$$\Sph^{\on{spec}}_{\cG,\ul{x}}\twoheadrightarrow \Sph^{\on{spec}}_{\cG,\on{temp},\ul{x}},$$
see \secref{sss:temp spec}. This follows from the fact that the action of $\Sph^{\on{spec}}_{\cG,\ul{x}}$ on $\QCoh(\LS_\cG)$ 
is given by t-exact functors, combined with the fact that the t-structure on $\QCoh(\LS_\cG)$ is separated.

\medskip

We do \emph{not} know how to formulate a parallel property for the action of 
$(\Sph^{\on{spec}}_{\cG,\Ran})^{\sotimes}$ on the category $\QCoh(\LS_\cG)\otimes \Dmod(\Ran)$,
see Remark \ref{r:no Ran temp}. 

\end{rem} 

\sssec{}

We now claim:

\begin{prop} \label{p:L coarse Sph compat}
The functor 
$$\BL_{G,\on{coarse}}\otimes \on{Id}: \Dmod_{\frac{1}{2}}(\Bun_G)\otimes \Dmod(\Ran)\to \QCoh(\LS_\cG)\otimes \Dmod(\Ran)$$
intertwines the $(\Sph_{G,\Ran})^{\sotimes}$ action on the left-hand side with the $(\Sph^{\on{spec}}_{\cG,\Ran})^{\sotimes}$-action
on the right-hand side via the functor 
$$(\Sph_{G,\Ran})^{\sotimes}\to (\Sph^{\on{spec}}_{\cG,\Ran})^{\sotimes},$$
induced by the factorization functor 
$$\on{Sat}_G:\Sph_G\to \Sph^{\on{spec}}_\cG.$$
\end{prop}

The rest of this subsection is devoted to the proof of \propref{p:L coarse Sph compat}.

\sssec{}

By Corollaries \ref{c:unital vs non-unital Ran actions} and \ref{c:integrated actions}, we can 
reformulate the assertion of the proposition as follows: the functor 
$$\BL_{G,\on{coarse}}: \Dmod_{\frac{1}{2}}(\Bun_G)\to \QCoh(\LS_\cG)$$
intertwines the actions of 
$$\Sph_{G,\Ran^{\on{untl}},\on{indep}}$$
on the left-hand side with the action of
$$\Sph^{\on{spec}}_{\cG,\Ran^{\on{untl}},\on{indep}}$$
on the right-hand side via the functor
$$\Sph_{G,\Ran^{\on{untl}},\on{indep}}\to \Sph^{\on{spec}}_{\cG,\Ran^{\on{untl}},\on{indep}},$$
induced by the factorization functor 
$$\on{Sat}_G:\Sph_G\to \Sph^{\on{spec}}_\cG.$$

\sssec{}

We start with the commutative diagram \eqref{e:coeff comp L coarse}
\begin{equation} \label{e:coeff comp L coarse again}
\CD
\Whit^!(G)_\Ran @>{\on{CS}_G}>> \Rep(\cG)_\Ran \\
@A{\on{coeff}_G[2\delta_{N_{\rho(\omega_X)}}]}AA @AA{\Gamma^{\on{spec}}_\cG}A \\
\Dmod_{\frac{1}{2}}(\Bun_G) @>{\BL_{G,\on{coarse}}}>> \QCoh(\LS_\cG),
\endCD
\end{equation} 
and note that the vertical arrows factor as
$$\Dmod_{\frac{1}{2}}(\Bun_G)  \to \Whit^!(G)_{\Ran^{\on{untl}},\on{indep}}\hookrightarrow \Whit^!(G)_\Ran$$
and 
$$\QCoh(\LS_\cG)\to \Rep(\cG)_{\Ran^{\on{untl}},\on{indep}}\hookrightarrow \Rep(\cG)_\Ran,$$
respectively, 
so that we obtain a commutative diagram
\begin{equation} \label{e:coeff comp L coarse indep}
\CD
\Whit^!(G)_{\Ran^{\on{untl}},\on{indep}} @>{\on{CS}_G}>> \Rep(\cG)_{\Ran^{\on{untl}},\on{indep}} \\
@A{\on{coeff}_{G,\Ran^{\on{untl}},\on{indep}}[2\delta_{N_{\rho(\omega_X)}}]}AA @AA{\Gamma^{\on{spec}}_{\cG,\Ran^{\on{untl}},\on{indep}}}A \\
\Dmod_{\frac{1}{2}}(\Bun_G) @>{\BL_{G,\on{coarse}}}>> \QCoh(\LS_\cG).
\endCD
\end{equation}

\sssec{}

Since the functor
$$\on{coeff}_{G,\Ran}:\Dmod_{\frac{1}{2}}(\Bun_G) \otimes \Dmod(\Ran) \to \Whit^!(G)_\Ran$$
is compatible with the action of $(\Sph_{G,\Ran})^{\sotimes}$, from Corollaries \ref{c:unital vs non-unital Ran actions} 
and \ref{c:integrated actions} we obtain that the functor
\begin{multline} \label{e:coeff indep}
\Dmod_{\frac{1}{2}}(\Bun_G) \overset{-\otimes \omega_{\Ran^{\on{untl}}}}\simeq \\
\simeq \Dmod_{\frac{1}{2}}(\Bun_G) \otimes \Dmod(\Ran^{\on{untl}})_{\on{indep}} \overset{\on{coeff}_{G,\Ran^{\on{untl}}}}\longrightarrow
\Whit^!(G)_{\Ran^{\on{untl}},\on{indep}},
\end{multline} 
appearing as the left vertical arrow in \eqref{e:coeff comp L coarse indep}, 
is compatible with the action of $\Sph_{G,\Ran^{\on{untl}},\on{indep}}$.

\medskip

Similarly, since the functor
$$\Gamma^{\on{spec}}_{\cG,\Ran}:\QCoh(\LS_\cG)\otimes \Dmod(\Ran) \to \Rep(\cG)_\Ran$$
is compatible with the actions of $(\Sph^{\on{spec}}_{\cG,\Ran})^{\sotimes}$, we obtain that the functor
\begin{equation} \label{e:Gamma indep}
\QCoh(\LS_\cG)\overset{-\otimes \omega_{\Ran^{\on{untl}}}}
\simeq \QCoh(\LS_\cG) \otimes \Dmod(\Ran^{\on{untl}})_{\on{indep}} 
\overset{\Gamma^{\on{spec}}_{\cG,\Ran^{\on{untl}}}}\longrightarrow \Rep(\cG)_{\Ran^{\on{untl}},\on{indep}},
\end{equation} 
appearing as the right vertical arrow in \eqref{e:coeff comp L coarse indep}, 
is compatible with the action of $\Sph^{\on{spec}}_{\cG,\Ran^{\on{untl}},\on{indep}}$. 

\sssec{}

Recall now that the functor
$$\Gamma^{\on{spec}}_\cG: \QCoh(\LS_\cG)\to \Rep(\cG)_\Ran$$
is fully faithful. By \propref{p:Ran vs indep}, this implies that the functor $\Gamma^{\on{spec}}_{\cG,\Ran^{\on{untl}},\on{indep}}$
is fully faithful.

\medskip

Hence, in order to equip $\BL_{G,\on{coarse}}$ with a datum of compatibility with respect to 
$$\Sph_{G,\Ran^{\on{untl}},\on{indep}}\overset{\Sat_G}\simeq \Sph^{\on{spec}}_{\cG,\Ran^{\on{untl}},\on{indep}}$$
it suffices to do so for the counter-clockwise composition in \eqref{e:coeff comp L coarse indep}.

\medskip

By the commutativity of \eqref{e:coeff comp L coarse indep}, this is equivalent to endowing 
the clockwise composition in \eqref{e:coeff comp L coarse indep} with a datum of compatibility 
with the above action. 

\medskip

However, this follows from the compatibility for \eqref{e:coeff indep} mentioned above, combined with
the compatibility of $\on{CS}_G$ with $\Sat_G$. 

\qed[\propref{p:L coarse Sph compat}]

\ssec{The actual Langlands functor} \label{ss:L functor}

\sssec{}

We now quote the following result established in \cite[Corollary 1.6.5]{GLC1}: 

\begin{thm} \label{t:true Langlands}
There exists a uniquely defined functor
$$\BL_G:\Dmod_{\frac{1}{2}}(\Bun_G)\to \IndCoh_\Nilp(\LS_\cG),$$
subject to the following conditions:

\begin{itemize}

\item $\Psi_{\LS_\cG}\circ \BL_G\simeq \BL_{G,\on{coarse}}$;

\item The functor $\BL_G$ sends compact objects in $\Dmod_{\frac{1}{2}}(\Bun_G)$ to
$$\IndCoh_\Nilp(\LS_\cG)^{>-\infty}\subset \IndCoh_\Nilp(\LS_\cG).$$

\end{itemize}

\end{thm} 

\sssec{}

Let
$$\Xi_{0,\Nilp}:\QCoh(\LS_\cG) \rightleftarrows \IndCoh_\Nilp(\LS_\cG):\Psi_{\Nilp,0}$$
and
$$\Xi_{\Nilp,\on{all}}:\IndCoh_\Nilp(\LS_\cG) \rightleftarrows \IndCoh(\LS_\cG):\Psi_{\on{all},\Nilp}$$
denote the resulting pairs of adjoint functors.

\medskip

By a slight abuse of notation we will denote by $\Gamma^{\on{spec},\IndCoh}_\cG$ the functor
$$\IndCoh_\Nilp(\LS_\cG)\overset{\Xi_{\Nilp,\on{all}}}\hookrightarrow \IndCoh(\LS_\cG)  
\overset{\Gamma^{\on{spec},\IndCoh}_\cG}\longrightarrow \Rep(\cG)_\Ran,$$
which is the same as
$$\IndCoh_\Nilp(\LS_\cG)\overset{\Psi_{\Nilp,0}}\twoheadrightarrow \QCoh(\LS_\cG)  
\overset{\Gamma^{\on{spec}}_\cG}\longrightarrow \Rep(\cG)_\Ran.$$

\sssec{}

From \propref{p:L coarse Whit compat} we formally obtain:

\begin{cor} \label{c:L Whit compat} 
The following diagram commutes:
\begin{equation} \label{e:coeff comp L}
\CD
\Whit^!(G)_\Ran @>{\on{CS}_G}>> \Rep(\cG)_\Ran \\
@A{\on{coeff}_G[2\delta_{N_{\rho(\omega_X)}}]}AA @AA{\Gamma^{\on{spec},\IndCoh}_\cG}A \\
\Dmod_{\frac{1}{2}}(\Bun_G) @>{\BL_G}>> \IndCoh_\Nilp(\LS_\cG).
\endCD
\end{equation} 
\end{cor} 

\sssec{}

Consider again the action of $(\Sph^{\on{spec}}_{\cG,\Ran})^{\sotimes}$ on $\IndCoh(\LS_\cG)\otimes \Dmod(\Ran)$. 
By the same mechanism as in \lemref{l:Sph spec acts on QCoh}, this action gives rise to an action of 
$(\Sph^{\on{spec}}_{\cG,\Ran})^{\sotimes}$ on the category $\IndCoh_\Nilp(\LS_\cG)\otimes \Dmod(\Ran)$, which is compatible with 
the functors
$$(\Xi_{0,\Nilp},\Psi_{\Nilp,0}) \text{ and } (\Xi_{\Nilp,\on{all}},\Psi_{\on{all},\Nilp}).$$

In particular, we obtain an action of 
$$(\Sph^{\on{spec}}_{\cG,\Ran})^\star\to \Sph^{\on{spec}}_{\cG,\Ran^{\on{untl}},\on{indep}}$$ on 
$\IndCoh_\Nilp(\LS_\cG)$.

\sssec{}

We now claim:

\begin{prop} \label{p:L Sph compat}
The functor $\BL_G$ intertwines the action of $(\Sph_{G,\Ran})^\star$ on the left-hand side and the action of 
$(\Sph^{\on{spec}}_{\cG,\Ran})^\star$ on the right-hand side via 
$$\on{Sat}_G:\Sph_G\to \Sph^{\on{spec}}_\cG.$$
\end{prop} 

\begin{proof}

Note that $(\Sph_{G,\Ran})^\star$ is compactly generated, and the subcategory 
$$((\Sph_{G,\Ran})^\star)^c\subset (\Sph_{G,\Ran})^\star$$
is closed under the monoidal operation, and its action on $\Dmod_{\frac{1}{2}}(\Bun_G)$ 
preserves the subcategory 
$$\Dmod_{\frac{1}{2}}(\Bun_G)^c\subset \Dmod_{\frac{1}{2}}(\Bun_G).$$

\medskip

Hence, in order to prove the proposition, it suffices to equip the functor 
$$\BL_G|_{\Dmod_{\frac{1}{2}}(\Bun_G)^c}:\Dmod_{\frac{1}{2}}(\Bun_G)^c\to \IndCoh_\Nilp(\LS_\cG)$$
with a datum of compatibility with respect to the action of $((\Sph_{G,\Ran})^\star)^c$.

\medskip

By the definition of the functor $\BL_G$, the restriction $\BL_G|_{\Dmod_{\frac{1}{2}}(\Bun_G)^c}$ factors as
$$\Dmod_{\frac{1}{2}}(\Bun_G)^c \overset{\BL_G}\to \IndCoh_\Nilp(\LS_\cG)^{>-\infty}\hookrightarrow \IndCoh_\Nilp(\LS_\cG).$$
where 
$$\IndCoh_\Nilp(\LS_\cG)^{>-\infty}\subset \IndCoh_\Nilp(\LS_\cG)$$
is also preserved by the action of
$$((\Sph_{G,\Ran})^\star)^c\simeq ((\Sph^{\on{spec}}_{\cG,\Ran})^\star)^c.$$

Hence, it suffices to endow the functor
$$\Dmod_{\frac{1}{2}}(\Bun_G)^c \overset{\BL_G}\to \IndCoh_\Nilp(\LS_\cG)^{>-\infty}$$
with a datum of compatibility with respect to the action of $((\Sph_{G,\Ran})^\star)^c$.

\medskip

Next, we note that the functor
$$\Psi_{\Nilp,0}|_{\IndCoh_\Nilp(\LS_\cG)^{>-\infty}}: \IndCoh_\Nilp(\LS_\cG)^{>-\infty}\to \QCoh(\LS_\cG)$$
is compatible with the action of $((\Sph^{\on{spec}}_{\cG,\Ran})^\star)^c$ and is fully faithful.

\medskip

Hence, it suffices to endow the composition
$$\Dmod_{\frac{1}{2}}(\Bun_G)^c \overset{\BL_G}\to \IndCoh_\Nilp(\LS_\cG)^{>-\infty}
\overset{\Psi_{\Nilp,0}}\to \QCoh(\LS_\cG)$$
with a datum of compatibility with respect to the action of $((\Sph_{G,\Ran})^\star)^c$.

\medskip

However, the latter composition is the functor 
$$\BL_{G,\on{coarse}}|_{\Dmod_{\frac{1}{2}}(\Bun_G)^c},$$
and the required datum is supplied by \propref{p:L coarse Sph compat}.

\end{proof}

\sssec{}

Combining \propref{p:L Sph compat} with Corollaries \ref{c:unital vs non-unital Ran actions}, \ref{c:integrated actions} and 
\ref{c:indep action from conv}, we obtain:

\begin{cor} \label{c:L Sph compat} \hfill 

\smallskip

\noindent{\em(a)} The functor $\BL_G$
intertwines the action of $\Sph_{G,\Ran,\on{indep}}$ on the left-hand side and the action of 
$\Sph^{\on{spec}}_{\cG,\Ran^{\on{untl}},\on{indep}}$ on the right-hand side.

\smallskip

\noindent{\em(b)} The functor 
$$\BL_G\otimes \on{Id}:\Dmod_{\frac{1}{2}}(\Bun_G) \otimes \Dmod(\Ran)\to 
\IndCoh_\Nilp(\LS_\cG)\otimes \Dmod(\Ran)$$
intertwines the action of $(\Sph_{G,\Ran})^{\sotimes}$ on the left-hand side and the action of 
$(\Sph^{\on{spec}}_{\cG,\Ran})^{\sotimes}$ on the right-hand side.

\smallskip

\noindent{\em(c)} The functor 
$$\BL_G\otimes \on{Id}:\Dmod_{\frac{1}{2}}(\Bun_G) \otimes \ul\Dmod(\Ran^{\on{untl}})\to 
\IndCoh_\Nilp(\LS_\cG)\otimes \ul\Dmod(\Ran^{\on{untl}})$$
intertwines the actions of $\Sph_G$ and $\Sph^{\on{spec}}_\cG$, viewed as crystals 
of monoidal categores over $\Ran^{\on{untl}}$.

\end{cor} 

\ssec{Critical localization and temperedness}

\sssec{}

Choose $\ul{x}\in \Ran$, and let
$$\Dmod_{\frac{1}{2}}(\Bun_G)_{\on{temp},\ul{x}}:=\Sph_{G,\on{temp},\ul{x}}\underset{\Sph_{G,\ul{x}}}\otimes \Dmod_{\frac{1}{2}}(\Bun_G).$$

The pair of adjoint functors
$$\Sph_{G,\on{temp},\ul{x}}\rightleftarrows \Sph_{G,\ul{x}}$$ 
allows us to view $\Dmod_{\frac{1}{2}}(\Bun_G)_{\on{temp},\ul{x}}$ as a colocalization of $\Dmod_{\frac{1}{2}}(\Bun_G)$. 

\medskip

According to \cite[Sect. 2.6.2]{FR}, this colocalization is actually independent of the choice of $\ul{x}$. So from now
on we will omit the subscript and denote the corresponding sub/quotient category by $\Dmod_{\frac{1}{2}}(\Bun_G)_{\on{temp}}$. 
Denote by
$$\bu:\Dmod_{\frac{1}{2}}(\Bun_G)_{\on{temp}}\rightleftarrows \Dmod_{\frac{1}{2}}(\Bun_G):\bu^R$$
the corresponding pair of adjoint functors.

\sssec{}

From \propref{p:L Sph compat} we obtain:

\begin{cor} \label{c:temp Langlands}
There exists a uniquely defined functor
$$\BL_{G,\on{temp}}:\Dmod_{\frac{1}{2}}(\Bun_G)_{\on{temp}}\to \QCoh(\LS_\cG),$$
which makes both squares in the next diagram commute:
$$
\CD 
\Dmod_{\frac{1}{2}}(\Bun_G) @>{\BL_G}>> \IndCoh_\Nilp(\LS_\cG) \\
@V{\bu^R}VV @VV{\Psi_{\Nilp,0}}V \\
\Dmod_{\frac{1}{2}}(\Bun_G)_{\on{temp}} @>{\BL_{G,\on{temp}}}>> \QCoh(\LS_\cG) \\
@V{\bu}VV @VV{\Xi_{0,\Nilp}}V \\
\Dmod_{\frac{1}{2}}(\Bun_G) @>{\BL_G}>> \IndCoh_\Nilp(\LS_\cG).
\endCD
$$
Furthermore,
$$\BL_{G,\on{temp}} \simeq \BL_{G,\on{coarse}}\circ \bu.$$ 
\end{cor}

\sssec{}

Let 
$$\Loc_G:\KL(G)_{\crit,\Ran} \to \Dmod_{\frac{1}{2}}(\Bun_G)$$
be as in \secref{sss:just Loc}. 

\medskip

The following assertion is a counterpart of \propref{p:Poinc to QCoh}:

\begin{prop} \label{p:Loc is temp}
The essential image of the functor 
$$\Loc_G:\KL(G)_{\crit,\Ran}\to \Dmod_{\frac{1}{2}}(\Bun_G)$$
lies in 
$$\Dmod_{\frac{1}{2}}(\Bun_G)_{\on{temp}} \subset \Dmod_{\frac{1}{2}}(\Bun_G).$$
\end{prop} 

\begin{proof}

Repeats the proof of \propref{p:Poinc to QCoh} using \propref{p:KL is temp}. 

\end{proof}

\ssec{Compatibility of the Langlands functor with critical localization}

\sssec{}

The following theorem expresses the compatibility of the Langlands functor with critical localization: 

\begin{thm} \label{t:Langlands critical compat}
The diagram
$$
\CD
\Dmod_{\frac{1}{2}}(\Bun_G) @>{\BL_G}>> \IndCoh_\Nilp(\LS_\cG) \\
@A{\Loc_G\otimes \fl^{\otimes \frac{1}{2}}_{G,N_{\rho(\omega_X)}}\otimes 
\fl^{\otimes -1}_{N_{\rho(\omega_X)}}[-\delta_{N_{\rho(\omega_X)}}]}AA 
@AA{\Poinc^{\on{spec}}_{\cG,*}}A \\
\KL(G)_{\crit,\Ran} @>{\FLE_{G,\crit}}>>  \IndCoh^*(\Op^{\on{mon-free}}_\cG)_\Ran
\endCD
$$
commutes, where the lines $\fl^{\otimes \frac{1}{2}}_{G,N_{\rho(\omega_X)}}$ and 
$\fl_{N_{\rho(\omega_X)}}$ are as in \eqref{e:fl N} and \eqref{e:fl G N}, respectively. 
\end{thm}

The rest of the subsection is devoted to the proof of \thmref{t:Langlands critical compat}. 

\sssec{}

First, by Propositions \ref{p:Loc is temp} and \ref{p:Poinc to QCoh}, 
the commutativity of the diagram in \thmref{t:Langlands critical compat} is equivalent to 
the commutativity of the following one:
$$
\CD
\Dmod_{\frac{1}{2}}(\Bun_G)_{\on{temp}} @>{\BL_{G,\on{temp}}}>> \QCoh(\LS_\cG) \\
@AAA
@AA{\Poinc^{\on{spec}}_{\cG,*}}A \\
\KL(G)_{\crit,\Ran} @>{\FLE_{G,\crit}}>>  \IndCoh^*(\Op^{\on{mon-free}}_\cG)_\Ran,
\endCD
$$
and is further equivalent to the commutativity of 
\begin{equation} \label{e:Langlands critical compat 1}
\CD
\Dmod_{\frac{1}{2}}(\Bun_G) @>{\BL_{G,\on{coarse}}}>> \QCoh(\LS_\cG) \\
@AAA
@AA{\Poinc^{\on{spec}}_{\cG,*}}A \\
\KL(G)_{\crit,\Ran} @>{\FLE_{G,\crit}}>>  \IndCoh^*(\Op^{\on{mon-free}}_\cG)_\Ran. 
\endCD
\end{equation} 

\sssec{}

Since the right vertical arrow in \eqref{e:coeff comp L coarse} is fully faithful, it suffices to show 
that the two circuits in \eqref{e:Langlands critical compat 1} become isomorphic after composing with
the functor $\Gamma^{\on{spec}}_\cG$.

\medskip

Since the diagram \eqref{e:coeff comp L coarse} is commutative, we obtain that it suffices to establish the commutativity
of the diagram
$$
\CD
\Whit^!(G)_\Ran @>{\on{CS}_G}>> \Rep(\cG)_\Ran \\
@A{\on{coeff}_G[2\delta_{N_{\rho(\omega_X)}}]}AA @AA{\Gamma^{\on{spec}}_\cG}A \\
\Dmod_{\frac{1}{2}}(\Bun_G) & & \QCoh(\LS_\cG) \\
@A{\Loc_G\otimes \fl^{\otimes \frac{1}{2}}_{G,N_{\rho(\omega_X)}}\otimes 
\fl^{\otimes -1}_{N_{\rho(\omega_X)}}[-\delta_{N_{\rho(\omega_X)}}]}AA 
@AA{\Poinc^{\on{spec}}_{\cG,*}}A \\
\KL(G)_{\crit,\Ran} @>{\FLE_{G,\crit}}>>  \IndCoh^*(\Op^{\on{mon-free}}_\cG)_\Ran,
\endCD
$$
or which is the same 
\begin{equation} \label{e:FCD}
\CD
\Whit^!(G)_\Ran @>{\on{CS}_G}>> \Rep(\cG)_\Ran \\
@A{\on{coeff}_G}AA @AA{\Gamma^{\on{spec}}_\cG}A \\
\Dmod_{\frac{1}{2}}(\Bun_G) & & \QCoh(\LS_\cG) \\
@A{\Loc_G\otimes \fl^{\otimes \frac{1}{2}}_{G,N_{\rho(\omega_X)}}\otimes \fl^{\otimes -1}_{N_{\rho(\omega_X)}}}A{[\delta_{N_{\rho(\omega_X)}}]}A 
@AA{\Poinc^{\on{spec}}_{\cG,*}}A \\
\KL(G)_{\crit,\Ran} @>{\FLE_{G,\crit}}>>  \IndCoh^*(\Op^{\on{mon-free}}_\cG)_\Ran.
\endCD
\end{equation} 

\sssec{}

Applying duality, we obtain that it suffices to show that the pairing
\begin{multline} \label{e:FCD side 1}
\KL(G)_{\crit,\Ran} \otimes \Whit_*(G)_{\Ran} 
\overset{\Loc_G\otimes \on{Id}}\longrightarrow  \Dmod_{\frac{1}{2}}(\Bun_G)\otimes \Whit_*(G)_{\Ran} 
\overset{\on{coeff}_G\otimes \on{Id}}\longrightarrow \\
\to \Whit^!(G)_{\Ran}\otimes \Whit_*(G)_{\Ran}\to 
\Vect\overset{(-\otimes \fl^{\otimes \frac{1}{2}}_{G,N_{\rho(\omega_X)}}\otimes \fl^{\otimes -1}_{N_{\rho(\omega_X)}})[\delta_{N_{\rho(\omega_X)}}]}
\longrightarrow \Vect 
\end{multline}
agrees under the FLE equivalences
$$\KL(G)_{\crit,\Ran} \overset{\FLE_{G,\crit}}\simeq \IndCoh^*(\Op^{\on{mon-free}}_\cG)_\Ran
\text{ and } \Rep(\cG)_\Ran \overset{\FLE_{\cG,\infty}}\simeq \Whit_*(G)$$
with 
\begin{multline}  \label{e:FCD side 2}
\IndCoh^*(\Op^{\on{mon-free}}_\cG)_\Ran \otimes \Rep(\cG)_\Ran
\overset{\Poinc^{\on{spec}}_{\cG,*}\otimes \on{Id}}\longrightarrow 
\QCoh(\LS_\cG)\otimes \Rep(\cG)_\Ran \to \\
\overset{\Gamma^{\on{spec}}_\cG\otimes \on{Id}} \longrightarrow \Rep(\cG)_\Ran\otimes \Rep(\cG)_\Ran \to \Vect.
\end{multline}


\sssec{}

By \thmref{t:Loc and coeff}, the functor \eqref{e:FCD side 1} identifies canonically with \eqref{e:Loc and coeff pairing 1}. 
By \thmref{t:spec Poinc and Gamma}, the functor \eqref{e:FCD side 2} identifies canonically with \eqref{e:spec Poinc and Gamma pairing 1}. 

\medskip

The desired assertion follows now from \corref{c:two pairings coarse}. 

\qed[\thmref{t:Langlands critical compat}]



%
%
%
%

\newpage

\appendix

\centerline{\bf Part III: Appendix}

\bigskip

\centerline{ By J.~Campbell, L.~Chen, D.~Gaitsgory, K.~Lin, S.~Raskin and N.~Rozenblyum} 

\bigskip

The main body of the paper relies on a lot of foundational material, which is developed in this Appendix.
The main points are:

\medskip

\begin{itemize}

\item Ind-coherent sheaves on algebrao-geometric objects of infinite type 
(our main, but \emph{by far, not only} example is $\Op^\mf_\cG$). This is developed
out in \secref{s:IndCoh inf type};

\item The notion of \emph{factorization category}, and associated objects (factorization module
categories, factorization algebras within a factorization category, etc.). This is developed in 
\secref{s:fact}. 

\smallskip

\item The notion \emph{unitality} in the factorization setting. This is developed in \ref{s:unit};

\item A result connecting the (pre dual of the) category of quasi-coherent sheaves on the loop space 
and the category of factorization modules over the corresponding factorization algebra \ref{s:chiral mods}; 

\smallskip

\item The definition of the \emph{spectral spherical category} in the factorization setting 
(the underlying algebro-geometric object is so unwieldy that one cannot algorithmically apply
a procedure from \secref{s:IndCoh inf type}). This is developed in \secref{s:spec Sph fact}. 

\end{itemize} 

The majority of this material is of \emph{local} nature, i.e., it is needed to set up the local
Langlands theory. That said, some sections in this Appendix (notably, Sects. \ref{s:hor sect D sch} ,
\ref{s:indep} and \ref{s:add unit colax}) consider local-to-global constructions. 

\bigskip

\section{Ind-coherent sheaves in infinite type} \label{s:IndCoh inf type}

This section is devoted to the development of the theory of ind-coherent sheaves on algebro-geometric
objects of infinite-type.

\medskip

Prior to doing so, we introduce another player, which in some sense lies in between $\QCoh(-)$ and $\IndCoh(-)$. 
This object is denoted by
$$\CY\in \on{PreStk} \, \rightsquigarrow\, \QCoh_{\on{co}}(\CY),$$
and it is defined by a \emph{colimit} procedure (unlike $\QCoh(\CY)$, which is defined as a \emph{limit}). The category
$\QCoh_{\on{co}}(\CY)$ is a \emph{predual} of $\QCoh(\CY)$.

\medskip

We now turn to the IndCoh theory. A priori, $\IndCoh(-)$ is defined for affine schemes almost of finite type,
and by the process of right Kan extension on all prestacks that are locally almost of finite type. When $S$
is an affine scheme that is not of finite type, one can approximate it by affine schemes of finite type $S_\alpha$,
but then one faces a choice: one can define $\IndCoh(S)$ either as the colimit of $\IndCoh(S_\alpha)$ with respect to !-pullbacks
and as a limit of $\IndCoh(S_\alpha)$ with respect to *-pullbacks. This leads to two different categories, denoted $\IndCoh^!(S)$
and $\IndCoh^*(S)$, respectively. In good cases (technically, when $S$ is \emph{placid}), the Serre duality
in finite type gives rise to a duality between $\IndCoh^!(S)$ and $\IndCoh^*(S)$. 

\medskip

The majority of this section is devoted to developing the $\IndCoh^!(-)$ and $\IndCoh^*(-)$ theories,
and their interactions with other actors.

\ssec{The category \texorpdfstring{$\QCoh_{\on{co}}(-)$}{QCohco}}

\sssec{}

In section we will work with \emph{all} affine schemes (i.e., ones not necessarily almost of finite type). 
We denote the corresponding category by $\affSch$. We denote by $\on{PreStk}$ the category of all functors
$$(\affSch)^{\on{op}}\to \on{Spc}.$$

\sssec{}

Consider the functor
\begin{equation} \label{e:QCoh direct image}
\affSch\to \DGCat, \quad S\mapsto \QCoh(S),\quad (S_1\overset{f}\to S_2) \, \rightsquigarrow\, \QCoh(S_1)\overset{f_*}\to \QCoh(S_2).
\end{equation} 

\medskip

Consider the left Kan extension of \eqref{e:QCoh direct image} along the fully faithful embedding
$$\affSch\hookrightarrow \on{PreStk};$$
this yields a functor
\begin{equation} \label{e:QCoh direct image preStk}
\on{PreStk}\to \DGCat.
\end{equation} 

We will denote the value of \eqref{e:QCoh direct image preStk} on a given prestack $\CY$ by
$$\QCoh_{\on{co}}(\CY).$$

Explicitly,
\begin{equation} \label{e:QCohco expl}
\QCoh_{\on{co}}(\CY)\simeq \underset{S\to \CY,S\in \affSch}{\on{colim}}\, \QCoh(S),
\end{equation} 
where the colimit is formed using the pushforward functors. 

\begin{rem}

The above definition turns out to be useful in many contexts, but in the present
paper its main application is the following (see \secref{ss:QCoh co IndSch}):

\medskip

Let $\CY$ be and ind-affine ind-scheme 
$$\CY=\underset{i}{``\on{colim}"}\, \Spec(R_i).$$
Let $R$ denote the topological ring $\underset{i}{\on{lim}}\, R_i$.

\medskip

Then 
$$\QCoh_{\on{co}}(\CY)\simeq \underset{i}{\on{colim}}\, R_i\mod,$$
i.e., this formalizes the notion of ``the category of discrete $R$-modules". 

\medskip

Note that the above definition is close, but not the same, as $\IndCoh(\CY)$.
The difference is two-fold:

\begin{itemize}

\item $\IndCoh(\CY)$ is a priori defined only when $\CY$ is locally
almost of finite type (however, we will generalize that in \secref{sss:IndCoh *} below),
while $\QCoh_{\on{co}}(\CY)$ does not require this assumption;

\medskip

\item When $\CY$ is an affine scheme $Y$, we have $\QCoh_{\on{co}}(Y)=\QCoh(Y)$,
i.e., there is no renormalization procedure involved. (The price we will have to pay
for this is that, even for ind-schemes locally almost of finite type, the category 
$\QCoh_{\on{co}}(\CY)$ is not necessarily dualizable.)

\end{itemize}

\end{rem}

\sssec{} \label{sss:functoriality of QCohco}

By construction, the assignment
$$\CY\mapsto \QCoh_{\on{co}}(\CY)$$
has a functoriality with respect to pushforwards, i.e., for a map $\CY_1\to \CY_2$
we have the functor
$$f_*: \QCoh_{\on{co}}(\CY_1)\to  \QCoh_{\on{co}}(\CY_2).$$

\sssec{} \label{sss:QCoh co mult}

The construction 
$$\CY \rightsquigarrow \QCoh_{\on{co}}(\CY)$$ has a natural multiplicativity property. Namely, for
a pair of prestacks $\CY_1$ and $\CY_2$, we have a naturally defined equivalence
\begin{equation} \label{e:QCoh co mult}
\QCoh_{\on{co}}(\CY_1)\otimes \QCoh_{\on{co}}(\CY_2)\simeq \QCoh_{\on{co}}(\CY_1\times \CY_2).
\end{equation}

Namely, we have, by definition:
\begin{multline*} 
\QCoh_{\on{co}}(\CY_1)\otimes \QCoh_{\on{co}}(\CY_2)\simeq
\underset{S_1\in \affSch_{/\CY_1},S_2\in \affSch_{/\CY_2}}{\on{colim}}\, \QCoh(S_1)\otimes \QCoh(S_2)\simeq \\
\simeq \underset{S_1\in \affSch_{/\CY_1},S_2\in \affSch_{/\CY_2}}{\on{colim}}\, \QCoh(S_1\times S_2),
\end{multline*} 
and
$$\QCoh_{\on{co}}(\CY_1\times \CY_2)\simeq \underset{S\in \affSch_{/\CY_1\times \CY_2}}{\on{colim}}\, \QCoh(S).$$

Now, the functor
$$\affSch_{/\CY_1}\times \affSch_{/\CY_2}\to \affSch_{/\CY_1\times \CY_2}, \quad S_1,S_2\mapsto S_1\times S_2$$
is cofinal.

\sssec{} \label{sss:co predual}

Note that there is no reason for the category $\QCoh_{\on{co}}(\CY)$ to be dualizable. However,
we claim that $\QCoh_{\on{co}}(\CY)$ is \emph{a pre-dual} of $\QCoh(\CY)$, i.e., we 
have a canonical identification
\begin{equation} \label{e:QCoh co as predual}
\on{Funct}_{\on{cont}}(\QCoh_{\on{co}}(\CY),\Vect)=:\QCoh_{\on{co}}(\CY)^\vee\simeq \QCoh(\CY).
\end{equation} 

Indeed, using \eqref{e:QCohco expl}, we have
\begin{multline*} 
\QCoh_{\on{co}}(\CY)^\vee \simeq \left(\underset{S\to \CY,S\in \affSch}{\on{colim}}\, \QCoh(S)\right)^\vee \simeq \\
\simeq \underset{S\to \CY,S\in \affSch}{\on{lim}}\, \QCoh(S)^\vee \simeq \underset{S\to \CY,S\in \affSch}{\on{lim}}\, \QCoh(S)=:
\QCoh(\CY),
\end{multline*}
where we recall that for $f:S_1\to S_2$, with respect to the self-dualities
$$\QCoh(S_i)^\vee\simeq \QCoh(S_i),\,\, i=1,2,$$
the dual of $f_*$ is $f^*$.

\sssec{}

We claim that there is a natural action of $\QCoh(\CY)$ on $\QCoh_{\on{co}}(\CY)$. Namely, in terms of
\eqref{e:QCohco expl}, an object $\CF\in \QCoh(\CY)$ gives rise to a compatible family of endofunctors of
$\QCoh(S)$ for $y:S\to \CY$, namely
$$S\rightsquigarrow y^*(\CF)\otimes (-).$$

\medskip

This action is compatible with the identification \eqref{e:QCoh co as predual}. 

\medskip

Furthermore, it satisfies the projection formula: for $f:\CY_1\to \CY_2$ we have
$$f_*(f^*(\CF_2)\otimes \CF_1)\simeq \CF_2\otimes f_*(\CF_1), \quad \CF_1\in \QCoh_{\on{co}}(\CY_1),\,\ \CF_2\in \QCoh(\CY_2).$$

\sssec{}

We can rewrite the canonical pairing
\begin{equation} \label{e:QCoh co pairing}
\QCoh(\CY)\otimes \QCoh_{\on{co}}(\CY)\to \Vect
\end{equation}
in terms of the above action of $\QCoh(\CY)$ on $\QCoh_{\on{co}}(\CY)$. 

\medskip

Namely, it is given by
$$\QCoh(\CY)\otimes \QCoh_{\on{co}}(\CY)\overset{\text{action}}\to \QCoh_{\on{co}}(\CY) \overset{\Gamma(\CY,-)}\longrightarrow \Vect.$$

\sssec{} \label{sss:* pullback QCoh co}

Let $f:\CY_1\to \CY_2$ be affine. In this case we claim that the functor
$$f_*:\QCoh_{\on{co}}(\CY_1)\to \QCoh_{\on{co}}(\CY_2)$$
admits a left adjoint, to be denoted $f^*$.

\medskip

Indeed, the functor
$$\affSch_{/\CY_2}\to \affSch_{/\CY_1}, \quad S\mapsto S\underset{\CY_2}\times \CY_1$$
is cofinal, so the functor
$$\underset{S\in \CY_2,S\in \affSch}{\on{colim}}\, \QCoh(S\underset{\CY_2}\times \CY_1)\to \QCoh_{\on{co}}(\CY_1).$$
is an equivalence.

\medskip

In terms of this identification, the functor $f^*$ is given by the compatible family of the pullback functors
$$\QCoh(S)\to \QCoh(S\underset{\CY_2}\times \CY_1).$$


\sssec{}

Let $\CY$ be a prestack over an affine scheme $S$. Let $f:S'\to S$ be a map between affine schemes;
set $\CY':=S'\underset{S}\times \CY$. By a slight abuse of notation we will denote by the same character
$f$ the resulting map $\CY'\to \CY$. 

\medskip

The category $\QCoh_{\on{co}}(\CY)$ (resp., $\QCoh_{\on{co}}(\CY')$) is naturally tensored over $\QCoh(S)$
(resp., $\QCoh(S')$, and the functor $f_*:\QCoh_{\on{co}}(\CY')\to \QCoh_{\on{co}}(\CY)$ is $\QCoh(S)$-linear.
Hence, so is its right adjoint $f^*$. From here we obtain a functor
\begin{equation} \label{e:QCoh co base change}
\QCoh(S')\underset{\QCoh(S)}\otimes \QCoh_{\on{co}}(\CY)\to \QCoh_{\on{co}}(\CY').
\end{equation}

We claim:

\begin{lem} \label{l:QCoh co base change}
The functor \eqref{e:QCoh co base change} is an equivalence.
\end{lem} 

\begin{proof}

Follows from the fact that the functor
$$\wt{S}\in \affSch_{/\CY} \mapsto S'\underset{S}\times \wt{S}\in \affSch_{/\CY'}$$
is cofinal.

\end{proof}

\sssec{}

Let $\CY$ have an affine diagonal. In this case we claim that there is a naturally defined functor 
\begin{equation} \label{e:from QCoh co to QCoh}
\Omega_\CY:\QCoh_{\on{co}}(\CY)\to \QCoh(\CY).
\end{equation}

Namely, in terms of \eqref{e:QCohco expl}, the functor \eqref{e:from QCoh co to QCoh} is given by the (compatible 
family) of direct image functors 
$$\QCoh(S)\to \QCoh(\CY),$$
which are well-defined, since the morphisms $S\to \CY$ are affine.

\medskip

Let $f:\CY_1\to \CY_2$ be a schematic map. Note that, by construction, the following diagram commutes:
$$
\CD
\QCoh_{\on{co}}(\CY_1) @>{f_*}>> \QCoh_{\on{co}}(\CY_2)  \\
@V{\Omega_{\CY_1}}VV @VV{\Omega_{\CY_2}}V \\
\QCoh(\CY_1) @>{f_*}>> \QCoh(\CY_2).
\endCD
$$

\sssec{}

The following assertion is established in \cite[Theorems 2.2.4 or 2.2.6 and Proposition 6.3.8]{Ga5}:

\begin{thm} \label{t:alg stack QCoh co} Let $\CY$ be an quasi-compact algebraic stack 
with an affine diagonal. Suppose that one of the following conditions holds:

\smallskip

\noindent{\em(i)} $\CY$ can be realized as a quotient of an algebraic space by an action of
a (finite-dimensional) algebraic group;

\smallskip

\noindent{\em(ii)} $\CY$ is eventually coconnective algebraic stack almost of finite type.

\medskip

Then the functor $\Omega_\CY$ of \eqref{e:from QCoh co to QCoh} is an equivalence.
\end{thm}

\sssec{} \label{sss:QCoh* on alg stacks}

Let $\CY$ be a (not necessarily quasi-compact) algebraic stack. Suppose that $\CY$
can be written as a union of quasi-compact open substacks $\CY_i$ that satisfy one
of the conditions in \thmref{t:alg stack QCoh co}. 

\medskip

\begin{cor} \label{c:QCoh co alg stacks}
Under the above circumstances, we have a canonical equivalence
$$\QCoh_{\on{co}}(\CY) \simeq \underset{i}{\on{colim}}\, \QCoh(\CY_i),$$
where the colimit is taken with respect to the pushforward functors.  
\end{cor}

\begin{proof}

Note that the map 
$$\underset{i}{\on{colim}}\, \CY_i \to \CY$$
is an isomorphism in $\on{PreStk}$. 

\medskip

Hence, the functor
$$\underset{i}{\on{colim}}\, \QCoh_{\on{co}}(\CY_i) \to \QCoh_{\on{co}}(\CY)$$
is an equivalence.

\medskip

Now the assertion follows from \thmref{t:alg stack QCoh co}.

\end{proof}

\ssec{The category \texorpdfstring{$\QCoh_{\on{co}}$}{QCohco} on ind-schemes} \label{ss:QCoh co IndSch}

\sssec{}  \label{sss:QCoh co IndSch}

Let $\CY$ be an ind-affine ind-scheme (see \cite{GaRo1}). According to Corollary 1.6.6 in {\it loc. cit.},
the map
$$\underset{S\,\text{closed in}\,\CY}{\on{colim}}\, S\to \CY$$
is an isomorphism in $\on{PreStk}$, where the index category is that of affine schemes, equipped with a closed embedding to $\CY$.

\medskip

Hence, in this case, we have
\begin{equation} \label{e:QCoh co IndSch}
\QCoh_{\on{co}}(\CY)\simeq \underset{S\,\text{closed in}\,\CY}{\on{colim}}\, \QCoh(S).
\end{equation}

\sssec{} 

Recall the following general paradigm:

\medskip

Let
$$i \mapsto \bC_i, \quad i\in I$$
be a diagram in $\DGCat$. Denote
$$\bC:=\underset{i\in I}{\on{colim}}\, \bC_i,$$
where, per our conventions, the colimit is taken in $\DGCat$ (i.e., the category of cocomplete 
DG categories and continuous functors). 

\medskip 

Let $\on{ins}_i:\bC_i\to \bC$ denote the tautological functors. 

\sssec{} \label{sss:t-structure colimit}

Suppose that each $\bC_i$ is equipped with a t-structure, compatible with filtered colimits, i.e., $\bC_i^{\geq 0}$
is closed under filtered colimits. And suppose that the transition functors
$$\bC_i \overset{F_{i,j}}\to \bC_j$$
are t-exact.

\medskip

We equip with $\bC$ with a t-structure by declaring that  $\bC^{\leq 0}$ is generated under colimits
by the essential images of $\on{ins}_i(\bC_i^{\leq 0})$.

\medskip

So by construction, the functors $\on{ins}_i$ are right t-exact. 

\sssec{}

We claim:

\begin{lem} \label{l:colimit t-structures} 
Assume that $I$ is filtered. Then:

\smallskip

\noindent{\em(a)} The t-structure on $\bC$ is compatible with filtered colimits.

\smallskip

\noindent{\em(b)} The functors $\on{ins}_i$ are t-exact.

\smallskip

\noindent{\em(c)} The category $\bC^{\geq 0}$ is generated under filtered colimits
by the essential images of $\bC_i^{\geq 0}$ along the functors $\on{ins}_i$.  

\end{lem}

\begin{proof}

The first two points are proved in \cite[Proposition C.3.3.5]{Lu3}. 

\medskip 

Namely, in the notation of loc. cit., the category $\on{Groth}^{\on{lex}}$ is equivalent to the category of presentable stable $\infty$-categories with right-complete t-structures and colimit preserving functors which are t-exact.  The first two points are equivalent to the assertion that $\on{Groth}^{\on{lex}}$ admits filtered colimits and the forgetful functor
$$ \on{Groth}^{\on{lex}} \to \on{Pr^L} $$
given by $(\bC, \bC^{\leq 0}) \mapsto \bC$ preserves filtered colimits.  By \cite[Proposition C.3.3.5]{Lu3}, $\on{Groth}^{\on{lex}}$ admits filtered colimits and the functor $(\bC, \bC^{\leq 0}) \mapsto \bC^{\leq 0}$ preserves filtered colimits.  The result now follows from the fact that for any $(\bC, \bC^{\leq 0}) \in \on{Groth}^{\on{lex}}$, we have
$$ \bC \simeq \on{Stab}(\bC^{\leq 0})$$
is given by the stabilization, and the functor $\on{Stab}: \on{Pr^L} \to \on{Pr^L}$ preserves colimits. 

\medskip

To prove the third point we note that
any $\bc\in \bC$ is canonically isomorphic to
$$\underset{i\in I}{\on{colim}}\, \on{ins}_i\circ \on{ins}_i^R(\bc).$$

If $\bc\in \bC^{\geq 0}$, then so are all $\on{ins}_i^R(\bc)$ (since the functors $\on{ins}_i^R$ are left t-exact,
being right adjoints of (right) t-exact functors. 

\end{proof} 

\begin{rem}

Note that in the situation of \lemref{l:colimit t-structures}, we have
$$\bC\simeq \underset{i\i I^{\on{op}}}{\on{lim}}\, \bC_i,$$
where:

\begin{itemize}

\item The limit is taken in the category of $\infty$-categories;

\item The functor $\to$ is given by the (compatible collection of) the functors $\on{ins}_i^R$. 

\end{itemize} 

Since the functors $\on{ins}_i^R$ send $\bC^{\geq 0}$ to $\bC_i^{\geq 0}$, we obtain they also induce
an equivalence
\begin{equation} \label{e:colimit coconn}
\bC^{\geq 0}\simeq \underset{i\i I^{\on{op}}}{\on{lim}}\, \bC^{\geq 0}_i.
\end{equation} 

\end{rem}

\begin{cor} \label{c:t-exact out of colimit}
Let $\Phi:\bC\to \bD$ be a continuous functor, where $\bD$ is also equipped with a t-structure.

\medskip

\noindent{\em(a)} The functor $\Phi$ is right t-exact if and only if each $\Phi\circ \on{ins}_i=:\Phi_i:\bC_i\to \bD$
is right t-exact.

\medskip

\noindent{\em(b)} If $\Phi$ is left t-exact, then so is each $\Phi_i$.

\medskip

\noindent{\em(b')} Suppose that the t-structure on $\bD$ is compatible with filtered colimits.
Then the assertion in (b) is ``if and only if".

\end{cor}

\sssec{} \label{sss:t on QCoh co}

Let $\CY$ be an ind-affine ind-scheme. We use the presentation \eqref{e:QCoh co IndSch} and the construction in \secref{sss:t-structure colimit}
to equip the category $\QCoh_{\on{co}}(\CY)$ with a t-structure:

\medskip

By definition, $\QCoh_{\on{co}}(\CY)^{\leq 0}$ is generated under colimits by the essential images of $\QCoh(S)^{\leq 0}$ for
$S\in \affSch_{/\CY}$. 

\sssec{}

From \lemref{l:colimit t-structures}, we obtain:

\begin{cor} \label{c:t-structure on QCoh co} \hfill

\smallskip

\noindent{\em(a)} The t-structure on $\QCoh_{\on{co}}(\CY)$ is compatible with filtered colimits.

\smallskip

\noindent{\em(b)} For every $S\in \affSch_{/\CY}$, the direct image functor $\QCoh(S)\to \QCoh_{\on{co}}(\CY)$
is t-exact. 

\end{cor}

\sssec{}

Let $f:\CY_1\to \CY_2$ be a map between ind-affine ind-schemes. It follows by definition that the functor
$$f_*:\QCoh_{\on{co}}(\CY_1)\to \QCoh_{\on{co}}(\CY_2)$$
is right t-exact. 

\bigskip

However, from Corollaries \ref{c:t-structure on QCoh co}(b) and \ref{c:t-exact out of colimit} we obtain:

\begin{cor}  \label{c:dir im QCoh co t exact}
The functor $f_*$ is t-exact.
\end{cor}

\sssec{}

Assume that $f$ is affine, in which case we have a well-defined functor
$$f^*:\QCoh_{\on{co}}(\CY_2)\to \QCoh_{\on{co}}(\CY_1).$$

\bigskip

From \corref{c:dir im QCoh co t exact} we obtain: 

\begin{cor}  \label{c:inv im QCoh co t right exact}
The functor $f^*$ is right t-exact.
\end{cor}

\bigskip

Finally, assume that $f$ is flat. In this case, unwinding the construction of $f^*$ in \secref{sss:* pullback QCoh co}
and using Corollaries \ref{c:t-structure on QCoh co}(b) and \ref{c:t-exact out of colimit}, we obtain: 
 
\begin{lem} \label{c:inv im QCoh co t exact}
For a flat map $f$ between ind-affine ind-schemes, 
the functor $f^*:\QCoh_{\on{co}}(\CY_2)\to \QCoh_{\on{co}}(\CY_1)$
is t-exact.
\end{lem}

\ssec{A descent property of \texorpdfstring{$\QCoh_{\on{co}}(-)$}{QCohco}}

In this subsection we will prove a certain technical assertion used in the main body of the text.

\sssec{}

Let $\CY$ be an ind-affine ind-scheme. Let $g:\wt\CY\to \CY$ be a map of prestacks that is an affine fpqc cover,
and let $\wt\CY^\bullet$ denote its \v{C}ech nerve.  

\medskip

Consider $\QCoh_{\on{co}}(\wt\CY^\bullet)$ as a cosimplicial category, equipped with an augmentation by $\QCoh_{\on{co}}(\CY)$
using *-pullbacks (they are well-defined since the maps involved are affine, see \secref{sss:* pullback QCoh co}).

\medskip

Thus, we obtain a functor
\begin{equation} \label{e:QCoh co descent}
\QCoh_{\on{co}}(\CY)\to \on{Tot}(\QCoh_{\on{co}}(\wt\CY^\bullet)).
\end{equation} 

\sssec{}

From now on we will perceive $\QCoh_{\on{co}}(\wt\CY^\bullet)$ as a \emph{semi-}cosimplicial category,
so that transition functors involved are t-exact (by the flatness assumption on $g$, see \lemref{c:inv im QCoh co t exact}).
Hence, the functor \eqref{e:QCoh co descent} induces a functor  
\begin{equation} \label{e:QCoh co descent plus}
\QCoh_{\on{co}}(\CY)^{>-\infty}\to \on{Tot}(\QCoh_{\on{co}}(\wt\CY^\bullet)^{>-\infty}).
\end{equation} 

\medskip

We will prove:

\begin{prop} \label{p:QCoh co descent}
Suppose that $\CY$ can be exhibited as a filtered colimit in $\on{PreStk}$ of 
affine schemes with transition maps that are almost finitely presented.\footnote{I.e., finitely presented after
each coconnective truncation.} Then the functor \eqref{e:QCoh co descent plus} is an equivalence.
\end{prop}

The rest of this subsection is devoted to the proof of this proposition. 

\sssec{}

It suffices to show that 
\begin{equation} \label{e:QCoh co descent >0}
\QCoh_{\on{co}}(\CY)^{\geq 0}\to \on{Tot}(\QCoh_{\on{co}}(\wt\CY^\bullet)^{\geq 0})
\end{equation} 
is an equivalence. 

\sssec{}

Let 
$$\CY=\underset{i\in I}{\on{colim}}\, Y_i$$
be the presentation of $\CY$ as in the statement of the proposition. Denote the map
$Y_i\to Y_j$ for $(i\to j)\in I$ by $f_{i,j}$. 

\medskip 

Set $$\wt{Y}^\bullet_i:=\wt\CY^\bullet\underset{\CY}\times Y_i.$$ 
Then for every $m$, we also have
$$\wt\CY^m\simeq \underset{i}{\on{colim}}\, \wt{Y}^m_i.$$
Denote the corresponding maps
$\wt{Y}^m_i\to \wt{Y}^m_j$ by $f^m_{i,j}$. 

\sssec{}

We have
$$\QCoh_{\on{co}}(\CY)\simeq \underset{i\in I}{\on{colim}}\, \QCoh(Y_i)$$
and 
$$\QCoh_{\on{co}}(\wt\CY^m)\simeq \underset{i\in I}{\on{colim}}\, \QCoh(\wt{Y}^m_i),$$
where in both cases the colimit is taken with respect to the pushforward functors (recall that by default,
colimits are taken in $\DGCat$, i.e., in the $\infty$-category of cocomplete DG categories and continuous functors). 

\medskip

We can rewrite the above colimits as \emph{limits} (in the category of DG categories with not necessarily continuous functors)
$$\QCoh_{\on{co}}(\CY)\simeq \underset{i\in I^{\on{op}}}{\on{lim}}\, \QCoh(Y_i)$$
and 
$$\QCoh_{\on{co}}(\wt\CY^m)\simeq \underset{i\in I^{\on{op}}}{\on{lim}}\, \QCoh(\wt{Y}^m_i),$$
where the transition functors are $f_{i,j}^!:=(f_{i,j})_*^R$ and $(f^m_{i,j})^!:=(f^m_{i,j})_*^R$,
respectively (note that these right adjoints are indeed in general \emph{discontinuous}). 

\medskip

From here we obtain:
$$\QCoh_{\on{co}}(\CY)^{\geq 0}\simeq \underset{i\in I^{\on{op}}}{\on{lim}}\, \QCoh(Y_i)^{\geq 0}$$
and 
$$\QCoh_{\on{co}}(\wt\CY^m)^{\geq 0}\simeq \underset{i\in I^{\on{op}}}{\on{lim}}\, \QCoh(\wt{Y}^m_i)^{\geq 0},$$
where limits are taken in the category of $\infty$-categories and all functors, see \eqref{e:colimit coconn}. 

\sssec{}

For every $\phi:[m']\to [m'']$ in\footnote{We remind that we only consider injective maps $\phi$.} 
$\bDelta$ denote by $g^\phi$ the corresponding map
$$\wt\CY^{m''}\to \wt\CY^{m'}.$$
For every index $i$, let $g^\phi_i$ denote the resulting map $\wt{Y}_i^{m''}\to \wt{Y}_i^{m'}$. 

\medskip

For every arrow $(i\to j)\in I$, we have a Cartesian diagram of schemes
$$
\CD
\wt{Y}_i^{m''} @>{f^{m''}_{i,j}}>> \wt{Y}_j^{m''} \\
@V{g^\phi_i}VV @VV{g^\phi_j}V \\
\wt{Y}_i^{m'} @>{f^{m'}_{i,j}}>> \wt{Y}_j^{m'},
\endCD
$$
which gives rise to a commutative diagram
$$
\CD
\QCoh(\wt{Y}_i^{m''}) @>{(f^{m''}_{i,j})_*}>> \QCoh(\wt{Y}_j^{m''}) \\
@A{(g^\phi_i)^*}AA @AA{(g^\phi_j)^*}A \\
\QCoh(\wt{Y}_i^{m'}) @>{(f^{m'}_{i,j})_*}>> \QCoh(\wt{Y}_j^{m'}). 
\endCD
$$

From here we obtain a natural transformation
\begin{equation} \label{e:BC indsch}
(g^\phi_i)^*\circ (f^{m'}_{i,j})^!\to (f^{m''}_{i,j})^! \circ (g^\phi_j)^*.
\end{equation}

Now, the assumption that the maps $f_{i,j}$ are almost of finite presentation and
the maps $g^\phi$ are flat implies that the natural transformations \eqref{e:BC indsch} are
isomorphisms when evaluated on $\QCoh(\wt{Y}_j^{m'})^{\geq 0}$. 

\medskip

Hence, we obtain a family of commutative diagrams
$$
\CD
\QCoh(\wt{Y}_i^{m''})^{\geq 0} @<{(f^{m''}_{i,j})^!}<< \QCoh(\wt{Y}_j^{m''})^{\geq 0} \\
@A{(g^\phi_i)^*}AA @AA{(g^\phi_j)^*}A \\
\QCoh(\wt{Y}_i^{m'})^{\geq 0} @<{(f^{m'}_{i,j})^!}<< \QCoh(\wt{Y}_j^{m'})^{\geq 0}. 
\endCD
$$

\sssec{}

Thus, we obtain a well-defined functor from $\bDelta\times I$ to the category of $\infty$-categories 
$$m,i\mapsto \QCoh(\wt{Y}_i^m)^{\geq 0},$$
and we can rewrite
$$\on{Tot}(\QCoh_{\on{co}}(\wt\CY^\bullet)^{\geq 0}):=
\underset{m\in \bDelta}{\on{lim}}\, \underset{i\in I^{\on{op}}}{\on{lim}}\, \QCoh(\wt{Y}_i^m)^{\geq 0}$$
as 
$$\underset{(m,i)\in \bDelta\times I^{\on{op}}}{\on{lim}}\,  \QCoh(\wt{Y}_i^m)^{\geq 0}$$
and further as
$$\underset{i\in I^{\on{op}}}{\on{lim}}\, 
\underset{m\in \bDelta}{\on{lim}}\, \QCoh(\wt{Y}_i^m)^{\geq 0}.$$

Unwinding the construction, we obtain that the following diagram commutes
\begin{equation} \label{e:compare Tot}
\CD
\on{Tot}(\QCoh_{\on{co}}(\wt\CY^\bullet)^{\geq 0}) @>{\sim}>> \underset{i\in I^{\on{op}}}{\on{lim}}\, 
\underset{m\in \bDelta}{\on{lim}}\, \QCoh(\wt{Y}_i^m)^{\geq 0} \\
@AAA @AAA \\
\QCoh_{\on{co}}(\wt\CY)^{\geq 0} @>{\sim}>> \underset{i\in I^{\on{op}}}{\on{lim}}\, \QCoh(Y_i)^{\geq 0},
\endCD
\end{equation} 
where the right vertical arrow is comprised of the functors
\begin{equation} \label{e:usual fpqc}
\QCoh(Y_i)^{\geq 0}\to \on{Tot}(\QCoh_{\on{co}}(\wt{Y}_i^\bullet)^{\geq 0}).
\end{equation} 

Now, the functors \eqref{e:usual fpqc} are equivalences by the usual fpqc descent. Hence, the right vertical
arrow in \eqref{e:compare Tot} is an equivalence. 

\medskip

Hence, the left vertical arrow is also an equivalence, as required. 

\qed[\propref{p:QCoh co descent}]

\ssec{The category \texorpdfstring{$\IndCoh^!(-)$}{IndCoh!}} \label{sss:IndCoh non-placid !}

\sssec{}

Let $^{\leq n}\!\affSch$ denote the category of $n$-coconnective
affine schemes, i.e., 
$$^{\leq n}\!\affSch=(\on{ComAlg}(\Vect^{\geq -n,\leq 0}))^{\on{op}}.$$

\medskip

We have 
$$\affSch\simeq \underset{n}{\on{lim}}\, ({}^{\leq n}\!\affSch).$$

\sssec{}

Let 
$$^{\leq n}\!\affSch_{\on{ft}}\subset {}^{\leq n}\!\affSch$$
be the full subcategory consisting of $n$-coconnective affine schemes \emph{of finite type}
(see \cite[Chapter 1, Sect. 1.5]{GaRo3}).

\medskip

Note that 
$$^{\leq n}\!\affSch\simeq \on{Pro}(^{\leq n}\!\affSch_{\on{ft}}).$$

\medskip

Let $\affSch_{\on{aft}}$ denote the full subcategory of $\affSch$ consisting of affine schemes 
almost of finite type, which is by definition
$$\underset{n}{\on{lim}}\, ({}^{\leq n}\!\affSch_{\on{ft}}).$$

\medskip

Let $\on{PreStk}_{\on{aft}}\subset \on{PreStk}$ the full subcategory consisting of prestacks 
locally almost of finite type. We have
$$\on{PreStk}_{\on{aft}}\simeq \underset{n}{\on{lim}}\, ({}^{\leq n}\!\on{PreStk}_{\on{lft}}),$$
where 
$$^{\leq n}\!\on{PreStk}_{\on{lft}}\simeq \on{Funct}(({}^{\leq n}\!\affSch_{\on{ft}})^{\on{op}},\inftygroup).$$

\sssec{} 

We define the functor 
$$^{\leq n}\IndCoh^!:({}^{\leq n}\!\affSch)^{\on{op}}\to \DGCat,$$
to be the left Kan extension of the functor
$$\IndCoh:({}^{\leq n}\!\affSch_{\on{ft}})^{\on{op}}\to \DGCat$$
along
$$({}^{\leq n}\!\affSch_{\on{ft}})^{\on{op}}\hookrightarrow ({}^{\leq n}\!\affSch)^{\on{op}}.$$

\medskip

Explicitly, for $S\in {}^{\leq n}\!\affSch$, we have
\begin{equation} \label{e:IndCoh! sch}
\IndCoh^!(S)=\underset{S\to S_0,S_0\in {}^{\leq n}\!\affSch_{\on{ft}}}{\on{colim}}\, \IndCoh(S_0),
\end{equation} 
where the transition functors are given by 
$$(S\to S'_0\overset{f}\to S''_0)\rightsquigarrow \IndCoh(S''_0)\overset{f^!}\to \IndCoh(S'_0).$$

\sssec{} \label{sss:IndCoh ! stacks}

It is easy to see that the natural transformation from
$$({}^{\leq n}\!\affSch)^{\on{op}} \overset{^{\leq n}\IndCoh^!}\longrightarrow \DGCat$$
to
$$({}^{\leq n}\!\affSch)^{\on{op}} \hookrightarrow ({}^{\leq n+1}\!\affSch)^{\on{op}}  
\overset{^{\leq n+1}\!\IndCoh^!}\longrightarrow \DGCat$$
is an isomorphism. Hence, we obtain a well-defined functor
\begin{equation} \label{e:IndCoh ! event coconn}
\IndCoh^!:({}^{<\infty}\!\affSch)^{\on{op}}\to \DGCat,
\end{equation}
where
$$^{<\infty} \!\affSch=\underset{n}{\on{colim}}\, {}^{\leq n}\!\affSch.$$

\sssec{}

We define the functor 
$$\IndCoh^!:(\on{PreStk})^{\on{op}}\to \DGCat$$
to be the right Kan extension of \eqref{e:IndCoh ! event coconn} along the embedding
$$({}^{<\infty}\!\affSch)^{\on{op}}\hookrightarrow (\on{PreStk})^{\on{op}}.$$

Explicitly, for $\CY\in \on{PreStk}$, we have
$$\IndCoh^!(\CY)=\underset{S\to \CY,S\in {}^{<\infty}\!\affSch}{\on{lim}}\, \IndCoh^!(S),$$
where the transition functors are given by
$$(S'\overset{f}\to S''\to \CY) \, \rightsquigarrow\, \IndCoh^!(S'')\overset{f^!}\to \IndCoh^!(S').$$

\sssec{}

Thus, by definition, for a map $f:\CY_1\to \CY_2$ in $\on{PreStk}$, we have a well-defined functor
$$f^!:\IndCoh^!(\CY_2)\to \IndCoh^!(\CY_1).$$

In particular, taking the projection $\CY\to \on{pt}$, we obtain that for any $\CY$, we have a well-defined 
object 
$$\omega_\CY\in \IndCoh^!(\CY).$$

\sssec{}

It follows from the convergence property of the usual $\IndCoh$ functor 
$$\IndCoh:(\on{PreStk}_{\on{laft}})^{\on{op}}\to \DGCat$$
(see \cite[Chapter 4, Prop. 6.4.3]{GaRo3})
that the natural transformation 
$$\IndCoh^!|_{\on{PreStk}_{\on{laft}}}\to \IndCoh$$ is
an isomorphism. 

\medskip

I.e., the value of $\IndCoh^!(\CY)$ on a prestack locally almost of finite type recovers
the usual $\IndCoh(\CY)$. 

\begin{rem}
The above construction gives a definition of the functor $\IndCoh^!(\CY)$ for a general prestack $\CY$.
But unless some conditions on $\CY$ are imposed, we will not be able to say much about the 
properties of this category. 

\medskip

For example, it is not even clear (and, probably, not true) whether for $S\in \affSch$, the category $\IndCoh(S)$ is dualizable. 

\medskip

A condition on $\CY$ that makes $\IndCoh^!(\CY)$ manageable is called ``placidity", to be discussed in
\secref{ss:placid}. 

\end{rem} 

\sssec{} \label{sss:IndCoh^! tensor up}

In the sequel we will need the following property of $\IndCoh^!$:

\medskip

Let $\CY$ be a prestack mapping to a \emph{smooth} affine scheme $S$ of finite type. 
Let $f:S'\to S$ be a map, where $S'\in \affSch_{\on{aft}}$. Denote
$$\CY':=S'\underset{S}\times \CY.$$
By a slight abuse of notation, we will denote by the same symbol $f$ the resulting map $\CY'\to \CY$. 

\medskip

The functor $f^!:\IndCoh^!(\CY)\to \IndCoh^!(\CY')$ extends to a functor
\begin{equation} \label{e:IndCoh^! tensor up}
\QCoh(S')\underset{\QCoh(S)}\otimes \IndCoh^!(\CY)\to \IndCoh^!(\CY').
\end{equation}

We claim:

\begin{lem}  \label{l:IndCoh^! tensor up}
The functor \eqref{e:IndCoh^! tensor up} is fully faithful. If $f$ is smooth, it is an equivalence.
\end{lem}

\begin{proof}

Follows from \cite[Propositions 4.4.2 and 7.5.7]{Ga7}. 

\end{proof}

\sssec{}

In the sequel we will need the following assertion about the behavior of $\IndCoh^!(-)$. Let
$$
\CD
\CY_1 @>{f_Y}>> \CY_2 \\
@VVV @VVV \\
S_1 @>{f_S}>> S_2
\endCD
$$ 
be a fiber square, where $S_i$ are affine schemes almost of finite type, and $f_S$ is a closed embedding of finite
Tor-dimension. Unwinding the definitions, we obtain:

\begin{lem}
Under the above circumstances, the functor $f_Y^!:\IndCoh^!(\CY_2)\to \IndCoh^!(\CY_1)$ admits a left adjoint,
to be denoted $(f_Y)^\IndCoh_*$. Furthermore, the functor $(f_Y)^\IndCoh_*$ satisfies base change for any fiber
square
$$
\CD
\CY'_1 @>{f'_Y}>> \CY'_2 \\
@VVV @VVV \\
\CY_1 @>{f_Y}>> \CY_2.
\endCD
$$ 
\end{lem}

\sssec{}

Let $S$ be an affine scheme almost of finite type, and let $S'\subset S$ we a Zariski-closed subset.
Let $S^\wedge$ denote the formal completion of $S$ along $S'$. Let $\CY$ be a prestack over $S$, and
set $\CY^\wedge:=S^\wedge\underset{S}\times \CY$:
$$
\CD
\CY^\wedge @>{i_Y}>> \CY \\
@VVV @VVV \\
S^\wedge @>{i_S}>> S.
\endCD
$$

We now claim: 

\begin{prop} \label{p:form compl IndCoh!}
Under the above circumstances, the functor $i_Y^!:\IndCoh^!(\CY)\to \IndCoh^!(\CY^\wedge)$ is a colocalization,
i.e., it admits a fully faithful left adjoint.  The essential image of this left adjoint
is the full subcategory of $\IndCoh^!(\CY)$ consisting of objects with set-theoretic support on $\CY':=S'\underset{S}\times \CY$.
The formation of the left adjoint satisfies base change for any fiber
square
$$
\CD
\wt\CY^\wedge_1 @>>> \wt\CY \\
@VVV @VVV \\
\wt\CY @>>> \CY.
\endCD
$$ 
\end{prop} 

\begin{proof}

By \cite[Proposition 6.7.4]{GaRo1}, we can write $S^\wedge$ as $\underset{i}{``\on{colim}"}\, S_i$, where $S_i$ are closed subschemes of
$S$, and the maps $S_i\to S$ are of finite Tor-dimension. 

\medskip

Unwinding the definition of $\IndCoh^!(-)$, we reduce the assertion to the case when $\CY$ is an eventually 
connective affine scheme $\wt{S}$, so that 
$$\wt{S}^\wedge\simeq \underset{i}{``\on{colim}"}\, \wt{S}_i, \quad \wt{S}_i:=S_i\underset{S}\times \wt{S};$$
note that all $\wt{S}_i$ are eventually coconnective. 

\medskip

Unwinding further, we reduce the assertion to the case when $\wt{S}$ is of finite type; in this case the assertion
follows from \cite[Proposition 7.4.5]{GaRo1}.

\end{proof}

\ssec{The category \texorpdfstring{$\IndCoh^*(-)$}{IndCoh*}} \label{sss:IndCoh non-placid *}

\sssec{}

We define the functor
$$^{\leq n}\IndCoh^*:{}^{\leq n}\!\affSch\to \DGCat,$$
to be 
$${}^{\leq n}\!\affSch \overset{({}^{\leq n}\IndCoh^!)^{\on{op}}}\longrightarrow (\DGCat)^{\on{op}}\to \DGCat,$$
where:

\begin{itemize}

\item The first arrow is the opposite of the functor
$$S\mapsto \IndCoh^!(S), \quad (S_1 \overset{f}\to S_2) \rightsquigarrow (\IndCoh^!(S_2)\overset{f^!}\to \IndCoh^!(S_1));$$

\item The second arrow is
$$\bD\mapsto \bD^\vee:=\on{Funct}_{\on{cont}}(\bD,\Vect)$$
(note that in the above formula, the DG category $\bD$ is not assumed dualizable).

\end{itemize}

\sssec{}

Explicitly, for $S\in {}^{\leq n}\!\affSch$, we have:
\begin{equation} \label{e:IndCoh* as limit}
^{\leq n}\IndCoh^*(S)=\underset{S\to S_0,S_0\in {}^{\leq n}\!\affSch_{\on{ft}}}{\on{lim}}\, \IndCoh(S_0),
\end{equation} 
where the transition functors are given by 
$$(S\to S'_0\overset{f}\to S''_0)\rightsquigarrow \IndCoh(S'_0)\overset{f^\IndCoh_*}\to \IndCoh(S''_0).$$

\sssec{}

As in \secref{sss:IndCoh ! stacks}, it is easy to see that the collection 
$$n\rightsquigarrow {}^{\leq n}\IndCoh^*$$
gives rise to a well-defined functor
\begin{equation} \label{e:IndCoh * event coconn}
\IndCoh^*:{}^{<\infty}\!\affSch\to \DGCat.
\end{equation}

\sssec{} \label{sss:! predual of *}

Note that, \emph{by construction}, for $S\in {}^{<\infty}\!\affSch$, the category 
$\IndCoh^!(S)$ is naturally a pre-dual of $\IndCoh^*(S)$. 

\medskip

This will be a perfect duality if $S$ is placid, see \secref{sss:! * duality placid} below.  

\sssec{} \label{sss:dual of ^!}

Let $f:S_1\to S_2$ be a morphism between eventually coconnective affine schemes. Unwinding the definitions,
we obtain that with respect to the identifications
$$\IndCoh^!(S_i)^\vee \simeq \IndCoh^*(S_i), \quad i=1,2,$$
we have
$$(f^!)^\vee \simeq f^\IndCoh_*.$$

\sssec{} \label{sss:IndCoh *}

Unlike $\IndCoh^!$, we do not even attempt to define $\IndCoh^*$ on all prestacks. Rather, we define
it on \emph{ind-affine ind-schemes} (see \cite[Chap. 3.1]{GaRo4}).

\medskip

Namely, we let the functor 
$$\IndCoh^*:\IndSch\to \DGCat$$
to be the left Kan extension of \eqref{e:IndCoh * event coconn} 
along the embedding
$${}^{<\infty}\!\affSch\hookrightarrow \IndSch^{\on{ind-aff}}.$$

\sssec{}

Explicitly, for $\CY\in \IndSch^{\on{ind-aff}}$, we have
\begin{equation} \label{e:IndCoh* indsch}
\IndCoh^*(\CY)=\underset{S}{\on{colim}}\, \IndCoh^*(S),
\end{equation} 
where:

\begin{itemize}

\item The index category is that of $S\in {}^{<\infty}\!\affSch$ equipped with a closed embedding
$S\to \CY$; 

\smallskip

\item The transition functors are given by
$$(S'\overset{f}\to S''\to \CY) \, \rightsquigarrow\, \IndCoh^*(S')\overset{f^\IndCoh_*}\to \IndCoh^*(S'').$$

\end{itemize} 

\sssec{}

Thus, by definition, for a map $f:\CY_1\to \CY_2$ in $\IndSch^{\on{ind-aff}}$, we have a well-defined functor
$$f^\IndCoh_*:\IndCoh^*(\CY_1)\to \IndCoh^*(\CY_2).$$

\medskip

In particular, taking the projection $\CY\to \on{pt}$, we obtain that for any $\CY\in \IndSch^{\on{ind-aff}}$
there is a well-defined functor 
$$\Gamma^\IndCoh(\CY,-):\IndCoh^*(\CY)\to \Vect.$$

\sssec{}

It follows from \cite[Sect. 2.4.2]{GaRo1} that if $\CY\in \IndSch^{\on{ind-aff}}_{\on{laft}}$, the naturally defined functor
$$\IndCoh(\CY)\to \IndCoh^*(\CY)$$
is an equivalence.

\ssec{The multiplicative structure}

\sssec{}

Note that since the index category in \eqref{e:IndCoh! sch} is filtered (and, in particular, sifted), for 
$S\in {}^{\leq n}\!\affSch$, the category $\IndCoh^!(S)$ carries a naturally defined monoidal structure. 

\medskip

Explicitly, the corresponding binary operation is given by
$$\IndCoh^!(S)\otimes \IndCoh^!(S)\to \IndCoh^!(S\times S) \overset{\Delta_S^!}\to \IndCoh^!(S).$$

\medskip

In other words, the functor
$$\IndCoh^!:({}^{<\infty}\!\affSch)^{\on{op}}\to \DGCat$$
lifts to a functor
$$({}^{<\infty}\!\affSch)^{\on{op}}\to \on{ComAlg}(\DGCat)=\DGCat^{\on{SymMon}}.$$

\sssec{}

By construction, we obtain that the functor 
$$\IndCoh^!: (\on{PreStk})^{\on{op}}\to \DGCat$$
also lifts to a functor
$$(\on{PreStk})^{\on{op}}\to \on{ComAlg}(\DGCat)=\DGCat^{\on{SymMon}},$$
i.e., for every $\CY\in \on{PreStk}$, the category $\IndCoh^!(\CY)$ has a naturally defined symmetric monoidal structure. 

\medskip

Namely, the corresponding binary operation is given by
$$\IndCoh^!(\CY)\otimes \IndCoh^!(\CY)\to \IndCoh^!(\CY\times \CY) \overset{\Delta_\CY^!}\to \IndCoh^!(\CY).$$

The unit for this symmetric monoidal structure is the object $\omega_\CY$. 

\sssec{} \label{sss:boxtimes IndCoh !}

Let $\CY_1$ and $\CY_2$ be a pair of prestacks. The operation of pullback and tensor product gives rise to a functor
\begin{equation} \label{e:boxtimes IndCoh !}
\IndCoh^!(\CY_1)\otimes \IndCoh^!(\CY_2)\to \IndCoh^!(\CY_1\times \CY_2).
\end{equation}

For general prestacks there is no reason for \eqref{e:boxtimes IndCoh !} to be an equivalence. 

\sssec{} \label{sss:boxtimes IndCoh *}

Let $S_1$ and $S_2$ be a pair of eventually coconnective schemes. Given an eventually coconnective scheme $S$
of finite type and a map $S_1\times S_2\to S$, the category of factorizations of this map as 
$$S_1\times S_2 \to S_{1,0}\times S_{2,0}\to S$$
is contractible, where:

\begin{itemize}

\item $S_{i,0}$ are eventually coconnective;

\item The first arrow comes from a pair of maps $S_i\to S_{i,0}$.

\end{itemize}

This implies that we have a well-defined functor 
$$\IndCoh^*(S_1)\otimes \IndCoh^*(S_2)\to \IndCoh(S).$$

Passing to the limit over $S$, we obtain a functor 
\begin{equation} \label{e:boxtimes IndCoh * aff}
\IndCoh^*(S_1)\otimes \IndCoh^*(S_2)\to \IndCoh^*(S_1\times S_2).
\end{equation}

Ind-extending, we obtain a functor
\begin{equation} \label{e:boxtimes IndCoh *}
\IndCoh^*(\CY_1)\otimes \IndCoh^*(\CY_2)\to \IndCoh^*(\CY_1\times \CY_2),
\end{equation}
where $\CY_1$ and $\CY_2$ are ind-schemes. 

\medskip

For general ind-schemes there is no reason for \eqref{e:boxtimes IndCoh *} to be an equivalence. 

\sssec{}

By a similar principle, we obtain a symmetric monoidal functor
$$\Upsilon_\CY:\QCoh(\CY)\to \IndCoh^!(\CY).$$

\sssec{} \label{sss:! acts on *}

Note also that by the definition of $\IndCoh^*$, for $S\in {}^{<\infty}\!\affSch$ we have a 
naturally defined action of $\IndCoh^!(S)$ on $\IndCoh^*(S)$.

\medskip

For a map $f:S_1\to S_2$, this action satisfies the projection formula
$$\CF_2\sotimes f^\IndCoh_*(\CF_1) \simeq f^\IndCoh_*(f^!(\CF_2)\sotimes \CF_1), \quad
\CF_1\in \IndCoh^*(S_1),\, \CF_2\in \IndCoh^!(S_2).$$

\medskip

This implies that for $\CY\in \IndSch^{\on{ind-aff}}_{\on{laft}}$, we also have a natural action of
$\IndCoh^!(\CY)$ on $\IndCoh^*(\CY)$, and the projection formula holds. 

\sssec{}

For an ind-scheme $\CY$ we have a canonically defined pairing:
\begin{equation} \label{e:IndCoh ! * pairing}
\IndCoh^*(\CY)\otimes \IndCoh^!(\CY) \overset{\on{action}}\to \IndCoh^*(\CY) \overset{\Gamma^\IndCoh(\CY,-)}\longrightarrow \Vect.
\end{equation}

Note, however, that unlike the case of schemes, we do not claim that the above pairing realizes 
$\IndCoh^!(\CY)$ as the predual of $\IndCoh^*(\CY)$. (It will, however, be a perfect duality, under
the placidity assumption.)

\ssec{Further properties of \texorpdfstring{$\IndCoh^*$}{IndCoh*}}

\sssec{The functor $\Psi$} 

Let $S$ be an eventually coconnective scheme. The presentation in \eqref{e:IndCoh* as limit} shows that
we have a canonically defined functor
$$\Psi_S:\IndCoh^*(S)\to \QCoh(S).$$

Indeed, it is given by the compatible family of functors 
$$\Psi_{S_0}:\IndCoh(S_0)\to \QCoh(S_0),$$
where we use the fact that for any affine scheme $S$, written as a limit of other affine schemes
$$S\simeq \underset{\alpha}{\on{lim}}\, S_\alpha,$$
the functor
$$\QCoh(S)\to \underset{\alpha}{\on{lim}}\, \QCoh(S_\alpha),$$
given by taking direct images along $S\to S_\alpha$, is an equivalence. 

\medskip

For a map $f:S_1\to S_2$, we have a commutative diagram
$$
\CD
\IndCoh^*(S_1) @>{\Psi_{S_1}}>> \QCoh(S_1) \\
@V{f^\IndCoh_*}VV @VV{f_*}V \\
\IndCoh^*(S_2) @>{\Psi_{S_2}}>> \QCoh(S_2).
\endCD
$$

\sssec{}

Unwinding, we obtain that with respect to the identification 
$$\IndCoh^*(S)\simeq \IndCoh^!(S)^\vee$$
of \secref{sss:! predual of *} and the canonical self-duality on $\QCoh(S)$, 
we have
$$\Psi_S \simeq (\Upsilon_S)^\vee$$

%
%
%
%
%
%

\sssec{} \label{sss:Psi for IndCoh* indsch}

For an ind-affine ind-scheme $\CY$, we have the functor
$$\Psi_\CY:\IndCoh^*(\CY)\to \QCoh_{\on{co}}(\CY)$$ defined in terms of the presentation \eqref{e:IndCoh* indsch} by
$$\IndCoh^*(\CY)\simeq \underset{S\to \CY}{\on{colim}}\, \IndCoh^*(S) \overset{\{\Psi_S\}}\longrightarrow
\underset{S\to \CY}{\on{colim}}\, \QCoh(S) \to \QCoh_{\on{co}}(\CY),$$
where the colimits are taken over the index category of eventually coconnective affine schemes equipped
with a closed embedding into $\CY$. 

\medskip

Note that the functors $\Psi_\CY$ and $\Upsilon_\CY$ are mutually dual in the sense that the following diagram
commutes
$$
\CD
\IndCoh^*(\CY) \otimes \QCoh(\CY) @>{\on{Id}\otimes \Upsilon_\CY}>> \IndCoh^*(\CY) \otimes \IndCoh^!(\CY)  \\
@V{\Psi_\CY\otimes \on{Id}}VV @VV{\text{\eqref{e:IndCoh ! * pairing}}}V \\
\QCoh_{\on{co}}(\CY) \otimes \QCoh(\CY) @>{\text{\eqref{e:QCoh co pairing}}}>> \Vect.
\endCD
$$

\sssec{} \label{sss:^* IndCoh *}

Let $f:S_1\to S_2$ be a map between eventually connective affine schemes. Assume that $S_1$
is finitely presented over $S_2$ (inside the category of $n$-coconnective schemes for some $n$)
and that $f$ is of finite Tor-dimension. 

\medskip

We claim that in this case the functor 
$$f^\IndCoh_*:\IndCoh^*(S_1)\to \IndCoh^*(S_2)$$
admits a left adjoint, to be denoted $f^{*,\IndCoh}$. 

\medskip

By Noetherian approximation (see \cite[Sect. 4.4.1]{Lu3}), the assumption on $f$ implies that we can write 
$S_2$ as $\underset{\alpha}{\on{lim}}\, S_{2,\alpha}$, where $S_{2,\alpha}\in {}^{\leq n}\!\affSch$, so that there exists a compatible 
family of Cartesian diagrams
$$
\CD
S_1 @>{f}>> S_2 \\
@VVV @VVV \\
S_{1,\alpha} @>{f_\alpha}>> S_{2,\alpha}
\endCD
$$
with $S_{1,\alpha}\in {}^{\leq n}\!\affSch$ and $S_1\simeq \underset{\alpha}{\on{lim}}\, S_{1,\alpha}$. 

\medskip

In this case we have
$$\IndCoh^*(S_2) \simeq \underset{\alpha}{\on{lim}}\, \IndCoh(S_{2,\alpha}) \text{ and } 
\IndCoh^*(S_1) \simeq \underset{\alpha}{\on{lim}}\, \IndCoh(S_{1,\alpha}),$$
and the functor $f^{*,\IndCoh}$ is given by the compatible family of functors
$$\IndCoh(S_{2,\alpha})  \overset{f_\alpha^{*,\IndCoh}}\longrightarrow \IndCoh(S_{1,\alpha}),$$
which exists thanks to \cite[Lemma 3.5.8]{Ga7}.

\sssec{}

Let 
$$
\CD
S'_1 @>{f'}>> S'_2 \\
@V{g_1}VV @VV{g_2}V \\
S_1 @>{f}>> S_2 
\endCD
$$
be a Cartesian diagram of eventually connective affine schemes, where the horizontal arrows
are of finite presentation and finite Tor-dimension. 

\medskip

We have the tautological isomorphism
$$(g_2)^\IndCoh_*\circ f'{}^{\IndCoh}_* \simeq f^{\IndCoh}_*\circ (g_1)^\IndCoh_*,$$
from which we obtain a natural transformation
$$f^{*,\IndCoh}\circ (g_2)^\IndCoh_*\to (g_1)^\IndCoh_*\circ f'{}^{*,\IndCoh}.$$

However, it is easy to see that the above natural transformation
is an isomorphism, see \cite[Lemma 3.6.9]{Ga7}.

\sssec{}

Let $f:\CY_1\to \CY_2$ be a a map between ind-affine ind-schemes, and assume that $f$ is
(i) affine, (ii) of finite presentation, (iii) of finite Tor-dimension (i.e., the above properties hold
after base change of $f$ by an affine scheme). 

\medskip

We claim that in this case the functor
$$f^\IndCoh_*:\IndCoh^*(\CY_1)\to \IndCoh^*(\CY_2)$$
admits a left adjoint (to be denoted $f^{*,\IndCoh}$). 

\medskip

Indeed, write
$$\IndCoh^*(\CY_2)\simeq \underset{S_{2,\alpha}\to \CY_2}{\on{colim}}\, \IndCoh^*(S_{2,\alpha}),$$
where the index category is that of $S_{2,\alpha}\in {}^{<\infty}\!\affSch$ equipped with a closed embedding
$S_{2,\alpha}\to \CY_2$.

\medskip

For $S_{2,\alpha}$ as above, set
$$S_{1,\alpha}:=\CY_1\underset{\CY_2}\times S_{2,\alpha}.$$

Then the family
$$S_{1,\alpha}\to \CY_1$$
is cofinal in the category of eventually coconnective affine schemes mapping to $\CY_1$, and hence we have
$$\IndCoh^*(\CY_1)\simeq \underset{S_{1,\alpha}\to \CY_1}{\on{colim}}\, \IndCoh^*(S_1).$$

\medskip

In terms of this presentation, the functor $f^{*,\IndCoh}$ is given by the (compatible) family of functors 
$$f_\alpha^{*,\IndCoh}:\IndCoh^*(S_{2,\alpha})\to \IndCoh^*(S_{1,\alpha}),$$
see \secref{sss:^* IndCoh *}. 

\sssec{} \label{sss:IndCoh^* tensor up}

The following is a counterpart of \secref{sss:IndCoh^! tensor up} for $\IndCoh^*$. Let
us be in the situation of {\it loc. cit.}, but let us assume that $\CY$ is an ind-scheme. 

\medskip

Since $f$ is finite Tor-dimension, we can consider the functor 
$f^{*,\IndCoh}:\IndCoh^*(\CY)\to \IndCoh^*(\CY')$, and it extends to a functor
\begin{equation} \label{e:IndCoh^* tensor up}
\QCoh(S')\underset{\QCoh(S)}\otimes \IndCoh^*(\CY)\to \IndCoh^*(\CY').
\end{equation}

We have:

\begin{lem}  \label{l:IndCoh^* tensor up}
The functor \eqref{e:IndCoh^* tensor up} is fully faithful. If $f$ is smooth, it is an equivalence.
\end{lem}

\ssec{The t-structure on \texorpdfstring{$\IndCoh^*$}{IndCoh*}} \label{ss:t-str IndCoh *}

\sssec{}

Let $S$ be an eventually coconnective affine scheme. The presentation in \eqref{e:IndCoh* as limit} endows the
category $\IndCoh^*(S)$ with a t-structure. It is uniquely characterized by the property that for a map
$$f:S\to S_0,\quad S_0\in {}^{<\infty}\!\affSch_{\on{ft}}$$
the functor
$$f_*:\IndCoh^*(S)\to \IndCoh(S_0)$$
is t-exact.

\medskip

This t-structure is compatible with filtered colimits, by construction. 

\sssec{}

For a map $S_1\to S_2$ between eventually coconnective affine schemes, the corresponding functor
$$f^\IndCoh_*:\IndCoh^*(S_1)\to \IndCoh^*(S_2)$$
is t-exact.

\sssec{}

By construction, the functor $\Psi_S$ is t-exact and induces an equivalence 
$$\IndCoh^*(S)^{>-\infty}\overset{\sim}\to \QCoh(S)^{>-\infty}.$$

\sssec{}

Let $\CY$ be an ind-affine ind-scheme. We use the presentation \eqref{e:IndCoh* indsch} and the
construction in \secref{sss:t-structure colimit} to equip 
$\IndCoh^*(\CY)$ with a t-structure:

\medskip

By definition, $\IndCoh^*(\CY)^{\leq 0}$ is generated under colimits by the essential images of
$\IndCoh^*(S)^{\leq 0}$ for $S\in {}^{<\infty}\!\affSch$ equipped with a closed embedding $S\to \CY$. 

\sssec{}

From \lemref{l:colimit t-structures}, we obtain: 

\begin{cor} \label{c:t-structure on IndCoh*} \hfill

\smallskip

\noindent{\em(a)} The t-structure on $\IndCoh^*(\CY)$ is compatible with filtered colimits.

\smallskip

\noindent{\em(b)} For every $S\in {}^{<\infty}\!\affSch$ equipped with a closed embedding $S\to \CY$, 
the direct image functor $\IndCoh^*(S)\to \IndCoh^*(\CY)$ is t-exact.

\end{cor}

\sssec{}

Let $f:\CY_1\to \CY_2$ be a map between ind-affine ind-schemes. As in \corref{c:dir im QCoh co t exact}, we have:

\begin{cor} \label{c:dir im IndCoh* t exact}
The functor
$$f^\IndCoh_*:\IndCoh^*(\CY_1)\to \IndCoh^*(\CY_2)$$
is t-exact.
\end{cor}

\sssec{} 

Recall the functor
$$\Psi_\CY:\IndCoh^*(\CY)\to \QCoh_{\on{co}}(\CY),$$
see \secref{sss:Psi for IndCoh* indsch}. We claim:

\begin{lem} \label{l:Psi IndCoh* IndSch}
The functor $\Psi_\CY$ is t-exact and induces an equivalence:
$$\IndCoh^*(\CY)^{>-\infty}\to \QCoh_{\on{co}}(\CY)^{>-\infty}.$$
\end{lem}

\begin{proof}

The fact that $\Psi$ is t-exact follows from from Corollaries 
\ref{c:t-structure on QCoh co}(b) and \ref{c:t-exact out of colimit}.

\medskip

To prove the equivalence statement, it suffices to show that for any $n$, the corresponding functor
$$\IndCoh^*(\CY)^{\geq 0,\leq n}\to \QCoh_{\on{co}}(\CY)^{\geq 0,\leq n}$$
is an equivalence.

\medskip

As in the proof of \lemref{l:colimit t-structures}, we can write 
$$\QCoh_{\on{co}}(\CY)^{\geq 0,\leq n}\simeq \underset{S\to \CY}{\on{colim}}\, \QCoh(S)^{\geq 0,\leq n}$$
and 
$$\IndCoh^*(\CY)^{\geq 0,\leq n}\simeq \underset{S'\to \CY}{\on{colim}}\, \IndCoh^*(S')^{\geq 0,\leq n},$$
where:

\begin{itemize}

\item The index $S$ runs over the category affine schemes equipped with 
a closed embedding into $\CY$; 

\item The index $S'$ runs over the category of eventually coconnective affine schemes equipped with 
a closed embedding into $\CY$.

\item Both colimits are taken in the $\infty$-category of categories closed under filtered colimits and
functors that preserve filtered colimits. 

\end{itemize}

Now, the assertion follows from the fact that that for $m\geq n$, the direct image functor
$$\QCoh({}^{\leq m}\!S)\to \QCoh(S)$$
induces an equivalence
$$\QCoh({}^{\leq m}\!S)^{\geq 0,\leq n}\to \QCoh(S)^{\geq 0,\leq n}.$$

\end{proof} 

\ssec{Placidity} \label{ss:placid}

\sssec{} \label{sss:placid}

An affine scheme $S$ is said to be \emph{placid} if it can be written as a limit
\begin{equation} \label{e:present placid}
S\simeq \underset{\alpha\in A}{\on{lim}}\, S_\alpha,
\end{equation} 
where:

\begin{itemize}

\item $S_\alpha\in \affSch_{\on{aft}}$;

\item The transition maps $f_{\beta,\alpha}:S_\beta\to S_\alpha$ are flat.

\item The category $A$ of indices is co-filtered (i.e., the opposite category is filtered). 

\end{itemize}

\begin{rem}

There are in fact two variants of the definition of placidity. The less restricted one is what we just
gave above. In the more restrictive one, one requires that the maps $f_{\beta,\alpha}$ be smooth.

\medskip

The flatness condition is sufficient for our purposes, which are to ensure the compact generation
of the categories $\IndCoh^*$ and $\IndCoh^!$ (see \secref{sss:IndCoh ! comp gen}). One needs smoothness when one
works with D-modules. 

\medskip

That said, in most examples of placid schemes that we will encounter, the smoothness condition 
is satisfied as well. 

\end{rem} 

\sssec{}

Note that if $S$ is placid, then so is any of its truncations: for a presentation \eqref{e:present placid}, we have
$$^{\leq n}\!S\simeq \underset{\alpha}{\on{lim}}\, {}^{\leq n}\!S_\alpha,$$
and the truncated maps $^{\leq n}\!S_\beta\to {}^{\leq n}\!S_\alpha$ are also flat; in fact the flatness of $f_{\beta,\alpha}$ implies that 
$$^{\leq n}\!S_\beta \simeq S_\beta\underset{S_\alpha}\times {}^{\leq n}\!S_\alpha.$$

\sssec{} \label{sss:afp}

Fix an integer $n$. Let $R\to R'$ be a map in $\on{ComAlg}(\Vect^{\leq 0,\geq -n})$. Recall 
that $R'$ is said to be \emph{finitely presented} as an $R$-algebra if $R'$ is compact as
an object of $\on{ComAlg}((R\mod)^{\leq 0,\geq -n})$.

\medskip

Let $R\to R'$ be a map in $\on{ComAlg}(\Vect^{\leq 0})$. We shall say that $R'$ 
\emph{almost finitely presented} as an $R$-algebra if for every $n$, the map
$\tau^{\geq -n}(R)\to \tau^{\geq -n}(R')$ realizes $\tau^{\geq -n}(R')$ as a finitely
presented $\tau^{\geq -n}(R)$-algebra.

\medskip

We shall say that a morphism of $n$-coconnective affine schemes (resp., affine schemes) 
$\Spec(R')=S'\to S=\Spec(R)$ is of finite presentation (resp., almost of finite presentation) 
if $R'$ is finitely presented (resp., almost finitely presented) as an $R$-algbera. 

\sssec{} \label{sss:lafp}

Let $\CY$ be an ind-affine ind-scheme mapping to an affine scheme $S$. We shall say that $\CY$ is
locally almost of finite presentation over $S$ if for every $n$, the truncation
$$^{\leq n}\CY\in  {}^{\leq n}\!\on{PreStk}$$
can be exhibited as a filtered colimit
$$^{\leq n}\CY\simeq \underset{i\in I}{\on{colim}}\, S_i, \quad S_i\in {}^{\leq n}\!\affSch,$$
such that the maps $S_i\to S$ are of finite presentation.

\medskip 

Let $\CY_1\to \CY_2$ be a map between ind-affine ind-schemes. We shall say that $f$ is locally almost of 
finite presentation if the base change of $f$ by any affine scheme $S$
yields an ind-affine ind-scheme locally almost of finite presentation over $S$. 

\sssec{}

We have the following hereditary property of placidity:

\begin{lem} \label{l:placid inherited}
Let $S'\to S$ be a map almost of finite presentation between affine schemes. Suppose that every coconnective
truncation of $S$ is placid. Then the same is true for $S'$.
\end{lem}

\begin{proof}

Fix $n$ and consider the corresponding map $^{\leq n}\!S'\to {}^{\leq n}\!S$. Write
$^{\leq n}\!S$ as 
$$\underset{\alpha\in A}{\on{lim}}\, S_\alpha, \quad S_\alpha\in {}^{\leq n}\!\affSch_{\on{ft}}$$
as in \eqref{e:present placid}.

\medskip

Then for some index $\alpha$, we have an affine scheme $S'_\alpha\in {}^{\leq n}\!\affSch_{\on{ft}}$
and a Cartesian square 
$$
\CD
^{\leq n}\!S' @>>>  {}^{\leq n}\!S \\
@VVV @VVV \\
S'_\alpha @>>> S_\alpha.
\endCD
$$

Consider the category $A_{/\alpha}$. Since $A$ is cofiltered, the category $A_{/\alpha}$ is also cofiltered and 
the opposite of the inclusion functor $A_{/\alpha}\to A$ is cofinal. 

\medskip

For any $\beta\in  A_{/\alpha}$ denote
$$S'_\beta:=S'_\alpha\underset{S_\alpha}\times S_\beta.$$ Then we have 
$$^{\leq n}\!S'  \simeq \underset{\beta\in A_{/\alpha}}{\on{lim}}\, S'_\beta,$$
and the maps $S'_{\beta_1}\to S'_{\beta_2}$ are flat.

\end{proof} 

\sssec{} \label{sss:ind-placid} 

Let $\CY$ be an ind-affine ind-scheme. We shall say that $\CY$ is \emph{ind-placid} if for every $n$,
the truncation 
$$^{\leq n}\CY\in {}^{\leq n}\!\on{PreStk}$$
can be exhibited as a filtered colimit
\begin{equation} \label{e:present placid ind-sch}
^{\leq n}\CY\simeq \underset{i\in I}{\on{colim}}\, S_i, \quad S_i\in {}^{\leq n}\!\affSch,
\end{equation} 
where:

\begin{itemize}

\item The affine schemes $S_i$ are placid;

\item The transition maps $S_i\to S_j$ are of finite presentation.

\end{itemize}

\sssec{}

From \lemref{l:placid inherited} we obtain:

\begin{cor} \label{c:placid inherited}
Let $\CY$ be an ind-placid ind-scheme and let $f:\CY'\to \CY$ be a map locally 
almost of finite presentation. Then $\CY'$ is also an ind-placid ind-scheme.
\end{cor}

\begin{proof}

Fix an integer $n$ and a presentation of $^{\leq n}\CY$ as in \eqref{e:present placid ind-sch}. 
For every index $i$, set
$$\CY'_i:={}^{\leq n}(\CY'\underset{\CY}\times S_i).$$

By assumption, this is an ind-affine ind-scheme of ind-finite presentation over $S_i$. 
Consider the category 
$$\sF_i:=\{S'_i\in {}^{\leq n}\!\affSch, \,\, S'_i\overset{\on{closed\, emb}}\hookrightarrow \CY'_i\}.$$

By assumption, this subcategory contains a full cofinal subcategory, denoted $\sF'_i$
consisting of those objects 
for which the map $S'_i\to S_i$ is of finite presentation. 

\medskip

The assignment $i\mapsto \sF_i$ extends to a co-Cartesian fibration
$$\sF\to I.$$

The assumption that the maps $S_i\to S_j$ are of finite presentation implies that 
the assignment $i\mapsto \sF'_i$ corresponds to a full cofinal subcategory $\sF'\subset \sF$,
and 
$$\sF'\to I$$
is also a co-Cartesian fibration. By cofinality, we obtain that the map
\begin{equation} \label{e:present placid ind-sch bis}
\underset{(i,S'_i)\in \sF}{\on{colim}}\, S'_i\to \CY'
\end{equation}
is an isomorphism in $^{\leq n}\!\on{PreStk}$. Moreover, for a map $(i,S'_i)\to (j,S'_j)$ in $\sF'$,
the corresponding map $S'_i\to S'_j$ is of finite presentation (because its composition with $S'_j\to S_j$ is). 

\medskip

Since $I$ is filtered and each $\sF'_i$ is filtered, we obtain that $\sF'$ is filtered. 

\medskip

Finally, by \lemref{l:placid inherited}, the affine schemes $S'_i$ are placid. Hence, \eqref{e:present placid ind-sch bis}
gives the desired presentation of $\CY'$.

\end{proof} 

\ssec{The categories \texorpdfstring{$\IndCoh^!(-)$}{IndCoh!} and \texorpdfstring{$\IndCoh^*(-)$}{IndCoh*} in the (ind)-placid case}

\sssec{} \label{sss:IndCoh ! comp gen}

Let $S\in {}^{\leq n}\!\affSch$ be placid. We claim that in this case, the category $\IndCoh^!(S)$ is compactly 
generated.

\medskip

Indeed, write $S$ as in \eqref{e:present placid}.  
Since the maps $f_{\beta,\alpha}$ are flat, the functors $f^!_{\beta,\alpha}$ preserves coherence,
and hence compactness. We obtain 
$$\IndCoh^!(S)\simeq \underset{\alpha}{\on{colim}}\, \IndCoh(S_\alpha),$$
where the terms are compactly generated, and the transition functors preserve compactness.

\medskip

Hence, the images of $\Coh(S_\alpha)\subset \IndCoh(S_\alpha)$ under the !-pullback functors
$$\IndCoh(S_\alpha)\to \IndCoh^!(S)$$
provide a set of compact generators of $\IndCoh^!(S)$.

\sssec{}  \label{sss:! * duality placid}

From the identification 
\begin{equation} \label{e:* dual of !}
\IndCoh^*(S) \simeq \IndCoh^!(S)^\vee
\end{equation} 
of \secref{sss:! predual of *}, we obtain that $\IndCoh^*(S)$ is also compactly generated, and 
\eqref{e:* dual of !} is a duality in $\DGCat$.

\sssec{}

Explicitly, the presentation 
$$\IndCoh^*(S) \simeq \underset{\alpha}{\on{lim}}\, \IndCoh(S_\alpha)$$
(with respect to the *-pushforward functors) implies that 
$$\IndCoh^*(S) \simeq \underset{\alpha}{\on{colim}}\, \IndCoh(S_\alpha),$$
with respect to the *-pullback functors (which are well-defined, due to the flatness
assumption).

\medskip

Thus, the compact generators of $\IndCoh^*(S)$ are the images of
$\Coh(S_\alpha)\subset \IndCoh(S_\alpha)$ along the *-pullback functors
$$\IndCoh(S_\alpha)\to \IndCoh^*(S).$$

\sssec{} \label{sss:_* preserve compactness}

Let $f:S_1\to S_2$ be a morphism almost of finite presentation between placid 
eventually coconnective affine schemes. Assume now that $f$ is a closed embedding.
We claim that in this case the functor
$$f^\IndCoh_*:\IndCoh^*(S_1)\to \IndCoh^*(S_2)$$
admits a continuous right adjoint (to be denoted $f^!$).

\medskip

Indeed, we claim that $f^\IndCoh_*$ preserves compactness. This follows
from the manipulation in the proof of \lemref{l:placid inherited} using the fact that
for a Cartesian diagram of affine schemes almost of finite type 
$$
\CD
S'_1 @>{f'}>> S'_2 \\
@V{g_1}VV @VV{g_2}V \\
S''_1 @>{f''}>> S''_2
\endCD
$$
with the maps $g_1,g_2$ flat, the natural transformation
$$g_2^{\IndCoh,*}\circ f''{}^\IndCoh_*\to f'{}^\IndCoh_*\circ g_1^{\IndCoh,*}$$
is an isomorphism.

\sssec{} \label{sss:BeckChev placid}

Moreover, in the above situation, for a Cartesian diagram 
$$
\CD
\wt{S}_1 @>{\wt{f}}>> \wt{S}_2 \\
@V{g_1}VV @VV{g_2}V \\
S_1 @>{f}>> S_2,
\endCD
$$
where:

\begin{itemize}

\item All affine schemes involved are eventually coconnective and placid;

\item The maps $f$ and $\wt{f}$ are closed embeddings of finite presentation,

\end{itemize}
the natural transformation
$$(g_1)^\IndCoh_*\circ \wt{f}^! \to f^!\circ (g_2)^\IndCoh_*, \quad \IndCoh^*(\wt{S}_2)\rightrightarrows \IndCoh^*(S_1),$$
obtained by adjunction from
$$f^\IndCoh_*\circ (g_1)^\IndCoh_*\to (g_2)^\IndCoh_*\circ \wt{f}^\IndCoh_*,$$
is an isomorphism.

\sssec{}

Recall (see \secref{sss:dual of ^!}) that that with respect to the dualities
$$\IndCoh^*(S_i) \simeq \IndCoh^!(S_i)^\vee, \quad i=1,2,$$
we have
$$f^\IndCoh_*\simeq (f^!)^\vee.$$

\medskip

Hence, the existence of a continuous right adjoint of 
$$f^\IndCoh_*:\IndCoh^*(S_1)\to \IndCoh^*(S_2)$$
implies the existence of a \emph{left} adjoint of the functor
$$f^!:\IndCoh^!(S_2)\to \IndCoh^!(S_1).$$

We will denote this left adjoint by
$$f^\IndCoh_*:\IndCoh^!(S_1)\to \IndCoh^!(S_2).$$

\sssec{} \label{sss:IndCoh* on ind-placid}

Let now $\CY$ be an ind-placid ind-scheme. We claim that in this case the 
category $\IndCoh^*(\CY)$ is compactly generated. 

\medskip

Indeed, writing 
$$\IndCoh^*(\CY) \simeq \underset{i}{\on{colim}}\, \IndCoh^*(S_i),$$
where $S_i$ are eventually coconnective placid affine schemes and the
transition maps $S_i\overset{f_{i,j}}\to S_j$ are closed embeddings almost of finite presentation,
we obtain that the transition functors 
$$\IndCoh^*(S_i)\to \IndCoh^*(S_j)$$
preserve compactness (see \secref{sss:_* preserve compactness} above).

\medskip

Similarly, the category $\IndCoh^!(\CY)$ can be written as 
$$\underset{i}{\on{colim}}\, \IndCoh^!(S_i),$$
where the transition functors are $(f_{i,j})^\IndCoh_*$. This implies that 
$\IndCoh^!(\CY)$ is also compactly generated by the essential images of
$\IndCoh^!(S_i)^c$.

\sssec{} \label{sss:IndCoh ! * duality}

Note also that it follows that if $\CY$ is placid, the pairing \eqref{e:IndCoh ! * pairing}
is a perfect duality.

\sssec{}

By a similar logic we obtain:

\begin{lem} \label{l:boxtimes IndCoh placid}
Let $\CY_1$ and $\CY_2$ be a pair of ind-placid ind-schemes. Then the functors
\eqref{e:boxtimes IndCoh !} and \eqref{e:boxtimes IndCoh *} are equivalences.
\end{lem}

\sssec{} \label{sss:closed emb placid}

Let $f:\CY_1\to \CY_2$ be a morphism between ind-placid ind-affine ind-schemes that
is locally almost of finite presentation. Assume that $f$ is an ind-closed embedding 
(i.e., for every closed embedding $S\to \CY_1$, the composite map $S\to \CY_1\to \CY_2$
is also a closed embedding).

\medskip

Similarly to the above, we obtain that in this case, the functor
$$f^\IndCoh_*:\IndCoh^*(\CY_1)\to  \IndCoh^*(\CY_2)$$
preserves compactness, and hence admits a continuous right adjoint, to be denoted $f^!$.

\medskip

By duality, the functor 
$$f^!:\IndCoh^!(\CY_2)\to \IndCoh^!(\CY_1)$$
admits a left adjoint, to be denoted
$$f^\IndCoh_*:\IndCoh^!(\CY_1)\to \IndCoh^!(\CY_2).$$

\sssec{} \label{sss:BeckChev placid ind}

Let
$$
\CD
\wt\CY_1 @>{\wt{f}}>> \wt\CY_2 \\
@V{g_1}VV @VV{g_2}V \\
\CY_1 @>{f}>> \CY_2
\endCD
$$
be a Cartesian diagram, where: 

\begin{itemize}

\item All objects are involved are ind-placid affine ind-schemes;

\item The maps $f$ and $\wt{f}$ are ind-closed embeddings of finite presentation.

\end{itemize}

\medskip

Unwinding, it follows from \secref{sss:BeckChev placid} that in this case,  the natural transformation
$$(g_1^\IndCoh)_*\circ \wt{f}^!\to  f^!\circ (g_2^\IndCoh)_*, \quad \IndCoh^*(\wt\CY_2)\rightrightarrows \IndCoh^*(\CY_1),
\quad \IndCoh^*(\wt\CY_2)\rightrightarrows \IndCoh^*(\CY_1),$$
obtained by adjunction from
$$f^\IndCoh_*\circ (g_1^\IndCoh)_*\to (g_2^\IndCoh)_*\circ \wt{f}^\IndCoh_*,$$
is an isomorphism.

\medskip

By duality, the natural transformation
$$\wt{f}^\IndCoh_*\circ g_1^! \to g_2^!\circ f^\IndCoh_*, \quad \IndCoh^*(\CY_1)\rightrightarrows \IndCoh^*(\wt\CY_2),$$
obtained by adjunction from
$$g_1^!\circ f^!\simeq \wt{f}^!\circ g_2^!,$$
is also an isomorphism.

\sssec{} \label{sss:* pullback on IndCoh* finite Tor dim}

Let now $f:\CY_1\to \CY_2$ be a morphism between ind-placid ind-schemes. Assume that $f$ is affine and of finite Tor-dimension
(but we are not assuming that $f$ be of finite presentation). 

\medskip

We claim, generalizing \secref{sss:^* IndCoh *}, 
that in this case the functor $$f^{*,\IndCoh}:\IndCoh^*(\CY_2)\to \IndCoh^*(\CY_1),$$ left adjoint to $f^\IndCoh_*$,
exists, and satisfies base change against functors $g^\IndCoh_*$ for $g:\CY'_1\to \CY_1$, where $g$ is another
ind-placid ind-scheme. 

\medskip

Indeed, by \secref{sss:closed emb placid}, the question reduces to the case when 
$\CY_2$ is an eventually coconnective affine scheme, to be denoted $S_2$. In this case, $\CY_1$ is also 
an eventually coconnective affine scheme, to be denoted $S_1$.

\medskip 

Furthermore, in this case we can further reduce to the case when $S_2$ is of finite type. Since $S_1$ is placid, we
can factor the morphism $f$ as
$$S_1\overset{h}\to S_{1,0} \overset{f_0}\to S_2,$$
where:

\begin{itemize} 

\item $S_{1,0}$ is an eventually coconnective scheme of finite type,

\item The morphism $h$ is flat;

\item The morphism $f_0$ is of finite Tor-dimension. 

\end{itemize}

This reduces the assertion to the case of eventually coconnective affine schemes of finite type,
where it follows from \cite[Lemma 3.5.8]{Ga7}.

\section{Factorization patterns} \label{s:fact}

The local Langlands theory considers various representation-theoretic categories $\bA_x$ attached to 
the group $G$ (or its dual $\cG$) and the formal disc $\cD_x$, attached to a point $x\in X$.
More generally, one is led to consider the multi-disc $\cD_{\ul{x}}$, $\ul{x}=x_1,...,x_n$; moreover
the points $x_1,...,x_n$ are allowed to move in families over $X$ and collide. In this case, we shall say
that $\ul{x}$ is a (scheme-theoretic) point of the \emph{Ran} space of $X$. 

\medskip

The datum of a 
\emph{factorization category} attached to such $\ul{x}$ a category $\bA_{\ul{x}}$, such that if 
$\ul{x}_1$ and $\ul{x}_2$ are disjoint, we are given an isomorphism
$$\bA_{\ul{x}_1\sqcup \ul{x}_2}\simeq \bA_{\ul{x}_1}\otimes \bA_{\ul{x}_2}.$$

\medskip

We develop the theory of factorization categories in this section, along with various adjoining
notions (factorization spaces, factorization algebras). Can can view this section as a natural development of 
the theory of chiral algebras, initiated in \cite{BD2}. The main difference with {\it loc. cit.} is that all
our constructions take place in the world of $\infty$-categories, whereas in \cite{BD2} one mainly worked
at the abelian level.

\medskip

In order to produce examples of factorization categories one often uses geometric objects associated
to the formal (resp., formal punctured) disc, such as arcs and loop spaces. Some of the work in this
section is devoted to the study of the relevant geometries.

\ssec{Factorization spaces} 

\sssec{}

The Ran space of $X$, denoted $\Ran$, is the prestack that assigns to an affine test scheme $S$ the set
of finite-non-empty subsets of $\Hom(S_{\on{red}},X)$. 

\medskip

Note that, by definition, the map $\Ran\to \Ran_{\dr}$ is an isomorphism.

\medskip

We denote $k$-points of $\Ran$ by $\ul{x}$. By definition, these are finite non-empty collections 
\begin{equation} \label{e:ul x}
\ul{x}=\{x_1,...,x_n\}
\end{equation}
of $k$-points of $X$. 

\sssec{}

In what follows we will use the following notations. Let $\ul{x}:S\to \Ran$ be a map corresponding to 
$I\subset \Hom(S_{\on{red}},X)$.

\medskip

\begin{itemize}

\item For $i\in I$, we will denote by $x_i$ the corresponding map $S_{\on{red}}\to X$;

\item We will denote by $\on{Graph}_{x_i}\subset S\times X$ the graph of $x_i$ (viewed as a closed \emph{subset},
i.e., we ignore its scheme-theoretic structure); 

\item We will denote by $\on{Graph}_{\ul{x}}$ the Zariski-closed subset of $S\times X$ equal to 
$\underset{i}\cup\,  \on{Graph}_{x_i}$.

\end{itemize} 

\sssec{} \label{sss:Ran as colim}

One can exhibit $\Ran$ explicitly as a colimit of de Rham spaces of schemes. Namely,
$$\Ran \simeq \underset{I\in (\on{fSet}^{\on{surj}})^{\on{op}}}{\on{colim}}\, X^I_\dr,$$
where $\on{fSet}^{\on{surj}}$ is the category of non-empty finite sets and surjective maps
(see \cite[Sect. 2]{Ro2} for a detailed discussion).

\sssec{}

The presentation \eqref{sss:Ran as colim} implies, in particular, that $X_\dr$ is locally almost of finite type as a prestack.
Hence, it is sufficient to probe it by eventually connective affine schemes of finite type. 

\medskip

Hence, in the discussion below, we will be tacitly assuming that schemes and prestacks mapping to $\Ran$
are laft (locally almost of finite type).  

\sssec{} \label{sss:disj}

Here are the two basic features of $\Ran$ that will be used in the sequel:

\medskip

\noindent(i) There is a canonically defined map 
$$\on{union}:\Ran\times \Ran\to \Ran,$$
given by the operation of \emph{union} of finite subsets.

\medskip

\noindent(ii) There exists an open subspace
$$(\Ran\times \Ran)_{\on{disj}} \subset \Ran\times \Ran,$$
corresponding to the condition that the two subsets are disjoint. Namely, for a affine test scheme $S$,
a pair of $S$-points $\ul{x}_1,\ul{x}_2$ of $\Ran$
maps to $(\Ran\times \Ran)_{\on{disj}}$ if 
$$\on{Graph}_{\ul{x}_1}\cap \on{Graph}_{\ul{x}_2}=\emptyset.$$

\medskip

Note that the restriction of the map $\on{union}$ to $(\Ran\times \Ran)_{\on{disj}}$ is \'etale. 
 
\sssec{} \label{sss:fact spaces}

By a factorization space $\CT$ over $X$ we will mean a prestack 
$$\CT_\Ran\to \Ran,$$ equipped with a \emph{factorization structure}, which is by definition the datum of an isomorphism
\begin{equation} \label{e:fact space}
\CT_\Ran\underset{\Ran,\on{union}}\times (\Ran\times \Ran)_{\on{disj}} \simeq 
(\CT_\Ran\times \CT_\Ran)\underset{\Ran\times \Ran}\times (\Ran\times \Ran)_{\on{disj}},
\end{equation} 
equipped with a homotopy-coherent data of associativity and commutativity (see \cite[Sect. 6]{Ra6}, where this is spelled out in detail). 

\sssec{}

Given a map $\CZ\to \Ran$, we will denote by $\CT_\CZ$ the base change
$$\CZ\underset{\Ran}\times \CT_\Ran.$$

\medskip

For $\CZ=\on{pt}$ so that $\CZ\to \Ran$ corresponds to $\ul{x}\in \Ran$, we will write $\CT_{\ul{x}}$ for the corresponding $\CZ$.
The factorization structure on $\CT$ implies that for $\ul{x}$ as in \eqref{e:ul x}, we have
$$\CT_{\ul{x}}\simeq \underset{i}\Pi\, \CT_{x_i}.$$

\sssec{} \label{sss:aff grass}

A basic example of a factorization space is the affine Grassmannian $\Gr_G$. Namely, for an affine test scheme $S$
and a map $\ul{x}:S\to \Ran$, its lift to $\Gr_{G,\Ran}$ is a datum of 
$$(\CP_G,\alpha),$$
where $\CP_G$ is a $G$-bundle on $S\times X$, and $\alpha$ is a trivialization of $\CP_G$ over the open
$$S\times X-\on{Graph}_{\ul{x}}.$$

\sssec{}  \label{sss:local prop fact}

Let $\CT$ be a factorization space. We can talk about its local properties, such as being a scheme, being an 
ind-scheme, being (ind)-placid, being formally smooth, etc. 

\medskip

By definition, this means that these properties hold for $\CT_S$ relatively to $S$ for every $S\in \affSch_{/\Ran}$. 

\medskip 

In a similar way, we can talk about local properties of a map between factorization spaces
(e.g., being flat or an fpqc cover). 

\ssec{Factorization \emph{module} spaces}

\sssec{}  \label{sss:Ran subseteq}

Let $\Ran^{\subseteq}$ be the subfunctor of $\Ran\times \Ran$, such that $\Maps(S,\Ran^{\subseteq})$
corresponds to pairs $\ul{x}\subseteq \ul{x}'$, as subsets of $\Hom(S_{\on{red}},X)$.  

\medskip

Denote by $\on{pr}_{\on{small}}$ and $\on{pr}_{\on{big}}$ the two projections $\Ran^{\subseteq}\rightrightarrows \Ran$
that send a pair $(\ul{x},\ul{x}')$ to $\ul{x}$ and $\ul{x}'$, respectively.

\sssec{} \label{sss:Z Ran}

Let $\CZ$ be a prestack equipped with a map to $\Ran$. Denote
$$\CZ^{\subseteq}:=\CZ\underset{\Ran,\on{pr}_{\on{small}}}\times \Ran^{\subseteq}.$$

For $\CZ=\on{pt}$, so that $\CZ\to \Ran$ corresponds to $\ul{x}\in \Ran$, we will write $\Ran_{\ul{x}}$
for the corresponding space $\CZ^{\subseteq}$. 

\sssec{}

Note that we have a variant of the map $\on{union}$:
$$\on{union}:\Ran\times \CZ^{\subseteq}\to \CZ^{\subseteq}, \quad (\ul{x},(z,\ul{x}'))\mapsto (z,\ul{x}\cup \ul{x}').$$

\medskip

Denote by
$$(\Ran\times \CZ^{\subseteq})_{\on{disj}}\subset \Ran \times \CZ^{\subseteq}$$
the open subfunctor equal to the preimage of $(\Ran\times \Ran)_{\on{disj}}$ under
$$\Ran\times \CZ^{\subseteq} \to \Ran\times \Ran^{\subseteq}\overset{\on{id}\times \on{pr}_{\on{big}}}\to
\Ran\times \Ran.$$

\sssec{} \label{sss:factorization module spaces}

Given a factorization space $\CT$, a factorization module space $\CT_m$ \emph{over} $\CT$ \emph{at} $\CZ$
is a prestack 
$$(\CT_m)_{\CZ^{\subseteq}}\to \CZ^{\subseteq},$$
equipped with a datum of factorization \emph{against} $\CT$:
\begin{equation} \label{e:fact mod space}
(\CT_m)_{\CZ^{\subseteq}}\underset{\CZ^{\subseteq},\on{union}}\times (\Ran\times \CZ^{\subseteq})_{\on{disj}}\simeq
(\CT_\Ran\times (\CT_m)_{\CZ^{\subseteq}})\underset{\Ran\times \CZ^{\subseteq}}\times (\Ran\times \CZ^{\subseteq})_{\on{disj}},
\end{equation} 
equipped with a homotopy-coherent data of associativity; see \cite[Sect. 6]{Ra6} for complete details. 

\sssec{}

For a factorization module space $\CT_m$ at $\CZ$, denote
$$(\CT_m)_\CZ:=\CZ\underset{\CZ^{\subseteq}}\times (\CT_m)_{\CZ^{\subseteq}},$$
where $\CZ\to \CZ^{\subseteq}$ is the map $\on{diag}_\CZ$ of \secref{sss:diag Z}.

\medskip

We will refer to $(\CT_m)_\CZ$ as the \emph{prestack underlying} the factorization module space $\CT_m$.

\sssec{} \label{sss:vac fact space}

A basic example of a factorization module space over $\CT$, defined for any $\CZ\to \Ran$, denoted
$\CT^{\on{fact}_\CZ}$, is constructed as follows:

\medskip

$$\CT^{\on{fact}_\CZ}_{\CZ^{\subseteq}}:=\CT_{\CZ^{\subseteq}},$$
where $\CZ^{\subseteq}\to \Ran$ is the map $\on{pr}_{\on{big},\CZ}$, 
with the datum of \eqref{e:fact mod space} being provided by the factorization structure on $\CT$ itself.

\medskip

We refer to $\CT^{\on{fact}_\CZ}$ as the \emph{vacuum} factorization module space over $\CT$ at $\CZ$. 

\sssec{}  \label{sss:Gr G level}

Here is an example of a factorization module space over $\Gr_G$, also defined for any $\CZ\to \Ran$, to be denoted 
$$\Gr_G^{\on{level}_\CZ}.$$

\medskip 

For $(z,\ul{x}):S\to \CZ^{\subseteq}$, a lift of this point to $\Gr_G^{\on{level}_\CZ}$ is a lift of $\ul{x}$ to a point
of $\Gr_{G,\Ran}$, and the trivialization of the restriction of the resulting $G$-bundle $\CP_G$ over $S\times X$ 
to $\wh\cD_{\ul{x}}$ (see \secref{sss:formal formal discs} below). 

\ssec{Digression: formal discs} \label{ss:formal discs}

\sssec{}  \label{sss:formal formal discs}

Fix a map $\ul{x}:S\to \Ran$. Let $\wh\cD_{\ul{x}}$ be the formal scheme equal to the formal completion of $S\times X$ along the closed subset
$$\on{Graph}_{\ul{x}}\subset S\times X,$$
i.e.,
$$\wh\cD_{\ul{x}}:=(S\times X)\underset{(S\times X)_\dr}\times (\on{Graph}_{\ul{x}})_\dr,$$
where $(\on{Graph}_{\ul{x}})_\dr\to (S\times X)_\dr$ is given by the embedding
$$(\on{Graph}_{\ul{x}})_{\on{red}}\hookrightarrow S\times X.$$

\medskip

Note that $\wh\cD_{\ul{x}}$ is naturally the pullback of a relative formal scheme, denoted  $\wh\cD_{\ul{x},\nabla}$,
over $S\times X_\dr$. Namely, 
$$\wh\cD_{\ul{x},\nabla}:=(S\times X_\dr)\underset{(S\times X)_\dr}\times (\on{Graph}_{\ul{x}})_\dr.$$

\sssec{} \label{sss:disj local}

Let $\on{PreStk}_{\on{disj-loc}}\subset \on{PreStk}$ be the full subcategory of disjoint-union-local prestacks, i.e.,
prestacks $\CZ$, for which for a disjoint union of
affine schemes 
$$S=S_1\sqcup S_2,$$
the map
\begin{equation} \label{e:mapping out of disj union}
\Maps(S,\CZ)\to \Maps(S_1,\CZ)\times \Maps(S_2,\CZ)
\end{equation}
is an isomorphism. 

\medskip

Most prestacks that one encounters in practice satisfy this condition. E.g., note that if $\CZ$ satisfies Zariski descent, then
it is disjoint-union-local. 

\medskip

Note that the Yoneda embedding
$$\affSch\to \on{PreStk}_{\on{disj-loc}}$$
commutes with coproducts (this would not be true for the original Yoneda embedding into $\on{PreStk}$).

\sssec{}

Let $\ul{x}_1$ and $\ul{x}_2$ be a pair of $S$-points of $\Ran$, such that $(\ul{x}_1,\ul{x}_2)$ 
lands in $(\Ran\times \Ran)_{\on{disj}}$. Set $\ul{x}:=\on{union}(\ul{x}_1,\ul{x}_2)$. 
In this case we have
\begin{equation} \label{e:splitting of formal disc} 
\wh\cD_{\ul{x}}\simeq \wh\cD_{\ul{x}_1}\sqcup \wh\cD_{\ul{x}_2} \text{ and }
\wh\cD_{\ul{x},\nabla}\simeq \wh\cD_{\ul{x}_1,\nabla}\sqcup \wh\cD_{\ul{x}_2,\nabla},
\end{equation} 
where $\sqcup$ is the coproduct taken in $\on{PreStk}_{\on{disj-loc}}$.

\sssec{} \label{sss:actual formal disc}

Note that, when viewed as an ind-scheme, $\wh\cD_{\ul{x}}$ is ind-affine, i.e., of the form
$$\underset{\alpha}{``\on{colim}"}\, \Spec(R_\alpha).$$

\medskip

Let $\cD_{\ul{x}}$ denote the affine \emph{scheme} equal to
$$\underset{\alpha}{\on{colim}}\, \Spec(R_\alpha),$$
where the colimit is taken in $\affSch$. I.e., 
$$\cD_{\ul{x}}=\Spec(R), \quad R=\underset{\alpha}{\on{lim}}\, R_\alpha.$$

\medskip

Now, by \cite[Theorem 1.1]{Bh}, the map 
\begin{equation} \label{e:formal disc projects to curve hat}
\wh\cD_{\ul{x}}\to X
\end{equation} 
canonically extends to a map 
\begin{equation} \label{e:formal disc projects to curve no hat}
\cD_{\ul{x}}\to X.
\end{equation} 

\begin{rem}

Here we have used \cite[Theorem 1.1]{Bh} in a very elementary situation, in which 
the required assertion can be handled explicitly:

\medskip

The observation is that maps from an affine scheme $\wt{S}$ to $X$ can be described as 
symmetric monoidal functors 
$$\on{Perf}(X)\to \on{Perf}(\wt{S}),$$
and hence this is true for $\wt{S}$ replaced by any prestack. 

\medskip

Now, the required assertion is that the restriction functor
$$\on{Perf}(\cD_{\ul{x}})\to \on{Perf}(\wh\cD_{\ul{x}})$$
is an equivalence, so the restriction map
$$\Maps(\cD_{\ul{x}},X)\to \Maps(\wh\cD_{\ul{x}},X)$$
is an isomorphism. 

\end{rem} 

\begin{rem}

Note that while the assignment
$$\ul{x}\rightsquigarrow \wh\cD_{\ul{x}}$$ is compatible with Zariski localization along $S$,
the formation of $\cD_{\ul{x}}$ is \emph{not}.

\medskip

I.e., for an open $S'\subset S$ and $\ul{x}':=\ul{x}|_{S'}$, the square
$$
\CD
\cD_{\ul{x}'} @>>> \cD_{\ul{x}} \\
@VVV @VVV \\
S' @>>> S
\endCD
$$
is \emph{not} Cartesian. 

\end{rem} 

\sssec{}

We have an ind-closed embedding 
$$\wh\cD_{\ul{x}}\to \cD_{\ul{x}}.$$

In particular, $\on{Graph}_{\ul{x}}$ is a Zariski-closed subset of $\cD_{\ul{x}}$. Set
$$\cD^\times_{\ul{x}}:=\cD_{\ul{x}}-\on{Graph}_{\ul{x}}.$$ 

\medskip

Composing with \eqref{e:formal disc projects to curve no hat}, we obtain a map
$$\cD^\times_{\ul{x}}\to X.$$

\sssec{}  \label{sss:descented formal disc}

The following material reproduces \cite[Construction A.1.3]{Bogd}. 

\medskip 

We claim that $\cD_{\ul{x}}$ descends to a relative affine scheme $\cD_{\ul{x},\nabla}$
over $X_\dr$. In order to construct this descent, it is enough to construct a version of $\cD_{\ul{x}}$ 
over the \v{C}ech nerve of the infinitesimal groupoid 
$$(X\times X)^\wedge$$ 
of $X$. 

\medskip

I.e., we have to construct a \emph{compatible} family of affine schemes
$$\cD'_{\ul{x}}\to X'$$
over each infinitesimal thickening $X'$ of the main diagonal in $X^{\times n}$ for all $n$.

\medskip

We let
$$\wh\cD'_{\ul{x}}:=X'\underset{X_\dr}\times \wh\cD_{\ul{x},\nabla},$$
viewed as an ind-affine ind-scheme, and we let $\cD'_{\ul{x}}$ be the colimit
of $\wh\cD'_{\ul{x}}$, taken in the category of affine schemes. 

\medskip

The map $\wh\cD'_{\ul{x}}\to X'$ extends to a map $\cD'_{\ul{x}}\to X'$
extends by the same principle as in the case of $\cD_{\ul{x}}$.

\sssec{}

This material repeats \cite[Lemma A.1.5]{Bogd}.

\medskip 

In order to establish the compatibility of the above construction, it suffices to show that
(with respect to any of the projections $X'\to X$), the map
\begin{equation} \label{e:compat descent}
\cD'_{\ul{x}}\to X'\underset{X}\times \cD_{\ul{x}}
\end{equation} 
is an isomorphism.

\medskip

We note that $X'\to X$ has the form $\Spec_X(A_0)$, where 
is given by $A_0\in \on{Perf}(X)$, and hence 
$X'\underset{X}\times \cD_{\ul{x}}$ has the form $\Spec_{\cD_{\ul{x}}}(A)$,
where $A\in \on{Perf}(\cD_{\ul{x}})$. 

\medskip

In other words, if $R$ is as in \secref{sss:actual formal disc}, then $X'\underset{X}\times \cD_{\ul{x}}\simeq \Spec(R')$,
where $R'$ is compact as an object of $R\mod$. 

\medskip

Let $R_\alpha$ be as in \secref{sss:actual formal disc}. Then
$$\wh\cD'_{\ul{x}} = \underset{\alpha}{``\on{colim}"}\, \Spec(R'\underset{R}\otimes R_\alpha).$$

The assertion that \eqref{e:compat descent} is an isomorphism is equivalent to saying that the map
$$R'\simeq R'\underset{R}\otimes (\underset{\alpha}{\on{lim}}\, R_\alpha)
\to  \underset{\alpha}{\on{lim}}\, (R'\underset{R}\otimes R_\alpha)$$
is an isomorphism. However, this follows from the fact that $R'$ is compact as an $R$-module.

\sssec{}  \label{sss:descented formal disc punct}

We let $\cD^\times_{\ul{x},\nabla}$ be the open sub-functor of $\cD_{\ul{x},\nabla}$,
equal to
$$\cD_{\ul{x},\nabla}-\on{Graph}_{\ul{x}}.$$

\medskip

It is easy to see that the prestacks
\begin{equation} \label{e:descented formal disc}
\cD_{\ul{x},\nabla} \text{ and } \cD^\times_{\ul{x},\nabla},
\end{equation}
when viewed as relative affine schemes over $X_\dr$, are independent of the choice of $X^o$. 

\sssec{} \label{sss:splitting punctured}

The prestacks \eqref{e:descented formal disc} satisfy a splitting property parallel to
\eqref{e:splitting of formal disc}. 

\ssec{Formation of (horizontal) arc and loop spaces} \label{ss:forming loops gen}

In this subsection we will discuss two ubiquitous sources of examples of factorization spaces. 

\sssec{} 

Let $\CY\to X$ be a D-prestack\footnote{A.k.a., crystal of prestacks.}
over $X$, i.e., $\CY$ is the pullback along $X\to X_\dr$ of a prestack 
$\CY_\nabla\to X_\dr$. 

\medskip

We will assume that $\CY_\nabla$ is disjoint-union-local (see \secref{sss:disj local}). 

\sssec{} \label{sss:forming arcs gen} 

We define the factorization prestack $\fL^+_\nabla(\CY)$ as follows. For an affine test scheme $S$ and a map
$\ul{x}:S\to \Ran$, a lift of this map to $\fL^+_\nabla(\CY)$ is an $X_\dr$-map
$$\wh\cD_{\ul{x},\nabla}\to \CY_\nabla.$$

\medskip

The factorization structure on $\fL^+_\nabla(\CY)$ follows from \eqref{e:splitting of formal disc} and \eqref{e:mapping out of disj union}. 

\begin{rem}

Note that the restriction $\fL^+_\nabla(\CY)|_X$ recovers the original $\CY$. In \secref{sss:another take on arcs}
we will upgrade the assignment
$$\CY\mapsto \fL^+_\nabla(\CY)$$
to an equivalence of categories. 

\end{rem} 

\sssec{}   \label{sss:aff D-sch}

For the rest of this subsection we will assume that $\CY$ is affine over $X$. I.e., $\CY$ is an affine D-scheme. 

\medskip

There is an (obvious) equivalence between the category of D-schemes 
and that of commutative algebras $A$ in $\Dmod(X)$ so that $\oblv^l(A)$ is connective. 

\medskip

The plain affine scheme underlying a given affine D-scheme is
$$\Spec_X(\oblv^l(A)).$$

We will write 
$$\Spec_X(A)$$
when we want to emphasize the D-structure. 

\sssec{} \label{sss:arcs via affine formal disc}

Note that the assumption that $\CY$ is affine over $X$ allows us to interpret the datum of a map
$$\wh\cD_{\ul{x},\nabla}\to \CY_\nabla$$
as 
$$\cD_{\ul{x},\nabla}\to \CY_\nabla,$$
where $\cD_{\ul{x},\nabla}$ is an in \secref{sss:descented formal disc}. 

\sssec{} \label{sss:forming loops aff}

We will now define another factorization  space, denoted $\fL_\nabla(\CY)$. 
A lift of this map to $\fL_\nabla(\CY)_\Ran$
is by definition a $X_\dR$ map
$$\cD^\times_{\ul{x},\nabla}\to \CY_\nabla.$$

\medskip

%

%


The factorization structure on $\fL_\nabla(\CY)_\Ran$ follows from the splitting property in 
\secref{sss:splitting punctured}. 

\sssec{}

Using \secref{sss:arcs via affine formal disc} and the open embeddings
$$\cD^\times_{\ul{x},\nabla}\hookrightarrow \cD_{\ul{x},\nabla}$$
we obtain a map of factorization spaces:
$$\iota:\fL^+_\nabla(\CY)\to \fL_\nabla(\CY).$$

We claim:

\begin{lem} \label{l:loops into affine} \hfill 

\smallskip

\noindent{\em(a)} $\fL^+_\nabla(\CY)$ is an factorization affine scheme.

\smallskip

\noindent{\em(b)} $\fL_\nabla(\CY)$ is a factorization ind-affine ind-scheme.

\smallskip

\noindent{\em(c)} The map $\iota$ is a closed embedding.\footnote{By a slight abuse of terminology, we call a map $Z\to \CZ$
from a scheme $Z$ to an ind-scheme $\CZ$ a ``closed embedding" if the map from $Z$ to some/any of the schemes that comprise
$\CZ$ is a closed embedding. Note that such a map is \emph{not} a closed embedding in the DAG sense, but rather an 
\emph{ind-closed embedding}.}  

\end{lem}

\sssec{Proof of \lemref{l:loops into affine}}

At any level of coconnective truncation, we can write $\CY$ as a finite product
of affine D-schemes that are (potentially infinite) products of affine D-schemes of 
the form 
\begin{equation} \label{e:loops into free}
\Spec_X(A), \quad A=\Sym^!(\ind^l(\CE)),
\end{equation}
where $\CE$ is a vector bundle on $X$. 

\medskip

Since (ind-)affine (ind-)schemes are closed under finite limits and products, and the functors
$$\CY\mapsto \fL^+_\nabla(\CY) \text{ and } \CY\mapsto \fL_\nabla(\CY)$$
map products to products, we are reduced to considering $\CY$ of the form specified in \eqref{e:loops into free}. 

\medskip

In the latter case, the assertions of the lemma
can be (easily) checked directly (see \secref{sss:jets into affine space} below).

\qed

\sssec{Example} \label{sss:Gr G disc}

Recall the factorization space $\Gr_G$. By Beauville-Laszlo theorem, we can rewrite it as follows: for $\ul{x}:S\to \Ran$,
a lift of this point to a point of $\Gr_{G,\Ran}$ is the datum of $G$-bundle on $\cD_{\ul{x}}$ 
(which is equivalent to that of a $G$-bundle on $\wh\cD_{\ul{x}}$) and the trivialization of its restriction to $\cD^\times_{\ul{x}}$.

\medskip

From this description, we obtain a canonical projection
$$\fL(G)\to \Gr_G.$$

We claim that this projection identifies $\Gr_G$ with the \'etale quotient $\fL(G)/\fL^+(G)$. 
Indeed, this follows from the fact that a $G$-bundle on $\wh\cD_{\ul{x}}$ can be trivialized
after an \'etale sheafification along $S$ (see \secref{sss:et top disc}). 

\medskip

A similar description applies to factorization $\Gr_G$-module space $\Gr^{\on{level}_\CZ}_G$,
see \secref{sss:Gr G level}.

\ssec{Digression: the jet construction}

\sssec{}  \label{sss:jet construction}

Let $\CY$ be a prestack over $X$. Its jets construction, denoted
$$\on{Jets}(\CY)\to X_\dr$$
is by definition restriction of scalars \`a la Weil of $\CY$ along the projection
$$X\to X_\dr.$$

\sssec{}

Explicitly for $x:S\to X_\dr$, a lift of this map to a map $S\to \on{Jets}(\CY)$ is an $X$-map 
$$\wh\cD_x\to \CY.$$

\sssec{} \label{sss:jets into affine}

Suppose for a moment that $\CY$ is affine over $X$, i.e., 
$$\CY=\Spec_X(A_0), \quad
A_0\in \on{ComAlg}(\QCoh(X)^{\leq 0}).$$ In this case
$$\on{Jets}(\CY)=\Spec_X(A),$$
where $A$ is the D-algebra obtained by applying to $A_0$ the left adjoint of the forgetful functor
$$\oblv^l:\on{ComAlg}(\Dmod(X))\to \on{ComAlg}(\QCoh(X)).$$

For example, when 
$$A_0=\on{Sym}(\CE), \quad \CE\in \QCoh(X),$$
the above left adjoint produces 
$$\on{Sym}^!(\ind^l(\CE)).$$

\sssec{}

Denote
$$\fL^+(\CY):=\fL^+_\nabla(\on{Jets}(\CY)).$$

One can tautologically rewrite the definition of $\fL^+(\CY)$ as follows. For $\ul{x}:S\to \Ran$,
its lift to $\fL^+(\CY)_\Ran$ is an $X$-map
$$\wh\cD_{\ul{x}}\to \CY.$$

\sssec{} \label{sss:just loops}

Assume again that $\CY$ is affine over $X$. Then $\on{Jets}(\CY)$ is affine over $X_\dr$
(see \secref{sss:jets into affine}), and one can consider 
$$\fL(\CY):=\fL_\nabla(\on{Jets}(\CY)).$$

Explicitly, for an affine test scheme $S$ and $\ul{x}:S\to \Ran$,
its lift to a map $\fL(\CY)_\Ran$ is an $X$-map
$$\cD^\times_{\ul{x}}\to \CY.$$

Note that in this case a lift to $\fL^+(\CY)_\Ran$ can also be described as an $X$-map
$$\cD_{\ul{x}}\to \CY.$$

\sssec{Example} \label{sss:jets into affine space} 

Let $\CY=\Spec_X(\Sym(\CE))$, where $\CE$ is a vector bundle on $X$. Then the spaces
$$\fL^+(\CY) \text{ and } \fL(\CY)$$
can be described explicitly as follows.

\medskip

Fix a map $\ul{x}:S\to \Ran$. Fix lifts of the maps $x_i:S_{\on{red}}\to X$ to maps $\wt{x}_i:S\to X$. 
Let $D$ be the divisor on $S\times X$ equal 
$$\underset{i}\Sigma\, n_i\cdot \on{Graph}_{\wt{x}_i}$$
for some/any choice of $n_i\geq 1$.

\medskip

Then a lift of $\ul{x}$ to a point of $\fL^+(\CY)$ is a point of
$$\underset{n}{\on{lim}}\, \Gamma(S\times X,\CO_X\otimes \CE^\vee/\CO_X\otimes \CE^\vee(-n\cdot D)).$$

A lift of $\ul{x}$ to a point of $\fL(\CY)$ is a point of
$$\underset{m}{\on{colim}}\, 
\underset{n}{\on{lim}}\, \Gamma(S\times X,\CO_X\otimes \CE^\vee(m\cdot D)/\CO_X\otimes \CE^\vee(-n\cdot D)).$$

\sssec{Example}

Let $H$ be a smooth group-scheme over $X$. On the one hand, we can consider the algebraic stack\footnote{Throughout the paper
we write $\on{pt}/H$, where we mean the \'etale sheafification of $B(H)$, where the latter is the quotient of the base
(in this case, $X$)  by the trivial action of $H$.} 
$\on{pt}/H$ over $X$,
and consider the corresponding D-prestack
$$\on{Jets}(\on{pt}/H).$$

\medskip

On the other hand, we can consider 
$$\on{pt}/\on{Jets}(H),$$
i.e., the \'etale sheafification of $B(\on{Jets}(H))$. 

\medskip

We have a tautological map
\begin{equation} \label{e:jets into quotient}
\on{pt}/\on{Jets}(H)\to \on{Jets}(\on{pt}/H).
\end{equation} 

We claim that \eqref{e:jets into quotient} is an isomorphism. Indeed, this follows from the fact that if
an $H$-bundle on $\wh\cD_{\ul{x}}$ is such that its restriction to $S$ is trivial, then it is itself trivial. 

\ssec{Digression: the notion of (almost) finite presentation in the D-sense} \label{ss:afp D}

\sssec{} \label{sss:afp D n}

Fix a natural number $n$, and consider the category $\on{ComAlg}(\Dmod(X))^{\leq 0,\geq -n}$ of connective $n$-coconnective
commutative algebras in $\Dmod(X)$, which is by definition
$$\on{ComAlg}(\Dmod(X))\underset{\QCoh(X)}\times \QCoh(X)^{\leq 0,\geq -n},$$
where the functor $\on{ComAlg}(\Dmod(X))\to \QCoh(X)$ is
$$\on{ComAlg}(\Dmod(X))\overset{\oblv_{\on{ComAlg}}}\longrightarrow \Dmod(X)\overset{\oblv^l}\to \QCoh(X).$$

\medskip

We shall say that $A\in \on{ComAlg}(\Dmod(X))^{\leq 0,\geq -n}$ is $n$-D-afp if it is \emph{compact} as an object of this category.

\sssec{}  \label{sss:afp D n retract}

Suppose that $A\in \on{ComAlg}(\Dmod(X))^{\leq 0,\geq -n}$ is isomorphic to the 
geometric realization of a simplicial object 
$A_\bullet$ in $\on{ComAlg}(\Dmod(X))^{\leq 0,\geq -n}$ with terms of the form
$$A_n=\Sym^!(\CM_n), \quad \CM_n\in \Dmod(X), \quad \oblv^l(\CM_n)\in \Dmod^{\heartsuit,\on{f.g.}}.$$

\medskip

It is clear that such $A$ is compact. 

\medskip

Furthermore, it easy to see that objects of this form generate
$\on{ComAlg}(\Dmod(X))^{\leq 0,\geq -n}$ under filtered colimits. Furthermore, for the generation 
statement, we can take $\CM_n$ to be locally free, i.e., of the form $\ind^l(\CE_n)$, where 
$\CE_n$ is a vector bundle on $X$. 

\medskip

From here it follows that every $n$-D-afp algebra can be written as a retract of such $|A_\bullet|$. 

\sssec{} \label{sss:afp D}

Consider now the category $\on{ComAlg}(\Dmod(X))^{\leq 0}$ of connective
commutative algebras in $\Dmod(X)$, i.e., 
$$\on{ComAlg}(\Dmod(X))\underset{\QCoh(X)}\times \QCoh(X)^{\leq 0}.$$

\medskip

We shall say that $A$ is D-afp if for every $n$, the truncation $\tau^{\geq -n}(A)\in \on{ComAlg}(\Dmod(X))^{\leq 0,\geq -n}$ 
is $n$-D-afp. 

\medskip

We shall say that $\CY=\Spec_X(A)$ is D-afp if $A$ is. 

\medskip

Here are some examples of D-afp algebras. 

\sssec{}

Let $\CY_0\to X$ be an affine scheme almost of finite type, and take $\CY:=\on{Jets}(\CY_0)$. 
Then $\CY$ is D-afp. 

\medskip

Indeed, let $\CY_0=\Spec_X(A_0)$, where $A_0\in \on{ComAlg}(\QCoh(X)^{\leq 0})$. 
We can write 
$$A_0=|A_{0,\bullet}|, \quad A_{0,n}=\Sym_{\CO_X}(\CE_n), \quad \CE_n\in \on{Perf}(X)^\heartsuit.$$

Then $\CY=\Spec_X(A)$ for
$$A=|A_\bullet|, \quad A_n=\Sym^!(\ind^l(\CE_n)).$$

\sssec{}

Let now $\CY$ be the constant affine D-scheme with fiber $\CY_0$, where $\CY_0$ is almost of finite type. 
We claim that it is D-afp.

\medskip

Indeed, write $\CY_0=\Spec(A_0)$, where
$$A_0=|A_{0,\bullet}|, \quad A_{0,\bullet}=\Sym(V_n), \quad V_n\in \Vect^{\heartsuit,\on{f.d.}}.$$

Then $\CY=\Spec(A)$ where 
$$A=|A_\bullet|, \quad A_n=\Sym(V_n)\otimes \omega_X[1]=\Sym^!(V_n\otimes \omega_X[1])$$
(recall that according to our conventions in \secref{sss:omega X}, the object $\omega_X[1]\in \Dmod(X)$ is
the dualizing sheaf on $X$). 

\sssec{}

Let $\CY\to \CY_0$ be a map of affine D-schemes. The notion of being D-afp has a straightforward analog 
in this situation: 

\medskip

If $\CY_0=\Spec_X(A_0)$ and $\CY=\Spec_X(A)$, so that $\CY\to \CY_0$ corresponds to a map
$A_0\to A$, it makes sense to talk about $A$ being finitely presented over $A_0$ in the D-sense. 

\medskip 

The notion of being D-afp is transitive: if
$$\CY''\to \CY'\to \CY$$
are maps of affine D-schemes, with $\CY'\to \CY$ and $\CY''\to \CY'$  D-afp, then so
is $\CY''\to \CY$.

%

\ssec{Local systems on the formal (punctured) disc}

\sssec{} \label{sss:LS loc}

Let $H$ be a (finite-dimensional) algebraic group. Consider the algebraic stack $\on{pt}/H$. 
We regard it as a constant D-prestack over $X$, 
i.e., 
$$\CY_\nabla:=\on{pt}/H\times X_\dr.$$

\medskip

We define the factorization space
$$\LS_H^\reg:=\fL^+_\nabla(\on{pt}/H).$$

\sssec{} \label{sss:et top disc}

Note that the natural map
\begin{equation} \label{e:arcs and B}
\on{pt}/\fL^+_\nabla(H)\to \fL^+_\nabla(\on{pt}/H)=:\LS_H^\reg
\end{equation}
is an isomorphism (as was mentioned earlier, the notation $\on{pt}/\fL^+_\nabla(H)$ means
the \'etale sheafification of $B(\fL^+_\nabla(H))$, the quotient of the base by the trivial action of the corresponding group-scheme).

\medskip

Indeed, this follows from the fact that for $\ul{x}:S\to \Ran$, the \'etale topology on $S$ generates
the \'etale topology on $\wh\cD_{\ul{x}}$. 

\sssec{}

First, we claim:

\begin{lem} \label{l:pt to pt/H}
The map $\on{pt}\to \LS^\reg_H$ is \'etale-locally surjective and is 
an fpqc cover.
\end{lem}

\begin{proof}

The \'etale surjectivity follows from the isomorphism \eqref{e:arcs and B}. 

\medskip

Hence, to prove the lemma it remains to show that
$$\on{pt}\underset{ \LS^\reg_H}\times \on{pt}\simeq \fL^+_\nabla(H)$$
is flat. I.e., we need to show that for every finite set $I$, the scheme $\fL^+_\nabla(H)_{X^I}$
is flat over $X^I$. 

\medskip 

The latter is the assertion that for any flat D-algebra $A$ on $X$, the factorization algebra
$\on{Fact}(A)$ (see \eqref{e:arcs as fact again}) has the property that for every $I$, the restriction 
$\oblv^l(\on{Fact}(A)_{X^I})$ is flat over $X^I$ (this is essentially \cite[Lemma 3.4.12]{BD2}).

\end{proof}

\sssec{}

One of the key facts about $\LS_H^\reg$ is that it can be accessed via \emph{gauge} forms. Namely,
consider the tautological map
$$\on{pt}/H\to \on{Jets}(\on{pt}/H)\simeq \on{pt}/\on{Jets}(H)$$

Note that the fiber product
$$\on{pt}/H\underset{\on{pt}/\on{Jets}(H)}\times \on{pt}$$
identifies with the D-scheme
$$\on{Jets}(\fh\otimes \omega_X)=:\on{Conn}(\fh)$$
of jets into $\fh\otimes \omega_X$ (i.e., the total space
of the corresponding vector bundle, viewed as an affine 
scheme over $X$).\footnote{Indeed, the map 
$\on{dlog}:\on{Jets}(H)\to \on{Jets}(\fh\otimes \omega_X)$ identies the target with $\on{Jets}(H)/H$.}
The resulting $\on{Jets}(H)$-action on $\on{Conn}(\fh)$ is called \emph{the gauge action}. 

\medskip

Hence, we can identify
$$\on{pt}/H \simeq \on{Conn}(\fh)/\on{Jets}(H)$$
as D-prestacks over $X$. 

\medskip

In particular, we obtain an action of $\fL^+(H)$ on $\fL^+_\nabla(\on{Conn}(\fh))$ as factorization spaces, and an
identification: 
\begin{equation} \label{e:LS via gauge}
\LS_H^\reg \simeq \fL^+_\nabla(\on{Conn}(\fh))/\fL^+(H).
\end{equation}

\sssec{}

We note:

\begin{lem} \label{l:LS formall smooth}
The factorization space $\LS^\reg_H$ is fomally smooth. 
\end{lem}

\begin{proof}

This follows from \eqref{e:LS via gauge}. 

\end{proof}

\begin{rem}

Note that the map 
$$\on{pt}\to \LS^\reg_H$$
is \emph{not} formally smooth. Indeed, the fiber product
$$\on{pt}\underset{\LS^\reg_H}\times \on{pt}$$
is the affine factorization scheme $\fL^+_\nabla(H)$, associated to the constant
D-scheme $H$ (see \secref{sss:fact schemes ass to aff}), and one can show that 
$\fL^+_\nabla(H)_{X^2}$ is \emph{not} formally smooth over $X^2$.

\medskip

Namely, the cotangent space to $\fL^+_\nabla(H)_{X^2}$ at the unit section 
is the object of $\QCoh(X^2)$ underlying the \emph{left} D-module 
\begin{equation} \label{e:cotan G}
\on{Fib}\left(j_*((\fh^*\oplus \fh^*)\otimes \CO_{X^2-\Delta(X)})\to \Delta_*(\fh^*\otimes \CO_X)\right),
\end{equation} 
where the map in \eqref{e:cotan G} is the composition
$$\on{Fib}\left(j_*((\fh^*\oplus \fh^*)\otimes \CO_{X^2-\Delta(X)})\to \Delta_*((\fh^*\oplus \fh^*)\otimes \CO_X)\right)\to 
\Delta_*(\fh^*\otimes \CO_X),$$
where the second arrow is induced by the addition map $\fh^*\oplus \fh^*\to \fh^*$. 

\medskip

The object \eqref{e:cotan G} of $\QCoh(X^2)$ is flat, but not projective, implying that $\fL^+_\nabla(H)$ is not formally smooth. 

\medskip

Note that this is not in contradiction with the conclusion of \corref{c:mf formally smooth}. Let us see that 
$T^*_{X^2}(\LS^\reg_H)$ does satisfy the infinitesimal condition for formal smoothness, i.e., that
$$\Hom(T^*_{X^2}(\LS^\reg_H),\CF)=0 \text{ if } \CF\in \QCoh(X^2)^{\leq -1}.$$

We have
$$T^*_{X^2}(\LS^\reg_H)\simeq T^*_{X^2}(\fL^+_\nabla(H))[-1].$$

Now, we claim for any flat \emph{countably generated} object $\CE\in \QCoh(X^2)$, we have
$$\Hom(\CE[-1],\CF)=0 \text{ if } \CF\in \QCoh(X^2)^{\leq -1}.$$
Indeed, any such $\CE$ can be written as a countable filtered colimit of projective modules $\CE_i$, and
$$\CHom(\CE[-1],\CF)=\on{lim.proj.}\, \CHom(\CE_i[-1],\CF).$$

Here all $\CHom(\CE_i[-1],\CF)$ live in cohomological degrees $\leq -2$, while $\on{lim.proj.}$
has amplitude $1$, due to the countability assumption. 

\end{rem} 

\sssec{}  \label{sss:LS punctured bad}

We now proceed to defining the factorization space $\LS_H^\mer$. Naively, one would want to apply the functor $\fL_\nabla(-)$
to the constant D-prestack $\on{pt}/H$. But we cannot quite do this, and that is for two reasons, which already occur
for $\on{pt}/\on{Jets}(H)$: 

\medskip

\noindent{(i)} When considering $\fL_\nabla(-)$, we only allow affine targets.

\medskip

\noindent{(ii)} Even over a fixed point $x=\ul{x}\in \Ran$ (but an arbitrary affine scheme $S=\Spec(R)$ of parameters), the space of maps 
$$\cD^\times_x\to \on{pt}/H$$
is the space of \'etale $H$-torsors over the affine scheme $\Spec(R\ppart)$, where $t$ is a local coordinate near $x$.
However, we do not want to consider the \'etale topology on $\Spec(R\ppart)$. Rather, we want to \'etale-localize with respect to
$S=\Spec(R)$ itself, i.e., we only want to consider covers of $\Spec(R\ppart)$ of the form $\Spec(\wt{R}\ppart)$, where
$\Spec(\wt{R})\to \Spec(R)$ is \'etale. 

\sssec{} \label{sss:LS punctured}

The gauge action of $\on{Jets}(H)$ on $\on{Conn}(\fh)$ gives rise to
an action of $\fL(H)$ as a factorization group ind-scheme on $\fL_\nabla(\on{Conn}(\fh))$.

\medskip

We define $\LS_H^\mer$ to be
$$\fL_\nabla(\on{Conn}(\fh))/\fL(H),$$
the \'etale sheafification of the non-sheafified quotient of
$\fL_\nabla(\on{Conn}(\fh))$ by the gauge action of $\fL(H)$. 

\sssec{} \label{sss:LS punctured bis}

We can rephrase the definition of $\LS_H^\mer$ as follows:

\medskip

Consider the \'etale quotient 
$$\fL_\nabla(\on{Conn}(\fh))/\fL^+(H).$$

It carries an action of the groupoid
$$\fL^+(H)\backslash \fL(H)/\fL^+(H).$$

Then $\LS_H^\mer$ is the quotient of $\fL_\nabla(\on{Conn}(\fh))/\fL^+(H)$ by $\fL^+(H)\backslash \fL(H)/\fL^+(H)$,
subsequent sheafified in the \'etale topology. 

\sssec{} \label{sss:LS reg to LS mer}

By construction, we have a naturally defined map
\begin{equation} \label{e:LS reg to mer}
\LS_H^\reg\to \LS_H^\mer.
\end{equation}

We claim:

\begin{lem} \label{l:LS reg to LS mer}
The map \eqref{e:LS reg to mer} is ind-affine locally almost of finite presentation\footnote{See \secref{sss:lafp} for what this means.}.
\end{lem}

\begin{proof}

By construction, it suffices to show that the fiber product
\begin{equation} \label{e:gauge those forms}
\LS_H^\reg\underset{\LS_H^\mer}\times \fL_\nabla(\on{Conn}(\fh))\to \fL_\nabla(\on{Conn}(\fh))
\end{equation}
is ind-schematic locally almost of finite presentation. 

\medskip

The left-hand side in \eqref{e:gauge those forms} is the space
$$\{g\in \fL(H),\alpha\in \fL_\nabla(\on{Conn}(\fh))\,|\, g\cdot \alpha\in \fL^+_\nabla(\on{Conn}(\fh))\}/\fL^+(H).$$

In other words, we can rewrite it as
$$(\Gr_H\times \fL_\nabla(\on{Conn}(\fh)))\underset{\fL_\nabla(\on{Conn}(\fh)/\fL^+(H)}\times \fL^+_\nabla(\on{Conn}(\fh))/\fL^+(H)$$
and its map to the right-hand side of \eqref{e:gauge those forms} is the composition
\begin{multline} \label{e:gauge those forms 1}
(\Gr_H\times \fL_\nabla(\on{Conn}(\fh)))\underset{\fL_\nabla(\on{Conn}(\fh)/\fL^+(H)}\times \fL^+_\nabla(\on{Conn}(\fh))/\fL^+(H)\to \\
\to \Gr_H\times \fL_\nabla(\on{Conn}(\fh))\to \fL_\nabla(\on{Conn}(\fh)).
\end{multline} 

Since the map $\fL^+_\nabla(\on{Conn}(\fh))\to \fL_\nabla(\on{Conn}(\fh))$ is an ind-closed embedding locally almost of finite presentation, we 
obtain that the first arrow in \eqref{e:gauge those forms 1} has this property. 

\medskip

The second arrow in \eqref{e:gauge those forms 1} is ind-schematic and locally almost of finite presentation since 
$\Gr_H$ is an ind-scheme locally almost of finite type. Hence \eqref{e:gauge those forms} is ind-schematic and locally almost of finite presentation.

\medskip

Finally, let us show that \eqref{e:LS reg to mer} is ind-affine. Let $H'$ denote the reductive quotient of $H$. Factor the map
\eqref{e:LS reg to mer} as
\begin{equation} \label{e:LS reg to mer 1}
\LS_H^\reg\to \LS_H^\mer\underset{\LS_{H'}^\mer}\times \LS_{H'}^\reg\to \LS_H^\mer,
\end{equation}
and it is enough to show that both arrows in \eqref{e:LS reg to mer 1} are ind-affine. 

\medskip
 
For the first arrow in \eqref{e:LS reg to mer 1}, after base-changing to the fpqc cover $\on{pt}\to \LS_{H'}^\reg$, it suffices to show that the map
$$\LS^\reg_{H''}\to \LS^\mer_{H''}$$
is ind-affine, where $H''$ is the unipotent radical of $H$. This follows from the fact that in this
case the second arrow \eqref{e:gauge those forms 1} is ind-affine, since the affine Grassmannian for a unipotent group
is ind-affine. 

\medskip 

For the second arrow in \eqref{e:LS reg to mer 1}, it suffices to show that the original map \eqref{e:LS reg to mer}
is ind-affine when $H$ is reductive. Note that in this case $\Gr_H$ is ind-proper, so the second arrow in \eqref{e:gauge those forms 1}
is ind-proper. Hence, the map \eqref{e:LS reg to mer} is ind-proper. However, since it is also injective at the level of $k$-points, we
obtain that it is ind-affine (indeed, a proper map between schemes that is injective at $k$-points is affine). 

\end{proof}

\sssec{} \label{sss:LS open curve}

We now discuss a global version of the factorization space $\LS_H^\mer$. 

\medskip

We define the
(non-factorization!) space
$$\LS^{\mer,\on{glob}}_{H,\Ran}\to \Ran$$
as follows:

\medskip

For an affine test scheme $S$ and a map $\ul{x}\to \Ran$, its lift to $\LS^{\mer,\on{glob}}_{H,\Ran}$ is a datum of
a map
\begin{equation} \label{e:LS open curve 1}
(S\times X_\dr-\on{Graph}_{\ul{x}}) \to \on{pt}/H
\end{equation} 
such that \'etale-locally on $S$, the composite map
$$(S\times X-\on{Graph}_{\ul{x}}) \to (S\times X_\dr-\on{Graph}_{\ul{x}}) \to \on{pt}/H$$
\emph{admits} an extension to a map 
\begin{equation} \label{e:LS open curve 2}
S\times X\to \on{pt}/H.
\end{equation} 

\sssec{} \label{sss:LS open curve bis}

We claim that we have a naturally defined evaluation map
$$\on{ev}_\Ran:\LS^{\mer,\on{glob}}_{H,\Ran}\to \LS^\mer_{H,\Ran}.$$

Indeed, for $\ul{x}:S\to \Ran$, 
given a map \eqref{e:LS open curve 1} and a lift \eqref{e:LS open curve 2}, the restriction 
of the connection form to $\cD_{\ul{x}}$ gives rise to a section of 
$$\fL_\nabla(\on{Conn}(\fh))_S/\fL^+(H)_S.$$

\medskip

A modification of \eqref{e:LS open curve 2} results in an action of the groupoid
$$\fL^+(H)_S\backslash \fL(H)_S/\fL^+(H)_S.$$

\ssec{Sheaves of categories over the Ran space} \label{ss:cat over Ran}

\sssec{}

Let $\CY$ be a prestack. When discussing sheaves of categories over $\CY$, we will
assume that $\CY$ is locally almost of finite type. In this context, when considering 
affine schemes $S$ or general prestacks $\CZ$ mapping to $\CY$, we will assume that
they are also locally almost of finite type. 

\sssec{} \label{sss:shvs of cat}

 A sheaf of categories $\ul\bC$ on $\CY$ is an assignment
$$(S,y)\in \affSch_{/\CY}\rightsquigarrow \bC_{S,y}\in \QCoh(S)\mmod,$$
equipped with identifications: 
$$S'\overset{f}\to S, \quad \bC_{S',y\circ f}\simeq \QCoh(S')\underset{\QCoh(S)}\otimes \bC_{S,y},$$
satisfying a homotopy-coherent system of compatibilities; we refer the reader to \cite{Ga5} for details. 

\medskip

Let $\on{ShvCat}(\CY)$ denote the (2-)-category of sheaves of categories over $\CY$.
It has a natural symmetric monoidal structure given by
$$(\bC_1\otimes \bC_2)_{S,y}:=\bC_{1,S,y}\underset{\QCoh(S)}\otimes \bC_{2,S,y}.$$

\sssec{}

A basic example of a sheaf of categories over $\CY$ is
$$S\rightsquigarrow \QCoh(S);$$
we denote it by $\ul\QCoh(\CY)$. 

\medskip

This is the unit for the above symmetric monoidal structure on $\on{ShvCat}(\CY)$. 

\sssec{}

Given a sheaf of categories $\ul\bC$ on $\CY$, the category of its global sections,
denoted $\Gamma(\CY,\ul\bC)$ is defined as
$$\underset{(S,y)\in \affSch_{/\CY}}{\on{lim}}\, \bC_{S,y}.$$

\medskip

This category is naturally a module over
$$\underset{(S,y)\in \affSch_{/\CY}}{\on{lim}}\, \QCoh(S)=:\QCoh(\CY).$$

\medskip

Thus, the functor 
$$\Gamma(\CY,-):\on{ShvCat}(\CY)\to \DGCat$$
upgrades to a functor
$$\Gamma(\CY,-)^{\on{enh}}:\on{ShvCat}(\CY)\to \QCoh(\CY)\mmod.$$

\sssec{Example}

We have
$$\Gamma(\CY,\ul\QCoh(\CY))\simeq \QCoh(\CY),$$
as a module over itself.

\sssec{}

A prestack $\CY$ is said to be \emph{1-affine} if the functor $\Gamma(\CY,-)^{\on{enh}}$
is an equivalence.

\medskip

In the paper \cite{Ga5} a number of results was proved, showing that various classes of 
prestacks are 1-affine.

\medskip

The most relevant for us are:

\begin{itemize}

\item The stack $\on{pt}/H$, where $H$ is a (finite-dimensional) algebraic group, is 1-affine.

\item The stack $Z_\dr$, where $Z$ is a scheme of finite type, is 1-affine. 

\end{itemize} 

\sssec{} \label{sss:pullback shf of cat}

Let $f:\CY_1\to \CY_2$ be a map between prestacks. We define a functor
$$f^*:\on{ShvCat}(\CY_2)\to \on{ShvCat}(\CY_1)$$
by sending $\ul\bC\in \on{ShvCat}(\CY_2)$ to $f^*(\ul\bC)\in  \on{ShvCat}(\CY_1)$ that assigns 
$$(S,y)\in \affSch_{/\CY_1} \rightsquigarrow \bC_{S,f\circ y}\in \QCoh(S)\mmod.$$

\medskip

The above functor $f^*$ has a right adjoint, denoted $f_*$. Explicitly, $f_*$ sends $\ul\bC\in \on{ShvCat}(\CY_1)$ to
$f_*(\ul\bC)\in  \on{ShvCat}(\CY_2)$ that assigns 
$$(S,y)\in \affSch_{/\CY_2} \rightsquigarrow \Gamma(S\underset{\CY_1}\times \CY_2,y'{}^*(\ul\bC))\in 
\QCoh(S\underset{\CY_1}\times \CY_2)\mmod \overset{f'{}^*}\to \QCoh(S)\mmod,$$
where $y'$ denotes the map
$$S\underset{\CY_1}\times \CY_2\to \CY_1$$
and $f'$ denotes the map
$$S\underset{\CY_1}\times \CY_2\to S.$$

\sssec{}

Let $\ul\bC$ be a sheaf of categories over $\CY$. We shall say that $\ul\bC$ is compactly generated
if for every $S\in \affSch_{/\CY}$, the category $\bC_S$ is compactly generated. 

\sssec{}

Suppose for a moment that $\ul\bC$ is pulled back from a sheaf of categories on $\CY_\dr$.
Hence, for $S\in \affSch_{/\CY}$ we have a well-defined category $\bC_{S_\dr}$. Note that if
$\bC_S$ is compactly generated, then so is $\bC_{S_\dr}$:  

\medskip

Indeed, with no restriction of generality we can assume that $S$ is eventually coconnective. By 
the 1-affineness of $S_\dr$, we have
$$\bC_S\simeq \QCoh(S)\underset{\Dmod(S)}\otimes  \bC_{S_\dr},$$
and the pair of ($\Dmod(S)$-linear) adjoint functors
$$\ind^l:\QCoh(S)\rightleftarrows \Dmod(S):\oblv^l$$
induces an adjoint pair
$$(\ind^l\otimes \on{Id}): \bC_S \rightleftarrows \bC_{S_\dr}: (\oblv^l\otimes \on{Id}),$$
where the essential image of the left adjoint generates the target category. 

\sssec{} 

Let $F:\ul\bC_1\to \ul\bC_2$ be a functor between sheaves of categories. Suppose that $\ul\bC_1$
is compactly generated (in the above sense). 

\medskip

We shall say that $F$ preserves compactness if for every $S\in \affSch_{/\Ran}$, the corresponding functor
$$F_S:\bC_{1,S}\to \bC_{2,S}$$
preserves compactness.

\medskip

Then the usual argument (using the fact that for an affine schemes $S$, the category $\QCoh(S)$ is rigid)
shows that in this case the functor $F$ admits a right adjoint, to be denoted $F^R$, as a functor between
sheaves of categories over $\Ran$.

\sssec{}   \label{sss:t-str on shf of cat}

Let $\ul\bC$ be a sheaf of categories over $\CY$. A t-structure on $\bC$ is a collection of t-structures on
$\bC_S$ for any $S\in \affSch_{/\CY}$ such that:

\medskip

For any $f:S'\to S$, the functor
$$f_*:\bC_{S'}\to \bC_S$$
is t-exact. 

\medskip

Note that this condition can be rewritten as saying that with respect to the identification
$$\bC_{S'}\simeq \QCoh(S')\underset{\QCoh(S)}\otimes \bC_S,$$
the t-structure on $\bC_{S'}$ is the \emph{tensor product} t-structure, i.e., $\left(\QCoh(S')\underset{\QCoh(S)}\otimes \bC_S\right)^{\leq 0}$
is generated under colimits by the essential image of
$$\QCoh(S')^{\leq 0}\times  \bC_S^{\leq 0}\to \QCoh(S')\underset{\QCoh(S)}\otimes \bC_S.$$

\sssec{}

Let $\ul\bC$ be equipped with a t-structure. For $S\in \affSch_{/\CY}$, consider the category
$$\IndCoh(S)\underset{\QCoh(S)}\otimes\bC_S.$$

\medskip 

We equip it with the tensor product t-structure. 

\medskip

Assume now that the t-structure on $\bC_S$ is right-complete and compatible with filtered colimits. Then 
a standard argument shows that the functor
$$\IndCoh(S)\underset{\QCoh(S)}\otimes\bC_S \overset{\Psi_S\otimes \on{Id}}\longrightarrow 
\QCoh(S)\underset{\QCoh(S)}\otimes\bC_S\simeq \bC_S$$
is t-exact and induces an equivalence of the eventually coconnective subcategories of the two sides
(see, e.g., \cite[Proposition C.4.6.1]{Lu3}). 

\sssec{} 

Suppose again that $\ul\bC$ is the pull back of a sheaf of categories on $\CY_\dr$. 
For $S$ as above consider the category $\bC_{S_\dr}$. Note that we can identify
$$\IndCoh(S)\underset{\QCoh(S)}\otimes\bC_S\simeq  \IndCoh(S)\underset{\QCoh(S)}\otimes  \QCoh(S) \underset{\Dmod(S)}\otimes
\bC_{S_\dr}\simeq \IndCoh(S)\underset{\Dmod(S)}\otimes \bC_{S_\dr}.$$

The pair of ($\Dmod(S)$-linear) adjoint functors
$$\ind^r:\IndCoh(S) \rightleftarrows \Dmod(S):\oblv^r$$
gives rise to an adjucntion
$$(\ind^r\otimes \on{Id}): \IndCoh(S)\underset{\QCoh(S)}\otimes\bC_S \rightleftarrows \bC_{S_\dr}:(\oblv^r\otimes \on{Id}).$$

We define a t-structure on $\bC_{S_\dr}$ by letting $\bC_{S_\dr}^{\leq 0}$ be generated under colimits
by the essential image under $(\ind^r\otimes \on{Id})$ of $(\IndCoh(S)\underset{\QCoh(S)}\otimes\bC_S)^{\leq 0}$. 

\medskip

Assume for a moment that $S$ is smooth. In this case the functor $\ind^r$ is t-exact, and hence the 
endofunctor $\oblv^r\circ \ind^r$ of $\IndCoh(S)$ is t-exact. This implies that the functor
$$(\oblv^r\otimes \on{Id}):\bC_{S_\dr}\to \IndCoh(S)\underset{\QCoh(S)}\otimes\bC_S$$
is t-exact. 

\sssec{}

We now specialize the case of $\CY=\Ran$. We claim:

\begin{lem}  \label{l:Ran is 1-affine}
The prestack $\Ran$ is 1-affine.
\end{lem}

\begin{proof}

The proof that we will give applies to any prestack $\CY$ that can be written as a colimit
$$\underset{i\in I}{\on{colim}}\, (Z_i)_\dr,$$
where $Z_i$ are schemes and the transition maps are proper, and for which the diagonal morphism
$$\CY\to \CY\times \CY$$
is closed (at the reduced level\footnote{I.e., the base change of this morphism by an affine scheme yields a morphism of
prestacks, such that the morphism of the underlying reduced prestacks is a closed embedding.}). 

\medskip

We can think of an object of $\on{ShvCat}(\CY)$ as a compatible collection of categories $\{\bC_i\}$
$$\bC_i \in \Dmod(Z_i), \quad (Z_i\overset{f_{i,j}}\to Z_j)\rightsquigarrow \bC_i\simeq \Dmod(Z_i)\underset{\Dmod(Z_j)}\otimes \bC_j.$$

The functor $\Gamma(\CY,-)$ sends such a collection to
$$\underset{i\in I^{\on{op}}}{\on{lim}}\, \bC_i,$$
viewed as a module over
$$\QCoh(\CY)=\Dmod(\CY)\simeq \underset{i}{\on{lim}}\, \Dmod(Z_i).$$

Given a $\Dmod(\CY)$-module category $\bD$, we attach to it an object of $\on{ShvCat}(\CY)$ by setting
$$\bC_i:=\on{Funct}_{\Dmod(\CY)}(\Dmod(Z_i),\bD).$$

For $(i\to j)\in I$, the corresponding transition functor
$$\on{Funct}_{\Dmod(\CY)}(\Dmod(Z_j),\bD)\to \on{Funct}_{\Dmod(\CY)}(\Dmod(Z_i),\bD)$$
is given by precomposition with $(f_{i,j})_!$. 

\medskip

For $(i\to j)\in I$, we have, tautologically,
$$\on{Funct}_{\Dmod(Z_j)}(\Dmod(Z_i),\on{Funct}_{\Dmod(\CY)}(\Dmod(Z_j),\bD)) \simeq \on{Funct}_{\Dmod(\CY)}(\Dmod(Z_i),\bD).$$

However, this implies that we also have 
$$\Dmod(Z_i)\underset{\Dmod(Z_j)}\otimes \on{Funct}_{\Dmod(\CY)}(\Dmod(Z_j),\bD)\simeq 
\on{Funct}_{\Dmod(\CY)}(\Dmod(Z_i),\bD),$$
since $\Dmod(Z_i)$ is self-dual as a $\Dmod(Z_j)$-module. 

\medskip

Let us establish the equivalence
\begin{equation} \label{e:Ran is 1-affine}
\underset{i\in I^{\on{op}}}{\on{lim}}\,  \on{Funct}_{\Dmod(\CY)}(\Dmod(Z_i),\bD) \simeq \bD.
\end{equation}

Indeed, we can rewrite the left-hand side as
$$\on{Funct}_{\Dmod(\CY)}(\underset{i\in I}{\on{colim}}\, \Dmod(Z_i),\bD),$$
where the colimit is taken with respect to the $(f_{i,j})_!$-functors.

\medskip

Now, we can rewrite
$$\underset{i\in I}{\on{colim}}\, \Dmod(Z_i)\simeq\underset{i\in I^{\on{op}}}{\on{lim}}\, \Dmod(Z_i),$$
where in the right-hand side the limit is taken with respect to the $f_{i,j}^!$ functors. Finally,
$$\underset{i\in I^{\on{op}}}{\on{colim}}\, \Dmod(Z_i)\simeq \Dmod(\CY),$$
by definition. Hence, the left-hand side in \eqref{e:Ran is 1-affine} identifies with
$$\on{Funct}_{\Dmod(\CY)}(\Dmod(\CY),\bD) \simeq \bD,$$
as required. 

\medskip

It is easy to see that in order to show that the above two functors are mutually inverse, it remains to show that
for $i_1,i_2\in I$, and the corresponding maps
$$Z_{i_1}\overset{f_{i_1}}\to \CY \overset{f_{i_2}}\leftarrow Z_{i_2}$$
the naturally defined functor
\begin{equation} \label{e:Ran is 1-affine bis}
\Dmod(Z_{i_1}\underset{\CY}\times Z_{i_2})\to \on{Funct}_{\Dmod(\CY)}(\Dmod(Z_{i_1}),\Dmod(Z_{i_2}))
\end{equation}
is an equivalence. 

\medskip

We calculate the right-hand side in \eqref{e:Ran is 1-affine bis} as the totalization of the cosimplicial
category with terms 
$$\Dmod(Z_{i_1}\times \CY^m\times Z_{i_2}),$$
with terms given by !-pushforward functors along the maps in the corresponding cosimplicial prestack.

\medskip

However, due to the assumption that $\CY$ has a closed diagonal, the face maps in this 
cosimplicial prestack are closed embeddings, and hence the corresponding !-pushforward functors
are fully faithful. Hence, the above totalization is the equalizer of 
$$\Dmod(Z_{i_1}\times Z_{i_2}) \rightrightarrows  \Dmod(Z_{i_1}\times \CY \times Z_{i_2}).$$

Objects of this equalizer are supported on
$$(Z_{i_1}\times Z_{i_2})\underset{\on{Graph}_{f_{i_1}}\times \on{id},
Z_{i_1}\times \CY \times Z_{i_2},\on{id}\times \on{Graph}_{f_{i_2}}}\times (Z_{i_1}\times Z_{i_2}).$$

However, the above fiber product identifies with 
$$Z_{i_1}\underset{\CY}\times Z_{i_2}.$$

\end{proof} 

\begin{rem}

Let $\ul\bC$ be a sheaf of categories over $\Ran$. In the course of the proof of \lemref{l:Ran is 1-affine}, we have encountered
another way of how one may think of the category $\Gamma(\Ran,\ul\bC)$. Namely, 
$$\Gamma(\Ran,\ul\bC) \simeq \underset{I\in (\on{fSet}^{\on{surj}})^{\on{op}}}{\on{colim}}\, \bC_{X^I_\dr},$$
where the colimit is formed using the !-pushforward functors, i.e., for $(I\overset{\phi}\to J)\in \on{fSet}^{\on{surj}}$
and the corresponding map 
$$\Delta_\phi:X^J\to X^I,$$
the functor in question is
$$\bC_{X^J_\dr}\simeq \Dmod(X^J)\underset{\Dmod(X^I)}\otimes \bC_{X^I_\dr}\overset{(\Delta_\phi)_!\otimes \on{Id}}\longrightarrow
\Dmod(X^I)\underset{\Dmod(X^I)}\otimes \bC_{X^I_\dr}\simeq  \bC_{X^I_\dr}.$$

\end{rem} 

\sssec{Notational convention} \label{sss:cristalline obj}

%
%
%


We use the term \emph{crystal of categories} for sheaves of categories on 
$\CY_\dr$ and denote $$\on{CrystCat}(\CY):= \on{ShvCat}(\CY_\dr).$$
From this perspective, we refer to crystalline objects of $\ul{\bC}$ in $\on{CrystCat}(\CY)$ to mean 
the objects of the category of global sections of $\ul\bC$ as a sheaf of categories on $\CY_\dr$.

\medskip

When we use such ``crystalline" terminology and are given $f:Z \to \CY$, we use the symbol $\bC_Z$ to denote the crystalline objects 
of $\ul\bC$ restricted to $Z$, that is,
$$\bC_\CZ:=\Gamma(\CZ_\dr,(f_\dr)^*(\ul\bC)).$$

\ssec{Factorization algebras and modules} \label{ss:fact alg}

\sssec{}

A factorization algebra $\CA$ on $X$ is an object 
$$\CA_\Ran\in \Dmod(\Ran),$$
equipped with a datum of factorization
\begin{equation} \label{e:fact alg}
\on{union}^!(\CA_\Ran)|_{(\Ran\times \Ran)_{\on{disj}}}\simeq \CA_\Ran\boxtimes \CA_\Ran|_{(\Ran\times \Ran)_{\on{disj}}}
\end{equation} 
and with a homotopy-coherent data of associativity and commutativity, see \cite[Sect. 6]{Ra6} for details. 

\medskip

Let $\on{FactAlg}(X)$ denote the category of factorization algebras on $X$. 

\medskip

For $\CA\in \on{FactAlg}(X)$ and $\CZ\to \CA$, we will denote by $\CA_\CZ\in \Dmod(\CZ)$ the pullback of $\CA$ to $\CZ$. 

\sssec{Example} \label{sss:unit fact alg nu}

We let $k$ denote the unit factorization algebra. I.e., the underlying object
$$k_\Ran\in  \Dmod(\Ran)$$
is $\omega_\Ran$, equipped with the natural factorization structure. 

\sssec{} \label{sss:fact modules}

Let $\CA$ be a factorization algebra on $X$. Let $\CZ$ be a prestack equipped with a map
$\CZ\to \Ran$. Recall the space $\CZ^{\subseteq}$, see \secref{sss:Z Ran}. 

\medskip

A factorization $\CA$-module $\CM$ \emph{at} $\CZ$ is an object $$\CM_{\CZ^{\subseteq}}\in \QCoh(\CZ^{\subseteq})$$
equipped with a datum of factorization against $\CA$:
\begin{equation} \label{e:fact mod}
\CM_{\CZ^{\subseteq}}|_{(\Ran\times \CZ^{\subseteq})_{\on{disj}}}\simeq
(\CA_\Ran \boxtimes \CM_{\CZ^{\subseteq}})|_{(\Ran\times  \CZ^{\subseteq})_{\on{disj}}}
\end{equation} 
and a homotopy-coherent data of associativity.

\medskip

We denote the category of factorization $\CA$-modules \emph{at} $\CZ$ by
$$\CA\mod^{\on{fact}}_\CZ.$$

This category is naturally tensored over $\QCoh(\CZ)$ via the projection
$$\CZ^{\subseteq}\overset{\on{pr}_{\on{small},\CZ}}\longrightarrow \CZ.$$

\begin{rem}

For $\CZ=X^I$, one can described the category $\CA\mod^{\on{fact}}_{X^I}$ explicitly as chiral modules, see
\cite[Sect. 3]{Ro1}.

\end{rem}

\sssec{} 

Let $\on{diag}_\CZ:\CZ\to \CZ^{\subseteq}$ be the diagonal map (see \secref{sss:diag Z}). For $\CM\in \CA\mod^{\on{fact}}_\CZ$
we will denote by $\CM_\CZ$ the object 
$$\on{diag}_\CZ^*(\CM_{\CZ^{\subseteq}})\in \QCoh(\CZ).$$

We will refer to $\CM_\CZ$ as the \emph{quasi-coherent sheaf on $\CZ$ underlying $\CM$}. 

\medskip

The resulting functor 
\begin{equation} \label{e:forgetful functor fact mod}
\oblv_\CA:\CA\mod^{\on{fact}}_\CZ\to \QCoh(\CZ)
\end{equation}
is conservative and compatible with colimits. (Sometimes instead of $\oblv_\CA$ we write $\oblv_{\CA,\CZ}$ in order
to emphasize the dependence on $\CZ$.) 

\medskip

We will think of the datum of $\CM$ as $\CM_\CZ\in \QCoh(\CZ)$, equipped with an additional
datum of factorization against $\CA$. 

\sssec{} \label{sss:fact mod over unit nu}

Take $\CA=k$ from \secref{sss:unit fact alg nu}. We claim that there is a naturally defined functor
\begin{equation} \label{e:fact mod over unit}
\QCoh(\CZ)\to k\mod^{\on{fact}}_\CZ.
\end{equation}

Namely, the functor \eqref{e:fact mod over unit} is given by pullback along the projection
$$\CZ^\subseteq\to \CZ.$$

\medskip

In \secref{sss:unital fact mod for k}, we will upgrade the functor \eqref{e:fact mod over unit} to an equivalence of
categories, once we replace the right-hand side by the category of \emph{unital} factorization modules. 

\sssec{} \label{sss:vacuum module}

A basic example of a factorization $\CA$-module is the \emph{vacuum} module, denoted
$$\CA^{\on{fact}_\CZ}\in \CA\mod^{\on{fact}}_\CZ.$$

Namely, the corresponding object
$$\CA^{\on{fact}_\CZ}_{\CZ^{\subseteq}}\in \QCoh(\CZ^{\subseteq})$$
is the pullback of $\CA_\Ran$ along the map
$$\CZ^{\subseteq} \overset{\on{pr}_{\on{big}}}\longrightarrow \Ran.$$

Note that 
$$\CA^{\on{fact}_\CZ}_\CZ\simeq \CA_\CZ.$$

For that reason, we will sometimes abuse the notation and write $\CA_\CZ$ instead of $\CA^{\on{fact}_\CZ}$. 
(This is similar to denoting the free object over an associative algebra $A$ 
by the same symbol $A$ as the underlying vector space.)

\sssec{} \label{sss:propagate modules}

Generalizing the above construction, given $\CZ\to \Ran$ and $f:\CZ'\to \CZ^{\subseteq}$, which we perceive as mapping to $\Ran$
by means of 
$$\CZ'\to \CZ^{\subseteq} \overset{\on{pr}_{\on{big}}}\to \Ran,$$
there is a naturally defined functor\footnote{See \secref{sss:untl str on fact mod} for more details.} 
\begin{equation} \label{e:propagate modules}
\CA\mod^{\on{fact}}_\CZ\to \CA\mod^{\on{fact}}_{\CZ'}
\end{equation}
that makes the following diagram commute
$$
\CD
\CA\mod^{\on{fact}}_{\CZ'} @>{\oblv_{\CA,\CZ'}}>> \QCoh(\CZ') \\
@AAA @AA{f^!}A \\
\CA\mod^{\on{fact}}_{\CZ} @>>> \QCoh(\CZ^{\subseteq}).
\endCD
$$

\medskip

We will denote the functor \eqref{e:propagate modules} by
$$\CM\mapsto \CM|_{\CZ'}.$$

\begin{rem}
The functor \eqref{e:propagate modules} reflects the \emph{unital structure} on the assugnment
$$\CZ\mapsto \CA\mod^{\on{fact}}_\CZ,$$
see \secref{sss:untl str on fact mod}.
\end{rem}

\sssec{}

Let $S$ be an affine scheme mapping to $\Ran$, and let $S'\overset{f}\to S$ be a map of affine schemes. 
Pullback along
$$f^{\subseteq}:S'{}^{\subseteq}\to S^{\subseteq}$$
defines a functor
\begin{equation} \label{e:pullback fact mod}
f^*:\CA\mod^{\on{fact}}_S\to \CA\mod^{\on{fact}}_{S'}.
\end{equation} 

This functor is $\QCoh(S)$-linear, and hence gives rise to a functor
\begin{equation} \label{e:pullback fact mod tensored up}
\QCoh(S')\underset{\QCoh(S)}\otimes \CA\mod^{\on{fact}}_S\to \CA\mod^{\on{fact}}_{S'}.
\end{equation} 

We claim:

\begin{lem} \label{l:pullback fact mod tensored up}
The functor \eqref{e:pullback fact mod tensored up} is an equivalence.
\end{lem} 

\begin{proof}

Pushforward along $f^{\subseteq}$
gives rise to a functor
$$f_*:\CA\mod^{\on{fact}}_{S'}\to \CA\mod^{\on{fact}}_S,$$
which is a right adjoint of \eqref{e:pullback fact mod}. It is easy to see that it is monadic. 

\medskip

The functor \eqref{e:pullback fact mod tensored up} gives rise to a map of monads,
and in order to prove the lemma, it suffices to show that this map of monads induces
an isomorphism of the underlying endofunctors. 

\medskip

However, the latter is obvious: both endofunctors are given by tensoring by 
$f_*(\CO_{S'})$, viewed as an algebra on $\QCoh(S)$. 

\end{proof} 

\sssec{}

From \lemref{l:pullback fact mod tensored up} we obtain that the assignment
$$(S\to \Ran)\rightsquigarrow \CA\mod^{\on{fact}}_S$$
forms a sheaf of categories over $\Ran$.

\medskip

We will denote this sheaf of categories by $\CA\ul\mod^{\on{fact}}$. For any
$f:\CZ\to \Ran$, we have
$$\CA\mod^{\on{fact}}_\CZ\simeq \Gamma(\CZ,f^*(\CA\ul\mod^{\on{fact}})).$$

\sssec{}

We shall say that a factorization algebra $\CA$ is \emph{connective} 
if $$\oblv^l(\CA_X)\in \QCoh(X)$$
is connective, where $X\to \Ran$ is the tautological map. 

\sssec{}

Recall that for any finite non-empty set, we have a tautological map
$$X^I_\dr\to \Ran.$$

We have (by the argument in \cite[Lemma 3.4.12]{BD2}):

\begin{lem} \label{l:ch Chevalley perv}
Suppose that $\CA$ is connective in the above sense, and \emph{unital}
(see \secref{sss:untl fact alg}) for what this means. Then for any $I$, the object
$$\oblv^l(\CA_{X^I})\in \QCoh(X^I)$$
is connective.
\end{lem} 

\begin{cor}
Under the assumptions of \lemref{l:ch Chevalley perv}, for any \emph{scheme} $Z$ equipped with a map to $\Ran$, the object 
$$\CA_Z\in \QCoh(Z)$$
is connective.
\end{cor}

\sssec{} \label{sss:t-structure fact mod}

We now claim:

\begin{prop} \label{p:t-structure fact mod}
Let $\CA$ be connective.\footnote{The unitality assumption on $\CA$ here is irrelevant, as one can add a unit
to it (and then consider the corresponding category of unital modules).}
Then for any scheme $Z$ mapping to $\Ran$, the category
$\CA\mod^{\on{fact}}_{Z}$
carries a t-structure, uniquely characterized by the property that the forgetful 
functor \eqref{e:forgetful functor fact mod} 
is t-exact. Furthermore,
$\CA\mod^{\on{fact}}_Z$ is left-complete in its t-structure.
\end{prop} 

\begin{proof}

This follows by interpreting factorization modules as chiral modules, see \cite{FraG}.
The connectivity assumption translates into the fact that the corresponding chiral
algebra is connective (for the \emph{right} t-structure on $\Dmod(X)$):

\medskip

Namely, we consider the category $\QCoh(Z^{\subseteq})$ as acted on via the map
$\on{union}$ by $\Dmod(\Ran)$, equipped with the \emph{chiral} symmetric monoidal
structure (see \cite[Sect. 3.4.10]{BD2} and/or \cite[Sect 2.2]{FraG}). 

\medskip

Let 
$$\CA^{\on{ch}}\in \Dmod(X)\subset \Dmod(\Ran)$$
be the chiral algebra, corresponding to $\CA$ (see \secref{sss:chiral algebras}); it has a structure of Lie algebra,
viewed as an object of $\Dmod(\Ran)$ in the chiral symmetric monoidal structure.

\medskip

We can interpret $\CA\mod^{\on{fact}}_{Z}$ as the full subcategory of 
$$\CA^{\on{ch}}\mod(\QCoh(Z^{\subseteq})),$$
consisting of objects with set-theoretic support on $Z\overset{\on{diag}_Z}\to Z^{\subseteq}$ 
(note that the map $\on{diag}_Z$ identifies $Z$ with its formal completion inside $Z^{\subseteq}$).

\medskip

Now the assertion of the proposition follows from \lemref{l:act t} below.

\end{proof}

\sssec{}

Let $\bC$ be a DG category with a t-structure compatible with filtered colimits,
and suppose we are given an embedding $\bC \subset \bC'$. Let $\bA$ be
a symmetric monoidal (resp., monoidal) category acting on $\bC'$; denote
the action functor by $\star$. 

\medskip

Let $L$ (resp., $A$) be a Lie algebra (resp., associative) in $\bA$ with the property that
for every $\bc_1 \in \bC^{\leq 0}$ and every $\bc_2 \in \bC^{>0}$,
$\Maps(L \star \bc_1,\bc_2) = 0$ (resp., $\Maps(A \star \bc_1,\bc_2) = 0$).

\begin{lem} \label{l:act t}
Under the above circumstances, the category $L\mod(\bC) := L\mod(\bC') \underset{\bC'}\times \bC$
(resp., $A\mod(\bC) := A\mod(\bC') \underset{\bC'}\times \bC$) admits a unique
t-structure for which the forgetful functor to $\bC$ is t-exact.  
\end{lem}

\begin{rem}

In the terminology of \secref{sss:t-str on shf of cat}, the assertion of \propref{p:t-structure fact mod} is that
if $\CA$ is connective, the sheaf of categories $\CA\ul\mod^{\on{fact}}$ carries a t-structure,
for which the forgetful functor $\oblv_CA$ is conservative. 

\end{rem} 

\begin{rem}

An assertion parallel to \propref{p:t-structure fact mod} holds for
the category 
$$\CA\mod^{\on{fact}}_{Z_\dr},$$
i.e., for crystalline factorization $\CA$-modules at $Z$. In this case, we consider 
$\QCoh(Z_\dr)\simeq \Dmod(Z)$ equipped with the \emph{right} t-structure, i.e.,
one for which the functor $\ind^r$ is t-exact.

\end{rem}

\begin{rem} \label{r:alg in t}
Assume for a moment that $\CA$ is such that $\oblv^l(\CA_X)$ belongs to $\QCoh(X)^\heartsuit$
and is flat. Then \cite[Lemma 3.4.12]{BD2} implies that the same is true for $\oblv^l(\CA_{X^I})$
for all finite sets $I$. 

\medskip

This in turn implies that for any $S\to \Ran$ with $S\in \affSch$, the object $\CA_S\in \QCoh(S)$ 
also belongs to $\QCoh(S)^\heartsuit$ and is flat. 

\end{rem} 

\sssec{} \label{sss:compatible maps of modules}

Let now $\phi:\CA_1\to \CA_2$ be a homomorphism of factorization algebras. 
Let 
$$\CM_i\in \CA_i\mod^{\on{fact}}_\CZ.$$

There is a natural notion of map $\CM_1\to \CM_2$, compatible with factorization. 
Namely, this is a map
$$\phi_m:\CM_{1,\CZ^{\subseteq}}\to \CM_{2,\CZ^{\subseteq}}$$
in $\QCoh(\CZ^{\subseteq})$, which makes the following diagram commute,
$$
\CD
\CM_{1,\CZ^{\subseteq}}|_{(\Ran\subset \CZ^{\subseteq})} @>{\sim}>>
(\CA_{1,\Ran} \boxtimes \CM_{1,\CZ^{\subseteq}})|_{(\Ran\subset \CZ^{\subseteq})} \\
@V{\phi_m}VV @VV{\phi \boxtimes \phi_m}V \\
\CM_{2,\CZ^{\subseteq}}|_{(\Ran\subset \CZ^{\subseteq})} @>{\sim}>>
(\CA_{2,\Ran} \boxtimes \CM_{2,\CZ^{\subseteq}})|_{(\Ran\subset \CZ^{\subseteq})},
\endCD
$$ 
along with a homotopy-coherent system of higher compatibilities.

\medskip

Denote the corresponding space of maps by
\begin{equation} \label{e:compatible maps of modules}
\Maps_{\CA_1\to \CA_2}(\CM_1,\CM_2).
\end{equation} 

\sssec{} \label{sss:univ property restr}

For a fixed $\CM_2\in \CA_2\mod^{\on{fact}}_\CZ$, we can consider the functor 
$$(\CA_1\mod^{\on{fact}}_\CZ)^{\on{op}}\to \Spc,$$
given by 
\begin{equation} \label{e:restr modules}
\CM_1\mapsto  \Maps_{\CA_1\to\CA_2}(\CM_1,\CM_2).
\end{equation}

This functor sends colimits to limits, and hence is representable. We denote the representing object 
by 
$$\on{Res}_\phi(\CM_2)\in \CA_1\mod^{\on{fact}}_\CZ.$$

We have a tautologically defined object in $\Maps_{\CA_1\to\CA_2}(\on{Res}_\phi(\CM_2),\CM_2)$, i.e., a map
\begin{equation} \label{e:univ restr map}
\on{Res}_\phi(\CM_2)\to \CM_2,
\end{equation}
compatible with factorization. 

\medskip

We claim:

\begin{lem} \label{l:fact restr}
The map \eqref{e:univ restr map} induces an isomorphism 
$$\on{Res}_\phi(\CM_2)_\CZ\to \CM_{2,\CZ}.$$
\end{lem}

\begin{proof}

This also follows by interpreting factorization modules as chiral modules. Under this
equivalence, restriction corresponds to restriction of modules along a homomorphism
of chiral algebras.

\end{proof} 

\begin{rem}
Let us emphasize that \lemref{l:fact restr} says that the operation of restriction acts as identity 
on the underlying object $\QCoh(\CZ)$, and its content is that one can restrict the factorization 
action on it of $\CA_2$ to obtain a factorization action of $\CA_1$.
\end{rem}

\sssec{Example} \label{sss:fact restr limit}

One can describe the object $\on{Res}_\phi(\CM_2)_{\CZ^{\subseteq}}$ explicitly. Let us consider the example
of $\CZ=\on{pt}$, where $\CZ\to \Ran$ corresponds to a singleton $x\in X$, so that $\CZ^{\subseteq}=\Ran_x$. 

\medskip

We have a natural map 
$$X\to \Ran_x, \quad x'\mapsto x\cup x',$$
and let us describe the restriction of $\on{Res}_\phi(\CM_2)$ along this map. 

\medskip 

Let $j$ denote the embedding $X-x\hookrightarrow X$. Then
$$\on{Res}_\phi(\CM_2)_X\simeq \CM_{2,X}\underset{j_*\circ j^*(\CA_2\boxtimes \CM_{2,x})}\times j_*\circ j^*(\CA_1\boxtimes \CM_{2,x}),$$
where the map 
$$\CM_{2,X}\to j_*\circ j^*(\CA_2\boxtimes \CM_{2,x})$$ is
$$\CM_{2,X}\overset{\text{adjunction}}\longrightarrow j_*\circ j^*(\CM_{2,X}) \overset{\on{factorization}}\simeq
j_*\circ j^*(\CA_2\boxtimes \CM_{2,x}).$$
 
\medskip

Using a description of $\on{Res}_\phi(\CM_2)_{\CZ^{\subseteq}}$ along these lines, one can prove \lemref{l:fact restr}
without resorting to chiral algebras.

\sssec{Example} \label{sss:restr vacuum}

Let $\phi:\CA_1\to \CA_2$ be a map of factorization algebras. For $\CZ\to \Ran$, consider the resulting map
$$\on{pr}_{\on{big}}^*(\CA_{1,\Ran})\to \on{pr}_{\on{big}}^*(\CA_{2,\Ran})$$
in $\QCoh(\CZ^{\subseteq})$. 

\medskip

It is easy to see that this map is compatible with factorization, i.e., gives rise to a map 
$$\phi:\CA_1^{\on{fact}_\CZ}\to \CA_2^{\on{fact}_\CZ},$$
compatible with factorization.

\medskip

In particular, we obtain a map
$$\CA_1^{\on{fact}_\CZ}\to \Res_\phi(\CA_2^{\on{fact}_\CZ})$$
in $\CA_1\mod^{\on{fact}}_\CZ$.

\ssec{Commutative factorization algebras} \label{ss:fact com}

\sssec{}   \label{sss:fact com}

The category $\on{FactAlg}(X)$ has an evident (pointwise) symmetric monoidal structure.
Consider the category 
$$\on{ComAlg}(\on{FactAlg}(X)).$$

\medskip

We will refer to its objects as \emph{commutative factorization algebras}. 

\medskip

We have a tautological
forgetful functor 
\begin{equation} \label{e:forget fact com}
\on{ComAlg}(\on{FactAlg}(X))\to \on{ComAlg}(\Dmod(\Ran))
\end{equation}

\sssec{} \label{sss:com fact vs Dmod com} 

Restriction along $X\to \Ran$ defines a functor
\begin{equation} \label{e:from com fact to fact in Dmod}
\on{ComAlg}(\on{FactAlg}(X))\to \on{ComAlg}(\Dmod(\Ran))\to \on{ComAlg}(\Dmod(X)).
\end{equation}

This functor has a right inverse, to be denoted 
$$A\mapsto \on{Fact}(A)$$. 

\begin{rem}
In \secref{ss:com untl fact alg} we will 
upgrade the functor $\on{Fact}(-)$ to an equivalence 
of categories. 
\end{rem} 

\sssec{} \label{sss:expl fact com}

The object
$$\on{Fact}(A)_\Ran\in \Dmod(\Ran)$$
can described by an explicit colimit procedure:

\medskip

Let $\on{TwArr}(\on{fSet}^{\on{surj}})$ be the twisted arrows category of $\on{fSet}^{\on{surj}}$. I.e.,
its objects are 
$$I\overset{\phi}\to J$$
and morphisms $(I_1\overset{\phi_1}\to J_1)\to (I_2\overset{\phi_2}\to J_2)$ are diagrams
\begin{equation} \label{e:morphism tw arr}
\CD
I_1 @>{\phi_1}>> J_1 \\
@V{\psi_I}VV @AA{\psi_J}A \\
I_2 @>{\phi_2}>> J_2.
\endCD
\end{equation} 

We define a functor 
$$A_{\on{TwArr}}:\on{TwArr}(\on{fSet}^{\on{surj}})\to \Dmod(\Ran)$$
be sending $I\overset{\phi}\to J$ to the direct image along 
$$\Delta_{X^J,\Ran}:X^J_\dr\to \Ran$$ of
$$\underset{j\in J}\boxtimes\, A^{\otimes \phi^{-1}(j)}.$$

For a morphism \eqref{e:morphism tw arr}, the corresponding map
$$(\Delta_{X^{J_1},\Ran})_!(\underset{j_1\in J_1}\boxtimes\, A^{\otimes \phi_1^{-1}(j_1)})\to
(\Delta_{X^{J_2},\Ran})_!(\underset{j_2\in J_2}\boxtimes\, A^{\otimes \phi_2^{-1}(j_2)})$$
is
\begin{multline*}
(\Delta_{X^{J_1},\Ran})_!(\underset{j_1\in J_1}\boxtimes\, A^{\otimes \phi_1^{-1}(j_1)})\simeq
(\Delta_{X^{J_2},\Ran})_! \circ (\Delta_{\psi_J})_!(\underset{j_1\in J_1}\boxtimes\, A^{\otimes \phi_1^{-1}(j_1)})\simeq \\
\simeq 
(\Delta_{X^{J_2},\Ran})_! \circ (\Delta_{\psi_J})_! \circ (\Delta_{\psi_J})^! 
(\underset{j_2\in J_2}\boxtimes\, A^{\otimes (\phi_2\circ \psi_I)^{-1}(j_2)})\to
(\Delta_{X^{J_2},\Ran})_! (\underset{j_2\in J_2}\boxtimes\, A^{\otimes (\phi_2\circ \psi_I)^{-1}(j_2)})\simeq \\
\simeq 
(\Delta_{X^{J_2},\Ran})_! (\underset{j_2\in J_2}\boxtimes\, \underset{i_2\in \phi_2^{-1}(j_2)}\otimes\, A^{\otimes \psi_I^{-1}(i_2)})\to 
(\Delta_{X^{J_2},\Ran})_! (\underset{j_2\in J_2}\boxtimes\, \underset{i_2\in \phi_2^{-1}(j_2)}\otimes\, A) \simeq \\
\simeq (\Delta_{X^{J_2},\Ran})_! (\underset{j_2\in J_2}\boxtimes\, A^{\otimes (\phi_2)^{-1}(j_2)}),
\end{multline*}
where the arrow in the third line is the tensor product of the maps
$$A^{\otimes \psi_I^{-1}(i_2)}\to A,$$
given by the commutative algebra structure on $A$. 

\medskip 

Then
$$\on{Fact}(A)_\Ran \simeq \underset{\on{TwArr}(\on{fSet}^{\on{surj}})}{\on{colim}}\, A_{\on{TwArr}}.$$

\sssec{} \label{sss:com mod as shf of cat}

Let $\CA$ be a commutative factorization algebra on $X$. Let $\CZ$ be a prestack mapping to $\Ran$. Denote
$$\CA\mod^{\on{com}}_\CZ:=\CA_\CZ\mod(\QCoh(\CZ)),$$
where we regard $\CA_\CZ$ as an object of $\on{ComAlg}(\QCoh(\CZ))$. 

\medskip

We will view the assignment 
$$\CZ\rightsquigarrow \CA\mod^{\on{com}}_\CZ$$
as a sheaf of categories over $\Ran$, which we will denote by $\CA\ul\mod^{\on{com}}$. 

\sssec{}

Let $\CA$ be of the form $\on{Fact}(A)$, for $A\in \on{ComAlg}(\Dmod(X)$. 
We claim that in this case there is a canonically defined functor of sheaves of categories,
$$\CA\ul\mod^{\on{com}}\to \CA\ul\mod^{\on{fact}},$$
i.e., a compatible collection of functors
\begin{equation} \label{e:com to fact}
\CA\mod^{\on{com}}_\CZ\to \CA\mod^{\on{fact}}_\CZ.
\end{equation}

\medskip

In order to construct \eqref{e:com to fact}, it suffices to exhibit an object of $\CA\mod^{\on{fact}}_\CZ$, equipped with an action of
$\CA_\CZ\in \on{ComAlg}(\QCoh(\CZ))$, where we view $\CA\mod^{\on{fact}}_\CZ$ as tensored over 
$\QCoh(\CZ)$. 

\medskip

The object in question is the vacuum module $\CA^{\on{fact}_\CZ}$, see \secref{sss:vacuum module}. Thus,
we need to define an action of $\CA_\CZ$ on $\CA^{\on{fact}_\CZ}$. To do so, we can consider the universal
case, i.e., $\CZ=\Ran$. 

\medskip

Thus, we need to define an action of 
$\on{pr}_{\on{small}}^*(\CA_\Ran)$ on $\on{pr}_{\on{big}}^*(\CA_\Ran)$,
compatible with the factorization. This map comes from the commutative algebra structure on $\CA_\Ran$,
and a homomorphism
$$\on{pr}_{\on{small}}^*(\CA_\Ran)\to \on{pr}_{\on{big}}^*(\CA_\Ran),$$
given by the \emph{unital structure} on $\CA$, see \secref{sss:com fact vs Dmod com untl} below. 

\medskip

Note that since 
$$\CA^{\on{fact}_\CZ}_\CZ\simeq \CA_\CZ,$$
the functor \eqref{e:com to fact} is compatible with the forgetful functors to $\QCoh(\CZ)$.

\sssec{} \label{sss:fact mods for com}

Let $\CA$ be a commutative factorization algebra. Restriction along the binary operation
$$\CA\sotimes \CA\to \CA$$
in the sense of \secref{sss:compatible maps of modules} defines on $\CA\mod^{\on{fact}}_\CZ$ a 
structure of (symmetric) \emph{pseudo-monoidal} category
(see \cite[Sect. 1.1]{BD2}) for what this means. 

\medskip

Concretely, for $\CM_1,\CM_2,\CN\in \CA\mod^{\on{fact}}_\CZ$
this means that we know what it means to map 
$$\CM_1``\otimes"\CM_2\to \CN.$$

Namely, by definition, the space of such maps is 
$$\Maps_{\CA\sotimes \CA\to \CA}(\CM_1\otimes \CM_2,\CN)$$
in the notation of \eqref{e:compatible maps of modules}, 
where we regard $\CM_1\otimes \CM_2$ as an object of $(\CA\sotimes \CA)\mod^{\on{fact}}_\CZ$. 

\medskip

In particular, it makes sense to talk about
(commutative) algebra objects in $\CA\mod^{\on{fact}}_\CZ$.

\medskip

Note that by construction, the functor \eqref{e:com to fact} is right-lax (pseudo)-monoidal.
In particular, it maps (commutative) algebras to (commutative) algebras. In addition, the
forgetful functor
$$\CA\mod^{\on{fact}}_\CZ\to \QCoh(\CZ^{\subseteq})$$
is also right-lax (pseudo)-monoidal; in particular, it maps (commutative) algebras to (commutative) algebras.

\ssec{Factorization categories}

\sssec{} \label{sss:fact cat}

A factorization category $\bA$ on $X$ is a sheaf/crystal of categories $\ul\bA$ over $\Ran$, equipped with a datum 
of factorization
\begin{equation} \label{e:fact category}
\on{union}^*(\ul\bA)|_{(\Ran\times \Ran)_{\on{disj}}}\simeq \ul\bA\boxtimes \ul\bA|_{(\Ran\times \Ran)_{\on{disj}}}
\end{equation} 
(here $-|_-$ denotes pullback of sheaves of categories along an open embedding) and with a homotopy-coherent
datum of associativity and commutativity; see \cite[Sect. 6]{Ra6}, where the definition is written out in detail. 

\medskip

For $f:\CZ\to \Ran$, we will denote
$$\bA_\CZ:=\Gamma(\CZ,f^*(\ul\bA)).$$

In particular, we denote
$$\bA_\Ran:=\Gamma(\Ran,\ul\bA);$$
this is a category tensored over $\Dmod(\Ran)$. 

\medskip

Let $\on{FactCat}(X)$ denote the (2-)category of factorization categories over $X$. This category
carries a naturally defined symmetric monoidal structure. 

\sssec{} \label{sss:unit fact cat}

The unit object in this symmetric monoidal category is the sheaf of categories $\ul\QCoh(\Ran)$. 

\medskip

By a slight abuse of notation, we will denote this factorization category by $\Vect$
(i.e., its fiber at any $\ul{x}\in \Ran$ is $\Vect\in \DGCat$). 

\sssec{}  \label{sss:fact glob sect}

Note that since $\Ran$ is 1-affine (by \lemref{l:Ran is 1-affine}), the datum of the sheaf of categories
$\ul\bA$ is equivalent to that of the category $\bA_\Ran$, equipped with an action of $\Dmod(\Ran)$. 

\medskip

Furthermore, the datum of factorization is equivalent to
\begin{multline} \label{e:fact category Sect}
\bA_\Ran\underset{\Dmod(\Ran),\on{union}^!}\otimes \Dmod((\Ran\times \Ran)_{\on{disj}})\simeq \\
\simeq (\bA_\Ran\otimes \bA_\Ran)\underset{\Dmod(\Ran)\otimes \Dmod(\Ran)}\otimes \Dmod((\Ran\times \Ran)_{\on{disj}}),
\end{multline}
equipped with a homotopy-coherent datum of associativity and commutativity.

\sssec{} \label{sss:fact alg in fact cat}

Let $\bA$ be a factorization category. A factorization algebra $\CA\in \bA$ is an object $\CA_\Ran\in \bA_\Ran$,
equipped with a datum of factorization
$$\on{union}^!(\CA_\Ran)|_{(\Ran\times \Ran)_{\on{disj}}}\simeq (\CA_\Ran\boxtimes \CA_\Ran)|_{(\Ran\times \Ran)_{\on{disj}}},$$
as objects in the two sides of \eqref{e:fact category}, and with a homotopy-coherent datum of associativity and commutativity.

\medskip

Let $\on{FactAlg}(X,\bA)$ denote the category of factorization algebras in $\bA$. 

\medskip

A factorization functor $\bA_1\to \bA_2$ induces a functor
$$\on{FactAlg}(X,\bA_1)\to \on{FactAlg}(X,\bA_2).$$

\sssec{}

Given $\CA\in \on{FactAlg}(X,\bA)$, parallel to \secref{sss:fact modules}, given $\CZ\to \Ran$, one defines 
the category of factorization $\CA$-modules \emph{in} $\bA$ \emph{at} $\CZ$, denoted
$$\CA\mod^{\on{fact}}(\bA)_\CZ.$$

\medskip

The assignment
\begin{equation} \label{e:fact A mod in self}
\CZ\to \CA\mod^{\on{fact}}(\bA)_\CZ
\end{equation} 
forms a sheaf of categories over $\Ran$ that we denote by
$$\CA\ul\mod^{\on{fact}}(\bA).$$

\sssec{} \label{sss:com fact cat}

We will see many examples of factorization categories in the sequel. However, one family of examples we can 
produce right away:

\medskip

Let $\ul\bA_X$ be a crystal of symmetric monoidal categories over $X$. To it we attach a (symmetric monoidal) factorization category,
denoted $\on{Fact}(\ul\bA_X)$, see \cite[Sect. 8.1]{GLys}.

\medskip 

Equivalently, let $\bA_X$ be a symmetric monoidal category tensored over $\Dmod(X)$. To it we can attach a DG category,
denoted $\on{Fact}(\bA_X)$, tensored over $\Dmod(\Ran)$ (in fact, a commutative algebra object in $\Dmod(\Ran)\mmod$),
and equipped with a factorization structure as in \secref{sss:fact glob sect}.

\medskip

This can be done by directly mimicking the procedure in \secref{sss:expl fact com}. 

\medskip

Similarly, if $A\in \bA_X$ is a commutative algebra object, it gives rise to a commutative algebra object
$$\on{Fact}(A)\in \on{FactAlg}(X,\on{Fact}(\ul\bA_X)).$$

\sssec{} \label{sss:com mod as fact cat}

Let $A$ be an object in $\on{ComAlg}(\Dmod(X))$. Denote $\CA:=\on{Fact}(A)\in \on{FactAlg}(X)$. It is easy to
see that the sheaf of categories $\CA\ul\mod^{\on{com}}$ (see \secref{sss:com mod as shf of cat}) has a natural factorization structure.
We denote the resulting factorization category by $\CA\mod^{\on{com}}$.

\medskip

Furthermore, we have:
$$\CA\mod^{\on{com}}\simeq \on{Fact}(A\mod(\Dmod(X))).$$

\sssec{Notational convention} \label{sss:const fact cat} 

In the special case when $\bA_X$ is constant, i.e., of the form $\bA\otimes \Dmod(X)$, where $\bA$ is a symmetric
monoidal category, we will denote the factorization category $\on{Fact}(\bA_X)$ simply by $\bA$. 

\medskip

For example, when $\bA=\Rep(\cG)$, we will use the symbol $\Rep(\cG)$ to denote the corresponding 
factorization category. 

\medskip

Note that this is in line with the notation for $\Vect$ in \secref{sss:unit fact cat}.

\sssec{}  \label{sss:comp gen fact}

Let $\bA$ be a factorization category over $X$. We shall say that $\bA$ is dualizable
if for every $S\in \affSch_{/\Ran}$, the category $\bA_S$ is dualizable. This is equivalent to
$\ul\bA$ being a dualizable object in $\on{ShvCat}(\Ran)$. 

\medskip

If $\bA$ is dualizable, then the dual $\ul\bA^\vee$ of $\ul\bA$ as a sheaf of categories 
admits a natural factorization structure; we will denote the resulting factorization category by
$\bA^\vee$.

\medskip

The datum of duality between two factorization categories $\bA_1$ and $\bA_2$ is equivalent to
that of factorization functors
$$\Vect\to \bA_1\otimes \bA_2 \text{ and } \bA_1\otimes \bA_2\to \Vect,$$
satisfying the usual axioms. 

\sssec{} \label{sss:pres compac fact}

Let $\bA$ be a factorization category over $X$. We shall say that $\bA$ is compactly generated
if $\ul\bA$ is compactly generated as a sheaf of categories over $\Ran$. 

\medskip

Let $\Phi:\bA_1\to \bA_2$ be a factorization functor between factorization categories. Suppose that $\bA_1$
is compactly generated. Suppose also that $\Phi$ preserves compactness, so that $\Phi$ admits a right
adjoint $\Phi^R$, as a functor between sheaves of categories. 

\medskip

It follows automatically that $\Phi^R$ carries a structure of compatibility with factorization. 

\sssec{}  \label{sss:t-str fact}

Let $\bA$ be a factorization category over $X$. A t-structure on $\bA$ is a t-structure on $\ul\bA$
as a sheaf of categories that is compatible with factorization in the following sense:

\medskip 

Let us be given a map 
$$(\ul{x}_1,\ul{x}_2):S\to (\Ran\times \Ran)_{\on{disj}},\quad S\in \affSch.$$ 
We need that the factorization equivalence
$$\bA_{S,\on{union}\circ (\ul{x}_1,\ul{x}_2)} \simeq \bA_{S,\ul{x_1}}\underset{\QCoh(S)}\otimes \bA_{S,\ul{x_2}}$$
be t-exact, where the right-hand side is equipped with the tensor product t-structure. 

\sssec{} \label{sss:lax fact}

A lax factorization category $\bA$ is a sheaf of categories $\bA_\Ran$ over $\Ran$, equipped with \emph{functors} 
\begin{equation} \label{e:lax fact category}
\ul\bA\boxtimes \ul\bA|_{(\Ran\times \Ran)_{\on{disj}}}\to \on{union}^*(\ul\bA)|_{(\Ran\times \Ran)_{\on{disj}}},
\end{equation}
equipped with a homotopy-coherent datum of associativity and commutativity.

\medskip

The entire discussion above equally applies to lax factorization categories. 

\sssec{}

In particular, we can talk about 
factorization algebras in a lax factorization category:

\medskip

We require that the corresponding functor 
\begin{multline} \label{e:lax fact category Sect}
(\bA_\Ran\otimes \bA_\Ran)\underset{\Dmod(\Ran)\otimes \Dmod(\Ran)}\otimes \Dmod((\Ran\times \Ran)_{\on{disj}})\to \\
\to \bA_\Ran\underset{\Dmod(\Ran),\on{union}^!}\otimes \Dmod((\Ran\times \Ran)_{\on{disj}}),
\end{multline}
maps the object 
$$(\CA_\Ran\boxtimes \CA_\Ran)|_{(\Ran\times \Ran)_{\on{disj}}}\in 
(\bA_\Ran\otimes \bA_\Ran)\underset{\Dmod(\Ran)\otimes \Dmod(\Ran)}\otimes \Dmod((\Ran\times \Ran)_{\on{disj}})$$
to the object
$$\on{union}^!(\CA_\Ran)|_{(\Ran\times \Ran)_{\on{disj}}}\in 
\bA_\Ran\underset{\Dmod(\Ran),\on{union}^!}\otimes \Dmod((\Ran\times \Ran)_{\on{disj}}).$$

\sssec{Example} \label{sss:fact alg lax fact}

Let $\CA$ be a factorization algebra on $X$. Note that the sheaf of categories $\CA\ul\mod^{\on{fact}}$
carries a natural lax factorization structure:

\medskip

Let us be given a pair of maps $\ul{x}_i:S_i\to \Ran$, $i=1,2$ such
that 
$$S_1\times S_2\overset{\ul{x}_1\times \ul{x}_2}\longrightarrow \Ran\times \Ran$$
lands in $(\Ran\times \Ran)_{\on{disj}}$. Let us be given a pair of objects
$\CM_i\in \CA\mod^{\on{fact}}_{S_i}$. We define the corresponding object
$$\CM_1\boxtimes \CM_2\in \CA\mod^{\on{fact}}_{S_1\times S_2}$$
as follows.

\medskip

Consider the prestack
$$(S_1\times S_2)^{\subseteq}=\{\ul{x'}\in \Ran,\, |\, \ul{x}_1\subseteq \ul{x}',\ul{x}_2\subseteq \ul{x}'\}.$$

For $i=1,2$, let $U_i\subset (S_1\times S_2)^{\subseteq}$ be the open sub-prestack corresponding to the condition that
$(\ul{x}_i,\ul{x}')$ belongs to $(\Ran\times \Ran)_{\on{disj}}$. The condition that $\ul{x}_1$ and $\ul{x}_2$ are disjoint
implies that 
$$U_1\cup U_2=(S_1\times S_2)^{\subseteq}.$$
Denote $U_{1,2}:=U_1\cap U_2$. 

\medskip

Let $\on{pr}_1$ (resp., $\on{pr}_2$)
denote the map $U_1\to S_1\times  S_2^{\subseteq}$ (resp., $U_2\to S_1^{\subseteq}\times S_2$) whose second (resp., first) 
component remembers the data of $\ul{x}$. Let $\on{pr}_{1,2}$ denote the map $U_{1,2}\to S_1\times S_2\times \Ran$,
whose last component remembers $\ul{x}'$. 

\medskip

We let
$$(\CM_1\boxtimes \CM_2)_{(S_1\times S_2)^{\subseteq}}|_{U_1} \simeq \on{pr}_1^*(\CM_{1,S_1}\boxtimes \CM_{2,S_2^{\subseteq}})$$
and 
$$(\CM_1\boxtimes \CM_2)_{(S_1\times S_2)^{\subseteq}}|_{U_2} \simeq \on{pr}_2^*(\CM_{1,S^{\subseteq}_1}\boxtimes \CM_{2,S_2}).$$

The factorization structures on $\CM_1$ and $\CM_2$ imply that
$$(\CM_1\boxtimes \CM_2)_{(S_1\times S_2)^{\subseteq}}|_{U_1}|_{U_{1,2}}\simeq
\on{pr}_{1,2}^*(\CM_{1,S_1}\boxtimes \CM_{2,S_2}\boxtimes \CA_\Ran) \simeq 
(\CM_1\boxtimes \CM_2)_{(S_1\times S_2)^{\subseteq}}|_{U_2}|_{U_{1,2}}.$$

\medskip 

This defines the object
$$(\CM_1\boxtimes \CM_2)_{(S_1\times S_2)^{\subseteq}}\in \QCoh((S_1\times S_2)^{\subseteq}).$$

The factorization structure on it against $\CA$ follows from the construction. 

\sssec{}  \label{sss:A mod is lax}

Let us denote the resulting lax factorization category by $\CA\mod^{\on{fact}}$. For a map $\phi:\CA_1\to \CA_2$
between the factorization algebras, the restriction functor
$$\Res_\phi: \CA_2\ul\mod^{\on{fact}}\to \CA_1\ul\mod^{\on{fact}}$$
upgrades to a factorization functor, denoted 
$$\Res_\phi: \CA_2\mod^{\on{fact}}\to \CA_1\mod^{\on{fact}}.$$

\medskip

Denote by $\oblv_\CA$ the forgetful functor $\CA\mod^{\on{fact}}\to \Vect$ (see \secref{sss:const fact cat} above 
for the convention regarding $\Vect$).  

\medskip

Note that the assignment 
$$\CZ\mapsto \CA^{\on{fact}_\CZ}$$
(see \secref{sss:vacuum module}) gives rise to a factorization algebra object in $\CA\mod^{\on{fact}}$, to be denoted $\CA^{\on{fact}}$,
so that $$\oblv_\CA(\CA^{\on{fact}})=\CA.$$ That said, we will sometimes abuse the notation
and instead of $\CA^{\on{fact}}$ simply write $\CA$. 

\medskip

Suppose for moment that $\CA$ is connective. Then the construction in \secref{sss:t-structure fact mod} equips $\CA\mod^{\on{fact}}$ with a t-structure
in the sense of \secref{sss:t-str fact}.\footnote{In the case of lax factorization categories, we require that the functor \eqref{e:lax fact category} be t-exact.}

\sssec{} \label{sss:A mod in A}

The example in \secref{sss:fact alg lax fact} generalizes to the situation when $\CA$ is a factorization algebra
in a given lax factorization category $\bA$. We obtain that the sheaf of categories 
$$\CA\ul\mod^{\on{fact}}(\bA)$$
carries a structure of lax factorization category. We denote it by $\CA\mod^{\on{fact}}(\bA)$.

\medskip

As in \secref{sss:A mod is lax}, we can consider $\CA^{\on{fact}}$
as an object of $\on{FactAlg}(X,\CA\ul\mod^{\on{fact}}(\bA))$. By a slight abuse of notation, 
we will sometimes denote this factorization algebra simply by $\CA$. 

\medskip

By definition, for any $\CZ\to \Ran$
$$(\CA^{\on{fact}})_\CZ\simeq \CA^{\on{fact}_\CZ}$$
as objects of $\CA\mod^{\on{fact}}(\bA)_\CZ$.

\medskip

If $\bA$ carries a t-structure (in the sense of \secref{sss:t-str fact})
and $\CA$ is connective (i.e., $\oblv^l(\CA_X)\in \bA_X$ is connective), 
then the construction in \secref{sss:t-structure fact mod} equips $\CA\mod^{\on{fact}}(\bA)$ with a t-structure. 

\sssec{} \label{sss:A and A' mod}

Let $\bA$ be a lax factorization category, and let $\phi:\CA\to \CA'$ be a homomorphism between factorization algebras
in $\bA$. Restriction along $\phi$ denotes a factorization functor
$$\Res_\phi:\CA'\mod^{\on{fact}}(\bA)\to \CA\mod^{\on{fact}}(\bA).$$

In particular, we can consider 
$$\Res_\phi(\CA'{}^{\on{fact}})\in \on{FactAlg}(X,\CA\mod^{\on{fact}}(\bA)).$$

We have
$$\oblv_\CA(\Res_\phi(\CA'{}^{\on{fact}}))=\CA'.$$

\medskip

We will sometimes abuse the notation, and instead of $\Res_\phi(\CA'{}^{\on{fact}})$ simply write $\CA'.$

\ssec{Factorization module categories} \label{ss:fact mod cat}

\sssec{} \label{sss:fact mod cat}

Let $\bA$ be a factorization category over $X$. Let $\CZ$ be a prestack mapping to $\Ran$.
A factorization module category $\bC$ \emph{over} $\bA$ \emph{at} $\CZ$ is a sheaf of categories $\ul\bC$
on $\CZ^{\subseteq}$, equipped with a factorization structure:
\begin{equation} \label{e:fact mod category}
\on{union}^*(\ul\bC)|_{(\Ran\times \CZ^{\subseteq})_{\on{disj}}}\simeq 
\ul\bA\boxtimes \ul\bC|_{(\Ran\times \CZ^{\subseteq})_{\on{disj}}}
\end{equation} 
and with a homotopy-coherent datum of associativity; see \cite[Sect. 6]{Ra6} for details. 

\medskip

Let $\bA\mmod^{\on{fact}}_\CZ$ denote the (2-)category of 
factorization module categories over $\bA$ at $\CZ$.  

\sssec{}

For an object $\bC\in \bA\mmod^{\on{fact}}_\CZ$ and $f:\CZ'\to \CZ^{\subseteq}$, we denote
$$\bC_{\CZ'}:=\Gamma(\CZ',f^*(\ul\bC)).$$

Taking $\CZ'=\CZ$ and $f=\on{diag}_\CZ$, we obtain the category
$\bC_\CZ$, tensored over $\QCoh(\CZ)$. We will refer to $\bC_\CZ$ as \emph{the category underlying} $\bC$. 

\sssec{}  \label{sss:propagate modules cat}

As in \secref{sss:propagate modules}, the above category $\bC_{\CZ'}$ is in fact the category underlying
an $\bA$-module category \emph{at} $\CZ'$, to be denoted $\bC|_{\CZ'}$. 

\sssec{Example} \label{sss:vac fact mod cat}

Repeating the construction in \secref{sss:vacuum module}, we obtain that for any $\CZ$ there exists a distinguished object
$$\bA^{\on{fact}_\CZ}\in \bA\mmod^{\on{fact}}_\CZ,$$
whose underlying category is $\bA_\CZ$.

\medskip

We will refer to $\bA^{\on{fact}_\CZ}$ as the \emph{vacuum} factorization module category at $\CZ$.  

\sssec{Example} \label{sss:cats to fact mods over Vect}

Take $\bA=\Vect$, and let $\ul\bC_0$ be a sheaf of categories over $\CZ$. We claim that it gives rise to
a factorization module category over $\Vect$ at $\CZ$, to be denoted $\bC$ (cf. \secref{sss:fact mod over unit nu}).

\medskip

Namely, the corresponding sheaf of categories $\ul\bC$ over $\CZ^{\subseteq}$ is the pullback of 
$\ul\bC_0$ along the projection $\on{pr}_{\on{small}}:\CZ^\subseteq\to \CZ$. 

\sssec{} \label{sss:modules inside modules}

Let $\CA$ be a factorization algebra in $\bA$, and let $\bC$ be an object of $\bA\mmod^{\on{fact}}_\CZ$.

\medskip

A factorization $\CA$-module $\CM$ \emph{in} $\bC$ is an object
$$\CM_{\CZ^{\subseteq}}\in \bC_{\CZ^{\subseteq}},$$
equipped with an isomorphism
\begin{equation} \label{e:fact mod in fact mod category}
\on{union}^*(\CM_{\CZ^{\subseteq}})|_{(\Ran\times \CZ^{\subseteq})_{\on{disj}}}\simeq 
(\CA_\Ran \boxtimes \CM_{\CZ^{\subseteq}})|_{(\Ran\times \CZ^{\subseteq})_{\on{disj}}},
\end{equation}
where:

\begin{itemize}

\item $\on{union}^*(\CM_{\CZ^{\subseteq}})$ is an object in
$$\Gamma\left((\Ran\times \CZ^{\subseteq})_{\on{disj}},\on{union}^*(\ul\bC)|_{(\Ran\times \CZ^{\subseteq})_{\on{disj}}}\right);$$

\item $\CA_\Ran \boxtimes \CM_{\CZ^{\subseteq}}$ is an object in
$$\Gamma\left((\Ran\times \CZ^{\subseteq})_{\on{disj}},\ul\bA\boxtimes \ul\bC|_{(\Ran\times \CZ^{\subseteq})_{\on{disj}}}\right);$$

\item The isomorphisms between the two sides is understood in the sense of the identification \eqref{e:fact mod category}. 

\end{itemize}

\medskip

The isomorphism \eqref{e:fact mod in fact mod category} is required to be equipped with a homotopy-coherent
datum of associativity. 

\sssec{}

We denote the category of factorization $\CA$-modules in $\bC$ by
$$\CA\mod^{\on{fact}}(\bC)_\CZ.$$

Note that when $\bC:=\bA^{\on{fact}_\CZ}$, recover the category \eqref{e:fact A mod in self}. 

\sssec{}

Let $\ul\bC_0$ be as in \secref{sss:cats to fact mods over Vect}. Denote
$$\bC_0:=\Gamma(\CZ,\ul\bC_0).$$

Let $\CA$ be a factorization algebra (i.e., a factorization algebra in $\Vect$, viewed as a factorization category). 
Then it makes sense to consider the category
$$\CA\mod^{\on{fact}}(\bC)_\CZ.$$

Note that we have a tautologically defined functor
\begin{equation} \label{e:tensor up fact mod}
\CA\mod^{\on{fact}}_\CZ\underset{\QCoh(\CZ)}\otimes \bC_0\to \CA\mod^{\on{fact}}(\bC)_\CZ.
\end{equation} 

We claim:

\begin{lem} \label{l:tensor up fact mod}
Assume that $\CZ=S$ is an affine scheme, and 
assume that $\ul\bC_0$ dualizable as a sheaf of categories. Then the
functor \eqref{e:tensor up fact mod} is an equivalence. 
\end{lem}

\begin{proof}

Let $\ul\bC_0^\vee$ be the dual sheaf of categories; note that
$\Gamma(\CZ,\ul\bC_0^\vee)$ identifies with the dual of $\bC_0$ as a $\QCoh(S)$-linear
category, and hence also as a plain DG category. 

\medskip

We have a naturally defined functor
$$\CA\mod^{\on{fact}}(\bC)_S\underset{\QCoh(S)}\otimes \bC_0^\vee\to
\CA\mod^{\on{fact}}(\bC\otimes \bC^\vee)_S,$$
where $\bC^\vee\in \Vect\mmod^{\on{fact}}_\CZ$ is attached to $\ul\bC_0^\vee$
by the procedure of \secref{sss:cats to fact mods over Vect}. 

\medskip

Composing with the evaluation map
$$\CA\mod^{\on{fact}}(\bC\otimes \bC^\vee)_S\to \CA\mod^{\on{fact}}(\Vect^{\on{fact}}_S)_S\simeq
\CA\mod^{\on{fact}}_S,$$
we obtain a functor
\begin{equation} \label{e:tensor up fact mod 1}
\CA\mod^{\on{fact}}(\bC)_S\underset{\QCoh(S)}\otimes \bC_0^\vee\to 
\CA\mod^{\on{fact}}_S.
\end{equation} 

\medskip

The sought-for inverse functor to \eqref{e:tensor up fact mod} is given by
$$\CA\mod^{\on{fact}}(\bC)_\CZ\overset{\on{Id}\otimes \on{co-eval}}\longrightarrow
\CA\mod^{\on{fact}}(\bC)_\CZ\underset{\QCoh(S)}\otimes \ul\bC_0^\vee\underset{\QCoh(S)}\otimes \ul\bC_0 
\overset{\text{\eqref{e:tensor up fact mod 1}}}\longrightarrow 
\CA\mod^{\on{fact}}_S \underset{\QCoh(S)}\otimes \ul\bC_0.$$

\end{proof} 

\sssec{} \label{sss:factorization restriction}

Let $\Phi:\bA_1\to \bA_2$ be a factorization functor between factorization categories. Let 
$\bC_1$ and $\bC_2$ be objects in $\bA_1\mmod^{\on{fact}}_\CZ$ and $\bA_2\mmod^{\on{fact}}_\CZ$,
respectively. 

\medskip

A functor $\Phi_m:\bC_1\to \bC_2$ compatible with factorization is a functor 
$$\Phi_m:\ul\bC_1\to \ul\bC_2$$
between sheaves of categories on $\CZ^{\subseteq}$ that makes the following diagram commute
$$
\CD
\on{union}^*(\ul\bC_1)|_{(\Ran\times \CZ^{\subseteq})_{\on{disj}}} @>{\sim}>>
\ul\bA_1\boxtimes \ul\bC_1|_{(\Ran\times \CZ^{\subseteq})_{\on{disj}}} \\
@V{\Phi_m}VV @VV{\Phi\boxtimes \Phi_m}V \\
\on{union}^*(\ul\bC_2)|_{(\Ran\times \CZ^{\subseteq})_{\on{disj}}} @>{\sim}>>
\ul\bA_2\boxtimes \ul\bC_2|_{(\Ran\times \CZ^{\subseteq})_{\on{disj}}}
\endCD
$$ 
along with a homotopy-coherent system of higher compatibilities.

\medskip

Let 
$$\on{Funct}_{\bA_1\to \bA_2}(\bC_1,\bC_2)$$
denote the category of such functors. When an ambiguity is likely to occur, we will use the notation 
$\on{Funct}_{\Phi:\bA_1\to \bA_2}(\bC_1,\bC_2)$, i.e., we will insert $\Phi$ in the subscript. 

\sssec{} \label{sss:univ property restr cat}

Given $\bC_2\in \bA_2\mmod^{\on{fact}}_\CZ$, one defines its restriction, 
$$\on{Res}_\Phi(\bC_2)\in \bA_1\mmod^{\on{fact}}_\CZ$$
by the universal property as in \secref{sss:univ property restr}, i.e.,
$$\on{Funct}_{\bA_1\mmod^{\on{fact}}_\CZ}(\bC_1,\on{Res}_\Phi(\bC_2))\simeq
\on{Funct}_{\bA_1\to \bA_2}(\bC_1,\bC_2).$$

\medskip

One can explicitly describe $\on{Res}_\Phi(\bC_2)$ by a limit procedure as in \secref{sss:fact restr limit}.
\footnote{This will be written out in detail in \cite{CFGY}.}

\medskip

We have a tautologically defined functor
\begin{equation} \label{e:restr fact cat}
\on{Res}_\Phi(\bC_2)\to \bC_2
\end{equation}
compatible with factorization. 

\medskip

Parallel to \lemref{l:fact restr}, we have:

\begin{lem} \label{l:fact restr cat}
The functor \eqref{e:restr fact cat} induces an equivalence of the underlying categories
$$\on{Res}_\Phi(\bC_2)_\CZ\to \bC_{2,\CZ}.$$
\end{lem} 

\sssec{}

Let $\Phi:\bA_1\to \bA_2$ be a factorization functor, and let $\CA_1$ be a factorization algebra in $\bA_1$. 
Note that $\Phi(\CA_1)$ has a natural structure of factorization algebra in $\bA_2$.

\medskip

Let $\bC_1$ and $\bC_2$ as in \secref{sss:factorization restriction}. 
Given a functor $\Phi_m:\bC_1\to \bC_2$ compatible with factorization, we obtain a naturally defined functor
\begin{equation} \label{e:pushforward fact modules}
\CA_1\mod^{\on{fact}}(\bC_1)_\CZ\to \Phi(\CA_1)\mod^{\on{fact}}(\bC_2)_\CZ.
\end{equation}

Parallel with \lemref{l:fact restr cat} one proves: 

\begin{lem} \label{l:modules for fact alg restr}
Assume that the functor $\Phi_m:\bC_1\to \bC_2$ induces an equivalence
$$\bC_1\to \on{Res}_\Phi(\bC_2)$$
as factorization module categories over $\bA_1$. Then the functor \eqref{e:pushforward fact modules}
is an equivalence.
\end{lem} 

\sssec{} \label{sss:restr pairs} 

As in \secref{sss:restr vacuum}, given a factorization functor
$\Phi:\bA_1\to \bA_2$, for any $\CZ\to \Ran$, we have a canonically defined functor
$$\bA_1^{\on{fact}_\CZ}\to \on{Res}_\Phi(\bA^{\on{fact}_\CZ}_2).$$

\medskip 

In particular, given a commutative diagram of factorization categories
$$
\CD
\bA'_1 @>>> \bA'_2 \\
@A{\Psi_1}AA @AA{\Psi_2}A \\
\bA_1 @>>> \bA_2.
\endCD
$$
we obtain a functor
$$\on{Res}_{\Psi_1}(\bA'_1{}^{\on{fact}_\CZ})\to \on{Res}_{\Psi_2}(\bA'_2{}^{\on{fact}_\CZ}),$$
compatible with factorization. 

\ssec{Factorization categories of algebro-geometric nature}

\sssec{} 

Let $\CZ_\Ran\to \Ran$ be a prestack. We attach to it a sheaf of categories $\ul\QCoh(\CZ)$ over $\Ran$, namely,
$$\ul\QCoh(\CZ):=\pi_*(\ul\QCoh(\CZ_\Ran)),$$
where $\pi$ denotes the projection $\CZ_\Ran\to \Ran$.

\medskip

Explicitly, for $S\in \affSch_{/\Ran}$, we have
$$\QCoh(\CZ)_S=\QCoh(S\underset{\Ran}\times \CZ_\Ran).$$

\sssec{} \label{sss:QCoh fact}

Let $\CT$ be a factorization space over $X$. Consider the corresponding sheaf of categories $\ul\QCoh(\CT)$, i.e.,
$$\QCoh(\CT)_S=\QCoh(\CT_S).$$
 
\medskip

The factorization structure on $\CT$ equips $\ul\QCoh(\CT)$ with a lax factorization structure. We denote the
corresponding lax factorization category by $\QCoh(\CT)$.

\begin{rem}

The reason that $\CT$ is a priori only lax is that for a pair of prestacks $\CY_1$ and $\CY_2$, the naturally defined functor
$$\QCoh(\CY_1)\otimes \QCoh(\CY_2)\to \QCoh(\CY_1\otimes \CY_2)$$
is not necessarily an equivalence. 

\medskip

This lax structure is strict, e.g., if $\ul\QCoh(\CT)$ is dualizable (as a sheaf of categories over $\Ran$, which is equivalent to
each $\QCoh(\CT)_S$ being dualizable). This happens, e.g., if $\CT$ is a factorization \emph{affine scheme}.

\end{rem}

\sssec{Example}

Let $\CY$ be an affine D-scheme over $X$. On the one hand, we can consider the factorization scheme $\fL^+_\nabla(\CY)$
and the corresponding factorization category $\QCoh(\fL^+_\nabla(\CY))$. As such, it is equipped with a symmetric monoidal
structure.

\medskip

On the other hand, we can consider the symmetric monoidal factorization category 
$$\on{Fact}(\QCoh(\CY)).$$

We claim that there is a canonical equivalence:
$$\on{Fact}(\QCoh(\CY))\simeq \QCoh(\fL^+_\nabla(\CY)).$$

Indeed, both sides identify, as factorization categories with $\CA\mod^{\on{com}}$ (see \secref{sss:com mod as fact cat}),
where 
$$\CA=\on{Fact}(A), \quad \CY=\Spec_X(A),$$
see \eqref{e:arcs as fact} below.

\sssec{} \label{sss:QCOh LS loc}

Let now $\CT:=\LS_H^\reg$. We claim that the lax factorization category $\QCoh(\LS_H^\reg)$ identifies with
$\Rep(H)$ (see \secref{sss:const fact cat} for the notational conventions); in particular $\QCoh(\LS_H^\reg)$
is a \emph{factorization} category. 

\medskip

Indeed, on the one hand, unwinding the definitions, we obtain that $\Rep(H)$, viewed as a factorization category
is the (factorization) category of comodules with respect to $\on{Fact}(\CO_H)$, viewed as a 
factorization \emph{coalgebra}.

\medskip

We can rewrite it as the totalization of the cosimplicial factorization category with terms
$$\on{Fact}(\CO_{H^\bullet})\mod^{\on{com}},$$
where $H^\bullet$ is the \v{C}ech nerve of $\on{pt}\to \on{pt}/H$.

\medskip

On the other hand since $\on{pt}\to \LS_H^\reg$ is an fpqc cover (see \lemref{l:pt to pt/H}),
$\QCoh(\LS_H^\reg)$ identifies with the totalization of the cosimplicial factorization category
with terms $\QCoh(-)$ of the \v{C}ech nerve of $\on{pt}\to \LS_H^\reg$, the latter being
$\fL^+_\nabla(H^\bullet)$.

\medskip

How the desired equivalence follows from the fact that
$$\CO_{\fL^+_\nabla(H^\bullet)}\simeq \on{Fact}(\CO_{H^\bullet}),$$
see \eqref{e:arcs as fact} below.

\begin{rem}
For completeness, we remark on the following comparison with \cite{Ra3}. 

\medskip

In \emph{loc. cit}., a different construction of the factorization category associated to a symmetric monoidal 
category was used; in particular, $\Rep(H)_\Ran$ has an a priori different meaning than how it is used in this paper. 
One can directly compare the two constructions, but rather than doing so here, we note that the above material combined 
with \cite[Lemma 9.8.1]{Ra4} allows us to indirectly deduce that the two constructions coincide.
\end{rem}

\sssec{} 

Let $\CZ_\Ran\to \Ran$ be a prestack. We attach to it a sheaf of categories $\ul\QCoh{}_{\on{co}}(\CZ)$ over $\Ran$
by setting for $S\in \affSch_{/\Ran}$
$$\QCoh_{\on{co}}(\CZ)_S:=\QCoh_{\on{co}}(S\underset{\Ran}\times \CZ_\Ran).$$

The sheaf of categories structure holds thanks to \lemref{l:QCoh co base change}. 

\sssec{} \label{sss:QCoh co fact}

Let $\CT$ be factorization space over $X$. Consider the corresponding sheaf of categories $\ul\QCoh{}_{\on{co}}(\CT)$, i.e., 
$$\QCoh{}_{\on{co}}(\CT)_S:=\QCoh_{\on{co}}(\CT_S).$$

\medskip

The factorization structure on $\CT$ induces a \emph{factorization} structure on $\ul\QCoh{}_{\on{co}}(\CT)$,
see \lemref{sss:QCoh co mult}. We will denote the resulting factorization category by $\QCoh_{\on{co}}(\CT)$.

\medskip

By \secref{sss:t on QCoh co}, if $\CT$ is a factorization ind-scheme, the factorization category $\QCoh_{\on{co}}(\CT)$
carries a naturally defined t-structure. 

\sssec{}  \label{sss:S alpha}

Let $\CZ_\Ran\to \Ran$ be a prestack. We will attach to it a sheaf of categories over $\Ran$, denoted 
$\ul\IndCoh^!(\CZ)$. Unlike the cases of $\QCoh(-)$ and $\QCoh_{\on{co}}(-)$, this will use some special
features of $\Ran$. 

\medskip

Namely, we will use the fact that $\Ran$ can be exhibited as a colimit of prestacks
$S_{\alpha,\dr}$, where $S_\alpha$ are \emph{smooth} schemes with transition 
$f_{\alpha,\beta}:S_\alpha\to S_\beta$ maps being closed embeddings. In practice,
$S_\alpha=X^I$ for finite non-empty sets $I$, see \secref{sss:Ran as colim}. 

\sssec{} \label{sss:IndCoh ! Ran 1}

For $S$ as above set
$$\IndCoh^!(\CZ)_{S_\alpha}:=\IndCoh(\CZ_{S_\alpha}).$$

For a map $S_\alpha\overset{f_{\alpha,\beta}}\to S_\beta$ denote by the same symbol $f_{\alpha,\beta}$
the corresponding map
$$\CZ_{S_\alpha}\to \CZ_{S_\beta}.$$

The functor
$$f^!_{\alpha,\beta}:\IndCoh^!(\CZ_{S_\beta})\to \IndCoh^!(\CZ_{S_\alpha})$$
gives rise to a functor 
\begin{equation} \label{e:IndCoh^! transition}
\QCoh(S_\alpha)\underset{\QCoh(S_\beta)}\otimes \IndCoh^!(\CZ_{S_\beta})\to \IndCoh^!(\CZ_{S_\alpha}).
\end{equation}

We claim:

\begin{lem} \label{l:IndCoh^! transition}
The functor \eqref{e:IndCoh^! transition} is an equivalence.
\end{lem} 

\begin{proof}

The question is Zariski-local, so we can assume that $S_\alpha$ and $S_\beta$ are affine. 
By \lemref{l:IndCoh^! tensor up}, the functor \eqref{e:IndCoh^! transition} is fully faithful. So we only have
to show that it is essentially surjective. 

\medskip

Let $S_\beta^\wedge$ be the formal completion of $S_\beta$ along $S_\alpha$. Consider the corresponding functor
\begin{equation} \label{e:IndCoh^! transition compl}
\QCoh(S^\wedge_\beta)\underset{\QCoh(S_\beta)}\otimes \IndCoh^!(\CZ_{S_\beta})\to \IndCoh^!(\CZ_{S^\wedge_\beta}).
\end{equation} 

We have a commutative diagram
$$
\CD
\QCoh(S^\wedge_\beta)\underset{\QCoh(S_\beta)}\otimes \IndCoh^!(\CZ_{S_\beta}) @>>>  \IndCoh^!(\CZ_{S^\wedge_\beta}) \\
@AAA @AAA \\
\QCoh(S_\beta)\underset{\QCoh(S_\beta)}\otimes \IndCoh^!(\CZ_{S_\beta}) @>{=}>> \IndCoh^!(\CZ_{S_\beta}). 
\endCD
$$

The left vertical arrow is a colocalization. We claim that the right vertical arrow is also a colocalization: indeed,
this follows from \propref{p:form compl IndCoh!}. Hence, we obtain that the 
functor \eqref{e:IndCoh^! transition compl} is also a colocalization.\footnote{In fact, the above argument shows that
\eqref{e:IndCoh^! transition compl} is an equivalence.} In particular, \eqref{e:IndCoh^! transition compl} is essentially
surjective.

\medskip

We have a commutative diagram
\begin{equation} \label{e:IndCoh^! transition compl and not}
\CD
\QCoh(S_\alpha)\underset{\QCoh(S_\beta)}\otimes \IndCoh^!(\CZ_{S_\beta}) @>>>  \IndCoh^!(\CZ_{S_\alpha}) \\
@AAA @AAA \\
\QCoh(S^\wedge_\beta)\underset{\QCoh(S_\beta)}\otimes \IndCoh^!(\CZ_{S_\beta}) @>>>  \IndCoh^!(\CZ_{S^\wedge_\beta}),
\endCD 
\end{equation} 
where:

\begin{itemize}

\item The right vertical arrow is given by !-pullback along $\CZ_{S_\alpha}\to \CZ_{S^\wedge_\beta}$;

\item The left vertical arrows is given by the *-pullback functor along $S_\alpha\to S_\beta^\wedge$
along the $\QCoh(-)$ factors. 

\end{itemize}

We wish to show that the top horizontal arrow in \eqref{e:IndCoh^! transition compl and not} is essentially
surjective. By the above, the bottom horizontal arrow is essentially surjective. Hence, it suffices to show that the functor
$$\IndCoh^!(\CZ_{S^\wedge_\beta})\to \IndCoh^!(\CZ_{S_\alpha})$$
is essentially surjective. 

\medskip

We claim, however, that the map $\CZ_{S_\alpha}\to \CZ_{S^\wedge_\beta}$ admits a retraction. Indeed, since
$S_\alpha$ is smooth, the embedding
$$S_\alpha\to S^\wedge_\beta$$ 
admits a retraction, denote it by $g_{\beta,\alpha}$.  Note that the two maps
$$S^\wedge_\beta\to \Ran \text{ and } S^\wedge_\beta\overset{g_{\beta,\alpha}}\to S_\alpha\to \Ran$$
agree on $(S^\wedge_\beta)_{\on{red}}\simeq S_\alpha$. Hence, since $\Ran\to \Ran_{\dr}$ is an isomorphism, we can identify
$$\CZ_{S^\wedge_\beta} \simeq S^\wedge_\beta\underset{g_{\beta,\alpha},S_\alpha}\times \CZ_{S_\alpha}.$$

In terms of this identification, the projection
$$S^\wedge_\beta\underset{g_{\beta,\alpha},S_\alpha}\times \CZ_{S_\alpha}\to  \CZ_{S_\alpha}$$
provides the sought-for retraction. 

\end{proof} 

\sssec{}

Let $S$ be an affine scheme mapping to $\Ran$. This map factors as
$$S\overset{f}\to S_\alpha \to \Ran$$
for some $\alpha$. 

\medskip

Set
$$\IndCoh^!(\CZ)_{S,f}:=\QCoh(S)\underset{\QCoh(S_\alpha)}\otimes \IndCoh^!(\CZ_{S_\alpha}).$$

In order to show that the assignment
$$S\in \affSch_{/\Ran} \rightsquigarrow \IndCoh^!(\CZ)_{S,f}$$
gives a well-defined sheaf of categories over $\Ran$, it remains to show that for two maps 
$f_1$ and $f_2$ as above, for which $f_1|_{S_{\on{red}}}\simeq f_2|_{S_{\on{red}}}$, we have a canonical
identification
$$\IndCoh^!(\CZ)_{S,f_1}\simeq \IndCoh^!(\CZ)_{S,f_2},$$
i.e., 
\begin{equation} \label{e:IndCoh^! transition red comparison}
\QCoh(S)\underset{f_1^*,\QCoh(S_\alpha)}\otimes \IndCoh^!(\CZ_{S_\alpha})\simeq
\QCoh(S)\underset{f_2^*,\QCoh(S_\alpha)}\otimes \IndCoh^!(\CZ_{S_\alpha}).
\end{equation} 

(In addition, one needs to show that these identifications satisfy a homotopy-coherent system
of compatibilities for multi-fold comparisons $f_1|_{S_{\on{red}}}\simeq f_2|_{S_{\on{red}}}\simeq...\simeq f_n|_{S_{\on{red}}}$,
but this will be automatic from the construction explained below.) 

\sssec{}

Let $(S_\alpha\times S_\alpha)^\wedge$ be the formal completion of the diagonal in $S_\alpha\times S_\alpha$.
Note that we have a well-defined map
$$(S_\alpha\times S_\alpha)^\wedge\to \Ran,$$
Consider the corresponding prestack $\CZ_{(S_\alpha\times S_\alpha)^\wedge}$ and the category
$$\IndCoh^!(\CZ_{(S_\alpha\times S_\alpha)^\wedge}).$$

\medskip

For $i=1,2$, let 
$$p_i:(S_\alpha\times S_\alpha)^\wedge\to S_\alpha$$
denote the corresponding projection. We will denote by the same symbol $p_i$ the corresponding map
$$\CZ_{(S_\alpha\times S_\alpha)^\wedge}\to \CZ_{S_\alpha}.$$

The functor
$$p_i^!:\IndCoh^!(\CZ_{S_\alpha})\to \IndCoh^!(\CZ_{(S_\alpha\times S_\alpha)^\wedge})$$
gives rise to a functor
\begin{equation} \label{e:form compl diag base change}
\QCoh((S_\alpha\times S_\alpha)^\wedge)\underset{p_i^*,\QCoh(S_\alpha)}\otimes 
\IndCoh^!(\CZ_{S_\alpha})\to \IndCoh^!(\CZ_{(S_\alpha\times S_\alpha)^\wedge}).
\end{equation} 

We claim:

\begin{lem} \label{l:form compl diag base change}
The functor \eqref{e:form compl diag base change} is an equivalence.
\end{lem}

\begin{proof}
Proceeds along the same lines as the proof of \lemref{l:IndCoh^! transition}.
\end{proof} 

\sssec{} \label{sss:IndCoh ! Ran 2}

Using \lemref{l:form compl diag base change}, we obtain 
\begin{multline*} 
\QCoh(S)\underset{f_i^*,\QCoh(S_\alpha)}\otimes \IndCoh^!(\CZ_{S_\alpha})=\\
=\QCoh(S)\underset{(f_1\times f_2)^*,\QCoh((S_\alpha\times S_\alpha)^\wedge)}\otimes 
\QCoh((S_\alpha\times S_\alpha)^\wedge) \underset{p_i^*,\QCoh(S_\alpha)}\otimes 
\IndCoh^!(\CZ_{S_\alpha})\simeq \\
\simeq \QCoh(S)\underset{(f_1\times f_2)^*,\QCoh((S_\alpha\times S_\alpha)^\wedge)}\otimes 
\IndCoh^!(\CZ_{(S_\alpha\times S_\alpha)^\wedge}),
\end{multline*} 
thereby establishing \eqref{e:IndCoh^! transition red comparison}. 

\medskip

This completes the construction of $\ul\IndCoh^!(\CZ)$ as a sheaf of categories over $\Ran$.

\sssec{} \label{sss:IndCoh ! Ran}

Let $\CT$ be a factorization space over $X$. The multiplicative structure in \secref{sss:boxtimes IndCoh !}
equips the corresponding sheaf of categories $\ul\IndCoh^!(\CT)$ with a structure of lax factorization
category; we will denote it by $\IndCoh^!(\CT)$. 

\medskip

Suppose now that $\CT$ is ind-placid. In this case, from \lemref{l:boxtimes IndCoh placid}, we obtain that the
lax factorization structure on $\IndCoh^!(\CT)$ is a \emph{factorization} structure. 

\sssec{}  \label{sss:IndCoh* over Ran}

Let now $\CZ_\Ran\to \Ran$ be a relative ind-placid ind-scheme. In this case, we are going to define the sheaf of categories
$\ul\IndCoh^*(\CZ)$.

\medskip

We proceed with the same recipe as in the case of $\ul\IndCoh^!(\CZ)$. For an index $\alpha$, set 
$$\IndCoh^*(\CZ)_{S_\alpha}:=\IndCoh^*(\CZ_{S_\alpha}).$$

For a map $f_{\alpha,\beta}$, we have a well-defined functor
$$f_{\alpha,\beta}^{*,\IndCoh}:\IndCoh^*(\CZ_{S_\beta})\to \IndCoh^*(\CZ_{S_\alpha})$$
(see \secref{sss:^* IndCoh *}). Consider the resulting functor
\begin{equation} \label{e:IndCoh^* transition}
\QCoh(S_\alpha)\underset{\QCoh(S_\beta)}\otimes \IndCoh^*(\CZ_{S_\beta})\to \IndCoh^*(\CZ_{S_\alpha}).
\end{equation}

We claim:

\begin{lem} \label{l:IndCoh^* transition}
The functor \eqref{e:IndCoh^* transition} is an equivalence.
\end{lem} 

\begin{proof}

We can reformulate the assertion of the lemma as saying that the right adjoint of $f_{\alpha,\beta}^{*,\IndCoh}$, i.e.,
the functor
$$(f_{\alpha,\beta})^\IndCoh_*:\IndCoh^*(\CZ_{S_\alpha})\to \IndCoh^*(\CZ_{S_\beta})$$
gives rise to an equivalence
\begin{equation} \label{e:IndCoh^* transition adj}
\IndCoh^*(\CZ_{S_\alpha})\to (f_{\alpha,\beta})_*(\CO_{S_\alpha})\mod(\IndCoh^*(\CZ_{S_\beta})).
\end{equation}

\medskip

Note that \lemref{l:IndCoh^! transition} can be reformulated as saying that the functor
$$f_{\alpha,\beta}^!:\IndCoh^!(\CZ_{S_\beta})\to \IndCoh^!(\CZ_{S_\alpha})$$
admits a right adjoint, and this right adjoint identifies
$$(f_{\alpha,\beta})_*(\CO_{S_\alpha})\mod(\IndCoh^!(\CZ_{S_\beta})) \simeq \IndCoh^!(\CZ_{S_\alpha}).$$

Passing to the duals, we obtain that the dual of $f_{\alpha,\beta}^!$ identifies
\begin{equation} \label{e:IndCoh^! transition dual}
(f_{\alpha,\beta})_*(\CO_{S_\alpha})\mod(\IndCoh^!(\CZ_{S_\beta})^\vee) \simeq \IndCoh^!(\CZ_{S_\alpha})^\vee.
\end{equation}

We will now use the identifications
$$\IndCoh^*(\CZ_{S_\alpha})\simeq \IndCoh^!(\CZ_{S_\alpha})^\vee \text{ and } 
\IndCoh^*(\CZ_{S_\beta})\simeq \IndCoh^!(\CZ_{S_\beta})^\vee,$$
see \secref{sss:IndCoh ! * duality}.

\medskip

Under these identifications 
$$(f_{\alpha,\beta}^!)^\vee \simeq (f_{\alpha,\beta})^\IndCoh_*.$$

Unwinding the definitions, it is easy to see that the resulting identification \eqref{e:IndCoh^! transition dual} is the 
same as \eqref{e:IndCoh^* transition adj}. 

\end{proof} 

\sssec{} \label{sss:IndCoh* over Ran transition}

Let $S$ be an affine scheme mapping to $\Ran$. This map factors as
$$S\overset{f}\to S_\alpha \to \Ran$$
for some $\alpha$. 

\medskip

Set
$$\IndCoh^*(\CZ)_{S,f}:=\QCoh(S)\underset{\QCoh(S_\alpha)}\otimes \IndCoh^*(\CZ_{S_\alpha}).$$

As in the case of $\IndCoh^!$, in order to complete the construction of $\ul\IndCoh^*(\CZ)$ as a sheaf
of categories over $\Ran$, it suffices to show that for a pair of maps 
$f_1$ and $f_2$ as above, for which $f_1|_{S_{\on{red}}}\simeq f_2|_{S_{\on{red}}}$, we have a canonical
identification
$$\IndCoh^*(\CZ)_{S,f_1}\simeq \IndCoh^*(\CZ)_{S,f_2},$$
i.e., 
\begin{equation} \label{e:IndCoh^* transition red comparison}
\QCoh(S)\underset{f_1^*,\QCoh(S_\alpha)}\otimes \IndCoh^*(\CZ_{S_\alpha})\simeq
\QCoh(S)\underset{f_2^*,\QCoh(S_\alpha)}\otimes \IndCoh^*(\CZ_{S_\alpha}).
\end{equation} 

This follows from the following: 

\begin{lem} \label{l:form compl diag base change *}
In the notations of \lemref{l:form compl diag base change}, the functor
$$\QCoh((S_\alpha\times S_\alpha)^\wedge)\underset{p_i^*,\QCoh(S_\alpha)}\otimes 
\IndCoh^*(\CZ_{S_\alpha})\to \IndCoh^*(\CZ_{(S_\alpha\times S_\alpha)^\wedge}).$$
 is an equivalence.
\end{lem}

The lemma follows by duality from \lemref{l:form compl diag base change *}. 

\sssec{} 

By construction, the sheaves of categories $\ul\IndCoh^*(\CZ)$ and $\ul\IndCoh^!(\CZ)$
are mutually dual. 

\medskip

By \secref{ss:t-str IndCoh *}, the sheaf of categories $\ul\IndCoh^*(\CZ)$ is equipped with a t-structure. 

\sssec{} \label{sss:IndCoh * Ran}

Let $\CT$ be a factorization ind-placid ind-scheme over $X$. The multiplicative structure in \secref{sss:boxtimes IndCoh *}
and \lemref{l:boxtimes IndCoh placid} imply that in this case the sheaf of categories $\ul\IndCoh^*(\CT)$ carries a factorization
structure. Denote the resulting factorization category by $\IndCoh^*(\CT)$. 

\medskip

By the construction of $\IndCoh^*(\CT)$ in \secref{sss:IndCoh* over Ran}
and \secref{ss:t-str IndCoh *}, the factorization category $\IndCoh^*(\CT)$
carries a naturally defined t-structure. 

\medskip

By construction, $\IndCoh^*(\CT)$ is dual to $\IndCoh^!(\CT)$ as a factorization category. 

\sssec{}

Let $\CT$ be a factorization space over $X$, and let $\CT_m$ be a factorization module
space over $\CT$ at some $\CZ\to \Ran$.

\medskip

Suppose that $\CT$ is such that the categories $\QCoh(\CT_S)$ for $S\in \affSch_{/\Ran}$ are dualizable. 
Then the sheaf of categories
$$\pi_*(\ul\QCoh((\CT_m)_{\CZ^{\subseteq}}))$$
on $\CZ^{\subseteq}$ (here $\pi$ denotes the structural map $\CT_{\CZ^{\subseteq}}\to \CZ^{\subseteq}$) 
admits a natural structure of factorization module category over $\QCoh(\CT)$ at $\CZ$.
We will denote it by $\QCoh(\CT_m)$. 

\medskip

For general $\CT$ and $\CT_m$, the assignment 
$$S\to \CZ^{\subseteq}, \rightsquigarrow \, \QCoh_{\on{co}}((\CT_m)_S)$$
is a sheaf of categories over $\CZ^{\subseteq}$, to be denoted  $\ul\QCoh{}_{\on{co}}(\CT_m)$. It has a natural 
structure of factorization module category over $\QCoh_{\on{co}}(\CT)$ at $\CZ$. We will denote it by $\QCoh_{\on{co}}(\CT_m)$. 

\medskip

Assume now that $\CT_m$ is an ind-placid ind-scheme relative to $\CZ^{\subseteq}$. In this case, we can consider the
sheaves of categories 
$$\ul\IndCoh^!(\CT_m) \text{ and } \ul\IndCoh^*(\CT_m)$$
and they have natural factorization module structures over $\IndCoh^!(\CT)$ and $\IndCoh^*(\CT)$, respectively.
We will denote the resulting factorization module categories by $\IndCoh^!(\CT_m)$ and $\IndCoh^*(\CT_m)$,
respectively. 

\ssec{Modules over Kac-Moody Lie algebras} \label{ss:KM fact}

In this subsection we will show how to adapt the theory developed in \cite{Ra5} to the factorization setting. 
We start be defining the factorization category $\Rep(\fL^+(G))$. 

\medskip

In this subsection $G$ can be arbitrary an algebraic group (i.e., not necessarily reductive). 

\sssec{}  \label{sss:Rep L^+G} 

For $S\in \affSch_{/\Ran}$
consider the group scheme $\fL^+(G)_S$ over $S$. It is pro-smooth, 
$$\fL^+(G)_S \simeq \underset{\alpha\in A^{\on{op}}}{\on{lim}}\, G^\alpha_S,$$
where $G^\alpha_S$ are smooth group-schemes over $S$ of finite type, and the
transition maps are smooth and surjective. 

\medskip

For every $\alpha$ we can consider the algebraic stack $\on{pt}/G^\alpha_S$. Set
$$\Rep(G^\alpha_S):=\QCoh(\on{pt}/G^\alpha_S).$$

Set
\begin{equation} \label{e:Rep L^+G}
\Rep(\fL^+(G))_S:=\underset{\alpha\in A}{\on{colim}}\, \Rep(G^\alpha_S),
\end{equation} 
where the transition functors are given by restriction:
$$G^\beta_S\twoheadrightarrow G^\alpha_S\,\, \, \rightsquigarrow \,\,\, \Rep(G^\alpha_S)\to \Rep(G^\beta_S).$$

Note that the above transition functors admit (continuous) right adjoints, given by 
$$\on{inv}_{\on{ker}(G^\beta_S\to G^\alpha_S)}.$$

\medskip

Hence, we can rewrite $\Rep(\fL^+(G))_S$ also as the limit 
$$\Rep(\fL^+(G))_S:=\underset{\alpha\in A^{\on{op}}}{\on{lim}}\, \Rep(G^\alpha_S),$$
with respect to the above right adjoints. 

\sssec{} \label{sss:Rep(G) comp gen alpha}

We claim that $\Rep(\fL^+(G))_S$ is compactly generated. In order to prove that, it
suffices to show that each $\Rep(G^\alpha_S)$ is compactly generated. This can be proved
on general grounds (the category of quasi-coherent sheaves on a smooth algebraic stack 
is compactly generated). What follows below is an explicit construction of compact generators. 

\medskip

First, it is easy to see that if an object $\CV\in \Rep(G^\alpha_S)$ is such that $\oblv_{G^\alpha_S}(\CV)\in \QCoh(S)$ is
compact, then $\CV$ itself is compact. Hence, in order prove the compact generation, it suffices to exhibit a generating
collection of compact objects. We do that as follows.

\medskip

We can assume that the category of indices $A$ has an initial element $\alpha_0$ for
which $G^{\alpha_0}_S$ is the following explicit group-scheme:

\medskip

Let $I$ be a finite set such that the map $S\to \Ran$ factors as
$$S\to X^I\to \Ran.$$

Let $\on{Graph}_I\subset X^I\times X$ be the incidence divisor. Let $G_{X^I}$ be the group-scheme over $X^I$ equal to 
the restriction of scalars \`a la Weil along 
$$\on{Graph}_I\subset X^I\to X^I$$
of the pullback of the constant group-scheme with fiber $G$ along
$$\on{Graph}_I\subset X^I\to X.$$

\medskip

Note that $\on{Graph}_I$ receives a map from the disjoint union of $I$ many copies of $X^I$ 
(i.e., the pairwise diagonals of the $i$th and the last coordinate in $X^I\times X$). In particular,
we obtain a map
\begin{equation} \label{e:G X^I}
G_{X^I}\to (G\times X)^I
\end{equation} 
as group-schemes over $X^I$.

\medskip

We take $G^{\alpha_0}_S$ to be the pullback of $G_{X^I}$ along $S\to X^I$. 

\medskip

Note that for every $\alpha$, the kernel of the projection
\begin{equation} \label{e:proj to constant group-sch}
G^\alpha_S\to G^{\alpha_0}_S
\end{equation} 
is unipotent (in fact, admits a filtration with subquotients isomorphic to the constant group-scheme with fiber $\BG_a$). 

\medskip

Hence, the essential image of the forgetful functor
$$\Rep(G^{\alpha_0}_S)\to \Rep(G^\alpha_S)$$
generates $\Rep(G^\alpha_S)$.

\medskip

The map \eqref{e:G X^I} induces a map
$$G^\alpha_S\to G^I.$$

In particular, we obtain a functor
\begin{equation} \label{e:G X^I S}
\Rep(G)^{\otimes I}\to \Rep(G^{\alpha_0}_S).
\end{equation}

\medskip

It is easy to see that the essential image of \eqref{e:G X^I S} generates $\Rep(G^{\alpha_0}_S)$. Hence, we obtain 
that the images of the compact objects in $\Rep(G)^{\otimes I}$ under the composition
\begin{equation} \label{e:G X^I S alpha}
\Rep(G)^{\otimes I}\to \Rep(G^{\alpha_0}_S)\to  \Rep(G^\alpha_S)
\end{equation}
provide a set of compact generators of $\Rep(G^\alpha_S)$.

\begin{rem} \label{r:Rep L^+G}
Note that the category $\Rep(\fL^+(G))_S$ is \emph{not} the same as representations of the group-scheme
$\fL^+(G)_S$, i.e.,
$$\CO_{\fL^+(G)_S}\comod \simeq \QCoh(\on{pt}/\fL^+(G)_S).$$
Rather, it is its renormalized version, in which we declare the compacts to be the images of the compacts in 
$\Rep(G^\alpha_S)$ under the restriction along $\fL^+(G)_S\to G^\alpha_S$. 

\medskip

In particular, the forgetful functor
$$\oblv_{\fL^+(G)_S}:\Rep(\fL^+(G))_S\to \QCoh(S)$$
is \emph{not} conservative.

\end{rem}

\sssec{} \label{sss:t-str Rep L+G}

The presentation \eqref{e:Rep L^+G} equips $\Rep(\fL^+(G))_S$ with a t-structure, for which the forgetful functor 
$\oblv_{\fL^+(G)_S}$ is t-exact. 

\sssec{}

Let $f:S'\to S$ be a map of affine schemes. Pullback along $f$ gives rise to a functor
$$f^*: \Rep(\fL^+(G))_S\to \Rep(\fL^+(G))_{S'},$$
which in turn gives rise to a functor
\begin{equation} \label{e:Rep L^+G base change}
\QCoh(S')\underset{\QCoh(S)}\otimes \Rep(\fL^+(G))_S\to \Rep(\fL^+(G))_{S'}.
\end{equation}

We claim that the functor \eqref{e:Rep L^+G base change} is an equivalence. Indeed,
this follows from the fact that for every $\alpha$, the corresponding functor
$$\QCoh(S')\underset{\QCoh(S)}\otimes \Rep(G^\alpha_S)\to \Rep(G^\alpha_{S'})$$
is an equivalence. 

\medskip

This defines on the assignment 
$$S\in \affSch_{/\Ran}\,\, \rightsquigarrow \,\, \Rep(\fL^+(G))_S$$
a structure of sheaf of categories over $\Ran$.

\medskip

We will denote this sheaf of categories by $\ul\Rep(\fL^+(G))$. 

\sssec{} \label{sss:fact for L+G}

Let $\ul{x}_i:S_i\to \Ran$, $i=1,2$ be points such that
$$S_1\times S_2 \overset{\ul{x}_1,\ul{x}_2}\longrightarrow \Ran\times \Ran$$ lands in $(\Ran\times \Ran)_{\on{disj}}$. 

\medskip

Note that we have 
$$\fL^+(G)_{S_1\times S_2} \simeq \fL^+(G)_{S_1}\times \fL^+(G)_{S_2}.$$

Tensor product of representations defines a functor
\begin{equation} \label{e:pre-fact structure on Rep L^+ G}
\Rep(\fL^+(G)_{S_1})\otimes \Rep(\fL^+(G)_{S_2})\to \Rep(\fL^+(G)_{S_1\times S_2}).
\end{equation}

We claim that \eqref{e:pre-fact structure on Rep L^+ G} is an equivalence.\footnote{This would not be the case
(at least, not obviously so) if instead of $\Rep(\fL^+(G))_{S_i}$ we used their naive versions, i.e., the categories
of representations of the group-schemes $\fL^+(G)_{S_i}$.} Indeed, we can write $\fL^+(G)_{S_1\times S_2}$ as 
$$\underset{(\alpha_1,\alpha_2)\in A_1^{\on{op}}\times A_2^{\on{op}}}{\on{lim}}\, G^{\alpha_1}_{S_1}\times G^{\alpha_2}_{S_2},$$
and for each pair of indices $\alpha_1,\alpha_2$, the functor
$$\Rep(G^{\alpha_1}_{S_1})\otimes \Rep(G^{\alpha_2}_{S_2})\to \Rep(G^{\alpha_1}_{S_1}\times G^{\alpha_2}_{S_2})$$
is an equivalence (e.g., by \cite[Cor. 10.3.6]{Ga7}). 

\medskip

This endows the sheaf of categories $\ul\Rep(\fL^+(G))$ with a factorization structure. We denote the resulting
factorization category by $\Rep(\fL^+(G))$. 

\medskip

The t-structures in \secref{sss:t-str Rep L+G} define a t-structure on $\Rep(\fL^+(G))$ as a factorization category. 

\sssec{}

Let $S\in \affSch_{/\Ran}$ be as above. Following \cite{Ra5}, one defines the (2-)category 
$$\fL^+(G)_S\mmod^{\on{weak}}$$
of $\QCoh(S)$-linear categories equipped with a \emph{weak} action of $\fL^+(G)_S$ to be equivalent to
$$\Rep(\fL^+(G))_S\mmod.$$

We consider the forgetful functor
\begin{equation} \label{e:forget weak action}
\oblv_{\fL^+(G)_S,\on{weak}}:\Rep(\fL^+(G))_S\mmod\to \QCoh(S)\mmod
\end{equation} 
given by
$$\bC\mapsto \bC\underset{\Rep(\fL^+(G))_S,\oblv_{\fL^+(G)_S}}\otimes \QCoh(S).$$

\begin{rem}
As in Remark \ref{r:Rep L^+G}, the 2-category $\Rep(\fL^+(G))_S\mmod$ is \emph{not} the same
as $\QCoh(S)$-linear categories equipped with a co-action of $\QCoh(\fL^+(G)_S)$. Indeed,
for the unit object 
$$\Rep(\fL^+(G))_S\in \Rep(\fL^+(G))_S\mmod,$$
its category of endofunctors in $\Rep(\fL^+(G))_S\mmod$ is $\Rep(\fL^+(G))_S$, while for the unit
object 
$$\QCoh(S)\in \QCoh(\fL^+(G)_S)\commod,$$
its category of endofunctors is the non-renormalized category 
$$\CO_{\fL^+(G)_S}\comod \simeq \QCoh(\on{pt}/\fL^+(G)_S).$$

\medskip

That said, for each individual $\alpha$, the stack $\on{pt}/G^\alpha_S$ is 1-affine, and hence the
functor
\begin{equation} \label{e:pt mod G alpha 1-aff}
\Rep(G^\alpha_S)\mmod \to \QCoh(G^\alpha_S)\commod, \quad \bC\mapsto 
\bC\underset{\Rep(G^\alpha_S),\oblv_{G^\alpha_S}}\otimes \QCoh(S)
\end{equation} 
is an equivalence

\end{rem}

\sssec{}

Parallel to \cite{Ra5}, one defines the (2-)category 
$$\fL^+(G)_S\mmod$$
of $\QCoh(S)$-linear categories equipped with a \emph{strong} action of $\fL^+(G)_S$. Namely,
this is the category of comodules (inside $\QCoh(S)\mmod$) for 
$$\Dmod_{\on{rel}/S}(\fL^+(G)_S),$$
where:

\begin{itemize}

\item $\Dmod_{\on{rel}/S}(\fL^+(G)_S):= \underset{\alpha\in A}{\on{colim}}\, \Dmod_{\on{rel}/S}(G^\alpha_S)$;

\item $\Dmod_{\on{rel}/S}(G^\alpha_S):=\QCoh(G^\alpha_{S_\dr}\underset{S_\dr}\times S)$. 

\end{itemize} 

One shows that any object 
$$\bC\in \fL^+(G)_S\mmod$$ 
can be \emph{canonically written} as
$$\underset{\alpha\in A}{\on{colim}}\, \bC^{\on{ker}(\fL^+(G)_S\to G^\alpha_S)},$$
where
$$\bC^{\on{ker}(\fL^+(G)_S\to G^\alpha_S)}\subset \bC$$
is the full subcategory of strong invariants with resect to the (pro-unipotent) group-scheme
$$\on{ker}(\fL^+(G)_S\to G^\alpha_S).$$

\sssec{}

We have a forgetful functor
$$\oblv^{\on{strong}}_{\on{weak}}:
\fL^+(G)_S\mmod\to \fL^+(G)_S\mmod^{\on{weak}}$$
that sends an object $\bC$ to
$$\underset{\alpha\in A}{\on{colim}}\, \oblv^{\on{strong}}_{\on{weak}}(\bC^{\on{ker}(\fL^+(G)_S\to G^\alpha_S)}),$$
where in the right-hand side $\oblv^{\on{strong}}_{\on{weak}}$ denotes the family of forgetful functors
$$\Dmod_{\on{rel}/S}(G^\alpha_S)\commod\overset{\oblv^l}\to\QCoh(G^\alpha_S)\commod \overset{\text{\eqref{e:pt mod G alpha 1-aff}}}
\simeq \Rep(G^\alpha_S)\mmod.$$

This functor intertwines the natural forgetful functor 
$$\oblv_{\fL^+(G)_S,\on{weak}}:\fL^+(G)_S\mmod\to \QCoh(S)\mmod$$
with the functor $\oblv_{\fL^+(G)_S,\on{weak}}$ of \eqref{e:forget weak action}. 

\sssec{} \label{sss:mods for Lie alg}

Our next goal is to define the category $\fL^+(\fg)\mod$ of modules for the arc Lie algebra $\fL^+(\fg)$. The definition
that we are about to give mimics the following finite-dimensional situation:

\medskip

Let $H$ be a finite-dimensional algebraic group. Then the category $\fh\mod$ of modules over its
Lie algebra has the following structures:

\medskip

\begin{itemize}

\item It carries a (strong) action of $H$;

\item It is equipped with a forgetful functor $\oblv_\fh:\fh\mod\to \Vect$;

\item The functor $\oblv_\fh$ is equipped with a structure of compatibility with 
the induced \emph{weak} action of $H$.

\end{itemize}

Moreover, the category $\fh\mod$ is universal with respect to the above pieces of structure. I.e.,
for a category $\bC$, equipped with a (strong) action of $H$, compositing with $\oblv_\fh$ defines
an equivalence between:

\begin{itemize} 

\item Functors $\bC\to \fh\mod$, compatible with (strong) actions of $H$;

\item Functors $\bC\to \Vect$, compatible with weak actions of $H$.

\end{itemize}

The above universal property can be established by realizing $\fh\mod$ as follows
\begin{equation} \label{e:mods for Lie alg}
\fh\mod\simeq \Dmod(H)^{H\on{-weak}}.
\end{equation} 

\sssec{}

We define the category 
$$\fL^+(\fg)\mod_S$$
by the universal property as in \secref{sss:mods for Lie alg}. I.e., this is a category, equipped with:

\begin{itemize}

\item A (strong) action of $\fL^+(G)_S$;

\item A forgetful functor $\oblv_{\fL^+(\fg)}:\fL^+(\fg)\mod_S\to \QCoh(S)$;

\item A datum of compatibility on $\oblv_{\fL^+(\fg)}$ with the weak action of $\fL^+(G)_S$.

\end{itemize}

Moreover, $\fL^+(\fg)\mod_S$ is universal with respect to the above pieces of structure. 

\medskip

As in \eqref{e:mods for Lie alg},we can explicitly realize $\fL^+(\fg)\mod_S$ as follows:
\begin{equation} \label{e:fL + mod Lie expl}
\fL^+(\fg)\mod_S\simeq \Dmod_{\on{rel}/S}(\fL^+(G)_S)^{\fL^+(G)_S\on{-weak}}.
\end{equation} 

\begin{rem}
A feature of $\fL^+(\fg)\mod_S$ that one has to keep in mind, and which distinguishes it from the finite-dimensional
situation, is that the functor
$$\oblv_{\fL^+(\fg)}:\fL^+(\fg)\mod_S\to \QCoh(S)$$
is \emph{not} conservative. 
\end{rem} 

\sssec{}

Note that, by definition, we have
\begin{multline} \label{e:equiv L^+G Lie mods}
\fL^+(\fg)\mod_S^{\fL^+(G)_S} \simeq
\on{Funct}_{\fL^+(G)_S\mmod}(\QCoh(S),\fL^+(\fg)\mod_S)\simeq  \\
\simeq \on{Funct}_{\fL^+(G)_S\mmod^{\on{weak}}}(\QCoh(S),\QCoh(S))\simeq
\Rep(\fL^+(G))_S.
\end{multline}

\sssec{}

By the universal property of $\fL^+(\fg)\mod_S$, we have naturally defined restriction functors
$$\on{Lie}(G^\alpha_S)\mod\to \fL^+(\fg)\mod_S.$$

Furthermore, the functor
\begin{equation} \label{e:fL + mod Lie expl as colim}
\underset{\alpha\in A}{\on{colim}}\, \on{Lie}(G^\alpha_S)\mod\to \fL^+(\fg)\mod_S
\end{equation}
is an equivalence. 

\medskip

Since the transition functors in the left-hand side of \eqref{e:fL + mod Lie expl} preserve compactness,
we obtain that $\fL^+(\fg)\mod_S$ is compactly generated. 

\sssec{} \label{sss:t-str on L^+G Lie}

The presentation of $\fL^+(\fg)\mod_S$ as in \eqref{e:fL + mod Lie expl as colim} equips it with a t-structure,
for which the forgetful functor $\oblv_{\fL^+(\fg)}$ is t-exact. 

\sssec{} \label{sss:L^G Lie shf of cat}

For $S'\to S$, the universal property of $\fL^+(\fg)\mod_S$ gives rise to a functor
$$\QCoh(S')\underset{\QCoh(S)}\otimes \fL^+(\fg)\mod_S\to \fL^+(\fg)\mod_{S'}.$$

One shows (e.g., using \eqref{e:fL + mod Lie expl}) that the above functor is equivalence.

\medskip

This endows the assignment 
$$S\,\, \rightsquigarrow \fL^+(\fg)\mod_S$$
with a structure of sheaf of categories over $\Ran$. We will denote it by $\fL^+(\fg)\ul\mod$.

\sssec{} \label{sss:L^G Lie fact}

Let $\ul{x}_i:S_i\to \Ran$ be as in \secref{sss:fact for L+G}. By the universal property of  $\fL^+(\fg)\mod_{S_1\times S_2}$,
we obtain a functor
$$\fL^+(\fg)\mod_{S_1}\otimes \fL^+(\fg)\mod_{S_2}\to \fL^+(\fg)\mod_{S_1\times S_2}.$$

One shows (e,g., using \eqref{e:fL + mod Lie expl}) that this functor is an equivalence. This endows 
the sheaf of categories $\fL^+(\fg)\ul\mod$ with a factorization structure. 

\medskip

We denote the resulting factorization category by $\fL^+(\fg)\mod$. The t-structures from \secref{sss:t-str on L^+G Lie} give
rise to a t-structure on $\fL^+(\fg)\mod$ as a factorization category. 

\sssec{} \label{sss:KL defn}

We finally consider representations of Kac-Moody algebras. First, fix $S\in \affSch_{/\Ran}$. Having 
the (2-)categories
$$\fL^+(G)_S\mmod^{\on{weak}} \text{ and } \fL^+(G)_S\mmod,$$
proceeding as in \cite[Sect. 7-8]{Ra5}, we define the (2-)categories
$$\fL(G)_S\mmod^{\on{weak}} \text{ and } \fL(G)_{\kappa,S}\mmod$$
of $\QCoh(S)$-linear categories, equipped with weak (resp., strong at level $\kappa$) actions of $\fL(G)_S$. 

\medskip

We define the category $\hg\mod_{\kappa,S}$ as the category, equipped with and universal with respect
to the following pieces of structure:

\begin{itemize}

\item A \emph{strong} action of $\fL(G)_S$ at level $\kappa$;

\item A functor to $\QCoh(S)$;

\item A datum of compatibility on the above functor with respect to the \emph{weak} action of $\fL(G)_S$.

\end{itemize}

We can explicitly realize the category $\hg\mod_{\kappa,S}$ as
$$\hg\mod_{\kappa,S}\simeq \Dmod_{\kappa,\on{rel}/S}(\fL(G)_S)^{\fL(G)_S\on{-weak}}.$$

\sssec{}

By the universal property of $\hg\mod_{\kappa,S}$, it comes equipped with a (conservative) forgetful functor
$$\oblv^{\hg_\kappa}_{\fL^+(\fg)}:\hg\mod_{\kappa,S}\to \fL^+(\fg)\mod_S.$$

As in \cite[§9.12, §11.10]{Ra5}, one shows that this functor admits a left adjoint, to be denoted
$$\ind^{\hg_\kappa}_{\fL^+(\fg)}:\fL^+(\fg)\mod_S\to \hg\mod_{\kappa,S}.$$

In particular, the fact that $\fL^+(\fg)\mod_S$ is compactly generated implies that so is
$\hg\mod_{\kappa,S}$.

\sssec{} \label{sss:t-str on LG Lie}

We note that the endofunctor of $\fL^+(\fg)\mod_S$ underlying the monad 
$$\oblv^{\hg_\kappa}_{\fL^+(\fg)}\circ \ind^{\hg_\kappa}_{\fL^+(\fg)}$$
is t-exact.

\medskip

This allows us to equip $\hg\mod_{\kappa,S}$ with a t-structure for which both 
functors $\oblv^{\hg_\kappa}_{\fL^+(\fg)}$ and $\ind^{\hg_\kappa}_{\fL^+(\fg)}$
are t-exact. 

\sssec{}

As in Sects. \secref{sss:L^G Lie shf of cat} one endows the assignment
$$S\,\,\rightsquigarrow \,\, \hg\mod_{\kappa,S}$$
with a structure of sheaf of categories over $\Ran$. We denote it by $\hg\ul\mod{}_\kappa$. 

\medskip

As in \secref{sss:L^G Lie fact}, one endows $\hg\ul\mod{}_\kappa$ with a factorization structure.
We denote the resulting factorization category by $\hg\mod_\kappa$. 

\medskip

The t-structures in \secref{sss:t-str on LG Lie} combine to a t-structure on $\hg\mod_\kappa$
as a factorization category. 

\sssec{}

For $S\in \affSch_{/\Ran}$ set
$$\KL(G)_{\kappa,S}:=(\hg\mod_{\kappa,S})^{\fL^+(G)_S}.$$

The adjoint functors
$$\ind^{\hg_\kappa}_{\fL^+(\fg)}:\fL^+(\fg)\mod_S\rightleftarrows \hg\mod_{\kappa,S}:\oblv^{\hg_\kappa}_{\fL^+(\fg)}$$
induce a pair of adjoint functors 
\begin{equation} \label{e:induce KL App}
\ind^{(\hg,\fL^+(G))_\kappa}_{\fL^+(G)}:\Rep(\fL^+(G))_S \rightleftarrows \KL(G)_{\kappa,S}: 
\oblv^{(\hg,\fL^+(G))_\kappa}_{\fL^+(G)}.
\end{equation}

In particular, we obtain that $\KL(G)_{\kappa,S}$ is compactly generated. 

\sssec{}

Let us describe in concrete terms a set of compact generators of $\KL(G)_{\kappa,S}$. Let $I$
be a finite set as in \secref{sss:Rep(G) comp gen alpha}.

\medskip

For a compact object $V\in \Rep(G)^{\otimes I}$, let us denote by a light abuse of notation
by the same character $V$ its image under \eqref{e:G X^I S alpha}.

\medskip

The objects
$$\ind^{(\hg,\fL^+(G))_\kappa}_{\fL^+(G)}(V)\in \KL(G)_{\kappa,S}$$
are called Weyl modules, and they compactly generate $\KL(G)_{\kappa,S}$.

\begin{rem}

Let $\ul\bC$ be a sheaf of categories on $S_\dR$, where $S$ is an affine scheme of finite type.
There is a stronger notion than compact generation, called \emph{ULA generation}: 

\medskip  

An object $\bc\in \Gamma(S_\dr,\ul\bC)$ is said to be ULA if its image in $\Gamma(S,\ul\bC)$
is compact. 

\medskip

We say that $\ul\bC$ is ULA-generated if it contains a collection of ULA objects, whose images
in $\Gamma(S,\ul\bC)$ generate this category. 

\medskip

We will say that a factorization category $\bC$ is ULA-generated if for every $S\to \Ran$,
the corresponding category $\bC_{S_\dr}$ is. This property is enough to check for $S=X^I$.

\medskip

Many factorization categories that appear in geometric representation theory have this 
property, e.g., $\bC=\Sph_G$. However, the category $\KL(G)_\kappa$ does not: namely, 
$\KL(G)_{\kappa,X^I}$ is not ULA-generated for $|I|\geq 2$.

\end{rem}

\sssec{} \label{sss:t-str on KL}

As in \secref{sss:t-str on LG Lie}, we obtain that $\KL(G)_{\kappa,S}$ carries a t-structure, for which both functors
\eqref{e:induce KL App} are t-exact.

\sssec{}

%

For $S'\to S$, the equivalence
$$\QCoh(S')\underset{\QCoh(S)}\otimes \hg\mod_{\kappa,S}\to \hg\mod_{\kappa,S'}$$
induces a functor
\begin{multline}  \label{e:base change KL}
\QCoh(S')\underset{\QCoh(S)}\otimes \KL(G)_{\kappa,S}=
\QCoh(S')\underset{\QCoh(S)}\otimes \hg\mod_{\kappa,S}^{\fL^+(G)_S}\to \\
\to (\QCoh(S')\underset{\QCoh(S)}\otimes \KL(G)_{\kappa,S})^{\fL^+(G)_S}\simeq
(\QCoh(S')\underset{\QCoh(S)}\otimes \KL(G)_{\kappa,S})^{\fL^+(G)_{S'}}\simeq \KL(G)_{\kappa,S'}.
\end{multline}

However, the functor 
$$\bC\to \bC^{\fL^+(G)_S}, \quad \fL^+(G)_S\mmod\to \QCoh(S)\mmod$$
is known to commute with colimits. Hence, the functor \eqref{e:base change KL} is an equivalence.

\medskip

This endows the assignment
$$S\in \affSch_{/\Ran},\,\, \rightsquigarrow \,\, \KL(G)_{\kappa,S}$$
with a structure of sheaf of categories. We denote it by $\ul\KL(G)_\kappa$.

\sssec{}

Similarly, in the situation \secref{sss:fact for L+G}, the equivalence
$$\hg\mod_{\kappa,S_1}\otimes \hg\mod_{\kappa,S_2}\to \hg\mod_{\kappa,S_1\times S_2}$$
induces an equivalence
\begin{multline*}
\KL(G)_{\kappa,S_1}\otimes \KL(G)_{\kappa,S_2}=
(\hg\mod_{\kappa,S_1})^{\fL^+(G)_{S_1}}\otimes (\hg\mod_{\kappa,S_1})^{\fL^+(G)_{S_2}}\simeq \\
\simeq \left(\hg\mod_{\kappa,S_1}\otimes \hg\mod_{\kappa,S_2}\right)^{\fL^+(G)_{S_1}\times \fL^+(G)_{S_2}}=
\KL(G)_{\kappa,S_1\times S_2}.
\end{multline*}

\medskip

This endows $\ul\KL(G)_\kappa$ with a factorization structure. The resulting factorization category, denoted
$\KL(G)_\kappa$ is the sought-for factorization incarnation of the Kazhdan-Lusztig category. 

\medskip

The t-structures in \secref{sss:t-str on KL} define a t-structure on $\KL(G)_\kappa$ as a factorization category. 

\ssec{Restriction of factorization module categories, continued} \label{ss:restr fact mod}

In this subsection we will discuss some additional aspects of the operation of restriction
of factorization module categories, introduced in \secref{sss:univ property restr cat}.

\medskip

We will omit most proofs (they are elaborations of the limit construction described in \secref{sss:fact restr limit});
a more detailed discussion will appear in \cite{CFGY}.

\sssec{}

Fix $\CZ\to \Ran$. Consider the totality of pairs $(\bA,\bC)$, where $\bA$ is a factorization category,
and $\bC\in \bA\mmod^{\on{fact}}_\CZ$. We view it as a 2-category, to be denoted $\on{FactCat-and-Mod}(X)$
with the above objects, and the categories of morphisms defined as follows:

\medskip

For a pair of objects $(\bA_1,\bC_1)$ and $(\bA_2,\bC_2)$ of $\on{FactCat-and-Mod}(X)$, the category
$$\Maps_{\on{FactCat-and-Mod}(X)}((\bA_1,\bC_1),(\bA_2,\bC_2))$$
consists of pairs of morphisms $(\Phi:\bA_1\to \bA_2,\Phi_m:\bC_1\to \bC_2)$ as in \secref{sss:univ property restr cat}.

\medskip

Consider the natural projection
\begin{equation}  \label{e:forget module}
\on{FactCat-and-Mod}(X)\to \on{FactCat}(X).
\end{equation} 

The following assertion encodes the functoriality of the assignment 
$$\bA\rightsquigarrow \bC\in \bA\mmod^{\on{fact}}_\CZ.$$

\begin{thm} \label{t:2-Cart fact mod}
The functor 
$$\on{FactCat-and-Mod}(X)^{2\on{-op}}\to \on{FactCat}(X)^{2\on{-op}},$$
induced by \eqref{e:forget module} is a 2-Cartesian fibration, where:

\begin{itemize}

\item The symbol $(-)^{2\on{-op}}$ refers to reversing 2-morphisms in a given 2-category;

\item The notion of 2-Cartesian fibration is an in \cite[Chapter 11, Sect. 1.1.]{GaRo3}.

\end{itemize}

\end{thm} 

\begin{rem}

Note that the situation in \thmref{t:2-Cart fact mod} is parallel to the following more familiar paradigm,
when instead of $\on{FactCat-and-Mod}(X)$ we consider the category of pairs $(\bA,\bC)$, where
$\bA$ is a monoidal category, and $\bC\in \bA\mmod$, and where the category of functors
$$(\bA_1,\bC_1)\to (\bA_2,\bC_2)$$ consists of right-lax monoidal functors between the monoidal
categories and compatible right-lax functors between the module categories.

\end{rem}

\sssec{}

The concrete meaning of this theorem is the following. It says that:

\begin{itemize}

\item For  
$$\bA_1\overset{\Phi_{1,2}}\to \bA_2 \overset{\Phi_{1,2}}\to \bA_3, \quad \Phi_{1,3}=\Phi_{2,3}\circ \Phi_{1,2},$$
and $\bC_3\in \bA_3\mmod^{\on{fact}}_\CZ$, the tautological functor 
$$\Res_{\Phi_{1,2}}\circ \Res_{\Phi_{2,3}}(\bC_3)\to \Res_{\Phi_{1,3}}$$
is an equivalence; 

\item For $\Phi:\bA_1\to \bA_2$, $\bC_i\in  \bA_i\mmod^{\on{fact}}_\CZ$, $i=1,2$, 
$\Phi_m\in \on{Funct}_{\Phi:\bA_1\to \bA_2}(\bC_1,\bC_2)$ and a natural transformation $\Phi\to \Phi'$, there exists
an object $\Phi'_m\in \on{Funct}_{\Phi':\bA_1\to \bA_2}(\bC_1,\bC_2)$ together with a compatible natural transformation
$\Phi'_m\to \Phi_m$, which is universal with respect to this property. 

\end{itemize} 

\sssec{}

According to \cite[Chapter 11, Theorem 1.1.8]{GaRo3}, we can interpret \thmref{t:2-Cart fact mod} as saying that
the assignment
$$\bA\mapsto \bA\mmod^{\on{fact}}_\CZ$$
extends to a functor between 2-categories
$$\on{FactCat}(X)^{1\on{-op},2\on{-op}}\to 2\on{-Cat},$$
where:

\begin{itemize}

\item The symbol $(-)^{1\on{-op},2\on{-op}}$ refers to reversing both 1-morphisms and 2-morphisms in a given 2-category;

\item $2\on{-Cat}$ is the totality of $(\infty,2)$-categories, viewed as a 2-category.

\end{itemize}

\sssec{}

We will need the following corollary of \thmref{t:2-Cart fact mod}. Let $\bA_1$ and $\bA_2$ be a pair of factorization
categories, and let 
$$\Phi:\bA_1\leftrightarrows \bA_2:\Phi^R$$
be a pair of factorization functors.

\medskip

We claim:

\begin{cor} \label{c:basic adj non-unital}
For $\bC_i\in \bA_i\mmod^{\on{fact}}_\CZ$, $i=1,2$ there is a canonical equivalence
$$\on{Funct}_{\bA_1\mmod^{\on{fact}}_\CZ}(\Res_\Phi(\bC_2),\bC_1) \simeq 
\on{Funct}_{\bA_2\mmod^{\on{fact}}_\CZ}(\bC_2,\Res_{\Phi^R}(\bC_1)).$$
\end{cor} 
 
This proposition can be reformulated as saying that the functor $\Res_\Phi$
is the \emph{left} adjoint of the functor $\Res_{\Phi^R}$. 
 
\begin{proof}

This is a formal corollary of having a 2-Cartesian fibration. To be explicit, let us exhibit
the unit and the counit of the adjunction.

\medskip

The unit is given by
$$\on{Id}=\Res_{\on{Id}} \overset{\Phi\circ \Phi^R\to \on{Id}}\longrightarrow  
\Res_{\Phi\circ \Phi^R} \simeq \Res_{\Phi^R}\circ \Res_\Phi.$$

The counit is given by 
$$\Res_\Phi\circ \Res_{\Phi^R}\simeq \Res_{\Phi^R\circ \Phi} 
\overset{\on{Id}\to \Phi^R\circ \Phi}\longrightarrow  \Res_{\on{Id}}=\on{Id}.$$

\end{proof} 

\sssec{}

We will now use \corref{c:basic adj non-unital} to prove the following partial converse to \lemref{l:fact restr cat}:

\begin{lem} \label{l:fact res crit}
Let $\Phi:\bA_1\to \bA_2$ be a factorization functor, and let 
$\bC_1\to \bC_2$ be a functor between objects of $\bA_i\mmod^{\on{fact}}_\CZ$, $i=1,2$, 
compatible with factorization. Assume that:

\smallskip

\noindent{\em(i)} The functor $\Phi:\ul\bA_1\to \ul\bA_2$ between sheaves of categories on $\Ran$ admits
a right adjoint;

\smallskip

\noindent{\em(ii)} The functor $\Phi_m:\ul\bC_1\to \ul\bC_2$ between sheaves of categories on $\CZ^{\subseteq}$ admits
a right adjoint;

\smallskip

\noindent{\em(iii)} The induced functor $\Phi_m:\bC_{1,\CZ}\to \bC_{2,\CZ}$ is an equivalence.

\medskip

Then the resulting functor
\begin{equation} \label{e:fact res crit}
\bC_1\to \on{Res}_\Phi(\bC_2)
\end{equation}
as module categories over $\bA_1$, is an equivalence. 

\end{lem} 

\begin{proof}

We claim that the functor \eqref{e:fact res crit} admits a right adjoint. Once we prove this, 
the lemma will follow, because a functor between sheaves of categories that admits a right adjoint
is an equivalence if and only if it is an equivalence strata-wise.

\medskip

The right adjoint of $\Phi_m:\ul\bC_1\to \ul\bC_2$ is a functor compatible with factorization against 
$\Phi^R:\ul\bA_1\to \ul\bA_2$. Hence, it gives rise to an object 
$$\Phi_m^R\in \on{Funct}_{\Phi^R:\bA_2\to \bA_1}(\bC_2,\bC_1)\simeq 
\on{Funct}_{\bA_2\mmod^{\on{fact}}_\CZ}(\bC_2,\Res_{\Phi^R}(\bC_1)).$$

(In fact, the pair $(\Phi^R,\Phi_m^R)$ is the right adjoint of $(\Phi,\Phi_m)$ as a 
1-morphism in $\on{FactCat-and-Mod}(X)$.) 

\medskip 

Using \corref{c:basic adj non-unital}, we identify the latter category with
$$\on{Funct}_{\bA_1\mmod^{\on{fact}}_\CZ}(\Res_\Phi(\bC_2),\bC_1).$$

Unwinding the constructions, we obtain that the resulting functor
$$\on{Res}_\Phi(\bC_2)\to \bC_1$$
is indeed the right adjoint of \eqref{e:fact res crit}.

\end{proof} 

\section{Unital structures} \label{s:unit}

In the previous section we introduced the notion of factorization category. In this section we will
describe an extra structure that factorization categories often carry: the \emph{unital} structure.

\medskip

The role that the unital structure plays is two-fold. For one thing, it enables various local-to-global
constructions (see \secref{s:unitality}). But it also leads to purely local constructions, which play
a key role in this paper: given a lax-unital factorization functor $F:\bA_1\to \bA_2$ between
unital factorization categories, the image $F(\one_{\bA_1})$ of the unite $\one_{\bA_1}\in \bA_2$
is a factorization algebra in $\bA_2$, and the functor $F$ \emph{enhances} to a functor
$$F^{\on{enh}}:\bA_1\to F(\one_{\bA_1})\mod^{\on{fact}}(\bA_2).$$

Now the functor is often, if not an equivalence, but is close to be such\footnote{E.g., in multiple instances, 
in the presence of a t-structure, it is an equivalence on the bounded below categories.}, and that allows
to understand the more complicated category $\bA_1$ in terms of $\bA_1$. 

\medskip

In order to talk about unitality we have to enlarge our world of algebro-geometric objects. Namely,
we normally work with prestacks (i.e., spaces in algebro-geometric sense, whose functor of points takes
place in ($\infty$)-groupoids). But in in order to talk about unitality, we need to work with 
\emph{categorical prestacks}; whose functors of points take place in ($\infty$)-categories. The rudiments
of categorical prestacks (D-modules and sheaves of categories on them) are also developed in this section. 

\ssec{Categorical prestacks}

\sssec{}

When discussing categorical prestacks, we will work in the locally almost of finite type (laft)
category. Accordingly, when we write $\affSch$ (resp., $\on{PreStk}$), we will mean affine 
schemes (resp., prestacks) locally almost of finite type.  

\medskip

By a categorical prestack we shall mean a functor
$$(\affSch)^{\on{op}}\to \inftyCat.$$

Given a categorical prestack $\CY$ we will denote by
$$\CY(S) \text{ or } \Maps(S,\CY)$$
the category of its values on $S\in \affSch$. 

\medskip

We let $\on{CatPreStk}$ denote the category of categorical prestacks. 

\sssec{}

Let $f:\CY_1\to \CY_2$ be a map between categorical prestacks. We shall say that some categorical property of 
$f$ (such as being Cartesian/co-Cartesian, cofinal, admitting an adjoint) holds \emph{value-wise} if this property holds
for the corresponding functor
$$\CY_1(S)\to \CY_2(S)$$
for every $S\in \affSch$.

\sssec{}

Let $\CY$ be a categorical prestack. To it we will associate two prestacks in groupoids,
by applying the right and left adjoints to the embedding
$$\inftygroup\hookrightarrow  \inftyCat,$$
respectively. 

\medskip

We let $\CY^{\on{grpd}}$ denote the prestack whose value on $S\in \affSch$ is the groupoid
underlying the category $\CY(S)$.

\medskip

We will denote by $\sft$ the tautological map
$$\CY^{\on{grpd}} \to \CY.$$

\medskip

We will denote by $\CY^{\on{strict}}$ the prestack in groupoids, whose value on $S\in \affSch$
is the groupoid obtained from $\CY(S)$ by inverting all 1-morphisms. 

\medskip

We will denote the tautological projection by
$$\on{strict}:\CY\to \CY^{\on{strict}}.$$ 

\sssec{} \label{sss:Dmods on cat prestack}

Given a categorical prestack $\CY$, we associate to it several DG categories of algebro-geometric nature. 

\medskip

Let $\CY_{\affSch}$ denote the Cartesian fibration over $\affSch$ that attaches to $S\in \affSch$
the category $\CY(S)$. 

\medskip

Let 
\begin{equation} \label{e:alg geom Cart}
\QCoh_{\affSch},\,\, \IndCoh_{\affSch} \text{ and } \Dmod_{\affSch}$$
denote the Cartesian fibrations over $\affSch$ that attach to $S$ the categories
$$\QCoh(S),\,\, \IndCoh(S) \text{ and }  \Dmod(S),
\end{equation}
respectively. 

\medskip

We define the categories 
\begin{equation} \label{e:alg geom Cat}
\QCoh(\CY),\,\,\IndCoh(\CY) \text{ and } \Dmod(\CY)
\end{equation} 
as functors from $\CY_{\affSch}$ to the categories in \eqref{e:alg geom Cat} that map 
arrows that are Cartesian over $\affSch$ to arrows with a similar property. See \cite[Sect. C.3]{Ro2}
for more details. 

\sssec{} \label{sss:pullback cat prestack}

For a map of prestacks $f:\CY_1\to \CY_2$, precomposition with
$$(\CY_1)_{\affSch}\to (\CY_2)_{\affSch}$$
gives rise to a functor
$$f^!:\Dmod(\CY_2)\to \Dmod(\CY_1),$$
and similarly for $\IndCoh(\CY)$ and $\Dmod(\CY)$. 

\sssec{}

Categorical prestacks form an $(\infty,2)$-category, so there is a natural notion of adjunction between morphisms.
Explicitly, a morphism 
$$f:\CY_1\rightleftarrows \CY_2$$
admits a right adjoint, if for every $S\in \affSch$, the corresponding functor
$$f:\CY_1(S)\to \CY_2(S)$$
admits a right adjoint, and for every $S'\to S$, the natural transformation
$$
\xy  
(0,20)*+{\CY_1(S')}="B";
(0,0)*+{\CY_1(S)}="A";
(20,20)*+{\CY_2(S')}="D";
(20,0)*+{\CY_2(S)}="C";
{\ar@{->} "A";"B"};
{\ar@{->} "C";"D"};
{\ar@{<-}^{f^R} "A";"C"};
{\ar@{<-}^{f^R} "B";"D"};
{\ar@{=>} "A";"D"};
\endxy
$$
arising by adjunction from the commutative diagram
$$
\xy
(0,20)*+{\CY_1(S')}="B";
(0,0)*+{\CY_1(S)}="A";
(20,20)*+{\CY_2(S')}="D";
(20,0)*+{\CY_2(S),}="C";
{\ar@{->} "A";"B"};
{\ar@{->} "C";"D"};
{\ar@{->}^{f} "A";"C"};
{\ar@{->}^{f} "B";"D"};
\endxy
$$
is an isomorphism. 

\medskip

We have the following useful observation:

\begin{lem} \label{l:adj pullback}
Let
$$f:\CY_1\rightleftarrows \CY_2:g$$ be 
mutually adjoint maps. Then the functors $(g^!,f^!)$ form an adjoint pair.
\end{lem}

\sssec{}

For a categorical prestack $\CY$, we let $\CY_\dR$ denote the categorical prestack defined by
$$\CY_\dr(S):=\CY(S_{\on{red}}).$$

A standard manipulation shows that 
$$\IndCoh(\CY_\dr)\simeq \Dmod(\CY).$$

\sssec{} \label{sss:expl Dmod cat prestack}

One can describe the categories \eqref{e:alg geom Cat} explicitly as follows. We will
do this for $\Dmod(\CY)$, while the other two cases are similar. 

\medskip

An object $\CF\in \Dmod(\CY)$ is an assignment:

\begin{itemize}

\item For every $y:S\to \CY$ of an object $\CF_{S,y}\in \Dmod(S)$;

\item For a map $(y_1\overset{\alpha}\to y_2)\in \CY(S)$ of a map $\CF_{S,y_1}\overset{\CF_\alpha}\to \CF_{S,y_2}$ in $\Dmod(S)$;

\item For $f:S'\to S$ and $y'=y\circ f$ of an isomorphism $\CF_{S',y'}\simeq f^!(\CF_{S,y})$ in $\Dmod(S')$;

\item The datum of commutativity for the diagram
$$
\CD
\CF_{S',y'_1} @>{\CF_{f^*(\alpha)}}>> \CF_{S',y'_2} \\
@V{\simeq}VV @VV{\simeq}V \\
f^!(\CF_{S,y_1}) @>{f^!(\CF_\alpha)}>> f^!(\CF_{S,y_2});
\endCD
$$

\item A homotopy-coherent system of compatibilities for compositions. 

\end{itemize} 

\sssec{} \label{sss:strict Dmods on categ prestack}

In addition to the categories \eqref{e:alg geom Cat}, one can consider their strict versions:
$$\QCoh(\CY)^{\on{strict}}:=\QCoh(\CY^{\on{strict}}),\,\,\IndCoh(\CY)^{\on{strict}}:=\IndCoh(\CY^{\on{strict}})$$
and 
$$\Dmod(\CY)^{\on{strict}}:= \Dmod(\CY^{\on{strict}}),$$
respectively. 

\medskip

Unwinding the definitions, we obtain that pullback along
$$\on{strict}:\CY\to \CY^{\on{strict}}$$
is a fully faithful embedding, with essential image described as follows:

\medskip

It consists of those objects $\CF$ in \secref{sss:expl Dmod cat prestack}, for which the maps
$$\CF_{S,y_1}\overset{\CF_\alpha}\to \CF_{S,y_2}$$
are isomorphisms.

\medskip

Note that we can describe $\Dmod(\CY)^{\on{strict}}$ also as 
$$\Dmod(\CY)^{\on{strict}}=\underset{S\in \affSch_{/\CY}}{\on{lim}}\, \Dmod(S).$$

\medskip

The same applies also to $\QCoh(\CY)^{\on{strict}}$ and $\IndCoh(\CY)^{\on{strict}}$.

\ssec{Crystals of categories on categorical prestacks} \label{sss:shvs-of-cats categ}

\sssec{}

Let $\CY$ be a categorical prestack. Combining the ideas of Sects. \ref{sss:shvs of cat} and \ref{sss:Dmods on cat prestack}
we obtain the notion of \emph{sheaf of categories} over $\CY$. 

\medskip 

Thus, a sheaf of categories $\ul\bC$ over $\CY$ is an assignment:

\medskip

\begin{itemize}

\item For every affine scheme $S$ and a map $y:S\to \CY$ of a category $\bC_{S,y}$ tensored over $\QCoh(S)$;

\item For a map $y_1\overset{\alpha}\to y_2$ in $\Maps(S,\CY)$ of a $\Dmod(S)$-linear functor
$\bC_{S,y_1}\overset{\bC_{S,\alpha}}\longrightarrow \bC_{S,y_2}$;

\item For $S'\overset{f}\to S$ and $y'=y\circ f$ of an identification $\bC_{S',y'}\simeq \QCoh(S')\underset{\QCoh(S)}\otimes \bC_{S,y}$;

\item For $\alpha'=\alpha \circ f$ of a datum of commutativity for
$$
\CD
\bC_{S,y_1} @>{\bC_{S,\alpha}}>> \bC_{S,y_2}  \\
@V{f^!}VV @VV{f^!}V \\
\bC_{S',y'_1} @>{\bC_{S',\alpha'}}>> \bC_{S',y'_2}.
\endCD
$$

\item A homotopy-coherent system of compatibilities for compositions.

\end{itemize}

\sssec{}

We shall say that $\ul\bC$ is strict if the functors $\bC_{S,\alpha}$ are equivalences. Note that $\ul\bC$ is strict if 
it is the pullback\footnote{In the sense of \secref{sss:pullback sheaf of cat}.}
along
$$\CY\to \CY^{\on{strict}}$$
of a crystal of categories on $\CY^{\on{strict}}$. 

\sssec{}

One can assign to $\ul\bC$ two categories, denoted
$$\Gamma^{\on{lax}}(\CY,\ul\bC) \text{ and } \Gamma^{\on{strict}}(\CY,\ul\bC),$$
respectively, defined as follows, the latter being a full subcategory of the former. 

\medskip

An object of $\Gamma^{\on{lax}}(\CY,\ul\bC)$ assigns to every affine scheme $S$ and a map $y:S\to \CY$ an 
object $\bc_{S,y}\in \bC_{S,y}$ together with the following data:

\begin{itemize}

\item For a map $y_1\overset{\alpha}\to y_2$ in $\Maps(S,\CY)$ a morphism 
\begin{equation} \label{e:transition map lax}
\bC_{S,\alpha}(\bc_{S,y_1})\to \bc_{S,y_2};
\end{equation}

\smallskip

\item For a map $f:S'\to S$ and $y'=y\circ f$ an isomorphism
$$f^!(\bc_{S,y})\simeq \bc_{S',y'}$$
as objects in $\bC_{S',y'}$.

\smallskip

\item A homotopy-coherent datum of compatibility for the above pieces of data.

\end{itemize}

The subcategory $\Gamma^{\on{strict}}(\CY,\ul\bC)$ consists of those assignments for 
which the maps \eqref{e:transition map lax} are isomorphisms.

\begin{rem}

We alert the reader to the discrepancy between the notations 
$$\Gamma^{\on{lax}}(-,-) \text{ and } \Gamma^{\on{strict}}(-,-)$$
introduced above and those used in \cite[Sect. 4]{Ra6}. 

\medskip

Namely, what we denote $\Gamma^{\on{lax}}(-,-)$ is denoted $\Gamma(-,-)$ in 
{\it loc. cit.}, and what we denote $\Gamma^{\on{strict}}(-,-)$ is denoted $\Gamma^{\on{naive}}(-,-)$ in 
{\it loc. cit.}. 

\medskip 

Similarly, the notion of \emph{functor} of sheaves of categories
considered in \cite{Ra6} corresponds to the notion of \emph{right-lax functor}
considered below.

\end{rem} 

\sssec{Example}

Let $\ul\bC$ be $\ul\QCoh(\CY)$, the unit crystal of categories, i.e., its value for $(S,y)\in \affSch_{/\CY}$ is
$\QCoh(S)$.

\medskip

Then 
$$\Gamma^{\on{lax}}(\CY,\ul\bC) =\QCoh(\CY),$$
see \secref{sss:Dmods on cat prestack}
and 
$$\Gamma^{\on{strict}}(\CY,\ul\bC) =\QCoh(\CY)^{\on{strict}}\simeq \QCoh(\CY^{\on{strict}})$$
(see \secref{sss:strict Dmods on categ prestack}).

\sssec{}

Let $\ul\bC'$ and $\ul\bC''$ be two crystals of categories on $\CY$. In this case there is an (evident) notion of functor
$$\ul\Phi:\ul\bC'\to \ul\bC''.$$

When an ambiguity is likely to occur, we will call such functors \emph{strict}.

\sssec{} \label{sss:ff strict functors}

For future reference, a (strict) functor $\ul\Phi$ is said to be fully faithful if for every $y:S\to \CY$, 
the resulting functor
$$\bC'_{S,y}\to \bC''_{S,y}$$
is fully faithful. 

\sssec{} \label{sss:lax vs strict functors}

In addition, there is a notion of \emph{right-lax} functor. A right-lax functor $\ul\Phi:\ul\bC'\to \ul\bC''$
is an assignment:

\begin{itemize}

\item For every $(S,y)$ of a functor $\bC'_{S,y} \overset{\Phi_{S,y}}\to \bC''_{S,y}$;

\item For every map $y_1\overset{\alpha}\to y_2$ in $\Maps(S,\CY)$ we have a \emph{natural transformation}
\begin{equation} \label{e:nat transf lax functor}
\bC''_{S,\alpha}\circ \Phi_{S,y_1}\to \Phi_{S,y_2}\circ \bC'_{S,\alpha}.
\end{equation}

\item For $f:\wt{S}\to S$ and $\wt{y}=y\circ f$ of an isomorphism $f^!\circ \Phi_{S,y}\simeq \Phi_{\wt{S},\wt{y}}\circ f^!$;

\smallskip

\item A homotopy-coherent system of compatibilities for the above data.

\end{itemize} 

\bigskip

By definition, a right-lax functor is \emph{strict} if the natural transformations \eqref{e:nat transf lax functor}
are isomorphisms. 

\medskip

We denote the categories of right-lax and strict functors by
$$\on{Funct}^{\on{lax}}_{\on{CrystCat}(\CY)}(\ul\bC',\ul\bC'')  \text{ and }
\on{Funct}^{\on{strict}}_{\on{CrystCat}(\CY)}(\ul\bC',\ul\bC''),$$
respectively.

\sssec{}

We will denote the (2-)category of sheaves of categories on $\CY$, with 1-morphisms being strict functors by 
$\on{ShvCat}^{\on{strict}}(\CY)$.

\medskip

We will denote the (2-)category of sheaves of categories on $\CY$, with 1-morphisms being right-lax functors by 
$\on{ShvCat}^{\on{lax}}(\CY)$.

\medskip

Sometimes we will simply write $\on{ShvCat}(\CY)$, when the discussion is applicable in both contexts. 

\medskip

Both $\on{ShvCat}^{\on{strict}}(\CY)$ and $\on{ShvCat}^{\on{lax}}(\CY)$ carry a natural symmetric monoidal structure
with the unit being $\ul\QCoh(\CY)$. 

\sssec{} \label{sss:cryst term}

We set
$$\on{CrystCat}^{\on{strict}}(\CY):=\on{ShvCat}^{\on{strict}}(\CY_\dr) \text{ and } \on{CrystCat}^{\on{lax}}(\CY):=\on{ShvCat}^{\on{lax}}(\CY_\dr)$$

\medskip

Terminologically, when we talk about $\ul\bC$ being a \emph{crystal of categories} over $\CY$, for $(S,y)\in \affSch_{/\CY}$, we will
denote by 
$$\bC_{S,y}\in \Dmod(S)\mmod$$ the corresponding category of crystalline sections. 

\medskip

We let $\ul\Dmod(\CY)$ denote the unit crystal of categories on $\CY$, i.e., its value for $(S,y)\in \affSch_{/\CY}$ is
$\Dmod(S)$.

\sssec{} \label{sss:right adj lax}

Let $\ul\bC_1$ and $\ul\bC_2$ be two crystals of categories on $\CY$, and let
$$\ul\Phi:\ul\bC_1\to \ul\bC_2$$
be a strict functor. 

\medskip

Assume that the induced functor
$$\sft^!(\ul\bC_1) \overset{\sft^!(\ul\Phi)}\to \sft^!(\ul\bC_2)$$
admits a right adjoint\footnote{As a functor between sheaves of categories on $\CY^{\on{grpd}}$, i.e., it admits a continuous right adjoint value-wise.}, 
to be denoted $(\sft^!(\ul\Phi))^R$. 

\medskip

In this case $(\sft^!(\ul\Phi))^R$ admits a natural extension, to be denoted $\ul\Phi^R$ to
a right-lax functor 
$$\ul\bC_2\to \ul\bC_1,$$
see, e.g., \cite[Lemma B.5.9]{AMR}.

\sssec{} \label{sss:dual cryst cat categ}

Let $\ul\bC$ be a crystal of categories over $\CY$. Assume that $\sft^!(\ul\bC)$ is dualizable. Assume moreover
that for every $(y_1\overset{\alpha}\to y_2)\in \Maps(S,\CY)$, the functor
\begin{equation} \label{e:arrow alpha again}
\bC_{S,y_1}\overset{\bC_{S,\alpha}}\to \bC_{S,y_2}
\end{equation} 
admits a right adjoint. 

\medskip

In this case, we can extend the dual $(\sft^!(\ul\bC))^\vee$ to a crystal of categories $\ul\bC^\vee$
over $\CY$ by letting
$$\bC^\vee_{S,y_1}\overset{\bC^\vee_{S,\alpha}}\to \bC^\vee_{S,y_2}$$
be the dual of the right adjoint of \eqref{e:arrow alpha again}. 

\medskip

Under the above circumstances we will say that $\ul\bC$ is dualizable, and we will refer to the above crystal of categories $\ul\bC^\vee$
as the dual of $\ul\bC$. 

\medskip

Note that we have the natural evaluation and coevaluation that are \emph{right-lax} functors
\begin{equation} \label{e:eval coeval cat}
\ul\Dmod(\CY)\overset{\ul{\on{co-eval}}}\to  \ul\bC\otimes \ul\bD \text{ and } \ul\bD\otimes \ul\bC \overset{\ul{\on{eval}}}\to  \ul\Dmod(\CY),
\end{equation} 
i.e., the duality between $\ul\bC$ and $\ul\bC^\vee$ takes place in the symmetric monoidal category $\on{ShvCat}^{\on{lax}}(\CY)$. 

\sssec{} \label{sss:dual cryst cat categ 2}

Vice versa, let us be given two crystal of categories $\ul\bC$ and $\ul\bD$ and right-lax functors as in \eqref{e:eval coeval cat}. 
Suppose that the following conditions hold:

\begin{itemize}

\item For every $S\overset{y}\to \CY$, the functors
$$\Dmod(S)\to \bC_{S,y}\underset{\Dmod(S)}\otimes \bD_{S,y} \text{ and } \bD_{S,y}\underset{\Dmod(S)}\otimes \bC_{S,y}\to \Dmod(S)$$
define a perfect pairing;

\smallskip

\item The identification of the pullback of the composition 
$$\ul\bC \overset{\on{co-eval}\otimes \on{Id}}\longrightarrow 
\ul\bC \otimes \ul\bD\otimes \ul\bC\overset{\on{Id}\otimes \on{eval}}\longrightarrow \ul\bC$$
along $\sft$ with the identity endofunctor of $\ul\bC|_{\CY^{\on{grpd}}}$ extends to $\CY$;

\smallskip

\item The identification of the pullback of the composition 
$$\ul\bD \overset{\on{Id}\otimes \on{eval}}\longrightarrow 
\ul\bD \otimes \ul\bC\otimes \ul\bD\overset{\on{co-eval}\otimes \on{Id}}\longrightarrow \ul\bD$$
along $\sft$ with the identity functor of $\ul\bD|_{\CY^{\on{grpd}}}$ extends to $\CY$.

\end{itemize} 
 
Then $\ul\bC$ is dualizable, and $\ul\bD$ identifies canonically with the dual of $\ul\bC$ 
in the sense of \secref{sss:dual cryst cat categ}, see \cite[Remark 11.11.17]{CF}. 

\sssec{} \label{sss:dual cryst cat categ 3}

Still equivalently, let us be given a pair of crystals of categories $\ul\bC$ and $\ul\bD$ and \emph{either} a right-lax functor
$$\ul\bD\otimes \ul\bC \overset{\ul{\on{eval}}}\to \ul\Dmod(\CY)$$
\emph{or}
a right-lax functor
$$\ul\Dmod(\CY)\overset{\ul{\on{co-eval}}}\to \ul\bC\otimes \ul\bD.$$

Suppose that:

\begin{itemize}

\item The pullback of $\ul{\on{eval}}$ (resp., $\ul{\on{co-eval}}$) along $\sft^!$ is a perfect pairing;

\item For every $(y_1\overset{\alpha}\to y_2)\in \Maps(S,\CY)$, the resulting natural transformation
$\on{Id}\to \bD_{S,\alpha}^\vee\circ \bC_{S,\alpha}$ (resp., $\bC_{S,\alpha}\circ \bD_{S,\alpha}^\vee\to \on{Id}$)
is the unit (resp., counit) of an adjunction.

\end{itemize}

Then this datum extends uniquely to a datum of duality between $\ul\bC$ and $\ul\bD$ as crystals
of categories on $\CY$.

\sssec{} \label{sss:functoriality of lax global sections}

For $\ul\bC'$ and $\ul\bC''$ as above, we observe that a right-lax functor $\ul\Phi:\ul\bC'\to \ul\bC''$
gives rise to a functor
$$\Phi:\Gamma^{\on{lax}}(\CY,\ul\bC') \to \Gamma^{\on{lax}}(\CY,\ul\bC'').$$

If $\ul\Phi$ is strict, then $\Phi$ induces a functor
$$\Gamma^{\on{strict}}(\CY,\ul\bC') \to \Gamma^{\on{strict}}(\CY,\ul\bC'').$$

\sssec{}

Note that the functors
$$\ul\bC\mapsto \Gamma^{\on{lax}}(\CY,\ul\bC) \text{ and } \ul\bC\mapsto \Gamma^{\on{strict}}(\CY,\ul\bC)$$
can be recovered as adjoints:

\medskip

For $\bD\in \DGCat$, we have
$$\on{Funct}_{\on{cont}}(\bD,\Gamma^{\on{lax}}(\CY,\ul\bC))\simeq 
\on{Funct}^{\on{lax}}_{\on{CrystCat}(\CY)}(\bD\otimes \ul\Dmod(\CY),\ul\bC)$$  and 
$$\on{Funct}_{\on{cont}}(\bD,\Gamma^{\on{strict}}(\CY,\ul\bC))\simeq 
\on{Funct}^{\on{strict}}_{\on{CrystCat}(\CY)}(\bD\otimes \ul\Dmod(\CY),\ul\bC),$$
respectively. 

\sssec{}  \label{sss:pullback sheaf of cat}

Let $f:\CY_1\to \CY_2$ be a map between categorical prestacks. As in Sects. \ref{sss:pullback shf of cat} and 
\ref{sss:pullback cat prestack}, there is a naturally
defined functor
$$f^*:\on{CrystCat} (\CY_1)\to \on{CrystCat} (\CY_2).$$

\sssec{} \label{sss:pullback sect sheaf of cat}

For a map $f:\CY_1\to \CY_2$ between categorical prestacks and a crystal of categories  $\ul\bC$ on $\CY_2$, 
we have a naturally defined functor
\begin{equation} \label{e:pullback on lax sections}
f^!:\Gamma^{\on{lax}}(\CY_2,\ul\bC)\to \Gamma^{\on{lax}}(\CY_1,f^*(\ul\bC)),
\end{equation}
which induces a functor
$$f^!:\Gamma^{\on{strict}}(\CY_2,\ul\bC)\to \Gamma^{\on{strict}}(\CY_1,f^*(\ul\bC)).$$

\ssec{Two notions of direct image of a crystal of categories } \label{ss:dir mage sheaf of cat}

\sssec{} \label{sss:Cart prestacks}

Let $f:\CY_1\to \CY_2$ be a map between categorical prestacks. We shall say that $f$ is 
a (co)Cartesian fibration, i.e.:

\begin{itemize}

\item $f$ is a value-wise (co)Cartesian fibration;

\item For $S'\to S$, the functor $\CY_1(S)\to \CY_1(S')$ sends arrows that are (co)Cartesian over $\CY_2(S)$
to arrows in $\CY_1(S')$ with the same property.

\end{itemize} 

\medskip 

Note that the second condition is automatic if $f$ is a value-wise (co)Cartesian fibration \emph{in groupoids}. 

\sssec{}

Let $f:\CY_1\to \CY_2$ be a Cartesian fibration.

\medskip

Let $\ul\bC_1$ be a crystal of categories  over $\CY_1$. In this case one can form two sheaves of categories,
denoted
$$f_*(\ul\bC_1) \text{ and } f_{*,\on{lax}}(\ul\bC_1)$$
on $\CY_2$, as follows. 

\sssec{}

For $y:S\to \CY_2$, the value of $f_{*,\on{lax}}(\ul\bC_1)$ on $(S,y_2)$ is
$$\Gamma^{\on{lax}}(S\underset{\CY_2}\times \CY_1,\ul\bC_1|_{S\underset{\CY_2}\times \CY_1})$$
and the value of $f_{*,\on{strict}}(\ul\bC_1)$ is 
$$\Gamma^{\on{strict}}(S\underset{\CY_2}\times \CY_1,\ul\bC_1|_{S\underset{\CY_2}\times \CY_1}).$$

\medskip

The data of crystal of categories  on $f_{*,\on{lax}}(\ul\bC_1)$ is defined as follows.

\sssec{}

For a map $y'\overset{\alpha}\to y''$ in $\Maps(S,\CY_2)$ we have a map 
$$S\underset{y'',\CY_2}\times \CY_1\overset{\alpha_{\CY_2}^*}\to S\underset{y',\CY_2}\times \CY_1.$$

The structure on $\ul\bC_1$ of crystal of categories  gives rise to a (strict) functor between crystals
of categories on $S\underset{y'',\CY_2}\times \CY_1$
$$(\alpha_{\CY_2}^*)^*(\ul\bC_1|_{S\underset{y',\CY_2}\times \CY_1})\to 
\ul\bC_1|_{S\underset{y'',\CY_2}\times \CY_1}.$$

\medskip

Hence, the constructions of Sects. \ref{sss:functoriality of lax global sections} and 
\ref{sss:pullback sect sheaf of cat} combine to gives rise to a functor
\begin{multline} \label{e:moving between fibers}
f_{*,\on{lax}}(\ul\bC_1)_{S,y'}:=
\Gamma^{\on{lax}}(S\underset{y',\CY_2}\times \CY_1,\ul\bC_1|_{S\underset{\CY_2}\times \CY_1})
\overset{(\alpha_{\CY_2}^*)^!}\longrightarrow \\
\to
\Gamma^{\on{lax}}(S\underset{y'',\CY_2}\times \CY_1,(\alpha_{\CY_2}^*)^*(\ul\bC_1|_{S\underset{y',\CY_2}\times \CY_1}))\to 
\Gamma^{\on{lax}}(S\underset{y'',\CY_2}\times \CY_1,\ul\bC_1|_{S\underset{\CY_2}\times \CY_1})=:f_{*,\on{lax}}(\ul\bC_1)_{S,y''}.
\end{multline}

This functor induces a functor
\begin{multline*}
f_{*,\on{strict}}(\ul\bC_1)_{S,y'}:=
\Gamma^{\on{strict}}(S\underset{y',\CY_2}\times \CY_1,\ul\bC_1|_{S\underset{\CY_2}\times \CY_1})
\overset{(\alpha_{\CY_2}^*)^!}\longrightarrow \\
\to
\Gamma^{\on{strict}}(S\underset{y'',\CY_2}\times \CY_1,(\alpha_{\CY_2}^*)^*(\ul\bC_1|_{S\underset{y',\CY_2}\times \CY_1}))\to 
\Gamma^{\on{strict}}(S\underset{y'',\CY_2}\times \CY_1,\ul\bC_1|_{S\underset{\CY_2}\times \CY_1})=:f_{*,\on{strict}}(\ul\bC_1)_{S,y''}.
\end{multline*}

\sssec{}

For $g:\wt{S}\to S$ and $\wt{y}=f\circ g$, we have a map
$$g:\wt{S}\underset{\CY_2}\times \CY_1\to S\underset{\CY_2}\times \CY_1,$$
and the construction of  \secref{sss:pullback sect sheaf of cat} gives rise to a functor
$$f_{*,\on{lax}}(\ul\bC_1)_{S,y}:=
\Gamma^{\on{lax}}(S\underset{\CY_2}\times \CY_1,\ul\bC_1|_{S\underset{\CY_2}\times \CY_1})
\overset{g^!}\to \Gamma^{\on{lax}}(\wt{S}\underset{\CY_2}\times \CY_1,\ul\bC_1|_{\wt{S}\underset{\CY_2}\times \CY_1})
=:f_{*,\on{lax}}(\ul\bC_1)_{\wt{S},\wt{y}}.$$

This functor induces a functor
$$f_{*,\on{strict}}(\ul\bC_1)_{S,y}:=
\Gamma^{\on{strict}}(S\underset{\CY_2}\times \CY_1,\ul\bC_1|_{S\underset{\CY_2}\times \CY_1})
\overset{g^!}\to \Gamma^{\on{strict}}(\wt{S}\underset{\CY_2}\times \CY_1,\ul\bC_1|_{\wt{S}\underset{\CY_2}\times \CY_1})
=:f_{*,\on{strict}}(\ul\bC_1)_{\wt{S},\wt{y}}.$$

\sssec{}

By construction, we have a strict functor
$$f_{*,\on{strict}}(\ul\bC_1)\to f_{*,\on{lax}}(\ul\bC_1),$$
which is a value-wise fully faithful embedding.

\sssec{} \label{sss:pushforward sheaf of cat funct}

The above construction is functorial in the following sense: for a lax functor $\ul\Phi:\ul\bC'_1\to \ul\bC''_1$
between crystals of categories on $\CY_1$ we obtain a lax functor
$$f_*(\ul\Phi):f_{*,\on{lax}}(\ul\bC'_1)\to f_{*,\on{lax}}(\ul\bC''_1).$$

If $\ul\Phi$ is strict, then so is $f_*(\ul\Phi)$, and it also induces a strict functor
$$f_*(\Phi):f_{*,\on{strict}}(\ul\bC'_1)\to f_{*,\on{strict}}(\ul\bC''_1).$$

\sssec{} \label{sss:push-pull categorical prestacks}

The operations
$$\ul\bC_1\mapsto f_{*,\on{strict}}(\ul\bC_1) \text{ and } \ul\bC_1\mapsto f_{*,\on{lax}}(\ul\bC_1)$$
can also be realized as right adjoints. Namely, for $\ul\bC_2\in \on{CrystCat}(\CY_2)$, we have
$$\on{Funct}^{\on{strict}}_{\on{CrystCat}(\CY_1)}(f^*(\ul\bC_2),\ul\bC_1)\simeq
\on{Funct}^{\on{strict}}_{\on{CrystCat}(\CY_2)}(\ul\bC_2,f_{*,\on{strict}}(\ul\bC_1))$$
and 
$$\on{Funct}^{\on{lax}}_{\on{CrystCat}(\CY_1)}(f^*(\ul\bC_2),\ul\bC_1)\simeq
\on{Funct}^{\on{lax}}_{\on{CrystCat}(\CY_2)}(\ul\bC_2,f_{*,\on{lax}}(\ul\bC_1)).$$

\medskip

In particular, we have a canonically defined (strict) functor 
\begin{equation} \label{e:push-pull categorical prestacks}
\ul\bC_2\to f_{*,\on{strict}}\circ f^*(\ul\bC_2).
\end{equation}

\sssec{}

Consider the forgeful functor
$$\on{CrystCat}^{\on{strict}}(\CY)\to \on{CrystCat}^{\on{lax}}(\CY).$$

We claim that it admits a right adjoint. Namely, consider the map
$$\on{pr}_{\on{source}}:\CY^\to \to \CY$$
is a Cartesian fibration.

\medskip

The above right adjoint is given by
$$(\on{pr}_{\on{source}})_{*,\on{lax}}\circ (\on{pr}_{\on{target}})^*.$$

In particular, for $\ul\bC_1,\ul\bC_2\in \on{CrystCat}^{\on{strict}}(\CY)$
we have a canonical identification
\begin{equation} \label{e:lax vs strict abs}
\on{Funct}^{\on{lax}}(\ul\bC_1,\ul\bC_2)\simeq 
\on{Funct}^{\on{strict}}(\ul\bC_1,(\on{pr}_{\on{source}})_{*,\on{lax}}\circ (\on{pr}_{\on{target}})^*(\bC_2)),
\end{equation}
with the inclusion
$$\on{Funct}^{\on{lax}}(\ul\bC_1,\ul\bC_2)\hookrightarrow \on{Funct}^{\on{strict}}(\ul\bC_1,\ul\bC_2)$$
corresponding to the strict functor
$$\ul\bC_2 \to (\on{pr}_{\on{source}})_{*,\on{strict}}\circ (\on{pr}_{\on{target}})^*(\bC_2)\to
(\on{pr}_{\on{source}})_{*,\on{lax}}\circ (\on{pr}_{\on{stirct}})^*(\bC_2),$$
where the first arrow corresponds by adjunction to the functor
$$(\on{pr}_{\on{source}})^*(\ul\bC_2)\to (\on{pr}_{\on{target}})^*(\ul\bC_2),$$
of \eqref{e:source to target}. 

\sssec{} \label{sss:left adj lax to strict}

We now explain an abstract framework for the construction in \secref{ss:constr int}. 

\medskip

Let $\CY$ be a categorical prestack,
and let $\ul\bC_1,\ul\bC_2$ be crystals of categories on it. 

\medskip

Suppose that the functor
\begin{equation} \label{e:pullback untl from Ran to arr abs}
\ul\bC_2\to (\on{pr}_{\on{souce}})_{*,\on{strict}}(\bC_2)\to (\on{pr}_{\on{souce}})_{*,\on{lax}}(\bC_2)
\end{equation}
admits a left adjoint that is a strict functor. 

\medskip

Then the functor
$$\on{Funct}^{\on{strict}}(\ul\bC_1,\ul\bC_2)\to \on{Funct}^{\on{lax}}(\ul\bC_1,\ul\bC_2)$$
admits a left adjoint, given, in terms of \eqref{e:lax vs strict abs}, by composing with the
left adjoint of \eqref{e:pullback untl from Ran to arr abs}.

\ssec{Pseudo-properness} \label{ss:pseudo-proper} 
 
\sssec{}

We introduce the category of \emph{pseudo-proper} categorical prestacks as 
$$\on{CatPreStk}_{\on{ps-proper}}:=\on{Funct}((\Sch^{\on{proper}})^{\on{op}},\inftyCat).$$

The embedding
$$\Sch^{\on{proper}}\hookrightarrow \on{PreStk}\hookrightarrow \on{CatPreStk}$$
uniquely extends to a colimit-preserving functor
\begin{equation} \label{e:forget ps-proper cat}
\on{CatPreStk}_{\on{ps-proper}}\to \on{CatPreStk}.
\end{equation}

Similarly, we define the category
$$\on{PreStk}_{\on{ps-proper}}:=\on{Funct}((\Sch^{\on{proper}})^{\on{op}},\inftygroup)$$
and the functor
\begin{equation} \label{e:forget ps-proper}
\on{PreStk}_{\on{ps-proper}}\to \on{PreStk}.
\end{equation}

\begin{rem}

Note that the functors \eqref{e:forget ps-proper cat} and \eqref{e:forget ps-proper} are \emph{not} fully faithful. 
I.e., with the above definition, pseudo-properness is not a property, but extra structure.

\medskip

However, since the functor $\Sch^{\on{proper}}\hookrightarrow \on{PreStk}$ preserves fiber products, so do
the functors \eqref{e:forget ps-proper cat} and \eqref{e:forget ps-proper}. 

\end{rem} 

\sssec{}

Concretely, pseudo-properness means the following: 

\medskip

A prestack $\CY$ is pseud-proper when it can written as 
$$\underset{i\in I}{\on{colim}}\, Z_i,$$
where:

\begin{itemize}

\item $Z_i$ are proper schemes;

\item The colimit is taken in $\on{PreStk}$.

\end{itemize}

 As morphisms
$$\underset{i\in I}{\on{colim}}\, Z_i\to \underset{i\in I}{\on{colim}}\, Z'_{i'}$$
we take 
$$\underset{i\in I^{\on{op}}}{\on{lim}}\, \underset{i'\in I'}{\on{colim}}\, \Maps(Z_i,Z'_{i'}).$$

\medskip

A categorical prestack $\CY$ is pseudo-proper if the prestacks $\on{Mor}^n(\CY)$ classifying $n$-fold composition of morphisms in $\CY$
are pseudo-proper, and for the maps $[n_1]\to [n_2]$ in $\bDelta$, the corresponding maps
$$\on{Mor}^{[n_2]}(\CY)\to \on{Mor}^{[n_1]}(\CY)$$
take place in $\on{PreStk}_{\on{ps-proper}}$.

\sssec{Example}

The prestack $\Ran$ and the categorical prestacks $\Ran^{\on{untl}}$ and $\Ran^{\on{untl},*}$
(see \secref{sss:Ranu ps-proper})
are pseudo-proper. 

\sssec{}

We define the functors $\Dmod(-)$ and $\on{CrystCat}(-)$ on $\on{CatPreStk}_{\on{ps-proper}}$ precomposing 
the same-named functors out of $\on{CatPreStk}$ with \eqref{e:forget ps-proper cat}.

\sssec{} \label{sss:proper approx}

For $\CY\in \on{CatPreStk}_{\on{ps-proper}}$, we can describe $\Dmod(\CY)$, $\on{CrystCat}(\CY)$ and
$$\Gamma^{\on{lax}}(\CY,\ul\bC), \quad \bC\in \on{CrystCat}(\CY)$$
in terms of proper schemes mapping to $\CY$. 

\medskip

I.e., in the appropriate definitions, we can replace 
$$\affSch_{/\CY} \rightsquigarrow \Sch^{\on{proper}}_{/\CY},$$
where in the right-hand side the morphisms take place in $\on{CatPreStk}_{\on{ps-proper}}$, 

\sssec{}

Let $\ul\bC$ be a crystal of categories over $\CY$, where $\CY$ is pseudo-proper. 
We claim:

\begin{lem} \label{l:limits ps-proper}
For $Z\overset{y}\to \CY$ with $Z$ proper and $y$ taking place in $\on{CatPreStk}_{\on{ps-proper}}$,
the functor of evaluation
$$\Gamma^{\on{lax}}(\CY,\ul\bC)\to \Gamma(Z,y^*(\ul\bC))=:\bC_{Z,y}$$
commutes with limits.
\end{lem} 

\begin{proof} 

We can describe the category 
$$\Gamma^{\on{lax}}(\CY,\ul\bC)$$ as a family of assignments
$$(Z\overset{y}\to \CY)\in \on{CatPreStk}_{\on{ps-proper}}\,\rightsquigarrow\, 
\bc_{Z,y}\in \bC_{Z,y},\quad Z \text{ is proper},$$
compatible under pullbacks:

\medskip

For $f:Z'\to Z$ we are given an isomorphism 
$$f^!(\bc_{Z,y})\simeq \bc_{Z',y\circ f}$$
in $\bC_{Z',y\circ f}$. 

\medskip

To prove the lemma, it suffices to show that the functors 
$$f^!:\bC_{Z,y}\to \bC_{Z',y\circ f}$$
commute with limits. Indeed, this would imply that limits in 
$\Gamma^{\on{lax}}(\CY,\ul\bC)$ are
computed component-wise in terms of $\{\bc_{Z,y}\}$.

\medskip

We claim that for any proper map $f:Z'\to Z$ and a crystal of categories $\ul\bD$ on $Z$,
the functor
$$f^!:\Gamma(Z,\ul\bD)\to \Gamma(Z',f^*(\ul\bD))$$
commutes with limits. 

\medskip

Indeed, since for a scheme its de Rham space is 1-affine, we can think of $\ul\bD$ as a 
$\Dmod(Z)$-linear category $\bD$, so that the functor $f^!$ is
\begin{equation} \label{e:^! crystals}
\Gamma(Z,\ul\bD)=\bD\simeq \Dmod(Z)\underset{\Dmod(Z)}\otimes \bD \overset{f^!\otimes \on{Id}}\to
\Dmod(Z')\underset{\Dmod(Z)}\otimes \bD=\Gamma(Z',f^*(\ul\bD)).
\end{equation} 

Now, since $f$ is proper, the functor $f^!:\Dmod(Z)\to \Dmod(Z')$ admits a left adjoint, namely,
$f_!$, which is automatically $\Dmod(Z)$-linear. Hence, the functor \eqref{e:^! crystals} admits a left adjoint, namely,
$$\Dmod(Z')\underset{\Dmod(Z)}\otimes \bD \overset{f_!\otimes \on{Id}}\to \Dmod(Z)\underset{\Dmod(Z)}\otimes \bD\simeq \bD.$$

This implies that $f^!$ commutes with limits. 

\end{proof}

\sssec{}

As a consequence of \lemref{l:limits ps-proper} we obtain:

\begin{cor} \label{c:!-dir image ps-proper}
Let $f:\CY_1\to \CY_2$ be a map between pseudo-proper categorical prestacks.
Then for $\ul\bC\in \on{CrystCat}(\CY_2)$, the functor
$$f^!:\Gamma^{\on{lax}}(\CY_2,\ul\bC)\to \Gamma^{\on{lax}}(\CY_1',f^*(\ul\bC))$$
admits a left adjoint (to be denoted $f_!$).
\end{cor}

\begin{proof}

It suffices to check that the functor $f^!$ commutes with limits. The latter follows from 
\lemref{l:limits ps-proper}.

\end{proof}

\begin{rem}

We warn the reader that although for a map $f$ between pseudo-proper prestacks,
the functor $f_!$ exists, it does not in general satisfy base change. (It does, however,
if $f$ is a value-wise co-Cartesian fibration, see \lemref{l:_! coCart} below.)

\end{rem} 

As a particular case of \corref{c:!-dir image ps-proper strict} we have:

\begin{cor} \label{c:integration ps-proper}
Let $\CY$ be pseudo-proper. Then the functor
$$\on{C}^\cdot_c(\CY,-):\Dmod(\CY)\to \Vect,$$ 
left adjoint to
$$\Vect \overset{\omega_\CY}\to \Dmod(\CY)$$
is well-defined.
\end{cor}

\sssec{}

Let $\CY$ be pseudo-proper. It follows formally that the prestack in groupoids $\CY^{\on{strict}}$ is also pseudo-proper.
Hence, from \corref{c:!-dir image ps-proper} we obtain:

\begin{cor} \label{c:!-dir image ps-proper strict}
Let $f:\CY_1\to \CY_2$ be a map between pseudo-proper categorical prestacks.
Then for a \emph{strict} $\ul\bC$, the functor 
$$f^!:\Gamma^{\on{strict}}(\CY_2,\ul\bC)\to \Gamma^{\on{strict}}(\CY_1,f^*(\ul\bC))$$
admits a left adjoint (to be denoted $f_!$).
\end{cor}

\begin{rem}

We warn the reader that the functors
$$f_!:\Gamma^{\on{lax}}(\CY_1,f^*(\ul\bC))\to \Gamma^{\on{lax}}(\CY_2,\ul\bC)$$
and 
$$f_!:\Gamma^{\on{strict}}(\CY_1,f^*(\ul\bC))\to \Gamma^{\on{strict}}(\CY_2,\ul\bC)$$
are in general \emph{incompatible} with the embeddings
$$\Gamma^{\on{strict}}(\CY_1,f^*(\ul\bC))\hookrightarrow \Gamma^{\on{lax}}(\CY_1,f^*(\ul\bC)) \text{ and }
\Gamma^{\on{strict}}(\CY_2,\ul\bC)\hookrightarrow \Gamma^{\on{lax}}(\CY_2,\ul\bC),$$
respectively.

\end{rem} 

\sssec{}

We will now show how to compute the functor $f_!$ more explicitly. First, unwinding the definitions, we obtain:

\begin{lem}  \label{l:_! coCart}
Let $f:\CY_1\to \CY_2$ be a map in $\on{CatPreStk}_{\on{ps-proper}}$ that is a 
co-Cartesian fibration.\footnote{As in \secref{sss:Cart prestacks}, but in the category $\on{Funct}((\Sch^{\on{proper}})^{\on{op}},\inftyCat)$.} 
Then for $\ul\bC\in \on{CrystCat}(\CY_2)$,
the functor $f_!$ satisfies base change, i.e., for a 
pullback diagram in $\on{CatPreStk}_{\on{ps-proper}}$
$$
\CD
\CY'_1 @>g_1>> \CY_1 \\
@V{f'}VV @VV{f}V \\
\CY'_2 @>g_2>> \CY_2,
\endCD
$$
the natural transformation 
$$f'_!\circ g_1^!\to g_2^!\circ f_!, \quad \Gamma^{\on{lax}}(\CY_1,f^*(\ul\bC))\rightrightarrows \Gamma^{\on{lax}}(\CY'_2,g_2^*(\ul\bC)),$$
obtained by adjunction from
$$g_1^!\circ f^!\simeq f'{}^!\circ g_2^!.$$
is an isomorphism.
\end{lem}

\begin{cor} \label{c:_! coCart}
Under the assumptions of \lemref{l:_! coCart}, for a proper scheme $Z$ equipped with a map $Z\overset{y_1}\to \CY_1$
in $\on{CatPreStk}_{\on{ps-proper}}$, the composition
$$y_1^!\circ f_!: \Gamma^{\on{lax}}(\CY_1,f^*(\ul\bC))\to \bC_{Z,y_1}$$
identifies canonically with 
$$\Gamma^{\on{lax}}(\CY_1,f^*(\ul\bC)) \overset{\on{pullback}}\longrightarrow 
\Gamma^{\on{lax}}(\CY_{1,Z},\ul\bC|_{\CY_{1,Z}}) \overset{(f_Z)_!}\to \bC_{Z,y_1},$$
where:

\begin{itemize}

\item $\CY_{1,Z}:=Z\underset{\CY_2}\times \CY_1$;

\item $f_Z$ is the map $\CY_{1,Z}\to Z$.

\end{itemize}

\end{cor} 

\sssec{}

Let now $f$ be an arbitrary map in $\on{CatPreStk}_{\on{ps-proper}}$. 
Denote $\CY_{2,f/}$ be the (pseudo-proper) categorical prestack given by the slice construction, i.e.,
$$\CY_{2,f/}(S)=\{y_1\in \CY_1(S),\,y_2\in \CY_2(S),\,f(y_1)\to y_2\}.$$

Let $\wt{f}$ and $\on{pr}_f$ denote the projections
$$\CY_{2,f/}\to \CY_2, \quad (y_1,y_2,f(y_1)\to y_2)\mapsto y_2, \quad 
\CY_{2,f/}\to \CY_1, \quad (y_1,y_2,f(y_1)\to y_2)\mapsto y_1,$$
respectively.

\medskip

Assume now that $\ul\bC$ is strict. In this case we have a canonical equivalence
$$\wt{f}^*(\ul\bC)\simeq \on{pr}_f^*\circ f^*(\ul\bC).$$

\medskip 

We claim:

\begin{lem}  \label{l:calc !}
Assume that $\ul\bC$ is strict. Then the functor $f_!$ identifies canonically with
$$\Gamma^{\on{lax}}(\CY_1,f^*(\ul\bC)) \overset{\on{pr}_f^!}\longrightarrow
\Gamma^{\on{lax}}(\CY_{2,f/},\on{pr}^*_f\circ f^*(\ul\bC)) \simeq \Gamma^{\on{lax}}(\CY_{2,f/},\wt{f}^*(\ul\bC))
\overset{\wt{f}_!}\longrightarrow \Gamma^{\on{lax}}(\CY_2,\ul\bC),$$
\end{lem}

\begin{proof} 

 Note that $f$ factors as
$$\on{diag}_f\circ \wt{f},$$
where 
$$\on{diag}_f:\CY_1\to \CY_{2,f/}, \quad y_1\mapsto (y_1,f(y_1),f(y_1)\overset{\on{id}}\to f(y_1)).$$

Hence, we obtain
$$f_!\simeq \wt{f}_!\circ (\on{diag}_f)_!.$$

\medskip

Now we claim that 
$$(\on{diag}_f)_!\simeq \on{pr}_f^!, \quad \Gamma^{\on{lax}}(\CY_1,f^*(\ul\bC)) \rightrightarrows \Gamma^{\on{lax}}(\CY_{2,f/},
\on{pr}^*_f\circ f^*(\ul\bC)).$$
Indeed, this follows from the fact that the morphisms
$(\on{diag}_f,\on{pr}_f)$ form an adjoint pair, cf. \lemref{l:adj pullback}. 

\end{proof}

\begin{cor} \label{c:calc !}
In the setting of \lemref{l:calc !}, for a proper scheme $Z$ equipped with a map $Z\overset{y_1}\to \CY_1$ in $\on{CatPreStk}_{\on{ps-proper}}$,
the composition
$$y_1^!\circ f_!: \Gamma^{\on{lax}}(\CY_1,f^*(\ul\bC))\to \bC_{Z,y_1}$$
identifies canonically with 
$$\Gamma^{\on{lax}}(\CY_1,f^*(\ul\bC)) \overset{\on{pullback}}\longrightarrow 
\Gamma^{\on{lax}}(\CY_{2,f/,Z},\ul\bC|_{\CY_{2,f/,Z}}) \overset{(\wt{f}_Z)_!}\longrightarrow \bC_{Z,y_1},$$
where:

\begin{itemize}

\item $\CY_{2,f/,Z}:=Z\underset{\CY_2}\times \CY_{2,f/}$;

\item $\wt{f}_Z$ is the map $\CY_{2,f/,Z}\to Z$.

\end{itemize}

\end{cor} 

\sssec{Example}

Let $\CY_1\to \CY_2$ be the map
$$\sft:\CY^{\on{grpd}}\to \CY.$$

The functor $\sft_!$ is the left adjoint of the forgetful functor, and the formula for it,
given by \propref{l:calc !}, coincides with that of \cite[Proposition 4.4.2]{Ga4}.

\sssec{}

Consider now the commutative square
\begin{equation} \label{e:op sq strict}
\CD
\CY^{\on{grpd}} @>{\sft}>> \CY \\
@V{\sft^{\on{op}}}VV  @VV{\on{strict}_\CY}V  \\
\CY^{\on{op}} @>{\on{strict}_{\CY^{\on{op}}}}>> \CY^{\on{strict}}.
\endCD
\end{equation}

We obtain a natural transformation
\begin{equation} \label{e:e:op sq strict 1}
(\sft^{\on{op}})_!\circ \sft^!\to  (\on{strict}_{\CY^{\on{op}}})^! \circ (\on{strict}_\CY)_!.
\end{equation}

Suppose now that $\CY$ has the property that \eqref{e:op sq strict} is Cartesian. In particular,
since the right vertical arrow is a co-Cartesian fibration, then so is the left vertical arrow. 
By \lemref{l:_! coCart}, we obtain that in this case \eqref{e:e:op sq strict 1} is an isomorphism.

\medskip

Applying this in the case $\CY=\Ran^{\on{untl}}$, this gives a conceptual explanation of the
commutativity of \eqref{e:ins unit int as functor untl compat}, at least in the particular case when $\ul\bC^{\on{loc}}=\ul\Dmod(\Ran^{\on{untl}})$,
$\bC^{\on{glob}}=\Vect$.

\sssec{} \label{sss:pseudo-proper rel} 

The above definitions and assertions admit a variant when we consider prestacks over a given affine base scheme $S$. In 
this case, one can talk about pseudo-properness relative to $S$, and the entire discussion applies.

\ssec{The unital Ran space} \label{ss:untl Ran}

\sssec{} 

There are two versions of the unital Ran space that we will consider: $\Ran^{\on{untl}}$ and
$\Ran^{\on{untl},*}$.

\medskip

For $S\in \affSch$, the category $\Ran^{\on{untl}}(S)$ is that of finite \emph{non-empty} subsets
in $\Maps(S_\dr,X)$, with the morphisms defined as follows:
$$\Maps_{\Maps(S_\dr,X)}(\ul{x}_1,\ul{x}_2)=
\begin{cases}
&\{*\} \text{ if } \ul{x}_1\subseteq \ul{x}_2,\\
&\emptyset \text{ otherwise}.
\end{cases}
$$

In the above formula, $\ul{x}_i$ denotes the finite subset of $\Hom(S_{\on{red}},X)$ corresponding to the same-named
$S\to \Ran$. 

\medskip

In the case of $\Ran^{\on{untl},*}$, we allow $\ul{x}$ to be empty, i.e., we add to $\Ran^{\on{untl}}$ a point $\{\emptyset\}$,
which is value-wise initial, corresponding to the empty set. 

\begin{rem}

There is a variant of the above definition, where the morphisms are defined by 
$$\Maps_{\Maps(S_\dr,X)}(\ul{x}_1,\ul{x}_2)=
\begin{cases}
&\{*\} \text{ if } \on{Graph}_{\ul{x}_1}\subseteq \on{Graph}_{\ul{x}_2},\\
&\emptyset \text{ otherwise}.
\end{cases}
$$
where in the formula $\subseteq$ means containment as closed subsets of $S\times X$. Denote the resulting categorical 
prestack by $'\Ran^{\on{untl}}$. 

\medskip

The two versions are equivalent for most practical purposes. Namely, we have a naturally defined map 
$\Ran^{\on{untl}}\to {}'\Ran^{\on{untl}}$, and we claim that it induces an equivalence between the
corresponding (2-)categories of crystals of categories. 

\medskip

This follows from the fact that the corresponding map
$$\on{Mor}^1(\Ran^{\on{untl}})\to \on{Mor}^1({}'\Ran^{\on{untl}})$$
becomes an isomorphism after sheafification in the Grothendieck topology generated by finite
surjective maps, while D-modules satisfy descent for this topology. 

\end{rem} 

\sssec{}

The above two versions of the unital Ran space, i.e., $\Ran^{\on{untl}}$ and $\Ran^{\on{untl},*}$,
are convenient in slightly different situations: the
$\Ran^{\on{untl},*}$ version is more convenient for discussing factorization, while the 
$\Ran^{\on{untl}}$ version is more convenient for the discussion of local-to-global functors. Yet,
the next assertion says that we could use $\Ran^{\on{untl},*}$ for the latter too.  

\medskip

\begin{prop} \label{p:local-to-global empty set}
Let $\ul\bC$ be a sheaf of categories over $\Ran^{\on{untl},*}$. Then for a DG category $\bD$,
pullback along $\Ran^{\on{untl}}\to \Ran^{\on{untl},*}$ gives rise to an an equivalence
\begin{multline*} 
\on{Funct}^{\on{strict}}_{\on{CrystCat}(\Ran^{\on{untl},*})}(\ul\bC,\bD\otimes \ul\Dmod(\Ran^{\on{untl},*}))\overset{\sim}\to \\
\overset{\sim}\to \on{Funct}^{\on{strict}}_{\on{CrystCat}(\Ran^{\on{untl}})}(\ul\bC|_{\Ran^{\on{untl}}},\bD\otimes \ul\Dmod(\Ran^{\on{untl}})).
\end{multline*}
\end{prop}

\begin{proof}

The assertion follows from the fact that the inclusion
$$\Ran^{\on{untl}}\to \Ran^{\on{untl},*}$$
is value-wise cofinal, and the following general claim:

\begin{lem} 
Let $f:\CY_1\to \CY_2$ be a value-wise cofinal morphism between categorical prestacks. Then for any pair of
crystals of categories $\ul\bC',\ul\bC''$ on $\CY_2$ with $\ul\bC''$ \emph{strict}, the functor
$$\on{Funct}^{\on{strict}}_{\CY_2}(\ul\bC',\ul\bC'')\to \on{Funct}^{\on{strict}}_{\CY_1}(f^*(\ul\bC'),f^*(\ul\bC'')).$$
\end{lem} 

\end{proof}

\sssec{} \label{sss:Ranu ps-proper}

Note that the presentation of $\Ran$ as in \secref{sss:Ran as colim} shows that it is pseudo-proper.

\medskip

We claim that $\Ran^{\on{untl}}$ is also pseudo-proper. We can write 
$$\on{Mor}^1(\Ran^{\on{untl}})=((\Ran^{\on{untl}})^\to)^{\on{grpd}} \simeq \Ran^{\subseteq}$$
as the colimit 
$$\underset{I_{\on{small}}\subseteq I_{\on{big}}}{\on{colim}}\, X^{I_{\on{big}}},$$
where the colimit is taken over the (opposite of the) category whose objects are pairs of non-empty finite sets
$I_1\to I_2$ and whose morphisms are commutative squares
$$
\CD
I _{\on{small}} @>>> I_{\on{big}} \\
@VVV @VVV \\
I'_{\on{small}} @>>> I'_{\on{big}}
\endCD
$$
with the vertical arrows surjective. 

\medskip

The maps $\on{pr}_{\on{big}}$ and $\on{pr}_{\on{small}}$ send the term corresponding to $I_{\on{small}}\subseteq I_{\on{big}}$
to 
$$X^{I_{\on{big}}}\to \Ran \text{ and } X^{I_{\on{big}}}\to X^{I_{\on{small}}}\to \Ran,$$
respectively. 

\medskip

The prestacks of higher-order compositions $\on{Mor}^n(\Ran^{\on{untl}})$ are described similarly. 

\sssec{}

By \corref{c:integration ps-proper}, the functors 
$$\on{C}^\cdot_c(\Ran^{\on{untl}},-):\Dmod(\Ran^{\on{untl}})\to \Vect \text{ and }
\on{C}^\cdot_c(\Ran,-):\Dmod(\Ran^{\on{untl}})\to \Vect,$$
left adjoint to 
$$k\mapsto \omega_{\Ran^{\on{untl}}} \text{ and } k\mapsto \omega_{\Ran},$$
respectively, are well-defined. 

\medskip

The same applies to $\Ran^{\on{untl},*}$.

\sssec{}

Being the left adjoint of a symmetric monoidal structure, the functor $\on{C}^\cdot_c(\Ran^{\on{untl}},-)$ carries
a naturally defined left-lax symmetric monoidal structure.

\medskip

We claim: 

\begin{lem}  \label{l:monoidal on Ran untl}
The left-lax monoidal structure on $\on{C}^\cdot_c(\Ran^{\on{untl}},-)$ is strict.
\end{lem}

\begin{proof}

We need to show that the natural transformation
$$\on{C}^\cdot_c(\Ran^{\on{untl}},-)\circ (\Delta_{\Ran^{\on{untl}}})^!\to
\on{C}^\cdot_c(\Ran^{\on{untl}}\times \Ran^{\on{untl}},-),$$
induced by the $((\Delta_{\Ran^{\on{untl}}})_!,(\Delta_{\Ran^{\on{untl}}})^!)$-adjunction,
is an isomorphism.

\medskip

This follows, however, from the fact that the morphism
$$\Delta_{\Ran^{\on{untl}}}:\Ran^{\on{untl}}\to \Ran^{\on{untl}}\times \Ran^{\on{untl}}$$
is value-wise cofinal. 

\end{proof}

\sssec{}

We now consider the relation between the functors
$$\on{C}^\cdot_c(\Ran^{\on{untl}},-) \text{ and } \on{C}^\cdot_c(\Ran,-).$$

We have the natural transformation
\begin{equation} \label{e:int over Ran and Ran untl}
\on{C}^\cdot_c(\Ran,-)\circ \sft^!\simeq \on{C}^\cdot_c(\Ran^{\on{untl}},-) \circ \sft_!\circ \sft^!\to \on{C}^\cdot_c(\Ran^{\on{untl}},-)
\end{equation}
as functors $\Dmod(\Ran^{\on{untl}})\to \Vect$. 

\sssec{}

We claim:

\begin{lem} \label{l:int over Ran and Ran untl}
The natural transformation transformation \eqref{e:int over Ran and Ran untl} is an isomorphism.
\end{lem}

This assertion is proved in \cite[Theorem 4.6.2]{Ga4}. We include the proof for completeness:

\medskip

Recall what it means for a morphism between categorical prestacks to be \emph{universally homologically cofinal},
see \cite[Sect. 3.5.1]{Ga4}. Using \cite[Corollary 3.5.12]{Ga4}, the assertion of \lemref{l:int over Ran and Ran untl}
follows from the next one: 

\begin{lem} \label{l:non untl to untl cofinal}
The map $\sft$ is universally homologically cofinal.
\end{lem} 

\begin{proof}[Proof of \lemref{l:non untl to untl cofinal}]

Let $S$ be an affine scheme and let us be given an $S$-point $\ul{x}$ of $\Ran$.
Consider the corresponding prestack
$$\Ran_{\ul{x}/},$$
see \cite[Sect. 3.5.1]{Ga4}. We need to show that it is \emph{universally homologically contractible}
over $S$ (see \cite[Sect. 2.5.1]{Ga4} for what this means).

\medskip

Note, however, that $\Ran_{\ul{x}/}$ is a prestack in groupoids isomorphic to $S^{\subseteq}_{\ul{x}}$,
and its projection to $S$ is pseudo-proper. Hence, it is enough to show that the fibers of the map
$$S^{\subseteq}_{\ul{x}}\to S$$
have trivial homology. The latter follows by the usual argument for the contractibility of the Ran space.

\end{proof} 

\begin{rem}
Statements parallel to Lemmas \ref{l:monoidal on Ran untl} and \ref{l:int over Ran and Ran untl}
hold for the $\Ran^{\on{untl},*}$ version of the unital Ran space, see \cite[Sect. 2.5]{Ro2}.
\end{rem} 

\sssec{} \label{sss:almost unital}

An analog of the assertion of \lemref{l:monoidal on Ran untl} would of course fail for the usual (i.e., non-unital)
Ran space. I.e., the natural transformation
\begin{equation} \label{e:left-lax nu}
\on{C}^\cdot_c(\Ran,-)\circ (\Delta_{\Ran})^!\simeq 
\on{C}^\cdot_c(\Ran\times \Ran,-)\circ (\Delta_{\Ran})_!\circ (\Delta_{\Ran})^!\to \on{C}^\cdot_c(\Ran\times \Ran,-)
\end{equation}
is not an isomorphism. 

\medskip

However, it admits the following variant:

\medskip

Let us denote by
$$\Dmod(\Ran)^{\on{almost-untl}}\subset \Dmod(\Ran)$$
the full subcategory generated by the essential image of the forgetful functor
$$\sft^!:\Dmod(\Ran^{\on{untl}})\to \Dmod(\Ran).$$

Note this subcategory is preserved by the monoidal operation.

\medskip

We claim the left-lax monoidal structure on $\on{C}^\cdot_c(\Ran,-)$, given by \eqref{e:left-lax nu}, 
becomes strict when restricted to 
to $\Dmod(\Ran)^{\on{almost-untl}}$. Indeed, this follows from \lemref{l:int over Ran and Ran untl}. 

\sssec{} \label{sss:Z subset untl}

By a similar token, for $\CZ\to \Ran$ one defines a unital version 
$$\CZ^{\subseteq,\on{untl}}:=\CZ\underset{\Ran^{\on{untl}}}\times (\Ran^{\on{untl}})^{\to}$$
of $\CZ^{\subseteq}$.

\medskip

The (categorical) prestacks $\CZ^{\subseteq}$ and $\CZ^{\subseteq,\on{untl}}$ are pseudo-proper
relative to $\CZ$. 

\medskip

The assertion of \lemref{l:int over Ran and Ran untl} renders to the present context, when instead of 
the functors
$$\on{C}^\cdot_c(\Ran,-) \text{ and } \on{C}^\cdot_c(\Ran^{\on{untl}},-)$$
we use the functors
$$(\on{pr}_{\on{small},\CZ})_!: \Dmod(\CZ^{\subseteq})\to \Dmod(\CZ)$$
and 
$$(\on{pr}^{\on{untl}}_{\on{small},\CZ})_!: \Dmod(\CZ^{\subseteq,\on{untl}})\to \Dmod(\CZ),$$
respectively. 

\ssec{Unital and counital factorization spaces}

\sssec{} \label{sss:unital spaces}

Let $\CZ_\Ran\to \Ran$ be a prestack. A unital structure on $\CZ_\Ran$ is its extension to a categorical prestack
\begin{equation} \label{e:unital prestack}
\CZ_{\Ran^{\on{untl},*}}\to \Ran^{\on{untl},*},
\end{equation}
such that \eqref{e:unital prestack} is a value-wise \emph{co-Cartesian fibration in groupoids}.

\medskip

Let $\on{PreStk}^{\on{untl}}_{/\Ran}$ denote the category of prestacks over $\Ran$, equipped with a unital
structure. 

\sssec{}

In concrete terms, an upgrade
$$\CZ_\Ran\rightsquigarrow \CZ_{\Ran^{\on{untl},*}}$$
means that for every $\ul{x}\subseteq \ul{x}'$
we give ourselves a map 
$$\on{ins.unit}_{\ul{x}\subseteq \ul{x}'}:\CZ_{\ul{x}}\to \CZ_{\ul{x}'},$$
in a way compatible with compositions.

\medskip

In addition, we give ourselves a space $\CZ_\emptyset$ and a system of maps
$$\on{ins.unit}_{\emptyset\subset \ul{x}}:\CZ_\emptyset\to \CZ_{\ul{x}}$$
equipped with identifications
$$\on{ins.unit}_{\ul{x}\subseteq \ul{x}'}\circ \on{ins.unit}_{\emptyset\subset \ul{x}}\simeq \on{ins.unit}_{\emptyset\subset \ul{x}'}.$$

\sssec{}

Here is an example of a prestack equipped with a unital structure (see \cite[Sect. 3.3]{Ro2}). Let $\CY$ be an affine D-scheme. Let
$$\on{Sect}_\nabla(X^{\on{gen}},\CY)_\Ran$$ 
be the space over $\Ran$ that attaches to $\ul{x}\in \Ran$ the space
$$\Sect_\nabla(X^{\on{gen}},\CY)_{\ul{x}}:=\Sect_\nabla(X-\ul{x},\CY).$$

This prestack has a natural unital structure: namely for $\ul{x}\subset \ul{x}'$, the corresponding map
$$\Sect_\nabla(X-\ul{x},\CY)\to \Sect_\nabla(X-\ul{x}',\CY)$$
is given by restriction.

\medskip

Note that 
$$\Sect_\nabla(X^{\on{gen}},\CY)_\emptyset\simeq \Sect_\nabla(X,\CY).$$

\sssec{} \label{sss:unital fact spaces}

Let $\CT$ be a factorization space over $X$. A unital structure on it is a unital structure on $\CT_\Ran$ (in the sense of \secref{sss:unital spaces})
and an extension of \eqref{e:fact space} to an isomorphism
\begin{multline} \label{e:fact space unital}
\CT_{\Ran^{\on{untl},*}}\underset{{\Ran^{\on{untl},*}},\on{union}}\times ({\Ran^{\on{untl},*}}\times {\Ran^{\on{untl},*}})_{\on{disj}} \simeq  \\
\simeq (\CT_{\Ran^{\on{untl},*}}\times \CT_{\Ran^{\on{untl},*}})\underset{{\Ran^{\on{untl},*}}\times {\Ran^{\on{untl},*}}}\times ({\Ran^{\on{untl},*}}\times {\Ran^{\on{untl},*}})_{\on{disj}},
\end{multline} 
equipped with a homotopy-coherent data of associativity and commutativity, where 
$$({\Ran^{\on{untl},*}}\times {\Ran^{\on{untl},*}})_{\on{disj}}\subset {\Ran^{\on{untl},*}}\times {\Ran^{\on{untl},*}}.$$
is the corresponding open subfunctor. 

\medskip

In addition, we stipulate that
\begin{equation} \label{e:fact space unital emptyset}
\CT_{\emptyset}\simeq \on{pt},
\end{equation} 
and this identification behaves (homotopically coherently) as a unit for the isomorphisms \eqref{e:fact space unital}, i.e., the map
$$\CT_{\Ran^{\on{untl},*}} \to \CT_{\emptyset}\times \CT_{\Ran^{\on{untl},*}},$$
obtained by base-changing \eqref{e:fact space unital} with respect to
$$\{\emptyset\} \times \Ran^{\on{untl},*}\to ({\Ran^{\on{untl},*}}\times {\Ran^{\on{untl},*}})_{\on{disj}}$$ 
identifies via \eqref{e:fact space unital emptyset}
with the identity map. 

\sssec{}

A typical example of a unital factorization space is $\Gr_G$. Namely, for $\ul{x}\subseteq \ul{x}'$ the corresponding map
$$\Gr_{G,\ul{x}}\to \Gr_{G,\ul{x}'}$$
is defined as follows. 

\medskip

Recall (see \secref{sss:Gr G disc})
that $\Gr_{G,\ul{x}}$ can be described as the space of $G$-bundles on $\cD_{\ul{x}}$ equipped with a trivialization on $\cD_{\ul{x}}-\ul{x}$.
We have:

\begin{lem} \label{l:glue discs}
The map 
$$\cD_{\ul{x}}\underset{\cD_{\ul{x}}-\ul{x}}\sqcup (\cD_{\ul{x}'}-\ul{x})\to \cD_{\ul{x}'}$$
is an isomorphism when the pushout is taken in the category of affine schemes. 
\end{lem} 

Using this lemma, we can interpret $\Gr_{G,\ul{x}}$ as the space of $G$-bundles on 
$\cD_{\ul{x}'}$ with a trivialization on $\cD_{\ul{x}'}-\ul{x}$. 

\medskip

The desired map
$$\on{ins.unit}_{\ul{x}\subseteq \ul{x}'}:\Gr_{G,\ul{x}}\to \Gr_{G,\ul{x}'}$$
is given by restricting the trivialization from $\cD_{\ul{x}'}-\ul{x}$ to $\cD_{\ul{x}'}-\ul{x}'$.

\sssec{}  \label{sss:counital spaces}

Let $\CZ_\Ran\to \Ran$ be a prestack. A \emph{counital} structure on $\CZ_\Ran$ is its extension to a categorical prestack
\begin{equation} \label{e:counital prestack}
\CZ_{\Ran^{\on{untl},*}}\to \Ran^{\on{untl},*},
\end{equation}
such that \eqref{e:counital prestack} is a value-wise \emph{Cartesian fibration in groupoids}.

\medskip

Let $\on{PreStk}^{\on{co-untl}}_{/\Ran}$ denote the category of prestacks over $\Ran$, equipped with a counital
structure. 

\sssec{}

In concrete terms, an upgrade
$$\CZ_\Ran\rightsquigarrow \CZ_{\Ran^{\on{untl},*}}$$
means that for every $\ul{x}\subseteq \ul{x}'$
we give ourselves a map 
$$\on{proj.counit}_{\ul{x}\subseteq \ul{x}'}:\CZ_{\ul{x}'}\to \CZ_{\ul{x}},$$
in a way compatible with compositions.

\medskip

In addition, we give ourselves a space $\CZ_\emptyset$ and a system of maps
$$\on{proj.counit}_{\emptyset\subset \ul{x}}:\CZ_{\ul{x}}\to \CZ_\emptyset$$
equipped with identifications
$$\on{proj.counit}_{\emptyset\subset \ul{x}}\circ \on{proj.counit}_{\ul{x}\subseteq \ul{x}'}\simeq
\on{proj.counit}_{\emptyset\subset \ul{x}'}.$$

\sssec{} 

Let $\CY$ be a D-prestack over $X$. Note that the arc space $\fL^+_\nabla(\CY)_\Ran$ has a natural counital structure:

\medskip

For $\ul{x}\subseteq \ul{x}'$, the corresponding maps
$$\fL^+_\nabla(\CY)_{\ul{x}'}\to \fL^+_\nabla(\CY)_{\ul{x}}$$
are given by restriction along $\wh\cD_{\ul{x}}\to \wh\cD_{\ul{x}'}$.

\sssec{}  \label{sss:another take on arcs}

We claim:

\begin{prop} \label{p:another take on arcs}
The functor $\CY\mapsto \fL^+_\nabla(\CY)_\Ran$ is the right adjoint to the functor
\begin{equation} \label{e:another take on arcs 0}
\on{PreStk}^{\on{co-untl}}_{/\Ran}\to \on{PreStk}_{/\Ran}\to \on{PreStk}_{X_\dr},
\end{equation}
where the last arrow is given by pullback along $X_\dr\to \Ran$.
\end{prop} 

\begin{proof} 

Let us construct the unit and counit map for the adjunction. The counit is easy, and it is actually an isomorphism:
we have
$$\fL^+_\nabla(\CY)_{X_\dr}\simeq \CY_\nabla.$$

\medskip

To construct the unit, for any $\ul{x}:S\to \Ran$ and $Z_\Ran\in \on{PreStk}^{\on{co-untl}}_{/\Ran}$, 
we need to define the map
$$\CZ_S \to \on{Weil-Res}^{\wh\cD_{\ul{x},\nabla}}_S(\wh\cD_{\ul{x},\nabla}\underset{X_\dr}\times \CZ_{X_\dr}),$$
where $\on{Weil-Res}$ is the functor of restriction of scalars \`a la Weil.

\medskip

By adjunction, the datum of the latter map is equivalent to that of a map
\begin{equation} \label{e:another take on arcs}
\wh\cD_{\ul{x},\nabla}\underset{S}\times \CZ_S \to \wh\cD_{\ul{x},\nabla}\underset{X_\dr}\times \CZ_{X_\dr}.
\end{equation}

Note that the two sides in \eqref{e:another take on arcs} are the pullbacks of $\CZ_\Ran$ along the following two maps
$$\wh\cD_{\ul{x},\nabla}\rightrightarrows \Ran.$$

One is
$$\wh\cD_{\ul{x},\nabla}\to S \overset{\ul{x}}\to \Ran,$$
and the other is
$$\wh\cD_{\ul{x},\nabla}\to X_\dr\to \Ran.$$

Now, by the definition of $\wh\cD_{\ul{x},\nabla}$, the there is a natural map from latter map to the former map 
inside the category
$$\Maps(\wh\cD_{\ul{x},\nabla},\Ran^{\on{untl},*}).$$

Hence, the required map is provided by the unital structure on $\CZ_\Ran$.

\medskip

The fact that the unit and counit maps constructed above satisfy the adjunction axioms is a straightforward verification. 

\end{proof} 

\sssec{} \label{sss:fact counital spaces}

Let $\CZ_\Ran$ be equipped with a counital structure. Note that this structure gives rise to a map
\begin{equation} \label{e:fact from counital prel} 
\CZ_\Ran\underset{\Ran,\on{union}}\times (\Ran\times \Ran) \to \CZ_\Ran\times \CZ_\Ran.
\end{equation} 

Base changing along
$$(\Ran\times \Ran)_{\on{disj}}\to \Ran\times \Ran,$$
we obtain a map 
\begin{equation} \label{e:fact from counital} 
\CZ_\Ran\underset{\Ran,\on{union}}\times (\Ran\times \Ran)_{\on{disj}} \to 
(\CZ_\Ran\times \CZ_\Ran)\underset{\Ran\times \Ran}\times (\Ran\times \Ran)_{\on{disj}}. 
\end{equation} 

We shall say that the counital structure is \emph{factorizable} if \eqref{e:fact from counital} is an isomorphism
and the diagonal map
$$\CZ_\emptyset\to \CZ_\emptyset\times \CZ_\emptyset$$
is also an isomorphism (implying that $\CZ_\emptyset\simeq \on{pt}$). 

\medskip

Note that the fact that the maps \eqref{e:fact from counital} are isomorphisms implies that the maps
\begin{multline} \label{e:fact space counital}
\CZ_{\Ran^{\on{untl},*}}\underset{{\Ran^{\on{untl},*}},\on{union}}\times ({\Ran^{\on{untl},*}}\times {\Ran^{\on{untl},*}})_{\on{disj}} \to \\
\to (\CZ_{\Ran^{\on{untl},*}}\times \CZ_{\Ran^{\on{untl},*}})\underset{{\Ran^{\on{untl},*}}\times {\Ran^{\on{untl},*}}}\times ({\Ran^{\on{untl},*}}\times {\Ran^{\on{untl},*}})_{\on{disj}},
\end{multline} 
are also isomorphisms.

\medskip

Let 
$$(\on{PreStk}^{\on{co-untl}}_{/\Ran})^{\on{factzbl}}\subset \on{PreStk}^{\on{co-untl}}_{/\Ran}$$
denote the full subcategory that consists of factorizable objects. 

\sssec{}

We have the following more precise version of \propref{p:another take on arcs}:

\begin{prop} \label{p:another take on arcs bis}
The functor 
\begin{equation} \label{e:another take on arcs bis}
\CY\mapsto \fL^+_\nabla(\CY)_\Ran
\end{equation} 
has an essential image in $(\on{PreStk}^{\on{co-untl}}_{/\Ran})^{\on{factzbl}}$. The resulting functor
$$\on{PreStk}_{/X_\dr}\to (\on{PreStk}^{\on{co-untl}}_{/\Ran})^{\on{factzbl}}$$
is fully faithful and defines an equivalence between the full subcategories consisting of objects that are affine 
over $X_\dr$ in the left-hand side and over $\Ran$ in the right-hand side. 
\end{prop} 

\begin{proof}

The fact that the essential image of the functor \eqref{e:another take on arcs bis}
lands in $(\on{PreStk}^{\on{co-untl}}_{/\Ran})^{\on{factzbl}}$ has been established in \secref{sss:forming arcs gen}.

\medskip

We have seen that the counit of the adjunction in \propref{p:another take on arcs} is an isomorphism.
Hence, the functor \eqref{e:another take on arcs bis} is fully faithful. 

\medskip

In order to prove the proposition, it remains to show that the functor \eqref{e:another take on arcs 0}
is conservative on objects in $(\on{PreStk}^{\on{co-untl}}_{/\Ran})^{\on{factzbl}}$ that are affine over $\Ran$. 

\medskip

Let $\phi:\CZ_{1,\Ran}\to \CZ_{2,\Ran}$ be a map between two objects in $(\on{PreStk}^{\on{co-untl}}_{/\Ran})^{\on{factzbl}}$,
such that the map 
$$\phi|_{X_\dr}:\CZ_{1,X_\dr}\to \CZ_{2,X_\dr}$$
is an isomorphism.  

\medskip

In order to check that $\phi$ is an isomorphism, it suffices to show that it is such when restricted to $X^I_\dr$ for every finite non-empty set $I$:
$$\CZ_{1,X^I_\dr}\to \CZ_{2,X^I_\dr}.$$

Since both prestacks are affine over $X^I_\dr$, the question of a map being an isomorphism can be checked strata-wise. 
Using the diagonal stratification of $X^I$, it suffices to show that that the further restriction to 
$$\overset{\circ}{X}{}^I_\dr\subset X^I_\dr$$
is an isomorphism. Since $\CZ_{1,\Ran}$ and $\CZ_{2,\Ran}$ are both factorizable, the latter map is the direct product of $I$ copies of the map
$\phi|_{X_\dr}$.

\end{proof} 

\sssec{} \label{sss:counital fact spaces}

As in \secref{sss:unital fact spaces}, given a factorization space, we can talk about a counital structure on it.

\medskip

It is easy to see, however, that if $\CT$ is a counital factorization space, then the corresponding prestack $\CT_{\Ran^{\on{untl},*}}$
is factorizable (in the sense of \secref{sss:fact counital spaces}). Vice versa, given an object
$$\CZ_{\Ran^{\on{untl},*}}\in (\on{PreStk}^{\on{co-untl}}_{/\Ran})^{\on{factzbl}},$$
the isomorphism \eqref{e:fact space counital} defines on it a factorization structure.

\medskip

So the categories of counital factorization spaces and factorizable counital prestacks over $\Ran$
are tautologically equivalent. 

\ssec{Unital factorization algebras}

\sssec{} \label{sss:untl fact alg}

Let $\CA$ be a factorization algebra on $X$. A unital structure on $\CA$ is an extension of $\CA_\Ran$ to
an object
$$\CA_{\Ran^{\on{untl},*}}\in \Dmod(\Ran^{\on{untl},*}),$$
and an extension of the isomorphism \eqref{e:fact alg} to an isomorphism 
\begin{equation} \label{e:untl fact alg}
\on{union}^!(\CA_{\Ran^{\on{untl},*}})|_{(\Ran^{\on{untl},*}\times \Ran^{\on{untl},*})_{\on{disj}}}\simeq 
\CA_{\Ran^{\on{untl},*}}\boxtimes \CA_{\Ran^{\on{untl},*}}|_{(\Ran^{\on{untl},*}\times \Ran^{\on{untl},*})_{\on{disj}}},
\end{equation} 
equipped with a homotopy-coherent data of associativity and commutativity.

\medskip

In addition, we stipulate that 
\begin{equation} \label{e:untl fact alg emptyset}
\CA_\emptyset\simeq k
\end{equation} 
and this isomorphism behaves (homotopically coherently) as a unit for the identifications \eqref{e:untl fact alg}, i.e., the map
$$\CA_{\Ran^{\on{untl},*}}\to \CA_\emptyset\otimes \CA_{\Ran^{\on{untl},*}},$$
obtained by restricting \eqref{e:untl fact alg} to 
$$\{\emptyset\}\times \Ran^{\on{untl},*}\to (\Ran^{\on{untl},*}\times \Ran^{\on{untl},*})_{\on{disj}}$$
identifies via \eqref{e:untl fact alg emptyset} with the identity map.

\medskip

Let $\on{FactAlg}^{\on{untl}}(X)$ denote the category of unital factorization algebras on $X$. 

\begin{rem}
Pullback along 
$$\ft:\Ran\to \Ran^{\on{untl},*}$$
gives rise to a functor 
\begin{equation} \label{e:forget unit on fact alg}
\on{FactAlg}^{\on{untl}}(X)\to \on{FactAlg}(X).
\end{equation} 

This functor is \emph{not} fully faithful. However, by analogy with the topological situation,
we expect that: 

\begin{itemize}

\item The functor \eqref{e:forget unit on fact alg} induces a \emph{monomorphism} on the mapping spaces;

\item The functor \eqref{e:forget unit on fact alg} induces an isomorphism on the union of the components of
the mapping spaces that correspond to \emph{isomorphisms}.

\end{itemize} 

The second property can be phrased as saying that for a factorization algebra, 
being unital is a property and not a structure. 

\medskip

We will not prove this in this paper. However, we will prove a result in this direction, 
see \propref{p:quasi-untl}. 

\end{rem} 

\sssec{} \label{sss:unit fact alg}

Let us denote by $k$ the unit factorization algebra, see \secref{sss:unit fact alg nu}. Note that it naturally upgrades
to a unital factorization algebra: namely, the corresponding object in $\Dmod(\Ran^{\on{untl},*})$
is $\omega_{\Ran^{\on{untl},*}}$. 

\medskip

Let $\CA$ be a unital factorization algebra on $X$. Note that the initial point $\{\emptyset\}\in \Ran^{\on{untl},*}$ gives rise to a map 
$$\omega_{\Ran^{\on{untl},*}}\to \CA_{\Ran^{\on{untl},*}}$$
in $\Dmod(\Ran^{\on{untl},*})$. It follows from the axioms that this map is compatible with factorization. 

\medskip

I.e., we obtain a map of unital factorization algebras
\begin{equation} \label{e:unit in fact alg}
\on{vac}_\CA:k\to \CA,
\end{equation} 
which we will refer to it as the \emph{vacuum} map for $\CA$. 

\sssec{} \label{sss:untl fact mod}

Let $\CZ$ be a prestack mapping to $\Ran$. Consider the corresponding categorical prestack $\CZ^{\subseteq,\on{untl}}$,
see \secref{sss:Z subset untl}.

\medskip

Let $\CA$ be a factorization algebra on $X$, and let $\CM$ be a factorization module $\CM$ over $\CA$ 
at $\CZ$. Let $\CA$ be equipped with a unital structure. A unital structure on $\CM$ is an extension of $\CM_{\CZ^{\subseteq}}$ 
to an object
$$\CM_{\CZ^{\subseteq}}\in \Dmod(\CZ^{\subseteq,\on{untl}})$$
and an extension of the isomorphism \eqref{e:fact mod} to an isomorphism 
\begin{equation} \label{e:fact mod untl}
\CM_{\CZ^{\subseteq,\on{untl}}}|_{(\Ran^{\on{untl},*}\times \CZ^{\subseteq,\on{untl}})_{\on{disj}}}\simeq
(\CA_{\Ran^{\on{untl},*}} \boxtimes \CM_{\CZ^{\subseteq,\on{untl}}})|_{(\Ran^{\on{untl},*}\times  \CZ^{\subseteq,\on{untl}})_{\on{disj}}}. 
\end{equation} 

In addition, we stipulate that the isomorphism
$$\CM_{\CZ^{\subseteq,\on{untl}}}\to \CA_{\emptyset}\otimes \CM_{\CZ^{\subseteq,\on{untl}}},$$
obtained by restricting \eqref{e:fact mod untl} along
$$\{\emptyset\} \times \CZ^{\subseteq,\on{untl}}\to (\Ran^{\on{untl},*}\times  \CZ^{\subseteq,\on{untl}})_{\on{disj}},$$
identifies via \eqref{e:untl fact alg emptyset} with the identity map.

\medskip

The contents of \secref{ss:fact alg} apply to unital factorization modules.

\sssec{Notational convention} 

When $\CA$ is a unital factorization algebra, we will denote by
$$\CA\mod^{\on{fact}}_\CZ$$
the category of \emph{unital} factorization modules over $\CA$ at $\CZ$.

\medskip

The category of modules over $\CA$ as a plain\footnote{I.e., non-unital.} factorization algebra will
be denoted by 
$$\CA\mod^{\on{fact-n.u.}}_\CZ.$$

\sssec{Example} \label{sss:vacuum module untl}

Let $\CA$ be a unital factorization algebra. For $\CZ\to \Ran$ consider the factorization module
$\CA^{\on{fact}_\CZ}$ from \secref{sss:vacuum module}.

\medskip

Unwinding the definitions, we obtain that $\CA^{\on{fact}_\CZ}$ carries a natural unital structure. 

\sssec{Example} \label{sss:unital fact mod for k}

Take $\CA=k$ from \secref{sss:unit fact alg}. Note that pullback along
$$\on{pr}^{\on{untl}}_{\on{small},\CZ}:\CZ^{\subseteq,\on{untl}}\to \CZ$$
gives rise to a functor
\begin{equation} \label{e:fact mod for unit}
\Dmod(\CZ)\to k\mod^{\on{fact}}_\CZ,
\end{equation}
cf. \secref{sss:fact mod over unit nu}.

\medskip

In other words, the functor \eqref{e:fact mod for unit} is given by tensoring (over $\Dmod(\CZ)$)
with the object
$$k^{\on{fact}_\CZ}\in  k\mod^{\on{fact}}_\CZ.$$

\medskip

We claim that the functor \eqref{e:fact mod for unit} is an equivalence, with the inverse functor
being
$$\CM\mapsto \CM_\CZ:=(\on{diag}^{\on{untl}}_\CZ)^!(\CM).$$

Indeed, the fact that the map
$$\on{diag}^{\on{untl}}_\CZ\to \CZ^{\subseteq,\on{untl}}$$
is initial relative to $\CZ$, implies that for any 
$\CM\in  k\mod^{\on{fact}}_\CZ$, we have a canonically defined map
$$(\on{pr}^{\on{untl}}_{\on{small},\CZ})^!\circ (\on{diag}^{\on{untl}}_\CZ)^!(\CM)\to \CM.$$

Now, the factorization condition implies that this map is actually an isomorphism. 

\sssec{}

Let $\CA$ be a unital factorization algebra on $X$.

\medskip 

Restriction along 
$$\sft:\CZ^{\subseteq}\to \CZ^{\subseteq,\on{untl}}$$
gives rise to a functor
\begin{equation} \label{e:unital vs non-unital modules}
\CA\mod^{\on{fact}}_\CZ\to \CA\mod^{\on{fact-n.u.}}_\CZ.
\end{equation}

\medskip

We have the following assertion, proved in \cite[Proposition 3.8.4]{CR}: 

\begin{prop} \label{p:unital vs non-unital modules}
The functor \eqref{e:unital vs non-unital modules} is fully faithful with essential image
$$\CA\mod^{\on{fact-n.u.}}_\CZ\underset{k\mod^{\on{fact-n.u.}}_\CZ}\times k\mod^{\on{fact}}_\CZ\simeq
\CA\mod^{\on{fact-n.u.}}_\CZ\underset{k\mod^{\on{fact-n.u.}}_\CZ}\times \Dmod(\CZ),$$
where 
$$\CA\mod^{\on{fact-n.u.}}_\CZ\to k\mod^{\on{fact-n.u.}}_\CZ$$
is the functor of restriction along \eqref{e:unit in fact alg}. 
\end{prop} 

\sssec{}

Let $\phi:\CA_1\to \CA_2$ be a unital map between unital factorization algebras. From 
\propref{p:unital vs non-unital modules} we obtain:

\begin{cor} \label{c:unital vs non-unital modules}
The restriction functor
$$\Res_\phi:\CA_2\mod^{\on{fact-n.u.}}_\CZ\to \CA_1\mod^{\on{fact-n.u.}}_\CZ$$
sends
$$\CA_2\mod^{\on{fact}}_\CZ\to \CA_1\mod^{\on{fact}}_\CZ.$$
\end{cor} 

\sssec{} \label{sss:q-untl}

Let $\on{FactAlg}^{\on{q-untl}}(X)$ be the category of pairs $(\CA,\on{vac}_\CA)$, where $\CA$ is a non-unital factorization algebra,
and $\on{vac}_\CA$ is a homomorphism $k\to \CA$, such that the object 
$$\Res_{\on{vac}_\CA}(\CA^{\on{fact}_\Ran})\in k\mod^{\on{fact-n.u.}}_\Ran$$
belongs to 
$$\Dmod(\Ran)\simeq k\mod^{\on{fact}}_\Ran\subset k\mod^{\on{fact-n.u.}}_\Ran.$$

We will call objects of $\on{FactAlg}^{\on{q-untl}}(X)$ ``quasi-unital factorization algebras". 

\medskip

We have a tautological functor 
\begin{equation} \label{e:untl to q-untl}
\on{FactAlg}^{\on{untl}}(X)\to \on{FactAlg}^{\on{q-untl}}(X).
\end{equation}

We claim:

\begin{prop} \label{p:quasi-untl}
The functor \eqref{e:untl to q-untl} is an equivalence. 
\end{prop}

The proof will be given in \secref{sss:proof quasi-unital}. 

\medskip

Note that the second assertion 
of the proposition says that a quasi-unital factorization algebra 
$$k \overset{\on{vac}_\CA}\to \CA$$
carries a canonical unital structure, for which $\on{vac}_\CA$ is the unit.

\ssec{Commutative \emph{unital} factorization algebras} \label{ss:com untl fact alg}

\sssec{} \label{sss:untl fact alg com} 

By the same token as in \secref{sss:fact com}, one can consider the category
$$\on{ComAlg}(\on{FactAlg}^{\on{untl}}(X)).$$

It is equipped with a tautological forgetful functor 
\begin{equation} \label{e:forget fact untl com}
\on{ComAlg}(\on{FactAlg}^{\on{untl}}(X))\to \on{ComAlg}(\Dmod(\Ran^{\on{untl},*}))
\end{equation}
and also with a functor
\begin{equation}  \label{e:forget unit on com fact alg}
\sft^!:\on{ComAlg}(\on{FactAlg}^{\on{untl}}(X))\to \on{ComAlg}(\on{FactAlg}(X)).
\end{equation}

\sssec{}

Let $\CA_{\Ran^{\on{untl},*}}$ be an object of $\on{ComAlg}(\Dmod(\Ran^{\on{untl},*}))$. Note that the unital structure on $\CA$ 
gives rise to the maps
$$(p_i)^!(\CA_{\Ran^{\on{untl},*}})\to \on{union}^!(\CA_{\Ran^{\on{untl},*}}),\quad i=1,2$$
where $p_1$ and $p_2$ are the two projections $\Ran^{\on{untl},*}\times \Ran^{\on{untl},*}\to \Ran^{\on{untl},*}$. 

\medskip

Since the coproduct in $\on{ComAlg}$ is the tensor product, we obtain a map
\begin{equation} \label{e:ten to union prel}
\CA_{\Ran^{\on{untl},*}}\boxtimes \CA_{\Ran^{\on{untl},*}}\to \on{union}^!(\CA_{\Ran^{\on{untl},*}}).
\end{equation}

The map \eqref{e:ten to union prel} gives rise to a map
\begin{equation} \label{e:ten to union}
\CA_{\Ran^{\on{untl},*}}\boxtimes \CA_{\Ran^{\on{untl},*}}|_{(\Ran^{\on{untl},*}\times 
\Ran^{\on{untl},*})_{\on{disj}}}\to \on{union}^!(\CA_{\Ran^{\on{untl},*}})|_{(\Ran^{\on{untl},*}\times \Ran^{\on{untl},*})_{\on{disj}}}.
\end{equation}

We shall say that $\CA_{\Ran^{\on{untl},*}}$ is \emph{factorizable} if the map \eqref{e:ten to union} is an isomorphism. 
Note that this automatically implies that $\CA_\emptyset\simeq k$.

\sssec{}

Let 
$$\on{ComAlg}(\Dmod(\Ran^{\on{untl},*}))^{\on{factzble}}\subset \on{ComAlg}(\Dmod(\Ran^{\on{untl},*}))$$
denote the full subcategory consisting of factorizable objects. 

\medskip

It follows from the axioms that the essential image of the functor \eqref{e:forget fact untl com}
lands in $$\on{ComAlg}(\Dmod(\Ran^{\on{untl},*}))^{\on{factzble}}.$$ Vice versa, for an object 
$$\CA_{\Ran^{\on{untl},*}}\in \on{ComAlg}(\Dmod(\Ran^{\on{untl},*}))^{\on{factzble}},$$
the isomorphism \eqref{e:ten to union} defines on $\CA_{\Ran^{\on{untl},*}}$ a factorization structure.

\medskip

It is easy to see that the resulting two functors
\begin{equation} \label{e:com fact untl two versions} 
\on{ComAlg}(\on{FactAlg}^{\on{untl}}(X))\leftrightarrow \on{ComAlg}(\Dmod(\Ran^{\on{untl},*}))^{\on{factzble}}
\end{equation} 
are mutually inverse.

\sssec{} \label{sss:com fact vs Dmod com untl}

Recall the functor
$$A\mapsto \on{Fact}(A),$$
see \secref{sss:com fact vs Dmod com}.

\medskip

Unwinding the construction, we obtain that $\on{Fact}(-)$ upgrades to a functor
$$\on{ComAlg}(\Dmod(X))\to \on{ComAlg}(\on{FactAlg}^{\on{untl}}(X)).$$

By a slight abuse of notation, we will use the same symbol $\on{Fact}(-)$ to denote the latter functor. 

\sssec{}

Consider now the functor
\begin{equation} \label{e:fact com untl to X}
\on{ComAlg}(\Dmod(\Ran^{\on{untl},*}))\to \on{ComAlg}(\Dmod(X)),
\end{equation} 
given by restriction along $\Delta_{X,\Ran^{\on{untl},*}}:X_\dr\to \Ran^{\on{untl},*}$. 

\medskip

We claim:

\begin{prop} \label{p:fact com untl to X} 
The functor \eqref{e:fact com untl to X} admits a left adjoint. Moreover, this left adjoint
is fully faithful and lands in $\on{ComAlg}(\Dmod(\Ran^{\on{untl},*}))^{\on{factzble}}$.
\end{prop}

\begin{proof}

It is enough to prove the assertion of the proposition on objects of $\on{ComAlg}(\Dmod(X))$
of the form 
$$\Sym^!(\CM), \quad \CM\in \Dmod(X).$$

The value of the left adjoint on such an object is
\begin{equation} \label{e:Sym dir im}
\Sym^!((\Delta_{X,\Ran^{\on{untl},*}})_!(\CM)).
\end{equation} 

To prove that this left adjoint is fully faithful, it is enough to show that the unit of the adjunction 
$$\CM\to (\Delta_{X,\Ran^{\on{untl},*}})^!\circ (\Delta_{X,\Ran^{\on{untl},*}})_!(\CM)$$
is an isomorphism. However, this follows from \corref{c:calc !}: indeed, the categorical prestack
$$X_\dr \underset{\Ran^{\on{untl},*}}\times (\Ran^{\on{untl},*})_{\Delta_{X,\Ran^{\on{untl},*}}/}$$
identifies with $X_\dr$. 

\medskip

In order to show that \eqref{e:Sym dir im} belongs to $\on{ComAlg}(\Dmod(\Ran^{\on{untl},*}))^{\on{factzble}}$,
it suffices to show that the canonical map
$$p_1^!((\Delta_{X,\Ran^{\on{untl},*}})_!(\CM))\oplus 
p_2^!((\Delta_{X,\Ran^{\on{untl},*}})_!(\CM))\to \on{union}^!((\Delta_{X,\Ran^{\on{untl},*}})_!(\CM))$$
becomes an isomorphism after restricting to $(\Ran^{\on{untl},*}\times \Ran^{\on{untl},*})_{\on{disj}}$. 
This is the content of \cite[Theorem 2.7.6]{Ro2}. We include the argument for completeness. 

\medskip

However, this follows again from \corref{c:calc !}: we have a canonical isomorphism
\begin{multline*} 
\left((\Ran^{\on{untl},*})_{\Delta_{X,\Ran^{\on{untl},*}}/} \times \Ran^{\on{untl},*} \sqcup 
\Ran^{\on{untl},*}\times (\Ran^{\on{untl},*})_{\Delta_{X,\Ran^{\on{untl},*}}/}\right)  
\underset{\Ran^{\on{untl},*}\times \Ran^{\on{untl},*}}\times \\
\underset{\Ran^{\on{untl},*}\times \Ran^{\on{untl},*}}\times (\Ran^{\on{untl},*}\times \Ran^{\on{untl},*})_{\on{disj}}\simeq \\
\simeq (\Ran^{\on{untl},*})_{\Delta_{X,\Ran^{\on{untl},*}}/}\underset{\Ran^{\on{untl},*},\on{union}}\times 
(\Ran^{\on{untl},*}\times \Ran^{\on{untl},*})_{\on{disj}}.
\end{multline*} 

Indeed, this is just the fact that for a disjoint pair $\ul{x}_1,\ul{x}_2$ of points of Ran, and a singleton $x$, 
$$x\subseteq \ul{x}_1\cup \ul{x}_2\,\,  \Leftrightarrow\,\, x\subseteq \ul{x}_1 \text{ or } x\subseteq \ul{x}_2.$$

\end{proof} 

\sssec{}

As a corollary, we obtain: 

\begin{cor} \label{c:fact com untl to X} 
The composite functor
$$\on{ComAlg}(\on{FactAlg}^{\on{untl}}(X)) \to \on{ComAlg}(\Dmod(\Ran^{\on{untl},*})) \overset{\Delta_{X,\Ran^{\on{untl},*}}^!}\longrightarrow
\on{ComAlg}(\Dmod(X))$$
is an equivalence, with the inverse given by $\on{Fact}(-)$.
\end{cor} 

\begin{proof}

Given \propref{p:fact com untl to X}, we only need to prove that the functor in \corref{c:fact com untl to X} is conservative.
But this is immediate from the factorization.

\end{proof}

\begin{cor} \label{c:fact com untl to X bis} 
The functor
$$\on{ComAlg}(\Dmod(X))\overset{\on{Fact}(-)}\longrightarrow \on{ComAlg}(\on{FactAlg}^{\on{untl}}(X)) \to  \on{ComAlg}(\Dmod(\Ran^{\on{untl},*}))$$
is the left adjoint of 
$$\on{ComAlg}(\Dmod(\Ran^{\on{untl},*})) \overset{\Delta_{X,\Ran^{\on{untl},*}}^!}\longrightarrow
\on{ComAlg}(\Dmod(X)).$$
\end{cor} 

\sssec{} \label{sss:com fact and arcs gen}

Let $\CY$ be a D-prestack over $X$. Suppose that the prestack $\fL_\nabla^+(\CY)_\Ran\to \Ran$ is such that the formation
of direct image of the structure sheaf along
$$\fL^+_\nabla(\CY)_S\to S, \quad S\in \affSch_{/\Ran}$$
is compatible with base change, and satisfies Kunneth formula. 

\medskip

This happens, e.g., when $\CY$ is affine over $X$, and hence $\fL^+_\nabla(\CY)_\Ran$ is affine over $\Ran$. 

\medskip

Taking the direct image of the structure sheaf along 
$$\fL^+_\nabla(\CY)_{\Ran^{\on{untl},*}}\to \Ran^{\on{untl},*},$$
we obtain an object in $\on{ComAlg}(\Dmod(\Ran^{\on{untl},*}))$, which by a slight abuse of notation we denote by
$\CO_{\fL^+_\nabla(\CY),\Ran^{\on{untl},*}}$. By Kunneth formula and factorization, $\CO_{\fL^+_\nabla(\CY),\Ran^{\on{untl},*}}$
has a natural structure of factorization algebra. 

\medskip

Denote the resulting object by
$$\CO_{\fL^+_\nabla(\CY)}\in \on{FactAlg}(X).$$

\medskip

The value of $\CO_{\fL^+_\nabla(\CY)}$ on $X$, i.e., the 
restriction of $\CO_{\fL^+_\nabla(\CY),\Ran^{\on{untl},*}}$ along $X\to \Ran^{\on{untl},*}$, is the direct image of the structure sheaf
along $\CY\to X$, which by a slight abuse of notation we denote by $\CO_\CY$. 

\medskip

Hence, by the equivalence of \corref{c:fact com untl to X}, we obtain that
\begin{equation} \label{e:arcs as fact}
\CO_{\fL^+_\nabla(\CY)} \simeq \on{Fact}(\CO_\CY).
\end{equation} 

\sssec{} \label{sss:com fact and arcs aff}

Assume that $\CY$ is affine over $X$, i.e., $\CY=\Spec_X(A)$ for 
$A\in \on{ComAlg}(\Dmod(X))$ with $\oblv^l(A)$ is connective. 

\medskip

Let $\CA:=\on{Fact}(A)$.  
Then \eqref{e:arcs as fact} says that for $S\to \Ran$,
\begin{equation} \label{e:arcs as fact again}
\fL^+_\nabla(\CY)_S\simeq \Spec_S(\CA_S).
\end{equation} 

\sssec{}

Note that the identification \eqref{e:arcs as fact again}
is in agreement with \propref{p:another take on arcs}. 

\medskip

Namely, let us be given a counital prestack 
$\CZ_\Ran$ over $\Ran$, such that the direct image of its structure sheaf satisfies base change.
Let $$\CB_{\Ran^{\on{untl},*}}\in\on{ComAlg}(\Dmod(\Ran^{\on{untl},*}))$$
denote the corresponding object. Set $\CB_X:=\CB_{\Ran^{\on{untl},*}}|_{X_\dr}$. 

\medskip

Then the following diagram commutes: 
$$
\CD
\Maps_{\on{PreStk}^{\on{co-untl}}_{/\Ran}}(\CZ_{\Ran^{\on{untl},*}},\fL^+_\nabla(\CY)_{\Ran^{\on{untl},*}}) 
@>{\text{\propref{p:another take on arcs}}}>{\sim}>
\Maps_{X_\dr}(X_\dr\underset{\Ran^{\on{untl},*}}\times \CZ_{\Ran^{\on{untl},*}},\CY_\nabla) \\
@V{\text{\eqref{e:arcs as fact again}}}V{\sim}V @VV{\sim}V \\
\Maps_{\on{ComAlg}(\Dmod(\Ran^{\on{untl},*}))}(\CA,\CB_{\Ran^{\on{untl},*}}) @>{\text{\propref{p:fact com untl to X}}}>> 
\Maps_{\on{ComAlg}(\Dmod(X))}(A,\CB_X).
\endCD
$$

\ssec{Factorization homology of commutative factorization algebras}

\sssec{}

Consider the functor
$$\on{ComAlg}(\Vect)\to \on{ComAlg}(\Dmod(X)), \quad R\mapsto R\otimes \CO_X,$$
where $\CO_X$ is perceived as a \emph{left} D-module.

\medskip

In this subsection, we will describe, following \cite[Sect. 4.6.1]{BD2}, its left adjoint.

\sssec{} \label{sss:left adj on untl alg}

From \lemref{l:monoidal on Ran untl} we obtain that the functor $\on{C}^\cdot_c(\Ran^{\on{untl}},-)$
gives rise to a functor
$$\on{ComAlg}(\Dmod(\Ran^{\on{untl}}))\to \on{ComAlg}(\Vect),$$
left adjoint to 
$$\on{ComAlg}(\Vect)\to \on{ComAlg}(\Dmod(\Ran^{\on{untl}})), \quad R\mapsto R\otimes \omega_{\Ran^{\on{untl}}}.$$

\begin{rem} \label{r:monoidal on Ran}

Note that by \secref{sss:almost unital}, the restriction of $\on{C}^\cdot_c(\Ran,-)$ to 
the subcategory 
$$\on{ComAlg}(\Dmod(\Ran)^{\on{almost-untl}})\subset \on{ComAlg}(\Dmod(\Ran))$$
defines a left adjoint to the functor 
$$\on{ComAlg}(\Vect)\to \on{ComAlg}(\Dmod(\Ran)^{\on{almost-untl}}), \quad R\mapsto R\otimes \omega_\Ran.$$

\end{rem}

\sssec{}

Recall now that according to \corref{c:fact com untl to X bis}, the functor 
$$\on{ComAlg}(\Dmod(X)) \overset{\on{Fact}}\to \on{ComAlg}(\on{FactAlg}^{\on{untl}}(X)) \to \on{ComAlg}(\Dmod(\Ran^{\on{untl},*}))$$
provides a left adjoint to the restriction functor
$$\on{ComAlg}(\Dmod(\Ran^{\on{untl},*}))\to \on{ComAlg}(\Dmod(X)).$$

Combined with \secref{sss:left adj on untl alg}, we obtain:

\begin{cor} \label{c:ch homology as left adj}
The functor 
\begin{multline*}
\on{ComAlg}(\Dmod(X)) \overset{\on{Fact}}\to \on{ComAlg}(\on{FactAlg}^{\on{untl}}(X)) \to \\
\to \on{ComAlg}(\Dmod(\Ran^{\on{untl},*}))
\overset{\on{C}^\cdot_c(\Ran^{\on{untl},*},-)}\longrightarrow \on{ComAlg}(\Vect)
\end{multline*} 
is the left adjoint of 
$$R\mapsto R\otimes \CO_X, \quad \on{ComAlg}(\Vect)\to\on{ComAlg}(\Dmod(X)).$$
\end{cor}

\sssec{} \label{sss:vac fact hom again}

For $\CA\in \on{FactAlg}^{\on{untl}}(X)$ recall the object 
$$\on{C}^{\on{fact}}_\cdot(X,\CA)=\on{C}^\cdot_c(\Ran,\CA_\Ran),$$
see \secref{sss:vac fact hom}. 

\medskip

Note that by \lemref{l:int over Ran and Ran untl}, we can rewrite this also as 
$$\on{C}^\cdot_c(\Ran^{\on{untl}},\CA_{\Ran^{\on{untl}}}).$$

\sssec{}

Thus, \corref{c:ch homology as left adj} says that the functor
$$A\mapsto \on{C}^{\on{fact}}_\cdot(X,\on{Fact}(A)), \quad \on{ComAlg}(\Dmod(X)) \to \on{ComAlg}(\Vect)$$
is the left adjoint of 
$$R\mapsto R\otimes \CO_X, \quad \on{ComAlg}(\Vect)\to\on{ComAlg}(\Dmod(X)).$$

\begin{rem}

Note that when we think of $\on{C}^{\on{fact}}_\cdot(X,\CA)$ as $\on{C}^\cdot_c(\Ran,\CA_\Ran)$, 
the commutative algebra structure on it follows from Remark \ref{r:monoidal on Ran}, since
$$\CA_\Ran\in \Dmod(\Dmod_\Ran)^{\on{almost-untl}}.$$
 
\end{rem} 

\sssec{}

Let $A\to B$ be a map in $\on{ComAlg}(\Dmod(X))$. Denote
$$\CA:=\on{Fact}(A),\,\, \CB:=\on{Fact}(B).$$

Let $R$ be an object of $\on{ComAlg}(\Vect)$, and fix a map
$$A\to R\otimes \CO_X,$$
or, equivalently by \corref{c:ch homology as left adj}, a map
$$\on{C}^{\on{fact}}_\cdot(X,\CA)\to R.$$

\sssec{}

Denote
$$B_R:=B\underset{A}\otimes (R\otimes \CO_X).$$

We can view $B_R$ as an object of 
$$\on{ComAlg}(\Dmod(X)\otimes R\mod),$$
i.e., as an $R$-linear object in $\on{ComAlg}(\Dmod(X))$. 

\sssec{}

Consider the corresponding object
$$\CB_R\in \on{ComAlg}(\on{FactAlg}(X)\otimes R\mod).$$

\medskip

We can apply the construction of factorization homology in the $R$-linear context, and 
thus form
$$\on{C}^{\on{fact}}_\cdot(X,\CB_R)\in \on{ComAlg}(R\mod).$$

\sssec{}

We have the naturally defined maps in $\on{ComAlg}(\Vect)$
$$R\to \on{C}^{\on{fact}}_\cdot(X,\CB_R) \leftarrow \on{C}^{\on{fact}}_\cdot(X,\CB),$$
which fit into the commutative diagram
$$
\CD
\on{C}^{\on{fact}}_\cdot(X,\CA) @>>> R \\
@VVV @VVV \\
\on{C}^{\on{fact}}_\cdot(X,\CB) @>>> \on{C}^{\on{fact}}_\cdot(X,\CB_R).
\endCD
$$

In particular, we obtain a map
\begin{equation} \label{e:rel fact homology}
R\underset{\on{C}^{\on{fact}}_\cdot(X,\CA)}\otimes \on{C}^{\on{fact}}_\cdot(X,\CB)\to \on{C}^{\on{fact}}_\cdot(X,\CB_R).
\end{equation}

\sssec{}

We claim:

\begin{lem} \label{l:rel fact homology}
The map \eqref{e:rel fact homology} is an isomorphism.
\end{lem}

\begin{proof}

Follows immediately from  \corref{c:ch homology as left adj}.

\end{proof} 

\ssec{Unitality in correspondenes} \label{ss:untl in corr}

\sssec{} \label{sss:defn corr}

Let $\Phi:\bC\to \bD$ be a functor between $\infty$-categories. We shall say that $\Phi$
is a \emph{fibration-in-correspondences}\footnote{Another name for this is 
``conservative exponentiable fibration", see \cite{AF}.} 
if the following two conditions hold:

\medskip

\begin{itemize}

\item For every $\bd\in \bD$, the fiber $\bC_\bd$ is a groupoid;

\item For every composable pair of arrows $\bd_0\overset{\alpha_{0,1}}\to \bd_1\overset{\alpha_{1,2}}\to \bd_2$ in $\bD$, the map
$$\bC_{\alpha_{0,1}}\underset{\bC_{\bd_1}}\times \bC_{\alpha_{1,2}}\to \bC_{\alpha_{1,2}\circ \alpha_{0,1}}$$
is an isomorphism, where for an arrow $\bd'\overset{\alpha}\to \bd''$, viewed as a functor $[0,1]\to \bD$, 
we denote by $\bC_\alpha$ the category of lifts of $\alpha:[0,1]\to \bD$ to a functor 
$[0,1]\to \bC$, i.e.,
$$\on{Funct}_{/\bD}([0,1],\bC).$$

\end{itemize}

\medskip

Note that the second condition can be reformulated as follows: given an arrow $\bc_0\overset{\beta_{0,2}}\to \bc_2$ in $\bC$ 
and a factorization of its image $\Phi(\beta_{0,2})$ as 
$$\Phi(\bd_0)=:\bd_0\overset{\alpha_{0,1}}\to \bd_1\overset{\alpha_{1,2}}\to \bd_2:=\Phi(\bd_2)$$
the space of the factorizations 
$$\beta_{0,2}=\beta_{1,2}\circ \beta_{0,1}, \quad \Phi(\beta_{0,1})=\alpha_{0,1},\,\, \Phi(\beta_{1,2})=\alpha_{1,2}$$
is contractible. 

\begin{rem}

There is a version of straightening construction that attaches to a fibration-in-correspondences $\Phi:\bC\to \bD$
a functor from $\bD$ to $\on{Corr}(\on{Grpd})$, i.e., the category, whose objects are groupoids and the morphisms are correspondences
between groupoids, see \cite[Theorem 0.8(2) and Theorem 0.10(1)]{AF}.  

\end{rem}

\sssec{}

Note that if $\Phi$ is either a Cartesian or a co-Cartesian fibration in groupoids, then it is a 
fibration-in-correspondences.

%
%
%
%
%

\sssec{} \label{sss:Z+}

Let $\Phi:\bC\to \bD$ be a fibration-in-correspondences. Suppose for a moment that $\bD$ contains an 
initial object $\{\emptyset\}$. Let $\bC^+$ denote the category 
$$\{\bc_\emptyset\in \bC_{\{\emptyset\}},\, \bc\in \bC, \,\bc_\emptyset\overset{\beta}\to \bc\}.$$

Note the functor
$$\Phi^+:\bC^+\to \bD, \quad (\bc_\emptyset,\bc,\beta)\mapsto \Phi(\bc)$$
is a Cartesian fibration in groupoids. 

\medskip

Indeed, given 
$$\bc_\emptyset\overset{\beta}\to \bc \text{ and } \bd'\to \Phi(\bc),$$
we let 
$$\bc_\emptyset\overset{\beta'}\to \bc'\to \bc$$
be its unique  factorization covering the canonical factorization
$$\{\emptyset\}\to \bd' \to \Phi(\bc).$$

\medskip

Note that if $\bC\to \bD$ is co-Cartesian fibration in groupoids, then $\bC^+\simeq \bC_{\emptyset}\times \bD$. 
If $\bC\to \bD$ is Cartesian fibration in groupoids, then $\bC^+\to \bC$ is an equivalence. 

\sssec{} 

Let $\CZ_\Ran\to \Ran$ be a prestack. A unital-in-correspondences structure on it is an extension of $\CZ_\Ran$ to 
a categorical prestack
$$\CZ_{\Ran^{\on{untl},*}}\to \Ran^{\on{untl},*},$$
which is a value-wise fibration-in-correspondences

\sssec{} \label{sss:rel unital}

Note that the construction in \secref{sss:Z+} associates to such $\CZ_\Ran$ a prestack $\CZ^+_\Ran$,
equipped with a counital structure. 

\medskip

In what follows we will say that $\CZ_\Ran$ admits a unital-in-correspondences structure \emph{relative} to $\CZ^+$. 
And we will refer to $\CZ^+_\Ran$ as the counital prestack \emph{underlying} $\CZ_\Ran$. 

\medskip

For an arrow in $\Ran^{\on{untl},*}$ given by $\ul{x}\subseteq \ul{x}'$ we will denote by
$$\CZ^{\on{all}\rightsquigarrow +}_{\ul{x}\subseteq \ul{x}'}$$ the prestack of its lifts to $\CZ_\Ran$
(see \secref{sss:defn corr} above). It is equipped with maps
$$\CZ_{\ul{x}} \overset{\on{pr}^\CZ_{\on{small}}}\longleftarrow \CZ^{\on{all}\rightsquigarrow +}_{\ul{x}\subseteq \ul{x}'} 
 \overset{\on{pr}^\CZ_{\on{big}}}\longrightarrow\CZ_{\ul{x}'}.$$

\sssec{}

Note that by \secref{sss:Z+}, a unital structure on $\CZ_\Ran$ gives rise to a unital-in-correspondences
structure, with $\CZ^+\simeq \CZ_\emptyset\times \Ran$. 

\medskip

A counital structure on $\CZ_\Ran$ gives rise to a unital-in-correspondences
structure with $\CZ^+_\Ran\to \CZ_\Ran$ being an isomorphism.

\sssec{}

Let $\CT$ be factorization space over $X$. There is a natural notion of unital-in-correspondences structure
on $\CT$ (i.e., $\CT_\Ran$ has a unital-in-correspondences structure, compatible with factorization, and we stipulate
$\CT_\emptyset\simeq \on{pt}$). 

\medskip

Let $\CT^+$ be the corresponding counital factorization space (see \secref{sss:Z+}); in this case we will
say that $\CT$ has a unital-in-correspondences structure 
\emph{relative} to $\CT^+$. We will refer to $\CT^+$ as the counital factorization space
\emph{underlying} $\CT$. 

\sssec{}

Let $\CT$ be a unital factorization space. Then $\CT$ acquires a natural unital-in-correspondences structure,
for which $\CT^+\to \on{pt}$ is an isomorphism. 

\medskip

According to \secref{sss:Z+}, if $\CT$ is counital as a factorization space, it acquires a natural structure of
unitality-in-correspondences with $\CT^+\to \CT$ being an isomorphism. 

\sssec{} \label{sss:unital loops}

Let $\CY$ be an affine D-scheme over $X$. We claim that the factorization space $\fL_\nabla(\CY)$ 
(see \secref{sss:forming loops aff}) possesses a natural unital-in-correspondences structure relative to
$\fL^+_\nabla(\CY)$.

\medskip

Namely, let $\ul{x}\subseteq \ul{x}':S\to \Ran$ be a pair of $S$-points of $\Ran$. Consider the
corresponding prestack $\cD_{\ul{x}',\nabla}$. It contains $\on{Graph}_{\ul{x}}$ as a closed subset.

\medskip

By definition, a lift of $\ul{x}\subseteq \ul{x}'$ to an $S$-point of 
$\fL_\nabla(\CY)^{\on{all}\rightsquigarrow +}_{\ul{x}\subseteq \ul{x}'}$ is a $X_\dr$-map
$$(\cD_{\ul{x}',\nabla}-\on{Graph}_{\ul{x}})\to \CY_\nabla.$$

\medskip

\noindent NB: in the particular case of $\CT=\fL_\nabla(\CY)$, we use the notation 
$\fL_\nabla^{\mer\rightsquigarrow \reg}(\CY)_{\ul{x}\subseteq \ul{x}'}$ instead of 
$\fL_\nabla(\CY)^{\on{all}\rightsquigarrow +}_{\ul{x}\subseteq \ul{x}'}$. 

\medskip

In order to define the composition of morphisms, we need to establish an isomorphism 
$$\fL_\nabla^{\mer\rightsquigarrow \reg}(\CY)_{\ul{x}_1\subseteq \ul{x}_2}
\underset{\fL_\nabla(\CY)^\mer_{\ul{x}_2}}\times \fL_\nabla^{\mer\rightsquigarrow \reg}(\CY)_{\ul{x}_2\subseteq \ul{x}_3}
\simeq \fL_\nabla^{\mer\rightsquigarrow \reg}(\CY)_{\ul{x}_1\subseteq \ul{x}_3}$$
for $\ul{x}_1\subseteq \ul{x}_2\subseteq \ul{x}_3$.  This follows from the isomorphism
$$(\cD_{\ul{x}_2,\nabla}-\ul{x}_1)\underset{\cD_{\ul{x}_2,\nabla}-\ul{x}_2}\sqcup (\cD_{\ul{x}_3,\nabla}-\ul{x}_2)\simeq (\cD_{\ul{x}_3,\nabla}-\ul{x}_1),$$
cf. \lemref{l:glue discs}.

%

%
%
%
%
%
%
%
%
%
%

\sssec{} \label{sss:LS partial puncture}

We will now define a unital-in-correspondences structure on the factorization space $\LS^\mer_H$
from \secref{sss:LS punctured}. 

\medskip

Let $\ul{x}\subseteq \ul{x}':S\to \Ran$ be as above. We consider $\fL_\nabla(\on{Conn}(\fh))_{\ul{x}'}$ as
acted on by $\fL(H)_{\ul{x}'}$. Similarly, $\fL^{\mer\rightsquigarrow \reg}_\nabla(\on{Conn}(\fh))_{\ul{x}\subseteq \ul{x}'}$
is acted on by $$\fL^{\mer\rightsquigarrow \reg}(H)_{\ul{x}\subseteq \ul{x}'}.$$

\medskip

Let $\LS^{\mer\rightsquigarrow \reg}_{H,\ul{x}\subseteq \ul{x}'}$
be the prestack equal to the \'etale sheafification of the (non-sheafified) quotient of 
$\fL^{\mer\rightsquigarrow \reg}_\nabla(\on{Conn}(\fh))_{\ul{x}\subseteq \ul{x}'}$
by $\fL^{\mer\rightsquigarrow \reg}(H)_{\ul{x}\subseteq \ul{x}'}$. 

\medskip

The composition of morphisms is defined as in \secref{sss:unital loops}.

\medskip

Note that the underlying counital factorization space of $\LS^\mer_H$ is $\LS^\reg_H$. 

\sssec{} \label{sss:unital on mf}

We now define a unital-in-correspondences structure on $\Op_\cG^\mf$. This is, however, automatic since
$$\Op_\cG^\mf:=\Op_\cG^\mer\underset{\LS_\cG^\mer}\times \LS_\sG^\reg,$$
so the unital-in-correspondences structure on $\Op_\cG^\mer$, $\LS_\cG^\mer$ and $\LS_\cG^\reg$
induces one on $\Op_\cG^\mf$.

\medskip

Note that by construction, the underlying counital factorization space of $\Op_\cG^\mf$ is $\Op^\reg_\cG$.

\ssec{Unital factorization categories}

\sssec{}  \label{sss:fact cat untl}

Let $\bA$ be a factorization category on $X$. A unital structure on $\bA$ is an extension of 
the crystal of categories $\ul\bA$ over $\Ran$ to a crystal of categories 
over $\Ran^{\on{untl},*}$ in a way compatible with factorization, i.e., we extend the identifications \eqref{e:fact category} to 
\begin{equation} \label{e:fact category untl}
\on{union}^*(\ul\bA)|_{(\Ran^{\on{untl},*}\times \Ran^{\on{untl},*})_{\on{disj}}}\simeq 
\ul\bA\boxtimes \ul\bA|_{(\Ran^{\on{untl},*}\times \Ran^{\on{untl},*})_{\on{disj}}}.
\end{equation} 

In addition, we stipulate
\begin{equation} \label{e:untl fact cat emptyset}
\bA_\emptyset\simeq \Vect,
\end{equation} 
so that this identification behaves (homotopically coherently) as a unit for the identifications \eqref{e:fact category untl}, i.e., the functor
$$\ul\bA\to \bA_\emptyset\otimes \ul\bA,$$
obtained by restricting \eqref{e:fact category untl} to 
$$\{\emptyset\}\times \Ran^{\on{untl},*}\to (\Ran^{\on{untl},*}\times \Ran^{\on{untl},*})_{\on{disj}}$$
identifies via \eqref{e:untl fact cat emptyset} with the identity functor.

\medskip

We let
$$\one_\bA\in \Gamma^{\on{strict}}(\Ran^{\on{untl},*},\ul\bA)$$
the canonical object, whose value at any $\ul{x}\in \Ran^{\on{untl},*}$
is 
$$\on{ins.unit}_{\emptyset\subset \ul{x}}(k),$$
where:
\begin{itemize}

\item $k\in \Vect\simeq \bA_\emptyset$;

\item $\on{ins.unit}_{\emptyset\subset \ul{x}}$ is the functor $\bA_\emptyset\to \bA_{\ul{x}}$
corresponding to the unique morphism $\{\emptyset\}\to \ul{x}$. 

\end{itemize} 

\sssec{} \label{sss:unit fact cat untl} 

Note that the factorization category $\Vect$ from \secref{sss:unit fact cat} admits a tautological unital structure. Namely,
the underlying crystal of categories on $\Ran^{\on{untl},*}$ is $\ul\Dmod(\Ran^{\on{untl},*})$. 

\sssec{} \label{sss:lax vs strict unital}

Given a pair of unital factorization categories, we can talk about
lax unital or strictly unital functors between them, compatible with factorization,
see \secref{sss:lax vs strict functors}. 

\medskip

We denote the resulting (2-)categories by
$$\on{FactCat}^{\on{untl,lax}}(X) \text{ and } \on{FactCat}^{\on{untl,strict}}(X),$$
respectively. 

\sssec{}

Pointwise tensor product defines a symmetric monoidal structure on the category of unital factorization categories
(for both variants: lax or strictly unital functors).

\medskip

The unit for the above symmetric monoidal structure is $\Vect$. 

\sssec{} \label{sss:duality lax untl}

Let $\bA$ be a unital factorization category. We shall say that $\bA$ is dualizable if it is a 
dualizable as an object of the above category, with \emph{lax} unital functors as morphisms.

\medskip 

In this case, the evaluation and the coevaluation functors
$$\bA\otimes \bA^\vee\to \Vect \text{ and } \Vect\to \bA^\vee\otimes \bA$$
carry lax unital structures.

\sssec{Example} \label{sss:com fact cat untl}

Let $\ul\bA_X$ be as in \secref{sss:com fact cat}. The explicit construction of the (symmetric monoidal)
factorization category $\on{Fact}(\ul\bA_X)$ shows that it admits a natural unital structure.

\medskip

We will denote by the same symbol $\on{Fact}(\ul\bA_X)$ the resulting (symmetric monoidal)
unital factorization category. 

\sssec{}

Given a unital factorization category $\bA$, we can talk about unital factorization algebras in it: by
definition, a unital factorization algebra $\CA$ in $\bA$ in it a lax unital factorization functor $\Vect\to \bA$. 

\medskip

Explicitly, the datum of $\CA$ is an object
$$\CA_{\Ran^{\on{untl},*}}\in \Gamma^{\on{lax}}(\Ran^{\on{untl},*},\ul\bA),$$
which is compatible with factorization in the natural sense (i.e., combine the ideas from
Sects. \ref{sss:untl fact alg} and \ref{sss:fact alg in fact cat}). 

\medskip

We denote by 
$$\on{FactAlg}^{\on{untl}}(X,\bA)$$
the category of unital factorization algebras in $\bA$. 

\medskip

The object $\one_\bA$ has a natural structure of unital factorization algebra in $\bA$. Furthermore, for any $\CA\in \on{FactAlg}^{\on{untl}}(X,\bA)$
we have a canonically defined map 
$$\on{vac}_\CA:\one_\bA\to \CA,$$
which we will refer to as the \emph{unit} or \emph{vacuum} for $\CA$. 

\sssec{} \label{sss:Vect to A}

If $\bA$ is a unital factorization category, it admits a canonically defined (strictly) unital functor
$$\on{Vac}_\bA:\Vect\to \bA.$$

By a slight abuse of notation, we will denote this functor by $\one_\bA$; the image of 
$$k\in \on{FactAlg}^{\on{untl}}(X,\Vect)$$ under $\on{Vac}_\bA$ is $\one_\bA$. 

\sssec{} \label{sss:fact mod in A untl}

Given $\CA\in \on{FactAlg}^{\on{untl}}(X,\bA)$ and $\CZ\to \Ran$, we can talk about unital factorization
$\CA$-modules \emph{at} $\CZ$. Denote this category by
$$\CA\mod^{\on{fact}}(\bA)_\CZ.$$

\medskip

When we talk about \emph{non-unital} factorization $\CA$-modules, we will denote the corresponding category 
by 
$$\CA\mod^{\on{fact-n.u.}}(\bA)_\CZ.$$

We have a forgetful functor 
$$\CA\mod^{\on{fact}}(\bA)_\CZ\to \CA\mod^{\on{fact-n.u.}}(\bA)_\CZ$$
and \propref{p:unital vs non-unital modules} applies in the present context as well. 

\medskip

The assignment 
$$\CZ\rightsquigarrow \CA\mod^{\on{fact}}(\bA)_\CZ$$
is a crystal of categories over $\Ran$ that we will denote by $\CA\ul\mod^{\on{fact}}(\bA)$. This crystal of categories 
has a natural lax factorization structure.  

\medskip

We will denote the resulting lax factorization category by $\CA\mod^{\on{fact}}(\bA)$. 

\sssec{} \label{sss:image of unit fact alg}

Let $\Phi:\bA_1\to \bA_2$ be a lax unital functor between unital factorization categories. 
Then it naturally gives rise to a functor
$$\Phi:\on{FactAlg}^{\on{untl}}(X,\bA_1)\to \on{FactAlg}^{\on{untl}}(X,\bA_2).$$

\medskip

In particular, $\Phi(\one_{\bA_1})$ has a natural structure of factorization algebra in $\bA_2$,
and we have a map of factorization algebras
\begin{equation} \label{e:unit to unit}
\one_{\bA_2}\to \Phi(\one_{\bA_1}).
\end{equation} 

\sssec{}

Let $\ul\bA_X$ be as in \secref{sss:com fact cat untl}. Let $A\in \bA_X$ be a commutative algebra object. 
Then the object 
$$\on{Fact}(A)\in \on{ComAlg}(\on{FactAlg}(X,\on{Fact}(\ul\bA_X)))$$
from \secref{sss:com fact cat} naturally lifts to an object of 
$$\on{ComAlg}(\on{FactAlg}^{\on{untl}}(X,\on{Fact}(\ul\bA_X))).$$

By a slight abuse of notation, we will denote it by the same symbol $\on{Fact}(A)$. 

\sssec{}

By a similar token, one defines the notion of unital structure on a \emph{lax} factorization 
category (see \secref{sss:lax fact}). The entire preceding discussion equally applies to
unital lax factorization categories.

\sssec{Example} \label{sss:untl str on fact mod}

Let $\CA$ be a (not necessarily unital) factorization algebra. Recall the lax factorization category
$\CA\mod^{\on{fact}}$ (see \secref{sss:A mod is lax}). We claim that it acquires a natural unital structure.

\medskip

In order to define it, we need to provide the following data: for any $\CZ\to \Ran$, we need to
define a functor
\begin{equation} \label{e:untl structure on fact mod}
\on{ins.unit}_{\CZ}:\CA\mod^{\on{fact}}_\CZ\to \CA\mod^{\on{fact}}_{\CZ^{\subseteq}},
\end{equation} 
equipped with an appropriate associativity structure. This construction was already mentioned in \secref{sss:propagate modules}: 

\medskip

Let  
$\CM$ be an object of $\CA\mod^{\on{fact}}_\CZ$. The corresponding object 
$\on{ins.unit}_{\CZ}(\CM)\in \CA\mod^{\on{fact}}_{\CZ^{\subseteq}}$ is constructed as follows.

\medskip

Note that there is a canonical projection
$$\CZ^{\subseteq^2}:=(\CZ^\subseteq)^{\subseteq}\overset{\on{pr}_{\on{comp},\CZ}}\longrightarrow
\CZ^{\subseteq}, \quad (z,\ul{x}_1\subseteq \ul{x}_2)\mapsto (z,\ul{x}_2).$$

We let 
$$(\on{ins.unit}_{\CZ}(\CM))_{\CZ^{\subseteq^2}}:=(\on{pr}_{\on{comp},\CZ})^!(\CM_{\CZ^{\subseteq}}).$$

The factorization structure on $(\on{ins.unit}_{\CZ}(\CM))_{\CZ^{\subseteq^2}}$ against $\CA$ is induced
by that on $\CM_{\CZ^{\subseteq}}$.

\medskip

The factorization unit in $\CA\mod^{\on{fact}}$ is the object $\CA^{\on{fact}}$ from \secref{sss:A mod is lax}.

\medskip

Let $\phi:\CA_1\to \CA_2$ be a map of (non-unital) factorization algebras. Then the functor
$$\Res_\phi:\CA_2\mod^{\on{fact}}\to \CA_1\mod^{\on{fact}}$$
carries a natural lax unital structure. 

\sssec{} \label{sss:unital on oblv}

Let now $\CA$ be a unital factorization algebra. The above construction applies verbatim to the category
of unital factorization $\CA$-modules. 

\medskip

Consider the forgetful functor
$$\oblv_\CA:\CA\mod^{\on{fact}}\to \Vect$$
(see \eqref{e:forgetful functor fact mod})
as a factorization functor (where the left-hand side is a lax factorization category). We 
claim that it carries a naturally defined lax unital structure.

\medskip

In order to define it, we need to provide the following data: for any $\CZ\to \Ran$ and
$\CM\in \CA\mod^{\on{fact}}_\CZ$ we need to define a map
$$\on{pr}_{\on{small},\CZ}^!(\CM_\CZ)\to (\on{ins.unit}_{\CZ}(\CM))_{\CZ^{\subseteq}}.$$

We note, however, that by construction
$$(\on{ins.unit}_{\CZ}(\CM))_{\CZ^{\subseteq}}:=\CM_{\CZ^{\subseteq}}.$$

Now, the sought-for map
$$\on{pr}_{\on{small},\CZ}^!(\CM_\CZ)\to \CM_{\CZ^{\subseteq}}$$
is exactly provided by the structure on $\CM$ of unital factorization $\CA$-module. 

\sssec{} \label{sss:fact A mod in A untl}

By a similar token, given a (non-unital) lax factorization category $\bA$ and a factorization algebra $\CA$ in it, the lax 
factorization category $\CA\mod^{\on{fact}}(\bA)$ (see \secref{sss:A mod in A}) acquires a naturally defined unital structure. 
The object 
$$\CA^{\on{fact}}\in \on{FactAlg}(X,\CA\mod^{\on{fact}}(\bA))$$
from \secref{sss:A mod in A} extends to an object of $\on{FactAlg}^{\on{untl}}(X,\CA\mod^{\on{fact}}(\bA))$, 
and equals in fact the unit in $\CA\mod^{\on{fact}}(\bA)$, i.e., the map
$$\one_{\CA\mod^{\on{fact}}}\to \CA^{\on{fact}}$$
is an isomorphism. 

\medskip

For a factorization functor $\Phi:\bA_1\to \bA_2$ between non-unital lax factorization categories and $\CA_1\in \on{FactAlg}(\bA_1,X)$,
the resulting functor
\begin{equation}  \label{e:fact A mod in A untl 1}
\Phi:\CA_1\mod^{\on{fact}}(\bA_1) \to \Phi(\CA_1)\mod^{\on{fact}}(\bA_2)
\end{equation}
has a natural (stictly) unital structure.  

\medskip

Similarly, if $\bA$ is unital and $\CA$ is unital, then the lax factorization category $\CA\mod^{\on{fact}}(\bA)$
(of \emph{unital} $\CA$-modules) acquires a naturally defined unital structure. Furthermore, in this case 
the forgetful functor
$$\oblv_\CA:\CA\mod^{\on{fact}}(\bA)\to \bA$$
carries a naturally defined lax unital structure.

\medskip

For a lax unital functor $\Phi:\bA_1\to \bA_2$ between unital lax factorization categories and $\CA_1\in \on{FactAlg}^{\on{untl}}(\bA_1,X)$,
the resulting functor
\begin{equation}  \label{e:fact A mod in A untl 2}
\Phi:\CA_1\mod^{\on{fact}}(\bA_1) \to \Phi(\CA_1)\mod^{\on{fact}}(\bA_2)
\end{equation}
has a natural (stictly) unital structure.  

\sssec{} \label{sss:A and A' mod bis} 

Let $\bA$ be a (non-unital) lax factorization category, and let $\phi:\CA\to \CA'$ be a homomorphism between factorization
(non-unital) algebras in it. 

\medskip

Since the functor
$$\Res_\phi:\CA'\mod^{\on{fact}}(\bA)\to \CA\mod^{\on{fact}}(\bA)$$
is lax unital and $\CA'{}^{\on{fact}}$ is the factorization unit in $\CA'\mod^{\on{fact}}(\bA)$, we obtain that the
factorization algebra 
$\Res_\phi(\CA'{}^{\on{fact}})\in \on{FactAlg}(X,\CA\mod^{\on{fact}}(\bA))$ acquires a natural unital structure. 

\medskip

Moreover, the functor \eqref{e:fact A mod in A untl 2} applied to $\bA_1=\CA'\mod^{\on{fact}}(\bA)$, $\bA_2=\CA\mod^{\on{fact}}(\bA)$
and $\Phi=\Res_\phi$ gives rise to a unital factorization functor
\begin{equation}  \label{e:A and A' mod bis}
\CA'\mod^{\on{fact}}(\bA)\to \Res_\phi(\CA'{}^{\on{fact}})\mod^{\on{fact-untl}}(\CA\mod^{\on{fact}}(\bA)). 
\end{equation}

Unwinding the definitions, we obtain: 

\begin{lem} \label{l:two algebras mods non-untl}
The functor \eqref{e:A and A' mod bis} is an equivalence. 
\end{lem}

\sssec{} \label{sss:A and A' mod untl} 

Let $\bA$ be a unital factorization category, and let $\phi:\CA\to \CA'$ be a homomorphism between 
unital factorization algebras in $\bA$. 

\medskip

Consider the unital lax factorization categories $\CA\mod^{\on{fact}}(\bA)$ and $\CA'\mod^{\on{fact}}(\bA)$,
and the restriction functor
$$\Res_\phi:\CA'\mod^{\on{fact}}(\bA)\to \CA\mod^{\on{fact}}(\bA).$$

This functor carries a natural lax unital structure. In particular, the object 
$$\Res_\phi(\CA'{}^{\on{fact}})\in \on{FactAlg}(\CA\mod^{\on{fact}}(\bA))$$
from \secref{sss:A and A' mod} lifts to an object of $\on{FactAlg}^{\on{untl}}(\CA\mod^{\on{fact}}(\bA))$.

\medskip

Note that the forgetful functor $\oblv_\CA:\CA\mod^{\on{fact}}(\bA)\to \bA$
sends $\Res_\phi(\CA'{}^{\on{fact}})\to \CA'$. In particular, it induces a functor
\begin{equation} \label{e:two algebras mods}
\Res_\phi(\CA'{}^{\on{fact}})\mod^{\on{fact}}(\CA\mod^{\on{fact}}(\bA))\to \CA'\mod^{\on{fact}}(\bA).
\end{equation} 

\medskip

The following is obtained by unwinding the definitions:

\begin{lem} \label{l:two algebras mods}
The functor \eqref{e:two algebras mods} is an equivalence.
\end{lem}

In fact, the inverse of the functor \eqref{e:two algebras mods} is given by the functor \eqref{e:fact A mod in A untl 2}
for $\bA_1=\CA'\mod^{\on{fact}}(\bA)$, $\bA_2=\CA\mod^{\on{fact}}(\bA)$, $\Phi=\Res_\phi$ and $\CA_1=\CA'{}^{\on{fact}}$.

\sssec{Proof of \propref{p:quasi-untl}} \label{sss:proof quasi-unital}

We will explicitly construct an inverse functor. 

\medskip

Let $k\overset{\on{vac}_\CA}\to\CA$ be a quasi-unital factorization algebra. Consider the unital
lax factorization category\footnote{The appearance of ``n.u." in the superscript in the next formula is
meant to emphasize that we are dealing with non-unital factorization modules, even though since $\CA$ is non-unital, 
we cannot even talk about unital modules.}
$\CA\mod^{\on{fact-n.u.}}$. It is equipped with a lax unital factorization functor
$$\Res_{\on{vac}_\CA}:\CA\mod^{\on{fact-n.u.}}\to k\mod^{\on{fact-n.u.}}.$$

Consider the fiber product
$$\CA\mod^{\on{fact-q.u.}}:=\Vect\underset{k\mod^{\on{fact-n.u.}}}\times \CA\mod^{\on{fact-n.u.}}\subset CA\mod^{\on{fact-n.u.}}.$$
It has a natural factorization structure, and the fact that 
$$\Res_{\on{vac}_\CA}(\CA^{\on{fact}_\Ran})\in k\mod^{\on{fact}}_\Ran\subset k\mod^{\on{fact-n.u.}}_\Ran$$
implies that the unital structure on $\CA\mod^{\on{fact-n.u.}}$ induces one on $\CA\mod^{\on{fact-q.u.}}$. The 
factorization unit $\one_{\CA\mod^{\on{fact-q.u.}}}$ in $\CA\mod^{\on{fact-q.u.}}$ is $\CA^{\on{fact}}$ (see \secref{sss:untl str on fact mod}). 

\medskip

Restricting $\Res_{\on{vac}_\CA}$ to $\CA\mod^{\on{fact-q.u.}}$, we obtain a lax unital factorization functor
\begin{equation} \label{e:quasi-untl to Vect}
\CA\mod^{\on{fact-q.u.}}\to \Vect.
\end{equation} 

The image of the factorization unit $\one_{\CA\mod^{\on{fact-q.u.}}}=\CA^{\on{fact}}$ under the above functor
is a unital factorization algebra (whose underlying plain factorization algebra is $\CA$ itself). This defines a functor
\begin{equation} \label{e:quasi-untl to untl}
\on{FactAlg}^{\on{q-untl}}(X)\to \on{FactAlg}^{\on{untl}}(X),
\end{equation} 
which commutes with a forgetful functor to $\on{FactAlg}(X)$.

\medskip

Let us show that the functors \eqref{e:untl to q-untl} and \eqref{e:quasi-untl to untl} are mutually inverse.
We first show that the composition \eqref{e:quasi-untl to untl}$\circ$\eqref{e:untl to q-untl} is isomorphic to
the identity functor.

\medskip

Indeed, when $\CA$ is unital, by \propref{p:unital vs non-unital modules}, the unital factorization category $\CA\mod^{\on{fact-q.u.}}$ identifies with 
$\CA\mod^{\on{fact}}$, and the functor \eqref{e:quasi-untl to untl} identifies with $\oblv_\CA$, equipped with its
natural lax unital structure. Hence, in this case, the image of $\CA^{\on{fact}}$ under \eqref{e:quasi-untl to untl}
identifies with $\CA$ as a unital factorization algebra. 

\medskip

Vice versa, let us start with a quasi-unital factorization algebra $k\overset{\on{vac}_\CA}\to\CA$. We have a commutative diagram
of unital lax factorization categories and lax unital functors, with the horizontal arrows being strict: 
$$
\CD
\CA\mod^{\on{fact-q.u.}} @>>> \CA\mod^{\on{fact-n.u.}} \\
@VVV @VVV \\
\Vect @>>> k\mod^{\on{fact-n.u.}}. 
\endCD
$$

Applying each circuit to $\one_{\CA\mod^{\on{fact-q.u.}}}$ we obtain a factorization algebra in $k\mod^{\on{fact-n.u.}}$,
equipped with a homomorphism from $k$. For the clockwise circuit, we obtain the original $k\overset{\on{vac}_\CA}\to\CA$.
For the anti-clockwise circuit, we obtain \eqref{e:quasi-untl to untl}$(k\overset{\on{vac}_\CA}\to\CA)$, equipped with a map
from $k$ to it, given by its unital structure.  

\qed[\propref{p:quasi-untl}]

\sssec{} \label{sss:adj of unital}

In the rest of this subsection we will focus on strict (i.e., non-lax) factorization categories.

\medskip

Let $\bA_1$ and $\bA_2$ be a pair of unital factorization categories, and let $\Phi:\bA_1\to \bA_2$
be a strictly unital functor. Assume that $\Phi$ admits a right adjoint $\Phi^R$ as a functor between the underlying 
crystals of categories on $\Ran$. 

\medskip

Then $\Phi^R$, viewed as a factorization functor between
factorization categories admits a natural extension to a lax unital functor between unital factorization
categories, see \secref{sss:right adj lax}. 

\sssec{}

We claim:

\begin{lem} \label{l:unit determines strict}
Let $\Phi:\bA_1\to \bA_2$ be a lax unital factorization functor between unital factorization categories.
Then $\Phi$ is strictly unital if and only if the map \eqref{e:unit to unit} is an isomorphism.
\end{lem} 

\begin{proof}

The ``only if" direction is tautological. Let us prove the ``if" direction, so let us assume that 
\eqref{e:unit to unit} is an isomorphism. 

\medskip

Let $\ul{x}\subseteq \ul{x}'$ be an arrow in $\Ran^{\on{untl},*}(S)$ for $S\in \affSch_{/\Ran}$. 
We need to show that the natural transformation
\begin{equation} \label{e:unit determines strict}
\on{ins.unit}_{\bA_2,\ul{x}\subseteq\ul{x}'}\circ \Phi_{\ul{x}}\to \Phi_{\ul{x}'}\circ \on{ins.unit}_{\bA_1,\ul{x}\subseteq\ul{x}'},
\end{equation} 
given by the lax unital structure on $\Phi$, is an isomorphism. 

\medskip

This asssertion can be checked strata-wise, so we can assume that $\ul{x}$ and $\ul{x}'$ are field-valued points.
Write 
$$\ul{x}'=\ul{x}\sqcup \ul{x}''.$$

We have
$$\bA_{i,\ul{x}'}\simeq \bA_{i,\ul{x}}\otimes \bA_{i,\ul{x}''}, \quad i=1,2$$
and the natural transformation \eqref{e:unit determines strict}
is the tensor product of the identity endomorphism 
of the functor
$$\Phi_{\ul{x}}:\bA_{1,\ul{x}}\to \bA_{2,\ul{x}}$$
and the natural transformation \eqref{e:unit determines strict} for $\emptyset \subset \ul{x}$. 

\medskip

However, the latter is exactly the map
$$\one_{\bA_2,\ul{x}}\to \Phi_{\ul{x}}(\one_{\bA_1,\ul{x}}).$$

\end{proof} 

\ssec{Examples of unital factorization categories arising from algebraic geometry}

\sssec{}

Let $\CT$ be a counital factorization space. We claim that the lax factorization category 
$\QCoh(\CT)$ admits a natural unital structure.

\medskip

Indeed, for a pair of $S$-points $\ul{x}\subseteq \ul{x}'$ of $\Ran^{\on{untl},*}$, the corresponding functor
$$\QCoh(\CT_{S,\ul{x}})\to \QCoh(\CT_{S,\ul{x}'})$$
is given by pullback along 
\begin{equation} \label{e:counital fibers}
\CT_{S,\ul{x}'}\to \CT_{S,\ul{x}}.
\end{equation} 

\medskip

In particular, the unit $\one_{\QCoh(\CT)}$ is the structure sheaf 
$$\CO_\CT\in \QCoh(\CT).$$

\sssec{}

By  a similar token, using \secref{sss:IndCoh ! Ran} the (lax) factorization category
$\IndCoh^!(\CT)$ admits a natural unital structure.

\medskip

Assume now that $\CT$ is affine and placid. Recall that according to \secref{sss:IndCoh * Ran}, we can consider the
factorization category $\IndCoh^*(\CT)$.

\medskip 

Assume now that the maps 
$$\CT_{S,X^{I_2}}\to \CT_{S,X^{I_1}}$$
for inclusions of finite sets $I_1\hookrightarrow I_2$ are of finite Tor-dimension (e.g., they are flat). In this 
case, by \secref{sss:* pullback on IndCoh* finite Tor dim}, the functors of *-pullback along the maps \eqref{e:counital fibers}
are defined for $\IndCoh^*(-)$
$$\IndCoh^*(\CT_{S,\ul{x}})\to \IndCoh^*(\CT_{S,\ul{x}'}).$$

\medskip

Hence, we obtain that in this case, $\IndCoh^*(\CT)$ also admits a natural unital structure.

\sssec{}

Let now $\CT$ be a unital factorization space. We claim that in this case the factorization category 
$\QCoh_{\on{co}}(\CT)$ has a natural unital structure. 

\medskip

Indeed, for a pair of $S$-points $\ul{x}\subseteq \ul{x}'$ of $\Ran^{\on{untl},*}$, the corresponding functor
$$\QCoh_{\on{co}}(\CT_{S,\ul{x}})\to \QCoh_{\on{co}}(\CT_{S,\ul{x}'})$$
is given by pushforward along 
\begin{equation} \label{e:unital fibers}
\CT_{S,\ul{x}}\to \CT_{S,\ul{x}'}.
\end{equation} 

\medskip

Assume now that $\CT$ is an ind-placid factorization ind-scheme, so that we can consider the factorization category
$\IndCoh^*(\CT)$.

\medskip

Taking the $\IndCoh$-pushforwards along \eqref{e:unital fibers} we obtain that $\IndCoh^*(\CT)$ 
acquires a natural unital structure. 

\sssec{Example}

We obtain that the factorization categories $\QCoh_{\on{co}}(\Gr_G)$ and $\IndCoh(\Gr_G)$ 
acquire a natural unital structure.

\medskip

By a similar mechanism, the factorization category $\Dmod(\Gr_G)$ also acquires a unital 
structure. 

\sssec{}

Let us now in addition assume that for an injection of finite sets $I_1\hookrightarrow I_2$,
the corresponding map
$$X^{I_2}\underset{X^{I_1}}\times \CT_{X^{I_1}}\to \CT_{X^{I_2}}$$
is an ind-closed embedding locally almost of finite-presentation.

\medskip

In this case, by \secref{sss:closed emb placid}, the $\IndCoh$-pushforward functors along \eqref{e:unital fibers}
are defined on $\IndCoh^!(-)$:
$$\IndCoh^!(\CT_{S,\ul{x}})\to \IndCoh^!(\CT_{S,\ul{x}'}).$$

Hence, we obtain that in this case, $\IndCoh^!(\CT)$ also acquires a unital structure.

\sssec{}

Let $\CT$ be a factorization space over $X$, equipped with a unital-in-correspondences structure. 
Assume that for a pair of $S$-points $\ul{x}\subseteq \ul{x}'$ of $\Ran^{\on{untl},*}$, the map
\begin{equation} \label{e:pr T small maps}
\CT^{\on{all}\rightsquigarrow +}_{\ul{x}\subseteq \ul{x}'}\overset{\on{pr}^\CT_{\on{small}}}\longrightarrow 
\CT_{\ul{x}}
\end{equation} 
is affine. 

\medskip

We claim that in this case, the factorization category $\QCoh_{\on{co}}(\CT)$ acquires a unital structure,
and the functor
$$\QCoh(\CT^+)=\QCoh_{\on{co}}(\CT^+)\to \QCoh_{\on{co}}(\CT),$$
given by direct image along 
$$\iota:\CT^+\to \CT$$ is strictly unital. In particular, the unit 
$\one_{\QCoh_{\on{co}}(\CT)}$ is the direct image of $\CO_{\CT^+}$ along $\iota$. 

\medskip

Namely, for a pair of $S$-points $\ul{x}\subseteq \ul{x}'$ of $\Ran^{\on{untl},*}$, the corresponding functor
$$\QCoh_{\on{co}}(\CT_{S,\ul{x}})\to \QCoh_{\on{co}}(\CT_{S,\ul{x}'})$$
is given by pull-push along
\begin{equation} \label{e:pr T small big maps}
\CT_{\ul{x}} \overset{\on{pr}^\CT_{\on{small}}}\longleftarrow 
\CT^{\on{all}\rightsquigarrow +}_{\ul{x}\subseteq \ul{x}'}\overset{\on{pr}^\CT_{\on{big}}}\longrightarrow 
\CT_{\ul{x}'}.
\end{equation} 

\sssec{Example}

Thus, we obtain that for an affine D-scheme $\CY$, the factorization category 
$$\QCoh_{\on{co}}(\fL_\nabla(\CY))$$ 
acquires a unital structure, with the unit $\one_{\QCoh_{\on{co}}(\fL_\nabla(\CY))}$ being 
the direct image of $\CO_{\fL^+_\nabla(\CY)}$ along 
$$\fL^+_\nabla(\CY)\overset{\iota}\to \fL_\nabla(\CY).$$

\medskip

As another example, we can take $\CT=\Op^\mf_\cG$, and we obtain that $\QCoh_{\on{co}}(\Op^\mf_\cG)$
acquires a unital structure. The unit $\one_{\QCoh_{\on{co}}(\Op^\mf_\cG)}$ is the direct image of
$\CO_{\Op^\reg_\cG}$ along
$$\Op^\reg_\cG \overset{\iota^{+,\mf}}\to \Op^\mf_\cG.$$

\sssec{} \label{sss:corr *}

Assume now that $\CT$ is an ind-placid ind-affine factorization ind-scheme. Assume that the maps
$$\CT^{\on{all}\rightsquigarrow +}_{I_1\subseteq I_2}\overset{\on{pr}^\CT_{\on{small}}}\longrightarrow \CT_{X^{I_1}}$$
for an injection of finite sets $I_1\subseteq I_2$ are affine and of finite Tor-dimension.

\medskip

In this case, by \secref{sss:* pullback on IndCoh* finite Tor dim}, the functors of *-pullback along
the maps \eqref{e:pr T small maps} are well-defined on $\IndCoh^*(-)$ and satisfy base change. 

\medskip

We define a unital structure on the factorization category 
$\IndCoh^*(\CT)$ by ($\IndCoh$,*)-pull followed by $\IndCoh$-pushforward along \eqref{e:pr T small big maps}.

\medskip

Note that by construction, the functor of $\IndCoh$-pushforward along $\iota$
$$\IndCoh^*(\CT^+)\to \IndCoh^*(\CT)$$
is (strictly) unital.

\medskip

In particular, the unit in $\IndCoh^*(\CT)$ is the direct image of $\CO_{\CT^+}\in \IndCoh^*(\CT^+)$ 
along $\iota$. 

\sssec{} \label{sss:corr !}

Let us continue to assume that $\CT$ is an ind-placid ind-affine factorization ind-scheme. Assume now that
for an inclusion of finite sets $I_1\subseteq I_2$, the map
$$\CT^{\on{all}\rightsquigarrow +}_{I_1\subseteq I_2}\overset{\on{pr}^\CT_{\on{big}}}\longrightarrow \CT_{X^{I_2}}$$
is an ind-closed embedding locally almost of finite-presentation. 

\medskip

In this case, \secref{sss:closed emb placid}, the $\IndCoh$-pushforward functors along 
$$\CT^{\on{all}\rightsquigarrow +}_{\ul{x}\subseteq \ul{x}'}\overset{\on{pr}^\CT_{\on{big}}}\longrightarrow 
\CT_{\ul{x}'}$$
are defined on $\IndCoh^!(-)$ and satisfy base change. 

\medskip

We define a unital structure on the factorization category 
$\IndCoh^!(\CT)$ by !-pull followed by $\IndCoh$-pushforward along \eqref{e:pr T small big maps}.

\medskip

Note that by construction, the functor of $\IndCoh$-pushforwardforward along $\iota$
$$\IndCoh^!(\CT^+)\to \IndCoh^!(\CT)$$
is (strictly) unital.

\medskip

In particular, the unit in $\IndCoh^!(\CT)$ is the direct image of $\omega_{\CT^+}\in \IndCoh^!(\CT^+)$ 
along $\iota$. 

\sssec{Example} \label{sss:placid for unital}

Let $\CY$ be an affine D-scheme, such that:

\begin{itemize}

\item $\fL_\nabla(\CY)$ is ind-placid;

\medskip

\item The maps 
$$\fL_\nabla^{\mer\rightsquigarrow \reg}(\CY)_{I_1\subseteq I_2}\overset{\on{pr}^\CY_{\on{small}}}\longrightarrow \fL_\nabla(\CY)_{X^{I_1}}$$
for $I_1\subseteq I_2$ are flat;

\medskip

\item The maps
$$\fL_\nabla^{\mer\rightsquigarrow \reg}(\CY)_{I_1\subseteq I_2}\overset{\on{pr}^\CY_{\on{big}}}\longrightarrow \fL_\nabla(\CY)_{X^{I_2}}$$
for $I_1\subseteq I_2$ are locally almost of finite-presentation. 

\end{itemize} 

\bigskip

This happens, e.g., for
$\CY=\on{Jets}(H)$, where $H$ is a smooth group-scheme over $X$. 

\medskip

We obtain that the factorization categories $\IndCoh^!(\fL_\nabla(\CY))$ and $\IndCoh^*(\fL_\nabla(\CY))$ 
acquire unital structures. 

\sssec{}

By a similar procedure, the factorization category $\Dmod(\fL(H))$ acquires a unital structure. 

\sssec{}

As yet another example, we obtain that the categories 
$$\IndCoh^!(\Op^\mf_\cG) \text{ and } \IndCoh^*(\Op^\mf_\cG)$$
acquire unital structures. 

\sssec{}

Let $\CT$ be an ind-placid ind-affine factorization ind-scheme. Assume that both additional conditions 
in Sects. \ref{sss:corr *} and \ref{sss:corr !} are satisfied, so both $\IndCoh^*(\CT)$ and 
$\IndCoh^!(\CT)$ acquire unital structures. 

\medskip

Note that the canonical pairing
\begin{equation} \label{e:dual cryst cat categ 3}
\IndCoh^!(\CT)\otimes \IndCoh^*(\CT)\to \Vect
\end{equation} 
(see \secref{sss:IndCoh ! * duality}) as factorization categories, 
admits a natural lax unital structure as a functor between 
unital factorization categories. 

\medskip

Moreover, unwinding the definitions we obtain that the condition from  
\secref{sss:dual cryst cat categ 3} is satisfied. 

\medskip

Hence, we obtain that \eqref{e:dual cryst cat categ 3} realizes $\IndCoh^!(\CT)$ and $\IndCoh^*(\CT)$
as each other duals as unital factorization categories. 

\ssec{Unital structure on Kac-Moody representations}

\sssec{}

Our current goal is to construct a unital structure on the factorization category $\hg\mod_\kappa$.
Let us be given a pair of $S$-points $\ul{x},\ul{x}'$ of $\Ran$ with $\ul{x}\subseteq \ul{x}'$. 
We need to construct a functor
\begin{equation} \label{e:unital structure on KM}
\hg\mod_{\kappa,S,\ul{x}}\to \hg\mod_{\kappa,S,\ul{x}'}.
\end{equation}

\medskip

The unital-in-correspondences structure on $\fL(G)$ gives rise to the following diagram of group ind-schemes over $S$:
$$\fL(G)_{\ul{x}} \overset{\on{pr}_{\on{small}}^G}\longleftarrow \fL^{\mer\rightsquigarrow \reg}(G)_{\ul{x}\subseteq \ul{x}'} 
\overset{\on{pr}_{\on{small}}^G}\longrightarrow \fL(G)_{\ul{x}'}.$$

Proceeding as in \secref{sss:KL defn}, we consider the corresponding categories
$$\hg\mod_{\kappa,S,\ul{x}},\,\, \hg\mod_{\kappa,S,\ul{x}'} \text{ and } \hg\mod_{\kappa,S,\ul{x}\subseteq\ul{x}'}.$$

The universal properties of these categories give rise to restriction functors
\begin{equation} \label{e:unital structure on KM 1}
\oblv^{\hg_{\ul{x}'}}_{\hg_{\ul{x}\subseteq\ul{x}'}}:\hg\mod_{\kappa,S,\ul{x}'} \to \hg\mod_{\kappa,S,\ul{x}\subseteq\ul{x}'}
\end{equation}
and 
\begin{equation} \label{e:unital structure on KM 2}
\oblv^{\hg_{\ul{x}}}_{\hg_{\ul{x}\subseteq\ul{x}'}}:\hg\mod_{\kappa,S,\ul{x}} \to \hg\mod_{\kappa,S,\ul{x}\subseteq\ul{x}'},
\end{equation}
both strongly compatible with the forgetful functors to $\Dmod(S)$ and the 
actions of $\fL^{\mer\rightsquigarrow \reg}(G)_{\ul{x}\subseteq \ul{x}'}$ at level $\kappa$. 

\sssec{}

One checks directly that the functor $\oblv^{\hg_{\ul{x}}}_{\hg_{\ul{x}\subseteq\ul{x}'}}$ admits a left adjoint.

\medskip

We define the functor
$$\on{ins.vac}_{\ul{x}\subseteq \ul{x}'}:\hg\mod_{\kappa,S,\ul{x}}\to \hg\mod_{\kappa,S,\ul{x}'}$$
as the composition 
\begin{equation} \label{e:unitality KM 0}
(\oblv^{\hg_{\ul{x}'}}_{\hg_{\ul{x}\subseteq\ul{x}'}})^L \circ \oblv^{\hg_{\ul{x}}}_{\hg_{\ul{x}\subseteq\ul{x}'}}.
\end{equation} 

\sssec{}

In order to upgrade the collection of functors $\on{ins.vac}_{\ul{x}\subseteq \ul{x}'}$ to a unital
structure, we need to construct associativity isomorphisms 
\begin{equation} \label{e:unitality KM 1}
\on{ins.vac}_{\ul{x}_3\subseteq \ul{x}_2}\circ \on{ins.vac}_{\ul{x}_1\subseteq \ul{x}_2}\simeq \on{ins.vac}_{\ul{x}_1\subseteq \ul{x}_3}
\end{equation} 
for $\ul{x}_1\subseteq \ul{x}_2 \subseteq \ul{x}_3$.

\medskip

Note that the inclusions
$$(\cD_{\ul{x}_2}-\ul{x}_1)\hookrightarrow (\cD_{\ul{x}_3}-\ul{x}_1) \hookleftarrow (\cD_{\ul{x}_3}-\ul{x}_2)$$
give rise to maps 
$$\fL^{\mer\rightsquigarrow \reg}(G)_{\ul{x}_1\subseteq \ul{x}_2}\leftarrow 
\fL^{\mer\rightsquigarrow \reg}(G)_{\ul{x}_1\subseteq \ul{x}_3}\to \fL^{\mer\rightsquigarrow \reg}(G)_{\ul{x}_3\subseteq \ul{x}_2}.$$

We have the corresponding functors
$$\hg\mod_{\kappa,S,\ul{x}_1\subseteq \ul{x}_2}
\overset{\oblv^{\hg_{\ul{x}_1\subseteq \ul{x}_2}}_{\hg_{\ul{x}_1\subseteq\ul{x}_3}}}\longrightarrow 
\hg\mod_{\kappa,S,\ul{x}_1\subseteq \ul{x}_3} \overset{\oblv^{\hg_{\ul{x}_2\subseteq \ul{x}_3}}_{\hg_{\ul{x}_1\subseteq\ul{x}_3}}}\longleftarrow 
\hg\mod_{\kappa,S,\ul{x}_2\subseteq \ul{x}_2}$$
and an isomorphism
$$\oblv^{\hg_{\ul{x}_1\subseteq \ul{x}_2}}_{\hg_{\ul{x}_1\subseteq\ul{x}_3}}\circ 
\oblv^{\hg_{\ul{x}_2}}_{\hg_{\ul{x}_1\subseteq\ul{x}_2}} \simeq 
\oblv^{\hg_{\ul{x}_2\subseteq \ul{x}_3}}_{\hg_{\ul{x}_1\subseteq\ul{x}_3}}\circ 
\oblv^{\hg_{\ul{x}_2}}_{\hg_{\ul{x}_2\subseteq\ul{x}_3}}.$$

From here we obtain a natural transformation
\begin{equation} \label{e:unitality KM 2}
(\oblv^{\hg_{\ul{x}_2\subseteq \ul{x}_3}}_{\hg_{\ul{x}_1\subseteq\ul{x}_3}})^L\circ 
\oblv^{\hg_{\ul{x}_1\subseteq \ul{x}_2}}_{\hg_{\ul{x}_1\subseteq\ul{x}_3}}\to
\oblv^{\hg_{\ul{x}_2}}_{\hg_{\ul{x}_2\subseteq\ul{x}_3}}\circ (\oblv^{\hg_{\ul{x}_2}}_{\hg_{\ul{x}_1\subseteq\ul{x}_2}})^L.
\end{equation} 

We claim that \eqref{e:unitality KM 2} is an isomorphism. Indeed, this follows from the fact that the diagram of Lie algebras
$$
\CD
\fL^{\mer\rightsquigarrow \reg}(\fg)_{\ul{x}_1\subseteq \ul{x}_3} @>>> \fL^{\mer\rightsquigarrow \reg}(\fg)_{\ul{x}_2\subseteq \ul{x}_3} \\
@VVV @VVV \\
\fL^{\mer\rightsquigarrow \reg}(\fg)_{\ul{x}_1\subseteq \ul{x}_2} @>>> \fL(\fg)_{\ul{x}_2}
\endCD
$$
is Cartesian. 

\medskip

Now, \eqref{e:unitality KM 1} follows
by precomposing both sides of \eqref{e:unitality KM 2} with $\oblv^{\hg_{\ul{x}_1}}_{\hg_{\ul{x}_1\subseteq\ul{x}_2}}$
and post-composing with $(\oblv^{\hg_{\ul{x}_3}}_{\hg_{\ul{x}_2\subseteq\ul{x}_3}})^L$. 

\sssec{}

Recall that the functor $\oblv^{\hg_{\ul{x}}}_{\hg_{\ul{x}\subseteq\ul{x}'}}$ is compatible with the action of
$\fL^{\mer\rightsquigarrow \reg}(G)_{\ul{x}\subseteq \ul{x}'}$. In particular, it is compatible with the action of
$\fL^+(G)_{\ul{x}'}$, which acts on $\hg\mod_{\kappa,S,\ul{x}}$ via $\fL^+(G)_{\ul{x}'}\to \fL^+(G)_{\ul{x}}$.
Hence, it induces a functor
$$\KL(G)_{\kappa,S,\ul{x}}=(\hg\mod_{\kappa,S,\ul{x}})^{\fL^+(G)_{\ul{x}}}
\to \left(\hg\mod_{\kappa,S,\ul{x}\subseteq\ul{x}'}\right)^{\fL^+(G)_{\ul{x}'}}.$$

\medskip

Since the functor $\oblv^{\hg_{\ul{x}'}}_{\hg_{\ul{x}\subseteq\ul{x}'}}$ is compatible with the action of 
$\fL^{\mer\rightsquigarrow \reg}(G)_{\ul{x}\subseteq \ul{x}'}$, then so is its left adjoint. In particular, 
we obtain that $(\oblv^{\hg_{\ul{x}'}}_{\hg_{\ul{x}\subseteq\ul{x}'}})^L $ is compatible with the action
of $\fL^+(G)_{\ul{x}'}$, and hence induces a functor
$$\left(\hg\mod_{\kappa,S,\ul{x}\subseteq\ul{x}'}\right)^{\fL^+(G)_{\ul{x}'}}\to
(\hg\mod_{\kappa,S,\ul{x}'})^{\fL^+(G)_{\ul{x}'}}=\KL(G)_{\kappa,S,\ul{x}'}.$$

\medskip

Hence, we obtain that the functor \eqref{e:unitality KM 0} induces a functor
\begin{equation} \label{e:unitality KM 3}
\on{ins.vac}_{\ul{x}\subseteq \ul{x}'}:\KL(G)_{\kappa,S,\ul{x}}\to \KL(G)_{\kappa,S,\ul{x}'}.
\end{equation} 

The functors \eqref{e:unitality KM 3} define a unital structure on $\KL(G)_\kappa$.

\sssec{}

By construction, the unit object $\one_{\hg\mod_\kappa}$ is the vacuum module $\on{Vac}(G)_\kappa$, i.e.,
$$\on{Vac}(G)_{\kappa,\ul{x}}:=(\oblv^{\hg_{\ul{x}}}_{\fL^+(\fg)_{\ul{x}}})^L\circ \oblv^0_{\fL^+(\fg)_{\ul{x}}}(k).$$
 
Furthermore, $\on{Vac}(G)_\kappa$ naturally upgrades to an object of $\KL(G)_\kappa$ and coincides with its unit.

\sssec{}

Let $\kappa'$ be as in \secref{ss:KL duality}.

\medskip

Recall that according to \cite[Sect. 9.16.11]{Ra5}, we have canonical pairings 
$$\hg\mod_{\kappa,\ul{x}}\otimes \hg\mod_{\kappa',\ul{x}}\to \Dmod(S)$$
making 
$$\hg\mod_\kappa \text{ and } \hg\mod_{\kappa'}$$
mutually dual as factorization categories.

\medskip

We claim that this duality extends to a duality between $\hg\mod_\kappa$ and $\hg\mod_{\kappa'}$
as unital factorization categories (see \secref{sss:dual cryst cat categ} for what this means). 

\medskip

Namely, imitating the construction in \cite[Sect. 9.16]{Ra5} we obtain a duality between
$$\hg\mod_{\kappa,S,\ul{x}\subseteq\ul{x}'} \text{ and }  \hg\mod_{\kappa',S,\ul{x}\subseteq\ul{x}'}.$$

Under this identification, the dual of the functor
$\oblv^{\hg_{\ul{x}}}_{\hg_{\ul{x}\subseteq\ul{x}'}}$ of \eqref{e:unital structure on KM 2} is the \emph{right adjoint} of 
$$\oblv^{\hg_{\ul{x}}}_{\hg_{\ul{x}\subseteq\ul{x}'}}:\hg\mod_{\kappa',S,\ul{x}} \to \hg\mod_{\kappa',S,\ul{x}\subseteq\ul{x}'},$$
and the dual of the functor $\oblv^{\hg_{\ul{x}'}}_{\hg_{\ul{x}\subseteq\ul{x}'}}$ of \eqref{e:unital structure on KM 1} is the 
left adjoint of the functor
$$\oblv^{\hg_{\ul{x}}}_{\hg_{\ul{x}\subseteq\ul{x}'}}:\hg\mod_{\kappa',S,\ul{x}'} \to \hg\mod_{\kappa',S,\ul{x}\subseteq\ul{x}'}.$$

From here we obtain the desired identification between the dual of 
$$\on{ins.vac}_{\ul{x}\subseteq \ul{x}'}:\hg\mod_{\kappa,S,\ul{x}}\to \hg\mod_{\kappa,S,\ul{x}'}$$
and the right adjoint of 
$$\on{ins.vac}_{\ul{x}\subseteq \ul{x}'}:\hg\mod_{\kappa',S,\ul{x}}\to \hg\mod_{\kappa',S,\ul{x}'}.$$

\sssec{}

In a similar fashion we obtain that the duality between
$$\KL(G)_\kappa \text{ and } \KL(G)_{\kappa'}$$
as factorization categories extends to a duality as unital factorization categories. 

\ssec{Unital factorization module categories}

\sssec{}

Let $\bA$ be a factorization category. Let $\CZ$ be a prestack mapping to $\Ran$, and let $\bC$
be a factorization module category over $\bA$ at $\CZ$.

\medskip

Suppose that $\bA$ is equipped with a unital structure. Combining the ideas of Sects. \ref{sss:fact cat untl},
\ref{sss:fact mod cat} and \ref{sss:untl fact mod}, we define the notion of unital structure on $\bC$. 

\sssec{}

Concretely, the unital structure on $\bC$ amounts to the following: given a pair of points
$(z,\ul{x})$ and $(z,\ul{x}')$ of $\CZ^{\subseteq}$ with $\ul{x}\subseteq \ul{x}'$, we must be given a functor
$$\on{ins.unit}_{\ul{x}\subseteq \ul{x}'}:\bC_{(z,\ul{x})}\to \bC_{(z,\ul{x}')},$$
compatible with factorization. 

\medskip

The latter compatibility means the following: for $\ul{x}'=\ul{x}\sqcup \ul{x}''$, with respect to the identification
$$\bC_{(z,\ul{x}')}\simeq \bA_{\ul{x}''}\otimes \bC_{(z,\ul{x})},$$
the functor $\on{ins.unit}_{\ul{x}\subseteq \ul{x}'}$ is
$$\one_{\bA_{\ul{x}''}}\otimes \on{Id}.$$

\medskip

In what follows we will denote by $\on{ins.unit}_\CZ$ the corresponding functor
\begin{equation} \label{e:ins unit module}
\bC_\CZ\to \bC_{\CZ^{\subseteq,\on{untl}}}.
\end{equation}

\sssec{}

For a unital factorization category $\bA$ and any $\CZ\to \Ran$, consider the factorization module
category $\bA^{\on{fact}_\CZ}$ from \secref{sss:vac fact mod cat}. 

\medskip

Unwinding the definitions, we obtain that $\bA^{\on{fact}_\CZ}$ carries a natural unital structure
(cf. \secref{sss:vacuum module untl}). 

\sssec{}

Given a pair of unital factorization module categories $\bC_1$ and $\bC_2$ over $\bA$ at $\CZ$
one can a priori talk about strictly unital or lax unital functors between them, compatible with factorization. However,
as in \lemref{l:unit determines strict} one shows that any lax unital functor between them is automatically
strictly unital.

\medskip

Thus, we can unambiguously talk about the (2-) category of unital factorization module categories over $\bA$ at $\CZ$.

\sssec{Notational convention} 

When $\bA$ is unital, we will denote the (2-) category
of unital factorization module categories over $\bA$ at $\CZ$ by
$$\bA\mmod^{\on{fact}}_\CZ.$$

\medskip

We will denote the category of plain (i.e., non-unital) factorization module categories over $\bA$ at $\CZ$ by
$$\bA\mmod^{\on{fact-n.u.}}_\CZ.$$

\begin{rem}
Unlike the case of modules over factorization algebras, the forgetful functor
\begin{equation} \label{e:oblv unit cat}
\bA\mmod^{\on{fact}}_\CZ\to \bA\mmod^{\on{fact-n.u.}}_\CZ
\end{equation}
is \emph{not} fully faithful.

\medskip

Indeed, take $\bA=\Vect$ and $\CZ=\on{pt}$ corresponding to a singleton $\{x\}=\ul{x}\in \Ran$. 
Take 
$$\bC=\Vect^{\on{fact}_x}\in \Vect\mmod^{\on{fact}}_x.$$

Then the category of endofunctors of $\bC$ as an object of $\Vect\mmod^{\on{fact}}_x$ is
$$k\mod^{\on{fact}}_x$$
(where $k$ is the unit factorization algebra), 
and the latter identifies with $\Vect$, see \secref{sss:unital fact mod for k}.

\medskip

By contrast, the category of endofunctors of $\bC$ as an object of $\Vect\mmod^{\on{fact-n.u.}}_x$ is
$$k\mod^{\on{fact-n.u.}}_x.$$

So, at the level of endofunctors of the above object, the forgetful functor \eqref{e:oblv unit cat}
is the forgetful functor 
$$k\mod^{\on{fact}}_x\to k\mod^{\on{fact-n.u.}}_x,$$
which is fully faithful, but not an equivalence. 

\end{rem} 

\sssec{Example}  \label{sss:untl fact mod Vect}

Recall the construction from \secref{sss:cats to fact mods over Vect}. It is easy to see that the resulting
$\Vect$-factorization module category $\bC$ is unital.

\medskip

In \secref{sss:tight vect} we will show that the functor
$$\on{CrystCat}(\CZ_0)\to \Vect\mmod^{\on{fact}}_{\CZ_0}$$
is fully faithful, and we will characterize its essential image.

\medskip

Note, however, that one categorical level down, the corresponding functor
was an equivalence, see \secref{sss:unital fact mod for k}.

\sssec{}

Let $\CA$ be a unital factorization algebra in a unital factorization category $\bA$. Let $\bC$ be a unital factorization
module category over $\bA$ at some $\CZ\to \Ran$. 

\medskip

Parallel to Sects. \ref{sss:modules inside modules} and \ref{sss:untl fact mod}, one defines the notion of 
unital factorization modules \emph{over} $\CA$ \emph{in} $\bC$. 

\medskip

We will denote the corresponding category by 
$$\CA\mod^{\on{fact}}(\bC)_\CZ.$$

By contrast, we will denote the category of plain (i.e., non-unital) $\CA$-modules in $\bC$ by
$$\CA\mod^{\on{fact-n.u.}}(\bC)_\CZ.$$

\sssec{}

We claim:

\begin{lem}  \label{l:unital modules for unit}
$\one_\bA\mod^{\on{fact}}(\bC)_\CZ\simeq \bC_\CZ$.
\end{lem}

\begin{proof}

The proof essentially repeats the contents of \secref{sss:unital fact mod for k}:

\medskip

Starting from an object $\CM'\in \bC_\CZ$, consider
$$\on{ins.unit}_\CZ(\CM')\in \bC_{\CZ^{\subseteq,\on{untl}}}.$$

This object has a tautological factorization structure against $\one_\bA$.

\medskip

Vice versa, starting from $\CM\in \one_\bA\mod^{\on{fact}}(\bC)_\CZ$, the unital structure on $\CM$
gives rise to a map
$$\on{ins.unit}_\CZ(\CM_\CZ)\to \CM_{\CZ^{\subseteq,\on{untl}}}$$
in $\bC_{\CZ^{\subseteq,\on{untl}}}$, compatible with factorization. 

\medskip

It suffices to show that the latter map is an isomorphism. This can be checked stratawise, in
which case it follows from factorization.

\end{proof} 

\sssec{} \label{sss:funct between pairs modules untl}

Let $\Phi:\bA_1\to \bA_2$ be a lax unital factorization functor between unital factorization categories.
Let $\bC_1$ and $\bC_2$ be unital factorization module categories over $\bA_1$ and $\bA_2$, respectively,
at some $\CZ\to \Ran$. 

\medskip

Mimicking \secref{sss:factorization restriction} we have the notion of lax unital functor
$$\Phi_m:\bC_1\to \bC_2,$$
compatible with factorization. Denote the category of such functors
$$\on{Funct}_{\bA_1\to \bA_2}(\bC_1,\bC_2).$$

\sssec{} \label{sss:funct between pairs modules untl alg}

Let $(\Phi,\Phi_m):(\bA_1,\bC_1)\to (\bA_2,\bC_2)$ be as above. Let $\CA_1\in \bA_1$ be a unital
factorization algebra, and consider its image $\Phi(\CA_1)\in \on{FactAlg}^{\on{untl}}(X,\bA_2)$. 

\medskip

Then the functor $\Phi_m$ induces a functor
\begin{equation} \label{e:funct between pairs modules untl alg}
\Phi_m:\CA_1\mod^{\on{fact}}(\bC_1)_\CZ\to \Phi(\CA_1)\mod^{\on{fact}}(\bC_2)_\CZ.
\end{equation} 

\sssec{}

As in \lemref{l:unit determines strict}, we have:

\begin{lem} \label{l:modules strict}
Suppose that $\Phi$ is strictly unital. Then the functor between crystals 
of categories on $\CZ^{\subseteq,\on{untl}}$ underlying every $\Phi_m\in \on{Funct}_{\bA_1\to \bA_2}(\bC_1,\bC_2)$
$$\ul\bC_1\to \ul\bC_2$$
is strict.
\end{lem} 

\sssec{} \label{sss:unital restr categ} 

Let $\Phi:\bA_1\to \bA_2$ be a strictly unital functor between unital factorization categories. 
Let $\bC_2$ be a unital module category over $\bA_2$ at some $\CZ\to \Ran$.

\medskip

In this case, it follows from the construction of the restriction functor $\Res_\Phi$ that 
$\Res_\Phi(\bC_2)$ possesses a natural unital structure, and the tautological functor
\begin{equation} \label{e:Res to original again}
\Res_\Phi(\bC_2)\to \bC_2
\end{equation}
admits a natural lax unital structure compatible with factorization.

\medskip 

Furthermore, the resulting object
$$\Res_\Phi(\bC_2)\in \bA_1\mmod^{\on{fact}}_\CZ$$
has a universal property parallel to that in the non-unital case:

\medskip

\begin{lem} \label{l:unital restr categ univ} 
For $\bC_1\in \bA_1\mmod^{\on{fact}}_\CZ$, composition with \eqref{e:Res to original again} defines 
an equivalence 
$$\on{Funct}_{\bA_1\mmod^{\on{fact}}_\CZ}(\bC_1,\Res_\Phi(\bC_2))\overset{\sim}\to \on{Funct}_{\bA_1\to \bA_2}(\bC_1,\bC_2).$$
\end{lem}

\sssec{}

Let $\CA_1\in \bA_1$ be a unital factorization algebra. As in \lemref{l:modules for fact alg restr}, we have:

\begin{lem} \label{l:modules for fact alg restr untl} 
The functor \eqref{e:Res to original again} induces an equivalence
$$\CA_1\mod^{\on{fact}}(\Res_\Phi(\bC_2))_\CZ\to \Phi(\CA_1)\mod^{\on{fact}}(\bC_2)_\CZ.$$
\end{lem}

\begin{rem}
The material in \secref{ss:restr fact mod} is applicable in the context of unital factorization categories
and strictly unital factorization functors between them.
\end{rem} 

\sssec{} \label{sss:unital restr categ pairs} 

Let us place ourselves momentarily in the context of \secref{sss:restr pairs}, where the factorization categories
in the diagram 
$$
\CD
\bA'_1 @>{\Phi'}>> \bA'_2 \\
@A{\Psi_1}AA @AA{\Psi_2}A \\
\bA_1 @>{\Phi}>> \bA_2.
\endCD
$$
are unital, all functors are lax unital, and the vertical arrows are strictly unital. 

\medskip

Unwinding the constructions, one obtains that in this case the resulting functor 
$$\on{Res}_{\Psi_1}(\bA'_1{}^{\on{fact}_\CZ})\to \on{Res}_{\Psi_2}(\bA'_2{}^{\on{fact}_\CZ})$$
viewed as a functor between unital module categories over $\bA_1$ and $\bA_2$, respectively,
possesses a natural lax unital structure compatible with factorization. 

\ssec{Restriction along lax unital functors}

In this subsection we will study the phenomenon of restriction with respect to a functor
$$\Phi:\bA_1\to \bA_2,$$
which is only lax unital. 

\sssec{} \label{sss:enh mod untl}

Consider the (unital) factorization algebra
$$\Phi(\one_{\bA_1})\in \on{FactAlg}^{\on{untl}}(X,\bA_2),$$
see \secref{sss:funct between pairs modules untl}. 

\medskip

Let $\CZ$ be a prestack mapping to $\Ran$, and let
$$\Phi_m:\bC_1\to \bC_2$$
be an object of $\on{Funct}_{\bA_1\to \bA_2}(\bC_1,\bC_2)$.

\medskip

Consider the induced functor
\begin{equation} \label{e:enhancement modules}
\Phi_m:\bC_{1,\CZ}\to \bC_{2,\CZ}
\end{equation}
between the underlying DG categories. 

\sssec{}

We claim:

\begin{lemconstr} \label{l:enhancement modules}
The functor \eqref{e:enhancement modules} naturally enhances to a functor
$$\Phi_m^{\on{enh}}:\bC_{1,\CZ}\to \Phi(\one_{\bA_1})\mod^{\on{fact}}(\bC_2)_\CZ.$$
\end{lemconstr}

\begin{proof} 

We rewrite
$$\bC_{1,\CZ}\simeq \one_{\bA_1}\mod^{\on{fact}}(\bC_1)_\CZ,$$
and now the required functor is a particular case of \eqref{e:funct between pairs modules untl alg}.

\end{proof}

\sssec{Variant} \label{sss:pairs ult}

Let us return for a moment to the setting of \secref{sss:unital restr categ pairs}. By \lemref{l:enhancement modules} for every
$\ul{x}:S\to \Ran$ we obtain that the functor
$$\Phi':\Res_{\Psi_1}(\bA'_1{}^{\on{fact}_{\ul{x}}})\to \Res_{\Psi_2}(\bA'_2{}^{\on{fact}_{\ul{x}}})$$
gives rise to a functor
\begin{multline} \label{e:pairs ult 1}
\Phi'{}^{\on{enh}}:\bA'_{1,\ul{x}}\to 
\Phi(\one_{\bA_1})\mod^{\on{fact}}(\Res_{\Psi_2}(\bA'_2{}^{\on{fact}_{\ul{x}}}))_{\ul{x}}
\overset{\text{\lemref{l:modules for fact alg restr untl}}}\simeq \\
\simeq (\Psi_2\circ \Phi)(\one_{\bA_1})\mod^{\on{fact}}(\bA'_2{}^{\on{fact}_{\ul{x}}})_{\ul{x}}.
\end{multline} 

Note that the left-hand side in \eqref{e:pairs ult 1} is the value at $\ul{x}$ of a unital 
factorization category, namely, $\bA'_1$ itself. 

\medskip

The right-hand side in \eqref{e:pairs ult 1} is the value at $\ul{x}$ of a unital 
\emph{lax} factorization category, namely,
$$(\Psi_2\circ \Phi)(\one_{\bA_1})\mod^{\on{fact}}(\bA'_2),$$
see \secref{sss:fact A mod in A untl}. 

\medskip

Unwinding the constructions, we obtain that \eqref{e:pairs ult 1} upgrades to a \emph{unital}
factorization functor
\begin{equation} \label{e:pairs ult 2}
\Phi'{}^{\on{enh}}:\bA'_1\to (\Psi_2\circ \Phi)(\one_{\bA_1})\mod^{\on{fact}}(\bA'_2).
\end{equation}  

\begin{rem}

Note that in the setting of \secref{sss:pairs ult}, we have
$$(\Psi_2\circ \Phi)(\one_{\bA_1})\simeq  (\Phi'\circ \Psi_1)(\one_{\bA_1})\simeq \Phi'(\one_{\bA'_1}).$$

So, the information contained by the functor \eqref{e:pairs ult 2} is completely captured by the case when
$\Psi_1$ and $\Psi_2$ are the identity functors. I.e., the claim is that the lax unital factorization functor
$$\Phi:\bA_1\to \bA_2$$
upgrades to a unital factorization functor
$$\Phi^{\on{enh}}:\bA_1\to \Phi(\one_{\bA_1})\mod^{\on{fact}}(\bA_2),$$
with the caveat that the right-hand side is a lax factorization category. 

\end{rem} 

\sssec{}

Let $\CZ$ and $\bC_2$ be as in \secref{sss:enh mod untl}. Consider the contravariant functor on $\bA_1\mmod^{\on{fact}}_\CZ$
that assigns to $\bC_1$ the category $\on{Funct}_{\bA_1\to \bA_2}(\bC_1,\bC_2)$. One shows
that this functor is representable, and let
$$\Res^{\on{untl}}_\Phi(\bC_2)\in \bA_1\mmod^{\on{fact}}_\CZ$$
denote the representing object. 

\medskip

By \lemref{l:enhancement modules}, the tautological functor
\begin{equation} \label{e:fact restr cat untl prel}
\Res^{\on{untl}}_\Phi(\bC_2)\to \bC_2
\end{equation} 
upgrades to a functor
\begin{equation} \label{e:fact restr cat untl}
\Res^{\on{untl}}_\Phi(\bC_2)_\CZ\to \Phi(\one_{\bA_1})\mod^{\on{fact}}(\bC_2)_\CZ.
\end{equation} 

\sssec{}

We have the following generalization of \lemref{l:fact restr cat}: 

\begin{lem} \label{l:fact restr cat untl}
The functor \eqref{e:fact restr cat untl} is an equivalence.
\end{lem} 

\begin{rem}
Note that when $\Phi$ is strictly unital, then by \lemref{l:unital restr categ univ}
$$\Res^{\on{untl}}_\Phi(\bC_2)\simeq \Res_\Phi(\bC_2),$$
and the assertion of \lemref{l:fact restr cat untl} coincides with that of
\lemref{l:fact restr cat}.
\end{rem} 

\sssec{Example} \label{sss:unital retsr expl}

Let $\CZ=\on{pt}$ and let $\CZ\to \Ran$ correspond to a singleton $\{x\}\in \Ran$. Recall the notations of 
\secref{sss:fact restr limit}. Let us give an explicit description of the category 
$$\Res^{\on{untl}}_\Phi(\bC_2)_X.$$

Namely,
\begin{multline*} 
\Res^{\on{untl}}_\Phi(\bC_2)_X \simeq \\
\simeq \Phi(\one_{\bA_1})\mod^{\on{fact}}((\bC_2)|_X)_X \underset{\Phi(\one_{\bA_1})\mod^{\on{fact}}((\bC_2)|_{X-x})_{X-x}}\times
\left(\bA_{1,X-x}\otimes \Phi(\one_{\bA_1})\mod^{\on{fact}}(\bC_2)_x\right),
\end{multline*} 
where:

\begin{itemize}

\item The notation $(\bC_2)|_X$ is as in \secref{sss:propagate modules cat}, i.e., we regard the pullback of $\ul\bC_2$
along $X\to \Ran_x$ as a module category over $\bA_2$ at $X$;

\item The functor $$\bA_{1,X-x}\otimes \Phi(\one_{\bA_1})\mod^{\on{fact}}(\bC_2)_x\to 
\Phi(\one_{\bA_1})\mod^{\on{fact}}((\bC_2)|_{X-x})_{X-x}$$
is the composition
\begin{multline*} 
\bA_{1,X-x}\otimes \Phi(\one_{\bA_1})\mod^{\on{fact}}(\bC_2)_x
\overset{\Phi^{\on{enh}}\otimes \on{Id}}\longrightarrow 
\Phi(\one_{\bA_1})\mod^{\on{fact}}(\bA_2)_{X-x} \otimes \Phi(\one_{\bA_1})\mod^{\on{fact}}(\bC_2)_x\to \\
\overset{\on{factorization\,of\,}\bC_2}\longrightarrow \Phi(\one_{\bA_1})\mod^{\on{fact}}((\bC_2)|_{X-x})_{X-x},
\end{multline*} 
where the second arrow is defined as in \secref{sss:fact alg lax fact}. 

\end{itemize}

\sssec{}

We have the following analog of \lemref{l:modules for fact alg restr untl}

\begin{lem} \label{l:modules for fact alg restr lax untl} 
The functor \eqref{e:fact restr cat untl prel} induces an equivalence
$$\CA_1\mod^{\on{fact}}(\Res^{\on{untl}}_\Phi(\bC_2))_\CZ\to \Phi(\CA_1)\mod^{\on{fact}}(\bC_2)_\CZ.$$
\end{lem}

\sssec{}

The material in \secref{ss:restr fact mod} is applicable in the context of unital factorization categories
and lax unital factorization functors between them. We will not need it in the full generality, except for
an analog of \corref{c:basic adj non-unital}, formulated as \lemref{l:basic adj} below.

\sssec{} \label{sss:basic adj}

Let
$$\Phi:\bA_1\to \bA_2$$
be a unital functor between unital factorization categories. 

\medskip

Suppose that $\Phi$ admits a right adjoint as a functor between the underlying crystals
of categories over $\Ran$. According to \secref{sss:adj of unital}, the right adjoint $\Phi^R$
of $\Phi$ admits a natural extension to a lax unital functor between unital factorization categories.

\medskip

Let $\bC_1$ (resp., $\bC_2$) be a unital module category over $\bA_1$ (resp., $\bA_2$) at some 
$\CZ\to \Ran$. 

\sssec{}

We claim:

\begin{lem} \label{l:basic adj}
There is a canonical equivalence
$$\on{Funct}_{\bA_1\mmod^{\on{fact}}_\CZ}(\Res_\Phi(\bC_2),\bC_1)\simeq
\on{Funct}_{\bA_2\mmod^{\on{fact}}_\CZ}(\bC_2,\Res_{\Phi^R}^{\on{untl}}(\bC_1)).$$
\end{lem} 

\sssec{}

We will not give a full proof of this lemma; rather we will sketch the construction of the maps
in both directions in the framework of \secref{sss:unital retsr expl}. 

\medskip

Namely, for $\bC_2\in \bA_2\mmod^{\on{fact}}_\CZ$ we will construct a functor
\begin{equation} \label{e:basic adj 1}
\bC_2\to \Res_{\Phi^R}^{\on{untl}}\circ \Res_\Phi(\bC_2) 
\end{equation} 
and for $\bC_1\in \bA_1\mmod^{\on{fact}}_\CZ$ we will construct a functor
\begin{equation} \label{e:basic adj 2}
\Res_\Phi \circ \Res_{\Phi^R}^{\on{untl}}(\bC_1)\to \bC_1. 
\end{equation} 

\sssec{}

By the universal property of $\Res_{\Phi^R}^{\on{untl}}$, the datum of \eqref{e:basic adj 1} is equivalent to the datum of a functor
\begin{equation} \label{e:basic adj 1 prime}
\bC_2\to \Res_\Phi(\bC_2) 
\end{equation} 
as module categories over $\bA_2$ and $\bA_1$, respectively, compatible with factorization against the functor $\Phi^R$. 

\medskip

We now specialize to the context of \secref{sss:unital retsr expl} and construct the corresponding functor
\begin{equation} \label{e:basic adj 1 prime X}
\bC_{2,X}\to  \Res_\Phi(\bC_2)_X.
\end{equation} 

\medskip

We write
\begin{equation} \label{e:basic adj 1 prime a}
\Res_\Phi(\bC_2)_X\simeq \bC_{2,X}\underset{\bC_{2,X-x}}\times (\bA_{1,X-x}\otimes \bC_{2,x}),
\end{equation} 
where 
$$\bA_{1,X-x}\otimes \bC_{2,x}\to \bC_{2,X-x}$$ is the functor
$$\bA_{1,X-x}\otimes \bC_{2,x}\overset{\Phi\otimes \on{Id}}\to
\bA_{2,X-x}\otimes \bC_{2,x}\overset{\on{fact}}\simeq \bC_{2,X-x}.$$
 
We write 
\begin{equation} \label{e:basic adj 1 prime b}
\bC_{2,X}\simeq \bC_{2,X}\underset{\bC_{2,X-x}}\times (\bA_{2,X-x}\otimes \bC_{2,x}),
\end{equation} 
where 
$$\bA_{2,X-x}\otimes \bC_{2,x}\overset{\sim}\to \bC_{2,X-x}$$
is the factorization equivalence. 

\medskip

We define the functor \eqref{e:basic adj 1 prime X} by sending an object $\bc_2\in \bC_{2,X}$ to 
$$\left(\bc_2\underset{j_*\circ j^*(\bc_2)}\times j_*(((\Phi\circ\Phi^R)\otimes \on{Id})(j^*(\bc_2))), (\Phi^R\otimes \on{Id})(j^*(\bc_2))\right)\in
\bC_{2,X}\underset{\bC_{2,X-x}}\times (\bA_{1,X-x}\otimes \bC_{2,x}),$$
where the map
$$((\Phi\circ\Phi^R)\otimes \on{Id})(j^*(\bc_2))\to j^*(\bc_2)$$
is the counit of the adjunction. 

\sssec{}

Write
\begin{multline} \label{e:basic adj 2 prime a}
\Res_{\Phi}\circ \Res^{\on{untl}}_{\Phi^R}(\bC_1)_X\simeq \\
\simeq (\Phi^R\circ \Phi)(\one_{\bA_1})\mod^{\on{fact}}((\bC_1)|_X)_X 
\underset{(\Phi^R\circ \Phi)(\one_{\bA_1})\mod^{\on{fact}}((\bC_1)|_{X-x})_{X-x}}\times \\
\underset{(\Phi^R\circ \Phi)(\one_{\bA_1})\mod^{\on{fact}}((\bC_1)|_{X-x})_{X-x}}\times \left(\bA_{1,X-x}\otimes ((\Phi^R\circ \Phi)(\one_{\bA_1})\mod^{\on{fact}}(\bC_1)_x)\right).
\end{multline} 

\medskip


The sought-for functor 
$$\Res_{\Phi}\circ \Res^{\on{untl}}_{\Phi^R}(\bC_1)_X\to \bC_{1,X}$$
sends an object $(\bc'_1,\bc''_1)$ in the right-hand side of \eqref{e:basic adj 2 prime a}, i.e.,
$$\bc'_1\in (\Phi^R\circ \Phi)(\one_{\bA_1})\mod^{\on{fact}}((\bC_1)|_X)_X$$
and
$$\bc''_1\in \bA_{1,X-x}\otimes ((\Phi^R\circ \Phi)(\one_{\bA_1})\mod^{\on{fact}}(\bC_1)_x)$$
to the object
$$\oblv_{(\Phi^R\circ \Phi)(\one_{\bA_1})}(\bc'_1)\underset{j_*\circ j^*\circ \oblv_{(\Phi^R\circ \Phi)(\one_{\bA_1})}(\bc'_1)}\times 
\left(j_*\circ (\on{Id}\otimes \oblv_{(\Phi^R\circ \Phi)(\one_{\bA_1})})(\bc''_1)\right),$$
where the map
$$(\on{Id}\otimes \oblv_{(\Phi^R\circ \Phi)(\one_{\bA_1})})(\bc''_1)\to j^*\circ \oblv_{(\Phi^R\circ \Phi)(\one_{\bA_1})}(\bc'_1)$$
is obtained by applying the functor $\oblv_{(\Phi^R\circ \Phi)(\one_{\bA_1})}$ to the isomorphism
$$((\Phi^R\circ \Phi)^{\on{enh}}\otimes \on{Id})(\bc''_1)\simeq j^*(\bc'_1),$$
precomposed with the unit of the $(\Phi,\Phi^R)$-adjunction
\begin{multline*}
(\on{Id}\otimes \oblv_{(\Phi^R\circ \Phi)(\one_{\bA_1})})(\bc''_1)
\to ((\Phi^R\circ \Phi)\otimes \oblv_{(\Phi^R\circ \Phi)(\one_{\bA_1})})(\bc''_1)\simeq \\
\simeq \oblv_{(\Phi^R\circ \Phi)(\one_{\bA_1})}\left(((\Phi^R\circ \Phi)^{\on{enh}}\otimes \on{Id})(\bc''_1)\right).
\end{multline*}

\ssec{Tightness}

\sssec{}

Let $\bA$ be a unital factorization category. We will say that $\bA$ is \emph{tight} if for every $\ul{x},\ul{x}':S\to \Ran$
with $\ul{x}\subseteq \ul{x}'$, the corresponding functor
$$\on{ins}_{\ul{x}\subseteq \ul{x}'}:\bA_{S,\ul{x}}\to \bA_{S,\ul{x}'}$$
admits a continuous right adjoint.

\sssec{}

Many of the  unital factorization categories we have introduced satisfy this property. This includes representation-theoretic
examples, e.g., 
$$\hg\mod_\kappa,\,\, \KL(G)_\kappa$$
and algebro-geometric examples:
$$\QCoh(\CT),$$
where $\CT$ is an affine counital factorization space, and 
$$\IndCoh^*(\CT) \text{ and } \IndCoh^!(\CT),$$
 where $\CT$ is a unital-in-correspondences ind-placid factorization ind-schemes, satisfying the conditions 
from Sects. \ref{sss:corr *} and \ref{sss:corr !}, respectively.

\sssec{}

Let $\bA$ be a tight unital factorization category, and let $\bC$ be a unital factorization module category 
over $\bC$ at some $\CZ$. 

\medskip

We shall say that $\bC$ is \emph{tight} if for every $(z,\ul{x}),(z,\ul{x}'):S\to \CZ^\subseteq$
with $\ul{x}\subseteq \ul{x}'$, the corresponding functor
$$\on{ins}_{\ul{x}\subseteq \ul{x}'}:\bC_{S,\ul{x}}\to \bC_{S,\ul{x}'}$$
admits a continuous right adjoint.

\sssec{}

The following is immediate:

\begin{lem}
Suppose that $\bA$ is tight. Then for any $\CZ\to \Ran$, the factorization module category 
$\bA^{\on{fact}_\CZ}$ is tight. 
\end{lem}

\sssec{} \label{sss:tight vect}

Take $\bA=\Vect$, and recall the construction from \secref{sss:untl fact mod Vect}
\begin{equation} \label{e:tight vect}
\ul\bC_0\in \on{CrysCat}(\CZ_0) \, \rightsquigarrow \, \bC\in \Vect\mmod^{\on{fact}}_\CZ.
\end{equation}

It is clear that the essential image of \eqref{e:tight vect} is contained in the full subcategory of 
$\Vect\mmod^{\on{fact}}_\CZ$ consisting of tight
unital factorization module categories. 

\medskip

We claim:

\begin{lem} \label{l:tight vect}
The functor \eqref{e:tight vect} is an equivalence onto the full subcategory of $\Vect\mmod^{\on{fact}}_\CZ$ 
consisting of tight objects. 
\end{lem}

\begin{proof}

The functor \eqref{e:tight vect} admits a retraction (i.e., a left inverse), given by restricting the crystal
of categories from $\CZ^\subseteq$ to $\CZ$ along $\on{diag}_\CZ$. We claim that this left inverse
is an actual inverse when applied to tight objects.

\medskip

Indeed, let $\bC'$ be a tight unital factorization module category at $\CZ$. Let $\ul\bC'$ be the 
corresponding crystal of categories over $\CZ^{\subseteq}$, and let $\ul\bC'_0$ denote the restriction of
$\ul\bC'$ along $\on{diag}_\CZ$.

\medskip

The unital structure on $\bC'$ gives rise to a functor
\begin{equation} \label{e:tight vect 1}
\on{pr}_{\on{small}}^*(\ul\bC'_0)\to \ul\bC.
\end{equation}

We have to show that \eqref{e:tight vect 1} is an equivalence. 

\medskip

By the tightness assumption, the functor \eqref{e:tight vect 1} admits a right adjoint.
Hence, in order to show that it is an equivalence, it suffices to show that it is an equivalence
strata-wise. However, this follows from factorization.

\end{proof} 

\sssec{}

Consider the following situation. Let $\bA$ be a tight unital factorization category, and let
$\bC$ be a tight unital factorization module category over it at some $\CZ$. 

\medskip

Consider the (strictly) unital factorization functor
$$\on{Vac}_\bA:\Vect\to \bA,$$ 
see \secref{sss:Vect to A}. 

\sssec{}

We claim: 

\begin{lem} \label{l:restr to Vect is tight}
The unital factorization module category
$$\Res_{\on{Vac}_\bA}(\bC)\in \Vect\mmod^{\on{fact}}_\CZ$$
is tight.
\end{lem}

\begin{proof} 

Set $\ul\bC'_0:=\ul\bC_\CZ$ be the sheaf of categories on $\CZ$ underlying $\bC$. Let $\bC'$ 
be the (tight, unital) factorization module over $\Vect$, attached $\ul\bC'_0$ by the functor \eqref{e:tight vect}.

\medskip

The unital structure on $\bC$ gives to a (strictly unital) functor
$$\bC'\to \bC,$$
compatible with factorization (in the sense of \secref{sss:funct between pairs modules untl}). Hence,
we obtain a functor
\begin{equation} \label{e:tight vect 2}
\bC'\to \Res_{\on{Vac}_\bA} (\bC).
\end{equation}

We claim that the functor \eqref{e:tight vect 2} is an equivalence. Indeed, this follows from
the assumptions by applying \lemref{l:fact res crit}.

\end{proof}

\sssec{}

We will use \lemref{l:restr to Vect is tight} as follows. Let $\bA$ be a tight unital factorization category,
and let $R$ be a factorization algebra (in $\Vect$). 

\medskip

Using the functor $\on{Vac}_\bA$, we can consider the factorization algebra 
$$\on{Vac}_\bA(R)\simeq R\otimes \one_\bA$$ 
in $\bA$, and consider the corresponding lax factorization category
$$R\mod^{\on{fact}}(\bA):=(R\otimes \one_\bA)\mod^{\on{fact}}(\bA).$$

\medskip

We have a naturally defined functor between lax factorization categories
\begin{equation} \label{e:R mod in A}
R\mod^{\on{fact}}\otimes \bA\to R\mod^{\on{fact}}(\bA),
\end{equation} 
see \eqref{e:tensor up fact mod}.

\medskip 

Combining Lemmas \ref{l:restr to Vect is tight} and \ref{l:tensor up fact mod}, we obtain:

\begin{cor} \label{c:R mod in A}
Assume that $\bA$ is dualizable as a factorization category. Then the functor \eqref{e:R mod in A}
is an equivalence.
\end{cor}

\section{Chiral modules} \label{s:chiral mods}

The main purpose of this section is to prove \thmref{t:from QCoh star to fact}, which gives a geometric description of modules 
over commutative factorization/chiral algebras. 

\medskip

To do so, we first develop a general theory describing modules over chiral algebras in terms of modules over topological algebras 
(although we do not write in these exact terms), with an especially explicit understanding for ``nice" Lie-* algebras. 

\medskip

This material largely consists of transporting \cite[Sect. 3.6]{BD2} into the derived setting. However, we will encounter a surprise:
a certain equivalence that always takes place at an abelian level, in order to hold at the derived level requires a finiteness
condition (see \secref{ss:QCoh* non-fp}). 

\ssec{A reminder: chiral vs factorization algebras} 

\sssec{} \label{sss:chiral algebras}

Recall that, according to \cite[Proposition 3.4.19]{BD2} (see \cite{FraG} for the derived version), 
factorization algebras are the same as chiral algebras. Given a factorization algebra $\CA$, 
the corresponding chiral algebra, thought of as a D-module on $X$, is given by
$$\CA^{\on{ch}}:=\CA_X[-1].$$

For example, the unit factorization algebra corresponds to the chiral algebra $\omega_X$
(see our conventions in \secref{sss:omega X}\footnote{This was one of the main reasons for this choice
of conventions, i.e., in order to be in line with \cite{BD2}.}).

\sssec{}

Generalizing \cite[Proposition 3.4.19]{BD2},
we have
\begin{equation} \label{e:ch vs fact}
\CA\mod^{\on{fact}}\simeq \CA^{\on{ch}}\mod^{\on{ch}},
\end{equation} 
where modules on both sides can be taken on any space that maps to $\Ran$. 

\sssec{} \label{sss:com ch}

Let us recall also how the bijection between between chiral and factorization algebras plays out in the 
commutative case.

\medskip

Let $A$ be a commutative algebra object in $\Dmod(X)$. Then the corresponding factorization algebra 
$\CA:=\on{Fact}(A)$ is such that
$$\CA_X=A.$$

And the corresponding chiral algebra $\CA^{\on{ch}}$ is $A[-1]$. Note that
$$\oblv^r(\CA^{\on{ch}})=\oblv^l(A)\otimes \omega_X.$$

\sssec{} \label{sss:com ch free}

For example, for the free commutative algebra $A=\Sym^!(\CM[1])$ for $\CM\in \Dmod(X)$, we have
$$\CA^{\on{ch}}=\Sym^!(\CM[1])[-1]\simeq U^{\on{ch}}(\CM),$$
where:

\begin{itemize}

\item In the left-hand side $\CM$ is considered as an abelian Lie-* algebra;

\item $U^{\on{ch}}$ is the functor of chiral envelope.

\end{itemize}

Note also that in this case.
$$\oblv^l(A)=\Sym_{\CO_X}(\oblv^l(\CM[1]))\simeq \Sym_{\CO_X}(\oblv^r(\CM)\otimes \omega_X^{\otimes -1}).$$

So, if $\CM=\ind^r(\CE)=\CE\otimes \on{D}_X$ for a classical locally free sheaf $\CE$ on $X$, then the corresponding D-scheme
$$\Spec_X(A)$$ is the scheme $\on{Jets}(\CE^\vee\otimes \omega_X)$
of jets into the vector bundle $\CE^\vee\otimes \omega_X$.

\sssec{Convention}

In what follows, by a slight abuse of notation, for a given factorization algebra $\CA$, we will use
the same symbol $\CA$ to denote the corresponding chiral algebra (i.e., we will not write $\CA^{\on{ch}}$). 

\ssec{The pro-projective generator for chiral modules} 

\sssec{}

Let $\CA$ be a unital chiral algebra on $X$. Let $\CA\mod_x^{\on{ch}}$ denote the category of \emph{unital} chiral $\CA$-modules at $x$. 
We let $\oblv_\CA$ denote the forgetful functor $\CA\mod_x^{\on{ch}}\to \Vect$. 

\medskip

For $\CM\in \CA\mod_x^{\on{ch}}$, we will consider the action map
$$j^*(\CA)\otimes \CM\overset{\text{action}}\longrightarrow i_*(\CM),$$
as a map of chiral $\CA$-modules on $X$, where:

\begin{itemize}


\item $j$ denotes the open embedding $X-x\hookrightarrow X$;

\item $i$ denotes the embedding of the point $x$ into $X$.

\end{itemize}

\sssec{}

We will also use a short-hand notation
$$\ul\CM:=\oblv_\CA(\CM), \quad \CM\in \CA\mod_x^{\on{ch}}.$$

In what follows we will take about ``elements" of $\ul\CM$: 

\medskip

For $V\in \Vect$, by an element 
$v\in V$ we mean a point of the space $\Maps_{\Vect}(k,V)$.

\sssec{}

In what follows we will denote by $\CA_x$ is the $[1]$-shifted !-fiber of $\CA$ at $x$ (this is the same as the !-fiber 
at $x\in \Ran$ of the factorization algebra
corresponding to $\CA$), viewed as an object\footnote{In fact, $\CA_x$ should more properly be denoted $\CA^{\on{fact}_x}$,
see \secref{sss:vacuum module}.}
of $\CA\mod^{\on{ch}}_x$.

\medskip

We let 
$$1_{\CA,x}\in \ul\CA_x$$
denote the vacuum vector, i.e., the image of $1\in k\simeq (\omega_X)_x$ under the unit map 
$$\on{vac}_\CA:\omega_X\to \CA.$$

\sssec{}

Consider the category $\on{Modif}(\CA)$ of unital chiral algebras $\CA'$ equipped
with an isomorphism 
$$\CA'|_{X-x}\simeq \CA|_{X-x}.$$

\medskip

This category has fiber products, and hence is 
cofiltered. 

\medskip

Note that the category $\CA\mod_x^{\on{ch}}$ only depends on $\CA|_{X-x}$, so for any $\CA'\in \on{Modif}(\CA)$ 
we have a canonical identification
$$\CA\mod_x^{\on{ch}}\simeq \CA'\mod_x^{\on{ch}}.$$

\sssec{}

Consider the functor
$$\on{Modif}(\CA)\to \CA\mod^{\on{ch}}_x, \quad \CA'\mapsto \CA'_x\in \CA'\mod^{\on{ch}}\simeq \CA\mod^{\on{ch}}.$$

\medskip

Set 
$$\sP_{\CA,x}:=\underset{\CA'\in \on{Modif}(\CA)}{``\on{lim}"}\, \CA'_x\in \on{Pro}(\CA\mod^{\on{ch}}_x).$$

The object
$$\oblv_\CA(\sP_{\CA,x})\in \on{Pro}(\Vect)$$
is equipped with a canonical vector $1_{\sP,x}$ comprised of the vacuum vectors
$1_{\CA',x}\in \ul\CA'_x$. 

\sssec{}

Evaluation on $1_{\sP,x}$ gives rise to a natural transformation
\begin{equation} \label{e:proj generator}
\CHom(\sP_{\CA,x},-)\to \oblv_\CA, \quad \CA\mod^{\on{ch}}_x\to \Vect.
\end{equation} 

The following assertion is a derived version of \cite[Proposition 3.6.16]{BD2}:

\begin{prop} \label{p:proj generator}
The natural transformation is an isomorphism.
\end{prop}

\sssec{} \label{sss:modif as pairs}

The proof of \propref{p:proj generator} is based on the following observation: we claim that the category $\on{Modif}(\CA)$ 
is equivalent to the category of pairs
$$(\CM\in \CA\mod_x^{\on{ch}},m\in \ul\CM).$$

Namely, in one direction, to $\CA'\in \on{Modif}(\CA)$ we attach the pair $(\CA'_x,1_{\CA',x})$.

\medskip

Vice versa, given $(\CM,m)$
we let $\CA'$ be the fiber of the map
\begin{equation} \label{e:construct modif}
j_*\circ j^*(\CA)\overset{\on{id}\otimes m}\longrightarrow j_*\circ j^*(\CA)\otimes \CM\overset{\text{action}}\longrightarrow i_*(\CM).
\end{equation} 

We claim that $\CA'$ has a natural structure of unital chiral algebra. This fits into the following general paradigm:

\medskip

Let $L$ be a Lie algebra in a symmetric pseudo-monoidal monoidal category $\bA$, and let $M$ be a module over it. Consider $L\oplus M$
as a split square-zero extension of $M$. For an element $m\in M$ (i.e., a map $\one_\bA\to M$), the action of $L$ on $M$
gives rise to an automorphism $\phi_m$ of $L\oplus M$. Then we can identify
$$\on{Fib}(L\overset{\text{act on }m}\longrightarrow M)$$
with 
$$L\underset{L\oplus M}\times L,$$
where the two maps $L\rightrightarrows L\oplus M$ are the compositions of the tautological embedding
with: (1) the identity map, and (2) $\phi_m$.

\medskip

We apply this to $\bA:=\Dmod(X)$, equipped with the chiral pseudo-monoidal monoidal structure, so 
that $\one_\bA$ is the ``constant sheaf" on $X$. We take $L=j_*\circ j^*(\CA)$ and $M=i_*(\CM)$. 
This endows $\CA'$ with a chiral
algebra structure, i.e., 
$$\CA'=j_*\circ j^*(\CA)\underset{j_*\circ j^*(\CA)\oplus i_*(\CM)}\times j_*\circ j^*(\CA).$$

In order to show that $\CA'$ is unital, by \propref{p:quasi-untl}, it is enough to equip it with a quasi-unital structure
(see \secref{sss:q-untl} for what this means).  The above fiber product presentation defines this structure on the nose.

\begin{proof}[Proof of \propref{p:proj generator}]

This is tautological from \secref{sss:modif as pairs}: the assertion of the proposition is just the fact that the map
$$\underset{(\CM,m\in \ul\CM)}{\on{colim}}\, \CHom(\CM,\CM')\to \ul\CM'$$
is an isomorphism. 

\end{proof} 

\sssec{}

Assume now that $j^*(\CA)$ is connective. Let 
\begin{equation} \label{e:conn modif}
\on{Modif}_{\on{conn}}(\CA)\subset \on{Modif}(\CA)
\end{equation} 
be the full subcategory consisting of those $\CA'$ that are connective.

\medskip

Truncation $\leq 0$ on chiral algebras defines a right adjoint to the above inclusion.
Hence, the opposite of \eqref{e:conn modif} is cofinal. 

\medskip

In particular, the object $\sP_{\CA,x}$ maps isomorphically to 
$$\underset{\CA'\in \on{Modif}_{\on{conn}}(\CA)}{``\on{lim}"}\, \CA'_x\in \on{Pro}(\CA\mod^{\on{ch}}_x).$$

\sssec{}

For an integer $m$, let 
$$\on{Modif}_{\on{conn,\geq -m}}(\CA)\subset \on{Modif}_{\on{conn}}(\CA)$$
be the full subcategory consisting of those objects $\CA'$, 
for which $\CA'_x\in \Vect^{\geq -m,\leq 0}$.

\medskip 

Let $\on{Modif}_{\on{conn,ev-c}}(\CA)$ be the full subcategory of $\on{Modif}_{\on{conn}}(\CA)$ consisting of those objects $\CA'$, 
for which $\CA'_x$ is eventually coconnective (as an object of $\Vect$). I.e.,
$$\on{Modif}_{\on{conn,ev-c}}(\CA)=\underset{m}{\on{colim}}\, \on{Modif}_{\on{conn,\geq -m}}(\CA).$$

\medskip

Note that if $j^*(\CA)$ is itself eventually coconnective, the above condition on $\CA'_x$ is equivalent to $\CA'$ being eventually coconnective.

\sssec{}

Set
$$\sP_{\CA,x,\on{ev-c}}:=\underset{\CA'\in \on{Modif}_{\on{conn,ev-c}}(\CA)}{``\on{lim}"}\, \CA'_x\in \on{Pro}(\CA\mod^{\on{ch}}_x).$$

We have a tautological map
\begin{equation} \label{e:P ev c}
\sP_{\CA,x}\to \sP_{\CA,x,\on{ev-c}}
\end{equation}
in $\on{Pro}(\CA\mod^{\on{ch}}_x)$. 

\begin{lem} \label{l:P ev c}
The map 
$$\CHom(\sP_{\CA,x,\on{ev-c}},-)\to \CHom(\sP_{\CA,x},-),$$
defined by \eqref{e:P ev c} is an isomorphism, when evaluated on $(\CA\mod_x^{\on{ch}})^{>-\infty}$.
\end{lem}

\begin{proof}

It suffices to show that if $\CM\in (\CA\mod_x^{\on{ch}})^{\geq -m}$ for some $m$, then 
$$\underset{\CA'\in \on{Modif}_{\on{conn,ev-c}}(\CA)}{\on{colim}}\, \Maps(\CA'_x,\CM)\to 
\underset{\CA'\in \on{Modif}_{\on{conn}}(\CA)}{\on{colim}}\, \Maps(\CA'_x,\CM)$$
is an isomorphism.

\medskip

However, it is clear that the map
$$\underset{\CA'\in \on{Modif}_{\on{conn,\geq -m'}}(\CA)}{\on{colim}}\, \Maps(\CA'_x,\CM)\to 
\underset{\CA'\in \on{Modif}_{\on{conn}}(\CA)}{\on{colim}}\, \Maps(\CA'_x,\CM)$$
is an isomorphism for every $m'\geq m$. 

\end{proof}

\begin{cor}  \label{c:P ev c}
The map 
$$\CHom(\sP_{\CA,x,\on{ev-c}},-)\to \oblv_\CA$$
is an isomorphism, when evaluated on $(\CA\mod_x^{\on{ch}})^{>-\infty}$.
\end{cor} 

\ssec{The case of Lie-* algebras} \label{ss:ass of top Lie}

In this subsection we will assume that the chiral algebra $\CA$ is the chiral universal envelope $U^{\on{ch}}(L)$  of a (connective)
Lie-* algebra $L$. Recall that we can identify
$$\CA\mod^{\on{ch}}_x\simeq L\mod_x^{\on{ch}}.$$

We will study how the object $\sP_{\CA,x}$ looks like in this case. 

\sssec{}

Consider the categories 
$$\on{Modif}(L) \text{ and } \on{Modif}_{\on{conn}}(L)$$
defined as in the case of chiral algebras, i.e., these are Lie-* algebras equipped with an isomorphism with $L$ over $X-x$. 

\medskip

We have two pairs of adjoint functors
$$U^{\on{ch}}:\on{Modif}(L)\rightleftarrows \on{Modif}(\CA):\oblv^{\on{ch}\to \on{Lie-}^*}$$ and  
$$U^{\on{ch}}:\on{Modif}_{\on{conn}}(L)\rightleftarrows \on{Modif}_{\on{conn}}(\CA):\oblv^{\on{ch}\to \on{Lie-}^*}.$$

In particular, the corresponding functors
$$\on{Modif}(L)^{\on{op}}\to \on{Modif}(\CA)^{\on{op}} \text{ and } \on{Modif}_{\on{conn}}(L)^{\on{op}}\to \on{Modif}(\CA)_{\on{conn}}^{\on{op}}$$
are cofinal.

\sssec{}

In particular, we obtain that we can write 
$$\sP_{\CA,x}\simeq \underset{L'\in \on{Modif}(L)}{``\on{lim}"}\, \ind^{L\mod^{\on{ch}}_x}_{L'\mod^{\on{Lie-}^*}_x}(k)$$
and when $L$ is connective also as 
$$\sP_{\CA,x}\simeq \underset{L'\in \on{Modif}_{\on{conn}}(L)}{``\on{lim}"}\, \ind^{L\mod^{\on{ch}}_x}_{L'\mod^{\on{Lie-}^*}_x}(k),$$
where
\begin{itemize}

\item $\ind^{L\mod^{\on{ch}}_x}_{L'\mod^{\on{Lie-}^*}_x}:L'\mod^{\on{Lie-}^*}_x\to L'\mod^{\on{ch}}_x\simeq L\mod^{\on{ch}}_x$
is the left adjoint of the restriction functor.

\end{itemize}

\sssec{} \label{sss:ass of top Lie gen good}

Assume now that $j^*(L)$ is classical (i.e., is in cohomological degree $0$ as a right D-module), and let
$$\on{Modif}_{\on{cl,flat}}(L)\subset \on{Modif}_{\on{conn}}(L)$$
be the full subcategory, consisting of those modifications $L'$ that are classical and flat (as $\CO_X$-modules). 

\medskip

Note that the functor of chiral universal envelope maps 
\begin{equation} \label{e:ch env good}
\on{Modif}_{\on{cl,flat}}(L)\to  \on{Modif}_{\on{conn,ev-c}}(\CA).
\end{equation}

\begin{prop} \label{p:Lie-star class}
Assume that $j^*(L)$ is finitely generated and locally free as a D-module. Then the (opposite of the) functor \eqref{e:ch env good} is cofinal. 
\end{prop} 

We note that the assertion of the proposition would be false without the finite generation assumption,
see \secref{ss:QCoh* non-fp}. 

\medskip

The proposition will be proved in \secref{ss:Lie-star class}. We will now consider some applications. 

\begin{cor}  \label{c:Lie-* class}
Under the assumptions of \propref{p:Lie-star class}, we have 
$$\sP_{\CA,x,\on{ev-c}}\simeq \underset{L'\in \on{Modif}_{\on{cl,flat}}(L)}{``\on{lim}"}\, \ind^{L\mod^{\on{ch}}_x}_{L'\mod^{\on{Lie-}^*}_x}(k).$$
\end{cor} 

\sssec{}  \label{sss:ass of top Lie}

Take 
$$L=L_\fg=\fg\otimes \on{D}_X,$$
where $\fg$ is a classical finite-dimensional Lie algebra (or a central extension of $L_\fg$). 

\medskip

In this case, the category $\on{Modif}_{\on{cl,flat}}(L)^{\on{op}}$ contains a cofinal family of objects of the
form
$$\fg\otimes \CO_X(-n\cdot x)\underset{\CO_X}\otimes \on{D}_X.$$
  
Note that the corresponding objects $\ind^{L\mod^{\on{ch}}_x}_{L'\mod^{\on{Lie-}^*}_x}(k)$ are the images
under \eqref{e:KM to Lie*} of 
$$\ind^{\hg}_{\fg_n}(k),$$
where $\fg_n\subset \hg$ is the $n$th congruence subalgebra. 

\medskip

Hence, combining Corollaries \ref{c:Lie-* class} and \ref{c:P ev c}, we obtain:

\begin{cor} \label{c:Lie-*}
The natural transformation
$$\underset{n}{\on{colim}}\, \CHom(\ind^{\hg}_{\fg_n}(k),-) \to \oblv_{L}, \quad (L\mod^{\on{ch}}_x)^{>-\infty}\to \Vect$$
is an isomorphism.
\end{cor} 

\sssec{} \label{sss:ass of top com}

Assume now that $\CM$ is abelian (and finitely generated and locally fee as a D-module). Note that in this case 
$$\CA=U^{ch}(\CM)\simeq \Sym^!(\CM[1])[1]$$
is a commutative chiral algebra (see \secref{sss:com ch free}). 

\medskip

In this case we can talk about \emph{commutative} chiral $\CA$-modules:
this is by definition the category of modules over the commutative algebra $\Sym(\CM_x)$, and it has a natural forgetful functor to 
$$\CA\mod_x^{\on{ch}}\simeq \CM\mod_x^{\on{ch}}.$$

\medskip

Denote
$$\on{Modif}_{\on{cl,flat,f.g.}}(\CM):=
\{\CM'\in \Dmod(X)^{\heartsuit,\on{flat,f.g.}},\,\, \CM'|_{X-x}\simeq j^*(\CM)\}.$$

\medskip

From Corollaries \ref{c:Lie-* class} and \ref{c:P ev c}, we obtain:
\begin{cor} \label{c:com ch}
The natural transformation
$$\underset{\CM'\in \on{Modif}_{\on{cl,flat,f.g.}}(\CM) }{\on{colim}}\, \CHom(\Sym(\CM'_x),-) \to \oblv_{\CM}, \quad (\CM\mod^{\on{ch}}_x)^{>-\infty}\to \Vect$$
is an isomorphism, where 
$$\Sym(\CM'_x)\in \Sym(\CM'_x)\mod \to \CM'\mod_x^{\on{ch}}\simeq  \CM\mod_x^{\on{ch}}.$$
\end{cor} 

\ssec{Proof of \propref{p:Lie-star class}}  \label{ss:Lie-star class}

\sssec{}

Let $\wt{L}$ be an eventually coconnective Lie-* algebra on $X$, and let us be given a map $$j^*(L)\overset{\alpha}\to j^*(\wt{L}).$$
Consider the category
$$\bC_{\wt{L}}:=\{L'\in \on{Modif}_{\on{cl,flat}}(L),\, L'\overset{\alpha'}\to \wt{L}\, |\, \alpha'|_{X-x}=\alpha\}.$$ 

We will prove:

\begin{prop} \label{p:Lie-star class nonempty}
Assume that $j^*(L)$ is finitely generated and locally free as a D-module. Then the category $\bC_{\wt{L}}$ is non-empty.
\end{prop}

Let us show how \propref{p:Lie-star class nonempty} implies \propref{p:Lie-star class}.

\begin{proof}[Proof of \propref{p:Lie-star class}]

By adjunction, it suffices to prove that in the setup of \propref{p:Lie-star class nonempty}, the category 
$\bC_{\wt{L}}$ is contractible. We will show that it is cofiltered. 

\medskip

Let
$$L'_I:I\to \bC_{\wt{L}}$$
be a finite diagram. We need to show that it can be extended to a diagram
$$L'_{I^{\triangleleft}}:I^{\triangleleft}\to \bC_{\wt{L}},$$
where $I^{\triangleleft}$ is a left cone over $I$. 

\medskip

Set 
$$\wt{L}_I:=\underset{I}{\on{lim}}\, L'_I,$$
where the limit is taken in the category of Lie-* algebras over $\wt{L}$.  Note that $\wt{L}_I$ is eventually coconnective. 
By construction,
we have a map
$$\alpha_I:j^*(L)\to j^*(\wt{L}_I).$$

The datum of $L'_{I^{\triangleleft}}$
is equivalent to finding $L'\in \on{Modif}_{\on{cl,flat}}(L)$ and a map
$$\alpha'_I:L'\to \wt{L}_I,$$
extending $\alpha_I$.

\medskip

However, the existence of $(L',\alpha'_I)$ is guaranteed by \propref{p:Lie-star class nonempty}

\end{proof}[\propref{p:Lie-star class}]

\sssec{}

The rest of this subsection is devoted to the proof of \propref{p:Lie-star class nonempty}. 
The starting point is the following observation (\cite[Lemma 2.5.13]{BD2}):

\begin{lem} \label{l:Lie* cofinal}
Assume that $j^*(L)$ is classical, and the underlying D-module is finitely generated and 
locally free.  
Then the tautological map
$$j_!\circ j^*(\oblv_{\on{Lie-}^*}(L))\to \underset{L'\in \on{Modif}_{\on{cl,flat}}(L)}{``\on{lim}"}\, \oblv_{\on{Lie-}^*}(L')$$
is an isomorphism in $\on{Pro}(\Dmod(X))$, where:

\medskip

\begin{itemize}

\item $\oblv_{\on{Lie-}^*}$ is the forgetful functor from the category of Lie-* algebras to the category of D-modules on $X$;

\item $j_!:\Dmod(X-x)\to \on{Pro}(\Dmod(X))$ is the pro-left adjoint of $j^*$.

\end{itemize}

\end{lem} 

\sssec{}

Let $\wt{L}$ be concentrated in degrees $[-n,0]$. We will argue by induction on $n$.

\medskip

Consider first the case $n=0$. 

\medskip

By \lemref{l:Lie* cofinal}, we can find $L'\in  \on{Modif}_{\on{cl,flat}}(L)$, so that
$\alpha$ extends to a map 
$$\alpha':L'\to \wt{L}$$
\emph{as plain D-modules}. 

\medskip

The obstruction to $\alpha'$ being a map of Lie-* algebras is a map
$$L'\boxtimes L'\to \Delta_*(\wt{L}),$$
which vanishes on $(X-x)\times (X-x)$. 

\medskip

The assumption that $j^*(L)$ is finitely generated implies that the naturally defined map
$$(j\times j)_!\circ (j\times j)^*(L\boxtimes L)\to j_!\circ j^*(L)\boxtimes j_!\circ j^*(L)$$
in $\on{Pro}(\Dmod(X\times X))$ is an isomorphism. Hence, again by \lemref{l:Lie* cofinal}, 
there exists an arrow $L''\to L'$ in $\on{Modif}_{\on{cl,flat}}(L)$ such that the composition 
$$L''\boxtimes L''\to L'\boxtimes L'\to \Delta_*(\wt{L})$$
vanishes.

\medskip

Hence, $L''\to L'\to \wt{L}$ provides the desired object of $\bC_{\wt{L}}$. 

\sssec{} 

We now perform the induction step. Suppose the assertion is valid for
$\tau^{\geq n-1}(\wt{L})$, i.e., that can find an object $L'\in \on{Modif}_{\on{cl,flat}}(L)$ and a lift of $\alpha$ to 
a map of Lie-* algebras
$$L'\to \tau^{\geq n-1}(\wt{L}).$$

Fix this map, and consider the fiber product
$$\wt{L}\underset{\tau^{\geq n-1}(\wt{L})}\times L'=:\wt{L}'.$$

We wish to find an arrow $L''\to L'$ in $\on{Modif}_{\on{cl,flat}}(L)$, so that the map
$$\wt{L}'\underset{L'}\times L''=:\wt{L}''\to L''$$
admits a left inverse. 

\medskip

By \lemref{l:Lie* cofinal}, after replacing $L'$, the extension
$$\wt{L}'\to L'$$
is given by an $(n+2)$-cocycle, which is a map  
$$(L')^{\boxtimes (n+2)}\to \Delta^n_*(\wt{L}),$$
which vanishes on $(X-x)^n$, where $\Delta^n$ denotes the main diagonal $X\to X^n$.

\medskip

Now, by the same argument as above, using the fact that $L$ is finitely generated, we can find an arrow 
$L''\to L'$ in $\on{Modif}_{\on{cl,flat}}(L)$ such that the composition 
$$(L'')^{\boxtimes n}\to (L')^{\boxtimes n}\to \Delta^n_*(\wt{L})$$
vanishes.

\medskip

Hence, the extension $\wt{L}''\to L''$ admits a splitting.

\qed[\propref{p:Lie-star class nonempty}]

\ssec{Proof of \thmref{t:from QCoh star to fact}} \label{ss:proof of fact over com}

We will prove the variant of the lemma with a fixed $\ul{x}=x\in \Ran$. The factorization
version is a variant of this in families.

\sssec{}

We start with the following observation: let us regard the assignment
\begin{equation} \label{e:QCoh star and lim}
\CY\rightsquigarrow \QCoh_{\on{co}}(\CY)^{>-\infty}, \quad (\CY_1\overset{f}\to \CY_2)\rightsquigarrow f_*
\end{equation} 
as a functor from the category of ind-affine ind-schemes to $\infty$-categories. 

\medskip

We claim:

\begin{prop} \label{p:QCoh star and lim}
The functor \eqref{e:QCoh star and lim} commutes with totalizations, in the sense that if $\CY^\bullet$
is a cosimplicial ind-affine ind-scheme and 
$$\CY\simeq \on{Tot}(\CF^\bullet),$$
where the limit is taken in $\on{PreStk}$, then the functor
$$\QCoh_{\on{co}}(\CY)^{>-\infty}\to  \on{Tot}(\QCoh_{\on{co}}(\CY^\bullet)^{>-\infty})$$
is an equivalence. 
\end{prop} 

The proposition will be proved in \secref{ss:QCoh star and lim}. We now proceed with the proof of 
\thmref{t:from QCoh star to fact}. 

\sssec{}

For an affine D-scheme $\CY$, let $A$ denote the corresponding commutative algebra in $\Dmod(X)$
so that $\oblv^l(A)$ is connective and $\CY=\Spec_X(A)$. 

\medskip

Let $\CA$ denote the corresponding (commutative) chiral 
algebra, see \secref{sss:com ch}, so that $\CA$ corresponds to the factorization algebra $\CO_\CY$, and 
$$\CO_\CY\mod^{\on{fact}}_x\simeq \CA\mod_x^{\on{ch}}.$$

\medskip

Over the next few subsections we will reduce the statement of \thmref{t:from QCoh star to fact} to the case
when $A=\Sym^!(\CM[1])$
for $\CM$ a classical locally free finitely generated D-module. 

\sssec{} \label{sss:reduction to free 1}

We interpret the functor $\Gamma(\fL_\nabla(\CY),-)^{\on{enh}}$ of \eqref{e:from QCoh star to fact}
as
$$\QCoh_{\on{co}}(\fL_\nabla(\CY))\to \CA\mod_x^{\on{ch}}.$$

This functor is t-exact and both categories are right complete in their respective t-structures. Hence, 
it is in enough to show that the functor $\Gamma(\fL_\nabla(\CY),-)^{\on{enh}}$ induces an equivalence
$$\QCoh_{\on{co}}(\fL_\nabla(\CY))^{\geq 0,\leq m}\to (\CA\mod_x^{\on{ch}})^{\geq 0,\leq m}$$
for every $m$. 

\medskip

Note now that if $A_1\to A_2$ is a map in $\on{ComAlg}(\Dmod(X))$, such that the induced map
$\tau^{\geq -m}(A_1)\to \tau^{\geq -m}(A_2)$ is an isomorphism, the corresponding functors
$$\QCoh_{\on{co}}(\fL_\nabla(\CY_2))\to \QCoh_{\on{co}}(\fL_\nabla(\CY_1)$$
and 
$$\CA_2\mod_x^{\on{ch}}\to \CA_1\mod_x^{\on{ch}}$$
induces equivalences on the corresponding $(-)^{\geq 0,\leq m}$ categories. 

\medskip

In particular, we obtain that it is enough to show that the functor 
$$\QCoh_{\on{co}}(\fL_\nabla({}^m\CY))^{\geq 0,\leq m}\to ({}^m\!\CA\mod_x^{\on{ch}})^{\geq 0,\leq m}$$
is an equivalence for $^m\!A:=\tau^{\geq -m}(A)$, and the corresponding $^m\CY$ and $^m\!\CA$. 

\sssec{} \label{sss:reduction to free 2}

By the assumption that $A$ is D-afp, we can find a simplicial object $A_\bullet$ in $\on{ComAlg}(\Dmod(X))^{\leq 0}$
with terms $A_n=\Sym^!(\CM_n[1])$, where $\CM_n$ is a classical\footnote{Recall that according to 
our conventions, this means that $\oblv^r(\CM_n)$ is classical, i.e., $\oblv^l(\CM_n)[1]$ is classical.}
locally free finitely generated D-module, such that $^m\!A$ is a retract of $\tau^{\geq -m}(|A_\bullet|)$,
see \secref{sss:afp D n retract}. 

\medskip

Hence, it is enough to prove that 
$$\QCoh_{\on{co}}(\fL_\nabla(\CY'))^{\geq 0,\leq m}\to (\CA'\mod_x^{\on{ch}})^{\geq 0,\leq m}$$
is an equivalence 
for $A':=\tau^{\geq -m}(|A_\bullet|)$ for $A_\bullet$ as above. 

\medskip

Applying \secref{sss:reduction to free 1} again, we obtain that it is enough to prove that 
$$\QCoh_{\on{co}}(\fL_\nabla(\CY''))^{\geq 0,\leq m}\to (\CA''\mod_x^{\on{ch}})^{\geq 0,\leq m}$$
is an equivalence for $A'':=|A_\bullet|$ for $A_\bullet$ as above. 

\medskip

Hence, we can obtain that it is enough to prove \thmref{t:from QCoh star to fact} for 
$A$ is of the form $|A_\bullet|$ for $A_\bullet$ as above. 

\sssec{} \label{sss:reduction to free 3}

For $A_\bullet$ as above set $\CY^n:=\Spec_X(A_n)$, and consider the corresponding simplicial
affine D-scheme $\CY^\bullet$, so that
$$\CY\simeq \on{Tot}(\CY^\bullet).$$

\medskip

It is clear that that the functor
$$\CY\mapsto \fL_\nabla(\CY)$$
preserves limits, so that
$$\fL_\nabla(\CY)\simeq |\fL_\nabla(\CY^\bullet)|.$$

Hence, by \propref{p:QCoh star and lim}, the functor
$$\QCoh(\fL_\nabla(\CY))^{>-\infty}\to \on{Tot}(\QCoh(\fL_\nabla(\CY^\bullet)))^{>-\infty}$$
is an equivalence.

\medskip

The functor
$$\CA\mod_x^{\on{ch}}\to\on{Tot}(\CA_n\mod_x^{\on{ch}})$$
is also an equivalence: indeed, this is obvious for non-unital modules 
(this is a general property of categories of modules over operad algebras),
and this property is inherited by unital modules by \cite[Proposition 3.8.4]{CR}. 

\medskip

Hence, it is enough to show that the functors
$$\QCoh_{\on{co}}(\fL_\nabla(\CY^n))^{>-\infty}\to (\CA\mod_x^{\on{ch}})^{>-\infty}$$
are equivalences for every $n$. 

\medskip

This reduces the assertion of \thmref{t:from QCoh star to fact} to the case when $A=\Sym^!(\CM[1])$
for $\CM$ a classical locally free finitely generated generated D-module. 

\sssec{}

Let $\CY$ be general an affine D-scheme. Consider the functor 
$$\Gamma(\fL_\nabla(\CY),-)^{\on{enh}}:\QCoh_{\on{co}}(\fL_\nabla(\CY))\to \CA\mod_x^{\on{ch}}$$
and its (a priori discontinuous) right adjoint. Tautologically, we have
$$\Gamma(\fL_\nabla(\CY),-) \simeq \oblv_\CA\circ \Gamma(\fL_\nabla(\CY),-)^{\on{enh}},$$
hence we obtain a natural transformation 
\begin{equation} \label{e:oblv and G R}
\Gamma(\fL_\nabla(\CY),-)  \circ (\Gamma(\fL_\nabla(\CY),-)^{\on{enh}})^R\to \oblv_\CA.
\end{equation}

We claim that it suffices to show that \eqref{e:oblv and G R} is an isomorphism when evaluated
on $(\CA\mod_x^{\on{ch}})^{>-\infty}$. 

\medskip

Indeed, this follows from the fact that both 
$\oblv_\CA$ and $\Gamma(\fL_\nabla(\CY),-)$ are conservative on the eventually coconnective subcategories
(see \secref{sss:Gamma cons}).  

\sssec{}

We will now specialize to the case when $\CA=\Sym^!(\CM[1])$
for $\CM\in \Dmod^l(X)^{\heartsuit,\on{loc.free,f.g.}}$, and prove that \eqref{e:oblv and G R} is an isomorphism 
on $(\CA\mod_x^{\on{ch}})^{>-\infty}$ by an explicit calculation. 

\sssec{} \label{sss:describe loops}

Let $\on{Modif}_{\on{cl,flat,f.g.}}(\CM)$ be the category 
$$\{\CM'\in  \Dmod^l(X)^{\heartsuit,\on{flat,f.g.}},\,\, \CM'|_{X-x}\simeq \CM|_{X-x}\}.$$

We claim that the ind-scheme $\fL_\nabla(\CY)$ identifies in this case with
$$\underset{\CM'\in \on{Modif}_{\on{cl,f.g.,flat}}}{``\on{colim}"}\, \Spec(\Sym(\CM'_x)).$$

Indeed, since the question is local around $x$, with no restriction of generality we can assume that $X$ is affine. Let $t$
be a uniformizer at $x$. Then for a connective commutative algebra $R$, we have by definition
$$\Maps(\Spec(R),\fL_\nabla(\CY))=\Maps_{\on{ComAlg}(\Dmod(X))}(\Sym_{\CO_X}(\CM),R\ppart)\simeq
\Maps_{\Dmod(X)}(\CM,R\ppart)$$

Since the map
$$j_!\circ j^*(\CM)\to \underset{\CM'\in \on{Modif}_{\on{cl,flat,f.g.}}}{``\on{lim}"}\, \CM'$$ 
is an isomorphism in $\on{Pro}(\Dmod(X))$, we have
\begin{multline*}
\Maps_{\Dmod(X)}(\CM,R\ppart)\simeq \underset{\CM'}{\on{colim}}\, \Maps_{\Dmod(X)}(\CM',R\qqart)\simeq 
\underset{\CM'}{\on{colim}}\,  \Maps_{\Vect}(\CM'_x,R)\simeq \\
\simeq \underset{\CM'}{\on{colim}}\, \Maps_{\on{ComAlg}(\Vect)}(\Sym(\CM'_x),R)=
\underset{\CM'}{\on{colim}}\, \Maps(\Spec(R),\Spec(\Sym(\CM'_x)))=\\
=:\Maps\left(\Spec(R),\underset{\CM'}{``\on{colim}"}\, \Spec(\Sym(\CM'_x))\right),
\end{multline*}
as desired. 

\sssec{} \label{sss:proof of QCoh * to fact for flat}

For $\CM'$ as above, let us denote by 
$$\CO_{(\CM'_x)^\vee}\in \QCoh_{\on{co}}(\fL_\nabla(\CY))$$
the direct images of the structure sheaf along the map,
$$\Spec(\Sym(\CM'_x)) \to \fL_\nabla(\CY).$$

\medskip

The above description of $\fL_\nabla(\CY)$ implies that the functor $\Gamma(\fL_\nabla(\CY),-)$ is isomorphic to 
$$\underset{\CM'}{\on{colim}}\, \CHom(\CO_{(\CM'_x)^\vee},-).$$

The functor
$$\Gamma(\fL_\nabla(\CY),-)^{\on{enh}}:\QCoh_{\on{co}}(\fL_\nabla(\CY))\to \CA\mod_x^{\on{ch}}$$
sends $\CO_{(\CM'_x)^\vee}$ to
$$\Sym(\CM'_x)\in \CA\mod_x^{\on{ch}}.$$

\medskip

Now, the required isomorphism follows from \corref{c:com ch}.

\qed[\thmref{t:from QCoh star to fact}]

\ssec{Failure of \thmref{t:from QCoh star to fact} in the non-finitely presented case} \label{ss:QCoh* non-fp}

In this subsection we will explain why \thmref{t:from QCoh star to fact} does \emph{not} hold when $\CY$
is \emph{not} almost finitely presented (in the D-sense).

\begin{rem}
One can show that the functor 
\begin{equation} \label{e:from QCoh star to fact again}
\QCoh_{\on{co}}(\fL_\nabla(\CY))^{>-\infty}\to (\CA\mod_x^{\on{ch}})^{>-\infty}
\end{equation}
induces an equivalence of the abelian categories
$$\QCoh_{\on{co}}(\fL_\nabla(\CY))^\heartsuit\to (\CA\mod_x^{\on{ch}})^\heartsuit.$$

So, the failure of \eqref{e:from QCoh star to fact again} to be 
an equivalence occurs at the derived level.
\end{rem}

\sssec{}

We will take $\CA=\Sym_{\CO_X}(\CM)$, where $\CM$ is an infinitely-generated D-module
(i.e., the direct sum of countably many copies of $\on{D}_X$). We will show that 
\eqref{e:from QCoh star to fact again} fails to be an equivalence in this case.

\medskip

Namely, we will construct two objects 
$\CF_1,\CF_2\in (\QCoh_{\on{co}}(\fL_\nabla(\CY))^\heartsuit$ with images in 
$(\CA\mod_x^{\on{ch}})^\heartsuit$ denoted $\CM_1,\CM_2$, respectively, and an element in 
$\Ext^2_{\CA\mod_x^{\on{ch}}}(\CM_1,\CM_2)$ 
that does not come from an element in $\Ext^2_{\QCoh_{\on{co}}(\fL_\nabla(\CY))}(\CF_1,\CF_2)$. 

\begin{rem}
It follows from the description of the category $\QCoh_{\on{co}}(\fL_\nabla(\CY))^\heartsuit$ in 
\secref{sss:QCoh* as comonad} that the category
$\QCoh_{\on{co}}(\fL_\nabla(\CY))^b$
is the bounded derived category of its heart. So the above inequality of the $\Ext^2$ spaces means
that $(\CA\mod_x^{\on{ch}})^b$ is \emph{not} the bounded derived category of its heart.
\end{rem} 

\sssec{}

Write 
$$\CM=\underset{i}{\on{colim}}\, \CM_i,$$
where $\CM_i$ are finitely generated.

\medskip

Let us call a modification $\CM'$ ``quasi-finitely generated" if all the intersections $\CM'\cap \CM_i$ are finitely generated.
As in \secref{sss:describe loops}, we have
$$\fL_\nabla(\CY)\simeq \underset{\CM'\in  \on{Modif}_{\on{cl,flat,q-f.g.}}}{``\on{colim}"}\, \Spec(\Sym(\CM'_x)).$$

By the same logic as in \secref{sss:proof of QCoh * to fact for flat}, the functor
$$\CF\mapsto \underset{\CM'\in  \on{Modif}_{\on{cl,flat,q-f.g.}}}{\on{colim}}\, 
\CHom(\CO_{(\CM'_x)^\vee},-), \quad \QCoh_{\on{co}}(\fL_\nabla(\CY))\to \Vect$$
identifies with $\Gamma(\fL_\nabla(\CY),-)$. In particular, it is t-exact. 

\medskip

Hence, if for some $\CM'$ and $\CF\in \QCoh_{\on{co}}(\fL_\nabla(\CY))^\heartsuit$ we have a class
$$\alpha'\in \Ext^2_{\QCoh_{\on{co}}(\fL_\nabla(\CY))}(\CO_{(\CM'_x)^\vee},\CF),$$
we can find $\CM''\subset \CM'$ such that the image $\alpha''$ of $\alpha'$ in 
$$\Ext^2_{\QCoh_{\on{co}}(\fL_\nabla(\CY))}(\CO_{(\CM''_x)^\vee},\CF)$$
vanishes.

\sssec{}

As in \secref{sss:proof of QCoh * to fact for flat}, the image of $\CO_{(\CM'_x)^\vee}$ in
$\CA\mod_x^{\on{ch}}$ is $\Sym(\CM'_x)\in \CA\mod_x^{\on{ch}}$.

\medskip

By adjunction for any $\CM\in \CA\mod_x^{\on{ch}}$,
we have
$$\CHom_{\CA\mod_x^{\on{ch}}}(\Sym(\CM'_x),\CM)\simeq \CHom_{\CM'\mod_x^{\on{Lie-}^*}}(k,\CM).$$

\sssec{}

Hence, it suffices to find $\CF\in \QCoh_{\on{co}}(\fL_\nabla(\CY))^\heartsuit$ with image
$\CM\in (\CA\mod_x^{\on{ch}})^\heartsuit$ and a class 
$$\beta'\in \Ext^2_{\CM'\mod_x^{\on{Lie-}^*}}(k,\CM)$$
such that for any $\CM''\subset \CM'$, the image $\beta''$ of $\beta'$ in 
$$\Ext^2_{\CM''\mod_x^{\on{Lie-}^*}}(k,\CM)$$
is non-zero. 

\medskip

We will take $\CF$ to be the sky-scraper at the origin of $\CY$, so that $\CM=k$, with the trivial 
chiral action of $\CM$.  

\sssec{}

We calculate
$$\CHom_{\CM'\mod_x^{\on{Lie-}^*}}(\CM_1,\CM_2)\simeq \underset{i}{\on{lim}}\, 
\CHom_{(\CM'\cap \CM_i)\mod_x^{\on{Lie-}^*}}(\CM_1,\CM_2).$$

For $\CM_1=\CM_2=k$, we obtain
\begin{multline*} 
\Ext^2_{\CM'\mod_x^{\on{Lie-}^*}}(k,k)\simeq \underset{i}{\on{lim}}\, \Lambda^2(\BD(\CM'\cap \CM_i)_x)\simeq \\
\simeq \underset{i}{\on{lim}}\, \Hom_{\Dmod(X\times X)^\heartsuit}^{\on{antisym}}((\CM'\cap \CM_i)\boxtimes (\CM'\cap \CM_i),\delta_{x,x})\simeq
\Hom^{\on{antisym}}_{\Dmod(X\times X)^\heartsuit}(\CM'\boxtimes \CM',\delta_{x,x}).
\end{multline*}

Now, it is clear that since $\CM'$ is infinitely generated, we can find an element
$$\beta'\in \Hom^{\on{antisym}}_{\Dmod(X\times X)^\heartsuit}(\CM'\boxtimes \CM',\delta_{x,x}),$$
such that for any $\CM''$, the restriction $\beta''$ of $\beta'$ to 
$$\Hom^{\on{antisym}}_{\Dmod(X\times X)^\heartsuit}(\CM''\boxtimes \CM'',\delta_{x,x})$$
is non-zero.

\begin{rem}
Note that the above counterexample does not work for $\Ext^1$ instead of $\Ext^2$
(as must be the case, since \eqref{e:from QCoh star to fact again} is an equivalence at the 
abelian level).  Indeed, for any 
$$\gamma'\in \Hom_{\Dmod(X)^\heartsuit}(\CM',\delta_{x}),$$
there exists $\CM''\subset \CM'$, such that the restriction $\gamma''$ of $\gamma'$ to
$$\Hom_{\Dmod(X)^\heartsuit}(\CM'',\delta_{x})$$
vanishes.
\end{rem}

\ssec{Proof of \propref{p:QCoh star and lim}: preparations} \label{ss:QCoh star and lim prep} 

\sssec{} \label{sss:Gamma cons}

First, we claim that if $\CY$ is an ind-affine ind-scheme, then the functor
$$\Gamma(\CY,-):\QCoh_{\on{co}}(\CY)\to \Vect$$
is conservative on $\QCoh_{\on{co}}(\CY)^{>-\infty}$. 

\medskip

Indeed, since $\Gamma(\CY,-)$ is t-exact and the t-structure on $\QCoh_{\on{co}}(\CY)$ is right-complete, 
it suffices to show that $\Gamma(\CY,-)$ does not annihilate objects from $\QCoh_{\on{co}}(\CY)^\heartsuit$. 

\medskip

Write $\CY$ is a filtered colimit of schemes $Y_\alpha$ under closed embeddings
$$f_{\alpha,\beta}:Y_\alpha\to Y_\beta.$$

An object $\CF\in \QCoh_{\on{co}}(\CY)^\heartsuit$ amounts to a collection
$$\{\CF_\alpha\in \QCoh(Y_\alpha)^\heartsuit, \CF_\alpha\simeq H^0(f_{\alpha,\beta}^!(\CF_\beta))\}.$$

In particular, the maps
$$\Gamma(\CY_\alpha,\CF_\alpha)\to \Gamma(\CY_\beta,\CF_\beta)$$ are
injective.  

\medskip

We have
$$\Gamma(\CY,\CF)\simeq \underset{\alpha}{\on{colim}}\, \Gamma(\CY_\alpha,\CF_\alpha)$$
and the statement is manifest. 

\sssec{} \label{sss:Gamma comonadic}

We claim:

\begin{prop} \label{p:Gamma comonadic}
The functor
$$\Gamma(\CY,-):\QCoh_{\on{co}}(\CY)^{\geq 0}\to \Vect^{\geq 0}$$
is comonadic. 
\end{prop} 

\medskip

Given the conservativity, the assertion of the proposition follows from the next general observation:

\medskip

\begin{lem} \label{l:comonadic}
Let $\bC,\bD$ be cocomplete DG categories, equipped with t-structures, compatible with filtered
colimits. Assume that $\bD$ is right-complete in its t-structure. Let $F:\bC\to \bD$ be a t-exact continuous functor. 
Assume that $F$ is conservative on $\bC^{>-\infty}$. Then the induced functor
$$\bC^{\geq 0}\to \bD^{\geq 0}$$
is comonadic.
\end{lem} 

\begin{proof}

By Barr-Beck-Lurie, suffices to show that $F$ preserves totalizations of cosimplicial objects in
$\bC^{\geq 0}$. Thus, let $\bc^\bullet$ be a cosimplicial object in $\bC$. We have to show that the map
$$F(\on{Tot}(\bc^\bullet))\to \on{Tot}(F(\bc^\bullet))$$
is an isomorphism. 

\medskip

Since $\bD$ is right-complete in its t-structure, it suffices to show that for every $n$,
$$\tau^{\leq n}(F(\on{Tot}(\bc^\bullet))) \to \tau^{\leq n}(\on{Tot}(F(\bc^\bullet)))$$
is an isomorphism.

\medskip

Let 
$$\on{Tot}^{\leq n}(\bc^\bullet) \text{ and } \on{Tot}^{\leq n}(F(\bc^\bullet))$$
be the totalizations of the corresponding $(n+1)$-skeleta.

\medskip

We have natural maps
$$\on{Tot}(\bc^\bullet)\to \on{Tot}^{\leq n}(\bc^\bullet)  \text{ and } 
\on{Tot}(F(\bc^\bullet))\to \on{Tot}^{\leq n}(F(\bc^\bullet)).$$

Since the terms of $\bc^\bullet$ and $F(\bc^\bullet)$ 
are in $\bC^{\geq 0}$, the above maps induce isomorphisms between the $\tau^{\leq n}$ truncations. 

\medskip

Hence, it suffices to show that
$$F(\on{Tot}^{\leq n}(\bc^\bullet))\to \on{Tot}^{\leq n}(F(\bc^\bullet))$$
is an isomorphism. 

\medskip

However, this is obvious, since the limit over $\bDelta_{\leq n}$ is a finite limit. 

\end{proof}

\sssec{} \label{sss:QCoh* as comonad}

Note that we have a canonical equivalence
$$\on{Pro}(\Vect)^{\on{op}}\simeq \on{Funct}_{\on{discnt}}(\Vect,\Vect),\quad \bV\mapsto \CHom(\bV,-).$$
where $\on{Funct}_{\on{discnt}}(-,-)$ denotes the category of exact $k$-linear functors that are not necessarily continuous.

\medskip

Under this equivalence, the monoidal structure on $\on{Funct}_{\on{discnt}}(-,-)$ given by composition 
corresponds to the $\arrowtimes$ monoidal structure on $\on{Pro}(\Vect)$ (see \cite{Bei}):

\medskip

For 
$$\bV=\underset{i\in I}{``\on{lim}"}\, V_i \text{ and } \bW=\underset{j\in j}{``\on{lim}"}\, W_j,$$
we have 
$$\bV\arrowtimes \bW=
\underset{j\in j}{\on{lim}}\,  \left( \underset{V_j^f\subset V_j}{\on{colim}}\, \left((\underset{i\in I}{``\on{lim}"}\, V_i)\otimes V_j^f\right)\right),$$
where:

\begin{itemize}

\item $V_j^f$ runs the category of compact objects mapping to $V_j$;

\item In the right-hand side, the outer limit and the inner colimit are taken in $\on{Pro}(\Vect)$.

\end{itemize} 

Thus, comonads on $\Vect$ correspond to algebra objects in $\on{Pro}(\Vect)$ with respect to $\arrowtimes$.
Comonads that are left t-exact correspond to algebra objects in $\on{Pro}(\Vect^{\leq 0})$. 

\sssec{} \label{sss:descr comonad}

Let $\sM_\CY$ denote the comonad on $\Vect$ corresponding to the functor $\Gamma(\CY,-)$, so that
$$\QCoh_{\on{co}}(\CY)^{\geq 0}\simeq \sM_\CY\comod(\Vect^{\geq 0}).$$

The object of $\on{Pro}(\Vect^{\leq 0})$ that corresponds to $\sM_\CY$ is described as follows.

\medskip

Let $\CO_\CY$ be the object of $\on{Pro}(\on{ComAlg}(\Vect^{\leq 0}))$ associated to the ind-affine ind-scheme $\CY$.
I.e., if 
$$\CY=\underset{i}{``\on{colim}"}\, Y_i,$$
then 
$$\CO_\CY:=\underset{i}{``\on{lim}"}\, \CO_{Y_i}.$$

Let 
$$\oblv^{\on{Pro}}_{\on{ComAlg}}:\on{Pro}(\on{ComAlg}(\Vect))\to \on{Pro}(\Vect)$$ be the pro-extension of the functor
$$\oblv_{\on{ComAlg}}:\on{ComAlg}(\Vect)\to \Vect.$$

\medskip

The endofunctor of $\Vect$ underlying the comonad $\sM_\CY$ is given by 
$$V\mapsto \CHom(\oblv^{\on{Pro}}_{\on{ComAlg}}(\CO_\CY),V).$$

\sssec{}

The proof of \propref{p:QCoh star and lim} will be based on the following two observations: 


\begin{lem} \label{l:O commute with limits b}
For every natural number $n$, the functor $\oblv^{\on{Pro}}_{\on{ComAlg}}:\on{Pro}(\on{ComAlg}(\Vect^{\leq 0}))\to \on{Pro}(\Vect^{\leq 0})$,
followed by the truncation 
$$\on{Pro}(\Vect^{\leq 0})\to \on{Pro}(\Vect^{\leq 0,\geq -n})$$
commutes with geometric realizations.
\end{lem} 

\begin{lem} \label{l:O commute with limits c}
For every natural number $n$ and for every natural number $m$, the functor
$$\bV\mapsto \bV^{\arrowtimes m}, \quad \on{Pro}(\Vect^{\leq 0})\to \on{Pro}(\Vect^{\leq 0}),$$
followed by the truncation 
$$\on{Pro}(\Vect^{\leq 0})\to \on{Pro}(\Vect^{\leq 0,\geq -n})$$
commutes with geometric realizations. 
\end{lem} 

Let us temporarily assume these lemmas and prove \propref{p:QCoh star and lim}. 

\ssec{Proof of \propref{p:QCoh star and lim}} \label{ss:QCoh star and lim} 

\sssec{}

It is enough to show that for every $n$, the functor
\begin{equation} \label{e:QCoh star and lim n}
\QCoh_{\on{co}}(\CY)^{\geq 0,\leq n}\to \on{Tot}(\QCoh_{\on{co}}(\CY^\bullet)^{\geq 0,\leq n})
\end{equation} 
is an equivalence.

\sssec{}

Note that $\sM$ is a left-exact endofunctor $\Vect$, we can create its truncation
$\sM^{\leq n}$ that acts on $\Vect^{\geq 0,\leq n}$, namely,
$$\sM^{\leq n}(V):=\tau^{\leq n}(\sM(V)).$$

This assignment is monoidal, so if $\sM$ is a comonad, then $\sM^{\leq n}$ inherits a natural comonad structure, and we have
$$\sM\comod(\Vect^{\geq 0})\underset{\Vect^{\geq 0}}\times \Vect^{\geq 0,\leq n}\simeq \sM^{\leq n}\comod(\Vect^{\geq 0,\leq n}).$$

\sssec{}

According to \propref{p:Gamma comonadic}, we have:
$$\QCoh_{\on{co}}(\CY)^{\geq 0}\simeq \sM_\CY\comod(\Vect^{\geq 0}) \text{ and }
\QCoh_{\on{co}}(\CY^\bullet)^{\geq 0}\simeq \sM_{\CY^\bullet}\comod(\Vect^{\geq 0}).$$

Hence, 
$$\QCoh_{\on{co}}(\CY)^{\geq 0,\leq n}\simeq \sM^{\leq n}_\CY\comod(\Vect^{\geq 0,\leq n}) 
\text{ and }
\QCoh_{\on{co}}(\CY^\bullet)^{\geq 0,\leq n}\simeq \sM^{\leq n}_{\CY^\bullet}\comod(\Vect^{\geq 0,\leq n}).$$

\sssec{}

We now observe:

\begin{lem} \label{l:inverse limit of comonads}
Let $I$ be an index category and let 
$$\sM_I:I\to \on{Comonad}(\bC), \quad i\mapsto \sM_i$$ be an $I$-diagram 
of comonads on a category $\bC$. Let $\sM$ be another comonad, equipped
with a compatible collection of maps $\sM\to \sM_i$. 
Suppose that for every natural number $m$, the map 
$$\sM^{\times m}\to \underset{i}{\on{lim}}\, \sM_i^{\times m}$$
is an isomorphism, where:

\begin{itemize}

\item The notation $\sM^{\times m}$ means an $m$-fold composition
of $\sM$, and similarly for $\sM_i$;

\item The limit in the right-hand side is taken in the category of endofunctors of $\bC$
(i.e., is computed value-wise).

\end{itemize}

Then the induced functor
$$\sM\comod(\bC)\to  \underset{i}{\on{lim}}\, \left(\sM_i\comod(\bC)\right)$$
is an equivalence.

\end{lem} 

\sssec{}

Hence, in order to prove that \eqref{e:QCoh star and lim n} is an equivalence, it suffices to show
that for any $m$, the map
\begin{equation} \label{e:tot of comonads}
(\sM^{\leq n}_\CY)^{\times m}\to \on{Tot}\left((\sM^{\leq n}_{\CY^\bullet})^{\times m}\right)
\end{equation} 
is an equivalence. 

\sssec{}

By \secref{sss:descr comonad}, we have
$$\sM_\CY\comod(\Vect^{\leq 0}) \simeq \CO_\CY\mod(\Vect^{\leq 0})$$ 
and 
$$\sM_{\CY^\bullet}\comod(\Vect^{\leq 0}) \simeq \CO_{\CY^\bullet}\mod(\Vect^{\leq 0}),$$
where we regard $\CO_\CY$ as an algebra object in $(\on{Pro}(\Vect),\arrowtimes)$ via the natural
forgetful functor 
$$\on{Pro}(\on{ComAlg}(\Vect))\to \on{AssocAlg}(\on{Pro}(\Vect),\arrowtimes).$$

Note now that in the situation of \propref{p:QCoh star and lim}, the map
$$\CO_\CY\to |\CO_{\CY^\bullet}|$$
is an isomorphism, where the geometric realization is taken in $\on{Pro}(\on{ComAlg}(\Vect))$.

\medskip 

Hence, \ref{l:O commute with limits b} and \ref{l:O commute with limits c} guarantee that the maps
\eqref{e:tot of comonads} are isomorphisms, as required. 

\qed[\propref{p:QCoh star and lim}]

\ssec{Proof of Lemmas \ref{l:O commute with limits b} and \ref{l:O commute with limits c}} 

\sssec{}

The key input is the following:

\begin{lem} \label{l:trunc geom real}
Let $\bC$ be an $n$-truncated category (i.e., the mapping spaces have homotopy groups
$\pi_m$ vanish for $m>n$). Then the map
$$|-|_{\leq n+1}\to |-|, \quad \on{Funct}(\bDelta^{\on{op}},\bC)\to \bC$$
is an isomorphism, where the left-hand side is the colimit over $\bDelta^{\on{op}}_{\leq n+1}$. 
\end{lem}

\sssec{Proof \lemref{l:O commute with limits c}}

Since $\bDelta^{\on{op}}$ is sifted, it suffices to show that the binary operation
\begin{equation} \label{e:sotimes in Pro}
\on{Pro}(\Vect^{\leq 0,\geq -n})\times \on{Pro}(\Vect^{\leq 0,\geq -n})\to \on{Pro}(\Vect^{\leq 0,\geq -n}), \quad
\bV,\bW\mapsto \bV\arrowtimes \bW
\end{equation} 
commutes with geometric realizations in each variable. 

\medskip

However, by \lemref{l:trunc geom real}, the functor of geometric realization in $\Vect^{\leq 0,\geq -n}$
is a \emph{finite colimit}, and it is clear that \eqref{e:sotimes in Pro} 
commutes with finite colimits in each variable. 

\qed[\lemref{l:O commute with limits c}] 

\sssec{Proof of \lemref{l:O commute with limits b}}

By \lemref{l:trunc geom real}, it suffices to show
that the functor 
$$\oblv^{\on{Pro}}_{\on{ComAlg}}:\on{Pro}(\on{ComAlg}(\Vect^{\leq 0,\geq -n}))\to \on{Pro}(\Vect^{\leq 0,\geq -n})$$
commutes with colimits over $\bDelta^{\on{op}}_{\leq n+1}$. 

\medskip

The rest of the argument essentially reproduces \cite[Proposition 6.1.5.3]{Lu2}: 

\medskip

Note that for any finite index category $I$
and a category $\bC$, the naturally defined functor
$$\on{Pro}\on{Funct}(I,\bC)\to \on{Funct}(I,\on{Pro}(\bC))$$
is an equivalence.

\medskip

Furthermore, for an object
$$\underset{\alpha}{``\on{lim}"}\, (i\mapsto \bc_{i,\alpha}) \in \on{Pro}(\on{Funct}(I,\bC))$$
and the resulting object
$$i\mapsto \underset{\alpha}{``\on{lim}"}\, \bc_{i,\alpha}\in \on{Funct}(I,\on{Pro}(\bC)),$$
we have
$$\underset{i\in I}{\on{colim}}\, \left(\underset{\alpha}{``\on{lim}"}\, \bc_{i,\alpha}\right) 
\simeq \underset{\alpha}{``\on{lim}"}\, \left(\underset{i\in I}{\on{colim}}\, \bc_{i,\alpha}\right).$$

We apply this observation to $I=\bDelta^{\on{op}}_{\leq n+1}$ and $\bC$ being
$$\on{ComAlg}(\Vect^{\leq 0,\geq -n}) \text{ and } \Vect^{\leq 0,\geq -n}.$$

Hence, in order to prove \lemref{l:O commute with limits b}, it suffices to show that the functor
$$\oblv_{\on{ComAlg}}:\on{ComAlg}(\Vect^{\leq 0,\geq -n}) \to \Vect^{\leq 0,\geq -n}$$
commutes with colimits over $\bDelta^{\on{op}}_{\leq n+1}$. 

\medskip

However, this follows from \lemref{l:trunc geom real} combined with the fact that the usual
geometric realization functor commutes with $\oblv_{\on{ComAlg}}$. 

\qed[\lemref{l:O commute with limits b}]

\section{The spectral spherical category} \label{s:spec Sph fact}

Throughout this section we let $H$ be an arbitrary finite-dimensional algebraic group. Our goal
is to define the factorization category
$$\IndCoh^*(\on{Hecke}_H^{\on{spec,loc}}).$$

\medskip

When $H=\cG$, the category $\IndCoh^*(\on{Hecke}_\cG^{\on{spec,loc}})=:\Spc_\cG^{\on{spec}}$
is the spectral counterpart of $\Sph_G$, and it acts by Hecke functors on the global spectral category. 
This action will play a key role in the sequel to this paper. 

\medskip

The difficulty we face is that we have not found a way to plug $\on{Hecke}_H^{\on{spec,loc}}$ into one of the 
previously discussed constructions, i.e.,  
$$\QCoh(-),\,\, \QCoh_{\on{co}}(-),\,\, \text{ or }\,\, \IndCoh^*(-)$$
to obtained the desired category. 

\medskip

Instead, we will define it as bi-coinvariants with respect to $\fL^+_\nabla(H)$ inside $\IndCoh^*(\fL_\nabla(H))$,
where the latter also requires some care, as the factorization scheme $\fL_\nabla(H)$ is \emph{not}
ind-placid. 

\ssec{A 1-affineness property of \texorpdfstring{$\LS_H^\reg$}{LSHreg}}

Throughout this subsection we fix an affine scheme $S$ and a map $\ul{x}:S\to \Ran$. 

\sssec{}

%

Consider $\QCoh(\fL^+_\nabla(H))_S$
as an $\QCoh(S)$-linear monoidal category with respect to convolution. 

\medskip

Note that we have:
$$\on{Funct}_{\QCoh(\fL^+_\nabla(H))_S\mmod}(\QCoh(S),\QCoh(S))\simeq \QCoh(\LS^\reg_{H,S})$$
as monoidal categories.

\medskip

This gives rise to a pair of adjoint functors 
\begin{equation} \label{e:1-aff rep H}
\QCoh(\LS^\reg_{H,S})\mmod \rightleftarrows \QCoh(\fL^+_\nabla(H))_S\mmod,
\end{equation}
$$\bC\mapsto \bC\underset{\QCoh(\LS^\reg_{H,S})}\otimes \QCoh(S), \quad \bC'\mapsto (\bC')^{\fL^+_\nabla(H)_S}.$$

\begin{rem}
Throughout this section, the symbols $(-)^{\fL^+_\nabla(H)_S}$ and $(-)_{\fL^+_\nabla(H)_S}$ indicate 
\emph{weak}\footnote{Here ``weak" is as opposed to ``strong". Note that we could not even talk about 
strong invariants/coinvariants, because we are talking about weak actions.} 
invariants/coinvariants of the group-scheme $\fL^+_\nabla(H)_S$ acting on a $\QCoh(S)$-linear
category. 

\end{rem} 

\sssec{}

We have the following assertion (\cite[Lemma 9.8.1]{Ra4}):

\begin{prop}  \label{p:1-aff rep H}
The functors \eqref{e:1-aff rep H} are mutually inverse equivalences. 
\end{prop}

\sssec{}

Let $\bC$ be a module category over $\QCoh(\fL^+_\nabla(H))_S$. The functor of $\fL^+_\nabla(H)_S$-averaging
$$\on{Av}^{\fL^+_\nabla(H)_S}_*:\bC\to \bC^{\fL^+_\nabla(H)_S}$$
naturally factors via a functor 
\begin{equation} \label{e:inv vs coinv}
\bC_{\fL^+_\nabla(H)_S}\to \bC^{\fL^+_\nabla(H)_S}.
\end{equation} 

We claim: 

\begin{cor} \label{c:inv vs coinv}
The functor \eqref{e:inv vs coinv} is an equivalence.
\end{cor}

\begin{proof}

\propref{p:1-aff rep H} implies that the functor
$$\bC\mapsto  \bC^{\fL^+_\nabla(H)_S}$$
commutes with colimits. 

\medskip

Hence, both sides in \eqref{e:inv vs coinv} commute with colimits. Any object 
in $\QCoh(\fL^+_\nabla(H))_S\mmod$ can be written as a colimit of objects of the form  
$\QCoh(\fL^+_\nabla(H))_S\otimes \bD$, where the module structure comes from the
first factor. Hence, we obtain that it is sufficient to prove that \eqref{e:inv vs coinv} 
is an equivalence for such objects. 

\medskip

However, the latter is obvious: the corresponding functor is the identity functor
$$\QCoh(S)\otimes \bD\simeq \left(\QCoh(\fL^+_\nabla(H))_S\otimes \bD\right)_{\fL^+_\nabla(H)_S}\to
\left(\QCoh(\fL^+_\nabla(H))_S\otimes \bD\right)^{\fL^+_\nabla(H)_S}\simeq \QCoh(S)\otimes \bD.$$

\end{proof} 

\sssec{} 

We now claim: 

\begin{cor} \label{c:QCoh co LS}
The functor $$\Omega_{\LS^\reg_{H,S}}:\QCoh_{\on{co}}(\LS^\reg_{H,S})\to \QCoh(\LS^\reg_{H,S})$$
is an equivalence.
\end{cor} 

\begin{proof}

The fact that $S\to \LS^\reg_{H,S}$ is an fpqc cover implies that the functor
$$\QCoh(S)\underset{\QCoh(\fL^+_\nabla(H))_S}\otimes \QCoh(S)\to \QCoh_{\on{co}}(\LS^\reg_{H,S})$$
is an equivalence (cf. \cite[Proposition 6.2.7]{Ga5}).

\medskip

Hence, it remains to show that the functor
$$\QCoh(S)\underset{\QCoh(\fL^+_\nabla(H))_S}\otimes \QCoh(S) \to \QCoh(\LS^\reg_{H,S})$$
is an equivalence.

\medskip

However, the latter functor is the functor \eqref{e:inv vs coinv} for $\bC=\QCoh(S)$.

\end{proof} 

\ssec{Definition of \texorpdfstring{$\Sph_H^{\on{spec}}$}{SphHspec}} \label{ss:defn Sph spec}

%

\sssec{}

Recall that the local spectral Hecke stack is by definition
$$\on{Hecke}_H^{\on{spec,loc}}:=\LS^\reg_H\underset{\LS^\mer_H}\times \LS^\reg_H.$$

\medskip

Our approach to the definition of $\IndCoh^*(\on{Hecke}_H^{\on{spec,loc}})$ is based on the following observation:

\begin{lem} \label{l:Hecke spec via loops}
The factorization prestack $\on{Hecke}_H^{\on{spec,loc}}$ identifies canonically with the
double quotient
$$\fL^+_\nabla(H)\backslash \fL_\nabla(H)/\fL^+_\nabla(H).$$
\end{lem}

\begin{proof}

By definition, the fiber product $\LS^\reg_H\underset{\LS^\mer_H}\times \LS^\reg_H$ identifies with
$$\fL^+_\nabla(H)\backslash \on{Stab}_{\fL_\nabla(\on{Jets}(H))}(0)/\fL^+_\nabla(H),$$
where:
\begin{itemize}

\item $0\in \fL_\nabla(\on{Conn}(\fh))$ is the trivial connection;

\item $\on{Stab}_{\fL_\nabla(\on{Jets}(H))}(0)$ denotes the stabilizer of $0$ with respect to the gauge action of 
$\fL_\nabla(\on{Jets}(H))\simeq \fL(H)$.

\end{itemize}

However, 
$$\on{Stab}_{\fL_\nabla(\on{Jets}(H))}(0)=\fL_\nabla(\on{Stab}_{\on{Jets}(H)}(0)),$$
while 
$$\on{Stab}_{\on{Jets}(H)}(0)\simeq H,$$
as a group D-scheme.

\end{proof} 

\sssec{}

We also note: 

\begin{lem} \label{l:const Gr is laft}
For $S\in \affSch_{\on{aft}}$, the quotient $(\fL_\nabla(H)/\fL^+_\nabla(H))_S$ is locally almost of finite type.
\end{lem} 

\begin{proof}

First, we note that the unit section
$$S\to (\fL_\nabla(H)/\fL^+_\nabla(H))_S$$
is an isomorphism at the classical level. (Indeed, for any affine $Y$, the map $\fL^+_\nabla(Y)\to \fL_\nabla(Y)$
is an isomorphism at the classical level.)

\medskip

Hence, by \cite[Chapter 1, Theorem 9.1.2]{GaRo4}, it suffices to show that the cotangent space to $(\fL_\nabla(H)/\fL^+_\nabla(H))_S$
at the unit section is laft (see \cite[Chapter 1, Sect. 3.4.1]{GaRo4} for what this means). However, this cotangent space is the dual of
$$(\fL_\nabla(\fh)/\fL^+_\nabla(\fh))_S,$$
which makes the assertion manifest.

\end{proof} 

\sssec{}

In what follows we will define the (monoidal) factorization category $\IndCoh^*(\fL_\nabla(H))$, equipped with an action of
the monoidal category $\QCoh(\fL^+_\nabla(H))$ on the two sides. We will then set
\begin{equation} \label{e:IndCoh* Hecke as inv}
\IndCoh^*(\on{Hecke}_H^{\on{spec,loc}}):=\left(\IndCoh^*(\fL_\nabla(H))\right)_{\fL^+_\nabla(H)\times \fL^+_\nabla(H)}.
\end{equation}

The caveat here is that $\fL^+_\nabla(H)$ is \emph{not} placid (and hence, $\fL_\nabla(H)$ is not ind-placid).
Yet, we will show that the construction of $\IndCoh^*(-)$ in \ref{sss:IndCoh* over Ran}-\ref{sss:IndCoh * Ran}
is applicable in this particular case.

\sssec{}

Let $S_\alpha$ be as in \secref{sss:S alpha}. Consider the relative affine scheme $\fL_\nabla^+(H)_{S_\alpha}$ and the
relative ind-affine ind-scheme $\fL_\nabla(H)_{S_\alpha}$. 

\medskip

First, we claim: 

\begin{lem} \label{l:Psi arcs H}
The functor 
$$\Psi_{\fL^+_\nabla(H)_{S_\alpha}}:\IndCoh^*(\fL^+_\nabla(H))_{S_\alpha}\to \QCoh(\fL^+_\nabla(H))_{S_\alpha}$$
is an equivalence.
\end{lem}

\begin{proof}

The assertion holds for any smooth target scheme $Y$. Indeed, one shows that for $S_\alpha=X^I$, the relative
affine scheme $\fL^+_\nabla(Y)_{X^I}$ is isomorphic to the limit of a sequence of affine blow-ups with smooth
centers, starting with $Y^I\to X^I$, see \secref{ss:arcs sm}. 

\medskip

In particular, $\fL^+_\nabla(Y)_{X^I}$ is isomorphic to a filtered limit of relative affine schemes $Y_n\to S_\alpha$
that are \emph{smooth}. 

\medskip

We have:
$$\IndCoh^*(\fL^+_\nabla(Y))_{S_\alpha}\simeq \underset{n}{\on{lim}}\, \IndCoh(Y_n)$$
(with respect to push-forwards) and
$$\QCoh(\fL^+_\nabla(Y))_{S_\alpha}\simeq \underset{n}{\on{lim}}\, \QCoh(Y_n),$$
and the functor $\Psi_{\fL^+_\nabla(H)_{S_\alpha}}$ corresponds to the compatible family of functors
$$\Psi_{Y_n}:\IndCoh(Y_n)\to \QCoh(Y_n),$$
all of which are equivalences, since $Y_n$ are smooth.

\end{proof}  

\sssec{}

According to \secref{sss:IndCoh non-placid *}, since $\fL_\nabla(H)_{S_\alpha}$ is an ind-affine ind-scheme, we have a well-defined category 
$\IndCoh^*(\fL_\nabla(H))_{S_\alpha}$. 

\medskip

The action of $\fL_\nabla^+(H)_{S_\alpha}\times \fL_\nabla^+(H)_{S_\alpha}$ on $\fL_\nabla(H)_{S_\alpha}$ defines on
$\IndCoh^*(\fL_\nabla(H))_{S_\alpha}$ a structure of bimodule with respect to $\IndCoh^*(\fL^+_\nabla(H))_{S_\alpha}$.

\medskip

Hence, thanks to \lemref{l:Psi arcs H}, we can think of $\IndCoh^*(\fL_\nabla(H))_{S_\alpha}$ as a bimodule with respect to
$\QCoh(\fL^+_\nabla(H))_{S_\alpha}$.

\sssec{} 

Consider $\IndCoh^*(\fL_\nabla(H))_{S_\alpha}$ as a module over $\QCoh(\fL^+_\nabla(H))_{S_\alpha}$ with respect to
the action on the right.

\medskip

Direct image with respect to the projection
\begin{equation} \label{e:hor loops to hor Gr}
\fL_\nabla(H)_{S_\alpha}\to (\fL_\nabla(H)/\fL_\nabla^+(H))_{S_\alpha}
\end{equation}
gives rise to a functor
\begin{multline} \label{e:IndCoh on Gr const}
\left(\IndCoh^*(\fL_\nabla(H))_{S_\alpha}\right)_{\fL^+_\nabla(H)_{S_\alpha}}\to \IndCoh^*((\fL_\nabla(H)/\fL_\nabla^+(H))_{S_\alpha})\simeq \\
\overset{\text{\lemref{l:const Gr is laft}}}\simeq \IndCoh((\fL_\nabla(H)/\fL_\nabla^+(H))_{S_\alpha}).
\end{multline}

\medskip

We claim:

\begin{lem} \label{l:IndCoh on Gr const}
The functor \eqref{e:IndCoh on Gr const} is an equivalence.
\end{lem}

\begin{proof}

Since $S_\alpha\to (\fL_\nabla(H)/\fL^+_\nabla(H))_{S_\alpha}$ is an isomorphism at the reduced level,
Zariski-locally on $S_\alpha$, the map \eqref{e:hor loops to hor Gr} splits as a product: indeed,
the restriction of the \'etale $\fL^+_\nabla(H)$-torsor \eqref{e:hor loops to hor Gr} 
to $S_\alpha$ is trivial, and hence over any open affine of $S^\circ_\alpha\subset S_\alpha$, the 
$\fL^+_\nabla(H)$-torsor \eqref{e:hor loops to hor Gr} itself is trivial. 

\medskip

Hence, by Zariski descent, it suffices
to show that for $Z_\alpha\to S_\alpha$, where $Z_\alpha$ is an ind-affine ind-scheme locally almost of finite type,
the functor
$$\left(\IndCoh^*(Z_\alpha \underset{S_\alpha}\times \fL^+_\nabla(H)_{S_\alpha})\right)_{\fL^+_\nabla(H)_{S_\alpha}}\to \IndCoh(Z_\alpha)$$
is an equivalence. 

\medskip

To prove this, it suffices to show that the functor 
\begin{equation} \label{e:IndCoh on Gr cons}
\IndCoh(Z_\alpha)\underset{\QCoh(S_\alpha)}\otimes \QCoh(\fL^+_\nabla(H))_{S_\alpha}\simeq
\IndCoh^*(Z_\alpha \underset{S_\alpha}\times \fL^+_\nabla(H)_{S_\alpha})
\end{equation}
is an equivalence. 

\medskip

To prove \eqref{e:IndCoh on Gr cons}, we can assume that $Z_\alpha$ is an affine scheme. 
Writing $\fL^+_\nabla(H)_{S_\alpha}$ as a limit of relative affine schemes $Y_n$ \emph{smooth} over $S_\alpha$ 
as in the proof of \lemref{l:Psi arcs H}, it suffices to show that each of the functors
\begin{equation} \label{e:IndCoh on Gr cons n}
\IndCoh(Z_\alpha)\underset{\QCoh(S_\alpha)}\otimes \QCoh(Y_n)\simeq
\IndCoh^*(Z_\alpha \underset{S_\alpha}\times Y_n)
\end{equation}
is an equivalence. 

\medskip

However, this follows from \lemref{l:IndCoh^! tensor up}. 

\end{proof}

\sssec{}

We are finally ready to define $\IndCoh^*(\fL_\nabla(H))$ as a factorization category. Proceeding as in 
\ref{sss:IndCoh* over Ran}-\ref{sss:IndCoh * Ran}, we need to show that Lemmas \ref{l:IndCoh^* transition} 
and \ref{l:form compl diag base change *} hold for $\IndCoh^*(\fL_\nabla(H))$. 

\medskip

We will prove \lemref{l:IndCoh^* transition}; \lemref{l:form compl diag base change *} is proved similarly.

\medskip

We need to show that the functor
\begin{equation} \label{e:IndCoh^* transition H}
\QCoh(S_\alpha)\underset{\QCoh(S_\beta)}\otimes \IndCoh^*(\fL_\nabla(H))_{S_\beta}\to \IndCoh^*(\fL_\nabla(H))_{S_\alpha}.
\end{equation} 
is an equivalence.

\medskip

We consider both sides as modules over
$$\QCoh(S_\alpha)\underset{\QCoh(S_\beta)}\otimes \QCoh(\fL^+_\nabla(H))_{S_\beta}\simeq \QCoh(\fL^+_\nabla(H))_{S_\alpha}.$$

By \propref{p:1-aff rep H}, it suffices to show that \eqref{e:IndCoh^* transition H} becomes an equivalence
after taking $\fL^+_\nabla(H)_{S_\alpha}$-invariants, or, equivalently, thanks to \lemref{c:inv vs coinv},
$\fL^+_\nabla(H)_{S_\alpha}$-coinvariants. 

\medskip

However, by \lemref{l:IndCoh on Gr const}, when we take $\fL^+_\nabla(H)_{S_\alpha}$-coinvariants in \eqref{e:IndCoh^* transition H},
the resulting functor identifies with 
\begin{equation} \label{e:IndCoh^* transition H 1}
\QCoh(S_\alpha)\underset{\QCoh(S_\beta)}\otimes \IndCoh^*(\fL_\nabla(H)/\fL^+_\nabla(H))_{S_\beta}\to 
\IndCoh^*(\fL_\nabla(H)/\fL^+_\nabla(H))_{S_\alpha}.
\end{equation} 

Now, \eqref{e:IndCoh^* transition H 1} is an equivalence by \lemref{l:IndCoh^* transition}, since
$\fL_\nabla(H)/\fL^+_\nabla(H)$ is locally almost of finite type (by \lemref{l:const Gr is laft}) and in particular is
placid. 

\sssec{} \label{sss:defn IndCoh* Hecke spec}

By construction, $\IndCoh^*(\fL_\nabla(H))$ is equipped, as a factorization category, with an action of
$\QCoh(\fL^+_\nabla(H))\otimes \QCoh(\fL^+_\nabla(H))$. 

\medskip

We define $\IndCoh^*(\on{Hecke}_H^{\on{spec,loc}})$
by formula \eqref{e:IndCoh* Hecke as inv}.

\medskip

By \propref{p:1-aff rep H}, we have
$$\IndCoh^*(\on{Hecke}_H^{\on{spec,loc}})\underset{\QCoh(\LS^\reg_H)\otimes \QCoh(\LS^\reg_H)}\otimes \Vect
\simeq \IndCoh^*(\fL_\nabla(H))$$
and 
$$\IndCoh^*(\on{Hecke}_H^{\on{spec,loc}})\underset{\QCoh(\LS^\reg_H)}\otimes \Vect
\simeq \IndCoh(\fL_\nabla(H)/\fL^+_\nabla(H)).$$

\sssec{}

The pair of adjoint functors
\begin{equation} \label{e:gen Hecke spec prel}
\iota^\IndCoh_*:\IndCoh^*(\fL^+_\nabla(H))\rightleftarrows \IndCoh^*(\fL_\nabla(H)):\iota^!
\end{equation} 
gives rise via 
\begin{multline*} 
\QCoh(\LS^\reg_H)\overset{\text{\corref{c:QCoh co LS}}}\simeq 
\QCoh_{\on{co}}(\LS^\reg_H)\simeq \Vect_{\fL^+_\nabla(H)}\simeq \\
\simeq \QCoh(\fL^+_\nabla(H))_{\fL^+_\nabla(H)\times \fL^+_\nabla(H)} \overset{\Psi_{\fL^+_\nabla(H)}}
\simeq \IndCoh^*(\fL^+_\nabla(H))_{\fL^+_\nabla(H)\times \fL^+_\nabla(H)})
\end{multline*} 
to an adjoint pair
\begin{equation} \label{e:gen Hecke spec}
\iota^\IndCoh_*:\QCoh(\LS^\reg_H) \rightleftarrows \IndCoh^*(\on{Hecke}_H^{\on{spec,loc}}):\iota^!.
\end{equation} 

Since the essential image of the left adjoint in \eqref{e:gen Hecke spec prel} generates the essential 
image, the same is true for \eqref{e:gen Hecke spec}.

\medskip

In particular, images of compact objects in $\QCoh(\LS^\reg_H)$ under $\iota^\IndCoh_*$ provide
compact generators of $\IndCoh^*(\on{Hecke}_H^{\on{spec,loc}})$. 

\ssec{Unital structure} \label{ss:unital Hecke spec}

\sssec{}

We claim that the factorization categories we defined above, namely,
$$\IndCoh^*(\fL_\nabla(H)) \text{ and } \IndCoh^*(\on{Hecke}_H^{\on{spec,loc}})$$
carry naturally defined unital structures.

\medskip

Let us carry out the construction for $\IndCoh^*(\fL_\nabla(H))$; the case of $\IndCoh^*(\on{Hecke}_H^{\on{spec,loc}})$
will follow by taking $\fL^+_\nabla(H)\times \fL^+_\nabla(H)$-coinvariants.

\sssec{}

By \secref{sss:unital loops}, the factorization space $\fL_\nabla(H)$ carries a natural unital-in-correspondences
structure. (Note, however, that we cannot deduce from there the unital structure on $\IndCoh^*(\fL_\nabla(H))$
by applying \secref{sss:corr *} directly because we are not in an ind-placid situation.) 

\medskip

For an injection of finite sets $I_1\subseteq I_2$ consider the corresponding diagram
\begin{equation}  \label{e:unital for H nabla 0}
\fL_\nabla(H)_{X^{I_1}}
\overset{\on{pr}^H_{\on{small}}}\longleftarrow 
\fL^{\mer\rightsquigarrow \reg}_\nabla(H)_{I_1\subseteq I_2} \overset{\on{pr}^H_{\on{big}}}\longrightarrow 
\fL_\nabla(H)_{X^{I_2}}.
\end{equation}

\sssec{} \label{sss:left adj loops H 1}

We claim that the functor
$$(\on{pr}^H_{\on{small}})^\IndCoh_*:\IndCoh^*(\fL^{\mer\rightsquigarrow \reg}_\nabla(H))_{I_1\subseteq I_2} \to
\IndCoh^*(\fL_\nabla(H))_{X^{I_1}}$$
admits a left adjoint, to be denoted $(\on{pr}^H_{\on{small}})^{*,\IndCoh}$. 

\medskip

Factor the map $\on{pr}^H_{\on{small}}$ as
$$\fL^{\mer\rightsquigarrow \reg}_\nabla(H)_{I_1\subseteq I_2}  \overset{'\!\on{pr}^H_{\on{small}}}\longrightarrow 
\fL_\nabla(H)_{X^{I_1}}\underset{X^{I_1}}\times X^{I_2}\to
\fL_\nabla(H)_{X^{I_1}},$$
and it is sufficient to prove the existence of the left adjoint for the first arrow, i.e., that the functor
\begin{equation} \label{e:unital for H nabla 1}
({}'\!\on{pr}^H_{\on{small}})^\IndCoh_*:
\IndCoh^*(\fL^{\mer\rightsquigarrow \reg}_\nabla(H)_{I_1\subseteq I_2})\to 
\IndCoh^*(\fL_\nabla(H)_{X^{I_1}}\underset{X^{I_1}}\times X^{I_2})
\end{equation}
admits a left adjoint.

\sssec{} \label{sss:left adj loops H 2}

We consider the two sides of \eqref{e:unital for H nabla 1} as acted on by
$$\fL^+_\nabla(H)_{X^{I_2}} \text{ and } \fL^+_\nabla(H)_{X^{I_1}}\underset{X^{I_1}}\times X^{I_2},$$
respectively. These actions are compatible via the the map
$$\fL^+_\nabla(H)_{X^{I_2}} \to \fL^+_\nabla(H)_{X^{I_1}}\underset{X^{I_1}}\times X^{I_2},$$
corresponding to the counital structure on $\fL^+_\nabla(H)$.

\medskip 

By \propref{p:1-aff rep H}, it suffices to show that the functor
\begin{equation} \label{e:unital for H nabla 2}
\left(\IndCoh^*(\fL^{\mer\rightsquigarrow \reg}_\nabla(H)_{I_1\subseteq I_2})\right)_{\fL^+_\nabla(H)_{X^{I_2}}}\to
\left(\IndCoh^*(\fL_\nabla(H)_{X^{I_1}}\underset{X^{I_1}}\times X^{I_2})\right)_{\fL^+_\nabla(H)_{X^{I_1}}\underset{X^{I_1}}\times X^{I_2}},
\end{equation}
induced by \eqref{e:unital for H nabla 1}, admits a left adjoint.

\medskip

However, by \lemref{l:IndCoh on Gr const}, the latter functor is the identity endofunctor of
$$\IndCoh\left((\fL_\nabla(H)/\fL^+_\nabla(H))_{X^{I_1}}\underset{X^{I_1}}\times  X^{I_2}\right).$$

\sssec{}

We define the functor 
$$\IndCoh^*(\fL_\nabla(H))_{X^{I_1}}\to \IndCoh^*(\fL_\nabla(H))_{X^{I_2}},$$
to be denoted $\on{ins.unit}_{I_1\subseteq I_2}$, to be 
$$(\on{pr}^H_{\on{big}})^\IndCoh_*\circ (\on{pr}^H_{\on{small}})^{*,\IndCoh}.$$

\sssec{}

In order to promote this to a unital structure on $\IndCoh^*(\fL_\nabla(H))$, we need to 
construct isomorphisms
$$\on{ins.unit}_{I_2\subseteq I_2}\circ \on{ins.unit}_{I_1\subseteq I_2}\simeq \on{ins.unit}_{I_1\subseteq I_3}$$
for $I_1\subseteq I_2\subseteq I_3$. 

\medskip

Denote the maps in \eqref{e:unital for H nabla 0} by 
$$\on{pr}^H_{\on{small},I_1\subseteq I_2} \text{ and } \on{pr}^H_{\on{big},I_1\subseteq I_2}$$
to indicate the dependence on the finite sets involved.

\medskip

Thus, we need to construct an isomorphism
\begin{multline} \label{e:unital for H nabla 3}
(\on{pr}^H_{\on{big},I_2\subseteq I_3})^\IndCoh_*\circ 
(\on{pr}^H_{\on{small},I_2\subseteq I_3})^{*,\IndCoh}\circ (\on{pr}^H_{\on{big},I_1\subseteq I_2})^\IndCoh_*
\circ (\on{pr}^H_{\on{small},I_1\subseteq I_1})^{*,\IndCoh} \overset{\sim}\to \\
\overset{\sim}\to (\on{pr}^H_{\on{big},I_1\subseteq I_3})^\IndCoh_*\circ 
(\on{pr}^H_{\on{small},I_1\subseteq I_3})^{*,\IndCoh}.
\end{multline}

\sssec{}

Note that we have a commutative diagram,
$$
\CD
\fL^{\mer\rightsquigarrow \reg}_\nabla(H)_{I_1\subseteq I_3}  @>{'\on{pr}^H_{\on{big},I_1\subseteq I_2}}>> 
\fL^{\mer\rightsquigarrow \reg}_\nabla(H)_{I_2\subseteq I_3} @>{\on{pr}^H_{\on{big},I_2\subseteq I_3}}>>
\fL_\nabla(H)_{X^{I_3}} \\
@V{'\on{pr}^H_{\on{small},I_2\subseteq I_3}}VV @VV{\on{pr}^H_{\on{small},I_2\subseteq I_3}}V \\
\fL^{\mer\rightsquigarrow \reg}_\nabla(H)_{I_1\subseteq I_2} @>{\on{pr}^H_{\on{big},I_1\subseteq I_2}}>> \fL_\nabla(H)_{X^{I_2}} \\
@V{\on{pr}^H_{\on{small},I_1\subseteq I_2}}VV \\
\fL_\nabla(H)_{X^{I_1}},
\endCD
$$
in which the inner square is Cartesian.

\medskip

We rewrite the right-hand side in \eqref{e:unital for H nabla 3} as
$$(\on{pr}^H_{\on{big},I_2\subseteq I_3})^\IndCoh_*\circ 
({}'\!\on{pr}^H_{\on{big},I_1\subseteq I_2})^\IndCoh_*\circ 
({}'\!\on{pr}^H_{\on{small},I_2\subseteq I_3})^{*,\IndCoh}\circ (\on{pr}^H_{\on{small},I_1\subseteq I_2})^{*,\IndCoh}.$$

\medskip

The isomorphism
$$(\on{pr}^H_{\on{big},I_1\subseteq I_2})^\IndCoh_* \circ ({}'\!\on{pr}^H_{\on{small},I_2\subseteq I_3})^\IndCoh_*\simeq 
(\on{pr}^H_{\on{small},I_2\subseteq I_3})^\IndCoh_*\circ 
({}'\!\on{pr}^H_{\on{big},I_1\subseteq I_2})^\IndCoh_*$$
induces a natural transformation
\begin{equation} \label{e:unital for H nabla 4}
(\on{pr}^H_{\on{small},I_2\subseteq I_3})^{*,\IndCoh}\circ (\on{pr}^H_{\on{big},I_1\subseteq I_2})^\IndCoh_*\to
({}'\!\on{pr}^H_{\on{big},I_1\subseteq I_2})^\IndCoh_*\circ 
({}'\!\on{pr}^H_{\on{small},I_2\subseteq I_3})^{*,\IndCoh}.
\end{equation} 

\sssec{}
 
We claim that \eqref{e:unital for H nabla 4} is an isomorphism. Indeed, this follows by the same argument as that proving
the existence of $(\on{pr}^H_{\on{small}})^{*,\IndCoh}$ in Sects. \ref{sss:left adj loops H 1}-\ref{sss:left adj loops H 2}.  

\sssec{}

Finally, we define the natural isomorphism in \eqref{e:unital for H nabla 3} by precomposing 
the isomorphism \eqref{e:unital for H nabla 4} with $(\on{pr}^H_{\on{small},I_1\subseteq I_2})^{*,\IndCoh}$ and
post-composing with $(\on{pr}^H_{\on{big},I_2\subseteq I_3})^\IndCoh_*$. 

\medskip

The higher compatibilities are constructed by a similar procedure. 

\sssec{}

By construction, the functor
$$\QCoh(\fL^+_\nabla(H)) \overset{\Psi_{\fL^+_\nabla(H))}}\simeq \IndCoh^*(\fL^+_\nabla(H)) \overset{\iota^\IndCoh_*}\longrightarrow
\IndCoh^*(\fL_\nabla(H))$$
is unital. 

\medskip

In particular, the object
$$\iota^\IndCoh_*(\CO_{\fL^+_\nabla(H)})\in \on{FactAlg}(X,\IndCoh^*(\fL_\nabla(H)))$$
is the factorization unit in $\IndCoh^*(\fL_\nabla(H))$.

\medskip

Similarly, the functor
$$\QCoh(\LS^\reg_H) \overset{\iota^\IndCoh_*}\longrightarrow  \IndCoh^*(\on{Hecke}_H^{\on{spec,loc}})$$
is unital, and
$$ \iota^\IndCoh_*(\CO_{\LS_H^\reg})\in \on{FactAlg}(X,\IndCoh^*(\on{Hecke}_H^{\on{spec,loc}}))$$
is the factorization unit in $\IndCoh^*(\on{Hecke}_H^{\on{spec,loc}})$.

\ssec{Duality}

\sssec{}

Let $S_\alpha$ be as in \secref{sss:S alpha}.

\medskip

Note that the same argument as in \lemref{l:Psi arcs H} shows that the functor
$$\Upsilon_{\fL^+_\nabla(H))_{S_\alpha}}:\QCoh(\fL^+_\nabla(H))_{S_\alpha}\to \IndCoh^!(\fL^+_\nabla(H))_{S_\alpha}$$
is an equivalence. 

\medskip

As a consequence, we obtain that the pairing 
$$\IndCoh^!(\fL^+_\nabla(H))_{S_\alpha}\otimes \IndCoh^*(\fL^+_\nabla(H))_{S_\alpha}\to \Vect$$
of \eqref{e:IndCoh ! * pairing} is perfect. 

\sssec{}

Similarly, an argument parallel to that in \lemref{l:IndCoh on Gr const} shows that the !-pullback functor along
$$\IndCoh(\fL_\nabla(H)/\fL^+_\nabla(H))_{S_\alpha}\to \IndCoh^!(\fL_\nabla(H))_{S_\alpha}$$
gives rise to an equivalence
$$\IndCoh(\fL_\nabla(H)/\fL^+_\nabla(H))_{S_\alpha}\simeq \left(\IndCoh^!(\fL_\nabla(H))_{S_\alpha}\right)^{\fL^+_\nabla(H)_{S_\alpha}}.$$

Combining with \propref{p:1-aff rep H}, we obtain that 
the pairing
$$\IndCoh^!(\fL_\nabla(H))_{S_\alpha}\otimes \IndCoh^*(\fL_\nabla(H))_{S_\alpha}\to \Vect$$
of \eqref{e:IndCoh ! * pairing} is perfect.  

\sssec{}

In particular, we obtain that Lemmas \ref{l:IndCoh^! transition} and \ref{l:form compl diag base change} 
hold for $\IndCoh^!(\fL_\nabla(H))$. I.e., the recipe
in Sects. \ref{sss:IndCoh ! Ran 1}-\ref{sss:IndCoh ! Ran 2} gives 
rise to a well-defined factorization category $\IndCoh^!(\fL_\nabla(H))$.

\medskip

Moreover, we obtain that \eqref{e:IndCoh ! * pairing} defines a perfect pairing between 
$$\IndCoh^!(\fL_\nabla(H)) \text{ and } \IndCoh^*(\fL_\nabla(H))$$
as factorization categories.

\sssec{}

By a similar logic as in \secref{sss:defn IndCoh* Hecke spec}, we obtain that the assignment
$$S_\alpha \rightsquigarrow \IndCoh^!(\on{Hecke}_H^{\on{spec,loc}})_{S_\alpha}$$
extends to a well-defined factorization category $\IndCoh^!(\on{Hecke}_H^{\on{spec,loc}})$.

\medskip

Moreover, we have
$$\IndCoh^!(\on{Hecke}_H^{\on{spec,loc}})\simeq \IndCoh^!(\fL_\nabla(H))^{\fL^+_\nabla(H)\times \fL^+_\nabla(H)}$$
and
$$\IndCoh^!(\on{Hecke}_H^{\on{spec,loc}})\underset{\QCoh(\LS^\reg_H)}\otimes \Vect \simeq \IndCoh^!(\fL_\nabla(H)/\fL^+_\nabla(H)).$$

We obtain that 
$$\IndCoh^*(\on{Hecke}_H^{\on{spec,loc}}) \text{ and } \IndCoh^!(\on{Hecke}_H^{\on{spec,loc}})$$
are mutually dual as factorization categories.

\sssec{}

Finally, a procedure dual to that in \secref{ss:unital Hecke spec} defines on 
$$\IndCoh^!(\fL_\nabla(H)) \text{ and } \IndCoh^*(\on{Hecke}_H^{\on{spec,loc}})$$
unital structures, and the identifications
$$\IndCoh^!(\fL_\nabla(H))^\vee\simeq \IndCoh^*(\fL_\nabla(H)) \text{ and } 
\IndCoh^!(\on{Hecke}_H^{\on{spec,loc}})^\vee \simeq \IndCoh^*(\on{Hecke}_H^{\on{spec,loc}})$$
extend to identifications of the corresponding unital factorization categories. 

\ssec{t-structures}

In this section we will discuss an alternative approach to the definition of $\IndCoh^*(\on{Hecke}_H^{\on{spec,loc}})$.

\medskip

Namely, we can start with (the more elementary) $\QCoh_{\on{co}}(\on{Hecke}_H^{\on{spec,loc}})$, and obtain from
it $\IndCoh^*(\on{Hecke}_H^{\on{spec,loc}})$ by a \emph{renormalization procedure} (i.e., ind-completion of a specified
small subcategory), see \secref{sss:recover Hecke spec from QCoh co} below. 

\sssec{}

As in \secref{sss:IndCoh * Ran}, the factorization category
$$\IndCoh^*(\fL_\nabla(H))$$
carries a naturally defined t-structure. 

\medskip

By construction, the functor
$$\Gamma^\IndCoh:\IndCoh^*(\fL_\nabla(H))\to \Vect$$
is t-exact, and conservative when restricted to the eventually coconnective subcategory.

\medskip

Moreover, the functor
$$\iota^\IndCoh_*:\QCoh(\fL^+_\nabla(H))\to \IndCoh^*(\fL_\nabla(H))$$
is t-exact.

\sssec{} \label{sss:t-str on Hecke spec}

We now define a t-structure on $\IndCoh^*(\on{Hecke}_H^{\on{spec,loc}})$. Namely, by \propref{p:1-aff rep H}
the projection
$$\IndCoh^*(\fL_\nabla(H))\to \IndCoh^*(\on{Hecke}_H^{\on{spec,loc}})$$
admits a left adjoint, which is comonadic.

\medskip

Moreover, the resulting comonad on $\IndCoh^*(\fL_\nabla(H))$ is t-exact. This implies that the category 
$\IndCoh^*(\on{Hecke}_H^{\on{spec,loc}})$
acquires a unique t-structure for which both functors
$$\IndCoh^*(\on{Hecke}_H^{\on{spec,loc}})\rightleftarrows \IndCoh^*(\fL_\nabla(H))$$
are t-exact. 

\sssec{} 

Consider the untal factorization categories
$$\QCoh_{\on{co}}(\fL_\nabla(H)),\,\, \QCoh_{\on{co}}(\fL_\nabla(H)/\fL^+_\nabla(H)) \text{ and } 
\QCoh_{\on{co}}(\on{Hecke}_H^{\on{spec,loc}}).$$

\medskip

The category $\QCoh_{\on{co}}(\fL_\nabla(H))$ carries a natural action of $\QCoh(\fL^+_\nabla(H))$, and it follows 
formally that the functors
\begin{equation} \label{e:coinv co LH nabla}
\left(\QCoh_{\on{co}}(\fL_\nabla(H))\right)_{\fL^+_\nabla(H)\times \fL^+_\nabla(H)}\to 
\QCoh_{\on{co}}(\on{Hecke}_H^{\on{spec,loc}})
\end{equation}  
and 
\begin{equation} \label{e:coinv co LH nabla bis}
\left(\QCoh_{\on{co}}(\fL_\nabla(H))\right)_{\fL^+_\nabla(H)}\to 
\QCoh_{\on{co}}(\fL_\nabla(H)/\fL^+_\nabla(H)),
\end{equation}  
induced by the direct image functors 
$$\QCoh_{\on{co}}(\fL_\nabla(H))\to \QCoh_{\on{co}}(\on{Hecke}_H^{\on{spec,loc}}) \text{ and }
\QCoh_{\on{co}}(\fL_\nabla(H))\to \QCoh_{\on{co}}(\fL_\nabla(H)/\fL^+_\nabla(H)),$$
respectively, 
are equivalences.

\medskip

In addition, we have the unital factorization functors
$$\QCoh(\fL^+_\nabla(H)) \overset{\iota_*}\to \QCoh_{\on{co}}(\fL_\nabla(H)) \text{ and }
\QCoh(\LS_H^\reg)\simeq \QCoh_{\on{co}}(\LS_H^\reg) \overset{\iota_*}\to \QCoh_{\on{co}}(\on{Hecke}_H^{\on{spec,loc}}).$$ 

\sssec{} \label{sss:IndCoh and QCoh co loops L H}

Let $S_\alpha$ be as in \secref{sss:S alpha}. Recall (see \lemref{l:Psi IndCoh* IndSch}) that the functor
$$\Psi_{\fL_\nabla(H)_{S_\alpha}}:\IndCoh^*(\fL_\nabla(H))_{S_\alpha}\to  \QCoh_{\on{co}}(\fL_\nabla(H))_{S_\alpha}$$
is t-exact, and induces an equivalence between the eventually coconnective subcategories on both sides.

\medskip

It follows from the definition of $\IndCoh^*(\fL_\nabla(H))$ as a factorization category that
the functors $\Psi_{\fL_\nabla(H)_{S_\alpha}}$ combine to give rise to a factorization functor
$$\Psi_{\fL_\nabla(H)}:\IndCoh^*(\fL_\nabla(H))\to  \QCoh_{\on{co}}(\fL_\nabla(H)).$$

\medskip

Moreover, the functor $\Psi_{\fL_\nabla(H)}$ is t-exact and induces an equivalences between the 
eventually coconnective subcategories on both sides.

\medskip

Furthermore, the functor $\Psi_{\fL_\nabla(H)}$ has a naturally defined unital structure. 

\sssec{} \label{sss:recover from QCoh co loops H}

Note that the contents of \secref{sss:IndCoh and QCoh co loops L H} allows us to recover $\IndCoh^*(\fL_\nabla(H))$,
as a unital factorization category, from $\QCoh_{\on{co}}(\fL_\nabla(H))$ with its t-structure. 

\medskip

Namely, for $S\to \Ran$, the category $\IndCoh^*(\fL_\nabla(H))_S$ identifies with the ind-completion of the
full subcategory of $\QCoh_{\on{co}}(\fL_\nabla(H))_S$, generated by \emph{finite} colimits by the essential image
of $\QCoh(\fL^+_\nabla(H))_S$ along $\iota_*$. 

\sssec{} \label{sss:recover Hecke spec from QCoh co}

It follows formally from \secref{sss:t-str on Hecke spec} that we have a naturally defined t-exact unital factorization
functor
$$\Psi_{\on{Hecke}_H^{\on{spec,loc}}}: \IndCoh^*(\on{Hecke}_H^{\on{spec,loc}})\to \QCoh_{\on{co}}(\on{Hecke}_H^{\on{spec,loc}}),$$
which induces an equivalences between the 
eventually coconnective subcategories on both sides.

\medskip

Furthermore, as in \secref{sss:recover from QCoh co loops H}, we can recover $\IndCoh^*(\on{Hecke}_H^{\on{spec,loc}})$
(as a unital factorization category) from
$\QCoh_{\on{co}}(\on{Hecke}_H^{\on{spec,loc}})$ with its t-structure.

\medskip 

Namely, for $S\to \Ran$, the category $\IndCoh^*(\on{Hecke}_H^{\on{spec,loc}})_S$ identifies with the ind-completion of the
full subcategory of $(\QCoh_{\on{co}}(\on{Hecke}_H^{\on{spec,loc}})_S)^{>-\infty}$, generated by \emph{finite} colimits by the essential image
of $\QCoh(\LS^\reg_H)_S$ along $\iota_*$. 

\ssec{The monoidal structure}

\sssec{}

Let $S_\alpha$ be as in \secref{sss:S alpha}. The group-scheme structure on $\fL_\nabla(H)_{S_\alpha}$ 
induces on the category $\IndCoh^*(\fL_\nabla(H))_{S_\alpha}$ a structure of monoidal category (under convolution).

\medskip

By the construction of $\IndCoh^*(\fL_\nabla(H))$ as a factorization category, we obtain that 
$\IndCoh^*(\fL_\nabla(H))$ acquires a structure of monoidal factorization category.  

\medskip

This structure is compatible with the unital structure on $\IndCoh^*(\fL_\nabla(H))$ in the sense that
the monoidal operation has a natural lax unital structure. I.e., $\IndCoh^*(\fL_\nabla(H))$ is an 
associative algebra object in the symmetric monoidal category of unital factorization categories
with lax unital functors as morphisms.

\sssec{}

Note that for $S_\alpha$ as above, the monoidal operation on $\IndCoh^*(\fL_\nabla(H))_{S_\alpha}$, viewed as a functor
$$\IndCoh^*(\fL_\nabla(H))_{S_\alpha}\underset{\QCoh(S_\alpha)}\otimes \IndCoh^*(\fL_\nabla(H))_{S_\alpha}\to
\IndCoh^*(\fL_\nabla(H))_{S_\alpha}$$
is t-exact.

\medskip

Hence, the monoidal operation on $\IndCoh^*(\fL_\nabla(H))$ is t-exact.

\sssec{}

We will now descend the above monoidal structure to one on $\IndCoh^*(\on{Hecke}_H^{\on{spec,loc}})$.

\medskip

Consider the correspondence
\begin{multline}  \label{e:corr Hecke 0}
(\LS_H^\reg\underset{\LS_H^\mer}\times \LS_H^\reg)\times  (\LS_H^\reg\underset{\LS_H^\mer}\times \LS_H^\reg) 
\overset{\Delta_{\LS_H^\reg}}\leftarrow \\
\leftarrow \LS_H^\reg\underset{\LS_H^\mer}\times \LS_H^\reg \underset{\LS_H^\mer}\times \LS_H^\reg\overset{\on{mult}}\to 
\LS_H^\reg \underset{\LS_H^\mer}\times \LS_H^\reg.
\end{multline} 

Note that by \lemref{l:Hecke spec via loops} we can think of this diagram also as
\begin{multline}  \label{e:corr Hecke}
(\fL^+_\nabla(H)\backslash \fL_\nabla(H)/\fL^+_\nabla(H)) \times 
(\fL^+_\nabla(H)\backslash \fL_\nabla(H)/\fL^+_\nabla(H)) \leftarrow \\
\leftarrow \fL^+_\nabla(H)\backslash \left(\underset{\fL^+_\nabla(H)}{\underline{\fL_\nabla(H)\times \fL_\nabla(H)}}\right)/\fL^+_\nabla(H))
\to \fL^+_\nabla(H)\backslash \fL_\nabla(H)/\fL^+_\nabla(H),
\end{multline} 
where:

\begin{itemize}

\item $\underset{\fL^+_\nabla(H)}{\underline{\fL_\nabla(H)\times \fL_\nabla(H)}}$ denotes the quotient with respect to
the diagonal action by right multiplication along the left factor and the left multiplication along the right factor;

\medskip

\item The arrow $\to$ is induced by the product map 
$$\underset{\fL^+_\nabla(H)}{\underline{\fL_\nabla(H)\times \fL_\nabla(H)}}\to \fL_\nabla(H).$$

\end{itemize} 

\sssec{} \label{sss:mon op Hecke spec}

As in \secref{ss:defn Sph spec}, we obtain a well-defined factorization category
$$\IndCoh^*(\LS_H^\reg\underset{\LS_H^\mer}\times \LS_H^\reg \underset{\LS_H^\mer}\times \LS_H^\reg),$$
equipped with factorization functors
\begin{multline}  \label{e:corr Hecke 1}
\IndCoh^*(\LS_H^\reg\underset{\LS_H^\mer}\times \LS_H^\reg)\otimes 
\IndCoh^*(\LS_H^\reg\underset{\LS_H^\mer}\times \LS_H^\reg) \overset{(\Delta_{\LS_H^\reg})^\IndCoh_*}\longleftarrow \\
\leftarrow \IndCoh^*(\LS_H^\reg\underset{\LS_H^\mer}\times \LS_H^\reg \underset{\LS_H^\mer}\times \LS_H^\reg)
\overset{(\on{mult})^\IndCoh_*}\longrightarrow \IndCoh^*(\LS_H^\reg\underset{\LS_H^\mer}\times \LS_H^\reg).
\end{multline} 

It is easy to see that the functor $(\Delta_{\LS_H^\reg})^\IndCoh_*$ admits a left adjoint, to be denoted 
$(\Delta_{\LS_H^\reg})^{*,\IndCoh}$.  We define the monoidal structure on $\IndCoh^*(\LS_H^\reg\underset{\LS_H^\mer}\times \LS_H^\reg)$
with the binary operation given by
$$(\on{mult})^\IndCoh_*\circ (\Delta_{\LS_H^\reg})^{*,\IndCoh}.$$

\medskip

One defines similarly $n$-fold compositions, and they form a compatible system thanks to the fact that
the functors $(\Delta_{\LS_H^\reg})^{*,\IndCoh}$ satisfy base change against $\IndCoh$-pushforwards. 

\medskip

This defines on 
$$\IndCoh^*(\LS_H^\reg\underset{\LS_H^\mer}\times \LS_H^\reg)=:\IndCoh^*(\on{Hecke}_H^{\on{spec,loc}})$$
a structure of monoidal factorization category. 

\medskip

As in the case of $\IndCoh^*(\fL_\nabla(H))$, it is easy to see that the monoidal operation on the category
$\IndCoh^*(\on{Hecke}_H^{\on{spec,loc}})$ is t-exact. 

\sssec{} \label{sss:Sph acts on Rep}

A similar procedure gives rise to an action of $\IndCoh^*(\on{Hecke}_H^{\on{spec,loc}})$ on $\QCoh(\LS^\reg_H)$. 

\sssec{}

By construction, the monoidal structure on $\IndCoh^*(\on{Hecke}_H^{\on{spec,loc}})$
is compatible with the unital structure, in the sense that $\IndCoh^*(\on{Hecke}_H^{\on{spec,loc}})$ is
an associative algebra object in the symmetric monoidal category of unital factorization categories
with lax unital functors as morphisms.

\medskip

However, we claim that, unlike $\IndCoh^*(\fL_\nabla(H))$, in the case of $\IndCoh^*(\on{Hecke}_H^{\on{spec,loc}})$
more is true: namely, the monoidal operation is strictly unital.

\medskip

Indeed, this follows from \lemref{l:unit determines strict}, since the factorization unit for $\IndCoh^*(\on{Hecke}_H^{\on{spec,loc}})$,
namely, $\iota^\IndCoh_*(\CO_{\LS^\reg_H})$ is also the monoidal unit.

\medskip

The same observation applies to the action of $\IndCoh^*(\on{Hecke}_H^{\on{spec,loc}})$ on $\QCoh(\LS^\reg_H)$. . 

\sssec{}

A similar procedure defines the monoidal factorization categories
\begin{equation} \label{e:monoidal Hecke spec QCoh co}
\QCoh_{\on{co}}(\fL_\nabla(H)) \text{ and } \QCoh_{\on{co}}(\on{Hecke}_H^{\on{spec,loc}}),
\end{equation} 
with t-exact monoidal operations.

\medskip

We note that the procedure \secref{sss:recover Hecke spec from QCoh co} allows us to recover 
\begin{equation} \label{e:monoidal Hecke spec IndCoh}
\IndCoh^*(\fL_\nabla(H)) \text{ and } \IndCoh^*(\on{Hecke}_H^{\on{spec,loc}})
\end{equation} 
as \emph{monoidal factorization categories} from those in \eqref{e:monoidal Hecke spec QCoh co}.

\medskip

Namely, these monoidal structures are uniquely determined by the requirement that
that the functors
$$\Psi_{\fL_\nabla(H)}:\IndCoh^*(\fL_\nabla(H))\to \QCoh_{\on{co}}(\fL_\nabla(H))$$
and 
$$\Psi_{\on{Hecke}_H^{\on{spec,loc}}}: \IndCoh^*(\on{Hecke}_H^{\on{spec,loc}})\to \QCoh_{\on{co}}(\on{Hecke}_H^{\on{spec,loc}})$$
are monoidal. 

\sssec{}

Note that the functor $(\on{mult})^\IndCoh_*$, involved in the definition of the monoidal structure on 
$\IndCoh^*(\on{Hecke}_H^{\on{spec,loc}})$ admits a (continuous) right adjoint, to be denoted $(\on{mult})^!$. 

\medskip

To prove this, it suffices to show that the corresponding functor
$$\IndCoh^*\left(\underset{\fL^+_\nabla(H)}{\underline{\fL_\nabla(H)\times \fL_\nabla(H)}}\right) \overset{(\on{mult})^\IndCoh_*}\longrightarrow
\IndCoh^*(\fL_\nabla(H))$$
admits a (continuous) right adjoint.

\medskip

However, we can isomorph the projection
$$\underset{\fL^+_\nabla(H)}{\underline{\fL_\nabla(H)\times \fL_\nabla(H)}}\overset{\on{mult}}\longrightarrow \fL_\nabla(H)$$
to the projection
$$\fL_\nabla(H)\times (\fL_\nabla(H)/\fL^+_\nabla(H))\to \fL_\nabla(H),$$
and the assertion follows from the fact that $\fL_\nabla(H)/\fL^+_\nabla(H)$ is ind-proper. 

\sssec{}

In particular, we obtain that the monoidal operation on $\IndCoh^*(\on{Hecke}_H^{\on{spec,loc}})$ admits a continuous right adjoint, namely,
$$(\Delta_{\LS_H^\reg})^\IndCoh_*\circ (\on{mult})^!.$$

\medskip

Furthermore, it easy to see that this right adjoint 
$$\IndCoh^*(\on{Hecke}_H^{\on{spec,loc}})\to \IndCoh^*(\on{Hecke}_H^{\on{spec,loc}})\otimes \IndCoh^*(\on{Hecke}_H^{\on{spec,loc}})$$
is compatible with the $\IndCoh^*(\on{Hecke}_H^{\on{spec,loc}})$-bimodule structure. 

\medskip

Since the monoidal unit in $\IndCoh^*(\on{Hecke}_H^{\on{spec,loc}})$ is compact, we obtain that $\IndCoh^*(\on{Hecke}_H^{\on{spec,loc}})$
is \emph{rigid} (see \cite[Chapter 1, Definition 9.1.2]{GaRo3} for what this means), i.e., for any $S\to \Ran$, the monoidal category
$$\IndCoh^*(\on{Hecke}_H^{\on{spec,loc}})_S$$
is rigid. 

\sssec{}

Passing to duals, the monoidal structure on $\IndCoh^*(\on{Hecke}_H^{\on{spec,loc}})$ induces a comonoidal
structure on its dual, i.e., $\IndCoh^!(\on{Hecke}_H^{\on{spec,loc}})$. Since $\IndCoh^*(\on{Hecke}_H^{\on{spec,loc}})$ is
rigid, the comonoidal operation on $\IndCoh^!(\on{Hecke}_H^{\on{spec,loc}})$ admits a left adjoint, i.e., 
$\IndCoh^!(\on{Hecke}_H^{\on{spec,loc}})$ is naturally a monoidal category. 

\medskip

Let us describe the monoidal operation on $\IndCoh^!(\on{Hecke}_H^{\on{spec,loc}})$ explicitly. In terms of \eqref{e:corr Hecke 0},
it is given by
$$(\on{mult})^\IndCoh_*\circ (\Delta_{\LS_H^\reg})^!,$$
where $(\on{mult})^\IndCoh_*$ is the left adjoint of 
$$\on{mult}^!:\IndCoh^!(\LS_H^\reg \underset{\LS_H^\mer}\times \LS_H^\reg)\to 
\IndCoh^!(\LS_H^\reg\underset{\LS_H^\mer}\times \LS_H^\reg \underset{\LS_H^\mer}\times \LS_H^\reg),$$
or, which is the same as the dual of
$$(\on{mult})^!: \IndCoh^*(\LS_H^\reg \underset{\LS_H^\mer}\times \LS_H^\reg)\to 
\IndCoh^*(\LS_H^\reg\underset{\LS_H^\mer}\times \LS_H^\reg \underset{\LS_H^\mer}\times \LS_H^\reg),$$
whose existence was proved above. 

\sssec{}

Since $\IndCoh^*(\on{Hecke}_H^{\on{spec,loc}})$ is rigid, by \cite[Sect. 9.2.1]{GaRo3}, a choice of ``right" or ``left" determines
an equivalence
\begin{equation} \label{e:rigid equiv Sph spec}
\IndCoh^*(\on{Hecke}_H^{\on{spec,loc}})\simeq \IndCoh^!(\on{Hecke}_H^{\on{spec,loc}})
\end{equation} 
as monoidal categories.

\ssec{Action on \texorpdfstring{$\IndCoh(\LS_H)$}{IndCohLSH}}

In this subsection we will define a (local action) of the monoidal factorization category 
$\IndCoh^*(\on{Hecke}_H^{\on{spec,loc}})$ on $\IndCoh(\LS_H)$
in the sense of \secref{sss:local actions}. 

\sssec{} \label{sss:Sph spec action on LS}

Our goal is to define an action of the monoidal category 
$(\IndCoh^*(\on{Hecke}_H^{\on{spec,loc}}))^{\sotimes}_\Ran$ 
(see \secref{sss:pointwise monoidal on Ran} for the notation) on $\IndCoh(\LS_H)\otimes \ul\Dmod(\Ran)$.
In other words, we need to define an action of $\IndCoh^*(\on{Hecke}_H^{\on{spec,loc}})_S$ on 
$\IndCoh(\LS_H)\otimes \QCoh(S)$ for any $S\to \Ran$.

\medskip

Let $\on{Hecke}_{\cG,S}^{\on{spec,glob}}$ denote the fiber product
$$(\LS_H\times S)\underset{\LS^{\mer,\on{glob}}_{H,S}}\times (\LS_H\times S),$$
where
$$\LS^{\mer,\on{glob}}_{H,S}:=S\underset{\Ran}\times \LS^{\mer,\on{glob}}_{H,\Ran}$$
for $\LS^{\mer,\on{glob}}_{H,\Ran}$ as in \secref{sss:LS open curve}.

\medskip

Restriction to the formal disc gives rise to vertical arrows in the following diagram, see \secref{sss:LS open curve bis}:
\begin{equation} \label{e:Hecke spec diag}
\CD
\LS_H\times S  @<{\hl^{\on{spec,glob}}}<< \on{Hecke}_{H,S}^{\on{spec,glob}} @>{\hr^{\on{spec,glob}}}>> \LS_H\times S  \\
@V{\on{ev}_S}VV @V{\on{ev}_S}VV @VV{\on{ev}_S}V \\
\LS^\reg_{H,S} @<{\hl^{\on{spec,loc}}}<< \on{Hecke}_{H,S}^{\on{spec,loc}}  @>{\hr^{\on{spec,loc}}}>> \LS^\reg_{H,S}. 
\endCD
\end{equation} 

\begin{lem}  \label{l:Hecke spec diag}
Both squares in \eqref{e:Hecke spec diag} are Cartesian. 
\end{lem}

\begin{proof}

First, we claim that both squares are Cartesian at the classical level. Indeed, classically,
the horizontal arrows in \eqref{e:Hecke spec diag} are isomorphisms.

\medskip

Hence, in order to prove the lemma, it remains to check the Cartesian property of the
tangent spaces on the unit section of $\on{Hecke}_{H,S}^{\on{spec,glob}}$. 

\medskip

Let $\sigma$ be a point of $\LS_H$. Then for any $s\in S$, the relative tangent space of
$$\on{Hecke}_{H,S}^{\on{spec,glob}} \to \LS_H\times S$$
at (the image along the unit section of) $(\sigma,s)$ identifies with
$$\on{Fib}\left(\on{C}^\cdot(X,\fh_\sigma)[1]\to \on{C}^\cdot(X-\ul{x},\fh_\sigma)[1]\right)\simeq
\on{coFib}\left(\on{C}^\cdot(X,\fh_\sigma)\to \on{C}^\cdot(X-\ul{x},\fh_\sigma)\right),$$
where $\ul{x}\in \Ran$ is the image of $s$. 

\medskip

The relative tangent space of 
$$\on{Hecke}_{H,S}^{\on{spec,loc}} \to \LS^\reg_{H,S}$$ at the image of this point identifies with 
$$(\fL_\nabla(\fh_\sigma)/\fL^+_\nabla(\fh_\sigma))_{\ul{x}}.$$

\medskip

The required Cartesian property of the tangent spaces follows from the fact that diagram 
$$
\CD
\on{C}^\cdot(X,\fh_\sigma) @>>> \on{C}^\cdot(X-\ul{x},\fh_\sigma) \\
@VVV @VVV  \\
\fL^+_\nabla(\fh_\sigma)_{\ul{x}} @>>> \fL_\nabla(\fh_\sigma)_{\ul{x}} 
\endCD
$$
is Cartesian. 

\end{proof} 

\sssec{}

For $S\to \Ran$ as above, denote
$$(\LS_H\times S)^{\on{level}}:=(\LS_H\times S)\underset{\LS^\reg_{H,S}}\times S,$$
where $S\to \LS^\reg_{H,S}$ is the unit point.

\medskip 

By construction, $(\LS_H\times S)^{\on{level}}$ is acted on by $\fL_\nabla^+(H)_S$, so that
$$(\LS_H\times S)^{\on{level}}/\fL_\nabla^+(H)_S\simeq \LS_H\times S.$$

\medskip

Note that \lemref{l:Hecke spec diag} can be reformulated as saying that the above action of $\fL^+(H)_S$
on $(\LS_H\times S)^{\on{level}}$ extends to an action of $\fL(H)_S$. 

\sssec{}

Let $S_\alpha$ be as in \secref{sss:S alpha}. Consider the category
$$\IndCoh^*((\LS_H\times S_\alpha)^{\on{level}}).$$

\medskip

The action of $\fL_\nabla(H)_{S_\alpha}$ on $(\LS_H\times S_\alpha)^{\on{level}}$ gives rise
to an action of the monoidal category $\IndCoh^*(\fL_\nabla(H)_{S_\alpha})$ on 
$\IndCoh^*((\LS_H\times S_\alpha)^{\on{level}})$. 

\medskip 

As in \lemref{l:Psi arcs H}, the $\IndCoh$-pushforward functor 
$$\IndCoh^*((\LS_H\times S_\alpha)^{\on{level}})\to \IndCoh(\LS_H\times S_\alpha)$$
gives rise to an equivalence
$$\left(\IndCoh^*((\LS_H\times S_\alpha)^{\on{level}})\right)_{\fL_\nabla^+(H)_{S_\alpha}}\overset{\sim}\to
\IndCoh(\LS_H\times S_\alpha).$$

\medskip

Then by the same mechanism as in \secref{sss:mon op Hecke spec}, we obtain an action of the monoidal 
category $\IndCoh^*(\on{Hecke}_H^{\on{spec,loc}})_{S_\alpha}$ on $\IndCoh(\LS_H\times S_\alpha)$. 

\sssec{}

Note that we can explicitly describe the action functor as follows:
\begin{multline} \label{e:Hecke spec 1}
\IndCoh^*(\on{Hecke}_H^{\on{spec,loc}})_{S_\alpha}\underset{\QCoh(S_\alpha)}\otimes \IndCoh(\LS_H\times S_\alpha) 
\overset{(\on{ev}_S\times \hr^{\on{spec,glob}})^{*,\IndCoh}}\longrightarrow \\
\to \IndCoh((\on{Hecke}_H^{\on{spec,glob}})_{S_\alpha})\overset{(\hl^{\on{spec,glob}})^\IndCoh_*}\longrightarrow
\IndCoh(\LS_H\times S_\alpha),
\end{multline}
where the first arrow is obtained by identifying
\begin{multline*} 
\IndCoh^*(\on{Hecke}_H^{\on{spec,loc}})_{S_\alpha}\underset{\QCoh(S_\alpha)}\otimes \IndCoh(\LS_H\times S_\alpha)
\simeq \\
\simeq \left(\IndCoh^*(\fL_\nabla(H)\underset{S_\alpha}\times (\LS_H\times S_\alpha)^{\on{level}})\right)_{(\fL^+(H)\times \fL^+(H)\times \fL^+(H))_{S_\alpha}}
\end{multline*}
and 
$$\IndCoh((\on{Hecke}_H^{\on{spec,glob}})_{S_\alpha}) 
\simeq \left(\IndCoh^*(\fL_\nabla(H)\underset{S_\alpha}\times (\LS_H\times S_\alpha)^{\on{level}})\right)_{(\fL^+(H)\times \fL^+(H))_{S_\alpha}},$$
and the functor $(\on{ev}_S\times \hr^{\on{spec,glob}})^{*,\IndCoh}$ is the left adjoint to the projection from 
$(\fL^+(H)\times \fL^+(H))_{S_\alpha}$-coinvariants
to $(\fL^+(H)\times \fL^+(H)\times \fL^+(H))_{S_\alpha}$-coinvariants. 

\sssec{}

Having defined the action of $\IndCoh^*(\on{Hecke}_H^{\on{spec,loc}})_{S_\alpha}$ on $\IndCoh(\LS_H\times S_\alpha)$,
the procedure in Sects. \ref{sss:IndCoh* over Ran}-\ref{sss:IndCoh* over Ran transition} defines an action of 
$\IndCoh^*(\on{Hecke}_H^{\on{spec,loc}})_S$ on $\IndCoh(\LS_H\times S)$ for any $S\to \Ran$. 

\medskip

Thus, we obtain the sought-for local action of $\IndCoh^*(\on{Hecke}_H^{\on{spec,loc}})$ on $\IndCoh(\LS_H)$. Furthermore,
unwinding the construction, we obtain that this action has a natural Ran-unital structure (see \secref{sss:local actions} for what this means). 

\sssec{}

Recall the functor
$$\Loc_{H,\Ran}^{\on{spec}}:\Rep(H)_\Ran \to \IndCoh(\LS_H)\otimes \ul\Dmod(\Ran).$$

We claim:

\begin{prop} \label{p:Loc spec comp Sph}
The functor $\Loc_{H,\Ran}^{\on{spec}}$ intertwines the actions of $\IndCoh^*(\on{Hecke}_H^{\on{spec,loc}})^{\sotimes}_\Ran$
on the two sides.
\end{prop}

\begin{proof}

Unwinding the construction, we need to construct the datum of compatibility for the functor
$$\Loc_{H,S_\alpha}^{\on{spec}}:\Rep(H)_{S_\alpha} \to \IndCoh(\LS_H)\otimes \QCoh(S_\alpha)$$
for $S_\alpha$ as in \secref{sss:S alpha}.

\medskip

We identify
$$\Rep(H)_{S_\alpha} \simeq \QCoh(\LS^\reg_{H,S_\alpha}),$$
so that the functor $\Loc_{H,S_\alpha}^{\on{spec}}$ identifies with the functor
$$\QCoh(\LS^\reg_{H,S_\alpha})\overset{(\on{ev}_{S_\alpha})^{*,\IndCoh}}\longrightarrow \IndCoh(\LS_H\otimes S_\alpha)\simeq 
 \IndCoh(\LS_H)\otimes \QCoh(S_\alpha).$$
 
\medskip
 
Now the assertion of the proposition follows by unwinding the constructions, using the fact that
diagram \eqref{e:Hecke spec diag} is Cartesian (see \lemref{l:Hecke spec diag}).
 
\end{proof} 

\sssec{}

Consider now the functor
$$\Gamma^{\on{spec},\IndCoh}_{H,\Ran}: \IndCoh(\LS_H)\otimes \Dmod(\Ran)\to \Rep(H)_\Ran,$$
right adjoint to $\Loc_{H,\Ran}^{\on{spec}}$.

\medskip

Since $\IndCoh^*(\on{Hecke}_H^{\on{spec,loc}})$ is rigid, from \propref{p:Loc spec comp Sph}, we obtain:

\begin{cor} \label{c:Gamma spec comp Sph}
The functor $\Gamma^{\on{spec},\IndCoh}_{H,\Ran}$ intertwines the actions of $\IndCoh^*(\on{Hecke}_H^{\on{spec,loc}})^{\sotimes}_\Ran$
on the two sides.
\end{cor}

\sssec{} \label{sss:flip factors glob}

The (local) action of $\IndCoh^*(\on{Hecke}_H^{\on{spec,loc}})$ on $\IndCoh(\LS_H)$ gives rise to a (local)
\emph{right} action of $\IndCoh^*(\on{Hecke}_H^{\on{spec,loc}})$ on the dual $\IndCoh(\LS_H)^\vee$ of $\IndCoh(\LS_H)$.

\medskip

We identify 
$$\IndCoh(\LS_H)^\vee\simeq \IndCoh(\LS_H)$$
by Serre duality. Thus, we obtain a new (local) right action of $\IndCoh^*(\on{Hecke}_H^{\on{spec,loc}})$ on $\IndCoh(\LS_H)$.

\medskip

Note now that we can pass between right and left modules over $\IndCoh^*(\on{Hecke}_H^{\on{spec,loc}})$ 
using the anti-involution $\sigma^{\on{spec}}$, induced by the inversion operation of $\fL_\nabla(H)$.

\medskip

Unwinding the construction, we obtain that the resulting new (local) action of $\IndCoh^*(\on{Hecke}_H^{\on{spec,loc}})$ 
on $\IndCoh(\LS_H)$ coincides with the original one. 

\sssec{} \label{sss:mon op on LS}

The (local) action of $\IndCoh^*(\on{Hecke}_H^{\on{spec,loc}})$ on $\IndCoh(\LS_H)$ gives rise to a (local) 
\emph{coaction} of the factorization comonidal category $\IndCoh^!(\on{Hecke}_H^{\on{spec,loc}})$ on $\IndCoh(\LS_H)$. 

\medskip

Since $\IndCoh^*(\on{Hecke}_H^{\on{spec,loc}})$ is rigid, the coaction
functor admits a left adjoint, so we obtain a (local) \emph{action} of $\IndCoh^!(\on{Hecke}_H^{\on{spec,loc}})$,
viewed as a factorization monidal category, on $\IndCoh(\LS_H)$. 

\medskip

The corresponding monoidal operation is described explicitly as follows.

\medskip

For $S_\alpha$ as in \secref{sss:S alpha}, the action functor 
$$\IndCoh^!(\on{Hecke}_H^{\on{spec,loc}})_{S_\alpha}\underset{\QCoh(S_\alpha)}\otimes 
\IndCoh(\LS_H\times S_\alpha) \to \IndCoh(\LS_H\times S_\alpha)$$
is given by
\begin{multline} \label{e:Hecke spec 2}
\IndCoh(\LS_H\times S_\alpha)\underset{\QCoh(S_\alpha)}\otimes 
\IndCoh^!(\on{Hecke}_H^{\on{spec,loc}})_{S_\alpha}\simeq \\
\simeq \IndCoh^!\left((\LS_H\times S_\alpha) \underset{S_\alpha}\times (\on{Hecke}_H^{\on{spec,loc}})_{S_\alpha}\right)
\overset{(\hl^{\on{spec,glob}}\times \on{ev})^!}\longrightarrow \\
\to \IndCoh((\on{Hecke}_H^{\on{spec,glob}})_{S_\alpha})\overset{(\hr^{\on{spec,glob}})^\IndCoh_*}\longrightarrow
\IndCoh(\LS_H\times S_\alpha).
\end{multline}

For an arbitrary $S\to \Ran$, this action is extended by the mechanism of Sects. \ref{sss:IndCoh ! Ran 1}-\ref{sss:IndCoh ! Ran 2}. 

\sssec{} \label{sss:rigid act on LS}

Recall now that according to \eqref{e:rigid equiv Sph spec}, we can identify $\IndCoh^!(\on{Hecke}_H^{\on{spec,loc}})$
as a monoidal category with $\IndCoh^*(\on{Hecke}_H^{\on{spec,loc}})$.

\medskip

It is a formal property of actions of rigid categories that with respect to this identification, the above action of 
$\IndCoh^!(\on{Hecke}_H^{\on{spec,loc}})$ on $\IndCoh(\LS_H)$ identifies with the right action of
$\IndCoh^*(\on{Hecke}_H^{\on{spec,loc}})$ on $\IndCoh(\LS_H)$ from \secref{sss:flip factors glob}. 

\ssec{Action on monodromy-free opers}  \label{ss:Sph spec acts on mon-free}

In this subsection we take $H=\cG$, the Langlands dual of a reductive group $G$. We will
construct a factorization version of the action of $\IndCoh^*(\on{Hecke}_\cG^{\on{spec,loc}})$ on
$\IndCoh^*(\Op_\cG^\mf)$.

\sssec{} \label{sss:Sph spec acts on mon-free}

Let $\on{Hecke}_\cG^{\on{spec},\Op^\mf_\cG}$ be the factorization ind-scheme, defined as 
$$\on{Hecke}_\cG^{\on{spec},\Op^\mf_\cG}:=\Op_\cG^\mer\underset{\LS_\cG^\mer}\times \on{Hecke}_\cG^{\on{spec,loc}}.$$

Note that we have a commutative diagram
\begin{equation} \label{e:Sph spec acts on mon-free}
\CD
\Op_\cG^\mf @<{\hl^{\on{spec},\Op}}<< \on{Hecke}_\cG^{\on{spec},\Op^\mf_\cG} @>{\hr^{\on{spec},\Op}}>> \Op_\cG^\mf  \\
@V{\fr}VV @VV{\fr}V @VV{\fr}V \\
\LS^\reg_\cG @<{\hl^{\on{spec,loc}}}<< \on{Hecke}_\cG^{\on{spec,loc}} @>{\hr^{\on{spec,loc}}}>> \LS^\reg_\cG,
\endCD
\end{equation} 
in which both arrows are Cartesian.

\medskip

From here, by the same mechanism as in \secref{sss:mon op Hecke spec}, we obtain an action
of the monoidal category 
$\IndCoh^*(\on{Hecke}_H^{\on{spec,loc}})$ on $\IndCoh^*(\Op_\cG^\mf)$. 

\medskip

We write the action functor symbolically as 
\begin{multline*}
\IndCoh^*(\on{Hecke}_H^{\on{spec,loc}})\otimes \IndCoh^*(\Op_\cG^\mf) \overset{(\fr\times \hr^{\on{spec},\Op})^{*,\IndCoh}}\longrightarrow \\
\to \IndCoh^*(\on{Hecke}_\cG^{\on{spec},\Op^\mf_\cG}) \overset{(\hl^{\on{spec},\Op})^\IndCoh_*}\longrightarrow 
\IndCoh^*(\Op_\cG^\mf),
\end{multline*} 
where the functor $(\fr\times \hr^{\on{spec},\Op})^{*,\IndCoh}$ is assigned a meaning as in \secref{sss:mon op on LS}. 

\begin{rem}

Note that $\on{Hecke}_\cG^{\on{spec},\Op^\mf_\cG}$ is ind-placid. Indeed, this follows from
the fact that $\Op_\cG^\mf$ is placid, combined with the fact that the map $\hl^{\on{spec,loc}}$
is locally almost of finite presentation (see \lemref{l:const Gr is laft}).

\end{rem}

\sssec{}

Recall the functor
$$\on{Poinc}^{\on{spec}}_{\cG,*,\Ran}:\IndCoh^*(\Op^\mf_{\cG,\Ran})\to \IndCoh(\LS_\cG)\otimes \Dmod(\Ran),$$
see \secref{sss:Poinc spec *}. 

\medskip

We claim:

\begin{prop} \label{p:Poinc * Sph compat}
The functor $\on{Poinc}^{\on{spec}}_{\cG,*,\Ran}$ intertwines the actions of $\IndCoh^*(\on{Hecke}_\cG^{\on{spec,loc}})$
on the two sides.
\end{prop}

\begin{proof}

Unwinding the construction, we need to show that the functor
$$\on{Poinc}^{\on{spec}}_{\cG,*,S_\alpha}:\IndCoh^*(\Op^\mf_{\cG})_S\to \IndCoh(\LS_\cG)\otimes \QCoh(S_\alpha)$$
is compatible with the action of $\IndCoh^*(\on{Hecke}_\cG^{\on{spec,loc}})_{S_\alpha}$ for $S_\alpha$ as in 
\secref{sss:S alpha}.

\medskip

Recall that the functor $\on{Poinc}^{\on{spec}}_{\cG,*,S_\alpha}$ is the composition of:

\begin{itemize}

\item *-pullback along
$$\Op^{\mf,\on{glob}}_{\cG,S_\alpha} \overset{\on{ev}_{S_\alpha}}\longrightarrow \Op^\mf_{\cG,S_{\alpha}};$$

\item $\IndCoh$-pushforward along $\Op^{\mf,\on{glob}}_{\cG,S_\alpha} \overset{\fr^{\on{glob}}}\longrightarrow \LS_\cG\times S_\alpha$. 

\end{itemize}

Denote
$$\on{Hecke}_{\cG,S_\alpha}^{\on{spec,glob},\Op^\mf_\cG}:=
\Op^{\mer,\on{glob}}_{\cG,S_\alpha}\underset{\LS^{\mer,\on{glob}}_{\cG,S_\alpha}}\times \on{Hecke}_{H,S_\alpha}^{\on{spec,glob}}.$$

The assertion of the proposition holds by unwinding the constructions from the fact that in the following 
diagrams both square are Cartesian
$$
\CD
\Op^{\mf,\on{glob}}_{\cG,S_\alpha} @<{\hl^{\on{spec},\Op}}<< \on{Hecke}_{\cG,S_\alpha}^{\on{spec,glob},\Op^\mf_\cG} 
@>{\hl^{\on{spec},\Op}}>> \Op^{\mf,\on{glob}}_{\cG,S_\alpha}  \\
@V{\on{ev}}VV @V{\on{ev}}VV @VV{\on{ev}}V \\
 \Op_{\cG,S_\alpha}^\mf @<{\hl^{\on{spec,glob},\Op}}<< \on{Hecke}_{\cG,S_\alpha}^{\on{spec},\Op^\mf_\cG} @>{\hr^{\on{spec,glob},\Op}}>> 
 \Op_{\cG,S_\alpha}^\mf 
 \endCD
 $$
and
$$
\CD
\Op^{\mf,\on{glob}}_{\cG,S_\alpha} @<{\hl^{\on{spec},\Op}}<< \on{Hecke}_{\cG,S_\alpha}^{\on{spec,glob},\Op^\mf_\cG} 
@>{\hl^{\on{spec},\Op}}>> \Op^{\mf,\on{glob}}_{\cG,S_\alpha}  \\
@V{\fr^{\on{glob}}}VV @V{\fr^{\on{glob}}}VV @VV{\fr^{\on{glob}}}V \\
\LS_\cG @<{\hl^{\on{spec,loc}}}<< \on{Hecke}_\cG^{\on{spec,glob}} @>{\hr^{\on{spec,loc}}}>> \LS_\cG.
\endCD
$$

\end{proof} 

\sssec{} \label{sss:Hecke act on Op !}

The action of $\IndCoh^*(\on{Hecke}_\cG^{\on{spec,loc}})$ on $\IndCoh^*(\Op_\cG^\mf)$ gives rise to a \emph{right}
action of $\IndCoh^*(\on{Hecke}_\cG^{\on{spec,loc}})$ on the dual of $\IndCoh^*(\Op_\cG^\mf)$. Identifying 
$$\IndCoh^*(\Op_\cG^\mf)^\vee \simeq \IndCoh^!(\Op_\cG^\mf),$$
we thus obtain a right action of $\IndCoh^*(\on{Hecke}_\cG^{\on{spec,loc}})$ on $\IndCoh^!(\Op_\cG^\mf)$.

\medskip

Applying the anti-involution $\sigma^{\on{spec}}$, we can turn this right action into a left action.

\sssec{}

Passing to dual functors, from the action of $\IndCoh^*(\on{Hecke}_\cG^{\on{spec,loc}})$ on $\IndCoh^*(\Op_\cG^\mf)$, we obtain 
a coaction of $\IndCoh^!(\on{Hecke}_\cG^{\on{spec,loc}})$ on $\IndCoh^!(\Op_\cG^\mf)$. 

\medskip

By rigidity, the above coaction functor admits
a left adjoint, i.e., we obtain an action of $\IndCoh^!(\on{Hecke}_\cG^{\on{spec,loc}})$, viewed as a monoidal category, on
$\IndCoh^!(\Op_\cG^\mf)$.  

\medskip

The corresponding action functor is explicitly given by
\begin{multline*}
\IndCoh^!(\Op_\cG^\mf) \otimes \IndCoh^!(\on{Hecke}_\cG^{\on{spec,loc}})
\overset{(\hl^{\on{spec},\Op}\times \fr)^{*,\IndCoh}}\longrightarrow \\
\to \IndCoh^!(\on{Hecke}_\cG^{\on{spec},\Op^\mf_\cG}) \overset{(\hr^{\on{spec},\Op})^\IndCoh_*}\longrightarrow 
\IndCoh^!(\Op_\cG^\mf). 
\end{multline*} 


\medskip

As in \secref{sss:rigid act on LS}, it follows formally that with respect to the identification \eqref{e:rigid equiv Sph spec},
the above action of $\IndCoh^!(\on{Hecke}_\cG^{\on{spec,loc}})$ on $\IndCoh^!(\Op_\cG^\mf)$ identifies with the 
action of $\IndCoh^*(\on{Hecke}\cG^{\on{spec,loc}})$ on $\IndCoh^!(\Op_\cG^\mf)$ from \secref{sss:Hecke act on Op !}. 

\sssec{}

In a way analogous to \propref{p:Poinc * Sph compat}, one proves:

\begin{prop} \label{p:Poinc ! Sph compat}
The functor $\on{Poinc}^{\on{spec}}_{\cG,!,\Ran}$ intertwines the coactions of $\IndCoh^!(\on{Hecke}_\cG^{\on{spec,loc}})$
on the two sides.
\end{prop}

\begin{cor} \label{c:Poinc ! Sph compat}
The functor $\on{Poinc}^{\on{spec}}_{\cG,!,\Ran}$ intertwines the actions of $\IndCoh^*(\on{Hecke}_\cG^{\on{spec,loc}})$
on the two sides, where the action on the left-hand side is as in \secref{sss:Hecke act on Op !}.  
\end{cor}

\sssec{}

Note that \lemref{l:Theta and Sph} adapts to the factorization setting as follows: 

\begin{lem} \label{l:Theta and Sph rigid}
The equivalence 
$$\IndCoh^!(\Op_{\cG,\Ran}^\mf) \overset{\Theta_{\Op^\mf_\cG}}\simeq \IndCoh^*(\Op_{\cG,\Ran}^\mf)$$
is compatible with the actions of $\IndCoh^*(\on{Hecke}_{\cG,\Ran}^{\on{spec,loc}})$ on the two sides, where:
the action on the left-hand side is one from \secref{sss:Hecke act on Op !}. 
\end{lem}

\sssec{}

Recall now that according to \thmref{t:Poinc spec * vs !}, we have a canonical isomorphism: 
\begin{equation} \label{e:Poinc spec * vs ! Sph}
\Poinc^{\on{spec}}_{\cG,!,\Ran}\otimes \fl_{\on{Kost}(\cG)}[-\delta_G]\simeq 
\Poinc^{\on{spec}}_{\cG,*,\Ran}\circ \Theta_{\Op^\mf_\cG}.
\end{equation}

\medskip

Unwinding the constructions, we obtain:

\begin{lem}
The commutative diagram 
$$
\CD
\IndCoh^!(\Op_{\cG,\Ran}^\mf) @>{\Theta_{\Op^\mf_\cG}}>>  \IndCoh^*(\Op_{\cG,\Ran}^\mf)  \\
@V{\Poinc^{\on{spec}}_{\cG,!,\Ran}}VV @VV{\Poinc^{\on{spec}}_{\cG,*,\Ran}}V \\
\IndCoh(\LS_\cG)\otimes \Dmod(\Ran) @>{\on{Id}}>> \IndCoh(\LS_\cG)\otimes \Dmod(\Ran) 
\endCD
$$
upgrades to a commutative diagram of categories equipped with actions of $\IndCoh^*(\on{Hecke}_{\cG,\Ran}^{\on{spec,loc}})$,
where: 

\begin{itemize}

\item The compatibility for the left vertical arrow is given by \corref{c:Poinc ! Sph compat};

\item The compatibility for the right vertical arrow is given by \propref{p:Poinc * Sph compat};

\item The compatibility for the top horizontal arrow with \eqref{e:rigid equiv Sph spec} is given by \lemref{l:Theta and Sph rigid};


\end{itemize} 

\end{lem} 

\ssec{An approach to \texorpdfstring{$\Sph_H^{\on{spec}}$}{SphH} via factorization modules}  

In this subsection we will review the connection between the definition of
$\IndCoh^*(\on{Hecke}_{H,\Ran}^{\on{spec,loc}})$
developed above, and one given in \cite{CR}.

\sssec{}

The projection
$$\on{Hecke}_{H}^{\on{spec,loc}}:=\LS_H^\reg\underset{\LS^\mer_H}\times \LS^\reg_H\overset{f}\to \LS^\reg_H\times \LS^\reg_H$$
gives rise to a lax unital functor
\begin{multline*} 
f^\IndCoh_*: \IndCoh^*(\on{Hecke}_H^{\on{spec,loc}})\to \IndCoh^*(\LS^\reg_H\times \LS^\reg_H)\overset{\Psi_{\LS^\reg_H\times \LS^\reg_H}}\simeq \\
\simeq \QCoh(\LS^\reg_H\times \LS^\reg_H)\simeq \Rep_{H\times H}\simeq \Rep_H\otimes \Rep_H,
\end{multline*} 
where the first arrow is obtained by applying the functor of $\fL^+_\nabla(H)\times \fL^+_\nabla(H)$-coinvariants to
$$\Gamma^\IndCoh(\fL_\nabla(H),-)\to \Vect.$$

\sssec{}

Note that the image of the factorization unit 
$$\one_{\IndCoh^*(\on{Hecke}_H^{\on{spec,loc}})}\simeq \iota^\IndCoh_*(\CO_{\LS^\reg_H})$$
identifies with
$$(\Delta_{\LS^\reg_H})_*(\CO_{\LS^\reg_H})\in \QCoh(\LS^\reg_H\times \LS^\reg_H) \, \Leftrightarrow \, 
R_H\in \Rep_{H\times H},$$
where $R_H$ denotes the regular representation, viewed as a commutative factorization algebra in $\Rep_{H\times H}$

\medskip

By \lemref{l:enhancement modules}, the functor $f_*$ enhances to a unital functor
\begin{equation} \label{e:enh Sph Spec}
(f^\IndCoh_*)^{\on{enh}}:\IndCoh^*(\on{Hecke}_H^{\on{spec,loc}})\to R_H\mod^{\on{fact}}(\Rep_H\otimes \Rep_H),
\end{equation}
where the right-hand side is viewed as a unital lax factorization category.

\sssec{}

We claim:

\begin{prop} \label{p:enh Sph Spec} \hfill

\smallskip

\noindent{\em(a)} 
The functor \eqref{e:enh Sph Spec} induces an equivalences between the eventually
coconnective subcategories of the two sides. 

\smallskip

\noindent{\em(b)} 
The essential image of $\IndCoh^*(\on{Hecke}_H^{\on{spec,loc}})^c\subset \IndCoh^*(\on{Hecke}_H^{\on{spec,loc}})$ 
under the functor \eqref{e:enh Sph Spec} is contained
in $\left(R_H\mod^{\on{fact}}(\Rep_H\otimes \Rep_H)\right)^{>-\infty}$.

\end{prop}

\begin{proof}

The proof proceeds along the same lines as that of \propref{p:IndCoh Op via fact almost}, with the following
difference:

\medskip

Instead of appealing to \propref{p:QCoh co descent}, we claim that $\QCoh_{\on{co}}(\on{Hecke}_H^{\on{spec,loc}})$ identifies
with the totalization of the cosimplicial category 
$$\QCoh_{\on{co}}(\on{Hecke}_{H}^{\on{spec,loc}} \underset{\LS^\reg_H\times \LS^\reg_H}\times \on{pt}^\bullet),$$
where $\on{pt}^\bullet$ is the \v{C}ech nerve of the projection
$$\on{pt}\to \LS^\reg_H\times \LS^\reg_H,$$
i.e., we claim that the functor 
\begin{equation}  \label{e:enh Sph Spec 1}
\QCoh_{\on{co}}(\on{Hecke}_{H}^{\on{spec,loc}})\to \QCoh_{\on{co}}(\fL_\nabla(H))^{\fL_\nabla^+(H)\times \fL_\nabla^+(H)}
\end{equation} 
is an equivalence. 

\medskip

Indeed, the precomposition of \eqref{e:enh Sph Spec 1} with \eqref{e:coinv co LH nabla} is the functor
$$\QCoh_{\on{co}}(\fL_\nabla(H))_{\fL_\nabla^+(H)\times \fL_\nabla^+(H)}\to \QCoh_{\on{co}}(\fL_\nabla(H))^{\fL_\nabla^+(H)\times \fL_\nabla^+(H)}$$
of \eqref{e:inv vs coinv}, which is an equivalence by \corref{c:inv vs coinv}.

\end{proof}

\sssec{} \label{sss:enh Sph Spec}

Note that we have a commutative diagram
$$
\CD 
\QCoh(\LS^\reg_H) @>{\iota^\IndCoh_*}>>  \IndCoh^*(\on{Hecke}_H^{\on{spec,loc}}) \\
@V{\sim}VV \\
\Rep(H) & & @VV{(f^\IndCoh_*)^{\on{enh}}}V \\
@V{\sim}VV \\
R_H\mod^{\on{com}}(\Rep(H)\otimes \Rep(H)) @>>> R_H\mod^{\on{fact}}(\Rep(H)\otimes \Rep(H)).
\endCD
$$ 

\medskip 

Recall now that $\IndCoh^*(\on{Hecke}_H^{\on{spec,loc}})^c$ is generated under finite colimits by the essential image
of $\QCoh(\LS^\reg_H)$ along $\iota^\IndCoh_*$. 

\medskip

Hence, from \propref{p:enh Sph Spec}, we obtain that the essential image of $\IndCoh^*(\on{Hecke}_H^{\on{spec,loc}})^c$ 
under \eqref{e:enh Sph Spec} 
is the full subcategory of $R_H\mod^{\on{fact}}(\Rep(H)\otimes \Rep(H))$ generated under finite colimits by the essential image
of $\Rep(H)^c$ under the functor
\begin{equation} \label{e:RH com}
\Rep(H)\simeq R_H\mod^{\on{com}}(\Rep(H)\otimes \Rep(H)) \to R_H\mod^{\on{fact}}(\Rep(H)\otimes \Rep(H)).
\end{equation}

This allows us to recover $\IndCoh^*(\on{Hecke}_H^{\on{spec,loc}})$ from $R_H\mod^{\on{fact}}(\Rep(H)\otimes \Rep(H))$
by an explicit procedure:

\medskip

Namely, $\IndCoh^*(\on{Hecke}_H^{\on{spec,loc}})$ identifies with the ind-completion of the full subcategory of 
$R_H\mod^{\on{fact}}(\Rep(H)\otimes \Rep(H))$
generated under finite colimits by the essential image of $\Rep(H)^c$ under the functor \eqref{e:RH com}. 

\medskip

In the rest of this subsection we will show how to recover various pieces of structure on $\IndCoh^*(\on{Hecke}_H^{\on{spec,loc}})$ 
from those on $R_H\mod^{\on{fact}}(\Rep(H)\otimes \Rep(H))$.

\sssec{}

Note that $\Rep(H)\otimes \Rep(H)$ is naturally a factorization monoidal category under \emph{convolution},
and $R_H$ is the monoidal unit.

\medskip

Hence, the monoidal operation
\begin{equation} \label{e:conv H}
(\Rep(H)\otimes \Rep(H)) \otimes (\Rep(H)\otimes \Rep(H))\overset{\star}\to \Rep(H)\otimes \Rep(H)
\end{equation}
sends the factorization algebra 
$$R_H\otimes R_H\in \on{FactAlg}^{\on{untl}}(X,(\Rep(H)\otimes \Rep(H)) \otimes (\Rep(H)\otimes \Rep(H)))$$
to
$$R_H\in \on{FactAlg}^{\on{untl}}(X,\Rep(H)\otimes \Rep(H)).$$

\medskip

From here we obtain that the functor \eqref{e:conv H} induces a functor
\begin{multline} \label{e:conv H R}
(R_H\mod^{\on{fact}}(\Rep(H)\otimes \Rep(H)))\otimes (R_H\mod^{\on{fact}}(\Rep(H)\otimes \Rep(H))) \overset{\star}\to \\
\to R_H\mod^{\on{fact}}(\Rep(H)\otimes \Rep(H)),
\end{multline}
which naturally extends to a monoidal structure on $R_H\mod^{\on{fact}}(\Rep(H)\otimes \Rep(H))$.

\sssec{} \label{sss:recover mon structure Sph}

Unwinding the constructions, we obtain that the functor \eqref{e:enh Sph Spec} is monoidal.

\medskip

The monoidal operation \eqref{e:conv H R} is t-exact. Hence, the procedure of recovering $\IndCoh^*(\on{Hecke}_H^{\on{spec,loc}})$
from $R_H\mod^{\on{fact}}(\Rep(H)\otimes \Rep(H))$ described in \secref{sss:enh Sph Spec} allows us also to recover its monoidal structure. 

\sssec{} \label{sss:bimod acting}

Let $\bA$ be a (unital) factorization category, equipped with a monoidal action of $\Rep(H)\otimes \Rep(H)$. Let
$\CA$ be \emph{any} (unital) factorization algebra in $\bA$. 

\medskip

The action functor 
\begin{equation} \label{e:conv H act}
(\Rep(H)\otimes \Rep(H)) \otimes \bA \overset{\star}\to \bA
\end{equation}
automatically sends the factorization algebra 
$$R_H\otimes \CA\in \on{FactAlg}^{(\on{untl})}(X,(\Rep(H)\otimes \Rep(H)) \otimes \bA)$$
to
$$\CA \in \on{FactAlg}^{(\on{untl})}(X,\Rep(H)\otimes \Rep(H)).$$

Hence, we obtain that the functor \eqref{e:conv H act} gives rise to a functor
\begin{equation} \label{e:conv H act R}
(R_H\mod^{\on{fact}}(\Rep(H)\otimes \Rep(H)))\otimes (\CA\mod^{\on{fact}}\mod(\bA)) \overset{\star}\to 
\CA\mod^{\on{fact}}\mod(\bA),
\end{equation}
which extends to a monoidal action of $R_H\mod^{\on{fact}}(\Rep(H)\otimes \Rep(H))$ on $\CA\mod^{\on{fact}}\mod(\bA)$.

\sssec{}

Take $\bA=\Rep(H)$ and $\CA=\one_{\Rep(H)}=\CO_{\Rep(H)}$, so that
$$\CA\mod^{\on{fact}}\mod(\bA)=\Rep(H).$$

\medskip

Hence, we obtain an action of $R_H\mod^{\on{fact}}(\Rep(H)\otimes \Rep(H))$ on $\Rep(H)$.

\medskip

Unwinding the definitions, we obtain that the functor \eqref{e:enh Sph Spec} intertwines the above action
with the action of $\IndCoh^*(\on{Hecke}_H^{\on{spec,loc}})$ on $\QCoh(\LS^\reg_H)$ from \secref{sss:Sph acts on Rep}. 

\medskip

Furthermore, as in \secref{sss:recover mon structure Sph}, this allows us to recover the latter
action from the action of $R_H\mod^{\on{fact}}(\Rep(H)\otimes \Rep(H))$ on $\Rep(H)$.

\sssec{}

Let us now take $H=\cG$. Let us take again $\bA=\Rep(\cG)$, but let us take $\CA:=R_{\cG,\Op}$
from \eqref{e:R G Op}. 

\medskip

We obtain an action of $R_\cG\mod^{\on{fact}}(\Rep(\cG)\otimes \Rep(\cG))$ on 
$R_{\cG,\Op}\mod^{\on{fact}}(\Rep(\cG))$. 

\medskip

Unwinding the constructions, we obtain that the functors \eqref{e:enh Sph Spec} and $(\fr^\IndCoh_*)^{\on{enh}}$ 
intertwine the above action with the action of $\IndCoh^*(\on{Hecke}_H^{\on{spec,loc}})$ on $\IndCoh^*(\Op^\mf_\cG)$
from \secref{sss:Sph spec acts on mon-free}. 

\medskip

Recall now that according to \propref{p:IndCoh Op via fact almost} the functor
$$(\fr^\IndCoh_*)^{\on{enh}}:\IndCoh^*(\Op^\mf_\cG)\to R_{\cG,\Op}\mod^{\on{fact}}(\Rep(\cG))$$
induces an equivalence between the eventually coconnective subcategories of the two sides
and sends compact objects to eventually coconnective ones. Combined with \propref{p:enh Sph Spec},
this allows us to recover the action of $\IndCoh^*(\on{Hecke}_H^{\on{spec,loc}})$ on $\IndCoh^*(\Op^\mf_\cG)$
from the above action of $R_\cG\mod^{\on{fact}}(\Rep(\cG)\otimes \Rep(\cG))$ on 
$R_{\cG,\Op}\mod^{\on{fact}}(\Rep(\cG))$. 

\ssec{Compatibility of the FLE with (derived) Satake} \label{ss:FLE and Sat}

In this subsection we continue to take $H=\cG$. Our goal is to prove \thmref{t:FLE and Sat} in the factorization
setting.

\sssec{}

As a first step, we recall the construction of the geometric equivalence functor $\on{Sat}_G$. 
Consider the factorization category 
$$\Whit^!(G)\otimes \Whit_*(G)$$
and note that it is naturally a bimodule\footnote{Recall that according to \secref{sss:left to right Sph} 
we freely pass between left and right modules for $\Sph_G$ using the anti-involution $\sigma$.} 
for the factorization monoidal category $\Sph_G$.

\medskip

We identify
$$\Whit^!(G) \overset{\on{CS}_G}\simeq \Rep(\cG),$$
and we identify
$$\Whit_*(G) \overset{\FLE_{\cG,\infty}}\simeq \Rep(\cG).$$

Thus, we obtain that $\Rep(\cG)\otimes \Rep(\cG)$ acquires a bimodule structure with respect to 
$\Sph_G$. The action on the factorization algebra object $R_\cG\in \Rep(G)\otimes \Rep(\cG)$
gives rise to a factorization functor
$$\on{pre-Sat}_G:\Sph_G\to \Rep(\cG)\otimes \Rep(\cG).$$

Since the monoidal unit $\delta_{1,\Gr_G}\in \Sph_G$ equals the factorization unit, the functor $\on{pre-Sat}_G$
sends $\one_{\Sph_G}$ to $R_\cG$.

\medskip

The above operations are compatible with unital structures. Hence, by \lemref{l:enhancement modules}, 
the functor $\on{pre-Sat}_G$ upgrades to a functor
$$(\on{pre-Sat}_G)^{\on{enh}}: \Sph_G\to R_\cG\mod^{\on{fact}}(\Rep(\cG)\otimes \Rep(\cG)).$$

Unwinding the constructions, we obtain that the above functor $(\on{pre-Sat}_G)^{\on{enh}}$ carries
a naturally defined monoidal structure. 

\sssec{} 

We observe:

\begin{lem}
The functor $(\on{pre-Sat}_G)^{\on{enh}}$ sends the subcategory $(\Sph_G)^c\subset \Sph_G$ to 
the full subcategory of $R_\cG\mod^{\on{fact}}(\Rep(\cG)\otimes \Rep(\cG))$ generated under finite colimits by the essential image
of $\Rep(\cG)^c$ under the functor
$$R_\cG\mod^{\on{com}}(\Rep(\cG)\otimes \Rep(\cG)) \to R_\cG\mod^{\on{fact}}(\Rep(\cG)\otimes \Rep(\cG)).$$
\end{lem}

Thus, combining with \secref{sss:enh Sph Spec}, we obtain that the functor $(\on{pre-Sat}_G)^{\on{enh}}$ can
be uniquely factored as
$$\Sph_G\to \Sph_\cG^{\on{spec}}:=\IndCoh^*(\on{Hecke}_\cG^{\on{spec,loc}})\overset{\text{\eqref{e:enh Sph Spec}}}\longrightarrow
R_\cG\mod^{\on{fact}}(\Rep_\cG\otimes \Rep_\cG),$$
where the first arrow preserves compactness.

\medskip

The resulting functor
$$\Sph_G\to \Sph_\cG^{\on{spec}}$$
is the functor $\Sat_G$, as it was constructed in \cite{CR}. 

\sssec{} \label{sss:cor CR}

For the proof of  \thmref{t:FLE and Sat}, we will need the following output of the above 
construction: 

\medskip

Consider the functor
\begin{equation} \label{e:functor to use}
\Sph_G \overset{\on{pre-Sat}_G}\longrightarrow \Rep(\cG)\otimes \Rep(\cG)
\overset{\on{Id}\otimes \FLE_{\cG,\infty}}\longrightarrow \Rep(\cG)\otimes \Whit_*(G).
\end{equation}

By construction, it has the following properties:

\begin{itemize}

\item It intertwines the action of $\Sph_G$ on itself by right multiplication with the natural action of
$\Sph_G$ on $\Whit_*(G)$; 

\item It intertwines the action of $\Sph_G$ on itself by left multiplication with the action of $\Sph_\cG^{\on{spec}}$ on 
$\Rep(\cG)$ via $\Sat_G$.

\end{itemize} 

\sssec{}

Let $\bA$ be a factorization category, equipped with a factorization action of $\fL(G)_{\rho(\omega_X)}$ at the critical level. We have a naturally defined
factorization functor 
$$\Dmod_\crit(\Gr_{G,\rho(\omega_X)}) \underset{\Sph_G}\otimes \Sph(\bA) \to \bA,$$
which gives rise to a factorization functor
\begin{equation} \label{e:Sph to Whit App}
\Whit_*(G)\underset{\Sph_G}\otimes \Sph(\bA) \to \Whit_*(\bA).
\end{equation} 

By duality, the functor \eqref{e:Sph to Whit App} gives rise to a functor
\begin{equation} \label{e:Sph to Whit App1}
\Sph(\bA) \to \Whit^!(G)\otimes \Whit_*(\bA),
\end{equation} 
compatible with the action of $\Sph_G$. 

\sssec{}

Composing with $\on{CS}_G$ along the first factor, from \eqref{e:Sph to Whit App1} we obtain a functor
\begin{equation} \label{e:Sph to Whit App2}
\Sph(\bA) \to \Rep(\cG)\otimes \Whit_*(\bA). 
\end{equation} 

We claim:

\begin{lem}   \label{l:Sph to Whit App}
The functor \eqref{e:Sph to Whit App2} intertwines the action of $\Sph_G$ on 
$\Sph(\bA)$ with the action of $\Sph_\cG^{\on{spec}}$ on $\Rep(\cG)$ via
$\Sat_G$.
\end{lem}

\begin{proof}

We can interpret the functor \eqref{e:Sph to Whit App2} as follows:
\begin{multline}
\Sph(\bA)   \simeq \Sph_G\underset{\Sph_G}\otimes \Sph(\bA) 
\overset{\text{\eqref{e:functor to use}}}\longrightarrow \left(\Rep(\cG)\otimes \Whit_*(G)\right) \underset{\Sph_G}\otimes \Sph(\bA)  \simeq \\
\simeq \Rep(\cG) \otimes (\Whit_*(G)\underset{\Sph_G}\otimes \Sph(\bA)) \overset{\on{Id}\otimes \text{\eqref{e:Sph to Whit App}}}
\longrightarrow  \Rep(\cG)\otimes \Whit_*(\bA).
\end{multline} 

Now the assertion follows from the second bullet point in \secref{sss:cor CR}. 

\end{proof} 

\sssec{} \label{sss:FLE and Sat App}

Consider the morphism
$$\Op^\mf_\cG \overset{\fr\times \iota^\mf}\to \LS^\reg_\cG\times \Op^\mer_\cG,$$
and the corresponding factorization functor
\begin{equation} \label{e:FLE and Sat App 1}
\IndCoh^*(\Op^\mf_\cG) \overset{(\fr\times \iota^\mf)^\IndCoh_*}\longrightarrow \IndCoh^*(\LS^\reg_\cG\times \Op^\mer_\cG)
\simeq \Rep(\cG)\otimes \IndCoh^*(\Op^\mer_\cG).
\end{equation}

Unwinding the construction, we obtain that \eqref{e:FLE and Sat App 1} is compatible with the actions
of $\Sph_\cG^{\on{spec}}$ on the two sides, where $\Sph_\cG^{\on{spec}}$ acts on the right-hand side via the $\LS^\reg_\cG$-factor. 

\medskip

Denote 
$$R^{\Rep}_{\cG,\Op}:=((\fr\times \iota^\mf)^\IndCoh_*)(\CO_{\Op^\reg_\cG})\in 
\on{FactAlg}^{\on{untl}}(X,\Rep(\cG)\otimes \IndCoh^*(\Op^\mer_\cG)).$$

By \lemref{l:enhancement modules}, the functor \eqref{e:FLE and Sat App 1} enhances to a functor
\begin{equation} \label{e:FLE and Sat App 2}
((\fr\times \iota^\mf)^\IndCoh_*)^{\on{enh}}:
\IndCoh^*(\Op^\mf_\cG)\to R^{\Rep}_{\cG,\Op}\mod^{\on{fact}}(\Rep(\cG)\otimes \IndCoh^*(\Op^\mer_\cG)).
\end{equation}

The functor \eqref{e:FLE and Sat App 2} intertwines the $\Sph_\cG^{\on{spec}}$-action on the left-hand side
and the action of $R_\cG\mod^{\on{fact}}(\Rep(\cG)\otimes \Rep(\cG))$ from \secref{sss:bimod acting} on the right-hand side
via the functor \eqref{e:enh Sph Spec}. 
 
\sssec{}

We will prove:

\begin{prop} \label{p:FLE and Sat App} \hfill

\smallskip

\noindent{\em(a)} The functor \eqref{e:FLE and Sat App 2} induces an equivalence between the eventually coconnective subcategories
of the two sides.

\smallskip

\noindent{\em(b)} The essential image of the subcategory of compact objects in $\IndCoh^*(\Op^\mf_\cG)$ under the functor
\eqref{e:FLE and Sat App 2} is contained in $\left(R^{\Rep}_{\cG,\Op}\mod^{\on{fact}}(\Rep(\cG)\otimes \IndCoh^*(\Op^\mer_\cG))\right)^{>-\infty}$.

\end{prop}

The proof will be given in \secref{ss:FLE and Sat App proof}. Let us accept this proposition temporarily and proceed
with the proof of \thmref{t:FLE and Sat}. 

\sssec{}

We now launch the proof of \thmref{t:FLE and Sat} proper.

\medskip

Consider the functor
\begin{equation} \label{e:FLE and Sat App 3}
\KL(G)_\crit \overset{\FLE_{G,\crit}}\longrightarrow \IndCoh^*(\Op^\mf_\cG)  \overset{\text{\eqref{e:FLE and Sat App 1}}}\longrightarrow
\Rep(\cG)\otimes \IndCoh^*(\Op^\mer_\cG).
\end{equation}

We consider its enhancement
\begin{equation} \label{e:FLE and Sat App 3 bis}
\KL(G)_\crit \overset{\FLE_{G,\crit}}\longrightarrow \IndCoh^*(\Op^\mf_\cG)  \overset{\text{\eqref{e:FLE and Sat App 2}}}\longrightarrow 
R^{\Rep}_{\cG,\Op}\mod^{\on{fact}}(\Rep(\cG)\otimes \IndCoh^*(\Op^\mer_\cG)).
\end{equation}

\medskip

By Propositions \ref{p:enh Sph Spec} and \ref{p:FLE and Sat App}, in order to show that the functor $\FLE_{G,\crit}$
intertwines the action of $\Sph_G$ on $\KL(G)_\crit$ with the action of $\Sph_\cG^{\on{spec}}$ on $\IndCoh^*(\Op^\mf_\cG)$, 
it suffices to show that the composite functor in \eqref{e:FLE and Sat App 3 bis}
intertwines the action of $\Sph_G$ on the left-hand side with the action 
of $R_\cG\mod^{\on{fact}}(\Rep(\cG)\otimes \Rep(\cG))$ from \secref{sss:bimod acting} on the right-hand side via
the functor $(\on{pre-Sat}_G)^{\on{enh}}$.

\sssec{}

By the construction of the action of $R_\cG\mod^{\on{fact}}(\Rep(\cG)\otimes \Rep(\cG))$ on the right-hand side
of \eqref{e:FLE and Sat App 3 bis} in \secref{sss:bimod acting}, it suffices to show that the original functor \eqref{e:FLE and Sat App 3}
intertwines the action of $\Sph_G$ on the left-hand side with the action of $R_\cG\mod^{\on{fact}}(\Rep(\cG)\otimes \Rep(\cG))$ 
on the right-hand side via the functor $(\on{pre-Sat}_G)^{\on{enh}}$.

\medskip

We will show that the functor \eqref{e:FLE and Sat App 3} intertwines the action of $\Sph_G$ on the left-hand side with the action of
$\Spc_\cG^{\on{spec}}$ on the right-hand side via $\Sat_G$.

\sssec{}

We claim that the functor \eqref{e:FLE and Sat App 3} identifies canonically with the functor
\begin{multline} \label{e:FLE and Sat App 4}
\KL(G)_\crit \overset{\alpha_{\rho(\omega_X),\on{taut}}} \longrightarrow
\KL(G)_{\crit,\rho(\omega_X)} = \Sph(\hg\mod_{\crit,\rho(\omega_X)}) \overset{\text{\eqref{e:Sph to Whit App2}}}\longrightarrow \\
\to \Rep(\cG) \otimes \Whit_*(\hg\mod_{\crit,\rho(\omega_X)}) 
\overset{\on{Id}\otimes \ol\DS^{\on{enh,rfnd}}}\longrightarrow \Rep(\cG) \otimes \IndCoh^*(\Op^\mer_\cG).
\end{multline} 

Indeed, by construction, both functors are $\Rep(\cG)$-linear. Hence, since $\Rep(\cG)$ is rigid as a monoidal category, 
it suffices to identify the compositions of \eqref{e:FLE and Sat App 3} and \eqref{e:FLE and Sat App 4} with the
functor
$$(\on{inv}_\cG\otimes \on{Id}):\Rep(\cG) \otimes \IndCoh^*(\Op^\mer_\cG)\to \IndCoh^*(\Op^\mer_\cG).$$

The composition involving \eqref{e:FLE and Sat App 3} becomes
$$\KL(G)_\crit \overset{\FLE_{G,\crit}}\longrightarrow \IndCoh^*(\Op^\mf_\cG)  \overset{(\iota^\mf)^\IndCoh_*}\longrightarrow 
\IndCoh^*(\Op^\mer_\cG),$$
which is by construction
\begin{multline} \label{e:FLE and Sat App 5}
\KL(G)_\crit \overset{\alpha_{\rho(\omega_X),\on{taut}}} \longrightarrow
\KL(G)_{\crit,\rho(\omega_X)} \to \hg\mod_{\crit,\rho(\omega_X)}\to \\
\to \Whit_*(\hg\mod_{\crit,\rho(\omega_X)}) 
\overset{\ol\DS^{\on{enh,rfnd}}}\longrightarrow  \IndCoh^*(\Op^\mer_\cG).
\end{multline} 

Unwinding, we obtain that the composition involving \eqref{e:FLE and Sat App 4} also 
identifies with \eqref{e:FLE and Sat App 5}. 

\sssec{}

Thus, we have to show that the functor \eqref{e:FLE and Sat App 4} intertwines the action of $\Sph_G$ on the left-hand side with the action of
$\Spc_\cG^{\on{spec}}$ on the right-hand side via $\Sat_G$.

\medskip

However, this follows from \lemref{l:Sph to Whit App}. 

\sssec{}

It remains to establish the commutativity of \eqref{e:FLE and Sat diag}. This is equivalent to the commutativity of the diagram
$$
\CD
\KL(G)_\crit  @>{\alpha_{\rho(\omega_X),\on{taut}}}>> \KL(G)_{\crit ,\rho(\omega_X)} @>>> \Whit^!(G)\otimes \Whit_*(\hg\mod_{\crit,\rho(\omega_X)})  \\
@V{\FLE_{G,\crit}}VV & & @VV{\on{CS}_G\otimes \ol\DS^{\on{enh,rfnd}}}V \\
\IndCoh^*(\Op^\mf_\cG) & @>{\text{\eqref{e:FLE and Sat App 1}}}>> & \Rep(\cG) \otimes \IndCoh^*(\Op^\mer_\cG),
\endCD
$$
compatibly with the action of $\Sph_G$ on the top row and the action of $\Sph_\cG^{\on{spec}}$ on the bottom row
via $\Sat_G$. 

\medskip

However, this amounts to the identification between \eqref{e:FLE and Sat App 3} and \eqref{e:FLE and Sat App 4} established above. 

\qed[\thmref{t:FLE and Sat}]

\ssec{Proof of \propref{p:FLE and Sat App}} \label{ss:FLE and Sat App proof}

\sssec{}

Consider the functor
\begin{equation} \label{e:rel ch alg}
\on{Id}\otimes \Gamma^\IndCoh(\Op^\mer_\cG,-): \Rep(\cG)\otimes \IndCoh^*(\Op^\mer_\cG)\to \Rep(\cG).
\end{equation}

It sends $R^{\Rep}_{\cG,\Op}\mapsto R_{\cG,\Op}$ and hence induces a functor
\begin{equation} \label{e:rel ch alg 1}
R^{\Rep}_{\cG,\Op}\mod^{\on{fact}}(\Rep(\cG)\otimes \IndCoh^*(\Op^\mer_\cG))\to
R_{\cG,\Op}\mod^{\on{fact}}(\Rep(\cG))
\end{equation}

The composition of the functor \eqref{e:FLE and Sat App 2} with \eqref{e:rel ch alg 1} is the functor
\eqref{e:Op to LS enh}. Hence, by \propref{p:IndCoh Op via fact almost}, it suffices to show that the functor 
\eqref{e:rel ch alg 1} induces an equivalence between the eventually coconnective subcategories
of the two sides.
 
\sssec{}
 
The functor
$\on{Id}\otimes \Gamma^\IndCoh(\Op^\mer_\cG,-)$ enhances to a functor
\begin{equation} \label{e:rel ch alg 2}
\Rep(\cG)\otimes \IndCoh^*(\Op^\mer_\cG)\to 
\CO_{\Op^\reg_\cG}\mod^{\on{fact}}(\Rep(\cG))
\end{equation}
where $\CO_{\Op^\reg_\cG}$ is regarded as a constant (commutative)
factorization algebra in $\Rep(\cG)$. 
 
 \medskip
 
The functor \eqref{e:rel ch alg 2} sends $R^{\Rep}_{\cG,\Op}\in \on{FactAlg}^{\on{untl}}(X,\Rep(\cG)\otimes \IndCoh^*(\Op^\mer_\cG))$
to $R_{\cG,\Op}$, regarded as a unital factorization algebra in $\CO_{\Op^\reg_\cG}\mod^{\on{fact}}(\Rep(\cG))$
via the homomorphism $\CO_{\Op^\reg_\cG}\to R_{\cG,\Op}$, see \secref{sss:A and A' mod untl}.

\medskip

Hence, \eqref{e:rel ch alg 2} induces a functor
\begin{equation} \label{e:rel ch alg 3}
R^{\Rep}_{\cG,\Op}\mod^{\on{fact}}(\Rep(\cG)\otimes \IndCoh^*(\Op^\mer_\cG))\to
R_{\cG,\Op}\mod^{\on{fact}}\left(\CO_{\Op^\reg_\cG}\mod^{\on{fact}}(\Rep(\cG))\right).
\end{equation}

\medskip

The functor \eqref{e:rel ch alg 1} is the composition of \eqref{e:rel ch alg 3} with the forgetful functor
\begin{equation} \label{e:rel ch alg 4}
\oblv_{\CO_{\Op^\reg_\cG}}:R_{\cG,\Op}\mod^{\on{fact}}\left(\CO_{\Op^\reg_\cG}\mod^{\on{fact}}(\Rep(\cG))\right)\to
R_{\cG,\Op}\mod^{\on{fact}}(\Rep(\cG)). 
\end{equation}

\sssec{}

Note now that by (a relative version of) \corref{c:IndCoh* Op bdd below}, the functor \eqref{e:rel ch alg 2}
induces an equivalence 
\begin{equation} \label{e:rel ch alg 2 +}
\left(\Rep(\cG)\otimes \IndCoh^*(\Op^\mer_\cG)\right)^{>-\infty}\overset{\sim}\to 
\left(\CO_{\Op^\reg_\cG}\mod^{\on{fact}}(\Rep(\cG))\right){}^{>-\infty}.
\end{equation}
 
\medskip

Note also that $R_{\cG,\Op}\in \on{FactAlg}^{\on{untl}}(X,\Rep(\cG))$ belongs to
$$\Rep(\cG)^\heartsuit\subset \Rep(\cG)^{>-\infty}.$$

Hence, it makes sense to consider the full subcategory
$$\left(R_{\cG,\Op}\mod^{\on{fact}}\left(\CO_{\Op^\reg_\cG}\mod^{\on{fact}}(\Rep(\cG))\right)\right)^{>-\infty}
\subset R_{\cG,\Op}\mod^{\on{fact}}\left(\CO_{\Op^\reg_\cG}\mod^{\on{fact}}(\Rep(\cG))\right).$$

\medskip

It follows formally from the equivalence \eqref{e:rel ch alg 2 +} that the functor \eqref{e:rel ch alg 3} induces an 
equivalence
\begin{multline} \label{e:rel ch alg 5}
\left(R^{\Rep}_{\cG,\Op}\mod^{\on{fact}}(\Rep(\cG)\otimes \IndCoh^*(\Op^\mer_\cG))\right)^{>-\infty}\to \\
\to \left(R_{\cG,\Op}\mod^{\on{fact}}\left(\CO_{\Op^\reg_\cG}\mod^{\on{fact}}(\Rep(\cG))\right)\right)^{>-\infty}.
\end{multline}

\sssec{}

Hence, it suffices to show that the functor \eqref{e:rel ch alg 4} induces an equivalence between the eventually
coconnective subcategories of the two sides. 

\medskip

However, we claim that the functor \eqref{e:rel ch alg 4} it itself 
an equivalence. Indeed, this is a particular case of \lemref{l:two algebras mods}. 

\qed[\propref{p:FLE and Sat App}]

\ssec{Arc spaces of smooth D-schemes} \label{ss:arcs sm}

In this subsection, we will prove the following result.

\begin{prop}
	\label{p:arc space for smooth}
	Let $\CY \to X$ be a smooth affine D-scheme over $X$. Then for any test scheme $S$ and a map $\ul{x}:S\to \Ran$, 
	the $S$-scheme $\fL_\nabla^+(\CY)_{\ul{x}}$ is isomorphic to the limit of a sequence of smooth affine $S$-schemes.
\end{prop}

\sssec{}
To prove \propref{p:arc space for smooth}, it is enough to treat the case for the canonical map $X^I \to \Ran$. 
So our goal will be to prove the following: 

\begin{prop}
	\label{p:arc space for smooth XI}
	Let $\CY \to X$ be a smooth D-scheme over $X$. Then $\fL_\nabla^+(\CY)_{X^I}$ is isomorphic to the limit of a sequence of 
	relative smooth affine schemes over $X^I$.

\end{prop}

We will now describe $\fL_\nabla^+(\CY)_{X^I}$ as a limit of affine blow-ups $\fL_\nabla^+(\CY)_{X^I}^n$ (see below) for \emph{any} $\CY$. 
When $\CY$ is smooth, we will show that all $\fL_\nabla^+(\CY)_{X^I}^n$ are smooth.

\sssec{}
\label{sss: affine blowup}
We first give a brief review of the classical theory of affine blow-ups (a.k.a. dilations). From now on, we only work in \emph{classical} algebraic 
geometry, i.e., schemes means classical schemes, and fiber products of schemes mean non-derived fiber products, etc.

\medskip

Let $S$ be a smooth scheme and $E$ be an effective Cartier divisor on $S$. Let $Z$ be a $E$-regular $S$-scheme, i.e., the closed subscheme 
$Z_E:=Z\underset{S}\times E$ is an effective Cartier divisor on $Z$. Let $V$ be a closed subscheme of $Z_E$. The \emph{affine blow-up} of $Z$ with center $V$ (with respect to $E$) is defined to be
\[
	\on{Dil}_V^E(Z) := \on{Bl}_V(Z) - \wt{Z_E},
\]
where $\on{Bl}_V(Z)$ is the blow-up of $Z$ with center $V$, and $\wt{Z_E}$ is the strict transform of $Z_E$.

\sssec{} \label{sss:blow-up expl}

More explicitly, let $\CO_S(-E) \subset \CO_S$ be the ideal sheaf defining $E$, and $\CI \subset \CO_Z$ be 
the ideal sheaf defining $V$. Note that $\CO_Z(-E):= \CO_Z \underset{\CO_S}\otimes \CO_S(-E)$ is a 
subsheaf of $\CI$. Consider the inductive colimit
\[
	\mathcal{D}il_\CI^E(\CO_Z):= \on{colim} ( \CO_Z \to \CI(E) \to \CI^2(2E) \to \cdots).
\]  
Note that the connecting morphisms are injective because $Z$ is $E$-regular. This is a quasi-coherent $\CO_Z$-algebra with 
multiplication defined in the obvious way. More explicitly, if $E$ is locally cut out by a function $f$ of $\CO_S$, then 
$\mathcal{D}il_\CI^E(\CO_Z)$ is obtained from $\CO_Z$ by adding local sections 
$$f^{-n}\cdot a_1\cdots a_n, \quad a_k\in \CI.$$ We have
\[
	\on{Dil}_V^E(Z) \simeq \on{Spec}_{Z}( \mathcal{D}il_\CI^E(\CO_Z) ).
\]

\sssec{}

We observe:

\begin{lem}
	\label{lem: unipro affine blowup}
	In the setting of Sect. \ref{sss: affine blowup}, let $\on{Sch}_{S,E\mathrm{-reg}}$ be the category of $E$-regular $S$-schemes. Then $\on{Dil}_V^E(Z)$ represents the functor 
	\[
		\on{Sch}_{S,E\mathrm{-reg}}^{\on{op}} \to \on{Set},\; W\mapsto  \{ f: W\to Z\,\vert\, f(W_E) \subseteq V\}.
	\]
	Here $f(W_E) \subseteq V$ means the restriction of $f|_{W_E}$ factors through $V$.
\end{lem}

\proof
	Follows by unwinding the definitions. 

\qed

\sssec{}

Let $S^\circ \subset S$ be an open, whose complement has codimension $\geq 2$. Denote $Z^\circ:=Z\underset{S}\times S^\circ$.
Observe:

\begin{lem}  \label{l:affine closure}
Assume that $Z$ equals the affine closure of $Z^\circ$. Assume moreover that $V$ is the closure of
$V^\circ:=V\underset{S}\times S^\circ$ in $Z$. Then $\on{Dil}_V^E(Z)$ equals the affine closure of
$$\on{Dil}_V^E(Z)^\circ:=\on{Dil}_V^E(Z)\underset{S}\times S^\circ.$$
\end{lem}

\sssec{}

We first construct the schemes $\fL_\nabla^+(\CY)_{X^I}^n$ for $I=\{1,2\}$. 

\medskip

Let $\Delta_X$ be the diagonal of $X^2$ and consider the divisor $\Delta_{X,n}:=n\cdot \Delta_X$. 
Note that the connection on $\CY$ relative to $X$ provides a closed subscheme $\Delta_{\CY,n}$ of $\CY^2$: 
\begin{equation} \label{e:thick diagonal}
\Delta_{\CY,n}:= \CY_\nabla \underset{X_\dR}\times \Delta_{x,n},
\end{equation} 
viewed as a closed subscheme of $\CY^2_\nabla \times_{X_\dR^2} X^2$. 

\medskip

Set
$$\fL_\nabla^+(\CY)_{X^2}^n:=\on{Dil}^{n\cdot \Delta_X}_{\Delta_{\CY,n}}(\CY^2).$$

We claim that we have a canonical isomorphism
$$\fL_\nabla^+(\CY)_{X^2} \simeq \underset{n}{\on{lim}}\, \fL_\nabla^+(\CY)_{X^2}^n.$$

First, by the universal property of \lemref{lem: unipro affine blowup}, we have a canonically defined map 
\begin{equation} \label{e:arcs to limit}
\fL_\nabla^+(\CY)_{X^2} \to \underset{n}{\on{lim}}\, \fL_\nabla^+(\CY)_{X^2}^n.
\end{equation}

Let us show that this map is an isomorphism. Let $A$ be the commutative algebra in left D-modules
on $X$ corresponding to $\CY$. Choose a local coordinate $t$ on $X$, so that $X^2$ has coordinates
$t_1,t_2$.

\medskip

By \secref{sss:blow-up expl}, the algebra of functions on the right-hand side
in \eqref{e:arcs to limit} is the submodule of
$$j_*(A\boxtimes A)$$
that consists of local sections of the form
$$(t_1-t_2)^{-n}\cdot a, \quad a|_{n\cdot \Delta}\in \on{ker}(A\underset{\CO_X}\otimes A|_{n\cdot \Delta_X}\to A|_{n\cdot \Delta_X}),$$

where:

\begin{itemize}

\item The connection on $A$ allows us to extend it and also $A\underset{\CO_X}\otimes A$ to quasi-coherent sheaves on $n\cdot \Delta_X$,
denoted $A\underset{\CO_X}\otimes A|_{n\cdot \Delta_X}$ and $A|_{n\cdot \Delta_X}$, respectively.

\item We identify $A\boxtimes A|_{n\cdot \Delta_X}\simeq A\underset{\CO_X}\otimes A|_{n\cdot \Delta_X}$.

\end{itemize} 

However, this description coincides with the description of 
$$\CA_{X^2}\subset j_*(A\boxtimes A), \quad \CA=\on{Fact}(\CA).$$

\sssec{}

We now construct $\fL_\nabla^+(\CY)_{X^I}^n$ for an arbitrary finite set $I$. Let $\Delta_{X,I}$ be the diagonal divisor in 
$X^I$, i.e., the sum of the pairwise diagonals. We define the subscheme
$$\Delta_{\CY,I,n}\subset \CY^I$$
over $n\cdot \Delta_{X,I}$
as follows:

\medskip

Let 
$$X^{I,\circ}\overset{j_{I,\circ}}\hookrightarrow X^I$$ be the
open corresponding to the condition that not more than two coordinates coincide (i.e., we remove diagonals
of codimension $\geq 2$). Denote
$$\CY^{I,\circ}:=\CY^I\underset{X^I}\times X^{I,\circ}.$$

\medskip

Formula \eqref{e:thick diagonal} defines a subscheme
$$\Delta^{\circ}_{\CY,I,n}\subset \CY^{I,\circ}.$$

We let $\Delta_{\CY,I,n}$ be the closure of $\Delta^{\circ}_{\CY,I,n}$ in $\CY^I$.

\sssec{}

Set
$$\fL_\nabla^+(\CY)_{X^I}^n:=\on{Dil}^{n\cdot \Delta_{X,I}}_{\Delta_{\CY,I,n}}(\CY^I).$$

By the the universal property of \lemref{lem: unipro affine blowup}, we have a map
\begin{equation} \label{e:arcs to limit I}
\fL_\nabla^+(\CY)_{X^I} \to \underset{n}{\on{lim}}\, \fL_\nabla^+(\CY)_{X^I}^n.
\end{equation}

We claim that \eqref{e:arcs to limit I} is an isomorphism. 

\medskip

Indeed, the fact that \eqref{e:arcs to limit I} is an isomorphism over $X^{I,\circ}$ follows from 
the case $I=\{1,2\}$, considered above. We now claim that both sides in \eqref{e:arcs to limit I}
are affine closures of their respective restrictions to $X^{I,\circ}$. 

\medskip

Indeed, for the right-hand side, this follows from \lemref{l:affine closure}. For the left-hand side, 
this follows from the fact that the map
$$\CA_{X^I}\to (j_{I,\circ})_*\circ (j_{I,\circ})^*(\CA_{X^I})$$
induces an isomorphism on $H^0$. 

\sssec{}

We will now show that if $\CY$ is smooth, then the schemes $\fL_\nabla^+(\CY)_{X^I}^n$ are smooth. 
Note that the assertion is invariant under \emph{formal isomorphisms}. This allows us to replace 
$\CY$ by $\BA^k$. 

\medskip

To prove the smoothness, we can assume $k=1$. Thus, from now on, we will consider the schemes 
$\fL_\nabla^+(\BA^1)_{X^I}^n$. 

\sssec{}

Let $\sH_I\subset X^I\times X$ be the incidence divisor. Consider the correspospondence
$$
\CD
n\cdot \sH_I @>>> X \\
@VVV \\
X^I. 
\endCD
$$

For an $X$-scheme $\CY$, let 
$$\on{Jets}^n_I(\CY)$$
be the scheme of $n$-truncated jets into $\CY$, i.e., the restriction of scalars \'a la Weil of $(n\cdot \sH_I)\underset{X}\times \CY$
along the map $n\cdot \sH_I\to X^I$.

\sssec{}

Consider the schemes 
$$\on{Jets}^n_I(\bA^1) \text{ and } \on{Jets}^{n-1}_I((\bA^1)_{\omega_X}),$$
where $(\bA^1)_{\omega_X}$ is the total space of the line bundle of 1-forms on $X$.

\medskip

Both schemes are vector bundles over $X^I$ of ranks $|I|\cdot n$ and $|I|\cdot (n-1)$, 
respectively. De Rham differential defines a map
\begin{equation} \label{e:jet 1 form} 
d_{\dr}:\on{Jets}^n_I(\bA^1)\to \on{Jets}^{n-1}_I((\bA^1)_{\omega_X}).
\end{equation} 

Set
$$'\fL_\nabla^+(\BA^1)_{X^I}^n:=\on{ker}(d_\dr).$$

Since $d_\dr$ is fiber-wise surjective as a map of vector bundles, $'\fL_\nabla^+(\BA^1)_{X^I}^n$ is
also a vector bundle. In particular, it is smooth over $X^I$. 

\sssec{}

We are going to prove:

\begin{prop} \label{p:acrs as jets n}
There exists a canonical isomorphism
$$'\fL_\nabla^+(\BA^1)_{X^I}^n \simeq \fL_\nabla^+(\BA^1)_{X^I}^n.$$
\end{prop}

\ssec{Proof of \propref{p:acrs as jets n}}

\sssec{}

First, we note that when $I=\{*\}$, there is an obvious isomorphism
\begin{multline*}
'\fL_\nabla^+(\BA^1)_{X}^n=
\on{ker}(\on{Jets}^n(\bA^1)\to \on{Jets}^{n-1}((\bA^1)_{\omega_X}) \overset{\sim}\leftarrow
\on{ker}(\on{Jets}(\bA^1)\to \on{Jets}((\bA^1)_{\omega_X}) \simeq \\
\simeq \fL_\nabla^+(\BA^1)_X=\BA^1\times X.
\end{multline*}

\sssec{}

Note that for a map of finite sets $I'\to I$, there is a canonically defined map 
$$\on{Jets}^n_I(\CY)\to \on{Jets}^n_{I'}(\CY),$$
covering $X^I\to X^{I'}$.

\medskip

In particular, we obtain a map
\begin{equation} \label{e:acrs as jets n step 0}
'\fL_\nabla^+(\BA^1)_{X^I}^n \to ({}'\fL_\nabla^+(\BA^1)_{X}^n)^I\simeq (\BA^1\times X)^I
\end{equation}
over $X^I$. This map is an isomorphism away from the diagonal divisor. 

\sssec{}

We will show that the map \eqref{e:acrs as jets n step 0} lifts to and defines an isomorphism
\begin{equation} \label{e:acrs as jets n step 1}
'\fL_\nabla^+(\BA^1)_{X^I}^n \to \fL_\nabla^+(\BA^1)_{X^I}^n
\end{equation}
\emph{over} $X^{I,\circ}$.

\medskip

This would imply \propref{p:acrs as jets n}, since both schemes in question are affine closures
of their respective restrictions to $X^{I,\circ}$.

\sssec{}

In order to prove \eqref{e:acrs as jets n step 1}, it suffices to consider the case $I=\{1,2\}$. The following
is obtained by a straightforward calculation:

\begin{lem}  \label{l:jets to arcs}
The map 
\begin{equation} \label{e:acrs as jets n step 0 sq}
'\fL_\nabla^+(\BA^1)_{X^2}^n \to (\BA^1\times X)^2
\end{equation}
has the following properties:

\smallskip

\noindent{\em(i)} Its restriction to $n\cdot \Delta_X$ maps to $\Delta_{\bA^1,n}$;

\smallskip

\noindent{\em(ii)} Its restriction to $(n+1)\cdot \Delta_X$ \emph{does not} map to $\Delta_{\bA^1,n+1}$.

\end{lem}

\sssec{}

From Lemmas \ref{l:jets to arcs}(i) and \eqref{lem: unipro affine blowup}, we obtain that \eqref{e:acrs as jets n step 0 sq} lifts to a map
\begin{equation} \label{e:acrs as jets n step 1 sq}
'\fL_\nabla^+(\BA^1)_{X^2}^n \to \fL_\nabla^+(\BA^2)_{X^I}^n.
\end{equation}

It remains to show that \eqref{e:acrs as jets n step 1 sq} is an isomorphism. 

\sssec{}

The map \eqref{e:acrs as jets n step 0 sq} respects the vector bundle structures. Hence, so does the map 
\eqref{e:acrs as jets n step 1 sq}. It is an isomorphism away from the diagonal, so it remains to show that it 
is surjective over the diagonal.

\sssec{}

It is clear that the image of 
$$'\fL_\nabla^+(\BA^1)_{X^2}^n|_{\Delta_X}\to  \fL_\nabla^+(\BA^2)_{X^I}^n|_{\Delta_X}\simeq \BA^1\times \BA^1$$
contains the diagonal copy of $\BA^1$.

\medskip

Furthermore, if the entire image \emph{had been} contained in the diagonal copy of $\BA^1$, it would have meant that
the map \eqref{e:acrs as jets n step 1 sq} lifts further to a map
$$'\fL_\nabla^+(\BA^1)_{X^2}^n \to \fL_\nabla^+(\BA^2)_{X^I}^{n+1}.$$

However, the latter contradicts point (ii) of \lemref{l:jets to arcs}.

\qed[\propref{p:acrs as jets n}]

\section{Horizontal sections of affine D-schemes} \label{s:hor sect D sch}

Let $\CY$ be an affine  D-prestack over $X$. Let $\on{Sect}_\nabla(X,\CY)$ denote the prestack of 
its horizontal sections, i.e.,
$$\Maps(\Spec(R),\on{Sect}_\nabla(X,\CY)):=
\Maps_{X,\nabla}(\Spec(R)\times X,\CY).$$

In this section we will describe, following \cite{BD2}, $\on{Sect}_\nabla(X,\CY)$ explicitly in terms of (vacuum) factorization
homology of the algebra of functions $\CO_\CY$ of $\CY$. We also describe spaces of sections of quasi-coherent sheaves on $\on{Sect}_\nabla(X,\CY)$,
in terms of factorization homology of $\CO_\CY$ with coefficients in corresponding modules. 

\medskip

We then generalize this discussion, when instead of $\on{Sect}_\nabla(X,\CY)$ we consider the space  $\on{Sect}_\nabla(X-\ul{x},\CY)$
of sections, where we allow punctures at $\ul{x}\subset X$. 

\ssec{Horizontal sections via factorization homology}

\sssec{} 

We start with an arbitrary D-prestack $\CY\to X$ (later on in this section, we will assume that $\CY$ is affine over $X$).

\medskip

Let $A\in \on{ComAlg}(\Dmod(X))$ be the direct image 
of the structure sheaf of $\CY$, and let $\CA$ be the corresponding object in $\on{ComAlg}(\on{FactAlg}(X))$,
i.e., $\CA=\on{Fact}(A)$. 

\medskip

Consider the evaluation map 
$$\on{Sect}_\nabla(X,\CY)\times X\to \CY.$$

It gives rise to a map 
\begin{equation} \label{e:alg to glob sect}
A \to \CO_{\on{Sect}_\nabla(X,\CY)}\otimes \CO_X
\end{equation}
in $\on{ComAlg}(\Dmod(X))$, where by a slight abuse of notation we denote by $\CO_{\on{Sect}_\nabla(X,\CY)}$
the algebra of global functions on $\on{Sect}_\nabla(X,\CY)$.

\medskip 

By the adjunction of \corref{c:fact com untl to X bis}, from \eqref{e:alg to glob sect} we obtain a map 
\begin{equation} \label{e:A on X}
\CA \to \CO_{\on{Sect}_\nabla(X,\CY)}\otimes \omega_{\Ran^{\on{untl},*}}
\end{equation}
in $\on{ComAlg}(\Dmod(\Ran^{\on{untl},*}))$, and further, by \secref{sss:left adj on untl alg}, a map
\begin{equation} \label{e:fact hom A on Ran}
\on{C}^{\on{fact}}_\cdot(X,\CA)\to \CO_{\on{Sect}_\nabla(X,\CY)}
\end{equation}
in $\on{ComAlg}(\Vect)$. 

\sssec{} 

Suppose that the prestack $\CY$ satisfies the assumption of \secref{sss:com fact and arcs gen}
(e.g., $\CY\to X$ is affine). In this case, the map \eqref{e:A on X} can also be interpreted as follows:

\medskip

Consider the map 
$$\on{Sect}_\nabla(X,\CY)\times \Ran \overset{\on{ev}}\to \fL_\nabla^+(\CY)_\Ran.$$

Pullback at the level of functions defines a map
\begin{equation} \label{e:A on Ran}
\CA \overset{\text{\eqref{e:arcs as fact}}}\simeq \CO_{\fL_\nabla^+(\CY),\Ran}\to \CO_{\on{Sect}_\nabla(X,\CY)}\otimes \omega_\Ran
\end{equation}
in $\on{ComAlg}(\Dmod(\Ran))$.

\medskip

It follows by unwinding the definition that the map \eqref{e:A on X} is the same as \eqref{e:A on Ran}. 

\sssec{} 

Suppose now that $\CY\to X$ is affine. In this case we claim:

\begin{prop} \label{p:hor sect}
The prestack $\on{Sect}_\nabla(X,\CY)$ is an affine scheme, and 
the map \eqref{e:fact hom A on Ran} is an isomorphism.
\end{prop}

\begin{proof}

The key fact is that for any factorization algebra $\CA$ such that
$\oblv^l(\CA_X)\in \QCoh(X)^{\leq 0}$, we have
$$\on{C}^{\on{fact}}_\cdot(X,\CA)\in \Vect^{\leq 0}.$$

This follows from the fact that 
$$\on{C}^{\on{fact}}_\cdot(X,\CA)\simeq \on{C}^\cdot_c(\Ran,\CA_\Ran)$$
can be written as a colimit with terms 
$$\on{C}^\cdot(X^I,\CA_{X^I})$$
for (non-empty) finite sets $I$, while each $\CA_{X^I}$ satisfies $\oblv^l(\CA_{X^I})\in \QCoh(X^I)^{\leq 0}$
(see \cite[Sect. 3.4.11]{BD2}) and hence $\on{C}^\cdot(X^I,\CA_{X^I})\in \Vect^{\leq 0}$. 

\medskip

In particular, in our case
$$\on{C}^{\on{fact}}_\cdot(X,\CA)\in \on{ComAlg}(\Vect^{\leq 0}).$$

Now, the assertion follows immediately from \corref{c:ch homology as left adj}: for $R\in \on{ComAlg}(\Vect^{\leq 0})$
we have
\begin{multline} \label{e:horiz sections functor}
\Maps(\Spec(R),\on{Sect}_\nabla(X,\CY))=\Maps_{X,\nabla}(\Spec(R)\times X,\CY)
\overset{\on{affineness\,of\,}\CY}\simeq \\
\simeq \Maps_{\on{ComAlg}(\Dmod(X))}(A,R\otimes \CO_X)\simeq
\Maps_{\on{ComAlg}(\Vect)}(\on{C}^{\on{fact}}_\cdot(X,\CA),R),
\end{multline}
so $\on{Sect}_\nabla(X,\CY)$ is the affine scheme $\Spec(\on{C}^{\on{fact}}_\cdot(X,\CA))$, and by construction
the map \eqref{e:fact hom A on Ran} is the map
$$\on{C}^{\on{fact}}_\cdot(X,\CA)\to \Gamma(\Spec(\on{C}^{\on{fact}}_\cdot(X,\CA)),\CO_{\Spec(\on{C}^{\on{fact}}_\cdot(X,\CA))})$$
resulting from \eqref{e:horiz sections functor}.

\end{proof} 

\sssec{}

Recall now that the functor 
$$\Gamma(\fL^+_\nabla(\CY)_{\ul{x}},-):\QCoh(\fL^+_\nabla(\CY)_{\ul{x}})\to \Vect$$
factors via an equivalence
$$\Gamma(\fL^+_\nabla(\CY)_{\ul{x}},-)^{\on{enh}}:\QCoh(\fL^+_\nabla(\CY)_{\ul{x}})\to \CA_{\ul{x}}\mod\simeq \CA\mod^{\on{com}}_x,$$
followed by the forgetful functor 
$$\oblv_{\CA,\ul{x}}:\CA\mod^{\on{com}}_{\ul{x}}\to \Vect.$$

\medskip

We claim: 

\begin{prop} \label{p:global section modules +}
There exists a canonical isomorphism between 
\begin{equation} \label{e:global section modules 5}
\QCoh(\fL^+_\nabla(\CY)_{\ul{x}}) \overset{\on{ev}_{\ul{x}}^*}\longrightarrow \QCoh(\on{Sect}_\nabla(X,\CY))
\overset{\Gamma(\on{Sect}_\nabla(X,\CY),-)}\longrightarrow \Vect
\end{equation} 
and
\begin{equation} \label{e:global section modules 6}
\QCoh(\fL^+_\nabla(\CY)_{\ul{x}}) \overset{\Gamma(\fL^+_\nabla(\CY)_{\ul{x}},-)^{\on{enh}}}\simeq
\CA\mod^{\on{com}}_{\ul{x}} \to \CA\mod^{\on{fact}}_{\ul{x}} \overset{\on{C}^{\on{fact}}_\cdot(X,\CA,-)_{\ul{x}}}\longrightarrow \Vect.
\end{equation} 
\end{prop} 

\begin{proof}

Since 
$$\QCoh(\fL^+_\nabla(\CY)_{\ul{x}}) \simeq \CA_{\ul{x}}\mod\simeq \CA\mod^{\on{com}}_{\ul{x}},$$
the assertion of the proposition amounts to the following:

\medskip

We have a canonical identification 
$$\Gamma(\on{Sect}_\nabla(X,\CY),\CO_{\on{Sect}_\nabla(X,\CY)}) \simeq \on{C}^{\on{fact}}_\cdot(X,\CA,\CA_{\ul{x}})_{\ul{x}}$$
as $\CA_{\ul{x}}$-modules, where:

\begin{itemize}

\item $\CA_{\ul{x}}$ acts on the left-hand side via
$$\CA_x\overset{\sim}\to \End(\CO_{\fL^+_\nabla(\CY)_{\ul{x}}}) \overset{\on{ev}_{\ul{x}}^*}\longrightarrow 
\End(\on{ev}_{\ul{x}}^*(\CO_{\fL^+_\nabla(\CY)_{\ul{x}}}))
\simeq \End(\CO_{\on{Sect}_\nabla(X,\CY)});$$

\item $\CA_{\ul{x}}$ acts on the right-hand side via 
$$\CA_x\overset{\sim}\to \End_{\CA\mod^{\on{com}}_{\ul{x}}}(\CA_x)\to \End_{\CA\mod^{\on{fact}}_{\ul{x}}}(\CA_{\ul{x}}).$$

\end{itemize}

\medskip

However, this follows from the fact that the isomorphism
$$\Gamma(\on{Sect}_\nabla(X,\CY),\CO_{\on{Sect}_\nabla(X,\CY)}) \overset{\text{\eqref{e:fact hom A on Ran}}}
\simeq \on{C}^{\on{fact}}_\cdot(X,\CA)\overset{\text{unitality}}\simeq  \on{C}^{\on{fact}}_\cdot(X,\CA,\CA_{\ul{x}})_{\ul{x}}$$
(as commutative algebras) is, by construction, compatible with the homomorphisms
$$\CA_{\ul{x}}\simeq \Gamma(\fL^+_\nabla(\CY)_{\ul{x}},\CO_{\fL^+_\nabla(\CY)_{\ul{x}}})\to 
\Gamma(\on{Sect}_\nabla(X,\CY),\CO_{\on{Sect}_\nabla(X,\CY)})$$
and 
$$\CA_{\ul{x}} \to \on{C}_c(\Ran^{\on{untl}}_{\ul{x}},\CA_{\Ran^{\on{untl}}_{\ul{x}}}) = \on{C}^{\on{fact}}_\cdot(X,\CA,\CA_{\ul{x}})_{\ul{x}}.$$

\end{proof} 

\begin{rem} \label{r:indep aff}

Recall that Remark \ref{r:indep com} says that for a commutative factorization algebra $\CA$, we have a canonical identification
$$\CA\mod^{\on{com}}_{\Ran,\on{indep}}\simeq \on{C}^\cdot(X,\CA)\mod.$$

Hence, by \propref{p:hor sect} we obtain:
$$\QCoh(\fL_\nabla^+(\CY))_{\Ran,\on{indep}}\simeq \QCoh(\on{Sect}_\nabla(X,\CY)).$$

We can interpret \propref{p:global section modules +} as saying that under the above equivalence, the corresponding functor
$$\QCoh(\fL_\nabla^+(\CY))\to  \QCoh(\on{Sect}_\nabla(X,\CY))\otimes \ul\Dmod(\Ran^{\on{untl}})$$
sends $\ul{x}\in \Ran$ to to the functor
$$\QCoh(\fL_\nabla^+(\CY)_{\ul{x}}) \overset{\on{ev}^*_{\ul{x}}}\longrightarrow \QCoh(\on{Sect}_\nabla(X,\CY)).$$

\end{rem} 

\begin{rem} 

The above remark applies to D-prestacks that are ``as good as affine", see \secref{ss:as good as affine}. For example, for
a unipotent group-scheme $N'$ over $X$ and $\CY=\on{pt}/\fL^+(N')$, we obtain an equivalence 
\begin{equation} \label{e:indep Rep L^+N}
\Rep(\on{pt}/\fL^+(N'))_{\Ran,\on{indep}}\simeq \QCoh(\Bun_{N'}),
\end{equation} 
where the composite functor
$$\Rep(\on{pt}/\fL^+(N'))_\Ran\twoheadrightarrow \Rep(\on{pt}/\fL^+(N'))_{\Ran,\on{indep}}\simeq \QCoh(\Bun_{N'})$$
is $\Loc^{\QCoh}_{N'}$.

\medskip

Similarly, if $N'$ is equipped with a connection (e.g., $N'$ is constant), we have
\begin{equation} \label{e:indep Rep N}
\Rep(N')_{\Ran,\on{indep}} \simeq \QCoh(\LS_{N'}),
\end{equation} 
where the composite functor 
$$\Rep(N')_\Ran \twoheadrightarrow \Rep(N')_{\Ran,\on{indep}} \simeq \QCoh(\LS_{N'})$$
is $\Loc^{\on{spec}}_{N'}$. 

\end{rem} 

\begin{rem}
We warn the reader that the equivalences \eqref{e:indep Rep L^+N} and \eqref{e:indep Rep N}
are a feature of unipotent group schemes. For a reductive $G$ (e.g., for $G=\BG_m$), 
the corresponding functors are far from being equivalences. 
\end{rem}

\ssec{Variant: allowing poles} \label{ss:sections with poles}

Let $A$ and $\CY$ be as in the previous subsection. For $\ul{x}\in \Ran$, consider the prestack
$$\on{Sect}_\nabla(X-\ul{x},\CY).$$

\medskip

We will now give an explicit description of $\on{Sect}_\nabla(X-\ul{x},\CY)$ in terms of factorization homology.
In particular, we will show that it is an ind-affine ind-scheme. 

\sssec{} \label{sss:glob poles as colimit}

Write 
$$\fL_\nabla(\CY)_{\ul{x}}\simeq \underset{R}{``\on{colim}"}\, \Spec(R'), \quad R'\in \on{ComAlg}(\Vect^{\leq 0}),$$
where each
\begin{equation} \label{e:local mr sect}
\Spec(R')\to \fL_\nabla(\CY)_{\ul{x}}
\end{equation} 
is a closed embedding.

\medskip

Consider the fiber product
$$\on{Sect}_\nabla(X-\ul{x},\CY)':=\Spec(R') \underset{\fL(\CY)_{\ul{x}}}\times \on{Sect}_\nabla(X-\ul{x},\CY).$$

We will show that $\on{Sect}_\nabla(X-\ul{x},\CY)'$ is a scheme and describe the algebra of functions on it
in terms of factorization homology. 

\sssec{}

First, we note that as in \secref{sss:modif as pairs} the datum of a closed embedding
\begin{equation} \label{e:map to mer Sect}
\Spec(R') \to  \fL_\nabla(\CY)_{\ul{x}}
\end{equation}
is equivalent to the datum of a modification $A'$ of $A$ at $\ul{x}$, i.e., $A'$ is an object of $\on{ComAlg}(\Dmod(X)^{\leq 0})$
equipped with an isomorphism 
$$A'|_{X-\ul{x}}\simeq A|_{X-\ul{x}}.$$

Indeed, given $A'$, set 
\begin{equation} \label{e:modified D-scheme}
\CY':=\Spec_X(A'), 
\end{equation}
and we recover $R'$ as $\CA'_{\ul{x}}$, where $\CA':=\on{Fact}(A')$, 
and \eqref{e:map to mer Sect} as 
$$\Spec(R')\simeq \fL^+_\nabla(\CY')\hookrightarrow \fL_\nabla(\CY')\simeq \fL_\nabla(\CY').$$

\sssec{}

Vice versa, given \eqref{e:map to mer Sect}, we interpret it as a map of commutative D-algebras 
$$j_*\circ j^*(A)\to R'\ppart$$
(here $t$ is a coordinate on the multidisc $\cD_{\ul{x}}$
and $j$ is the open embedding $X-\ul{x}\hookrightarrow X$), where the 
condition that \eqref{e:map to mer Sect} is a closed embedding corresponds to the condition that the map of
D-modules
\begin{equation} \label{e:recover A'}
j_*\circ j^*(A)\to R'\ppart\to R'\ppart/R'\qqart
\end{equation}
is surjective on $H^0$. 

\medskip

We recover $A'$ as the fiber product
$$j_*\circ j^*(A)\underset{R'\ppart}\times R\qqart,$$
and the surjectivity of \eqref{e:recover A'} is equivalent to the condition that $A'$ is connective. 

\sssec{}

Let $R'$ be as in \eqref{e:map to mer Sect}, and let $A'$ be the corresponding modification of $A$. 
Unwinding the definitions, we obtain:
$$\Spec(R') \underset{\fL(\CY)_{\ul{x}}}\times \on{Sect}_\nabla(X-\ul{x},\CY) \simeq
\on{Sect}_\nabla(X-x,\CY').$$

Hence, by \propref{p:hor sect}, we obtain:

\begin{cor} \label{c:hor sect mer}
The prestack $\Spec(R') \underset{\fL(\CY)_{\ul{x}}}\times \on{Sect}_\nabla(X-x,\CY)$ is affine and the naturally defined map
$$\on{C}_\cdot^{\on{fact}}(X,\CA')\to \CO_{\Spec(R') \underset{\fL(\CY)_{\ul{x}}}\times \on{Sect}_\nabla(X-\ul{x},\CY)}$$
is an isomorphism, where $\CA':=\on{Fact}(A')$.
\end{cor}
   
\sssec{} \label{sss:modified vacuum}

We will now rewrite $\on{C}_\cdot^{\on{fact}}(X,\CA')$ in terms of factorization homology of $\CA$ itself. 

\medskip

We can view $R'$ as an object of $\CA\mod^{\on{fact}}_{\ul{x}}$ via
$$R'\in \CA'\mod^{\on{com}}_{\ul{x}}\to \CA'\mod^{\on{fact}}_{\ul{x}}\simeq \CA\mod^{\on{fact}}_{\ul{x}}.$$

\medskip

When viewed as such, we will denote it by $R'_\CA$. Let $R'_{\CA,\Ran^{\on{untl}}_{\ul{x}}}$ be the corresponding object of 
$\Dmod(\Ran^{\on{untl}}_{\ul{x}})$. 

\medskip

Recall (see \secref{sss:fact mods for com}) that since $\CA$ is a \emph{commutative} factorization algebra, 
the category $\CA\mod^{\on{fact}}_{\ul{x}}$ has a natural symmetric pseudo-monoidal structure. The commutative algebra structure on $R'$ endows 
$R'_\CA$ with a structure of commutative algebra object in $\CA\mod^{\on{fact}}_{\ul{x}}$. In particular, 
$R'_{\CA,\Ran^{\on{untl}}_{\ul{x}}}$
acquires a structure of commutative algebra object in $\Dmod(\Ran^{\on{untl}}_{\ul{x}})$.

\sssec{}

The assertion of \lemref{l:monoidal on Ran untl} is valid for $\Ran^{\on{untl}}_{\ul{x}}$, hence
$$\on{C}^{\on{fact}}_\cdot(X,\CA,R'_\CA)_{\ul{x}}:=
\on{C}^\cdot_c(\Ran_{\ul{x}}^{\on{untl}},R'_{\CA,\Ran^{\on{untl}}_{\ul{x}}})$$
acquires a structure of commutative algebra in $\Vect$.

\medskip

We claim:

\begin{lem}  \label{l:fact hom R}
There is a canonical isomorphism
\begin{equation} \label{e:fact hom R}
\on{C}^{\on{fact}}_\cdot(X,\CA,R'_\CA)_{\ul{x}} \simeq \on{C}_\cdot^{\on{fact}}(X,\CA')
\end{equation}
as objects of $\on{ComAlg}(\Vect)$. 
\end{lem} 

\begin{proof} 

Note that with respect to the (symmetric monoidal) equivalence
$$\CA\mod^{\on{fact}}_{\ul{x}}\simeq \CA'\mod^{\on{fact}}_{\ul{x}},$$
the object
$$R'_\CA\in \CA\mod^{\on{fact}}_{\ul{x}}$$
corresponds to the vacuum object 
$$(\CA')^{\on{fact}_{\ul{x}}}\in \CA'\mod^{\on{fact}}_{\ul{x}},$$
where we recall that $(\CA')^{\on{fact}_{\ul{x}}}$ denotes the vacuum factorization module at $\ul{x}$.

\medskip

Moreover, this isomorphism respects the structure of commutative algebra object on both sides. 
In particular, 
$$R'_{\CA,\Ran^{\on{untl}}_{\ul{x}}}\simeq \CA'_{\Ran^{\on{untl}}_{\ul{x}}}$$
as objects of $\on{ComAlg}(\Dmod(\Ran^{\on{untl}}_{\ul{x}}))$. 

\medskip

Now, \eqref{e:fact hom R} follows by concatenating 
$$\on{C}^{\on{fact}}_\cdot(X,\CA,R'_\CA)_{\ul{x}} \simeq \on{C}^{\on{fact}}_\cdot(X,\CA',(\CA')^{\on{fact}_{\ul{x}}})_{\ul{x}},$$
and 
$$\on{C}^{\on{fact}}_\cdot(X,\CA',(\CA')^{\on{fact}_{\ul{x}}})_{\ul{x}}\overset{\text{unitality}}\simeq \on{C}^{\on{fact}}_\cdot(X,\CA').$$

\end{proof} 

\sssec{}

Thus, $\Spec(\on{C}^{\on{fact}}_\cdot(X,\CA,R'_\CA)_{\ul{x}})$ gives the desired expression for 
$$\Spec(R') \underset{\fL(\CY)_{\ul{x}}}\times \on{Sect}_\nabla(X-x,\CY)$$
in terms of factorization homology of $\CA$. 

\sssec{}

Let 
$$R'_{\CA,\Ran_{\ul{x}}}\in \on{ComAlg}(\Dmod(\Ran_{\ul{x}}))$$
be the restriction of $R'_{\CA,\Ran^{\on{untl}}_{\ul{x}}}$ to the non-unital Ran space.

\medskip

By construction, it belongs to $\on{ComAlg}(\Dmod(\Ran_{\ul{x}})^{\on{almost-untl}})$.
Hence, by \secref{sss:almost unital},
$$\on{C}^\cdot_c(\Ran_{\ul{x}},R'_{\CA,\Ran_{\ul{x}}})$$
acquires a commutative algebra structure.

\medskip

By \lemref{l:int over Ran and Ran untl} we have an isomorphism
$$\on{C}^\cdot_c(\Ran_{\ul{x}},R'_{\CA,\Ran_{\ul{x}}})\simeq \on{C}^\cdot_c(\Ran^{\on{untl}}_{\ul{x}},R'_{\CA,\Ran^{\on{untl}}_{\ul{x}}})=
\on{C}^{\on{fact}}_\cdot(X,\CA,R')_{\ul{x}}$$
as commutative algebras. 

\ssec{Factorization with poles}

In his subsection we will add a slightly different spin to the discussion in \secref{ss:sections with poles}. 

\sssec{}

As in \cite[Definition A.2.1]{Bogd}, let
 $\fL^{\mer\rightsquigarrow\reg}_\nabla(\CY)_{\Ran^{\subseteq}}$ be the factorization ind-scheme over $\Ran^{\subseteq}$
that attaches to 
$$(\ul{x}\subseteq \ul{x}')$$
the space
$$\fL^{\mer\rightsquigarrow\reg}_\nabla(\CY)_{\ul{x}\subseteq \ul{x}'}:=\on{Sect}_\nabla(\cD_{\ul{x}'}-\ul{x},\CY).$$

We have the projections
\begin{equation} \label{e:Op big and small}
\fL_\nabla(\CY) \overset{\on{pr}^\CY_{\on{small}}}\longleftarrow \fL^{\mer\rightsquigarrow\reg}_\nabla(\CY)_{\Ran^{\subseteq}}
\overset{\on{pr}^\CY_{\on{big}}}\longrightarrow \fL_\nabla(\CY)
\end{equation} 
given by restrictions along
$$\cD_{\ul{x}}-\ul{x}\hookrightarrow \cD_{\ul{x}'}-\ul{x}\hookleftarrow \cD_{\ul{x}'}-\ul{x}',$$
respectively.

\sssec{Example}

When 
$$\ul{x}'=\ul{x}\sqcup \ul{x}'',$$
we have 
$$\fL^{\mer\rightsquigarrow\reg}_\nabla(\CY)_{\ul{x}\subseteq \ul{x}'}\simeq \fL_\nabla(\CY)_{\ul{x}}\times \fL^+_\nabla(\CY)_{\ul{x}''}.$$

\sssec{} \label{sss:Y mer to reg}

Fix $\ul{x}$, and consider the space
$$\fL^{\mer\rightsquigarrow\reg}_\nabla(\CY)_{\Ran_{\ul{x}}}:=\fL^{\mer\rightsquigarrow\reg}_\nabla(\CY)_{\Ran^{\subseteq}}
\underset{\Ran^{\subseteq}}\times \Ran_{\ul{x}}.$$

It has a natural structure of factorization module space with respect to $\fL^+_\nabla(\CY)$, with the underlying space being
$\fL_\nabla(\CY)_{\ul{x}}$. 

\medskip

Moreover, viewed as such, $\fL^{\mer\rightsquigarrow\reg}_\nabla(\CY)_{\Ran_{\ul{x}}}$
has a natural \emph{counital} structure (see \secref{sss:counital spaces} for what this means).

\medskip

In particular, we have the map 
$$\on{pr}^\CY_{\on{small},\ul{x}}:\fL^{\mer\rightsquigarrow\reg}_\nabla(\CY)_{\Ran_{\ul{x}}}\to \fL_\nabla(\CY)_{\ul{x}}.$$

\sssec{}

Fix a closed embedding $\Spec(R')\to \fL_\nabla(\CY)_{\ul{x}}$, and consider the fiber product
$$\Spec(R')\underset{\fL_\nabla(\CY)_{\ul{x}}}\times \fL^{\mer\rightsquigarrow\reg}_\nabla(\CY)_{\Ran_{\ul{x}}}.$$

Unwinding the definitions, we obtain: 

\begin{lem} \label{l:mixed with poles}
We have a canonical isomorphism
$$\Spec(R')\underset{\fL_\nabla(\CY)_{\ul{x}}}\times \fL^{\mer\rightsquigarrow\reg}_\nabla(\CY)_{\Ran_{\ul{x}}}\simeq
\fL^+_\nabla(\CY')_{\Ran_{\ul{x}}},$$ where 
$\CY'$ is as in \eqref{e:modified D-scheme}.
\end{lem}

\ssec{Sections of quasi-coherent sheaves, meromorphic version}

\sssec{}

Recall (see \secref{sss:QCohco to fact mod}) that the functor 
$$\Gamma(\fL_\nabla(\CY)_{\ul{x}},-):\QCoh_{\on{co}}(\fL_\nabla(\CY)_{\ul{x}})\to \Vect$$
admits an an enhancement to a functor
$$\Gamma(\fL_\nabla(\CY)_{\ul{x}},-)^{\on{enh}}:\QCoh_{\on{co}}(\fL_\nabla(\CY)_{\ul{x}})\to \CA\mod^{\on{fact}}_{\ul{x}}.$$

\sssec{}

Note also that since the morphism
$$\on{ev}_{\ul{x}}:\on{Sect}_\nabla(X-\ul{x},\CY)\to \fL_\nabla(\CY)_{\ul{x}}$$
is schematic, we have a well-defined functor
$$\on{ev}_{\ul{x}}^*:\QCoh_{\on{co}}(\fL_\nabla(\CY)_{\ul{x}})\to \QCoh_{\on{co}}(\on{Sect}_\nabla(X-\ul{x},\CY)).$$

\sssec{}

The goal of this subsection is to prove the following:

\begin{prop} \label{p:global section modules}
There exists a canonical isomorphism between 
\begin{equation} \label{e:global section modules 1}
\QCoh_{\on{co}}(\fL_\nabla(\CY)_{\ul{x}}) \overset{\on{ev}_{\ul{x}}^*}\longrightarrow \QCoh_{\on{co}}(\on{Sect}_\nabla(X-\ul{x},\CY))
\overset{\Gamma(\on{Sect}_\nabla(X-\ul{x},\CY),-)}\longrightarrow \Vect
\end{equation} 
and
\begin{equation} \label{e:global section modules 2}
\QCoh_{\on{co}}(\fL_\nabla(\CY)_{\ul{x}}) \overset{\Gamma(\fL_\nabla(\CY)_{\ul{x}},-)^{\on{enh}}}\longrightarrow
\CA\mod^{\on{fact}}_{\ul{x}} \overset{\on{C}^{\on{fact}}_\cdot(X,\CA,-)_{\ul{x}}}\longrightarrow \Vect.
\end{equation} 
\end{prop} 

The rest of this subsection is devoted to the proof of \propref{p:global section modules}.

\sssec{}

First, we construct a natural transformation
$$\text{\eqref{e:global section modules 2}} \to \text{\eqref{e:global section modules 1}}.$$

Note that the functor
$$\on{ins.unit}_{\ul{x}}:\QCoh_{\on{co}}(\fL_\nabla(\CY)_{\ul{x}})\to \QCoh_{\on{co}}(\fL_\nabla(\CY)_{\Ran_{\ul{x}}})$$
is given by 
$$(\on{pr}^\CY_{\on{big},\ul{x}})_*\circ (\on{pr}^\CY_{\on{small},\ul{x}})^*,$$
where
$$\fL_\nabla(\CY)_{\ul{x}} \overset{\on{pr}^\CY_{\on{small},\ul{x}}}\longleftarrow \fL^{\mer\rightsquigarrow\reg}_\nabla(\CY)_{\Ran_{\ul{x}}}
\overset{\on{pr}^\CY_{\on{big},\ul{x}}}\longrightarrow \fL_\nabla(\CY)_{\Ran_{\ul{x}}}.$$

Hence, we obtain that the functor
$$\QCoh_{\on{co}}(\fL_\nabla(\CY)_{\ul{x}}) \overset{\Gamma(\fL_\nabla(\CY)_{\ul{x}},-)^{\on{enh}}}\longrightarrow
\CA\mod^{\on{fact}}_{\ul{x}} \to \Dmod(\Ran_{\ul{x}})$$
is given by
$$(p_{\Ran_{\ul{x}}})_*\circ (\on{pr}^\CY_{\on{small},\ul{x}})^*,$$
where 
$$p:\fL^{\mer\rightsquigarrow\reg}_\nabla(\CY)_{\Ran_{\ul{x}}}\to \Ran_{\ul{x}}.$$

\medskip

From here we obtain a natural transformation from 
\begin{equation} \label{e:global section modules 3}
\QCoh_{\on{co}}(\fL_\nabla(\CY)_{\ul{x}}) \overset{\Gamma(\fL_\nabla(\CY)_{\ul{x}},-)^{\on{enh}}}\longrightarrow
\CA\mod^{\on{fact}}_{\ul{x}} \overset{\on{ins.vac}_{\ul{x}}}\longrightarrow \CA\mod^{\on{fact}}_{\Ran_{\ul{x}}} 
\overset{\oblv_{\CA,\Ran_{\ul{x}}}}\longrightarrow \Dmod(\Ran_{\ul{x}})
\end{equation}
to
\begin{multline} \label{e:global section modules 4}
\QCoh_{\on{co}}(\fL_\nabla(\CY)_{\ul{x}}) \overset{(\on{pr}^\CY_{\on{small},\ul{x}})^*}\longrightarrow
\QCoh_{\on{co}}(\fL^{\mer\rightsquigarrow\reg}_\nabla(\CY)_{\Ran_{\ul{x}}}) \overset{(\on{ev}_{\Ran_{\ul{x}}})^*}\longrightarrow  \\
\to \QCoh_{\on{co}}(\on{Sect}_\nabla(X-\ul{x},\CY))\otimes 
\Dmod(\Ran_{\ul{x}})\overset{\Gamma(\on{Sect}_\nabla(X-\ul{x},\CY),-)\otimes \on{Id}}\longrightarrow \Dmod(\Ran_{\ul{x}}).
\end{multline}

\sssec{}

Note, however, that the map
$$\on{Sect}_\nabla(X-\ul{x},\CY)\times \Ran_{\ul{x}} \overset{\on{ev}_{\Ran_{\ul{x}}}}\longrightarrow 
\fL^{\mer\rightsquigarrow\reg}_\nabla(\CY)_{\Ran_{\ul{x}}} \overset{\on{pr}^\CY_{\on{small},\ul{x}}}\longrightarrow \fL_\nabla(\CY)_{\ul{x}}$$
equals
$$\on{Sect}_\nabla(X-\ul{x},\CY)\times \Ran_{\ul{x}} \to \on{Sect}_\nabla(X-\ul{x},\CY)\overset{\on{ev}_{\ul{x}}}\to \fL_\nabla(\CY)_{\ul{x}}.$$

Hence, we obtain that \eqref{e:global section modules 4} can be rewritten as
\begin{multline} \label{e:global section modules 4.5}
\QCoh_{\on{co}}(\fL_\nabla(\CY)_{\ul{x}}) \overset{\on{ev}_{\ul{x}}^*}\longrightarrow \QCoh_{\on{co}}(\on{Sect}_\nabla(X-\ul{x},\CY))
\overset{\on{Id}\otimes \omega_{\Ran{\ul{x}}}}\longrightarrow \\
\to \QCoh_{\on{co}}(\on{Sect}_\nabla(X-\ul{x},\CY)) \otimes \Dmod(\Ran_{\ul{x}})
\overset{\Gamma(\on{Sect}_\nabla(X-\ul{x},\CY),-)\otimes \on{Id}}\longrightarrow \Dmod(\Ran_{\ul{x}}),
\end{multline}
which is the same as 
\begin{multline} \label{e:global section modules 4.75}
\QCoh_{\on{co}}(\fL_\nabla(\CY)_{\ul{x}}) \overset{\on{ev}_{\ul{x}}^*}\longrightarrow \QCoh_{\on{co}}(\on{Sect}_\nabla(X-\ul{x},\CY))
\overset{\Gamma(\on{Sect}_\nabla(X-\ul{x},\CY),-)}\longrightarrow \Vect \overset{\omega_{\Ran{\ul{x}}}} \longrightarrow \Dmod(\Ran_{\ul{x}}).
\end{multline}

Thus, we obtain a natural transformation from \eqref{e:global section modules 3} to \eqref{e:global section modules 4.75}. By adjunction,
this gives rise to the desired natural transformation from \eqref{e:global section modules 2} to \eqref{e:global section modules 1}.

\sssec{}

We will now prove that the above natural transformation $\text{\eqref{e:global section modules 2}} \to \text{\eqref{e:global section modules 1}}$
is an isomorphism.

\medskip

Write
$$\fL_\nabla(\CY)=\underset{R'}{``\on{colim}"}\, \Spec(R')=\fL^+_\nabla(\CY')$$
as in \secref{sss:glob poles as colimit} and $\CY'$ as an \eqref{e:modified D-scheme}, so that 
$$\on{Sect}_\nabla(X-\ul{x},\CY)=\underset{R'}{``\on{colim}"}\, \on{Sect}_\nabla(X-\ul{x},\CY').$$

\medskip

Thus, we obtain that 
$$\QCoh_{\on{co}}(\fL_\nabla(\CY)) \simeq \underset{R'}{\on{colim}}\, \QCoh(\fL^+_\nabla(\CY')),$$
and we obtain in order to show that $\text{\eqref{e:global section modules 2}} \to \text{\eqref{e:global section modules 1}}$
is an isomorphism, it suffices to show that it becomes such if we precompose both functors with 
$$\QCoh(\fL^+_\nabla(\CY')) \overset{\iota'_*}\to \QCoh_{\on{co}}(\fL_\nabla(\CY))$$
for every $R'$ as above, where $\iota'$ denotes the corresponding map $\iota':\fL^+_\nabla(\CY')\to \fL_\nabla(\CY)$. 

\sssec{}

Note that for every $R'$ as above, the functor
\begin{equation} \label{e:global section modules 2.5}
\QCoh(\fL^+_\nabla(\CY')) \overset{\iota'_*}\to \QCoh_{\on{co}}(\fL_\nabla(\CY)) 
\overset{\Gamma(\fL_\nabla(\CY)_{\ul{x}},-)^{\on{enh}}} \longrightarrow \CA\mod^{\on{fact}}_{\ul{x}}
\end{equation} 
identifies with
\begin{equation} \label{e:global section modules 2.75}
\QCoh(\fL^+_\nabla(\CY')) \overset{\Gamma(\fL^+_\nabla(\CY')_{\ul{x}},-)}\longrightarrow 
\CA'\mod^{\on{com}}_{\ul{x}}\to \CA'\mod^{\on{fact}}_{\ul{x}}\simeq \CA\mod^{\on{fact}}_{\ul{x}},
\end{equation}
and the functor
\begin{equation} \label{e:global section modules 1.5}
\QCoh(\fL^+_\nabla(\CY'))  \overset{\iota'_*}\to \QCoh_{\on{co}}(\fL_\nabla(\CY)) 
\overset{\on{ev}_{\ul{x}}^*}\longrightarrow \QCoh_{\on{co}}(\on{Sect}_\nabla(X-\ul{x},\CY))
\overset{\Gamma(\on{Sect}_\nabla(X-\ul{x},\CY),-)}\longrightarrow \Vect
\end{equation}
identifies with
\begin{equation} \label{e:global section modules 1.75}
\QCoh(\fL^+_\nabla(\CY'))  \overset{\on{ev}_{\ul{x}}^*}\longrightarrow \QCoh_{\on{co}}(\on{Sect}_\nabla(X,\CY'))
\overset{\Gamma(\on{Sect}_\nabla(X,\CY'),-)}\longrightarrow \Vect.
\end{equation}

Composing \eqref{e:global section modules 2.5} and \eqref{e:global section modules 2.75} with the functor
$\on{C}^{\on{fact}}_\cdot(X,\CA,-)_{\ul{x}}$, we obtain that the natural transformation 
$\text{\eqref{e:global section modules 2}} \to \text{\eqref{e:global section modules 1}}$, constructed above, gives
rise to a natural transformation from
$$\QCoh(\fL^+_\nabla(\CY')) \overset{\Gamma(\fL^+_\nabla(\CY')_{\ul{x}},-)}\longrightarrow 
\CA'\mod^{\on{com}}_{\ul{x}}\to \CA'\mod^{\on{fact}}_{\ul{x}}\simeq \CA\mod^{\on{fact}}_{\ul{x}}
\overset{\on{C}^{\on{fact}}_\cdot(X,\CA,-)_{\ul{x}}}\longrightarrow \Vect,$$
which is the same as
\begin{equation} \label{e:global section modules 2.875}
\QCoh(\fL^+_\nabla(\CY')) \overset{\Gamma(\fL^+_\nabla(\CY')_{\ul{x}},-)}\longrightarrow 
\CA'\mod^{\on{com}}_{\ul{x}}\to \CA'\mod^{\on{fact}}_{\ul{x}} \overset{\on{C}^{\on{fact}}_\cdot(X,\CA',-)_{\ul{x}}}\longrightarrow \Vect,
\end{equation}
to \eqref{e:global section modules 1.75}.

\medskip

However, unwinding the definitions, we obtain that the resulting natural transformation
$$\text{\eqref{e:global section modules 2.875}}\to \text{\eqref{e:global section modules 1.75}},$$
coincides with the natural isomorphism of \propref{p:global section modules +}.

\qed[\propref{p:global section modules}]

\ssec{Interpretation as factorization restriction}

\sssec{} \label{sss:spread Y placid}

In this section we will assume that $\fL_\nabla(\CY)$ is ind-placid, and that the embeddings
$$\iota:\fL^+_\nabla(\CY)\to \fL_\nabla(\CY) \text{ and } 
\fL^{\mer\rightsquigarrow\reg}_\nabla(\CY)_{\Ran_{\ul{x}}}\overset{\iota^{\mer\rightsquigarrow\reg}}\longrightarrow \fL_\nabla(\CY)_{\Ran_{\ul{x}}}$$
are locally almost of finite presentation, where
$$\iota^{\mer\rightsquigarrow\reg}=\on{pr}^\CY_{\on{big},\ul{x}}.$$

This assumption implies that $\fL^+_\nabla(\CY)$ (resp., $\fL^{\mer\rightsquigarrow\reg}_\nabla(\CY)_{\Ran_{\ul{x}}}$)
is placid (resp., ind-placid). 

\sssec{Example}

The above assumptions hold for $\CY=\on{Jets}(\CE)$, where $\CE$ is the total space of a vector bundle on $X$.

\medskip

In particular, they hold for $\CY=\Op_\cG$. 

\sssec{}

We can consider the factorization categories
$$\IndCoh^*(\fL^+_\nabla(\CY)) \text{ and } \IndCoh^*(\fL_\nabla(\CY))$$
and  the factorization functor
\begin{equation} \label{e:iota Y}
\iota^\IndCoh_*:\IndCoh^*(\fL^+_\nabla(\CY)) \to \IndCoh^*(\fL_\nabla(\CY)).
\end{equation}

\sssec{} \label{sss:spread Y x}

Consider the assignment
$$(\CZ\to \Ran_{\ul{x}}) \mapsto 
\IndCoh^*(\fL^{\mer\rightsquigarrow\reg}_\nabla(\CY)_{\Ran_{\ul{x}}}\underset{\Ran_{\ul{x}}}\times \CZ)$$
as a crystal of categories over $\Ran_{\ul{x}}$. 

\medskip

It has a natural structure of factorization module category 
with respect to $\IndCoh^*(\fL^+_\nabla(\CY))$; when viewed in this capacity we will denote it by
$\IndCoh^*(\fL^{\mer\rightsquigarrow\reg}_\nabla(\CY))^{\on{fact}_{\ul{x}}}$.

\sssec{}

The map $\iota^{\mer\rightsquigarrow\reg}$ gives rise to a functor
\begin{equation} \label{e:iota Y mer to reg}
(\iota^{\mer\rightsquigarrow\reg})^\IndCoh_*: 
\IndCoh^*(\fL^{\mer\rightsquigarrow\reg}_\nabla(\CY))^{\on{fact}_{\ul{x}}}\to
\IndCoh^*(\fL_\nabla(\CY))^{\on{fact}_{\ul{x}}}
\end{equation}
as factorization module categories, compatible with the factorization functor \eqref{e:iota Y}.

\medskip

Hence, by \secref{sss:univ property restr}, we obtain a functor
\begin{equation} \label{e:iota Y mer to reg restr} 
\IndCoh^*(\fL^{\mer\rightsquigarrow\reg}_\nabla(\CY))^{\on{fact}_{\ul{x}}}\to
\Res_{\iota^\IndCoh_*}(\IndCoh^*(\fL_\nabla(\CY))^{\on{fact}_{\ul{x}}}). 
\end{equation}

\sssec{}

We claim:

\begin{lem} \label{l:iota Y mer to reg restr}
The functor \eqref{e:iota Y mer to reg restr} is an equivalence.
\end{lem}

\begin{proof}

Repeats that of \lemref{l:reg to mf fact module comparison}.

\end{proof}

\section{From module categories over \texorpdfstring{$\QCoh(\LS^\mer_\cG)$}{QCohLS} 
to factorization module categories over \texorpdfstring{$\Rep(\cG)$}{RepG}} \label{s:from LS}

This section is not logically necessary for the rest of the paper. Here we explain a procedure that associates to a module category over
$\QCoh(\LS^\mer_\cG)$ a factorization module category over $\Rep(\cG)$, and show that this functor is fully faithful (on a certain
subcategory). 

\medskip

Using this, we deduce an alternative proof of \propref{p:structure over LS compat}. 

\ssec{Creating factorization modules categories}

Throughout this section we will work at a fixed point $x\in X$. In this subsection we let $H$ be an arbitrary algebraic group.

\sssec{}

Consider the space $\LS^\mer_{H,x}$ and the monoidal category $\QCoh(\LS^\mer_{H,x})$. 
Let us recall the construction of a functor
\begin{equation} \label{e:from QCoh mmod to fact mmod}
\QCoh(\LS^\mer_{H,x})\mmod\to \Rep(H)\mmod^{\on{fact}}_x, \quad \bC\mapsto \bC^{\on{fact}_x,\Rep(H)}.
\end{equation} 

\sssec{}

Namely, we will construct an object 
$$\QCoh(\LS^{\mer\rightsquigarrow\reg}_H)^{\on{fact}_x,\Rep(H)}\in \Rep(H)\mmod^{\on{fact}}_x$$
that carries a commuting action of $\QCoh(\LS^\mer_{H,x})$. The functor \eqref{e:from QCoh mmod to fact mmod}
will then be given by
$$\QCoh(\LS^{\mer\rightsquigarrow\reg}_H)^{\on{fact}_x,\Rep(H)}\underset{\QCoh(\LS^\mer_{H,x})}\otimes -.$$

\sssec{}

The object $\QCoh(\LS^{\mer\rightsquigarrow\reg}_H)^{\on{fact}_x,\Rep(H)}$ will have the feature that its underlying DG category, equipped with an action of
$\QCoh(\LS^\mer_{H,x})$, identifies with $\QCoh(\LS^\mer_{H,x})$ itself.

\medskip

This will imply that the functor \eqref{e:from QCoh mmod to fact mmod} has the feature that for $\bC\in \QCoh(\LS^\mer_{H,x})\mmod$,
the category underlying $\bC^{\on{fact}_x,\Rep(H)}$ identifies with the original $\bC$. 

\sssec{}

The object $\QCoh(\LS^{\mer\rightsquigarrow\reg}_H)^{\on{fact}_x,\Rep(H)}$ is constructed as follows. 

\medskip 

For our fixed point $x\in X$ and a finite subset $x\in \ul{x}\subset X$, consider the multi-disc $\cD_{\ul{x}}$, and set
$$(\QCoh(\LS^{\mer\rightsquigarrow\reg}_H)^{\on{fact}_x,\Rep(H)})_{\ul{x}}:=
\QCoh(\LS^{\mer\rightsquigarrow \reg}_{H,x\subseteq \ul{x}}),$$
see \secref{sss:LS partial puncture}.

\begin{rem}

The construction above is a cousin of the following construction for affine D-schemes: starting from a
module category over $\QCoh(\fL_\nabla(\CY))$ we can associate to it a factorization module category
over $\QCoh(\fL^+_\nabla(\CY))$ by the operation of tensoring with 
$$\QCoh_{\on{co}}(\fL^{\mer\rightsquigarrow\reg}_\nabla(\CY))^{\on{fact}_{\ul{x}}},$$
see \secref{sss:spread Y x}. 

\end{rem}

\begin{rem} \label{r:compare with Gr}

This construction also be viewed as an analog of a construction in \cite{CFGY} that associates to a $\fL(G)_x$-module
category a factorization module category over $\Dmod(\Gr_G)$; in fact a version of the latter construction
has appeared in \secref{sss:from LG to Gr_G}. 

\end{rem}

\sssec{}

In what follows we will need a variant of the above construction, when instead of 
$\LS^\mer_{H,x}$ we use its formal completion $(\LS^\mer_{H,x})^\wedge_\reg$
around $\LS^\reg_{H,x}$.

\medskip

The corresponding functor 
\begin{equation} \label{e:from QCoh mmod to fact mmod form}
\QCoh((\LS^\mer_{H,x})^\wedge_\reg)\mmod\to \Rep(H)\mmod^{\on{fact}}_x, \quad \bC\mapsto \bC^{\on{fact}_x,\Rep(H)}
\end{equation} 
is constructed using the prestack
$$(x\subseteq \ul{x})\rightsquigarrow (\LS^{\mer\rightsquigarrow \reg}_{H,x\subseteq \ul{x}})^\wedge_\reg,$$
i.e., the formal completion of $\LS^{\mer\rightsquigarrow \reg}_{H,x\subseteq \ul{x}}$ along $(\LS^\reg_H)^{\on{fact}_x}_{\ul{x}}$.

\sssec{Example} \label{sss:Rep(G) fact}

Consider $\QCoh(\LS^\reg_{H,x})\simeq \Rep(H)$ as an object of $\QCoh((\LS^\mer_{H,x})^\wedge_\reg)\mmod$, via the restriction functor
$$\QCoh((\LS^\mer_{H,x})^\wedge_\reg)\to \QCoh(\LS^\reg_{H,x}).$$

We claim that 
\begin{equation} \label{e:Rep(G) fact}
\QCoh(\LS^\reg_{H,x})^{\on{fact}_x,\Rep(H)}\simeq \Rep(H)^{\on{fact}_x},
\end{equation}
i.e., the vacuum object of  $\Rep(H)\mmod^{\on{fact}}_x$.

\medskip

Indeed, this follows from \cite[Chapter 3, Proposition 3.5.3]{GaRo3} using the fact that
$$(\LS^{\mer\rightsquigarrow \reg}_{H,\ul{x}\subseteq \ul{x}'})^\wedge_\reg\underset{(\LS^\mer_{H,x})^\wedge_\reg}\times
\LS^\reg_{H,x} \simeq (\LS^\reg_H)^{\on{fact}_x}_{\ul{x}}.$$

When applying \cite[Chapter 3, Proposition 3.5.3]{GaRo3} we use the fact that we can identify
$$(\LS^\mer_{H,x})^\wedge_\reg\simeq \fh^\wedge_0/\on{Ad}(H),$$
so this prestack is passable. 

\sssec{}

We claim: 

\begin{thm} \label{t:create fact ff}
The functor \eqref{e:from QCoh mmod to fact mmod form} is fully faithful.
\end{thm}

This theorem is proved in \cite[Theorem 9.13.1]{Ra4}. We will supply a proof for completeness. 

\begin{rem}
We conjecture that the functor \eqref{e:from QCoh mmod to fact mmod} is itself fully faithful. A partial result in this direction has recently been established
in \cite{Bogd}: the composition \eqref{e:from QCoh mmod to fact mmod} with the restriction functor
$$\QCoh(\LS^{\on{restr}}_{H,x})\mmod\to \QCoh(\LS^\mer_{H,x})\mmod$$
is fully faithful, where
$$\LS^{\on{restr}}_{H,x}\subset \LS^\mer_{H,x}$$
is the stack of local systems with restricted variation (see \cite[Sect. 1.4]{AGKRRV}).
\end{rem}

\begin{rem}
The reason the discussion in this section is for a fixed point in the Ran space is that we do not know how to prove 
the analog of \thmref{t:create fact ff} in the factorization setting (and are not confident in its validity).
\end{rem}

\ssec{Proof of \thmref{t:create fact ff}}

The proof is a ``baby version" of the argument proving the corresponding assertion in \cite{CFGY},
see Remark \ref{r:compare with Gr}. 

\sssec{Reduction steps}

First, we note that the functor \eqref{e:from QCoh mmod to fact mmod form} preserves both limits and colimits. 

\medskip

\medskip

Second, the (symmetric monoidal) restriction functor
$$\QCoh((\LS^\mer_{H,x})^\wedge_\reg) \to \QCoh(\LS^\reg_{H,x})$$
is comonadic, and its right adjoint is $\QCoh((\LS^\mer_{H,x})^\wedge_\reg) $-linear. 

\medskip

It follows that the
2-category $\QCoh((\LS^\mer_{H,x})^\wedge_\reg)\mmod$ is generated under colimits
by the essential image of
\begin{equation}\label{e:restr from qcoh on reg opers compl}
\QCoh(\LS^\reg_{H,x})\mmod\to \QCoh((\LS^\mer_{H,x})^\wedge_\reg)\mmod,
\end{equation}
given by restriction along
$$\LS^\reg_{H,x}\to (\LS^\mer_{H,x})^\wedge_\reg.$$

Moreover, by passing to right adjoints, the 2-category $\QCoh((\LS^\mer_{H,x})^\wedge_\reg)\mmod$ is also generated under \emph{limits} by the essential image of \eqref{e:restr from qcoh on reg opers compl}.

\medskip

Third, $\QCoh(\LS^\reg_{H,x})\mmod\simeq \Rep(H)\mmod$ is generated under colimits (and separately, limits) by
objects of the form $\bC\otimes \Rep(H)$, for $\bC\in \on{DGCat}$.

\medskip

And fourth, since
$$ \Rep(H) \simeq \Vect \underset{\QCoh(H)}{\otimes} \Vect,$$
we obtain that the object $\Rep(H)\in \Rep(H)\mmod$ is a colimit of objects on which the action of $\Rep(H)$ is trivial, 
i.e., factors via the augmentation functor
$$\Rep(H)\to \Vect,$$
corresponding to 
$$\on{pt}\to \LS^\reg_{H,x}.$$

\medskip

Combining these observations, we obtain that the
2-category $\QCoh((\LS^\mer_{H,x})^\wedge_\reg)\mmod$ is 
generated under colimits by trivial modules and under limits by objects of the form $\bC_0\otimes \Rep(H)$
for a DG category $\bC_0$ with a trivial action of $\Rep(H)$. 

\medskip

Hence, it is enough to show that the functor
\begin{multline} \label{e:ff Rep G to prove}
\on{Funct}_{\QCoh((\LS^\mer_{H,x})^\wedge_\reg)\mmod}(\Vect,\bC_0\otimes \Rep(H))\to \\
\to \on{Funct}_{\Rep(H)\mmod^{\on{fact}}_x}\left(\Vect^{\on{fact}_x,\Rep(H)},\bC_0\otimes \Rep(H)^{\on{fact}_x}\right)
\end{multline}
is an equivalence.

%

\sssec{}

Recall that for a (unital) factorization functor $\Phi:\bA_1\to \bA_2$ and an object $\bC_2\in \bA_2\mmod^{\on{fact}}_x$, we denote
by $\Res_\Phi(\bC_2)$ its restriction along $\Phi$.

\medskip

Consider the tautological object $\Vect^{\on{fact}_x}\in \Vect\mmod^{\on{fact}}_x$, and the corresponding object 
$$\Res_{\oblv_H}(\Vect^{\on{fact}_x})\in \Rep(H)\mmod^{\on{fact}}_x,$$
where $\oblv_H:\Rep(H)\to \Vect$ is the forgetful functor, viewed as a (strictly unital) factorization functor.

\medskip

Pullback along the unit section of $(\LS^{\mer\rightsquigarrow \reg}_{H,\ul{x}\subseteq \ul{x}'})^\wedge_\reg$
gives rise to a functor
\begin{equation} \label{e:restr of Vect prel}
\Vect^{\on{fact}_x,\Rep(H)} \to \Vect^{\on{fact}_x}
\end{equation}
compatible with factorization (in the sense of \secref{sss:funct between pairs modules untl}). Hence, we obtain a morphism
\begin{equation} \label{e:restr of Vect}
\Vect^{\on{fact}_x,\Rep(H)} \to \Res_{\oblv_H}(\Vect^{\on{fact}_x})
\end{equation}
as unital factorization modules over $\Rep(H)$. 

\sssec{}

We claim:

\begin{lem}  \label{l:restr of Vect}
The functor \eqref{e:restr of Vect} is an equivalence.
\end{lem}

\begin{proof}

Follows from \lemref{l:fact res crit}.

\end{proof} 

\sssec{}

Recall now the paradigm of \lemref{l:basic adj}. We apply it to $\oblv_H:\Rep(H)\to \Vect$. Combined with
\lemref{l:restr of Vect} above, we obtain that for $\bC\in \Rep(H)\mmod^{\on{fact}}_x$
we have a canonical equivalence
\begin{equation} \label{e:basic adj Rep G}
\on{Funct}_{\Rep(H)\mmod^{\on{fact}}_x}\left(\Vect^{\on{fact}_x,\Rep(H)},\bC\right)\simeq 
R_H\mod^{\on{fact}}(\bC_x).
\end{equation}

\sssec{}

Taking $\bC:=\Rep(H)^{\on{fact}_x}\otimes \bC_0$ in \eqref{e:basic adj Rep G} and using \eqref{e:Rep(G) fact}, 
we obtain that \eqref{e:ff Rep G to prove} reduces to showing that the resulting functor
\begin{equation} \label{m:ff Rep G to prove q C0}
\on{Funct}_{\QCoh((\LS^\mer_{H,x})^\wedge_\reg)\mmod}(\Vect,\Rep(H)\otimes \bC_0)\to R_H\mod^{\on{fact}}(\Rep(H)\otimes \bC_0)
\end{equation}
is an equivalence.

\medskip

It is clear that the functor
$$\on{Funct}_{\QCoh((\LS^\mer_{H,x})^\wedge_\reg)\mmod}(\Vect,\Rep(H))\otimes \bC_0\to
\on{Funct}_{\QCoh((\LS^\mer_{H,x})^\wedge_\reg)\mmod}(\Vect,\Rep(H)\otimes \bC_0)$$
is an equivalence (since $\QCoh((\LS^\mer_{H,x})^\wedge_\reg$ is a semi-rigid monoidal category,
see \cite[Sect. C]{AGKRRV}). 

\medskip

Now the fact that $R_H$ is holonomic implies that
$$R_H\mod^{\on{fact}}(\Rep(H))\otimes \bC_0\to R_H\mod^{\on{fact}}(\Rep(H)\otimes \bC_0)$$
is also an equivalence (see \cite[Theorem 8.13.1]{Ra4}). 

\medskip

Hence, it is enough to show that the functor
\begin{equation} \label{m:ff Rep G to prove q}
\on{Funct}_{\QCoh((\LS^\mer_{H,x})^\wedge_\reg)\mmod}(\Vect,\Rep(H))\to R_H\mod^{\on{fact}}(\Rep(H))
\end{equation}
is an equivalence. 

\sssec{}

We rewrite the left-hand side in \eqref{m:ff Rep G to prove q} as 
$$\QCoh(\on{pt}\underset{(\LS^\mer_{H,x})^\wedge_\reg}\times \LS^\reg_{H,x})\simeq 
\QCoh(\on{pt}\underset{\LS^\mer_{H,x}}\times \LS^\reg_{H,x}).$$

\sssec{}

By \secref{sss:rel fact flat}, the category $(R_H\mod^{\on{fact}}(\Rep(H)))^{>-\infty}$ identifies with 
$$\QCoh_{\on{co}}\Bigl(\on{pt}\underset{\LS^\mer_{H,x}}\times \LS^\reg_{H,x}\Bigr)^{>-\infty}.$$

Since $R_H\mod^{\on{fact}}(\Rep(H))$ is left-complete (see \propref{p:t-structure fact mod}), we obtain that it identifies
with the left completion of
$$\QCoh_{\on{co}}(\on{pt}\underset{\LS^\mer_{H,x}}\times \LS^\reg_{H,x}).$$

\medskip

Since $\on{pt}\underset{\LS^\mer_{H,x}}\times \LS^\reg_{H,x}$ is almost of finite type, we have an equivalence
$$\Psi_{\on{pt}\underset{\LS^\mer_{H,x}}\times \LS^\reg_{H,x}}:
\QCoh_{\on{co}}\Bigl(\on{pt}\underset{\LS^\mer_{H,x}}\times \LS^\reg_{H,x}\Bigr)^{>-\infty}\overset{\sim}
\to \IndCoh\Bigl(\on{pt}\underset{\LS^\mer_{H,x}}\times \LS^\reg_{H,x}\Bigr)^{>-\infty}.$$

Hence, the above left-completion identifies with 
$$\QCoh(\on{pt}\underset{\LS^\mer_{H,x}}\times \LS^\reg_{H,x}).$$

\sssec{}

Unwinding the constructions, we obtain that the endo-functor of $\QCoh(\on{pt}\underset{\LS^\mer_{H,x}}\times \LS^\reg_{H,x})$,
induced by \eqref{m:ff Rep G to prove q} and the above two identifications is the identity.

\qed[\thmref{t:create fact ff}]

\ssec{Factorization module categories attached to affine D-schemes} \label{ss:fact mod sch}

\sssec{}

Let $\CY$ be an affine D-scheme over $X$, equipped with a map 
$$\CY\to \on{pt}/H.$$

Consider the corresponding factorization spaces
\begin{equation} \label{e:loops into Y}
\fL^+_\nabla(\CY) \overset{\iota}\hookrightarrow \fL_\nabla(\CY)
\end{equation}
and the commutative (but \emph{not necessarily} Cartesian) diagram
$$
\CD
\fL^+_\nabla(\CY)  @>>> \fL_\nabla(\CY) \\
@V{\fr^\reg}VV @VV{\fr}V \\
\LS^\reg_H @>>> \LS^{\on{mer}}_H
\endCD
$$

\sssec{}

On the one hand, we can consider $\QCoh_{\on{co}}(\fL_\nabla(\CY)_x)$ as an object of 
$\QCoh(\LS^\mer_{H,x})\mmod$. Consider the resulting object
\begin{equation} \label{e:fact mod abs 1} 
\QCoh_{\on{co}}(\fL_\nabla(\CY)_x)^{\on{fact}_x,\Rep(H)}\in \Rep(H)\mmod^{\on{fact}}.
\end{equation}

\sssec{} \label{sss:fact mod abs 2 QCoh co}

On the other hand, consider the $\QCoh(\fL_\nabla^+(\CY))$-factorization category
$$\QCoh_{\on{co}}(\fL^{\mer\rightsquigarrow\reg}_\nabla(\CY))^{\on{fact}_x},$$
defined as in \secref{sss:spread Y x}.

\medskip

Pullback along 
$$\fr^\reg:\fL^+_\nabla(\CY)\to \LS^\reg_H$$ defines a 
(strictly unital) factorizaton functor
$$\Rep(H)\simeq \QCoh(\LS^\reg_H) \overset{(\fr^\reg)^*}\longrightarrow \QCoh(\fL_\nabla^+(\CY)).$$

\medskip

Consider the resulting object 
\begin{equation} \label{e:fact mod abs 2 QCo co} 
\Res_{(\fr^\reg)^*}\left(\QCoh_{\on{co}}(\fL^{\mer\rightsquigarrow\reg}_\nabla(\CY))^{\on{fact}_x}\right)\in  \Rep(H)\mmod_x^{\on{fact}}.
\end{equation}

\sssec{}

Note that the natural morphism
\begin{equation} \label{e:Y LS spread} 
\fL^{\mer\rightsquigarrow\reg}_\nabla(\CY)_{x\subseteq \ul{x}}\to 
\LS^{\mer\rightsquigarrow \reg}_{H,x\subseteq \ul{x}}\underset{\LS^\mer_{H,x}}\times \fL_\nabla(\CY)_x
\end{equation} 
is affine. 

\medskip

Hence, pullback along \eqref{e:Y LS spread} gives rise to a functor
\begin{equation} \label{e:restriction and fact prel}
\QCoh_{\on{co}}(\fL_\nabla(\CY)_x)^{\on{fact}_x,\Rep(H)}\to  
\QCoh_{\on{co}}(\fL^{\mer\rightsquigarrow\reg}_\nabla(\CY))^{\on{fact}_x}
\end{equation}
compatible with factorization. 

\medskip

Hence, we obtain a functor
\begin{equation} \label{e:restriction and fact}
\QCoh_{\on{co}}(\fL_\nabla(\CY)_x)^{\on{fact}_x,\Rep(H)}\to 
\Res_{(\fr^\reg)^*}\left(\QCoh_{\on{co}}(\fL^{\mer\rightsquigarrow\reg}_\nabla(\CY))^{\on{fact}_x}\right)
\end{equation}
in $\Rep(H)\mmod^{\on{fact}}$. 

\sssec{Variant}

Let $\fL_\nabla(\CY)^\wedge_\mf$ denote the formal completion of $\fL_\nabla(\CY)$ along 
$$\fL^\mf_\nabla(\CY):=\fL_\nabla(\CY)\underset{\LS^\mer_H}\times \LS^\reg_H\to \fL_\nabla(\CY),$$
or which is the same
$$\fL_\nabla(\CY)\underset{\LS^\mer_H}\times (\LS^\mer_H)^\wedge_\reg.$$

\medskip

Similar to the above, we can consider
$$\QCoh_{\on{co}}((\fL_\nabla(\CY)^\wedge_\mf)_x)^{\on{fact}_x,\Rep(H)}\in \Rep(H)\mmod^{\on{fact}}$$
and
$$\QCoh_{\on{co}}(\fL^{\mer\rightsquigarrow\reg}_\nabla(\CY)^\wedge_\mf)^{\on{fact}_x},$$
and we obtain a functor
\begin{equation} \label{e:restriction and fact form compl}
\QCoh_{\on{co}}((\fL_\nabla(\CY)^\wedge_\mf)_x)^{\on{fact}_x,\Rep(H)}\to 
\Res_{(\fr^\reg)^*}\left(\QCoh_{\on{co}}(\fL^{\mer\rightsquigarrow\reg}_\nabla(\CY)^\wedge_\mf)^{\on{fact}_x}\right)
\end{equation} 
in $\Rep(H)\mmod^{\on{fact}}$. 

\sssec{Variant}

Let us assume now that $\fL_\nabla(\CY)$ and also $\on{pt}\underset{\LS^\mer_H}\times \fL_\nabla(\CY)$
satisfy the assumptions of \secref{sss:placid for unital}.

\medskip

This implies that the morphism
$$\LS^{\on{fact}_x}_{H,x\subseteq \ul{x}}\underset{\LS^{\mer\rightsquigarrow \reg}_{H,x\subseteq \ul{x}}}\times
\fL^{\mer\rightsquigarrow\reg}_\nabla(\CY)_{x\subseteq \ul{x}}\to
\LS^{\on{fact}_x}_{H,x\subseteq \ul{x}}\underset{\LS^\reg_{H,x}}\times \LS^\reg_{H,x}\underset{\LS^\mer_{H,x}}\times \fL_\nabla(\CY)_x
= \LS^{\on{fact}_x}_{H,x\subseteq \ul{x}}\underset{\LS^\mer_{H,x}}\times \fL_\nabla(\CY)_x
$$
is of finite Tor-dimension. 

\medskip

Hence, so is the morphism
\begin{equation} \label{e:Y LS spread form} 
\left(\fL^{\mer\rightsquigarrow\reg}_\nabla(\CY)_{x\subseteq \ul{x}}\right)^\wedge_\mf \to 
\left(\LS^{\mer\rightsquigarrow \reg}_{H,x\subseteq \ul{x}}\right)^\wedge_\mf
\underset{(\LS^\mer_{H,x})^\wedge_\reg}\times (\fL_\nabla(\CY)^\wedge_\mf)_x.
\end{equation} 

Hence, by \secref{sss:* pullback on IndCoh* finite Tor dim}, 
the $(\IndCoh,*)$-pullback along \eqref{e:Y LS spread form} gives rise to a functor
\begin{equation} \label{e:restriction and fact form IndCoh prel}
\IndCoh^*\left((\fL_\nabla(\CY)^\wedge_\mf)_x\right)^{\on{fact}_x,\Rep(H)}\to
\IndCoh^*\left(\left(\fL^{\mer\rightsquigarrow\reg}_\nabla(\CY)\right)^\wedge_\mf\right)^{\on{fact}_x}
\end{equation} 
compatible with factorization. 

\medskip

Hence, we obtain a functor
\begin{equation} \label{e:restriction and fact form IndCoh}
\IndCoh^*\left((\fL_\nabla(\CY)^\wedge_\mf)_x\right)^{\on{fact}_x,\Rep(H)}\to
\Res_{(\fr^\reg)^*}\left(\IndCoh^*\left(\left(\fL^{\mer\rightsquigarrow\reg}_\nabla(\CY)\right)^\wedge_\mf\right)^{\on{fact}_x}\right)
\end{equation} 
in $\Rep(H)\mmod^{\on{fact}}$.

\sssec{}

We claim:

\begin{lem} \label{l:restriction and fact form IndCoh}
The functor \eqref{e:restriction and fact form IndCoh} is an equivalence.
\end{lem}

\begin{proof}

Follows from \lemref{l:fact res crit}.

\end{proof} 

\sssec{}

Let $\iota^{+,\mf}$ denote the morphism
$$\fL^+_\nabla(\CY)\to \fL_\nabla(\CY)^\wedge_\mf.$$

The operation of $\IndCoh$-pushforward along $\iota^{+,\mf}$ gives rise to a (strictly) unital factorization functor 
$$(\iota^{+,\mf})^\IndCoh_*:\IndCoh^*(\fL^+_\nabla(\CY))\to \IndCoh^*(\fL^+_\nabla(\CY)^\wedge_\mf).$$

It follows from \lemref{l:iota Y mer to reg restr} that the naturally defined functor
$$\IndCoh^*\left(\left(\fL^{\mer\rightsquigarrow\reg}_\nabla(\CY)\right)^\wedge_\mf\right)^{\on{fact}_x}\to
\Res_{(\iota^{+,\mf})^\IndCoh_*}\left((\IndCoh^*(\fL^+_\nabla(\CY)^\wedge_\mf))^{\on{fact}_x}\right)$$
is an equivalence. 

\medskip

Hence, combining with \lemref{l:restriction and fact form IndCoh}, we obtain:

\begin{cor} \label{c:restriction and fact form IndCoh}
There is a canonical equivalence 
$$\IndCoh^*\left((\fL_\nabla(\CY)^\wedge_\mf)_x\right)^{\on{fact}_x,\Rep(H)}\simeq
\Res_{(\iota^{+,\mf})^\IndCoh_*\circ (\fr^\reg)^*}\left((\IndCoh^*(\fL^+_\nabla(\CY)^\wedge_\mf))^{\on{fact}_x}\right)$$
in $\Rep(H)\mmod^{\on{fact}}$.
\end{cor}

%
%

\ssec{Proof of \propref{p:structure over LS compat} via factorization}

\sssec{}

Consider the following (strictly unital) factorization functor, to be denoted $\Phi_1$
\begin{multline} \label{e:Rep to Whit KM}
\Rep(\cG) \overset{\FLE_{\cG,\infty}}\simeq \Whit_*(G)\overset{\on{Id}\otimes \on{Vac}(G)_{\crit,\rho(\omega_X)}}\longrightarrow \\
\to \Whit_*(G)\otimes \KL(G)_{\crit,\rho(\omega_X)} \to \Whit_*(\hg\mod^{\on{Sph-gen}}_{\crit,\rho(\omega_X)}).
\end{multline}

Restricting the vacuum factorization module category
$$\Whit_*(\hg\mod^{\on{Sph-gen}}_{\crit,\rho(\omega_X)})^{\on{fact}_x}\in 
\Whit_*(\hg\mod^{\on{Sph-gen}}_{\crit,\rho(\omega_X)})\mmod^{\on{fact}}_x$$
along \eqref{e:Rep to Whit KM}, we obtain an object 
\begin{equation} \label{e:Res Whit to Rep}
\Res_{\Phi_1}(\Whit_*(\hg\mod^{\on{Sph-gen}}_{\crit,\rho(\omega_X)})^{\on{fact}_x})\in \Rep(\cG)\mmod^{\on{fact}}_x.
\end{equation} 

\sssec{}

We claim:

\begin{lem} \label{l:ident Whit as fact mod}
The object \eqref{e:Res Whit to Rep} identifies canonically with 
$$\Whit_*(\hg\mod^{\on{Sph-gen}}_{\crit,\rho(\omega_X)})^{\on{fact}_x,\Rep(\cG)},$$
where we regard $\Whit_*(\hg\mod^{\on{Sph-gen}}_{\crit,\rho(\omega_X)})$ as an object of
$\QCoh((\LS^\mer_{\cG,x})^\wedge_\reg)\mmod$ by the recipe of \secref{sss:Whit ten over LS}. 
\end{lem} 

\begin{proof}

We rewrite $\Phi_1$ as
\begin{multline} \label{e:Rep to Whit KM 0}
\Rep(\cG) \overset{\FLE_{\cG,\infty}}\simeq \Whit_*(G)\overset{\on{Id}\otimes \on{Vac}(G)_{\crit,\rho(\omega_X)}}\longrightarrow 
\Whit_*(G)\otimes \KL(G)_{\crit,\rho(\omega_X)} \to \\
\to \Whit_*(G)\underset{\Sph_G}\otimes \KL(G)_{\crit,\rho(\omega_X)} \simeq \Whit_*(\hg\mod^{\on{Sph-gen}}_{\crit,\rho(\omega_X)}).
\end{multline}

\medskip

Therefore, we can reinterpret the object \eqref{e:Res Whit to Rep} as follows. Consider
$$\Res_{\FLE_{\cG,\infty}}(\Whit_*(G)^{\on{fact}_x})\in \Rep(\cG)\mmod^{\on{fact}}_x$$
as an object equipped with a commuting action of $\Sph_{G,x}$. Then \eqref{e:Res Whit to Rep}
identifies with
$$\Res_{\FLE_{\cG,\infty}}(\Whit_*(G)^{\on{fact}_x})\underset{\Sph_{G,x}}\otimes \KL(G)_{\crit,\rho(\omega_X),x}.$$ 

Comparing with the definition of 
$$\Whit_*(\hg\mod^{\on{Sph-gen}}_{\crit,\rho(\omega_X)})\in \QCoh((\LS^\mer_{\cG,x})^\wedge_\reg)\mmod$$
(see \secref{sss:Whit ten over LS}), the assertion of the lemma follows from the identification \eqref{e:Rep(G) fact}. 

\end{proof} 

\sssec{}

Consider the following (strictly unital) factorization functor, to be denoted $\Phi_2$
\begin{equation} \label{e:Rep to Op}
\Rep(\cG) \overset{(\fr^\reg)^*}\longrightarrow \QCoh(\Op^\reg_\cG) \overset{(\iota^\mf)^\IndCoh_*}\longrightarrow 
\IndCoh^*(\Op^\mer_\cG)_\mf.
\end{equation}

Restricting the vacuum factorization module category
$$(\IndCoh^*(\Op^\mer_\cG)_\mf)^{\on{fact}_x}\in \IndCoh^*(\Op^\mer_\cG)_\mf\mmod^{\on{fact}}_x$$
along \eqref{e:Rep to Op}, we obtain an object
\begin{equation} \label{e:Res Op to Rep}
\Res_{\Phi_2}((\IndCoh^*(\Op^\mer_\cG)_\mf)^{\on{fact}_x})\in \Rep(\cG)\mmod^{\on{fact}}_x.
\end{equation} 

\sssec{}

We claim:

\begin{lem} \label{l:ident Op as fact mod}
The object \eqref{e:Res Op to Rep} identifies canonically with 
$$(\IndCoh^*(\Op^\mer_\cG)_\mf)^{\on{fact}_x,\Rep(\cG)},$$
where we regard $\IndCoh^*(\Op^\mer_\cG)_\mf$ as an object of
$\QCoh((\LS^\mer_{\cG,x})^\wedge_\reg)\mmod$ via
$$\fr:(\Op^\mer_\cG)^\wedge_\mf\to (\LS^\mer_{\cG,x})^\wedge_\reg.$$
\end{lem}

\begin{proof}

This is a particular case of \lemref{l:restriction and fact form IndCoh}.

\end{proof}

\sssec{}

We now claim:

\begin{prop} \label{p:match fact structures}
We have a canonical isomorphism of (strictly unital) factorization functors
$$\ol\DS^{\on{enh,rfnd}}\circ \Phi_1\simeq \Phi_2.$$
\end{prop}

Let us accept this proposition for a moment and finish the proof of \propref{p:structure over LS compat}. 

\sssec{}

By \thmref{t:create fact ff}, combined with Lemmas \ref{l:ident Whit as fact mod} and \ref{l:ident Op as fact mod}, 
it suffices to show that the functor $\ol\DS^{\on{enh,rfnd}}$ appearing in \propref{p:structure over LS compat}
can be realized as the fiber at $x$ of a functor between 
$$\Res_{\Phi_1}(\Whit_*(\hg\mod^{\on{Sph-gen}}_{\crit,\rho(\omega_X)})^{\on{fact}_x}) \text{ and }
\Res_{\Phi_2}((\IndCoh^*(\Op^\mer_\cG)_\mf)^{\on{fact}_x}),$$
viewed as objects of $\Rep(\cG)\mmod^{\on{fact}}_x$. 

\medskip

The latter structure is supplied by \propref{p:match fact structures}. 

\qed[\propref{p:structure over LS compat}]

\ssec{Proof of \propref{p:match fact structures}}

\sssec{}

It is enough to show that the two functors match after we compose them with the fully faithful embedding
$$\IndCoh^*(\Op^\mer_\cG)_\mf\hookrightarrow \IndCoh^*(\Op^\mer_\cG).$$

\medskip

We will first establish an isomorphism between the compositions of the two functors in question with
$$\Gamma^\IndCoh(\Op^\mer_\cG,-)^{\on{enh}}:\IndCoh^*(\Op^\mer_\cG)\to \CO_{\Op^\reg_\cG}\mod^{\on{fact}}.$$

\sssec{}

We start by rewriting the corresponding composition for $\Phi_1$. It is equal to 
\begin{multline} \label{e:Rep to Whit KM 1}
\Rep(\cG) \overset{\FLE_{\cG,\infty}}\simeq \Whit_*(G)\overset{\on{Id}\otimes \on{Vac}(G)_{\crit,\rho(\omega_X)}}\longrightarrow 
\Whit_*(G)\otimes \KL(G)_{\crit,\rho(\omega_X)}\to \\
\to \Whit_*(G)\underset{\Sph_G}\otimes \KL(G)_{\crit,\rho(\omega_X)}\to \Whit_*(\hg\mod_{\crit,\rho(\omega_X)}) \overset{\ol\DS^{\on{enh}}}\longrightarrow 
\CO_{\Op^\reg_\cG}\mod^{\on{fact}}.
\end{multline} 

By Remark \ref{r:CS again}, we rewrite this as
\begin{multline} \label{e:Rep to Whit KM 2}
\Rep(\cG) \overset{\Sat_G^{\on{nv}}}\longrightarrow \Sph_G \overset{\sigma}\to \Sph_G
\overset{-\star \one_{\Whit_*(G)}}\longrightarrow \Whit_*(G)\overset{\on{Id}\otimes \on{Vac}(G)_{\crit,\rho(\omega_X)}}\longrightarrow \\
\to \Whit_*(G)\otimes \KL(G)_{\crit,\rho(\omega_X)}\to \Whit_*(G)\underset{\Sph_G}\otimes \KL(G)_{\crit,\rho(\omega_X)}\to \\
\to \Whit_*(\hg\mod_{\crit,\rho(\omega_X)}) \overset{\ol\DS^{\on{enh}}}\longrightarrow 
\CO_{\Op^\reg_\cG}\mod^{\on{fact}}.
\end{multline} 

We can rewrite the composition in the first two lines in \eqref{e:Rep to Whit KM 2} as
\begin{multline}  \label{e:Rep to Whit KM 2'}
\Rep(\cG) \overset{\Sat_G^{\on{nv}}}\longrightarrow \Sph_G \overset{-\star \on{Vac}(G)_{\crit,\rho(\omega_X)}}\longrightarrow 
 \KL(G)_{\crit,\rho(\omega_X)} \overset{\one_{\Whit_*(G)}\otimes \on{Id}}\longrightarrow \\
\to \Whit_*(G)\otimes \KL(G)_{\crit,\rho(\omega_X)}\to \Whit_*(G)\underset{\Sph_G}\otimes \KL(G)_{\crit,\rho(\omega_X)}.
\end{multline} 

Hence, we can rewrite \eqref{e:Rep to Whit KM 2} as
\begin{multline}  \label{e:Rep to Whit KM 3}
\Rep(\cG) \overset{\Sat_G^{\on{nv}}}\longrightarrow \Sph_G \overset{-\star \on{Vac}(G)_{\crit,\rho(\omega_X)}}\longrightarrow \\
\to \KL(G)_{\crit,\rho(\omega_X)}
\to \hg\mod_{\crit,\rho(\omega_X)} \overset{\DS^{\on{enh}}}\longrightarrow \CO_{\Op^\reg_\cG}\mod^{\on{fact}}.
\end{multline}  

\sssec{}

We now apply \thmref{t:birth}. It says that the functor
$$\Rep(\cG) \overset{\Sat_G^{\on{nv}}}\longrightarrow \Sph_G \overset{-\star \on{Vac}(G)_{\crit,\rho(\omega_X)}}\longrightarrow 
\KL(G)_{\crit,\rho(\omega_X)}$$
is isomorphic to
\begin{equation}  \label{e:Rep to Whit KM 4}
\Rep(\cG) \overset{(\fr^\reg)^*}\longrightarrow \QCoh(\Op^\reg_\cG) \overset{\Gamma(\Op^\reg_\cG,-)}\longrightarrow 
\CO_{\Op^\reg_\cG}\mod^{\on{com}} \overset{-\underset{\CO_{\Op^\reg_\cG}}\otimes \on{Vac}(G)_{\crit,\rho(\omega_X)}}\longrightarrow \KL(G)_{\crit,\rho(\omega_X)}.
\end{equation}  

Hence, we can rewrite \eqref{e:Rep to Whit KM 2} as 

\begin{equation}  \label{e:Rep to Whit KM 5}
\Rep(\cG) \overset{(\fr^\reg)^*}\longrightarrow \QCoh(\Op^\reg_\cG) \overset{\Gamma(\Op^\reg_\cG,-)}\longrightarrow 
\CO_{\Op^\reg_\cG}\mod^{\on{com}}  \overset{-\underset{\CO_{\Op^\reg_\cG}}\otimes \DS^{\on{enh}}(\on{Vac}(G)_{\crit,\rho(\omega_X)})}\longrightarrow 
\CO_{\Op^\reg_\cG}\mod^{\on{fact}}.
\end{equation}

\sssec{}

We now use \thmref{t:center as DS}, which says that the natural map
$$\CO_{\Op^\reg_\cG}\to \DS^{\on{enh}}(\on{Vac}(G)_{\crit,\rho(\omega_X)})$$
is an isomorphism.

\medskip

This allows us to rewrite \eqref{e:Rep to Whit KM 5} as
\begin{equation}  \label{e:Rep to Whit KM 6}
\Rep(\cG) \overset{(\fr^\reg)^*}\longrightarrow \QCoh(\Op^\reg_\cG) \overset{\Gamma(\Op^\reg_\cG,-)}\longrightarrow 
\CO_{\Op^\reg_\cG}\mod^{\on{com}} \to \CO_{\Op^\reg_\cG}\mod^{\on{fact}}.
\end{equation}

\sssec{}

We now consider the composition involving $\Phi_2$:
\begin{multline}  \label{e:Rep to Whit KM 7}
\Rep(\cG) \overset{(\fr^\reg)^*}\longrightarrow \QCoh(\Op^\reg_\cG) \overset{(\iota^\mf)^\IndCoh_*}\longrightarrow 
\IndCoh^*(\Op^\mer_\cG)\overset{\Gamma^\IndCoh(\Op^\mer_\cG,-)^{\on{enh}}}\longrightarrow \CO_{\Op^\reg_\cG}\mod^{\on{fact}}.
\end{multline}

This functor identifies with \eqref{e:Rep to Whit KM 6} on the nose. 

\sssec{}

Thus, we have identified the compositions of the two functors
\begin{equation} \label{e:F1 and F2}
\Rep(\cG)\rightrightarrows \IndCoh^*(\Op^\mer_\cG)
\end{equation}
with
$$\IndCoh^*(\Op^\mer_\cG)\overset{\Gamma^\IndCoh(\Op^\mer_\cG,-)^{\on{enh}}}\longrightarrow \CO_{\Op^\reg_\cG}\mod^{\on{fact}}.$$

Let us show how to upgrade this isomorphism to an isomorphism between the two functors in \eqref{e:F1 and F2}.

\sssec{}

It is enough to establish the isomorphism between the two functors in question on the compact generators of $\Rep(\cG)$.
These generators can be taken to be eventually coconnective. Hence, by \corref{c:IndCoh* Op bdd below}(a),
it is enough to show that both functors are t-exact (in fact, it is sufficient to know that they are t-exact on $\Rep(\cG)^c$). 

\medskip

The t-exactness of the composition involving $\Phi_2$ is clear. For $\Phi_1$, we interpret it as
\begin{multline}  \label{e:Rep to Whit KM 8}
\Rep(\cG) \overset{(\fr^\reg)^*}\longrightarrow \QCoh(\Op^\reg_\cG) \overset{\Gamma(\Op^\reg_\cG,-)}\longrightarrow 
\CO_{\Op^\reg_\cG}\mod^{\on{com}} \overset{-\underset{\CO_{\Op^\reg_\cG}}\otimes \on{Vac}(G)_{\crit,\rho(\omega_X)}}\longrightarrow \KL(G)_{\crit,\rho(\omega_X)}\to \\
\to \hg\mod_{\crit,\rho(\omega_X)} \overset{\DS^{\on{enh,rfnd}}}\longrightarrow \IndCoh^*(\Op^\mer_\cG). 
\end{multline}  

\medskip

As was remarked already, it is enough to show that this functor is t-exact on $\Rep(\cG)^c$. The composition in the first line of 
\eqref{e:Rep to Whit KM 8} is t-exact.  Hence, by the construction of $\DS^{\on{enh,rfnd}}$, 
the corresponding functor 
$$\Rep(\cG)^c\to \IndCoh^*(\Op^\mer_\cG)$$
maps to $\IndCoh^*(\Op^\mer_\cG)^{>-\infty}$, and its t-exactness follows from the t-exactness of its composition with
$\Gamma^\IndCoh(\Op^\mer_\cG,-)^{\on{enh}}$, while the latter is the functor \eqref{e:Rep to Whit KM 6}, which is evidently
t-exact.   

\section{Unital local-to-global functors and monoidal actions} \label{s:indep}

This sections serves as a complement to \secref{s:unitality}. Here we express the notion 
of a (strictly) unital local-to-global functor as a functor of what we call the \emph{independent}
category.

\medskip

This will allow us to study the interaction between various global monoidal categories attached
to a local unital monoidal categories. These various variants are handy when studying the pattern
of the Hecke action. 

\ssec{The ``independent" category} \label{ss:indep}

\sssec{}

Let $\ul\bC^{\on{loc,untl}}$ be a crystal of categories over $\Ran^{\on{untl}}$. Denote
$$\bC^{\on{loc}}_{\Ran^{\on{untl}}}:=\Gamma^{\on{lax}}(\Ran^{\on{untl}},\ul\bC^{\on{loc,untl}}).$$

For example, when $\ul\bC^{\on{loc,untl}}$ is the unit sheaf of categories, i.e., $\ul\Dmod(\Ran^{\on{untl}})$,
the above category is $\Dmod(\Ran^{\on{untl}})$. 

\sssec{} 

Let $\bC^{\on{glob}}$ be a target DG category. On the one hand we can consider the category
$\on{Funct}^{\on{loc}\to\on{glob},\on{lax-untl}}(\ul\bC^{\on{loc}},\bC^{\on{glob}})$ of 
lax unital local-to-global functors, i.e., 
right-lax functors
$$\ul\sF^{\on{untl}}:\ul\bC^{\on{loc,untl}}\to \bC^{\on{glob}}\otimes \ul\Dmod(\Ran^{\on{untl}})$$
between sheaves of categories, see \secref{sss:loc-to-glob notation}.

\medskip

On the other hand, we can consider the category $\on{Funct}(\bC^{\on{loc}}_{\Ran^{\on{untl}}},\bC^{\on{glob}})$
of (continuous) functors
$$\sF^{\on{untl}}:\bC^{\on{loc}}_{\Ran^{\on{untl}}}\to \bC^{\on{glob}}.$$

There is a naturally defined functor
\begin{equation} \label{e:F ul lax untl to F}
\on{Funct}^{\on{loc}\to\on{glob},\on{lax-untl}}(\ul\bC^{\on{loc}},\bC^{\on{glob}})\to
\on{Funct}(\bC^{\on{loc}}_{\Ran^{\on{untl}}},\bC^{\on{glob}}).
\end{equation} 

Namely, given $\ul\sF^{\on{untl}}$, we construct $\sF^{\on{untl}}$ by applying the functor 
$\Gamma^{\on{lax}}(\Ran^{\on{untl}},-)$, followed by
$$\bC^{\on{glob}}\otimes \Dmod(\Ran^{\on{untl}})\overset{\on{Id}\otimes \on{C}^\cdot_c(\Ran^{\on{untl}},-)}\longrightarrow
\bC^{\on{glob}}.$$

\begin{rem}

Note that unlike the case of the usual Ran space, the functor \eqref{e:F ul lax untl to F} is \emph{not} an equivalence.

\medskip

For example, for $$\ul\bC^{\on{loc,untl}}=\ul\Dmod(\Ran^{\on{untl}})
\text{ and }\bC^{\on{glob}}=\Vect,$$ we have
$$\on{Funct}^{\on{loc}\to\on{glob},\on{lax-untl}}(\ul\bC^{\on{loc}},\bC^{\on{glob}})\simeq \Dmod(\Ran^{\on{untl}}),$$
while $\on{Funct}(\bC_{\Ran^{\on{untl}}},\bC^{\on{glob}})$ identifies with $\Dmod(\Ran^{\on{untl}})^\vee$. In terms
of this identification, the functor \eqref{e:F ul untl to F} is the pairing
$$\Dmod(\Ran^{\on{untl}})\otimes \Dmod(\Ran^{\on{untl}})\to \Vect, \quad \CF_1,\CF_2\mapsto
\on{C}^\cdot_c(\Ran^{\on{untl}},\CF_1\sotimes \CF_2).$$

However, this pairing is not perfect. The reason for this is that, although the functor
$(\Delta_{\Ran^{\on{untl}}})_!$ is defined, it does \emph{not} satisfy the projection formula.

\medskip

More generally, the category $\on{Funct}^{\on{loc}\to\on{glob},\on{lax-untl}}(\ul\bC^{\on{loc}},\bC^{\on{glob}})$ can be described
explicitly as lax sections of another crystal of categories on $\Ran^{\on{untl}}$, see Remark \ref{r:dual of lax sections}.

\end{rem} 

\sssec{} \label{sss:right adj categ prestack}

From now on we will make the following assumption on $\ul\bC^{\on{loc,untl}}$: 

\medskip

For every $S\in \affSch$ and a map $\ul{x}_1\overset{\alpha}\to \ul{x}_2$ in $\Maps(S,\Ran^{\on{untl}})$,
the corresponding functor
$$\bC^{\on{loc}}_{S,\ul{x}_1} \overset{\on{ins.unit}_{\ul{x}_1\subseteq \ul{x}_2}}\to \bC^{\on{loc}}_{S,\ul{x}_2}$$
admits a \emph{continuous} right adjoint.

\medskip

Note that in this case, this right adjoint is automatically $\Dmod(S)$-linear.

\begin{rem} \label{r:shvs of cat on Ran via proper}
This assumption is made in order to simplify the exposition. One can make do without it,
but in what follows one will have to describe $\ul\bC^{\on{loc,untl}}$ using \emph{proper}
(rather than affine) schemes mapping to $\Ran$, see \secref{sss:proper approx}. 

\medskip

The properness assumption would guarantee that for map $f:Z'\to Z$ 
(of proper schemes) and 
$$(\ul{x}_1\to \ul{x}_2)\in \Maps(Z,\Ran^{\on{untl}}),$$
the diagram
$$
\CD
\bC_{Z',\ul{x_1}\circ f}  @<{\on{ins.unit}_{\ul{x_1}\circ f\subseteq \ul{x_2}\circ f}^R}<< \bC_{Z',\ul{x_2}\circ f}  \\
@A{f^!}AA @AA{f^!}A \\
\bC_{Z,\ul{x_1}}  @<{\on{ins.unit}_{\ul{x_1}\subseteq \ul{x_2}}^R}<< \bC_{Z,\ul{x_2}} 
\endCD
$$
obtained by passing to right adjoints in the commutative diagram
$$
\CD
\bC_{Z',\ul{x_1}\circ f}  @<{\on{ins.unit}_{\ul{x_1}\circ f\subseteq \ul{x_2}\circ f}}<< \bC_{Z',\ul{x_2}\circ f}  \\
@A{f_!}AA @AA{f_!}A \\
\bC_{Z,\ul{x_1}}  @<{\on{ins.unit}_{\ul{x_1}\subseteq \ul{x_2}}}<< \bC_{Z,\ul{x_2}}.
\endCD
$$
commutes. 
 
\end{rem} 

\sssec{} \label{sss:indep via duality}

Under the above assumption, passage to right adjoints defines a sheaf of categories
$$(\ul\bC^{\on{loc,untl}})^{\on{op}}$$
over $(\Ran^{\on{untl}})^{\on{op}}$.

\medskip

Note that we can tautologically identify
$\bC^{\on{loc}}_{\Ran^{\on{untl}}}:=\Gamma^{\on{lax}}(\Ran^{\on{untl}},\ul\bC^{\on{loc,untl}})$
with
$$\Gamma^{\on{lax}}((\Ran^{\on{untl}})^{\on{op}},(\ul\bC^{\on{loc,untl}})^{\on{op}}).$$

\medskip

Set
$$\bC^{\on{loc}}_{\Ran^{\on{untl}},\on{indep}}:=\Gamma^{\on{strict}}((\Ran^{\on{untl}})^{\on{op}},(\ul\bC^{\on{loc,untl}})^{\on{op}})\subset
\Gamma^{\on{lax}}((\Ran^{\on{untl}})^{\on{op}},(\ul\bC^{\on{loc,untl}})^{\on{op}}).$$

\sssec{}

Let us describe $\bC^{\on{loc}}_{\Ran^{\on{untl}},\on{indep}}$ explicitly as a full subcategory of $\bC^{\on{loc}}_{\Ran^{\on{untl}}}$:

\medskip

An object
$$(S\overset{\ul{x}}\to \Ran)\mapsto \bc_{S,\ul{x}}\in \bC^{\on{loc}}_{S,\ul{x}}$$
belongs to $\bC^{\on{loc}}_{\Ran^{\on{untl}},\on{indep}}$ if for every 
$$(\ul{x}_1\overset{\alpha}\to \ul{x}_2) \in \Maps(S,\Ran^{\on{untl}})$$
the map
$$\bc_{S,\ul{x}_1}\to (\on{ins.unit}_{\ul{x}_1\subseteq \ul{x}_2})^R(\bc_{S,\ul{x}_2}),$$
obtained by adjunction from
$$\on{ins.unit}_{\ul{x}_1\subseteq \ul{x}_2}(\ul{x}_1)\to \ul{x}_2,$$
is an isomorphism.

\begin{rem} \label{r:dual of lax sections}

Let assume in addition that $\ul\bC^{\on{loc,untl}}$ is value-wise dualizable. Then passing to duals in  
$(\ul\bC^{\on{loc,untl}})^{\on{op}}$, we obtain a crystal of categories $((\ul\bC^{\on{loc,untl}})^{\on{op}})^\vee$
on $\Ran^{\on{untl}}$.

\medskip

It is easy to see that the category 
$$\on{Funct}^{\on{loc}\to\on{glob},\on{lax-untl}}(\ul\bC^{\on{loc}},\bC^{\on{glob}})$$
identifies with
$$((\bC^{\on{op}})^\vee)_{\Ran^{\on{untl}}}:=
\Gamma^{\on{lax}}(\Ran^{\on{untl}},((\ul\bC^{\on{loc,untl}})^{\on{op}})^\vee\otimes \bC^{\on{glob}}).$$

\end{rem}

\sssec{}

We claim:

\begin{lem} \label{l:adj of lax sections}
The embedding
$$\on{emb.indep}:\bC^{\on{loc}}_{\Ran^{\on{untl}},\on{indep}}\hookrightarrow \bC^{\on{loc}}_{\Ran^{\on{untl}}}$$
admits a left adjoint. 
\end{lem}

\begin{proof}

It is sufficient to show that 
$$\bC^{\on{loc}}_{\Ran^{\on{untl}},\on{indep}}\subset \bC^{\on{loc}}_{\Ran^{\on{untl}}}$$
is closed under limits.

\medskip

However, this follows from the fact that limits in $\bC_{\Ran^{\on{untl}}}$
exists and have the property that they commute with evaluation on every
\emph{proper} $Z$ mapping to $\Ran^{\on{untl}}$, see Remark 
\ref{r:shvs of cat on Ran via proper}.

\end{proof} 

\sssec{}

Let $\on{emb.indep}^L$ denote the left adjoint of $\on{emb.indep}$. Thus, we can view 
$\bC^{\on{loc}}_{\Ran^{\on{untl}},\on{indep}}$ as a localization of $\bC^{\on{loc}}_{\Ran^{\on{untl}}}$.

\medskip

Let us view the category $\on{Funct}(\bC^{\on{loc}}_{\Ran^{\on{untl}},\on{indep}},\bC^{\on{glob}})$
as a full subcategory of $\on{Funct}(\bC^{\on{loc}}_{\Ran^{\on{untl}}},\bC^{\on{glob}})$ via precomposition
with $\on{emb.indep}^L$. 

\sssec{}

We claim:

\begin{lem} \label{l:dual of lax sections}
The functor \eqref{e:F ul lax untl to F} sends
$$\on{Funct}^{\on{loc}\to\on{glob},\on{untl}}(\ul\bC^{\on{loc}},\bC^{\on{glob}})\subset 
\on{Funct}^{\on{loc}\to\on{glob},\on{lax-untl}}(\ul\bC^{\on{loc}},\bC^{\on{glob}})$$
to
$$\on{Funct}(\bC^{\on{loc}}_{\Ran^{\on{untl}},\on{indep}},\bC^{\on{glob}}) \subset 
\on{Funct}(\bC^{\on{loc}}_{\Ran^{\on{untl}}},\bC^{\on{glob}}),$$
and the resulting functor 
\begin{equation} \label{e:F ul untl to F}
\on{Funct}^{\on{loc}\to\on{glob},\on{untl}}(\ul\bC^{\on{loc}},\bC^{\on{glob}})\to
\on{Funct}(\bC^{\on{loc}}_{\Ran^{\on{untl}},\on{indep}},\bC^{\on{glob}})
\end{equation} 
is an equivalence.
\end{lem}

\begin{proof}

We will access $\Ran^{\on{untl}}$ via proper schemes mapping to it, see \secref{sss:proper approx}.

\medskip

Let us be given an object $\ul\sF\in \on{Funct}^{\on{loc}\to\on{glob}}(\ul\bC^{\on{loc}},\bC^{\on{glob}})$.
For every proper $Z$ equipped with a map $\ul{x}:Z\to \Ran$, consider the corresponding functor
$$\sF_{\int_Z,\ul{x}}:\bC^{\on{loc}}_{Z,\ul{x}}\to \bC^{\on{glob}}$$
and its (not necessarily continuous) right adjoint
$$(\sF_{\int_Z,\ul{x}})^R: \bC^{\on{glob}}\to \bC^{\on{loc}}_{Z,\ul{x}}.$$

\medskip

For $Z'\overset{g}\to Z$, the diagram 
\begin{equation} \label{e:proper diagram 1}
\CD
\bC^{\on{loc}}_{Z',\ul{x}\circ g} @>{\sF_{Z',\ul{x}\circ g}}>> \bC^{\on{glob}}\otimes \Dmod(Z') \\
@A{g^!}AA @AA{\on{id}\otimes g^!}A \\
\bC^{\on{loc}}_{Z,\ul{x}} @>{\sF_{Z,\ul{x}}}>> \bC^{\on{glob}} \otimes \Dmod(Z)
\endCD
\end{equation}
is equipped with a datum of commutativity. Since $g$ is proper, the diagram
\begin{equation} \label{e:proper diagram 2}
\CD
\bC^{\on{loc}}_{Z',\ul{x}\circ g} @>{\sF_{Z',\ul{x}\circ g}}>> \bC^{\on{glob}}\otimes \Dmod(Z')  \\
@V{g_!}VV @VV{\on{id}\otimes g_!}V \\
\bC^{\on{loc}}_{Z,\ul{x}} @>{\sF_{Z,\ul{x}}}>> \bC^{\on{glob}} \otimes \Dmod(Z'),
\endCD
\end{equation}
obtained from \eqref{e:proper diagram 1} by passing to left adjoints along the vertical arrows,
also commutes.

\medskip

From \eqref{e:proper diagram 2}, we obtain a datum of commutativity for the diagram
\begin{equation} \label{e:proper diagram 2'}
\CD
\bC^{\on{loc}}_{Z',\ul{x}\circ g} @>{\sF_{\int_{Z'},\ul{x}\circ g}}>> \bC^{\on{glob}} \\
@V{g_!}VV @VV{\on{id}}V \\
\bC^{\on{loc}}_{Z,\ul{x}} @>{\sF_{\int_{Z},\ul{x}}}>> \bC^{\on{glob}}.
\endCD
\end{equation}

Finally, by passing to right adjoints in \eqref{e:proper diagram 2'}, we obtain a 
datum of commutativity for the diagram
\begin{equation} \label{e:proper diagram 3}
\CD
\bC^{\on{loc}}_{Z',\ul{x}\circ g} @<{(\sF_{\int_{Z'},\ul{x}\circ g})^R}<< \bC^{\on{glob}} \\
@A{g^!}AA @AA{\on{id}}A \\
\bC^{\on{loc}}_{Z,\ul{x}} @<{(\sF_{\int_{Z},\ul{x}}})^R<< \bC^{\on{glob}},
\endCD
\end{equation}
where the functors are $(\sF_{\int_Z,\ul{x}})^R$ are not necessarily continuous,
but limit-preserving. 

\medskip

The functor $\ul\sF$ is determined by the data of $\sF_{\int_Z,\ul{x}}$ plus the 
data of commutativity for the diagrams \eqref{e:proper diagram 2'}, which is equivalent to
having the data of limit-preserving functors $(\sF_{\int_Z,\ul{x}})^R$ plus the data 
of commutativity for the diagrams \eqref{e:proper diagram 3}. 

\medskip

Suppose now that $\ul\sF$ is upgraded to an object $\ul\sF^{\on{untl}}$ of
$\on{Funct}^{\on{loc}\to\on{glob},\on{lax-untl}}(\ul\bC^{\on{loc}},\bC^{\on{glob}})$.

\medskip

Let us be given a map $\ul{x}_1\to \ul{x}_2$ in $\Maps(Z,\Ran^{\on{untl}})$. The unital structure on
$\ul\bC^{\on{loc}}$ gives rise to a functor
$$\on{ins.unit}_{\ul{x}_1\subseteq \ul{x}_2}:\bC^{\on{loc}}_{Z,\ul{x}_1}\to \bC^{\on{loc}}_{Z,\ul{x}_2},$$
and the lax unital structure $\ul\sF^{\on{untl}}$ on $\ul\sF$ gives rise to a natural transformation
$$\sF_{\int_Z,\ul{x}_1}\to \sF_{\int_Z,\ul{x}_2}\circ \on{ins.unit}_{\ul{x}_1\subseteq \ul{x}_2}.$$

By definition, this natural transformation is an isomorphism if and only if the above lax unital structure on 
$\ul\sF$ is strict. By adjunction, we obtain a natural transformation
$$(\on{ins.unit}_{\ul{x}_1\subseteq \ul{x}_2})^R\circ (\sF_{\int_Z,\ul{x}_2})^R\to (\sF_{\int_Z,\ul{x}_1})^R,$$
which is an isomorphism if and only if $\ul\sF^{\on{untl}}$ is strictly unital. 

\medskip

Since for every $(Z,\ul{x})$, the evaluation functor 
$$\bC^{\on{loc}}_{\Ran^{\on{untl}},\on{indep}}\hookrightarrow \bC^{\on{loc}}_{\Ran^{\on{untl}}}\to \bC^{\on{loc}}_{Z,\ul{x}}$$
commutes with limits, we obtain that the datum of a strictly unital object $\ul\sF^{\on{untl}}$ is equivalent to
that of a limit-preserving functor
$$\bC^{\on{glob}}\to \bC^{\on{loc}}_{\Ran^{\on{untl}},\on{indep}},$$
which is equivalent to that of a continuous functor
$$\bC^{\on{loc}}_{\Ran^{\on{untl}},\on{indep}}\to \bC^{\on{glob}}.$$

Thus, we have constructed an equivalence
$$\on{Funct}^{\on{loc}\to\on{glob},\on{untl}}(\ul\bC^{\on{loc}},\bC^{\on{glob}})\simeq 
\on{Funct}(\bC^{\on{loc}}_{\Ran^{\on{untl}},\on{indep}},\bC^{\on{glob}}).$$

Unwinding the construction, it is easy to see that the diagram
$$
\CD
\on{Funct}^{\on{loc}\to\on{glob},\on{untl}}(\ul\bC^{\on{loc}},\bC^{\on{glob}}) @>>> 
\on{Funct}^{\on{loc}\to\on{glob},\on{lax-untl}}(\ul\bC^{\on{loc}},\bC^{\on{glob}}) \\
@V{\sim}VV @VVV \\
\on{Funct}(\bC^{\on{loc}}_{\Ran^{\on{untl}},\on{indep}},\bC^{\on{glob}}) @>>> \on{Funct}(\bC^{\on{loc}}_{\Ran^{\on{untl}}},\bC^{\on{glob}})
\endCD
$$
commutes.

\end{proof} 

\begin{rem}
The discussion in this subsection is not specific to $\Ran^{\on{untl}}$. It applies to any pseudo-proper 
categorical prestack. 
\end{rem}

\ssec{The calculation of the independent category in the vacuum case}

\sssec{}

Here is a sample calculation of the category $\bC^{\on{loc}}_{\Ran^{\on{untl}},\on{indep}}$. 
Take $\ul\bC^{\on{loc,untl}}:=\ul\Dmod(\Ran^{\on{untl}})$. Denote the corresponding independent
category by
$$\Vect_{\Ran^{\on{untl}},\on{indep}}.$$

\sssec{}

Take $\bC^{\on{glob}}=\Vect$. Using \lemref{l:dual of lax sections}, from the identity functor
$$\ul\Dmod(\Ran^{\on{untl}})\to \ul\Dmod(\Ran^{\on{untl}}),$$
we obtain a functor

\begin{equation} \label{e:indep Vect}
\Vect_{\Ran^{\on{untl}},\on{indep}}\to \Vect. 
\end{equation} 

We claim:

\begin{prop} \label{p:indep Vect}
The functor \eqref{e:indep Vect} is an equivalence.
\end{prop}

\begin{rem} \label{r:indep Vect}

Note that by \secref{sss:indep via duality} can be reformulated as saying that the functor
\begin{equation} \label{e:indep Vect 1}
\Vect\overset{k\mapsto \omega_{\Ran^{\on{untl}}}}\longrightarrow 
\Gamma^{\on{strict}}(\Ran^{\on{untl}},\ul\Dmod(\Ran^{\on{untl}}))
\end{equation}
is an equivalence. 

\end{rem} 

\begin{rem}

Note that, unlike the fact that $\Ran$ is contractible, which requires $X$ to be \emph{connected}, the assertion
of \propref{p:indep Vect} is valid for any $X$ that is non-empty.

\end{rem} 

\begin{proof}[Proof of \propref{p:indep Vect}]

We will show that the functor \eqref{e:indep Vect 1} is an equivalence. The right adjoint of this functor
is given by the restriction of $\on{C}^\cdot_c(\Ran^{\on{untl}},-)$ to 
$$\Gamma^{\on{strict}}(\Ran^{\on{untl}},\ul\Dmod(\Ran^{\on{untl}}))\subset 
\Gamma^{\on{lax}}(\Ran^{\on{untl}},\ul\Dmod(\Ran^{\on{untl}}))=\Dmod(\Ran^{\on{untl}}).$$

\medskip

It suffices to show that for $\CF\in \Gamma^{\on{strict}}(\Ran^{\on{untl}},\ul\Dmod(\Ran^{\on{untl}}))$, the map
$$\CF\to \on{C}^\cdot_c(\Ran^{\on{untl}},\CF)\otimes \omega_{\Ran^{\on{untl}}}$$
is an isomorphism.

\medskip

For that it suffices to show that for any $\ul{x}\in \Ran$, the map
\begin{equation} \label{e:indep Vect 2}
\CF_{\ul{x}}\to  \on{C}^\cdot_c(\Ran^{\on{untl}},\CF)
\end{equation}
corresponding to 
\begin{equation} \label{e:indep Vect 3}
\on{pt} \overset{\ul{x}}\longrightarrow \Ran^{\on{untl}},
\end{equation}
is an isomorphism.

\medskip

We factor \eqref{e:indep Vect 2} as
\begin{equation} \label{e:indep Vect 4}
\on{pt}\to \Ran^{\on{untl}}_{\ul{x}}\to \Ran^{\on{untl}},
\end{equation}
and hence \eqref{e:indep Vect 3} as 
\begin{equation} \label{e:indep Vect 5}
\CF_{\ul{x}}\to   \on{C}^\cdot_c(\Ran^{\on{untl}}_{\ul{x}},\CF|_{\Ran^{\on{untl}}_{\ul{x}}})\to   \on{C}^\cdot_c(\Ran^{\on{untl}},\CF).
\end{equation}

\medskip

We claim that both maps in \eqref{e:indep Vect 5} are isomorphisms. Indeed, the first map is an isomorphism, becauase
$$ \CF|_{\Ran^{\on{untl}}_{\ul{x}}}\in \Gamma^{\on{strict}}(\Ran^{\on{untl}},\ul\Dmod(\Ran_{\ul{x}}^{\on{untl}})),$$
and $\ul{x}$ is the (value-wise) initial point of $\Ran^{\on{untl}}_{\ul{x}}$.

\medskip

The second arrow in \eqref{e:indep Vect 5} is an isomorphism because the map
$$\Ran^{\on{untl}}_{\ul{x}}\to \Ran^{\on{untl}}$$
is value-wise cofinal: its value-wise left adjoint is given by $\ul{x}'\mapsto \ul{x}\cup \ul{x}'$. 

\end{proof}

\sssec{} \label{sss:extension to the emptyset}

Let us rerturn to \propref{p:local-to-global empty set}, and explain its meaning in terms of 
local-to-global functors.

\medskip

Let us be given a strictly unital functor
$$\ul\sF^{\on{untl}}:\ul\bC^{\on{loc,untl}}\to \bC^{\on{glob}}\otimes \ul\Dmod(\Ran^{\on{untl}}).$$

Assume now that $\ul\bC^{\on{loc,untl}}$ is the restriction of a sheaf of categories $\bC^{\on{loc,untl},*}$
on $\Ran^{\on{untl},*}$ (as is the case in all our examples). 

\medskip

Le $\bC^{\on{loc}}_\emptyset$ be the value of this extension on the initial point. The entire datum of
the extension
$$\ul\bC^{\on{loc,untl}}\rightsquigarrow \bC^{\on{loc,untl},*}$$
is equivalent to the datum of a functor
\begin{equation} \label{e:unit object abs}
\bC^{\on{loc}}_\emptyset \overset{\on{ins.unit}_\emptyset}\longrightarrow \Gamma^{\on{strict}}(\Ran^{\on{untl}},\ul\bC^{\on{loc,untl}})\hookrightarrow
\Gamma^{\on{lax}}(\Ran^{\on{untl}},\ul\bC^{\on{loc,untl}})=\bC^{\on{loc}}_{\Ran^{\on{untl}}}.
\end{equation}

\begin{rem}

Note that in most examples, $\bC^{\on{loc}}_\emptyset\simeq \Vect$, so that datum of \eqref{e:unit object abs}
is that of an object
$$\one_{\ul\bC^{\on{loc}}}\in \Gamma^{\on{strict}}(\Ran^{\on{untl}},\ul\bC^{\on{loc,untl}}).$$

I.e., for every $\ul{x}\in \Ran$, we have an object $\one_{\ul\bC^{\on{loc}},\ul{x}}\in \bC^{\on{loc}}_{\ul{x}}$
and for every $\ul{x}\subseteq \ul{x}'$ we have an isomorphism
$$\on{ins.unit}_{\ul{x}\subseteq \ul{x}'}(\one_{\ul\bC^{\on{loc}},\ul{x}}\simeq \one_{\ul\bC^{\on{loc}},\ul{x}'}.$$

\end{rem}

\sssec{}

Consider the composition
$$\ul\sF^{\on{untl}}\circ \on{ins.unit}_\emptyset: \bC^{\on{loc}}_\emptyset\to \bC^{\on{glob}}\otimes \Dmod(\Ran^{\on{untl}}).$$

The claim is that it factors canonically as
$$\bC^{\on{loc}}_\emptyset\overset{\sF_\emptyset}\to \bC^{\on{glob}} \overset{\on{Id}\otimes \omega_{\Ran^{\on{untl}}}}\longrightarrow 
\bC^{\on{glob}}\otimes \Dmod(\Ran^{\on{untl}}).$$

Indeed, the functor $\ul\sF^{\on{untl}}\circ \on{ins.unit}_\emptyset$ factors as
\begin{multline*} 
\bC^{\on{loc}}_\emptyset\to 
\bC^{\on{glob}}\otimes \Gamma^{\on{strict}}(\Ran^{\on{untl}},\ul\Dmod(\Ran^{\on{untl}}))\to \\
\to \bC^{\on{glob}}\otimes \Gamma^{\on{lax}}(\Ran^{\on{untl}},\ul\Dmod(\Ran^{\on{lax}}))=\bC^{\on{glob}}\otimes \Dmod(\Ran^{\on{untl}}),
\end{multline*} 
and according to Remark \ref{r:indep Vect}, the functor
$$\bC^{\on{glob}}\to \bC^{\on{glob}}\otimes \Gamma^{\on{strict}}(\Ran^{\on{untl}},\ul\Dmod(\Ran^{\on{untl}}))$$
is an equivalence.

\ssec{Non-unitality vs independence}

Let $\ul\bC^{\on{loc,untl}}$ be as in the previous subsection. 

\medskip

In this subsection we will utilize the contractibility of the Ran space to explain the relation between 
the ``indepndent" category $\bC^{\on{loc}}_{\Ran^{\on{untl}},\on{indep}}$ and the non-unital version $\bC^{\on{loc}}_\Ran$.

\sssec{}

We claim:

\begin{prop} \label{p:Ran vs indep}
The composite functor
$$\bC^{\on{loc}}_{\Ran^{\on{untl}},\on{indep}}\overset{\on{emb.indep}}\hookrightarrow \bC^{\on{loc}}_{\Ran^{\on{untl}}}\overset{\sft^!}\to
\bC^{\on{loc}}_\Ran$$
is fully faithful.
\end{prop}

\begin{proof}

The assertion of the proposition amounts to the following. Let us be given two objects
$$\bc',\bc''\in \bC^{\on{loc}}_{\Ran^{\on{untl}},\on{indep}}\subset \bC^{\on{loc}}_{\Ran^{\on{untl}}},$$
and a map 
$$\sft^!(\bc')\overset{\phi}\to \sft^!(\bc'').$$

We need to show that this map can be uniquely upgraded to a map
$$\bc'\overset{\phi_{\on{untl}}}\to \bc''.$$

This amounts to the following. Fix $S\in \affSch$ and let us be given a map
$$(\ul{x}_1\to  \ul{x}_2)\in \Maps(S,\Ran^{\on{untl}}).$$

We need to equip the following diagram (taking place in in $\bC^{\on{loc}}_{S,\ul{x}_2}$)
$$
\CD
\bc'_{S,\ul{x}_2} @>{\phi_{S,\ul{x}_2}}>> \bc''_{S,\ul{x}_2} \\
@AAA @AAA \\
\on{ins.unit}_{\ul{x}_1\subseteq \ul{x}_2}(\bc'_{S,\ul{x}_1}) @>{\on{ins.unit}_{\ul{x}_1\subseteq \ul{x}_2}(\phi_{S,\ul{x}_1})}>> 
\on{ins.unit}_{\ul{x}_1\subseteq \ul{x}_2}(\bc''_{S,\ul{x}_1})
\endCD
$$
(where the vertical arrows are given by the structure of objects of $\bC^{\on{loc}}_{\Ran^{\on{untl}}}$
on $\bc'$ and $\bc''$, respectively) with a datum of commutativity.

\medskip

The datum of commutativity of the above diagram is equivalent to the datum of commutativity
of the following diagram (taking place in in $\bC_{S,\ul{x}_1}$):
\begin{equation} \label{e:indep diag}
\CD
\on{ins.unit}_{\ul{x}_1\subseteq \ul{x}_2}^R(\bc'_{S,\ul{x}_2}) @>{\on{ins.unit}_{\ul{x}_1\subseteq \ul{x}_2}^R(\phi_{S,\ul{x}_2})}>> 
\on{ins.unit}_{\ul{x}_1\subseteq \ul{x}_2}^R(\bc''_{S,\ul{x}_2}) \\
@AAA @AAA \\
\bc'_{S,\ul{x}_1} @>{\phi_{S,\ul{x}_1}}>> \bc''_{S,\ul{x}_1}.
\endCD
\end{equation} 

Note the vertical arrows in \eqref{e:indep diag} are isomorphisms, by the assumption that
$\bc',\bc''\in \bC^{\on{loc}}_{\Ran^{\on{untl}},\on{indep}}$. Thus, we can view both circuits in \eqref{e:indep diag}
as an object of $\CHom_{\bC^{\on{loc}}_{S,\ul{x}_1}}(\bc'_{S,\ul{x}_1},\bc''_{S,\ul{x}_1})$. 

\medskip

Thus, letting $\ul{x}_2$ vary, we can view the clockwise circuit in \eqref{e:indep diag} as a map
\begin{equation} \label{e:indep diag 1}
(\on{Id}\otimes \on{pr}_{\on{small},S}^!)(\bc'_{S,\ul{x}_1})\to (\on{Id}\otimes \on{pr}_{\on{small},S})^!(\bc'_{S,\ul{x}_1})
\end{equation} 
in 
$$\bC^{\on{loc}}_{S,\ul{x}_1}\underset{\Dmod(S)}\otimes \Dmod(S_{\ul{x}_1}^{\subseteq}).$$

\medskip

Now, the universal homological contractibility of the map
$$\on{pr}_{\on{small},S}:S_{\ul{x}_1}^{\subseteq}\to S$$
implies that the functor 
$$\on{pr}_{\on{small},S}^!:\Dmod(S)\to \Dmod(S_{\ul{x}_1}^{\subseteq})$$
is fully faithful, and hence so is
$$\on{Id}\otimes \on{pr}_{\on{small},S}^!:\bC^{\on{loc}}_{S,\ul{x}_1}\to
\bC^{\on{loc}}_{S,\ul{x}_1}\underset{\Dmod(S)}\otimes \Dmod(S_{\ul{x}_1}^{\subseteq}).$$

\medskip

Hence, the space 
$$\CHom_{\bC^{\on{loc}}_{S,\ul{x}_1}\underset{\Dmod(S)}\otimes \Dmod(S_{\ul{x}_1}^{\subseteq})}
((\on{Id}\otimes \on{pr}_{\on{small},S}^!)(\bc'_{S,\ul{x}_1}),(\on{Id}\otimes \on{pr}_{\on{small},S}^!)(\bc''_{S,\ul{x}_1}))$$
is isomorphic to $\CHom^{\on{loc}}_{\bC_{S,\ul{x}_1}}(\bc'_{S,\ul{x}_1},\bc''_{S,\ul{x}_1})$, with the mutually inverse isomorphisms given by
the functors $\on{Id}\otimes \on{pr}_{\on{small},S}^!$ and $\on{Id}\otimes \on{diag}_S^!$, respectively. 

\medskip

This equips the family of diagrams \eqref{e:indep diag} with a unique datum of commutativity as 
$\ul{x}_2$ varies over $S_{\ul{x}_1}^{\subseteq}$. 

\end{proof} 

\sssec{} \label{sss:unitality as a property Take 2}

Combining \propref{p:Ran vs indep} and \lemref{l:dual of lax sections}, we obtain:  

\begin{cor} \label{c:Ran vs indep}
The functor from the category of \emph{strict} functors of crystals of categories
$$\ul\sF^{\on{untl}}:\ul\bC^{\on{loc,untl}}\to \bC^{\on{glob}}\otimes \ul\Dmod(\Ran^{\on{untl}})$$
to the category of functors
$$\ul\sF:\ul\bC^{\on{loc}}\to \bC^{\on{glob}}\otimes \ul\Dmod(\Ran),$$
given by restriction along $\Ran\to \Ran^{\on{untl}}$, is fully faithful.
\end{cor}

Unwinding the definitions, we obtain that the essential image of the above fully faithful functor 
$$\on{Funct}^{\on{loc}\to \on{glob},\on{untl}}(\ul\bC^{\on{loc}},\bC^{\on{glob}})\to
\on{Funct}^{\on{loc}\to \on{glob},\on{lax-untl}}(\ul\bC^{\on{loc}},\bC^{\on{glob}})\to 
\on{Funct}^{\on{loc}\to \on{glob}}(\ul\bC^{\on{loc}},\bC^{\on{glob}})$$
consists of objects that have a global unitality property.

\medskip

I.e., this proves \propref{p:unitality as a property}. 

\sssec{}

As another immediate corollary of \propref{p:Ran vs indep} we obtain:

\begin{cor} \label{c:Ran vs indep 1}
The natural transformation
$$\on{emb.indep}^L\circ \sft_!\circ \sft^!\to \on{emb.indep}^L$$
is an isomorphism.
\end{cor} 

\sssec{}

Finally, we record:

\begin{cor} \label{c:Ran vs indep 2}
The composite functor
$$\bC_\Ran \overset{\sft_!}\to \bC_{\Ran^{\on{untl}}} \overset{\on{emb.indep}^L}\to
\bC_{\Ran^{\on{untl}},\on{indep}}$$
is a localization.
\end{cor}  

\ssec{Sheaves of monoidal categories on the unital Ran space} \label{ss:monoidal on untl Ran}

Let $\bA^{\on{loc,untl}}$ be a sheaf of unital monoidal categories over $\Ran^{\on{untl}}$.
We will assume that $\bA^{\on{loc,untl}}$ satisfies the condition from \secref{sss:right adj categ prestack}.

\medskip

We will also assume that the monoidal operation
admits a right adjoint, which is a strict functor between sheaves of categories. 

\medskip

In this subsection we will study the categories 
\begin{equation} \label{e:monoidal Ran}
\bA_{\Ran^{\on{untl}}}:=\Gamma^{\on{lax}}(\Ran^{\on{untl}},\bA^{\on{loc,untl}}) 
 \text{ and } \bA_{\Ran^{\on{untl,indep}}}.
\end{equation}

We will equip them with monoidal structures, and study the interactions between them. 

\sssec{Example}

An example of such $\bA^{\on{loc,untl}}$ is the sheaf of unital monoidal categories
over $\Ran^{\on{untl}}$ attached to a unital monoidal factorization category $\bA$.

\medskip

An important example of such an $\bA$ is $\Sph_G$. 

\sssec{} \label{sss:constant monoidal}

Another example is the \emph{symmetric} monoidal factorization category $\bA$ associated to a 
crystal of symmetric monoidal categories over $X$ (this is a categorical counterpart of the procedure from \secref{ss:com untl fact alg}). 

\medskip

An example of this is the constant crystal of symmetric monoidal categories over $X$ with fiber $\Rep(\cG)$. 

\begin{rem}

Note that being a sheaf of \emph{unital} monoidal categories over $\Ran^{\on{untl}}$ automatically
imposes a condition on the compatibility between the monoidal unit and a structure of sheaf
of categories: 

\medskip

The monoidal unit
$$\one_{\bA,\Ran}\in \Gamma^{\on{lax}}(\Ran^{\on{untl}}, \bA^{\on{loc,untl}})=\bA_{\Ran^{\on{untl}}}$$
belongs to 
$$\Gamma^{\on{strict}}(\Ran^{\on{untl}}, \bA^{\on{loc,untl}})\subset
\Gamma^{\on{lax}}(\Ran^{\on{untl}}, \bA^{\on{loc,untl}}).$$

\end{rem}

\sssec{}

Consider the map
$$\on{union}:\Ran^{\on{untl}}\times \Ran^{\on{untl}}\to \Ran^{\on{untl}}.$$

We have the 1-morphisms
$$p_1\to \on{union} \leftarrow p_2$$
in the category of maps from $\Ran^{\on{untl}}\times \Ran^{\on{untl}}$ to $\Ran^{\on{untl}}$. From here we obtain functors
$$p_1^*(\bA^{\on{loc,untl}})\to \on{union}^*(\bA^{\on{loc,untl}})\leftarrow p_2^*(\bA^{\on{loc,untl}})$$
as crystals of categories over $\Ran^{\on{untl}}\times \Ran^{\on{untl}}$, and hence a functor
\begin{equation} \label{e:conv op unital -2}
\bA^{\on{loc,untl}}\boxtimes \bA^{\on{loc,untl}} \to \on{union}^*(\bA^{\on{loc,untl}}\otimes \bA^{\on{loc,untl}}).
\end{equation} 

Combining with the monoidal operation on $\bA^{\on{loc,untl}}$, we obtain a functor
\begin{equation} \label{e:conv op unital -1}
\bA^{\on{loc,untl}}\boxtimes \bA^{\on{loc,untl}}\to \on{union}^*(\bA^{\on{loc,untl}}).
\end{equation} 

Applying $\Gamma^{\on{lax}}(\Ran^{\on{untl}}\times \Ran^{\on{untl}},-)$, from \eqref{e:conv op unital -1} we obtain a functor
\begin{equation} \label{e:conv op unital 0}
\bA_{\Ran^{\on{untl}}}\otimes \bA_{\Ran^{\on{untl}}}\to 
\Gamma^{\on{lax}}(\Ran^{\on{untl}}\times \Ran^{\on{untl}},\on{union}^*(\bA^{\on{loc,untl}})).
\end{equation}

Finally, composing \eqref{e:conv op unital 0} with the functor 
$$\on{union}_!:\Gamma^{\on{lax}}(\Ran^{\on{untl}}\times \Ran^{\on{untl}},\on{union}^*(\bA^{\on{loc,untl}}))\to
\Gamma^{\on{lax}}(\Ran^{\on{untl}},\bA^{\on{loc,untl}}),$$
left adjoint to $\on{union}^!$ (it exists thanks to \corref{c:!-dir image ps-proper}), we obtain a functor 
\begin{equation} \label{e:conv op unital}
\bA_{\Ran^{\on{untl}}}\otimes \bA_{\Ran^{\on{untl}}}\to \bA_{\Ran^{\on{untl}}}.
\end{equation}

We will denote the resulting binary operation on $\bA_{\Ran^{\on{untl}}}$ by $\star$, and will refer to it
as ``convolution".

\sssec{}

The above binary operation extends to a monoidal structure on $\bA_{\Ran^{\on{untl}}}$,
which we will refer to us the \emph{convolution} monoidal structure. We will denote $\bA_{\Ran^{\on{untl}}}$,
viewed as a monoidal category equipped with the convolution structure by $\bA^\star_{\Ran^{\on{untl}}}$.

\sssec{} \label{sss:Ran unital pointwise}

Note now that the monoidal operation on $\bA^{\on{loc,untl}}$ defines a \emph{pointwise} monoidal
structure on $\bA_{\Ran^{\on{untl}}}$. 

\medskip

We will denote $\bA_{\Ran^{\on{untl}}}$,
viewed as a monoidal category equipped with the pointwise structure by $\bA^{\sotimes}_{\Ran^{\on{untl}}}$.

\medskip

Note that $\bA^{\sotimes}_{\Ran^{\on{untl}}}$ is unital: its unit is the object $\one_{\Ran^{\on{untl}}}$. 

\sssec{}

Unwinding the definitions, one obtains: 

\begin{lem} \label{l:mon str indep}
The pointwise monoidal structure on $\bA_{\Ran^{\on{untl}}}$ descends to the quotient 
$$\bA_{\Ran^{\on{untl}}}\twoheadrightarrow \bA_{\Ran^{\on{untl}},\on{indep}}.$$
\end{lem}

In what follows we will consider $\bA_{\Ran^{\on{untl}},\on{indep}}$ as a monoidal category, with the monoidal
structure furnished by \lemref{l:mon str indep}. 

\medskip

Since $\bA^{\sotimes}_{\Ran^{\on{untl}}}$ is unital, we obtain that so is $\bA_{\Ran^{\on{untl}},\on{indep}}$.

\sssec{}

Note now that the natural transformation
$$(\Delta_{\Ran^{\on{untl}}})^!\simeq \on{union}_!\circ (\Delta_{\Ran^{\on{untl}}})_!\circ (\Delta_{\Ran^{\on{untl}}})^!\to
\on{union}_!$$
defines on the identity functor on $\bA_{\Ran^{\on{untl}}}$ a structure of \emph{right-lax} monoidal functor
\begin{equation} \label{e:conv vs pointwise}
\bA^\star_{\Ran^{\on{untl}}}\to \bA^{\sotimes}_{\Ran^{\on{untl}}}.
\end{equation}

\sssec{}

However, we claim:

\begin{lem}  \label{l:conv vs pointwise}
The right-lax monoidal structure on the functor \eqref{e:conv vs pointwise} is strict.
\end{lem}

\begin{proof} 

Let $\ul\bC$ denote the crystal of categories 
$$\on{union}^*(\bA^{\on{loc,untl}}\otimes \bA^{\on{loc,untl}})$$
on $\Ran^{\on{untl}}\times \Ran^{\on{untl}}$. 

\medskip

Let $\on{ins.union}$ denote the functor \eqref{e:conv op unital -2}; we will use the same notation for 
the induced functor
$$\bA_{\Ran^{\on{untl}}}\otimes \bA_{\Ran^{\on{untl}}}\simeq 
\Gamma^{\on{lax}}(\Ran^{\on{untl}}\times \Ran^{\on{untl}},\bA^{\on{loc,untl}}\boxtimes \bA^{\on{loc,untl}})\to
\Gamma^{\on{lax}}(\Ran^{\on{untl}}\times \Ran^{\on{untl}},\ul\bC).$$

\medskip

We have to show that for $\CF_1,\CF_2,\CF\in \bA^\star_{\Ran^{\on{untl}}}$, the map
\begin{multline} \label{e:conv vs pointwise 1}
\CHom_{\Gamma^{\on{lax}}(\Ran^{\on{untl}}\times \Ran^{\on{untl}},\ul\bC)}
(\on{ins.union}(\CF_1\boxtimes \CF_2),\on{union}^!(\on{mult}^R(\CF))) \to \\
\to 
\CHom_{\Gamma^{\on{lax}}(\Ran^{\on{untl}},\bA^{\on{loc,untl}})}(\on{mult}(\CF_1\sotimes \CF_2),\CF),
\end{multline} 
is an isomorphism, where: 

\begin{itemize}

\item $\on{mult}$ denotes the functor 
$$\Gamma^{\on{lax}}(\Ran^{\on{untl}},\bA^{\on{loc,untl}}\otimes \bA^{\on{loc,untl}})\to \Gamma^{\on{lax}}(\Ran^{\on{untl}},\bA^{\on{loc,untl}})$$
induced by the monoidal operation
$$\bA^{\on{loc,untl}} \otimes \bA^{\on{loc,untl}}\overset{\on{mult}}\to \bA^{\on{loc,untl}} ;$$

\item $\on{mult}^R$ denotes the functor 
$$\Gamma^{\on{lax}}(\Ran^{\on{untl}},\bA^{\on{loc,untl}})\to 
\Gamma^{\on{lax}}(\Ran^{\on{untl}},\bA^{\on{loc,untl}}\otimes \bA^{\on{loc,untl}})$$
induced by the functor 
$$\bA^{\on{loc,untl}} \overset{\on{mult}^R}\to \bA^{\on{loc,untl}} \otimes \bA^{\on{loc,untl}},$$
right adjoint to the monoidal operation;

\smallskip

\item $\on{union}^!$ denotes the functor
$$\Gamma^{\on{lax}}(\Ran^{\on{untl}},\bA^{\on{loc,untl}}\otimes \bA^{\on{loc,untl}})\to
\Gamma^{\on{lax}}(\Ran^{\on{untl}}\times \Ran^{\on{untl}},\ul\bC);$$

\item The map \eqref{e:conv vs pointwise 1} is given by the composition
\begin{multline*}
\CHom_{\Gamma^{\on{lax}}(\Ran^{\on{untl}}\times \Ran^{\on{untl}},\bC)}
(\on{ins.union}(\CF_1\boxtimes \CF_2),\on{union}^!(\on{mult}^R(\CF))) \overset{\Delta^!_{\Ran^{\on{untl}}}}\longrightarrow \\
\CHom_{\Gamma^{\on{lax}}(\Ran^{\on{untl}},\Delta^*_{\Ran^{\on{untl}}}(\ul\bC))}
(\Delta^!_{\Ran^{\on{untl}}}(\on{ins.union}(\CF_1\boxtimes \CF_2)),\Delta^!_{\Ran^{\on{untl}}}(\on{union}^!(\on{mult}^R(\CF)))) \simeq \\
\simeq 
\CHom_{\Gamma^{\on{lax}}(\Ran^{\on{untl}},\bA^{\on{loc,untl}}\otimes \bA^{\on{loc,untl}})}
(\CF_1\sotimes \CF_2,\on{mult}^R(\CF)) \simeq \\
\simeq \CHom_{\Gamma^{\on{lax}}(\Ran^{\on{untl}},\bA^{\on{loc,untl}})}(\on{mult}(\CF_1\sotimes \CF_2),\CF). 
\end{multline*}

\end{itemize}

\medskip

Now the isomorphism assertion holds for any crystal of categories $\ul\bC$ on  $\Ran^{\on{untl}}\times \Ran^{\on{untl}}$, since the map
$$\Delta_{\Ran^{\on{untl}}}:\Ran^{\on{untl}}\to \Ran^{\on{untl}}\times \Ran^{\on{untl}}$$
is cofinal.

\end{proof} 

\sssec{}

Thanks to \lemref{l:conv vs pointwise}, we do not need to distinguish between the two monoidal structures 
on $\bA_{\Ran^{\on{untl}}}$. Thus we will use the symbol $\bA_{\Ran^{\on{untl}}}$ unambiguously to refer
to a tensor structure on $\bA_{\Ran^{\on{untl}}}$.

\medskip

Yet we will sometimes use the symbols $\bA^\star_{\Ran^{\on{untl}}}$ or $\bA^{\sotimes}_{\Ran^{\on{untl}}}$
to emphasize that we are thinking of the monoidal structure as convolution or pointwise tensor product,
respectively. 

\sssec{} \label{sss:integrated va Tw arr}

Let $\bA^{\on{loc,untl}}$ be as in \secref{sss:constant monoidal}, i.e., it is
associated to a crystal $\ul\bA_X$ of symmetric monoidal categories on $X$. 
In this case, one can describe 
the corresponding (symmetric) monoidal category $\bA_{\Ran^{\on{untl}},\on{indep}}$ explicitly.

\medskip

Namely, as a DG category it is isomorphic to the colimit over the twisted arrows category of
the category of finite non-empty sets (and arbitrary maps) of the functor that associates to 
$$I_1\overset{\phi}\to I_2$$
the category
\begin{equation} \label{e:Tw arr functor}
\Gamma(X^{I_2},\underset{i\in I_2}\boxtimes\, \ul\bA_X^{\otimes \phi^{-1}(i_2)}).
\end{equation} 

The symmetric monoidal structure is given by the operation of disjoint union of finite sets,
see \cite[Sect. 2.2.1]{FraG}. 

\ssec{Sheaves of monoidal categories on the \emph{non}-unital Ran space}

We now consider the usual (i.e., non-unital) Ran space. 
For $\bA^{\on{loc,untl}}$ as above, let $\bA^{\on{loc}}$ denote its restriction along the map $\sft:\Ran\to \Ran^{\on{untl}}$. 

\medskip

Denote 
$$\bA_\Ran:=\Gamma(\Ran,\bA^{\on{loc}}).$$

In this subsection we will endow $\bA_\Ran$ with monoidal structure(s) and study its
interactions with the unital counterparts. 

\sssec{}

By a slight abuse of notation, we will
use the same symbol $\on{union}$ to denote the corresponding map
$$\Ran\times \Ran\to \Ran.$$

\medskip

Restricting \eqref{e:conv op unital -1} along
$$\sft:\Ran\to \Ran_{\on{untl}}$$
we obtain a map of crystals of categories
\begin{equation} \label{e:conv op unital -1 nu}
\bA^{\on{loc}}\boxtimes \bA^{\on{loc}}\to \on{union}^*(\bA^{\on{loc}})
\end{equation} 
on $\Ran\times \Ran$. 

\medskip

Since the map $\on{union}$ is pseudo-proper, the functor \eqref{e:conv op unital -1 nu} induces a functor 
$$\bA_\Ran\otimes \bA_\Ran \overset{\star}\to \bA_\Ran.$$

\sssec{} \label{sss:conv non-unital}

The above binary operation extends to a monoidal structure on $\bA_\Ran$,
which we will refer to us the \emph{convolution} monoidal structure. We will denote $\bA_\Ran$,
viewed as a monoidal category equipped with the convolution monoidal structure by $\bA^\star_\Ran$.

\sssec{} \label{sss:indep monoidal loc non-un}

By construction, the functor
$$\sft_!:\bA_\Ran\to \bA_{\Ran^{\on{untl}}}$$
has a monoidal structure, when we consider both as equipped with the convolution monoidal structure. 

\medskip

In particular, we obtain that the functor
\begin{equation} \label{e:from nu Ran}
\on{emb.indep}^L\circ \sft_!:\bA^\star_\Ran\to \bA_{\Ran^{\on{untl}},\on{indep}}
\end{equation}
acquires a monoidal structure. 

\medskip

Combining with \corref{c:Ran vs indep 2}, we obtain that the functor 
\eqref{e:from nu Ran} is a monoidal localization.

\sssec{}

Consider now the functor
$$\sft^!:\bA_{\Ran^{\on{untl}}}\to \bA_\Ran.$$

Being the right adjoint of a monoidal functor, the functor $\sft^!$ acquires a right-lax monoidal structure as a functor
$$\bA^\star_{\Ran^{\on{untl}}}\to \bA^\star_\Ran.$$

\medskip

Note that the natural transformation
\begin{equation} \label{e:t back and forth}
\on{emb.indep}^L\circ \sft_!\circ \sft^!\to \on{emb.indep}^L
\end{equation}
has a natural right-lax monoidal structure. However, from \corref{c:Ran vs indep 1} we obtain 
that the natural transformation (as right-lax monoidal functors) is an isomorphism. 

\medskip

In particular, 
we obtain that the right-lax monoidal structure on \eqref{e:t back and forth} is strict. 

\sssec{} \label{sss:pointwise monoidal on Ran}

As in \secref{sss:Ran unital pointwise}, we can also consider the $\sotimes$-monoidal structure
on $\bA_\Ran$. We will denote $\bA_\Ran$,
viewed as a monoidal category equipped with the pointwise structure by $\bA^{\sotimes}_\Ran$.

\medskip

As in \secref{sss:Ran unital pointwise}, the identity functor on $\bA_\Ran$ has a right-lax monoidal
structure, when viewed as a functor
$$\bA^\star_\Ran\to \bA^{\sotimes}_\Ran.$$

\sssec{}

The functor $\sft^!$ is monoidal, when viewed as a functor
$$\bA^{\sotimes}_{\Ran^{\on{untl}}}\to \bA^{\sotimes}_\Ran.$$

Hence, the functor 
$$\sft_!:\bA^{\sotimes}_\Ran\to \bA^{\sotimes}_{\Ran^{\on{untl}}}$$
acquires a left-lax monoidal structure, when viewed as a functor
$$\bA^{\sotimes}_\Ran\to \bA^{\sotimes}_{\Ran^{\on{untl}}}.$$

\sssec{}

Let 
$$\bA^{\on{almost-untl}}_\Ran\subset \bA_\Ran$$
be the full subcategory generated by the essential image of the functor $\sft^!$.

\medskip

It is easy to see that it is preserved by the $\sotimes$ monoidal operation, and
hence it acquires a monoidal structure.

\sssec{}

We claim:

\begin{lem} \label{l:! vs star indep 1}
The left-lax monoidal structure on the functor
$$\bA^{\sotimes}_\Ran \overset{\sft_!}\to \bA^{\sotimes}_{\Ran^{\on{untl}}} 
\overset{\on{emb.indep}^L}\to \bA_{\Ran^{\on{untl}},\on{indep}}$$
becomes strict when restricted to $\bA^{\on{almost-untl}}_\Ran$. 
\end{lem} 

\begin{proof}

It suffices to show that the left-lax monoidal structure on the functor 
$$\bA^{\sotimes}_{\Ran^{\on{untl}}} \overset{\sft^!}\to \bA^{\sotimes}_\Ran 
\overset{\sft_!}\to \bA^{\sotimes}_{\Ran^{\on{untl}}} 
\overset{\on{emb.indep}^L}\to \bA_{\Ran^{\on{untl}},\on{indep}}$$
is strict.

\medskip

The assertion follows now from \corref{c:Ran vs indep 2}, which implies that the above composition is
isomorphic to 
$$\bA^{\sotimes}_{\Ran^{\on{untl}}} 
\overset{\on{emb.indep}^L}\to \bA_{\Ran^{\on{untl}},\on{indep}}$$
as a left-lax monoidal functor. 

\end{proof}

\sssec{}

Let $\bA$ be as in \secref{sss:integrated va Tw arr}. In this case, we can also describe the category $\bA^\star_\Ran$
explicitly.

\medskip

It is given by the colimit of the functor \eqref{e:Tw arr functor}, with the only difference that we take the twisted
arrows category of the category of finite non-empty sets and \emph{surjective} maps.

\medskip

The functor
$$\sft_!:\bA^\star_\Ran\to \bA^\star_{\Ran^{\on{untl}}}$$
is given by embedding the index categories one into the other. 

\ssec{Local and integrated monoidal actions} \label{ss:actions on Ran}

\sssec{} \label{sss:local actions}

Let $\bA^{\on{loc,untl}}$ be as above. Let $\bA^{\on{loc}}$ denote the restriction of $\bA^{\on{loc,untl}}$
along $\sft:\Ran\to \Ran^{\on{untl}}$. 

\medskip

Let $\bD$ be a DG category. We give the following definitions:

\medskip

\begin{itemize}

\item A local action of $\bA^{\on{loc}}$ on $\bD$ is a (unital) 
action of the crystal of monoidal categories $\bA^{\on{loc}}$ on $\bD\otimes \ul\Dmod(\Ran)$;

\medskip

\item A local lax-Ran-unital action of $\bA^{\on{loc,untl}}$ on $\bD$ is a (unital) action of the crystal of monoidal categories $\bA^{\on{loc,untl}}$ on $\bD\otimes \ul\Dmod(\Ran^{\on{untl}})$, in the 2-category of 
crystals of categories and right-lax functors between them; 

\medskip

\item A local Ran-unital action of $\bA^{\on{loc,untl}}$ on $\bD$ is a (unital) action of the crystal of monoidal categories $\bA^{\on{loc,untl}}$ on $\bD\otimes \ul\Dmod(\Ran^{\on{untl}})$, in the 2-category of 
crystals of categories and strict functors between them.

\end{itemize}

\sssec{}

At the poinwtise level, a local action of $\bA^{\on{loc}}$ on $\bD$ yields an action of the 
monoidal category $\bA_{\ul{x}}$ on $\bD$ for every $\ul{x}\in \Ran$, denoted
$$\ba_{\ul{x}},\bd\mapsto \ba\cdot \bd, \quad \ba_{\ul{x}}\in \bA_{\ul{x}},\, \bd\in \bD.$$

\medskip

A lax-Ran-unital structure on such an action is a natural transformation
$$\ba_{\ul{x}_1}\cdot \bd\to 
\on{ins.unit}_{\ul{x}_1\subseteq \ul{x}_2}(\ba_{\ul{x}_1})\cdot \bd, \quad 
\ba_{\ul{x}_1}\in \bA_{\ul{x}_1},\, \bd\in \bD.$$

A lax-Ran-unital structure is strict if the above natural transformation is an isomorphism. 

\sssec{}

Thus, we obtain the 2-categories
$$\bA^{\on{loc}}\mmod,\,\, \bA^{\on{loc},\on{untl}}\mmod^{\on{lax}} \text{ and }
\bA^{\on{loc},\on{untl}}\mmod.$$

We have a fully faithful functor
$$\bA^{\on{loc},\on{untl}}\mmod\hookrightarrow \bA^{\on{loc},\on{untl}}\mmod^{\on{lax}}$$
and a forgetful functor
$$ \bA^{\on{loc},\on{untl}}\mmod^{\on{lax}}\to \bA^{\on{loc}}\mmod.$$

\medskip

Note, however, that from \corref{c:Ran vs indep}, we obtain:

\begin{cor} \label{c:unital vs non-unital Ran actions}
The composite functor
\begin{equation} \label{e:unital vs non-unital Ran actions}
\bA^{\on{loc},\on{untl}}\mmod\hookrightarrow \bA^{\on{loc},\on{untl}}\mmod^{\on{lax}}\to 
\bA^{\on{loc}}\mmod
\end{equation}
is fully faithful.
\end{cor}

\sssec{} \label{sss:integrated actions}

Consider now the (unital) monoidal category $\bA^{\sotimes}_{\Ran^{\on{untl}}}$.
We claim that given a local lax-Ran-unital action of $\bA^{\on{loc,untl}}$ on $\bD$, we can construct
a (unital) action of $\bA^{\sotimes}_{\Ran^{\on{untl}}}$ on $\bD$.

\medskip

Indeed, this follows from the fact that the functor 
$$\on{C}^\cdot_c(\Ran^{\on{untl}},-):\Dmod(\Ran^{\on{untl}})\to \Vect$$
is symmetric monoidal. 

\medskip

Explicitly, given $\ba\in \bA_{\Ran^{\on{untl}}}$, the action is given by
\begin{multline*}
\bD \overset{-\otimes \omega_{\Ran^{\on{untl}}}}\longrightarrow \bD\otimes \Dmod(\Ran^{\on{untl}})^{\on{strict}}\hookrightarrow \\
\to \bD\otimes \Dmod(\Ran^{\on{untl}})\overset{\ba\cdot -}\to \bD\otimes \Dmod(\Ran^{\on{untl}})
\overset{\on{Id}\otimes \on{C}^\cdot_c(\Ran^{\on{untl}},-)}\longrightarrow \bD.
\end{multline*}

\sssec{}

From \lemref{l:dual of lax sections} we obtain:

\begin{cor} \label{c:integrated actions}
The composite functor
$$\bA^{\on{loc},\on{untl}}\mmod \hookrightarrow \bA^{\on{loc},\on{untl}}\mmod^{\on{lax}} \to 
\bA^{\sotimes}_{\Ran^{\on{untl}}}\mmod$$
is fully faithful with essential image being
$$\bA_{\Ran^{\on{untl}},\on{indep}}\mmod\subset 
\bA^{\sotimes}_{\Ran^{\on{untl}}}\mmod.$$
\end{cor}

\sssec{} \label{sss:action of conv}

Precomposing the construction in \secref{sss:integrated actions} with the monoidal functor
$$\bA^\star_\Ran\overset{\sft_!}\to \bA^\star_{\Ran^{\on{untl}}}\simeq \bA^{\sotimes}_{\Ran^{\on{untl}}}$$
we obtain that for
$$\bD\in \bA^{\on{loc},\on{untl}}\mmod^{\on{lax}}$$
we have an action of $\bA^\star_\Ran$ on $\bD$.

\begin{rem}

Unwinding the definitions, for $\ba\in \bA_\Ran$, its action on $\bD$ is given by the composition
$$\bD \overset{-\otimes \omega_{\Ran}}\longrightarrow \bD\otimes \Dmod(\Ran)
\overset{\ba\cdot -}\to \bD\otimes \Dmod(\Ran)
\overset{\on{Id}\otimes \on{C}^\cdot_c(\Ran,-)}\longrightarrow \bD.$$

I.e., the binary operation only depends on structure on $\bD$ of object of 
$\bA^{\on{loc}}\mmod$. 

\medskip

However, one would not be able to define either left-lax
or right-lax structure on this binary operation without $\bD$ being extended 
to an object of $\bA^{\on{loc},\on{untl}}\mmod^{\on{lax}}$. 

\end{rem}

\sssec{}

From \secref{sss:indep monoidal loc non-un} and \corref{c:integrated actions} we obtain:

\begin{cor} \label{c:indep action from conv}
The composite functor
$$\bA^{\on{loc},\on{untl}}\mmod \hookrightarrow \bA^{\on{loc},\on{untl}}\mmod^{\on{lax}} \to 
\bA^{\sotimes}_{\Ran^{\on{untl}}}\mmod\simeq \bA^\star_{\Ran^{\on{untl}}}\mmod\to \bA^\star_\Ran\mmod$$
is fully faithful with the essential image being
$$\bA_{\Ran^{\on{untl}},\on{indep}}\mmod\subset \bA^\star_\Ran\mmod.$$
\end{cor}

\ssec{Local actions with parameters} \label{ss:loc act param}

\sssec{}

Let now $S$ be an affine scheme mapping to $\Ran$. The discussion in 
Sects. \ref{ss:monoidal on untl Ran}-\ref{ss:actions on Ran} applies when we replace
$\Ran$ (resp., $\Ran^{\on{untl}}$) by $S^{\subseteq}$ (resp., $S^{\subseteq,\on{untl}}$).

\medskip

Note that $S^{\subseteq}$ and $S^{\subseteq,\on{untl}}$ are pseudo-proper \emph{relative}
to $S$, so \secref{sss:pseudo-proper rel} applies. 

\sssec{}

In particular, one can consider the $\Dmod(S)$-linear monoidal categories
$$\bA^\star_{S^{\subseteq,\on{untl}}}\,\, \bA^{\sotimes}_{S^{\subseteq,\on{untl}}},\,\,
\bA^\star_{S^{\subseteq}},\,\, \bA^{\sotimes}_{S^{\subseteq}} \text{ and } \bA_{S^{\subseteq,\on{untl}},\on{indep}},$$
equipped with the monoidal functors
$$\bA^\star_{S^{\subseteq,\on{untl}}}\overset{\sim}\to \bA^{\sotimes}_{S^{\subseteq,\on{untl}}},$$
$$\bA^\star_{S^{\subseteq}} \overset{\sft_!}\to \bA^\star_{S^{\subseteq,\on{untl}}},$$
$$\bA^{\sotimes}_{S^{\subseteq,\on{untl}}} \overset{\sft^!}\to \bA^{\sotimes}_{S^{\subseteq,\on{untl}}},$$
$$\bA^\star_{S^{\subseteq,\on{untl}}} \overset{\on{embed.indep}^L}\twoheadrightarrow \bA_{S^{\subseteq,\on{untl}},\on{indep}},$$
$$\bA^\star_{S^{\subseteq}} 
\overset{\on{embed.indep}^L\circ \sft_!}\twoheadrightarrow \bA_{S^{\subseteq,\on{untl}},\on{indep}},$$
where the last two functors are monoidal localizations. 

\sssec{} \label{sss:local actions params}

We also have the corresponding notions of local action, so we have the 2-categories
$$\bA_S^{\on{loc}}\mmod,\,\, \bA_S^{\on{loc},\on{untl}}\mmod^{\on{lax}} \text{ and }
\bA_S^{\on{loc},\on{untl}}\mmod,$$
the fully faithful embedding  
$$\bA_S^{\on{loc},\on{untl}}\mmod\hookrightarrow \bA_S^{\on{loc},\on{untl}}\mmod^{\on{lax}}$$
and a forgetful functor
$$\bA_S^{\on{loc},\on{untl}}\mmod^{\on{lax}}\to \bA_S^{\on{loc}}\mmod,$$
so that composition
$$\bA_S^{\on{loc},\on{untl}}\mmod\hookrightarrow \bA_S^{\on{loc},\on{untl}}\mmod^{\on{lax}}\to \bA_S^{\on{loc}}\mmod$$
is fully faithful.

\sssec{}

In addition, we have the commutative diagrams
$$
\CD
\bA_S^{\on{loc},\on{untl}}\mmod @>>> \bA_S^{\on{loc},\on{untl}}\mmod^{\on{lax}} \\
@V{\sim}VV @VVV \\
\bA_{S^{\subseteq,\on{untl}},\on{indep}}\mmod @>>> \bA^\star_{S^{\subseteq,\on{untl}}}\mmod 
\endCD
$$

\sssec{}

We now claim:

\begin{prop} \label{p:get rid of parameters}
The functor
$$\Dmod(S)\otimes \bA^{\sotimes}_{\Ran^{\on{untl}},\on{indep}}\overset{(\on{pr}^{\on{untl}}_{\on{small},S})^!\otimes (\on{pr}^{\on{untl}}_{\on{big}})^!}
\longrightarrow \bA^{\sotimes}_{S^{\subseteq,\on{untl}},\on{indep}}$$
 is an equivalence.
\end{prop}

\begin{proof}

Note that operation of union defines a map 
$$S\times \Ran^{\on{untl}} \overset{\on{union}_S}\to S^{\subseteq,\on{untl}}.$$

Hence, we obtain a functor
$$\bA^{\sotimes}_{S^{\subseteq,\on{untl}}}=\Gamma^{\on{lax}}(S^{\subseteq,\on{untl}},\bA^{\on{loc,untl}})
\overset{\on{union}_S^!}\longrightarrow
\Gamma^{\on{lax}}(S\times \Ran^{\on{untl}},\bA^{\on{loc,untl}})\simeq \Dmod(S)\otimes \bA^{\sotimes}_{\Ran^{\on{untl}}},$$
and it is easy to see that it sends 
$$\bA^{\sotimes}_{S^{\subseteq,\on{untl}},\on{indep}}\to 
\Dmod(S)\otimes \bA^{\sotimes}_{\Ran^{\on{untl}},\on{indep}}.$$

The composition
$$S^{\subseteq,\on{untl}} \overset{\on{pr}^{\on{untl}}_{\on{small},S}\times \on{pr}^{\on{untl}}_{\on{big}}}\longrightarrow 
S\times \Ran^{\on{untl}} \overset{\on{union}_S}\to S^{\subseteq,\on{untl}}$$
is the identity map. Hence, the composition
$$\bA^{\sotimes}_{S^{\subseteq,\on{untl}},\on{indep}}\overset{\on{union}_S^!}\longrightarrow 
\Dmod(S)\otimes \bA^{\sotimes}_{\Ran^{\on{untl}},\on{indep}} \overset{(\on{pr}^{\on{untl}}_{\on{small},S})^!\otimes (\on{pr}^{\on{untl}}_{\on{big}})^!}
\longrightarrow \bA^{\sotimes}_{S^{\subseteq,\on{untl}},\on{indep}}$$
is the identity functor.

\medskip

The composition 
$$S\times \Ran^{\on{untl}} \overset{\on{union}_S}\to S^{\subseteq,\on{untl}} 
\overset{\on{pr}^{\on{untl}}_{\on{small},S}\times \on{pr}^{\on{untl}}_{\on{big}}}\longrightarrow 
S\times \Ran^{\on{untl}}$$
receives a 1-morphism from the identity map. This 1-morphism defines a natural transformation from the identity endofunctor on
$\Dmod(S)\otimes \bA^{\sotimes}_{\Ran^{\on{untl}}}$ to the composition
$$\Dmod(S)\otimes \bA^{\sotimes}_{\Ran^{\on{untl}}}\overset{(\on{pr}^{\on{untl}}_{\on{small},S}\times \on{pr}^{\on{untl}}_{\on{big}})^!}
\longrightarrow \bA^{\sotimes}_{\Ran^{\on{untl}}} \overset{\on{union}_S^!}\longrightarrow 
\Dmod(S)\otimes \bA^{\sotimes}_{\Ran^{\on{untl}}}.$$

However, this natural transformation is an isomorphism when restricted to 
$\Dmod(S)\otimes \bA^{\sotimes}_{\Ran^{\on{untl}}}$, by the definition of this subcategory.

\end{proof}

\sssec{}

From \propref{p:get rid of parameters} we obtain: 

\begin{cor} \label{c:get rid of parameters}
For a $\Dmod(S)$-module category $\bD$, pullback along $\on{pr}_{\on{big}}:S^{\subseteq,\on{untl}}\to \Ran^{\on{untl}}$
defines an equivalence between the following data:

\medskip

\noindent{\em(i)} A local Ran-unital action of $\bA^{\on{loc},\on{untl}}_S$ on 
$\ul\Dmod(S^{\subseteq,\on{untl}})\underset{\Dmod(S)}\otimes \bD$;

\medskip

\noindent{\em(ii)} A local Ran-unital action of $\bA^{\on{loc},\on{untl}}$ on 
$\bD$, compatible with the $\Dmod(S)$-action. 

\medskip

\noindent{\em(ii')} An action on $\bD$ of the monoidal category $\bA^{\sotimes}_{\Ran^{\on{untl}},\on{indep}}$,
compatible with the $\Dmod(S)$-action. 

\end{cor}

%
%
%
%
%
%
%
%
%
%
%
%
%

\section{The integration functor in the non-unital setting} \label{s:add unit colax}

In this appendix, we give a categorical meaning to the functor \eqref{e:ins unit int as functor}
\[
	\int\on{ins.unit}: \on{Funct}^{\on{loc}\to \on{glob}}(\ul\bC^{\on{loc}},\bC^{\on{glob}})\to \on{Funct}^{\on{loc}\to \on{glob}}(\ul\bC^{\on{loc}},\bC^{\on{glob}}).
\]
using the notion of \emph{left-lax} (a.k.a. left-lax) unital structure on a local-to-global functor $\ul\sF$. As a byproduct, we provide a proof to \propref{p:unitality as a property}.

\ssec{A left-lax unital structure on a local-to-global functor}
\label{ss: left-lax intro}

\sssec{}
	In this subsection, we describe the notion of a left-lax unital structure on a local-to-global functor $\ul\sF$ in concrete words. The precise definition will be given in Sect. \ref{ss:left-lax}. Also, we explain why left-lax unital structures provide categorical meaning to the functor $\int\on{ins.unit}$.

\sssec{}
	Let $\CZ$ be a space, and let $\ul{x}_1\overset{\alpha}\to \ul{x}_2$ be a morphism in the category $\Maps(\CZ,\Ran^{\on{untl}})$. Recall (see Sect. \ref{sss:unitality concrete}) that a \emph{lax} unital structure on $\ul\sF$ means there is a natural transformation
	\[
		\sF_{\CZ,\ul{x}_1}\to \sF_{\CZ,\ul{x}_2}\circ \bC^{\on{loc}}_\alpha
	\]
	as functors $\bC^{\on{loc}}_{\CZ,\ul{x}_1}\to \bC^{\on{glob}}\otimes \Dmod(\CZ)$. Then a \emph{left-lax} unital structure on $\ul\sF$ means there is a natural transformation of the opposite direction, i.e.,
	\[
		\sF_{\CZ,\ul{x}_2}\circ \bC^{\on{loc}}_\alpha \to \sF_{\CZ,\ul{x}_1}.
	\]

\sssec{}
	As in Sect. \ref{sss:F subset Z}, using the prestack $\CZ^\subseteq$, we can rewrite the datum of these natural transformations as a natural transformation 
	\begin{equation}
		\label{e:left-lax unitality transform original}
		\sF_{\CZ^\subseteq} \circ \on{ins.unit}_\CZ \to
		(\on{Id}\otimes (\on{pr}_{\on{small},\CZ})^!) \circ \sF_\CZ
	\end{equation}
	as functors
	$$\bC^{\on{loc}}_\CZ\to \bC^{\on{glob}}\otimes \Dmod(\CZ^{\subseteq})$$
	which is supposed to be equipped with a datum of associativity.

\sssec{}
	Using the adjunction $((\on{pr}_{\on{small},\CZ})_!, (\on{pr}_{\on{small},\CZ})^!)$, knowing \eqref{e:left-lax unitality transform original} is equivalent to knowing a natural transformation
	\[
		(\on{Id}\otimes (\on{pr}_{\on{small},\CZ})_!) \circ \sF_{\CZ^\subseteq} \circ \on{ins.unit}_\CZ  \to \sF_\CZ.
	\]
	Note that the LHS is exactly the functor $\sF^{\int \on{ins.unit}}_\CZ$ (see \eqref{e:ins unit int as functor}). Hence a left-lax unital structure on $\ul\sF$ means a natural transformation
	\[
		\ul\sF^{\int \on{ins.unit}} \to \ul\sF
	\]
	equipped with a datum of associativity. As we will see in the proof of \propref{prop: all to left-lax}, the endofunctor
	\[
		\int\on{ins.unit}: \on{Funct}^{\on{loc}\to \on{glob}}(\ul\bC^{\on{loc}},\bC^{\on{glob}})\to \on{Funct}^{\on{loc}\to \on{glob}}(\ul\bC^{\on{loc}},\bC^{\on{glob}})
	\]
	has a natural monad structure, and this datum of associativity says exactly that $\ul\sF$ is a module for this monad. 

\sssec{}
	In particular, for \emph{any} local-to-global functor $\ul\sF$, the functor $\ul\sF^{\int \on{ins.unit}}$ is an induced (a.k.a. free) module for this monad. Therefore we obtain a left-lax unital structure on $\ul\sF^{\int \on{ins.unit}}$, and this left-lax unital structure satisfies the following universal property. For any local-to-global functor $\ul\sG$ equipped with a left-lax unital structure, knowing a (plain) natural transformation $\ul\sF \to \ul\sG$ is equivalent to knowing a natural transformation $\ul\sF^{\int \on{ins.unit}}\to \ul\sG$ compatible with the left-lax unital structure.

\sssec{}
	Let 
	\[
		\on{Funct}^{\on{loc}\to \on{glob},\on{left-lax-untl}}(\ul\bC^{\on{loc}},\bC^{\on{glob}})
	\]
	be the category of left-lax unital local-to-global functors. The above observation suggests the following result, which will be proved in \ref{ss: proof prop: all to left-lax} (after we give a precise definition to left-lax unital functors).

	\begin{prop}
	\label{prop: all to left-lax}
	The forgetful functor
	\[
		\iota_{\on{left-lax \to all}}:\on{Funct}^{\on{loc}\to \on{glob},\on{left-lax-untl}}(\ul\bC^{\on{loc}},\bC^{\on{glob}})
		\to
			\on{Funct}^{\on{loc}\to \on{glob}}(\ul\bC^{\on{loc}},\bC^{\on{glob}})
	\]
	has a left adjoint 
	\[
		\iota_{\on{left-lax \to all}}^L: \on{Funct}^{\on{loc}\to \on{glob}}(\ul\bC^{\on{loc}},\bC^{\on{glob}}) \to \on{Funct}^{\on{loc}\to \on{glob},\on{left-lax-untl}}(\ul\bC^{\on{loc}},\bC^{\on{glob}})
	\]
	that sends 
	\[
		\ul\sF \mapsto \ul\sF^{\int \on{ins.unit}}.
	\]
	In particular, the latter has a natural left-lax unital structure.
\end{prop}

\ssec{Comparison with the integration functor in the lax unital setting}

\sssec{}
	Tautologically, there is a commutative square
	\begin{equation}
		\label{e: lax colax sq general local global}
		\xymatrix{
			\on{Funct}^{\on{loc}\to \on{glob},\on{untl}}(\ul\bC^{\on{loc}},\bC^{\on{glob}})
			\ar[r]_{\iota_{\on{st\to lax}}}^-\subset \ar[d]_-{\iota_{\on{st\to left-lax}}}^-\subset 
			& \on{Funct}^{\on{loc}\to \on{glob},\on{lax-untl}}(\ul\bC^{\on{loc}},\bC^{\on{glob}}) \ar[d]^-{\iota_{\on{lax \to all}}} \\
			\on{Funct}^{\on{loc}\to \on{glob},\on{left-lax-untl}}(\ul\bC^{\on{loc}},\bC^{\on{glob}}) \ar[r]_-{\iota_{\on{left-lax \to all}}} &
			\on{Funct}^{\on{loc}\to \on{glob}}(\ul\bC^{\on{loc}},\bC^{\on{glob}})
		}
	\end{equation}
	where each $\iota_?$ is a forgetful functor from a category of local-to-global functors equipped with certain unital structures to a category of such functors equipped with coarser structures. Here only $\iota_{\on{st\to lax}}$ and $\iota_{\on{st\to left-lax}}$ are fully faithful.

\sssec{}
	Recall in Sect. \ref{ss:constr int}, we also constructed an adjoint pair:
	\[
		\iota_{\on{st \to lax}}^L: \on{Funct}^{\on{loc}\to \on{glob},\on{lax-untl}}(\ul\bC^{\on{loc}},\bC^{\on{glob}}) \rightleftarrows \on{Funct}^{\on{loc}\to \on{glob},\on{untl}}(\ul\bC^{\on{loc}},\bC^{\on{glob}}) :\iota_{\on{st \to lax}}
	\]
	such that the left adjoint sends 
	\[
		\ul\sF^{\on{untl}} \mapsto \ul\sF^{\on{untl},\int \on{ins.unit}}.
	\]
	The following result, which will be proved in \ref{ss: proof prop: all to left-lax adjointable}, says this adjoint pair is compatible with that in \propref{prop: all to left-lax}.

	\begin{prop}
		\label{prop: all to left-lax adjointable}
		The commutative square \eqref{e: lax colax sq general local global} is left adjointable along the horizontal direction, i.e., the Bech--Chevalley natural transformation from the clockwise circuit in the diagram below to the counter-clockwise circuit is invertible:
	\begin{equation}
		\label{e: lax colax sq general local global left adj}
		\xymatrix{
			\on{Funct}^{\on{loc}\to \on{glob},\on{untl}}(\ul\bC^{\on{loc}},\bC^{\on{glob}})
			 \ar[d]_-{\iota_{\on{st\to left-lax}}}^-\subset 
			& \on{Funct}^{\on{loc}\to \on{glob},\on{lax-untl}}(\ul\bC^{\on{loc}},\bC^{\on{glob}}) \ar[d]^-{\iota_{\on{lax \to all}}} 
				\ar[l]_{\iota^L_{\on{st\to lax}}}
			\\
			\on{Funct}^{\on{loc}\to \on{glob},\on{left-lax-untl}}(\ul\bC^{\on{loc}},\bC^{\on{glob}})  &
			\on{Funct}^{\on{loc}\to \on{glob}}(\ul\bC^{\on{loc}},\bC^{\on{glob}})
			\ar[l]_-{\iota^L_{\on{left-lax \to all}}}.
		}
	\end{equation}
	Moreover, 
	\begin{itemize}
		\item 
			The monad $\iota_{\on{st\to lax}}\circ \iota^L_{\on{st\to lax}}$ can be identified with \eqref{e:ins unit int as functor untl}.
		\item
			The underlying endofunctor of the monad $\iota_{\on{left-lax\to all}}\circ \iota^L_{\on{left-lax\to all}}$ can be identified with \eqref{e:ins unit int as functor}, and the unit of the monad can be identified with \eqref{e:correlator F}.
		\item
			The combination of \eqref{e: lax colax sq general local global left adj} and \eqref{e: lax colax sq general local global} gives the commutative diagram \eqref{e:ins unit int as functor untl compat}.			
	\end{itemize}
	\end{prop}

\ssec{Proof of \propref{p:unitality as a property}}
	
	In this subsection, we deduce \propref{p:unitality as a property} from \propref{prop: all to left-lax}. This is essentially a formal diagram chase. 

\sssec{}
	We need to show the functor
	\[
		\iota_{\on{st\to all}}:= \iota_{\on{lax\to all}} \circ \iota_{\on{st\to lax}} \simeq  \iota_{\on{left-lax\to all}} \circ \iota_{\on{st\to left-lax}}
	\]
	is fully faithful and identify its essential image. 

	\medskip
	Given 
	\[
		\ul\sF^{\on{untl}},\ul\sG^{\on{untl}}\in \on{Funct}^{\on{loc}\to \on{glob},\on{untl}}(\ul\bC^{\on{loc}},\bC^{\on{glob}}),
	\]
	we have
	\begin{eqnarray*}
		& & \Maps(\iota_{\on{st\to all}}( \ul\sF^{\on{untl}} ), \iota_{\on{st\to all}}( \ul\sG^{\on{untl}} )) \simeq\\
		&\simeq& \Maps(\iota_{\on{st\to all}}( \ul\sF^{\on{untl}} ), \iota_{\on{left-lax\to all}}\circ \iota_{\on{st\to left-lax}}( \ul\sG^{\on{untl}} )) \simeq\\
		&\simeq& \Maps( \iota^L_{\on{left-lax\to all}}\circ \iota_{\on{st\to all}}( \ul\sF^{\on{untl}} ), \iota_{\on{st\to left-lax}}( \ul\sG^{\on{untl}} ) ) \simeq\\
	        &\simeq& \Maps(\iota^L_{\on{left-lax\to all}}\circ \iota_{\on{lax\to all}}\circ \iota_{\on{st\to lax}}( \ul\sF^{\on{untl}}  ),\iota_{\on{st\to left-lax}}  
	        (\ul\sG^{\on{untl}} ) ) \simeq \\
	        &\simeq& \Maps(\iota_{\on{st\to left-lax}}\circ \iota_{\on{st\to lax}}^L\circ \iota_{\on{st\to lax}}( \ul\sF^{\on{untl}}  ),\iota_{\on{st\to left-lax}} 
	        (\ul\sG^{\on{untl}} ) ) \simeq \\
	        &\simeq& \Maps(\iota_{\on{st\to lax}}^L\circ \iota_{\on{st\to lax}}( \ul\sF^{\on{untl}}  ),
	        (\ul\sG^{\on{untl}} ) ) \simeq \\
	        &\simeq& \Maps(\iota_{\on{st\to left-lax}}( \ul\sF^{\on{untl}}  ),\iota_{\on{st\to left-lax}}( \ul\sG^{\on{untl}} ) ) \simeq \\
		&\simeq&  \Maps( \ul\sF^{\on{untl}}, \ul\sG^{\on{untl}}  ),
	\end{eqnarray*}
	where 
	\begin{itemize}
		\item The fourth equivalence is due to \propref{prop: all to left-lax};
		\item The sixth equivalence is because $\iota_{\on{st\to lax}}$ is fully faithful;
		\item The seventh equivalence is because $\iota_{\on{st\to left-lax}}$ is fully faithful.
	\end{itemize}
	We leave it to the reader to check the resulting composition is the inverse to the obvious map from the RHS to the LHS. 
	This proves $\iota_{\on{st\to all}}$ is fully faithful.

\sssec{}
	\label{sss: proof p:unitality as a property colax}
	Let $\ul\sF$ be a local-to-global functor. By definition, the Global Unitality Axiom (i) means exactly the unit adjunction
	\[
		\ul\sF \to \iota_{\on{left-lax\to all}} \circ \iota^L_{\on{left-lax\to all}}(\ul\sF)
	\]
	is invertible. In particular, $\ul\sF$ can be naturally lifted to the object
	\[
		\iota^L_{\on{left-lax\to all}}(\ul\sF)\in \on{Funct}^{\on{loc}\to \on{glob},\on{left-lax-untl}}(\ul\bC^{\on{loc}},\bC^{\on{glob}}).
	\]
	By the proof of \propref{prop: all to left-lax}, the Global Unitality Axiom (ii) means exactly the left-lax structure on $\iota^L_{\on{left-lax\to all}}(\ul\sF)$ 
	is strict (see Remark \ref{rem: left-lax GUAii} below). It follows that $\ul\sF$ is contained in the essential image of $\iota_{\on{st\to all}}$ if and only if 
	it satisfies the Global Unitality Axioms.

\qed[\propref{p:unitality as a property}]

\ssec{Definition of left-lax unital structures}
	\label{ss:left-lax}
	To give a homotopy-coherent definition of a left-lax unital structure on a local-to-global functor, we need some higher algebra. This will also complete the omitted higher datum in Sect. \ref{sss:unital insertion}, e.t.c..

\sssec{}
	For any $[n]\in \Delta^{\on{op}}$, let $\Ran^{\subseteq,[n]}$ be the moduli space of chains $\ul{x}_0\subseteq \ul{x}_1 \subseteq \dots \subseteq \ul{x}_n$. In particular, we have
	\[
		\Ran^{\subseteq,[0]} = \Ran,\; \Ran^{\subseteq,[1]} = \Ran^{\subseteq}.
	\]
	We obtain a simplicial prestack $\Ran^{\subseteq,\bullet}$ which is a \emph{($\infty$-)categocial object} in $\on{PreStk}$. Indeed, $\Ran^{\subseteq,[0]} = \Ran$ is the ``prestack of objects'' and $\Ran^{\subseteq,[1]} = \Ran^{\subseteq}$ is the ``prestack of $1$-morphisms''. The projections	
	\[
		\on{pr}_{\on{small}},\,\on{pr}_{\on{big}}: \Ran^{\subseteq} \to \Ran
	\]
	remember respectively the source and the target of a $1$-morphism, while
	\[
		\on{diag}: \Ran\to \Ran^{\subseteq}
	\]
	sends an object to the identity $1$-morphism at it. The higher categorical structure on $\Ran^{\subseteq,\bullet}$ is provided by its simplicial structure.

\sssec{}
	Let $\CY^\bullet$ be any categorical object in $\on{PreStk}$. Write $\mathcal{O}bj:=\CY^0$, $\mathcal{M}or_1:=\CY^1$ and $\mathcal{M}or_2:=\CY^2$. In general, there is a \emph{comonad} on the 2-category
	\[
		\on{CrystCat}(\mathcal{O}bj)
	\]
	of sheaves of categories on $\mathcal{O}bj$, whose underlying endofunctor is the composition
	\[
		\on{CrystCat}(\mathcal{O}bj) \xrightarrow{p_t^*} \on{CrystCat}(\mathcal{M}or_1) \xrightarrow{p_{s,*}} \on{CrystCat}(\mathcal{O}bj),
	\]
	where $p_s,p_t: \mathcal{M}or_1  \to \mathcal{O}bj $ are the projections that remember respectively the sources and the targets.

\sssec{}
	The operation of composition on $p_{s,*}\circ p_t^*$ can be obtained as follows. View
	\begin{equation}
		\label{e: corr Mor1}
		\mathcal{O}bj \xleftarrow{p_s} \mathcal{M}or_1\xrightarrow{p_t}  \mathcal{O}bj 
	\end{equation}
	as an endomorphism $\phi$ on $\mathcal{O}bj$ in the 2-category $\mathbf{Corr}(\on{PreStk})_{\on{all},\on{all}}^{\on{all}}$ of correspondences. The categorical structure gives a monad structure on this endomorphism. Namely, $\phi\circ \phi$ can be identified with the correspondence
	\begin{equation}
		\label{e: corr Mor2}
		\mathcal{O}bj \xleftarrow{p_s} \mathcal{M}or_2\xrightarrow{p_t}  \mathcal{O}bj 
	\end{equation}
	Then $\phi\circ \phi\to \phi$ is induced by the 2-morphism from \eqref{e: corr Mor2} to \eqref{e: corr Mor1} given by the projection
	\[
		p_{02}: \mathcal{M}or_2 \to  \mathcal{M}or_1.
	\]
	that corresponds to the map $[1] \to [2]$, $0\mapsto 0, 1\mapsto 2$. We have a functor between 
	2-categories\footnote{See \cite[Chapter 7]{GaRo3} for the notation.}
	\[
		\on{CrystCat}: \mathbf{Corr}(\on{PreStk})_{\on{all},\on{all}}^{\on{all},2-\on{op}} \to \mathbf{2-Cat}
	\]
	that sends the endomorphism $\phi$ to the endomorphism $p_{s,*}\circ p_t^*$. It follows that the monad structure on $\phi$ gives a comonad structure on $p_{s,*}\circ p_t^*$.

\sssec{}
	Applying to $\Ran^{\subseteq,\bullet}$, we obtain a comonad 
	\[
		\sP : \on{CrystCat}(\Ran) \to \on{CrystCat}(\Ran)
	\]
	whose underlying endofunctor is $(\on{pr}_{\on{small}})_* \circ (\on{pr}_{\on{big}})^*$. In other words, we have (see Sect. \ref{sss:unital insertion})
	\[
		\sP(\ul\bC^{\on{loc}}) \simeq \bC^{\on{loc},\subseteq}.
	\]
	As explained in \emph{loc.cit.}, we have

\begin{lem}
	A local unital structure on $\ul\bC^{\on{loc}}$ is the same as a $\sP$-comodule structure on it, where
	\[
		\sP=(\on{pr}_{\on{small}})_* \circ (\on{pr}_{\on{big}})^* : \on{CrystCat}(\Ran) \to \on{CrystCat}(\Ran)
	\]
	is a comonad acting on $\on{CrystCat}(\Ran)$.
\end{lem}

\sssec{}
	From now on, whenever $\ul\bC^{\on{loc}}$ is equipped with a local unital structure, we view it as a $\sP$-comodule via the above lemma. Note that the coaction functor fits into the following commutative diagram
	\[
		\xymatrix{
			\ul\bC^{\on{loc}} \ar[r]^-{\on{coact}} \ar[d]_{\on{Id}}
			& \sP(\ul\bC^{\on{loc}}) \ar[d]^-\simeq \\
			\ul\bC^{\on{loc}} \ar[r]^-{\on{ins.unit}} 
			& \bC^{\on{loc},\subseteq}.
		}
	\]

\sssec{Example}
	\label{sss: exam const sheaf cat}
	The constant sheaf of categories
	\[
		\bD^{\on{cons}}:=\bC^{\on{glob}}\otimes \ul\Dmod(\Ran) \in \on{CrystCat}(\Ran)
	\]
	has an obvious local unital structure. Note that 
	\[
		\sP(\bC^{\on{glob}}\otimes \ul\Dmod(\Ran) ) \simeq \bC^{\on{glob}}\otimes \ul\Dmod(\Ran^\subseteq),
	\]
	where the RHS is viewed as a sheaf of categories over $\Ran$ via the small projection $\on{pr}_{\on{small}}: \Ran^\subseteq\to \Ran$. The corresponding $\sP$-comodule structure is given by the functor
	\[
		\on{Id} \otimes (\on{pr}_{\on{small}})^!: \bC^{\on{glob}}\otimes \ul\Dmod(\Ran) \to \bC^{\on{glob}}\otimes \ul\Dmod(\Ran^\subseteq).
	\]

	\medskip
	Note that we have an adjunction 
	\[
		\on{coact}^L: \sP(\bD^{\on{cons}}) \rightleftarrows \bD^{\on{cons}} : \on{coact}
	\]
	in the 2-category $\on{CrystCat}(\Ran)$, where the left adjoint $\on{coact}^L$ is given by the functor $\on{Id} \otimes (\on{pr}_{\on{small}})_!$. Also, the right adjoint $\on{coact}$ is fully faithful, i.e. 
	\[
		\on{coact}^L\circ \on{coact} \xrightarrow{\sim} \on{Id}.
	\]
	Indeed, this follows from the contractibility of the map $\on{pr}_{\on{small}}: \Ran^\subseteq \to \Ran$.

\sssec{A general paradigm}
	\label{sss: paradigm colax}
	Let $\sP$ be a comonad acting on a $2$-category $\BS$. Let $\bc,\bd$ be two $\sP$-comodules. For any morphism $f:\bc\to \bd$ in $\BS$, we can talk about (co)lax $\sP$-linear structures on $f$. Namely, a lax $\sP$-linear structure on $f$ is a 2-morphism
	\[
		\xymatrix{
			\bc \ar[r]^-{\on{coact}} \ar[d]_-f
			& \sP(\bc) \ar[d]^-{\sP(f)} \\
			\bd \ar[r]_-{\on{coact}}  \ar@{=>}[ru]
			& \sP(\bd),
		}
	\]
	i.e.,
	\[
		\alpha:\on{coact}\circ f \to \sP(f) \circ \on{coact},
	\]
	equipped with a datum of associativity.

	\medskip
	Recall any comonad $\sP$ has a counit natural transformation $\epsilon:\sP \to \on{Id}$, and for any comodule $\bc$, the composition $\bc\xrightarrow{\on{coact}}\sP(\bc) \xrightarrow{\epsilon}\bc$ is isomorphic to the identity morphism. Then the above datum of associativity in particular says the outer square in the following diagram commutes:
	\[
		\xymatrix{
			\bc \ar[r]^-{\on{coact}} \ar[d]_-f
			& \sP(\bc) \ar[d]^-{\sP(f)} \ar[r]^-\epsilon
			& \bc \ar[d]^-f \\
			\bd \ar[r]_-{\on{coact}}  \ar@{=>}[ru]
			& \sP(\bd) \ar[r]_-\epsilon
			& \bd.
		}
	\]
	In other words, $\alpha$ becomes invertible after composing with the counit $\sP(\bd) \xrightarrow{\epsilon} \bd$.

	\medskip
	Dually, a colax $\sP$-linear structure on $f$ is a 2-morphism
	\[
		\beta:\sP(f) \circ \on{coact} \to \on{coact}\circ f 
	\]
	equipped with a datum of associativity. Also, $\beta$ becomes invertible after composing with $\epsilon$.

	\medskip
	Given a (co)lax $\sP$-linear structure on $f$, we say it is \emph{strict}, or equivalently $f$ is $\sP$-linear, if the above 2-morphism $\alpha$ (resp. $\beta$) is invertible.

\sssec{}
	\label{sssec: def left-lax}
	Now for $\BS:=\on{CrystCat}(\Ran)$, $\bc:=\ul\bC^{\on{loc}}$ and $\bd:=\bD^{\on{const}}:=\bC^{\on{glob}}\otimes \ul\Dmod(\Ran)$, a morphism $f:\bc\to \bd$ in $\BS$ is just a local-to-global functor
	\[
		\ul\sF: \ul\bC^{\on{loc}}\to \bC^{\on{glob}}\otimes \ul\Dmod(\Ran).
	\]	

	\medskip
	Suppose $\ul\bC^{\on{loc}}$ is equipped with a local unital structure, i.e., $\bc$ is equipped with a $\sP$-comodule structure, where recall $\sP(\ul\bC^{\on{loc}})\simeq \bC^{\on{loc},\subseteq}$. As explained in Sect. \ref{sss:F subset Z}, we have

	\begin{lem}
		In the above notations, knowing a lax unital structure on $\ul\sF$ is equivalent to knowing a lax $\sP$-linear structure on $f$. Via this correspondence, the natural transformation
		\[
			\alpha:\on{coact}\circ f \to \sP(f) \circ \on{coact},
		\]
		is given by \eqref{e:unitality transform original}.

	\end{lem}

\sssec{}
	Now we \emph{define} a left-lax unital structure on $\ul\sF$ to be a left-lax $\sP$-linear structre on $f$. This is the homotopically sound definition promised in Sect. \ref{ss: left-lax intro}.

\ssec{Proof of \propref{prop: all to left-lax}}
\label{ss: proof prop: all to left-lax}

\sssec{}
	Using the notations in Sect. \ref{sss: paradigm colax}, the forgetful functor $\iota_{\on{left-lax \to all}}$ is given by
	\[
		\iota_{\on{left-lax \to all}}:\on{Funct}_{\on{left-lax-}\sP}(\bc,\bd) \to \on{Funct}(\bc,\bd)
	\]
	which sends a left-lax $\sP$-linear morphism $f:\bc\to \bd$ to its underlying morphism. Now \propref{prop: all to left-lax} is a particular case of the following general result.

\begin{lem}
	\label{lem: all to left-lax general}
	Let $(\BS,\sP,\bc,\bd)$ be as in Sect. \ref{sss: paradigm colax}. Suppose:
	\begin{itemize}
		\item[(*)]
			The morphism $\on{coact}:\bd\to \sP(\bd)$ has a left adjoint $\on{coact}^L: \sP(\bd)\to \bd$.
	\end{itemize}
	Then the forgetful functor
	\begin{equation}
		\label{e: oblv left-lax to all}
		\on{Funct}_{\on{left-lax-}\sP}(\bc,\bd) \to \on{Funct}(\bc,\bd)
	\end{equation}
	has a left adjoint given by
	\begin{equation}
		\label{e: int all to left-lax}
		\on{Funct}(\bc,\bd)\to \on{Funct}_{\on{left-lax-}\sP}(\bc,\bd),\; f\mapsto \on{coact}^L \circ \sP(f)\circ \on{coact}.
	\end{equation}
\end{lem}

\sssec{Proof}
	The rest of this subsection is devoted to the proof of the lemma. We first define the desired left-lax $\sP$-linear structure on the composition
		\begin{equation}
			\label{e: int all to left-lax value}
			\bc \xrightarrow{\on{coact}} \sP(\bc) \xrightarrow{\sP(f)} \sP(\bd) \xrightarrow{\on{coact}^L} \bd.
		\end{equation}
	It is the composition of the following three left-lax $\sP$-linear structures:
		\begin{itemize}
			\item The coaction morphism $\bc \xrightarrow{\on{coact}} \sP(\bc) $ always has a natural strict $\sP$-linear structure;
			\item The morphism $\sP(f): \sP(\bc) \to \sP(\bd)$ is a strictly $\sP$-linear morphism between cofree $\sP$-comodules;
			\item As the left adjoint of the strictly $\sP$-linear morphism $\bd \xrightarrow{\on{coact}} \sP(\bd) $, the morphism $\sP(\bd) \xrightarrow{\on{coact}^L} \bd$ has a natural left-lax $\sP$-linear structure.
		\end{itemize}

\sssec{}
	To show \eqref{e: int all to left-lax} is indeed left adjoint to \eqref{e: oblv left-lax to all}, we provide the unit and counit natural transformations for this adjunction. 

\sssec{Unit natural transformation}
	\label{sss: left-lax unit natural trans}
				Given $f\in \on{Funct}(\bc,\bd)$, we have the following commutative diagram
				\[
					\xymatrix{
						\bc \ar[r]^-{\on{coact}}
						& \sP(\bc) \ar[d]^-{\sP(f)} \ar[r]^-\epsilon
						& \bc \ar[d]^-f \\
						\bd \ar[r]^-{\on{coact}}
						& \sP(\bd )  \ar[r]^-\epsilon
						& \bd.
					}
				\]
				Then the desired morphism $f\to  \on{coact}^L \circ \sP(f)\circ \on{coact}$ is obtained by passing to left adjoints along the bottom line and using the fact that $\epsilon \circ \on{coact} \simeq \on{Id}$. More precisely, we have a morphism
				\[
					\epsilon \to \epsilon \circ \on{coact} \circ \on{coact}^L \simeq \on{coact}^L,
				\]
				where the first morphism is given by the unit adjunction of $(\on{coact}^L,\on{coact})$. Hence we obtain a morphism
				\begin{equation}
					\label{eqn: unit adjunction comodule}
					f \simeq f\circ \epsilon \circ \on{coact} \simeq \epsilon \circ \sP(f) \circ \on{coact} \to  \on{coact}^L \circ \sP(f)\circ \on{coact}.
				\end{equation}
				which is the value of the desired unit natural transformation at $f$.

\sssec{Counit natural transformation}
	Given $g\in \on{Funct}_{\on{left-lax-}\sP}(\bc,\bd)$, the left-lax $\sP$-linear structure on it provides a morphism in $\on{Funct}_{\on{left-lax-}\sP}(\bc,\bd)$:
				\[
					\sP(g)\circ \on{coact} \to \on{coact} \circ g.
				\]
				Using the adjunction $(\on{coact}^L,\on{coact})$, we obtain a morphism
				\[
					\on{coact}^L\circ \sP(g)\circ \on{coact} \to g
				\]
				which is the value of the desired counit natural transformation at $g$.

\sssec{}
	It is a routine exercise to verify these natural transformations indeed satisfy the axioms of an adjunction. We leave it to the readers.
		
\qed

\sssec{Remark}
	\label{rem: left-lax GUAii}
	Note that \propref{prop: all to left-lax} is indeed a particular case of \lemref{lem: all to left-lax general} because of Sect. \ref{sss: exam const sheaf cat}. Also, unwinding the definitions, for a local-to-global functor $\ul\sF$ and the corresponding $f:\bc\to \bd$, the unit adjunction (\ref{eqn: unit adjunction comodule}) is exactly \eqref{e:correlator F}.

	\medskip
	Moreover, by definition, $\ul\sF$ satisfies the Global Unitality Axioms iff the corresponding $f$ satisfies
	\begin{itemize}
		\item[(i)]
			The unit adjunction (\ref{eqn: unit adjunction comodule}) is invertible;
		\item[(ii)]
			The natural transformation
			\begin{equation}
				\label{e: GUAii comodule}
				\sP(f)\circ \on{coact} \to \on{coact}\circ f
			\end{equation}
			obtained from the inverse of (\ref{eqn: unit adjunction comodule}) is invertible.
	\end{itemize}
	Now (i) implies $f$ has a natural left-lax $\sP$-linear structure given by that of $\on{coact}^L \circ \sP(f)\circ \on{coact}$. A direct diagram chasing shows this left-lax $\sP$-linear structure is exhibited by the natural transformation \eqref{e: GUAii comodule}. 
	Hence (ii) says this left-lax $\sP$-linear structure on $f$ is strict. This was used in the proof of \propref{p:unitality as a property} in 
	Sect. \ref{sss: proof p:unitality as a property colax}.

\ssec{Proof of \propref{prop: all to left-lax adjointable}}
\label{ss: proof prop: all to left-lax adjointable}

\sssec{}
	Unlike \propref{prop: all to left-lax}, the proof of \propref{prop: all to left-lax adjointable} is \emph{not} completely formal. Instead, we need the following particular feature of the setting of local-to-global functors.

\begin{lem}
	\label{lem: lax to st}
	The composition
		\[
			\iota^L_{\on{left-lax \to all}}\circ \iota_{\on{lax\to all}}:
			\on{Funct}^{\on{loc}\to \on{glob},\on{lax-untl}}(\ul\bC^{\on{loc}},\bC^{\on{glob}}) \to \on{Funct}^{\on{loc}\to \on{glob},\on{left-lax-untl}}(\ul\bC^{\on{loc}},\bC^{\on{glob}})
		\]
		takes image in $\on{Funct}^{\on{loc}\to \on{glob},\on{untl}}(\ul\bC^{\on{loc}},\bC^{\on{glob}})$.
\end{lem}

\sssec{Remark}
	The claim of the lemma is a priori \emph{stronger} than the results listed in Sect. \ref{ss:constr int}. Namely, given a local-to-global functor $\ul{\sF}$ equipped with a lax unital structure, Sect. \ref{ss:constr int} says there is \emph{a} strictly unital structure on $\ul{\sF}^{\int \on{ins.unit}}$, while the lemma says the natural left-lax unital structure on it (provided by \propref{prop: all to left-lax}) is strict. Nevertheless, it is easy to see the proof in loc.cit. actually implies this stronger claim. For completeness, we repeat this proof.

\proof
	We will deduce the claim from Lemma \ref{l:pr small untl}. Let $\ul\sF$ be a local-to-global functor equipped with a lax unital structure, and $f:\bc\to \bd$ be the corresponding lax $\sP$-linear morphism. We need to show that the left-lax $\sP$-linear structure on $\on{coact}^L \circ \sP(f)\circ \on{coact}$ is strict. Recall that $\sP(f)\circ \on{coact}$, i.e.,
	\begin{equation}
		\label{e: Pfcoact explicit}
		\bC_\Ran^{\on{loc}} \xrightarrow{\on{ins.unit}} \bC_\Ran^{\on{loc},\subseteq} \xrightarrow{\ul\sF^\subseteq} \bC^{\on{glob}} \otimes \ul{\Dmod}(\Ran^{\subseteq}) 
	\end{equation}
	is naturally $\sP$-linear. On the other hand, $\on{coact}^L: \sP(\bd) \to \bd$, i.e.,
	\begin{equation}
		\label{e: coactL explicit}
		\on{Id} \otimes (\on{pr}_{\on{small}})_!: \bC^{\on{glob}} \otimes \ul{\Dmod}(\Ran^{\subseteq}) \to  \bC^{\on{glob}} \otimes \ul{\Dmod}(\Ran).
	\end{equation}
	has a natural left-lax $\sP$-linear structure because its right adjoint is naturally $\sP$-linear.

	\medskip
	Consider the categorical prestack
	\[
		\Ran \underset{\Ran^{\on{untl}}}\times \Ran^{\subseteq,\on{untl}}
	\]
	where $\Ran^{\subseteq,\on{untl}}\to \Ran^{\on{untl}}$ is the small projection $\on{pr}_{\on{small}}^{\on{untl}}$. Viewing it as a categorical prestack over $\Ran$, we obtain a crystal of categories  over $\Ran$
	\[
		\ul{\Dmod}( \Ran \underset{\Ran^{\on{untl}}}\times \Ran^{\subseteq,\on{untl}} )
	\]
	and a functor
	\begin{equation}
		\label{e: oblv Ransubset untl to Ransubset}
		\mathsf{t}^!:\ul{\Dmod}( \Ran \underset{\Ran^{\on{untl}}}\times \Ran^{\subseteq,\on{untl}} ) \to \ul{\Dmod}(\Ran^{\subseteq})
	\end{equation}
	given by pullback along $\mathsf{t}:\Ran^{\subseteq}\to  \Ran \underset{\Ran^{\on{untl}}}\times \Ran^{\subseteq,\on{untl}} $. Note that both sides of \eqref{e: oblv Ransubset untl to Ransubset} have natural local unital structures and the functor is obviously unital. Now \lemref{l:pr small untl}, combined with a variant of \lemref{l:int over Ran and Ran untl}, says the composition 
	\[
		\ul{\Dmod}( \Ran \underset{\Ran^{\on{untl}}}\times \Ran^{\subseteq,\on{untl}} ) \to \ul{\Dmod}(\Ran^{\subseteq}) \xrightarrow{(\on{pr}_{\on{small}})_!} \ul{\Dmod}(\Ran),
	\]
	which is a priori left-lax unital, is strictly unital. The statement remains true if we tensor it with the DG category $ \bC^{\on{glob}}$. Then we finish the proof because the composition \eqref{e: Pfcoact explicit} factors through \eqref{e: oblv Ransubset untl to Ransubset} by the lax unital structure on $\ul\sF$.

\qed

\sssec{}
	By \lemref{lem: lax to st}, there is a unique functor $\tau_{\on{lax \to st}}$ making the following diagram commute
	\begin{equation}
		\label{e: lax colax sq general local global left adj proof}
		\xymatrix{
			\on{Funct}^{\on{loc}\to \on{glob},\on{untl}}(\ul\bC^{\on{loc}},\bC^{\on{glob}})
			 \ar[d]_-{\iota_{\on{st\to left-lax}}}^-\subset 
			& \on{Funct}^{\on{loc}\to \on{glob},\on{lax-untl}}(\ul\bC^{\on{loc}},\bC^{\on{glob}}) \ar[d]^-{\iota_{\on{lax \to all}}} 
				\ar[l]_{\tau_{\on{lax \to st}}}
			\\
			\on{Funct}^{\on{loc}\to \on{glob},\on{left-lax-untl}}(\ul\bC^{\on{loc}},\bC^{\on{glob}})  &
			\on{Funct}^{\on{loc}\to \on{glob}}(\ul\bC^{\on{loc}},\bC^{\on{glob}})
			\ar[l]_-{\iota^L_{\on{left-lax \to all}}}.
		}
	\end{equation}
	In particular, we obtained:
	\begin{itemize}
		\item 
			A lifting of the endofunctor 
			\[
				\iota_{\on{left-lax \to all}}\circ\iota^L_{\on{left-lax \to all}}: \on{Funct}^{\on{loc}\to \on{glob}}(\ul\bC^{\on{loc}},\bC^{\on{glob}}) \to \on{Funct}^{\on{loc}\to \on{glob}}(\ul\bC^{\on{loc}},\bC^{\on{glob}})
			\]
			to an endofunctor 
			\[
				\iota_{\on{st\to lax}}\circ \tau_{\on{lax \to st}}: \on{Funct}^{\on{loc}\to \on{glob},\on{lax-untl}}(\ul\bC^{\on{loc}},\bC^{\on{glob}})  \to \on{Funct}^{\on{loc}\to \on{glob},\on{lax-untl}}(\ul\bC^{\on{loc}},\bC^{\on{glob}});
			\]
		\item 
			A lifting of the endofunctor 
			\[
				\iota^L_{\on{left-lax \to all}}\circ\iota_{\on{left-lax \to all}}: \on{Funct}^{\on{loc}\to \on{glob},\on{left-lax-untl}}(\ul\bC^{\on{loc}},\bC^{\on{glob}})\to\on{Funct}^{\on{loc}\to \on{glob},\on{left-lax-untl}}(\ul\bC^{\on{loc}},\bC^{\on{glob}}))
			\]
			to an endofunctor 
			\[
				 \tau_{\on{lax \to st}}\circ \iota_{\on{st\to lax}}:  \on{Funct}^{\on{loc}\to \on{glob},\on{untl}}(\ul\bC^{\on{loc}},\bC^{\on{glob}}) \to \on{Funct}^{\on{loc}\to \on{glob},\on{untl}}(\ul\bC^{\on{loc}},\bC^{\on{glob}});
			\]
	\end{itemize}

\sssec{}
	To finish the proof, we only need to lift the unit and counit adjunctions
	\[
		\on{Id} \to \iota_{\on{left-lax \to all}}\circ\iota^L_{\on{left-lax \to all}},\;  \iota^L_{\on{left-lax \to all}}\circ\iota_{\on{left-lax \to all}} \to \on{Id}
	\]
	to 
	\[
		\on{Id} \to \iota_{\on{st\to lax}}\circ \tau_{\on{lax \to st}},\; \tau_{\on{lax \to st}}\circ \iota_{\on{st\to lax}}\to \on{Id}
	\]
	and verify they satisfy the axioms of an adjunction. Indeed, these will induce an equivalence $ \tau_{\on{lax \to st}}\simeq  \iota^L_{\on{st\to lax}}$ such that \eqref{e: lax colax sq general local global left adj proof} can be identified with the Bech--Chevalley natural transformation \eqref{e: lax colax sq general local global left adj}. Then the other claims in \propref{prop: all to left-lax adjointable} follow from definitions.

\sssec{}
	We will lift the counit and unit natural transformations, and leave it to the readers to verify they satisfy the axioms of an adjunction.

\sssec{Lift the counit}
	The lifting of the counit is obvious because the forgetful functor
	\[
		 \on{Funct}^{\on{loc}\to \on{glob},\on{untl}}(\ul\bC^{\on{loc}},\bC^{\on{glob}})\to \on{Funct}^{\on{loc}\to \on{glob},\on{left-lax-untl}}(\ul\bC^{\on{loc}},\bC^{\on{glob}})
	\]
	is fully faithful and the any morphism in the RHS (such as the counit $\iota^L_{\on{left-lax \to all}}\circ\iota_{\on{left-lax \to all}} \to \on{Id}$) has a unique lifting to the LHS as long as its source and target are contained in the LHS.

\sssec{Lift the unit}
	The rest of this subsection is devoted to lift the unit. Recall its definition in Sect. \ref{sss: left-lax unit natural trans}. By definition, the morphism $\epsilon:\sP(\bd) \to \bd$ is
	\[
		\on{Id} \otimes \on{diag}^!: \bC^{\on{glob}} \otimes \ul{\Dmod}(\Ran^{\subseteq}) \to \bC^{\on{glob}} \otimes\ul{\Dmod}(\Ran)
	\]
	and the morphism 
	\[
		\epsilon \to \on{coact}^L
	\]
	is induced by the natural transformation
	\begin{equation}
		\label{e: epsilon to coactL}
		\on{diag}^! \to \on{diag}^!\circ \on{pr}_{\on{small}}^! \circ (\on{pr}_{\on{small}})_! \simeq(\on{pr}_{\on{small}})_!
	\end{equation}
	provided by the isomorphism $\on{pr}_{\on{small}}\circ \on{diag} \simeq \on{Id}$.

\sssec{}
	\label{sss: lax str on epsilon to coactL}
	Recall the restriction functor
	\[
		\mathsf{t}^!:\ul{\Dmod}( \Ran \underset{\Ran^{\on{untl}}}\times \Ran^{\subseteq,\on{untl}} ) \to \ul{\Dmod}(\Ran^{\subseteq}).
	\]
	We claim the (horizontal) composition of \eqref{e: epsilon to coactL} with $\mathsf{t}^!$, i.e.,
	\begin{equation}
		\label{e: epsilon to coactL rest}
		\on{diag}^!\circ \mathsf{t}^! \to  (\on{pr}_{\on{small}})_!\circ \mathsf{t}^!
	\end{equation}
	can be naturally lifted to a natural transformation between lax unital functors
	\[
		\ul{\Dmod}( \Ran \underset{\Ran^{\on{untl}}}\times \Ran^{\subseteq,\on{untl}} )  \to \ul{\Dmod}(\Ran).
	\]
	Here $\ul{\Dmod}( \Ran \underset{\Ran^{\on{untl}}}\times \Ran^{\subseteq,\on{untl}} )$ and $\ul{\Dmod}(\Ran)$ are equipped with the local unital structures given respectively by $\Ran^{\subseteq,\on{untl}}$ and $\Ran^{\on{untl}}$.

\sssec{}
	Consider the map 
	\[
		\on{diag}^{\on{untl}}: \Ran^{\on{untl}}\to \Ran^{\subseteq,\on{untl}}
	\]
	and its left inverse
	\[
		\on{pr}_{\on{small}}^{\on{untl}}:  \Ran^{\subseteq,\on{untl}}\to \Ran^{\on{untl}}.
	\]
	Similar to \eqref{e: epsilon to coactL}, we have a natural transformation
	\begin{equation}
		\label{e: epsilon to coactL untl}
		(\on{diag}^{\on{untl}})^! \to (\on{diag}^{\on{untl}})^! \circ (\on{pr}_{\on{small}}^{\on{untl}})^!\circ (\on{pr}_{\on{small}}^{\on{untl}})_!\simeq (\on{pr}_{\on{small}}^{\on{untl}})_!.
	\end{equation}
	Note that $(\on{diag}^{\on{untl}})^!$ has an obvious lax unital structure, i.e., is a lax functor between sheaves of categories on $\Ran^{\on{untl}}$, while $(\on{pr}_{\on{small}}^{\on{untl}})_!$ has a strictly unital structure by Lemma \ref{l:pr small untl}. The composition \eqref{e: epsilon to coactL untl} is compatible with the lax unital structures on both sides because each natural transformation is.

	\medskip
	By (a variant of) \lemref{l:int over Ran and Ran untl}, when restricted along $\Ran \to \Ran^{\on{untl}}$, \eqref{e: epsilon to coactL untl} gives exactly \eqref{e: epsilon to coactL rest}. In other words, we have proved the claim in Sect. \ref{sss: lax str on epsilon to coactL}.

\sssec{}
	Recall that $\sP(f)\circ \on{coact}$ (which is just \eqref{e: Pfcoact explicit}) factors through $\on{Id}\otimes\mathsf{t}^!$. It follows that \eqref{e: epsilon to coactL rest} induces a natural transformation
	\[
		\epsilon \circ \sP(f)\circ \on{coact} \to \on{coact}^L\circ \sP(f)\circ \on{coact}
	\]
	compatible with the lax unital structures on both sides.

\sssec{}
	On the other hand, we have
	\[
		f \simeq f\circ \epsilon \circ \on{coact} \simeq 
		\epsilon \circ \sP(f)\circ \on{coact} 
	\]
	because $\epsilon$ is the counit of the comonad $\sP$. Moreover, this isomorphism is obviously compatible with the lax unital structures on both sides.

\sssec{}
	Combining the above two subections, we obtain that
	\[
		f \to  \on{coact}^L\circ \sP(f)\circ \on{coact}
	\]
	is naturally compatible with the lax unital structures on both sides. In other words, we have found the desired lifting of the unit adjunction.

\sssec{}
	We leave it to the readers to check the above liftings indeed satisfy the axioms of an adjunction.

\qed[\propref{prop: all to left-lax}]

\section{A homotopical device for coaction} \label{s:device}

The goal of this section is to introduce a homotopical device that will help us carry out the 
constructions in Sects. \ref{ss:action of center} and \ref{ss:z on KL} up to \emph{coherent} homotopy.

\ssec{Associative algebras via mock-simplicial sets} \label{ss:mon cat}

\sssec{} 

Let us recall the following device of encoding the structure of associative algebra (resp., module over a given associative algebra) 
in a monoidal category (see \cite[Sect. 2.2.4]{Lu2}). 

\medskip

Let $\bDelta^{\on{op,mock}}$ be the category of (possibly empty) finite ordered sets. The operation of (ordered) union 
defines on $\bDelta^{\on{op,mock}}$ a structure of monoidal category. Its monoidal unit is $\emptyset$.

\medskip

In what follows we will denote\footnote{The category $\bDelta^{\on{mock}}$ is equivalent to the subcategory of $\bDelta$ consisting of \emph{active} morphisms, i.e. those that preserve the maximal and minimal element.  
For our purposes, the nonstandard ``mock" terminology will be more convenient. We caution the reader that
the usual simplex category, denoted $\bDelta$, is that of non-empty finite ordered sets. Note that arrows in $\bDelta$ and $\bDelta^{\on{mock}}$
go in opposite directions.}
$$\bDelta^{\on{mock}}:=(\bDelta^{\on{op,mock}})^{\on{op}}.$$

\sssec{}

The datum of a unital associative algebra in a monoidal category $\fA$ is equivalent to that of a monoidal functor
$$\bDelta^{\on{op,mock}}\to \fA.$$

Under this correspondence, for a given functor $F:\bDelta^{\on{op,mock}}\to \fA$, the corresponding algebra object $\fa\in \fA$
is $F(\{*\})$. The unit in $\fa$ is given by the map
$$\one_\fA=F(\emptyset)\to F(\{*\})=\fa,$$
corresponding to the (unique) map $\emptyset\to \{*\}$. 

\medskip

The binary operation on $\fa$ corresponds to the map
$$\fa\otimes \fa=F(\{*\})\otimes F(\{*\})\simeq F(\{1,2\})\to F(\{*\})=\fa,$$
where the arrows is given by the (unique) map $\{1,2\}\to \{*\}$ in $\bDelta^{\on{op,mock}}$. 

\sssec{} \label{sss:Delta *}

Let $\bDelta_*^{\on{op,mock}}$ be the category of non-empty finite ordered sets, pointed by their maximal element.
The category $\bDelta_*^{\on{op,mock}}$ is naturally a module over $\bDelta^{\on{op,mock}}$. 

\sssec{}

Let us be given a monoidal category $\fA$. Let $\fa$ be an associative algebra in $\fA$, 
thought of as a monoidal functor
$$F_a:\bDelta^{\on{op,mock}}\to \fA.$$

\medskip

Then the category  
$$\fa\mod(\fA)$$
of $\fa$-modules in $\fA$ is equivalent to that of functors of left $\bDelta^{\on{op,mock}}$-module
categories
$$\bDelta_*^{\on{op,mock}}\to \fA,$$
where $\fA$ is a left $\bDelta^{\on{op,mock}}$-module via $F_a$.

\sssec{}

Under this correspondence, given a functor $F_m:\bDelta_*^{\on{op,mock}}\to \fA$, the object
of $\fA$ underlying the corresponding $\fa$-module is
$$\fm:=F_m(\{*\}).$$

The action map $\fa\otimes M\to M$ is given by
$$\fa\otimes \fm\simeq F_a(\{0\})\otimes F_m(\{*\}) \simeq F_m(\{0,*\})\to F_m(\{*\})\simeq \fm.$$

\sssec{}

As above, in what follows we will denote
$$\bDelta_*^{\on{mock}}:=(\bDelta_*^{\on{op,mock}})^{\on{op}}.$$

We will refer to functors $\bDelta^{\on{op,mock}}\to \fA$ (resp., $\bDelta_*^{\on{op,mock}}\to \fA$, $\bDelta^{\on{mock}}\to \fA$,
$\bDelta_*^{\on{mock}}\to \fA$) as \emph{mock-simplicial} (resp., pointed mock-simplicial, mock cosimplicial, pointed
mock cosimplicial) objects of $\fA$. 

\ssec{Factorization categories attached to \emph{associative} factorization algebras} \label{ss:mon cat to ass fact alg}

Recall that our goal is to carry out the constructions in Sects. \ref{ss:action of center} and \ref{ss:z on KL}.
We will achieve this by introducing appropriate objects to feed into the machine in \secref{ss:mon cat}. 

\sssec{} \label{sss:mon cat to ass fact alg}

Let $\fA_1$ be a monoidal category and let $\fC_1$ be a module category over it. Note that given two such pairs,
we can talk about \emph{strictly} and \emph{right-lax} monoidal functors of pairs
$$(F_a,F_c):(\fA_1,\fC_1)\to (\fA_1,\fC_2),$$
i.e., $F_a:\fA_1\to \fA_2$ is righ-lax monoidal, and $F_c:\fC_1\to \fC_2$ is right-lax compatible with the actions.

\medskip

We introduce several such pairs now. 

\sssec{}

Take 
$$\fA_1=\bDelta^{\on{mock}},\,\, \fA_2=(\on{FactAlg}^{\on{untl}}(X))^{\on{op}},\,\, \fA_3=\on{FactCat}_{\on{lax}}^{\on{untl,lax}}(X),$$
where the subscript ``lax" means that we are considering lax factorization categories. 

\sssec{}

By \secref{ss:mon cat}, an object 
\begin{equation} \label{e:ass fact alg}
\CR\in \on{AssAlg}(\on{FactAlg}^{\on{untl}}(X))
\end{equation}
can be viewed as a \emph{stricty} monoidal functor
\begin{equation} \label{e:ass fact alg bis}
\bDelta^{\on{op,mock}}\to \on{FactAlg}^{\on{untl}}(X),
\end{equation}
and hence also as a \emph{stricty} monoidal functor 
\begin{equation} \label{e:ass fact alg bis op}
\bDelta^{\on{mock}}\to (\on{FactAlg}^{\on{untl}}(X))^{\on{op}}.
\end{equation}

\sssec{}

In addition, we have a naturally defined \emph{right-lax} monoidal functor 
\begin{equation}  \label{e:from fact alg to fact cat again}
(\on{FactAlg}^{\on{untl}})^{\on{op}}\to \on{FactCat}_{\on{lax}}^{\on{untl,lax}}(X), \quad  \CR\to \CR\mod^{\on{fact}}, \quad
(\CR_1\overset{\phi}\to \CR_2)\rightsquigarrow \Res_\phi.
\end{equation} 

\medskip

Composing, for $\CR$ as above, we obtain a \emph{right-lax} monoidal functor
\begin{equation} \label{e:Delta to FactCat}
\bDelta^{\on{mock}}\to \on{FactCat}_{\on{lax}}^{\on{untl,lax}}(X), \quad n\mapsto \CR^{\otimes n}\mod^{\on{fact}}.
\end{equation} 

\sssec{} \label{sss:mon cat to ass fact mod}

Take 
$$\fC_1=\bDelta_*^{\on{mock}},\,\,\fC_2=\on{FactCat}_{\on{lax}}^{\on{untl,lax}}(X),\,\, 
\fC_3=\on{FactCat}_{\on{lax}}^{\on{untl,lax}}(X),$$
where:

\begin{itemize}

\item $\bDelta^{\on{mock}}$ acts on $\bDelta_*^{\on{mock}}$ as in \secref{sss:Delta *};

\medskip

\item $\on{FactCat}_{\on{lax}}^{\on{untl,lax}}(X)$ acts on itself via the (symmetric) monoidal structure
on $\on{FactCat}_{\on{lax}}^{\on{untl,lax}}(X)$;

\medskip

\item $(\on{FactAlg}^{\on{untl}}(X))^{\on{op}}$ acts on $\on{FactCat}_{\on{lax}}^{\on{untl,lax}}(X)$ by
$$\CR,\bA\mapsto \CR\mod^{\on{fact}}(\bA)\simeq (\CR\otimes \one_\bA)\mod^{\on{fact}}(\bA).$$

\end{itemize}

\medskip

We note that the above action of $(\on{FactAlg}^{\on{untl}})^{\on{op}}$ on $\on{FactCat}_{\on{lax}}^{\on{untl,lax}}(X)$
is monoidal thanks to \lemref{l:two algebras mods}. 

\sssec{} \label{sss:abs R acts on vacuum}

Let $\CR$ be as in \eqref{e:ass fact alg} and suppose that it acts on the unit $\one_\bC\in \bC$ for some
$\bC\in \on{FactCat}_{\on{lax}}^{\on{untl,lax}}(X)$.

\medskip

Then the functor \eqref{e:ass fact alg bis} extends to a \emph{strictly} monoidal functor of pairs 
\begin{equation}  \label{e:ass fact alg mod}
(\bDelta^{\on{op,mock}},\bDelta_*^{\on{op,mock}})\to 
(\on{FactAlg}^{\on{untl}}(X),(\on{FactCat}_{\on{lax}}^{\on{untl,lax}}(X))^{\on{op}}), \quad *\mapsto \bC.
\end{equation}

\medskip

Hence, we obtain a \emph{strictly} monoidal functor of pairs 
\begin{equation}  \label{e:ass fact alg mod op}
(\bDelta^{\on{mock}},\bDelta_*^{\on{mock}})\to 
((\on{FactAlg}^{\on{untl}}(X))^{\on{op}},\on{FactCat}_{\on{lax}}^{\on{untl,lax}}(X)), \quad *\mapsto \bC.
\end{equation}

\sssec{}

In addition, the functor \eqref{e:from fact alg to fact cat again} extends to a \emph{right-lax} monoidal functor
\begin{equation}  \label{e:from fact alg to fact cat mod}
((\on{FactAlg}^{\on{untl}}(X))^{\on{op}},\on{FactCat}_{\on{lax}}^{\on{untl,lax}}(X))\to
(\on{FactCat}_{\on{lax}}^{\on{untl,lax}}(X),\on{FactCat}_{\on{lax}}^{\on{untl,lax}}(X)),
\end{equation}
which acts as identity on the module component.

\sssec{} \label{sss:abs Delta module outcome}

Combining, we obtain that for $(\CR,\bC)$ as in \secref{sss:abs R acts on vacuum} we obtain that
\eqref{e:Delta to FactCat} extends to a \emph{right-lax} monoidal functor of pairs 
\begin{equation} \label{e:Delta * to FactCat}
(\bDelta_*^{\on{mock}},\bDelta_*^{\on{mock}})\to (\on{FactCat}_{\on{lax}}^{\on{untl,lax}}(X),\on{FactCat}_{\on{lax}}^{\on{untl,lax}}(X))
\end{equation}
that sends 
$$(*\in \bDelta_*^{\on{mock}})\,\rightsquigarrow \, \bC.$$

\ssec{The renormalization step}

\sssec{}

We apply \secref{sss:abs Delta module outcome} to 
$$\CR=\fz_\fg \text{ and } \bC=\hg\mod_\crit,$$
and consider the resulting right-lax monoidal functor of pairs \eqref{e:Delta * to FactCat}; denote it 
\begin{equation} \label{e:Delta module outcome KM}
(F_{\fz_\fg},F_{\hg\mod_\crit}).
\end{equation}

\sssec{}

We now introduce several more (symmetric) monoidal categories. Let 
$$\fA=\fA_3:=\on{FactCat}_{\on{lax}}^{\on{untl,lax}}(X)$$
be as above.

\medskip

We let $\fA'$ be the following 1-full subcategory of $\fA$, to be denoted 
$$\on{FactCat}_{\on{lax}}^{\on{untl,lax}}(X)^{\on{t-str}}.$$

Its objects are 
lax unital factorization categories, equipped with a t-structure.  For a pair of objects $\bC_1,\bC_2$, we let
$$\Maps_{\on{FactCat}_{\on{lax}}^{\on{untl,lax}}(X)^{\on{t-str}}}(\bC_1,\bC_2)\subset 
\Maps_{\on{FactCat}_{\on{lax}}^{\on{untl,lax}}(X)}(\bC_1,\bC_2)$$
be the full subcategory consisting of \emph{left t-exact} functors.  

\sssec{}

We let $\fA''$ be the full subcategory of the category of arrows in $\fA'$, whose objects are those pairs
$$\Phi:\bC^{\on{ren}}\to \bC,$$
for which:

\begin{itemize}

\item $\Phi$ induces an equivalence between the eventually coconnective subcategories;

\item $\bC^{\on{ren}}$ is compactly generated by objects that are eventually coconnective. 

\end{itemize}

Note that the tautological forgetful functor
$$\fA''\to \fA', \quad (\Phi:\bC^{\on{ren}}\to \bC)\mapsto \bC$$
is 1-fully faithful, i.e., induces a fully faithful functor on spaces of morphisms. 

\medskip

Explicitly, given two pairs $\bC_1^{\on{ren}}\overset{\Phi_1}\to \bC_1$ and $\bC_2^{\on{ren}}\overset{\Phi_1}\to \bC_2$,
a functor $\Psi:\bC_1\to \bC_2$ in $\fA'$ lifts to $\fA''$ if and only if the ind-extension of
$$(\bC_1^{\on{ren}})^c \subset (\bC_1^{\on{ren}})^{>-\infty}\simeq \bC_1^{>-\infty} \overset{\Psi}\to 
\bC_2^{>-\infty} \simeq (\bC_2^{\on{ren}})^{>-\infty}$$
is left t-exact (equivalently, has a bounded cohomological amplitude on the left). 

\sssec{}

Finally, we let $\fA'''$ be again $\on{FactCat}_{\on{lax}}^{\on{untl,lax}}(X)$, and
we consider the forgetful functor
\begin{equation} \label{e:forget aux}
\fA''\to \fA''', \quad (\Phi:\bC^{\on{ren}}\to \bC)\mapsto \bC^{\on{ren}}.
\end{equation} 

\sssec{}

Note now that given a lax monoidal functor
$$(\bDelta^{\on{mock}},\bDelta_*^{\on{mock}})\to (\fA,\fA),$$
in order to lift to a lax monoidal functor 
$$(\bDelta^{\on{mock}},\bDelta_*^{\on{mock}})\to (\fA'',\fA''),$$
it suffices to do so at the level of objects and 1-morphisms, i.e., it is sufficient to
do so \emph{at the homotopy level} (moreover, the lift at the level of objects defines
it completely). 

\medskip

We start with the functor $(F_{\fz_\fg},F_{\hg\mod_\crit})$ of \eqref{e:Delta module outcome KM}. We lift to a functor
$$(F_{\fz_\fg},F_{\hg\mod_\crit})'':(\bDelta^{\on{mock}},\bDelta_*^{\on{mock}})\to (\fA'',\fA'')$$
at the level of objects by sending
$$(n\in \bDelta^{\on{mock}})\,\, \rightsquigarrow \,\, 
\left(\IndCoh^*(``\Spec"(\fZ_\fg))^{\otimes n}\to \fz_\fg^{\otimes n}\mod^{\on{fact}}\right)$$
$$(\{n,*\}\in \bDelta_*^{\on{mock}}) \,\,\rightsquigarrow \,\, \left(\IndCoh^*(``\Spec"(\fZ_\fg))^{\otimes n}\otimes \hg\mod_\crit\to
\fz_\fg^{\otimes n}\mod(\hg\mod_\crit)\right).$$

The existence of the lift at the level of 1-morphisms is guaranteed by \lemref{l:coaction functor term n}.

\sssec{} \label{sss:final lift}

We compose $(F_{\fz_\fg},F_{\hg\mod_\crit})''$ with \eqref{e:forget aux} to obtain a right-lax monoidal functor of pairs 
$$(F_{\fz_\fg},F_{\hg\mod_\crit})''':
(\bDelta^{\on{mock}},\bDelta_*^{\on{mock}})\to (\on{FactCat}_{\on{lax}}^{\on{untl,lax}}(X),\on{FactCat}_{\on{lax}}^{\on{untl,lax}}(X)).$$

However, by construction, the latter functor is strictly monoidal. I.e., the functor
$$F'''_{\fz_\fg}:\bDelta^{\on{mock}}\to \on{FactCat}_{\on{lax}}^{\on{untl,lax}}(X)$$
is strictly monoidal, and 
$$F'''_{\hg\mod_\crit}:\bDelta_*^{\on{mock}}\to \on{FactCat}_{\on{lax}}^{\on{untl,lax}}(X)$$
is a functor between $\bDelta^{\on{mock}}$-module categories. 

\medskip

Applying \secref{ss:mon cat} to $(F_{\fz_\fg},F_{\hg\mod_\crit})'''$ we obtain the desired coaction of 
$\IndCoh^*(``\Spec"(\fZ_\fg))$ on $\hg\mod_\crit$.

\medskip

This completes the construction from \secref{ss:action of center}. 

\ssec{Adding another monoidal category} \label{ss:hom device 2}

In order to carry out the construction in \secref{ss:z on KL}, we need to enhance the setting 
of \secref{ss:mon cat to ass fact alg}.

\sssec{}

We modify the setting of \secref{sss:mon cat to ass fact alg}, and we now take $\fA_2$ to be the category, denoted
$$(\on{FactAlg}^{\on{op}}\in \on{FactCat}_{\on{lax}})^{\on{untl}}(X),$$
whose objects are pairs $(\bA,\CR)$, where:

\medskip

\begin{itemize}

\item $\bA\in \on{FactCat}_{\on{lax}}^{\on{untl}}(X)$;

\item $\CR\in \on{FactAlg}^{\on{untl}}(X,\bA)$.

\end{itemize} 

\medskip

The space of morphisms 
$$\Maps_{(\on{FactAlg}^{\on{op}}\in\on{FactCat}_{\on{lax}})^{\on{untl}}(X)}((\bA_1,\CR_1),(\bA_2,\CR_2))$$
consists of pairs: $(\Phi,\phi)$, where:

\medskip

\begin{itemize}

\item $\Phi$ is a lax unital functor $\bA_1\to \bA_2$;

\medskip

\item $\phi$ is a map of unital factorization algebras in $\bA_2$
$$\CR_2\to \Phi(\CR_1).$$

\end{itemize} 

\begin{rem} \label{r:straighten fact cat}

Note that we can interpret $(\on{FactAlg}^{\on{op}}\in \on{FactCat}_{\on{lax}})^{\on{untl}}(X)$, equipped with the
forgetful functor to $\on{FactCat}_{\on{lax}}^{\on{untl}}(X)$ as the co-Cartesian fibration corresponding to the functor
$$\on{FactCat}_{\on{lax}}^{\on{untl}}(X)\to \inftyCat, \quad \bA\mapsto (\on{FactAlg}^{\on{untl}}(X,\bA))^{\on{op}}.$$

\medskip

The monoidal structure on $(\on{FactAlg}^{\on{op}}\in \on{FactCat}_{\on{lax}})^{\on{untl}}(X)$ corresponds to the
right-lax monoidal structure on the above functor. 

\end{rem}

\sssec{} \label{sss:2 to 3 dec}

We let $\fA_3$ be the same as  in \secref{sss:mon cat to ass fact mod}, i.e., 
$\on{FactCat}_{\on{lax}}^{\on{untl,lax}}(X)$. 

\medskip

Note that we have a lax monoidal functor
\begin{equation}  \label{e:from fact alg to fact cat dec}
(\on{FactAlg}^{\on{op}}\in \on{FactCat}_{\on{lax}})^{\on{untl}}(X)\to \on{FactCat}_{\on{lax}}^{\on{untl,lax}}(X), \quad (\bA,\CR)\mapsto \CR\mod^{\on{fact}}(\bA).
\end{equation}

\sssec{}

We take $\fC_2=\fA_2$ and $\fC_3=\fA_3$. The functor \eqref{e:from fact alg to fact cat dec} gives rise to a right-lax monoidal functor of pairs 
\begin{multline}  \label{e:from fact alg to fact cat mod dec}
\left((\on{FactAlg}^{\on{op}}\in \on{FactCat}_{\on{lax}})^{\on{untl}}(X),(\on{FactAlg}^{\on{op}}\in \on{FactCat}_{\on{lax}})^{\on{untl}}(X)\right)\to \\
\to (\on{FactCat}_{\on{lax}}^{\on{untl,lax}}(X),\on{FactCat}_{\on{lax}}^{\on{untl,lax}}(X)). 
\end{multline}

\sssec{}

We now explain a procedure that gives rise to strictly monoidal functors
$$\bDelta^{\on{mock}}\to (\on{FactAlg}^{\on{op}}\in \on{FactCat}_{\on{lax}})^{\on{untl}}(X)$$
and strictly monoidal functors of pairs
$$(\bDelta^{\on{mock}},\bDelta_*^{\on{mock}})\to 
\left((\on{FactAlg}^{\on{op}}\in \on{FactCat}_{\on{lax}})^{\on{untl}}(X),(\on{FactAlg}^{\on{op}}\in \on{FactCat}_{\on{lax}})^{\on{untl}}(X)\right).$$

\sssec{}

Let $(\on{FactAlg}\in \on{FactCat}_{\on{lax}})^{\on{untl}}(X)$ denote the category defined as follows. 
It objects are pairs $(\bA,\CR)$, where:

\medskip

\begin{itemize}

\item $\bA\in \on{FactCat}_{\on{lax}}^{\on{untl}}(X)$;

\item $\CR\in \on{FactAlg}^{\on{untl}}(X,\bA)$.

\end{itemize} 

\medskip

The space of morphisms 
$$\Maps_{(\on{FactAlg}\in \on{FactCat}_{\on{lax}})^{\on{untl}}(X)}((\bA_1,\CR_1),(\bA_2,\CR_2))$$
consists of pairs: $(\Phi,\phi)$, where:

\begin{itemize}

\item $\Phi$ is a strictly unital functor $\bA_1\to \bA_2$;

\item $\phi$ is a map of unital factorization algebras in $\bA_2$
$$\Phi(\CR_1)\to \CR_2.$$

\end{itemize} 

\begin{rem}

As in Remark \ref{r:straighten fact cat}, the category $(\on{FactAlg}\in \on{FactCat}_{\on{lax}})^{\on{untl}}(X)$, equipped with the
forgetful functor to $\on{FactCat}_{\on{lax}}^{\on{untl}}(X)$, is the co-Cartesian fibration corresponding to the functor
$$\on{FactCat}_{\on{lax}}^{\on{untl}}(X)\to \inftyCat, \quad \bA\mapsto \on{FactAlg}^{\on{untl}}(X,\bA).$$

\end{rem}

\sssec{}  \label{sss:create functor to pairs dec}

Let $\bA$ a monoidal unital lax factorization category $\bA$ and let $\CR\in \on{FactAlg}^{\on{untl}}(X,\bA)$ be an object, 
equipped with a structure of \emph{associative algebra}, in the sense of the monoidal structure on $\bA$.

\medskip

We can think of $(\bA,\CR)$ as an associative algebra object in $(\on{FactAlg}\in \on{FactCat}_{\on{lax}})^{\on{untl}}(X)$. 
Hence, by \secref{ss:mon cat}, it gives rise to a (strictly) monoidal functor
$$F_{\bA,\CR}: \bDelta^{\on{op,mock}} \to (\on{FactAlg}\in \on{FactCat}_{\on{lax}})^{\on{untl}}(X).$$

\medskip

Let $\bC$ be a unital lax factorization category, equipped with an action of $\bA$ as a monoidal factorization category. 
Suppose, moreover, that $\CR$ acts on $\one_\bA$ in the sense of the action of $\bA$ on $\bC$. 

\medskip

Then we can consider $(\bC,\one_\bC)$ as a module over $(\bA,\CR)$ in $(\on{FactAlg}\in \on{FactCat}_{\on{lax}})^{\on{untl}}(X)$.
Hence, by \secref{ss:mon cat}, the functor $F_{\bA,\CR}$ extends to a (strictly) monoidal functor of pairs
$$(F_{\bA,\CR},F_\bC):(\bDelta^{\on{op,mock}},\bDelta_*^{\on{op,mock}})\to
\left((\on{FactAlg}\in \on{FactCat}_{\on{lax}})^{\on{untl}}(X),(\on{FactAlg}\in \on{FactCat}_{\on{lax}})^{\on{untl}}(X)\right).$$

\sssec{}

Let 
$$((\on{FactAlg}\in \on{FactCat}_{\on{lax}})^{\on{untl}}(X))_{\on{adj}} \subset (\on{FactAlg}\in \on{FactCat}_{\on{lax}})^{\on{untl}}(X)$$
be a 1-full subcategory, where we take the same objects, but as 1-morphisms we let
\begin{multline*} 
\Maps_{((\on{FactAlg}\in \on{FactCat}_{\on{lax}})^{\on{untl}}(X))_{\on{adj}}}((\bA_1,\CR_1),(\bA_2,\CR_2))\subset \\
\subset \Maps_{(\on{FactAlg}\in \on{FactCat}_{\on{lax}})^{\on{untl}}(X)}((\bA_1,\CR_1),(\bA_2,\CR_2))
\end{multline*} 
be the full subcategory consisting of those pairs $(\Phi,\phi)$, for which $\Phi$ admits a factorization right adjoint,
which then automatically acquires a lax unital structure (see \secref{sss:adj of unital}). 

\sssec{}

Note that the operation of passage to the right adjoint functor defines a (1-fully faithful, symmetric) monoidal functor
\begin{equation} \label{e:passage to right adj}
((\on{FactAlg}\in \on{FactCat}_{\on{lax}})^{\on{untl}}(X))_{\on{adj}} \to
((\on{FactAlg}^{\on{op}}\in \on{FactCat}_{\on{lax}})^{\on{untl}}(X))^{\on{op}},
\end{equation} 
which acts as identity on objects, and sends
$$(\Phi:\bA_1\to \bA_2,\phi(\CR_1)\to \CR_2) \,\rightsquigarrow \, 
(\Phi^R:\bA_2\to \bA_1, \psi:\CR_1\to \Phi^R(\CR_2)),$$
where $\psi$ is obtained from $\phi$ by adjunction. 

\sssec{} \label{sss:pass to right adjoints}

Therefore, given a (strictly) monoidal functor of pairs
$$(\bDelta^{\on{op,mock}},\bDelta_*^{\on{op,mock}})\to 
\left((\on{FactAlg}\in \on{FactCat}_{\on{lax}})^{\on{untl}}(X),(\on{FactAlg}\in \on{FactCat}_{\on{lax}})^{\on{untl}}(X)\right),$$
which at the level of 1-morphisms lands in 
\begin{equation} \label{e:adjntbl subcat}
\left(((\on{FactAlg}\in \on{FactCat}_{\on{lax}})^{\on{untl}}(X))_{\on{adj}},
((\on{FactAlg}\in \on{FactCat}_{\on{lax}})^{\on{untl}}(X))_{\on{adj}}\right),
\end{equation} 
by passing to right adjoints, we can create from it a (strictly) monoidal functor of pairs
$$(\bDelta^{\on{mock}},\bDelta_*^{\on{mock}})\to 
\left((\on{FactAlg}^{\on{op}}\in \on{FactCat}_{\on{lax}})^{\on{untl}}(X),(\on{FactAlg}^{\on{op}}\in \on{FactCat}_{\on{lax}})^{\on{untl}}(X)\right).$$

\ssec{Applying the paradigm} \label{ss:coaction functor term KL coh hom}

We are finally ready to complete the construction from \secref{ss:z on KL}. 

\sssec{}

In the context of \secref{sss:create functor to pairs dec}, we take 
$$\bA:=\Rep(\cG),\,\,  \bC:=\KL(G)_\crit,\,\, \CR=R_{\cG,\Op}.$$

The data of action from \secref{sss:create functor to pairs dec} is provided by \secref{sss:R Op acts on Vac}. 

\medskip

Denote the resulting monoidal functor of pairs by
\begin{multline*} 
(F_{\Rep(\cG),R_{\cG,\Op}},F_{\KL(G)_\crit}):
(\bDelta^{\on{op,mock}},\bDelta_*^{\on{op,mock}})\to \\
\to \left((\on{FactAlg}\in \on{FactCat}_{\on{lax}})^{\on{untl}}(X),(\on{FactAlg}\in \on{FactCat}_{\on{lax}})^{\on{untl}}(X)\right).
\end{multline*} 

\sssec{}

Note now that since $\Rep(\cG)$ is rigid, the above functor $(F_{\Rep(\cG),R_{\cG,\Op}},F_{\KL(G)_\crit})$ lands in 
the subcategory \eqref{e:adjntbl subcat}. 

\medskip

Hence, by \secref{sss:pass to right adjoints}, we can produce from it a (strictly) monoidal functor of pairs
\begin{multline*} 
(F_{\Rep(\cG),R_{\cG,\Op}},F_{\KL(G)_\crit})_{\on{adj}}:
(\bDelta^{\on{mock}},\bDelta_*^{\on{mock}})\to  \\
\to \left((\on{FactAlg}^{\on{op}}\in \on{FactCat}_{\on{lax}})^{\on{untl}}(X),(\on{FactAlg}^{\on{op}}\in \on{FactCat}_{\on{lax}})^{\on{untl}}(X)\right).
\end{multline*} 

\sssec{}

Composing with \eqref{e:from fact alg to fact cat mod dec}, we obtain a right-lax monoidal functor of pairs
$$(F_{\Rep(\cG),R_{\cG,\Op}},F_{\KL(G)_\crit})_{\on{adj}}^\sim:
(\bDelta^{\on{mock}},\bDelta_*^{\on{mock}})
\to (\on{FactCat}_{\on{lax}}^{\on{untl,lax}}(X),\on{FactCat}_{\on{lax}}^{\on{untl,lax}}(X)).$$

\sssec{}

The final step consists of lifting the functor $(F_{\Rep(\cG),R_{\cG,\Op}},F_{\KL(G)_\crit})_{\on{adj}}^\sim$ to a 
(strictly) monoidal functor of pairs
$$(F_{\Rep(\cG),R_{\cG,\Op}},F_{\KL(G)_\crit})'_{\on{adj}}:
(\bDelta^{\on{mock}},\bDelta_*^{\on{mock}})
\to (\on{FactCat}_{\on{lax}}^{\on{untl,lax}}(X),\on{FactCat}_{\on{lax}}^{\on{untl,lax}}(X)),$$
which at the level of objects sends
$$n\mapsto \IndCoh^*(\Op^\mf_\cG)^{\otimes n}, \quad
(n\sqcup *)\mapsto \IndCoh^*(\Op^\mf_\cG)^{\otimes n}\otimes \KL(G)_\crit.$$

\medskip

This is achieved by repeating the procedure in \secref{sss:final lift} using \lemref{l:coaction functor term KL n}.


\newpage 


\begin{thebibliography}{99}

\bibitem[AF]{AF} D.~Ayala and J.~Francis, {\it Fibrations of $\infty$-categories}, Higher Structures {\bf 4} (2020), 168--265. 

\bibitem[AG]{AG} D.~Arinkin and D.~Gaitsgory, {\it Singular support of coherent sheaves, and the geometric Langlands conjecture}, 
Selecta Math. N.S. {\bf 21} (2015), 1--199.

\bibitem[AGKRRV]{AGKRRV} D.~Arinkin, D.Gaitsgory, D.~Kazhdan, S.~Raskin, N.~Rozenblyum and Y.~Varshavsky, \newline
{\em The stack of local systems with restricted variation and geometric Langlands theory with nilpotent singular support},
arXiv:2010.01906. 

\bibitem[AMR]{AMR} D.~Ayala, A.~Mazel-Gee and N.~Rozenblyum, {\it Stratified non-commutative geometry}, arXiv:1910.14602. 

\bibitem[BD1]{BD1} A.~Beilinson and V.~Drinfeld, {\it Quantization of Hitchin's integrable system and Hecke eigensheaves}, 
available at http://people.math.harvard.edu/$\sim$gaitsgde/grad$\underline{\text{\,\,\,}}$2009/

\bibitem[BD2]{BD2} A.~Beilinson and V.~Drinfeld, {\it Chiral algebras}, American Mathematical Soc. {\bf 51} (2004).

\bibitem[Be1]{Be1} D.~Beraldo, {\it On the extended Whittaker category}, Selecta Mathematica {\bf 25} (2019). 

\bibitem[Be2]{Be2} D.~Beraldo, {\it On the geometric Ramanujan conjecture}, arXiv:2103.17211.

\bibitem[Bei]{Bei} A.~Beilinson, {\it Remarks on topological algebras}, arXiv:0711.2527.

\bibitem[Bh]{Bh} B.~Bhatt, {\it Algebraization and Tannaka duality}, Cambridge Journal of Math {\bf 4} (2016), 403-461. 


\bibitem[Bogd]{Bogd} E.~Bogdabova, {\it in preparation}. 

\bibitem[CR]{CR} J.~Campbell and S.~Raskin, {\it Langlands duality on the Beilinson-Drinfeld Grassmannian}, arXiv:2310.19734. 

\bibitem[Ch1]{Ch1} L.~Chen, {\it Nearby cycles on Drinfeld-Gaitsgory-Vinberg interpolation Grassmannian and long intertwining functor},
Duke Math. Jour. {\bf 172} (2023), 447--553. 

\bibitem[Ch2]{Ch2} L.~Chen, {\it Nearby cycles and dualities in geometric Langlands program}, PhD Thesis (Harvard University 2021), \newline
available at: https://windshower.github.io/linchen/papers/dissertation.pdf

\bibitem[CF]{CF} L.~Chen and Y.~Fu, {\it An extension of Kazhdan-Lusztig equivalence}, arXiv:2111.14606

\bibitem[CFGY]{CFGY} L.~Chen, Y.~Fu, D.~Gaitsgory and D.~Yang, {\it Loop group actions via factorization}, preprint.

\bibitem[Dh]{Dh} G.~Dhillon, {\it Semi-infinite cohomology and the linkage principle for W-algebras}, Advances in Mathematics {\bf 381} (2021). 

\bibitem[Dr]{Dr} V.~Drinfeld, {\it Two-dimensional $\ell$--adic representations of the fundamental group of a curve over a finite field
and automorphic forms on $GL(2)$}, Amer. Jour. of Math. {\bf 105} (1983), 85--114. 

\bibitem[DG1]{DrGa1} V.~Drinfeld and D.~Gaitsgory, {\it On some finiteness questions for algebraic stacks}, Geometric and Functional Analysis. {\bf 23} (2013), 149--294.

\bibitem[DrGa2]{DrGa2} V.~Drinfeld and D.~Gaitsgory, {\it Compact generation of the category of D-modules on the stack of G-bundles on a curve}
Cambridge Math Journal, {\bf 3} (2015), 19--125. 

\bibitem[FR]{FR} J.~Faergeman and S.~Raskin, {\it Non-vanishing of geometric Whittaker coefficients for reductive groups}, \newline
arXiv:2207.02955.

\bibitem[FraG]{FraG} J.~Francis and D.~Gaitsgory, {\it Chiral Koszul duality}, Selecta Math. (N.S.) {\bf 18} (2012), 27--87.

\bibitem[FF]{FF} B.~Feigin and E.~Frenkel, {\it Affine Kac-Moody algebras at the critical level and Gelfand-Dikii algebras}, International Journal of Modern Physics A. {\bf 7} (1992), 197--215.

\bibitem[FG1]{FG1} E.~Frenkel and D.~Gaitsgory, {\it D-modules on the affine Grassmannian and representations of affine Kac-Moody algebras}, 
Duke Math.Jour {\bf 152} (2004), 279--327. 

\bibitem[FG2]{FG2} E.~Frenkel and D.~Gaitsgory, {\it Local geometric Langlands correspondence: the spherical case},
Algebraic Analysis and around, {\bf 54} (2009), 167--186.

\bibitem[FG3]{FG3} E.~Frenkel and D.~Gaitsgory, {\it Weyl modules and opers without monodromy}, 
Arithmetic and geometry around quantization, (2010), 101--121. 

\bibitem[FG4]{FG4} E.~Frenkel and D.~Gaitsgory, {\it Local geometric Langlands correspondence and affine Kac-Moody algebras}, 
in: Algebraic Geometry and Number Theory: In Honor of Vladimir Drinfeld's 50th birthday, (2006), 69--260.

\bibitem[FG5]{FG5} E.~Frenkel and D.~Gaitsgory, {\it Fusion and convolution: applications to affine Kac-Moody algebras at the
critical level}, Pure Appl. Math. Q. {\bf 2} (2006), 1255--1312.

\bibitem[FG6]{FG6} E.~Frenkel and D.~Gaitsgory,  {\it D-modules on the affine flag variety and representations of 
affine Kac-Moody algebras}, Represent. Theory {\bf 13} (2009), 470--608.

\bibitem[FGV]{FGV} E.~Frenkel, D.~Gaitsgory and K.~Vilonen, 
{\it Whittaker patterns in the geometry of moduli spaces of bundles on curves},
Annals of Math. {\bf 153} (2001), 699--748.

\bibitem[FLMM1]{FLMM1} G.~Fortuna, D.~Lombardo, A.~Maffei, V.~Melani,
{\it Local opers with two singularities: the case of $\mathfrak{sl}_2$}, Comm.Math.Phys {\bf 394}, (2022), 303--1360.

\bibitem[FLMM2]{FLMM2} G.~Fortuna, D.~Lombardo, A.~Maffei, V.~Melani,
{\it The semi-infinite cohomology of Weyl modules with two singular points}, PAMQ {\bf 20} (2024), 1251--1284. 

\bibitem[Ga1]{Ga1} D.~Gaitsgory, {\it Outline of the proof of the geometric Langlands conjecture for $GL(2)$}, 
Ast\'erisque {\bf 370} (2015), 1--112. 

\bibitem[Ga2]{Ga2} D.~Gaitsgory, {\it A conjectural extension of the Kazhdan-Lusztig equivalence}, Publications of RIMS
{\bf 57} (2021), 1227--1376. 

\bibitem[Ga3]{Ga3} D.~Gaitsgory, {\it A "strange" functional equation for Eisenstein series and miraculous
duality on the moduli stack of bundles}, Annales Scientifiques de l'ENS {\bf 50} (2017), 1123--1162. 

\bibitem[Ga4]{Ga4} D.~Gaitsgory, {\it The Atiyah-Bott formula for the cohomology of the moduli space of bundles on a curve}, arXiv:1505.02331 (2015).

\bibitem[Ga5]{Ga5} D.~Gaitsgory, {\it Sheaves of categories and the notion of 1-affineness}, Stacks and categories in geometry, topology, and algebra, {\bf 643} (2015), 127--225.

\bibitem[Ga6]{Ga6} D.~Gaitsgory, {\it The local and global versions of the Whittaker category}, 
PAMQ {\bf 16}, (2020), 775--904. 

\bibitem[Ga7]{Ga7} D.~Gaitsgory, {\it Ind-coherent sheaves}, Moscow Mathematical Journal, {\bf 13} (2013), 399-528.

\bibitem[GaLu]{GaLu} D.~Gaitsgory and J.~Lurie, {\it Weil's conjecture for function fields. Volume 1}, Annals of Mathematics Studies (2014).

\bibitem[GLys]{GLys} D.~Gaitsgory and S.~Lysenko, {\it Parameters and duality for the metaplectic geometric Langlands theory}, 
Selecta Math. New Ser. {\bf 24} (2018), 227--301. Also arXiv: 1608.00284

\bibitem[GaRo1]{GaRo1}  D.~Gaitsgory and N.~Rozenblyum, {\it DG indschemes}, Perspectives in representation theory, {\bf 610} (2014), 139--251.

\bibitem[GaRo2]{GaRo2}  D.~Gaitsgory and N.~Rozenblyum, {\it Crystals and D-modules}, Pure and Applied Mathematics Quarterly, {\bf 10} (2014), 57--154.

\bibitem[GaRo3]{GaRo3}  D.~Gaitsgory and N.~Rozenblyum, {\it A study in derived algebraic geometry, Vol. 1: Correspondences and Duality}, 
Mathematical surveys and monographs {\bf 221} (2017), AMS, Providence, RI.

\bibitem[GaRo4]{GaRo4}  D.~Gaitsgory and N.~Rozenblyum, 
{\it A study in derived algebraic geometry, Vol. 2: Deformtions, Lie theory and formal geometry}, 
Mathematical surveys and monographs {\bf 221} (2017), AMS, Providence, RI.

\bibitem[GLC1]{GLC1} D.~Gaitsgory and S.~Raskin, {\it Proof of the geometric Langlands conjecture I: construction of the functor}, 
arXiv:2405.03599

\bibitem[GLC4]{GLC4} D.~Arinkin, D.~Beraldo, L.~Chen, J.~Faergeman, D.~Gaitsgory, K.~Lin, S.~Raskin and N.~Rozenblyum, \hfill \newline
{\it Proof of the geometric Langlands conjecture IV: ambidexterity},  \hfill \newline
available at: https://people.mpim-bonn.mpg.de/gaitsgde/GLC/

\bibitem[HR]{HR} J.~Hilburn and S.~Raskin, {\it Tate's thesis in the de Rham setting}, Journal of the American Mathematical Society, {\bf 36} (2023), 917--1001.

\bibitem[Lin]{Lin} K.~Lin, {\it Poincar\'e series and miraculous duality}, arXiv:2211.05282.

\bibitem[Lu1]{Lu1} J.~Lurie, {\it Higher topos theory}, Princeton University Press (2009).

\bibitem[Lu2]{Lu2} J.~Lurie, {\it Higher Algebra} (2017). Available at: https://www.math.ias.edu/~lurie/papers/HA.pdf.

\bibitem[Lu3]{Lu3} J.~Lurie, {\it Spectral algebraic geometry} (2018). Available at: https://www.math.ias.edu/~lurie/papers/SAG-rootfile.pdf.

\bibitem[Ra1]{Ra1} S.~Raskin, {\it A geometric proof of the Feigin-Frenkel theorem}, Representation Theory {\bf 16} (2012), 489--512. 

\bibitem[Ra2]{Ra2} S.~Raskin, {\it W-algebras and Whittaker categories}, Selecta Mathematica {\bf 27} (2021), 46.

\bibitem[Ra3]{Ra3} S.~Raskin, {\it Chiral principal series categories I: finite dimensional calculations}, Adv. in Math. {\bf 388} (2021). 

\bibitem[Ra4]{Ra4} S.~Raskin, {\it Chiral principal series categories II: the factorizable Whittaker category}, \newline
available at: https://gauss.math.yale.edu/~sr2532/cpsii.pdf.

\bibitem[Ra5]{Ra5} S.~Raskin, {\it Homological methods in semi-infinite contexts}, arXiv:2002.01395 (2020).

\bibitem[Ra6]{Ra6} S.~Raskin, {\it Chiral categories} (2015), available at: https://gauss.math.yale.edu/~sr2532/chiralcats.pdf.

\bibitem[Ro1]{Ro1} N.~Rozenblyum, {\it Modules over a chiral algebra}, arXiv:1010.1998 (2010). 

\bibitem[Ro2]{Ro2} N.~Rozenblyum, {\it Connections on moduli spaces and infinitesimal Hecke modifications}, arXiv:2108.07745 (2021). 


\end{thebibliography}
\end{document}